\documentclass[11pt,oneside,openany]{report}
\usepackage[utf8]{inputenc}
\usepackage[T1]{fontenc}
\usepackage{lmodern,amsmath,amssymb,amsthm,mathtools,mathrsfs,booktabs,array,tabularx,enumitem,microtype}
\usepackage[a4paper,top=24mm,bottom=23mm,left=26mm,right=26mm,headheight=15pt]{geometry}
\usepackage{xcolor}
\usepackage{tikz}
\definecolor{linkblue}{HTML}{16465B}
\usepackage[unicode,colorlinks=true,linkcolor=linkblue,citecolor=linkblue,urlcolor=linkblue,bookmarksnumbered=true]{hyperref}
\usepackage{bookmark}
\hypersetup{pdftitle={K\"ahler–Einstein metrics on quartic hypersurfaces},pdfauthor={Kefeng Liu and Hao Xu},pdfsubject={Uniform K-stability of smooth complex Fano hypersurfaces}}
\newtheorem{theorem}{Theorem}[chapter]
\newtheorem{lemma}[theorem]{Lemma}\newtheorem{proposition}[theorem]{Proposition}
\newtheorem{corollary}[theorem]{Corollary}
\usepackage{titlesec}
\titleformat{\chapter}[display]{\normalfont\huge\bfseries}
  {\normalfont\small\scshape Chapter \thechapter}{7pt}{}
\titlespacing*{\chapter}{0pt}{5pt}{18pt}
\title{K\"ahler–Einstein metrics on quartic hypersurfaces}
\author{Kefeng Liu \and Hao Xu}\date{}
\begin{document}
\pagenumbering{arabic}
\thispagestyle{plain}
\begin{center}
  \vspace*{12mm}
  {\Huge\bfseries K\"ahler–Einstein metrics on\\[4pt] quartic hypersurfaces\par}
  \vspace{5mm}
  \begin{tabular}{@{}c@{\qquad\qquad}c@{}}
    {\large Kefeng Liu} & {\large Hao Xu}\\[2pt]
    \href{mailto:kefengliu@cqut.edu.cn}{\nolinkurl{kefengliu@cqut.edu.cn}} &
    \href{mailto:haoxu@cqut.edu.cn}{\nolinkurl{haoxu@cqut.edu.cn}}
  \end{tabular}\par
  \vspace{2mm}
  Chongqing University of Technology\par
\end{center}
\vspace{6mm}
\begin{center}\textbf{Abstract}\end{center}
\noindent We prove K-stability, and hence the existence of K\"ahler--Einstein metrics, for all smooth quartic Fano hypersurfaces.
\par\medskip
\begingroup
  \let\chapter\section
  \tableofcontents
\endgroup
\clearpage
\chapter{Introduction}\label{ch:intro}

Shing-Tung Yau's proof of the Calabi conjecture laid the foundations for the
modern existence theory of K\"ahler--Einstein metrics \cite{yau1978}. In
particular, on a compact K\"ahler manifold with vanishing real first Chern
class, every K\"ahler class contains a unique Ricci-flat K\"ahler metric.
For Fano manifolds, Yau was the first to propose that an algebraic
stability condition should characterize the existence of a K\"ahler--Einstein
metric \cite[Problem~65]{yau1993}. The Yau--Tian--Donaldson theorem
identifies this condition as K-polystability \cite{ar1205.6214,kdCDS,kdTian}.
We use this characterization to study smooth quartic Fano hypersurfaces.

\section{The hypersurface problem}
Let \(X=X_d^n\) be a smooth complex hypersurface of degree \(d\) in
\(\mathbf P^{n+1}\), and write \(H=\mathcal O_X(1)\). Adjunction gives
\(-K_X=(n+2-d)H\). Thus \(X\) is Fano exactly when \(d\le n+1\).
We ask when \emph{every} smooth member of a fixed degree and dimension
admits a K\"ahler form \(\omega\in 2\pi c_1(X)\) with
\(\operatorname{Ric}(\omega)=\omega\). The cases \(d=1,2\) are homogeneous:
projective space and smooth quadrics carry invariant K\"ahler--Einstein
metrics. For \(d\ge3\) finite automorphism groups allow us to work with
the strict stability threshold.

\section{K-stability and the stability threshold}
\label{sec:kd-definitions}

Throughout this section $X$ is a smooth complex projective Fano
variety of dimension $n$, polarized by $L=-K_X$. We use the
following conventions for test configurations
\cite{kdDonaldson,kdBHJ}. A normal ample test configuration of
exponent $e\ge 1$ consists of a flat projective morphism
$\pi\colon\mathcal X\to\mathbf A^1$ with normal total space, a
relatively ample line bundle $\mathcal L$, and a $\mathbf C^*$-action
on $(\mathcal X,\mathcal L)$ lifting multiplication on $\mathbf A^1$.
It includes a polarized identification of the fiber at $1$ with
$(X,L^{\otimes e})$. Transport by the action then identifies the
restriction over $\mathbf C^*$ equivariantly with
$(X,L^{\otimes e})\times\mathbf C^*$, where the action in this
trivialization is on the second factor.

Write $N_m$ for the dimension of
$H^0(\mathcal X_0,\mathcal L_0^{\otimes m})$ and $w_m$ for the
sum of its weights, with a character $t\mapsto t^\lambda$
assigned weight $\lambda$. Here the action on sections is
$(t\cdot s)(x)=t\bigl(s(t^{-1}x)\bigr)$. In large degree write
\[
 N_m=a_0m^n+a_1m^{n-1}+O(m^{n-2}),\qquad
 w_m=b_0m^{n+1}+b_1m^n+O(m^{n-1}).
\]
Our normalization of the Donaldson--Futaki invariant is
\[
 \operatorname{DF}(\mathcal X,\mathcal L)
       =\frac{2(a_1b_0-a_0b_1)}{a_0^2}.
\]
A character twist of the linearization adds $cmN_m$ to $w_m$
and leaves this expression invariant. Replacing $\mathcal L$ by a
positive tensor power also preserves it, as one sees by replacing
$m$ by the corresponding multiple in both expansions.

A product configuration has family $X\times\mathbf A^1$ and line
bundle $\operatorname{pr}_1^*L^{\otimes e}$. Its action has the form
$(x,z)\mapsto(\rho(t)x,tz)$ for a one-parameter subgroup $\rho$ of
polarized automorphisms of $X$, with a lift to the line bundle.
These descriptions are taken up to equivariant polarized
isomorphism and character twist. The product is called trivial
when $\rho$ acts trivially on $X$. The variety $X$ is
\emph{K-semistable} if every normal ample test configuration has
$\operatorname{DF}\ge0$. It is \emph{K-polystable} if in addition
every configuration with $\operatorname{DF}=0$ is a product, and
\emph{K-stable} if every nontrivial configuration has
$\operatorname{DF}>0$.

For the uniform version let $\lambda_{\max,m}$ be the greatest
weight on the same space of sections. The non-Archimedean
$J$-functional in units of $L$ is
\[
 J^{\mathrm{NA}}(\mathcal X,\mathcal L)
   =\lim_{m\to\infty}\frac{1}{em}
            \left(\lambda_{\max,m}-\frac{w_m}{N_m}\right).
\]
It is nonnegative and, for normal ample test configurations, vanishes
precisely on trivial configurations \cite{kdBHJ}. Uniform K-stability
means that there is an $\varepsilon>0$ such that
$\operatorname{DF}\ge\varepsilon J^{\mathrm{NA}}$ for every such
configuration.

\paragraph{The valuative criterion.}
For a prime divisor $E$ on a smooth birational model
$\mu\colon Y\to X$ and an ample line bundle $M$ on $X$, put
\[
 A_X(E)=1+\operatorname{ord}_E(K_Y-\mu^*K_X),\qquad
 S_M(E)=\frac{1}{M^n}\int_0^\infty
               \operatorname{vol}_Y(\mu^*M-tE)\,dt .
\]
The stability threshold is $\delta(X,M)=\inf_E A_X(E)/S_M(E)$,
where $E$ ranges over prime divisors on all such models.
The valuative criterion gives
\begin{align*}
 X\ \text{is K-semistable}
      &\ \Longleftrightarrow\ \delta(X,-K_X)\ge1,\\
 X\ \text{is uniformly K-stable}
      &\ \Longleftrightarrow\ \delta(X,-K_X)>1 .
\end{align*}
See Fujita and Blum--Jonsson \cite{kdFujita,ar1706.04548}. For Fano varieties
K-stability also implies uniform K-stability
\cite[Theorem~1.6]{ar2102.09405}. We use the hyperplane polarization
$H=\mathcal O_X(1)$ on a degree-$d$ hypersurface
$X\subset\mathbf P^{n+1}$. Adjunction gives $-K_X=rH$ with
$r=n+2-d$. Homogeneity of volume and the substitution $t=rs$
give $S_{rH}(E)=rS_H(E)$, hence
\[
 \delta(X,-K_X)=\frac{\delta(X,H)}{r}.
\]
In particular our strict target in hyperplane units is
$\delta(X,H)>r$.

\paragraph{K\"ahler--Einstein metrics and automorphisms.}
The Yau--Tian--Donaldson theorem identifies K-polystability of a
smooth complex Fano variety with the existence of a K\"ahler metric
satisfying $\operatorname{Ric}(\omega)=\omega$ in
$2\pi c_1(X)$ \cite{ar1205.6214,kdCDS,kdTian}.
For a smooth Fano hypersurface of degree $d\ge3$ the automorphism
group is finite. In dimension $n\ge3$ the Picard group is generated
by $H$ by Lefschetz, so each automorphism preserves the hyperplane
embedding. In dimension two the only such hypersurface is the cubic
surface, whose embedding is anticanonical. Finiteness now follows
from the finiteness of the projective linear stabilizer
\cite{MM}. Consequently its one-parameter automorphisms are
trivial. For these hypersurfaces K-polystability, K-stability,
uniform K-stability and existence of a K\"ahler--Einstein metric
are equivalent.

\section{Established cases}\label{sec:established-cases}
Let $X_d^n\subset\mathbf P^{n+1}_{\mathbf C}$ be a smooth hypersurface of
degree $3\le d\le n+1$, and let $H=\mathcal O_X(1)$. Adjunction gives
$-K_X=rH$ with $r=n+2-d$. Xu and Zhuang ask whether every such
hypersurface is K-stable \cite[Problem~6]{ar2601.15576}. For a complex Fano
manifold, K-stability is equivalent to uniform K-stability, and
K-polystability is equivalent to the existence of a K\"ahler--Einstein
metric \cite[Theorem~1.6]{ar2102.09405}. All stability statements below use the
anticanonical polarization.

\begin{samepage}
\textbf{Results for every smooth member.}
Throughout the table \(d\ge3\) and \(r=n+2-d\).
\begin{center}
\begin{tabularx}{\linewidth}{@{}>{\raggedright\arraybackslash}p{.23\linewidth}>{\raggedright\arraybackslash}X@{}}
\toprule
Hypothesis & Conclusion and source \\
\midrule
\(d=3,\ n\ge2\) &
The next table gives the cases \(2\le n\le4\). In dimensions \(n\ge5\),
the manuscripts of Kong--Shen--Zhao--Zhao \cite[Theorem~1.1 and
Corollary~1.2]{keKongCubics} and Pan--Wu \cite[Theorem~1.1 and
Remark~1.2]{kePanWuCubics} identify K-stability with GIT stability.
Every smooth cubic is GIT-stable.\\[3pt]
\(r=1\), i.e.\ \(d=n+1\) &
Fujita proves K-stability and K\"ahler--Einstein existence \cite[Corollary~1.4]{ar1606.08261}.
Abban and Zhuang prove uniform K-stability \cite[Corollary~4.9(2)]{ar2003.13788}.\\[3pt]
\(r=2\), i.e.\ \(d=n\) &
Abban and Zhuang prove uniform K-stability \cite[Theorem~1.1, restated as Theorem~4.21]{ar2003.13788}.\\[3pt]
\(r\ge3\) and \(n\ge r^3\) &
Abban and Zhuang prove uniform K-stability \cite[Theorem~B, restated as Theorem~4.1]{ar2101.09246}.\\
\bottomrule
\end{tabularx}
\end{center}
\end{samepage}

\begin{samepage}
\textbf{Every smooth member in dimensions \(n\le4\).}
For \(d\ge3\) the Fano condition leaves exactly the following cases. In each, \(X\) is K-stable and admits a K\"ahler--Einstein metric.
\begin{center}
\begin{tabularx}{\linewidth}{@{}>{\raggedright\arraybackslash}p{.23\linewidth}>{\raggedright\arraybackslash}X@{}}
\toprule
Hypersurface & Source for this case \\
\midrule
Cubic surface &
Tian proves K\"ahler--Einstein existence \cite[Main theorem, p.~101]{keTian}.
Abban and Zhuang compute \(\delta(X)=3/2\) if \(X\) has an Eckardt point and \(27/17\) otherwise \cite[Corollary~4.15]{ar2003.13788}.\\[3pt]
Cubic threefold &
Liu and Xu \cite[Theorem~1.1 and Corollary~1.2(1)]{ar1706.01933}.\\[3pt]
Quartic threefold &
Fujita \cite[Corollary~1.4]{ar1606.08261}, with \(n=3\) and \(d=4\).\\[3pt]
Cubic fourfold &
Liu \cite[Corollary~1.2(1)]{ar2007.14320}.\\[3pt]
Quartic fourfold &
Abban and Zhuang \cite[Theorem~1.1]{ar2003.13788}, with \(n=d=4\).\\[3pt]
Quintic fourfold &
Fujita \cite[Corollary~1.4]{ar1606.08261}, with \(n=4\) and \(d=5\).\\
\bottomrule
\end{tabularx}
\end{center}
\end{samepage}

\paragraph{General members and explicit equations.}
For each $3\le d\le n+1$, the Fermat equation
$x_0^d+\cdots+x_{n+1}^d=0$ defines a K-stable hypersurface
\cite[Corollary~1.4(1), restated as Corollary~4.17(1)]{ar2004.09413}.
The smooth locus of the projective space of degree-$d$ equations is
irreducible. Its universal family is a smooth family of Fano manifolds,
so openness of K-stability gives a dense Zariski-open subset of
K-stable members containing the Fermat point
\cite[Theorem~4.5]{ar1907.02408}. Campo, DeVleming and Duarte Guerreiro give
the same conclusion for hypersurfaces in a preprint that also treats
general Fano complete intersections \cite[Theorem~1.3 and
Corollary~5.5]{ar2608.09293}. For individual smooth hypersurfaces, a related
criterion allows adjoining a pure power. If $Y=(f=0)\subset\mathbf P^m$
is a smooth K-polystable Fano hypersurface of degree $3\le d\le m$,
then $(f+x_{m+1}^d=0)\subset\mathbf P^{m+1}$ is K-polystable by the
cyclic-cover theorem applied to the degree-$d$ cover of $\mathbf P^m$
branched along $Y$ \cite[Theorem~5.1]{ar2001.10557}. The hypersurface
formulation gives K-stability \cite[Corollary~5.3]{ar2608.09293}.

\section{Results of our research}\label{sec:results}
Write \(X=X_d^n\subset\mathbf P^{n+1}\) for a smooth complex hypersurface,
\(H=\mathcal O_X(1)\) and \(-K_X=rH\), where \(r=n+2-d\).

\begin{theorem}[Ranges for every smooth member]\label{thm:range}
The arguments in the chapters below give uniform K-stability, and hence a
positive K\"ahler--Einstein metric, for every smooth complex hypersurface in
the following ranges.
\[
\begin{array}{c|l}
  d & n\\ \hline
  3 & n\ge2\\
  4 & n\ge3\\
  5 & 4\le n\le11\\
  d\ge6 & d-1\le n\le2d .
\end{array}
\]
\end{theorem}
For cubics, the result in every dimension follows from the works cited in
Section~\ref{sec:established-cases}. Chapter~\ref{ch:cubics} gives a proof
for \(n\ge5\) using the limiting hyperplane system and its section ring.
For quartics, Section~\ref{sec:established-cases} covers \(n=3,4\).
The valuation estimates in Chapters~\ref{ch:fivefold}--\ref{ch:eight} treat \(5\le n\le8\), and
Chapter~\ref{ch:nov30-replacements} treats every \(n\ge9\) by
K-polystable replacement.

For \(d\ge5\), the common range is \(d-1\le n\le2d\).
Since \(d=n+2-r\), the upper bound \(n\le2d\) is equivalent to
\(n\ge2r-4\). Theorem~\ref{thm:nov27-second-row} proves the endpoint
\(n=2d\). The quintic elevenfold of Theorem~\ref{thm:nov28-quintic-eleven}
lies one dimension beyond this interval.

For an integral curve \(D\subset X\), write \(\delta_D(X,H)\) for the
infimum of \(A_X(v)/S_H(v)\) over divisorial valuations \(v\) whose
centers on \(X\) contain \(D\). In each case below, \(M\) is a common lower bound for
\(\delta_D(X,H)\) over all lines and integral conics \(D\subset X\).
Together with the point estimates in the indicated chapter, the condition
in the last column gives \(\delta(X,H)>r\).
\begin{center}
\begin{tabular}{@{}lccc@{}}
\toprule
\(X\) & \(r\) & Lower bound \(M\) & Sufficient condition on \(M\)\\
\midrule
Quartic sevenfold (Ch.~\ref{ch:boundary}) & 5 & \(16/3\) & \(M\ge16/3\)\\
Quartic eightfold (Ch.~\ref{ch:eight}) & 6 & \(198/31\) & \(M>51/8\)\\
Quintic ninefold (Ch.~\ref{ch:quintic}) & 6 & \(58/9\) & \(M\ge58/9\)\\
Sextic elevenfold (Ch.~\ref{ch:sextic}) & 7 & \(2448/329\) & \(M\ge7433/1000\)\\
\bottomrule
\end{tabular}
\end{center}
For the quintic ninefold, counting jets bounds the mean of the lowest
fifth of the valuation weights at a closed point. For the sextic
elevenfold, the weakest line estimate exceeds the sufficient bound by
\(2448/329-7433/1000=2543/329000\).

\subsection*{Main points of the proof for quartics}

\paragraph{Dimensions three through eight.}
The threefold and fourfold cases follow from Fujita and from Abban and
Zhuang, respectively (Section~\ref{sec:established-cases}). For a quartic
fivefold the valuative criterion reduces the proof to estimates at
tangent hyperplane sections of multiplicity two, three or four.
In dimensions six through eight, bounds along lines and conics exclude
positive-dimensional minimizing centers. It remains to control valuations
centered at a point. On the eightfold, for instance, a free line on an
auxiliary quartic sixfold gives a point estimate at cubic tangent contact.
The resulting criterion requires a common curve bound greater than
\(51/8\), and the bound \(198/31\) exceeds it by \(3/248\)
(Theorem~\ref{thm:nov28-quartic-eightfold}).

\paragraph{Dimension nine.}
Specialize a family of K-stable quartics to the given smooth quartic
and let \(Z\) be its K-polystable replacement, with limiting hyperplane
class \(L\). If \(Z\) differs from the given ninefold, the image of
\(|L|\) must be the seven-dimensional cone over \(\operatorname{Gr}(2,5)\)
(Corollary~\ref{cor:dec1-root-grassmannian-only}). The five Pl\"ucker
relations and the quartic Hilbert function force new quadratic sections.
These vary on the general hyperplane fiber, so the image of \(|2L|\)
has dimension eight or nine. Theorem~\ref{thm:next-q9-no-eight-quadratic}
excludes dimension eight. In dimension nine a general surface fiber
maps birationally to a plane. On the inverse image of a general line,
the pullback of \(L\) has degree \(a<2/3\). The pole at infinity
accounts for at least \(1/2\) of this degree, and the local Cartier
index bound makes each further contribution at least \(1/6\).
Thus a single exceptional divisor contributes to adjunction on the
line. Its ratio of log discrepancy to order along \(|L|\) is
\(7-2/a<4\), contradicting
\(\operatorname{lct}(Z;|L|)>4\)
(Proposition~\ref{prop:next-q9-no-nine-quadratic}).

\paragraph{Dimensions ten and higher.}
Here the hyperplane system itself excludes a distinct replacement.
After the comparison of fillings and the classification of the possible
hyperplane images, the remaining case has a smooth del Pezzo base of
dimension \(b\le6\). The first-jet slope estimate forces
\(3b\ge 2n-1\ge19\), whereas \(b\le6\) gives \(3b\le18\)
(Lemma~\ref{lem:dec1-hyperplane-atiyah-slope} and
Theorem~\ref{thm:dec1-quartic-high-dimensions}). Thus in every
dimension \(n\ge9\) the smooth quartic is its own K-polystable
replacement. Its finite automorphism group gives K-stability, and the
Yau--Tian--Donaldson theorem gives a K\"ahler--Einstein metric.

\subsection{Reading the proofs}\label{sec:reading-the-proofs}

\paragraph{The quartic theorem.}
Chapters~\ref{ch:fivefold}, \ref{ch:sixfold}, \ref{ch:boundary}
and~\ref{ch:eight} contain the proofs for dimensions five through eight
in order. In Chapter~\ref{ch:nov30-replacements}, the ninefold case is
Theorem~\ref{thm:next-all-smooth-quartic-ninefolds} and the result for
\(n\ge10\) is Theorem~\ref{thm:dec1-quartic-high-dimensions}.

\paragraph{Shared estimates.}
Chapter~\ref{ch:cubics} develops the limiting hyperplane argument for cubics.
Chapter~\ref{ch:nov30-filtrations} gives the results on filtered section
rings used in the replacement argument. Chapter~\ref{ch:foundations} develops
the point estimates, Chapter~\ref{ch:global} treats positive-dimensional
centers and proves the range \(n\ge2r-1\), and Chapter~\ref{ch:models}
gives the local estimates along lines and conics. The eightfold proof
also uses the interpolation estimate in Section~\ref{sec:nov28-direct-point}.

\paragraph{Other degrees and the parameter space.}
Chapters~\ref{ch:quintic} and~\ref{ch:sextic} treat the quintic ninefold
and the sextic elevenfold. For \(n=2d-1\), Chapter~\ref{ch:higher} gives
the global criterion for \(d\ge14\) and Chapter~\ref{ch:fullrow}
completes the proof for \(d\ge7\). Chapter~\ref{ch:secondrow} treats
\(n=2d\) for \(d\ge5\). Chapter~\ref{ch:beyond} gives sufficient
conditions beyond \(n=2d\) and proves the quintic elevenfold case.
Chapter~\ref{ch:families} constructs explicit K-stable families.
For \(d\ge4\) and \(n\ge d-1\), Chapter~\ref{ch:new-counterexamples}
bounds the codimension of the complement of the uniformly K-stable klt
locus in the full parameter space by \(\binom{2d-2}{d-1}/d\) from below
(Theorem~\ref{thm:new-catalan}). For quartic ninefolds the smooth
exceptional locus is empty. In the full space the codimension lies
between \(9\) and \(211\). For the upper bound, Section~\ref{sec:dec1-unstable-flags}
exhibits a family of quartics with a single double point that fail uniform K-stability. Chapter~\ref{ch:checks}
explains how to obtain and run the arithmetic certificates.

\begin{figure}[!ht]
\centering
\begin{tikzpicture}[
  x=1mm,y=1mm,
  case/.style={draw=linkblue!65,fill=linkblue!3,line width=.45pt,
    anchor=north,text width=46mm,minimum height=36mm,
    inner xsep=1.6mm,inner ysep=2mm,align=center,font=\small},
  conclusion/.style={draw=linkblue,line width=.65pt,fill=linkblue!7,
    anchor=north,inner xsep=4mm,inner ysep=2.4mm,align=center,font=\small}
]
  \node[case] (published) at (-53,0) {
    \textbf{Published cases}\\[3pt]
    Section~\ref{sec:established-cases}\\[4pt]
    \begin{tabular}{@{}l@{\quad}l@{}}
      \(n=3\) & Fujita\\
      \(n=4\) & Abban--Zhuang
    \end{tabular}};
  \node[case] (valuative) at (0,0) {
    \textbf{Point and curve estimates}\\[4pt]
    \begin{tabular}{@{}l@{\quad}l@{}}
      \(n=5\) & Chapter~\ref{ch:fivefold}\\
      \(n=6\) & Theorem~\ref{sixfold:thm}\\
      \(n=7\) & Theorem~\ref{boundary:main}\\
      \(n=8\) & Theorem~\ref{thm:nov28-quartic-eightfold}
    \end{tabular}};
  \node[case] (replacement) at (53,0) {
    \textbf{K-polystable replacements}\\[4pt]
    \(n=9\)\quad Theorem~\ref{thm:next-all-smooth-quartic-ninefolds}\\
    Quadratic system \(|2L|\)\\[3pt]
    \(n\ge10\)\quad Theorem~\ref{thm:dec1-quartic-high-dimensions}\\
    Hyperplane system \(|L|\)};
  \coordinate (join) at (0,-41);
  \node[conclusion] (quartics) at (0,-47) {
    \textbf{Every smooth quartic \(X_4^n\subset\mathbf P^{n+1}\), \(n\ge3\), is uniformly K-stable}\\[2pt]
    and admits a K\"ahler--Einstein metric.};
  \draw[draw=linkblue,line width=.6pt] (published.south) |- (join)
    -| (replacement.south);
  \draw[draw=linkblue,line width=.6pt,->,>=stealth] (valuative.south) -- (quartics.north);
\end{tikzpicture}
\caption{The quartic case of Theorem~\ref{thm:range} by dimension.
Here \(L\) is the limiting hyperplane class on a K-polystable replacement.}
\label{fig:quartic-proof-by-dimension}
\end{figure}
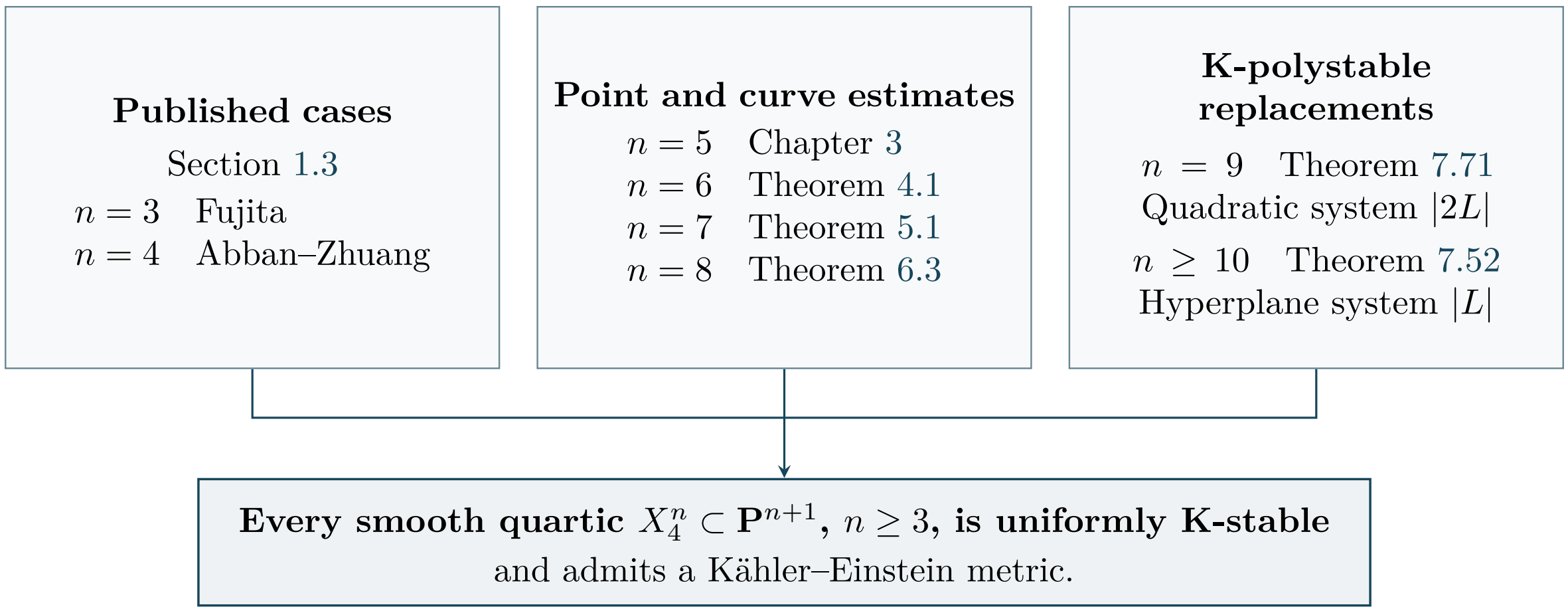

\begin{table}[!ht]
\centering
\caption{Results from later chapters used in the quartic proofs.}
\label{tab:quartic-proof-dependencies}
\smallskip
\begin{tabularx}{\linewidth}{@{}>{\raggedright\arraybackslash}p{0.23\linewidth}
  >{\raggedright\arraybackslash}X@{}}
\toprule
Quartic case & Results used in the proof \\
\midrule
\(n=5,6\)\newline Chapters~\ref{ch:fivefold}--\ref{ch:sixfold}
& Both chapters derive their point and curve estimates using the preceding
material and the cited literature. \\[5pt]
\(n=7\)\newline Chapter~\ref{ch:boundary}
& Equations \eqref{models:eq-mixed-monomial} and \eqref{models:eq-pure-powers}
in Section~\ref{models:sec:local-line-models} give the line threshold \(16/3\).
The proof in the same chapter of the general bound for \(n=2d-2\)
also uses the estimates of Chapter~\ref{ch:global}. \\[5pt]
\(n=8\)\newline Chapter~\ref{ch:eight}
& The local calculations in Chapters~\ref{ch:models}, \ref{ch:fullrow}
and~\ref{ch:secondrow} give a common lower bound of \(198/31\) for lines
and conics. Section~\ref{sec:nov27-quintic-ten} gives the Hessian
rank estimate used for rank-two lines. The interpolation argument in Section~\ref{sec:nov28-direct-point}
uses the partial-basis estimate of Chapter~\ref{ch:global} and shows that
a bound greater than \(51/8\) suffices. Here \(198/31-51/8=3/248\). \\[5pt]
\(n=9\) and \(n\ge10\)\newline Chapter~\ref{ch:nov30-replacements}
& The comparison of fillings in Theorem~\ref{thm:nov30-all-ci-fillings}
uses Chapter~\ref{ch:nov30-filtrations} on filtered section rings and
the CM degree. The analysis of the hyperplane image also uses the
determinant slope inequality \eqref{global:eq:optimizer-covering-slope}
from Chapter~\ref{ch:global}. \\
\bottomrule
\end{tabularx}
\end{table}

\

\noindent{\bf Acknowledgements} This work was supported by the
National Key Research and Development Program of China (NKPs) [Grant
Number 2024YFA1013201], the Natural Science Foundation of Chongqing
(NSFCQ) [Grant Number CSTB2024NSCQ-LZX0040], the Special Project of
Chongqing Municipal Science and Technology Bureau [Grant Number
2025CCZ015].

\

\noindent{\bf Disclosure of AI assistance.}

The authors took the lead in developing the mathematical ideas and proof
strategies, with assistance from AI tools in working main parts of the proofs. AI tools drafted
the text, which the authors revised and reviewed in full.

% ===== Cubics and their limiting hyperplane systems =====
\chapter{Cubics and their limiting hyperplane systems}\label{ch:cubics}
\begin{theorem}\label{cubics:main}
Let $X\subset\mathbf{P}^{n+1}_{\mathbf{C}}$ be a smooth cubic hypersurface with $n\ge2$. Then $(X,-K_X)$ is K-stable and admits a positive K\"ahler--Einstein metric.
\end{theorem}

Adjunction gives $-K_X=(n-1)\mathcal O_X(1)$. The low-dimensional cases are due to Tian for surfaces \cite{keTian}, Liu--Xu for threefolds \cite{ar1706.01933}, and Liu for fourfolds \cite{ar2007.14320}. Sections~\ref{cubics:setup}--\ref{cubics:moduli} prove the result for $n\ge6$. Section~\ref{cubics:rankfour} treats the remaining rank-four quadric when $n=5$. The all-dimensional cubic conclusion also follows directly from \cite[Theorem~1.1 and Corollary~1.2]{keKongCubics} and \cite[Theorem~1.1 and Remark~1.2]{kePanWuCubics}.

\paragraph{Comparison with the two 2026 proofs.}
After completing the proof in this chapter, we learned of the work of
Kong--Shen--Zhao--Zhao~\cite{keKongCubics} and Pan--Wu~\cite{kePanWuCubics}.
All three proofs recover the cubic embedding from the limiting hyperplane
class \(L\) and use GIT to identify the smooth cubic. For a lower-dimensional
image, our minimal-degree and ruling arguments are close to those of
Kong--Shen--Zhao--Zhao, while Pan--Wu analyze the general fibers
\cite[Propositions~3.1 and~3.4]{keKongCubics} \cite[\S6]{kePanWuCubics}.
In the full-dimensional case, Kong--Shen--Zhao--Zhao exclude a quadric
image using threshold and multiplicity estimates on a local index cover
of the limit, and Pan--Wu
prove Cartierness using mixed intersections and the codimension of the base
locus \cite[Theorem~4.1]{keKongCubics} \cite[Propositions~4.9 and~5.5]{kePanWuCubics}.
When the image of \(|L|\) is a \textbf{full-dimensional quadric}, our proof
takes a different route by deriving a cubic relation from the first three
graded pieces of the hyperplane section ring.

\paragraph{Our section-ring strategy for the quadric image.}
Let \(S=\bigoplus_{m\ge0}H^0(Z,mL)\) be the hyperplane section ring of a
K-polystable limit \(Z\), and let \(T\subset S\) be generated by its linear
sections. Once the image of \(|L|\) has dimension \(n\), a cubic image gives
\(T=S\) by equality of Hilbert series. For a quadric image \(q=0\), we
recover equations for the limit by comparing degrees two and three.
The quotient \(S_2/T_2\) is spanned by the class of one extra section \(y\).
In degree three, multiplication by \(y\) maps the \(n+2\) linear sections
to the \((n+1)\)-dimensional space \(S_3/T_3\). A nonzero linear form
\(\ell\) in the kernel therefore gives a cubic relation \(\ell y=g_3(x)\).
When \(\ell,g_3\) form a regular sequence in \(T\), this relation and the
Hilbert series give the presentation
\[
 S\simeq\mathbf C[x_0,\ldots,x_{n+1},Y]\big/(q(x),\,\ell(x)Y-g_3(x)),
 \qquad \deg Y=2 .
\]
The equations exhibit a singular complete-intersection double cover at the
weight-two point. K-semistability forces its normalized volume above Liu's
upper bound for such a singularity (Section~\ref{cubics:full}). In the remaining
rank-four fivefold case, we use a Hilbert--Burch matrix and its bigraded
Hilbert series to find a ruling divisor \(F\) with
\(4S_L(F)\ge19/18>1=A_Z(F)\), again contradicting K-semistability
(Section~\ref{cubics:rankfour}).

\section{The limiting system and the tangent sheaf}\label{cubics:setup}
Fix $n\ge5$. Let $Z$ be the K-polystable central fiber of a $\mathbf{Q}$-Gorenstein family over a smooth pointed curve with smooth K-stable cubic general fibers. After base change, Liu's Theorem~3.1 provides the hyperplane Weil class $L$ on $Z$. It is ample and $\mathbf{Q}$-Cartier and satisfies
\[
  -K_Z\sim(n-1)L,\qquad L^n=3,\qquad 2L\text{ is Cartier},\qquad
  h^0(Z,mL)=h^0(X,mH)\quad(m\ge0),
\]
where $X$ is a smooth cubic of the same dimension and $H=\mathcal O_X(1)$ \cite{ar2007.14320}. In particular $h^0(Z,L)=n+2$. Write $W\subset\mathbf{P}^{n+1}$ for the closure of the rational image of the complete system $|L|$. It is integral and nondegenerate. The K-polystable $\mathbf{Q}$-Fano variety $Z$ admits a weak K\"ahler--Einstein metric \cite{ar1411.1725}. Hence its reflexive tangent sheaf $T_Z$ is slope semistable with respect to $L$, by Druel--Guenancia--P\u{a}un \cite[Theorem A]{ar2008.05352}. The degree of a Weil divisor below means intersection with $L^{n-1}$, defined via the Cartier multiple $2L$.

\begin{lemma}\label{cubics:slope}
For a dominant rational map $h:Z\dashrightarrow Y$ to a smooth projective $s$-fold, $1\le s<n$, take a smooth resolution $\pi:T\to Z$, $p:T\to Y$ of the map. Then
\[
   (\pi^*L)^{n-1}\cdot p^*K_Y\ \le\
       -\frac{s(n-1)}{n}L^n.
\]
In particular, for $Y=\mathbf{P}^s$ and $V=p^*\mathcal O(1)$,
\begin{equation}
   \frac{(\pi^*L)^{n-1}\cdot V}{L^n}\ \ge\
      \frac{s(n-1)}{n(s+1)}. \label{cubics:projective-slope}
\end{equation}
\end{lemma}
\begin{proof}
On the big smooth open set where $h$ is regular, the differential gives a generically injective map $h^*\Omega^1_Y\to\Omega^1_Z$. Take the saturated rank-$s$ subsheaf in $\Omega_Z^{[1]}$ with this generic fiber. Its determinant contains $h^*K_Y$ in codimension one. The dual of the slope-semistable $T_Z$ is semistable and has slope $-(n-1)L^n/n$. The degree of the saturated subsheaf is therefore at most $-s(n-1)L^n/n$. Pushing down $p^*K_Y$ computes precisely the determinant class of $h^*K_Y$ on this open set. The projection formula gives the assertion. This is the differential calculation used for rational hyperplane maps by H\"oring and Peternell \cite[proof of Theorem 1.2]{ar2403.05250}.
\end{proof}

\section{A lower-dimensional image has minimal degree}\label{cubics:minimal}
Suppose $r=\dim W<n$. Nondegeneracy gives $r\ge1$. Resolve the normalization of $W$ by a smooth projective $Y$ and write $M$ for the pullback of $\mathcal O_W(1)$. Thus $M$ is big, nef and basepoint-free, with $M^r=\deg W$. Use $\pi:T\to Z$ and $p:T\to Y$ as in Lemma~\ref{cubics:slope}, resolving also the hyperplane system. Put $P=\pi^*L$. There is an effective $\mathbf{Q}$-divisor linearly equivalent over $\mathbf{Q}$ to $P-p^*M$. Indeed, the products of pairs of sections of $L$ define a subsystem of $|2L|$ whose moving part after resolution is $2p^*M$.

In fact $h^0(Y,M)=n+2$. Fix a coordinate section $q_0$ of $M$ and the corresponding $s_0\in H^0(Z,L)$. The fixed-part computation gives $\pi^*\operatorname{div}(s_0)-p^*\operatorname{div}(q_0)\ge0$. Thus for $q\in H^0(Y,M)$ the rational section $s_0 h^*(q/q_0)$ has effective Weil divisor on $Z$. This embeds the complete space $H^0(Y,M)$ into $H^0(Z,L)$. Its image contains all $n+2$ independent coordinate sections.

There is no nonzero section of $K_Y+(r-1)M$. If one existed, its effective pullback, the inequality $P-p^*M\ge0$ up to $\mathbf{Q}$-linear equivalence, and nefness of $P$ would give
\[
    P^{n-1}p^*K_Y\ge -(r-1)P^{n-1}p^*M
                    \ge -(r-1)L^n .
\]
This contradicts Lemma~\ref{cubics:slope}, since $r(n-1)/n>r-1$.

We now compute the sectional genus. The irregularity of $Y$ is zero. Indeed $Z$ has rational singularities and $H^1(Z,\mathcal O_Z)=0$ by Fano vanishing. Hence $H^1(T,\mathcal O_T)=0$, and pullback of holomorphic one-forms along the dominant map $T\to Y$ is injective. For $r=1$ it follows that $Y\simeq\mathbf{P}^1$ and $M^r=h^0(Y,M)-1$. Suppose $r\ge2$ and cut by $r-1$ general members of $|M|$ to a smooth curve $C$. On any rung $V$ of dimension $j\ge2$, $M|_V$ stays big and basepoint-free. Bertini gives smoothness of the next cut $D$. Its connectedness follows from $H^1(V,\mathcal O_V(-M))=0$, by Serre duality and Kawamata--Viehweg vanishing. For $j\ge3$, $H^2(V,\mathcal O_V(-M))=0$ as well, so $H^1(D,\mathcal O_D)=H^1(V,\mathcal O_V)=0$ down to the surface.

On a rung of dimension $j\ge3$, the adjunction restriction
\[
 H^0(V,K_V+(j-1)M)\longrightarrow H^0(D,K_D+(j-2)M|_D)
\]
is surjective by $H^1(V,K_V+(j-2)M)=0$. On the surface the obstruction group is $H^1(V,K_V)=0$. It follows that $H^0(C,K_C)=0$, so $C\simeq\mathbf{P}^1$. Also, $H^0(V,M)\to H^0(D,M|_D)$ is surjective with one-dimensional kernel at each step, by irregularity zero. Therefore
\[
 n+2=h^0(Y,M)=r-1+h^0(C,M|_C)=r+M^r .
\]
Together with the curve case this gives $\deg W=n+2-r=\operatorname{codim}W+1$.

\section{Rulings and the valuative criterion}\label{cubics:rulings}
A nondegenerate integral variety of minimal degree and codimension at least two is a rational normal scroll, possibly a cone, or the Veronese surface in $\mathbf{P}^5$ or a cone over it \cite{EH}. On the Veronese surface a double line gives a hyperplane section containing a moving prime divisor twice. On a scroll with a summand $\mathcal O_{\mathbf{P}^1}(a)$ of degree $a\ge2$, use a binary polynomial in that summand with a double zero at a general base point, with all other components zero. Its hyperplane section contains the ruling fiber there twice. These multiplicities persist on the complement of any cone vertex. On the dense open of $Z$ where the map to the relevant smooth locus is dominant and smooth, a general such divisor pulls back with a prime component. Some member of $|L|$ therefore has coefficient at least two along a prime divisor. This would imply
\[
 \alpha(Z)\le\frac{1}{2(n-1)}<\frac1{n+1},
\]
contradicting the alpha bound for K-semistable Fano varieties \cite[Theorem 2.5]{ar2007.14320}.

The remaining scroll has only summands of degrees zero and one. Its image is a cone, with possibly empty vertex, over $\mathbf{P}^1\times\mathbf{P}^s$ in its Segre embedding. Here $s+1=\deg W\ge3$. The following computation also applies to the cone over $\mathbf{P}^1\times\mathbf{P}^1$ when considering a full-dimensional rank-four quadric.

On a common resolution of the projections to the two factors let $U,V$ be the pullbacks of their hyperplane classes. Then $P-U-V$ is $\mathbf{Q}$-linearly equivalent to an effective divisor. Indeed the subsystem of matrix-coordinate hyperplanes off the vertex defines the Segre morphism with moving class $U+V$. As before one can take the system of all pairwise products in $|2L|$ and divide its fixed part by two. A dense open of $Z$ maps dominantly to $\mathbf{P}^1$, so a general fiber has a prime component $D'$ mapping onto a prime divisor $D$ on $Z$. Since $U-D'$ is effective up to linear equivalence, $P-D'-V$ is effective up to $\mathbf{Q}$-linear equivalence. Nefness gives, for $0\le t\le1$,
\[
 \operatorname{vol}(P-tD')\ge ((1-t)P+tV)^n
 \ge (1-t)^n L^n+nt(1-t)^{n-1}P^{n-1}V .
\]
The expected vanishing $S_L(D)$ is the integral of $\operatorname{vol}(P-tD')/L^n$ over $t\ge0$. The prime divisor $D$ on $Z$ has log discrepancy one. K-semistability and \eqref{cubics:projective-slope} now give
\[
 1\ \ge\ (n-1)S_L(D)\ \ge\ \frac{n-1}{n+1}
       \left(1+\frac{s(n-1)}{n(s+1)}\right).
\]
Equivalently $s(n-1)^2\le 2n(s+1)$. For $s\ge2$ and $n\ge5$ this is impossible, since $(n-1)^2>3n$. This excludes every lower-dimensional image. The same calculation with $s=1$ excludes a full-dimensional rank-four quadric for $n\ge6$, since $(n-1)^2>4n$.

\section{Recovering the cubic equation}\label{cubics:full}

\begin{proposition}
Let \(Z\) be a K-polystable \(\mathbb Q\)-Gorenstein limit of smooth cubic
\(n\)-folds as in \cite[Theorem~3.1]{ar2007.14320}, and let \(L\) be its limiting
hyperplane class. If \(n\ge5\) and the rational map defined by
\(\lvert L\rvert\) has \(n\)-dimensional image \(W\), then \(Z\) is a cubic
hypersurface and \(L\) is its hyperplane line bundle.
\end{proposition}

\begin{proof}
Put \(S=\bigoplus_{j\ge0}H^0(Z,\mathcal O_Z(jL))\) and let \(T\subset S\)
be the subalgebra generated by \(S_1\). By Liu's theorem,
\(-K_Z\sim(n-1)L\), \(2L\) is Cartier, and
\[
       \operatorname{Hilb}_S(z)=\frac{1-z^3}{(1-z)^{n+2}}.
\]
The image \(W=\operatorname{Proj}T\) is an integral nondegenerate
hypersurface of degree \(d\) in \(\mathbb P^{n+1}\). Comparison of the
leading Hilbert coefficients of \(T\subset S\) gives \(2\le d\le3\).
If \(d=3\), the two Hilbert series agree, so \(S=T\). The ampleness
of \(2L\) identifies \(Z\) with \(\operatorname{Proj}S\).
The degree-one map is then the cubic embedding. Indeed on each open
set where a degree-one section \(s\) is nonzero, its square trivializes
\(\mathcal O_Z(2L)\), and \(s\) trivializes \(\mathcal O_Z(L)\).
It remains to exclude \(d=2\).

Write \(T=\mathbb C[x_0,\ldots,x_{n+1}]/(q)\), with \(q\) an irreducible
quadratic form. Hilbert functions give
\(\dim S_2/T_2=1\) and \(\dim S_3/T_3=n+1\).
Choose \(y\in S_2\setminus T_2\). The map
\(T_1\to S_3/T_3,\ x\mapsto xy\), has a nonzero kernel, so
\(\ell y=g\) for some \(0\ne\ell\in T_1\) and \(g\in T_3\).
We have \(g\notin\ell T\) since \(S\) is a domain.

\emph{A regular sequence.}
Suppose first that \(\ell,g\) is a regular sequence on \(T\).
The map \(T[Y]/(\ell Y-g)\to S\) taking \(Y\) to \(y\) is injective.
It is injective after inverting \(\ell\). Moreover, reducing a relation
\(\ell f=(\ell Y-g)h\) modulo \(\ell\) shows that \(\ell\) divides \(h\),
so \(\ell\) is a nonzerodivisor on the source. With \(\deg Y=2\), the source has the Hilbert series of \(S\),
so this map is an isomorphism degree by degree. Consequently \(Z\) is
the weighted complete intersection
\[
  q(x)=0,\qquad \ell(x)Y=g(x)
  \quad\hbox{in }\mathbb P(1^{\,n+2},2).
\]
The usual double chart at the weight-two point \(p\) has equations
\(q(x)=\ell(x)-g(x)=0\) in \(\mathbb A^{n+2}\), and the involution
sends \(x\) to \(-x\). The chart is finite over its quotient, so its
two equations have height two and define a complete intersection of dimension \(n\).
Only the second equation has a linear term, so its vertex \(p'\) is
singular. The chart is \'etale over its normal quotient away from its
vertex. It is normal also at the vertex by Serre's criterion, since a
complete intersection is Cohen--Macaulay. Thus this is a quasi-\'etale
map of degree two from a klt germ with zero boundary.
The finite-degree formula, the local-to-global inequality for the
K-semistable \(Z\), and \(L^n=3\) give
\[
 \widehat{\operatorname{vol}}(p')
 =2\widehat{\operatorname{vol}}(p,Z)
 \ \ge\ 6\left(\frac{n(n-1)}{n+1}\right)^n
 \ >\ 2(n-1)^n .
\]
The strict inequality follows from \((1+1/n)^n<e<3\).
This contradicts Liu's ordinary-double-point upper bound for
singular local complete intersections
\cite[Theorems~2.12--2.14]{ar2007.14320}.
This excludes every \(q\) of rank at least five, since restriction to
any hyperplane lowers rank by at most two and \(T/\ell T\) is then a
domain.

\emph{Rank three.}
Write \(q=ac-b^2\). The projection of this quadric away from its
vertex is the smooth conic \([a:b:c]\). A tangent line at its general
point cuts twice that point, so the corresponding hyperplane divisor
on \(W\) is twice a ruling divisor. On a dense open of the smooth
Cartier locus of \(Z\), the composite map to the conic is defined and
dominant. A general fiber there has a prime divisorial component \(F\)
on taking its closure in \(Z\). The corresponding section of \(L\)
vanishes along \(F\) to order at least two, as one sees by trivializing
a nonvanishing coordinate of the conic there. Hence
\(\alpha(Z)\le1/(2(n-1))<1/(n+1)\). This contradicts the alpha bound for
K-semistable Fano varieties \cite[Theorem~2.5]{ar2007.14320}.

The remaining quadric has rank four, hence is a cone over the Segre quadric surface. Section~\ref{cubics:rulings} rules it out when $n\ge6$. Section~\ref{cubics:rankfour} excludes it at $n=5$. This proves the proposition.
\end{proof}

\section{From limits to every smooth cubic}\label{cubics:moduli}
Let $X$ be an arbitrary smooth cubic $n$-fold. For $n\ge6$ we have proved that a K-polystable replacement of a family of smooth K-stable cubics is a cubic. Section~\ref{cubics:rankfour} gives the same conclusion for $n=5$. We apply the moduli argument of Spotti and Sun \cite[\S5.2]{SS} and Liu \cite{ar2007.14320}, spelling out the passage to a specified $X$.

The smooth cubic parameter space is irreducible. K-stability of the Fermat cubic \cite[Corollary~1.4]{ar2004.09413} and openness \cite[Theorem~4.5]{ar1907.02408} give a nonempty open set of K-stable cubics. Choose a smooth pointed curve through $[X]$ with punctured family in this open set, shrinking the curve if necessary. K-moduli properness gives a K-polystable replacement $Z$ after finite base change \cite{proper,ar2102.09405}. By Sections~\ref{cubics:minimal}--\ref{cubics:full}, the limiting class $L$ on $Z$ is a very ample Cartier hyperplane class and $(Z,L)$ is a cubic.

Liu constructs a relative rank-one reflexive sheaf lifting $L$, flat over the curve with the required fiberwise restrictions \cite[proof of Theorem 3.1]{ar2007.14320}. By Nakayama's lemma a local generator on the central fiber lifts, and the resulting surjection from the local ring is an isomorphism by torsion-freeness and rank one. Properness allows shrinking the base to obtain a line bundle. Since $h^0(L)=n+2$ and $h^1(L)=0$ on all fibers, cohomology and base change, relative very ampleness and a trivialization of the pushforward give a family of embedded cubics.

K-semistability of $Z$ implies GIT semistability of its cubic equation by the CM-line comparison \cite[\S5.2]{SS}. The original family and its replacement have the same punctured GIT quotient map, hence the same limiting quotient point by separatedness. A smooth cubic is GIT-stable. The fiber of a stable point in the semistable locus is a single orbit. Thus $Z$ is projectively isomorphic to $X$, proving K-polystability of $X$. Smooth cubics have finite polarized automorphism group \cite{MM}, so K-polystability gives K-stability. The K\"ahler--Einstein conclusion follows from the Yau--Tian--Donaldson theorem \cite{ar1411.1725}. Together with the cited low-dimensional results, this proves Theorem~\ref{cubics:main}.

\paragraph{The monomial calculation in dimension five.}
Consider the tangent model in which the quadratic and cubic initial forms
are \(Q=xy\) and \(C=s^3\) in \(n\) tangent coordinates \(y_1,\ldots,y_n\),
with \(x,y,s\) distinct. Its basis in section degree \(m\) consists of
\(T^k Z^i U^g M\), where \(M\) is a tangent monomial divisible by neither
\(Q\) nor \(C\), and \(k+i+2g+\deg M=m\). Indeed, in tangent degree \(j\)
the model uses \(((Q,C)^{\max(0,j-m)})_j\). Every such monomial factors
uniquely as \(Q^iC^g M\) by removing the largest powers of the two coprime
monomials. Its least section degree is \(j-i-g\), since \(Z\) represents
\(Q\) in section degree one and \(U\) represents \(C\) in degree two.
The remaining section degree is the exponent of \(T\).

Fix a coordinate \(y_h\), and put \(a=\operatorname{ord}_{y_h}(Q)\),
\(b=\operatorname{ord}_{y_h}(C)\). In the generating series that counts
this basis by section degree \(t\) and \(y_h\)-exponent \(u\), the
restriction on \(M\) contributes the numerator
\((1-t^2u^a)(1-t^3u^b)\). The factors for \(y_h,Z,U\) contribute
\((1-tu)^{-1}(1-tu^a)^{-1}(1-t^2u^b)^{-1}\).
If \(H(t,u)\) denotes the resulting series, logarithmic differentiation gives
\[
 \left.\frac{u\partial_u H}{H}\right|_{u=1}
 =\frac{(1+a)t}{1-t}+\frac{(b-a)t^2}{1-t^2}
     -\frac{bt^3}{1-t^3}.
\]
Here \(H(t,1)=(1-t^3)/(1-t)^{n+2}\) has a pole of order \(n+1\) at
\(t=1\). Multiplication by \(1/(1-t)\) therefore multiplies its
degree-\(m\) coefficient asymptotically by \(m/(n+1)\). Applying this
to the simple pole of the logarithmic derivative gives the limiting
mean \(y_h\)-exponent divided by \(m\),
\[
                  B_h=\frac{1+a/2+b/6}{n+1}.
\]
Thus each of \(x,y,s\) has coefficient \(3/[2(n+1)]\) and each other
coordinate has coefficient \(1/(n+1)\). The coordinate hyperplanes have
simple normal crossings, so their basis divisor has threshold
\(1/\max_h B_h=2(n+1)/3\). This equals the cubic index \(n-1\) for
\(n=5\) and is smaller for \(n\ge6\).

\paragraph{Scope for singular cubics.}
For $n\ge3$ the projective cone over a smooth cubic $(n-1)$-fold is a klt cubic Fano variety. At its vertex the exceptional divisor $E$ of the blowup has $A(E)=n-2$ by hypersurface adjunction. Write its equation as $f(x_0,\ldots,x_n)=0$ with one extra vertex coordinate $z$. A section in degree $m$ is the direct sum of terms $z^{m-j}$ times base sections of degree $j$, of order $j$ at the vertex. The base Hilbert function has degree $n-1$, giving mean order divided by $m$ tending to $n/(n+1)$. Thus $A(E)-S_{-K}(E)=-2/(n+1)<0$.

\section{The rank-four quadric in dimension five}\label{cubics:rankfour}
It remains to exclude a rank-four quadric as the hyperplane image of a fivefold replacement. Retain the notation of Section~\ref{cubics:full} and suppose the image is a rank-four quadric. We keep $n$ in the Hilbert series to show the dimension dependence.

If \(q|_{\ell=0}\) is irreducible, the regular-sequence argument applies.
Otherwise a change of coordinates gives
\[
         q=ad-bc,\qquad \ell=a ,
\]
with \(n-2\) further independent variables. Indeed such a hyperplane
contains the radical of \(q\) and corresponds to a tangent plane of the
smooth quadric surface on the nondegenerate quotient. Thus
\(T/aT=\mathbb C[b,c,d,\ldots]/(bc)\) is reduced.
The sole remaining case has \(g\) a zero divisor modulo \(a\).
Interchanging \(b,c\) if necessary, write \(g=ap+bh\) in \(T\)
with \(p,h\) quadratic. After replacing \(y\) by \(y-p\), we have
\[
     y=(b/a)h,\qquad (b/a)c=d,\qquad h\notin (a,c)T.
\]
The last assertion follows because \(b(a,c)/a\subset T\) and
\(y\notin T\).

The subalgebra \(B=T[y]\subset S\) is presented by the maximal minors
of the matrix
\[
      \begin{pmatrix} a&c&h\\ b&d&Y\end{pmatrix}.
\]
Here the ambient polynomial ring has \(n+2\) degree-one variables and
the extra variable \(Y\) of degree two. On \(a\ne0\) or \(c\ne0\)
these minors give the graph of \(y\), which is integral of dimension
\(n+1\). On \(a=c=0\) the remaining equations are \(bh=dh=0\).
Since \(h\) restricts to a nonzero polynomial on this subspace, this
closed locus has dimension at most \(n\). The ideal of minors has
height two and is perfect by the Hilbert--Burch theorem. It has no
embedded primes and is generically reduced along its unique component,
so it is the prime graph ideal.

The ring \(B\) is normal. The open set where \(a\ne0\) or \(c\ne0\)
is smooth, and the locus \(a=b=c=d=0\) has dimension \(n-1\), hence
codimension two. In the chart \(b\ne0\) we can eliminate \(c=ad/b\).
The remaining equation is
\[
   aY-bh(a,b,ad/b,d,\ldots)=0.
\]
Each component of \(a=0\) in this chart has \(Y\) as a free variable,
since its equation in the remaining variables is \(h(0,b,0,d,\ldots)=0\),
with the latter polynomial nonzero. The derivative of the chart equation
with respect to \(a\) has the form \(Y\) minus a function independent of
\(Y\), and hence is nonzero generically on each component. The chart
\(d\ne0\) has the same argument with \(a=bc/d\) and equation
\(cY-dh=0\). This proves regularity in codimension one, and perfectness
already gives the Cohen--Macaulay condition in Serre's criterion.

We have \(S\subset\operatorname{Frac}(T)=\operatorname{Frac}(B)\).
Indeed a homogeneous \(s\in S\) outside \(\operatorname{Frac}(T)\)
would give a direct sum \(T\oplus sT\subset S\), of leading Hilbert
coefficient \(4/n!\), whereas the coefficient of \(S\) is \(3/n!\).
The ring \(S\) is finitely generated because \(2L\) is ample Cartier
(the odd graded summands form a finite section module over the even
Veronese ring). The denominators of finitely many generators therefore
give a finite set of height-one primes of \(B\) outside which every
element of \(S\) is regular. A general member
\[
      e=\lambda(b-t_0a)+\mu(d-t_0c)
\]
avoids all of these primes because the common zero locus of the entire family,
as \(t_0,\lambda,\mu\) vary, is \(a=b=c=d=0\), of codimension two.
Normality, or the intersection description of \(B\) by its height-one
valuation rings, now gives, for \(0\le j\le m\),
\begin{equation}
     e^j S_{m-j}\cap B_m=e^j B_{m-j}.                   \label{cubics:pole-equality}
\end{equation}

Consider \(t=b/a\) as a rational function on \(Z\). On the open
Cartier locus where \(a\ne0\) this function is defined and dominant
onto an open subset of the line. For a general value \(t_0\), choose
a prime divisor \(F\) in the closure of its fiber. Replace the variables
\(b,d,Y\) by \(b-t_0a,d-t_0c,Y-t_0h\) and use the same letters for
the new variables. This invertible triangular change preserves the form
of the matrix after expressing \(h\) in the changed variables. All three
new sections vanish along \(F\) with order at least one. Give these
three generators weight one and give the remaining generators weight
zero. The filtration induced on \(B\) from the span of monomials of
weight at least \(j\) is bounded above by section order along \(F\).

Let \(h_0\) be \(h\) with the new \(b,d\) set to zero. The lowest-weight
forms of the three minors are those of the same matrix with \(h_0\) in
place of \(h\). Their ideal still has height two, since \(ad-bc\) is prime
and does not divide \(aY-bh_0\) and maximal minors of this size have height at most two. By Hilbert--Burch these initial minors
have the same ordinary Hilbert series as the original minors. Since
the initial ideal contains them and preserves ordinary Hilbert series
in each finite degree, they generate the entire initial ideal. Its
bigraded Hilbert series, with \(z\) for degree and \(u\) for weight, is
\[
 H(z,u)=
 \frac{1-u z^2-2u z^3+(u+u^2)z^4}
 {(1-z)^n(1-zu)^2(1-z^2u)}.
\]
The syzygies of the minors have degree four and respective weights one
and two, as given by the two rows of the matrix.

The ordinary series \(H(z,1)\) has leading term
\(\tfrac52(1-z)^{-n-1}\). Its numerator equals
\((1-z)^2(1+2z+2z^2)\) and the \(u\)-derivative of the numerator at
\(u=1\) is \(-z^2-2z^3+3z^4=-4(1-z)+O((1-z)^2)\).
The logarithmic derivative of the inverse denominator contributes
\((2+\tfrac12)/(1-z)+O(1)\), giving
\[
   \partial_u H(z,1)
   =\frac{17}{4}(1-z)^{-n-2}
      +O((1-z)^{-n-1})\qquad (z\to1).
\]
The only other possible pole, at \(z=-1\), has order at most two.
It follows that an adapted basis of \(B_m\) has total weight
\(\tfrac{17}{4}m^{n+1}/(n+1)!+O(m^n)\). Write \(G^jB_m\) for this
monomial filtration, and \(\mathcal F^jS_m\) for the filtration by
vanishing to order at least \(j\) along \(F\). The general \(t_0\) may
be taken as in~\eqref{cubics:pole-equality}, so \(e\) also vanishes along \(F\). Hence
\(G^jB_m+e^j S_{m-j}\subset\mathcal F^j S_m\) for \(1\le j\le m\).
By~\eqref{cubics:pole-equality} this yields
\[
 \dim\mathcal F^j S_m\ \ge\ \dim G^j B_m+
     \dim S_{m-j}-\dim B_{m-j}.
\]
For all \(j>m\) we still have \(G^j B_m\subset\mathcal F^j S_m\).
The difference of the two Hilbert functions has leading term
\(\tfrac12 m^n/n!\). Summing over \(j\ge1\) thus improves the total
weight lower bound on \(S_m\) to
\(\tfrac{19}{4}m^{n+1}/(n+1)!+O(m^n)\). Since
\(\dim S_m=3m^n/n!+O(m^{n-1})\), dividing by \(m\dim S_m\) and
passing to the limit over even \(m\) gives
\[
        S_L(F)\ \ge\ \frac{19}{12(n+1)} .
\]
Here \(S_L(F)\) is the limiting maximal mean vanishing of a basis of
\(\lvert mL\rvert\), normalized by \(m\) \cite[Corollary 3.6]{ar1706.04548}. As \(F\) is a prime divisor
on \(Z\), its log discrepancy is one. For \(n\ge5\),
\((n-1)19/[12(n+1)]>1\), contradicting the valuative inequality
\((n-1)S_L(F)\le A_Z(F)\) for K-semistability. This eliminates the
last quadric and proves the proposition.
\paragraph{Uniformity.} For smooth Fano varieties, K-stability implies
uniform K-stability by \cite[Theorem~1.6]{ar2102.09405}. This gives the
cubic row in the introduction.

% ===== The quartic fivefold =====
\chapter{The quartic fivefold}\label{ch:fivefold}
\section{A tangent-space criterion for a quartic fivefold}

For a smooth quartic \(X=(h=0)\subset\mathbf P^6\), write \(b(p)\) for the multiplicity at \(p\) of the tangent hyperplane section \(X\cap T_pX\). The following sufficient condition uses just derivatives of \(h\).
\[
 \boxed{\quad b(p)\in\{2,4\}\quad\hbox{for every }p\in X
       \quad\Longrightarrow\quad X\ \hbox{is K-stable}.\quad}
\]
Equivalently, wherever the quadratic part of \(h|_{T_pX}\) vanishes, its cubic part must also vanish. We give the argument because it isolates a specific set of points at which an individual equation can escape this test.

For a subvariety \(Z\), let \(\delta_Z(X,H)\) be the infimum of \(A_X(v)/S_H(v)\) over valuations whose center contains \(Z\) \cite[Definition~2.5]{ar2003.13788}. First suppose \(b(p)=2\). A general linear surface section \(S\) through \(p\) is smooth, and \(B=S\cap T_pX\) is an integral plane quartic with multiplicity two at \(p\). Indeed the Gauss map is finite, since it pulls the dual hyperplane bundle back to the ample bundle \(3H\). Hence \(X\cap T_pX\) has isolated singularities and is integral. Projection of this tangent section from \(p\) is generically finite because \(b(p)<4\). Successive applications of Bertini irreducibility \cite[Theorem~6.3(4)]{Jou83} to hyperplanes pulled back from the projection give the integral curve \(B\). Its multiplicity is two since a general tangent plane retains a nonzero quadratic term. On \(S\), every other curve \(C\) through \(p\) satisfies \(H\cdot C=B\cdot C\ge2\operatorname{mult}_p C\), and \(B\) itself has ratio \(4/2=2\). Thus the Seshadri constant is exactly two. Abban and Zhuang's surface-section inequality \cite[Theorem~3.1]{ar2101.09246} gives \(\delta_p(X,H)\ge(6/4)2=3\). In case of equality, their alternatives~(2) and~(3) require that a computing valuation have center of dimension at least three or be induced by a prime divisor on \(X\) \cite[Theorem~3.1(2)--(3)]{ar2101.09246}.

If \(b(p)=4\), the tangent section is a cone over a smooth quartic threefold. Smoothness of the base follows again from isolation of the singularities of the tangent section. The base is an index-one K-stable hypersurface by Fujita's theorem \cite[Corollary~1.4]{ar1606.08261}. The generalized Eckardt point formula of Abban and Zhuang \cite[Theorem~4.18]{ar2003.13788} then gives \(\delta_p(X,H)=30/8>3\) for the hyperplane polarization. Moreover the locus where \(b(p)\ge3\) is finite. Indeed the differential of the Gauss map vanishes there, and its restriction to any curve cannot vanish identically in characteristic zero, since the Gauss map is finite. Together these arguments exclude \(A_X(v)\le3S_H(v)\) for centers of dimension zero, one, or two under the displayed hypothesis.

For completeness take an irreducible center \(Z\) of dimension at least three. We show that three general points of \(Z\) lie on a smooth plane quartic section of \(X\). Write \(L=\langle Z\rangle=\mathbf P^l\). A smooth quartic fivefold contains no \(\mathbf P^3\), since its three normal cubic derivatives would have a common zero there. Thus \(4\le l\le6\). For \(l<6\), \(G=X\cap L\) has singular locus of dimension at most \(5-l\), since it consists of points with tangent hyperplane containing \(L\) and the Gauss map is finite. For \(l=6\) take \(G=X\).

The planes of triples from \(Z\) fill \(L\). Indeed joining a proper secant variety to \(Z\) increases its dimension by at least two or fills \(L\). Otherwise the join would equal the cone over the previous secant variety with vertex \(p\) for every general \(p\in Z\), and its vertex would span \(L\). On the open set of triples \((p,q,r)\) with independent representatives, send coefficients \([a:b:c]\) to \([ap+bq+cr]\). On \(abc\ne0\) the image of the differential is the span of the three projective tangent spaces to \(Z\), independently of the nonzero coefficients. Dominance and generic smoothness \cite[Corollary~III.10.7]{Hart77} show that the differential is surjective on this whole coefficient torus for a general triple. The inverse image of \(G_{\rm sm}\) there is smooth. Generic smoothness of its projection to the triples gives smooth general plane sections on that torus. The locus mapping to \(G_{\rm sing}\) has codimension at least three, hence misses a general two-dimensional coefficient fiber.

For the sides of the triangle, a general chord \(pq\) is not contained in \(G\), since otherwise \(Z\subset T_pG\). When \(G_{\rm sing}\) is nonempty, the chords fill \(L\), because then \(l\le5\) and \(\dim Z\ge3\). The same differential and dimension argument for two points shows that a general chord misses \(G_{\rm sing}\). A general \(r\) avoids the finitely many tangent hyperplanes to \(G\) at \(G\cap pq\). Consequently the plane section is smooth along that side, and symmetry treats the other two sides. This proves the assertion about the plane quartic.

Take general nested hyperplanes through this plane to form a smooth flag. The plane quartic \(B\) is smooth, so transversality holds along \(B\), and Bertini applies elsewhere. Hyperplane restriction is surjective in every nonnegative degree by the hypersurface restriction sequence. For four successive orders \(a_1,\ldots,a_4\), the terminal series is \(H^0(B,kH)\), where \(k=m-\sum a_i\), of dimension \(4k+O(1)\). The limiting density on the simplex of normalized orders is proportional to \(1-\sum_i a_i/m\). Its mean is \(1/6\) for each \(a_i/m\) and \(1/3\) for \(k/m\), giving terminal mean degree \(4/3\). At a fourth very general point of that curve, a basis compatible with its order uses asymptotically half the degree by Riemann--Roch \cite[Theorem~IV.1.3]{Hart77}. The three marked points can therefore share at most \(4/6\) of the normalized order, so at least one has coefficient at most \(2/9\) in the limit, as in Hu's curve estimate \cite[Lemma~3.5]{ar2610.03311}. Multipoint adjunction \cite[Theorem~2.26]{ar2610.03311} gives \(\delta_Z(X,H)\ge\min\{6,9/2\}=9/2\). This excludes the remaining centers. Since \(-K_X=3H\), the valuative criterion \cite[Theorem--Definition~2.3(2)]{ar2007.14320} gives the claimed K-stability.

We give a proof of K-stability for every smooth quartic fivefold \(X\subset\mathbf P^6\). The preceding tangent-space criterion covers quartics with tangent multiplicities two and four. For an arbitrary smooth quartic fivefold, a valuation with positive-dimensional center of dimension at most two can still be tested at a point of tangent multiplicity two because the locus of multiplicity at least three is finite. The proof of the same criterion treats all centers of dimension at least three without any tangent multiplicity hypothesis. Thus arbitrary smooth quartics require just one more calculation, at a closed point of tangent multiplicity three.

\section{Order-three tangent sections}

\textbf{Local claim.}
Let \(X\subset\mathbb P^6\) be a smooth quartic, let \(H=\mathcal O_X(1)\),
and let \(p\) be a point at which the tangent hyperplane section has
multiplicity three. The argument gives
\(A_X(v)>3S_H(v)\) for every divisorial valuation with center exactly
\(p\). Here \(3H=-K_X\).

Write the projective coordinates as \([T:u_0:\cdots:u_4:z]\), with
\(p=[1:0:\cdots:0]\) and tangent hyperplane \(z=0\). The equation is
\[
 zB(T,u,z)+T f(u)+g(u)=0,\qquad B(1,0,0)=1,
\]
where \(B,f\) are cubic and \(g\) is quartic. The tangent section has
finite singular locus since its singular points lie in a fiber of the
Gauss morphism of \(X\), whose pullback of the hyperplane bundle is
\(\mathcal O_X(3)\). In particular \(f,g\) are coprime. A common
factor would give intersecting components, or a multiple component, in
the tangent section and a positive-dimensional singular locus.

For an extraction \(\pi:Y\to X\) with unique exceptional divisor \(E\)
and \((Y,E)\) purely log terminal, write \(W_{m,J}\) for the image of
sections of \(mH\) of \(E\)-order at least \(J\) in the quotient by
those of order at least \(J+1\). Its elements give sections on \(E\).
We count every slice with its dimension and divide all sums of
vanishings by \(m N_m\), where \(N_m=h^0(X,mH)\). The resulting
basis divisors on \(E\) are those of Abban and Zhuang
\cite[Definitions~2.18--2.19]{ar2003.13788}. With this normalization \(S_H(E)\)
is the limiting mean of \(J/m\) \cite[Lemma~2.21]{ar2003.13788}.
We apply their filtered adjunction inequality
\cite[proof of Theorem~3.3, eqs.~(3.5)--(3.6)]{ar2003.13788}.
If the asymptotic threshold on \((E,\operatorname{Diff}_E)\) of \(W\)
is at least \(3\), then for a divisorial \(v\) centered at \(p\),
\[
 A_X(v)-3S_H(v)\ \ge\
       \bigl(A_X(E)-3S_H(E)\bigr)v(E).
\]
Indeed choose a basis simultaneously adapted to \(v\) and to \(E\)
\cite[Lemma~3.1]{ar2003.13788}.
Its pullback splits into its mean \(E\)-order times \(E\) and a
residual divisor whose restriction is a basis divisor of \(W\).
For any fixed coefficient less than three, inversion of adjunction
bounds the value of this residual divisor by
\(A_Y(v)-v(E)\), after multiplication by that coefficient and passage
to high levels \cite[proof of Theorem~3.3, following eq.~(3.5)]{ar2003.13788}.
Use \(A_X(v)=A_Y(v)+(A_X(E)-1)v(E)\) and take the
coefficient up to three. Every such \(v\) has \(v(E)>0\).

\subsubsection*{The ordinary blowup and its residual spaces}

For the ordinary point extraction the initial sections \(T,u,z\)
have orders \(0,1,3\), respectively. Write \(B=T^3+T^2 B_1+T B_2+B_3\),
where \(B_i\) is homogeneous of degree \(i\) in \(u,z\).
The section \(w=(T^2+T B_1+B_2)z+f\) has degree three.
The equation becomes \(Tw+B_3z+g=0\).
The initial relations are \(T^2z+f=0\) and \(Tw+g=0\).
In the affine chart \(T=1\), \(w\) has order four since \(B_3z\)
has order at least six. Consequently the associated graded ring is
\[
 \mathbb C[T,u_0,\ldots,u_4,Z,U]/(T^2Z+f,\ TU+g).
\]
Here \(Z,U\) have section degrees \(1,3\) and point orders \(3,4\).
To verify the presentation, the two equations and \(T\) form a regular
sequence, in either order, since modulo \(T\) the forms \(f,g\) are
coprime. After inverting \(T\) the graph is a Laurent polynomial
ring, and its map to the corresponding associated graded ring is
injective because \(u/T\) are regular parameters at \(p\).
Regularity of \(T\) gives injectivity before localization. The
section-degree Hilbert series is \((1-t^4)/(1-t)^7\), the series of
\(X\), which gives surjectivity in every degree.

Put \(R=\mathbb C[u_0,\ldots,u_4]\), with \(R_b=0\) for \(b<0\).
On \(E=\mathbb P^4\), the sections in fixed bidegree are precisely
\[
 W_{m,J}=\sum_{2i+j\ge J-m} f^i g^j R_{J-3i-4j}.
\]
For this description take \(T=1\) and use nonnegative exponents of
\(T,Z,U\) in the graph presentation. The section ring has a basis
\(T^a Z^i U^j P_l\), where \(P_l\) runs through a homogeneous
complement to \((f,g)\) in \(R_l\).
Indeed graded lifts of a homogeneous basis of \(R/(f,g)\) give
a free basis of \(R\) over \(\mathbb C[f,g]\), by degree induction
and the regular-sequence Hilbert series. In the graph, \(f,g\)
can be replaced by \(-T^2Z,-TU\). Thus
\[
 m=a+i+3j+l,\qquad J=3i+4j+l.
\]
The Hilbert series of \(R/(f,g)\) is a polynomial divided by
\((1-t)^3\). Hence \(l\) counts three asymptotically free
degree-one indices. The normalized indices \(a,i,3j\) and these
three indices have the uniform measure on the simplex with total
one in the limit. This follows by ordinary lattice-point Riemann
sums, or by writing the Hilbert series with nonnegative numerator
\((1+t+t^2)(1+t+t^2+t^3)\). Each simplex coordinate has mean \(1/6\).
It follows that \(S_H(E)=11/9\), whereas \(A_X(E)=5\).

We will specialize the pair by rescaling coordinates. Whenever the
initial pair stays coprime, its graph presentation has the same
bigraded Hilbert series. The initial subspaces contain the products
defining the specialized spaces and equality follows from their
dimensions. In each bidegree a basis may be made simultaneously
adapted to a fixed target valuation and to the weight filtration
\cite[Lemma~3.1]{ar2003.13788}.
Its independent initial forms are a basis of the specialized space.
Lower semicontinuity of log canonical thresholds therefore transfers
uniform lower bounds from the specialized basis divisors back to the
original spaces \cite[Main Theorem~0.2(3)]{armath-9910118}. For successive degenerations apply this statement
one weight at a time. We use the same argument locally for an action
fixing a point where the threshold is being tested.

\subsubsection*{A bound away from the cubic}

Fix a point \(x\in\mathbb P^4\) with \(f(x)\ne0\). A general line
through \(x\) misses \(V(f,g)\), since this base locus has dimension
two. Let the line be the binary coordinate space and let
\(v_1,v_2,v_3\) be complementary linear coordinates. Scaling all
three to zero fixes \(x\) and gives coprime binary \(f,g\). The three
complementary powers can now be split by a diagonal torus, giving
three coordinate hyperplanes of mean coefficient \(1/6\) each.
For fixed complementary exponents of total \(k\), set
\(D=J-k\) and \(h=\max(0,J-m)\). The binary space is
\[
 \sum_{2i+j\ge h} f^i g^j \mathbb C[s,t]_{D-3i-4j}.
\]
In any nonzero such space the least possible \(f\)-exponent is
\(b=\max(0,\lceil(4h-D)/5\rceil)\).
For a given \(i\) one needs
\(3i+4\max(0,h-2i)\le D\), which gives this formula and the
usual feasibility condition at \(i=\lceil h/2\rceil\).

After removing \(f^b\) this space has codimension at most three in
the full binary space of degree \(D-3b\). Filter the binary forms
successively by \(f^i\) from \(i=b\) through \(\lceil h/2\rceil\).
Before the last index, the coefficients modulo \(f\) contain
\(g^{h-2i}\) times the full degree \(D-3i-4(h-2i)\) space.
Multiplication by \(g\) is invertible on the length-three scheme
\(f=0\) after local trivialization, and degrees at least two
surject onto this scheme. The coefficient degree increases by
five at each step. Only the first layer can therefore miss up to
three dimensions. At \(i=\lceil h/2\rceil\) all binary multiples
of \(f^i\) are allowed. This argument includes repeated binary
roots and the case in which the last layer is zero.

In the binary free basis, the residual degree in a complement to
\((f,g)\) is bounded. Thus \(D=3i+4j+O(1)\) and
\(b=(i-4a/5)_++O(1)\) in the preceding simplex notation, with
\(t_+=\max(0,t)\). The normalized means of \(D,b\) are \(13/18,5/54\).
For the second, if \(I,A\) are independent unit exponentials, then
\(\mathbb E (I-4A/5)_+=5/9\). Normalizing six such exponentials
to sum one gives the stated division by six, by independence of
the total and the normalized vector. At \(x\) the factor \(f^b\)
is a unit. A codimension-at-most-three subspace of binary degree
\(D-3b\) has mean vanishing at a point at most
\((D-3b+3)/2\), by counting the orders of a basis. Its limiting
contribution is at most
\[
 \frac12\left(\frac{13}{18}-3\frac{5}{54}\right)=\frac29.
\]
The binary divisors through \(x\) all have the same linear
equation there. Together with the three complementary coordinate
hyperplanes they form a simple normal crossings divisor locally.
Consequently the localized asymptotic threshold off \(V(f)\) is at
least \(9/2\). In particular, in any coprime pair a prime coordinate
hyperplane with generic point outside \(V(f)\) has mean coefficient
at most \(2/9+o(1)\), uniformly over bases. One can obtain this
last conclusion by testing its valuation, or by taking a general closed point on it. For the latter argument, for each trial valuation choose a basis compatible with the first weight filtration. On the specialized binary model use the valuation and the diagonal action separating the three complementary exponents. In each resulting block every divisor has the stated binary degree after removal of \(f^b\). The basis average at the chosen point is bounded by the codimension estimate, independently of the choice of basis. The two specializations may thus be applied successively at each level.

\subsubsection*{The flex degeneration for an irreducible cubic}

Suppose the cubic \(f\) is irreducible. We prove that the residual series
\[
 W_{m,J}=\sum_{2i+j\geq J-m}f^ig^jR_{J-3i-4j}
\]
on \(\mathbf P^4\) satisfies \(\delta(\mathbf P^4;W)\geq3\).
All averages below count the individual basis elements of all the
\(W_{m,J}\) at section degree \(m\).

Choose a general plane \(\Pi\subset\mathbf P^4\). Its cubic section
\(B\) is integral by Bertini's theorem
\cite[Tags 0G4F and 0FD6]{Stacks}. The complete intersection \((f=g=0)\) has
dimension two, so its intersection with \(\Pi\) is finite. Scaling the
two coordinates complementary to \(\Pi\) gives as initial forms the
plane restrictions of \(f\) and \(g\), which are coprime.

Every integral complex plane cubic has a smooth flex. In fact, after
choosing a smooth base point, the smooth locus is the degree-zero
generalized Jacobian. In the smooth, nodal and cuspidal cases this group
is an elliptic curve, the multiplicative group and the additive group,
respectively \cite[Ch.~II, \S\S1, 3, in particular pp.~56--58]{Milne06}.
One can see the identification with \(\operatorname{Pic}^0(B)\) by
twisting a degree-zero line bundle by the chosen base point.
Riemann--Roch and \(\omega_B\simeq\mathcal O_B\) give the resulting
degree-one bundle a unique section up to scalar. Its zero divisor is
a smooth point, since an effective Cartier divisor of length one has
regular local ring in dimension one. Multiplication by three is
surjective in each of the three groups. Hence
there is a smooth point \(p\) with
\(\mathcal O_B(3p)\simeq\mathcal O_B(1)\). Since the linear forms on
\(\Pi\) restrict isomorphically to \(H^0(B,\mathcal O_B(1))\), some
line cuts out \(3p\).

Choose plane coordinates \(u,y,z\) with \(p=[0:0:1]\) and flex tangent
\(y=0\). The coefficient of \(yz^2\) in \(f\) is nonzero by smoothness.
The restriction \(f(u,0,z)\) is a nonzero multiple of \(u^3\). Thus,
after rescaling, the initial form for weights \((1,3,0)\) is
\(f_0=yz^2-u^3\).
We can arrange that the initial quartic is a monomial. First observe that
replacing \(g\) by \(g+fL\) for a linear form \(L\) preserves every
\(W_{m,J}\). In the binomial expansion each replacement increases
\(2i+j\) by one and preserves ordinary degree, and the inverse
substitution gives the reverse inclusion. Now evaluation at
\((t,t^3,1)\) on homogeneous quartics has kernel \(f_0R_1\). It sends
every weight-\(r\) monomial to \(t^r\). If the lowest weight part of
\(g\) belongs to \(f_0R_1\), subtract the corresponding \(fL\) and
repeat. Its weight strictly increases, and \(g\) is nonzero modulo
\(f\), so the process stops by weight twelve. In the first nonzero
weight space modulo \(f_0\), subtracting one further \(fL\) if necessary
leaves a scalar multiple of a single monomial \(M\).

The limiting pair \((f_0,M)\) is coprime. At each of the two
specializations the graph presentation therefore retains its bigraded
Hilbert function. The initial space contains all the products in the
limiting generators that define its counterpart \(W_{m,J}\), and equality
follows by dimension. For each valuation and each \(W_{m,J}\), choose
a basis adapted simultaneously to the valuation and the weight
filtration. Its initial forms constitute a basis of the limiting space.
Lower semicontinuity of log canonical thresholds for the resulting
divisors on projective space \cite[Main Theorem~0.2(3)]{armath-9910118} transfers a lower bound
on all limiting basis thresholds to the original series. We apply this argument at each
specialization.

Interchange \(y\) and \(z\), so that the limiting cubic is \(y^2z-u^3\).
Use the weights \((2,3,0)\) on \(u,y,z\) and scale the two complementary
coordinates \(v,w\) independently. For fixed \(m\), a generic integral
combination of these actions separates the finitely many distinct
characters in each piece. The simultaneous basis argument reduces the
threshold estimate to bases of torus eigenvectors. Such a vector is
\(v^aw^b\) times a homogeneous plane form of ordinary degree \(D\) and
weight \(q\). Differences of exponent vectors with the same degree and
weight belong to \(\mathbf Z(3,-2,-1)\). Factoring a binary polynomial
in \(u^3\) and \(y^2z\) therefore factors the plane form as
\[
 u^\alpha y^\beta z^\gamma
 \prod_{j=1}^{e}(y^2z-c_j u^3),\qquad c_j\in\mathbf C^* .
\]

Put \(K=\operatorname{wt}(M)\) in the new weights, so that
\(0\leq K\leq12\). Write \(\mathbb E_m\) for the average over the
degree-\(m\) basis. The graph character gives
\[
 \frac{\mathbb E_m q}{m}\longrightarrow\frac32+\frac{K}{72},
 \qquad
 \frac{\mathbb E_m a}{m},\frac{\mathbb E_m b}{m}\longrightarrow\frac16,
 \qquad
 \frac{\mathbb E_m D}{m}\longrightarrow\frac89.
\]
Indeed the pairs consisting of section degree and weight for the graph
variables are \((1,2),(1,3)\), four copies of \((1,0)\), and
\((1,6),(3,K)\). The relations have pairs \((3,6),(4,K)\).
Logarithmic differentiation of the character in the weight variable at
one gives a factor whose leading term at section parameter \(t=1\) is
\[
 \frac{2+3+6+K/3-6/3-K/4}{1-t}
 =\frac{9+K/12}{1-t}.
\]
The unweighted Hilbert series has a pole of order six at \(t=1\).
Comparing coefficients gives the first limit. The character factors
with a root of unity different from one contribute only lower-order
terms. Each independent complementary coordinate contributes \(1/6\)
by the same calculation. Subtracting the two contributions from the
normalized mean total degree \(11/9\) gives \(8/9\).

Each individual cubic factor has weight six. Its coefficient in the
normalized basis divisor is therefore at most
\(\mathbb E_m q/(6m)\leq5/18+o(1)\), uniformly in its parameter and in
the basis. The three coordinate hyperplanes \(u=0,y=0,z=0\) have
generic points outside the limiting cubic. The preceding off-cubic
estimate applied to their divisorial valuations bounds each coefficient
by \(2/9+o(1)\), including when \(M\) contains that coordinate.

On the plane, distinct members \(y^2z=c u^3\) meet only at
\([0:0:1]\) and \([0:1:0]\), and they meet the coordinate axes only
at those points. Away from them the cubic components are smooth and
the coordinate axes have simple normal crossings. At \([0:0:1]\) the
weighted blowup with weights \((2,3)\) in \(u/z,y/z\) has log
discrepancy \(2+3=5\), by the toric discrepancy formula
\cite[Proposition~7.2]{ar1706.04548}. The exceptional multiplicity of the normalized
divisor tends to \(3/2+K/72\leq5/3\). Each distinct strict transform of
a cubic meets the exceptional curve transversely at its own torus point,
as seen from the coordinate \((y/z)^2/(u/z)^3\). The strict coordinate
axes meet the other two points transversely on the smooth local
uniformizing charts, which are \'etale in codimension one.

At \([0:1:0]\) use the weighted blowup with weights \((1,3)\) in
\(u/y,z/y\). For a plane form as above its local weight is \(3D-q\).
The exceptional multiplicity of the normalized divisor thus tends to
\[
 3\cdot\frac89-\left(\frac32+\frac K{72}\right)
 =\frac76-\frac{K}{72}<\frac43.
\]
The log discrepancy here is four. The strict transforms of the curves
\(z/y=c(u/y)^3\) again meet the exceptional curve transversely at
distinct ordinary points and the coordinate axes meet its two remaining
coordinate points. These charts prove that for every fixed
\(\lambda<3\) all the resulting normalized plane basis divisors are
log canonical with coefficient \(\lambda\) for sufficiently large \(m\).

In \(\mathbf P^4\) the plane forms define cones with common vertex
\(u=y=z=0\). Blowing up this line has log discrepancy three. The
exceptional multiplicity tends to \(8/9\), and the induced divisor on
each exceptional plane is precisely the plane divisor treated above.
In local blowup charts the radial parameter and the coordinate along
the vertex are independent of the two direction coordinates. At most
one of \(v,w\) vanishes at a given vertex point, and its normalized
coefficient tends to \(1/6\). These charts, together with the projection
to \(\Pi\) away from the vertex, give the same lower bound on
\(\mathbf P^4\). Passing to the limit and reversing the specializations
proves \(\delta(\mathbf P^4;W)\geq3\).

For the ordinary blowup we have \(A_X(E)-3S_H(E)=5-3(11/9)=4/3\). The retained
extraction inequality gives the stated strict bound for irreducible cubics.

\subsubsection*{A linear factor in the cubic tangent form}

Suppose \(f=L^s Q(y,L)\), with \(s\in\{1,2,3\}\) maximal, and
\(y=(y_0,\ldots,y_3)\). The quartic \(G(y)=g(y,0)\) is nonzero
by coprimality. Its surface in \(\mathbb P^3\) has only isolated
singularities. To see this, suppose an irreducible complete curve
\(C\) lies in its singular locus. At \([y:0]\) for \(y\in C\), the
gradients of \(f,g\) have only their \(L\)-component possibly nonzero.
If \(\partial_L f\) is generically nonzero there, choose \(T\) at
each such point to cancel \(T\partial_L f+\partial_Lg\). This gives
a curve in the singular locus of the tangent section. If
\(\partial_Lf\) vanishes identically on \(C\), a zero of the cubic
\(\partial_Lg\) on \(C\) gives an entire singular line as \(T\) varies.
Both contradict finiteness. Bertini's theorem gives a smooth quartic
curve on a general plane \cite[Tag 0FD6]{Stacks}. Scaling a complementary coordinate
degenerates \(G\) to the projective cone over this smooth section.
Blowing up the cone vertex gives threshold \(3/4\) for the special
fiber. Lower semicontinuity \cite[Main Theorem~0.2(3)]{armath-9910118} gives
\(\operatorname{lct}(\mathbb P^3,G)\ge 3/4\).

Give the tangent parameters \(y\) weight \(s\) and \(L\) weight
\(s+1\). The weighted blowup is a plt extraction of discrepancy
\(A_X(E)=4s+(s+1)=5s+1\) by the toric discrepancy formula
\cite[Proposition~7.2]{ar1706.04548}. The plt property follows on the weighted
blowup charts: the smooth covers have a single coordinate divisor
with coefficient one, and their maps are crepant for the pair
\cite[Proposition~5.20]{KM98}. The equation gives \(z\) weight \(4s\) and
initial homogeneous relation
\[
 T^3 Z+T L^s Q(y,0)+G(y)=0 .
\]
Indeed \(zB\) contributes \(T^3 Z\), and the other terms of \(B\)
involving \(u,z\) have positive weight. The displayed polynomial
is primitive linear in \(Z\), hence irreducible, because \(G\ne0\)
prevents \(T\) from dividing its constant term. After inverting
\(T\) its map to the associated graded section ring is injective
by the tangent parameters. The section-degree series is the same as
the original quartic, proving the presentation as before. For the
seven section-degree-one variables \(T,y,L,Z\), the sum of their
weights is \(9s+1\), and the degree-four relation has weight \(4s\).
Differentiating the Hilbert character gives
\[
 S_H(E)=\frac{9s+1-4s/4}{6}=\frac{8s+1}{6},
 \qquad A_X(E)>3S_H(E).
\]

The exceptional divisor has weighted projective coordinates
\(\mathbb P(s,s,s,s,s+1)\). Scaling \(L\) to zero in the polynomial
representing \(Z\) leaves \(G\). The initial relation \(T^3Z+G\)
is still irreducible and retains all bidegree dimensions. Split the
limiting pieces by the exponent \(\beta\) of \(L\). Fixing
\((m,J,\beta)\) leaves on the \(y\) variables the space
\[
 G^b\,\mathbb C[y]_{D-4b},\qquad
 D=\frac{J-(s+1)\beta}{s},\qquad
 b=\max\{0,\lceil(D+\beta-m)/3\rceil\}.
\]
Only nonzero pieces are counted. Indeed a monomial
\(T^a L^\beta Z^\gamma P_l(y)\), for \(P_l\) in a homogeneous
complement to \((G)\), has
\(m=a+\beta+\gamma+l\) and \(D=4\gamma+l\).
At fixed \((m,J,\beta)\) the minimum possible \(\gamma\) is \(b\)
and its multiples fill the indicated complete polynomial space.
In the free basis \(b=(\gamma-a/3)_++O(1)\).

The six asymptotic simplex parameters are \(a,\beta,\gamma\) and
the three unrestricted indices in the Hilbert series of
\(\mathbb C[y]/(G)\). Consequently the limiting means divided by
\(m\) of \(\beta,D,b,D-4b\) are respectively
\[
 \frac16,\qquad \frac76,\qquad \frac18,\qquad \frac23 .
\]
For the third value the positive part of the difference of independent
unit exponentials \(\Gamma-A/3\) has mean \(3/4\). For a base valuation
\(v\) on \(\mathbb P^3\), the fixed \(G\) factors contribute at most
\((1/8)(4/3)A(v)\) by the preceding threshold. The complete systems
contribute at most \((2/3)A(v)/4\), giving a total of \(A(v)/3\).
For the standard three-simplex, the barycenter formula of Blum and
Jonsson gives \(\delta(\mathbb P^3,\mathcal O(1))=4\)
\cite[Corollary~7.16, eq.~(7.3)]{ar1706.04548}. Their uniform bound for
finite-level expected vanishing applies to every moving degree
sufficiently large \cite[Corollary~2.10]{ar1706.04548}.
Explicitly, for large moving degrees \(b'\) the maximal normalized
mean for any \(v\) is at most \((1+\epsilon)A(v)/4\).
Bounded \(b'\) contribute \(O(1/m)A(v)\) after our normalization,
using the elementary bound \(\operatorname{ord}_v(P)\le b'A(v)\)
for any degree-\(b'\) polynomial on projective space. This also
follows from \(\alpha(\mathbf P^3,\mathcal O(1))=1\)
\cite[Corollary~7.16, eq.~(7.2)]{ar1706.04548}.

On a chart \(y_i\ne0\), the uniformizing cover is affine four-space,
with three affine base coordinates and the fiber coordinate coming
from \(L\). The cyclic group \(\mu_s\) acts on this fiber coordinate.
On the cover the canonical divisor is the pullback of
\(K_E+\operatorname{Diff}_E\), by adjunction on the usual weighted
blowup charts. Equivalently the different on the coarse space has
coefficient \(1-1/s\) along \(L=0\). The limiting bases split by
the \(L\) exponent into a fiber factor of mean \(1/6\) and a pullback
of the basis divisor on the base just estimated. A base log resolution
and the transverse fiber coordinate give threshold at least three
asymptotically. At the remaining point \(L\ne0,\ y=0\), take the cover
with \(L=1\). The basis polynomials are now homogeneous in the four
\(y\)'s with mean degree \(7/6\), and their direction divisors on
\(\mathbb P^3\) are the same basis divisors used above. Blow up the
origin on this cover. The radial log discrepancy is four and
\(3(7/6)<4\). The total transforms away from the radial component
are products over the direction space, so the base bound proves
the required threshold at this chart as well. In particular, the coefficient bound for the residual divisor, on the base and on both smooth covers, holds for every fixed \(\lambda<3\) at all sufficiently large levels, uniformly over the adapted bases. The finite quotient
and its different preserve log canonicity
\cite[Proposition~5.20]{KM98}. The specialization
argument then gives the threshold at least three for the original
residual spaces. In every linear-factor case the retained extraction term is
\[
 \bigl(A_X(E)-3S_H(E)\bigr)v(E)
 =\left(5s+1-\frac{8s+1}{2}\right)v(E)
 =\left(s+\frac12\right)v(E)>0
\]
for every valuation centered at \(p\). This gives the strict point-centered inequality.

\enlargethispage{2\baselineskip}
\section{Consequence for dimension five}

For every prime divisor over \(X\), the preceding arguments therefore give \(A_X(v)>3S_H(v)\). If \(\delta(X)\le1\), the optimal destabilization theorem of Liu, Xu and Zhuang \cite[Theorem~1.2]{ar2102.09405} would give a prime divisor computing this threshold, with \(A_X(v)\le S_{-K_X}(v)=3S_H(v)\). Thus the argument gives \(\delta(X)>1\), hence uniform K-stability \cite[Theorem~B(ii)]{ar1706.04548} and a K\"ahler--Einstein metric \cite[Theorem~1.6]{ar2102.09405}.

% ===== The quartic sixfold =====
\chapter{The quartic sixfold}\label{ch:sixfold}
\section{Statement and proof strategy}
\begin{theorem}\label{sixfold:thm}
For every smooth quartic $X\subset\mathbf{P}^7_{\mathbf{C}}$, one has $\delta(X,-K_X)>1$. In particular $X$ is uniformly K-stable.
\end{theorem}
Write $H=\mathcal{O}_X(1)$, so that $-K_X=4H$. For a nonzero divisorial valuation $v$ write $A(v)$ for its log discrepancy, $S(v)=S_H(v)$ for its expected vanishing, and $\beta=\delta(X,H)=\inf_v A(v)/S(v)$. For an irreducible $C\subset X$, $\delta_C(X,H)$ tests all valuations whose center \emph{contains} $C$. The following bounds are in $H$-units.
\begin{center}
\begin{tabular}{llc}
\toprule
Curve & Geometry along the curve & Lower bound for $\delta_C$\\
\midrule
Line & normal cubics span at least three dimensions & $14/3$\\
Line & normal cubics span a pencil & $35/8$\\
Conic & proper plane section, distinct residual divisor & $28/5$\\
Conic & proper plane section, doubled conic & $21/4$\\
Conic & contained plane, polar span at least three & $56/13$\\
Conic & contained plane, polar pencil & $35/8$\\
\bottomrule
\end{tabular}
\end{center}
These six cases exhaust integral curves of degree at most two. The weakest bound exceeds $17/4$ by $3/52$. If $\beta\le4$, a divisor computes $\beta$ by \cite[Theorem 1.2]{ar2102.09405}. Section~\ref{sixfold:positive} produces a line or conic in its center if that center has positive dimension. For a point center $p$, normalize $v(\mathfrak m_p)=1$. Section~\ref{sixfold:points} shows that either $A(v)>4S(v)$ or $A(v)\ge56/5$. Section~\ref{sixfold:nadel} uses the table to deduce $1\ge A(v)/16+5/14$ in the latter case. At $A(v)=56/5$ the right hand side is $37/35$.

\paragraph{Basis adjunction conventions.}
Write $R_m=H^0(X,mH)$ and $N_m=\dim R_m$. An $m$-basis divisor is the sum of the divisors of a basis, divided by $mN_m$. Its maximum $v$-order is $S_m(v)\to S(v)$, and $S_m(u)\le(1+o(1))S(u)$ uniformly in $u$ by \cite[Corollary 3.6]{ar1706.04548}. Two filtrations of a finite-dimensional vector space admit a common adapted basis. For an extraction with exceptional divisor $E$, start with a basis adapted to $v$ and to $E$. Its pullback has coefficient $S_m(E)$ at $E$; the residual restricts to bases of the order quotients on $E$. If all these residual basis divisors admit coefficient $c$ on $E$ with its different, inversion of adjunction and the discrepancy formula give in the limit
\[
 A_X(v)\ \ge\ c S(v)+(A_X(E)-cS(E))v(E).
\]
Smaller coefficients and passage to the limit suffice. This is the proof of \cite[Theorem 3.3, (3.6)]{ar2003.13788}. Over a curve $C$ it suffices to work over its generic point and bound the residual basis divisor on the \emph{whole} geometric generic exceptional fiber. Shrink the curve to spread a log resolution; inversion of adjunction and properness give log canonicity near its generic point downstairs. This also tests every larger center containing $C$.

In using a torus specialization to bound arbitrary bases of a residual series, fix a valuation on that series first. On each piece choose a basis adapted to it and to the torus filtration, keeping the coefficient spaces over $\mathbf{C}$, even on a generic fiber. Its valuation sum is maximal. Its initial forms give bases of the limiting weight pieces. A uniform lower bound on the log canonical thresholds of these initial divisors passes to the adapted divisors by lower semicontinuity on the proper fiber. Thus the bound controls every incoming basis. These changes of basis never require simultaneous adaptation to three filtrations.

For a smooth hyperplane flag the refined residual pieces on the next rung are the full spaces of degrees $m-i$ by the hypersurface restriction sequence. In ambient dimension $j$ the mean hyperplane coefficient tends to $1/(j+1)$ and the mean normalized residual degree tends to $j/(j+1)$. Thus a lower threshold $c$ on the next rung lifts to $\min\{j+1,(j+1)c/j\}$; see \cite[Lemma 4.6]{ar2003.13788}. A chain may be used just near the generic point under consideration when the flag members are normal and smooth there.
\section{Lines and conics}
\paragraph{Lines with at least three independent normal polars.}
Along the line use binary coordinates $u,v$ and six linear normal coordinates $z_i$. The first normal part of the quartic is $\sum z_iQ_i(u,v)$, with basepoint-free cubics $Q_i$. In
\[
 0\longrightarrow\mathcal{O}_{\mathbf{P}^1}(-3)\xrightarrow{(Q_i)}
 \mathcal{O}_{\mathbf{P}^1}^{\,6}\longrightarrow V=\bigoplus_{i=1}^5\mathcal{O}(a_i)\longrightarrow0
\]
the $a_i$ are nonnegative and sum to three. The ordinary blowup of the line has discrepancy five. Its order-graded section ring is cut out by $\sum z_iQ_i$: that equation is primitive and prime, and after localization at the generic point of the line its symbols are the ordinary normal parameters with this one linear relation. The Hilbert series $(1-t^4w)/((1-t)^2(1-tw)^6)$ gives ordinary expected order $23/28$.

In degree $m$, order $j\le m$, the symbol space maps into $H^0(\mathbf{P}^1,\operatorname{Sym}^jV(m-j))$, and the map onto this last space is surjective for $m-j\ge2$. Indeed the kernel bundle in the presentation by $\operatorname{Sym}^j\mathcal{O}^6(m-j)$ is $\operatorname{Sym}^{j-1}\mathcal{O}^6(m-j-3)$, whose $H^1$ then vanishes. In a split frame the monomial $z^\gamma$, $|\gamma|=j$, has coefficient multiplicity $m-j+1+\sum a_i\gamma_i$. The last two order slices have total dimension $O(m^5)$ and normalized angular contribution $O(1/m)$ per unit discrepancy, since a fiber polynomial of degree $j$ has order at most $jA$ on projective space. A diagonal specialization in the split generic frame gives coordinate hyperplanes. Their limiting coefficients are $(4+a_i)/28$. Indeed on the simplex with coordinates $(m-j,\gamma_1,\ldots,\gamma_5)/m$ the density is proportional to $x_0+\sum_i a_i x_i$; the second moments of the uniform simplex give the stated barycenters. When the polar span is at least three, $\max a_i\le2$: the number of trivial summands equals the dimension of the constant relations between the $Q_i$, by duality of the displayed sequence. Adjunction gives $\delta_C\ge\min\{140/23,28/6\}=14/3$.

\paragraph{Exceptional lines on a quartic sixfold.}
Let $X=X_4^6$ be smooth, and suppose that along a line $L$ the normal
cubic polars span a basepoint-free pencil. Then
\[
                       \delta_{L}(X,H)\ge\frac{35}{8}.
\]
Here the local threshold ranges over valuations whose centers contain
$L$. Use coordinates $u,v$ along $L$, four light normal coordinates
$x_i$, and two more normals $z_1,z_2$ so the first normal part is
$f(u,v)z_1+g(u,v)z_2$ for coprime binary cubics. On the five-plane
$z_1=z_2=0$ the first nonzero part is $G(u,v,x)$ of normal degree
$2\le h\le4$. The five-plane cannot be contained in a smooth quartic
sixfold (its two transverse cubic polars would have a common zero).
At the generic point of $L$, $x_1,\ldots,x_4,z_1$ are regular
parameters and the weighted extraction with weights $(1,1,1,1,h)$
has log discrepancy $4+h$. The initial ring of the complete series is
\[
 R'=\mathbb C[u,v,x_1,\ldots,x_4,z_1,z_2]/(fz_1+gz_2+G).
\]
The equation is linear and primitive in the heavy variables since
$(f,g)$ has no common factor. The displayed hypersurface ring injects
into the polynomial ring over $K=\mathbb C(L)$ in the independent
normal parameters after retaining the homogeneous base frame (equivalently,
localize in $u,v$ first). The map sending the displayed generators
to the weighted symbols in the section ring is therefore injective
and its ordinary Hilbert series $(1-t^4)/(1-t)^8$ gives equality.
The mean exceptional order is $(4+2h-h/4)/7$.

Work on the generic exceptional weighted projective space
$P=\mathbb P_K(1,1,1,1,h)$, with quotient charts (or its smooth stack)
for testing thresholds. Choose a line in the $x$-space where $G$
restricts to $a(u,v)x^h$, $a\ne0$ of base degree $4-h$. Denote the
three complementary coordinates by $x'$. Give them distinct positive
weights and the binary coordinates $x,z_1,z_2$ weight zero, choosing
weights separating the $x'$ exponents in the finite degree range at
hand. The limiting initial equation is $fz_1+gz_2+a x^h$, still
primitive. Viewing $z_1,x,x'$ as the independent coordinates on $P$
makes this a geometric diagonal degeneration: $z_2$ is replaced by
$-(fz_1+G)/g$ before taking the leading term. The injection and
ordinary Hilbert-series arguments just given show that the initial
spaces are exactly the corresponding pieces of the specialized
ring. In the specialized bases the $x'$ factors split, with mean
exponent $1/7$ for each. The mean weighted degree in the remaining
binary variables is
\[
               B=\frac{1+2h-h/4}{7}=\frac17+\frac h4 .
\]
These means, and the earlier exceptional order, can be checked by
logarithmic differentiation of the graded Hilbert series. A free
generator of ordinary degree one contributes asymptotic mean $1/7$,
and the quartic equation contributes minus one fourth of its weight
divided by seven.

For clarity the binary calculation retains complex coefficient
dimensions even though the threshold is tested over $K$. Put $Q=x^h$
and write each binary piece as $x^s$ times a heavy binary system,
$0\le s<h$. Indeed, adjoining $x$ with $x^h=Q$ is finite free
over
\[
 T=\mathbb C[u,v,x',z_1,z_2,Q]/(fz_1+gz_2+aQ),
\]
with basis $1,x,\ldots,x^{h-1}$. The ordinary degrees of $z_1,z_2,Q$
are $1,1,h$, their binary heavy degrees are all one, and the binary
weighted degree in the old grading is $hk+s$ on heavy degree $k$.
Since $a\ne0$, $z_1,z_2$ form independent homogeneous coordinates on
the coarse binary projective line over $K$; the distinguished point
at infinity is $Q=0$. At any geometric point different from it,
the factor $x^s$ is a unit. In the direct sum of these heavy systems,
the limiting mean $k/m$ is $B/h$.

Choose the heavy coordinates so that $f$ is coprime to $a$, which is
possible using the basepoint-free pencil. In each heavy binary system,
scale $z_2$ against $z_1$ on the coarse binary line. The least-weight
symbol of $Q=-(fz_1+gz_2)/a$ is then $-fz_1/a$; the new ring has
equation $fz_1+aQ$, is primitive, and injects after localization at
$K$. Ordinary Hilbert series, with the $x^s$ decomposition or with
each extra grading, again give all initial spaces. The variable $z_2$
is free, so its limiting mean exponent in the full normalized basis
(summing over all the blocks) is $1/7$. The complementary exponent
of the surviving binary coordinate $z_1$ has limiting mean
\[
              B/h-1/7=(3/4+1/h)/7\le 5/28 .
\]
The diagonal initial basis on the binary line therefore cuts out
only the two reduced coordinate points with these mean coefficients.
For any chosen geometric binary point, take complex bases in the
blocks simultaneously adapted to its order filtration and to this
diagonal filtration. The initials span their respective complex
coefficient spaces; each nonzero coefficient is a unit in the
geometric generic fiber. Semicontinuity of the log canonical
threshold on the binary projective line bounds the maximum normalized
order at the chosen point by $5/28+o(1)$. The bounds by the common
finite-level diagonal coefficients are uniform in that point. This
uses the coarse line only for finite points, where its parameter
is an unramified parameter on the weighted projective stack.

At infinity compute the upstairs order along $x=0$ directly in the
specialized $x'$-split ring, before the second diagonal degeneration.
The quotient by $x$ has equation $fz_1+gz_2$, a primitive prime
polynomial. Its map to the generic binary fiber, with the extra
gradings retained, is injective. Consequently divisibility by $x^i$
in the section ring equals vanishing to order $i$ along $x=0$ in
the affine normal ring over $K$, by induction using the quotient
by $x$. On the stack chart $z_1\ne0$ extracting the $h$-th root
of $z_1$ gives exactly that order at the point at infinity. The
normalized sum of its maxima is the mean of the powers of a free
filtration parameter of ordinary degree one, hence $1/7+o(1)$.
Explicitly $\sum_{i=1}^m\dim(x^i R'_{m-i})=\sum_{i\ge1}\dim R'_{m-i}$
in the specialized ring, and one may sum this identity over the
angular gradings and $x'$ characters first. This explains why
multiplying the coarse-line bound by $h$ is unnecessary.

It remains to check simultaneous intersections in $P$. Over the
algebraic closure a homogeneous binary form splits into powers of
$x$ and of $z_1-\lambda x^h$. Distinct finite-point divisors and
the infinity divisor meet only on the binary axis $x=z_1=0$.
This axis lies on the smooth charts of $P$. Its weighted blowup
in the normal parameters $(x,z_1)$ with weights $(1,h)$ has log
discrepancy $1+h$ and total mean order $B$. On its quotient charts
the strict transforms of all the distinct binary factors are
disjoint, smooth, and transverse to the exceptional divisor and
to the independent coordinate divisors of the $x'$ variables.
For example the first chart has $z_1=x^h T$ and the second chart
has $z_1=s^h,\ x=sU$ before quotient by $\mu_h$. On the second
chart only $U=0$ goes through the infinity point; the finite roots
there have nonzero $U$ and are distinct on the coarse line.
Off the axis this same transversality holds, including the ambient
cyclic point with $z_1\ne0$, as checked on its quasi-\'etale chart.
Thus the limiting threshold of the specialized angular series is
at least
\[
          \min\{7,28/5,(1+h)/B\}\ \ge\ 35/8.
\]
To apply the first semicontinuity step for a test valuation, choose
a complex basis simultaneously adapted to that valuation and the
split-$x'$ filtration (with order slices fixed). Its specialized
binary factors have coefficients bounded by the just computed
maxima at every point, regardless of the choice of the specialized
bases. The resulting lower bound for these adapted basis divisors
therefore also holds prior to degeneration. Filtered adjunction
over the generic point of $L$ on the weighted extraction now takes
the minimum with $7(4+h)/(4+7h/4)$, which is at least $56/11$.
This proves the asserted bound.

\paragraph{Conics with distinct residual plane divisor.}
Let $C\subset\Pi$ be a smooth conic with $X|_\Pi=C+R$, where $R$ is a possibly nonreduced conic divisor having no component $C$. A general linear flag through $\Pi$ gives a quartic surface $S$. It is smooth along $C$: at smooth points of $C+R$ transversality is automatic, and a general remaining direction avoids tangency at the finitely many other points of $C$. For a reduced plane section the same argument applies everywhere. If $R$ is a doubled line, singularities can occur only at finitely many other points on that line, by Bertini and smoothness of $X$. The surface is normal and integral, since a reducible hypersurface surface is singular along a curve. Take the pullback of $H$ to a resolution. The strict transform of $C$ has $C^2=-2$ by adjunction in its smooth neighborhood. Write $R=H-C$ also for the residual effective divisor on this resolution. Its exceptional components are disjoint from $C$. For each other component $B$ of $R$ one has $C\cdot B=2\deg B$ by plane B\'ezout; the intersection scheme is the same on the surface and in the plane.

It follows that $H-tC$ is nef on $[0,1/2]$; use $(1-t)H+tR$ on curves outside $R$. Its volume there is $4-4t-2t^2$. For $t\in[1/2,1]$ the big nef class $2H-C$ is orthogonal to $R$, and the positive part of $H-tC$ is $(1-t)(2H-C)$, with effective negative part $(2t-1)R$. Equivalently concavity of the square root of volume gives the same upper bound $6(1-t)^2$. Intersection against $2H-C$ excludes bigness past one. Consequently
\[
 S_{H|S}(C)\le\frac14\left(\int_0^{1/2}(4-4t-2t^2)\,dt+
                  \int_{1/2}^1 6(1-t)^2\,dt\right)=\frac5{12}.
\]
This controls the only divisor order centered at the generic point of $C$ on the surface. The four hyperplane refinements give $\delta_C(X,H)\ge(7/3)(12/5)=28/5$. The general intermediate threefold and higher rungs are normal and smooth near the generic point of $C$ by the same Bertini and transversality argument.
\paragraph{A plane section equal to $2C$.}

Take a general \(\mathbf P^4\) through \(\Pi\) and let \(Y\) be the induced quartic threefold. The normal polar cubics of \(X\) along \(\Pi\) have no common zero on \(C\), so two general linear combinations have none on this curve. Consequently \(Y\) is smooth along all of \(C\). Every larger flag member containing the chosen four-space is also transverse to \(X\) there. Off \(C\) smoothness of the general flag follows from Bertini. With \(u_0,u_1,u_2\) homogeneous coordinates on \(\Pi\), and \(x,y\) transverse coordinates on the chosen four-space, write the equation of \(Y\) as \(g(u)^2+xA(u)+yB(u)+O((x,y)^2)\). Here \(g=0\) cuts out \(C\), and \(A,B\) are cubic forms without a common zero on \(C\).

At the generic point of \(C\) give \(g\) weight one and the transverse parameter independent of \(g\) weight two. This yields the weighted blowup with weights \((1,2)\) in two smooth transverse parameters, with exceptional divisor \(E\) of log discrepancy three. The initial algebra, where \(G\) represents \(g(u)\), is
\[
 T=\mathbf C[u_0,u_1,u_2,x,y,G]/(g(u),G^2+A(u)x+B(u)y).
\]
Weights are \(0,0,0,2,2,1\), while section degrees are \(1,1,1,1,1,2\). To check the presentation, the displayed equations are the initials of \(g(u)-G\) and \(G^2+xA+yB+O((x,y)^2)\). The quotient is free with basis \(1,G\) over \(R(C)[x,y]\) and is a domain since the nonzero form \(Ax+By\) is nonsquare over \(\mathbf C(C)[x,y]\), also after adjoining a variable recording ordinary section degree. Its ordinary Hilbert series is that of \(Y\). Localizing at a generic frame for the conic plane coordinates identifies the displayed initials with the monomial valuation on the local ring: on a chart with \(A\ne0\), solve for \(x\), and use \(g,y\) as parameters. Thus this presentation computes every section slice.

The coarse exceptional fiber is \(\mathbf P^1_{\mathbf C(C)}\) with coordinates \([x:y]\). The different equals \(\frac12 D\) for the point \(D=(Ax+By=0)\): in the weighted presentation its equation is \(G=0\), with a generic stabilizer of order two. Locally on \(C\) this is the standard \((1,2)\) blowup quotient chart. Differentiating the bigraded Hilbert series gives \(S_Y(E)=[2+2+1/2-2/4]/4=1\). Here the relation \(g=0\) contributes zero to the secondary weight and the relation of degree four has secondary weight two.

For section degree \(m\) and weight \(2l+\epsilon\) with \(\epsilon\in\{0,1\}\), the description as a free module identifies the slice with \(G^\epsilon\) times
\[
 H^0\!\left(C,\mathcal O_C(kH)\right)\otimes H^0(\mathbf P^1,\mathcal O(l)),
 \qquad k=m-l-2\epsilon\ge0.
\]
The conic is projectively normal and \(\mathcal O_C(H)=\mathcal O_{\mathbf P^1}(2)\). The fixed factor \(G^\epsilon\) contributes at most \(1/(2m)\) at \(D\) to normalized basis divisors and zero elsewhere.

Every nontrivial divisorial valuation on the regular curve \(\mathbf P^1_{\mathbf C(C)}\) is proportional to the order at a closed point. Its closure in \(C\times\mathbf P^1\) is an effective prime Cartier divisor of bidegree \((b,a)\), with \(b\ge0,\ a\ge1\). Conic base and exceptional fiber are ordered first and second, respectively. With \(h(t)=\max(t+1,0)\), the complex dimension of the moving part of this slice with order at least \(s\) equals
\[
 h(2k-sb)\,h(l-sa).
\]
This follows by subtracting \(s\) times the closure as a Cartier divisor. There are no additional vertical conditions because a rational section of the original bidegree is already regular everywhere on the product; division gives precisely the regular sections of the subtracted line bundle. In particular this is at most \(h(2k)h(l-s)\), the value for a constant fiber point. Upon summing over \(l,\epsilon\), the last contribution is exactly \(h^0(Y,(m-s)H)\) by reindexing \(l'=l-s\). Summing over \(s\ge1\) and using the cubic-growth Hilbert function gives normalized expected order at most \(1/4+o(1)\), uniformly in the closed point.

For \(D\) its closure is the graph divisor \(Ax+By=0\) of bidegree \((6,1)\), since \(A,B\) have no common zero on the conic. In this case the sum of order-\(s\) dimensions is exactly \(h^0(Y,(m-4s)H)\), by the substitution \(l'=l-s,\ k'=k-3s\). The normalized expected order tends to \(1/16\), including the vanishing contribution of the bounded parity factor. The log discrepancy at this point is \(1/2\); the other closed points have discrepancy one. Therefore the residual series at the generic fiber has threshold at least four. Abban--Zhuang, Theorem 3.3, at the generic point of \(C\), gives \(\delta_C(Y,H)\ge\min(3,4)=3\). Local weighted extraction near the generic point suffices: the restrictions of the absolute complete section series used in the refined filtration are exactly the spaces above. The first-polynomial basis argument in the proof of their theorem applies to every center containing \(C\) by localization.

For each larger hyperplane section, apply the generic-point form of the same theorem to the Cartier hyperplane through \(C\). Its own ratio in ambient dimension \(j\) is \(j+1\), and the complete restricted series carries the scale \(j/(j+1)\). Iteration from the threefold (dimension three) through the sixfold gives \(\delta_C(X,H)\ge\min(7,\frac74\delta_C(Y,H))\ge21/4\). The factor is \(7/4\), the product \((5/4)(6/5)(7/6)\) over the three hyperplane cuts.

\paragraph{The normal polar series of a conic on a contained plane.}
Let $C$ be a smooth conic in a plane $\Pi\subset X_4^6$ and
$H=\mathcal O_X(1)$. The five normal cubic polars on $\Pi$ restrict
to a basepoint-free system of sextics on $C\simeq\mathbf P^1$.
Define $E$ on $\mathbf P^1$ by the exact sequence
\[
 0\longrightarrow\mathcal O(-6)\xrightarrow{\text{normal polars}}
 \mathcal O^{\oplus5}\longrightarrow E\longrightarrow0 ,
 \qquad E=\bigoplus_{i=1}^4\mathcal O(a_i).
\]
Thus $a_i\ge0$ and $\sum a_i=6$. For divisorial valuations $v$
whose center contains $C$ one has
\[
 \frac{A_X(v)}{S_H(v)}
 \ \ge\ \min\left\{\frac{20}{3},\frac{56}{8+\max_i a_i}\right\}.
\]
When the polars span dimension at least three on $C$, this gives
$56/13>17/4$.

Write $y$ for the three plane variables, $x_1,\ldots,x_5$ for the
normal variables, and $g(y)$ for the conic equation. The ideal of
$C$ in the polynomial ring of the ambient affine cone is the
complete intersection $(x_1,\ldots,x_5,g)$. Its associated graded
ring is the symmetric algebra over $\mathbf C[y]/(g)$, with
symbols $x_i,w=\operatorname{in}(g)$ of degree/order $(1,1)$ and
$(2,1)$ respectively. The initial equation of $X$ is the linear
form $\sum x_i A_i(y)$ in the cubic polars. It cuts out an integral
scheme in this algebra. Off the vertex of $\operatorname{Spec}
\mathbf C[y]/(g)$ it is a vector bundle by the absence of common
polar zeros. That locus is dense in the zero scheme since the
fiber over the vertex has dimension six and the zero scheme is
equidimensional of dimension seven. The hypersurface quotient of
the Cohen--Macaulay ring has no embedded primes, giving the
assertion. The initial form of the equation is a nonzerodivisor,
so the resulting quotient is the full associated graded ring for
the ideal filtration; its integrality identifies this order with
the ordinary order at the generic point of $C$.

Fix the $w$ exponent $s$ and the combined $x$ exponent $k$.
In section degree $m$ put $p=m-k-2s\ge0$. The polynomial
coefficient space in this component is the image of
\[
 H^0(\mathcal O(2p))\otimes\operatorname{Sym}^k\mathbf C^5
      \longrightarrow H^0(\operatorname{Sym}^k E(2p)).
\]
For $p\ge3$ this map is surjective, by the symmetric-power
sequence whose kernel bundle is the direct sum of copies of
$\mathcal O(2p-6)$. For a split monomial indexed by
$\beta\in\mathbf Z_{\ge0}^4$, $|\beta|=k$, the coefficient
dimension is $2p+1+a\cdot\beta$. The slices with $p=0,1,2$
have total dimension $O(m^5)$ by the same upper bound, and their
normalized contribution to any fixed-coordinate mean is negligible.

The simplex variables are
$(\beta_1,\ldots,\beta_4,2s,p)/m$, with sum one and density
proportional to $a\cdot\beta/m+2p/m$. Its six vertex coefficients
are $(a_1,a_2,a_3,a_4,0,2)$ of sum eight. The elementary Dirichlet
moments give $\overline{\beta_i}/m=(8+a_i)/56$ and
$\bar s/m=1/14$. Thus the mean ordinary extraction coefficient is
$3/4$, and its codimension-five discrepancy is five. On the
generic exceptional $\mathbf P^4$ a basis can be chosen
simultaneously for $v$ and the lexicographic refinement of the
normal degree by the $w$ and split-frame monomials. Over
$\mathbf C(C)$ specialize those initial forms by a diagonal torus
weight representing the finite lexicographic comparisons. The
divisor coefficients of the resulting coordinate hyperplanes have
the displayed means, and the residual threshold is at least
$56/(8+\max a_i)$ by semicontinuity. The omitted slices have
degree at most $m$ and their proportion is $O(m^{-1})$.
Inversion of adjunction with the exceptional divisor gives the
stated minimum. It applies also when the center of the valuation
strictly contains $C$: for any smaller coefficient, every
valuation-adapted basis has log canonical pullback near the entire
fiber over the geometric generic point of $C$, after increasing
the coefficient of the exceptional divisor to one. Indeed the
angular estimates are uniform for bases adapted simultaneously to
the valuation and the extraction filtration. Inversion of adjunction
and properness imply log canonicity near the generic point of
$C$ downstairs. Thus a larger center containing $C$ is tested too.

Finally $\max a_i=6$ implies three trivial summands. Taking
sections of the dual kernel sequence identifies that number
with the dimension of constant linear relations on the five
polars. Consequently the exceptional case is precisely a pencil
of polar sextics.

\paragraph{The polar-pencil conic on a contained plane.}
Suppose the five polar cubics of the plane restrict to a pencil
on the smooth conic $C\subset\Pi\subset X_4^6$. The weighted-line
calculation gives the same angular bound $35/8$ for every
divisorial valuation whose center contains $C$. We record the changes
to the weighted presentation and the binary calculation.

Write $u_0,u_1,u_2$ for the plane variables and $g(u)$ for the
conic. Choose five linear normal coordinates consisting of light
variables $x_1,x_2,x_3$ and two heavy variables $z_1,z_2$, so
the polars paired with the $x_i$ vanish on $C$, while those for
the $z_i$ give the basepoint-free sextic pencil $(f_1,f_2)$.
Let $h$ be the order of the nonzero quartic
$F|_{z_1=z_2=0}$ along the complete-intersection ideal
$(x_1,x_2,x_3,g)$. The restriction is nonzero because a smooth
quartic sixfold contains no five-plane, and its first jets vanish
on $C$ since the plane is contained and the three specified
polars vanish on $C$. Hence $2\le h\le4$. For the upper bound,
decompose by the monomials $x^\alpha$: the order of a nonzero
coefficient times $x^\alpha$ is
$|\alpha|+\operatorname{ord}_g(F_\alpha)$, at most four by its
section degree. The symbol of the degree-two generator can be
taken to be
\[
                 w=\operatorname{in}\bigl(g+\sum_i L_i(u)x_i\bigr),
\]
where $L_i$ are arbitrary linear forms. This changes the original
symbol by a shear over the conic, preserving the same homogeneous
generators and weights. After this choice and a constant linear
change of the three $x_i$, the order-$h$ normal form $G$ has a
nonzero restriction $a(u)x_1^h$ on $w=x_2=x_3=0$. Indeed at the
generic point $G$ is a nonzero polynomial in four coordinates,
and some constant direction $\xi\in\mathbf C^3\setminus\{0\}$
gives $G(\xi,W)$ not identically zero as a polynomial over
$\mathbf C(C)$. Choose a general linear form for $\sum\xi_iL_i$;
its value in a fixed trivialization can be varied through infinitely
many distinct elements of $\mathbf C(C)$, so it avoids its finitely
many possible roots. Put this direction on the $x_1$ axis.
The coefficient $a$ is a nonzero section of $\mathcal O_C(4-h)$.

Give all four light generators order one and $z_1,z_2$ order $h$.
With $R_C=\mathbf C[u_0,u_1,u_2]/(g)$ the weighted initial ring is
\[
 R_C[x_1,x_2,x_3,w,z_1,z_2]\big/
            (f_1 z_1+f_2 z_2+G(x_1,x_2,x_3,w)).
\]
The section degrees of the three $x_i$ and the two $z_i$ are one,
and the section degree of $w$ is two. The initial presentation
before the hypersurface relation follows from the regular sequence
of the five linear normal generators and $g$ in the polynomial
ring, or directly by adjoining the symbol of $g$ and then taking
the initial equation $g(u)=0$. This weighted hypersurface quotient
is a domain. Off the vertex of the conic cone the two coefficients
$f_1,f_2$ generate the unit ideal locally, giving an integral
affine-space fibration there. The vertex fiber has dimension at
most six and the quotient has pure dimension seven and is
Cohen--Macaulay. Thus there are no further components or
embedded primes. In particular localization at the generic
point of $C$ identifies its order with the weighted divisorial
extraction of the four light variables and one heavy parameter
left after eliminating the other. Its discrepancy is $h+4$.

Specialize by splitting off $x_2,x_3,w$ to leave the relation
$f_1 z_1+f_2 z_2+a x_1^h$.
The limiting relation is prime by the same argument, and the
Hilbert functions of the weighted pieces agree. The free angular
variables $x_2,x_3$ now have mean exponent $1/7$ each, and
$w$ has mean exponent $1/14$. Indeed a free generator of
section degree $j$ has mean exponent $1/(7j)$; the relation
for $g(u)$ takes weight zero. The extraction order has mean
\[
   s=\frac{3+1/2+2h-h/4}{7},\qquad \frac{h+4}{s}\ge\frac{16}{3}.
\]
The remaining binary block $(x_1,z_1,z_2)$ has mean weighted
order $B=(1+7h/4)/7$, unchanged from the line calculation. For
the finite chart on $\mathbf P(1,h)$ apply that calculation over
$K=\mathbf C(C)$ with $Q=x_1^h$. In its polynomial degeneration
to $f_1z_1+aQ$ one can first choose $f_1$ general in the pencil,
so it misses the finitely many zeros of $a$ on $C$.
This proves primitivity on the punctured cone. The possible
vertex fiber has codimension at least one, and the same
Cohen--Macaulay argument establishes primeness, including for
the conic ring which need not be factorial. Thus $z_2$ has mean
exponent $1/7$ in that last limit, giving the same finite-chart
residual coefficient $B/h-1/7=(1/h+3/4)/7$. At $x_1=0$ the
relation modulo $x_1$ is $f_1z_1+f_2z_2$, prime as before.
The same exact divisibility computation therefore bounds the
orbifold coefficient there by $1/7$. On the opposite charts
the binary-axis weighted blowup gives discrepancy $1+h$ over
radial mean $B$, and $(1+h)/B\ge35/8$ for $2\le h\le4$.
The disjointness of binary directions, product charts and
weighted inversion of adjunction are exactly those of the line
calculation. All free angular coefficients are at most $1/7$.
Together with $(1/h+3/4)/7\le5/28$ this proves the assertion.
The angular threshold holds over the whole geometric generic
exceptional fiber. For every smaller coefficient, log
canonicity for valuation-adapted bases follows near that whole
fiber by weighted inversion of adjunction, and near the generic
point of $C$ downstairs by properness. This applies to each
valuation whose center contains $C$.

\section{Point discrepancies}\label{sixfold:points}
Let $p\in X$, let $T$ be its projective tangent hyperplane section, and let $b=\operatorname{mult}_p T\in\{2,3,4\}$. The Gauss map is finite, as the pullback of the dual hyperplane class is $3H$. Thus $T$ has isolated singularities and is integral. For $b=2,3$ projection from $p$ on $T$ is generically finite (on a general direction the residual polynomial has degree $4-b>0$). General linear cuts through $p$ down to a surface $S\subset X$ are smooth, and $C=T|_S$ is an integral plane quartic of multiplicity $b$. Indeed these cuts on $T$ are pullbacks of hyperplanes of the projection, whose system has image dimension at least two at every Bertini irreducibility step, including the last. General tangent directions preserve the leading form. In particular, if $b=3$ and the cubic tangent form is reduced, its restriction to the tangent plane of $S$ consists of three distinct directions. If $b=4$, $T$ is a cone over a smooth quartic fourfold (a singularity of the base would give a line of singularities of $T$). That fourfold is K-stable by \cite[Theorem 4.21]{ar2003.13788}; the generalized Eckardt point formula \cite[Corollary 4.19]{ar2003.13788} gives $\delta_p(X,H)=14/3$.

\paragraph{Retaining the pointed hyperplanes.}
Suppose $\delta_p(S,H_S)\ge3/B$ with $B>1$. For the four general hyperplane cuts through $p$, include on $X$ each of their defining divisors $D_i$ with coefficient $1-1/B$. They are transverse near $p$. Iterated basis adjunction gives
\[
 A_X(v)\ \ge\ (7/B)S(v)+(1-1/B)\sum_{i=1}^4v(D_i).
\]
Indeed on a $j$-dimensional rung a multiplier $(j+1)/B$ gives limiting cutting-divisor coefficient $1/B$ and restricts to multiplier $j/B$ in the complete residual series. Increasing the cut coefficient to one and using inversion of adjunction proves the assertion inductively. At finite levels take both the multipliers and the added coefficients slightly below their limiting values; the uniform complete-series comparison controls the large residual degrees, and bounded residual degrees have vanishing normalized contribution. For a valuation centered at $p$ normalize $v(\mathfrak m_p)=1$. Each cut has order at least one. Writing $A=A_X(v)$ yields
\begin{equation}
                    7S(v)\le B(A-4)+4 .\label{sixfold:signed}
\end{equation}
When $b=2$, the Seshadri constant on $S$ equals two. For a curve $D\ne C$ through $p$, $H_SD=C D\ge2\operatorname{mult}_pD$, and $C$ itself has ratio two. Hence $\delta_p(S,H_S)\ge3/2$ by \cite[Theorem 3.1]{ar2101.09246}. Taking $B=2$ gives $A\ge16$ whenever $A\le4S(v)$. The same theorem directly on $X$ gives $\delta_p(X,H)\ge7/2$ at every point with $b=2$, also for centers larger than the point.

\paragraph{The reduced cubic tangent form.}
On $S$ the integral $C\sim H_S$ now has an ordinary triple point. We prove $\delta_p(S,H_S)\ge18/13$. In the degree-$m$ complete series, split by ordinary order $j$ at $p$. A section of order $j$ has at least
\[
 a_j=\left\lfloor(3j-4m)_+/5\right\rfloor
\]
copies of $C$: after removing the actual $b'$ copies, intersection with $C$ gives $4(m-b')\ge3(j-3b')$. Also $j\le3m$ since $b'\le m$. Divide these forced copies out of the order-$j$ piece. The space of its leading binary forms has dimension $r_j\le q_j+1$, $q_j=j-3a_j$. Up to bounded errors $q_j=m g(j/m)$, where $g(t)=t$ for $0\le t\le4/3$ and $g(t)=4(3-t)/5$ for $4/3\le t\le3$. Since $\int g=2$ and $h^0(mH_S)=2m^2+O(m)$, the total dimension deficit $\sum_j(q_j+1-r_j)$ is $O(m)$.

For a basis adapted to ordinary order, remove $a_j C$ from each section of order $j$, giving a normalized residual divisor $R_m$ and a removed coefficient $b_m$. On the exceptional line of the point blowup the sum of the vanishing orders at any one point of $r_j$ independent degree-$q_j$ forms is at most $r_jq_j-r_j(r_j-1)/2$. Replacing the $r_j$ by $q_j+1$ in the normalized sums costs $o(1)$, uniformly in that point. Integration gives
\[
 \operatorname{mult}_p R_m\to\tfrac12\int_0^3 g^2=\tfrac89,\qquad
 \sup_{x\in E}\operatorname{mult}_x(R'_m|_E)\le\tfrac14\int_0^3g^2+o(1)=\tfrac49+o(1),
 \qquad b_m\to\tfrac5{27}.
\]
For the last identity integrate $g(t)(3t-4)_+/10$. Inversion of adjunction on $E$ therefore gives $u(R_m)\le(4/9+o(1))A_S(u)$ for any $u$ whose center contains $p$. An ordinary triple point has local threshold $2/3$ by its single blowup. Taking a simultaneous basis for $u$ and ordinary order gives
\[
 S_{H_S}(u)\le \left(\frac49+\frac5{27}\frac32\right)A_S(u)
               =\frac{13}{18}A_S(u).
\]
With $B=13/6$, inequality \eqref{sixfold:signed} gives $A\ge56/5$ whenever $A\le4S(v)$.

\paragraph{The nonreduced cubic tangent form.}
Every nonreduced cubic has a linear factor. We include the weighted calculation in this case, which also works for any cubic with a linear factor. This calculation uses the tangent-form method of Chapter~\ref{ch:fivefold}.

Choose coordinates \([T:y_1:\cdots:y_5:L:z]\) with \(p=[1:0:\cdots:0]\), \(z=0\) tangent and with equation
\[
 zB(T,y,L,z)+T L^s Q(y,L)+g(y,L)=0 .
\]
Here \(B(p)=1\), \(1\le s\le3\), \(Q(y,0)\ne0\), and \(g\) is quartic. The Gauss map is finite because it pulls back the hyperplane class to \(3H\). Consequently the tangent section has isolated singularities, and its hypersurface equation is integral. In particular \(G(y)=g(y,0)\) is nonzero. The hypersurface \(G=0\) in \(\mathbf{P}^4\) also has isolated singularities. Indeed on a hypothetical singular curve the derivatives of \(L^sQ\) and \(g\) at \(L=0\) can be nonzero only in the \(L\) direction. If \(\partial_L(L^sQ)\) is generically nonzero there, choose \(T\) to cancel the remaining derivative, giving a curve singular in the tangent section. If it vanishes identically, choose a zero on this projective curve of the cubic \(\partial_Lg\). The corresponding entire line of \(T\)'s is then singular. This rules out both possibilities.

The pair \((\mathbf{P}^4,G)\) is log canonical. A general \(\mathbf{P}^3\) avoids the isolated singularities and cuts out a smooth quartic surface. Scaling a complementary coordinate degenerates \(G\) to the cone over that surface. On blowing up the vertex the log discrepancy is \(4\) and the multiplicity is \(4\); the strict transform is smooth and transverse to the exceptional divisor. Semicontinuity of log canonical thresholds proves the claim.

In the affine chart \(T=1\) give the local parameters \(y_i,L\) weights \(s,s+1\), respectively. Let \(E\) denote the exceptional divisor. The weight of \(z\) is \(4s\), and in the section ring the initial relation is
\[
 T^3 Z+T L^s Q(y,0)+G(y)=0 .
\]
This polynomial is irreducible, since it is primitive and linear in \(Z\). After inverting \(T\), the initial functions of the local parameters are algebraically independent. Thus the displayed quotient injects into the associated graded ring, and its ordinary Hilbert series, the Hilbert series of a quartic, gives equality in each degree. Differentiating its weight character gives
\[
       S_H(E)=\frac{9s+1}{7},\qquad A_X(E)=6s+1,\qquad
       A_X(E)-4S_H(E)=\frac{6s+3}{7}.
\]
Indeed the sum of variable weights is \(10s+1\) and the degree-four relation subtracts \(4s/4\) from this sum.

On \(E=\mathbf{P}(s,s,s,s,s,s+1)\), scale \(L\) to zero in the polynomial representing \(Z\). The relation becomes \(T^3Z+G\), still prime with the same two-graded Hilbert function. The limiting spaces are the full initial spaces: products in the limiting generators occur as initial sections, and the dimensions agree. Splitting them further by the exponent \(\beta\) of \(L\), a piece of section degree \(m\) and extraction order \(J\) is precisely
\[
 L^\beta G^b\,\mathbf{C}[y]_{D-4b},\qquad
 D=\frac{J-(s+1)\beta}{s},\qquad
 b=\max\{0,\lceil(D+\beta-m)/3\rceil\}.
\]
Only nonzero pieces occur. In fact a monomial \(T^aL^\beta Z^\gamma P_l\) has \(m=a+\beta+\gamma+l,\ D=4\gamma+l\); eliminating \(Z\) in a chart proves the assertion by allowing all multiples with \(\gamma\ge b\).

Choose \(P_l\) from a graded complement to \((G)\). The resulting monomials form a basis by replacing \(G\) with \(-T^3Z\). The Hilbert series of \(\mathbf{C}[y]/(G)\) is \((1+t+t^2+t^3)/(1-t)^4\). The seven normalized indices \(a,\beta,\gamma\) and the four unrestricted indices of this last ring are therefore asymptotically uniform on the simplex of sum one, with indices divided by \(m\). The bounded degree shifts in the numerator do not affect their means. In particular,
\[
 \overline{\beta}=\frac17,\qquad
 \overline D=\frac87,\qquad
 \overline b=\frac3{28},\qquad
 \overline{D-4b}=\frac57 .
\]
Here \(b=(\gamma-a/3)_++O(1)\). The positive part of the difference of two independent mean-one exponential variables with coefficients \(1,-1/3\) has mean \(3/4\); normalizing their seven-variable vector by its sum gives the stated \(3/28\).

For any valuation \(u\) on the base \(\mathbf{P}^4_y\), the fixed \(G\)-factors contribute at most \(3A(u)/28\) by log canonicity. The complete moving systems contribute at most \((5/7)A(u)/5\), since \(\delta(\mathbf{P}^4,\mathcal{O}(1))=5\) \cite{ar1706.04548}. Their total is at most \(A(u)/4\). The estimates are uniform at large section degree for any fixed coefficient below four: the finite-level upper bound for \(S_j(u)/S(u)\) tends uniformly to one \cite[Corollary 3.6]{ar1706.04548}, while bounded moving degrees contribute \(O(A(u)/m)\) by \(\alpha(\mathbf{P}^4,\mathcal{O}(1))=1\).

On a uniformizing chart \(y_i\ne0\), the specialized basis divisor is a pullback from this base plus the fiber coordinate divisor of mean \(1/7\). These meet as a product on a base log resolution. The cyclic chart cover is crepant for the different of the weighted blowup (on the coarse space it has coefficient \(1-1/s\) along \(L=0\)). On the chart \(L=1\) the polynomials in the five \(y\)'s are homogeneous. Blow up its origin; the radial discrepancy is \(5\) and \(4\overline D=32/7<5\). The strict transforms are pullbacks of the same direction divisors on \(\mathbf{P}^4\). These two chart calculations prove the asymptotic residual threshold is at least four on the specialized model. To transfer it back, for each valuation choose bases in each filtration piece simultaneously adapted to it and to the specialization weights. Their initial forms are independent homogeneous bases. Semicontinuity for the resulting basis divisors, applied one specialization at a time for each fixed coefficient below four, gives the same residual bound on the original \(E\).

Finally use filtered adjunction \cite[proof of Theorem 3.3]{ar2003.13788}. For a basis adapted to \(v\) and \(\operatorname{ord}_E\), subtract its mean \(E\)-coefficient from the pullback. Inversion of adjunction bounds the value of the residual divisor by \((A_{\widetilde X}(v)-v(E))/c\) for every fixed \(c<4\) at sufficiently high levels near \(E\). Passing to the limit and then taking \(c\uparrow4\) gives
\[
             A_X(v)-4S_H(v)\ \ge\ \frac{6s+3}{7}\,v(E)>0,
\]
as \(v(E)>0\) for center exactly \(p\).

\section{Positive-dimensional minimizing centers}\label{sixfold:positive}
Suppose $\beta\le4$ and let a divisor computing it have positive-dimensional center $Z$, of dimension $s$. First, $\beta\ge7/2$. In fact the points of $X$ at which the differential of the Gauss map is zero form a finite set. A positive-dimensional irreducible component of that locus would map to an image of the same positive dimension by finiteness; in characteristic zero the restriction has positive differential rank at a general smooth point. The differential of the Gauss map is the second fundamental form, so at a general $p\in Z$ the tangent multiplicity is two. The point estimate $\delta_p(X,H)\ge7/2$ of Section~\ref{sixfold:points} applies to every valuation with center containing $p$, in particular to this divisor.

\paragraph{The map to the center.}
By \cite[Theorems 4.5 and 4.7]{ar1907.05399} the divisor induces a special degeneration to a klt $\mathbf Q$-Fano sixfold $Y$ with $\delta(Y,-K_Y)=\beta/4$. Put $L=-K_Y/4$. For sufficiently divisible degrees, the section ring on $Y$ is the associated graded of the filtration of the section ring on $X$ by the given divisor. Shift the linearization so that weights are orders of section divisors. Then the zero-weight subalgebra is the homogeneous coordinate ring of $Z$ in those degrees: positive order is precisely vanishing at its generic point. All high-degree sections on $Z$ lift by Serre vanishing. This gives a dominant rational map $Y\dashrightarrow Z$ from the subsystem in $|aL|$ corresponding to $|a(H|_Z)|$ for $a$ sufficiently divisible and large. Let $W\to Z$ be a smooth projective resolution and resolve the rational map to obtain $\pi:\widetilde Y\to Y$, $h:\widetilde Y\to W$. If $P=\pi^*L$ and $H_W$ is the pullback of $H|_Z$, then
\[
                         P-h^*H_W\sim_{\mathbf Q} E_0\ge0 .
\]
Indeed the difference multiplied by $a$ is the fixed divisor on a resolution of the subsystem.

If $s=1$, write $e=\deg_H Z$. At a general smooth point $q\in Z$, a section of $kH|_Z$ vanishes to order at least $ke-p_a(Z)$, by the Hilbert polynomial and jet evaluation. A component of a general fiber of the rational map is a prime divisor $F$ on $Y$ after closure: on the domain of the rational map the fiber is a nonempty divisor and meets the regular locus. Pulling the section back by the degree-zero-weight subsystem, with its possible effective fixed divisor, shows that the maximum asymptotic $F$-order satisfies $T_L(F)\ge e$. Since $A_Y(F)=1$ and $S_L(F)\ge T_L(F)/7$ by \cite[Proposition 3.11]{ar1706.04548}, one has $\beta=\delta(Y,L)\le7/e$. Hence $e\le2$ and $Z$ itself is a line or conic, contradicting the table.

\paragraph{A quotient slope on the central fiber.}
For $s\ge2$ we use
\begin{equation}
 -P^5\cdot h^*K_W\ \ge\ \frac{s\beta}{6}\,P^6 .\label{sixfold:slope}
\end{equation}
We give the metric input, since the central fiber may be singular and may have threshold below one. Choose a rational $0<t<\beta/4$ and a smooth positive reference curvature form $\omega_0\in c_1(-K_Y)$ on $Y$, in the sense of restriction from a projective embedding. Let $\mu_{\omega_0}$ be the adapted anticanonical volume measure associated to the reference metric. By \cite[Corollary 1.5 and equation (9)]{ar2608.23505}, with boundary zero, model potential zero and class $c_1(-K_Y)$, the equation
\[
      (\omega_0+dd^c\varphi)^6 = c\,e^{-t\varphi}\mu_{\omega_0}
\]
has a solution with minimal singularities, hence bounded since the class is ample. Here $c>0$ is a normalization constant. Its metric is smooth on the regular locus by regularity for the smooth-twist equation, as in \cite[\S2.1.2]{ar2008.05352} and \cite[Appendix B]{BBEGZ}. Taking $-dd^c\log$ shows that $\operatorname{Ric}(\omega)=t\omega+(1-t)\omega_0$ there.

Apply \cite[\S2.2, inequality (18)]{ar2008.05352} with $\gamma=1-t$ on the resolution. The estimate up to (18) uses the one fixed twisted solution; its later limit $\gamma\to0$ is unnecessary here. With slopes in $\pi^*c_1(-K_Y)$ it gives
\[
 \mu(\mathcal F)\le [1+(1-t)(6/k-1)]\mu(T_{\widetilde Y})
\]
for any rank-$k$ saturated subsheaf $\mathcal F\subset T_{\widetilde Y}$. Exceptional discrepancy divisors have zero degree against $P^5$. Subtracting the degree of the kernel, the slope of every rank-$s$ torsion-free quotient of $T_{\widetilde Y}$ is at least $t\mu(T_{\widetilde Y})$. Apply this to the image of the differential in $h^*T_W$; its determinant has degree at most that of $h^*T_W$. In the polarization $P$ this gives $-P^5h^*K_W\ge (4ts/6)P^6$. Let $t\uparrow\beta/4$ to get \eqref{sixfold:slope}.

A general complete intersection curve of high multiples of $L$ on $Y$ avoids the exceptional images of codimension at least two. Its strict transform has class proportional to $P^5$. By \eqref{sixfold:slope} its image on $W$ is a curve on which $K_W$ has negative degree. For such a curve, \cite[Theorem 1.1]{ar2509.08065} gives through a general point a rational curve of $H_W$-degree at most $(s+1)$ times the ratio of its $H_W$-degree to its anticanonical degree. Passing from the image cycle to its reduced support preserves this ratio. Choose the general point in the isomorphism locus of $W\to Z$. Since $P^5h^*H_W\le P^6$, there is an integral rational curve in $Z$ of positive $H$-degree
\[
           e\le\frac{6(s+1)}{s\beta}\le\frac{9}{\beta}
               \le\frac{18}{7}<3 .
\]
The table applies because the center contains this curve, contradicting $\beta\le4$. This excludes all positive-dimensional minimizing centers.

\section{The isolated point and Nadel vanishing}\label{sixfold:nadel}
Take a divisorial minimizer $v$, now centered at $p$. Normalize as in Section~\ref{sixfold:points}, and suppose $\beta\le4$. Thus $A(v)=\beta S(v)$ and $A(v)\ge56/5$. Fix $T_*=17/4$. Choose a hyperplane $H_0$ avoiding $p$ and, in every large degree, a basis of $R_m$ adapted to both $v$ and $\operatorname{ord}_{H_0}$. Remove its $H_0$-coefficient $a_m$ from its basis divisor $D_m$, obtaining $\Gamma_m=D_m-a_mH_0$. Then $a_m\to1/7$, $v(\Gamma_m)=S_m(v)\to S(v)$, and $u(\Gamma_m)\le S_m(u)$ for every $u$.

We claim $(X,T_*\Gamma_m)$ is klt outside finitely many points for all sufficiently large $m$. Otherwise let $t_m\le T_*$ be the first coefficient admitting a positive-dimensional non-klt center. It is attained on a log resolution, by taking the minimum over divisors with positive-dimensional image. At this coefficient the non-lc locus is finite. The uniform upper estimate for $S_m$ gives $t_m\ge\beta-o(1)$. Add a general hyperplane through $p$ with a rational coefficient $\epsilon_m>0$ tending to zero, chosen so the pair is non-lc at $p$. This is possible because $A(v)-t_m v(\Gamma_m)\le o(1)$ and the hyperplane has $v$-order at least one. On the punctured variety the system is basepoint-free. Bertini on a log resolution preserves log canonicity there wherever the pair was lc and preserves its positive-dimensional lc centers at their generic points. It creates no new lc place dominating a new positive-dimensional center: on the log resolution the general divisor is transverse to the coefficient-one strata and has coefficient less than one. The new non-lc locus $N$ is thus finite, nonempty. For the perturbed boundary $\Delta_m$ write
\[
 -(K_X+\Delta_m)\equiv\eta_m H,\qquad
 \eta_m=4-t_m(1-a_m)-\epsilon_m\ge5/14-o(1)>1/3 .
\]

Choose a minimal positive-dimensional lc center $Z$ under inclusion. For completeness, apply the quasi-log structure of the effective normal pair
$[X,\omega=K_X+\Delta_m]$ as in \cite[Example 4.11]{ar2102.11986}. Its non-qlc scheme has support $N$ and its qlc centers are the lc centers. By quasi-log adjunction \cite[Theorem 4.6(i)]{ar2102.11986}, $X'=Z\cup N$ inherits this structure, with canonical class $\omega|_{X'}$. The vanishing in Theorem 4.6(ii) for the zero Cartier divisor on $X$ gives $H^1(X,\mathcal I_{X'})=0$. Thus constants surject onto $H^0(X',\mathcal{O}_{X'})$, so $X'$ is connected and $N\subset Z$ set-theoretically. By minimality the non-qklt locus of $X'$ is a finite nonempty set. A complete curve in $Z$ meets it, for instance one through $p$. Fujino--Hashizume's length-one theorem \cite[Theorem 1.8]{ar2103.01417}, for $X'\to\operatorname{Spec}\mathbf{C}$, supplies an integral rational curve $C\subset Z$ with $0<\eta_m H\cdot C\le1$. In particular $C$ is a line or a smooth conic.

At the generic point of $Z$ the pointed hyperplane contributes no coefficient. There is an lc place $w$ of the original pair with center $Z$, hence
\[
  A(w)=t_m w(\Gamma_m)\ \le\ T_* (1+o(1))S(w).
\]
This contradicts the uniform table bound $A(w)\ge(56/13) S(w)$, because $C\subset Z$ and $17/4<56/13$. The claim follows.

Finally choose rational $0\le\lambda_m<1$ tending to $5/14$, with $T_*(1-a_m)+\lambda_m<4$, and take a general hyperplane $H_p$ through $p$. For $\Xi_m=T_*\Gamma_m+\lambda_m H_p$ the multiplier ideal has finite cosupport by the claim and Bertini. Since $-(K_X+\Xi_m)$ is ample, Nadel vanishing gives a surjection of constants onto $H^0(\mathcal{O}_X/\mathcal J(\Xi_m))$. The quotient has length at most one, so $\mathfrak m_p\subset\mathcal J(\Xi_m)_p$. The valuative membership criterion gives $v(\Xi_m)<A(v)+1$. Passage to the limit, using $S(v)\ge A(v)/4$, gives
\[
 1\ \ge\ T_* S(v)-A(v)+\frac5{14}
 \ \ge\ \frac{A(v)}{16}+\frac5{14}
 \ \ge\ \frac{56}{80}+\frac5{14}=\frac{37}{35}.
\]
This contradiction proves Theorem~\ref{sixfold:thm}.\qed

% ===== The quartic sevenfold and the boundary \(n=2d-2\) =====
\chapter{\texorpdfstring{The quartic sevenfold and the boundary \(n=2d-2\)}{The quartic sevenfold and the boundary n = 2d - 2}}\label{ch:boundary}
\section{Statement of the dimension bound}
\begin{theorem}\label{boundary:main}
Let \(H=\mathcal{O}_X(1)\) and \(d\ge4\). If \(d-1\le n\le2d-2\), or if
\((d,n)=(4,7)\), then \(\delta(X,H)>r\).
In particular \(\delta(X,-K_X)>1\) and \(X\) is uniformly K-stable.
\end{theorem}
The equivalence \(n\ge2r-2\Longleftrightarrow n\le2d-2\) follows
by substituting \(r=n+2-d\); the Fano condition is \(n\ge d-1\).
The added infinite family is \(d=r,\ n=2d-2\). The quartic sixfold
on that boundary is proved in Chapter~\ref{ch:sixfold}. The first further
cases on the boundary are the quintic eightfold, the sextic tenfold,
and the septic twelvefold. The isolated quartic sevenfold has
\((r,n)=(5,7)\) and lies on \(n=2r-3\).

\paragraph{Inputs and scope.}
The theorem for \(n\ge2r-1\), the special-degeneration slope
inequalities, and the principal-factor estimates for point valuations
are deductions in Chapter~\ref{ch:global}.
Its sections on the dimension and rational image of a minimizing
center give, for a klt boundary of class \(\lambda H\) with
\(\beta'=\delta(X,B;H)\le r-\lambda\), the inequality
\(s+1\ge r-\lambda-\tfrac n{n+1}\beta'\) for a divisorial minimizing
center of dimension \(s\). For \(s=1\), or when \(s>0\) and
\((s+1)\beta'\ge n+1\), that center contains a rational curve of
degree at most \((n+1)/\beta'\). Its point estimates are used explicitly
in Section~\ref{boundary:numbers}. The quartic sixfold proof is Chapter~\ref{ch:sixfold}.
Smooth cubics have a separate all-dimensional proof in
Chapter~\ref{ch:cubics}. The cubic surface, threefold and fourfold cases already
have published proofs \cite{keTian,ar1706.01933,ar2007.14320}. This manuscript does not require
the general cubic theorem.
The present upper endpoint in dimension is a sufficient condition;
optimality is open in this argument.

\paragraph{Basis adjunction and local thresholds.}

Write \(A_X(v)\) for the log discrepancy, and \(S(v)=S_H(v)\) for
expected vanishing in the complete section ring of \(H\). We use
\(\delta(X,H)=\inf_v A_X(v)/S(v)\) and
\(\delta_C(X,H)=\inf_{C\subseteq c_X(v)} A_X(v)/S(v)\), with
divisorial valuations in both infima. Scaling \(H\) by \(r\) divides
the threshold by \(r\).

Write $R_m=H^0(X,mH)$ and $N_m=\dim R_m$. An $m$-basis divisor is the sum of the divisors of a basis, divided by $mN_m$. Its maximum $v$-order is $S_m(v)\to S(v)$, and $S_m(u)\le(1+o(1))S(u)$ uniformly in $u$ by \cite[Corollary 3.6]{ar1706.04548}. Two filtrations of a finite-dimensional vector space admit a common adapted basis. For an extraction with exceptional divisor $E$, start with a basis adapted to $v$ and to $E$. Its pullback has coefficient $S_m(E)$ at $E$; the residual restricts to bases of the order quotients on $E$. If all these residual basis divisors admit coefficient $c$ on $E$ with its different, inversion of adjunction and the discrepancy formula give in the limit
\[
 A_X(v)\ \ge\ c S(v)+(A_X(E)-cS(E))v(E).
\]
Smaller coefficients and passage to the limit suffice. This is the proof of \cite[Theorem 3.3, (3.6)]{ar2003.13788}. Over a curve $C$ it suffices to work over its generic point and bound the residual basis divisor on the \emph{whole} geometric generic exceptional fiber. Shrink the curve to spread a log resolution; inversion of adjunction and properness give log canonicity near its generic point downstairs. This also tests every larger center containing $C$.

In using a torus specialization to bound arbitrary bases of a residual series, fix a valuation on that series first. On each piece choose a basis adapted to it and to the torus filtration, keeping the coefficient spaces over $\mathbf{C}$, even on a generic fiber. Its valuation sum is maximal. Its initial forms give bases of the limiting weight pieces. A uniform lower bound on the log canonical thresholds of these initial divisors passes to the adapted divisors by lower semicontinuity on the proper fiber. Thus the bound controls every incoming basis. These changes of basis never require simultaneous adaptation to three filtrations.

For a smooth hyperplane flag the refined residual pieces on the next rung are the full spaces of degrees $m-i$ by the hypersurface restriction sequence. In ambient dimension $j$ the mean hyperplane coefficient tends to $1/(j+1)$ and the mean normalized residual degree tends to $j/(j+1)$. Thus a lower threshold $c$ on the next rung lifts to $\min\{j+1,(j+1)c/j\}$; see \cite[Lemma 4.6]{ar2003.13788}. A chain may be used just near the generic point under consideration when the flag members are normal and smooth there.

\paragraph{Curve inputs.} Chapter~\ref{ch:models} proves the common
threshold \((n+1)(d+1)/(2d)\) used below for lines and conics,
together with its refinements in small degree.
\section{From the curve bound to uniform K-stability}\label{boundary:numbers}
Assume \(r=d\ge5,\ n=2r-2,\ N=n+1,\ u=n/N\). The preceding sections give
\[
  L_C=\frac{N(d+1)}{2d}=r+\frac12-\frac1{2r}
   \le\delta_C(X,H)
\]
for every integral curve of degree at most two. An integral conic
over \(\mathbb C\) is a smooth plane conic.
Suppose \(\beta=\delta(X,H)\le r\). By \cite[Theorem 1.2]{ar2102.09405} a
divisorial valuation \(v\) computes \(\beta\); finite generation makes
its threshold rational. Take any point \(p\) in its center and a
general smooth hyperplane section \(Y\) through \(p\).
It has dimension \(n-1=2(r-1)-1\) and index \(r-1\).
Chapter~\ref{ch:global} therefore gives \(\delta(Y,H|_Y)>r-1\).
The hyperplane adjunction bound gives
\[
 \beta\ \ge\ \delta_p(X,H)\ \ge\
 \min\{N,(N/n)\delta(Y,H|_Y)\}\
 >\frac{N(r-1)}{n}=\frac N2 .
\]
In particular \(\beta>r-1\). This applies the dimension bound of Chapter~\ref{ch:global} to a
single smooth pointed section, and deals with every possible low
value of the original threshold.

Let $r-1<\beta\le r$, set $b=r/N$, and fix $\epsilon=r^{-4}$.
Use the upper fraction $1-1/d$ of a basis adapted to a minimizer and
a free hyperplane that avoids its generic center. Remove its copies
of that hyperplane, and multiply the resulting divisor by
$\lambda=(d-1)x$, where $x=(\beta-N/3-3\epsilon)/d$.
For the resulting boundary $B$ of class $\lambda cH$ one has
$c<1-[1-\log(d)/(d-1)]/N$. Indeed the full hyperplane coefficient has survival distribution
\((1-t)^n\) on \([0,1]\). Retaining any fraction \(1-f\) removes
at least \(\int_0^1((1-t)^n-f)_+\,dt/(1-f)\) in the limit.
Taking \(f=1/d\) and \(d^{-1/n}>1-\log(d)/n\) gives the displayed
strict estimate. There is a choice at a sufficiently
large level with $B(w)\le kS(w)$ for all $w$,
$k=\lambda+x+\epsilon$. This follows by extending the chosen basis vectors to a complete basis and
using uniform finite-level comparison, since \(\lambda/(1-1/d)=\lambda+x\).
The polynomial subseries of any general finite
projection has fraction $1/d$ and its mean at the original minimizer
is at most $bS(v)$: the projection is \'etale at the generic center
and \(\delta(\mathbf{P}^n,\mathcal{O}(1))=N\). In each degree the mean of the lowest
weights of the full space is at most the mean of a subspace of the
same dimension (compare the dimensions of all filtered quotients).
The gain on its complement is consequently at least
\((1-b)/(d-1)\) times \(S(v)\), up to a term tending to zero. Hence
\[
 N/3<\beta-k\le\beta'=\delta(X,B;H)
 \le\beta-\lambda-x(1-b)+\epsilon.
\]
The upper bound is strictly less than the index $r-\lambda c$ of the pair, since the gain $x(1-b)$ exceeds the finite-level allowance in the indices using this tilt. This strict range supplies a dreamy minimizing divisor for the tilted
log Fano pair. Put $x_r=(r-N/3-3\epsilon)/d$ and set
\[
 U=r+bx_r+2\epsilon,\qquad
 g=u(1-b)-\log d/N,\qquad z=b+x_r g-u\epsilon .
\]
Then $A_X(w)\le U S(w)$ at a tilted optimizer, and
$r-\lambda c-u\beta'\ge z$. For the latter use
$r-u\beta+xg-u\epsilon$, a decreasing affine function of $\beta$
since $0<g<ud$. The curve-image and canonical-extension inequalities
give a line or conic in every positive center since $\beta'>N/3$.
If $U<L_C$ this excludes positive centers.

At a point center normalize $w(\mathfrak m_p)=1$, and write
$A_X(w)=n+q$. It suffices to find $\Delta>0$ and an established
lower bound $Q\le q$ with
\[
 U+\Delta<L_C,\qquad z-u\Delta>1/3,\qquad
 z+\Delta((n+Q)/U-u)>1 .                                      \tag{*}
\]
Indeed the first positive-dimensional non-klt center of a residual
basis divisor at coefficient at most $\beta'+\Delta$ would contain
a curve of degree at most two by the length-one theorem and the
tiny pointed perturbation. Its place would violate the first
inequality. At the isolated point Nadel and the remaining ample slack
contradict the last inequality.

\paragraph{The multiplier ideal at the optimizing point.}
Here are details of the passage used in \((*)\). The pair \((X,B)\)
is fixed at a sufficiently high basis level. Choose a divisor \(w\)
computing \(\beta'\) and, if it has point center \(p\), normalize
\(w(\mathfrak m_p)=1\). In a new degree \(m\), remove the copies of a
hyperplane avoiding \(p\) from a \(w\)-adapted, hyperplane-adapted
basis divisor. The residual \(\Gamma_m\) has class tending to \(uH\),
\(w(\Gamma_m)\to S_H(w)\), and \(\nu(\Gamma_m)\le(1+o(1))S_H(\nu)\)
uniformly. Put \(t=\beta'+\Delta\). If positive-dimensional
non-klt loci occur at \(t\), take the first positive-center threshold
\(t_m\le t\). It satisfies \(t_m\ge\beta'/(1+o(1))\).
Adding a general pointed hyperplane with coefficient \(o(1)\) if
needed makes the pair non-lc at \(p\) and preserves the positive
centers at their generic points. The non-lc scheme is finite.
The ample class left has coefficient at least \(z-u\Delta-o(1)>1/3\).
Quasi-log adjunction and vanishing on the union of a minimal positive
lc center with the non-lc scheme make that union connected, with
finite nonempty non-qklt locus meeting the positive component.
The length-one theorem \cite[Theorem 1.8]{ar2103.01417} gives a rational curve
in that center of degree at most two (see the short-curve proof in
Chapter~\ref{ch:global}). Its original lc place has \(A_X/S_H\le
k+t+o(1)\le U+\Delta+o(1)<L_C\), impossible. Thus the multiplier
ideal for \(B+t\Gamma_m\) has finite cosupport for large \(m\).
Now add a general hyperplane through \(p\) of coefficient approaching
\(\eta=z-u\Delta\) from below, taking first \(m\to\infty\).
All displayed choices have \(0<\eta<1\).
The augmented pair remains klt outside finitely many points and its
anticanonical remainder is ample. Nadel vanishing makes constants
surject onto the multiplier-ideal quotient, so its length is at
most one. The ideal therefore contains \(\mathfrak m_p\) locally.
The valuative multiplier-ideal inequality gives
\(1\ge\Delta S_H(w)+\eta\) after passing to the limits, contrary to
the final inequality in \((*)\).

For $r\ge1024$ write $y=\sqrt r$ and
$E=0.30+0.18\log r$. One has $U<r+0.168$ and, more precisely,
$U-r<1/6+0.26/r$ by $x_r<(r+1)/(3r)$.
Since $x_r\ge1/3$, $u(1-b)\ge1/2-1/N$ and $b>1/2$, one has
$z\ge2/3-E/r$. Take $\Delta=1/3-(E+0.1)/r$; small rational
perturbations are harmless. As $L_C=r+1/2-1/(2r)$ the first
two tests in $(*)$ are strict. The simple universal point envelope
$S\le J+Gq/n$ satisfies $J<1.515$ and $G<1.56y$ in this range.
For example use
\[
 J\le \max\!\left\{\frac{2\sqrt{2r}+n-2}{N},
                   \frac{n-2+r+\sqrt{r/2}}{N}\right\};
 \quad G\le\max\{\sqrt{2r}+4,\sqrt{r/2}+\log(12r^2)\}
\]
and the decreasing quotients $\log(12r^2)/\sqrt r$ and $4/\sqrt r$.
Unless the point is excluded outright, $n-UJ>0.482r$ and
$UG/n-1<0.782y$ imply $q>0.615y$. Therefore $(n+q)/U>2+0.54/y$,
since $0.075y>2.336+0.091/y$. The last test in $(*)$ follows from
$0.18y>2E+0.1+0.54(E+0.1)/y$, which holds at $y=32$ and persists by
monotonicity after division by $y$.

For $56\le r<1024$ the rational simpler envelope verifies $(*)$
with $\Delta=\min\{L_C-U,(z-1/3)/u\}-\epsilon$. For
$15\le r\le55$ the principal-bin point bound of Chapter~\ref{ch:global}
does the same. For $7\le r\le14$ choose a fraction $f$ and a share
$a$ from the exact witnesses. Write
\[
 k_0\ge n(1-(1-f)^{1/n}),\quad
 \gamma\le(1-b)\frac f{1-f}\min(1,\log(1/f)/\log d),\quad
 \alpha\le\frac{1-\frac f{1-f}n(1-f^{1/n})}{N}.
\]
Use $\epsilon=10^{-8}$ and
$\lambda=a(\beta-N/3-3\epsilon)/(1+k_0)$ on the retained
upper $1-f$ fraction after free-hyperplane subtraction. The
beta-law tail, projection-mean and free-hyperplane bounds give
$k=\lambda(1+k_0)+\epsilon$, $\beta'>N/3$,
$\beta'\le\beta-\lambda(1+\gamma)+\epsilon$. Thus at the
worst endpoint $\beta=r$ one may take
\[
 U=r+\lambda_r(k_0-\gamma)+2\epsilon,\qquad
 z=b+\lambda_r[u(1+\gamma)-1+\alpha]-u\epsilon .
\]
The expressions before taking the endpoint are respectively
increasing and decreasing, as is checked directly for the finite
witnesses. For the pairs $(1000f,100a)$ in increasing order of $r$ use
$(143,70),(123,75),(111,80),(95,85),(91,85),(83,85),(75,90),(71,90)$.
For these smaller tilts the bounds on \(k_0,\gamma,\alpha\) are those
of the top-fraction argument in Chapter~\ref{ch:global}. The upper tail mean for
any valuation is at most \(1+n(1-(1-f)^{1/n})\) times its full mean
by concavity of the \(n\)-th root of the survival distribution. For
\(f\le1/d\) the projection mean bounds the lowest \(f\)-fraction.
For \(f>1/d\) the distribution function is log-concave by
Brunn--Minkowski on the Okounkov body; its inverse at \(e^t\) and
its bottom mean at \(e^t\) are convex in \(t\). Interpolation between
\(1/d\) and \(1\) gives the displayed logarithmic factor in the gain.
Hyperplane subtraction gives \(\alpha\) by the integral above. These
estimates apply to finite levels with the indicated positive allowances.
The included checker verifies by rational arithmetic that $(*)$ follows
from the principal bins in all eight cases.

\paragraph{The two small indices on $n=2r-2$.}
Keep the notation of the finite-fraction proof, including the weaker
curve bound $L_C=N(d+1)/(2d)$. At $r=5$ take $(1000f,100a)=(198,53)$
and $Q=31/10$. At $r=6$ take $(166,63)$ and $Q=67/20$. The exact
verifier gives respectively $U<5.09395$, $U<6.11141$, and the
lower bounds on $q$ required by $(*)$ are $<3.095$ and $<3.347$.
We show that a point minimizer with $A_X(w)\le U S_H(w)$ has $q>Q$.
We normalize $w(\mathfrak m_p)=1$ and write $q=A_X(w)-n$.

Let $b$ be the integer multiplicity of the tangent hyperplane $T$ at
$p$. If $b<d$, general pointed cuts to a smooth surface $S$ give an
integral plane curve $C=T|_S$ of multiplicity $b$. All tangent-plane
contact computations may be performed in that plane: projection of
the smooth surface from a transverse coordinate is a local analytic
isomorphism, identifying its tangent-plane section with $C$. One has
$\epsilon_p(H_S)\ge\min(b,d/b)$ by intersection with $C$.
If $b^2>d$, equality holds, and $\tau_p(H_X)=b$. Indeed the
restriction of any effective divisor of degree $m$ to the general
pointed surface has order at most $bm$: after removing all copies
of $C$ its intersection with $C$ bounds the remaining order by
$(d/b)$ times the remaining degree. Equality is attained by $T$.

For $r=5$, the case $b=2$ admits $B=5/2$ in
$NS\le B(q+2)+n-2$ by the Seshadri estimate. For $b=3,4$ use the
surface argument by splitting off the forced copies of $C$.
In degree $m$ and ordinary order $j$ at $p$, there are at least
$a_j=\lfloor(bj-dm)_+/(b^2-d)\rfloor$ copies. After removing them,
the leading binary forms have degree $j-ba_j$. Its scaled envelope
is the tent $g(t)=t$ for $0\le t\le d/b$, and
$g(t)=d(b-t)/(b^2-d)$ for $d/b\le t\le b$. Its integral is $d/2$,
so the total deficit from full binary spaces is $O(m)$ in a degree-$m$
series of dimension $dm^2/2+O(m)$. As in the quartic calculation,
the residual basis divisor has local $w$-order at most
$(d/(3b)+o(1))A_S(w)$: its ordinary coefficient tends to
$(2/d)\int g^2=2d/(3b)$ and its angular coefficient at any point to
at most $d/(3b)$. The average removed coefficient tends to
$(2/d)\int g(t)(bt-d)_+/(b^2-d)\,dt=(1-d/b^2)/3$. Thus if
$(S,cC)$ is lc near $p$, the signed flag may use
\[
              B=d/b+(1-d/b^2)/c .
\]

For a plane quintic of multiplicity $b=3$, one may take $c=4/9$.
Blow up $p$, with exceptional $E$ and strict transform $C'$.
At $x\in E\cap C'$ let $\ell$ be the line in direction $x$.
As $C$ is integral, its intersection with $\ell$ is at most five,
so $\operatorname{mult}_x C'\le 5-3=2$. Hence at any valuation centered at $x$,
$w(C')\le2A_{\widetilde S}(w)$ by the smooth multiplicity bound.
Also $w(C')\le3(A_{\widetilde S}(w)-w(E))$ by inversion of adjunction,
since $C'\cdot E=3$. Set $h=w(E)/A_{\widetilde S}(w)$. Then
\[
 \frac{w(C)}{A_S(w)}
 \le\frac{3h+\min(2,3(1-h))}{1+h}\le\frac94 .
\]
The exceptional divisor itself satisfies the same bound.
Thus $B=8/3$. If $b=4$ and no irreducible factor of the tangent
quartic has multiplicity at least three, every root on a general
projective direction line has multiplicity at most two. On the blowup
of $S$ the boundary $\tfrac12 C'+E$ is lc along $E$ by inversion of
adjunction, so $c=1/2$ and $B=21/8$. If a factor occurs at least
three times, its degree is one; write its multiplicity $s=3,4$.
The weighted linear-factor calculation for a tangent form of degree
$d-1$ gives
\[
         A_X(w)-d S_H(w)\ge
         \frac{(d-1)((d-2)s+1)}{N(s+1)}
\]
at such normalized valuations. With $A_X(w)\le U S_H(w)$, the
right side times $U/(U-d)$ is a lower bound for $A_X(w)$.
In every quintic case these bounds imply $q>31/10$ by the verifier.

For the weighted inequality just used, write the equation on the quintic
as \(zB+T L^sQ+g=0\), with \(T=1\) at the point, \(B(p)=1\),
\(g\) of degree five, and \(y_1,\ldots,y_7,L\) regular parameters.
Here \(G=g(y,0)\) is nonzero and has only isolated projective
singularities. Otherwise on a positive-dimensional singular locus
take a zero of \(\partial_Lg\): since \(s\ge3\), its entire line of
\(T\)'s is singular in the tangent section, contradicting the finite
Gauss map. The same reasoning excludes \(G=0\). A general hyperplane
of \(\mathbf{P}^6_y\) avoids these singularities and cuts a smooth quintic;
degenerating to its cone shows \((\mathbf{P}^6,G)\) is lc (vertex discrepancy
six, multiplicity five). Weight \(y,L,z\) by \(s,s+1,5s\).
The initial relation \(T^4Z+TL^sQ(y,0)+G\) is prime by primitivity,
and gives the whole initial ring by generic localization and ordinary
Hilbert functions. It has \(A(E)=8s+1,\ S(E)=(12s+1)/9\).
In the angular series, scaling \(L\) leaves \(T^4Z+G\), again prime.
In degree \(m\), the fixed-\(L^k\), fixed-extraction-order pieces are
\(L^kG^j\mathbf{C}[y]_{D-5j}\), \(j=\max(0,\lceil(D+k-m)/4\rceil)\), with
\(D\) the resulting \(y\)-degree. Use basis monomials
\(T^aL^k Z^h P_l\) for a graded complement \(P_l\) to \((G)\).
The Hilbert series of \(\mathbf{C}[y]/G\) has six free indices and bounded
shifts. Including \(a,k,h\), this is the nine-variable simplex.
Here \(D=5h+l,\ j=(h-a/4)_++O(1)\), so the normalized means are
\(\bar k=1/9,\ \bar D=11/9,\ \bar j=4/45,\ \overline{D-5j}=7/9\).
For the positive part use independent exponentials for the simplex
coordinates. Fixed factors and complete plane systems together give
mean base order at most \((4/45+(7/9)/7)A=A/5\) on \(\mathbf{P}^6\).
On the chart of a nonzero \(y_i\), the remaining fiber coordinate
contributes the independent divisor of coefficient \(1/9\). On
\(L=1\), the radial blowup of the \(y\)-origin has discrepancy seven
and \(5\bar D=55/9<7\). Its directions have the preceding base bound.
These are the cyclic uniformizing charts with their induced different.
Semicontinuity and extraction adjunction thus give
\[
A_X(w)-5S_H(w)\ge (A(E)-5S(E))w(E)
 =\frac{12s+4}{9} w(E),\qquad w(E)\ge\frac1{s+1}.
\]
The last bound follows by taking a Cartier multiple of the weighted
exceptional ideal: all monomials of weight at least \(k\) have ordinary
order at least \(k/(s+1)\).

For $r=6$ and $b=2,3$, the pointed surface has
$\epsilon_p(H_S)\ge2$ and $\tau_p(H_S)\le3$. The latter again
follows by intersecting with the tangent plane curve after removing
copies of it. The restricted angular estimate on the exceptional line
of the point blowup in the proof of Abban--Zhuang,
\emph{Seshadri constants and K-stability}, Lemma~3.2, applies at
$c=3(2)/6=1$. Write $\sigma=S_{H_S}(E)$. Along the ordinary
exceptional valuation the normalized order law on $[0,2]$ has density
$2t/6$, since $\pi^*H_S-tE$ is nef there and its volume is $6-t^2$.
Its remaining mass is $1/3$ and has support in $[2,3]$. Therefore
\[
     \sigma \le \int_0^2\frac{t^2}{3}\,dt+3\left(1-\frac46\right)
                  =\frac{17}{9}.
\]
In the same surface adjunction take the additional ideal
$\mathfrak m_p^{\,2-\sigma}$: its pullback contributes precisely its
exponent along $E$ and has no residual divisor. More explicitly, for
all coefficients $c'<1$, the angular threshold bound for
order-adapted bases applies asymptotically, and their ordinary
coefficient tends to $\sigma$. The discrepancy formula and inversion
of adjunction then allow the additional ideal exponent strictly below
$2-c'\sigma$. Pass to the limit. Since the restriction of
$\mathfrak m_{p,X}$ to every general pointed flag member is its
maximal ideal there, the signed flag with $B=3$ retains this same
additional exponent at every step. Hence
\[
       NS_H(w)\le 3q+n-2+3\sigma\le 3q+n-2+17/3 .
\]
This gives $q>67/20$. For $b=4,5$ the radial threshold equals $b$
and is computed by $T$. Apply the earlier principal-divisor bin
estimate at this exact value rather than an interval of possible
thresholds. The checker uses $K=20$ for each and verifies both affine
tests at $q=0,67/20$ with rational arithmetic.

Finally if $b=d$, the smooth $(n-2)$-dimensional degree-$d$ variety
of lines at the generalized Eckardt point has index $r-2$ and lies
in the established range $n-2\ge2(r-2)-1$. The formula of
Abban--Zhuang, \emph{K-stability via admissible flags}, Cor.~4.19,
gives $\delta_p(X,H)=nN/(n-1+d)>U$ for both indices.

\section{The quartic sevenfold}

\paragraph{The sevenfold line bound.}
For every line $L$ on a smooth quartic sevenfold,
$\delta_L(X,H)\ge16/3$. Use the weighted mixed-monomial and two-power
calculations \eqref{models:eq-mixed-monomial} and \eqref{models:eq-pure-powers}, with $N=n+1=8$ and $\alpha=3/4$.
For three or more independent normal polars, the positive degrees
of $N^*_{L/X}(1)$ sum to three and there are at least two of them.
The ordinary splitting calculation therefore has largest degree at
most two, and yields $\delta_L/N\ge4/(4+2)=2/3$ (the radial entry
is larger). For a polar pencil there are five light variables.
At contact two, (L1) gives $\delta_L/N\ge2/3$ already with a pure
light monomial. At contact three, write $G=uP+vQ$ for the transverse
leading coefficient, for cubics $P,Q$ on the light $\mathbf P^4$.
If some common zero $e$ is not a common cone vertex, restrict to a
general line through $e$. Its contact order is one or two, yielding
mixed exponents $(1,2)$ after a diagonal degeneration. Thus (L1)
again gives $2/3$. Otherwise the common vertex has codimension at
most two by the dimension theorem applied to $P,Q$. The quartic
term has a nonzero value in some direction $y$ of this vertex:
otherwise its join with $L$ gives at least a four-plane in $X$,
impossible since the four normal cubic derivatives would have a
common zero there. (In ambient $\mathbf P^8$ a four-plane has
four normal directions.) Retain such a $y$ and a transverse $x$
with $G(x)\ne0$. After the weights $(1/3,1/4,1)$ the two surviving
pure powers are $x^3,y^4$ as in \eqref{models:eq-pure-powers}. The two entries of \eqref{models:eq-pure-powers} are
$12/13$ and $19/28$, both $>2/3$.

At contact four the transverse form $G$ defines a smooth quartic
threefold in $\mathbf P^4$, since a singular point would join with
$L$ to give a singular plane in the fixed section by the two heavy
coordinates, impossible by finite Gauss. At a general point its
second fundamental form is nonzero because the finite Gauss map
is separable; equivalently its differential is nonzero. A tangent
direction where this quadratic form does not vanish gives a line
of contact exactly two. Restricting to its two-dimensional vector
space and degenerating the binary form by the order at that point
gives $x^2y^2$. (L1) now gives the same $2/3$. This proves the
assertion for all ambient valuations with center containing $L$.

\paragraph{Conics on contained planes of a quartic sevenfold.}
Let $C\subset\Pi\subset X_4^7\subset\mathbf P^8$ be a smooth conic
on a contained plane of a smooth quartic. Then
\[
                         \delta_C(X,H)\ge16/3 .
\]
The threshold tests all valuations with center containing $C$. All
mean coefficients in the proof are multiplied by eight, so the
upper bound on order per unit log discrepancy is $3/2$.

The six plane polars restrict on $C\simeq\mathbf P^1$ to a
basepoint-free system of sextics. For the quotient of $\mathcal
O^{6}$ by their copy of $\mathcal O(-6)$ let $a_i$ be the five
nonnegative split degrees. They sum to six. The ordinary conic
extraction gives normalized angular bound $1+\max a_i/8$, and
radial bound $(25/4)/6<3/2$. Indeed the complete split normal
monomial with exponent $\beta$ and ordinary remaining base cost
$p$ has coefficient dimension $2p+1+a\cdot\beta$. The conic
equation supplies an independent degree-two normal generator.
The simplex barycenters give normalized split means $1+a_i/8$
and mean $1/2$ for the quadratic generator, as in the ordinary
conic computation. Only the pencil with positive splitting $(6)$
or the net with positive splitting $(1,5)$ require more, since
the number of positive summands is the polar dimension minus one.

For a pencil use the weighted binary calculation with two heavy
normal variables and light variables comprising the others and the
conic equation. The first nonzero light part has order $h=2$:
the two heavy cubic polars have a common zero on $\Pi$, and some
light polar must be nonzero there by smoothness. It vanishes on
$C$ and hence equals the conic equation times a nonzero linear
form. The binary-axis ratio is consequently
$(1+7h/4)/(1+h)=3/2$. The other binary bounds are $1$ and
$3/4+1/h=5/4$, and the free normalized exponents are at most
one. For the binary calculation, write the core equation as
$f_1z_1+f_2z_2+a x^h=0$, with the two sextic polars $f_1,f_2$
basepoint-free. The choice and splitting of $x$ are exactly the light
axis construction given below. The graded ring over the conic has
one quartic relation, so normalized mean binary weight is
$B=1+2h-h/4$. In the coarse variables $Q=x^h,z_1,z_2$ the block
is a binary linear system, up to the bounded powers $x^i$, $i<h$.
Take $f_1$ general in the pencil, avoiding all roots of $a$, then
diagonally scale $z_2$ out of the relation. The remaining relation
$f_1z_1+aQ=0$ is primitive over the punctured conic cone and its
ordinary Hilbert series is unchanged. The free coordinate $z_2$
has mean exponent one, the other coarse coordinate has mean
$B/h-1$. These bound all finite binary-point orders by
semicontinuity after choosing bases simultaneously adapted to the
point order and the diagonal weight. At infinity the quotient by
$x$ is prime because $f_1,f_2$ have no common zeros. Thus the
actual upstairs order equals divisibility by the section-cost-one
parameter $x$, giving normalized coefficient one. The weighted
blowup of the binary axis has discrepancy $1+h$ and normalized
order $B$; the binary factors separate on its covering charts,
transversely to the independent free coordinates. This proves the
pencil case and the same calculation used for $h=2$ in the net.

For the net write its polars $(t\ell,-s\ell,b_0)$ after a
constant coordinate change, where $[s:t]$ parameterizes $C$,
$\deg\ell=5$, $\deg b_0=6$, and they have no common zero. This
is precisely the linear syzygy $(s,t,0)$. Give the three heavy
normals $z_1,z_2,z_3$ weight $h$, the first order of the quartic
in the light variables on their common five-plane. This five-plane
cannot be contained in $X$ since its three normal cubics have
common zero there. Thus $2\le h\le4$. There are three light
linear normals and one light conic equation. After shearing the
quadratic symbol by linear combinations of the light normals with
linear base coefficients, the first nonzero form has nonzero
restriction $a x^h$ on a light linear axis. Split off the remaining
light coordinates as in the binary calculation. The exact ternary
core of the initial ring over the conic cone has equation
\[
              a x^h+\ell T+b_0z_3=0,\qquad T=t z_1-s z_2 .
\]
Here $\deg a=2(4-h)$ on $C$. The initial presentations follow by
the absence of common polar zeros on the punctured conic cone,
the lower dimension of the vertex fiber, and ordinary Hilbert
series. If $h=2$, split off one general heavy coordinate, leaving
two basepoint-free polars; the same binary bounds apply and the
free heavy normalized exponent equals one. Assume $h=3$ or $4$.

On the coarse weighted plane the heavy coordinates $z_i$ are
independent; its boundary is $(1-1/h)L$ for the line
$L=(\ell T+b_0z_3=0)$. Put $Q=x^h=-(\ell T+b_0z_3)/a$. Removing
the bounded parity $x^s$, $0\le s<h$, leaves the coarse systems
generated by the $z_i,Q$, all of coarse degree one, of section
costs $1,1,1,h$. In a block of coarse degree $k$ and section
degree $p+k$, $Q$ has exponent at most $\lfloor p/(h-1)\rfloor$.
Multiplying by that power of $a$ embeds this block in the full
product of base degree $\lfloor2rp\rfloor$ and plane degree $k$.
For $h=3$ take $r=3/2$, and for $h=4$ take $r=1$. Base factors
are units over the generic fiber. The total enlarged dimension
divided by the true total tends to $R=rh/4$, since the true heavy
ring adds the degree-$h$ generator and one quartic relation to
the polynomial ring on the $z_i$ over the conic, whereas the enlarged
base degrees multiply its leading coefficient by $r$. Free and
bounded-power factors are common to both sums. For each valuation,
extend each original basis to the enlarged block. It suffices to
bound the independently normalized enlarged orders multiplied by $R$.

Off $L$ complete plane series give bound $R$ since the normalized
mean $k$ is three. The divisor $L$ has section cost $D=1+3/r=h$,
so its normalized order is $1/D$ and its discrepancy is $1/h$.
On $L$ the residual system has binary splitting $(1,5)$, with
$z_3$ generically the high-degree coordinate. For $P=[s:t:0]$
and $P'\in L\setminus\{P\}$ dehomogenize $z_3$ to one. The order
at $P'$ is bounded by the degree in the low binary coordinate.
The complete quotient blocks, after peeling off $L$, have density
$i+5j+2rp$ on the simplex of the two exponents, the base cost,
and the peeled line cost. Thus the normalized low mean is
$1+1/(2r+6)$. This uses the symmetric-power quotient of the
basepoint-free sextic net, complete outside a bounded strip in
base degree by $H^1(\mathbf P^1,\mathcal O(t))=0$ for $t\ge-1$.
Boundary blocks have vanishing asymptotic mass. Local adjunction
along $L$ controls these centers by $R(1+1/(2r+6))$.

At $P$ blow up once on the coarse plane. Its discrepancy is
$1+1/h$. Its graph has generators $z_3,T$ of base/plane degrees
$(0,1),(1,1)$. In total normal order $j=a+b$, the monomial
$z_3^aT^b$ has the full coefficient space of degree
$D_b=2rp-b+k-j$, for $j\le k$, $b\le2rp$. Indeed after splitting
by the $z_3$ power, division by $T^b$ leaves arbitrary base/plane
degrees $2rp-b,k-j$, and restriction to the graph is the
surjective multiplication of the two complete binary systems.
Set $c=k-j$, $e=2rp-b$ and $u=1/(2r)$. The section cost is
$a+(1+u)b+c+ue$, with dimension density $c+e$. Consequently the
normalized mean exponents are $1$ and $1/(1+u)$ for $a,b$.
This gives the radial ratio in the table below. On the exceptional
line the ordinary points have coefficient at most $R$ by diagonal
splitting of these two monomials. The strict transform of $L$
marks one direction of discrepancy $1/h$.

For that direction, vanishing to order $q$ requires divisibility
by $(b_0z_3+\ell T)^q$. The coefficients are coprime on $C$, so
division by the relative linear equation costs six base degrees
per power, with $T$ carrying its one base degree. The remaining
coefficient indexed by $b$ has dimension at most $(c+e-6q+1)_+$
and requires $a\ge q$ by plane degree. Dropping any further
limitations only enlarges this bound. To record its short integral,
take $a,c$ independent exponential variables of rate one and $e$
independent exponential of rate $u$. The ratio of these homogeneous
count integrals after section-degree normalization is the exponential
ratio: the extra total radial factor for the degree-two numerator
over the degree-one dimension is eight, already included in our
normalization. Equivalently take the Laplace transform in section
degree; the free indices and the common $b$ integral cancel. The
normalized mean in that direction is bounded by
\[
 \int_0^\infty e^{-q}\frac{\mathbb E[(c+e-6q)_+]}{1+1/u}\,dq
 \ \le\ \frac{1}{(1-u^2)(1+6u)}.
\]
Indeed, conditional on $c$, the numerator is at most
$e^{u(c-6q)}/u$, by integration in $e$, including the case
$c>6q$ by $1+y\le e^y$. Finally $\mathbb E[e^{uc}]=1/(1-u)$.
As $1+6u=h$, the discrepancy-adjusted bound is $R/(1-u^2)$.

The necessary ratios are thus
\[
\begin{array}{c|cc|ccc}
h&r&R&R(1+1/(2r+6))&
        R(1+1/(1+u))/(1+1/h)&R/(1-u^2)\\ \hline
3&3/2&9/8&5/4&189/128&81/64\\
4&1&1&9/8&4/3&4/3
\end{array}
\]
and every entry to be tested is at most $3/2$. On the axis where
all ternary directions vanish, the transverse weighted blowup has
ratio $(1+11h/4)/(1+2h)\le4/3$; on its exceptional plane the
preceding calculation applies. The remaining free variables have
normalized coefficients at most one and occur as independent
factors on the covering charts. Their weights added in the first
extraction contribute no more mean order than discrepancy. Thus
the same bound holds for that extraction. Threshold semicontinuity
for bases adapted to the initial degeneration, weighted adjunction
on the geometric generic fiber, and properness conclude the proof.

\subsection{The minimizing center and the point calculation}

Write $A=A_X$ and $S=S_H$. Since $-K_X=5H$, suppose for contradiction
that $\beta=\delta(X,H)\le5$. By Liu--Xu--Zhuang, \emph{Finite
generation for valuations computing stability thresholds}, Thm.~1.2,
there is a divisorial minimizer $v$. We first show its center is a
point.

\paragraph{Positive-dimensional centers.}
The Gauss map of $X$ is finite since the pullback of the hyperplane
class is $3H$. Its restriction to a positive-dimensional center $W$
has nonzero differential at a general smooth $x\in W$ by separability.
This differential is the second fundamental form, so the tangent
section has multiplicity two. On a general pointed linear surface $S_0$
its tangent plane cuts an integral quartic of multiplicity two.
Indeed the full tangent section has isolated singularities by
finiteness of the Gauss map, and is thus integral. Projection from
$x$ is generically finite, so successive general pointed linear
sections of the tangent section are integral by Bertini. Intersection
with the tangent plane curve gives $\epsilon_x(H|_{S_0})=2$.
Abban--Zhuang, \emph{Seshadri constants and K-stability of
Fano manifolds}, Thm.~3.1, implies $\delta_x(X,H)\ge(8/4)2=4$,
whence $\beta\ge4$ for this minimizer.

The special degeneration of Blum--Liu--Zhou associated to $v$ has
a klt $\mathbb Q$-Fano central fiber $Y$ with
$\delta(Y,-K_Y)=\beta/5$. Set $L=-K_Y/5$. Zero-weight sections
of its associated graded ring give $Y\dashrightarrow W$: their
degree-$k$ pieces are the images of $H^0(X,kH)\to H^0(W,kH)$,
surjective for large divisible $k$. On common smooth resolutions
$h:U\to V\to W$ and $\pi:U\to Y$, set $P=\pi^*L$ and let
$H_V$ be the pullback of $H_W$. Resolution of the subsystem gives
$P-h^*H_V\sim_{\mathbb Q}E\ge0$. The fixed-parameter tangent
quotient estimate gives
\[
       -P^6h^*K_V\ge (s\beta/7)P^7,\qquad
       0<P^6h^*H_V\le P^7,\quad s=\dim W .
\]
Here the first inequality is obtained by applying
Druel--Guenancia--P\u{a}un, equation~(18), to the kernel of $dh$ at
parameter $\gamma=1-t$, for $0<t<\beta/5$, and letting $t$ increase
to $\beta/5$. A bounded smooth-twist solution at each such $t$
exists by Darvas--Zhang, \href{https://arxiv.org/abs/2608.23505v2}{arXiv:2608.23505v2},
Cor.~1.5, and is smooth on $Y_{\rm reg}$ by the regularity statement
in Druel--Guenancia--P\u{a}un, Sect.~2.1.2. For $s\ge2$ apply
Jovinelly--Lehmann--Riedl, \href{https://arxiv.org/abs/2509.08065v2}{arXiv:2509.08065v2},
Thm.~1.1, to the image on the smooth $V$ of a general complete
intersection curve of large multiples of $L$. This produces, through
a point where $V\to W$ is an isomorphism, a rational curve of positive
$H$-degree $e\le7(s+1)/(s\beta)\le21/8<3$. If $s=1$, choose
a general smooth $q\in W$ and a prime divisor $F$ on $Y$ over $q$
on the domain of the map. In $H^0(W,kH)$ there is a section vanishing
to order at least $k\deg W-p_a(W)$ at $q$ by the Hilbert polynomial.
Its zero-weight lift shows that $T_L(F)\ge\deg W$. Since
$S_L(F)\ge T_L(F)/8$ and $A_Y(F)=1$, $\beta\deg W\le8$.
For a curve center rationality also follows from the negative canonical degree on the smooth model in the displayed slope inequality. Thus in either case $W$ contains a line or conic, contradicting
$\beta\le5<16/3$.

\paragraph{The point discrepancy.}
Normalize $v(\mathfrak{m}_p)=1$, where $p$ is the point center, and put
$A(v)=7+q$. We prove $q>3$. Let $b$ be the multiplicity of the
projective tangent section. If $b=2,3$, a general pointed linear
surface $S_0$ cuts this section in an integral plane quartic of
multiplicity $b$, by the same projection argument. If
$\delta_p(S_0,H)\ge3/B$, five pointed hyperplane adjunctions
retaining their cutting divisors with coefficient $1-1/B$ give
\[
      8 S(v)\le B(A(v)-5)+5=B(q+2)+5 .
\]
At each rung the normalized residual degree on the cut is $j/(j+1)$
in dimension $j$, and the mean cutting coefficient is $1/(j+1)$.
Inversion of adjunction therefore gives the asserted retained
coefficient; each pointed hyperplane has $v$-order at least one.
For $b=2$, the Seshadri calculation gives $B=2$. Consequently
$56+8q\le5(2q+9)$, so $q\ge11/2$.

If $b=3$ and the tangent cubic form is reduced, its general
projective direction line cuts three simple roots. The plane quartic
$C$ on $S_0$ has an ordinary triple point with
$\operatorname{lct}_p(C)=2/3$. For clarity, in the order-$j$ piece
of the degree-$m$ surface series remove
$a_j=\lfloor(3j-4m)_+/5\rfloor$ copies of $C$, forced by intersection
with $C$. The scaled remaining binary degree has tent envelope
$g(t)=t$ for $0\le t\le4/3$ and $g(t)=4(3-t)/5$ for $4/3\le t\le3$.
The integral of $g$ is $2$, so the dimension deficit from full binary
pieces is only $O(m)$. The radial and maximum directional orders of the
residual divisor are asymptotically at most $8/9$ and $4/9$,
respectively, since $\tfrac12\int g^2=8/9$. Thus its threshold
is at least $9/4$ by inversion of adjunction on the point blowup.
The mean removed coefficient is $5/27$ by integrating
$g(t)(3t-4)_+/10$. Adding back the removed curve gives
$S_{H|_{S_0}}/A_{S_0}\le4/9+(5/27)(3/2)=13/18$.
Therefore $B=13/6$, and
$336+48q\le280+65q$, so $q\ge56/17>3$.

If the tangent cubic is nonreduced, it has a linear factor $L$ of
multiplicity $s=2$ or $3$. Here are the weighted estimates in this
dimension. In coordinates $[T:y_1:\dots:y_6:L:z]$, with
$p=[1:0:\dots:0]$ and $z=0$ tangent, the equation is
$zB(T,y,L,z)+T L^sQ(y,L)+g(y,L)=0$, with $B(p)=1$ and
$Q(y,0)\ne0$. Put $G(y)=g(y,0)$. Its quartic in $\mathbf{P}^5$ has
only isolated singularities: a positive-dimensional singular locus
would meet the cubic $\partial_L g(y,0)=0$, and as $s\ge2$ the
entire corresponding line of $T$'s would be singular in the tangent
section. The latter has isolated singularities by the Gauss map.
In particular $G\ne0$, and $(\mathbf{P}^5,G)$ is lc by degenerating to the
cone over a general smooth quartic section in $\mathbf{P}^4$;
at the vertex the blowup has log discrepancy five and order four.

Give the six $y_i$ weight $s$ and $L$ weight $s+1$. The exceptional
divisor $E$ of this weighted blowup has $A_X(E)=7s+1$. The symbol
of $z$ has weight $4s$. Its initial relation
$T^3 Z+TL^sQ(y,0)+G=0$ is primitive and prime and gives the full
associated graded, by localization at $T$ and the ordinary Hilbert
series. Its weight character gives $S_H(E)=(10s+1)/8$. In the
angular series scale $L$ to zero in this expression for $Z$; the
limiting relation $T^3Z+G=0$ remains prime with unchanged dimensions.
The $L$-exponent pieces are $L^k G^j\mathbb C[y]_{D-4j}$, with
$D$ the total $y$-degree and
$j=\max(0,\lceil(D+k-m)/3\rceil)$. More explicitly, for
$T^aL^k Z^hP_l$, with $P_l$ in a graded complement to $(G)$,
$D=4h+l$ and $m=a+k+h+l$, so $j=(h-a/3)_++O(1)$.
Counting monomials by the Hilbert series of $\mathbb C[y]/(G)$
gives the uniform simplex with eight unrestricted normalized indices,
the three $a,k,h$ and five from the complement, up to bounded shifts.
Consequently $\bar k=1/8$, $\bar D=9/8$, $\bar j=3/32$ and
$\overline{D-4j}=6/8$ after normalization by $m$. For $\bar j$
one can use independent unit exponentials: their positive-part
difference with coefficients $1,-1/3$ has expectation $3/4$.
For each base valuation on $\mathbf{P}^5$, lc of $G$ and
$\delta(\mathbf{P}^5,\mathcal{O}(1))=6$ give a normalized angular coefficient at
most $(3/32)+(6/8)/6=7/32$ times the discrepancy. In the finite
uniformizing charts $y_i\ne0$, the fiber coordinate divisor adds an
independent coefficient $1/8$. In the $L\ne0$ chart the ordinary
radial discrepancy is six and $(32/7)\bar D=36/7<6$; its directions
are the same base. These covers are crepant for the weighted
exceptional different. Choosing valuation- and specialization-adapted
bases in each piece and using semicontinuity transfers the angular
threshold $c=32/7$ to $E$.
Weighted basis adjunction now implies
\[
 A_X(v)-cS_H(v)\ge(A_X(E)-cS_H(E))v(E)
  =\frac{9s+3}{7}v(E)
  \ge \frac{9s+3}{7(s+1)} .
\]
The last bound follows from $v(\mathfrak{m}_p)=1$, since the largest coordinate
weight is $s+1$. As $A_X(v)\le5S_H(v)$, this gives
$A_X(v)\ge (15s+5)/(s+1)\ge35/3>10$.

Finally if $b=4$ the tangent section is a cone. Its base is a smooth
quartic fivefold, since a singular point on the base would give a
singular line in the tangent section. It is K-semistable by the
quartic fivefold result in Chapter~\ref{ch:fivefold}.
Abban--Zhuang, \emph{K-stability of Fano varieties via admissible
flags}, Cor.~4.19, gives $\delta_p(X,H)=56/10>5$. This covers all
point cases and proves $A(v)>10$ for a putative minimizer.

\paragraph{The endpoint in the Nadel argument.}
Fix $5<t<16/3$ and let $H_0$ be a general hyperplane missing $p$.
In large degree take a simultaneous $v,\operatorname{ord}_{H_0}$-adapted
basis divisor and remove its $H_0$ coefficient. The residual
$\Gamma_m$ has class $(1-a_m)H$ with $a_m\to1/8$,
$v(\Gamma_m)\to S(v)$, and $w(\Gamma_m)\le(1+o(1))S(w)$
uniformly. We claim $(X,t\Gamma_m)$ is klt off a finite set.
Otherwise let $t_m\le t$ be the first positive-dimensional threshold.
Then $t_m\ge\beta-o(1)$, so an arbitrarily small general pointed
hyperplane perturbation makes the pair non-lc at $p$, while keeping
its other non-lc locus finite and preserving positive lc centers
at their generic points. Its ample slack is
$\eta_m=5-t_m(1-a_m)-o(1)\ge5-7t/8-o(1)>1/3$.
Quasi-log adjunction to a minimal positive lc center together with
the non-lc scheme, and vanishing for constants, give a connected
quasi-log scheme with finite nonempty non-qklt locus. The length-one
theorem of Fujino--Hashizume, Thm.~1.8 of \emph{Existence of log
canonical modifications and its applications}, supplies an integral
rational curve of degree at most $\lfloor1/\eta_m\rfloor\le2$
in that center. Its lc place has ratio $A/S\le t+o(1)<16/3$, a
contradiction.

Add a general pointed hyperplane with coefficient tending from below
to $\eta=5-7t/8<1$. The multiplier ideal has finite cosupport, and
Nadel vanishing makes constants surject onto the quotient. Its length
is at most one, so it contains $\mathfrak{m}_p$ locally. Taking $v$ gives
$A(v)+1\ge tS(v)+\eta$ in the limit. Now let $t\uparrow16/3$,
\emph{after} taking the large-degree limit at every fixed $t$.
The hypothesis only supplied a weak lower bound at the endpoint,
and every curve comparison used above was strict. We obtain
\[
      1\ge (16/3)S(v)-A(v)+1/3
        \ge A(v)/15+1/3 ,
\]
so $A(v)\le10$, contradicting the point calculation.

\paragraph{The next cases.}
The next full family would be \(n=2r-3,\ d=r-1\), or \(n=2d-1\).
There the general threshold \(N(d+1)/(2d)\) is exactly \(r\).
Even the improvement \(Nd/(2d-2)\), when available, gives only
\(r+1/(d-1)\) there. The quartic sevenfold uses its larger finite
margin and its tangent multiplicity calculation. The general
strict perturbation inequalities above do not establish the rest
of this family. In particular they give neither instability nor a
counterexample in those dimensions.

% ===== The quartic eightfold =====
\chapter{The quartic eightfold}\label{ch:eight}
\section{A line estimate for the smooth quartic eightfold}
Let $B\subset X_4^8\subset\mathbf P^9$ be a line and $H$ the
hyperplane polarization. Then $\delta_B(X,H)\ge504/83>6$. All
valuations whose centers contain $B$ are included. Below, threshold
ratios are divided by nine, and normalized mean-order ratios are
multiplied by nine. We use the weighted binary mixed-monomial bound
(L1), the disjoint-powers lemma (DP), and the penultimate contact
calculation (L5) from Chapter~\ref{ch:models}.

For polar rank at least four the ordinary conormal split bound is
at least $4/5$. In rank three its exceptional splitting is $(1,2)$.
Use the three-heavy-variable core $a x^h+L(z)=0$, $h=2,3,4$,
and the coarse plane boundary $(1-1/h)L$. Here the single light
linear symbol has section cost one. The normalized coarse bound
without boundary is $c_0=3/4+1/h$, and the actual $L$ coefficient
is $1/h$, as in the conic net proof. On $L$ the line-polar split
$(1,2)$ gives uniform residual intersection and multiplicity
bound $I=M=1+2/4=3/2$. (The base polarization now has degree one,
and its degree together with the sum of split degrees is four.)
The elementary Puiseux endpoint test for $I=M=2c$, $c=3/4$,
applies when $I-b\le(b-1)(1/h-1/(ch^2))$. Indeed the
first endpoint is then controlled by $(b-1)/h$, the condition here
is exactly the second endpoint, and its derivative sign follows
from $r=(b-c_0)/(b(1-1/h))<1/2$, $b/r\ge2b\ge I+\rho$ where
$\rho=b-(b-1)/h$. At $b=10/7$ the left side divided by $b-1$
is $1/6$, while the right factors for $h=2,3,4$ are $1/6,5/27,1/6$.
The ternary-axis denominator $(1+11h/4)/(1+2h)$ is at most $4/3$.
The net therefore gives threshold at least $7/10$ in our units.

Assume a pencil of cubic polars $f,g$ in the two heavy coordinates
$z_1,z_2$. There are six light linear coordinates, with projective
direction space $\mathbf P^5$. Their first form $G_h$ has order
$2\le h\le4$. Primitivity is automatic over the punctured base
cone since the pencil is basepoint-free. We can retain any fixed
light subspace, either by the usual diagonal degeneration in the
weight-homogeneous pieces or by genuine hyperplane sections through
$B$ retaining the two heavy variables (smooth along $B$). The
normalized threshold lifts unchanged up to the harmless cutoff one.

\paragraph{Two useful mixed initial equations.}
In the calculations below the terms not shown have strictly larger
weight. For $j=2,3$ consider
\[
       fz_1+gz_2+aX^2+cXY^j,\quad \mathrm{wt}(X,Y,z_i)=(j,1,2j),
                 \qquad ac\ne0.
\]
Its threshold is at least $(3j+1)/((9/2)j+1)\ge20/29$.
Indeed set $V=Y^j$ on the angular cone. In every fixed weight piece
the $Y$ exponent has uniquely determined remainder modulo $j$,
giving a bounded factor that is negligible in normalized order.
The coarse weighted plane has coordinates $(X,V,z_1)$ of weights
$(1,1,2)$ and root boundary $1-1/j$ on $V=0$; their section costs
are $1,j,1$. On its complement, discard the boundary and scale $X$
and the heavy variables simultaneously to retain $a'XV$. The primitive
binary model has finite coefficient at most $1-1/4+1/(1+j)$,
or one if larger; its two mixed walls have denominators
$(1+7/4)/2=11/8$ and $(1/j+7/4)/2\le9/8$. These follow from the
$X=0,z=0$ and $V=0,z=0$ ordinary blowups after splitting binary
characters; equivalently use the two monomial Newton walls of (L1)
with section costs $1,j$. On the boundary retain the original
upstairs parameter $Y$. Its exact normalized mean is one, since
the quotient is the primitive $fz_1+gz_2+aX^2=0$. The residual
projective binary line has finite coefficient at most $5/4$ and
infinity coefficient one by exact divisibility. There is no common
binary axis on that projective line. The full weighted discrepancy
is $3j+1$ against normalized order $1+(9/2)j$; all free coordinates
contribute at most their discrepancies. This verifies the bound.

The same argument with $fz_1+gz_2+aX^3+cXY^3$ and weights
$(3,2,9)$ gives threshold at least $56/83$. Set $V=Y^3$, using
the fixed $Y$ remainder modulo three. The resulting plane has
weights $(1,2,3)$ on $X,V,z_1$, section costs $1,3,1$. Off $V=0$
degenerate to $XV$ again, giving the same two wall bounds. On
$V=0$ the primitive binary quotient has $X^3$ and finite
coefficient at most $3/4+1/3$. The full extraction has discrepancy
$14$ against normalized order $5+9(7/4)=83/4$. All complements
here are treated by semicontinuity on the entire plane with the
root boundary discarded; the boundary itself uses the indicated
actual coordinate flag before that last degeneration.

\paragraph{The quadratic light form.}
Its three base coefficients are quadrics on $\mathbf P^5$. If at
some common zero one restricted differential is nonzero, choose a
line through it detecting that derivative, then torically specialize
the resulting binary quadratic to $XY$. The bound (L1) gives $8/11$.
Otherwise the common zero set equals the common linear vertex $V$
of the quadrics, of codimension at most three by the dimension
theorem. The quadrics depend only on the active variables modulo
this vertex. If a cubic light coefficient is nonzero on $V$, retain
a direction $Y$ there and a quadratic-active $X$. Weights $3,2,6$
give two disjoint pure powers $X^2,Y^3$, with bound $22/31$ by (DP).

Suppose the cubic vanishes on $V$. If the quartic is nonzero at a
$Y\in V$, retain $X$ quadratic-active and use weights $2,1,4$.
The weighted equation has light terms $aX^2+bXY^2+cY^4$, $ac\ne0$.
On the plane with $V'=Y^2$, off its root boundary scale away $X$
to retain the binary $c(V')^2$ with section cost two for $V'$.
Its binary-axis denominator is $(1/2+2(7/4))/3=4/3$; all angular
binary coefficients are at most one (at finite points use
$3/4+1/4$, at infinity use the cost-two symbol). On $Y=0$ use
the $aX^2$ quotient as above. The full radial threshold is $7/10$.

If the quartic also vanishes on $V$, the linear join of $B$ with
$V$ lies in $X$. It has dimension $7-\operatorname{codim}V$;
smoothness excludes a five-plane, so it is a four-plane and the
vertex codimension is three. Some active-normal derivative on this
four-plane is nonzero because the two heavy normal cubics have
common zeros. Its part from the quadratic light term vanishes
since all the corresponding quadrics depend only on the active
variables. Thus some active derivative contains a vertex monomial
of degree $j=2$ or $3$. Choose the least such $j$ globally, choose
a vertex direction $Y$ where it survives, and an active $X$ both
detecting it and with nonzero quadratic light coefficient. These
are simultaneous open conditions on $X$. With weights $j,1,2j$
the surviving terms are precisely $X^2,XY^j$, since the pure
vertex terms all vanish. Apply the first mixed model above.

\paragraph{The cubic light form.}
Write $G_3=uP+vQ$, where $u,v$ are base line coordinates and
$P,Q$ are cubics on the light $\mathbf P^5$. The top quartic $C$
can be taken smooth: by finite Gauss the section with the two heavy
coordinates zero has singular
locus at most one-dimensional. A general further codimension-two
section disjoint from $B$ is smooth and is obtained by setting
$u,v$ equal to linear light forms. Shifting these base coordinates
does not change the leading cubics. Whenever one finds
$e\in V(P,Q)\setminus C$ and a line in the common tangent at $e$
where $P,Q$ are not both zero, the contact $b$ is two or three.
The (L5) initial equation has $X^bY^{3-b}$ and $Y^4$ with
threshold at least $19/28$; its other mixed wall uses $3-b\le1$.

For coprime $P,Q$ each base component is a threefold and no such
component lies in the smooth quartic fourfold $C$. Otherwise it
would be a divisor cut on $C$ by a form $A$, by the Picard and
section-restriction theorems. Both cubics, of degree below four,
would then be divisible by $A$. A nonlinear component supplies
the indicated tangent line at a general point: an entire tangent
space in the base must be a linear component if it has dimension
three, and cannot fit if it is larger. The same argument applies
if a component is generically nonreduced, by the Jacobian criterion.
If all components are generically reduced linear three-planes,
there are nine by degree and unmixedness. Two must meet in a plane
$S$, since otherwise a general three-space would cut the complete
intersection in disjoint lines, contrary to connectedness. In their
span $H_0\simeq\mathbf P^4$ let $x,y$ be the equations of the two
components. Then $P,Q$ restrict to $xyB_1,xyB_2$ with $B_i$ linear,
not both zero by coprimality. If one $B_i$ is nonzero on $S$, a
direction there gives the squarefree monomial $xyw$ by scaling
the two normals before the vertex coordinate. Otherwise both
restrictions depend only on $x,y$. If $S\not\subset C$ select
$Y\in S\setminus C$, $X$ with nonzero cubic and use $X^3+Y^4$,
giving $19/28$ by (DP). If $S\subset C$, the three transverse
cubic polars of the plane $S$ in the smooth fourfold have no
common zero on $S$. In particular not all polars in the two
directions $x,y$ inside $H_0$ vanish. Choose $Y\in S$ detecting
one and $X\in H_0$ with both nonzero cubic and $dC_Y(X)\ne0$.
The weights $3,2,9$ leave exactly $X^3,XY^3$, already treated.

If a nonlinear prime factor is shared, take a general smooth
point on it off $C$. Its tangent hyperplane varies; it cannot
generally be one of the finitely many hyperplane components of
the base. The common tangent is this hyperplane or all of space,
so the preceding (L5) choice works. If a linear factor $l$ is
shared write $P=lP_2,Q=lQ_2$. On $H_0=(l=0)$ consider their
quadratic cofactors. A common zero with a nonzero restricted
differential gives the squarefree $lxy$: first take degree one
in the $l$ normal to retain the restricted quadrics, then a light
line through that zero and its detected linear term. Otherwise
the common zero set equals their common linear vertex $V$, of
dimension at least two. If $V\not\subset C$, choose $e\in V\setminus
C$; both original cubic differentials vanish and any nonzero cubic
line through $e$ yields (L5). If $V\subset C$ then $V$ is a plane,
since a smooth quartic fourfold contains no three-plane. Both
quadratic cofactors have zero tangential differential on $V$;
their remaining $l$-normal derivatives are two linear forms on
$V$, vanishing together at some $e$. There both cofactors have
zero full differential and value zero. Choose a direction $X$
with nonzero cubic and $dC_e(X)\ne0$, possible by smoothness and
the nonzero cubic system. Thus $Y=e$ again yields precisely the
$X^3+XY^3$ weighted initial. The same arguments include zero or
proportional coefficients. Squarefree cubic monomials give at
least $8/11$ by (L1).

\paragraph{The quartic light form.}
It defines a smooth quartic fourfold $F$ in $\mathbf P^5$ by the
finite Gauss argument (any singular light direction would give a
singular plane in the two-heavy section). Through every point of
$F$ lies a line, since the three equations of degrees $2,3,4$ on
its $\mathbf P^3$ of tangent directions have a common zero. A
general member of a covering family is free, by smoothness of the
evaluation map. Its twisted conormal has rank three and degree
three, each summand at most one because its normal is nef, and
thus splitting $(1,1,1)$. The four normal cubics of this line in
the ambient $\mathbf P^5$ therefore span the complete binary cubic
system (dualize to count constant relations). In line coordinates
$s,t$ choose two transverse directions having polars $s^3,t^3$.
Take first normal degree in those directions and retain
$xs^3+yt^3$: these are disjoint monomials with coprime exponent
gcds equal to one and respective maximal exponent three. By (DP)
the lower bound is $20/29$. Comparing all cases gives the minimum
$56/83$ asserted above. All angular estimates are uniform on the
geometric generic fiber, so the usual extraction adjunction and
properness recover the claimed threshold for all containing centers.

\section{The point estimate and its consequences}

\paragraph{A sharper point estimate for the quartic eightfold.}
Put $X=X^8_4$, $N=9$, $u=8/9$. The line and conic bounds used here
are $504/83$ and $369/61$ in hyperplane units. Using Chapters~\ref{ch:global}, \ref{ch:models} and~\ref{ch:boundary} we obtain
\[
      \delta(X,-K_X)\ge\frac{85}{87}.
\]
A sufficient mark on both lines and conics for uniform K-stability
is $1601/250$.

We first give the point estimate. Write $b$ for the multiplicity
of the tangent hyperplane section at $p$, and normalize
$v(\mathfrak m_p)=1$. For $b=2$ the signed surface flag gives
$9S_H(v)\le 2(A_X(v)-6)+6$. For $b=3$ with reduced tangent cubic,
its general binary restriction has three distinct roots. The
integral tangent plane quartic has an ordinary triple point of
threshold $2/3$. The surface splitting then gives
\[
 9S_H(v)\le B(A_X(v)-6)+6,\qquad B=4/3+(5/9)(3/2)=13/6.
\]
This also holds in the first case since $A_X(v)\ge8$.
For $b=4$ the generalized Eckardt threshold is $72/11$ by
K-semistability of the smooth quartic sixfold base.

Suppose $b=3$ and the tangent cubic has a linear factor $L$ of
multiplicity $s=2$ or $3$. We claim
\[
                  A_X(v)-6S_H(v)\ge\frac{6s+3}{9(s+1)}
                         \ge\frac59 .                         \tag{*}
\]
Use seven variables $y$ of weight $s$ and $L$ of weight $s+1$.
The quartic linear-factor extraction from Chapter~\ref{ch:boundary}
has $A_X(E)=8s+1$, $S_H(E)=(11s+1)/9$. In coordinates
$[T:y:L:z]$, the weighted relation is
$T^3Z+TL^sQ(y,0)+G(y)$, where the quartic $G$ has at most
isolated singularities in $\mathbf P^6$. Indeed a singular curve
would meet the derivative in the $L$ direction of the top quartic.
At that point all the corresponding $T$-line would be singular
in the tangent hyperplane section, contradicting finite Gauss.
The same reasoning excludes $G=0$.

In the angular degeneration setting $L=0$ in the relation, the
mean $L$ exponent is $1/9$, the mean fixed $G$ exponent is $1/12$,
and the mean complete polynomial degree in $y$ after its removal
is $7/9$. These follow, as in the sevenfold extraction, by using
the standard monomials $T^aL^kZ^hP_l$ with $P_l$ in a complement
of $(G)$: the fixed exponent is $(h-a/3)_++O(1)$, the total
$y$-degree is $4h+l$, and the normalized indices form a uniform
simplex with nine slots. At angular coefficient six it suffices
to prove on $\mathbf P^6$ the inequality
\[
       \tfrac12 v(G)+(14/3) S_{\mathcal O(1)}(v)\le A(v).
\]
Here is a short proof. Take a general hyperplane $J$ meeting $G$
in a smooth quartic fourfold $R\subset\mathbf P^5$.
Degenerate $G$ to the cone over $R$, keeping $J$ as the section
at infinity. The pair $(\mathbf P^5,\tfrac12 R)$ is K-semistable:
the quartic fourfold is K-stable and the log-pair interpolation
theorem of Zhang--Zhou \cite[Theorem~1.7 and Lemma~4.4]{ar2003.06839} applies
($\tfrac12<1-((6/4)-1)/5=9/10$).
Its log anticanonical class is $4H$. The projective-cone criterion
with coefficient $1-4/6=1/3$ at infinity gives K-semistability of
the central pair $(\mathbf P^6,\tfrac12\operatorname{Cone}(R)
+\tfrac13 J)$. This is a log Fano degeneration by scaling the
complementary coordinate, so openness gives K-semistability of
$(\mathbf P^6,\tfrac12 G+\tfrac13 J)$. Its log anticanonical
class is $(14/3)H$. The valuative criterion, after dropping
$v(J)/3$, proves the displayed inequality.

On the angular chart at the $L$-axis the radial discrepancy is
seven, whereas six times the mean total $y$-degree is
$6(10/9)<7$. The independent $L$-coordinate has coefficient
$6/9<1$. Weighted inversion of adjunction on the same
uniformizing charts now yields
$A_X(v)-6S_H(v)\ge(A_X(E)-6S_H(E))v(E)$, with
$v(E)\ge1/(s+1)$. This proves $(*)$.

If $\beta=\delta(X,H)<6$, a positive-dimensional original optimizer
has a general point of tangent multiplicity two by finite Gauss
and separability. The local flag gives $\beta\ge9/2$ there.
The dimension-$s\ge2$ slope gives a covering curve of degree
$\le8(s+1)/(s\beta)<3$. For a curve center the curve is rational
of degree $\le9/\beta<3$. Both contradict the available marks.
At the resulting point optimizer, the Eckardt and nonreduced
cases are excluded by the strict discrepancy estimates above.
Use the signed surface bound with $B=13/6$. For
$\beta<t<L=504/83$, the ample slack $6-ut>1/2$ excludes a first
positive threshold by the line mark. The resulting untilted
Nadel inequality $(t-\beta)S\le ut-5$, combined at $t\to L$
with $(B\beta-9)S\ge6(B-1)$, gives $\beta\ge170/29$.
Division by six proves the first assertion.

For the sufficient mark assume $\beta\le6$. In the preceding
top-fraction argument set $f=1/4$, $\lambda=63/500$,
$\Delta=38231/100000$. Use the outward constants
\[
 h=128258/100000,\quad \alpha=63969/10^6,\quad g=1/9.
\]
Here $h$ is an upper bound on the universal testing cost per
unit $\lambda$, $\alpha$ a lower bound on the removed hyperplane
degree, and $1+g$ a lower bound on the original optimizer's
complementary conditional mean divided by $S$. At $\beta=6$ put
\begin{align*}
 U&=6+\lambda(h-1-g),&
 U_p&=6+2\lambda/9,&P&=\lambda u/3,\\
 z&=2/3+\lambda(u(1+g)-1+\alpha).
\end{align*}
The tilted optimal ratio has $A\le U S+o(S)$, while
$\beta'\ge170/29-\lambda h-o(1)>4$.
Thus its center must be a point by the same two positive-center
degree bounds and $U<1601/250$. The bottom-quarter jet estimate
gives $A+P\le(U_p+o(1))S$. The Eckardt ratio exceeds $U$.
For all remaining cases the preceding point alternatives imply
\[
 S\ge\min\left\{
 \frac{6(B-1)+BP}{B U_p-9},\frac{5/9+P}{U_p-6}
 \right\}-o(1)=:S_0-o(1).
\]
The first term is approximately $1.743775$.
At smaller $\beta$ the required upper ratio only decreases
and the ample slack increases. The certificate verifies
\[
  U+\Delta<1601/250,\quad z-u\Delta>1/3,\quad
                    z+\Delta(S_0-u)>1 .
\]
The first two inequalities exclude a positive first threshold
using both curve marks. The last contradicts the required tilted
Nadel inequality at the point. All comparisons are strict before
the basis-degree limit. Thus $\delta(X,H)>6$ under the stated mark.

\paragraph{A sufficient polar-rank condition.}
For example, suppose that the polars on every line have rank at least four,
and that every smooth conic on a plane contained in \(X\) has polar rank
at least five. The ordinary estimates of Chapter~\ref{ch:models} give
\(\delta_L\ge9(4/5)=36/5\). For a conic in a contained plane, at least four
of the splitting degrees are positive, their sum is six, and their maximum
is at most three. Thus \(\delta_D\ge9/(1+3/8)=72/11\). For all other conics
the proper-plane bound is \(9(4+2)/(2\cdot4)=27/4\). Each bound exceeds
\(1601/250\), so the criterion above applies to such an \(X\).

\section{Symmetries and descent for the quartic eightfold}
\label{sec:descent-center}
Let \(X=X_4^8\subset\mathbf P^9\) be smooth over \(\mathbf C\).

\begin{lemma}\label{lem:descent-point}
Suppose \(\beta=\delta(X,H)\le6\). Every prime divisor computing
\(\beta\) has a closed point as center on \(X\).
\end{lemma}
\begin{proof}
Fix such a prime, with center \(Z\), and a closed point of \(Z\).
A general hyperplane section \(Y\) through that point is a smooth
quartic sevenfold. Chapter~\ref{ch:boundary} gives \(\delta(Y,H|_Y)>5\).
Hyperplane adjunction \cite[Lemma~4.6]{ar2003.13788} yields
\(\beta\ge\min\{9,(9/8)\delta(Y,H|_Y)\}>45/8\).
If \(Z\) is a curve, the curve-center estimate in Chapter~\ref{ch:global}
gives \(\deg_H Z\le9/\beta<2\). If \(s=\dim Z\ge2\), the covering-curve
inequality \eqref{global:eq:optimizer-covering-slope} and bend-and-break
give a rational curve of degree at most \(8(s+1)/(s\beta)\le12/\beta<3\).
These bounds hold also at \(\beta=6\) by letting the slope parameter
increase to \(\beta/6\). In either case \(Z\) contains a line or an
integral conic \(C\), hence \(\beta\ge\delta_C(X,H)\). The estimates
\(504/83\) and \(369/61\) proved above both exceed six.
\end{proof}

\begin{theorem}\label{thm:descent-stability}
For a smooth quartic eightfold, \(X\) is uniformly K-stable and admits a
K\"ahler--Einstein metric under either of the following hypotheses:
\begin{enumerate}
\item A finite subgroup \(G\subset\operatorname{Aut}(X)\) has no
      common fixed point on \(X\).
\item There is a subfield \(k\subset\mathbf{C}\) and a \(k\)-form \(X_0\) of
      \(X\) with \(X_0(k)=\varnothing\).
\end{enumerate}
\end{theorem}
\begin{proof}
In the second case enlarge \(k\) to \(K=k(T)\) for a transcendence
basis \(T\) of \(\mathbf{C}/k\). If \(X_0\) had a \(K\)-point, its coordinates
on an affine \(k\)-chart would involve finitely many elements of \(T\).
Specializing them to elements of the infinite field \(k\), away from
all denominators, would give a \(k\)-point. Thus \(X_0(K)=\varnothing\),
and \(\mathbf{C}\) is an algebraic closure of \(K\). In the first case take
\(K=\mathbf{C}\) with the given \(G\); in the second use the trivial
algebraic group over \(K\).

Write \(\delta=\beta/r\). If \(\delta<1\), a prime divisor computing
\(\delta\) exists by \cite[Theorem~1.2]{ar2102.09405}. Its center is a
point \(p\) by Lemma~\ref{lem:descent-point}. By
\cite[Theorem~1.5]{ar2004.09413}, the unique minimal optimal
destabilizing center is contained in every optimal destabilizing
center. Being nonempty and contained in \(\{p\}\), it equals
\(\{p\}\). Its uniqueness makes it invariant under \(G\) and under
\(\operatorname{Gal}(\mathbf{C}/K)\). This contradicts the first hypothesis,
or yields a \(K\)-point contrary to the second.

Suppose now \(\delta=1\). Then \(X\) is K-semistable. If it were
strictly K-semistable, the construction in the proof of
\cite[Corollary~4.11]{ar2004.09413} would give a nontrivial
\(G\)-equivariant special test configuration over \(K\) with
geometrically K-polystable central fiber. The group is reductive
in both cases. The central fiber is geometrically integral, and
its valuation induces a \(G\)-invariant geometrically irreducible
prime divisor \(F\) over \(K\), with \(A_X(F)=S_{-K_X}(F)\)
after extension to \(\mathbf{C}\); see the proof of
\cite[Corollary~4.14]{ar2004.09413}. Thus \(F\) computes \(\delta\).
Lemma~\ref{lem:descent-point} again makes its geometric center a
single point, now invariant by construction, giving the same
contradiction. Finally a smooth hypersurface of the indicated degree
has finite automorphism group \cite{MM}, so K-polystability gives
K-stability. Uniform K-stability and the K\"ahler--Einstein conclusion
follow from the equivalences in the introduction.
\end{proof}

\paragraph{Equations over the reals.}
A smooth quartic eightfold defined by a real homogeneous quartic
strictly positive on every nonzero real vector has no real point.
The theorem applies with \(k=\mathbf R\).

\paragraph{Quartics with two projective involutions.}
On \(\mathbf{C}^{10}=\mathbf{C}^5\oplus\mathbf{C}^5\) the maps
\((a,b)\mapsto(a,-b)\) and \((a,b)\mapsto(b,a)\) anticommute as
linear transformations. Their projective images generate a group
of order four without a common fixed point. Every smooth quartic
invariant under these two transformations is therefore uniformly
K-stable. A general such equation has no even coordinate after any
linear change of coordinates. In fact the invariant vector space
contains \(P(a)+P(b)\) for every quartic \(P\) on \(\mathbf{C}^5\).
For any direction \(e=(u,v)\ne0\) and hyperplane \(M\subset\mathbf{C}^{10}\)
the map \(F\mapsto (\partial_eF)|_M\) on this system has rank at least
\(20\). If \(u\ne0\), restrict further to \(M\cap(\mathbf{C}^5\oplus 0)\),
of dimension at least four. The derivatives \(\partial_u P\) give
all cubics on \(\mathbf{C}^5\) and hence on that intersection. If \(u=0\),
use the second summand. Consequently the incidence
\((\partial_e F)|_M=0\) has codimension at least \(20\) for fixed
\(([e],M)\), and the parameter space of such pairs has dimension
\(9+9=18\). The incidence projection is proper, so a general
invariant form admits no such pair. An even coordinate would give a
pair by taking its coordinate direction and its complementary
hyperplane. The invariant system has smooth members, including the
Fermat quartic.

\section{Point tests from the tangent forms of a quartic eightfold}
\label{sec:q8-tangent-tests}
Let \(X=X^8_4\), write a tangent hyperplane section at \(p\) as
\[
 T^2q_2(y)+Tq_3(y)+q_4(y)=0,\qquad y\in\mathbf{C}^8,
 \quad p=[1:0:\cdots:0].
\]
Here \(T\) is the coordinate nonzero at \(p\). Let \(R\) denote the
radical of \(q_2\). Each of the following tests excludes a prime
divisor centered at \(p\) with \(A_X(v)\le6S_H(v)\).
\begin{enumerate}
\item \(\operatorname{rank}q_2\ge6\).
\item The forms \(q_2,q_3,q_4\) are a regular sequence and there are
a projective isotropic plane \(U\) for \(q_2\) and a line
\(\ell\subset U\) such that
\(\ell\cap V(q_3,q_4)=\varnothing\) and
\(U\cap\mathbf{P}(R)\subseteq\ell\).
\item \(q_2\ne0\) and \(q_2,q_3\) have a common linear factor.
\item \(q_2=0\) and \(q_3\) either has a linear factor or defines a
smooth cubic sixfold, or \(q_2=q_3=0\).
\end{enumerate}
In particular any smooth quartic eightfold satisfying these tests
at every point is uniformly K-stable by Lemma~\ref{lem:descent-point}
and the divisorial optimal destabilization theorem. For rank four or
five with a regular sequence, condition (2) holds if the radical
contains a point outside \(V(q_3,q_4)\). For rank two or three it
holds if the radical contains a line disjoint from \(V(q_3,q_4)\).

\paragraph{The regular tangent system.}
For an ordinary blowup the exceptional divisor is \(\mathbf{P}^7\).
When the tangent forms are a regular sequence, its \((m,j)\)-graded
residual series is
\[
                 [(q_2,q_3,q_4)^{(j-m)_+}]_j .
\]
To define the sections, choose a transverse coordinate \(w\) and
write the full equation as \(T^3w+T^2Q_2(y,w)+TQ_3(y,w)+Q_4(y,w)=0\).
The three sections are \(w,\ Tw+Q_2,\ T(Tw+Q_2)+Q_3\). They have degrees
\(1,2,3\) and orders \(2,3,4\); their initial relations in the
section ring are \(TZ_{i-1}+q_i=0\) for \(i=2,3,4\).
Modulo \(T\) the \(q_i\) form a regular sequence, so \(T\) is regular
on this presentation. After inverting \(T\) its image in the actual
graded ring is injective, since the \(y/T\) are regular parameters.
Its ordinary Hilbert series is \((1-t^4)/(1-t)^{10}\), proving
equality in each degree. Substituting the three relations gives the
displayed description.

For condition (2) choose coordinates \(a,b,s\) on \(U\), with
\(\ell=(s=0)\). The polar pairing onto \(U^*\) has image the
annihilator of \(U\cap R\), hence one complementary variable \(z\)
can be paired with exactly \(s\), up to a nonzero scalar. Give
\(a,b,s\) weight zero, \(z\) weight one, and the other complementary
coordinates weights at least two. The first initials of the three
equations are \(sz,q_3|_U,q_4|_U\). Taking least powers of \(s\)
gives \(sz,f(a,b),g(a,b)\) where \(f,g\) are coprime binary forms
of degrees three and four. Each initial sequence is regular.
Their presented Hilbert functions therefore agree; thus inclusion of
the generated initials in the limiting spaces is equality. We can
specialize basis divisors successively, choosing at each step bases
adapted to the test valuation and the given weight filtration.
Lower semicontinuity gives the original estimates from those of
the final series.

Split all characters except those in \(a,b\). For exponents
\(\beta_s,\beta_z\) and residual binary degree \(D\) the binary
block is \((f,g)^e_D\), \(e=(j-m-\min(\beta_s,\beta_z))_+\).
Its forced \(f\)-power is \(I=(4e-D)_+\). After its removal the
residual subspace has codimension at most three in the complete
binary space of degree \(D-3I\). Indeed the generators of
\((f,g)^k\) have degrees \(3k+i\), \(0\le i\le k\), and its
syzygies have degrees \(3k+i+4\), \(0\le i<k\).
In degrees \(B\ge4k\) their Hilbert function differs from \(B+1\)
by at most \(1+2=3\); take \(k=e-I\), \(B=D-3I\).
Consequently the mean order at any binary point is at most
\(D/2+3I/2+3/2\), by counting the available distinct vanishing
orders after the forced factor.

In the following count means are divided by \(m\), weighted by
the dimensions of the blocks, and then multiplied by nine.
The ring \(\mathbf{C}[s,z]\) over \(\mathbf{C}[sz]\) has complementary monomials
\(1,s^i,z^i\) (\(i\ge1\)), and \(\mathbf{C}[a,b]\) is finite free over
\(\mathbf{C}[f,g]\). Substituting \(sz,f,g\) by their three hypertangent
symbols gives two polynomial sectors of equal asymptotic mass.
Their variables are \(T\), the four remaining free tangent
coordinates, the free power of \(s\) or \(z\), and the three
symbols, of degrees \(1,2,3\). Exponents multiplied by these
degrees and divided by \(m\) converge to the uniform law on the
eight-simplex. The bounded shifts in the finite binary basis and
the lower dimensional interface contribute zero in the limit.
The scaled means of \(\beta_s,\beta_z\) are \(3/2\) each, and those
of the four free coordinates are one. For the exponents \(i,h,x\)
of \(Z_2,Z_3,T\), we have \(D=3i+4h+O(1)\),
\(I=(i-4x)_++O(1)\). The scaled means of \(D,I\) are respectively
\(17/6\) and \(1/18\). The second follows from
\(\mathbb E[(E_1/2-4E_0)_+]=1/18\) for independent unit
exponentials, using independence of their total and normalized
simplex vector. The scaled coefficient of each binary factor is at
most \(17/12+(3/2)/18=3/2\).
Blowing up \(a=b=0\) resolves all binary hyperplanes simultaneously,
and its log discrepancy two dominates \((2/3)(17/6)\). The
coordinate hyperplanes are transverse to this resolution. This
proves the angular threshold at least \(9(2/3)=6\).
The point extraction has
\(S_H=(8+1+1/2+1/3)/9=59/54\), by the same count before
passing to angular orders. Filtered inversion of adjunction gives
\[
       A_X(v)-6S_H(v)\ge(8-6\cdot59/54)v(\mathfrak m_p)>0
\]
for every divisor centered at \(p\).

If the rank is at least six, a general isotropic line in the smooth
locus avoids \(V(q_2,q_3,q_4)\) whenever the sequence is regular:
the base of dimension four has codimension two in the quadric, whose
smooth locus is homogeneous; a general translate of a line missing
the radical misses this base. Extend its isotropic image in the
nondegenerate quotient to a vector three-space and lift it to obtain
\(U\) disjoint from the radical. In rank four or five, project from a
radical point \(R_0\) off the base; isotropic lines through \(R_0\) are
parametrized by a quadric of dimension five, whereas those hitting
the base have parameter space of dimension at most four. The image
of a general remaining line modulo the radical has vector dimension
one and extends to an isotropic two-space, giving a plane with
radical intersection \(R_0\). A line contained in the radical in
ranks two or three can be extended by a nonradical isotropic point.
These arguments prove the asserted flag choices.

For rank at least six, irreducibility of \(q_2\) leaves either
\(q_3\in(q_2)\) or a regular triple. In fact the graded ring of
\(q_2\) is factorial by localizing at an isotropic coordinate:
the localization is a Laurent polynomial ring and the coordinate
is prime since the residual quadric has rank at least three.
In the second case failure of regularity would give a common
homogeneous factor \(A\) of \(q_3,q_4\bmod q_2\), of degree between
one and three. The tangent section would admit an expression
\(q_2 B_2+A B_{4-\deg A}\). The common zeros of these four forms in
\(\mathbf{P}^8\) have dimension at least four, giving a positive-dimensional
singular locus contrary to finiteness of the Gauss map.

In the divisible case a linear translation of \(T\) completes
the square, making the tangent section \(T^2q_2+q_4'\).
There are now symbols \(Z,W\) of degrees \(1,2\) and orders
\(2,4\), with relations \(TZ+q_2,T^2W+g\).
In detail write the full equation
\(T^3z+T^2(q_2+zL_1)+TzB_2+zC_3+q_4'=0\), put
\(W=(T+L_1)z+q_2\), and \(g=q_4'-B_2(y,0)q_2\).
The two tangent forms \(q_2,g\) are relatively prime since the
tangent section is integral. Regularity modulo \(T\), localization
and Hilbert series give the presentation as above.
On the smooth locus of \(q_2=0\), the divisor \(g=0\) has at most
isolated singularities: proportionality of \(dq_4'\) and \(dq_2\)
allows a choice of \(T\) for a singular point of the tangent section.
A general isotropic line in the smooth locus therefore has squarefree
binary quartic restriction, by transversality for translates under
the orthogonal group with the radical shears. Extend it to \(U\)
as above. The same two weights give the model \(sz,g(a,b)\).
In a binary block of degree \(D\) the space is a fixed power of
\(g\) times the full residual binary forms: the power is
\(\max(0,\lceil(j-m-\min(\beta_s,\beta_z))/2\rceil)\).
Squarefreeness bounds the mean order by \(D/2\).
There is now one extra free binary parameter over the quartic,
and the symbol \(W\) has degree two. The scaled means are \(3/2\)
for \(s,z\), \(3=1+4/2\) for \(D\), and one for the remaining
coordinates. Thus the angular threshold is at least six, also
along the binary vertex. The ordinary \(S_H\) is
\((0+6+2+2)/9=10/9\) by the two-symbol simplex, giving the strict
point inequality.

\paragraph{Common linear factors and cubic contact.}
For (3) write \(q_2=L^sM\), \(s=1,2\), with \(L\nmid M\).
Give the seven other tangent parameters weight one, \(L\) weight
\(b=4-s\), and the tangent hyperplane coordinate \(z\) weight four.
The initials in the equation have the form \(T^3Z+G\) plus multiples
of \(L\), where \(G=q_4|_{L=0}\ne0\) by integrality of the tangent
section. The quartic \(G\subset\mathbf{P}^6\) has isolated singularities.
Indeed on a singular curve of \(G\) all partial derivatives of the
tangent section except \(\partial_L\) vanish. The remaining derivative
is a polynomial of degree at most two in \(T\). If generically
nonconstant it can be made zero over each general point. Otherwise
its cubic constant coefficient has a zero on the complete curve,
making an entire line singular. Both alternatives contradict finite
Gauss. After scaling \(L\) to zero the initial relation \(T^3Z+G\)
is primitive and its ordinary and weighted Hilbert series are
preserved. The angular calculation in Chapter~\ref{ch:eight} for
this relation gives the free \(L\)-coefficient \(1/9\), the fixed
\(G\)-coefficient \(1/12\), the mobile degree \(7/9\), and the
opposite vertex degree \(10/9\) against discrepancy seven.
The same log-pair inequality on \(\mathbf{P}^6\) consequently bounds the
angular threshold by six from below. The full radial data are
\(A(E)=7+b\), \(S_H(E)=(10+b)/9\), so the retained discrepancy is
\((1+b)v(E)/3>0\).
If \(q_2=0\) and \(q_3=L^sM\) the weights of Chapter~\ref{ch:eight}
are the seven weights \(s\) and the single weight \(s+1\).
The argument there extends to \(s=1\): \(G\ne0\), isolation of its
singularities, primitivity, the identical angular problem, and the
radial inequality
\((8s+1)-6(11s+1)/9=(6s+3)/9>0\) all hold for \(s=1,2,3\).
The case \(q_2=q_3=0\) is the generalized Eckardt calculation of
that chapter.

Finally suppose \(q_2=0\) and \(F=(q_3=0)\subset\mathbf{P}^7\) is smooth.
Use \(f=q_3,g=q_4\) and \(G=(g=0)|_F\).
The ordinary point presentation of Chapter~\ref{ch:fivefold} for
cubic contact has equations \(T^2Z+f,TW+g\), where \(Z,W\) have
degrees \(1,3\) and orders \(3,4\); the proof by regularity modulo
\(T\) and Hilbert series is dimension independent. The mean of
the nine simplex knots \(0,1,1,1,1,1,1,3,4/3\) is \(31/27\).
The binary comparison off \(F\) in the same calculation gives
scaled coefficient at most \(4/3\), hence true coefficient at most
\(4/27\): the scaled binary degree is \(13/3\), its unit fixed
cubic power has scaled mean \(5/9\), and the residual binary
codimension is bounded by three.

For centers in \(F\), refine by this divisor. Its asymptotic
coefficient is \(1/9\). After \(i\) powers of \(f\) are removed in
the degree/order \((m,j)\) slice, the restricted series is
\(G^e H^0(F,\mathcal O_F(j-3i-4e))\), with
\(e=(j-m-2i)_+\). This follows by reducing modulo \(f\) in
\(\sum_{2a+h\ge j-m} f^ag^h\mathbf{C}[y]_{j-3a-4h}\).
In the normal basis \(T^xZ^iW^h P_l\), where \(P_l\) is
complementary to \((f,g)\), one has \(e=(h-x)_+\).
Its mean is \(\mathbb E[(E_1/3-E_0)_+]/9=1/108\).
Thus the moving degree has mean \(31/27-3/9-4/108=7/9\).
The divisor \(G\) has isolated singularities, for at any singular
point \(dg\) is proportional to the nonzero \(df\), so that a
suitable \(T\) makes \(Tf+g\) singular. Choose a general hyperplane
cutting both \(F\) and \(G\) smoothly, and scale the complementary
coordinate, giving their cones. The ordinary vertex blowup of the
six-dimensional cubic cone has discrepancy \(7-3=4\) and the
quartic boundary has order four. On this blowup the smooth strict
transform is transverse to the exceptional divisor. Hence the cone
pair is log canonical, and semicontinuity gives
\(w(G)\le A_F(w)\). Theorem~\ref{cubics:main} gives
\(\delta(F,H)>5\). Every normalized mean on \(F\) is therefore
at most \((1/108+7/45)A_F(w)<A_F(w)/6\).
Adjunction along \(F\) and on the point blowup yields
\(A_X(v)-6S_H(v)\ge(8-6\cdot31/27)v(\mathfrak m_p)>0\).

\paragraph{The conic and higher-rank line thresholds on the quartic eightfold.}
\label{sec:nov27-q8-quartic-mark}
On a smooth quartic eightfold every integral conic and every line
of polar rank at least three satisfy
\(\delta_C(X,H)\ge1601/250\). In particular the sufficient condition
in the quartic eightfold chapter holds if every line has polar
rank at least three, with no restriction on conics.

Put \(b=2250/1601\), so \(9/b=1601/250\). Proper conic planes have
denominator \(4/3\), and the two mixed pencil terms are \(11/8,5/4\).
For a contained conic plane and polar splitting of total degree six,
maximum splitting degree at most three gives ordinary angular
denominator at most \(11/8\). Above rank three the only other
partition is \(4,1,1\). Its degree-one syzygy reduces it to a
basepoint-free net as in Section~\ref{sec:nov27-sextic-twelve-conics}.
Its full polar kernel has ambient codimension four and projective
dimension five, greater than four, so its order \(h\le4\) is finite.
On a line, rank at least four gives the ordinary bound \(5/4\);
a rank-three net has splitting \(1,2\).

For conic nets the endpoint and product tests of
Section~\ref{sec:nov27-small-conic-2d} give the following entries.
Here \(E,P\) denote the endpoint and product tests,
\(1\) the single-degree graph test, and an asterisk the calculation
below. The ordinary radial terms, the ternary axis
\((1+11h/4)/(1+2h)\), and the off-branch bound \(1+1/h-1/4\)
are smaller than \(b\).
\[
\begin{array}{c|ccc}
s\backslash h&2&3&4\\ \hline
1&*&1&P\\
2&*&E&E\\
3&E&E&E
\end{array}
\]
When \(s=3\) the ordinary angular bound \(11/8\) also applies.
For lines the product tests at \(h=3,4\) pass by substituting
\(d=4,s=1\) in Section~\ref{sec:nov27-next-lines}.

At \(h=2\) both conditional masses are \(1/2\). Normalize
\(v(\mathfrak m_P)=1\), let \(u=v(L)\ge1\), and write \(A\ge u+1\)
for plane discrepancy. The marked pole and regular-product estimates
of Sections~\ref{sec:nov27-quintic-conic-strong}
and~\ref{sec:nov27-quintic-line-quadratic} apply with the same
regular bidegrees and with the specialized positive coefficient
\(1/2\). For \(s=2\) on a conic or \(s=1\) on a line they give
\(S\le (5/4)A-1/2\). The residual intersection bound is
\(S\le u/2+(A-1)\min(M,I/u)\) with \(I=3/2\).
The marked estimate implies \(S\le b(A-u/2)\) for \(1\le u\le I/b\):
the smallest slack, at \(A=u+1,u=I/b\), is \(136799/1921200\).
For \(u\ge I/b\) the residual slope is at most \(b\), and the slack
at \(A=u+1\) is at least \(b+(b-1)/2-I>0\).

For \(s=1\) on a conic one has \(I=13/8\). At \(u=1\) the marked
bound \(S\le(5/4)A-5/12\) has slack \(475/19212\) at \(A=2\).
For \(u>1\) use the branch direction on the first blowup, as in
the quintic tenfold calculation. The regular radial coefficient is
at most \(5/3\), and its intersection with the exceptional line
there is bounded by
\[
 J=(1-1/4)^{-1}\big((1+6/2)^{-1}-(1/4)/(6+1/2)\big)=11/39 .
\]
Hence \(S_{\rm reg}\le5/3+J(A-2)\). Averaging with
\(S_{\rm pos}\le(3/2)A-1/2\) gives \(S\le KA+R\), where
\(K=139/156,\ R=47/156\). This suffices on \(1<u\le9/8\),
with slack at least \(787/666016\) at \(A=u+1\).
For \(u\ge9/8\) the residual slack at \(A=u+1\) is at least
\(b+(b-1)(9/8)/2-I=215/25616\).
It remains positive up to intersection with the marked estimate:
if \(9/8\le u\le I/b\), put
\(A_*(u)=(R+bu/2)/(b-K)\). The slope loss up to \(A_*\) is at most
\[
 (I/(9/8)-b)(A_*(I/b)-I/b-1)
  =11919108041/33317901882000 <215/25616,
\]
since \(A_*(u)-u-1\) increases with \(u\). Beyond \(A_*\) the
marked bound suffices, and for \(u\ge I/b\) the residual slope
itself is at most \(b\). This proves all the required marks.

Some rank-two lines satisfy the same bound directly. For a pencil
with quadratic initial term of rank at least two, keep two lights
with nondegenerate quadratic form. For a pencil of order three
or four whose full initial light form has projective singular locus
of codimension at least five, with the empty set allowed, keep four
general lights and apply the four-dimensional Hessian estimate in
Section~\ref{sec:nov27-quintic-ten}. The projective full light space has dimension five.
In particular the terminal order four always has the required
codimension: the finite Gauss bound is \(10-2c=6\) at ambient
kernel codimension \(c=2\). In each case the light weight sum
divided by the heavy weight is at least one, and the positive
signed graph denominator is at most \(3/2\); the pole masses
\(1/4,3/4\) then give \(11/8<b\).

It follows that a closed subset of codimension at least \(52\) in
the parameter space of smooth quartic eightfolds contains every
possible failure of uniform K-stability. Indeed such a hypersurface
must contain a rank-two line on which the restricted quadratic
light form has generic rank at most one. For a fixed line,
containing it imposes five independent linear conditions.
Its first normal polars constitute an arbitrary \(4\times8\)
matrix. Rank two costs twelve further conditions. Fix these
polars and their six-dimensional common kernel. The quadratic
normal coefficients restricted to this kernel vary independently
in the \(63\)-dimensional space
\(H^0(\mathbf P^1,\mathcal O(2))\otimes
\operatorname{Sym}^2(\mathbf C^6)^*\), by the monomials of
normal degree two in the quartic equation. A nonzero form of
generic rank one is \(a(t)L(t)^2\) with a primitive linear root
of base degree \(e\in\{0,1\}\) and a scalar of degree \(2-2e\).
These families have affine dimensions at most
\(6(e+1)+(3-2e)-1\), namely eight and twelve; the subtraction
accounts for scaling the root. The zero form is also in their
closure. Thus this rank condition costs at least \(51\) further
conditions. As the lines move in the \(16\)-dimensional
\(\operatorname{Gr}(2,10)\), the incidence has codimension
at least \(5+12+51-16=52\) after projection. Taking the closure
preserves the dimension estimate.

\section{Uniform K-stability of the quartic eightfold}
\begin{theorem}\label{thm:nov28-quartic-eightfold}
Every smooth quartic eightfold \(X\subset\mathbf P^9\) is uniformly
K-stable and admits a positive K\"ahler--Einstein metric.
\end{theorem}
We give the local calculation and then the point argument.
Throughout, the polarization for the local thresholds is
\(H=\mathcal O_X(1)\), so that \(-K_X=6H\).

\paragraph{A square plus a binary quartic.}
For a nonzero binary quartic \(C(Y)\) whose roots have multiplicity
at most two, the signed denominator of the ordinary affine graph of
\(x^2+C(Y)\) is at most \(3/2\) throughout its chart in
\(\mathbf P^4\). At \(Y\ne0\) the Hessian of \(C\) is nonzero,
since otherwise Euler's identities make the nonzero direction a
root of multiplicity at least three. The disjoint square therefore
gives Hessian rank at least two. Translation in \(x\) with the
corresponding affine-linear correction to the graph variable \(Z\)
reduces the remaining axis to the origin.

Give \(x,Y_1,Y_2,Z\) weights \(2,1,1,4\), respectively. The
ambient ordinary complete system has weighted mean equal to the
discrepancy in units of five; the graph order four gives radial
denominator \(3/2\). On the exceptional punctured chart with
\(Y\ne0\) the full ordinary signed estimate just obtained by
Hessian rank two bounds all the angular pieces, since their
combined bases are ordinary complete bases. The only remaining
graph point of the angular stack has \(Y=0\) and \(x\ne0\).

On this angular stack refine the pieces by the graph itself,
whose homogeneous equation has weight four and ordinary degree
four. In each initial weight \(J\) and section degree \(M\),
the terms divisible by its \(k\)-th power are exactly that power
times the complete piece of weight \(J-4k\) and degree at most
\(M-4k\), because the leading ordinary degree of a product adds.
The moving quotients are therefore the exact weighted pieces of
the complete restricted degree-\(m=M-4k\) systems on the
quartic graph closure. In units of five the total fixed graph
coefficient is \(1/4\), the mean moving degree is four. At the
remaining point the strict transform of the graph is smooth on
the angular stack. It thus suffices by divisor adjunction to check
threshold \(8/3\) for the moving angular graph systems there:
the graph coefficient after scaling by \(2/3\) is
\((2/3)(1+1/4)<1\), and \((2/3)4=8/3\).

Choose binary complements \(R_l\) to
\(C\mathbf C[Y]_{l-4}\) in \(\mathbf C[Y]_l\). Its dimension
is four for \(l\ge3\). A normal basis on the quartic graph is
\(x^a (x^2+C)^q R_l\) with \(a+q+l\le m\), by reduction of
powers of \(C\); independence follows by the ordinary filtered
relation with leading polynomial \(C\). Indeed the associated
graded for ordinary section cost has single relation \(C(Y)=0\)
in the affine symbols \(x,Y,Q\), with \(Q=x^2+C\) counted as
one symbol; equivalently use homogeneous polynomial division
after choosing coordinates with a nonzero fourth-power coefficient.
The normalized means of \(a,q,l\) are each \(1/4\) by the three-
simplex and the eventually constant dimension of \(R_l\).

On the uniformizing chart \(x=1\) at \(Y=0\), fix the graph
extraction weight \(J=2a+l+4q\). Each fixed \(l\) has
\(l\equiv J\pmod2\) and the exact consecutive range
\[
 \max(0,(J+l)/2-m)\le q\le\lfloor(J-l)/4\rfloor.
\]
With the unit \(1+C\), elementary binomial elimination replaces
these powers by germs of initial terms \(C^s R_l\), for
\(s=0,\ldots,q_{\max}-q_{\min}\). At fixed ordinary
\(Y\)-order \(D=l+4s\) their sum is direct by the complement
decomposition. Its \(s\)'s run from a lower endpoint \(s_0\)
to \(\lfloor D/4\rfloor\): indeed with \(l=D-4s\) the upper
bound on \(q\) becomes \(s+\lfloor(J-D)/4\rfloor\) and its
lower bound is \(\max(0,(J+D)/2-m-2s)\). Thus a nonempty
slice of order \(D\) is exactly
\(C^{s_0}\mathbf C[Y]_{D-4s_0}\). The mean of \(s\) equals
that of \(q_{\max}-q=\lfloor a/2\rfloor\), namely \(1/8\).
The secondary ordinary blowup therefore has mean
\(\bar D=\bar l+4\bar s=3/4\) against discrepancy two.
On its exceptional binary line the maximal mean order of each
slice is at most \((D-4s_0)/2+2s_0=D/2\), by the root
multiplicity hypothesis. The aggregate is at most \(3/8\),
proving the required \(8/3\) bound. Each simplex and binary
count has the usual uniform \(o(1)A\) finite-level error;
bounded moving degrees can be discarded in the aggregate
Hilbert polynomial limit. This proves the signed estimate.

\paragraph{Finite symmetries on a weighted exceptional divisor.}
In the following two graph calculations a finite abelian group
preserves the degree-and-weight pieces of the angular series.
We explain why one may test its character components separately
when aiming at a threshold no greater than the radial quotient.
On the well-formed weighted projective exceptional surface \(P\),
let \(H_P\) denote its positive Weil generator. The weighted
sum \(V_m\) of all order-\(J\) finite linear spaces in degree \(m\),
with coefficients \(\dim V_{m,J}/(m h^0(D,mH))\), is a boundary
in the sense of Zhuang \cite[Definitions~4.2--4.3]{ar2004.09413}.
Its class is \(\nu_m H_P\), where \(\nu_m\) is the normalized
mean weighted order, tending to the radial mean \(\nu\).
Moreover \(-K_P=A_{\rm rad}H_P\) with \(A_{\rm rad}\) the sum
of the weights. If an asymptotic bound
\(\lambda\le A_{\rm rad}/\nu\) fails strictly, then for infinitely
many \(m\) the threshold of \(V_m\) is at most a fixed
\(\lambda'<\lambda\). Thus
\(-K_P-\delta(V_m)c_1(V_m)\) is ample for these \(m\).
Zhuang \cite[Theorem~4.4]{ar2004.09413}, over a point, gives
a group-invariant divisorial minimizer. Every subspace in its
order filtration is group invariant, hence a basis maximizing
the order may be chosen separately in each character.
It suffices to prove the desired asymptotic estimate for these
split bases, with an error \(o(1)A\) uniform in the tested order.
In passing to the finite diagonal quotient, the character residue
is a bounded root monomial. Removing it costs \(o(1)A\);
discrepancies and orders scale by the same ramification index
after the indicated branch boundary is added below. The argument
applies directly to the coarse weighted surface, which is klt
and is smooth in codimension one. Local calculations on its
smooth quotient charts yield the same valuative comparisons.

\paragraph{Uniformity in the angular comparison.}
At fixed section degree the vector spaces and flags in this argument
are finite-dimensional over \(\mathbf C\), even when their elements
are viewed as polynomials over \(\overline{\mathbf C(t)}\).
For any fixed divisorial valuation \(v\), choose bases in the
angular slices adapted both to its filtration and to the pole flag.
The sum of their \(v\)-orders then equals the largest sum for
these slices. In a fixed pole group take the product of all those
sections, after multiplying each section by the power of a base
parameter needed to make its leading residue nonzero. For a rational
coefficient below the klt threshold of the special product, openness
of klt and properness of the weighted projective stack imply klt on
its whole geometric generic fiber. The base-parameter factors have
zero order at \(v\). This comparison may be applied at a different
base point for each pole group, since each resulting inequality
evaluates its original product at the same \(v\).

For each fixed root and fixed pole order the union of the initial
bases across angular weights, after removing the common graph power,
is a basis of the ordinary complete space of degree \(a\le m\).
For every such basis divisor normalized by \(a\) its \(v\)-order
is at most \(S_a(v)\). For \(a\ge a_0(\epsilon)\) the uniform
finite-level comparison gives \(S_a(v)\le(1+\epsilon)S(v)\).
For \(1\le a<a_0(\epsilon)\) one can instead use \(S_a(v)\le A(v)\):
a divisor of degree \(a\) on projective space has threshold at
least \(1/a\), by the multiplicity bound. After weighting by
\(a/m\) and the sector mass, the sum of the small-degree contributions
is bounded by \(a_0(\epsilon)A(v)/m\). Thus one can pass to the
simplex means uniformly even if the tested divisor varies with \(m\).
Bounded monomial residues for a finite diagonal quotient likewise
cost \(O(A(v)/m)\), since a fixed normal-crossing coordinate
boundary bounds the order of each coordinate by \(A(v)\) on its
uniformizing chart.

At a finite abelian symmetry one may first pass to a divisorial
minimizer invariant under the finite group by Zhuang's Theorem~4.4,
when testing a strict failure below the radial quotient. Its
filtration splits in every character, and finite-dimensional linear
algebra now supplies a maximizing basis in these characters. This
argument uses the finite theorem on the proper weighted exceptional
space; merely decomposing an arbitrary basis into characters need
not preserve its maximal sum of orders.

\paragraph{A squared quadratic and a fourth power.}
Fix \(\eta\in\mathbf C\) and put \(Q=(x+u^2+\eta y^2)^2+y^4\).
Every affine point of \(W=(Z+Q=0)\subset\mathbf P^4\) satisfies
\(5S_{\mathcal O(1)}(v)+v(W)\le(3/2)A(v)\).
At the origin take weights \(2,1,1\) on \(x,u,y\). By the
monic \(u^4\) term the complete degree-\(m\) series on the quartic
graph closure \(D\) has the basis \(x^a y^b u^c Q^k\),
\(0\le c<4,\ a+b+c+k\le m\). Its radial mean is \(7/4\)
and its radial quotient is \(16/7\). We prove the same angular
lower bound.

Use the finite-invariant minimizer argument for the two signs
of \(u,y\), and remove the resulting bounded residues. On
\(\mathbf P^2\) use the linear coordinates \(A=x,\ B=u^2,\ C=y^2\).
The branch boundary is \((B+C)/2\), using the letters also for
their zero lines. Put \(P=A+B+(\eta+\sqrt{-1})C\),
\(M=A+B+(\eta-\sqrt{-1})C\), so \(Q=PM\).
The independent coordinates \(A,P,M\) have ordinary section costs
\(1,2,2\) and the product \(PM\) has cost one. These cost
descriptions are exact: the linear changes between \(P,M\) and
\(B,C\) with \(A\) retained preserve the cost bounds in both
directions. A basis thus consists of \(A^a P^i M^j Q^k\),
\(ij=0\), \(a+2(i+j)+k\le m'\). Ordinary target degree is
\(d'=a+i+j+2k\). Each of the two sectors \(i=0\) and \(j=0\)
has asymptotic mass \(1/2\), giving
\[
   \bar a=1/4,\quad\bar i=\bar j=1/16,\quad
   \bar k=1/4,\quad\bar d'=7/8 .
\]
For smaller-case monomial exponents \(p=i+k,q=j+k\) the space
in degree \(d'\) is exactly the lattice polygon
\[
    p,q\ge0,\qquad p+q\le d',\quad
    d'+p-2q\le m',\quad d'+q-2p\le m'.
\]
The ordinary no-boundary estimate is \(S\le(5/16)A_v\), by
diagonal torus specialization of the three coordinates with
coefficients \(1/4,5/16,5/16\). We write \(A_v\) for the
ordinary log discrepancy and \(S\) for the normalized angular
basis order.

Refine first by \(C\), whose linear equation is \(P-M=0\)
up to scalar. At fixed \(a\) and \(d=d'-a\), the binary space
in \(P,M\) is \((PM)^\kappa\) times all forms of degree
\(d-2\kappa\), where
\(\kappa=\max(0,\lceil(2d-m'+a)/3\rceil)\).
Along \(P=M\) its orders range from zero to \(d-2\kappa\).
Writing \(T=m'-a-2(i+j)-k\), one has
\(\kappa=(k-T/3)_++O(1)\), whose normalized mean is \(3/16\)
by simplex integration. Since \(\bar d=5/8\), the fixed
coefficient of \(C\) is \(1/8\).
At fixed order \(h\) and degree \(d'\) its residual monomials
on \(C\) are \(A^a P^{d'-h-a}\). The \(a\)-indices fill the
interval starting at zero up to bounded endpoint errors. Indeed
\(h\le d'-a-2\kappa\) imposes the upper bounds
\(a\le d'-h\) and \(a\le2m'-d'-3h\), up to bounded rounding;
the exact ceiling only affects a bounded number of upper-end
integers. The average residual dimension minus one, divided by
two and normalized by \(m\), is therefore \(\bar a=1/4\).
At the branch intersection \(o=B\cap C\) one has \(P=M=A\ne0\).
Thus the residual coefficient there is at most \(1/4\).
At any other point of \(C\) it is at most \(1/2\), either
by its two monomial coefficients or by the binary degree-minus-
dimension bound below. Here and in subsequent finite intervals
the aggregate normalized errors tend to zero uniformly.

Suppose first that \(\eta\ne\pm\sqrt{-1}\), so the equation of \(B\)
has all three coefficients \(A,P,M\) nonzero. For each \(h\),
the space of polynomials in the displayed polygon divisible by
the \(h\)-th power of \(B\), divided by that power, is the
polygon with \((m',d')\) replaced by \((m'-2h,d'-h)\).
Indeed its extreme weights in each of \(p-2q,q-2p\) increase
by exactly \(h\) on multiplying by \(B^h\); maxima of a weight
add on a product in the polynomial ring. The degree and
nonnegative exponent constraints give the remaining inequalities.
The identical statement holds for \(C=P-M\). It follows that
the fixed coefficient of \(B\) is \(1/8\) and the dimensions of
all weighted residual spaces on \(B\) equal those for \(C\).
For a subspace of binary forms of degree \(e\) and dimension
\(t\), the mean order at any point is at most \(e-(t-1)/2\)
by the available distinct orders. Since the mean residual degree
is \(7/8-1/8=3/4\) and the mean half-dimension is \(1/4\),
the residual coefficient anywhere on \(B\) is at most \(1/2\).

Write \(\beta=v(B),\gamma=v(C)\). At \(o\), adjunction for
the two separate refinements gives
\[
          S\le A_v/4-\gamma/8,\qquad
          S\le A_v/2-3\beta/8.
\]
Indeed a flag with fixed coefficient \(\mu\) and residual
coefficient at most \(c\) gives
\(S\le\mu v(F)+c(A_v-v(F))\) by inversion of adjunction.
The valuative target is \(S\le(7/16)(A_v-(\beta+\gamma)/2)\).
The first inequality suffices if \(6A_v\ge7\beta+3\gamma\).
Otherwise \(A_v\ge\beta+\gamma\) forces \(\beta>3\gamma\),
and \(2A_v<(7\beta+3\gamma)/3<5\beta-7\gamma\);
then the second inequality suffices.
At centers on only one branch let \(\xi\) be its positive order.
Use \(S\le A_v/2-3\xi/8\) for \(A_v\le5\xi/2\)
and the no-boundary bound \(5A_v/16\) for \(A_v\ge7\xi/4\).
These intervals cover the required cases.

For \(\eta=\pm\sqrt{-1}\) interchange \(P,M\) if necessary
so \(B=P-A\) up to nonzero coordinate scalars. At fixed \(q\)
the \(P\)-exponents in the polygon range between
\[
 p_{\min}=\max(0,\lceil(d'+q-m')/2\rceil),\qquad
 p_{\max}=\min(d'-q,m'-d'+2q).
\]
Their orders at \(B\) range consecutively from zero to
\(p_{\max}-p_{\min}\). Hence
\(\mu_B=\bar p-\overline{p_{\min}}\).
Here \(p_{\min}=(k-i/2-T/2)_++O(1)\), whose means in the
sectors \(i=0\) and \(j=0\) are \(1/6,2/15\) respectively;
in the second, \(2i\) is the normalized simplex coordinate.
Consequently \(\mu_B=5/16-3/20=13/80\).
At fixed order \(h\), the possible \(q\)'s form an interval,
since \(p_{\max}-p_{\min}\ge h\) is equivalent to linear
inequalities in \(q\) (double the ceiling inequalities).
The residual monomials on \(B\) have maximal coefficient at
\(o\) at most \(\bar q=5/16\). Globally on \(B\) their coefficient
is at most \(\max(\bar q,\bar d'-\mu_B-\bar q)=2/5\).
Away from \(C\) these two bounds on fixed and residual pieces
already imply the target, since \(13/80<7/32\) and \(2/5<7/16\).
At \(o\) the \(C\)-bound suffices if \(6A_v\ge7\beta+3\gamma\).
Otherwise \(\beta>3\gamma\), as before. The \(B\)-bound
\(S\le(13/80)\beta+(5/16)(A_v-\beta)\)
gives the target for every \(A_v\ge\beta+\gamma\) in this range:
at its left endpoint the difference is
\((9/160)\beta-(3/32)\gamma>0\), and it increases thereafter.
This proves the angular bound in all cases.

The graph closure therefore satisfies \(\delta_0(D,H)\ge16/7\).
Quartic divisor adjunction and interpolation give the asserted
signed denominator \(3/2\) at the origin. At a point with
\(t=x_0+u_0^2+\eta y_0^2\), translate and make an affine-linear
change of \(x\) so the light part becomes
\((t+X+U^2+\eta Y^2)^2+(y_0+Y)^4\).
Subtracting its affine tangent gives quadratic part
\(X^2+2tU^2+(2\eta t+6y_0^2)Y^2\).
If the rank is one then \(t=y_0=0\) and this is the original
model at the origin. Otherwise the rank-two Hessian test applies.
The corresponding face \((x+uv)^2+y^4\) with independent \(u,v,y\)
follows as well: after translation and a linear shear it has
quadratic rank at least two unless \(t=y_0=0\); in that case
keep \(u=v\) by a linear flag and apply the result with \(\eta=0\).

\subsection{Every polar pencil at the quartic-eightfold cusp mark}
\label{sec:nov28-q8-cusp-mark}
Let \(X_4^8\subset\mathbf P^9\) be smooth. Every line \(L\subset X\)
whose first normal polars span a pencil satisfies
\(\delta_L(X,H)\ge198/31\). In particular every line and integral
conic now satisfies this bound, since the conic and higher-rank
line mark \(1601/250\) is larger. The proof uses the sharp signed
estimate \(17/11\) for the graph of a square plus a cube in
Section~\ref{sec:nov27-quintic-ten}, the
square-plus-binary-quartic estimate, and the completed-square
estimate for \((x+u^2+b y^2)^2+y^4\).

We use the ordinary linear flags through the line and the
weighted pencil flag in Chapter~\ref{ch:fullrow}, in the proof
under ``Ordinary polynomial series from the pole blocks.'' For
\(k\) retained lights the eliminated heavy coordinate has a
simple pole at each of three general roots of its polar. At
fixed pole power \(j\) the sum of the angular weight pieces is
the complete space of ordinary polynomials in the other heavy
coordinate and the \(k\) lights of degree at most
\(M-j+O(1)\), times the \(j\)-th graph power; the regular sector
is the same without the graph. Thus the four sectors have equal
asymptotic mass by summing \(\binom{M-j+O(1)+k+1}{k+1}\).
In units of \(k+3\), their moving degree is \(k+2\) and the
positive fixed exponent is one by the same simplex means.
The resulting ordinary affine signed inequality is exactly the
one needed in weighted adjunction; equality in the limiting
comparisons is tested through smaller rational coefficients.
For this mark the two-heavy pencil needs a sum of light weights
divided by the heavy weight of at least \(5/6\), and a signed
denominator at most \(17/11\) on every positive punctured pole
graph. The radial bound is then at most
\(1+(3/4)/(1+5/6)=31/22\); the regular and positive pole masses
\(1/4,3/4\) give the same angular bound
\(1/4+(3/4)(17/11)=31/22\). Since \(9/(31/22)=198/31\), it
suffices to establish these two estimates. At an ordinary graph
point, Hessian rank at least two gives signed denominator \(3/2\).
If the Hessian has rank one and the third Taylor term is nonzero
on its kernel, a general kernel direction and a quadratic direction
give a square plus a cube after subtracting the affine tangent and
using weights three and two, with graph weight six. Its radial
light ratio is \(5/6\) and its sharp signed bound is \(17/11\).
All constant choices in the initial flags are made over the base
line and specialize to general poles as in the pencil argument.
Subsequent ordinary graph calculations are made at each pole.

\paragraph{A terminal order-four test.}
Suppose that the preceding constant linear cuts give a terminal
quartic, and that the join of the line with its projective light
space has ambient codimension \(c\le3\). Its dimension is greater
than four, so it cannot lie in \(X\). The quartic is base
independent, and the finite Gauss map gives projective singular
codimension at least \(10-2c\) in the light space, as in the
preceding terminal calculation. Thus for \(c=3\) its singular
locus in \(\mathbf P^4\) is at most finite; for \(c=2\) it is
empty in \(\mathbf P^5\). Four general lights cut a smooth
quartic surface, and the outer radial ratio is one.

For any smooth quartic \(P\) in four lights, its graph has signed
denominator at most \(17/11\) at every punctured point. Its affine
Hessian cannot vanish at a nonzero vector by smoothness and Euler.
At a projective zero its gradient is nonzero, so its Hessian has
rank at least two: an isotropic radial vector of a rank-one symmetric
form would lie in its radical, contrary to Euler differentiation.
At a rank-one point off the zero surface, use the radial vector
as coordinate \(T\) and its gradient kernel for the other three
variables \(Y\). The form becomes
\(aT^4+T P_3(Y)+P_4(Y)\), \(a\ne0\). If \(P_3\ne0\), use the
ordinary square-plus-cube test just stated. Otherwise \(P_4\)
defines a smooth plane quartic, because a projective singular
point of \(P_4\) would give one of \(P\) at \(T=0\).
Retain two general \(Y\)'s where it is squarefree. Around the
chosen \(T\ne0\), after subtracting the affine tangent the weights
two on its variation and one on both other lights give a square
plus that binary quartic. The inner radial ratio is one and the
square-plus-binary-quartic estimate gives \(3/2\) on this face.
This proves the terminal test, also covering initial order four.

\paragraph{Quadratic order.}
For quadratic rank at least two, take a nondegenerate section on
two lights. In rank one make the quadratic root constant by at
most one preliminary linear cut, as in the quintic-tenfold proof.
Its primitive root has base degree at most one. Keep the remaining
quadratic direction \(x\) and the \(4\)- or \(5\)-dimensional
common kernel \(U\) of the coefficients of the primitive root.
If the restricted cubic on \(U\) is nonzero, one general direction
there and the quadratic direction give the pole face \(ax^2+bu^3\),
\(ab\ne0\), so the required radial and signed bounds both hold.

If the cubic on \(U\) is identically zero, the terminal quartic
\(C\) there is nonzero and reduced. Indeed the join of the line
with \(\mathbf P(U)\) has ambient codimension \(c\le4\), dimension
at least five, and the terminal Gauss bound gives projective
singular codimension at least \(10-2c\ge2\). Keep two general
directions \(Y=(u,v)\) in \(U\), so that \(C(Y)\) is squarefree.
The weights \((2,1,1)\) with heavy weight four give radial ratio
one and pole polynomial
\[
                  a x^2+x B(Y)+C(Y),\quad a\ne0.
\]
Here \(B\) is binary quadratic. At \(Y=0,x\ne0\) use the Hessian
rank-two test when \(B\ne0\), or the square-plus-binary-quartic
estimate otherwise. Put \(D=C-B^2/(4a)\). At a nonzero direction
\(Y=p\) with rank-one Hessian, its cubic Taylor term on the Hessian
kernel is the directional derivative \(\partial_pD\), up to a
nonzero constant and projection of that kernel to the \(Y\)-plane.
Indeed write the polynomial as \(a s^2+D\), \(s=x+B/(2a)\).
Use only the affine-linear part of \(s-s(x,p)\) as new centered
coordinate, so the remaining quadratic variation of \(s\) is
\(B(Y-p)/(2a)\). In these coordinates the quadratic part is a
nonzero square plus a quadratic in \(Y-p\); rank one makes the
second quadratic zero. The term of degree three in the expansion
of \(as^2\) vanishes on the kernel of that square. A binary
quartic has vanishing directional derivative in a nonzero
direction only if it is zero or a pure fourth power. Therefore
for nonpure \(D\) the ordinary square-plus-cube test suffices.
The case \(D=0\) contradicts squarefreeness of \(C\). For pure
nonzero \(D\), rescale to write the whole polynomial as
\((x+B'(Y))^2+L(Y)^4\). Squarefreeness of \(C=B'^2+L^4\)
forces \(B'|_{L=0}\ne0\). With \(y=L\), a binary linear change
completes a square in the other light while retaining \(y\),
giving \((x+u^2+b y^2)^2+y^4\). The completed-square estimate
gives signed denominator at most \(3/2\). This covers all quadratic
pencils.

\paragraph{Cubic order.}
Write \(G_t\) for the nonzero initial cubic on the full six-light
space over the geometric generic base. A general ternary section
of an irreducible \(G_t\) is an integral plane cubic by Bertini.
At every nonzero vector the graph of an integral ternary cubic
has signed denominator at most \(17/11\). Indeed it has no triple
point: a ternary cubic of multiplicity three at a projective point
depends on two variables and factors. At a double point, Hessian
rank one has a unique tangent line, and the cubic restriction to
that line is nonzero by integrality. Smooth zeros have Hessian
rank at least two by Euler. At a rank-one point off the cubic,
the radial and gradient-kernel coordinates express it as
\(aT^3+P_3(Y)\) with \(a\ne0\); the residual cubic must be nonzero
since the whole form is not a cube. The same proof works for a
reduced plane cubic with only ordinary nodes: at its singular
points the Hessian has rank two. The outer radial ratio for the
three ternary lights of ordinary weight one is one.

A constant hyperplane factor of the full cubic gives the terminal
test at \(c=3\), by cutting that hyperplane. In particular, if
the geometric generic cubic is nonreduced then its unique repeated
hyperplane descends to the base field, and is constant: twice the
primitive base degree of its linear root is at most the available
base degree one. Suppose the cubic is a line times an irreducible
quadric and the unique linear factor moves. Its primitive degree
is then exactly one and its quadric factor is constant up to scale.
The restriction of this quadric to the moving hyperplane has
generic rank at least two. Indeed restricting to a hyperplane
lowers rank by at most two. For rank three a rank-one restriction
requires that the hyperplane contain the radical and its projected
line be tangent to the nonsingular conic in the quotient plane.
The moving hyperplanes trace a projective line in the dual space,
which cannot lie in the dual conic. Thus three general constant
lights give a smooth conic and a line meeting it transversely,
and the preceding nodal plane-cubic test applies.

Finally let \(G_t\) have three distinct geometric linear factors.
If their span has dimension three, the general plane section is
a triangle, hence nodal. If it has dimension two, their span is
defined over the base field as the essential normal space, the
span of the second partial derivatives of \(G_t\). It extends to
a rank-two subbundle \(E\) of the constant dual light bundle by
the projective Grassmannian over the base curve. Write
\(E=\mathcal O(-i)\oplus\mathcal O(-j)\), \(0\le i\le j\).
The cubic belongs to \(H^0(\operatorname{Sym}^3E(1))\), since
membership holds generically and the quotient by this subbundle
is locally free. Hence \(3i\le1\), so \(i=0\). For \(j\ge1\)
the terms containing the second summand twice or more have negative
degree, which would force the first linear factor to divide the
cubic twice. Thus \(j=0\) and the two normal directions are
constant. Write them \(x,z\). On their four-dimensional common
kernel the terminal quartic is again nonzero and reduced by the
dimension and finite Gauss estimates at \(c=4\). Keep two general
kernel lights \(u,v\) where it is squarefree. The weights
\((4,4,3,3)\), with heavy weight twelve, give the face
\(G_3(x,z)+C_4(u,v)\), both binary forms squarefree.
At a point where both blocks are nonzero, each has nonzero
Hessian, so the total rank is at least two. On the quartic-block
axis, a rank-one Hessian has in its kernel the two cubic lights,
where the cubic term is nonzero. On the cubic-block axis a
rank-one Hessian at a nonzero binary vector is possible only
off its squarefree zero set, and its expression there is
\(aT^3+bY^3\), \(ab\ne0\), by the same radial argument. Thus
the ordinary square-plus-cube test applies on both axes. The
outer radial ratio is \(14/12\). This finishes the proof for
all orders of all polar pencils.

\section{Cubic contact on the quartic eightfold and the mark \texorpdfstring{\(51/8\)}{51/8}}
\label{sec:nov28-reverse-polar}
For a smooth quartic eightfold \(X\), suppose there is a number
\(M>51/8\) such that every line and every integral conic has
\(\delta_C(X,H)\ge M\). Then \(X\) is uniformly K-stable. In
particular \(M=198/31\) suffices, since \(198/31-51/8=3/248\).

We first prove a point estimate. Suppose in appropriate homogeneous
coordinates that the tangent hyperplane section at \(p\) has equation
\(T f(y)+g(y)=0\), with eight tangent coordinates and nonzero cubic
\(f\), quartic \(g\). If \(f\) is irreducible, then every point
valuation normalized by \(v(\mathfrak m_p)=1\), \(A_X(v)=8+q\), satisfies
\[
                  9S_H(v)\le 31/3+(17/8)q .                 \tag{2}
\]

\paragraph{A free line on a smooth quartic.}
Changing \(T\) by a linear form replaces \(g\) by \(g+Lf\).
A general such quartic \(G\subset\mathbf P^7\) is smooth.
Off \((f=g=0)\) this is Bertini applied to the linear system spanned
by \(g\) and the \(Lf\), on the open chart with \(g\)-coefficient one.
At a base point where \(df=0\), necessarily \(dg\ne0\), since
otherwise its whole \(T\)-line is singular on the tangent hyperplane
section. The latter has isolated singularities by the finite Gauss map.
At a base point with \(df\ne0\), a singularity of a member requires
\(dg\) proportional to \(df\), whereupon a unique \(T\) gives a
singular point of \(Tf+g=0\) above it. There are only finitely
many such base points, and the corresponding value at each is
avoided by a general \(L\).

The divisor \(D=(f=0)|_G\) is integral. The smooth quartic sixfold
has Picard group generated by \(H\), so any splitting of its
effective divisor of class \(3H\), including a repeated component,
expresses \(f|_G\) as a product of sections of smaller positive
degrees. These lift to ambient polynomials by projective normality.
The homogeneous ideal of \(G\) has zero degree-three part, so this
would factor \(f\). The same reasoning applies to any reduced
\(f\) by treating its irreducible factors separately.

Lines cover \(G\): through any point of \(G\) the remaining
line-direction equations on its projectivized tangent space
\(\mathbf P^5\) have degrees two, three and four. They have a
common zero by the projective dimension theorem. The projective
scheme of lines has finitely many irreducible components, so one
has a dominating universal family. The general line in this family
is free. Indeed generic smoothness of the evaluation on the smooth
locus of the reduced family makes the differential surjective
at a general point of a general line. In particular evaluations of
global infinitesimal deformations span the normal bundle there.
A negative summand in its splitting would have no global section,
so every summand is nonnegative. The Hilbert scheme is smooth there
by vanishing of \(H^1\) of this normal bundle, and evaluation is
smooth along the entire line by its global generation. Shrink to
the resulting open family of free lines. The inverse image of
the reduced divisor \(D\) by smooth evaluation is reduced.
Its map to the family is generically finite of degree three after
discarding lines contained in \(D\). Generic smoothness, and
removal of nondominating components, give a line \(\ell\subset G\)
with three distinct roots of \(f|_\ell\).

Since \(G\) is smooth along \(\ell\), its plane-normal cubic polars
form a basepoint-free system on this line. Take a normal direction
whose polar \(P\) avoids all three roots of \(f|_\ell\), and let its
span with \(\ell\) be a plane. Take binary coordinates on the line
and a transverse plane coordinate \(z\).

\paragraph{The angular degree and its fixed factors.}
Use the ordinary point presentation of Section~\ref{sec:q8-tangent-tests}
with the two symbols \(Z,W\) of section degrees \(1,3\) and orders
\(3,4\). Its relations are \(T^2Z+f,TW+g\). First specialize the
two tangent forms to their plane restrictions, then take least
\(z\)-orders. The resulting angular initials are the relatively
prime forms \(f_0=f|_\ell\) and \(g_0=zP\), with both \(f_0,P\)
binary cubics. At both stages the forms are a regular sequence.
Equivalently, in each section degree \(m\) and ordinary order \(j\)
the angular piece is
\(\sum_{2i+h\ge j-m}f^ig^h\mathbf C[y]_{j-3i-4h}\).
For a regular sequence, a homogeneous vector-space complement modulo
\((f,g)\) freely supplies the remainder in the unique expansion
in powers of \(f,g\). The resulting dimension of this piece depends
only on degrees. The specialized pieces contain those generated by
the two initials and their dimensions agree, proving equality.
Choosing compatible bases for each character block and the next
weight filtration gives the threshold comparison by specialization
as in Section~\ref{sec:q8-tangent-tests}.

Split the \(z\) exponent \(s\) and the exponents of the five
other independent linear variables outside the plane. At binary
degree \(D_0\) put \(l=j-m\). The resulting binary subspace is
\[
       \sum_{h\le s,\ 2i+h\ge l}
                f_0^i P^h\mathbf C[x,y]_{D_0-3i-3h},\qquad i,h\ge0 .
\]
There is a fixed power \(f_0^I\), where
\(I=\max(0,\lceil(l-s)/2\rceil,\lceil l-D_0/3\rceil)\).
After its removal the codimension in the complete binary space of
degree \(D_0-3I\) is bounded independently of all indices. Indeed
the binary ring is finite free over \(\mathbf C[f_0,P]\) with a
homogeneous remainder basis of degrees \(b\le4\). For a monomial
\(f_0^aP^c R_b\) of total degree \(D_0=3(a+c)+b\), membership
is precisely \(2a+\min(c,s)\ge l\), equivalently
\(a\ge(l-s)/2\) and \(a\ge l-(D_0-b)/3\). Only a bounded number
of endpoint monomials per residue can be missing. On the binary
projective line the mean of a codimension-\(O(1)\) subspace in
degree \(D_0-3I\) is at most \((D_0-3I)/2+O(1)\) at each point,
by the dimensions of the ordinary vanishing filtration. Since
\(f_0\) is squarefree its fixed power contributes at most \(I\).
Thus the whole block has mean at most \(D_0/2+O(1)\).

Scale the limiting means of exponents by nine and normalize by
\(m\). In the normal basis \(T^a Z^i W^h R\) the means of \(i,h\)
are \(1,1/3\), by their ordinary section degrees. The residual
ternary complete-intersection ring \(\mathbf C[x,y,z]/(f_0,zP)\)
has degree twelve. Its \(z\)-axis contributes degree nine: away from
the origin on that axis the transverse ring is the binary complete
intersection \((f_0,P)\) of length nine. The remaining degree three
is on \(z=0\), at the three directions of \(f_0\).
Thus the single free ternary simplex index in this ring contributes
\(3/4\) to the mean \(z\) exponent and \(1/4\) to binary degree.
For a direct verification, filter the ring by powers of \(z\).
In \(z\)-degree zero the ring is \(\mathbf C[x,y]/(f_0)\), with
constant Hilbert function three in large degrees. In each positive
\(z\)-degree the binary part is \(\mathbf C[x,y]/(f_0,P)\),
of length nine with bounded binary degrees. These are exactly the
two simplex contributions above. Including \(f_0^i(zP)^h\) gives
scaled mean \(1/3+3/4=13/12\) for \(z\), and
\(3+3/3+1/4=17/4\) for \(D_0\). Each other independent linear
coordinate has mean one. The radial ratio along the binary vertex
is \((17/4)/2=17/8\); on its exceptional binary line the same
upper bound follows from the marked-order calculation. After the
binary vertex blowup these strata and the independent coordinate
divisors have normal crossings. This bounds the angular series by
\(17/8\) in the scaled units. Ordinary point adjunction and the
known \(9S_H(\operatorname{ord}_p)=31/3\) give (2).

\paragraph{The resulting point comparison.}
Assume \(\beta=\delta(X,H)\le6\). Its original optimizing center
is a point by the line and conic bounds. At an optimizing point of
irreducible cubic contact, (2) and \(8+q\le6S_H(v)\) yield
\(q\ge8/3\) and \(S_H(v)\ge16/9\). At quadratic contact the
signed surface bound gives \(S_H(v)\ge2\). Reducible cubic forms
have a linear factor and their optimizing points are excluded by
the weighted extraction in Section~\ref{sec:q8-tangent-tests}.
The generalized Eckardt threshold exceeds six as well.

In the interpolation of Section~\ref{sec:nov28-direct-point} take
\(f=1/4,\ t=51/8-(10/9)\lambda\) for \(\lambda>0\) arbitrarily
small. Since \(1+g=10/9\), the additional point coefficient is
\(a=t+(1+g)\lambda-6=3/8\). At \(\lambda=0\) one has
\(\eta=1/3\) and \(a(16/9)+\eta=1\). Its derivative is
\(80/81-1+\alpha>0\), since the lower bound \(\alpha>0.063\)
comes from the explicit bottom hyperplane integral. Thus both
strict inequalities hold at positive \(\lambda\), and
\(t+h\lambda\to51/8<M\). This proves the asserted curve criterion.

\begin{proof}[Proof of Theorem~\ref{thm:nov28-quartic-eightfold}]
Section~\ref{sec:nov28-q8-cusp-mark} gives
\(\delta_C(X,H)\ge198/31\) for every line and integral conic.
Section~\ref{sec:nov28-reverse-polar} proves that
\(198/31>51/8\) forces \(\delta(X,H)>6\).
Thus \(\delta(X,-K_X)=\delta(X,H)/6>1\). The valuative criterion
and the Yau--Tian--Donaldson theorem give the assertions.
\end{proof}

\chapter{Quartic hypersurfaces of dimension nine and higher}\label{ch:nov30-replacements}
Let \(X^n\) be a smooth quartic hypersurface and \(L\) the limiting
hyperplane class on a K-polystable replacement. The two proofs in this
chapter use the same comparison of fillings
(Theorem~\ref{thm:nov30-all-ci-fillings}).
For \(n\ge10\), the geometry of the image of \(|L|\) excludes a distinct
replacement (Theorem~\ref{thm:dec1-quartic-high-dimensions},
page~\pageref{thm:dec1-quartic-high-dimensions}). For \(n=9\), the
hyperplane calculation leaves the cone over \(\operatorname{Gr}(2,5)\).
Section~\ref{sec:quartic-nine-higher-images} excludes both possible dimensions of
the image of \(|2L|\), completing the proof in
Theorem~\ref{thm:next-all-smooth-quartic-ninefolds}
(page~\pageref{thm:next-all-smooth-quartic-ninefolds}).
Chapter~\ref{ch:nov30-filtrations} proves the general results on filtered
hyperplane section rings used here.

\section{What a K-polystable replacement can look like}
\label{sec:nov30-special-fibers}
The argument with cubic limits also gives an obstruction in all
degrees. Let \(Z\) be a K-semistable \(n\)-dimensional
\(\mathbf Q\)-Fano variety with an ample integral
\(\mathbf Q\)-Cartier Weil class \(L\), \(-K_Z\sim rL\),
\(L^n=d\), and the section Hilbert function of a degree-\(d\)
hypersurface in \(\mathbf P^{n+1}\). Here \(3\le d\le n+1\)
and \(r=n+2-d\). Suppose the rational image of \(|L|\) has
dimension \(n\).

\begin{proposition}\label{prop:nov30-limit-degree-minus-one}
If this image is a hypersurface of degree \(d-1\), it has a
hyperplane section which is reducible or nonreduced.
In particular an image with factorial homogeneous coordinate ring
cannot occur in degree \(d-1\).
\end{proposition}
\begin{proof}
Write \(e=d-1\), \(S=\bigoplus_{k\ge0}H^0(Z,kL)\), and let
\(T\) be its degree-one generated subring. The full-dimensional
image has homogeneous coordinate ring
\(T=\mathbf C[x_0,\ldots,x_{n+1}]/(q_e)\).
Its Hilbert series is \((1-z^e)/(1-z)^{n+2}\), while that of
\(S\) is \((1-z^d)/(1-z)^{n+2}\).
Thus \(S_e/T_e\) has dimension one and
\(S_{e+1}/T_{e+1}\) has dimension \(n+1\).
Choose \(y\in S_e\setminus T_e\); multiplication by it from
\(T_1\) to the latter quotient has a nonzero kernel. Hence
\(\ell y=g\) for some nonzero linear \(\ell\in T_1\) and
\(g\in T_d\). We have \(g\notin\ell T\) by cancellation in \(S\).

If the section \((\ell=0)\) is integral, then \(T/\ell T\)
is a domain: the homogeneous restriction of \(q_e\) to that
hyperplane defines the section scheme-theoretically.
Consequently \(g\) is a nonzerodivisor modulo \(\ell\).
As in the degree-three calculation, the map
\(T[Y]/(\ell Y-g)\to S\), \(Y\mapsto y\), is injective:
it is injective after inverting \(\ell\), and reducing
\(\ell f=(\ell Y-g)h\) modulo \(\ell\) forces \(\ell\mid h\).
Giving \(Y\) degree \(e\), the source has exactly the Hilbert
series of \(S\), so the map is an isomorphism.

At the weight-\(e\) coordinate point of \(\operatorname{Proj}S\),
the cover of the weighted chart has equations
\[
              q_e(x)=0,\qquad \ell(x)-g_d(x)=0
              \quad\text{in }\mathbf A^{n+2},
\]
where \(g_d\) is a homogeneous lift of \(g\).
Its scalar \(\mu_e\)-action is free off the origin. The chart
has dimension \(n\), since it is finite over its quotient.
It is a complete intersection, normal by Serre's criterion:
off the origin it is \'etale over a normal variety, and
Cohen--Macaulayness gives the depth condition at the origin.
The finite quasi-\'etale cover of the klt germ is klt and has
degree \(e\). Eliminating the smooth second equation gives a
hypersurface germ of multiplicity exactly \(e\) in a smooth
\((n+1)\)-fold. Indeed its tangent equation is
\(q_e|_{\ell=0}\ne0\).

For such a germ \(V\), inversion of adjunction in its smooth
ambient space gives
\(\operatorname{lct}(V;\mathfrak m)\le n+1-e\): the ordinary
ambient blowup has log discrepancy \(n+1\), takes value \(e\)
on the hypersurface and value one on the maximal ideal.
The ideal formula for normalized volume, applied to the
maximal ideal of multiplicity \(e\), now gives
\(\widehat{\operatorname{vol}}(0,V)\le e(n+1-e)^n=er^n\).
On the other hand the finite-degree formula and the
local-to-global volume inequality on \(Z\) give
\[
 \widehat{\operatorname{vol}}(0,V)\ge
             ed\left(\frac{nr}{n+1}\right)^n>er^n ,
\]
since \((1+1/n)^n<\exp(1)<3\le d\).
The three volume results used here are
\cite[Theorems~2.7, 2.13, 2.14]{ar2007.14320}.
This rules out integrality of the indicated hyperplane section.
If \(T\) is factorial, every nonzero linear form in it is prime
by its degree, giving the last assertion.
\end{proof}

For \(n\ge3\) this forces a singular locus of dimension at least
\(n-3\) on the image. Indeed, in the indicated hyperplane
\(\mathbf P^n\) write the restricted form as \(AB\), with both
factors of positive degree. At \(A=B=0\) all tangential
derivatives vanish. The one normal derivative has positive degree,
so the projective dimension theorem gives a common zero locus of
dimension at least \(n-3\). If a factor is repeated the same
argument along that component gives dimension at least \(n-2\).
In particular a cubic image of a quartic limit must contain a
linear \(\mathbf P^{n-1}\), and isolated singularities are
insufficient in dimensions \(n\ge4\).

\begin{proposition}\label{prop:nov30-linear-image-component}
Suppose that $Z,L,d,r$ satisfy the hypotheses of
Proposition~\ref{prop:nov30-limit-degree-minus-one}, and that the
degree-$(d-1)$ image $W$ has a hyperplane section with a linear
component. Then
\[
                    r\left(2-\frac{2}{d}\right)\le n+1 .
\]
In particular, for $d=4$ and $n\ge9$ the hyperplane image cannot
be a cubic hypersurface, regardless of its singularities.
\end{proposition}
\begin{proof}
Put $e=d-1$. The linear component is a plane
$\Lambda=(a=b=0)\simeq\mathbf P^{n-1}$ for independent linear
forms $a,b$. The irreducible equation $q$ lies in $(a,b)$.
Use the full section ring $S$ and its degree-one subring
$T=\mathbf C[x_0,\ldots,x_{n+1}]/(q)$.
We have $S_i=T_i$ for $i<e$. The ring $S$ is integrally closed:
on writing $L$ as a Weil divisor on the normal $Z$, its homogeneous
pieces consist of rational sections with nonnegative order at
each prime divisor. Equivalently its embedding in $K(Z)[s]$ is
the intersection with the valuation rings given by the homogeneous
extensions of these orders.

First $q\notin(a,b)^2$. Otherwise write
$q=a^2U+abV+b^2W_0$, with $U,V,W_0$ homogeneous of degree $e-2$.
We have $a\nmid W_0$ by irreducibility. The homogeneous fraction
$h=bW_0/a$ satisfies $h^2+Vh+UW_0=0$, so it belongs to
$S_{e-2}=T_{e-2}$. Lifting it to a polynomial $h_0$ would give
$bW_0=ah_0$ as a polynomial identity, since both sides have degree
$e-1<e$. This contradicts $a\nmid W_0$.

Write $q=aD_0-bC_0$ with $C_0,D_0$ of degree $e-1$.
The homogeneous rational section $C_0/a=D_0/b$ of $(e-2)L$
does not belong to $S_{e-2}=T_{e-2}$, since the equality
$C_0=ah_0$ in degree $e-1$ would imply divisibility of the
polynomial $q$ by $a$. Thus it has a pole at a prime divisor $E$
of $Z$. Both $a$ and $b$ vanish along $E$ to order at least one,
as their products with this rational section are regular sections.

On a common smooth resolution of $Z\dashrightarrow W$ and the
pencil $[a:b]$, denote by $P$ the pullback of $L$, by $J$ the
pullback of $\mathcal O_W(1)$, and by $M$ the pullback of
$\mathcal O_{\mathbf P^1}(1)$. All three are nef. If $E'$ is the
strict transform of $E$, then $P-J$ and $P-M-E'$ are effective
up to $\mathbf Q$-linear equivalence. For example use all degree-$q_0$
products in the hyperplanes, respectively in the pencil, for a
positive integer $q_0$ making $q_0L$ Cartier. On the resolution
their moving classes are $q_0J$ and $q_0M$, and the pencil products
have $E'$ as a fixed component of multiplicity at least $q_0$.

For general $t$, the intersection of $W$ with $b=ta$ consists
of $\Lambda$ and a residual hypersurface of degree $e-1$ in that
hyperplane. Indeed $q|_{b=ta}$ is divisible by $a$ exactly once
for general $t$ since $q\notin(a,b)^2$. On the complement of
$\Lambda$ the residual hypersurface is the pencil fiber with its
scheme multiplicities; it has no component in $\Lambda$.
The image of the corresponding fiber cycle of the resolved pencil
therefore dominates that residual cycle with at least its
multiplicities, by the projection formula for the generically
finite map to $W$. Intersecting with $n-1$ general hyperplanes
gives $J^{n-1}M\ge e-1$. As $P-J$ is effective and $P,J,M$
are nef, expansion of $P^{n-1}-J^{n-1}$ gives
$P^{n-1}M\ge e-1$ as well.

For $0\le t\le1$, monotonicity of volume gives
\[
 \operatorname{vol}(P-tE')\ge ((1-t)P+tM)^n
  \ge (1-t)^n d + nt(1-t)^{n-1}(e-1).
\]
Integration yields $(n+1)S_L(E)\ge 1+(e-1)/d=2-2/d$.
The valuative inequality $rS_L(E)\le A_Z(E)=1$ proves the bound.
For a cubic image, Proposition~\ref{prop:nov30-limit-degree-minus-one}
provides a reducible or nonreduced hyperplane section. Every
nonintegral cubic hypersurface in a hyperplane has a linear
component. Substituting $r=n-2$ gives $3(n-2)/2\le n+1$, or
$n\le8$.
\end{proof}

\begin{proposition}\label{prop:nov30-quartic-birational-quadric}
Let $(Z,L)$ satisfy the hypotheses of
Proposition~\ref{prop:nov30-limit-degree-minus-one} with $d=4$ and
with $Z$ K-polystable. If $n\ge9$ and the image of $|L|$ is a quadric
of rank at least seven, the hyperplane map is a finite flat double
cover, and $L$ is the pullback of its hyperplane line bundle.
\end{proposition}
\begin{proof}
Write $S=R(Z,L)$ and $T=\mathbf C[x_0,\ldots,x_{n+1}]/(q)$ for
the degree-one image ring. The map degree is at most two, since the
moving hyperplane class on a resolution has top self-intersection
at most $L^n=4$. Suppose it is birational. Homogeneous elements of
$S$ then belong to $\operatorname{Frac}(T)$, by trivializing their
degree with a power of any nonzero linear section. The quadric ring
$T$ is factorial. Indeed, for a quadric of rank at least five use
coordinates $q=x_0x_1+q'(x_2,\dots)$; $x_0$ is prime as $q'$
has rank at least three, and $T[x_0^{-1}]$ is factorial. The same
criterion applies to $T/h$ for every linear form $h$, since its rank
is at least five.

The extra section in $S_2/T_2$ has a coprime homogeneous expression
$y=g/h$ in $\operatorname{Frac}(T)$, with $k=\deg h\ge1$ and
$\deg g=k+2$. Since $T$ is factorial,
$T[y]=T[Y]/(hY-g)$ with $\deg Y=2$. Its leading Hilbert multiplicity
is $k+2$, and inclusion into $S$ forces $k\le2$.

We record the test configuration used twice below. Suppose a normal
complete-intersection graded domain $B$ generated in positive degrees
and homogeneous $h,f\in B$ give the full ring
$S=B[v]/(hv-f)$, where $\deg h=k$, $\deg v=j$ and $h,f$ form a
regular sequence in $B$. With $\deg u=k$ take the relative Proj over
$\mathbf A^1_t$ of
\[
              B[t,u,v]/(h-tu,\ uv-f).
\]
On points the action has weights $1,-1,1$ on $t,u,v$ and fixes
$B$. The fiber at $1$ is $\operatorname{Proj}S$, and a divisible
Veronese polarizes this family by a multiple of $-K_Z$.
The total ring is a complete intersection: off $u=0$ it is
$B[u^{\pm1}]$ and on $u=0$ its equations require $h=f=0$ in
$B$, leaving $t,v$ free. Thus no component lies in $u=0$; the
ring is a domain, Cohen--Macaulay, and flat over $\mathbf C[t]$.
The locus $t=u=0$ has codimension two in the total ring. The
complement is normal, being covered by $u\ne0$ and the trivial
punctured family $t\ne0$. Serre's criterion gives normality.

On sections of the central fiber $u,v$ have weights $1,-1$.
Differentiating its character series at the trivial character gives
the total-weight series
\[
 H_S(z)\left(\frac{z^k}{1-z^k}-\frac{z^j}{1-z^j}\right),
 \qquad H_S(z)=\frac{1-z^4}{(1-z)^{n+2}}.
\]
Indeed the central presentation has coefficient ring $B/h$ of
weight zero and the relation $uv-f$ of weight zero. Degree pieces
compute section weights in a sufficiently divisible Veronese; the
central ring is Cohen--Macaulay of dimension $n+1$. Put
$\Delta=1/k-1/j$. At $z=e^{-s}$ the factor in parentheses is
$\Delta/s+O(s)$, and
$H_S(e^{-s})=4s^{-n-1}(1+rs/2+O(s^2))$. Other roots of unity
give at most simple poles in the weight series. Thus the Hilbert
coefficients are $a_0=4/n!$, $a_1=nr a_0/2$, and the two weight
coefficients are $b_0=4\Delta/(n+1)!$ and
$b_1=4r\Delta/(2n!)$. The Donaldson--Futaki invariant in the
introductory convention is $-r\Delta/(n+1)$. In particular $k<j$
contradicts K-semistability.

If $k=2$ for $y=g/h$, the Hilbert series of $T[y]$ equals that of
$S$, so $S=T[y]$. Apply this construction with $B=T$, $k=j=2$;
the Futaki invariant is zero. This is not a product configuration.
If the special fiber is not integral and normal it cannot be a
product. Otherwise suppose it is a product. The Weil root $L$
extends uniquely up to linear equivalence on the product once its
punctured-family extension is given, since the integral central
fiber is principal. In a product trivialization the one-parameter
action preserves the class of $L$: for instance the difference is
torsion killed by $r$ since $rL=-K_Z$, and the connected group acts
trivially on this finite ambiguity. Equivalently, the class group
of $Z\times\mathbf G_m$ is pulled back from $Z$ by polynomial
homotopy invariance for the class group of a normal variety and
localization. Hence the root on the product restricts to the same
class on its two fibers. The tautological rank-one sheaf in this
family gives precisely this extension in codimension one; on each
fiber it is invertible where some $x_i$ is nonzero, and the missing
locus consists of at most two points (on the special fiber,
$uv=0$ when all $x_i=0$). On a normal integral special fiber
its complete space of degree-one sections is spanned by the $x_i$,
by the Cohen--Macaulay presentation. Its image is $(q=h=0)$ of
dimension $n-1$, whereas the general fiber has image $(q=0)$ of
dimension $n$. This excludes a product and contradicts K-polystability.

There remains $k=1$. First suppose $g\bmod h$ has a linear factor
$b$ in $T/h$. Write $g=hp+bC$, where $b,h$ are independent linear
forms and $p,C$ are quadratic. The rational section
$C/h=(y-p)/b$ of $L$ has a pole on some prime divisor $E$ of $Z$:
if it belonged to $S_1=T_1$, we would have $y\in T_2$.
Both $h,b$ vanish along $E$. Apply the nef-pencil computation of
Proposition~\ref{prop:nov30-linear-image-component} to $[h:b]$.
Here the quadric does not contain $(h=b=0)$, since restriction
lowers rank by at most four. Thus each general pencil hyperplane
has residual degree two, with no fixed component on the image.
The resulting inequality $(n+1)S_L(E)\ge1+2/4$ contradicts
$(n-2)S_L(E)\le1$ when $n\ge9$.

Otherwise $g\bmod h$ is prime in the factorial ring $T/h$, as
it has degree three and has no linear factor. The ring $B=T[y]$
is normal: it is a hypersurface over $T$, is normal away from
$h=0$, and $B/h=(T/(h,g))[y]$ is a domain. Its remaining
codimension-one localization has maximal ideal generated by $h$.
Its Hilbert series equals $(1-z^3)/(1-z)^{n+2}$. Thus
$\dim S_3/B_3=1$ and $\dim S_4/B_4=n+1$; choose
$w\in S_3\setminus B_3$. Multiplication from $B_1$ to $S_4/B_4$
gives a relation $\ell w=f$ for nonzero linear $\ell$ and
$f\in B_4\setminus\ell B$. We claim $B/\ell$ is a domain.
This was shown when $\ell$ is proportional to $h$. Otherwise
$T/\ell$ and $T/(h,\ell)$ are domains by the rank bounds.
Moreover $g\notin(h,\ell)T$, since otherwise $g\bmod h$ would
have the linear factor $\ell$. Thus $h,g$ form a regular sequence
in $T/\ell$, and $(T/\ell)[Y]/(hY-g)$ embeds in its localization
at $h$ by cancellation. This proves the claim. Consequently
$\ell,f$ form a regular sequence in $B$, so likewise
$B[w]=B[V]/(\ell V-f)$. Its ordinary Hilbert series equals
$H_S$, whence $S=B[w]$. The preceding test configuration with
$(k,j)=(1,3)$ has negative Futaki invariant, the final contradiction
in the birational case.

The hyperplane map therefore has degree two. In this case the
extra section $y\in S_2\setminus T_2$ cannot belong to
$\operatorname{Frac}(T)$. Indeed the same denominator computation
would give $B=T[y]$ of multiplicity at least three. Choosing a
homogeneous $s\in S\setminus\operatorname{Frac}(T)$ then yields
$B\oplus sB\subset S$ of multiplicity at least six, exceeding four.
Thus $T\oplus yT\subset S$ is a direct sum, with the full Hilbert
series $(1+z^2)H_T=H_S$. Equality in every degree provides a monic
relation $y^2+b_2(x)y+b_4(x)=0$. Completing the square identifies
$S$ with $T[y]/(y^2-f_4(x))$ by the Hilbert series, proving the
finite flat double-cover description. The presentation avoids the
weight-two vertex, so its hyperplane class $L$ is Cartier.
\end{proof}

For comparison, a quadric of rank three cannot be the full-dimensional
image for $n\ge9$. Its hyperplane tangent along a general ruling cuts
that ruling twice, producing a section vanishing at least twice at a
prime divisor on $Z$ above the general ruling. This violates the
K-semistable bound $\alpha(Z)\ge1/(n+1)$ since $2(n-2)>n+1$.
A rank-four quadric cannot be a degree-two image for $n\ge9$ either.
It is the cone over $\mathbf P^1\times\mathbf P^1$. With $U,V$ the
two pulled-back ruling pencils and $J$ the resolved hyperplane class,
$P-U-V$ is effective up to $\mathbf Q$-linear equivalence, as in
Section~\ref{cubics:rulings}. Here $P^{n-1}V\ge J^{n-1}V\ge2$
because the map has degree two and a general ruling on the quadric is
a linear divisor. A prime component $E$ of a general first ruling
fiber on $Z$ then gives $(n+1)S_L(E)\ge1+2/4$ by the same volume
integral. This contradicts K-semistability.
In the birational rank-four case the same argument, with the
tangent-sheaf slope bound $P^{n-1}V\ge4(n-2)/(2n)$ from the
proof of Lemma~\ref{cubics:slope} (now $-K_Z=(n-2)L$), yields
$(n-2)^2\le6n$. Thus a rank-four image is also excluded
for $n\ge10$ in the birational case.

\subsection{The double-quadric boundary and its quartic GIT point}
\label{sec:nov30-double-quadric-limits}
\begin{proposition}\label{prop:nov30-double-quadric-limits}
For \(n\ge5\), let \(Q=(q=0)\subset\mathbf P^{n+1}\) be a smooth
quadric and let \(Z\to Q\) be the double cover branched along a
general quartic section. This smooth uniformly K-stable Fano
\(n\)-fold is a polarized limit of smooth K-stable quartic
hypersurfaces. For every such smoothing of \((Z,H)\), where
\(H\) is the pullback of the hyperplane class, the quartic GIT
quotient limit is the double quadric \([q^2]\).
In particular this boundary type cannot occur over a smooth
quartic GIT point.
\end{proposition}
\begin{proof}
The branch \(B\) is a general complete intersection of type
\((2,4)\), Fano of index \(n-4\) and K-stable by
\cite[Theorem~1.3]{ar2608.09293}. This applies with \(Q\)
fixed since smooth quadrics form one projective orbit. The
quadric is K-polystable, \(-K_Q=nH_Q\), and \(B\sim4H_Q\).
Zhang--Zhou's log interpolation theorem
\cite[Theorem~1.7]{ar2003.06839} gives K-polystability of
\((Q,\frac12 B)\) because
\(\frac12<1-(n/4-1)/n=\frac34+1/n\).
The double cover is crepant for this pair, hence is K-polystable
by \cite[Theorem~1.2(2)]{ar2001.10557}. Its connected
automorphisms preserve \(H\), since \(H^1(Z,\mathcal O_Z)=0\),
and \(|H|\) is the double cover followed by the quadric embedding.
Thus they preserve \((Q,B)\). The latter group is finite: its
action on \(B\) is faithful, since an irreducible spanning
variety cannot be pointwise fixed by a nonscalar projectivity,
and \(B\) is K-stable. Automorphisms over \(Q\) form the
finite deck group. It follows that \(Z\) is uniformly K-stable.

In \(\mathbf P(1^{n+2},2)\), with weight-two coordinate \(y\),
consider \(q(x)=ty,\ y^2=f_4(x)\), where \(f_4|_Q\) defines
\(B\). These complete intersections give a flat family around
zero. The central fiber is \((Z,H)\); it avoids the weighted
vertex, and is smooth. The nearby fibers are smooth by
openness of smoothness and properness, and uniformly K-stable
by openness in the smooth Fano family. Eliminating \(y\)
gives their quartic equations \(q^2-t^2 f_4=0\).

For an arbitrary polarized smoothing over a discrete valuation
ring with uniformizer \(t\), let \(A_j\) be the relative spaces
of sections of \(jH\). They are free and commute with base
change by Kodaira vanishing and cohomology and base change
for the smooth Fano fibers. Choose a basis \(x\) of \(A_1\).
The special fiber section algebra is
\[
        \mathbf C[x,y]/(q(x),y^2-f_4(x)),
        \qquad \deg x_i=1,\quad\deg y=2 .
\]
Indeed the pushforward from the double cover splits as
\(\mathcal O_Q\oplus\mathcal O_Q(-2)\). Hence
\(\operatorname{Sym}^2 A_1\to A_2\) is generically an
isomorphism between free modules of the same rank, with
one-dimensional kernel and cokernel on the special fiber.
Smith normal form yields a primitive quadratic \(q_R(x)\)
reducing to \(q\) up to scalar, a section \(y_R\in A_2\)
generating the central cokernel, and an integer \(e\ge1\)
such that \(q_R(x)=t^e y_R\).

The central algebra is also generated by the reductions of
\(x,y_R\), since \(\bar y_R\) differs from a nonzero scalar
multiple of \(y\) by a quadric in \(x\). In degree four,
the central kernel therefore contains a relation monic
quadratic in \(\bar y_R\). The degree-four map from the
weighted polynomial ring onto \(A_4\) is surjective by
Nakayama, and its kernel reduces onto the central kernel
by freeness of \(A_4\). Rescaling the leading unit, lift
this relation to \(y_R^2+b_2(x)y_R+b_4(x)=0\). It follows
that the primitive degree-four equation of the generic
embedded hypersurface is
\[
               q_R(x)^2+t^e b_2(x)q_R(x)+t^{2e}b_4(x)=0 .
\]
Its reduction is \(q^2\). The square of a nondegenerate
quadric is GIT-polystable: for example its affine
\(\mathrm{SL}_{n+2}\)-orbit is closed since the affine orbit of
the nondegenerate quadric with fixed determinant is closed
and the squaring map on the space of quadratic forms is
finite onto its image. The GIT quotient limit is thus
\([q^2]\), independently of the family or basis.
Smooth quartic equations are GIT-stable and their quotient
points are distinct from this closed orbit.
\end{proof}
The properness argument for cubics that recovers an embedded
hypersurface on every K-polystable limit therefore requires a
different conclusion in degree four. This proposition excludes
one smooth special-fiber type over smooth GIT points; other
section algebras and non-Cartier limits still require analysis.

\begin{corollary}\label{cor:nov30-smooth-quartic-replacement}
Let \(n\ge4\). A smooth K-polystable replacement of a family of
smooth quartic \(n\)-folds whose quartic GIT limit is a smooth
quartic is isomorphic to that quartic. Thus any smooth quartic
that fails K-stability would have to admit a \emph{singular}
K-polystable replacement when approached by K-stable quartics.
\end{corollary}
\begin{proof}
After finite base change extend the hyperplane class from the
punctured family by closing a divisor. The replacement is
smooth near the special fiber by flatness and smoothness of
the fibers, so the extension is Cartier. Denote its restriction
by \(H\). Then \(-K=(n-2)H\) and \(H^n=4\). Smooth proper
deformation invariance of integral cohomology and the exponential
sequence, using Fano vanishing, give \(\operatorname{Pic}=\mathbf Z H\)
as on the quartic fibers. Thus this is a smooth prime Fano variety
of coindex three and genus three. By Mukai's classification
(S.~Mukai, \emph{Biregular classification of Fano 3-folds and
Fano manifolds of coindex 3}, Proc. Natl. Acad. Sci. USA
\textbf{86} (1989), 3000--3002; see also E.~Jovinelly and
F.~Okamura, \emph{Rational Curves on Coindex 3 Fano Varieties},
\href{https://arxiv.org/abs/2409.00834}{arXiv:2409.00834},
Theorem~2.6), its fundamental model is either an embedded
quartic or a complete intersection of degrees two and four
in \(\mathbf P(1^{n+2},2)\).

In the latter presentation a nonzero \(y\)-coefficient in the
quadratic equation already eliminates \(y\), giving the first
case. Otherwise the equation is \(q(x)=0\). The coefficient
of \(y^2\) in the other equation must be nonzero. If it
vanished the weighted vertex would belong to the intersection:
on its index-two chart both equations have no linear term,
so the covering complete intersection is singular at the origin.
It is normal, since it is smooth away from the origin by the
quasi-\'etale map to the assumed smooth special fiber, and is
Cohen--Macaulay. Purity of the branch locus would make this
finite map \'etale over the smooth fiber at the origin as well,
a contradiction. Completing the square now makes the second
equation \(y^2=f_4(x)\). At any point over a singular point of
\((q=0)\) the Jacobian of these two equations would have rank
at most one, so the quadric is smooth. The arbitrary-smoothing
calculation above then gives the GIT limit \([q^2]\), even
when the smooth branch is special. This contradicts the smooth
GIT limit. For the embedded quartic central fiber the line
bundle and its complete relative system give a quartic family;
separatedness of the GIT quotient identifies its stable
smooth limit with the prescribed one. Smooth quartics have
finite automorphism group, so their K-polystability implies
uniform K-stability.
\end{proof}

\paragraph{A singular double cover at a quadric vertex.}
The conclusion that the GIT limit cannot be smooth also holds
for a polarized double cover \(y^2=f_4(x)\) of an integral
quadric of corank one if its branch avoids the quadric vertex.
Here the naive central equation \(q^2\) is GIT-unstable.
For a polarized smoothing to quartics assume that the relative
degree-one, degree-two and degree-four section spaces are free and commute
with base change, as holds in a klt Fano family with relatively
ample Cartier hyperplane root by vanishing. The central section
ring has the same two-generator presentation by the finite
double-cover construction. Apply the Smith normal form argument
above, choosing on the central fiber the trace-zero generator
and lifting a relation monic in the weight-two generator.
Completing the square in the lifted monic relation yields the
generic quartic equation
\[
                     q_R(x)^2+t^{2e}g_R(x)=0 ,
         \qquad  q_R\bmod t=q,\qquad g_R(\operatorname{vertex}(q))
                                    \not\equiv0\pmod t .
\]
Signs or invertible constants in \(g_R\) do not affect the last
condition. Over the completed DVR, splitting off the unimodular
part of the quadratic form and taking roots of units gives
\(q_R=q_A(x')+t^b c z^2\), with \(q_A\) nondegenerate on
\(n+1\) variables, \(b\ge1\), and \(c\) a unit; allow the last
term to vanish, with \(b=\infty\). After a finite base change
rescale \(z\) by \(t^{-\min(b,e)/2}\). Each coefficient remains
integral. If \(b<e\) the central quartic is the square of a
nondegenerate quadric. If \(b=e\) or \(b>e\), only the
\(z^4\) term of the rescaled \(t^{2e}g_R\) survives; its
coefficient is nonzero. In all cases the limit is
\[
              (q_A+\lambda z^2)(q_A+\mu z^2),
              \qquad (\lambda,\mu)\ne(0,0).
\]

This reducible quartic is GIT-semistable. If both factors are
nondegenerate, the product of their determinants gives the
equal-bidegree invariant test. If just one is nondegenerate,
write the factors as \(Q_1=q_A\) and \(Q_2=q_A+\mu z^2\).
The invariant pencil determinant
\(\det(uQ_1+vQ_2)=c'v(u+v)^{n+1}\) has nonzero coefficients
of \(u v^{n+1}\) and \(u^{n+1}v\). Their product is an
\(\mathrm{SL}_{n+2}\)-invariant of equal positive bidegree in
\(Q_1,Q_2\), nonzero at this pair. In either case, if a
one-parameter subgroup had strictly positive minimum weight
on \(Q_1Q_2\), the sum of its minimum weights on the two
factors would be strictly positive, forcing every such
equal-bidegree invariant to vanish. Hilbert--Mumford excludes
instability. The resulting semistable quotient point cannot be
that of a smooth stable quartic, since the latter has a single
orbit as its entire semistable fiber. This excludes the stated
singular double covers without requiring GIT semistability of
their original hyperplane cycle.

\begin{proposition}\label{prop:nov30-factorial-quartic-replacement}
Let \(n\ge4\), and let \(Z\) be a terminal factorial
\(\mathbf Q\)-Gorenstein Fano limit of smooth quartic
\(n\)-folds. Its hyperplane root is an ample, globally generated
Cartier divisor \(L\), and \(\operatorname{Pic}(Z)=\mathbf Z L\).
Either \((Z,L)\) is an embedded quartic, or \(|L|\) makes
\(Z\) a finite flat double cover of an integral quadric
\(Q\subset\mathbf P^{n+1}\), with section algebra
\[
        \mathbf C[x_0,\ldots,x_{n+1},y]/
                         (q(x),y^2-f_4(x)),\qquad \deg y=2 .
\]
If \(Z\) is K-semistable and the quartic GIT limit is smooth,
the embedded case gives exactly that smooth quartic. In the
double-cover case \(Q\) must be singular, and if it has corank
one, its branch contains the vertex.
\end{proposition}
\begin{proof}
Work over a complete discrete valuation ring after a finite base
change trivializing the generic hyperplane class. As in the
proof of \cite[Theorem~3.1(1)--(3)]{ar2007.14320}, close a
generic hyperplane to a Weil divisor \(\mathcal L\) on the
total space. The relation \(-K\sim(n-2)\mathcal L\) holds
relatively, since the special fiber is reduced, integral and
principal. Hence \(\mathcal L\) is \(\mathbf Q\)-Cartier.
The total space is klt near its klt special fiber by inversion
of adjunction, so its rank-one sheaves
\(\mathcal O(j\mathcal L)\) are Cohen--Macaulay.
Their restrictions to the special fiber are Cohen--Macaulay
and rank one; in codimension one they are the usual
restrictions of line bundles, since the family is smooth
at smooth points of the central fiber. Thus they are the
reflexive sheaves \(\mathcal O_Z(jL)\) of a Weil divisor
\(L\). Factoriality makes \(L\) Cartier. A local generator
of \(\mathcal O_Z(L)\) lifts by Nakayama to a generator of
\(\mathcal O(\mathcal L)\), which is then free of rank one
by torsion-freeness. Properness allows shrinking around
the central fiber. In particular \(L\) extends the generic
hyperplane line bundle, is ample since \(-K_Z=(n-2)L\),
and has \(L^n=4\).

Every line bundle on \(Z\) extends uniquely up to isomorphism
to the formal family: across each infinitesimal thickening the
obstructions and ambiguities lie respectively in
\(H^2(Z,\mathcal O_Z)\) and \(H^1(Z,\mathcal O_Z)\), which
vanish for a klt Fano variety. Grothendieck existence for the
proper family algebraizes these compatible invertible
sheaves, and their algebraizations are invertible by the local
criterion for flatness and properness. Restriction from the
normal total space to the generic fiber is injective on
Picard groups. Indeed a generically trivial line bundle has
a rational section with vertical Cartier divisor, whose Weil
divisor is an integral multiple of the principal special fiber.
Restriction to the geometric generic fiber is also injective
by Hilbert's theorem 90. The geometric generic quartic has
Picard group generated by its hyperplane class by Lefschetz.
Since \(L\) extends this generator, \(\operatorname{Pic}(Z)=
\mathbf Z L\).

Global generation in every dimension used here is precisely
A.~Bayer, A.~Kuznetsov and E.~Macrì, \emph{Mukai models
of Fano varieties}, \href{https://arxiv.org/abs/2501.16157}
{arXiv:2501.16157}, Lemma~6.1(a). The lemma assumes dimension
at least three, factorial terminal singularities,
\(\operatorname{Pic}=\mathbf Z L\), \(-K=(n-2)L\), and
\(L^n=2g-2\), and gives basepoint freeness for \(g\ge3\).
Here \(g=3\). The proof of that lemma reduces the
higher-dimensional assertion to threefolds by its part (b).

Vanishing for all nonnegative multiples of \(L\) and base
change give the quartic Hilbert series. Since \(L\) is
ample and globally generated, its complete linear system
gives a finite map to a nondegenerate hypersurface
\(W\subset\mathbf P^{n+1}\). The degree of this map times
\(\deg W\) is four and \(\deg W\ge2\). If \(\deg W=4\),
the subalgebra generated by degree one already has the
Hilbert series of the whole section ring, so the two agree.
Otherwise let \(T=\mathbf C[x]/(q)\) be the coordinate ring
of the integral quadric image. This ring is normal by
Serre's criterion (an integral quadric has rank at least
three). Its finite extension \(S\), the section ring of
\(Z\), has Hilbert series \((1+t^2)\operatorname{Hilb}_T(t)\).
Finiteness follows also by pushing the structure sheaf forward
under the finite hyperplane map. Choose \(y\in S_2\setminus
T_2\); integrality and normality of \(T\) give \(y\notin
\operatorname{Frac}(T)\). Thus \(T\oplus Ty\subset S\), and
the Hilbert-series identity makes this equality in every
degree. Expressing \(y^2\) in this free module and completing
the square gives the asserted equation. Sheafifying
\(T\oplus T(-2)\) gives finite flatness over \(Q\).

We give the GIT comparison also for a singular embedded
K-semistable quartic. In general let a one-parameter subgroup
of \(\mathrm{SL}_{n+2}\) have minimum monomial weight \(w\)
on an embedded degree-\(d\) equation, using weights on
coordinate sections with sum zero. The corresponding orbit
test configuration for the hyperplane polarization has central
coordinate ring the quotient by the initial equation of
weight \(w\). Its Hilbert function and degree-\(k\) weight
for \(k\ge d\) are
\[
 \binom{k+n+1}{n+1}-\binom{k-d+n+1}{n+1},
 \qquad -w\binom{k-d+n+1}{n+1}.
\]
If their first two coefficients are \(a_0,a_1\) and
\(b_0,b_1\), expansion gives
\[
       b_0a_1/a_0-b_1
              =-\frac{(d-1)(n+2)}{2(n+1)n!}\,w .
\]
Normalization of an ample test configuration cannot increase
its Donaldson--Futaki invariant (S.~Boucksom, T.~Hisamoto and
M.~Jonsson, \emph{Uniform K-stability, Duistermaat--Heckman
measures and singularities of pairs}, Proposition~3.15 and
Corollary~3.16,
\href{https://arxiv.org/abs/1504.06568}{arXiv:1504.06568}).
K-semistability therefore excludes \(w>0\). Hilbert--Mumford
gives GIT semistability of the equation. Since the complete
relative hyperplane system embeds the family in the embedded
case, separatedness of the GIT quotient and the single-orbit
fiber at a smooth stable equation identify \(Z\) with that
quartic. In the double-cover case the two preceding GIT
calculations exclude a smooth quadric or an unbranched
corank-one vertex. This proves the assertion.
\end{proof}

\begin{theorem}\label{thm:nov30-quartic-cartier-k3}
Let $Z$ be a complex K-semistable $n$-dimensional $\mathbf Q$-Fano
variety, $n\ge7$, with an ample Cartier divisor $L$ such that
$-K_Z\sim(n-2)L$ and the section Hilbert function of $(Z,L)$ is
that of a quartic hypersurface. Then $L$ is globally generated.
Its complete linear system identifies $Z$ with an embedded quartic
or with a finite flat double cover of an integral quadric branched
along a quartic section.
\end{theorem}
\begin{proof}
We first describe what a base point would imply without using
K-semistability. Since $K_Z$ is Cartier and $Z$ is klt, $Z$ is
Gorenstein canonical. Take a ladder of general members of $|L|$
down to a surface. More precisely, on each intermediate variety
$Y$ of dimension $j\ge3$, adjunction gives $-K_Y\sim(j-2)L|_Y$.
The Ambro--Kawamata general-member theorem in
\cite[Theorem~4.1]{ar2007.14320} applies: $L|_Y$ is nef Cartier,
$L|_Y-K_Y$ is ample, and the coefficient $j-2$ is strictly greater
than $j-3\ge0$. Hence a general member gives a plt pair. The member
is normal and klt by adjunction, is Gorenstein, and is connected
by ampleness and vanishing; in particular it is integral. The map
on hyperplane sections at each step is surjective by
$H^1(Y,\mathcal O_Y)=0$, and its kernel is one-dimensional.

The resulting surface $T$ has Du Val singularities, $K_T\sim0$,
$L_T^2=4$ and $h^0(T,L_T)=4$. Also $H^1(T,\mathcal O_T)=0$:
in the last threefold $Y$ the sequence cutting $T$ has kernel
$\mathcal O_Y(K_Y)$, whose second cohomology is dual to
$H^1(Y,\mathcal O_Y)=0$. Its minimal resolution
$\rho:\widetilde T\to T$ is a K3 surface, since Du Val singularities
are rational and crepant. The complete system of
$A=\rho^*L_T$ is the pullback of the complete system on $T$.
Here $A$ is nef and big.

We recall the precise K3 argument. A complete linear system on a
smooth K3 surface has no base points outside its fixed part
(B.~Saint-Donat, \emph{Projective models of $K$-3 surfaces},
Amer. J. Math. \textbf{96} (1974), 602--639, Corollary~3.2).
Write $A=M+F$ for its moving and fixed parts. If its image has
dimension two, $M$ is big and nef. Riemann--Roch and vanishing give
$2+A^2/2=h^0(A)=h^0(M)=2+M^2/2$. Thus
$(A+M)F=0$, and nefness gives $AF=MF=0$, hence $F^2=0$.
The Hodge index theorem applied to $A$ excludes a nonzero effective
$F$ with these two intersections. In this case $|A|$ is base point
free. If its image has dimension one, the base-point-free moving
part gives a morphism. Its Stein factor has smooth rational base,
since $H^1(\widetilde T,\mathcal O)=0$. Its moving line bundle is
pulled back from $\mathcal O_{\mathbf P^1}(k)$, with $k+1=h^0(M)=4$.
Thus its image in $\mathbf P^3$ is a twisted cubic. The image cannot
be a point, since the four sections are linearly independent.

Write $W$ for the rational image of $|L|$ and $s=\dim W$. The
successive section surjections show that our $n-2$ general divisors
are cuts by general hyperplanes of the target and that $T$ meets
the domain of the rational map. Its image is dense in the general
linear section of $W$ of codimension $n-2$: on a dense open of
$W$ the rational map is a surjective morphism from an open of $Z$,
and general linear sections of positive dimension meet that open
dense subset densely. If $s<n-2$ the general linear space misses
$W$; if $s=n-2$ it cuts a finite set. These cases contradict the
positive-dimensional image of $T$. Hence $s=n-1$ or $s=n$.
When $s=n$, $T$ has two-dimensional image, and global generation
on $\widetilde T$ implies global generation on $T$. Every original
base point would belong to $T$ and remain a base point there by
surjectivity of restriction. Consequently $|L|$ is globally generated.

If $s=n-1$, a general curve section of $W$ is the twisted cubic
just described; it is integral and generically reduced by Bertini,
so $\deg W=3$. A codimension-two variety of minimal degree is a
rational normal scroll, possibly a cone \cite{EH}; the Veronese
surface and its cones have codimension three. If a scroll summand
has degree at least two, a hyperplane divisor contains a general
ruling at least twice: choose a binary polynomial with a double
zero in that summand, and zero polynomials in the others. A general
ruling meets the open where the rational map from $Z$ is dominant
and smooth over a smooth open of $W$. Its inverse image has a prime
component, so some section of $L$ vanishes there at least twice.
This gives $\alpha(Z)\le1/(2(n-2))<1/(n+1)$, contrary to
K-semistability \cite[Theorem~2.5]{ar2007.14320}.

The other scroll is a cone over $\mathbf P^1\times\mathbf P^2$.
On a common resolution of $Z$, the map to $W$, and its two Segre
projections, let $P$ be the pullback of $L$, $J$ the pulled-back
hyperplane class of $W$, and $U,V$ the two pulled-back hyperplane
classes of $\mathbf P^1,\mathbf P^2$. Their nef classes satisfy
$P-J\ge0$ and $P-U-V\ge0$ up to linear equivalence: the first
difference is the fixed part of the resolved system, and the
Segre matrix-coordinate subsystem has moving class $U+V$.
Let $c$ be the $P$-degree of the general fiber cycle over a
resolution of $W$. This is a positive integer, since $L$ is ample
Cartier and the rational map has relative dimension one. Projection
formula and nefness yield
\[
 P^{n-1}V\ge P J^{n-2}V=2c,\qquad
 P^{n-2}V^2\ge P J^{n-3}V^2=c.
\]
Indeed pulling back one or two general hyperplanes of $\mathbf P^2$
on the scroll gives cones over $\mathbf P^1\times\mathbf P^1$ or
$\mathbf P^1\times\{\mathrm{pt}\}$, of degree two or one. A prime
component $E$ of a general first-ruling fiber on $Z$ has strict
transform $E'$ with $P-E'-V$ effective up to linear equivalence.
For $0\le t\le1$, integrate
$\operatorname{vol}(P-tE')\ge((1-t)P+tV)^n$. Each mixed
intersection in the binomial expansion has integrated coefficient
$1/(n+1)$. Thus
\[
 (n+1)S_L(E)\ge1+\frac{2c+c}{4}\ge\frac74.
\]
As $(7/4)(n-2)>n+1$ for $n\ge7$, this contradicts the valuative
K-semistability inequality. Hence only $s=n$ occurs.

Finally global generation and ampleness give a finite morphism,
with image degree times mapping degree equal to four. The image is
a nondegenerate integral hypersurface, so the possibilities are
an embedded quartic or a map of degree two onto a quadric. In the
first case the section ring agrees with the quartic coordinate
ring degree by degree by the Hilbert series. In the second, let
$B$ denote the quadric ring and $S=R(Z,L)$. The finite morphism
makes $S$ integral over the normal ring $B$ (every integral quadric
has rank at least three). A section $y\in S_2\setminus B_2$ must
lie outside $\operatorname{Frac}(B)$. Hence $B\oplus By\subset S$
is a direct sum with the full Hilbert series $(1+z^2)H_B=H_S$.
Completing the square in its monic quadratic relation gives
$S=B[y]/(y^2-f_4)$; this is the asserted flat double cover.
\end{proof}

\begin{corollary}\label{cor:nov30-gorenstein-quartic-index}
Let $n\ge9$ and $\gcd(n-2,30)=1$. The conclusion of
Theorem~\ref{thm:nov30-quartic-cartier-k3} holds for any
Gorenstein K-semistable quartic replacement with its limiting
integral $\mathbf Q$-Cartier Weil class $L$. In particular it
holds for every Gorenstein K-polystable quartic ninefold replacement.
\end{corollary}
\begin{proof}
At a point of local index $I$ for $L$, its degree-$I$ cyclic
quasi-\'etale index cover has normalized volume at most $n^n$.
The finite-degree formula and the local-to-global lower bound
\cite[Theorems~2.11(2), 2.13, 2.14]{ar2007.14320} therefore give
\[
    I\le\frac14\left(\frac{n+1}{n-2}\right)^n
      \le\frac{10^9}{4\cdot7^9}<7.
\]
For the middle inequality write $t=3/(n-2)$. The derivative of
$n\log((n+1)/(n-2))$ with respect to $n$ is negative because
$\log(1+t)<t(1+2t/3)/(1+t)$ for $t>0$; the difference on the
right minus the left has value zero at zero and derivative
$(t+2t^2)/(3(1+t)^2)>0$. Since $-K_Z\sim(n-2)L$ is Cartier,
$I$ divides $n-2$. Coprimality forces $I=1$ at every point, and
the theorem applies.
\end{proof}

\begin{corollary}\label{cor:nov30-optimal-quartic-cartier}
If a smooth quartic $X$ of dimension $n\ge9$ fails uniform
K-stability, the limiting hyperplane Weil divisor on the central
fiber of every optimal destabilizing special test configuration
is non-Cartier.
\end{corollary}
\begin{proof}
Put $N=n+1$ and $\beta=\delta(X,H)$. Such optimal configurations
exist by divisorial attainment \cite{ar2102.09405} and
\cite[Theorem~1.1]{ar1907.05399}. The latter theorem gives
$\delta(Z,L)=\beta$ on their central fibers. The pointed complete
hyperplane adjunction inequality applied repeatedly down to a
smooth quartic eightfold gives $\beta\ge 2N/3$. Indeed each
smooth pointed hyperplane section from dimension $j$ to $j-1$
lifts the local threshold bound $b$ to
$\min\{j+1,(j+1)b/j\}$; at dimension eight the hyperplane
threshold is at least six by Theorem~\ref{thm:nov28-quartic-eightfold}.
Such smooth flags exist through every point by Bertini. Taking the
infimum over points gives the stated non-strict bound, sufficient
for what follows.

The proof of Theorem~\ref{thm:nov30-quartic-cartier-k3} eliminates
the Cartier basepoint case whenever $\delta(Z,L)>4N/7$: a doubled
prime ruling gives $S_L(E)\ge2/N$ by its powers, and the Segre
scroll gives $S_L(E)\ge7/(4N)$. Both have $A_Z(E)=1$.
Thus a Cartier $L$ here would be globally generated, since
$2/3>4/7$. The same Hilbert-series argument makes its section
ring an embedded quartic ring or a weighted $(2,4)$
complete-intersection ring. For the hyperplane filtration of
the fixed smooth $X$ this is its associated graded ring. By
Corollary~\ref{cor:nov30-all-ci-hyperplane-degenerations} the
nontrivial configuration would then have positive
Donaldson--Futaki invariant, contradicting optimality.
\end{proof}

\subsection{Hyperplane images of optimal quartic degenerations}

\begin{corollary}\label{cor:nov30-optimal-image-ranks}
Suppose that a smooth quartic $X$ of dimension $n\ge9$ fails
uniform K-stability and that the image of the hyperplane system
on an optimal destabilizing special fiber $(Z,L)$ has dimension
$n$. Then this image is a quadric of rank four, five or six, and
the hyperplane map is birational.
\end{corollary}
\begin{proof}
In the notation of Corollary~\ref{cor:nov30-optimal-quartic-cartier},
one has the strict bound $\beta>2(n+1)/3$. Indeed take a closed
point in the center of a divisorial minimizer on $X$ and a finite
smooth pointed hyperplane flag to a quartic eightfold. The local
threshold at this point is at most the value of the minimizing
divisor. The eightfold hyperplane threshold is strictly greater
than six and the flag estimate gives the strict inequality at
this point. Write $N=n+1$ below. A complete-intersection full
central section ring is impossible by
Corollary~\ref{cor:nov30-all-ci-hyperplane-degenerations}.

The degree of the integral nondegenerate hypersurface image is two,
three or four by its leading Hilbert coefficient. In degree four
equality of Hilbert series makes the ring a hypersurface ring. In
degree three follow Proposition~\ref{prop:nov30-limit-degree-minus-one}
up to the presentation by $\ell y=g$. If $\ell,g$ is a regular
sequence, that proof gives a complete intersection; otherwise the
hyperplane section is nonintegral and has a linear component.
Proposition~\ref{prop:nov30-linear-image-component} constructs
a prime divisor with $N S_L(E)\ge3/2$ without using semistability.
This contradicts $\delta(Z,L)=\beta>2N/3$.

In degree two, if the generic map has degree two, let $J$ be its
resolved moving hyperplane class and $P$ the pullback of $L$.
They are big and nef with $P-J$ effective up to $\mathbf Q$-linear
equivalence and $P^n=J^n=4$. All intervening mixed intersections
therefore equal four. The proportionality theorem for big and nef
classes (B.~F.~Teissier's equality criterion for the
Khovanskii--Teissier inequalities) gives $P\equiv J$. The
effective difference is numerically trivial and hence zero.
Consequently the products of all sufficiently divisible degrees
in the linear sections have no base locus. At each point some
power of a single degree-one section is nonvanishing in a Cartier
multiple of $L$; its Weil zero divisor is then empty locally, so
it trivializes $L$ there. This contradicts
Corollary~\ref{cor:nov30-optimal-quartic-cartier}.

For a birational quadric of rank at least seven repeat the ring
calculation of Proposition~\ref{prop:nov30-quartic-birational-quadric}.
Both the degree-two denominator case and the regular final
degree-one denominator case give complete-intersection rings.
In the other case the fixed prime of the linear pencil satisfies
$N S_L(E)\ge3/2$, again impossible. Finally rank three gives a
doubled general ruling, with $N S_L(E)\ge2$, equally impossible.
\end{proof}

\begin{proposition}\label{prop:nov30-optimal-minimal-image}
Let $(Z,L)$ be an optimal special degeneration of a smooth quartic
with $n\ge9$ and $\delta(X,-K_X)\le1$. A lower-dimensional image
of $|L|$ cannot have minimal degree. The same conclusion holds for
a K-polystable quartic replacement. Neither assertion assumes $L$
Cartier.
\end{proposition}
\begin{proof}
Put $P=\pi^*L$ on smooth resolutions and $N=n+1$. Use
$\beta=\delta(Z,L)>2N/3$ for the optimal degeneration, as in
Corollary~\ref{cor:nov30-optimal-image-ranks}; in the polystable
case K-semistability gives $\delta(Z,L)\ge n-2$, so we may take
$\beta=n-2>2N/3$. In both cases $\beta S_L(E)\le A_Z(E)$.
The determinant slope inequality
\eqref{global:eq:optimizer-covering-slope} on a dominant map to
$\mathbf P^s$ bounds its nef hyperplane pullback $V$ by
\[
             \frac{P^{n-1}V}{P^n}\ge\frac{s\beta}{(s+1)n}.
\]
For an optimal degeneration take $t$ increasing to $\beta/(n-2)$.
For K-polystability this is the tangent-sheaf semistability argument
of Lemma~\ref{cubics:slope}, with index $n-2$.

Here is a second intersection estimate for $s\ge2$. Project to
$\mathbf P^2$, resolving the projection; let $V'$ be its moving
hyperplane class. The difference $V-V'$ is effective up to linear
equivalence, so $P^{n-2}V^2\ge P^{n-2}(V')^2$. Take one connected
component $F$ of a general fiber over $\mathbf P^2$. Generic
smoothness makes $F$ smooth with trivial rank-two conormal bundle.
It meets the isomorphism locus of the resolution of $Z$ densely,
since the general fiber of each exceptional divisor has strictly
smaller dimension. The exceptional divisor of its blowup thus has
discrepancy two over $Z$. Put $x=P^{n-2}F/P^n>0$. For divisible
$m$, vanishing to order $k$ along $F$ cuts out at most
$\binom{k+1}{2}h^0(F,mP)$ independent conditions: the first $k$
conormal layers have ranks $1,2,\ldots,k$ and are trivial before
tensoring by $mP$. The restriction $P|_F$ is big and semiample, so
$h^0(F,mP)=(P^{n-2}F)m^{n-2}/(n-2)!+O(m^{n-3})$. With $k/m\to u$
this gives the normalized volume lower bound
$1-\binom n2 xu^2$ for the vanishing filtration. Integration up to
its first zero gives $S_L(\operatorname{ord}_F)\ge
2/(3\sqrt{\binom n2 x})$. The discrepancy inequality now yields
\[
       \frac{P^{n-2}V^2}{P^n}\ge x\ge\frac{\beta^2}{9\binom n2}.
\]

The classification of varieties of minimal degree in codimension
at least two gives scrolls and cones over scrolls or the Veronese
surface \cite{EH}. In the Veronese case or for a scroll with a
summand of degree at least two, a hyperplane contains a general
ruling divisor (respectively the lift of a line on the plane)
twice. As in Section~\ref{cubics:rulings}, pullback on a dominant
smooth open gives a prime divisor $E$ with $NS_L(E)\ge2$,
contradicting $\beta>N/2$.

The remaining scroll is a cone over $\mathbf P^1\times\mathbf P^s$,
where $s+1$ is its degree and $s\ge2$. With $U,V$ the two ruling
classes on a resolution, $P-U-V$ is effective up to
$\mathbf Q$-linear equivalence by the Segre coordinate subsystem.
A general first-ruling fiber gives a prime component $E$ on $Z$.
The standard volume integral with $P-E'-V$ effective gives
\[
 NS_L(E)\ge1+\frac{P^{n-1}V+P^{n-2}V^2}{P^n}
 \ge1+\frac{2\beta}{3n}+\frac{\beta^2}{9\binom n2}
 >1+\frac49+\frac8{81}>\frac32.
\]
This contradicts $\beta S_L(E)\le1$.
\end{proof}

In particular a quartic optimal degeneration as above has hyperplane
image of dimension at least four, and a quartic K-polystable
replacement as above has image of dimension at least $\lceil n/2\rceil$.
Indeed, if its dimension is $a<n$ with $a\beta/n>a-1$, the proof in
Section~\ref{cubics:minimal} applies verbatim with this determinant
slope in place of $a(n-1)/n$. A resolution $Y$ of the image, with
hyperplane pullback $M$, has no section of $K_Y+(a-1)M$. Restriction
to successive smooth hyperplane sections and vanishing give a rational
general curve section and $\deg W=\operatorname{codim}W+1$, the
excluded minimal-degree case. The numerical inequality holds for
$a\le3$ when $\beta>2(n+1)/3$, and for $a<n/2$ when $\beta=n-2$.

\begin{proposition}\label{prop:nov30-optimal-quadratic-image}
In the optimal-destabilization setting of
Corollary~\ref{cor:nov30-optimal-image-ranks}, suppose the
degree-two section system $|2L|$ has $n$-dimensional image. Then
the degree-one image cannot be a complete intersection of $c\ge2$
quadrics if every nonzero quadric in their span has rank at least
$2c+5$, or if that complete intersection is smooth of dimension
at least five. In particular the smooth case applies to both
$c=2$ and $c=3$ in dimension $n=9$.
\end{proposition}
\begin{proof}
Let $T$ be the coordinate ring of the proposed image, and $S$ the
full hyperplane section ring. We will need that $T$, and $T/h$ for
each nonzero linear form $h$, are factorial, and that $T/(h,\ell)$
is a domain for independent linear forms $h,\ell$. Here is a rank
check. Restrict the quadrics to a vector subspace of codimension
$t\le2$ in the $m=n+2$ variables. Every nonzero combination has
rank at least $R-2t$, where $R=2c+5$. The locus of dependent
gradients in its projective space lies in the union of the
projectivized kernels of this $(c-1)$-dimensional system, so has
dimension at most $m+t-R+c-2$. This is smaller than the expected
intersection dimension $m-t-c-1$. A component exceeding that
expected dimension would lie in the dependent-gradient locus;
thus these sections are complete intersections. Their singular
loci have codimension at least $R-2t-2c+1$. The affine cones are
normal also at the vertex by Cohen--Macaulayness and Serre's
criterion; equivalent graded connectedness gives domains. For
$t=0,1$ the projective intersections are regular in codimension
three, hence locally factorial by complete-intersection
parafactoriality. Their Picard groups are generated by the
hyperplane class by complete-intersection Lefschetz, so their
homogeneous coordinate rings are factorial by the cone class-group
formula. The dimension requirements hold since $m\ge R$.

If the complete intersection $W$ is smooth, the same conclusions
follow without a rank assumption. On $W\times\mathbf P^{c-1}$
the gradients of a linear combination of the quadrics define
an everywhere-defined morphism to the dual projective space,
since the defining gradients are independent at every point of
$W$. Its hyperplane pullback is
$\mathcal O_W(1)\boxtimes\mathcal O_{\mathbf P^{c-1}}(1)$, which
is ample. The morphism therefore has finite fibers. A hyperplane
section is singular only in the contact locus of that hyperplane,
which is finite. On a section of codimension two, the dependent
gradients can occur only in the union of contact loci for the
pencil of hyperplanes through the given linear subspace. This
union has dimension at most one. The first section is a normal
integral complete intersection spanning its hyperplane, and the
second is therefore a complete intersection as well. Since
$\dim W\ge5$, the first section has singular locus of codimension
at least four and the second of codimension at least two. The
same factoriality and domain conclusions follow.

There are $c$ independent extra degree-two sections modulo $T_2$.
The ratio of elements of $|2L|$ recovers the ratios of all original
coordinates by taking products with a fixed linear section. Hence
the full-dimensionality assumption gives $c-1$ degree-two extras
algebraically independent over $T$. Adjoin them to form a polynomial
extension $B=T[y_1,\ldots,y_{c-1}]$ with each $\deg y_i=2$.
Then
\[
                H_B(z)=\frac{1-z^2}{(1-z)^{n+2}},
             \qquad H_S(z)=(1+z^2)H_B(z).
\]
The ring $B$ and its linear sections satisfy the factoriality and
domain conditions just established for $T$. There is one further
degree-two generator $v\in S_2\setminus B_2$. If
$v\notin\operatorname{Frac}(B)$ then $B\oplus Bv=S$ by Hilbert
series, and the monic quadratic relation of $v$ gives a complete
intersection, which the optimal configuration excludes.

For $v\in\operatorname{Frac}(B)$ write $v=g/h$ homogeneously in
coprime terms. As in Proposition~\ref{prop:nov30-quartic-birational-quadric},
the multiplicity of $B[v]=B[V]/(hV-g)$ is $2+\deg h$, forcing
$\deg h=1$ or $2$. Each ring here is built over the original
complete intersection $T$. When $\deg h=2$, its Hilbert series
equals $H_S$, yielding an excluded complete intersection. For
$\deg h=1$, if $g\bmod h$ has no linear factor it is prime in
$B/h$, as its degree is three. The proof for $k=1$ in the cited
proposition now applies without change. Explicitly $B[v]$ has the
Hilbert series of a cubic. An additional generator in degree three
satisfies $\ell w=f$ with $f\notin\ell B[v]$. The ring
$B[v]/\ell$ is a domain, using $B/h$ factorial when $\ell$
is proportional to $h$, and $B/(h,\ell)$ a domain and
$g\notin(h,\ell)$ otherwise. The regular-sequence presentation
of $B[v,w]$ therefore has Hilbert series $H_S$, making $S$ again
a complete intersection.

It remains that $g=hp+bC$ with $b,h$ independent linear forms.
Then $C/h=(v-p)/b$ is a rational section of $L$ outside
$S_1=B_1$, so $h,b$ vanish along a common prime divisor $E$.
To measure their pencil, resolve the dominant rational map from
$Z$ to $\operatorname{Proj}B$. Denote the pullback of its rational
hyperplane class by $J$ and that of the pencil by $M$. The class
$2J$ is globally generated by the degree-two pieces of $B$ since
the algebra has generators of degrees one and two. Thus with
$P=\pi^*L$ the difference $P-J$ is effective up to
$\mathbf Q$-linear equivalence. On $\operatorname{Proj}B$ the
linear forms $h,b$ have no common divisorial component, since
$B/(h,b)$ has codimension two. Consequently the degree of a
general pencil fiber against the $(n-1)$st power of the rational
hyperplane class is the volume of that hyperplane class, namely
two by $H_B$. The projection formula for the generically finite
map therefore gives $P^{n-1}M\ge J^{n-1}M\ge2$. Integrating
the same volume inequality as for the linear component in
Proposition~\ref{prop:nov30-linear-image-component} gives
$(n+1)S_L(E)\ge1+2/4$, contradicting $\beta>2(n+1)/3$.

\end{proof}

In dimension nine there is no degree-two-system qualification for
an image with finite polarized automorphism group. Indeed every
divisorial minimizer on the smooth quartic ninefold has point
center by the descent argument following
Section~\ref{sec:nov27-third-descent}. The valuation weights on
$H^0(X,H)$ are nonconstant: a hyperplane missing the center has
weight zero and a hyperplane through it has positive weight.
The associated projective one-parameter action on the hyperplane
image is therefore nontrivial, since a nonscalar projectivity
cannot fix an integral nondegenerate projective variety pointwise.
Thus its polarized stabilizer contains a one-dimensional torus.
In particular no smooth complete intersection of codimension at
least two and dimension at least three can be that image, by
\cite[Theorem~2.3]{ar2608.09293} (the finiteness theorem of
Libgober--Wood and Benoist, as cited there).

\section{Discriminants and comparison of Fano fillings}

\subsection{The discriminant of a successive monic presentation}

\begin{lemma}\label{lem:nov30-tower-discriminant}
Put $N=n+1$ and $D=(N+1)(d-1)^N$, the degree of the discriminant
$\Delta_d$ of degree-$d$ forms in $N+1$ variables. Consider the
tower with $2\le a_1\mid a_2\mid\cdots\mid a_s\mid a_{s+1}=d$
and all successive degrees strictly increasing. In variables
$x$ of degree one and $y_i$ of degree $a_i$, let
\[
  \tau_i y_i=P_i(x,y_1,\ldots,y_{i-1})\quad(1\le i\le s),
  \qquad P_{s+1}(x,y_1,\ldots,y_s)=0,
\]
where $P_i$ has section degree $a_i$, $P_1=q(x)$, and $P_i$
for $i\ge2$ is monic in $y_{i-1}$. Let $F(x)$ be the eliminated
degree-$d$ form, clearing denominators by
$\prod_i\tau_i^{d/a_i}$. In the universal polynomial coefficient
ring one has
\[
        \prod_{i=1}^s\tau_i^{\mu_i}\ \mid\ \Delta_d(F),
 \qquad
 \mu_i=\frac d{a_i}\left[D-
      \sum_{j=0}^{N}(d-1)^j(a_i-1)^{N-j}\right].
\]
In particular, after arbitrary specialization of both the
coefficients and the parameters to a complex discrete valuation
ring, the discriminant order is at least the sum of the
$\mu_i\operatorname{ord}(\tau_i)$ when all parameters are nonzero.
\end{lemma}
\begin{proof}
We first treat one parameter $t$, putting $a=a_1$, $m=d/a$.
The universally eliminated form is
\[
      F_t=q^m+\sum_{j=1}^m t^j f_{ja}(x)q^{m-j},
\]
where $f_{ja}$ has degree $ja$. The order of $\Delta_d(F_t)$
at the generic values of the coefficients is the minimum exponent
of $t$ in the universal polynomial. To compute it choose general
complex coefficients. In $\mathbf P(1^{N+1},a)$ the hypersurface
\[
     Y:\quad y^m+\sum_{j=1}^m f_{ja}(x)y^{m-j}=0
\]
is smooth, avoids the weighted vertex by monicity, and has
dimension $N$. Write $H=\mathcal O_Y(1)$. The line bundle $aH$
is very ample, since the $a$th Veronese of the ambient weighted
polynomial ring is generated by its degree-$a$ piece. All degree
$a$ sections are represented by $q'(x)+c y$; restriction from
the ambient space is surjective by the hypersurface sequence
(or the Cohen--Macaulay section ring). Consider the pencil
$q-t y$. Its endpoints $q=0$ and $y=0$ give smooth divisors
for general coefficients; for the first use Bertini on the
basepoint-free pullbacks from $\mathbf P^N$, and for the second
use smoothness of the general ordinary hypersurface $f_d=0$.
The pencil is a Lefschetz pencil for general coefficients. In
detail, a general pencil in $|aH|$ on a general $Y$ may be written
as $\langle q(x),y+q'(x)\rangle$: its underlying two-plane meets
the hyperplane of pullback sections in a line. The weighted
coordinate change $y\mapsto y+q'(x)$ preserves the parameter
space of monic equations and converts this to the indicated
form. The Lefschetz condition is open in this irreducible
parameter space and is attained by this change of variables.

The number of singular members of a general such pencil is
$M_Y=\int_Y c_N(J^1\mathcal O_Y(a))$. Indeed the universal
first jet on $Y\times\mathbf P^1$ counts the isolated singular
points of the pencil. By the characteristic-zero Lefschetz
pencil theorem each is an ordinary double point in a single
transverse singular member. The displayed Chern number is
positive, so the dual is a hypersurface. The Euler sequence on
the smooth ambient locus and the normal bundle sequence give
$c(T_Y)=(1+H)^{N+1}(1+aH)/(1+dH)$ and $H^N=m$. The first jet
sequence therefore gives
\[
 M_Y=m[z^N]\frac{(1-z)^{N+1}}{(1-dz)(1-az)}
    =m\sum_{j=0}^{N}(d-1)^j(a-1)^{N-j}.
\]
For the last identity use the binomial theorem in the form
$\sum_{k\ge0}t^k[z^k](1-z)^k A(z)=A(t/(1+t))/(1+t)$.
Every finite nonzero singular parameter of the pencil corresponds
to a simple zero of the ordinary discriminant: eliminating $y$
is an isomorphism of families there, and at a unique ordinary
double point a transverse smoothing has a simple discriminant
zero. Also $\Delta_d(F_t)$ has degree exactly $mD$ in $t$, since
its top coefficient is $\Delta_d(f_d)\ne0$. Neither endpoint of
the pencil is singular on $Y$. Subtracting the $M_Y$ finite
nonzero zeros gives the order $mD-M_Y$ at $t=0$. Consequently
the universal polynomial is divisible by that power of $t$,
with no smoothness requirement on specialized coefficients.

For several parameters write
$\sigma_i=\prod_{j\le i}\tau_j^{a_i/a_j}$. Substituting the
earlier expressions $y_j=A_j/\sigma_j$ into
$(\sigma_i/\tau_i)P_i$ yields a polynomial $A_i$, as each
monomial has section degree $a_i$. The same observation at
degree $d$ shows that $F$ is polynomial in all the coefficients
and parameters. Its specialization at all $\tau_i=0$ is
$q^{d/a_1}$ by the successive monic terms, so this denominator
clearing is primitive when the reduction of $q$ is nonzero.
Now fix $i$ and invert all parameters other than $\tau_i$.
Eliminating the earlier variables leaves an equation
$\tau_i y_i=q_i(x)$, where $q_i$ has degree $a_i$ and is
independent of $\tau_i$. Eliminating the later variables
leaves a degree-$d$ equation in $x,y_i$, independent of
$\tau_i$. Its coefficient of $y_i^{d/a_i}$ is a scalar in
the field of the other parameters and universal coefficients;
it is nonzero by specializing the later equations to pure
monic powers. Work over that field. Dividing by the
leading coefficient gives the one-parameter calculation of
degree $a_i$, up to a scalar independent of $\tau_i$. It
follows that $\tau_i^{\mu_i}$ divides the universal discriminant
in this one-variable polynomial ring. Its coefficients of lower
$\tau_i$-degree therefore vanish identically in the other
parameters. Since the $\tau_i$ are independent prime elements,
the product divides the discriminant in the full polynomial
ring. Specialization gives the claimed valuation bound.
\end{proof}

\subsection{Comparison with a smooth hypersurface filling}
\label{sec:nov30-finite-cm-model}
\begin{theorem}\label{thm:nov30-finite-cm-model}
Consider two polarized Fano fillings over a smooth complex pointed
curve, identified off the marked point. The first is a smooth
family of degree-\(d\) hypersurfaces at that point,
\(3\le d\le n+1\). Suppose the special fiber of the second is
K-semistable. If its hyperplane root \(L\) is Cartier and
globally generated, the second special fiber is the same polarized
hypersurface as the first. In particular a distinct K-polystable
replacement over a smooth hypersurface GIT point has a base
point of its limiting degree-one system.
\end{theorem}
\begin{proof}
Pass to the local discrete valuation ring \(R\) with parameter
\(t\), allowing finite base change. The hyperplane class extends
on the second model as an ample Cartier class in a neighborhood
of this fiber, by the reflexive restriction and lifting argument
in Proposition~\ref{prop:nov30-factorial-quartic-replacement};
there only Cartierness of the restricted class is needed for
the lifting step. Write \(S_R\) for its section algebra. Its
graded pieces are free and commute with base change by Fano
vanishing. The central algebra satisfies the hypotheses of
Theorem~\ref{thm:nov30-finite-tower}, as in
Corollary~\ref{cor:nov30-basepoint-necessary}.

Suppose first its tower is nontrivial. Write \(a_1,\ldots,a_s\)
for the degrees of the additional generators in increasing order,
\(a_{s+1}=d\) and \(m_i=a_{i+1}/a_i\). Lift these generators
\(y_i\) and the degree-one basis \(z\) to the relative algebra.
They generate in each degree by Nakayama. The kernels of the
relative presentation specialize onto the central kernels, by
freeness of every section space. Thus one can lift the \(s+1\)
central equations. They generate the relative ideal by graded
Nakayama in each degree. The equation of degree \(a_i\)
only involves \(z,y_1,\ldots,y_i\) and the term in \(y_i\) is
linear with a coefficient zero modulo \(t\). For \(i>1\) the
coefficient of \(y_{i-1}^{m_{i-1}}\) is a unit by the central
tower. In degree \(d\) the coefficient of \(y_s^{m_s}\) is
likewise a unit. Rescale the equations by these units to write
\[
     P_i(z,y_{<i})=\tau_i y_i\quad(1\le i\le s),
     \qquad P_{s+1}(z,y_1,\ldots,y_s)=0 .
\]
Here \(P_i\) is monic in \(y_{i-1}\) for \(i>1\); \(P_1\)
reduces to the irreducible degree-\(a_1\) image equation \(q\).
All \(\tau_i\) are nonzero: on the generic fiber degree one
generates the algebra. In degree \(a_i\), modulo the products
of lower-degree equations, only the indicated equation can
eliminate the independent linear monomial \(y_i\). Put
\(j_i=\operatorname{ord}_t(\tau_i)>0\).

Let \(F_B\in R[z]\) be the primitive degree-\(d\) polynomial
obtained by eliminating the \(y_i\). In detail put
\(\sigma_i=\prod_{j\le i}\tau_j^{a_i/a_j}\) and recursively
\[
    A_i=(\sigma_i/\tau_i)
           P_i(z,A_1/\sigma_1,\ldots,A_{i-1}/\sigma_{i-1}).
\]
Then \(A_i=\sigma_i y_i\) in \(S_R\). All these are polynomials
in \(z\) since each monomial has the prescribed section degree.
Their reductions are \(q^{a_i/a_1}\), by the successive
monic equations. The same recursion at degree \(d\), multiplying
\(P_{s+1}\) by \(\prod_j\tau_j^{d/a_j}\), gives \(F_B\) with
reduction \(q^{d/a_1}\), in particular primitive.

Write \(m=n+2\), \(N=n+1\), \(D=m(d-1)^N\), and \(\Delta\)
for the hypersurface discriminant. The discriminant multiplicity
estimate of Lemma~\ref{lem:nov30-tower-discriminant} gives
\begin{equation}
  \operatorname{ord}_t\Delta(F_B)\ \ge\
       \sum_{i=1}^s j_i\,\mu(a_i,d),\qquad
  \mu(a_i,d)> dN(d-1)^{N-1}(d/a_i-1).
  \label{eq:nov30-cm-discriminant}
\end{equation}
For the strict inequality, set \(a=a_i\), \(h=d/a\). Subtracting
\(dN(d-1)^{N-1}(h-1)\) from the formula for \(\mu_i\) in that
lemma gives exactly
\[
  h(d-1)^{N-1}(a-1)
   \left[N-\sum_{\ell=0}^{N-1}
       \left(\frac{a-1}{d-1}\right)^\ell\right]>0 ,
\]
since \(2\le a<d\) and \(N\ge2\).
We compare section lattices explicitly. The embedded algebra
\(E_R=R[z]/(F_B)\) injects into \(S_R\). Choose \(n+1\)
linear parameters \(u_0,\ldots,u_n\) so that the remaining
coordinate \(x\) gives \(q\) with nonzero coefficient of
\(x^{a_1}\). Over \(U_R=R[u_0,\ldots,u_n]\), the sections
\[
       x^k\prod_i y_i^{h_i},
       \qquad 0\le k<a_1,\quad 0\le h_i<m_i ,
\]
give a graded free basis of \(S_R\): this is the monic tower
basis on the special fiber, and in each section degree the
lifts are a basis by Nakayama and equality of ranks. Replacing
each \(y_i\) by \(A_i\) gives a free basis of \(E_R\) over
\(U_R\). Indeed their degrees range once over the integers
from zero through \(d-1\); on reduction the coefficient at
the corresponding highest power of \(x\) is nonzero. The
ordinary powers \(1,x,\ldots,x^{d-1}\) are already a basis
by the unit leading coefficient of \(F_B\). Thus the inclusion
of the two bases is diagonal, with entries
\(\prod_i\sigma_i^{h_i}\).
The mean length per section in degree \(k'\) of the quotient
lattice \((S_R)_{k'}/(E_R)_{k'}\) tends as \(k'\to\infty\) to
\[
 \frac12\sum_i (m_i-1)\operatorname{ord}_t\sigma_i
    =\frac12\sum_i (d/a_i-1)j_i .
\]
Each of the \(d\) basis summands has asymptotically the same
dimension in degree \(k'\). The equality telescopes using
\(\sum_{i\ge j}(m_i-1)a_i/a_j=d/a_j-1\).

It remains to compare the two embedded algebras. Smith bases
in degree one write \(z_i=t^{-w_i}x_i'\), where \(x_i'\)
are the first model's coordinates. A common shift allows
\(\min w_i=0\); put \(W=\sum_i w_i\). Up to a unit
\(F_B(z)=t^{-b}F_A(t^{w_i}z_i)\) for the primitive equation
\(F_A\) of the first family. Its discriminant is a unit by
smoothness. Hence the determinant character of the discriminant
gives \(\operatorname{ord}_t\Delta(F_B)=D(dW/m-b)\).
With \(C_k=\binom{k+m-1}{m-1}\), \(h_k=C_k-C_{k-d}\), the
relative determinant degree of the embedded section lattices is
\[
       (kW/m)h_k+(dW/m-b)C_{k-d}.
\]
This follows by subtracting the lattice change for the multiples
of their primitive equations from the lattice change in the
polynomial spaces. The doubled, normalized CM degree difference
is minus twice the constant term after division by \(h_k\), by
the Hilbert determinant formula for the CM line. As in
Section~\ref{sec:nov30-ambient-futaki}, this yields for the
second filling minus the first
\[
 \frac{(d-1)m}{dND}\operatorname{ord}_t\Delta(F_B)
            -\sum_i(d/a_i-1)j_i .
\]
A common tensor power giving the anticanonical polarization
does not change this sign. By \eqref{eq:nov30-cm-discriminant}
the difference is strictly positive. A filling with K-semistable
special fiber minimizes the CM degree among Fano fillings with
the prescribed punctured family, a contradiction
(C.~Xu, \emph{K-stability of Fano varieties: an algebro-geometric
approach}, EMS Surv. Math. Sci. \textbf{8} (2021), 265--354,
Remark~7.9, \href{https://doi.org/10.4171/EMSS/51}{doi:10.4171/EMSS/51}).
For the curve formulation one first completes the common family
after finite base change by C.~Li and C.~Xu,
\emph{Special test configuration and K-stability of Fano varieties},
Ann. of Math. \textbf{180} (2014), 197--232, Theorem~1,
\href{https://doi.org/10.4007/annals.2014.180.1.4}{doi:10.4007/annals.2014.180.1.4}.
Glue in the second filling only at the chosen point of the
normalized base curve. The resulting complete families coincide
elsewhere, and the local determinant difference is multiplied
by the positive ramification index.

If the tower is trivial, the relative hyperplane system embeds
the second family as degree-\(d\) hypersurfaces. The equation of
its K-semistable fiber is GIT-semistable by the ambient weight
calculation in the proof of
Proposition~\ref{prop:nov30-factorial-quartic-replacement},
which applies in every degree. The smooth first fiber is
GIT-stable. Separatedness of the GIT quotient identifies the
two points, and the fiber of the quotient at a stable point is
a single orbit. This proves the result.
\end{proof}

\subsection{Discriminant multiplicities without monic equations}

\begin{lemma}\label{lem:nov30-nonmonic-discriminant}
Write $N=n+1$, $D=(N+1)(d-1)^N$. Given $2\le a\le d$ set
$m=\lfloor d/a\rfloor$, $h=d-ma$, and
$B_N(u,v)=\sum_{i=0}^N(u-1)^i(v-1)^{N-i}$. For arbitrary forms
$q_a,f_{d-ia}$ of the indicated degrees, the polynomial
\[
   F_t=\sum_{i=0}^m t^{m-i}q_a^i f_{d-ia}
\]
satisfies $t^{\mu(a)}\mid\Delta_d(F_t)$ universally, with
\[
 \mu(a)=mD-M(a),\qquad
 M(a)=\frac d a B_N(d,a)-\frac h a B_N(h,a).
\]
When $h=0$ interpret the second term as zero. For $a>d$ set
$m=0$, $M(a)=\mu(a)=0$. For any triangular system
$\tau_i y_i=P_i(x,y_{<i})$, $H_d(x,y)=0$ with
$\deg y_i=\deg P_i=a_i\ge2$, allow all the coefficients to vary.
Its eliminated form $F$ after multiplication by
$\prod_i\tau_i^{\lfloor d/a_i\rfloor}$ is polynomial and obeys
$\prod_i\tau_i^{\mu(a_i)}\mid\Delta_d(F)$ in the universal
coefficient ring. Monicity and divisibilities between the $a_i$
are unnecessary.
\end{lemma}
\begin{proof}
For $h=0$ this is the one-parameter computation of
Lemma~\ref{lem:nov30-tower-discriminant}, since the coefficient
of the highest power is generically a unit if $a<d$.
For $a=d$ the formula gives $\mu(d)=0$ directly and the
claim is automatic. Suppose $h>0$, and take general complex
coefficients; the generic order is the minimum exponent in
the universal polynomial. Compactify the total space of
$\mathcal O_{\mathbf P^N}(a)$ to its projective line bundle.
Write $I$ for the infinity divisor and $H$ for the base hyperplane
class, and write $z,y$ for the canonical relative coordinates of
classes $I,I+aH$, respectively; $I(I+aH)=0$. The homogenized
spectral hypersurface
\[
      Y:\quad\sum_{i=0}^m f_{d-ia}(x)y^i z^{m-i}=0
\]
is smooth of class $mI+dH$. Indeed this system is basepoint-free
including at infinity, where its coefficient has positive degree
$h$. The pencil $q_a z-ty$ on $Y$ has smooth base locus
$Q_a\cap D_d$ in the zero section, where $Q_a=(q_a=0)$ and
$D_k$ denotes the smooth ordinary hypersurface cut by the
corresponding general form of degree $k$. Its infinity member
$y=0$ is smooth. Its zero member is the transverse union of the
infinity divisor of $Y$ and the pullback of $Q_a$, with intersection
$Q_a\cap D_h$ disjoint from the pencil base, smooth by Bertini.
The two components are smooth by the same generality and meet
transversely, since $D_h,Q_a$ meet transversely on the infinity
section. After blowing up the pencil base, the generic fiber is
a smooth ordinary $D_d$, and the Euler characteristic of the zero
fiber exceeds that of $D_d$ by $\chi(Q_a\cap D_h)$, by the
normal-crossing smoothing formula. The base blowup does not
change the individual pencil members because its ideal is Cartier
in each member near the base.

Away from zero and infinity the only singular fibers for general
data have a single ordinary double point and smooth transversely
in $t$. Here is a parameter count validating this claim for the
specified pencils. In the chart $ty=q_a$ elimination gives
ordinary degree-$d$ forms with freely variable $f_d$.
Their first jets at two distinct points and their second jet at
one point are surjective, so simultaneous singularities or a
degenerate Hessian impose at least two conditions after allowing
the marked points to vary at fixed $t$. The total pencil family
is smooth after blowing up its smooth base on the smooth $Y$;
therefore a remaining isolated nondegenerate critical point is
a transverse one-parameter smoothing. Each such member gives
one simple discriminant root. The affine polynomial in $t$ has
discriminant of degree $mD$, as its leading coefficient is
$\Delta_d(f_d)\ne0$. The number of its nonzero roots is therefore
\[
 M=(-1)^N\bigl[\chi(Y)+\chi(Q_a\cap D_d)-2\chi(D_d)
                         -\chi(Q_a\cap D_h)\bigr],
\]
using additivity of Euler characteristic on the pencil and the
Euler jump $(-1)^N$ for an ordinary double point.

In the projective bundle, integration of a class $g(I)$ down to
the base is $(g(0)-g(-aH))/(aH)$, and the relative tangent
Chern class is $(1+I)(1+I+aH)$. Hence $\chi(Y)$ is the degree-$N$
coefficient in $(1+H)^{N+1}$ times
\[
   \frac {d(1+aH)}{a(1+dH)}-\frac {h(1-aH)}{a(1+hH)}.
\]
For an intersection $Q_a\cap D_k$, the corresponding factor is
$akH^2/((1+aH)(1+kH))$; for $D_d$ it is $dH/(1+dH)$.
Substitution and cancellation give
\[
 M=m[z^N]\frac{(1-z)^{N+1}}{(1-az)(1-dz)(1-hz)}.
\]
The binomial-series identity used in
Lemma~\ref{lem:nov30-tower-discriminant} and partial fractions
in $d-1,h-1$ give exactly the asserted $M(a)$. Thus all lower
powers of $t$ vanish universally, not just for general coefficients.

For the multivariable claim consider any fully expanded monomial
after successive graph substitutions. Each occurrence of a
particular $y_i$ consumes section degree $a_i$, and none of its
descendants in the substitution contains $y_i$ again, by
triangularity. Distinct occurrences thus consume disjoint portions
of the original degree $d$. No denominator can contain more
than $\lfloor d/a_i\rfloor$ powers of $\tau_i$, proving
polynomiality after the prescribed clearing. Invert all parameters
other than $\tau_i$. Eliminating earlier variables leaves
$\tau_i y_i=q(x)$ and eliminating the later variables leaves a
degree-$d$ polynomial in $y_i,x$, independent of $\tau_i$.
The one-parameter universal result applies also if some of its
coefficients vanish. Thus every indicated prime power divides
in the universal polynomial ring, hence so does their product.
Specialization allows arbitrary varying coefficients over a DVR.
\end{proof}

\subsection{Complete intersections with base points in degree one}
\label{sec:nov30-ci-fillings}
\begin{theorem}\label{thm:nov30-all-ci-fillings}
In the setting of Theorem~\ref{thm:nov30-finite-cm-model}, the
same conclusion holds if the graded hyperplane section ring of
the K-semistable special fiber is a complete intersection, allowing
a non-Cartier hyperplane root and base points. Thus any distinct
K-polystable replacement over a smooth degree-\(d\) hypersurface
has a non-complete-intersection graded hyperplane section ring.
\end{theorem}
\begin{proof}
Use \(q=n+2,\ N=n+1,\ D=q(d-1)^N\) and the degree-one coordinates
\(z\) of the second model. The relative algebra \(S_R\) is defined
by the reflexive powers of the hyperplane root. Its section spaces
are free and commute with base change by the reflexive restriction
and vanishing argument in Corollary~\ref{cor:nov30-basepoint-necessary}.
Thus its central Hilbert series is that of a hypersurface. Choose
minimal homogeneous algebra generators lifting those of the central
ring, with additional variables \(y_i\) of degrees \(a_i\ge2\).
As in Corollary~\ref{cor:nov30-all-ci-hyperplane-degenerations},
the central regular sequence has relation degrees
\(\{a_i\}\cup\{d\}\), counted with multiplicities. Lift its
elements to relative equations. They generate the relative ideal
in every degree by Nakayama and flatness of the section spaces.

In each degree \(a\), the scalar coefficient matrix of the
degree-\(a\) variables in the degree-\(a\) equations has full
column rank over the fraction field, since degree one generates
the generic hypersurface ring. Its reduction modulo \(t\) is
zero by minimality. Smith normal form, using invertible changes
of basis among these rows and among the variables of this degree,
therefore puts the relations in the form
\[
   \tau_i y_i=P_i(z,y_{a_j<a_i}),\qquad H(z,y_{a_j<d})=0,
   \qquad j_i=\operatorname{ord}_t\tau_i>0 .
\]
There is exactly one extra row in degree \(d\); after row
reduction it has no term linear in any degree-\(d\) variable.
These changes preserve the central regular sequence. The relative
equations are a regular sequence as well, by flatness and the
regular special fiber sequence. In this presentation the
polynomials \(P_i,H\) need not be monic.

Let \(m_i=\lfloor d/a_i\rfloor\) for \(a_i<d\), and \(m_i=0\)
otherwise. Eliminate all \(y_i\) into \(H\), multiplying by
\(\prod\tau_i^{m_i}\) to obtain a polynomial \(F_*\in R[z]\).
Indeed in the full expansion a denominator \(\tau_i\) arises
at most \(\lfloor d/a_i\rfloor\) times, by ordinary degree
and by eliminating in decreasing degree order. Write
\(F_*=t^\kappa F_B\) with \(F_B\) primitive. Its degree is
\(d\) since it generates the generic hypersurface ideal.
Let \(E_R=R[z]/(F_B)\); it injects into \(S_R\), with identical
section spaces over the fraction field.

Here is the relative determinant calculation without finiteness
over degree one. Put
\[
  \sigma_i=\tau_i\prod_{a_j<a_i}\tau_j^{\lfloor a_i/a_j\rfloor},
  \quad U_i=\sigma_i y_i,\quad M=\sum_i m_i j_i .
\]
Present \(E_R\) by the variables \(z,U_i\) and the equations
\[
 U_i-(\sigma_i/\tau_i)P_i(z,U_j/\sigma_j)=0,\qquad F_B(z)=0 .
\]
The displayed polynomials have integral coefficients. In a
monomial of \(P_i\), the power of any particular \(\tau_j\) in
the substituted denominators is at most \(\lfloor a_i/a_j\rfloor\)
by addition of section degrees. Eliminating the \(U_i\) therefore
recovers precisely \(E_R\); the graph and hypersurface equations
are a regular sequence. Over the fraction field the map in
these two polynomial presentations scales \(U_i\) by \(\sigma_i\)
and its graph relation by \(\sigma_i/\tau_i\). Modulo graph
relations the last row scales by \(t^{-\kappa}\prod_i\tau_i^{m_i}\).
Extra off-diagonal multiples of graph equations in this last
row have determinant one. More explicitly the corresponding
change of basis on each graded Koszul term is triangular after
ordering its exterior basis by the presence of the last row.

For clarity write \(h(z)=(1-z^d)/(1-z)^q\). Taking relative
determinant degrees on the free graded Koszul complexes gives
the following series for the lengths of the lattice inclusions:
\[
 \sum_{k\ge0}\operatorname{length}_R
 \bigl((S_R)_k/(E_R)_k\bigr)z^k
 =h(z)\left(\sum_i\frac{j_i z^{a_i}}{1-z^{a_i}}
               -\frac{(M-\kappa)z^d}{1-z^d}\right).
\]
Indeed scale each polynomial generator and each exterior equation
basis by its indicated power of \(t\). The sum of exponents
in a basis of the quotient is the alternating sum on the graded
Koszul terms, by exactness over \(R\) for both presentations.
For a polynomial variable of degree \(a\) the sum of monomial
exponents contributes \(z^a/(1-z^a)\) times the ordinary series;
a relation of the same degree contributes its negative. Subtracting
\(\operatorname{ord}_t(\sigma_i/\tau_i)\) from
\(\operatorname{ord}_t\sigma_i\) gives \(j_i\).

Combining with the embedded determinant comparison of
Theorem~\ref{thm:nov30-finite-cm-model}, and expanding at one
as in \eqref{eq:nov30-triangular-futaki}, yields the following
value of \(N\) times the normalized relative CM degree of the
second filling minus the first:
\[
 \frac{q(d-1)}{dD}\operatorname{ord}_t\Delta_d(F_*)
    +\sum_i\left(N-\frac r{a_i}-\frac{q(d-1)m_i}{d}\right)j_i .
\]
The exponent \(\kappa\) cancels because the discriminant is
homogeneous of degree \(D\). Divisible Veronese section degrees
compute the CM line also for the Weil root.

For each \(a_i<d\), apply
Lemma~\ref{lem:nov30-nonmonic-discriminant}, with \(M_i=M(a_i)\)
in its notation. Graph variables of degree at least \(d\) can be
omitted in this elimination since none enters \(H\) or a
lower-degree row. The contribution of an index \(a_i<d\) to
this comparison is at least
\[
 j_i\left(N-\frac r{a_i}
        -\frac{M_i}{d(d-1)^{N-1}}\right)>0 .
\]
Indeed with \(a=a_i\), the positive geometric sum in the
formula for \(M_i\) gives
\[
 \frac{M_i}{d(d-1)^{N-1}}
 \ \le\ \frac{(d-1)^2}{a(d-a)}
 \ \le\ d-1
 \ <\ N-r/a ,
\]
using \(a(d-a)\ge d-1\) for \(1\le a<d\). For \(a_i\ge d\)
the positive term \(j_i(N-r/a_i)\) remains. With any
additional generator the comparison is therefore strictly
positive, contrary to CM minimization. In the absence of
additional generators the embedded GIT argument of the same
theorem identifies the special fibers.
\end{proof}

\section{Possible non-Cartier loci on a replacement}

\begin{proposition}\label{prop:nov30-quartic-nine-cartier-codim-two}
Let $Z$ be a K-polystable $\mathbf Q$-Gorenstein smoothable quartic
ninefold limit with integral ample $\mathbf Q$-Cartier hyperplane
root $L$, $-K_Z\sim7L$, $L^9=4$, and $h^0(Z,L)=11$. Then $L$ is
Cartier away from a set of codimension at least three. Thus a
distinct K-polystable replacement above a smooth quartic ninefold
requires non-Cartier points in codimension at least three.
\end{proposition}
\begin{proof}
Write $s$ for the dimension of the complete hyperplane image in
$\mathbf P^{10}$. This image is nondegenerate. If $s=1$, at a
general smooth point of the image its order-nine jet conditions
impose at most ten linear conditions on the eleven hyperplanes.
A nonzero hyperplane section therefore vanishes there to order
at least ten. At a general fiber divisor of the rational map a
local parameter of the image pulls back with positive order, so
$T_L(E)\ge10$. This contradicts $1=A_Z(E)\ge 7 S_L(E)\ge
7 T_L(E)/10$. A zero-dimensional image is impossible by linear
independence. Consequently $s\ge2$.

On a smooth model resolving $|L|$ and any given prime valuation
$v=\operatorname{ord}_E$, put $P=\pi^*L$ and let $J$ be the
moving hyperplane class. Use products in a Cartier multiple to
define this moving part also when $L$ is Weil. Write
$a=\min_{0\ne u\in H^0(L)}v(\operatorname{div}u)$. By construction
$P-aE-J$ is effective up to $\mathbf Q$-linear equivalence. A
general projection of the $s$-dimensional image to $\mathbf P^s$
is finite and pulls the hyperplane back to its hyperplane class.
Apply \eqref{global:eq:optimizer-covering-slope} with $t=1$ to
this dominant map to $\mathbf P^s$. Since $-K_{\mathbf P^s}=
(s+1)\mathcal O(1)$ it gives
\[
     P^8J/P^9\ \ge\ 7s/(9(s+1))\ \ge\ 14/27.
\]
When $a>0$, the volume of $P-taE$ for $0\le t\le1$ is at least
that of the nef class $(1-t)P+tJ$. Integrate its first two binomial
terms in $t$, obtaining
\[
   A_Z(v)\ \ge\ 7 S_L(v)\ \ge\ \frac7{10}
         \left(1+\frac{14}{27}\right)a>a.
\tag{QC2}
\]

At a general point of a codimension-two singular stratum, a klt
variety is analytically the product of a smooth germ and a quotient
surface $\mathbf A^2/G$ with $G$ finite and free off zero. Here
one uses the usual quotient charts for klt singularities in
codimension two, restricted away from smaller strata, and finite
group linearization at a general point of the stratum. The
normalized-volume finite-degree formula on its smooth cover and
the local-to-global K-semistable bound give
\[
    \frac{9^9}{|G|}\ \ge\ 9^9\,4(7/10)^9,
    \qquad |G|\le6.
\]
A general transverse surface germ there admits a $\mathbf Q$-Gorenstein
smoothing: cut the smoothing family by seven sufficiently general
local hyperplanes lifting transverse parameters on the special
stratum. They form a relative regular sequence by flatness and
Cohen--Macaulayness, preserve the relative $\mathbf Q$-Gorenstein
property by adjunction (equivalently cut also on the index-one
cover), and the nearby surface is smooth by Bertini. Such a
quotient surface is either Du Val or cyclic of type
$\tfrac1{b j^2}(1,bjc-1)$ with $j\ge2$, $b\ge1$, $(j,c)=1$,
by J.~Koll\'ar and N.~I.~Shepherd-Barron, \emph{Threefolds and
deformations of surface singularities}, Invent. Math. \textbf{91}
(1988), 299--338, \href{https://doi.org/10.1007/BF01389370}{Section 3}.

In the Du Val case the canonical class is locally Cartier, hence
$7L$ is locally trivial. The quotient local divisor class group
is killed by $|G|$, by pulling to its smooth cover and taking the
norm. As $|G|\le6$, it has no $7$-torsion and $L$ is Cartier.
In the non-Du Val case $bj^2\le6$ forces $j=2$, $b=1$, so the
quotient is $\tfrac14(1,1)$. Its local divisor class group is
$\mathbf Z/4$, and the anticanonical character has order two.
The relation $7L=-K_Z$ uniquely makes $L$ the class of character
two. Trivialize its pullback on the smooth cover; its sections
are represented by regular functions transforming by that
character. Every such function vanishes to at least order two in
the two transverse coordinates, because the group acts scalarly
on them and trivially along the stratum. The quotient of the
ordinary transverse blowup gives a prime valuation with log
discrepancy $2/4$ and order at least $2/4$ on every section
divisor: the cover blowup has discrepancy two and ramification
index four along its exceptional divisor. This valuation may
equally be taken over the generic point of the stratum by the
exceptional divisor of the normalized blowup of the stratum;
on the cover its ideal is the fourth power of the transverse
maximal ideal and the ordinary blowup gives the same divisor.
Thus $A_Z(v)\le a$, contradicting (QC2). This leaves Cartier
roots at the generic point of every codimension-two stratum.
The final assertion uses the Cartier-root basepoint-free result
and the relative CM comparison for finite hyperplane maps.
\end{proof}

\begin{proposition}\label{prop:nov30-general-cartier-codimension-two}
Let $Z$ be a K-polystable $\mathbf Q$-Gorenstein smoothable
degree-$d$ Fano hypersurface limit, $3\le d\le n+1$. Suppose
the integral ample $\mathbf Q$-Cartier hyperplane root $L$ on
$Z$ has $h^0(Z,L)=n+2$, $-K_Z\sim rL$, $r=n+2-d$, and comes
from a $\mathbf Q$-Cartier extension of the hyperplane class in
the smoothing. If
\[
                    2r^2>3n(d-1),
\]
then $L$ is Cartier in codimension two. This applies to
quartics in every dimension $n\ge9$.
\end{proposition}
\begin{proof}
Put $N=n+1$ and write $s$ for the dimension of the complete
hyperplane image. It is nondegenerate. If $s=1$, the $n+2$
sections contain a nonzero hyperplane vanishing to order at
least $n+1$ at a general smooth point of the image, by imposing
vanishing of its first $n+1$ Taylor coefficients. Pull it back
at a general fiber divisor $E$. Then $T_L(E)\ge n+1$ and
$1\ge rS_L(E)\ge rT_L(E)/N\ge r$, contrary to the hypothesis,
which implies $r>1$. Thus $s\ge2$.

Given a divisorial valuation $v=\operatorname{ord}_E$, set
$a=v(|L|)$, the minimum of its orders on effective section
divisors. On a common resolution $P$ denotes the pullback of
$L$ and $J$ the moving hyperplane class, so $P-aE-J$ is
effective up to $\mathbf Q$-linear equivalence. This follows
using hyperplane products in a Cartier multiple of $L$. A
general projection of the image to $\mathbf P^s$ pulls back its
hyperplane class to $J$. Applying the tangent-slope inequality
\eqref{global:eq:optimizer-covering-slope} at $t=1$ to the
resulting dominant map gives $P^{n-1}J/P^n\ge sr/(n(s+1))$.
For $0\le t\le1$, $P-taE$ dominates the nef class
$(1-t)P+tJ$. Integration of its first two binomial terms gives
\[
    A_Z(v)\ge rS_L(v)\ge
    \frac{r}{N}\left(1+\frac{sr}{n(s+1)}\right)a
    \ge \frac{r}{N}\left(1+\frac{2r}{3n}\right)a.
\]
The final coefficient exceeds one by the displayed hypothesis.
Thus the linear system $|L|$ has log canonical threshold above
one, with its base orders interpreted using products in a
Cartier multiple. This also holds analytically in every local
quotient chart, by restriction of a log resolution.

Take a general point of a codimension-two singular stratum. The
usual quotient charts for klt varieties in codimension two
identify its germ analytically with a quotient surface times a
smooth germ: the finite group linearizes at a general point,
fixes the stratum directions and acts freely off zero in the
two transverse directions. The transverse quotient surface $U$
admits a $\mathbf Q$-Gorenstein smoothing: cut the given family
by general local hyperplanes transverse to the stratum. These
form a relative regular sequence, the nearby surface is smooth
by Bertini, and adjunction on the index-one cover gives the
induced $\mathbf Q$-Gorenstein property. Consequently $U$ is
either Du Val or of type $\tfrac1{bj^2}(1,bjc-1)$, where
$b\ge1$, $j\ge2$ and $(j,c)=1$. This is the classification of
$T$-singularities in J.~Koll\'ar and N.~I.~Shepherd-Barron,
\emph{Threefolds and deformations of surface singularities},
Invent. Math. \textbf{91} (1988), 299--338,
\href{https://doi.org/10.1007/BF01389370}{Section 3}.

In both cases the local index-one cover of $U$ is Du Val. The
associated cover of its transverse smoothing is a smoothing of
this Du Val surface, by the local index-one-cover description
of $\mathbf Q$-Gorenstein deformations. Its three-dimensional
total germ is smooth or an isolated hypersurface singularity:
the central Du Val local ring is a hypersurface, a lift modulo
the nonzerodivisor $t$ adds at most one to its embedding
dimension, and the total space is smooth off its central
singular point. The punctured total germ has torsion-free Picard
group. For example the link of an isolated three-dimensional
hypersurface singularity is simply connected by the connectivity
theorem for hypersurface links, so any cyclic cover attached
to a torsion line bundle on the punctured germ splits. The
smooth case follows in the same way. This is also the local
Grothendieck--Lefschetz argument used in Y.~Liu and C.~Xu,
\emph{K-stability of cubic threefolds}, proof of Lemma~3.16.

The given $\mathbf Q$-Cartier extension of $L$ is a line bundle
away from the singular central point on the transverse slice,
with a power trivial locally. Its pullback to the punctured
cover total space is therefore trivial. Restricting to the
punctured central cover shows that the pullback of the transverse
divisor class of $L$ is trivial there. This proves the claim on
a Du Val stratum, where the index-one cover is the identity.

In the cyclic case the subgroup defining the canonical cover
is the determinant kernel $\mu_{bj}\subset\mu_{bj^2}$: the
determinant has character $bjc$. Thus the character of $L$ is
divisible by $bj$. If it is nonzero, it is neither the trivial
character nor either coordinate character, since $1$ and $bjc-1$
are units modulo $bj$. On the smooth uniformizing cover every
local section of $L$ is a semi-invariant function; it has no
constant term or linear terms in the two transverse coordinates.
The transverse ordinary blowup upstairs has log discrepancy
two and gives every section divisor order at least two. This
contradicts the threshold above one: the uniformizing map is
quasi-\'etale, so it preserves the log canonical threshold for
the pulled-back linear system by the finite-cover discrepancy
formula. Equivalently both the discrepancy and the orders
are divided by the same divisorial ramification index upon
descending this blowup valuation. Hence $L$ is Cartier at the
generic point of every codimension-two stratum.
\end{proof}

\begin{proposition}\label{prop:nov30-quartic-nine-index-cover}
In the quartic ninefold setting of
Proposition~\ref{prop:nov30-general-cartier-codimension-two},
the hyperplane root is Cartier away from a locus of codimension
at least four.
\end{proposition}
\begin{proof}
Suppose the non-Cartier locus has a component of codimension
three. The local index $j$ of $L$ at a general complex point is at
most six, by the finite-degree formula for its cyclic index
cover, the smooth-point upper bound for normalized volume, and
the local-to-global lower bound $9^9\,4(7/10)^9$. As $7$ is
prime to $j$, the relation $-K_Z\sim7L$ makes its cover the
canonical index-one cover. Shrink the codimension-three stratum
away from its smaller strata. On a general transverse threefold
slice its index cover has group $\mu_j$, with action free away
from the origin because $L$ is Cartier there. Cutting the
original smoothing transversely gives a $\mathbf Q$-Gorenstein
smoothing of the slice, as in the preceding proposition. The
induced deformation of its canonical cover has free group action
on every smooth nearby fiber, since the cover is \'etale over
the Gorenstein locus of the quotient family.

The global hyperplane-system log canonical threshold exceeds
one by the proof of
Proposition~\ref{prop:nov30-general-cartier-codimension-two}.
This remains true on the quasi-\'etale canonical cover near
the given stratum. The pullback of $L$ there is Cartier; in
its character trivialization the hyperplane sections are all
semi-invariants of nontrivial character and vanish on the
lifted stratum. Their ideal is therefore contained in that
stratum ideal, whose log canonical threshold also exceeds one.
Intersect the codimension-three stratum with six general
ambient hyperplanes meeting it transversely at general points.
The same strict threshold holds on the transverse threefold
cover $U$ for the restricted stratum ideal, which is the
maximal ideal at the point. Indeed on a simultaneous log
resolution these general cuts are smooth, meet the fixed
normal-crossings strata transversely, and preserve both the
orders and the discrepancies of all the divisors that meet
the resulting slice by adjunction. The restricted pulled-back
ideal is thus resolved by the sliced resolution. In particular
$U$ is klt and Gorenstein by Bertini and adjunction.

A general local section from the maximal ideal of $U$ through
the point consequently gives a plt pair in a neighborhood of
that point: the fixed exceptional orders on the resolution
are strictly below the log discrepancies, and the general free
strict transform is smooth transverse to the exceptional divisors.
This is also the local argument of Pan and Wu
\cite[Lemma~4.5]{kePanWuCubics}. The resulting surface is klt by plt adjunction and Gorenstein
by the Cartier adjunction formula. It is smooth or Du Val, of
embedding dimension at most three. Lifting its local generators
and adjoining the cut parameter gives embedding dimension at
most four for $U$. Therefore $U$ is smooth or a hypersurface.

If the cover germ were smooth, the equivariant submersion to the
base would retain a fixed point on every nearby fiber (choose
linearized local coordinates including the invariant parameter).
Suppose the transverse cover is a singular hypersurface. Lift an
eigenbasis of its four-dimensional cotangent space to the total
space of the smoothing. By the local embedding criterion the
total family is a hypersurface $F(x,t)=0$ in
$\mathbf C^4\times\mathbf C$, with linear diagonal action on
$x$ and trivial action on $t$. The equation can be chosen
semi-invariant by projection to the equation character followed
by Nakayama. It is invariant: otherwise $F(0,t)$ would vanish
identically, giving fixed points on the nearby fibers. Moreover
every nonidentity group element has zero fixed space in
$\mathbf C^4$. A fixed space of dimension at least two would
meet $F(x,0)=0$ in positive dimension. On a one-dimensional fixed
space, its restriction at $t=0$ has a zero of positive finite
order; the Weierstrass preparation theorem preserves a zero on
nearby fibers. Both alternatives contradict freeness. Hence the
four variable characters are all units modulo $j$.

The canonical character on this hypersurface is the determinant
character of the four variables, up to sign, by the residue
formula and invariance of the equation. It must have order $j$.
For even $j$ the sum of four units is even, a contradiction.
For odd $j\ge3$, the nondegenerate quotient of the transverse
quadratic part is a self-dual representation with no trivial
character. Its rank is even because the nontrivial characters
come in distinct inverse pairs. Rank four would give determinant
one and hence trivial canonical character. Its rank is therefore
at most two at a general point of the stratum.

Return to the full nine-dimensional canonical cover at such a
point. Six parameters $u$ along the stratum, together with four
generators of the transverse cotangent space, give a smooth
ten-dimensional ambient hypersurface chart. Indeed lifting the
generators adds at most six to the transverse embedding dimension,
and the nine-dimensional cover is Cohen--Macaulay, with the
transverse parameters a regular sequence. On a sufficiently small
open piece of the stratum the lifted stratum is smooth and the
parameters $u$ are coordinates there; the reduced lift maps
locally isomorphically after shrinking since the stabilizer is
the full local index group. Choose transverse coordinates $x$
vanishing on this lifted stratum. The local equation belongs to
$(x)^2$, since the transverse cover is singular at general
points throughout it. After shrinking to the constant rank locus,
the quadratic matrix has an analytic kernel subbundle of rank at
least two. A change of $x$ depending on $u$ therefore makes the
entire transverse quadratic part a form in $x_1,x_2$. Give
$x_1,x_2$ weights nine, $x_3,x_4$ weights six and the six base
parameters weight four. The minimum monomial weight $b$ of the
equation is at least eighteen.

For a klt hypersurface germ of dimension $n$ with positive weights
$w$ on smooth ambient parameters, the maximal-ideal normalized
volume admits the bound
\[
 \widehat{\operatorname{vol}}\ \le\
      \frac{b(\sum w_i-b)^n}{\prod w_i},
       \qquad b=\operatorname{wt}_{\min}(F).
\]
Indeed the induced weighted-order ideals $I_k$ on the
hypersurface have asymptotic multiplicity
$\lim e(I_k)/k^n=b/\prod w_i$, from the Hilbert series of the
ambient weighted polynomial ring modulo the initial form of $F$.
For instance $I_{p(k+\max w_i)}\subset I_k^p\subset I_{pk}$
by splitting ambient monomials into pieces. Weighted blowup
and inversion of adjunction give
$\operatorname{lct}(I_k)\le(\sum w_i-b)/k$; the ideal formula
for normalized volume proves the bound without an irreducibility
assumption on the initial form. In the current weights $\sum w_i=54$.
For $b\ge18$ with $b<54$ the expression $b(54-b)^9$ decreases.
Consequently
\[
  \frac{\widehat{\operatorname{vol}}(\text{cover})}{9^9}
    \le \frac{18\,4^9}{9^2 6^2 4^6}=\frac{32}{81}.
\]
Finite degree and K-semistability instead give a lower bound
$j\cdot4(7/10)^9\ge12(7/10)^9>32/81$. This contradiction
excludes the hypersurface cover.
\end{proof}

\begin{corollary}\label{cor:nov30-quartic-full-image-index}
In the setting of Proposition~\ref{prop:nov30-quartic-nine-index-cover},
suppose the hyperplane image has dimension nine. Then the
non-Cartier locus of $L$ has codimension at least seven. If it
has a codimension-seven component, its general local index is
three. On the transverse hypersurface cover the characters in
eight ambient coordinates are six copies of two and two copies
of one, after choosing the generator for which sections of $L$
have character two.
\end{corollary}
\begin{proof}
Use the notation $P,J,a,v$ of the preceding propositions.
For the mixed-intersection threshold estimate in the cubic setting,
see Kong--Shen--Zhao--Zhao \cite[Proposition~3.2]{keKongCubics}
and Pan--Wu \cite[Lemma~5.3]{kePanWuCubics}.
The full-dimensional nondegenerate image in $\mathbf P^{10}$
has degree at least two, whence $J^9\ge2$. Since $P^9=4$,
log-concavity of mixed nef intersections gives
$P^{9-i}J^i/P^9\ge2^{-i/9}$. For $0\le t\le1$ use
$P-taE\ge(1-t)P+tJ$ and integrate every term of the nef
binomial expansion. K-semistability yields
\[
    A_Z(v)\ge7S_L(v)\ge \frac7{10}
       \left(\sum_{i=0}^{9}2^{-i/9}\right)a>5a.
\]
For example $2^{-1/9}>37/40$, and
$(7/10)\sum_{i=0}^9(37/40)^i>5$ by direct evaluation.

At a general non-Cartier component of codimension $c$ the same
argument as in Proposition~\ref{prop:nov30-quartic-nine-index-cover}
transfers this strict threshold to the maximal ideal of the
transverse Gorenstein index cover $U$. If $3\le c\le7$, we
may therefore cut through the marked point with $c-2$ successive
general elements of its maximal ideal to reach a Du Val or smooth
surface. At each cut, the exceptional log discrepancies on a
resolution of the maximal ideal decrease by its fixed orders,
the residual free divisor meets the resolution transversely,
and plt adjunction gives a Gorenstein klt section. The strict
initial threshold above $c-2$ permits all the cuts. It follows
that $U$ has embedding dimension at most $3+(c-2)=c+1$. A
smooth transverse cover is excluded by the fixed-point argument
for equivariant smoothings. For a singular hypersurface of
multiplicity $b\ge2$, the ordinary ambient blowup and inversion
of adjunction give $\operatorname{lct}(U;\mathfrak m)\le c+1-b$.
This is at most five when $c\le6$, proving the first assertion.

For $c=7$ the transverse cover must be a hypersurface. As
before its smoothing forces an invariant equation and primitive
characters on all eight ambient variables. The index is at most
six, and an even index is impossible since the sum of eight odd
characters cannot have full even order. For index five the
normalized-volume upper bound on the full index cover is at
most $9^9\,2(6/7)^7$. Indeed choose transverse parameters of
weight seven and the two stratum parameters of weight six.
The transverse hypersurface equation has weight at least
fourteen since it vanishes quadratically along the lifted
stratum. The same weighted-ideal bound as in the preceding
proof is maximal at weight fourteen, giving exactly the stated
ratio. This contradicts the lower bound $9^9\,20(7/10)^9$.

It remains that the index is three. Choose the group generator
so lifted sections of $L$ transform by character two. Let $h$
be the number of transverse weight-two variables among the eight.
The character of anticanonical sections is the determinant
character of the variables, since the equation is invariant.
Thus $8+h\equiv7\cdot2\equiv2\pmod3$ and $h=0,3,6$.
Use ambient positive weights one and two according to these
characters. The invariant equation through the origin has weight
at least three, and all lifted sections of $L$ have weight at
least two. On the transverse slice the section ideal still has
threshold above five by restriction of the same log resolution.
The weighted ambient blowup bounds that threshold by
$(8+h-3)/2$ via inversion of adjunction. This is at most four
when $h=0$ or $3$, leaving only $h=6$.
\end{proof}

\begin{corollary}\label{cor:nov30-general-threefold-index}
Under the hypotheses of Proposition~\ref{prop:nov30-general-cartier-codimension-two},
including \(2r^2>3n(d-1)\), put \(\beta=d(r/(n+1))^n\).
At the general point of any
codimension-three component of the non-Cartier locus the local
index \(j>1\) satisfies
\[
 j\beta\le16/27,\qquad 2\mid j\ \Longrightarrow\ 2\mid r,
 \qquad j\ge3,\ j\nmid r\ \Longrightarrow\ j\beta\le32/81 .
\]
For quartics in odd dimension \(n\ge9\), any such component
would have \(j=3\mid r\) and the lifted transverse section ideal
would have log canonical threshold at most \(3/2\).
\end{corollary}
\begin{proof}
Use the cyclic cover of \(L\) in the transverse argument of
Proposition~\ref{prop:nov30-quartic-nine-index-cover}. To see
that the cover extends in the smoothing without a coprimality
condition, let \(\mathcal L\) be the relative hyperplane class.
Its reflexive powers are Cohen--Macaulay and restrict to the
corresponding powers of \(L\), as in the reflexive restriction
argument of Proposition~\ref{prop:nov30-factorial-quartic-replacement}.
At a general point of the chosen component, a local generator of
\(\mathcal O_Z(jL)\) lifts to a generator of the torsion-free
rank-one module \(\mathcal O(j\mathcal L)\), so the latter is
free there. Choose its trivialization and form the usual cyclic
algebra \(\bigoplus_{a=0}^{j-1}\mathcal O(-a\mathcal L)\).
It restricts to the central index cover of \(L\) and is a
\(\mu_j\)-torsor on the smooth generic fiber since the generic
hyperplane class is Cartier. The central cover is Gorenstein
by the integral linear equivalence \(-K_Z\sim rL\). Its
canonical character can have order smaller than \(j\), whereas
the character of \(L\) is faithful. Localizing the family over
the geometric generic point of the stratum in a separating
parameter space, or slicing a log resolution with general
transverse parameters, preserves the lower threshold bound
for the lifted section ideal. The strict threshold
therefore still gives a cDV transverse
hypersurface, its equivariant smoothing has invariant equation,
and all four transverse coordinate characters are units modulo
\(j\). The determinant character represents the canonical class
up to sign. For even \(j\) the sum of four units is even, so
\(r\) must be even. For \(j\ge3\) the primitive characters
are non-self-dual. Their invariant quadratic form has even rank
and has rank four only if its four characters are two opposite
pairs, in which case \(j\mid r\).

To apply the closed-point volume bound allow \(n-3\) parameters
along the stratum, of weight \(2/3\) each. Giving the four transverse
variables weight one gives equation weight at least two. The weighted
hypersurface bound in Proposition~\ref{prop:nov30-quartic-nine-index-cover}
then gives \(16/27\) as an upper bound for the normalized density
(the expression is decreasing from equation weight two).
See also Y.~Liu,
\emph{Optimal bounds for local volumes of threefold singularities},
Proposition~4.1, \href{https://arxiv.org/abs/2512.05429}
{arXiv:2512.05429}. For rank at most two, use parameters in
which the transverse quadratic part is supported on two variables,
as in the previous proof. Give those two variables
weight \(3/2\) and the remaining two weight one. The equation
has weight at least three. At weight three the bound on density
is \(3(2/3)^3/(3/2)^2=32/81\), and the expression decreases
thereafter. The density is at least \(j\beta\)
by the finite-degree and local-to-global formulas.

For quartics the sequence \(4((n-2)/(n+1))^n\) is increasing
for \(n\ge9\). Indeed for real \(x>2\) the logarithmic derivative
of \(((x-2)/(x+1))^x\) equals
\(\log(1-u)+u(1-u/3)/(1-u)\), \(u=3/(x+1)\); the coefficient
of \(u^k\) in this expression is \(2/3-1/k>0\) for \(k\ge2\).
At \(n=9\) already \(4\beta>16/27\) and \(3\beta>32/81\).
For odd \(n\), the only remaining possibility is \(j=3\mid r\)
and a full-rank quadratic. Then the representation has two
variables of each nontrivial character. Assign weight two to the
two coordinates of the section character of \(L\), and weight one
to the others. All lifted sections have weight at least two.
After equivariant linear changes, the weighted initial equation is
\(x_1y_1+x_2y_2+g_3(y_1,y_2)\): nondegeneracy gives the two
independent bilinear terms. This initial equation is irreducible
since the two linear coefficients \(y_1,y_2\) are coprime.
The weighted exceptional divisor therefore gives a valuation on
the cover of log discrepancy \(2+2+1+1-3=3\), by adjunction.
Its value on the section ideal is at least two, so the transverse
threshold is at most \(3/2\).
\end{proof}

\begin{corollary}\label{cor:q9-all-odd-transverse-cartier}
The limiting hyperplane class of a K-polystable replacement of smooth
quartic $n$-folds is Cartier through codimension three for every odd $n\ge9$.
\end{corollary}
\begin{proof}
Corollary~\ref{cor:nov30-general-threefold-index} leaves only the
possibility $3\mid n-2$, and thus $n\ge11$ in the stated odd range.
We sharpen the threshold estimate in
Proposition~\ref{prop:nov30-general-cartier-codimension-two}.
The strict bound above one excludes every fixed prime divisor of
$|L|$: such a divisor would have order at least one and log
discrepancy one. On a resolution $\pi$ of the hyperplane map, let
$P=\pi^*L$ and let $J$ be the basepoint-free moving hyperplane
class. Pushing down $J$ gives the class of $L$, since every fixed
component upstairs is exceptional. Consequently
$P^{n-1}J=P^n$. Moreover $P^{n-2}J^2>0$. Indeed the image
has dimension at least two, so two general image hyperplanes
cut a codimension-two cycle with a component meeting the open
set where $\pi$ is an isomorphism: choose the hyperplanes
generally so their intersection meets the open image of the
locus where the hyperplane map is smooth onto its image.
Its nonzero effective pushforward has positive degree against
the ample class $L$.

For $a=v(|L|)>0$, retain the first three terms in the integral
of $((1-t)P+tJ)^n$, exactly as in the proof of
Proposition~\ref{prop:nov30-general-cartier-codimension-two}.
Each binomial coefficient times the corresponding integral
$\int_0^1t^i(1-t)^{n-i}\,dt$ equals $1/(n+1)$.
K-semistability gives
\[
  \frac{A_Z(v)}{a}\ge
  \frac{n-2}{n+1}\left(2+\frac{P^{n-2}J^2}{P^n}\right)
   >\frac{2(n-2)}{n+1}.
\]
The mixed intersection is unchanged upon further resolution,
so the strict improvement is uniform in $v$.
The residual index-three transverse model in
Corollary~\ref{cor:nov30-general-threefold-index} has threshold
at most $3/2$ for the lifted system. The displayed global
lower bound passes to the transverse cover as in that corollary.
For every $n\ge11$ we have $2(n-2)/(n+1)\ge3/2$.
The strict improvement contradicts the transverse upper bound.
\end{proof}

\subsection{Two additional generators and a quadratic relation}

\begin{proposition}\label{prop:nov30-quartic-two-pfaffian}
Let $(Z,L)$ have the polarized section Hilbert function of a
degree-$d$ Fano hypersurface of dimension $n$, $3\le d\le n+1$,
with $Z$ K-semistable klt, $L$ integral ample $\mathbf Q$-Cartier,
and $-K_Z\sim rL$, $r=n+2-d$. If the full section ring $S$ has
exactly two additional minimal homogeneous generators beyond its
degree-one part and the smaller generator degree is two, then
$S$ is a complete intersection. Consequently these rings cannot
occur as distinct K-polystable replacements above smooth
degree-$d$ hypersurfaces. In particular a distinct K-polystable
replacement over a smooth quartic ninefold with full-dimensional
hyperplane image requires at least three additional generators.
\end{proposition}
\begin{proof}
The normal section ring $S$ is Gorenstein with canonical module
$S(-r)$. We recall the ingredients when $L$ is Weil. The
intermediate cohomology of $\mathcal O_Z(kL)$ vanishes in all
degrees: the Weil form of Kawamata--Viehweg applies for
$k\ge0$, and Serre duality and vanishing for $K_Z-kL$ apply
for $k<0$. The rank-one sheaves are Cohen--Macaulay, as follows
locally by the finite cyclic index cover and the character
summands on its klt cover. The local cohomology description of
the graded section ring gives the vertex depth. Off the vertex
it is locally the finite cyclic quotient of a line-bundle torsor
over that index cover, hence is Cohen--Macaulay there also.
The canonical-module formula for a normal section ring of an
integral Weil divisor gives $\bigoplus_k H^0(Z,K_Z+kL)=S(-r)$;
equivalently compute on the line-bundle torsor over the Cartier
locus and extend reflexively over codimension two. The cone
projection has one-dimensional fibers off its vertex, so that
complement still has codimension at least two.

There are exactly two additional minimal variables. Write $y,z$ for
them, of degrees $2,b$ with $b\ge2$. As $d\ge3$, comparison of the
degree-two Hilbert functions supplies at least one quadratic
relation, also when both variables have degree two. Every such
relation lies in the $x$ variables by minimality. Every nonzero
one is irreducible, since the central ideal is prime and contains
no linear relation. The minimal presentation over
$P=\mathbf C[x_0,\ldots,x_{n+1},y,z]$ has a prime height-three
Gorenstein ideal. Suppose it is not a complete intersection. By
D.~Buchsbaum and D.~Eisenbud, \emph{Algebra structures for finite
free resolutions, and some structure theorems for ideals of
codimension 3}, Amer. J. Math. \textbf{99} (1977), 447--485,
\href{https://doi.org/10.2307/2373926}{Theorem 2.1}, its minimal relations
are the submaximal Pfaffians of a minimal odd alternating
matrix. The matrix has size five: a Pfaffian from a matrix
of size at least seven has terms involving at least three
positive-degree entries, and the ideal contains a minimal
quadratic relation. The last free-module shift in the resolution
is $s=n+2+2+b-r=b+d+2$. If the five generator degrees are $g_i$,
the syzygy degrees are $s-g_i$ and the $(i,j)$ matrix entry has
degree $s-g_i-g_j$. A degree-zero or negative entry is zero by
minimality.

The Hilbert numerator over $P$ is
$(1-t^2)(1-t^b)(1-t^d)$. Comparing it with the self-dual
resolution shows that the five degrees form the multiset
\[
                       2,\ b,\ d,\ k,\ s-k .
\]
Indeed the excess of the multiplicity of any degree over its
complementary degree is the same as in the triple $2,b,d$.
Two members of that triple cannot be complementary since the
remaining degree is positive. If one is $s/2$ then its
multiplicity in the five-element list is odd: the parity from
all nonmidpoint pairs is fixed by their differences. Thus the
triple occurs with multiplicities, leaving one complementary
pair (possibly two midpoint entries).

The Pfaffian of degree two is a nonzero irreducible quadratic
form in the $x_i$ by the preceding prime-ideal argument. Every nonzero matching
in its four remaining rows must pair them into two entries of
degree one. Their corresponding degrees must each sum to
$s-1$. Pairing $b$ with $d$ fails since $b+d=s-2$. A single
remaining matching would give a product of two linear forms.
Thus both other matchings are required, forcing
$b+k=b+(s-k)=s-1$ and $d+k=d+(s-k)=s-1$.
Consequently $b=d$, $k=d+1$, and the five relation degrees are
$2,d,d,d+1,d+1$.

In that order write the upper triangular entries of the matrix
as $(A,B,C,D)$, $(E,f,g)$, $(h,k')$, $(0)$. Their degrees are
$(d,d,d-1,d-1)$, $(2,1,1)$, $(1,1)$ and zero. At least one of
$A,B$ involves $z$ with nonzero coefficient. Otherwise $z$ is
a free variable in $S$; the product configuration scaling $z$
has total-weight series $\operatorname{Hilb}_S(t)t^d/(1-t^d)$,
giving $N\operatorname{DF}=N-r/d>0$ for $N=n+1$ in the
normalization of Section~\ref{sec:kd-definitions}. Its inverse has negative invariant,
contradicting K-semistability. Interchange rows and columns two
and three to assume $A$ has that coefficient.

The point where only $z$ is nonzero belongs to $Z$ since none
of the Pfaffians contains a pure power of $z$. On its $z=1$
uniformizing chart the point is the origin, $A$ is a unit, and
the equations are the three Pfaffians containing the pivot
$A$, that is,
\[
       A h-B f+C E=0,\qquad A k'-B g+D E=0,
            \qquad D f-C g=0.
\]
The pivot identity for Pfaffians gives all the others. This is
a complete intersection of dimension $n$ in $n+3$ variables.
The uniformizing chart is normal: adjoining $u$ with $u^d=z$
to $S[z^{-1}]$ is \'etale, and the resulting scheme is the
product of the $z=1$ chart with $\mathbf G_m$, by scaling
coordinates according to their section degrees. The germ of
the chart maps to $Z$ quasi-\'etale with degree $d$. Indeed
the group acts with weight one on the $n+2$ coordinates $x_i$ and weight two on
$y$, freely away from $x=0$, and that fixed-set union has
dimension at most one. This has codimension at least two for
$n\ge3$; when $n=2$ one has $d=3$ and the action is free on
the punctured weight-two $y$-axis as well. The origin is the unique preimage of
the indicated point. K-semistability and the finite-degree
formula therefore require normalized volume of this chart
germ at least $d^2(nr/(n+1))^n$.

Put $\lambda=B(0)/A(0)$. The linear parts of the first two
equations are nonzero scalar multiples of $h-\lambda f$ and
$k'-\lambda g$. These two linear forms are independent, since
the irreducible quadratic Pfaffian is
$f k'-g h=f(k'-\lambda g)-g(h-\lambda f)$; dependence would
make it zero or a product of linear forms. Eliminate two
$x$ variables by the implicit function theorem. Give the
remaining $n$ variables $x$ weight one and give $y$ weight two.
The substituted expressions for the eliminated variables have
weight at least one. Thus $C,D$ retain weights at least $d-1$,
and the residual hypersurface equation has weight $u\ge d$.
The weighted hypersurface volume bound in the proof of
Proposition~\ref{prop:nov30-quartic-nine-index-cover} gives
\[
     \widehat{\operatorname{vol}}
       \ \le\ \frac{u(n+2-u)^n}{2}\ \le\ \frac{d r^n}{2}.
\]
Indeed $u<n+2$ by klt and the expression decreases for
$u\ge d\ge3$. But $d r^n/2<d^2(nr/(n+1))^n$ since
$(1+1/n)^n<\exp(1)<2d$. This proves the complete-intersection
assertion. The relative CM comparison for complete-intersection
fillings gives the claim about replacements.

In the quartic ninefold case a full-dimensional hyperplane
image has degree between two and four. Degree four gives the
whole embedded quartic ring by Hilbert functions. Degree three
is excluded by Proposition~\ref{prop:nov30-linear-image-component}.
For a quadric image the quartic Hilbert function gives precisely
one additional generator in degree two. The one-additional-
generator result and the proposition now exclude fewer than
three additional generators in a distinct replacement.
\end{proof}

\paragraph{A Cartier-index obstruction for a quartic ninefold limit.}
\begin{proposition}\label{prop:nov30-quartic-volume-obstruction}
The local-to-global normalized-volume bound for a K-semistable
quartic ninefold limit, the local relation \(-K\sim7L\), and
\(\mathbf Q\)-Gorenstein smoothability together permit a local
Cartier index of two for \(L\).
\end{proposition}
\begin{proof}
In a K-semistable limit \(Z\) of smooth quartic ninefolds one
has \((-K_Z)^9=4\cdot7^9\). The local-to-global bound
\cite[Theorem~2.14]{ar2007.14320} requires at every closed point
\[
            \widehat{\mathrm{vol}}(p,Z)
                     \ge 9^9\cdot4(7/10)^9 .
\]
At the origin of the quotient germ
\[
       W=(\mathbf A^2/\mu_4(1,1))\times\mathbf A^7
\]
the degree-four quotient from \(\mathbf A^9\) is
quasi-\'etale. The finite-degree formula and the smooth-point
volume \cite[Theorems~2.11(2) and~2.13]{ar2007.14320} give
\(\widehat{\mathrm{vol}}(0,W)=9^9/4\).
This is strictly above the required lower bound, since
\(16\cdot7^9=645657712<10^9\).

The local class group of the two-dimensional cyclic quotient
is its character group \(\mathbf Z/4\), by character descent
from the factorial cover free in codimension one. The determinant
of the action has character two, so the canonical and
anticanonical classes both represent its nonzero element of
order two. Equivalently one can read this off from the action
on a volume form of \(\mathbf A^9\). Put \(L=-K_W\).
Then \(7L\sim -K_W\) in the local class group and \(L\)
has Cartier index two. Taking the quotient by the subgroup
of order two first gives \(xy-z^2=0\) with the remaining
involution acting by simultaneous negation of \(x,y,z\).
The family
\[
       \bigl(\{xy-z^2=t\}/\mu_2\bigr)\times\mathbf A^7
                 \longrightarrow \mathbf A^1_t
\]
is a \(\mathbf Q\)-Gorenstein smoothing: its cover has smooth
total space, the quotient is quasi-\'etale, and the involution
acts freely on every fiber with \(t\ne0\). The indicated
class extends as the relative anticanonical class and is
Cartier on the smooth fibers.
\end{proof}
This example is local. No global K-semistable quartic limit
with this germ is asserted here, and it supplies no counterexample
to K-stability of a smooth quartic ninefold.

\paragraph{An isolated non-Gorenstein germ meeting the volume bound.}
The same numerical limitation persists when the non-Cartier locus
is a single point. In \(\mathbf A^{11}\) choose distinct numbers
\(c_i\), put \(q_1=\sum x_i^2,\ q_2=\sum c_i x_i^2\), and set
\[
 V=(q_1=q_2=0),\qquad W'=V/\{1,-1\},
\]
where \(-1\) negates all eleven coordinates. The projective
intersection \(B\subset\mathbf P^{10}\) is smooth by the
Jacobian criterion. Adjunction gives \(-K_B=7H\) and \(H^8=4\).
The degree valuation on its cone \(V\) has discrepancy seven
and volume four. It minimizes normalized volume since \(B\)
is K-semistable: smooth intersections of two quadrics admit
K\"ahler--Einstein metrics by \cite{SS}, and the cone criterion
is C.~Li and Y.~Liu, \emph{K\"ahler--Einstein metrics and volume
minimization}, Adv. Math. \textbf{341} (2019), 440--492,
\href{https://doi.org/10.1016/j.aim.2018.10.038}{doi:10.1016/j.aim.2018.10.038}.
The sign quotient is free away from the vertex,
so the finite-degree formula gives
\[
 \widehat{\operatorname{vol}}(0,W')=2\cdot7^9
             >4(9\cdot7/10)^9
\]
because \(10^9>2\cdot9^9\). The residue of the ambient
eleven-form by \(dq_1,dq_2\) changes sign. Thus the quotient
has canonical index two, and \(L=-K_{W'}\) has index two
and satisfies \(7L\sim-K_{W'}\) locally.

Choose \(c\) distinct from every \(c_i\) and vary the equations
as \(q_1=t,\ q_2=ct\). The Jacobian criterion gives smooth fibers
for \(t\ne0\), on which the sign action is free. The total cover
is the hypersurface \(q_2-cq_1=0\) after eliminating \(t\);
it is Gorenstein and the quotient is quasi-étale. Taking the
quotient thus gives a \(\mathbf Q\)-Gorenstein smoothing, with the
relative anticanonical class extending \(L\).
This local model tests the volume and Cartier-index obstructions
after exclusion of codimension-two strata. Realization on a
K-polystable replacement over a smooth GIT point would require
additional global information.

\subsection{K-polystable limits with two isolated non-Cartier points}
\label{sec:nov30-two-vertices}
The non-Cartier alternative has genuine projective examples on
the boundary of the hypersurface family. Their GIT limits can be
identified explicitly.
\begin{proposition}\label{prop:nov30-two-vertices}
Let \(d=2a\), \(a\ge2\), and \(n\ge2d-1\).
In \(n+2\) ordinary coordinates \(x\) choose general forms
\(Q_1,Q_2,G\) of degrees \(a,a,d\), respectively. Then
\[
 Z=(Q_1=Q_2=0,\ uv=G(x))
     \ \subset\ \mathbf P(1^{n+2},a,a)
\]
is a K-polystable limit of smooth K-stable degree-\(d\)
hypersurfaces. Its hyperplane root has local Cartier index
\(a\) at two isolated points and is invertible elsewhere. The
family in the proof has reducible GIT limit \([Q_1Q_2]\).
For a quartic ninefold both points are non-Gorenstein of index
two.
\end{proposition}
\begin{proof}
Put \(B=(Q_1=Q_2=0)\subset\mathbf P^{n+1}\),
\(D=(G|_B=0)\), \(r=n+2-d\), and \(N=n+1\).
Write \(H_B=\mathcal O_B(1)\) and \(L=\mathcal O_Z(1)\).
The complete intersections \(B,D\) are smooth Fano of indices
\(r,r-d\), respectively, and general such complete intersections
are K-stable by \cite[Theorem~1.3]{ar2608.09293}.
Over \(B\) the locus of \(Z\) where \(x\ne0\) is locally
\(uv=g\) in two line-bundle coordinates, with \(g\) a local
equation for \(D\). It is smooth, including where \(u=v=g=0\)
by smoothness of \(D\). The charts of the two points with
\(x=0\) are quotients by \(\mu_a\) of the affine cone over
\((B,H_B)\): put \(u=1\) or \(v=1\) on the cover and use
\(uv=G\). The scalar action is free away from the vertex.
The cone is klt since blowing up its vertex gives the line
bundle \(\mathcal O_B(-1)\), with exceptional discrepancy
coefficient \(r-1\). Thus \(Z\) is klt. The degree-one class
has faithful character on each vertex cover and its \(a\)th
multiple is Cartier. A local generator of a smaller Cartier
multiple would give a nonvanishing semi-invariant of nontrivial
character there, which is impossible at the fixed point.
Adjunction gives \(-K_Z\sim rL\). For \(a=2,n=9\) this has
local index two.

Act by a one-dimensional torus with opposite weights on \(u,v\).
We compute its invariant divisorial tests, allowing positive
scalings of orders. Let \(\nu\) be such a valuation and \(w\)
its restriction to \(\mathbf C(B)\). Choose a hyperplane
section \(x_0\) nonvanishing at the center of \(w\) on \(B\)
and put \(U=u/x_0^a\), \(b=\nu(\operatorname{div}_Z x_0)\).
The function field is \(\mathbf C(B)(U)\).
Invariance gives
\[
 \nu\!\left(\sum f_iU^i\right)=\min_i\{w(f_i)+i\nu(U)\}.
\]
Indeed each filtration subspace is torus invariant and the
distinct characters split, as is also seen by a Vandermonde
matrix of translates. The restriction \(w\) is zero or a
multiple of a divisorial valuation: its residue transcendence
degree drops by at most one on restricting across the
one-variable field extension, and its value group is discrete.

We claim, with terms involving the trivial valuation taken as zero,
\[
      A_Z(\nu)=rb+A_B(w),\qquad
      S_L(\nu)= b+\frac{n S_{H_B}(w)+w(D)/(2a)}{N}.
 \tag{V1}
\]
For the first identity take a rational top form \(\eta\) on \(B\)
with divisor \(-r(x_0=0)\). The top form
\(\eta\wedge dU/U\) on \(Z\) has divisor \(-r\operatorname{div}_Z x_0\).
At a generic point over \(D\) one can use the invertible \(u\)
or \(v\), since their logarithmic differentials sum to a
differential from \(B\) after trivialization. Over the other
codimension-one points this is the relative logarithmic form on
a torus; the two omitted vertices have codimension at least two.
To compute its order on \(\nu\), work over a model extracting
the prime of \(w\) and compactify the \(U\)-line by \(0,\infty\).
Near the generic prime the minimum formula makes \(\nu\)
monomial in a local parameter and \(U\) or \(U^{-1}\), with
weights given by \(w\) and \(|\nu(U)|\). The logarithmic
differential has simple poles at \(0,\infty\). Thus one plus the
order of the top form for the primitive valuation is \(A_B(w)\),
also when \(w=0\); the identity is homogeneous under rescaling.
Adding the pullback of \(r\operatorname{div}_Z x_0\) proves
the discrepancy formula.

For the second identity use the degree-\(k\) decomposition
of the section ring from \(uv=G\). In addition to the single
unshifted space \(H^0(B,kH_B)\), each \(1\le j\le k/a\)
gives two copies of \(H^0(B,(k-aj)H_B)\), multiplied by \(u^j\)
and \(v^j\). These are the full sections by the normal
Cohen--Macaulay cone presentation. The two opposite-character
summands can be based independently by the minimum formula.
After trivializing by \(x_0^k\), their shifts sum to \(j w(D)\).
The common factor contributes \(kb\). The unshifted summand has
asymptotic proportion zero. With \(z=aj/k\), the limiting
density of the other summands is proportional to
\((1-z)^{n-1}\) on \([0,1]\), by the Hilbert polynomial on \(B\).
Its means of \(1-z\) and \(z\) are \(n/N\) and \(1/N\).
Using convergence of the complete-series expected orders on
\(B\) gives the claimed formula, first in divisible section
degrees. The endpoint where the residual degree is bounded
has vanishing contribution by the same Hilbert estimate
and linear boundedness.

Set \(c=r/(dN)\). The pair \((B,cD)\) is K-polystable by
\cite[Theorem~1.7]{ar2003.06839}: \(D\sim(d/r)(-K_B)\)
and \(0<c<1-(r/d-1)/(n-1)\). Its automorphism group is
finite since \(B\) is K-stable, so the log pair is uniformly
K-stable. As \(-K_B-cD\sim rnH_B/N\), (V1) gives
\[
      A_Z(\nu)-rS_L(\nu)
        =A_B(w)-c w(D)-\frac{rn}{N}S_{H_B}(w).
\]
This is positive for \(w\ne0\). For \(w=0\) the minimum formula
shows the entire filtration is given by the indicated torus
character on the section ring, up to a constant times the
section degree; its configuration is a product. The torus
equivariant valuative criterion for special configurations and
equivariant K-polystability therefore give K-polystability of
\(Z\), by \cite[Theorem~1.1 and Corollary~4.11]{ar2004.09413}.

Finally replace the first two equations by
\(Q_1=tu,\ Q_2=tv\), retaining \(uv=G\). This is a flat
weighted complete-intersection family near zero by the regular
sequences in its fibers; relative adjunction gives a
\(\mathbf Q\)-Gorenstein family on the coarse space, since
the uniformizing charts are free in codimension one. For
\(t\ne0\) elimination gives \(Q_1Q_2-t^2G=0\). For the general
choice of the three forms its generic member is a smooth
K-stable hypersurface: already at \(t=1\) the quartic, or
more generally the degree-\(d\) equation, is general when \(G\)
varies. The central embedded equation \(Q_1Q_2\) is
GIT-semistable, since both factors are smooth GIT-semistable
forms and minimum weights add under multiplication for every
one-parameter subgroup. Its quotient point differs from every
smooth GIT-stable orbit. This gives the stated boundary example.
\end{proof}

\begin{corollary}\label{cor:nov30-composite-polystable-boundary}
For every composite \(d=ka\), \(k,a\ge2\), and every \(n\ge d-1\),
there are projective K-polystable limits of smooth K-stable
hypersurfaces with hyperplane Cartier index \(a\). One family is
\[
    Q_1(x)=\cdots=Q_k(x)=0,\qquad
    u_1\cdots u_k=G_d(x),\qquad \deg u_i=\deg Q_i=a .
\]
Here all forms are general, with \(n+2\) coordinates of degree one.
\end{corollary}
\begin{proof}
Let \(B=(Q_1=\cdots=Q_k=0)\), \(D=(G_d|_B=0)\) and
\(r=n+2-d\). Then \(B\) is smooth K-stable by
\cite[Theorem~1.3]{ar2608.09293}, \(D\) is smooth, and
\(-K_B=rH_B\). Its locus with \(x\ne0\) is smooth by the
Jacobian criterion and the smoothness of \(D\). At a weight-\(a\) chart put
one \(u_i=1\) on the local cover. It is the hypersurface
\(\prod z_j=G_d(x)\) over the affine cone of \(B\).
Blow up the vertex of the base cone. On its smooth resolution
the pullback of \(G_d\) is the zero section with multiplicity
\(d\) together with the strict transform of the smooth divisor
\(D\), meeting transversely. Consequently the fiber product is
locally given, up to smooth factors and units, by
\(\prod z_j=s^d w\) or \(\prod z_j=s^d\). Each is an irreducible
normal Gorenstein toric hypersurface: irreducibility follows
by solving on the torus, and the Jacobian criterion gives
regularity in codimension one. Such a toric singularity is
canonical, since its nonzero interior lattice weights have
positive integral canonical height. The original chart is also
normal by the complete-intersection depth condition and
regularity in codimension one; the locus over the vertex has
codimension at least two. Adjunction for its one equation
shows that the discrepancy of this partial resolution is the
pullback of \((r-1)\) times the exceptional divisor on the base
cone resolution, which is effective. Thus the chart is klt.
The finite scalar quotient is quasi-étale, with fixed locus
where \(x=0\). The same character argument as in the proposition
gives hyperplane index \(a\) there.

Use the torus scaling the \(u_i\) with product of characters
one. Its invariant valuations are given by the same minimum
formula, now in \(k-1\) Laurent variables over \(\mathbf C(B)\).
The residue transcendence-degree inequality again makes the
restriction \(w\) zero or divisorial up to scaling. Wedging
the base form with the \(k-1\) logarithmic differentials
gives \(A_Z(\nu)=rb+A_B(w)\): on a resolution over the
generic base prime, the order formula is the toric one with
the whole boundary of the torus compactification included.
Over the generic point of \(D\) any vanishing \(u_i\) can
be eliminated using the other invertible \(u_j\), so the
form has no extra divisor there.

A monomial basis over the base ring is indexed by exponent
vectors \(\gamma\ge0\) with \(\min\gamma_i=0\). They have
distinct torus characters. The number with \(\sum\gamma_i=l\)
has leading order a positive multiple of \(l^{k-2}\), by
counting the \(k\) coordinate facets; their overlaps have
lower order. For \(k=2\) the positive levels have two terms.
Thus in section degree \(m'\), the variable \(z=al/m'\)
has limiting density proportional to
\(z^{k-2}(1-z)^{n+1-k}\) on \([0,1]\). Symmetry gives the
same mean for all the exponents, and the sum of the shifts
\(\nu(u_i/x_0^a)\) is \(w(D)\). This gives
\[
 S_L(\nu)=b+
  \frac{(n+2-k)S_{H_B}(w)+(k-1)w(D)/d}{n+1}.
\]
Take \(c=r(k-1)/(d(n+1))<1\). When \(r>d\), the smooth
Fano \(D\) is a general complete intersection and K-stable,
and \(c<1-(r/d-1)/\dim B\) since \(kr<d(n+1)\).
Interpolation gives uniform K-stability of \((B,cD)\) as
in the proposition. When \(r\le d\), for nontrivial \(w\)
use \(w(D)\le A_B(w)\) by smoothness, and
\(A_B(w)>rS_{H_B}(w)\); then
\(A_B(w)-cw(D)>(r-cd)S_{H_B}(w)\). In both cases the
resulting sign is strict when \(w\ne0\). Valuations trivial
on the base give product torus configurations since their
weights are torus characters up to a shift by section degree.
The same reductive equivariant criterion proves K-polystability.
Finally \(Q_i=tu_i\) gives the flat smoothing whose ordinary
equation is \(\prod Q_i-t^kG_d=0\). Its general fibers are
smooth K-stable and the product has semistable GIT class, by
the same argument as in the proposition.
\end{proof}

\section{The hyperplane image of a quartic ninefold replacement}\label{sec:dec1-ninefold-replacements}

\begin{proposition}\label{prop:dec1-factorial-quadric-image}
For a K-polystable replacement $(Z,L)$ of smooth Fano quartic
$n$-folds, the hyperplane map cannot be birational onto a quadric
of rank at least five, in every dimension $n\ge3$.
\end{proposition}
\begin{proof}
Suppose the hyperplane map is birational onto such a quadric $W$.
Write $N=n+1$ and $r=n-2=N-3$, so $-K_Z\sim rL$.
Write $p:Y\to Z$, $g:Y\to W$ for a smooth resolution of its graph,
$H$ for the hyperplane class on $W$, $P=p^*L$ and $J=g^*H$.
Choose the representatives corresponding to the same general
linear section. Then $E=P-J$ is the effective fixed divisor.
Put $D=g_*E=\sum_i b_i\Delta_i$ with $b_i>0$. The quadric cone
coordinate ring is factorial: in the presentation
$x_0x_1+q'(x_2,\ldots,x_{n+1})$, the residual form has rank at least
three, so $x_0$ generates a prime ideal, and inverting $x_0$ gives a
Laurent polynomial ring. Hence $\Delta_i\sim m_iH$ for integers
$m_i\ge1$. Write $q=\sum_i m_i b_i$ and $a_j=P^{n-j}J^j/4$.
The projection formula gives $a_n=1/2$, $a_{n-1}=(1+q)/2$.
Nefness and effectivity of $E$ give $1=a_0\ge\cdots\ge a_n$, so
$q\le1$. Put $\Sigma=\sum_{j=0}^n a_j$ and $s=N-\Sigma$.
The mixed volume integral and the local cyclic-cover threshold
estimate of \cite[Lemmas~5.1, 5.3, 5.4]{kePanWuCubics} give
\[
               A_Z(v)\ge\theta\,v(|L|),\qquad
               \theta=r\Sigma/N.
\]
Since $P^n>J^n$, the base is nonempty. For its codimension $k$,
the same estimates give $k\ge\theta$ and $a_0=\cdots=a_{k-1}=1$.
Put $h=N-k$. If $h\ge6$, then $s\le h/2$ since $a_j\ge1/2$.
Thus $\theta=N-3-(r/N)s>N-3-h/2\ge k$, a contradiction. If $h=5$,
interpolation by log-concavity from $a_{k-1}=1$ to $a_n=1/2$ gives
$s\le\sum_{j=1}^5(1-2^{-j/5})<2$. Indeed $2^{-j/5}$ exceeds
$4/5,7/10,3/5,1/2$ for $j=1,2,3,4$ respectively, and equals $1/2$
for $j=5$. Again $\theta>N-5=k$. Therefore $h\le4$. At most three
entries preceding $a_n$ differ from one, all at least $a_{n-1}$.
Consequently
\[
             0<s\le 2-3q/2,\qquad qs\le q(2-3q/2)\le2/3.
\]

The quartic Hilbert function gives a quadratic section $y$ outside
the degree-two part of the homogeneous coordinate ring of $W$.
Birationality identifies $y$ with a rational section of $2H$; on
$Y$ its divisor relative to $2P$ is effective. Its pole divisor on
$W$ relative to $2H$ is consequently bounded by $\lfloor 2D\rfloor$.
Since $W$ is normal and projectively normal, $y$ must have a pole,
so some $b_j\ge1/2$.

With compatible canonical divisors set $F=g_*(K_Y-p^*K_Z)$.
Its support is contained in $\operatorname{Supp}D$. Indeed, off the
geometric base locus the linear system gives a morphism and any
$p$-exceptional divisor centered there has $g$-image of codimension
at least two. At the generic point of a center in the base the
fixed order is positive, as seen using the ideal generated by the
products of a Cartier-multiple number of degree-one sections in a
local trivialization. Divisors which are not $p$-exceptional have
discrepancy coefficient zero. Let $f_i$ be the coefficient of $F$
at $\Delta_i$, let $v_i$ be its divisorial valuation, and put
$k_i=rb_i-f_i$. Since $K_Z+rL$ is principal,
\[
 rE-(K_Y-p^*K_Z)=p^*(K_Z+rL)-rg^*H-K_Y
\]
is integral. Pushing down gives
\[
            rD-F=\sum_i k_i\Delta_i\sim -rH-K_W=2H,
            \qquad \sum_i m_i k_i=2,\quad k_i\in\mathbb Z.
\]
Since $A_Z(v_i)=1+f_i\ge\theta b_i$ and $b_i\le q$, we have
$k_i\le1+(r/N)sb_i<2$. Thus $k_i\le1$.

If $k_i=1$, then $A_Z(v_i)=r b_i$, and K-semistability gives
$S_L(v_i)\le b_i$. Denote the strict transform of $\Delta_i$ by
$B_i$. For $0\le t\le b_i$, the divisor $P-tB_i$ dominates
the nef divisor $(1-t/b_i)P+(t/b_i)J$ since $E\ge b_iB_i$.
For $b_i\le t\le b_i+1/m_i$ it dominates
$(1-m_i(t-b_i))J$ up to $\mathbb R$-linear equivalence, since
$g^*\Delta_i\sim m_iJ$ contains $B_i$. Integrating both lower
bounds for volume gives
\[
        NS_L(v_i)\ge b_i\Sigma+\frac1{2m_i}.
\]
Thus $m_i b_i s\ge1/2$. Since $qs<1$, at most one $i$ can have
$k_i=1$. The equality
$\sum m_i k_i=2$ forces this index to exist and $m_i\ge2$.
The quadratic section's pole with $b_j\ge1/2$ must occur for this
same index: otherwise $m_i b_i\le q-1/2\le q/2$, which gives
$m_i b_i s<1/2$. Consequently $m_i=2$, $b_i=1/2$,
$D=\Delta_i/2$, and $q=1$. Now $a_{n-1}=1$ and $\Sigma=N-1/2$.
The displayed lower bound is $NS_L(v_i)\ge N/2$, giving equality
$S_L(v_i)=1/2$ and $A_Z(v_i)=r/2$.

The divisor $v_i$ therefore computes $\delta(Z,-K_Z)=1$.
Blum--Liu--Zhou \cite[Theorem~1.1(1)]{ar1907.05399} yield a
nontrivial special test configuration induced by this divisor with
zero Futaki invariant. By K-polystability it is a product, induced
by a one-parameter subgroup of $\operatorname{Aut}(Z)$. Its
valuation on any eigenfunction is the corresponding character,
up to a common nonzero factor and choice of sign. The torus
preserves the Weil class $L$: apply the action map on
$Z\times\mathbf G_m$ and use
$\operatorname{Cl}(Z\times\mathbf G_m)=\operatorname{Cl}(Z)$.
It therefore acts projectively on $H^0(Z,L)$, with a linear lift
after taking a power. Every nonzero linear section has order exactly
$b_i$ at $v_i$, since the prime divisor $\Delta_i\sim2H$ cannot
be contained in any hyperplane section. Every ratio of two linear
sections has order zero. In particular every ratio of two
eigensections has weight zero, so the projective action is trivial.
Birationality of the hyperplane map makes the action trivial on the
function field. This contradicts the nonzero inducing valuation.
\end{proof}

\begin{theorem}\label{thm:dec1-general-degree-gap}
Let $(Z,L)$ be a K-semistable replacement of smooth degree-$d$
Fano hypersurfaces in $\mathbf P^{n+1}$ with $L^n=d$ and
$-K_Z\sim(n+2-d)L$. Write $N=n+1$ and $r=n+2-d$. Suppose that
the full-dimensional hyperplane image $W$ is a normal hypersurface
of degree $e<d$ whose homogeneous coordinate ring is factorial.
Set $\lambda=1-e/d$, $\tau=e/d$. If $e>d/2$, then necessarily
\[
    \min_{0\le u\le\lambda}
    \big[(d-2)u^2-\{\tau^2+(d-3)\lambda\}u+1-2\lambda\big]
                   \le0.
\]
In particular the degree $e=d-2$ is impossible for every $d\ge5$.
More generally, for $h=d-e\ge1$ the same conclusion holds whenever
\[
         d>2h,\qquad
         (h+1)^2\le 4(d-2)(1-2h/d).
\]
\end{theorem}
\begin{proof}
The map to $W$ is birational because its degree times $e$ is at
most $d$ and $e>d/2$. On a resolution $p:Y\to Z$, $g:Y\to W$,
take $P=p^*L$, $J=g^*H$ and representatives with $E=P-J\ge0$.
Write $g_*E=\sum_i b_i\Delta_i$, $\Delta_i\sim m_iH$ with
$m_i\ge1$ integers, $q=\sum_i m_i b_i$, $a_j=P^{n-j}J^j/d$,
$s=N-\sum_{j=0}^n a_j$. The intersection numbers are nonincreasing,
and $a_n=\tau$, $a_{n-1}=\tau(1+q)$. Define $u=1-a_{n-1}$.
Then $0\le u\le\lambda$ and $q=(\lambda-u)/\tau$.

The base locus is nonempty because $P^n>J^n$. Write $k$ for its
codimension. By the same mixed-intersection volume integral and
local cyclic-cover threshold estimate used in
Proposition~\ref{prop:dec1-factorial-quadric-image},
$k\ge r(1-s/N)$ and $a_0=\cdots=a_{k-1}=1$. Put $m=N-k$.
Since $N-r=d-1$ and $r/N<1$, we have
\[
      m\le d-1+(r/N)s<d-1+s,\qquad
      s\le\lambda+(m-1)u.
\]
Here $s>0$ since $a_n=\tau<1$. It follows that
\[
              s\le\frac{\lambda+(d-2)u}{1-u}.
\]

The discrepancy comparison in
Proposition~\ref{prop:dec1-factorial-quadric-image} works for this
normal birational image. With $F=g_*(K_Y-p^*K_Z)$, it gives
\[
           r g_*E-F=\sum_i k_i\Delta_i\sim -rH-K_W=(d-e)H,
           \qquad k_i\in\mathbb Z.
\]
The support and integrality follow from exactly the same local
argument there, using $K_Z+rL\sim0$. In particular some $k_i\ge1$
since $\sum m_i k_i=d-e>0$. For $v_i=\operatorname{ord}_{\Delta_i}$
we have $A_Z(v_i)=1+rb_i-k_i\le rb_i$ and $S_L(v_i)\le b_i$
by semistability. Its two-interval volume integral gives
\[
         N S_L(v_i)\ge b_i(N-s)+\frac{\tau}{m_i}.
\]
Indeed on $[0,b_i]$ one integrates the mixed nef lower bound
between $P$ and $J$, and on $[b_i,b_i+1/m_i]$ the nef multiple
$(1-m_i(t-b_i))J$ using $\Delta_i\sim m_iH$. Thus
$qs\ge m_i b_i s\ge\tau$. Combining with the preceding upper bound
on $s$ gives
\[
    \tau^2(1-u)\le(\lambda-u)(\lambda+(d-2)u).
\]
After rearrangement this is the claimed quadratic inequality.

For $h=d-e$ we have
$\tau^2+(d-3)\lambda< h+1$. The simpler sufficient condition in
the statement therefore makes the discriminant strictly negative.
For $h=2$ it applies for $d\ge8$, since $4(d-2)(1-4/d)\ge12>9$.
For $d=5,6,7$ one computes the discriminant directly: after
multiplication by $d^4$ the integer is respectively $-659$,
$-4208$ and $-14019$. The quadratic has positive leading
coefficient and no real roots. This finishes the proof.
\end{proof}

\begin{theorem}\label{thm:dec1-normal-degree-gap}
Let $(Z,L)$ be a K-semistable replacement of smooth Fano
hypersurfaces of degree $d\ge5$ in $\mathbf P^{n+1}$ as in
Theorem~\ref{thm:dec1-general-degree-gap}. Its complete hyperplane
system cannot have full-dimensional image of degree $d-2$.
\end{theorem}
\begin{proof}
Suppose the image $W$ has degree $e=d-2$, and let
$\nu:\overline W\to W$ be its normalization. Since $e>d/2$,
the hyperplane map is birational. The resolution $g:Y\to W$
factors through $\overline W$. Use $H=\nu^*\mathcal O_W(1)$ on
$\overline W$ and write $\Delta_i$ for prime divisors there.
Set $D_i=H^{n-1}\Delta_i$ and $m_i=D_i/e$. With
$\lambda=2/d$, $\tau=e/d$, the intersection and base-codimension
computations of Theorem~\ref{thm:dec1-general-degree-gap} still
apply without any target divisor-class assumption. In particular
$q=\sum m_i b_i=(\lambda-u)/\tau$, with the same $s$ and $u$.
Keep $m=N-k$ from the preceding theorem (it is distinct from the
rational $m_i$). If $m=1$, then $s=\lambda$, $u=0$, and the
bound $qs=\lambda^2/\tau<1/d$ holds for $d\ge9$. Suppose $m\ge2$. Interpolating
by log-concavity between $a_{k-1}=1$ and $a_{n-1}=1-u$ gives
\[
 s\le\lambda+\sum_{j=1}^{m-1}\bigl[1-(1-u)^{j/(m-1)}\bigr]
 \le\lambda+\frac{m u}{2(1-u)}.
\]
For the last bound use $1-(1-u)^x\le-x\log(1-u)\le xu/(1-u)$.
Since $m<d-1+s$ by the preceding theorem, we obtain
\[
          s\le\frac{2\lambda+(d-1-2\lambda)u}{2-3u}.
\]
For $d\ge9$ assuming $qs\ge1/d$ would give
$(\lambda-u)s\ge(d-2)/d^2$, hence
\[
       d(d^2-d-4)u^2-(2d^2-3d-14)u+2(d-6)\le0.
\]
Its discriminant $(2d^2-3d-14)^2-8d(d^2-d-4)(d-6)$ is strictly
negative for $d\ge9$. Indeed at $d=9$ its value and its successive
forward differences of orders one through four are
$-47,-3137,-2294,-744,-96$. The binomial-basis expansion at nine
therefore has only negative coefficients. Consequently $qs<1/d$ for $d\ge9$.

Taylor's theorem, for $0\le x\le1$, gives
\[
       1-(1-u)^x\le xu+\frac{x(1-x)u^2}{2(1-u)^2}.
\]
Sum this over $x=j/(m-1)$ to obtain
$s\le \lambda+m u(1/2+u/[12(1-u)^2])$. This bound also holds
for $m=1$ because then $u=0$ and $s=\lambda$. For $d=6,7,8$,
use $0\le u\le\lambda$ and put
$c=1/2+\lambda/[12(1-\lambda)^2]$. For these degrees $c$ equals
$9/16,41/75,29/54$, respectively. As $m<d-1+s$ we have
\[
    s\le\frac{\lambda+(d-1)c\lambda}{1-c\lambda}<2,
    \qquad m\le d.
\]
Consequently
\[
 qs\le\max_{u\in\mathbb R}
      \frac{(\lambda-u)(\lambda+dcu)}{\tau}
    =\frac{\lambda^2(1+dc)^2}{4dc\tau}.
\]
Substitution gives upper bounds less than $24/100,18/100,14/100$
respectively, each strictly smaller than $(2-1/e)/d$.

For $d=5$, the Taylor factor is at most $17/28$, since $u\le2/5$.
Putting $c=17/28$ in the preceding $s$ estimate again gives
$s<2$, so $m\le5$. If $m\le4$ the same quadratic maximum, now
using $4c$ in place of $dc$, is less than $1/3$. If $m=5$,
we must have $s>1$. Then $u>1/5$: otherwise the Taylor factor
is at most $1/2+(1/5)/[12(4/5)^2]<10/17$ and
$s\le 2/5+5(1/5)(10/17)<1$.
For $1/5<u\le1/4$ the factor is at most $11/20$; the quadratic
$(2/5-u)(2/5+5(11/20)u)/(3/5)$ decreases there and at $u=1/5$
is less than $1/3$. For $1/4\le u\le2/5$ the analogous quadratic
with $c=17/28$ decreases and its value at $1/4$ is less than
$1/3$. Hence $qs<1/3=(2-1/e)/d$, as required.

The estimates proved so far imply $e q s<1$ in every degree.
For $d=5$ this is $3qs<1$, for $d=6,7,8$ use the displayed
numerical bounds, and for $d\ge9$ use $e<d$. On $\overline W$
the conductor divisor $C=\sum c_i\Delta_i$ is effective and
satisfies $K_{\overline W}+C\sim\nu^*K_W$. This normalization
formula holds for the integral hypersurface $W$, whose dualizing
sheaf is invertible: by finite duality
$\nu_*\omega_{\overline W}=\mathcal Hom_W(\nu_*\mathcal O_{\overline W},\omega_W)$,
and $\mathcal Hom_W(\nu_*\mathcal O_{\overline W},\mathcal O_W)$
is the conductor ideal. Take the divisorial part on the normalization.

Push $rE-(K_Y-p^*K_Z)$ onto $\overline W$, as in
Proposition~\ref{prop:dec1-factorial-quadric-image}, and denote
its integral coefficients by $k_i$. Its divisor class is
$-rH-K_{\overline W}\sim (d-e)H+C$. Hence
$\sum_i D_i(k_i-c_i)=e(d-e)>0$, including in the finite sum the
conductor components outside the fixed support with $k_i=b_i=0$.
The threshold estimate $A_Z(v_i)\ge r(1-s/N)b_i$ gives
$k_i\le1+(r/N)sb_i$. Since $D_i\ge1$ and $b_i\le e q$, it follows
that $k_i<2$. Thus there is an $i$ with $k_i=1$, $c_i=0$. Its
center on $W$ lies generically in the normal locus. Indeed if the
local ring at the image prime were not normal, its conductor would
be a proper ideal of that local ring. It would be contained in its
maximal ideal, hence in every maximal ideal of its semilocal
normalization, giving strictly positive conductor coefficient on
each branch. Thus $D_i$ equals the ordinary projective degree of
the image of $\Delta_i$. For $v_i=\operatorname{ord}_{\Delta_i}$,
we have $A_Z(v_i)=rb_i$ and $S_L(v_i)\le b_i$.
We show the lower bound $qs\ge(2-1/e)/d$. The complete
system of degree-$D_i$ hypersurfaces containing the image of $\Delta_i$ has
no other fixed prime divisor on $W$, and its generic fixed order
along $\Delta_i$ is one. Here the generic normality just proved ensures that its image has a unique divisorial preimage. For the first assertion, at any given
point $x$ outside $\Delta_i$ choose a projection center outside
the cone over $\Delta_i$ with vertex $x$; the pullback of the
projected equation misses $x$. For the second, take a general
smooth point $z$ of both $W$ and $\Delta_i$ and project from a
center outside $T_zW$ and outside the secant cone joining $z$ to
the remaining points of $\Delta_i$. The linear projection is
finite on $\Delta_i$ and is an immersion at $z$, with singleton
fiber there. Its image is smooth at the projected point. The
pullback of the local hypersurface equation restricts to order
one on $W$ at $z$: its tangent hyperplane is spanned by the
center and $T_z\Delta_i$, so differs from $T_z W$. All centers
excluded here lie in a proper closed subset of $\mathbf P^{n+1}$.
Multiplying by forms missing the indicated point if needed gives
exactly degree $D_i$ in both verifications.

Write $M$ for the resulting moving linear system on the common
resolution and $B_i$ for the strict transform of $\Delta_i$. Then
$D_iJ-B_i-M$ is effective up to linear equivalence and
$J^{n-1}M=eD_i-D_i$, since the other fixed components are
exceptional over $W$. For $0\le v\le1/D_i$ use the nef lower
bound $(1-D_iv)J+vM$ for $J-vB_i$. Its first two mixed terms
in the tail integral for $P-(b_i+v)B_i$ give
\[
    N S_L(v_i)\ge b_i(N-s)+\frac{\tau}{D_i}(2-1/e),\qquad
    q s\ge(2-1/e)/d.
\]

For $d\ge9$ this contradicts $qs<1/d$. For $5\le d\le8$ it
contradicts the corresponding stronger numerical bounds proved
above, which finishes the proof.
\end{proof}

\begin{corollary}\label{cor:dec1-unconditional-degree-minus-one}
For every $d\ge3$, the hyperplane image of a K-semistable
replacement in Theorem~\ref{thm:dec1-general-degree-gap} cannot
be an $n$-dimensional hypersurface of degree $d-1$.
\end{corollary}
\begin{proof}
Take $e=d-1$, $\lambda=1/d$ and $\tau=1-\lambda$.
The image map would be birational. The elementary estimate
$s\le[\lambda+(d-2)u]/(1-u)$ in
Theorem~\ref{thm:dec1-general-degree-gap} uses only the mixed
intersections and the base-codimension bound and applies on the
normalization. Using $1-u\ge\tau$ and
$u(\lambda-u)\le\lambda^2/4$ gives
\[
 qs\le\frac{(\lambda-u)[\lambda+(d-2)u]}{\tau(1-u)}
   \le\frac{(d+2)\lambda^2}{4\tau^2}
   =\frac{d+2}{4(d-1)^2}<\frac1d.
\]
The final inequality is $3d^2-10d+4>0$ for $d\ge3$. In particular
$e q s<1$. The normalization, conductor and integral discrepancy
argument in Theorem~\ref{thm:dec1-normal-degree-gap} gives a prime
with $k_i=1$, $c_i=0$: the right-hand side of the degree relation
is now $e(d-e)=e>0$. For this prime, a degree-$D_i$ hypersurface
containing its image provides the tail estimate $qs\ge1/d$,
already a contradiction.
\end{proof}

\begin{proposition}\label{prop:dec1-optimal-factorial-quadric}
Let $(Z,L)$ be an optimal destabilizing special fiber of a smooth
quartic ninefold $X$ with $\delta(X,-K_X)\le1$. If the hyperplane
image of $Z$ is nine-dimensional, the hyperplane map is birational
onto a quadric of rank four.
\end{proposition}
\begin{proof}
By Corollary~\ref{cor:nov30-optimal-image-ranks}, it suffices to
exclude a birational quadric image of rank at least five. Write
$\beta=\delta(Z,L)$; as in that corollary $20/3<\beta\le7$.
Use the notation of Proposition~\ref{prop:dec1-factorial-quadric-image}.
The same first volume integral applied to any valuation with fixed
order $b$ gives $A_Z\ge\beta b\Sigma/10$. On the local cyclic
cover this gives the base-ideal threshold bound
$\theta=\beta\Sigma/10$ by the finite-cover discrepancy formula
in the proof of \cite[Lemma~5.3]{kePanWuCubics}. If the base is
nonempty of codimension $k$, then $k\ge\theta$ and $a_0=\cdots=
a_{k-1}=1$ as in Proposition~\ref{prop:dec1-factorial-quadric-image}.
The base cannot be empty since $P^9>J^9$.
Initially $a_j\ge2^{-j/9}$ gives $\theta>4$: the first five
terms exceed $2^{-1/2}$ and the last five are at least $1/2$.
Thus $k\ge5$. Interpolation from $a_4=1$ gives
$\Sigma\ge 5+2^{-1/5}+2^{-2/5}+2^{-3/5}+2^{-4/5}+1/2>15/2$.
Thus $\theta>5$, $k\ge6$, and $a_0=\cdots=a_5=1$.

Put $t=70/\beta-\Sigma>0$ so that $7-\theta=\beta t/10$.
From $a_6,a_7\ge a_8=(1+q)/2$ one gets
$t<5/2-3q/2$. For $q\le2/3$ this implies $qt<1$. When $q\ge2/3$,
$\Sigma\ge 6+3(1+q)/2+1/2\ge9$ gives $\theta>6$ and $a_6=1$.
Consequently $t<2-q$, again yielding $qt<1$ since $q\le1$.
The integral coefficients $k_i$ in the discrepancy comparison
satisfy $k_i\le1+(7-\theta)b_i<1+\beta/10<2$, because
$b_i\le q$. Each index with $k_i=1$ has $A_Z(v_i)=7b_i$.
Using $\beta S_L(v_i)\le A_Z(v_i)$ in the two-interval volume
inequality gives $q_i t\ge1/2$ for $q_i=m_i b_i$. There can be
at most one such index, by $qt<1$. As $\sum m_i k_i=2$, the
unique positive index has $m_i\ge2$.

The extra quadratic section $y$ has a pole at some $\Delta_j$
with $b_j\ge1/2$. If $j\ne i$, then $q_i\le q-1/2\le q/2$,
contrary to $q_i t\ge1/2$. Thus $m_i=2$, $b_i=1/2$, $q=1$
and $D=\Delta_i/2$. Let $h$ be a quadratic equation of $\Delta_i$
in the factorial ring $T$. The pole bound for $y$ shows that
$g=hy$ is a regular quartic section on $W$, hence $g\in T_4$.
The prime $h$ does not divide $g$ in $T$ since $y\notin T_2$.
Thus $T[y]=T[Y]/(hY-g)$, with $Y$ of degree two. Its Hilbert
series equals the quartic hypersurface series; this gives the
entire ring $R(Z,L)$ and makes it a complete intersection. This
contradicts Corollary~\ref{cor:nov30-all-ci-hyperplane-degenerations}.
\end{proof}

\begin{proposition}\label{prop:dec1-rank-four-pencils}
Let $Z$ be a normal projective klt variety of dimension $n\ge2$ and
$L$ an ample $\mathbf Q$-Cartier integral Weil divisor. If the
kernel of $\operatorname{Sym}^2 H^0(Z,L)\to H^0(Z,2L)$ contains
a rank-four quadric, then $\delta(Z,L)\le2(n+1)/3$. In fact a prime
divisor $B$ on $Z$ satisfies $(n+1)S_L(B)\ge3/2$.
\end{proposition}
\begin{proof}
Write a rank-four relation $s_{00}s_{11}=s_{01}s_{10}$ using four
linearly independent sections. The rational map they define takes
values in the Segre surface $\mathbf P^1\times\mathbf P^1$.
Each projection is nonconstant, since constancy would give a
linear relation between the sections. Resolve both maps and the
base ideal of the subsystem on a common smooth projective model
$p:Y\to Z$. Let $M_1,M_2$ be the pullbacks of the two pencil classes
and $P=p^*L$. Then
\[
       P\sim_{\mathbf Q} M_1+M_2+G,\qquad G\ge0.
\]
Indeed, take an integer $k>0$ for which $kL$ is Cartier. The degree-$k$
products of the four sections define the Veronese of the original
rational map, with free class $k(M_1+M_2)$. Divide its fixed divisor
by $k$ to obtain $G$. At the strict transform of any prime divisor
on $Z$ its coefficient is the minimum of the four vanishing orders,
which is integral.

Put $V=P^n$, $a_i=P^{n-1}M_i$, and $N=n+1$. We use two volume
inequalities. For a prime $B$ on $Z$ with strict transform $B'$,
if $P-B'$ dominates a nef divisor $M$ up to $\mathbf Q$-linear
equivalence, then for $0\le t\le1$ the volume of $P-tB'$ is at
least $((1-t)P+tM)^n$. Integration of the first two terms gives
\[
                  NS_L(B)\ge1+P^{n-1}M/V.
\]
Also, if $B'$ is a component of a general member of $|M_i|$,
then $M_i-B'$ is effective up to linear equivalence. The algebraic
Morse inequality for the nef divisors $P$ and $tM_i$ gives
$\operatorname{vol}(P-tB')\ge V-nt a_i$. Integration up to
$V/(n a_i)$ gives $NS_L(B)\ge NV/(2na_i)$. A nonexceptional
component of such a general member always exists: a nonconstant
rational pencil has a divisorial general fiber meeting the open
set where $p$ is an isomorphism. In particular $a_i>0$.

If $a_1+a_2<V/2$, choose $i$ with $a_i<V/4$ and apply the second
inequality. It gives $NS_L(B)>2N/n>3/2$. Suppose henceforth that
$a_1+a_2\ge V/2$. If $G$ contains a nonexceptional prime $B'$, its
coefficient is at least one. The first inequality with $M=M_1+M_2$
gives $NS_L(B)\ge3/2$. Otherwise the projection formula gives
$P^{n-1}G=0$, whence $a_1+a_2=V$. Relabel so $a_2\ge V/2$ and
choose a nonexceptional prime component $B'$ of a general fiber in
$|M_1|$. Now $P-B'$ dominates $M_2$, giving $NS_L(B)\ge3/2$ once
again. As $A_Z(B)=1$, this proves the threshold bound.
\end{proof}

The constant is sharp in every dimension. Let $Q^n\subset
\mathbf P^{n+1}$ be the rank-four quadric, the cone over
$\mathbf P^1\times\mathbf P^1$ with $l=n-2$ independent vertex
coordinates, and put $H=\mathcal O_Q(1)$. Its toric prime divisors
are the four ruling divisors and, if $l>0$, one divisor for each
vertex coordinate. In degree $m$, a toric monomial basis decomposes
by the Segre degree $k$: its summands have $(k+1)^2$ choices of
a Segre monomial and $\binom{m-k+l-1}{l-1}$ choices in the vertex
coordinates. The expected exponent of one of the two variables
of a Segre factor is $k/2$. Dividing by $m$, the limiting mean
of $k$ is $3/(l+3)$ by the beta integral with density proportional
to $x^2(1-x)^{l-1}$, or is one directly if $l=0$. By symmetry
each vertex coordinate has limiting mean exponent $1/(l+3)$.
Consequently the four ruling divisors have expected order
$3/(2(n+1))$ and all remaining toric prime divisors have expected
order $1/(n+1)$. The toric threshold formula (equivalently the
minimum of the barycentric ratios over the facets of the moment
polytope) gives $\delta(Q^n,H)=2(n+1)/3$.

In particular, for a smooth quartic ninefold with $\delta(X,-K_X)\le1$,
the hyperplane image on every optimal destabilizing special fiber has
dimension at most eight. Indeed its hyperplane threshold is strictly
greater than $20/3$ by Corollary~\ref{cor:nov30-optimal-image-ranks},
whereas Proposition~\ref{prop:dec1-optimal-factorial-quadric} would
force a rank-four quadric relation if the image had dimension nine.
More generally the same rank-four obstruction applies to every
optimal quartic fiber in dimension $n\ge9$ and to K-polystable
quartic replacements in these dimensions, independently of image
dimension.

\begin{proposition}\label{prop:dec1-lower-five-six}
On a K-polystable replacement of smooth quartic ninefolds the image
of the limiting hyperplane class cannot have dimension five or six.
\end{proposition}
\begin{proof}
Let $a\in\{5,6\}$ be the image dimension, $W\subset\mathbf P^{10}$
the image, $B\to W$ a smooth projective resolution, and $H$ the
hyperplane pullback on $B$. Resolve $Z\dashrightarrow B$ by
$g:Y\to B$ and $p:Y\to Z$; set $P=p^*L$, $J=g^*H$.
We have $P-J\ge0$ up to rational linear equivalence, and
$h^0(B,H)=11$ by the section-inclusion argument in
Section~\ref{cubics:minimal}. The degree $e=H^a$ exceeds the
minimal value $11-a$, by Proposition~\ref{prop:nov30-optimal-minimal-image}.

The slope bound for $Z\dashrightarrow B$ and the bound for a map
to $\mathbf P^1$ give, with $A=K_B+(a-1)H$,
\[
 \frac{P^8g^* A}{4}\le a-1-\frac{7a}{9}=\frac{2a-9}{9},
 \qquad \frac{P^8V}{4}\ge\frac7{18}
\]
for the resolved moving fiber class $V$ of any pencil on $Z$.
If $h^0(B,A)\ge2$, a pencil in this system pulls back to $Z$;
on a higher common resolution its moving divisor is dominated
by $g^*A$ up to linear equivalence. The two bounds contradict
$(2a-9)/9<7/18$. Thus $h^0(B,A)\le1$.

Cut $B$ by general members of $|H|$ down to a smooth curve $C$.
On each intermediate $j$-fold ($j\ge2$) the polarization remains
big, nef and free; connectedness and Bertini give integral smooth
sections. Irregularity is zero on $B$ since it is dominated by
a resolution of the klt Fano $Z$, and remains zero down to the
surface by Kawamata--Viehweg vanishing and Serre duality as in
Section~\ref{cubics:minimal}. The adjunction sequences, with
Kawamata--Viehweg vanishing when $j\ge3$ and surface irregularity
when $j=2$, give a surjection $H^0(B,A)\to H^0(C,K_C)$.
Thus $g(C)\le1$. The restriction maps on $H^0(H)$ are also
surjective down the ladder, each with one-dimensional kernel by
the irregularity vanishing. Consequently $h^0(C,H|_C)=12-a$.
If $C$ were rational this would give $e=11-a$, already excluded.
Thus $C$ is elliptic, and Riemann--Roch gives $e=12-a$. The
restriction of the hyperplane system embeds $C$ as an elliptic
normal curve in its linear span.

We use the following consequence of the proof of
Proposition~\ref{prop:nov30-optimal-minimal-image}. No subsystem
of $|L|$ can have a minimal-degree image of codimension at least
two in its span, of positive dimension less than nine. Indeed,
the classification supplies a scroll or the cone over the Veronese
surface. In the double hyperplane-section cases, pulling the
hyperplane back gives $10 S_L(E)\ge2$ for a prime divisor $E$
on $Z$. In the remaining case, the image is a cone over
$\mathbf P^1\times\mathbf P^t$, $t\ge2$. Its Segre coordinates
give $P-U-V$ effective up to rational linear equivalence, with
$U,V$ the resolved pullbacks from the factors. If $E'$ is a prime
component in a general first-pencil fiber mapping onto a divisor
$E$ of $Z$, the standard first mixed term of the volume integral
gives
\[
  10S_L(E)\ge1+\frac{P^8V}{4}
       \ge1+\frac{7t}{9(t+1)}\ge 1+\frac{14}{27}>\frac{10}{7}.
\]
Both alternatives contradict $7S_L(E)\le A_Z(E)=1$. These
arguments use only sections in the subsystem and effectivity
of $P$ minus their resolved moving divisor.

Remove all linear vertex coordinates of $W$, giving a nonconical
embedded base $W_0$ of dimension $b$, degree $e$ and codimension
$c=e-2$. Projection from the vertex is given by a subsystem of
$|L|$. If $W_0$ has a singular point $z$, inner projection from
$z$ is generically finite: a positive-dimensional general fiber
would be a line through $z$ contained in $W_0$, making $z$ a
vertex. Projection therefore gives a nondegenerate image $W'$ of
codimension $c-1$ satisfying
\[
     (\deg\text{ projection})\,\deg W'
             =e-\operatorname{mult}_z W_0\le e-2=c.
\]
This is the usual projection formula, obtained by intersecting with
$b$ general hyperplanes through $z$ and subtracting the local
multiplicity. Nondegeneracy gives $\deg W'\ge c$, so $W'$ has
minimal degree and codimension $e-3\ge3$. Its coordinates once
again give a subsystem of $|L|$, which the preceding paragraph
excludes. Therefore $W_0$ is smooth.

The base $W_0$ is isomorphic, with its polarization $H_0$, to a
general linear section complementary to the vertex, so its general
curve section is elliptic. If $b=1$, $W_0$ is an elliptic curve
dominated by $Z$, impossible by the irregularity vanishing. If
$b\ge2$, the same adjunction ladder gives a nonzero section of
$K_{W_0}+(b-1)H_0$. This divisor has $H_0^{b-1}$-degree zero by
adjunction on the elliptic curve. Ampleness of $H_0$ forces
$-K_{W_0}\sim(b-1)H_0$. Thus $W_0$ is a smooth, nonconical
del Pezzo variety of degree six or seven.

The rational map $Z\dashrightarrow W_0$ gives one more slope bound.
For its hyperplane moving class $M$ on the common resolution, we
have $P-M$ effective up to rational linear equivalence, while the
identity $-K_{W_0}=(b-1)H_0$ gives
\[
       1\ge \frac{P^8M}{4}\ge \frac{7b}{9(b-1)}.
\]
Thus $b\ge5$. The nonconical smooth del Pezzo varieties of degree
six have dimension at most four, and in degree seven at most three,
by Kuznetsov--Prokhorov, \emph{On higher-dimensional del Pezzo
varieties}, Theorem~1.2(vi)--(vii),
\href{https://arxiv.org/abs/2206.01549}{arXiv:2206.01549}.
The hyperplane embedding of the base is complete: the same
irregularity ladder gives $h^0(W_0,H_0)=e+b-1$, equal to the
number of ambient coordinates after removing the vertex. Thus
nonconicality in the classification is exactly the absence of a
linear vertex used above. This contradiction proves the result.
\end{proof}

\begin{proposition}\label{prop:dec1-lower-seven}
If the limiting hyperplane image of a K-polystable replacement of
smooth quartic ninefolds has dimension seven, it has degree five
and is a cone over $\operatorname{Gr}(2,5)$ in its Pl\"ucker
embedding or over its smooth hyperplane section.
\end{proposition}
\begin{proof}
We first record a consequence of the slope and volume bounds on
$Z$. A prime divisor $E$ cannot be contained in three independent
members of $|L|$, meaning that three independent linear sections
all vanish along $E$. Indeed, resolve the subsystem spanned by the
three sections and let $V$ be its moving class. The effective
divisor $P-V$ contains the strict transform $E'$. If the image
has dimension two, the slope inequality gives $P^8V/4\ge14/27$.
If it has dimension one, its normalization is rational and its
degree is at least two by nondegeneracy in the plane, so
$P^8V/4\ge 2(7/18)$. On integrating the nef lower bound
$(1-t)P+tV$ for $P-tE'$, $0\le t\le1$, one obtains
$10S_L(E)\ge1+P^8V/4>10/7$, a contradiction.

For the adjunction computation use the notation of
Proposition~\ref{prop:dec1-lower-five-six}. In fact for $a=7,8$
the restriction $H^0(B,A)\to H^0(C,K_C)$ with $A=K_B+(a-1)H$
is an isomorphism. The slope bound excludes an effective
$K_B+(a-2)H$, since $7a/9>a-2$. Write $B_j$ for the smooth
$j$-dimensional sections in the ladder. The absence of sections of
$K_{B_j}+(j-2)H$ descends down to $j=2$: restriction is surjective
from each $j\ge4$ by Kawamata--Viehweg, and from $j=3$ by
$H^1(K_{B_3})=0$. For the latter use $H^2(B_3,\mathcal O)=0$:
the smooth projective $B$ is rationally connected since dominated
by a resolution of the klt Fano $Z$, and $H^i(\mathcal O_{B_j})=0$
for $0<i<j$ follows down the ladder by duality and
Kawamata--Viehweg. Hence all the kernels in the restrictions of
the adjoint system to the respective adjoint systems are zero.
Surjectivity was proved by the same ladder in the preceding
proposition; this proves the isomorphism.

Now take $a=7$, and let $g_C$ be the sectional genus. The
nonminimal degree gives $g_C\ge1$ as before. On a resolution,
the normalized mixed degree of the moving adjoint class is at
most $5/9$ by the determinant slope as above. If its image has
dimension $t\ge1$, project this image generally and dominantly to
$\mathbf P^t$. The moving class of this projection is dominated
by the original adjoint class, so $7t/[9(t+1)]\le5/9$ and
$t\le2$. If $t=1$, its normalized image curve is rational, of
hyperplane degree at most $(5/9)/(7/18)<2$. Thus the image is
a line and $g_C=h^0(B,A)=2$.
If $t=2$, let $S$ resolve the adjoint image and $M$ be the
hyperplane pullback. Were $K_S+2M$ effective, the determinant
slope for $Z\dashrightarrow S$ would give normalized moving
degree at least $7/9$, a contradiction. The restriction to a
smooth curve $D\in|M|$ gives a surjection onto $H^0(K_D+M|_D)$
by Kawamata--Viehweg. Its vanishing forces $D$ rational and
$\deg M|_D=1$. Thus the image is a plane and $g_C=3$.

In the pencil case, use the morphism to the smooth Stein curve
of a resolution of the adjoint map; it is a rational curve. A
general prime fiber in $B$ has $H^6$-degree at most $2$, since
the complete adjoint restricts to $K_C$ of degree two and its
moving fiber occurs with positive coefficient. This fiber meets
the birational locus over $W$, and its image has dimension six.
An integral variety of degree at most two and dimension six spans
at most $\mathbf P^7$, by the minimal-degree inequality. There
are consequently at least three independent hyperplanes containing
this general fiber image. A prime component of its pullback on
the dominant regular locus of $Z$ is a prime divisor in all three
hyperplane zero divisors, already excluded.

In the plane case, a general moving member of the adjoint net on
$B$ maps to an integral divisor of $W$ by Bertini irreducibility
and birationality. It has $H^6$-degree at most $\deg K_C=4$, so
spans at most $\mathbf P^9$. A linear form in $\mathbf P^{10}$
therefore vanishes on it. On $Z$ the corresponding general pullback
of a line in the adjoint plane is again an integral moving divisor
on the regular locus, by Bertini. Denote its Weil divisor class
by $A_Z$. The linear form just found yields an effective
representative $G$ of $L-A_Z$. If $G$ contains a prime, multiplying
its defining section with the three independent sections of the
adjoint net in $A_Z$ gives the prohibited three hyperplanes through
that prime. If $G=0$, the moving net class on the resolution has
normalized mixed degree $P^8V/4=1$ by pushing down to $A_Z\sim L$.
This contradicts the upper bound $5/9$. Thus only $g_C=1$ remains.

The elliptic curve computation gives $e=12-7=5$. Removing the
linear vertex and projecting from any remaining singular point as
in Proposition~\ref{prop:dec1-lower-five-six} would give a forbidden
minimal-degree subsystem in codimension $e-3=2$. Thus the base is
smooth. Its dimension $b$ is at least five by the same computation
$1\ge7b/[9(b-1)]$, and its polarization is the complete fundamental
divisor of a degree-five del Pezzo variety. Theorem~1.2(v) of
Kuznetsov--Prokhorov, cited in the preceding proposition, gives
the two stated linear sections of $\operatorname{Gr}(2,5)$.
\end{proof}

\begin{lemma}\label{lem:dec1-hyperplane-atiyah-slope}
Let $Z$ be a K-polystable $n$-dimensional $\mathbb Q$-Fano with
$-K_Z\sim rL$, for an ample $\mathbb Q$-Cartier Weil divisor $L$.
Suppose a subsystem of $|L|$ has no fixed prime divisor and gives
a dominant rational map to a smooth resolution $B$ of its image.
Let $a=\dim B>0$, $H$ be the hyperplane pullback on $B$, and
$p:Y\to Z$, $g:Y\to B$ a smooth resolution of the map. Then
\[
   \frac{(p^*L)^{n-1}g^*K_B}{L^n}\le-\frac{r(a+1)}{n+1}.
\]
\end{lemma}
\begin{proof}
On a smooth variety let $\mathcal A(M)$ denote the normalized
first-jet extension for a line bundle $M$, with exact sequence
$0\to\Omega^1\to\mathcal A(M)\to\mathcal O\to0$ and extension
class $c_1(M)$ up to sign. Equivalently use
$J^1(M)\otimes M^{-1}$. On a big open $U\subset Z_{\rm reg}$
the map to the proper $B$ is defined and $h^*H=L|_U$: remove
the indeterminacy loci and the common zero locus of the subsystem,
which have codimension at least two. Pullback of first jets gives
\[
             h^*\mathcal A(H)\longrightarrow\mathcal A(L|_U),
\]
inducing the differential on the left terms and the identity on
the right terms. The map is generically injective by dominance
in characteristic zero. Since $c_1(-K_Z)=r c_1(L)$, scalar
change on one end of the extension identifies the target with the
dual of the canonical extension of $T_Z$ on $U$; common sign
conventions do not change its isomorphism class. The reflexive
extension to $Z$ is slope semistable with respect to $L$ by
Druel--Guenancia--P\u{a}un, \emph{A decomposition theorem for
$\mathbb Q$-Fano K\"ahler--Einstein varieties}, Theorem~B,
\href{https://arxiv.org/abs/2008.05352}{arXiv:2008.05352}.
The rank-$(a+1)$ source has determinant $h^*K_B$; saturating
the image can only increase its determinant degree. The target
has rank $n+1$ and first Chern class $K_Z$. Comparing slopes and
using the projection formula proves the inequality.
\end{proof}

\begin{corollary}\label{cor:dec1-only-grassmannian}
In Proposition~\ref{prop:dec1-lower-seven} the base is necessarily
the entire Grassmannian $\operatorname{Gr}(2,5)$. In particular
its seven-dimensional cone has a point as vertex.
\end{corollary}
\begin{proof}
The subsystem projecting to the smooth vertex-reduced base $W_0$
has at least three independent linear sections, so has no fixed
prime by the first paragraph of Proposition~\ref{prop:dec1-lower-seven}.
Let $b=\dim W_0$. The polarization on $W_0$ satisfies
$-K_{W_0}=(b-1)H_0$ and its pullback moving class on a resolution
has normalized mixed degree one by the projection formula.
Lemma~\ref{lem:dec1-hyperplane-atiyah-slope} gives
$b-1\ge7(b+1)/10$, so $b\ge6$. The smooth hyperplane-section
alternative has $b=5$, proving the assertion.
\end{proof}

\begin{proposition}\label{prop:dec1-eight-adjoint}
Suppose the hyperplane image of a K-polystable quartic-ninefold
replacement has dimension eight, and let $B,H,C$ denote its smooth
resolution, hyperplane class and smooth sectional curve. Put
$A=K_B+7H$ and $g_C=g(C)$. Then $g_C\ge1$, $g_C\ne2$, and if
$g_C\ge3$ the image of $|A|$ fills its projective space
$\mathbf P^{g_C-1}$. In that case $g_C\le9$.
\end{proposition}
\begin{proof}
As in Proposition~\ref{prop:dec1-lower-seven}, $h^0(B,A)=g_C$;
minimal degree is excluded so $g_C\ge1$. The original hyperplane
system has no fixed prime by the three-section criterion. Thus
Lemma~\ref{lem:dec1-hyperplane-atiyah-slope} gives
\[
             P^8g^*A/4\le7-63/10=7/10.
\]
Suppose $g_C\ge2$ and let $S$ resolve the adjoint image, $t=\dim S$,
and $M$ be the adjoint hyperplane pullback. The adjoint moving
system (on a higher common resolution) has normalized mixed degree
at most $7/10$. If $K_S+tM$ is effective, tangent-sheaf slope for
the map from $Z$ to $S$ gives normalized mixed degree at least
$7/9$, a contradiction. For $t=1$, the base is rational by
dominance, and this rules out $\deg M\ge2$. For $t\ge2$ restrict
$K_S+tM$ to a general hyperplane ladder to a smooth curve $D$.
The restrictions are surjective by Kawamata--Viehweg since at each
step from dimension $j\ge2$ their kernel sheaf is $K_{S_j}+(j-1)M$.
Hence $H^0(D,K_D+M|_D)=0$. Riemann--Roch on $D$ implies that
it is rational and $\deg M|_D=1$. Thus the adjoint image is
nondegenerate of degree one and fills $\mathbf P^{g_C-1}$.
Its dimension is at most eight.

It remains to exclude $g_C=2$. Here the adjoint pencil maps to a
rational curve. Pass to the rational Stein curve on $B$ and then,
if necessary, on a resolution of $Z$. A general connected fiber
of $B$ over its Stein curve maps to a moving prime divisor on
$W$ of $H^7$-degree at most two, by restriction of the adjoint
pencil to $C$. Its seven-dimensional image spans at most
$\mathbf P^8$ in $\mathbf P^{10}$. On $Z$ let $D$ be the moving
prime divisor from a general fiber over the final Stein curve.
Then $h^0(Z,D)\ge2$. At least two independent hyperplanes vanish
on $D$, by the bound on the span of its image, giving
$h^0(Z,L-D)\ge2$.

Choose a pencil in each of $D$ and $L-D$ and resolve them with
moving fiber classes $U,V$. Thus $P-U-V$ is effective up to
rational linear equivalence. If its fixed part includes a prime
$E'$ dominating a divisor $E$ of $Z$, its coefficient is at least
one (orders of Weil sections along a prime on $Z$ are integers).
The volume bound for $P-tE'$ using $U+V$ gives
$10S_L(E)\ge1+P^8(U+V)/4\ge1+7/9>10/7$, using slope for both
pencils. Otherwise $P^8(U+V)=4$ by projection. One of the
pencils, say $V$, has $P^8V/4\ge1/2$. Choose a prime component
of a general $U$-fiber mapping onto a prime divisor $E$ on $Z$.
The corresponding bound gives $10S_L(E)\ge1+P^8V/4\ge3/2$.
Both cases contradict K-semistability.
\end{proof}

\begin{proposition}\label{prop:dec1-eight-first-reduction}
In Proposition~\ref{prop:dec1-eight-adjoint}, the sectional genus
and degree of the eight-dimensional hyperplane image $W\subset
\mathbf P^{10}$ are either $(1,4)$ or $(3,6)$. In the first case
$W$ is a complete intersection of two quadrics, every quadric in
their pencil having rank at least five. In the second case $W$ is
covered by a two-parameter family of linear $\mathbf P^6$'s.
\end{proposition}
\begin{proof}
By the slope estimate $K_B+6H$ is not pseudoeffective: the linear
functional given by pullback followed by intersection with the nef
class $P^8$ is nonnegative on the pseudoeffective cone, and is
strictly negative here. Apply M.~Andreatta, \emph{Minimal model
program with scaling and adjunction theory}, Theorem~5.7 and
Proposition~3.5(A),(B), \href{https://arxiv.org/abs/1107.4878}
{arXiv:1107.4878}, to the smooth quasi-polarized $(B,H)$. Its
first reduction $(B'',H'')$ is reached by a polarization-trivial
minimal model program followed by contractions of exceptional
divisors to finitely many points; at those last steps the new
polarization pulls back to the old one plus an effective exceptional
divisor, by his Propositions~4.2 and~5.4. Pushing down effective
representatives gives an injection
\[
       H^0(B,K_B+7H)\ \lhook\joinrel\longrightarrow\
                         H^0(B'',K_{B''}+7H'').
\]
Indeed, at a divisorial contraction a section remains regular in
codimension one on the target, and a flip is an isomorphism in
codimension one. This also identifies ratios of these sections in
the common function field. The argument works with reflexive
canonical sheaves on the terminal models.

In dimension eight the first-reduction list comprises projective
space, a quadric, the generalized cone over the Veronese surface,
a projective-space bundle over a curve, a del Pezzo variety with
$-K_{B''}\sim_{\mathbb Q}7H''$, a quadric fibration over a curve,
or a $\mathbf P^6$-bundle over a normal surface. The three
low-dimensional exceptional polarizations in part~(B) do not
occur. Projective space and the quadric have $K_{B''}+7H''$
antiample. The generalized Veronese cone has sectional genus zero
by Andreatta's Corollary~6.1; equivalently its general curve section
is the rational conic hyperplane section of the Veronese surface.
Its adjoint $K_{B''}+7H''$ has negative $H''^7$-degree by the
sectional-genus formula, so has no sections. In the del Pezzo
case the rank-one sheaf of $K_{B''}+7H''$ has at most one section:
any effective representative of a rationally trivial divisor is
zero. The projective-space bundle over a curve has no adjoint
sections on a general fiber $\mathbf P^7$. For the quadric
fibration the restriction of the adjoint to the general quadric
sevenfold is trivial, and for the surface bundle its restriction
to the general $\mathbf P^6$ fiber is trivial. Thus section ratios
are constant on general fibers and the dimension of the adjoint
image is at most one and two respectively. In conjunction with
Proposition~\ref{prop:dec1-eight-adjoint}, the only possible moving
adjoint is a net filling $\mathbf P^2$, giving $g_C=3$ and the
surface-bundle alternative. The remaining case is $g_C=1$.

The sectional curve has complete hyperplane system of dimension
three. For $g_C=1$ this gives $e=4$. When $g_C=3$, Riemann--Roch
gives $e\le6$; if $H|_C$ were special, Clifford's inequality
$4\le1+e/2$ would force $e\ge6$, whereas the nonzero $h^1$
would give $e<6$. Thus $e=6$. On a general $\mathbf P^6$ fiber
of the first reduction the original linear series restricts to
a subsystem of $\mathcal O(1)$: on a common resolution the
polarizations for the zero-reduction agree, and for the remaining
contractions the general fiber avoids their finitely many point
images. Since the map to $W$ is birational on a dense open,
the subsystem on a general fiber maps it birationally onto its
six-dimensional image. A linear system in $\mathcal O_{\mathbf P^6}(1)$
with that image dimension has full rank seven. This produces the
required family of linear spaces.

For $g_C=1$, the elliptic normal sectional quartic is a complete
intersection of two quadrics. Both quadrics lift to equations of
$W$ along the flag. Indeed every embedded positive-dimensional
member $W_j$ of the flag is linearly normal: the pullback of
$H^0(W_j,\mathcal O(1))$ injects into the complete hyperplane
series on its smooth model $B_j$, and all sections there are
ambient linear forms by the irregularity ladder. Therefore
$H^1(\mathcal I_{W_j}(1))=0$, exactly the obstruction to lifting
an equation of degree two from a hyperplane section. General
embedded hyperplane sections here are integral by Bertini, and
the smooth-model ladder maps onto them. The two lifted quadrics
have no common factor since their restrictions cut the elliptic
quartic in codimension two. Their complete intersection contains
the integral $W$ of degree four, hence equals it scheme-theoretically
by the degree and unmixedness of a complete intersection.

Finally suppose a quadric relation has rank at most four. Rank one
or two contradicts integrality and linear nondegeneracy. A rank
three relation gives on $Z$ a rational map to a conic and a
hyperplane section with twice a general ruling prime, contrary to
$7S_L\le A_Z$. A rank four relation gives a rational map to the
Segre quadric surface. Its four coordinate sections are linearly
independent. The images under each ruling projection are dominant,
since otherwise the image in $\mathbf P^3$ would be linearly
degenerate. The product structure gives two moving pencil classes
$U,V$ with $P-U-V$ effective up to rational linear equivalence.
The two-pencil contradiction in Proposition~\ref{prop:dec1-eight-adjoint}
applies. This proves the rank bound for every quadric in the pencil.
\end{proof}

\begin{proposition}\label{prop:dec1-eight-vertical-action}
In the $(g_C,e)=(1,4)$ case of
Proposition~\ref{prop:dec1-eight-first-reduction}, a central
one-dimensional torus in $\operatorname{Aut}^0(Z)$ fixes the
hyperplane image pointwise. A smooth connected component of a
general resolved hyperplane fiber is a rational curve.
\end{proposition}
\begin{proof}
The complete intersection of two quadrics $W$ is normal. Indeed
every quadric of the pencil has rank at least five, and the locus
of linearly dependent gradients in $\mathbf P^{10}$ lies in the
union of the projective vertices of the pencil. Their dimension
is at most five each, giving a union of dimension at most six.
Serre's criterion now applies since $\dim W=8$ and the complete
intersection is Cohen--Macaulay. We have $K_W\sim-7H$.

In fact $W$ has canonical singularities. On a smooth resolution
$\pi:B\to W$ the divisor $K_B+7H_B$ has a nonzero section by
the genus-one adjunction calculation. With canonical divisors
chosen compatibly, this produces an effective divisor linearly
equivalent to $K_B-\pi^*K_W$. Its pushforward is effective and
linearly trivial on the projective $W$, hence zero. The two
exceptional divisors must coincide: their difference is principal
and the inducing rational function has divisor zero on $W$,
hence is constant. Thus $K_B-\pi^*K_W\ge0$.

On a big smooth open $U\subset Z$, the limiting hyperplane map
$h:U\to W$ is regular and identifies $h^*H$ with $L|_U$. The
relative canonical calculation yields a generically nonzero map
\[
                    h^*\omega_W\longrightarrow\Omega_U^8.
\]
For example lift to $B$ on a smaller big open, use the regular
inclusion $\pi^*\omega_W\to\omega_B$ just proved, then take the
ordinary wedge of the differential. The map to the proper $B$
is defined in codimension one by the valuative criterion, so this
construction loses no prime divisor of $Z$. As
$h^*\omega_W\simeq\omega_Z|_U$, the natural identification
$\Omega^8_U\otimes\omega_U^{-1}\simeq T_U$ gives a nonzero
global reflexive vector field $v\in H^0(Z,T_Z)$ tangent to
the general hyperplane fiber.

The line spanned by $v$ is preserved by $\operatorname{Aut}^0(Z)$.
Indeed this group preserves the class of $L$ (as in the product
argument of Proposition~\ref{prop:dec1-factorial-quadric-image}),
and hence acts on the degree-one hyperplane map and its image.
Pullback of top forms is functorial on the generic smooth locus.
The chosen identification of $h^*\omega_W$ with $\omega_U$ is
unique up to scalar, since invertible regular functions on a big
open of the projective normal $Z$ are constant. Likewise the
identification $\omega_W\simeq\mathcal O_W(-7)$ is unique up to
scalar. These observations give invariance of the line.

The automorphism group of a K-polystable Fano has reductive
identity component. An invariant line in its Lie algebra is a
one-dimensional ideal, hence lies in the center. Its elements
belong to the Lie algebra of the central torus. The algebraic
subtorus generated by $v$ fixes the hyperplane image pointwise,
since the Lie algebra element is vertical and the kernel of the
projective action on $W$ is an algebraic subgroup. Its action has
generic orbit dimension one: it acts faithfully on $Z$, and the
general hyperplane fibers are curves. A faithful torus acting on
an integral variety has generic stabilizer of dimension zero
(diagonalize on an equivariant projective embedding). Thus this
subtorus has dimension one. Its general orbit is open in its
one-dimensional fiber component. The smooth proper model of the
orbit closure is $\mathbf P^1$, giving the assertion.
\end{proof}

\begin{corollary}\label{cor:dec1-optimal-lower-images}
Suppose that $X$ is a smooth quartic ninefold with
$\delta(X,-K_X)\le1$ and that $(Z,L)$ is an optimal destabilizing
special fiber. The image of $|L|$ has dimension six, seven or eight.
If $\delta(X,-K_X)>90/91$, its dimension is seven or eight.
\end{corollary}
\begin{proof}
Write $\beta=\delta(Z,L)=7\delta(X,-K_X)>20/3$ as in
Corollary~\ref{cor:nov30-optimal-image-ranks}. Dimension nine is
excluded following Proposition~\ref{prop:dec1-rank-four-pencils},
and dimension at most three is excluded following
Proposition~\ref{prop:nov30-optimal-minimal-image}. Consider an
image of dimension $a=4,5$, or $a=6$ with $\beta>90/13$. Put
$A=K_B+(a-1)H$ on a smooth resolution $B$ of the image. The
determinant slope inequality for an optimal fiber gives upper
bound $a-1-a\beta/9$ for its normalized mixed degree, and lower
bound $\beta/18$ for the moving class of any pencil. The former is
strictly smaller than the latter: this asks for
$\beta>18(a-1)/(2a+1)$, respectively $6$, $72/11$, and $90/13$.
Hence $h^0(B,A)\le1$. Just as in
Proposition~\ref{prop:dec1-lower-five-six}, the adjunction ladder
gives sectional genus at most one. The minimal-degree case is
excluded even for a subsystem. Indeed a doubled moving divisor
violates $\beta>N/2$ for $N=10$; the remaining Segre scroll has
a rank-four quadratic relation, violating
Proposition~\ref{prop:dec1-rank-four-pencils}. The classification
of minimal degree gives precisely these alternatives.

The same inner-projection and adjunction argument now shows that the
vertex-reduced base of the image is a smooth del Pezzo variety of
degree $12-a$ in its complete embedding. In degree eight the only
possibility in dimension at least three is $(\mathbf P^3,2H)$ by
Kuznetsov--Prokhorov, Theorem~1.2(viii). The smooth surfaces are
the plane blown up at one point and $\mathbf P^1\times\mathbf P^1$.
Their anticanonical classes contain twice a general line of the
plane, respectively twice a member of a ruling, with effective
residual divisor. In all three cases a hyperplane pulls back to
at least twice a prime divisor on $Z$, again contradicting
$\beta>N/2$. In degree seven the threefold is the blowup of
$\mathbf P^3$ at a point and its polarization
$2h-E=2(h-E)+E$ contains a doubled moving divisor.
On the degree-seven surface, the anticanonical plane model
$3h-E_1-E_2$ contains twice a general line with residual the
line through the two blown-up points. On the degree-six surface
$3h-E_1-E_2-E_3$ contains twice a general member of
$|h-E_1|$, with residual $(h-E_2-E_3)+E_1$ effective.
The higher-dimensional degree-six cases have a rank-four
relation. For the Segre fourfold and the triple product this is a
$2$-by-$2$ minor; on a smooth hyperplane section of
$\mathbf P^2\times\mathbf P^2$ one can take such a minor using
four coordinates independent modulo the hyperplane. Indeed its
defining bilinear form has rank three by smoothness, whereas any
combination of the four entries in a $2$-by-$2$ submatrix has
rank at most two. These relations also persist on the cone over
the respective varieties. This excludes all the listed dimensions.
\end{proof}

\begin{corollary}\label{cor:dec1-ninefold-alternatives}
Let \(X\subset\mathbf P^{10}\) be a smooth quartic ninefold.
If a K-polystable replacement over its GIT point differs from \(X\),
the image of its hyperplane system has dimension seven or eight.
In dimension seven it is the projective cone over
\(\operatorname{Gr}(2,5)\subset\mathbf P^9\). In dimension eight
it is either a complete intersection of two quadrics, each nonzero
member of their pencil having rank at least five, or a degree-six
variety of sectional genus three covered by a two-parameter family
of linear \(\mathbf P^6\)'s.
\end{corollary}
\begin{proof}
Write \(W\) for the image. If \(\dim W=9\), the moving intersection
number bounds the product of its degree and the mapping degree
by four. Nondegeneracy gives \(\deg W\ge2\). In degree four
the Hilbert series identifies the full ring with the hypersurface
ring. In degree three Proposition~\ref{prop:nov30-linear-image-component}
excludes \(n=9\). In degree two, a generically two-to-one map gives
equality of the moving and total volumes. The proportionality
and base-point argument in Corollary~\ref{cor:nov30-optimal-image-ranks}
then gives a Cartier globally generated root. Both embedded
and globally generated cases recover \(X\) by
Theorem~\ref{thm:nov30-finite-cm-model}. For a birational quadric
the rank is at least three. A rank-three equation has a hyperplane
section twice a moving ruling divisor, giving \(10S_L(E)\ge2\)
for a prime \(E\) over a general ruling, whereas \(7S_L(E)\le1\).
Rank four is excluded by Proposition~\ref{prop:dec1-rank-four-pencils}
and ranks at least five by Proposition~\ref{prop:dec1-factorial-quadric-image}.

If \(\dim W<9\), the slope and minimal-degree argument following
Proposition~\ref{prop:nov30-optimal-minimal-image} gives \(\dim W\ge5\).
Proposition~\ref{prop:dec1-lower-five-six} eliminates five and six.
Corollary~\ref{cor:dec1-only-grassmannian} and
Proposition~\ref{prop:dec1-eight-first-reduction} give the remaining assertions.
\end{proof}
An optimal special test configuration of \(X\) with
\(\delta(X,-K_X)\le1\) also has image of dimension at most eight,
by Propositions~\ref{prop:dec1-optimal-factorial-quadric}
and~\ref{prop:dec1-rank-four-pencils}. A distinct K-polystable
replacement is a different moduli-theoretic object. The alternatives
in the corollary are necessary conditions; their realizability over
a smooth GIT point remains open.

\begin{proposition}\label{prop:dec1-stein-hyperplane-fibers}
In each of the three cases of
Corollary~\ref{cor:dec1-ninefold-alternatives}, the general fiber of
a resolution of the hyperplane map $Z\dashrightarrow W$ is connected.
In the intersection-of-two-quadrics case, the central torus
$G\simeq\mathbf G_m$ of Proposition~\ref{prop:dec1-eight-vertical-action}
satisfies $\mathbb C(Z)^G=\mathbb C(W)$.
\end{proposition}
\begin{proof}
Let $a=\dim W$ and let $B$ resolve the normalization of $W$ in the
relative algebraic closure of $\mathbb C(W)$ inside $\mathbb C(Z)$.
Write $p$ for the degree of this finite field extension and $H$ for
the pullback of the hyperplane class to $B$. On a common resolution
$\pi:Y\to Z$, $h:Y\to B$ set $P=\pi^*L$ and $J=h^*H$. We have
$H^a=p\deg W$. The eleven hyperplane sections on $Z$ have no common
prime divisor by Proposition~\ref{prop:dec1-lower-seven}. On the big
regular open where $Z\dashrightarrow B$ is defined their line bundle
therefore equals $h^*H$. Thus every section of $H$ gives a section
of $L$ by extension across codimension two. Pullback of the ambient
coordinates gives the opposite dimension bound, so $h^0(B,H)=11$.

Cut $B$ successively with $|H|$ down to a smooth curve $C$. The
class is big, nef and free. The sections are connected by vanishing
and integral by Bertini. The irregularity vanishing and the
adjunction restriction argument of Proposition~\ref{prop:dec1-lower-seven}
apply to this ladder: $B$ is dominated by a resolution of the klt
Fano $Z$ and hence is rationally connected. In particular
\[
   h^0(C,H|_C)=12-a,\qquad
   g(C)\le h^0(B,A),\qquad A=K_B+(a-1)H.
\]
Only the surjectivity of adjoint restriction is needed here.

The first-jet proof of Lemma~\ref{lem:dec1-hyperplane-atiyah-slope}
applies verbatim to $Z\dashrightarrow B$: its inputs are dominance
and the equality $h^*H=L$ on a big open. Therefore
\[
              P^8h^*A/4\le c_a:=a-1-7(a+1)/10.
\]
Both $c_7=2/5$ and $c_8=7/10$ are strictly smaller than $7/9$.
When $|A|$ moves, let $S$ resolve its image, of dimension $t$, and
$M$ the hyperplane pullback to $S$. Its normalized mixed degree
on $Y$ is at most $c_a$, by subtracting the fixed divisor on a
common resolution. The ordinary tangent-sheaf slope for
$Z\dashrightarrow S$ shows $K_S+tM$ cannot be effective: otherwise
that moving degree would be at least $7/9$. The restriction and
vanishing argument in Proposition~\ref{prop:dec1-eight-adjoint}
now gives a linear adjoint image, so $h^0(B,A)=t+1$. In dimension
seven the slope for a general projection to $\mathbf P^t$ gives
$7t/[9(t+1)]\le2/5$, hence $t\le1$. In dimension eight the same
first-jet bound gives $P^8h^*(K_B+6H)/4\le6-63/10<0$.
Consequently $K_B+6H$ is not pseudoeffective, and Andreatta's
first-reduction list for this smooth quasi-polarized pair shows
$t\le2$, exactly as in the first two paragraphs of
Proposition~\ref{prop:dec1-eight-first-reduction}. Those steps
bound adjoint sections on the first reduction by restriction to
its general Mori fibers and apply to a big, nef polarization.
Thus $g(C)\le2$ for $a=7$ and $g(C)\le3$ for $a=8$, also when
$h^0(B,A)\le1$.

Riemann--Roch gives
\[
   p\deg W=\deg H|_C\le h^0(C,H|_C)+g(C)-1\le6.
\]
The three image degrees are $5,4,6$, forcing $p=1$. In
characteristic zero the corresponding function field extension
is regular, and generic smoothness and Stein factorization give
connected general resolved fibers. In the case with the central
torus, the general orbit has dimension one and the quotient
function field has transcendence degree eight. It contains
$\mathbb C(W)$, so relative algebraic closedness proves equality.
For the corresponding section-ring statement let $T$ be the ring
generated by $H^0(Z,L)$ inside $S=R(Z,L)$. Choose the torus lift
(after an isogeny if necessary) with its common weight on $S_1$
normalized to zero. A nonzero $x\in S_1$ trivializes the grading:
$\operatorname{Frac}(S)=\mathbb C(Z)(x)$ and
$\operatorname{Frac}(T)=\mathbb C(W)(x)$. Here $x$ is transcendental
over the respective degree-zero fields, and ampleness of $L$ gives
the first degree-zero field by a very ample Cartier multiple. As
$x$ is invariant, $\operatorname{Frac}(S)^G=\mathbb C(Z)^G(x)$:
in a reduced fraction in $x$ with monic denominator, invariance
forces every coefficient to be invariant. Hence
$\operatorname{Frac}(S)^G=\operatorname{Frac}(T)$.
\end{proof}

\begin{proposition}\label{prop:dec1-factorial-ci-image}
In the quartic-ninefold replacement problem, suppose the hyperplane
image is a complete intersection $W$ of two quadrics and
that its homogeneous coordinate ring $T$ is factorial. Then the
full section ring is a complete intersection
\[
       R(Z,L)=T[y_+,y_-]/(y_+y_--f_4),\qquad \deg y_\pm=2,
              \quad 0\ne f_4\in T_4.
\]
Thus such an image is excluded for a distinct K-polystable
replacement over a smooth quartic GIT point by
Theorem~\ref{thm:nov30-all-ci-fillings}. In particular the image
cannot be a smooth complete intersection of two quadrics.
\end{proposition}
\begin{proof}
Let $R=R(Z,L)$, $\mathcal X=\operatorname{Spec}R$ and
$\mathcal W=\operatorname{Spec}T$. Write $G$ for the vertical torus
of Proposition~\ref{prop:dec1-eight-vertical-action}. Use the lift
fixing every degree-one section, after an isogeny and an ordinary
degree character twist. By Proposition~\ref{prop:dec1-stein-hyperplane-fibers},
$\operatorname{Frac}(R)^G=\operatorname{Frac}(T)$.
The projective vector field $v$ of $G$ corresponds to the pullback
of top differential forms from $W$ under the identification
$h^*\omega_W\simeq\omega_Z$ on a big smooth open of $Z$. Since
$v$ is fixed and $G$ acts identically on $W$, this identification
is $G$-equivariant for the natural linearizations, as may be
checked where the map is smooth and the differential has full rank.

No prime divisor on $Z$ is pointwise fixed. At a general smooth
point of any proposed fixed divisor the map $h:Z\dashrightarrow W$
is regular and the identification just obtained gives trivial
character on the canonical line. But the tangent representation
has eight trivial characters tangent to the divisor, and its one
normal character must be nonzero. Indeed a torus acting trivially
on the tangent space at a fixed smooth point acts trivially on
the formal local ring (choose equivariant formal parameters),
and hence on the variety. The determinant character gives a
contradiction. Thus $v$ is nonzero along the generic point of
every prime divisor.

Every irreducible component of $\operatorname{Sing}W$ is contained
in a linear $\mathbf P^7$. To see this let $s\ge5$ be the generic
rank of the symmetric pencil $A:\mathcal O^{11}_{\mathbf P^1}
\to\mathcal O_{\mathbf P^1}(1)^{11}$. All its members have rank
at least five by Proposition~\ref{prop:dec1-eight-first-reduction}.
Its kernel bundle splits as $K=\bigoplus_{i=1}^{11-s}\mathcal
O(-a_i)$, $a_i\ge0$; its quotient $E$ is locally free since
the image has no torsion. Symmetry gives a generically invertible
map $E\to E^*(1)$. The determinant is a nonzero section of a line
bundle of degree $s-2\sum a_i$, so $\sum a_i\le s/2$. The columns
realizing $K\to\mathcal O^{11}$ are vectors of binary forms of
degree $a_i$. The span of the general vertex spaces therefore has
vector dimension at most $11-s+\sum a_i\le8$. At the finitely
many rank drops the projective vertex has dimension at most five.
The Jacobian criterion puts the singular locus of $W$ in the
union of these linear spaces. The image of a prime divisor of
$Z$ cannot be contained in any $\mathbf P^7$, as that would give
three independent linear sections vanishing along it, contrary
to Proposition~\ref{prop:dec1-lower-seven}.

No prime divisor of $Z$ maps into a subset of $W$ of codimension
at least two. Indeed its image meets the smooth locus by the
preceding paragraph. At a general point of the prime divisor the
differential of $h$ would have rank at most seven, since the image
of the tangent space of the divisor has dimension at most six.
The pullback of a generator of $\omega_W$ would therefore vanish
along this divisor, contrary to the nonvanishing of $v$ there.

The same absence of such contractions holds for the affine cone
map $\mathcal X\to\mathcal W$. Indeed the locus where all the
degree-one sections vanish has codimension at least two in
$\mathcal X$, by the absence of fixed primes in the full linear
system. To see this directly, components of the zero locus of
this homogeneous ideal are stable under the ordinary grading
torus; a prime component of height one would project under
$\operatorname{Proj}$ to a fixed prime on $Z$. Over a chart with
some $x\in R_1$ invertible, the cone map is the product of the
projective map with the identity of $\mathbf G_m$ in the common
coordinate $x$. A divisor in this product either dominates the
projective source, or is the full product over a divisor there.
In the first case its image dominates the projective target, so
has codimension at most one on the target cone; the second case
follows from the projective assertion. The same charts show that
no cone divisor is fixed by $G$, using that $x$ is $G$-invariant.

Every base prime divisor is dominated by a prime on $\mathcal X$:
its equation $f$ in the factorial $T$ remains a nonunit of the
connected positively graded domain $R$, since the only units are
the nonzero scalars. Any component of $(f=0)$ has height one and
must dominate its target divisor by the preceding paragraph.
Consequently $R^G=T$. Indeed an invariant element of $R$ lies in
$\operatorname{Frac}(T)$, and a pole at a base prime would remain
a pole along the dominating divisor. Normality of $T$ concludes.

Put $K=\operatorname{Frac}(T)$ and $R_K=R\otimes_T K$. This is
a normal affine $K$-curve with a split torus action, degree-zero
character space $K$, and fraction field of transcendence degree
one over $K$. The occurring character spaces are one-dimensional
over $K$ since ratios are invariant. Their weights generate a
nonzero subgroup of $\mathbb Z$. Both signs occur: if the weights
were one-sided, their positive-degree ideal, after choosing sign,
would define a $K$-rational fixed point of this affine curve.
The closure in $\mathcal X$ would be a horizontal fixed divisor
by the dimension formula, which was excluded. In the presence of
both signs the weight semigroup is a group, since a positive and
a negative integer admit a positive integral relation. The generic
curve is consequently $\operatorname{Spec} K[z,z^{-1}]$ for an
eigengenerator of the group of occurring weights. Thus every
nonzero rational eigenfunction is a unit on the generic affine
fiber.

For each occurring character $w$ choose an ordinarily homogeneous
eigenfunction $h_w\in R$ of that weight. Every rational eigenfunction
of weight $w$ is $h_w f$ for $f\in K$. Normality of $R$ tests its
membership in $R$ by nonnegative order along all prime divisors.
The horizontal divisors give no conditions by the generic-fiber
description. Every remaining prime dominates a unique prime $D$
of $T$, and its order on $f$ is a positive integral multiple of
$\operatorname{ord}_D(f)$. Thus membership is equivalent to one
integral lower bound on $\operatorname{ord}_D(f)$ for each $D$,
taking the maximum of the rounded bounds at the primes above it.
There are finitely many such primes above each $D$, as they are
components of the pullback of its equation. These lower bounds
are zero away from the finite image of the divisor support of
$h_w$. Therefore the weight-$w$ module is a copy of a divisorial
fractional ideal of $T$. Factoriality makes it free of rank one,
with an ordinarily homogeneous generator: the divisor prescribing
the bounds is stable under the ordinary grading torus, so a
principal generator is homogeneous (its translates differ by units).

Write $d_w$ for the ordinary degree of this generator. Dividing
the Hilbert series gives
\[
         \sum_w t^{d_w}=H_R(t)/H_T(t)
                       =\frac{1+t^2}{1-t^2}.
\]
Hence exactly $2j+1$ different weights occur by ordinary degree
$2j$. The two generators in ordinary degree two have distinct
nonzero weights $a,b$. If they have the same sign, change signs
so $0<a<b$. The products of the degree-two sections already
supply at least $2j+1$ weights by degree $2j$: use the distinct
weights $0,a,\ldots,ja$ and $(j-i)a+ib$ for $1\le i\le j$.
A negative-weight section exists and for $j$ sufficiently large
gives a further weight, a contradiction.

Write the opposite weights as $a>0$ and $-b<0$. The $2j+1$
weights $0,a,\ldots,ja,-b,\ldots,-jb$ exhaust the bound. Their
degree-four product of weight $a-b$ must already occur in that
list for $j=2$. If $a>b$, then $0<a-b<a$, and if $b>a$, then
$-b<a-b<0$, both impossible. Thus $a=b$ and $y_+y_-=f_4\in T_4$
is nonzero. The homomorphism from the presented algebra to $R$
is injective: reduce a polynomial by $y_+y_-=f_4$ to a sum of
a constant and pure powers of $y_+$ or $y_-$ with coefficients
in $T$. These terms have distinct characters and each nonzero
coefficient gives a nonzero term in the domain. The Hilbert series
of this subalgebra is $H_T(t)(1-t^4)/(1-t^2)^2=H_R(t)$,
so it is all of $R$. For a smooth $W$, the usual Lefschetz
theorem gives $\operatorname{Pic}(W)=\mathbb ZH$; the complete
intersection is projectively normal, and the cone class-group
formula gives factoriality of $T$.
\end{proof}

\begin{lemma}\label{lem:dec1-ci-no-projective-seven-plane}
An intersection-of-two-quadrics hyperplane image in the replacement
problem cannot contain a linear $\mathbf P^7$.
\end{lemma}
\begin{proof}
Suppose such a plane is given by $x_1=x_2=x_3=0$. Write the
quadrics as $x\cdot y$ and $x\cdot z$, where $x,y,z$ are triples
of linear forms. On the function field the equations imply
$y\times z=u x$ for a rational section $u$ of the hyperplane
class. Pull this identity back to $Z$. For every prime divisor
some $x_i$ is nonvanishing at its generic point by the
three-section criterion. Hence the corresponding identity
$(y\times z)_i=u x_i$ shows that $u$ is a regular Weil section
of $L$ there. Normality puts $u$ in $H^0(Z,L)$, so it is the
pullback of a linear form. The identities therefore hold in the
coordinate ring of $W$. Restricting them to the plane makes
$y\times z$ identically zero there. But $W$ is normal, so smooth
at the generic point of this prime divisor; at that point the
two defining gradients are independent and have coefficient rows
$y,z$ in the three conormal directions $dx_i$. Their cross
product is nonzero, giving a contradiction.
\end{proof}

\subsection{Excluding the intersection of two quadrics}
\begin{theorem}\label{thm:dec1-root-no-two-quadrics}
Let \(X\subset\mathbf P^{10}\) be a smooth quartic ninefold.
The hyperplane image of a distinct K-polystable replacement over
its GIT point cannot be an eight-dimensional intersection of two
quadrics.
\end{theorem}
\begin{proof}
Suppose \(W=(Q_0=Q_1=0)\) is such an image. Every pencil member
has rank at least five, \(W\) is canonical by
Proposition~\ref{prop:dec1-eight-vertical-action}, and
Lemma~\ref{lem:dec1-ci-no-projective-seven-plane} excludes a linear
\(\mathbf P^7\) in \(W\). Write \(s\) for the generic pencil rank.

First we show \(s\ge7\). Let \(k\) be the dimension of the common
vector-space kernel, and let \(W_0\) be the vertex-reduced
intersection of quadrics of dimension \(b=8-k\). Its normality
and canonical singularities follow from the affine-space bundle
over \(W_0\) on the complement of the vertex in \(W\). By adjunction
\(-K_{W_0}=(b-1)H_0\). The rank bound on the pencil ensures the
subsystem projecting to \(W_0\) has at least five independent
linear sections. It has no fixed prime on \(Z\) by
Proposition~\ref{prop:dec1-lower-seven}. On a smooth resolution of
\(W_0\) its canonical divisor dominates the pullback of \(K_{W_0}\).
The normalized mixed degree of the moving hyperplane pullback is one,
by projection to \(Z\) and absence of a fixed prime. Applying
Lemma~\ref{lem:dec1-hyperplane-atiyah-slope} therefore gives
\[
                   b-1\ge 7(b+1)/10 .
\]
In particular \(b\ge6\) and \(k\le2\).

Use the kernel bundle \(K=\bigoplus_i\mathcal O(-a_i)\) and
\(E=\mathcal O^{11}/K\) of Proposition~\ref{prop:dec1-factorial-ci-image}.
The common kernel has dimension \(h^0(K)\), so
\[
   k\ge 11-s-\sum_i a_i\ge 11-s-\lfloor s/2\rfloor.
\]
If \(s\le6\), then \(s=6,\ k=2,\ \deg E=3\). The nonzero
determinant of \(E\to E^*(1)\) is now a section of degree zero,
so the map is an isomorphism. As \(E\) is globally generated,
its splitting summands have degrees zero and one. Each occurs
three times. The annihilator of the span of the general pencil
kernels is \(H^0(E^*)\subset(\mathbf C^{11})^*\), of dimension
three. The span thus has dimension eight. For two different
general parameters \(t,t'\) and vectors \(v_t,v_{t'}\) in their
respective kernels, both bilinear forms \(Q_t,Q_{t'}\) vanish on
the pair. Consequently both \(Q_0,Q_1\) vanish on it. Closure in
the product of the kernel bundles gives the same conclusion at
equal parameters. The span is totally isotropic for the pencil,
giving the prohibited \(\mathbf P^7\).

It remains to prove factoriality of the homogeneous coordinate
ring \(T\) when \(s\ge7\). Take a smooth point \(p\in W\) and
write \(Q_i=x_0 l_i+q_i\), where \(x_0\) is the coordinate of
\(p\) and \(l_i,q_i\) use the other ten coordinates. The forms
\(l_0,l_1\) are independent by smoothness. On \(l_0\ne0\) use
\(t=l_1/l_0\), eliminate \(x_0\) by \(Q_0\), and retain the
single equation \(q_1-tq_0=0\) on \(l_1=t l_0\). For every
parameter \(t\), including the analogous fiber at infinity, the
resulting quadratic form is the form induced by the pencil member
on \(p^\perp/\langle p\rangle\). Its rank is the rank of that
member minus two since \(p\) is isotropic and nonradical for
every member. Thus each finite fiber in this chart is a nonempty
open affine-quadric-cone chart of rank at least three and is
geometrically integral. In particular \(t-c\) is prime in
\(T[1/l_0]\) for every \(c\in\mathbf C\).

Over \(\mathbf C(t)\) the induced form has rank at least five.
Its nondegenerate part has a nonzero isotropic vector by Tsen's
theorem. Split off a hyperbolic plane to write the equation as
\(uv+q'=0\) with \(q'\) of rank at least three. The class of
\(u\) is prime, and inverting it gives a Laurent polynomial ring.
Nagata's criterion gives factoriality of this quadric ring and
of its localization at \(l_0\). Applying the same criterion to
the primes \(t-c\) then gives factoriality of \(T[1/l_0]\).

Finally \(l_0=0\) is integral on \(l_1\ne0\) by the fiber-at-infinity
calculation. An additional divisorial component would have to lie
in the linear \(\mathbf P^8\) given by \(l_0=l_1=0\). The two
quadrics restricted to this space would have a common factor.
A linear factor would give a \(\mathbf P^7\) in \(W\). A common
quadratic factor, or an identically zero restriction, would give
a pencil member containing \(\mathbf P^8\), whose rank is at most
four. Thus there is no such component. The hypersurface section of
the complete intersection is Cohen--Macaulay, and its unique
component is generically reduced on \(l_1\ne0\). Hence \(l_0\)
is prime in \(T\). A final application of Nagata's criterion gives
factoriality of \(T\). Proposition~\ref{prop:dec1-factorial-ci-image}
now gives the contradiction.
\end{proof}

\begin{corollary}\label{cor:dec1-root-final-ninefold-images}
For a distinct K-polystable replacement over the GIT point of a
smooth quartic ninefold, the limiting hyperplane image is either
the seven-dimensional cone over \(\operatorname{Gr}(2,5)\), or
an eight-dimensional degree-six variety of sectional genus three
covered by a two-parameter family of linear \(\mathbf P^6\)'s.
In both cases the general resolved hyperplane fiber is connected.
\end{corollary}
\begin{proof}
Combine Corollary~\ref{cor:dec1-ninefold-alternatives},
Proposition~\ref{prop:dec1-stein-hyperplane-fibers} and
Theorem~\ref{thm:dec1-root-no-two-quadrics}.
\end{proof}

\subsection{Vertical vector fields and the adjoint system}
\begin{proposition}\label{prop:dec1-root-adjoint-vector-fields}
Let \(Z\) be a K-polystable \(n\)-dimensional \(\mathbf Q\)-Fano
variety with \(-K_Z\sim rL\) for an ample \(\mathbf Q\)-Cartier
Weil divisor. Suppose the complete system of \(L\) has no fixed
prime divisor and has image \(W\) of dimension \(n-1\).
On a smooth resolution \(B\) of \(W\), let \(H\) be the pullback
of the hyperplane bundle. Then
\[
                       h^0(B,K_B+rH)\le1.
\]
\end{proposition}
\begin{proof}
Let \(G=\operatorname{Aut}^0(Z)\). The connected automorphism
group preserves \(L\): the class of \(g^*L-L\) is \(r\)-torsion
since \(g^*K_Z\sim K_Z\). Equivalently one can check preservation
on the one-parameter subgroups generating the connected reductive
group, using invariance of the divisor class group on a product
with \(\mathbf A^1\) or \(\mathbf G_m\). Thus the map to \(W\)
is equivariant for its induced projective action. Choose a
functorial resolution of \(W\) in characteristic zero, so that
this action lifts to \(B\). The space of adjoint sections is the
same on any smooth resolution, by birational invariance of the
pushforward of the canonical sheaf and the projection formula.

On a smooth open \(U\subset Z\) with complement of codimension
at least two the rational map \(h:Z\dashrightarrow B\) is regular
and \(h^*H=L|_U\). This follows from properness of \(B\) and
the absence of a fixed prime. Pull back an adjoint section using
the top exterior power of the differential and this identification.
Since \(K_Z+rL\sim0\), a choice of trivialization gives
\[
 H^0(B,K_B+rH)\ \lhook\joinrel\longrightarrow\
 H^0(U,\Omega_U^{n-1}\otimes\omega_U^{-1})
       =H^0(Z,T_Z).
\]
Dominance and generic separability give injectivity. The image
is a \(G\)-stable subspace. Indeed pullback of differentials
and the hyperplane map are functorial, and a change in a chosen
line-bundle identification on a big open multiplies it by a
scalar, since \(Z\) is normal and projective.

Each resulting vector field is tangent to the general fiber of
\(h\), by contraction of a local volume form. The quotient
of any two nonzero adjoint sections is rational on \(B\).
Consequently on the generic smooth locus all the corresponding
vector fields are multiples of one vertical field by functions
pulled back from \(B\). Their pairwise Lie brackets vanish there
and hence globally. Their span is an abelian ideal of the reductive
Lie algebra of \(G\), so it lies in the center. It also lies in
the Lie algebra of the kernel of the action on \(W\).
The identity component of the intersection of the central torus
with that kernel is a torus whose dimension is at least the
dimension of the span. Its general orbit on \(Z\) has dimension
equal to its rank. Indeed its action is faithful, and in a sufficiently
ample equivariant projective embedding the generic stabilizer
is the common kernel of the differences of the occurring weights.
The orbit lies in a fiber of dimension one. This bounds the
adjoint space by one.
\end{proof}

\begin{corollary}\label{cor:dec1-root-grassmannian-only}
If a K-polystable replacement over the GIT point of a smooth
quartic ninefold differs from that ninefold, its hyperplane
image is the seven-dimensional projective cone with point vertex
over the Pl\"ucker embedding of \(\operatorname{Gr}(2,5)\).
The general resolved hyperplane fiber is a connected surface.
\end{corollary}
\begin{proof}
In the degree-six alternative of
Corollary~\ref{cor:dec1-root-final-ninefold-images} the adjoint
system on a resolution of the eightfold has
\(h^0(B,K_B+7H)=3\) by Proposition~\ref{prop:dec1-eight-adjoint}.
This contradicts Proposition~\ref{prop:dec1-root-adjoint-vector-fields}.
Connectedness in the remaining case is
Proposition~\ref{prop:dec1-stein-hyperplane-fibers}.
\end{proof}

\begin{proposition}\label{prop:dec1-grassmannian-fibers}
Suppose that the limiting hyperplane map of a K-polystable quartic
ninefold replacement has image the projective cone $W$ over
$G=\operatorname{Gr}(2,5)$. On a simultaneous smooth resolution
of the maps to $W$ and $G$, the general fiber over $W$ is a rational
surface and the general fiber over $G$ is a rationally connected
threefold. For a general complex point of $G$, the induced fibration
of the resolved threefold over the ruling $\mathbf P^1$ has a
section. Each of its fibers has a component of multiplicity one.
\end{proposition}
\begin{proof}
Let $F_2\subset F_3\subset T_Z$ be the saturated foliations tangent
to the general fibers of the two rational maps, and use normalized
degree $\deg_L D=L^8D/L^9$. The vertex does not contain the image
of any prime divisor of $Z$ by the three-section argument. On
a big open of $Z$, the pullbacks of both hyperplane classes thus
equal $L$ and both targets are smooth there. Their canonical
classes are $K_G=-5H$ and $K_{W_{\rm reg}}=-6H$: away from the
vertex the cone is the total space of $\mathcal O_G(1)$, and the
canonical bundle formula for that line bundle gives the latter
identity. Saturate the horizontal cotangents. Their determinants
are the respective pullbacks of the canonical bundles plus effective
divisors $R_j$ of divisorial differential degeneracy. Taking the
determinant in the tangent sequence on the big smooth open gives
\[
                 c_1(F_2)=L+R_2,\qquad c_1(F_3)=2L+R_3.
\]
In particular $\deg_L F_2\ge1$ and $\deg_L F_3\ge2$, where the
degree of a sheaf means the degree of its determinant. The
first-jet inequality also gives $\deg_L R_2\le 2/5$ and
$\deg_L R_3\le1/10$: saturating the first-jet map adds the
same divisorial orders as saturating the differential. Indeed
with local splittings the jet map has identity on its quotient
line and the differential on its kernel, so their maximal-minor
ideals agree.

Choose a sufficiently high general complete-intersection curve
$\Gamma$ of a Cartier multiple of $L$, in the smooth locus and
where the foliations are regular. By slope semistability and the
Mehta--Ramanathan restriction theorem, $T_Z|_\Gamma$ is
semistable of normalized slope $7/9$, with curve degrees divided
by the same positive factor throughout. For a torsion-free
quotient of $F_j|_\Gamma$ of rank $q$, its kernel has degree at
most $(j-q)7/9$. Hence its normalized slope is at least
$(j-1-(j-q)7/9)/q$. Taking $1\le q\le j$ gives lower bounds
$2/9$ for $j=2$ and $4/9$ for $j=3$. In particular both
restrictions are ample. The foliation theorem of
Bogomolov--McQuillan (see Kebekus--Sol\'a Conde--Toma,
\emph{Rationally connected foliations after Bogomolov and
McQuillan}, Theorem~1, \href{https://doi.org/10.1090/S1056-3911-06-00435-8}{Contemp. Math. 397 (2006)})
gives rational connectedness of the normalizations of general leaf
closures. Equivalently use a sufficiently general $\Gamma$ and
openness of rational connectedness on the smooth family of
resolved fibers. Proposition~\ref{prop:dec1-stein-hyperplane-fibers}
gives connected general fibers over $W$. The same holds over $G$
since $\mathbb C(W)=\mathbb C(G)(t)$ and $\mathbb C(W)$ is
algebraically closed in $\mathbb C(Z)$. These fibers are birational
to the leaf closures. A smooth rationally connected complex surface
is rational by the classification of surfaces.

Resolve also the projection of the cone from its vertex, whose
standard resolution is a $\mathbf P^1$-bundle over $G$. For a
general point of $G$ this gives a smooth projective threefold
mapping onto its $\mathbf P^1$, with general fiber a smooth
rational surface. The section theorem of Graber--Harris--Starr,
\emph{Families of rationally connected varieties}, Theorem~1.1,
\href{https://arxiv.org/abs/math/0109220}{arXiv:math/0109220},
provides a rational section, which extends across the curve by
properness. Its intersection with each fiber divisor has total
degree one. Thus each fiber contains a component of multiplicity one.
\end{proof}

\begin{corollary}\label{cor:dec1-root-grassmannian-volume}
Let \(S\subset Z\) be the closure of a general fiber of the
rational map to the Grassmannian cone, and \(V\subset Z\) the
closure of a general fiber of the projection to the Grassmannian.
Their degrees with respect to the ample Weil divisor \(L\) satisfy
\[
                    0<L^2\cdot S\le 4/5,\qquad
                    0<L^3\cdot V\le 4/5 .
\]
In particular each general surface \(S\) meets the non-Cartier
locus of \(L\).
\end{corollary}
\begin{proof}
For each map take a resolution with \(P\) the pullback of \(L\)
and \(J\) the moving hyperplane pullback. Both divisors are nef,
and \(P-J\) is effective up to rational linear equivalence.
Write \(a=7\) for the cone and \(a=6\) for the Grassmannian.
Intersect \(a\) general hyperplanes on the target. In both cases
one obtains five reduced general points, because the Pl\"ucker
degree of \(\operatorname{Gr}(2,5)\) is five and coning preserves
degree. The general resolved fibers are integral and map
birationally to their closures on \(Z\). Generic flatness and
the projection formula therefore give
\[
  5(L^{9-a}\cdot S_a)=P^{9-a}J^a\le P^9=4 ,
\]
where \(S_7=S\) and \(S_6=V\). The inequality follows by
successive replacement of \(J\) by \(P\) against nef products,
and strict positivity follows from ampleness. If \(L\) were
Cartier near \(S\), the degree \(L^2\cdot S\) would be a
positive integer, contradicting the bound.
\end{proof}

\subsection{Quartic replacements in every dimension at least nine}
\begin{theorem}\label{thm:dec1-root-quartic-relative-dimension}
Let \(n\ge9\). On a distinct K-polystable replacement over the GIT
point of a smooth quartic \(X\subset\mathbf P^{n+1}\), the image of
the limiting hyperplane system has dimension at most \(n-2\).
\end{theorem}
\begin{proof}
Put \(N=n+1,\ r=n-2\), and use \(P\) for the pullback of the
limiting class \(L\) on a common resolution. First suppose the
image has dimension \(n\). Its degree times the generic mapping
degree is at most four. Degree four gives equality of the graded
rings by their Hilbert series, and degree three is excluded by
Corollary~\ref{cor:dec1-unconditional-degree-minus-one}. A double
map to a quadric has \(P^n=J^n=4\) for the moving class \(J\).
Proportionality, as in Corollary~\ref{cor:nov30-optimal-image-ranks},
gives \(P=J\) and a globally generated Cartier root. The embedded
and globally generated cases both recover the original hypersurface
by Theorem~\ref{thm:nov30-finite-cm-model}. For a birational
quadric, ranks at least five are excluded by
Proposition~\ref{prop:dec1-factorial-quadric-image}. Rank four
contradicts \(r>2N/3\) by
Proposition~\ref{prop:dec1-rank-four-pencils}. In rank three a
doubled moving ruling gives \(NS_L(E)\ge2\) for a prime on \(Z\),
contrary to \(rS_L(E)\le1\). Ranks at most two contradict
integrality and nondegeneracy of the image.

Now assume the image \(W\) has dimension \(a=n-1\). No prime
divisor of \(Z\) belongs to three independent hyperplane sections.
Indeed let \(M\) be the moving class of the resulting net on a
resolution. If it dominates the plane, the tangent-slope estimate
gives \(P^{n-1}M/4\ge2r/(3n)\). If the image is a curve, it is
rational and of degree at least two, giving the stronger
\(P^{n-1}M/4\ge r/n\). The two leading terms of the volume
integral on \([0,1]\) for a prime shared by the net would give
\[
     N/r\ \ge\ N S_L(E)\ \ge\ 1+2r/(3n).
\]
The strict reverse inequality is equivalent to
\(2(n-2)^2>9n\), valid for \(n\ge9\).

Resolve the normalization of \(W\) in the relative algebraic
closure of \(\mathbf C(W)\) inside \(\mathbf C(Z)\), and denote
the resulting smooth model by \(B\), its hyperplane pullback by
\(H\), and the finite field degree by \(p\). The pullback identity
\(h^*H=L\) holds on a big open of \(Z\), so \(h^0(B,H)=n+2\).
The group action lifts to this normalization since it preserves
the field inclusion and its relative algebraic closure. It lifts
also to a functorial resolution. Thus the proof of
Proposition~\ref{prop:dec1-root-adjoint-vector-fields} gives
\(h^0(B,K_B+rH)\le1\) on this model too. Cut by a general
hyperplane ladder down to the smooth sectional curve \(C\).
The irregularity vanishing and the adjunction restriction
surjections of Proposition~\ref{prop:dec1-lower-five-six} use
only a big, free polarization and rational connectedness of
\(B\), which holds by dominance from \(Z\). Because \(r=a-1\),
\[
          g(C)\le1,\qquad h^0(C,H|_C)=4.
\]
Thus \(p\deg W=\deg H|_C\) is three or four. The image has
codimension two and hence degree at least three, so \(p=1\).

The degree-three case is a variety of minimal degree and codimension
two. It is a rational normal scroll or a cone over one. In the
scroll matrix presentation a block of length at least two gives a
minor of rank three or four, while if every block has length one,
two blocks give a rank-four minor. These relations persist on
cones. Both ranks contradict the ruling bounds used above.
In degree four the sectional curve is an elliptic normal quartic.
The two quadrics lift through the ladder by linear normality,
exactly as in Proposition~\ref{prop:dec1-eight-first-reduction}.
Their complete intersection has degree four and is \(W\). Every
member of their pencil has rank at least five. Its singular locus
lies in the union of their projective vertices, of dimension at
most \(n-3\), so \(W\) is normal. It has \(K_W=-rH\).
Since \(g(C)=1\), the adjunction surjection gives a nonzero section
of \(K_B+rH\). By the exceptional-divisor argument of
Proposition~\ref{prop:dec1-eight-vertical-action} this implies
that \(W\) is canonical.

We detail the dimensional inputs in the factoriality argument.
A linear \(\mathbf P^{n-2}\) in \(W\) is a prime divisor cut out
by three linear equations in the ambient space. The cross product
argument of Lemma~\ref{lem:dec1-ci-no-projective-seven-plane}
excludes it, since no three independent sections share a divisor
on \(Z\). Let \(s\) be the generic pencil rank and \(k\) the
dimension of the common kernel. The split kernel bundle gives
\(k\ge n+2-s-\lfloor s/2\rfloor\). The vertex-free base has
dimension \(b=n-1-k\), and the jet inequality on its canonical
resolution gives \(b-1\ge r(b+1)/N\), that is
\(3b\ge2n-1\). These arguments use exactly the canonical
product chart and determinant of the kernel bundle in
Theorem~\ref{thm:dec1-root-no-two-quadrics}. If \(s\le6\)
then \(b\le6\), impossible for \(n\ge10\).
The borderline \(n=9\) was excluded in that theorem. Thus \(s\ge7\).

The pencil calculation at a smooth point in the same theorem now
proves factoriality of the homogeneous ring \(T\) in every such
dimension. Its finite fibers are quadric charts of rank at least
three, and its generic quadric over \(\mathbf C(t)\) has rank at
least five. At infinity a common linear factor in the two
restrictions would give the prohibited \(\mathbf P^{n-2}\).
Proportional restrictions would give a quadric containing
\(\mathbf P^{n-1}\), of rank at most four. These are precisely the
hypotheses for the two localizations and Nagata's criterion.

Finally, the vertical vector field associated to \(K_W=-rH\)
gives a central torus fixing \(W\), by
Proposition~\ref{prop:dec1-eight-vertical-action}. The field
degree \(p=1\) gives the invariant fraction field as in
Proposition~\ref{prop:dec1-stein-hyperplane-fibers}. For the
divisorial argument of Proposition~\ref{prop:dec1-factorial-ci-image},
the generic vertices span a vector space of dimension at most
\(n+2-s+\lfloor s/2\rfloor\le n-1\); each special vertex is
smaller than \(\mathbf P^{n-2}\). Thus a divisor on \(Z\) cannot
map into the singular locus. The nonzero vertical field has no
fixed divisor, by its canonical character, and the differential
excludes a divisor mapping to codimension at least two in the
smooth locus. The resulting weight modules are free rank one over
the factorial \(T\). The Hilbert-series quotient in all dimensions
is \((1-t^4)/(1-t^2)^2\). The count of positive and negative
weights in that proposition gives the complete-intersection
presentation \(T[y_+,y_-]/(y_+y_--f_4)\) with quadratic generators.
Theorem~\ref{thm:nov30-all-ci-fillings} excludes this presentation
over a distinct smooth quartic GIT point.
\end{proof}

\begin{theorem}\label{thm:dec1-quartic-high-dimensions}
Every smooth complex quartic hypersurface of dimension $n\ge10$
is uniformly K-stable and admits a K\"ahler--Einstein metric.
\end{theorem}
\begin{proof}
Write $N=n+1$ and $r=n-2$. Suppose a smooth quartic $X$ has a
distinct K-polystable replacement $(Z,L)$. Theorem~\ref{thm:dec1-root-quartic-relative-dimension}
gives $a:=\dim W\le n-2$ for the hyperplane image $W$. No
subsystem of hyperplanes on $Z$ can have a minimal-degree image
of codimension at least two in its span. Indeed, the proof of
Proposition~\ref{prop:nov30-optimal-minimal-image} uses just the
sections of that subsystem: a doubled ruling contradicts $r>N/2$,
and for a cone over $\mathbf P^1\times\mathbf P^t$, $t\ge2$,
the pullbacks of the two factor classes are dominated in sum by
$L$ on a resolution. The slope bound and first volume term at a
prime over a general ruling give $N S_L(E)\ge1+2r/(3n)>N/r$,
as in the three-section test of the preceding theorem. In what
follows we also use that test itself: three independent linear
sections cannot vanish along a common prime on $Z$.

Let $B$ resolve $W$, $H$ the pullback of its hyperplane, and let
$\pi:Y\to Z$, $g:Y\to B$ resolve the map, with $P=\pi^*L$.
Put $\mu(D)=P^{n-1}g^*D/4$ and use the same mixed-degree notation
after resolving further maps. Since the complete hyperplane
system has no fixed prime, $\mu(H)=1$. For $A=K_B+(a-1)H$ the
first-jet estimate of Lemma~\ref{lem:dec1-hyperplane-atiyah-slope}
gives
\[
   \mu(A)\le c:=a-1-\frac{r(a+1)}N
   =\frac{3a-2n+1}{N}\le\frac{n-5}{N}<\frac{r}{n}<1 .
\]
In particular $K_B+(a-2)H$ is not pseudoeffective, since
$\mu$ is nonnegative on pseudoeffective classes and $\mu(H)=1$.
The lower bound on dimension following
Proposition~\ref{prop:nov30-optimal-minimal-image} gives
$a\ge\lceil n/2\rceil\ge5$.
Apply Andreatta's Theorem~5.7 and Proposition~3.5 to $(B,H)$ as in
Proposition~\ref{prop:dec1-eight-first-reduction}. Pushforward
embeds the adjoint sections into the adjoint sections of the first
reduction and preserves their ratios. In dimensions $a\ge5$ the
three- and four-dimensional exceptional cases do not occur.
Projective space and the quadric have negative adjoint. For the
generalized Veronese cone its sectional genus zero gives negative
adjoint degree. In the del Pezzo case the adjoint is rationally
trivial. In a projective-space bundle over a curve its restriction
to the general fiber is negative. In a quadric fibration over
a curve and a projective-space bundle over a surface this
restriction is trivial. Consequently, if the adjoint system moves,
its image has dimension $t\le2$.

In fact this image is linear. Let $S$ resolve it and let $M$ be
the hyperplane pullback. The moving degree $\mu(M)$ on a higher
resolution is at most $c$. If $K_S+tM$ were effective, the
tangent-sheaf inequality for the map from $Z$ to $S$ would give
$\mu(M)\ge r/n$. This contradicts the displayed bound. For $t=1$
the curve is rational and this vanishing forces $\deg M=1$.
For $t=2$, take a general smooth curve $D\in|M|$.
Kawamata--Viehweg vanishing gives surjectivity from
$H^0(S,K_S+2M)$ to $H^0(D,K_D+M|_D)$, so the latter vanishes.
Riemann--Roch gives $g(D)=0$ and $M^2=1$. Nondegeneracy of
the image then implies it is a plane. The moving possibilities
are therefore a line or a plane.

By the vanishing ladder of Section~\ref{cubics:minimal}, a general
smooth sectional curve $C$ on $B$ has $h^0(C,H)=n+3-a$ and
$g(C)=h^0(B,A)$. For equality use the ladder of
Proposition~\ref{prop:dec1-lower-seven}: the kernels vanish since
$K_B+(a-2)H$ is not pseudoeffective, and this absence of sections
descends by vanishing, including $H^1(K_{B_3})=0$ on the
three-dimensional rung. Here $h^0(B,H)=n+2$ by pullback to $Z$,
and the irregularity is zero on the smooth ladder down to the
surface. The dimension-one image case already gives a
minimal-degree rational image by the same argument, so we may
assume $a\ge2$. We have $g(C)\le3$. Genus zero would give a
minimal-degree $W$ and is excluded. In genus two the moving
adjoint maps to a line and a general prime fiber divisor on $B$
has $H^{a-1}$-degree at most two, by restriction to $C$ where the
adjoint is $K_C$ of degree two. Use the rational Stein curve if
the pencil does not have connected fibers. The image of a general
such divisor on $W$ has dimension $a-1$, positive degree at most
two, and therefore spans at most $\mathbf P^a$. At least three
independent hyperplanes contain it since $a\le n-2$. Pulling back
a general such divisor on a dense smooth locus gives a prime
divisor on $Z$ in their common zero locus, a contradiction.

In genus three a general member of the moving adjoint net is an
integral divisor on the regular locus of the map to the plane,
by Bertini irreducibility. Its closure on $W$ has hyperplane
degree at most four, since the adjoint restricts to $K_C$. Its
dimension is $a-1$, so it spans at most $\mathbf P^{a+2}$ and
is contained in an ambient hyperplane because $a\le n-2$.
On $Z$ the pullback of a general line in the adjoint plane is
again generically reduced and integral on the regular locus.
Let $D_Z$ be its moving Weil class, with at least three
independent sections. The containing hyperplane supplies an
effective representative of $L-D_Z$. Any prime in it would be
common to the three independent products in $H^0(Z,L)$, which
the three-section estimate excludes. If the representative is
zero, then $D_Z\sim L$ and the moving degree is one by the
projection formula, contrary to $c<1$. This excludes genus three.

In genus one, $W$ has degree $e=n+3-a\ge5$. Remove its entire
linear vertex and call the base $W_0$. If $W_0$ had a singular
point, inner projection from it would give a minimal-degree image
of codimension $e-3\ge2$: the point is not a vertex, and
$e-\operatorname{mult}\le e-2$ equals the nondegeneracy lower
bound, as in Proposition~\ref{prop:dec1-lower-five-six}. The
projected image comes from a subsystem of $|L|$, already excluded.
Thus $W_0$ is smooth. Its general curve section is elliptic,
since $W_0$ is a general linear section complementary to the
vertex. A curve base would be elliptic and cannot be dominated
by $Z$. For $b=\dim W_0\ge2$ the same ladder makes $K_{W_0}+
(b-1)H_0$ effective. It has $H_0^{b-1}$-degree zero, so
$-K_{W_0}=(b-1)H_0$. The ladder gives
$h^0(W_0,H_0)=e+b-1$, hence its embedding is complete. The smooth
nonconical del Pezzo classification, Kuznetsov--Prokhorov,
Theorem~1.2(v)--(viii), \href{https://arxiv.org/abs/2206.01549}
{arXiv:2206.01549}, gives $b\le6$ in every degree $e\ge5$:
degree five yields linear sections of $\operatorname{Gr}(2,5)$,
and the higher-degree list has dimension at most four. For
$b=2$ the dimension bound is automatic. On the other hand the
subsystem projecting to $W_0$ has more than two independent
sections and hence no fixed prime. The first-jet bound gives
\[
                 b-1\ge\frac{r(b+1)}N,
                    \qquad 3b\ge2n-1>18,
\]
which contradicts $b\le6$. All distinct replacements are excluded.

Finally approach $X$ along a smooth pointed curve of quartics
with K-stable punctured family. Such a curve exists because the
Fermat quartic is K-stable \cite[Corollary~1.4]{ar2004.09413},
the smooth parameter space is irreducible, and uniform K-stability
is open. K-moduli properness supplies a K-polystable special fiber
after finite base change. Its limiting integral Weil root and its
graded dimensions are obtained by taking the closure of the
punctured hyperplane and restricting its reflexive powers, as in
the proof of Proposition~\ref{prop:nov30-factorial-quartic-replacement}
before its Cartier lifting step. More explicitly, the integral
special fiber is principal over the local base, so that
$-K_{\mathcal Z/B}\sim r\mathcal L$ relatively. The total space
is klt by inversion of adjunction. Its $\mathbf Q$-Cartier
divisorial sheaves $\mathcal O(j\mathcal L)$ are Cohen--Macaulay
and flat over the curve. Restriction to the normal special fiber
is Cohen--Macaulay and agrees with $\mathcal O_Z(jL)$ in
codimension one, hence everywhere. For $j\ge0$,
Kawamata--Viehweg vanishing and flatness give precisely the
quartic Hilbert function, as in the proof of
\cite[Theorem~3.1(1)--(3)]{ar2007.14320}. The preceding proof identifies that
fiber with $X$. Smooth quartic hypersurfaces have finite
automorphism group, so the equivalence between uniform K-stability
and K-stability for Fano manifolds and the Yau--Tian--Donaldson
theorem conclude the proof as in Section~\ref{sec:established-cases}.
\end{proof}

\paragraph{The cyclic suspension of a quartic ninefold.}
Let $Y=(f_4=0)\subset\mathbf P^{10}$ be smooth. The hypersurface
$X=(y^4+f_4=0)\subset\mathbf P^{11}$ is smooth by the Jacobian
criterion. Its projection $\pi:X\to\mathbf P^{10}$ is the quotient
by $y\mapsto\zeta_4 y$, and
$K_X=\pi^*(K_{\mathbf P^{10}}+\tfrac34Y)$.
Theorem~\ref{thm:dec1-quartic-high-dimensions} gives uniform
K-stability of the tenfold $X$. Hence the quotient pair
$(\mathbf P^{10},\tfrac34Y)$ is K-polystable by
\cite[Theorem~1.2]{ar2001.10557}. Liu and Zhu's branch-divisor
criterion in the other direction assumes that $Y$ is K-polystable
\cite[Theorem~5.1]{ar2001.10557}.

The coefficient needed for the usual degeneration to the cone over
$Y$ can be computed directly. For the valuation
$v=\operatorname{ord}_Y$ and $0\le a<1$, one has
\[
 A_{\mathbf P^{10},aY}(v)=1-a,\qquad
 S_{-(K_{\mathbf P^{10}}+aY)}(v)=\frac{11-4a}{44}.
\]
Indeed the polarization is $(11-4a)H$ and $Y\sim4H$, so the
second formula follows by integrating
$((11-4a)-4t)^{10}/(11-4a)^{10}$ up to $t=(11-4a)/4$.
Their difference vanishes at $a=33/40$ and equals $3/44$ at
$a=3/4$. At $33/40$, K-semistability of the ambient pair would
imply K-semistability of the branch: the special degeneration
associated with $v$ has zero Futaki invariant, hence a
K-semistable projective cone, whose delta invariant at this
coefficient is $\min\{1,\delta(Y)\}$ by
\cite[Lemmas~5.3--5.5]{ar2003.06839}. The cyclic quotient gives
the smaller coefficient $3/4$. Extending its semistable interval
to $33/40$ requires a further argument.

\section{Restrictions on a possible quartic ninefold replacement}

\subsection{A tangent pencil on the Grassmannian cone}
\begin{proposition}\label{prop:dec2-ninefold-tangent-pencil}
Let $(Z,L)$ be a quartic ninefold replacement as in
Corollary~\ref{cor:dec1-root-grassmannian-only}. Let $S$ be the
closure on $Z$ of a general surface fiber over the cone. Then
\[
       \frac{4}{75}<L^2 S\le\frac45,\qquad
       \operatorname{lct}(Z;|L|)>\frac{173}{50},\qquad
       \operatorname{codim}_Z\operatorname{Bs}|L|\ge4 .
\]
Here the threshold uses the minimum order of the degree-one
section divisors at each valuation. On an index cover it is the
threshold of the ideal of the lifted sections.
\end{proposition}
\begin{proof}
On a smooth resolution write $P$ for the pullback of $L$ and
$J$ for the moving hyperplane class. Put $a_i=P^{9-i}J^i/4$.
Over a general smooth point $w$ of the seven-dimensional cone
the resolved morphism is smooth with connected integral fiber
$F$. Its conormal bundle is trivial of rank seven. The
exceptional divisor $E$ of the blowup of $F$ has $A_Z(E)=7$,
because the generic point of $F$ is in the isomorphism locus
of the resolution. Use the same notation for further pullbacks.
We have $a_7=5P^2F/4$ by the degree of the cone.
Write $V(t)=\operatorname{vol}(P-tE)/4$. The layers of the
ideal of $F$ show that imposing order $\ell$ costs at most
$\binom{\ell+6}{7}h^0(F,mP|_F)$ conditions on $H^0(mP)$
for divisible $m$. Since $P|_F$ is big, nef and semiample,
$h^0(F,mP|_F)=m^2P^2F/2+O(m)$. Taking $\ell/m\to t$ gives
\[
                 V(t)\ge1-Ct^7,\qquad C=36a_7/5 .
\]
The valuative inequality for this replacement gives
$\int_0^\infty V(t)\,dt=S_L(E)\le1$.

We sharpen the integral on $1\le t\le2$. Hyperplanes containing
the embedded tangent space of the cone at $w$ form a projective
plane, since that tangent space is a projective seven-plane in
$\mathbf P^{10}$. Choose two independent members. Their lifted
sections have order at least two along $F$, and their ratio is
nonconstant by nondegeneracy. Choose a prime divisor $D$ on $Z$
where the ratio has a zero or pole. Such a divisor exists by
normality and projectivity. For the flag defining an Okounkov
body choose the strict transform of $D$ as first member and a
smooth flag at a general point of that transform. The first
coordinate of a section is then its order at $D$. We may use
divisible multiples throughout, with the same trivialization
for all section spaces.

Let $K$ be the body of the sections with order at least their
degree along $E$, inside the full body $\Delta$ of $P$. Normalize
Lebesgue measure by $\mu(\Delta)=1$. The volume theorem for
Okounkov bodies gives $\mu(K)=V(1)$, denoted by $v$.
Indeed multiplication by the canonical section of $mE$ identifies
$H^0(m(P-E))$ with this filtered subspace, up to translation of
its body. The first coordinate $x_1$ is nonnegative on $\Delta$.
Slicing at $x_1=t$ and integrating the upper-level volumes gives
\[
        \int_\Delta x_1\,d\mu=S_L(D)\le1/7 ,
\]
because $A_Z(D)=1$ and $-K_Z\sim7L$. We use the volume and
slicing results of \cite{ar0805.4559}, or equivalently their
section-counting interpretation for a flag starting with $D$.

The two tangent sections give a segment $I$ in the body
of the series with order at least twice the degree. For a
non-Cartier $L$ take equal Cartier powers and divide the
flag coordinates by that power. The displacement vector $b$
of the segment, oriented in the direction of increasing first
coordinate, has $b_1\ge1$ since the ratio has integral order
along $D$. Slice $K$ by lines parallel to $b$ and measure
length $\lambda$ in units of $b$. On the quotient space of
lines take the measure $d\sigma$ for which $d\mu=d\lambda\,d\sigma$.
Write $p$ for the measure of the projection of $K$. A chord of
length $\lambda$ contributes at least $\lambda^2/2$ to its
first moment, since $x_1$ starts nonnegative and increases by
at least one per unit length. Hence
\[
 v^2=\left(\int\lambda\,d\sigma\right)^2
 \le p\int\lambda^2\,d\sigma
 \le 2p\int_K x_1\,d\mu\le 2p/7 .
\]
For $0<u<1$ the filtered body at $1+u$ contains
$(1-u)K+uI$. This follows by multiplying a section with
normalized order at least one by products of the two sections
with normalized order at least two, in the proportions
$1-u$ and $u$. Rational proportions suffice by continuity
of divisor volume. Addition of the segment extends each chord
of $(1-u)K$ by $u$, giving
\[
     V(1+u)\ge (1-u)^9v+u(1-u)^8p
     \ge (1-u)^9v+\tfrac72u(1-u)^8v^2 .
\]
If $a_7\le1/15$, then $C\le12/25$ and $v\ge1-C>0$. Integrating
the two intervals would give
\[
 1\ge 1-\frac{C}{8}+\frac{1-C}{10}
               +\frac{7(1-C)^2}{180}
 \ge 1+\frac{283}{112500}.
\]
This proves $a_7>1/15$ and the stated lower bound for $L^2S$.
The upper bound is Corollary~\ref{cor:dec1-root-grassmannian-volume}.

For completeness put $\theta=(7/10)\sum_{i=0}^7a_i$ and let
$k$ be the base codimension of the complete eleven-section
system. The first-interval mixed-volume estimate, including
its local index-cover form, gives
$\operatorname{lct}(Z;|L|)\ge\theta$, $k\ge\theta$, and
$a_i=1$ for $i<k$ \cite[Lemmas~5.1, 5.3 and 5.4]{kePanWuCubics}.
There is no fixed divisor, so $a_1=1$. Since $a_7>2^{-6}$,
log-concavity between $a_1$ and $a_7$ gives
$\theta>(7/10)(2+\sum_{j=1}^6 2^{-j})>2$. Now $a_2=1$.
Since $(29/50)^5<1/15$, interpolation from $a_2$ gives
$\theta>(7/10)(3+\sum_{j=1}^5(29/50)^j)>3$. Thus $k\ge4$
and $a_3=1$. Finally $(101/200)^4<1/15$, so that
\[
   \theta>\frac7{10}\left(4+\sum_{j=1}^4(101/200)^j\right)
          >\frac{173}{50}.
\]
\end{proof}

\subsection{Cartierness through codimension five}

\begin{corollary}\label{cor:dec2-grass-local-index-five}
Suppose a quartic ninefold replacement has the Grassmannian-cone
image, and the mixed-intersection lower bound \(\theta\) for the
threshold of the degree-one system satisfies \(\theta>173/50\).
The hyperplane root \(L\) is Cartier in codimension four.
At the geometric generic point of each codimension-four base
component the ambient transverse germ is smooth. In suitable
regular parameters its degree-one ideal has one of the forms
\[
 (x_1,x_2,x_3,x_4),\qquad (x_1,x_2,x_3,x_4^2).
\]
\end{corollary}
\begin{proof}
Consider a codimension-four component of the non-Cartier locus
and a geometric generic transverse germ at a general point. The
local index is at most six by the volume estimate in
Proposition~\ref{prop:nov30-quartic-nine-index-cover}. It is
coprime to seven, so the cyclic root cover is the canonical cover.
As in that proposition, intersect the stratum by five general
ambient hyperplanes transverse at a general complex point. Their
pullbacks meet the strata of a log resolution transversely.
Adjunction preserves the log discrepancies on the cut, and the
restricted ideal has the same exceptional orders. Hence the
lifted ideal on the four-dimensional Gorenstein klt cover has
threshold at least \(\theta\), as in
\cite[Lemmas~5.3, 5.4]{kePanWuCubics}. Its degree-one generators
vanish at the marked point, since they have nontrivial character
and the full covering group fixes the point.

Three general combinations of these generators cut out a smooth
curve at the marked point. Indeed the base codimension on the
cover is at least four: a positive-dimensional component on a
transverse slice would spread out to a base component of codimension
at most three downstairs, whereas base codimension is at least
\(\theta>3\). On a log resolution the sum of the three fixed orders
is strictly smaller than the corresponding log discrepancy. The
three general moving members meet each other and every exceptional
stratum transversely by Bertini. Successive adjunction therefore
gives a klt Cartier complete intersection curve, which is regular.
Equivalently at each cut the coefficient-one strata of the resolved
pair are exactly the intersections of the moving members, and the
strict exceptional inequalities give klt singularities on their
common intersection. The embedding dimension of the transverse
cover is at most \(1+3=4\), so the cover is smooth.

Cut the given smoothing by the same transverse hyperplanes. Its
cyclic cover is the canonical cover of the sliced family, with
the character unchanged since all cutting functions are invariant.
The action on the smooth nearby fibers is free, as in the proof of
Proposition~\ref{prop:nov30-quartic-nine-index-cover}. At a smooth
point of the central cover its deformation is smooth and the
base parameter is invariant. Averaging lifts of an eigenbasis
of the relative cotangent space gives formal parameters on which
the finite cyclic group acts diagonally, together with the invariant
base parameter. Setting the relative parameters to zero gives a
fixed formal section, contrary to generic freeness. Thus there
is no codimension-four non-Cartier component. Lower codimensions
are covered by Proposition~\ref{prop:nov30-quartic-nine-index-cover}.

At a codimension-four component of the base locus the same
three-cut argument, now on the identity cover, gives a regular
local ring of dimension four. Let \(h\) be the dimension of the
span of the linear parts of the ideal \(I\) in that ring. Use
independent generators as the first \(h\) parameters. Weights two
on these parameters and one on the others give every element of
\(I\) weight at least two, so
\(\operatorname{lct}(I)\le(4+h)/2\) for \(h<4\).
Consequently \(h\ge3\). For \(h=4\) the ideal is the maximal
ideal. For \(h=3\), its image in the regular local curve cut
out by the first three parameters is \((x_4^j)\), with \(j\ge2\)
finite since \(I\) is primary there. Thus
\(I=(x_1,x_2,x_3,x_4^j)\). The weight valuation with weights
\((j,j,j,1)\) gives threshold at most \(3+1/j\).
Since \(3+1/3<173/50\), necessarily \(j=2\).
\end{proof}

\begin{proposition}\label{prop:dec2-index-codim-five}
Suppose the threshold of the degree-one system on the local index
covers of a quartic ninefold replacement is strictly greater than
three. If the limiting hyperplane root \(L\) is Cartier through
codimension four, it is Cartier through codimension five.
\end{proposition}
\begin{proof}
Suppose a component of the non-Cartier locus has codimension five.
At a general point of this component let \(I\) be the local index
of \(L\), and let \(Y\) be its local cyclic index cover. Set
\(D=4(7/10)^9\). The finite-degree formula and the local-to-global
inequality bound the normalized volume of the ninefold cover below
by \(9^9ID\). The smooth-point upper bound implies \(I\le6\).
Thus the relation \(-K_Z\sim7L\) identifies \(Y\) with the canonical
index cover. Its group fixes the reduced lift of a dense open subset
of the component pointwise, and acts freely off that lift near its
general point. Indeed outside the component \(L\) is Cartier after
shrinking about its general point, while on the component the local
index is \(I\). Choose four invariant parameters along the component.
General transverse slices, taken after shrinking this open set, give
a five-dimensional Gorenstein klt germ \(T\) on the cover. The
threshold bound holds on the transverse slice by restricting a log
resolution smooth over a dense open subset of the component, as in
Proposition~\ref{prop:nov30-quartic-nine-index-cover}. The restricted
degree-one sections have nontrivial character, so vanish at the
marked point.

The Gorenstein threshold criterion of
\cite[Lemma~4.5]{kePanWuCubics} now applies with \(c=5\).
Indeed three general combinations of these functions cut out a klt
surface at that point. The base ideal on the slice has
codimension at least four by the threshold bound and the
base-codimension inequality for the lifted system. In the
Cohen--Macaulay ring of \(T\) the three general functions thus
form a regular sequence. On a log resolution of the ideal, Bertini
makes the three general moving divisors smooth and transverse to one
another and to all exceptional strata. For each exceptional divisor
the sum of the three fixed orders is strictly smaller than its log
discrepancy. Successive plt Bertini and Cartier adjunction give
normal klt intersections at each cut: on the resolved intersection
the remaining exceptional coefficients are still less than one,
and its intersections with the fixed strata are transverse. The
base has codimension at least two before even the last cut, so no
intermediate component is forced into the base. The resulting klt
surface is Gorenstein by adjunction. It is smooth or Du Val and has
embedding dimension at most three. Lifting generators through the
three equations shows that \(T\) has embedding dimension at most six.

Restrict the relative canonical cover of the smoothing to these
transverse parameters, lifting them invariantly as in
Proposition~\ref{prop:nov30-quartic-nine-index-cover}. This gives
an equivariant deformation whose general fiber has free group
action. A smooth central transverse cover would give a fixed
section by diagonalizing formal parameters together with the
invariant deformation parameter. Consequently \(T\) is a singular
hypersurface in six variables. The full local cover is likewise a
singular hypersurface, of dimension nine: the four independent
stratum parameters together with the transverse generators bound
its embedding dimension by ten. A smooth ninefold cover would
admit a fixed section in its own equivariant smoothing. The
complete-intersection ordinary-double-point volume bound
\cite[Theorem~2.12]{ar2007.14320} now gives
\[
             I D \le 2(8/9)^9,\qquad I\le4 .
\]

In a diagonal minimal transverse embedding, the defining equation
is invariant. In fact it lifts to a semi-invariant equation of the
deformation, and a nontrivial equation character would force the
origin to remain a fixed point on nearby fibers. All six coordinate
characters are units modulo \(I\). For any nontrivial subgroup a
fixed ambient subspace of dimension at least two would meet the
central hypersurface away from the origin. A fixed line on which the
equation vanishes identically gives the same contradiction. If its
restriction has finite positive order, Weierstrass preparation gives
fixed points on the general deformation fiber. These contradictions
use freeness off the marked point centrally and freeness everywhere
on the general fiber.

Normalize characters so that the functions representing sections of
\(L\) have weight one. Write \(w_i\) for the weights of the coordinate
functions, under this same action on functions. For the invariant
equation, the residue form has weight \(\sum w_i\), and an invariant
canonical form equals this residue times a function of weight
\(-\sum w_i\). Thus functions representing anticanonical sections
have weight \(+\sum w_i\). Since \(-K_Z\sim7L\), their weight is
also \(7\) modulo \(I\). If \(I\) were even, the sum of the six
units \(w_i\) would be even. Hence \(I=3\). If \(h\) variables have
character one and the rest have character two, the character equation
reads \(12-h\equiv 7\pmod 3\). Thus \(h=2\) or \(h=5\).

For \(h=2\), give these \(h\) variables weight two and the remaining
variables weight one in a weighted blowup. All section functions
have weight at least two. The invariant equation, which has no
linear term, has weight at least three because no monomial involving
just two of the character-two variables is invariant. Inversion of
adjunction from the smooth ambient space and the log discrepancy of
its weighted blowup give a threshold at most
\((6+h-3)/2=5/2\), a contradiction.

For \(h=5\), all invariant transverse quadratic monomials pair the
unique character-two variable with a character-one variable, so
the quadratic rank is at most two. Work in an analytic chart of the
full hypersurface cover at a general closed point of the stratum,
with stratum parameters \(u_1,\ldots,u_4\) and transverse variables
\(x_1,\ldots,x_6\) vanishing on the lifted stratum. Its equation
belongs to \((x)^2\), since the cover is transversely singular all
along a dense open part of the stratum. The transverse quadratic
rank is at most two at each general point by the preceding
character argument. Shrink to its constant-rank locus. Its kernel
is a vector subbundle, so a change of transverse coordinates
depending on \(u\) makes the entire quadratic part involve only
\(x_1,x_2\).

Give \(x_1,x_2\) weight fifteen, \(x_3,\ldots,x_6\) weight ten,
and \(u_1,\ldots,u_4\) weight eight. The equation has minimum
weight \(b\ge30\). The weighted hypersurface estimate proved in
Proposition~\ref{prop:nov30-quartic-nine-index-cover} bounds the
volume by \(b(102-b)^9/(15^2 10^4 8^4)\). It applies for
arbitrary initial form. Here \(b<102\) by klt and inversion of
adjunction, and \(b(102-b)^9\) is decreasing for \(b\ge30\).
Consequently
\[
  \frac{\widehat{\operatorname{vol}}(Y)}{9^9}
    \le \frac{30\cdot 8^9}{15^2 10^4 8^4}
    =\frac{4096}{9375}
    <12(7/10)^9=3D .
\]
This contradicts the normalized-volume lower bound on the full
cover and excludes the proposed component.
\end{proof}

\subsection{The integral Schubert classes on a fivefold slice}
Write $B$ for the base of the ten sections projecting the cone from
its vertex, including the non-Cartier locus. In the replacement
setting it has codimension at least three. To see this, resolve
the projection to $G$ and let $F$ be its general resolved threefold
fiber. Write $P$ for the hyperplane pullback from $Z$, $J$ for the
moving Pl\"ucker class and $x=P^3F/4$. The count of conormal layers
in Proposition~\ref{prop:dec2-ninefold-tangent-pencil}, now with
trivial conormal of rank six, gives
\[
  6/7\ge S_L(\operatorname{ord}_F)
       \ge \int_0^\infty \max(0,1-84xt^6)\,dt
       = (6/7)(84x)^{-1/6}.
\]
Here $84=\binom96$ and the discrepancy of the blowup of $F$ is
six. Thus $a_6=P^3J^6/4=5x\ge5/84$. There is no fixed prime by
the three-section criterion, so $a_0=a_1=1$. Apply the nef
segment and the local threshold comparison of
\cite[Lemmas~5.1, 5.3 and 5.4]{kePanWuCubics} to this subsystem.
The proofs apply to its moving pullback $J$ since the base order
is the corresponding coefficient of $P-J$ on a resolving model.
As $5/84>(11/20)^5$, log-concavity gives
\[
  \operatorname{codim}_Z B\ge\frac7{10}\sum_{i=0}^6 a_i
           >\frac7{10}\left(2+\sum_{j=1}^5(11/20)^j\right)>2 .
\]

\begin{proposition}\label{prop:dec2-grass-singular-base}
Let $(Z,L)$ be a polarized smoothing limit of smooth quartic
ninefolds with $-K_Z\sim 7L$ and hyperplane image the cone over
$G=\operatorname{Gr}(2,5)$. Let $B$ be the complement of the open
locus on which $L$ is invertible and the ten Pl\"ucker sections
generate it. Suppose $\operatorname{codim}_Z B\ge3$. Then
\[
                 \operatorname{codim}_Z(B\cap\operatorname{Sing}Z)\le5 .
\]
\end{proposition}
\begin{proof}
Suppose $\dim(B\cap\operatorname{Sing}Z)\le3$. Cut by four general
members of a sufficiently high relative very ample Cartier multiple
of the anticanonical class. The central slice $V$ is an integral
projective fivefold by Bertini. The relative slice is flat near the
central fiber by the fiberwise regular sequence criterion (the
central klt variety is Cohen--Macaulay). The nearby fiber $V_t$ is
a smooth complete intersection of four ample divisors in a smooth
quartic ninefold. The closed subset $K=B\cap V$ has dimension at
most two and $V$ is smooth near $K$, by general choice and Bertini
on the smooth locus of $Z$. In particular $L$ is a line bundle
near $V$: at points of $K$ this follows from smoothness of $Z$, and
on $V\setminus K$ it follows from the definition of $B$. Put
$h=c_1(L|_V)$. Excision at a smooth neighborhood of $K$ and
duality with the Borel--Moore homology of $K$ give
$H^j_K(V,\mathbf Z)=0$ for $j=4,5$. Thus every integral fourth
cohomology class on $V\setminus K$ extends uniquely over $V$.

The map on $Z\setminus B$ to $G$ pulls back the hyperplane line
bundle to $L$. Let $a_1,a_2$ be the unique classes in
$H^4(V,\mathbf Z)$ extending the pullbacks of $c_2(\mathcal S^*)$
and $c_2(\mathcal Q)$, where $\mathcal S,\mathcal Q$ are the
tautological bundles on $G$. Their defining sequence gives
$a_1+a_2=h^2$. Each class has strictly positive degree against $h^3$.
To verify this, pull back general translates of the corresponding
Schubert cycles first on $Z\setminus B$ and then cut to $V$. Over
the dominant smooth part of the map these inverse images are
nonempty of codimension two; generically transverse pullback holds
there by the transitive group action. More explicitly, take three
additional general members of the same very ample system cutting
$V$ to an integral surface $D$. This surface avoids $K$. The smooth
locus of $Z$ meets $D$ densely, and its complement on $D$ is finite
by normality. General translates avoid the images of these finitely
many points, and also meet $D$ transversely off them by generic
transversality. Nonemptiness follows since each inverse image on
the smooth dominant part of $Z$ has a seven-dimensional closure,
and seven general very ample sections meet this closure away from
its proper boundary. Thus $\int_D a_i>0$. Since the cutting class
on $V$ is a positive multiple of $h$, the asserted positivity follows.

Use the integral specialization in ordinary cohomology from the
central fiber to a nearby smooth fiber. It preserves products and
is compatible with specialization of the fiber fundamental class,
so $\int_V\gamma=\int_{V_t}\operatorname{sp}(\gamma)$ for each
integral class of top degree. Equivalently one may restrict classes
from a sufficiently small inverse image of a closed disk, which
retracts onto $V$ by stratified trivialization. The equality of
the two fiber classes there follows by joining their values by an
arc: the regular fiber is reduced and the flat central fiber has
multiplicity one at its generic point. The multiple $7h$ equals
the restriction of the relative anticanonical class near $V$, or
one may take a further Cartier multiple of this class. Integral
weak Lefschetz gives $H^2(V_t,\mathbf Z)=\mathbf Z h_t$ and
$H^4(V_t,\mathbf Z)=\mathbf Z h_t^2$ for the primitive hyperplane
class $h_t$. Therefore $\operatorname{sp}(h)=h_t$ and
$\operatorname{sp}(a_i)=b_i h_t^2$ for integers $b_i$. Their
positivity follows from the degree equality and $\int_{V_t} h_t^5>0$.
Applying specialization to $a_1+a_2=h^2$ gives $b_1+b_2=1$,
a contradiction.
\end{proof}

\paragraph{A transverse node where the Schubert classes fail to extend.}
The codimension-five case in Proposition~\ref{prop:dec2-grass-singular-base}
has a local model that also realizes dominance over the cone.
Take a four-dimensional vector space $E$, set $V=E\oplus\mathbf C$,
and for $b\in E^*$ let $i_b:E\hookrightarrow V$ be the inclusion
of the graph of $b$. Choose a volume form on $E$. Apart from a
nonzero scalar the Pfaffian $q(x)$ of $x\in\bigwedge^2 E$ is
the nondegenerate quadratic form $x\wedge x$. On the affine ninefold
\[
       Y=\{(b,x)\in E^*\times\textstyle\bigwedge^2E:q(x)=0\}
\]
the ten entries of $y=(\bigwedge^2 i_b)(x)$ and the extra section
$z=1$ define a morphism to the chart $z=1$ of the projective
Grassmannian cone. In this chart a nonzero point is a decomposable
bivector $y\in\bigwedge^2 V$. For general $y$, its support
two-plane projects isomorphically to a two-plane in $E$.
For the image of $i_b$ to contain the support specifies the two
values of $b$ on a basis of that projected plane. Once $b$ is
chosen, $x$ is uniquely determined by projection. Thus the map
is dominant and its general fiber is an affine plane. The
ten-section base and the singular locus of $Y$ both equal
$E^*\times\{0\}$, since the projection of $y$ to $\bigwedge^2E$
equals $x$. The transverse germ is an ordinary double point of
dimension five and $q(x)=t$ smooths it.

For fixed $b$, the punctured transverse germ has link the circle
bundle of $\mathcal O(-1)$ on the smooth quadric
$Q^4=\operatorname{Gr}(2,E)$. The Gysin sequence identifies the
fourth cohomology of the link with
\[
   H^4(Q^4,\mathbf Z)\big/\big(h\,H^2(Q^4,\mathbf Z)\big)
      =(\mathbf Z\sigma_2\oplus\mathbf Z\sigma_{1,1})
             \big/\mathbf Z(\sigma_2+\sigma_{1,1}).
\]
Here $h$ is the Pl\"ucker class, $H^3(Q^4,\mathbf Z)=0$ by the
Schubert cell decomposition, and $h^2=\sigma_2+\sigma_{1,1}$
by the tautological sequence. The two Grassmannian Chern classes
restrict to these two Schubert classes: the tautological subbundle
restricts to that on $\operatorname{Gr}(2,E)$ and the quotient
acquires one constant summand. Their images in the displayed
group are opposite nonzero generators. Neither extends across the
contractible local cone. Producing a projective K-polystable
replacement over a smooth quartic GIT point would require
additional global geometry.

\subsection{Divisors and quadratic sections on the cone}
\begin{proposition}\label{prop:dec2-grass-all-primes}
Every prime divisor on the Grassmannian cone $W$ is dominated by
a prime divisor on the domain of the hyperplane map from $Z$.
For the subring $T$ of the hyperplane section ring $S=R(Z,L)$
generated in degree one, one has $S\cap\operatorname{Frac}(T)=T$.
\end{proposition}
\begin{proof}
The punctured affine Pl\"ucker cone is the complement of the zero
section of $\mathcal O_G(-1)$. Its class group is
$\operatorname{Pic}(G)/\mathbf Z[\mathcal O_G(1)]=0$ by localization
on this line bundle. The Pl\"ucker ring is normal by projective
normality, and its vertex has codimension greater than one. Hence
the ring is factorial, as is its polynomial extension $T$ by
the vertex coordinate. A prime divisor $\Delta$ on $W$ therefore
has an irreducible homogeneous equation $s$ of degree $m>0$.

Suppose no prime on $Z$ dominates $\Delta$. Write
$\operatorname{div}_Z(s)=\sum e_jD_j\sim mL$. Each $D_j$ maps
into a subset of $\Delta$ of codimension at least two in $W$.
Its generic image is away from the cone vertex by the
three-section criterion. At a general smooth point of $D_j$
the restricted map on $D_j$ has differential of rank at most
five. Adding the normal direction gives rank at most six. Thus
the divisorial degeneracy $R_2$ of the map to the sevenfold $W$
contains $\sum D_j$, and
$L^8R_2/L^9\le2/5$ by Proposition~\ref{prop:dec1-grassmannian-fibers}.
On the other hand for each prime divisor $D_j$ on $Z$,
\[
      e_j/m\le T_L(D_j)\le10S_L(D_j)\le10/7 .
\]
The middle inequality follows by integrating the concave ninth
root of divisor volume up to the pseudoeffective threshold.
Consequently
\[
     m=\sum_j e_j\, L^8D_j/L^9
         \le (10m/7)\,L^8R_2/L^9\le 4m/7 ,
\]
which is impossible. A homogeneous element of $S_k$ lying in
$\operatorname{Frac}(T)$ gives a rational section of
$\mathcal O_W(k)$. At the generic point of each prime of $W$ a
dominating divisor on $Z$ detects a pole, with order multiplied
by the positive ramification index. The polarizations agree at
the generic point of that divisor. Thus the section has no poles,
and normality and projective normality give membership in $T_k$.
The grading action preserves the intersection. Applying it at
distinct nonzero scalars and taking linear combinations recovers
the homogeneous components of each element, giving the assertion
for the entire ring.
\end{proof}

\subsection{The first nonzero terms of the Pl\"ucker Pfaffians}
\begin{proposition}\label{prop:dec2-first-pfaffians}
Suppose a pointed family with generic fiber a smooth quartic
ninefold has the polarized K-polystable replacement $(Z,L)$ of
Corollary~\ref{cor:dec1-root-grassmannian-only}. Write $S$ for its
hyperplane section ring and $T\subset S$ for the ring generated by
$S_1$. One can lift the ten Pl\"ucker coordinates to the relative
degree-one lattice so that the first nonzero special coefficients
of their five quadratic Pfaffians span a subspace of dimension at
least three in $S_2/T_2$. More precisely, if $A$ is the lifted skew
matrix and $m$ is the minimum of the orders of the five Pfaffians
in the relative degree-two lattice, take the reductions of
$t^{-m}\operatorname{Pf}_i(A)$ for a uniformizer $t$ and the usual
alternating signs on the maximal Pfaffians.
\end{proposition}
\begin{proof}
Work over the complete discrete valuation ring after base change.
The relative graded pieces of the hyperplane root are free and
commute with restriction as in the proof of
Theorem~\ref{thm:dec1-quartic-high-dimensions}. Start with any
relative lift of an adapted basis of $S_1$. Regard the reduction
$A_0$ as a bivector in a five-dimensional vector space $V$, and
use a volume form to identify $\bigwedge^4V$ with $V^*$. On the
special fiber $A_0$ is decomposable and the extra independent
linear coordinate is denoted by $z$. For the minimum $m$ above,
the reduction $h$ of $t^{-m}(A\wedge A)/2$ is nonzero and
satisfies $A_0 h=0$, by the Pfaffian identity in odd size.

First suppose every entry of $h$ belongs to $T_2$. We claim that
$h=A_0\wedge B$ for some skew matrix $B$ of linear forms in the
eleven coordinates. On $G=\operatorname{Gr}(2,V)$ write $U$ and
$Q$ for the tautological subbundle and quotient. At a point of
$G$ the kernel of $A_0:V^*\to V$ is $Q^*$. Since the cone ring
is $\mathbb C[G][z]$ in homogeneous coordinates, the claim reduces
to surjectivity of
\[
\bigwedge^2 V\otimes H^0(G,\mathcal O_G(j-1))
\longrightarrow H^0(G,Q^*(j))
\]
for $j=1,2$, and to $H^0(G,Q^*)=0$ for the $z^2$ coefficient.
The displayed map is exterior multiplication by the universal
bivector. The vanishing in degree zero also follows by intersecting
the annihilators of all two-planes in $V$. Since
$Q^*(1)=\bigwedge^2Q$ up to a constant determinant factor, the
two targets in degrees one and two are irreducible
$\operatorname{SL}(V)$-representations by the Borel--Weil theorem
on a Grassmannian. In the quotient-bundle convention the dominant
partitions are $(1,1,0,0,0)$ and $(2,2,1,0,0)$, respectively.
The equivariant maps are nonzero by evaluation at a plane and
choice of two independent vectors in its quotient. This proves
surjectivity. Projective normality of the Pl\"ucker embedding
identifies the coefficient polynomials used here with global
sections.

Replace $A$ by $A-t^m\widetilde B$, lifting the coefficients of
$B$. The new matrix has the same reduction, still occurs among
the coordinates of a basis, and its minimum Pfaffian order has
strictly increased. This procedure terminates. In fact the
relative multiplication map from the symmetric square of the
degree-one lattice to the degree-two lattice is an isomorphism
over the fraction field: the generic quartic has no quadratic
equation. The largest Smith exponent of this map bounds the
divisibility in the target of the image of every primitive
element in its source. The quadratic Pfaffian polynomials in any
of our bases are primitive, giving the required bound.

Let $d$ be the positive dimension of the span of the entries of
$h$ modulo $T_2$. If $d\le2$, write
$h=q+\sum_{j=1}^d b_j y_j$, with $q\in T_2\otimes V^*$,
linearly independent constant covectors $b_j$, and the classes of
the $y_j$ independent modulo $T_2$. At a general two-plane in $V$
the restrictions of these $d$ covectors are linearly independent.
Contraction with its nonzero decomposable bivector therefore makes
the columns $A_0 b_j$ independent over $\operatorname{Frac}(T)$.
Solving $A_0h=0$ over that field gives $y_j\in\operatorname{Frac}(T)$
for every $j$. Proposition~\ref{prop:dec2-grass-all-primes} forces
$y_j\in T_2$, a contradiction.
\end{proof}
Each quadratic section representing a nonzero one of these classes varies on the
geometric generic surface fiber after division by the square of a
general coordinate. Indeed $\mathbf C(W)$ is relatively algebraically
closed in $\mathbf C(Z)$ by
Proposition~\ref{prop:dec1-stein-hyperplane-fibers}. A section
constant on the geometric generic fiber would therefore give an
element of $\operatorname{Frac}(T)$ and would belong to $T_2$.
This measures independence modulo $T_2$ over the constants.
The dimension of the rational image of the moving quadratic
system on the surface still requires a separate argument.

\subsection{Excluding a quartic section of the Pfaffian incidence variety}
\begin{proposition}\label{prop:dec2-spinor-incidence-quartic}
In the Grassmannian-cone case, the graded hyperplane section ring
requires a generator of degree at least three. In particular the ring
cannot have the presentation
\[
  \mathbf C[A_{ij},z,y_i]\big/
   \big(\operatorname{Pf}_i(A),\ (Ay)_i,\ H(A,z,y)\big),
\]
where \(A\) is a skew \(5\times5\) matrix with ten independent
entries of degree one, \(z\) has degree one, the five entries of
\(y\) have degree two, \(\operatorname{Pf}_i(A)\) are the maximal
Pfaffians, and \(H\) has weighted degree four.
\end{proposition}
\begin{proof}
Suppose the ring is generated in degrees one and two.
Choose generators \(y_i\) of degree two, since
\(\dim S_2/T_2=5\). Write the five first Pfaffians
as \(By+c(A,z)\), where \(B\) is a constant matrix of rank
\(\nu\ge3\) by Proposition~\ref{prop:dec2-first-pfaffians},
and \(c\) is a vector of quadratic forms.
Their syzygy gives \(A(By+c)=0\) for the central matrix \(A\).
Let \(R\)
be the reduced coordinate ring of the affine incidence variety
\[
 \{(A,y,z): \operatorname{rank}A\le2,\ A(By+c(A,z))=0\},
 \qquad \deg A_{ij}=\deg z=1,\quad \deg y_i=2 .
\]
It surjects onto \(S\) because \(S\) is a domain and all its
generators satisfy these conditions. To compute its Hilbert
series write \(A=p\wedge q\)
for two vectors in a fixed five-space. In the sixteen variables
\((p,q,y,z)\) impose the two equations
\((By+c)(p)=(By+c)(q)=0\), with \(A=p\wedge q\) in \(c\).
Over the open set of \((p,q,z)\) where the two restrictions to
\(\operatorname{im}B\) are independent, their solution is an
affine three-space bundle in \(y\), smooth and irreducible.
In the base of this projection, the locus of restriction rank
one has codimension \(\nu-1\ge2\), since it is the rank-one
locus of a general two-by-\(\nu\) coordinate matrix. Each such
fiber has dimension at most four. Rank zero costs \(2\nu\) in
the base and adds at most two to the fiber dimension. Thus
the complement of the indicated smooth open has dimension at
most thirteen. Every component has dimension at least fourteen
by the two equations in sixteen variables. It follows that these
equations define an irreducible, generically reduced complete
intersection. Cohen--Macaulayness implies the scheme is reduced.
Its invariants for \(\operatorname{SL}_2\) acting on the two
vectors form a domain, generated by \(p\wedge q\) and
\(y\) by the first fundamental theorem and exactness of invariants.
The image is the entire incidence variety: when \(A\ne0\),
its plane is annihilated by \(By+c\), and when \(A=0\) take
\(p=q=0\). Thus this invariant ring is \(R\).
Its Hilbert function is independent of \(B,c\) for the ranks
under consideration. Indeed assign weights one to \(p,q\), weight
two to \(z\), and weight four to \(y\), twice the intended degrees
of the invariant generators. The two equations of the preceding
complete intersection always form a copy of the standard
\(\operatorname{SL}_2\) representation in degree five. Their
equivariant Koszul resolution thus has the same graded characters
as for \(B=1,c=0\), and taking invariants in each degree is exact.
For \(B=1,c=0\), before adjoining \(z\), the components of bidegree \((i,j)\) in \(A,y\)
are the coordinates of the incidence flag
\(\operatorname{Fl}(2,4;5)\) in the product of its two Pl\"ucker
spaces. By Borel--Weil these are the irreducible representations
of highest weight \(i\omega_2+j\omega_4\), up to duality: the
evaluation of the bihomogeneous polynomials on this flag is an
equivariant nonzero map to that irreducible space of sections.
Here \(i+2j\) is the weighted degree. The Littlewood identity for
the symmetric algebra of an exterior square gives
\[
 \operatorname{Sym}^k(\textstyle\bigwedge^2 \mathbf C^5)
  =\bigoplus_{i+2j=k}
          \mathbf S_{(i+j,i+j,j,j)}\mathbf C^5 .
\]
Indeed the Schur functions in the classical identity
\(\prod_{a<b}(1-x_ax_b)^{-1}=\sum_{\lambda'\ {\rm even}}s_\lambda(x)\)
have even column lengths. With five variables the columns have
length two or four, of multiplicities \(i,j\). These are exactly
the flag representations above. Hence the Hilbert series after
adjoining \(z\) is \(H_R(t)=(1-t)^{-11}\). The Hilbert series
of the quartic section ring is \((1-t^4)H_R(t)\), so the kernel
of \(R\to S\) has Hilbert series \(t^4H_R(t)\). Choose a nonzero
element \(H\) of its one-dimensional degree-four piece. Since
\(R\) is a domain, the principal ideal \((H)\) already has the
full Hilbert series of the kernel. Thus \(S=R/(H)\).
For the explicit presentation in the statement
one has the same quotient: the reduced incidence ideal already
vanishes in the domain \(S\), so passing from the displayed
Pfaffian equations to their radical preserves the quotient.

Put \(q(y)=H(0,0,y)\), a quadratic form. The base of the full
hyperplane system has support
\(\{A=z=0,\ q(y)=0\}\) in the projective four-space of the
degree-two variables. We describe the local cover near its general
point. On a chart normalize one nonzero entry of \(y\) to one,
leaving a group of order two which acts by sign on \(A,z\) and
fixes \(y\). At a general point of each base component
\(By\ne0\), as the projective kernel of \(B\) has dimension at
most one. Near that point choose a change of basis depending
analytically on the covector \(By+c\) which makes it a multiple
of a fixed coordinate covector. The corresponding congruence
transformation on \(A\), leaving \(y,z\) fixed, is a local
coordinate change: at \(A=z=0\) its derivative on \(A\) is
invertible. It preserves the sign action. The equations
\(A(By+c)=0\) then eliminate the entries of the transformed \(A\)
in the row and column of the fixed covector. On the resulting
four-dimensional complement the remaining skew matrix has six
entries \(x\). Its single four-by-four Pfaffian is a nondegenerate
quadric \(Q(x)\). All the original Pfaffian equations are
equivalent on this chart to \(Q=0\), by the invertible change of
basis sending \(y\) to a coordinate covector. Thus the chart
cover is the complete intersection \(Q(x)=H'=0\) in the smooth
ambient space with six variables \(x\), one variable \(z\), and
four affine \(y\)-coordinates. The degree-one ideal on this
cover is \((x,z)\). This indeed computes the local root cover:
on the even chart the displayed cover is given by trivializing
\(\mathcal O(2)\), with the two choices of its square root.
It is free away from \(\{x=z=0\}\), of codimension at least
five. The chart ring is a complete intersection, normal because
it is Cohen--Macaulay and regular in codimension one by this
\'etale description over the normal quotient. The odd eigensheaf
therefore agrees by reflexive extension with the hyperplane
root. The threshold bounds on the root cover apply.

If \(q=0\), fix the four \(y\)-coordinates at a general point of
the base. This gives a transverse fivefold complete intersection
in seven variables. Both equations have order at least two:
\(H'\) has no pure \(y\)-term, and invariance under the sign action
excludes terms of first degree in the odd variables. The ordinary
ambient blowup and inversion of adjunction for a codimension-two
local complete intersection bound the threshold of \((x,z)\) by
\(7-2\cdot2=3\). Explicitly the exceptional divisor has log
discrepancy seven in the ambient space, order at least two on
the defining ideal and order one on \((x,z)\). This contradicts
Proposition~\ref{prop:dec2-ninefold-tangent-pencil}.

Assume \(q\ne0\). Over a general point \([A_0:z_0]\) of the
cone, \(A_0B\) has rank two and its kernel is a three-space \(K\).
The closure of the surface fiber is described by
\((A,z)=s(A_0,z_0)\), \(y\in s^2 y_0+K\), and \(H=0\);
here choose \(y_0\) solving
\(A_0(By_0+c(A_0,z_0))=0\), possible since \(A_0B\) maps onto
\(\operatorname{im}A_0\). The parameter space is
\(\mathbf P(1,2,2,2)\) with coordinates \(s,y-s^2 y_0\). Indeed on
\(s\ne0\) this is exactly the fiber. The second Veronese of this
weighted space is \(\mathbf P^3\), with linear coordinates
\(s^2,y-s^2y_0\). The equation becomes a quadratic with restriction
\(q|_K\ne0\) at infinity. Indeed a general \(y\) lies in one such
\(K\): choose a general two-plane in \(\ker(By)\). The restrictions
of \(\operatorname{im}B\) to this plane span its dual since at
least two independent covectors remain on \(\ker(By)\). Thus
no nonzero quadratic vanishes on every \(K\). The
quadric has no component supported at infinity. It is generically
reduced and integral by the geometric integrality of the general
fiber. The indicated weighted coordinates give a closed immersion
of this fiber closure into the presented model by their linear
substitutions. The restriction of \(2L\) is the usual hyperplane
class of the quadric, hence \(L^2\cdot S=1/2\).

The base \(q=0\) now has codimension six, even when \(q\) is a
square. For \(a_i\) as in
Proposition~\ref{prop:dec2-ninefold-tangent-pencil}, we have
\(a_0=\cdots=a_5=1,\ a_7=5/8\), and \(a_6\ge\sqrt{5/8}\) by
log-concavity. The mixed threshold bound on the local cover is
therefore at least
\[
                 \frac7{10}\left(6+\sqrt{5/8}+5/8\right)>5 .
\]
At a general point of its reduced base cut three independent
parameters along the base. There remain seven odd variables and
one even parameter \(b\) transverse to the reduced quadric in
the \(y\)-space. If \(q\) is reduced, then \(H'\) has a nonzero
linear \(b\)-term. Eliminate \(b\). The transverse cover is
\(Q(x)=0\) in the seven odd variables \(x,z\), and the
degree-one ideal is their maximal ideal. The ambient blowup and
inversion of adjunction give threshold at most \(7-2=5\).
If \(q\) is a square, both \(Q\) and \(H'\) have order at
least two in the eight transverse variables: the pure even
term vanishes to order two and every term involving odd variables
has even total odd degree. The same codimension-two
inversion of adjunction bounds the threshold of \((x,z)\) by
\(8-2\cdot2=4\). The transverse threshold retains the mixed
lower bound by a general slice along the stratum, as in
Proposition~\ref{prop:nov30-quartic-nine-index-cover}. Both cases
give the contradiction.
\end{proof}

\begin{corollary}\label{cor:dec2-incidence-quadratic-base}
Let \(S\) be the section ring in
Proposition~\ref{prop:dec2-spinor-incidence-quartic}. Then \(S\)
is not finite over the subalgebra generated by \(S_1,S_2\). In
particular the complete system \(|2L|\) has a nonempty geometric
base locus.
\end{corollary}
\begin{proof}
Use \(B,c,\nu,R\) from that proposition and write \(S'\) for
the subalgebra generated by \(S_1,S_2\), always a quotient of
\(R\). We first show that the rank-drop locus in the construction
of \(R\) has codimension at least two in its solution variety,
also for \(\nu=3\). Only the rank-one stratum needs attention,
as the preceding dimension count gives the other cases. Choose
coordinates so that \(\operatorname{im} B\) consists of the
first three covectors. In the polynomial ring
\(P=\mathbf C[p,q,z]\) form the two-by-three matrix \(M\) with
rows the first three entries of \(p,q\). The equations in the
three effective \(y\)-variables have right-hand side
\[
                 b=-\big(c(p\wedge q,z)(p),\
                            c(p\wedge q,z)(q)\big)^{\rm t}.
\]
We claim they are incompatible at the generic point of the rank-one
stratum. The minors of \(M\) generate the prime ideal of its
codimension-two determinantal locus. The cokernel of \(M\) has no
embedded associated prime: its free resolution starts with the
column of its three signed minors and continues with \(M\).
Indeed that column spans the kernel over the fraction field,
and its entries have no common factor. The resolution has length
two, so Auslander--Buchsbaum at every prime of height greater
than two gives positive depth. At the generic point of the
determinantal locus a unit entry reduces its cokernel to the
residue field, since the two residual entries generate its minor
ideal. Compatibility there would therefore force the class of
\(b\) in the cokernel to vanish everywhere.

A polynomial solution lifting \(b\) can be taken homogeneous of
weight four for the grading assigning weight one to \(p,q\)
and two to \(z\). Its three entries can also be chosen
\(\operatorname{SL}_2\)-invariant: the map and datum transform by
the standard two-dimensional representation, and taking invariants
in the finite-dimensional homogeneous representation is exact.
Equivalently apply the Reynolds operator to a solution under the
natural action on the equation and unknowns. These invariant
entries are quadratics in \(A=p\wedge q,z\).
After subtracting them, the three effective independent quadratic
generators on \(Z\) would form a vector in the kernel of \(M\) at
its generic plane. They would therefore be a common scalar \(g\)
times the three signed minors. The proportionality is intrinsic
in the plane, since the minors are Pl\"ucker coordinates, so
\(g\) is a rational section of \(L\) on \(Z\). The three minors
are independent linear sections and share no prime divisor on
\(Z\) by the three-section criterion. The three regular products
thus show that \(g\) is regular in codimension one, hence
\(g\in S_1\) by reflexivity. This puts the three effective
quadratic generators in \(T_2\), contradicting their independent
classes. This proves generic incompatibility.

Compatibility on the rank-one stratum is a closed condition
there. Thus for \(\nu=3\) it costs at least one further base
condition, beyond the codimension two of rank drop, whereas
the fiber dimension increases by at most one. For \(\nu\ge4\)
the rank-one condition alone already costs at least three. Hence
the complete intersection in \((p,q,y,z)\) used to construct
\(R\) is smooth in codimension one in every case. It is normal
by Serre's criterion, so its invariant ring \(R\) is normal too.
On \(A\ne0\), choosing a basis \(p,q\) with prescribed exterior
product has fibers of constant dimension three, the dimension
of \(\operatorname{SL}_2\). Thus the same codimension estimate
applies on this part of \(\operatorname{Spec}R\). The part
\(A=0\) has dimension six, compared with dimension eleven for
\(\operatorname{Spec}R\).

On the resulting complement of codimension at least two in
\(\operatorname{Spec}R\), the map to \((A,z)\) is an affine
three-space bundle over the open where \(\operatorname{rank}A
=\operatorname{rank}(AB)=2\). Its base has trivial class group:
the punctured Pl\"ucker cone is the complement of the zero section
in \(\mathcal O_G(-1)\), and \(\operatorname{Pic}(G)=\mathbf Z H\).
Adjoining \(z\) and passing to the indicated open preserve
triviality of the class group. The total affine bundle also has
trivial class group, by restricting a divisor to its generic
affine-space fiber and then subtracting the pullbacks of base
primes. Normality and the codimension estimate prove that \(R\)
is factorial.

If \(S\) were finite over \(S'\), the two rings would have the
same dimension ten. The homogeneous kernel of \(R\to S'\) would
have height one and hence a single homogeneous generator of degree
\(k\). The degree-one and degree-two spaces of \(S'\) equal those
of \(S\) and have the same dimensions as those of \(R\); thus
\(k\ge3\). Let \(p\) be the generic rank of \(S\) over \(S'\).
Comparing the leading coefficients of their Hilbert series at
one using \(H_R=(1-t)^{-11}\) gives \(4=pk\). Indeed a finite
graded torsion-free module contains a graded free submodule of
the same generic rank with lower-dimensional quotient. It follows
that \(k=4\) and the Hilbert series of \(S'\) and \(S\) agree
in every degree, forcing equality. This contradicts
Proposition~\ref{prop:dec2-spinor-incidence-quartic}.
Finally, if the degree-two sections have no geometric base point,
their homogeneous ideal in \(S\) has radical \(S_+\). The quotient
has only finitely many nonzero graded pieces, so graded induction
makes \(S\) finite even over the algebra generated by \(S_2\).
\end{proof}

If the map given by \(|2L|\) has nine-dimensional image, it is
birational and its graded image ring \(S'\) has the form \(R/(H_3)\)
for a homogeneous cubic. Indeed \(\operatorname{Frac}(S)\) is then
a finite extension of \(\operatorname{Frac}(S')\). The same
height-one argument gives a single relation of degree \(k\ge3\).
Even without module finiteness, the generic field degree \(p\)
satisfies \(pk\le4\): one can choose \(p\) homogeneous elements
of \(S\) linearly independent over \(\operatorname{Frac}(S')\),
since the homogeneous elements span \(S\) and localization at
the nonzero elements of \(S'\) is a finite-dimensional domain
over \(\operatorname{Frac}(S')\), hence the full field.
The free graded submodule on these elements has leading Hilbert
coefficient \(pk\). As before \(k=4\) would give \(S'=S\).
Thus \(p=1\), \(k=3\), and
\[
                         H_S-H_{S'}=\frac{t^3}{(1-t)^{10}} .
\]
There is exactly one new cubic generator \(g\) modulo \(S'\),
and it satisfies \(\ell g\in S'_4\) for a nonzero linear form
\(\ell\). Indeed multiplication by \(g\) sends the
eleven-dimensional \(S_1\) to the ten-dimensional \(S_4/S'_4\).
On a
general affine Pl\"ucker chart over \(\mathbf C(W)\) the cubic
equation is affine linear in the three fiber variables. Its linear
part is nonzero: otherwise the equation would be a function of
the base alone, either zero in \(R\) or preventing dominance over
\(W\). Thus the generic fiber on this affine chart is an affine
plane.

\section{Excluding the quadratic images for quartic ninefolds}\label{sec:quartic-nine-higher-images}

\subsection{Excluding a three-dimensional span of the first Pfaffians}
\begin{proposition}\label{prop:next-quartic-pfaffian-rank-four}
For a quartic ninefold replacement with Grassmannian-cone image,
the span of the first Pfaffian coefficients in the choice of
Proposition~\ref{prop:dec2-first-pfaffians} has dimension at least
four in $S_2/T_2$.
\end{proposition}
\begin{proof}
Assume its dimension is three. In the notation of
Proposition~\ref{prop:dec2-spinor-incidence-quartic}, let $V$ be
the five-space, choose independent $b_i\in V^*$ for $1\le i\le3$,
and write the first Pfaffian vector as $q+\sum b_i y_i$, with
$q$ a vector of quadratics in the degree-one coordinates.
Let $D=T[y_1,y_2,y_3]\subset S$, and put $M=\operatorname{Proj}D$.
Remove the vertex and the rank-drop locus on $W$ where the $b_i$
span a space of dimension at most one on the Pl\"ucker plane.
On the resulting smooth open $W^\circ$, the equation
$A(q+\sum b_iY_i)=0$ is an affine system of rank two in the three
coordinates $Y_i$ of $\mathcal O_W(2)^{\oplus3}$. Its affine
solution space $M^\circ\to W^\circ$ is a torsor for $\mathcal O_W(1)$.
Indeed the kernel direction is given by the triple $c$ of signed
linear Pl\"ucker minors of the restrictions of $b_1,b_2,b_3$.
They do not vanish simultaneously there. Evaluating the $y_i$
dominates this torsor, since at least one first coefficient is
transcendental over the linear coordinate field by
Proposition~\ref{prop:dec2-grass-all-primes} and
Proposition~\ref{prop:dec1-stein-hyperplane-fibers}.

Consequently $D$ is the reduced affine incidence algebra for the
three effective unknowns. Indeed the incidence variety for the
five unknowns used in Proposition~\ref{prop:dec2-spinor-incidence-quartic}
is irreducible, with dense full-rank locus, and the last two
unused quadratic variables are free after a constant change of
basis. Its ring $R$ there equals $D[Y_4,Y_5]$, with the two
indeterminates of degree two. By the factoriality argument in
Corollary~\ref{cor:dec2-incidence-quadratic-base}, $D$ is therefore
normal and factorial, and
\[
 H_D(t)=\frac{(1-t^2)^2}{(1-t)^{11}},\qquad
 H_S(t)/H_D(t)=\frac{1+t^2}{1-t^2}.
\]
In particular $M^\circ$ is open in $M$. The three independent
linear sections in $c$ have no common prime divisor on $Z$ by
Proposition~\ref{prop:dec1-lower-seven}. Thus on a big smooth open
$U\subset Z$ the evaluations give a morphism $p:U\to M^\circ$.
The canonical bundle of $M^\circ$ is the pullback of
$\omega_{W^\circ}(-1)=\mathcal O_{W^\circ}(-7)$. The equality
$p^*\omega_{M^\circ}\simeq\omega_U$ and the eighth wedge of the
differential give a nonzero reflexive vector field $v$ on $Z$.
It annihilates the function field of $M$.

Let $K$ be the connected kernel of the projective action of
$\operatorname{Aut}^0 Z$ on $S_1$. This connected normal subgroup
is reductive. Lift its action to the section ring with $S_1$
fixed pointwise. To obtain this lift, take the scalar extension
of the group of automorphisms preserving the divisorial sheaf of
$L$ by its choices of sheaf isomorphism. The action on $S_1$ is
scalar on this extension of $K$, and twisting the lift by the
inverse scalar character to the power of the section degree makes
the grading action descend to $K$. For $g\in K$ the differences
$g y_i-y_i$ have the form $c_i a_g$ since they solve the homogeneous
rank-two system. The rational degree-one element $a_g$ is regular
in codimension one on $Z$ by the absence of any common prime of
the minors. Thus $a_g\in S_1$. The map $g\mapsto a_g$ is algebraic,
by multiplying by any nonzero $c_i$ and inverting the resulting
linear inclusion of $S_1$ into $S_2$. It is a homomorphism to the
additive group of $S_1$ since $S_1$ is fixed. Reductivity forces
it to vanish. Hence $K$ fixes $M^\circ$ pointwise.

The vector $v$ is in $\operatorname{Lie}(K)$ and its line is
$K$-invariant. Indeed the map to $M^\circ$ is $K$-invariant,
the differential construction is functorial, and two
identifications $p^*\omega_{M^\circ}\simeq\omega_U$ differ by a
scalar on a big open of the normal projective $Z$. An invariant
line in a reductive Lie algebra is central. The smallest algebraic
subtorus $G$ of the center whose Lie algebra contains $v$ has
dimension one: it acts faithfully on $Z$ and the general fibers
over $M$ have dimension one. A faithful torus acting on an
integral variety has trivial generic stabilizer, as follows by
diagonalizing its action in an equivariant projective embedding.

No prime divisor of $Z$ is pointwise fixed by $G$. In fact the
isomorphism $p^*\omega_{M^\circ}\simeq\omega_U$ is $G$-equivariant.
Both the vector $v$ and the differential giving it are invariant,
and the nonzero differential on the dense maximal-rank locus
determines the action on the isomorphism. At a general smooth
point of a putative fixed divisor the canonical fiber would thus
have weight zero. Its tangent representation has eight zero
weights and a nonzero normal weight, a contradiction. The normal
weight is nonzero since trivial tangent action of a torus at a
smooth fixed point implies trivial action by equivariant formal
parameters. The vector field $v$ is thus nonvanishing at every
codimension-one generic point, so $p$ has differential of full
rank eight at each such point. In particular, if an upstairs
prime maps into a divisor on the smooth $M^\circ$, it dominates
that divisor and has ramification index one. The first assertion
follows by restricting the differential to the prime, dropping
its rank by at most one. The second follows by taking the pullback
of a smooth local parameter at the target generic point.

Let $\widetilde M$ be the normalization of $M$ in
$\mathbf C(Z)^G$, and denote its finite degree by $e$. This field
is the relative algebraic closure of $\mathbf C(M)$ in
$\mathbf C(Z)$. Indeed the invariant field has transcendence
degree eight and is relatively algebraically closed for the
connected torus action. A resolution of the graph factors through
$\widetilde M$ and its general geometric fiber there is integral,
with an open free $G$-orbit. Write $\Gamma$ for its closure on $Z$.
Let $P$ be the pullback of $L$ and $J$ half the pullback of
$\mathcal O_M(2)$ to this resolution. The latter sheaf is very
ample since the second Veronese of $D$ is generated in degree one.
Both divisors are nef and $P-J$ is effective up to rational
linear equivalence. For this last assertion one can take all
products of the quadratic subsystem in a Cartier multiple of
$2L$. The Hilbert series gives $(\mathcal O_M(2))^8=2^8\,4$.
By the projection formula,
\[
                  4=P^9\ \ge P J^8=4e (L\cdot\Gamma).
                  \tag{*}
\]

For a general orbit the two endpoints of its normalized orbit
closure $\mathbf P^1$ both belong to the base of the degree-one
system. Otherwise one endpoint gives a fixed point where the
subsystem $D_2$ is generated, since it contains the squares of
the linear sections. The rational map to $M$ is defined there
and has the same value as on the general orbit by separatedness.
Here nonvanishing of a reflexive linear section means its Cartier
power is nonzero, which also makes $L$ invertible locally.
The fixed locus has dimension at most seven by the absence of
fixed divisors, and so the images of its points on the domain of
this map cannot contain a general point of the eightfold $M$.
This proves the endpoint assertion.
The local Cartier indices of $L$ are at most six by the
normalized-volume bound in
Proposition~\ref{prop:nov30-quartic-nine-index-cover}, which
bounds the degree of the local cyclic cover at every closed point.
Take a linear section nonzero at the generic point of $\Gamma$.
On each of the two branches at the endpoints its effective
rational section divisor has order at least $1/6$: raise it to
the local Cartier index, giving a regular section vanishing at
that point. Positivity and summation on the normalization give
$L\cdot\Gamma\ge1/3$. It follows from $(*)$ that $e\le3$.

Every prime $\Delta$ in $M^\circ$ is dominated by a prime on $Z$.
Indeed its closure is cut out by a positive-degree homogeneous
prime element of $D$ by factoriality. Its section divisor on $Z$
is nonempty by ampleness, and its components map into $\Delta$
on the big open $U$. Full differential rank makes each such
component dominant over $\Delta$. If $e=2$, ramification over any
divisor of $M^\circ$ would give a unique prime over it in the
finite quadratic extension and that prime would have ramification
index two. The index-one domination just proved rules this out.
The finite map is thus \'etale over $M^\circ$ by purity of the
branch locus. But $M^\circ$ is simply connected: $W^\circ$
is the total space of $\mathcal O_G(1)$ restricted to the
Grassmannian minus a Schubert variety of codimension two. Removing
that locus preserves the trivial topological fundamental group
of the smooth Grassmannian by general position. Both the line
bundle and the affine-line torsor have contractible fibers and
are locally trivial in the analytic topology. The connected cover
over $M^\circ$ therefore has degree one, excluding $e=2$.

If $e=3$, equality holds in $(*)$. Then the geometric base of
the subsystem $D_2$ on $Z$ is finite. Here are details of this
mixed-intersection consequence. Take the ordinary blowup of the
ideal generated by all products of a sufficiently divisible
number of elements of $D_2$ in a Cartier trivialization. On the
normalized blowup $2m(P-J)$ is the effective exceptional divisor
for the ideal of the $m$-fold products, and $2mJ$ is relatively
ample over $Z$. The numbers $P^{9-i}J^i$ are nonincreasing.
If a base component has codimension $c\le8$, an exceptional
divisor above its generic point has general fiber dimension
$c-1$ over it. Indeed the pullback ideal is an invertible
nontrivial ideal there, and its zero divisor covers the center,
by properness and the definition of the blowup. One of its
components dominates the base component. Relative ampleness gives
positive top degree of $J$ on its general fiber, and ampleness
of $L$ gives positive degree of the center. Thus
$P^{9-c}J^{c-1}(P-J)>0$. All other exceptional contributions
are nonnegative. This contradicts equality of the mixed
intersections through $i=8$ and proves finiteness.

Every prime of $\widetilde M$ above $M^\circ$ is now dominated
by a prime on $Z$. To see this, take the normalized graph of the
subsystem, which agrees with $Z$ off its finite base, and use its
proper morphism to $\widetilde M$. Localize this family over the
discrete valuation ring at the proposed prime. The family is
flat since the total space is integral and dominant. Its fibers
are projective curves and the degree of a fixed Cartier multiple
of $P$ on the generic fiber is positive. The degree on the
special fiber cycle is the same by flatness (or specialization of
intersections with Cartier divisors). Thus some special-fiber
component has positive $P$-degree. Its image in $Z$ is a curve
and cannot be contained in the finite base. The closure of this
component in the graph is a prime divisor dominating the proposed
prime. It meets the isomorphism locus of the graph over $Z$, so
it gives a divisor on $Z$ with the same property. Full differential
rank on $U$ now excludes ramification at every prime of
$\widetilde M$ above $M^\circ$. Purity and simple connectedness
give the contradiction also for $e=3$.

It remains to consider $e=1$. On the affine section cones the
invariant fraction field is then $\operatorname{Frac}(D)$:
adjoin the same invariant degree-one coordinate to the projective
fields as in Proposition~\ref{prop:dec1-stein-hyperplane-fibers}.
No affine prime divisor is fixed or contracts to codimension two.
Indeed the complement of the charts given by nonzero degree-one
coordinates has codimension at least two on the source and,
on these charts, the maps and actions are the products of the
projective maps and actions with an invariant multiplicative
coordinate. Every prime divisor on the affine target is dominated:
its prime equation in $D$ remains a nonunit of $S$, and any
component of its zero divisor must dominate the given prime.
Here all units of the connected positively graded domain $S$
are constants. Normality then gives $S^G=D$ by absence of poles
in $\operatorname{Frac}(D)$.

The weight-module argument in
Proposition~\ref{prop:dec1-factorial-ci-image} now applies with
base ring $D$. More explicitly, over $\operatorname{Frac}(D)$
each torus weight space has dimension at most one. Both signs
occur since a one-sided semigroup would have a fixed point on
the affine generic curve, whose closure is a fixed divisor.
Consequently this curve is a multiplicative group and its
eigenfunctions have no horizontal zeros or poles. For a rational
eigenfunction $s$ in any occurring weight, membership of $sf$ in
$S$ for $f\in\operatorname{Frac}(D)$ is tested by the orders at
primes above the base primes. All primes are dominated, so this
gives a divisorial fractional ideal of $D$. Factoriality makes
each weight module free of rank one over $D$ with a homogeneous
generator. The quotient $H_S/H_D=(1+t^2)/(1-t^2)$ has two terms
in each strictly positive even degree. The same elementary
weight count in Proposition~\ref{prop:dec1-factorial-ci-image}
therefore gives $S=D[y_+,y_-]/(y_+y_--f_4)$ with
$\deg y_\pm=2$. In detail the two weights occurring first in
degree two cannot have the same sign: $0,a,\ldots,ja$ and
$(j-i)a+ib$, $1\le i\le j$, already give $2j+1$ weights for
$0<a<b$, leaving no room for a negative weight for large $j$.
For opposite weights $a,-b$ their pure powers exhaust the
$2j+1$ available weights, and their product in degree four
forces $a=b$. Powers of the two sections then give an injective
map from the indicated presentation, whose Hilbert series is the
full Hilbert series of $S$. This makes $S$ generated in degrees
at most two, contrary to Proposition~\ref{prop:dec2-spinor-incidence-quartic}.
This final contradiction proves the claim.
\end{proof}

\subsection{Threefold fibers over the Grassmannian}

\begin{proposition}\label{prop:dec2-grass-projected-degree}
Let \(V\) be the closure of a general threefold fiber of the
projection to \(G\). Then \(L^3\cdot V<27/80\). On a resolution
of \(V\), every subsystem of the restricted degree-two sections
with three-dimensional image has mapping degree times image
degree at most two.
\end{proposition}
\begin{proof}
Use \(b_i=P^{9-i}J^i/4\) for the ten-section projection to \(G\),
with \(P\) the pullback of \(L\) and \(J\) the moving class.
Then \(b_6=5L^3\cdot V/4\). By the mixed bound for this system
the lifted local thresholds and its base codimension are at least
\(\theta_P=(7/10)\sum_{i=0}^6 b_i\).
We must have \(\theta_P\le4\). Otherwise the base would have
codimension at least five, while
Proposition~\ref{prop:dec2-grass-singular-base} would supply a
codimension-five base component in the singular locus. At a
general point of that component \(L\) is Cartier by
Proposition~\ref{prop:dec2-index-codim-five}. The transverse
fivefold is therefore Gorenstein, and the threshold on this
transverse slice still exceeds four by restriction of a general
log resolution along the stratum. Four general members cut out
a klt curve, by the successive adjunction argument in
\cite[Lemma~4.5]{kePanWuCubics} with four cuts and strict
threshold. This curve is regular. Lifting a parameter of the
curve and the four equations makes the transverse embedding
dimension at most five, contradicting singularity at a general
point of the component.

The projected base already has codimension at least three, so
\(b_0=b_1=b_2=1\). If \(b_6\ge27/64\), log-concavity gives
\(b_{2+j}\ge(3/4)^j\) for \(1\le j\le4\), since
\((3/4)^4<27/64\). Thus
\((7/10)(3+\sum_{j=1}^4(3/4)^j)>3\) and \(b_3=1\).
Log-concavity from \(b_3\) would now give
\[
        \theta_P\ge(7/10)(4+3/4+9/16+27/64)>4 ,
\]
a contradiction. Finally the moving class of a subsystem of
the restricted quadrics is bounded above by \(2P|_V\) on a
common resolution, with effective difference, both classes nef.
Its top degree is an integer at most \(8 L^3\cdot V<27/10\).
\end{proof}

Suppose in addition that \(|2L|\) has nine-dimensional image, and
let \(H_3\) be its cubic equation in the incidence ring \(R\) of
Corollary~\ref{cor:dec2-incidence-quadratic-base}. Then
\[
                           H_3(0,z,y)=0 .
\]
In fact fix a general plane of \(G\) with representative \(A_0\).
Put \(A=sA_0\); after discarding the common scalar \(s\), the
incidence equations give \(A_0(By+c(sA_0,z))=0\). Since
\(A_0B\) has rank two, these equations express two \(y\)-variables
as quadratics in \(s,z\) and linear combinations of the
remaining three. The parameter space of the resulting system is
\(\mathbf P(1,1,2,2,2)\), and its second Veronese embeds it as
the quadric specified by \((sz)^2=(s^2)(z^2)\). All these
Veronese coordinates are obtained from the restrictions of
\(S_2\), because the three remaining \(y\)-coordinates and
the binary quadratics in \(s,z\) are their restrictions.
On \(s\ne0\) the cubic cuts out the generic threefold fiber image.
The geometric closure of that fiber is irreducible and has
multiplicity one along its generic point, since the quadratic
map is birational and the image dominates \(G\) separably.
If the substituted cubic were not divisible by \(s\), there
would be no component at \(s=0\), so its Veronese image would
have degree \(3\): the degree is
\(2^3\cdot3/(1\cdot1\cdot2^3)=3\) from the weighted Hilbert
series. Proposition~\ref{prop:dec2-grass-projected-degree} excludes
this. Thus the cubic vanishes at \(s=0\) for all
\((z,y)\) satisfying \(A_0(By+c(0,z))=0\) for a general \(A_0\).
These conditions range over a dense subset of \((z,y)\): with
\(z,y\) general, a general two-plane in the kernel of the indicated
covector has independent restrictions from \(\operatorname{im}B\).
More formally the incidence with the general plane is an affine
bundle in the \(y\)-coordinates over the open plane locus and
the \(z\)-line, and its image in \((z,y)\) contains this dense
set. The polynomial vanishing thus holds identically.

\subsection{The quadratic curve over the general point of the cone}

\begin{proposition}\label{prop:next-q9-one-place}
Suppose that the quadratic image of a quartic ninefold replacement
with Grassmannian-cone hyperplane image $W$ has dimension eight.
Over the geometric generic point of $W$, the affine curve described
by the quadratic coordinates has normalization $\mathbf A^1$.
Its projective normalization has exactly one place above infinity.
If its projective degree in the affine quadratic coordinates and a
constant coordinate is $m$, and $e$ is the degree of the relative
algebraic closure of the quadratic image field in $\mathbf C(Z)$,
then $em\le6$.
\end{proposition}
\begin{proof}
Choose one nonzero linear coordinate to trivialize the quadratic
coordinates over $\mathbf C(W)$. They define an integral affine
curve whose function field is the subfield of $\mathbf C(Z)$
generated over $\mathbf C(W)$ by the quadratic section ratios.
Since the resolved general surface fiber over $W$ is integral,
$\mathbf C(W)$ is relatively algebraically closed also in this
intermediate field. Thus the curve is geometrically integral.
The smooth projective model of its geometric generic fiber has
genus zero: the general curve is dominated by the rational general
surface fiber in Proposition~\ref{prop:dec1-grassmannian-fibers},
and genus is constant on the smooth proper family of resolved
curves over a sufficiently small open on the base.

Let $\widetilde W$ be the usual resolution of the vertex by a
projective-line bundle over $\operatorname{Gr}(2,5)$. To compactify
the affine curve family, use projective completion of the vector
bundle of the five quadratic coordinates on $\widetilde W$, each
in the square of the hyperplane line bundle, and take the closure
of the generic image. Resolve this projective closure and the
inverse image of infinity to get a smooth projective eightfold
$f:Y\to\widetilde W$. Let $D$ be the sum with coefficient one
of all horizontal components of the inverse image of infinity.
It has simple normal crossings on the resolution. On a general
fiber it cuts precisely the distinct punctures of the affine
normalization. This holds since a smooth proper birational model
of an integral curve is its projective normalization.

Write $P$ for the pullback of $L$ on a common smooth resolution
of $Z\dashrightarrow Y$, with projection $p$ to $Y$. Every
prime of $Z$ has generic image away from the cone vertex by
Proposition~\ref{prop:dec1-lower-seven}. On a chart around that
generic point with a nonzero linear coordinate $x$, the ratios of
the quadratic sections by $x^2$ are regular. The image in the
projective completion lies in the affine vector bundle there.
Consequently $p^*D$ pushes forward to zero on $Z$ and
$P^8p^*D=0$. Likewise the pullback of $K_{\widetilde W}$ has
normalized mixed degree $-6$, since $K_W=-6H$ off the vertex
and every source prime avoids it.

Suppose the projective normalization has at least two punctures.
Then $K_{Y_w}+D|_{Y_w}$ has nonnegative degree and a nonzero
section on the general rational fiber. We use the relative
pseudo-effectivity conclusion of F.~Campana and M.~P\u{a}un,
\emph{Foliations with positive slopes and birational stability of
orbifold cotangent bundles}, statement and proof of Theorem~3.4
in \href{https://arxiv.org/abs/1508.02456v4}{arXiv:1508.02456v4}.
Here is how its global hypothesis is obtained. The direct image
$f_*\mathcal O_Y(K_Y+D)$ is a nonzero coherent sheaf, by the
generic-fiber section just found. After twisting by a sufficiently
ample line bundle on $\widetilde W$ it has a nonzero global
section. Represent the twist by a general very ample divisor $A$.
Then $K_Y+D+f^*A$ is pseudo-effective and $D+f^*A$ is reduced
with simple normal crossings by Bertini. The horizontal part of
this boundary is $D$. The cited theorem gives pseudo-effectivity
of $K_{Y/\widetilde W}+D$ after subtraction of the effective
divisorial ramification term, hence also before that subtraction.

Intersect its pullback with the nef product $P^8$. It follows
that $P^8p^*K_Y/P^9\ge-6$. On the other hand the saturated
horizontal cotangent of $Z\dashrightarrow Y$ has rank eight.
Taking its determinant and using semistability of the reflexive
cotangent of the K-polystable ninefold gives
\[
                  P^8p^*K_Y/P^9\le -8(7/9)=-56/9.
\]
Indeed the actual differential adds an effective divisorial
degeneracy term to the pullback determinant before semistability
is applied. Since $-56/9<-6$, there cannot be two punctures.
A nonconstant affine curve has at least one point at infinity,
so its normalized geometric affine curve is $\mathbf A^1$.

Replace the image family by its normalization in the relative
algebraic closure of its field in $\mathbf C(Z)$ and resolve.
The same argument applies to this model: the geometric generic
curve is integral since the constant base field is still relatively
algebraically closed, and the map from the resolved rational
surface dominates it. All components above infinity still project
to infinity, so they make zero contribution after pullback to
$Z$. Thus also on this finite degree-$e$ cover the geometric
generic projective normalization has a single point above infinity.

Now take a general linear combination $s$ of the quadratic
coordinates and quadratic monomials in linear coordinates. On the
generic projective image curve its zero divisor consists of $m$
distinct affine points in the smooth locus, by Bertini and the
definition of the projective degree. On the finite cover it
therefore consists of $em$ distinct points. Include in the
resolution the reduced sum $D_0$ of the horizontal closures of
these points, arranging simple normal crossings with $D$. The
pullback of $D_0$ has normalized mixed degree at most two.
In fact on a chart $x\ne0$ the rational function $s/x^2$ is
regular on the chart of the original image. Its pullback to the
cover vanishes along $D_0$, and no prime of $Z$ sees a pole in
this trivialization. Hence after pushforward to $Z$ the pullback
of $D_0$ is bounded by the section divisor of $s$ in $|2L|$.

On the general projective rational curve, the boundary
$D+(em)^{-1}D_0$ now has degree two. After clearing denominators
the previous twisting and log pseudo-effectivity argument applies
to this rational boundary, yielding
\[
                  P^8p^*K_Y/P^9\ \ge\ -6-\frac{2}{em}.
\]
The normalized first-jet argument of
Lemma~\ref{lem:dec1-hyperplane-atiyah-slope} gives the upper bound
$-63/10$ for the left-hand side. Indeed on the smooth target take
the pullback of the hyperplane bundle from $\widetilde W$, whose
pullback on a big open of $Z$ equals $L$. Pullback of normalized
first jets is generically injective by dominance, has rank nine
and determinant $p^*K_Y$, exactly as in that lemma. These bounds
give $em\le20/3$, hence $em\le6$.
\end{proof}

\begin{proposition}\label{prop:next-q9-exclude-quadratic-line}
Under the assumptions of Proposition~\ref{prop:next-q9-one-place},
the affine quadratic image over the geometric generic point of $W$
cannot be a line. Thus its projective degree $m$ satisfies $2\le m\le6$.
\end{proposition}
\begin{proof}
Suppose the image is a line. Take a basis $y_1,\ldots,y_5$ of the
five new quadrics modulo $T_2$. Its direction is a point of
$\mathbf P^4(\mathbf C(W))$, since the image line is defined over
$\mathbf C(W)$. Write its coordinates as forms $c_i\in T_k$ of
equal degree with no common factor, using the factorial graded
ring $T$. The five forms are linearly independent over $\mathbf C$:
a constant relation among them would make the corresponding
combination of the $y_i$, divided by the square of a linear
trivialization, constant on the line. This combination would
belong to $S\cap\operatorname{Frac}(T)=T$ by
Proposition~\ref{prop:dec2-grass-all-primes}. In particular $k\ge1$.

Write $R$ for the effective divisorial differential degeneracy of
the map to $W$ in Proposition~\ref{prop:dec1-grassmannian-fibers},
and put $\rho=L^8R/L^9\le2/5$. On $Z$ write $F$ for the common
fixed divisor of the $c_i$, and $x=k-L^8F/L^9$ for the normalized
degree of their moving part. Every component of $F$ maps into
codimension at least two on $W$, since the $c_i$ have no common
prime there. Hence each such component occurs with coefficient at
least one in $R$, by the differential rank argument of
Proposition~\ref{prop:dec2-grass-all-primes}.

For the foliation $\mathcal F$ tangent to the map to $W$ one
has normalized determinant degree $1+\rho$, by
Proposition~\ref{prop:dec1-grassmannian-fibers}. Differentiation
along $\mathcal F$ gives a map into $\mathcal O_Z(2L)^{\oplus5}$
on a big smooth open. In local linear trivializations pulled back
from $W$ it differentiates the five functions $y_i$, and changes of
trivialization are constant in the foliation directions. The image
has rank one and its generic direction is $[c_1:\cdots:c_5]$.
Its saturation in $\mathcal O_Z(2L)^{\oplus5}$ has normalized
degree $2-x$: locally at each prime one divides the common
coordinate vector by its common vanishing factor. The rank-one
kernel in $\mathcal F$ has degree at most $7/9$ by semistability
of $T_Z$. Taking normalized degrees therefore gives
\[
                         2-x\ \ge\ 2/9+\rho .
\]

Let $E$ have fixed multiplicity $a>0$ in $F$. On a resolution let
$P$ be the hyperplane pullback and $J$ the free moving class of
the $c_i$. Thus $P^8J/P^9=x$. On $0\le t\le a$, $kP-tE$
dominates the nef class $(1-t/a)kP+(t/a)J$. Integrating the
first two binomial terms gives
\[
             10k/7\ \ge\ 10S_{kL}(E)\ \ge\ a(1+x/k)
\]
by K-semistability and $A_Z(E)=1$. Summing the resulting bounds
on the fixed multiplicities against prime degrees gives
\[
       k-x\le \frac{10k}{7(1+x/k)}\rho,\qquad
       x^2\ge (3/7)k^2 .
\]
Since $x\le16/9$, this leaves $k\le2$. If $k=2$, the inequality
$x\ge2\sqrt{3/7}>6/7$ gives $a<2$ for every fixed multiplicity.
They are integers, so $F\le R$. But $2-x=\deg_L F\le\rho$
contradicts $2-x\ge2/9+\rho$.

For $k=1$ write the first Pfaffian vector as $q+By$, with $B$ a
constant matrix. Differentiation in the line direction yields
$A_0 Bc=0$. Every degree-one solution of $A_0s=0$ in $T$ is
$A_0\wedge b$ for a fixed $b\in\bigwedge^2V$, as in the
Borel--Weil calculation of Proposition~\ref{prop:dec2-first-pfaffians}.
Indeed the coefficient of the vertex variable would annihilate
all two-planes and must be zero. On the Grassmannian the map
$\bigwedge^2V\to H^0(G,Q^*(1))$ is surjective, where $Q$ is its
universal quotient bundle. The span dimension of the entries of
$Bc$ is $\operatorname{rank}B\ge4$ by
Proposition~\ref{prop:next-quartic-pfaffian-rank-four} and
independence of the $c_i$. A decomposable $b$ gives only three
independent entries in $A_0\wedge b$. A rank-four bivector gives
five: for $b=e_1\wedge e_2+e_3\wedge e_4$ the four entries
containing $e_5$ are independent and the remaining entry is
nonzero. Thus $B$ is invertible and $b$ has rank four.

The common zero locus of $c$ has codimension three in the affine
cone. In fact on $G$ a plane $U$ satisfies $(\bigwedge^2U)\wedge
b=0$ precisely when it is a Lagrangian plane in the support
four-space $E_4$ of $b$. A plane not contained in $E_4$ fails the
equation since $E_4\to\bigwedge^3E_4$, $u\mapsto u\wedge b$,
is injective. For planes in $E_4$ the equation is the symplectic
isotropy condition. The Lagrangian Grassmannian has dimension
three inside the six-dimensional $G$, and adjoining the vertex
variable preserves the codimension.

Along the affine image line all coordinates $c_i y_j-c_j y_i$
are constant after a common linear trivialization. They are
therefore polynomials in $T_3$ by the intersection property of
Proposition~\ref{prop:dec2-grass-all-primes}. Their two-form
$c\wedge y$ satisfies $c\wedge(c\wedge y)=0$. The ideal of
the five $c_i$ has grade three in the Cohen--Macaulay Pl\"ucker
cone ring $T$. The Koszul cochain complex with differential
$c\wedge$ is consequently exact in degrees below three
by the grade sensitivity of Koszul cohomology. It follows that
$c\wedge y=c\wedge q'$ for a vector $q'$ of homogeneous
quadratics in $T$. Thus $y_i-q'_i=c_i u$ for a rational
homogeneous $u$ of degree one on $Z$. The five independent linear
sections $c_i$ share no prime on $Z$ by the three-section
criterion, so $u$ is regular in codimension one and belongs to
$S_1=T_1$. This would put all $y_i$ in $T_2$, a contradiction.
\end{proof}

\begin{corollary}\label{cor:next-q9-conic-or-cubic}
With the same assumptions, the quadratic image field is relatively
algebraically closed in $\mathbf C(Z)$. Over a general point of $W$
the quadratic image curve has projective degree two or three.
\end{corollary}
\begin{proof}
For a general threefold fiber closure $V$ over the Grassmannian,
the estimate in Proposition~\ref{prop:dec2-grass-projected-degree}
in fact gives $L^3V<1/3$. Write $b_i$ for the normalized mixed
intersections of the ten-section subsystem as there. The proof
gives $b_0=b_1=b_2=1$, $b_6=5 L^3V/4$, and
$\frac{7}{10}\sum_{i=0}^6 b_i\le4$. If $b_6\ge5/12$, then
$b_6>(3/4)^4$ and the first log-concavity step in that proof
gives $b_3=1$. Set $t=373/500$. Since $t^3<5/12$,
log-concavity now bounds $b_{3+j}\ge t^j$ for $1\le j\le3$.
But $t+t^2+t^3>12/7$ contradicts the mixed-intersection upper
bound. This proves the sharper estimate.

On $V$, the quadratic sections have surface image $\Sigma$,
projectively compactified by the full restricted subsystem.
Indeed their image maps onto the ruling $\mathbf P^1$ of $W$
and has curve fibers. For two distinct general points on that
ruling, multiply the two linear equations cutting out the points
on the ruling. Their product defines a hyperplane section of
$\Sigma$, whose support contains both degree-$m$ curve fibers.
The section is nonzero on a general fiber. Thus
$\deg\Sigma\ge2m$.
On a common resolution of $V$ the moving pullback $J$ of the
quadratics and the pullback $P$ of $L$ are nef and $2P-J$ is
effective up to rational linear equivalence. The Stein degree is
$e$ as in Proposition~\ref{prop:next-q9-one-place}. In fact on a
dense open of the global eightfold image there are $e$ geometric
fiber components, permuted transitively, and a general point of
the slice over the Grassmannian lies in this open. At the generic
point of each such component the resolution of $V$ maps
birationally to its curve closure on $Z$. Denote this curve
class numerically by $C$ for the degree computation.

One has $L\cdot C\ge1/6$. Choose a linear section nonzero at the
generic point and restrict its rational Cartier section divisor
to the normalization of $C$. This gives an effective divisor of
strictly positive degree by ampleness. Each positive coefficient
is at least $1/6$ since the local Cartier index is at most six.
The projection formula and the mixed bounds therefore give
\[
       \frac83>8L^3V\ge (2P)J^2
          =2e(\deg\Sigma)(L\cdot C)\ge \frac23em.
\]
Hence $em\le3$. The line exclusion gives $m\ge2$, forcing
$e=1$ and $m\in\{2,3\}$.
\end{proof}

\begin{theorem}\label{thm:next-q9-no-eight-quadratic}
If the limiting hyperplane image of a K-polystable quartic ninefold
replacement is the Grassmannian cone, its complete quadratic system
has nine-dimensional image.
\end{theorem}
\begin{proof}
Suppose the quadratic image has dimension eight. Let $Y$ be its
smooth projective model in Proposition~\ref{prop:next-q9-one-place}
and $m$ the projective degree there. By
Corollary~\ref{cor:next-q9-conic-or-cubic} the degree $e$ is one
and $m\in\{2,3\}$. On a smooth resolution $\pi:\widehat Z\to Z$ of
the map $p:\widehat Z\to Y$ the general fiber $C$ is therefore a
smooth integral curve. Put $P=\pi^*L$. In the notation of that
corollary, the estimate $8/3>2(\deg\Sigma)(P\cdot C)$ and
$\deg\Sigma\ge2m$ give
\[
                      0<P\cdot C<2/(3m)\le1/3 .
\]
Choose a general linear section $s$ nonzero on $C$. Its rational
Cartier section divisor pulled back to the smooth curve is
effective, of degree $P\cdot C$. At each point of the curve the
positive coefficients are integral multiples of the reciprocal
of the local Cartier index at its image, which is at most six.
There is consequently precisely one positive coefficient, equal
to $1/j$ with $4\le j\le6$. Indeed the next possible numerator
at any length at most six gives a coefficient at least $1/3$.
All nonzero general linear sections have proportional restrictions
on $C$, since their ratios come from $Y$.

Take the resolution to principalize the products of linear and
quadratic sections in Cartier multiples and to have simple normal
crossings exceptional set. Over the regular locus of $Z$ where
the degree-one system is generated, the map to the projective
quadratic image is regular. Any exceptional divisor centered off
the degree-one base cannot meet the general $C$: the center has
dimension at most seven, and the map to the quadratic image is
defined generically there. This argument applies also to a
singular center outside the base since a nonvanishing local
Cartier power of a linear section trivializes $L$ near the point.
Additional modifications to lift the map to $Y$ preserve the
claim over its general point.

Outside the exceptional divisors the divisor of $s$ does not
meet the general $C$: its non-base part maps to the hyperplane
cut out by $s$ in $W$. The hyperplane system has no fixed primes
on $Z$. The sole point in its restriction to $C$ consequently
lies on an exceptional divisor $E$ over the base. Every
exceptional divisor dominating $Y$ has positive fixed order
there, by taking the products in a local Cartier multiple.
Intersecting the general curve with the exceptional divisors
is transverse by generic smoothness in characteristic zero.
Two exceptional divisors cannot meet at a point of the general
fiber since their intersection has dimension at most seven.
Thus $E$ is the unique exceptional divisor meeting $C$, it meets
$C$ once, and $a=\operatorname{ord}_E(|L|)=1/j$. Adjunction on
the smooth fiber gives
\[
             2g(C)-2=K_{\widehat Z}\cdot C
                      = -7/j+A_Z(E)-1 .
\]
If $g(C)=0$, then $A_Z(E)/a=7-j\le3$. This contradicts the
strict hyperplane threshold in
Proposition~\ref{prop:dec2-ninefold-tangent-pencil}. Hence $g(C)\ge1$.

Let $\mathcal F$ be the saturated rank-one foliation on $Z$ tangent
to the fibers. Its Weil canonical class $K_{\mathcal F}$ is
pseudo-effective on a small projective $\mathbf Q$-factorialization
$Z'$ of $Z$, identifying divisors in codimension one. We recall
the precise positivity input. For an algebraic foliation on
a smooth projective variety whose smooth general fiber has
pseudo-effective canonical class, the determinant of the
dual foliation is pseudo-effective. This follows by taking a neat
model of the map, applying the relative pseudo-effectivity result
to the relative canonical class minus the multiple-fiber divisor,
and pushing down. See Campana--P\u{a}un,
\emph{Foliations with positive slopes and birational stability of
orbifold cotangent bundles}, Remark~3.3, Remark~3.6 and the proof
of Theorem~4.7, in particular formula~(39) of
\href{https://arxiv.org/abs/1508.02456v4}{arXiv:1508.02456v4}.
This argument uses precisely the pseudo-effectivity on the general
fiber, here furnished by $g(C)\ge1$. Apply the result on a smooth
resolution of $Z'$ and push down. The pushforward agrees in
codimension one with $K_{\mathcal F}$, and preserves
pseudo-effectivity on the $\mathbf Q$-factorial target.

On the big open of $Z$ where the rational map is regular and the
relevant sheaves are locally free, saturating the pullback of the
canonical bundle of $Y$ in the eighth exterior cotangent gives
the relative determinant formula
\[
 K_{\mathcal F}=K_Z-\pi_*p^*K_Y-R_p
\]
with $R_p$ an effective integral divisor of differential degeneracy.
Let $\deg_L A=L^8A/L^9$ on $Z$, and put $k_Y=P^8p^*K_Y/P^9$.
Proposition~\ref{prop:next-q9-one-place} gives $k_Y\ge-6-2/m$
by its boundary argument. Therefore
\[
        0\le\deg_L K_{\mathcal F}
          =-7-k_Y-\deg_L R_p
          \le -1+2/m\le0 .
\]
The left-hand inequality follows also on $Z'$ by intersecting
the pseudo-effective class with the eighth power of the nef
pullback of $L$. Thus $\deg_L K_{\mathcal F}=0$.

The small $\mathbf Q$-factorialization $Z'$ is of Fano type.
Indeed it is klt and its anticanonical divisor is the nef and big
pullback of $-K_Z$. By the Mori dream space consequence of
Birkar--Cascini--Hacon--McKernan, \emph{Existence of minimal models
for varieties of log general type}, Corollary~1.3.2, its
pseudo-effective cone is generated by finitely many effective
$\mathbf Q$-divisor classes, and numerical and $\mathbf Q$-linear
equivalence agree for $\mathbf Q$-Cartier divisors. Every nonzero
effective divisor on $Z'$ has strictly positive $\deg_L$, because
the morphism is small and $L$ is ample on $Z$. Hence the
pseudo-effective class of degree zero just found is
$\mathbf Q$-linearly trivial. Pushing down gives
$K_{\mathcal F}\sim_{\mathbf Q}0$ as a Weil divisor on $Z$.

Take the cyclic cover $\tau:\overline Z\to Z$ trivializing this
torsion Weil class, with a connected component selected if necessary.
It is finite and \'etale in codimension one by the construction
from a trivialization of a reflexive power. Thus $\overline Z$
is again klt with ample rational anticanonical divisor.
On a big smooth open the pullback of the foliation line is
trivial, so its inclusion in the tangent sheaf gives a nonzero
global reflexive vector field on $\overline Z$. It kills the
pulled-back function field of $Y$. The identity component of the
automorphism group of a $\mathbf Q$-Fano is linear algebraic, by
its faithful action in a sufficiently divisible complete
anticanonical embedding. The Zariski closure of the analytic
one-parameter group generated by this vector field fixes the map
to $Y$: the field generators are constant along its flow wherever
defined, and the stabilizer of each rational function is
algebraically closed. Its general orbit has dimension one since
it acts nontrivially in fibers of dimension one. Such an orbit is
unirational, being the image of a connected linear algebraic group
over $\mathbf C$, and its smooth projective model is rational.
The finite map to $Z$ sends a general orbit closure nontrivially
onto the corresponding general curve fiber closure. Thus it
dominates its smooth projective model $C$. This is impossible
for $g(C)\ge1$. The contradiction excludes dimension eight.
\end{proof}

\subsection{The birational quadratic image and its cubic generator}

\begin{proposition}\label{prop:dec2-quadratic-full-classgroup}
Suppose the degree-two system of a Grassmannian-cone replacement has
nine-dimensional image. Let $W'=\operatorname{Proj}S'$, where $S'$
is generated by $S_1,S_2$, and let $H$ be the hyperplane Weil class
of degree one. Then $W'$ is nonnormal or
$\dim_{\mathbf Q}(\operatorname{Cl}(W')\otimes\mathbf Q)\ge2$.
In particular $S'$ is not factorial.
\end{proposition}
\begin{proof}
By Corollary~\ref{cor:dec2-incidence-quadratic-base} and the
subsequent calculation the map to $W'$ is birational and its ring
is $R/(H_3)$ with Hilbert series $(1-t^3)/(1-t)^{11}$. The degree
two Veronese is generated in degree one, so $2H$ is very ample.
The locus where the eleven degree-one coordinates vanish is
contained in the projective four-space of the $y$'s, of codimension
at least five. The class $H$ is thus integral in codimension one.

Suppose $W'$ is normal. We first compute $K_{W'}H^8=-8H^9$ and
$H^9=3$. The Hilbert polynomial of the quadratic embedding, by the
displayed series and eventual equality with the section Hilbert
function, is
\[
       \binom{2m+10}{10}-\binom{2m+7}{10}.
\]
Eight general members of $|2H|$ have smooth integral curve
intersection $C$ in the smooth locus of $W'$. Bertini gives this
since normality makes the singular locus of dimension at most
seven. The Hilbert polynomial records $(2H)^9$ and $\deg K_C$:
take eight successive finite differences to obtain the Hilbert
polynomial of $C$ and apply curve Riemann--Roch. This is the
same polynomial as for the eight quadratic sections on a smooth
cubic ninefold. There adjunction gives $K_C=(-8H+16H)|_C$.
Adjunction along the smooth locus of $W'$ therefore gives the
two asserted numerical identities.

If $\operatorname{Cl}(W')\otimes\mathbf Q=\mathbf Q H$, every
prime $\Delta_i$ below is $\mathbf Q$-linearly equivalent to
$m_iH$, where $m_i>0$ by ampleness. On a common resolution
$p:Y\to Z$, $g:Y\to W'$, write $P=p^*L$, $J=g^*H$ and arrange
$E=P-J\ge0$ using half of the fixed divisor of the quadratic
system. Put $a_j=P^{9-j}J^j/4$, $s=10-\sum_{j=0}^9a_j$.
The absence of a fixed prime downstairs gives $a_0=a_1=1$, because
the squares of the complete degree-one system occur among the
quadratic sections. The numbers $a_j$ decrease and $a_9=3/4$.
Hence $0<s\le2$. Write $g_*E=\sum_i b_i\Delta_i$ and
$q=\sum_i m_i b_i$. Since $a_8=(3/4)(1+q)\le1$, we have
$0\le q\le1/3$.

Push $7E-(K_Y-p^*K_Z)$ onto $W'$. Its support is among the
$\Delta_i$: over the domain of the quadratic morphism any
$p$-exceptional divisor has image of codimension at least two,
and every exceptional divisor over its base has positive fixed
order. Its coefficient $k_i$ at $\Delta_i$ is integral, since
the pushforward of $p^*(K_Z+7L)-K_Y-7g^*H$ is integral at
each generic point, using $K_Z+7L\sim0$ and the codimension-one
integrality of $H$. The pushforward is linearly equivalent to
$-K_{W'}-7H$, so $\sum_i k_i m_i=1$ by the intersection with
$H^8$. Thus some $k_i\ge1$. Its valuation $v_i$ on $Z$ has
$A_Z(v_i)=1+7b_i-k_i\le7b_i$, and semistability gives
$S_L(v_i)\le b_i$.

Let $D_i$ be the strict transform of $\Delta_i$ on $Y$. On
$0\le t\le b_i$, the class $P-tD_i$ dominates the nef segment
from $P$ to $J$, with parameter $t/b_i$. On
$b_i\le t\le b_i+1/m_i$ it dominates $(1-m_i(t-b_i))J$ up to
real linear equivalence, since $g^*\Delta_i\ge D_i$ and
$\Delta_i\sim_{\mathbf Q}m_iH$. Integration of these two
volume bounds as in Theorem~\ref{thm:dec1-general-degree-gap}
gives
\[
          10b_i\ge10S_L(v_i)\ge b_i(10-s)+\frac3{4m_i}.
\]
It follows that $qs\ge m_i b_i s\ge3/4$, whereas $qs\le2/3$.
Finally, factoriality of $S'$ would imply normality of $W'$ and
that each prime divisor of $W'$ has a homogeneous equation,
putting its class in $\mathbf Z H$. This proves the assertion.
\end{proof}

\begin{corollary}\label{cor:dec2-quadratic-vertical-prime}
In the situation of Proposition~\ref{prop:dec2-quadratic-full-classgroup},
the projection from $W'$ to the Grassmannian cone $W$ maps some
prime divisor into a subset of codimension at least two in $W$.
\end{corollary}
\begin{proof}
Use the incidence description $R/(H_3)$ above. On the open part
of $W$ where $A\ne0$ and $AB$ has rank two, the projection from
$\operatorname{Proj}R$ is an affine three-space bundle, in which
$H_3$ is affine linear on the fibers. The complement of this
part has codimension at least two since the rank-one condition on
the general two-by-$\nu$ restriction matrix has codimension
$\nu-1\ge2$. At the generic point of each prime divisor of $W$,
the affine linear polynomial cuts out a nonempty proper locus.
Indeed emptiness would contradict dominance of that prime by a
prime on $Z$ in Proposition~\ref{prop:dec2-grass-all-primes}:
on the domain of the hyperplane morphism its lift to the quadratic
image is regular since the quadratic system contains squares of
linear sections. If the polynomial vanished identically there,
the inverse image of that divisor in the affine three-space bundle
would have the same dimension as $W'$, contradicting its integrality
and dominance over $W$. Thus there is a big open subset $W^\circ$
of $W$, smooth with complement of codimension at least two, over
which $W'$ is an affine two-space bundle. Namely remove the rank
drop of the affine linear part, and solve locally for one variable
where a coefficient is a unit.

If the conclusion failed, the complement of this affine bundle
in $W'$ would have codimension at least two. Indeed the base of
the projection in $W'$ has $A=0$, hence is contained in weighted
projective space of the six remaining coordinates $z,y_1,\ldots,y_5$
of dimension five. Let $\widetilde W'$ be the normalization. Its
complement of the same affine bundle also has codimension at least
two by finiteness. Its class group is therefore generated by $H$:
the base $W^\circ$ has class group $\mathbf Z H$ by factoriality
of the homogeneous ring of $W$, and a divisor on the affine bundle
can be made vertical by subtracting the divisor of a rational
function on its generic affine-plane fiber. Every vertical prime
in this bundle is the reduced inverse image of a base prime.
The class of $H$ agrees with the pullback there.

The normalization is an isomorphism in codimension one since
the affine bundle is smooth. Thus eight general divisors of $|2H|$
cut the same smooth curve on $W'$ and its normalization. In
particular the canonical intersection and volume from the preceding
proof hold on the normalization. The remaining discrepancy and
volume argument there applies verbatim with this normal target,
since the generically birational quadratic morphism on a common
resolution factors through the normalization. Its conclusion
contradicts generation of the rational class group by $H$.
\end{proof}

\begin{proposition}\label{prop:next-cubic-image}
Assume the quadratic image of the replacement has dimension nine,
and write \(A'=R/(H_3)\) for its graded coordinate algebra in
Corollary~\ref{cor:dec2-incidence-quadratic-base}.
Let \(g\) represent the new cubic and choose a nonzero linear form
\(\ell\) such that \(P_4=\ell g\in A'\). Then \(P_4\) is a
zero divisor in \(A'/(\ell)\).
The hyperplane \(\ell=0\) on the quadratic image contains a
prime divisor mapping to codimension at least two in \(W\).
\end{proposition}
\begin{proof}
Suppose \(P_4\) were regular modulo \(\ell\). Introduce
a variable \(u\) of degree three. The algebra
\[
              A'[u]/(\ell u-P_4)
\]
injects into \(A'[\ell^{-1}]\), since its defining principal
ideal is saturated with respect to \(\ell\): reduce a proposed
multiple of the relation modulo \(\ell\), where \(P_4\) is a
nonzerodivisor, and cancel \(\ell\). Its image under \(u\mapsto g\)
is therefore a graded subalgebra of \(S\). Its Hilbert series is
\((1-t^4)(1-t)^{-11}\) since \(H_{A'}=(1-t^3)(1-t)^{-11}\).
Hence it equals \(S\).

The affine cyclic chart at \(u=1\) is consequently
\[
                   U=\operatorname{Spec} A'/(\ell-P_4).
\]
Its origin maps to \(o\in Z\); the residual action of
\(\boldsymbol\mu_3\) has characters one on \(A,z\) and two on
\(y\), hence is free off the origin. This is the local cover
trivializing the Cartier multiple \(3L\) on this open set.
It is normal near the origin. Indeed localizing the normal ring
\(S\) at \(u\) and adjoining an invertible cube root of \(u\)
is an \'etale extension isomorphic to the product of \(U\) with a
one-dimensional torus after rescaling the coordinates by that
root. At the origin the entire group fixes the point, so its
stabilizer preserves the unique local component. Its local degree
over \(Z\) is three. By finite degree and K-semistability,
\[
              \widehat{\operatorname{vol}}(0,U)/9^9
                    \ \ge\ 3\cdot4(7/10)^9 .
\]

We bound this volume from above using the grading on \(A'\).
The seven-dimensional base \(W\) is factorial in homogeneous
coordinates, so \(D=(\ell=0)\) there is integral. The locus on
which \(\operatorname{rank}(AB)<2\) has codimension at least
two in \(W\) by the Schubert rank calculation. Thus over a
general point of \(D\) the incidence space is an affine
three-space bundle in the \(y\)-coordinates on an affine
Pl\"ucker chart. The cubic \(H_3\) is affine linear in these
coordinates. Its restriction to a generic fiber over \(D\)
is a nonconstant affine-linear function. If it were identically
zero, the divisor of \(\ell\) in \(\operatorname{Spec}R\)
would divide \(H_3\), since \(R\) is factorial and its affine
bundle open has complement of codimension at least two. This
contradicts primeness of \((H_3)\). A nonzero constant would
leave the generic fiber over \(D\) empty, contrary to
Proposition~\ref{prop:dec2-grass-all-primes}, because the
quadratic map is regular on the hyperplane map's domain. Thus
\(\operatorname{Proj}(A')\) is smooth over a general point of
\(D\), with its hyperplane \(\ell=0\) reduced there. Let \(E\)
be this component over \(D\) in the weighted blowup of \(U\).
The associated graded ring for the quotient weight filtration
of \(U\) at its origin is \(A'/(\ell)\), since the lowest
piece of \(\ell-P_4\) is the nonzerodivisor \(\ell\) in the
domain \(A'\). In particular the component just found is a
generically smooth exceptional prime of multiplicity one.
Its valuation \(v\) dominates weighted order, with a degree-one
chart parameter having order one. Consequently its volume is
at most the leading Hilbert multiplicity of \(A'/(\ell)\),
namely three by \(H_{A'}=(1-t^3)(1-t)^{-11}\). Indeed the
quotients by its valuation ideals have length bounded above by
the sums of the graded dimensions in weights below the cutoff.

We compute \(A_U(v)=7\). On the affine incidence space where
\(A=p\wedge q\ne0\), the two equations in the auxiliary
coordinates \(p,q,y,z\) are a regular sequence, and passing
from the basis \(p,q\) to \(A\) is an \(\operatorname{SL}_2\)-torsor.
The residue of the ambient top differential along the two
equations, contracted with a determinant-one basis of the
three infinitesimal group directions, gives a trivializing
dualizing form on this incidence open. Using doubled weights
\(1,1,4,2\) on the groups of variables \(p,q,y,z\), its
weight is \(10+20+2-2\cdot5=22\), since the group directions
have weight zero. Thus its weight for the actual grading of
\(R\) is eleven. The residue along \(H_3\) gives a dualizing
form of weight eight on \(A'\) where \(A\ne0\). In the chart
\(U\) its further residue along \(\ell-P_4\) has no divisor:
these adjunctions trivialize the dualizing sheaf on the local
complete-intersection locus \(A\ne0\), and the remaining locus
has dimension at most six in the ninefold \(U\).
Near a general point of \(E\), choose a nonzero degree-one
Pl\"ucker coordinate as radial variable \(s\) on the weighted
blowup. In the rescaled coordinates the ambient form of weight
eight has factor \(s^7\,ds\), with unit nonradial coefficient
at that generic point. The equation has factor \(s\), leaving
the rescaled equation \(\ell_{\mathrm{chart}}-s^3
P_{4,\mathrm{chart}}\). Its differential in a direction along
the exceptional divisor is nonzero by smoothness of the
hyperplane cut there. The residue on the strict transform
therefore has order six, giving the claimed log discrepancy.
This yields
\[
            \widehat{\operatorname{vol}}(0,U)/9^9
               \le 3(7/9)^9 < 12(7/10)^9 ,
\]
the desired contradiction.

For the geometric assertion, \(R\) is Cohen--Macaulay. In fact
take \(\operatorname{SL}_2\)-invariants of the Koszul resolution of
the two auxiliary equations in \(p,q,y,z\) used in
Proposition~\ref{prop:dec2-spinor-incidence-quartic}. Its outer
terms are graded shifts of the polynomial extension of the
Pl\"ucker ring by \(y,z\). Its middle term is the same polynomial
extension of the standard-covariant module of the two-by-five
matrix. This module is maximal Cohen--Macaulay over the
Pl\"ucker ring. Here are details. The Cauchy decomposition of
\(\mathbf C[p,q]\) identifies the multiplicity space of the
standard \(\operatorname{SL}_2\)-representation with
\(\bigoplus_{j\ge0}\mathbf S_{(j+1,j)}(\mathbf C^5)^*\),
the full section module of the tautological dual
\(\mathcal S^*\) on \(G\) by Borel--Weil. The module identification
respects the Pl\"ucker action, by restricting covariants on the
rank-two matrix locus to frame bases of each plane. Complements
of this locus have codimension at least two, so regular
covariants extend. The intermediate cohomology of
\(\mathcal S^*(j)\) vanishes for every integer \(j\):
the shifted weights for Bott's theorem are
\((j+5,j+3,2,1,0)\). For \(j\ge0\) they are dominant, for
\(-5\le j\le-1\) they are singular, and for \(j\le-6\)
there are six inversions. The same conclusion for the Pl\"ucker
ring itself follows with the first two coordinates \(j+4,j+3\).
By the section-module local cohomology criterion both modules
are maximal Cohen--Macaulay. Thus the three terms of the
invariant Koszul resolution have depth thirteen after adjoining
the six \(y,z\)-variables. The depth lemma gives depth at least
eleven for \(R\), its dimension. This graded Cohen--Macaulay
property at the homogeneous maximal ideal implies the
Cohen--Macaulay property of the ring.

The quotient \(A'/(\ell)=R/(H_3,\ell)\) is therefore a complete
intersection in a Cohen--Macaulay ring and has no embedded
associated primes. If the component dominating \(D\) described
above were the sole component, generic reducedness would make
the quotient a domain. The proved zero-divisor property would
then force \(P_4\in(\ell)\) in \(A'\), which gives \(g\in A'\),
a contradiction. Thus another component exists. It does not
lie in \(A=0\), of dimension at most six on the affine cone.
Since the general fiber over \(D\) is an affine plane,
this other prime maps into a proper subvariety of \(D\).
\end{proof}

When \(\nu=\operatorname{rank}B\ge4\), every linear denominator
\(\ell\) as above belongs to the ten-dimensional Pl\"ucker system.
To prove this, work away from the vertex of \(W\), dividing by a
nonzero Pl\"ucker coordinate. Group the equations of the quadratic
image as an affine-linear system in the five \(y\)-variables. Its
coefficients depend only on the plane in \(G\): the first two rows
come from \(AB\), and the last row has coefficients linear in \(A\)
by \(H_3(0,z,y)=0\). Write \(k\) for its rank, generically three.
A nonempty fiber of this affine-linear system has dimension \(5-k\).
At rank three the fibers have dimension two. The locus of rank at
most two in the coefficient matrix has codimension at least two
in \(G\). Indeed if a divisor occurred, then \(AB\) would still
have rank two at its generic point, since its own rank-drop locus
has codimension \(\nu-1\ge3\) if nonempty. On the cone of that
divisor the cubic would have no linear part after solving the
two incidence equations. This was excluded in the preceding proof
for an arbitrary prime on \(W\), using its dominance by a prime
of \(Z\), factoriality of \(R\), and the fact that the cubic
generates a prime ideal of degree three. In this argument a common
divisor in \(R\) of degree three would identify \(H_3\) with its
pullback from \(W\), impossible by dominance over \(W\).

The prime furnished by the proposition has dimension eight. If
the matrix has generic rank two at its center, this center must
have dimension at least five. The rank-two locus on \(W\) is
pulled back from \(G\) and has dimension at most five. Thus the
center contains a general ruling of the cone over each general
point of its projection to \(G\). At rank one the same argument
applies: the locus is contained in the rank-one locus of \(AB\)
and hence has codimension at least three in \(W\). The fiber
dimension four requires the center to have dimension at least
four, which is the upper bound. Rank zero cannot occur because
\(\nu\ge4\) and the two-plane would have to annihilate
\(\operatorname{im}B\). Since the center belongs to
\(\{\ell=0\}\) and contains an open subset of a ruling, the
coefficient of the vertex coordinate \(z\) in \(\ell\) is zero.

In the nine-dimensional quadratic-image case, the resulting image
of a general threefold fiber \(V\) over \(G\) is the rank-three
quadric threefold, with its usual polarization. In the coordinates
of the proof of Proposition~\ref{prop:dec2-grass-projected-degree},
the substituted cubic for \(A=sA_0\) has a factor \(s\).
Its quotient is a nonzero constant linear combination of the three
remaining weight-two variables plus a binary quadratic in \(s,z\):
the coefficient of \(y\) in \(H_3\) is linear in \(A\) alone,
and its restriction to the kernel is nonzero by dominance and
the affine-plane calculation. Eliminating one weight-two variable
leaves \(\mathbf P(1,1,2,2)\) with its second Veronese system,
whose presentation has five coordinates and the single rank-three
equation \((sz)^2=s^2z^2\). Over a general point of \(W\), fixing
\([s:z]\) gives \(\mathbf P(1,2,2)\) in its second Veronese,
a plane. Thus the two general fiber closures satisfy
\[
             L^3\cdot V\ge1/4,\qquad L^2\cdot S_{\rm fib}\ge1/4
\]
by the moving-intersection bound for the quadrics on each fiber.
In particular \(a_7=5L^2\cdot S_{\rm fib}/4\ge5/16\). The
complete linear base has codimension at least four and \(a_3=1\)
by Proposition~\ref{prop:dec2-ninefold-tangent-pencil}. Since
\((7/10)^4<5/16\), log-concavity gives
\[
    \theta\ge(7/10)\big(4+\textstyle\sum_{j=1}^4(7/10)^j\big)>4 .
\]
Thus the complete base has codimension at least five, and its
codimension-five components have smooth generic transverse germs
by Cartierness of \(L\) and the four-cut argument in
Proposition~\ref{prop:dec2-grass-projected-degree}. Consequently
the codimension-at-most-five component of the singular Pl\"ucker
base in Proposition~\ref{prop:dec2-grass-singular-base} has
nonzero vertex coordinate at its general point.

In fact that singular intersection has codimension exactly five.
The threefold degree gives \(b_6=5L^3\cdot V/4\ge5/16\).
Interpolation from \(b_2=1\) and \((7/10)^4<5/16\) first gives
\(\theta_P>(7/10)(3+\sum_{j=1}^4(7/10)^j)>3\). Hence the
projected base has codimension at least four, and at a general
point of any codimension-four component three general members
cut a regular curve. As \(L\) is Cartier there, the transverse
fourfold is regular by the resulting embedding dimension bound.
Thus its singular intersection has codimension at least five.
The opposite bound is Proposition~\ref{prop:dec2-grass-singular-base}.

At the general point of a codimension-five component of this
intersection, the transverse germ is a hypersurface double point.
Indeed the stricter interpolation from \(b_3=1\) gives
\(\theta_P>(7/10)(4+\sum_{j=1}^3(27/40)^j)>19/5\).
Here \((27/40)^3<5/16\). The transverse germ is Gorenstein
because \(L\) is Cartier, and is a hypersurface by
\cite[Lemma~4.5]{kePanWuCubics}. Since it is singular it has six
ambient coordinates, and ordinary blowup and inversion of
adjunction bound the threshold of its maximal ideal by six
minus its multiplicity. The Pl\"ucker ideal is contained in the
maximal ideal, giving multiplicity two. At least four independent
linear parts occur in the Pl\"ucker functions: if there are
\(h\) such parts, choose those functions as the first \(h\)
smooth ambient parameters. Weights two on these \(h\) parameters
and one on the others give the ideal weight at least two, and
the hypersurface equation has weight at least two. The resulting
threshold is at most \((4+h)/2\), which excludes \(h\le3\).

Write \(A_1\) for the matrix of these linear parts in the six
ambient coordinates. Each of its five Pfaffians is a scalar
multiple of the initial quadratic of the transverse equation,
since the five full Pfaffians vanish in the hypersurface ring.
If a multiple is nonzero, the exterior-square vector of \(A_1\)
is proportional to a fixed nonzero covector of the five-space.
The Pfaffian identity then shows that every \(A_1\) takes values
in \(\bigwedge^2 E\) for the four-space \(E\) annihilated by
that covector. The initial quadratic is the pullback of the
nondegenerate Pfaffian on this six-space, of rank at least
\(2h-6\) on an \(h\)-dimensional linear image. In particular
six independent linear parts give an ordinary double point.
In this case there is also an exact analytic description. Lift the
skew matrix to the regular six-dimensional ambient germ and write
its Pfaffian vector as \(f b\), where \(f\) is a local defining
equation of the transverse hypersurface. A nonzero quadratic
Pfaffian initial means one component of \(b\) is a unit.
The Pfaffian syzygy \(A f b=0\) in the ambient domain gives
\(A b=0\). A change of basis depending analytically on \(b\)
therefore puts \(A\) into a four-by-four skew block and a zero
last row and column. Its six block entries define a map to
\(\bigwedge^2\mathbf C^4\) of differential rank \(h\), and
\(f\) is a unit times their nondegenerate Pfaffian. Each of
these six entries lies in the ideal of the original ten entries
in the regular ambient germ. Assign weight two to \(h\) of the
original entries with independent linear parts, completing
those to regular ambient parameters of weight one. Then every
original entry has weight at least two, and the Pfaffian in the
new entries has weight at least four. Thus the threshold is at
most \((6+h-4)/2\); as it exceeds \(19/5\) one must have
\(h=6\). This branch is an ordinary double point, with the
Pl\"ucker ideal equal to its maximal ideal by Nakayama.

It remains to describe the other local possibility. If all
quadratic Pfaffian initials vanish, the \(h\)-dimensional linear
image consists entirely of decomposable bivectors. A linear
space of dimension at least four of such bivectors in five
variables is \(e_0\wedge\mathbf C^4\) in a suitable basis.
Indeed intersecting pairs of distinct support planes meet in a
line, and any further plane either contains that common line or
lies in their three-dimensional span. A space in the latter
three-space has dimension at most three. Thus \(h=4\). Write
the four components involving \(e_0\) as \(x_i\), and denote
the four-by-four block on the complement by \(B\). The \(x_i\)
lift to part of a system of regular ambient parameters; the
initial quadric is nonzero on the two remaining parameters,
since otherwise the same weight bound gives threshold at most
\(7/2\). Consequently the \(x_i\) form a regular sequence in
the transverse hypersurface. The Pfaffian relation
\(x\wedge B=0\) and exactness of the Koszul complex in exterior
degree two give \(B=x\wedge b\) for a regular vector \(b\).
Locally the Pl\"ucker matrix is therefore decomposable as
\((e_0-b)\wedge x\). On its domain the corresponding planes
all contain the regular nowhere-vanishing vector \(e_0-b\).
The pulled-back tautological rank-two bundle has a nowhere-zero
section, so its second Chern class vanishes on this local domain.
As \(L\) is trivialized by \(z\), the same holds for the second
Chern class of the tautological quotient.

At least one codimension-five component has the ordinary
double-point type. Otherwise cut by four general high relative
very ample Cartier divisors as in
Proposition~\ref{prop:dec2-grass-singular-base}. Write \(V_5\)
for the fivefold slice and \(K\) for its intersection with the
projected base. The slice meets the singular base in finitely
many general points \(p_i\), all of the four-linear-part type.
On \(V_5\setminus\{p_i\}\) the two integral Schubert classes on
the complement of \(K\) extend uniquely: the slice is smooth
near \(K\setminus\{p_i\}\) and this closed subset of the punctured
slice has codimension at least four. Excision and duality with
Borel--Moore homology give vanishing of support cohomology in
degrees four and five. In a small punctured neighborhood of each
\(p_i\), the same argument gives injectivity on \(H^4\) on removing
\(K\). Our local factorization applies to the transverse slice
at that point and makes both Schubert classes zero off \(K\)
there. Thus the extended classes vanish on each punctured
neighborhood. The Mayer--Vietoris sequence glues them with zero
classes across the \(p_i\) to integral classes on \(V_5\).
Their sum can be arranged to be \(c_1(L)^2\): take any extension
of the first, and take \(c_1(L)^2\) minus this extension for
the second. The root is a line bundle on this slice since its
non-Cartier locus has codimension at least six. The same positive
degrees and integral specialization as in
Proposition~\ref{prop:dec2-grass-singular-base} would then express
the primitive hyperplane square on the smooth slice as a sum of
two strictly positive integral multiples of itself. This is
impossible, proving the assertion.

Let \(C_j\) be the closures of all the four-dimensional components
which have the ordinary double-point type. Then
\[
                0<\sum_j L^4\cdot C_j
                        \le 2-\frac{\sqrt5}{2}<1 .
\]
Indeed over a general point of \(C_j\) the projected base
blowup is the blowup of the maximal ideal of the transverse node.
The exceptional divisor there has coefficient one in \(P-J\)
and its geometric fiber is the smooth quadric fourfold polarized
by \(J\) as its usual hyperplane class. The top \(J\)-degree
of this fiber is two. Intersecting \(P^4J^4(P-J)\) and using
the projection formula gives \(4(b_4-b_5)\ge2\sum_j L^4 C_j\).
Every further exceptional contribution is nonnegative by nefness.
Since \(b_5\ge\sqrt{b_4b_6}\) and
\(5/16\le b_6\le b_4\le1\), the right bound follows by
maximizing \(b_4-\sqrt{5b_4/16}\) on this interval.
In particular every \(C_j\) meets the non-Cartier locus of
\(L\): if \(L\) were Cartier in a neighborhood of \(C_j\), the
positive degree \(L^4C_j\) would be integral. Any such point
lies on \(z=0\), since a nonzero value of a root section in
the geometric sense trivializes the root there.

If \(\nu\ge4\), the quadratic image of each such node fourfold
has dimension at most three. Work on \(z=1\), where this image
lies in the affine five-space with coordinates \(y\). Put
\(m(y)=By+c(0,1)\) and let \(\omega(y)\) be the derivative of
\(H_3\) in the ten entries of \(A\) at \((0,1,y)\), considered
as an alternating covariant two-form. This is affine linear in
\(y\). At a general node point the possible transverse linear
parts of \(A\) comprise \(\bigwedge^2 E\) with \(E=\ker b\)
for a nonzero covector \(b\). Differentiating \(A(By+c)=0\)
gives \(m(y)\in\mathbf C b\). Differentiating the cubic gives
\(\omega(y)|_{\bigwedge^2E}=0\), using \(H_3(0,1,y)=0\)
identically to remove derivatives in \(y\). Hence the quadratic
vector polynomial \(m(y)\wedge\omega(y)\) vanishes there.

Suppose its zero locus contained a hypersurface in the five-space.
Homogenize the two affine linear maps with a variable \(t\);
their homogeneous versions are denoted \(m,\omega\). If their
wedge is nonzero, its components have a common homogeneous
irreducible factor \(f\) of degree one or two. The entries of
\(m\) include at least four independent linear forms. In the
polynomial ring in \(y,t\) the Koszul complex for exterior
multiplication by \(m\) is exact in exterior degrees smaller
than four. This follows by choosing the independent linear forms
as coordinates, reducing the redundant entries to zero by a
constant change of basis, and using the Koszul complex of a
regular sequence. A degree-two factor would make
\((m\wedge\omega)/f\) a constant three-form killed by \(m\);
exactness and homogeneity force it to vanish. A linear factor
instead gives
\[
               \omega=f\,\eta +m\wedge\beta
\]
for a constant two-form \(\eta\) and a constant covector \(\beta\).
Indeed first apply exactness in exterior degree three to the
linear three-form \((m\wedge\omega)/f\), then apply it in
exterior degree two to \(\omega-f\eta\). If \(m\wedge\omega\)
is identically zero we obtain the same expression with
\(\eta=0\).

The \(y\)-linear part of \(\omega(y)\) records the coefficients
of the terms of \(H_3\) linear in both \(A\) and \(y\), and there
are no \(zy\)-terms in this cubic. For any general \(A\), pairing
\(A\) with \((By)\wedge\beta\) vanishes in the kernel directions
\(AB y=0\) of the affine incidence fiber. If \(\eta=0\), the
cubic consequently has zero linear part in that general fiber,
contrary to its affine-plane description. If \(\eta\ne0\), that
linear part vanishes on the fiber over the general point of the
Pl\"ucker hyperplane \(\langle\eta,A\rangle=0\), where \(AB\)
still has rank two. The prime-dominance and irreducibility
argument in the proof above rules this out. Thus the zero locus
of \(m\wedge\omega\) has codimension at least two, as asserted.

In this case the section ring also requires a generator
of degree at least four. Indeed set $A'=R/(H_3)$ and take the
cubic $g$ and its denominator $\ell$ from Proposition~\ref{prop:next-cubic-image}.
Put $P_4=\ell g$. The Hilbert series of
$A'[u]/(\ell u-P_4)$, with $u$ of degree three, equals the Hilbert
series of $S$. Choose a homogeneous $a\notin\ell A'$ with $aP_4=\ell b$
in $A'$, using that proposition. The relation $au-b$ maps to zero
in $A'[g]$ and is nonzero modulo $\ell u-P_4$, since its leading
coefficient is outside $\ell A'$ and $A'$ is a domain. Thus
$A'[g]$ has smaller dimension than $S$ in some degree.
All elements through degree three belong to $A'[g]$, proving the claim.

\subsection{Adjunction on the general quadratic plane and the ninefold theorem}

\begin{proposition}\label{prop:next-q9-no-nine-quadratic}
On a K-polystable quartic ninefold replacement whose limiting
hyperplane image is the Grassmannian cone, the complete quadratic
system cannot have nine-dimensional image.
\end{proposition}
\begin{proof}
Suppose it does. On a smooth common resolution $p:Y\to Z$, resolve
the maps to the Grassmannian $G$, the cone $W$ (also its bundle
resolution over $G$), and the quadratic image. Write $P=p^*L$,
let $J$ be the moving class for $|L|$ and $M$ the moving class for
$|2L|$. Thus $P,J,M$ are nef, $P-J$ and $2P-M$ are effective up
to rational linear equivalence, using products in divisible Cartier
multiples. All generic fibers below may be taken over general
complex points by generic smoothness. Let $V$ be a smooth resolved
threefold fiber over $G$, and $F$ a smooth resolved surface fiber
inside it over a general point $w$ of the ruling of $W$. On $V$ the
class $J$ pulls back $\mathcal O_{\mathbf P^1}(1)$ from this ruling,
so that $F\sim J|_V$.

The system $M|_F$ gives a birational morphism to a linear plane.
In fact on the affine chart of the general $w$ trivialized by a
nonzero degree-one coordinate $x$, the quadratic map has coordinates
the constant degree-one products and the five affine coordinates
$y_i/x^2$. The incidence equations cut these last coordinates by
two independent affine linear equations. The additional cubic
equation, on this chart, is one further independent affine linear
equation, as shown following
Corollary~\ref{cor:dec2-incidence-quadratic-base}. Thus the closure
in the ordinary quadratic-coordinate projective space is a plane
with a line at infinity. The full quadratic map is birational, and
the function fields of these planes over $W$ agree with the field of
$Z$ over $W$. Hence the morphism on a general $F$ is birational.
Take as $C\subset F$ the inverse image of a general line of that
plane. Such a line avoids the finitely many exceptional images of
the birational surface morphism, so $C\simeq\mathbf P^1$,
$C^2=1$ on $F$, and $K_F\cdot C=-3$ by adjunction. As $F$ is
a smooth fiber over the smooth locus of $W$, its normal bundle is
trivial and $K_Y\cdot C=-3$ as well.

Put $a=P\cdot C$. The mixed intersections on $V$ give
\[
       a=P M J\cdot V\le 2P^3\cdot V<2/3 .
\]
Indeed the replacements of $M,J$ by $2P,P$ are valid against nef
products on a general fiber, and the strict bound $L^3p_*V<1/3$
is the first estimate in the proof of
Corollary~\ref{cor:next-q9-conic-or-cubic}; that estimate for the
projection to $G$ is independent of the quadratic image dimension.

For a general linear section $x$ nonzero at $w$, pull back its
effective rational Cartier section divisor to $C$. It has degree
$a$ and every positive coefficient is at least $1/6$, since the
local Cartier index at every point of $Z$ is at most six. At the
point $q$ where the image line meets infinity some ratio $y_i/x^2$
has a simple pole: on the image plane these are the ordinary affine
coordinates relative to the constant coordinate, and the line is
transverse to infinity. Since the section divisor of $y_i$ pulls
back effectively on $C$, the coefficient of $p^*\operatorname{div}(x)$
at $q$ is at least $1/2$. The degree bound $a<2/3$ precludes another
positive coefficient. The restriction is therefore precisely $a q$.

We record which exceptional divisors can enter adjunction. Outside
the geometric degree-one base both the cone map and the quadratic
map are regular, by taking Cartier powers. If a $p$-exceptional
divisor has generic center there, its image in the simultaneous
graph of these maps has dimension at most seven, the dimension
bound on its center in $Z$. If that image dominates $W$, its fiber
over a general $w$ consists of finitely many points in the quadratic
plane; a general line avoids them. The other such images miss the
general $w$ entirely. Thus only exceptional divisors centered in
the geometric degree-one base can meet $C$, and each has strictly
positive order on the general linear section $x$. The strict
transform of its divisor does not meet $F$, since it maps to a
hyperplane section missing $w$. Consequently exactly one exceptional
divisor $E$ meets $C$, and it meets $C$ transversely at $q$.
For these genericity assertions take the exceptional set to have
simple normal crossings. Pairwise intersections have dimension at
most seven and give finitely many points on the general $F$ when
dominant. Each horizontal exceptional prime cuts the general $F$
in smooth curves, by generic smoothness. Their images which are
points are avoided by the general line; the rest meet the line
transversely on the isomorphism locus of the surface morphism.
In particular $\operatorname{ord}_E(|L|)=a$.

The discrepancy identity on $Y$ now gives
\[
        -3=K_Y\cdot C=-7a+A_Z(E)-1,
        \qquad A_Z(E)/a=7-2/a<4 .
\]
But the mixed-volume lower bound for the complete linear system is
strictly greater than four in this case. Indeed $M^2F=1$, so
$P^2F\ge1/4$. In the notation of
Proposition~\ref{prop:dec2-ninefold-tangent-pencil}, we have
$a_7=5P^2F/4\ge5/16$ and $a_0=a_1=a_2=a_3=1$ from the
codimension-four base bound. Log-concavity and $(7/10)^4<5/16$
give
\[
  \operatorname{lct}(Z;|L|)\ge\frac7{10}\sum_{i=0}^7a_i
   >\frac7{10}\left(4+\sum_{i=1}^4(7/10)^i\right)>4.
\]
This contradicts the inequality supplied by $E$.
\end{proof}

\begin{theorem}\label{thm:next-all-smooth-quartic-ninefolds}
Every smooth complex quartic ninefold $X\subset\mathbf P^{10}$
is uniformly K-stable and admits a K\"ahler--Einstein metric.
\end{theorem}
\begin{proof}
Choose a pointed curve through the parameter of $X$ in the space of
smooth quartic forms such that its punctured family is uniformly
K-stable. Such a curve exists by irreducibility, openness of
uniform K-stability, and K-stability of the Fermat quartic
\cite[Corollary~1.4]{ar2004.09413}. K-moduli properness, after
finite base change, gives a K-polystable special fiber $Z$ of
this family. The construction in the last paragraph of the proof of
Theorem~\ref{thm:dec1-quartic-high-dimensions} gives its limiting
ample Weil class $L$, with $-K_Z\sim7L$, $L^9=4$ and the
quartic Hilbert function. The smooth hypersurface filling over
the original curve has the GIT point of $X$.

If this replacement differs from $X$,
Corollary~\ref{cor:dec1-root-grassmannian-only} identifies the
hyperplane image $W$ with the projective cone over
$\operatorname{Gr}(2,5)$. The function field of the quadratic
image contains $\mathbf C(W)$, by taking ratios of products of
linear sections with a fixed nonzero linear section.
There is a quadratic section outside the degree-two part of the
ring generated by the linear sections, by the five independent
Pl\"ucker equations and the quartic Hilbert function. Divide such
a section by the square of a linear coordinate. The resulting
function lies outside $\mathbf C(W)$ by
Proposition~\ref{prop:dec2-grass-all-primes}, and is transcendental
over it by Proposition~\ref{prop:dec1-stein-hyperplane-fibers}.
Thus the quadratic image has dimension at least eight.
Theorem~\ref{thm:next-q9-no-eight-quadratic} excludes dimension
eight. Proposition~\ref{prop:next-q9-no-nine-quadratic} excludes
dimension nine. This contradicts the existence of a distinct
replacement.

Consequently $X$ itself is K-polystable. The automorphism group
of a smooth quartic ninefold is finite, giving K-stability.
The equivalence of K-stability and uniform K-stability for Fano
manifolds and the Yau--Tian--Donaldson theorem give the assertions,
as in Section~\ref{sec:established-cases}.
\end{proof}

Together with the quartic results through dimension eight and
Theorem~\ref{thm:dec1-quartic-high-dimensions}, this gives the
quartic row of Theorem~\ref{thm:range} for every $n\ge3$.

\section{A numerical obstruction to full-dimensional hyperplane images}

\begin{proposition}\label{prop:dec1-general-conductor-test}
Let $(Z,L)$ be a K-semistable replacement of smooth Fano hypersurfaces
of degree $d$, with $-K_Z\sim rL$, $L^n=d$ and $r=n+2-d$.
Suppose its complete hyperplane system has full-dimensional image $W$
of degree $d/2<e<d$. With $q$ and $s$ computed on the normalization
as in Theorem~\ref{thm:dec1-normal-degree-gap}, one has
\[
                      qs\ \ge\ \frac{2-1/e}{d}.
\]
\end{proposition}
\begin{proof}
Put $h=d-e$ and $N=n+1$. The map onto $W$ is birational. Retain
the primes $\Delta_i$ on the normalization, their hyperplane
degrees $D_i$, the fixed orders $b_i$, the conductor coefficients
$c_i$ and the integral coefficients $k_i$ of the preceding proof.
The normalization formula and the discrepancy bound give
\[
  \sum_iD_i(k_i-c_i)=eh,\qquad
  k_i\le 1+(r/N)s b_i,\qquad \sum_iD_i b_i=e q .
\]
Suppose $qs<(2-1/e)/d$. If $k_i\le0$ whenever $c_i=0$, then
$k_i-c_i\le(r/N)s b_i$ at every index, since each positive
conductor coefficient is at least one. Summing would give
$h\le(r/N)qs<(2-1/e)/d<1$, contrary to $h\ge1$.
Thus some prime satisfies $k_i\ge1$ and $c_i=0$. Its image is
generically normal on $W$, and
$A_Z(\operatorname{ord}_{\Delta_i})=1+rb_i-k_i\le rb_i$.
The moving system of degree-$D_i$ hypersurfaces through its image
and the corresponding two intervals of the volume integral in
Theorem~\ref{thm:dec1-normal-degree-gap} apply to this prime:
\[
 b_i N\ge N S_L(\operatorname{ord}_{\Delta_i})
       \ge b_i(N-s)+(e/d)(2-1/e)/D_i .
\]
Using $e q\ge D_i b_i$ contradicts the assumed strict bound on $qs$.
\end{proof}

\begin{theorem}\label{thm:dec1-unconditional-degree-minus-three}
For every $d\ge14$, the complete hyperplane image of a K-semistable
replacement as in Proposition~\ref{prop:dec1-general-conductor-test}
cannot be an $n$-dimensional hypersurface of degree $d-3$.
\end{theorem}
\begin{proof}
Assume $e=d-3$ and put $\lambda=3/d$, $\tau=e/d$ and
$T=(2-1/e)/d$. Write $u=1-a_{n-1}$ for the penultimate normalized
mixed intersection, $m=n+1-k$ for the integer determined by the
base codimension. As in Theorem~\ref{thm:dec1-normal-degree-gap},
$0\le u\le\lambda$, $q=(\lambda-u)/\tau$, $m<d-1+s$, and
interpolation between $a_{k-1}=1$ and $a_{n-1}=1-u$ gives
\[
 s\le\frac{2\lambda+(d-1-2\lambda)u}{2-3u},
 \qquad
 s\le\lambda+mu\left(\frac12+\frac{u}{12(1-u)^2}\right).
\]
For $m=1$ these inequalities hold with $s=\lambda$ and $u=0$.
By Proposition~\ref{prop:dec1-general-conductor-test}, $qs\ge T$.
Substitute the first upper estimate for $s$ and clear the positive
denominators. The resulting inequality is
\[
 d(d^2-d-6)u^2-(3d^2-3d-39)u+4d-32\ \le\ 0 .
\]
Its discriminant is
$(3d^2-3d-39)^2-8d(d^2-d-6)(2d-16)$.
At $d=16$ its value and its first four successive forward
differences are $-15471,-31968,-11952,-2184,-168$.
Its expansion in the binomial basis at $16$ is strictly negative
on every integer $d\ge16$. This excludes those degrees.

For $d=14,15$ set $c=1/2+\lambda/[12(1-\lambda)^2]$.
Using $m<d-1+s$ in the second interpolation bound yields
\[
    s\le\frac{\lambda+(d-1)c\lambda}{1-c\lambda}<2 .
\]
Indeed the right side is $2859/1502$ for $d=14$ and $1606/859$
for $d=15$. In particular $m\le d$. For $d=15$ the resulting
quadratic bound gives
\[
 qs\le\frac{\lambda^2(1+dc)^2}{4dc\tau}
     =\frac{323761}{2585600}<\frac{23}{180}=T.
\]
For $d=14$, first suppose $u\le1/8$. The Taylor factor is at most
$1/2+(1/8)/[12(7/8)^2]=151/294$. The same maximum with this
factor in place of $c$ is $11094/81389<3/22=T$.
For $1/8\le u\le\lambda$, the quadratic
$(\lambda-u)(\lambda+14cu)/\tau$, with $c=64/121$, decreases
throughout the interval. Its value at $1/8$ is
$9655/74536<3/22$. This gives the contradiction in both degrees.
\end{proof}

\begin{corollary}\label{cor:dec1-degree-thirteen-image-ten}
For $12\le n\le236$, a K-semistable replacement of smooth
degree-thirteen hypersurfaces in $\mathbf P^{n+1}$ cannot have a
full-dimensional hyperplane image of degree ten.
\end{corollary}
\begin{proof}
Here $\lambda=3/13$, $\tau=10/13$, and the necessary lower bound
for $qs$ is $T=19/130$. With $c=213/400$ the same Taylor estimate
gives $s\le8868/4561<2$, so $m\le13$. If $m\le12$ and
$u\le1/8$, use $151/294$ for the Taylor factor: the resulting
quadratic maximum is $85293/591920<T$. For $u\ge1/8$ the
quadratic $(\lambda-u)(\lambda+12cu)/\tau$ is decreasing and
its value at $1/8$ is $117777/832000<T$.

It remains to treat $m=13$. The base-codimension inequality gives
$s\ge N/r\ge79/75$ because $N=n+1\le237$ and $r=N-12$.
Set $z=(1-u)^{1/12}$, $S(z)=\lambda+12-\sum_{j=1}^{12}z^j$.
Log-concavity gives $s\le S(z)$. Put $z_0=39561/40000$.
As $S$ decreases and $S(z_0)<1053298/10^6<79/75$, we must have
$z<z_0$. Define $P(z)=(z^{12}-\tau)S(z)$. The ratio
\[
 \frac{P'(z)}{12z^{11}}
    =S(z)-\frac{z-\tau/z^{11}}{12}\sum_{j=1}^{12}jz^{j-1}
\]
decreases for $z\ge\tau^{1/12}$, where both factors of the
subtracted term are nonnegative and increasing. Substitution
at $z_0$ gives $P'(z_0)>7/2$ and
$P(z_0)/\tau<146152/10^6<19/130$. Hence $qs\le P(z)/\tau<T$,
which contradicts Proposition~\ref{prop:dec1-general-conductor-test}.
\end{proof}

\chapter{Point valuations and filtrations of the hyperplane section ring}\label{ch:nov30-filtrations}

\section{A divisor in the exceptional projective tangent space}
\label{sec:angular-divisor-bound}
Let \(p\) be a point of a smooth \(X=X_d^n\), \(3\le d\le n+1\).
Write \(E=\mathbf P(T_pX)\) for the exceptional divisor of the blowup
\(\pi:Y\to X\), \(N=n+1\), \(r=n+2-d\), and
\(\mu=S_H(\operatorname{ord}_p)\). The following test applies to any
prime divisor over \(X\) whose center on \(Y\) is a prime divisor on
\(E\). In particular the divisor over \(X\) may be obtained after
further blowups.

\begin{theorem}\label{thm:angular-divisor-bound}
Let \(v\) be such a divisorial valuation. Suppose its center on \(Y\)
is an integral hypersurface \(D\subset E\) of degree \(e\). Set
\(a=v(E)\) and \(\gamma=A_X(v)-na\). Then \(a>0,\ \gamma\ge0\), and
\[
 S_H(v)\le a\mu+
             \frac{\gamma}{e}\left(\mu-\frac{n-1}{N}\right).
 \tag{AD1}
\]
Consequently \(A_X(v)>rS_H(v)\) whenever \(Ne\ge rd\).
In particular this excludes every \(e\ge d\) in every Fano dimension.
For quartic ninefolds it excludes all \(e\ge3\).
\end{theorem}

\begin{proof}
Trivialize \(H\) by a section nonzero at \(p\). For a nonzero
section \(s\) of \(mH\) write \(j=\operatorname{ord}_p(s)\) and let
\(F\) be its degree-\(j\) initial form on \(T_pX\). Write
\(\ell=\operatorname{ord}_D F\) and let \(B\) be the divisor of the
strict residual section after dividing \(\pi^*s\) by \(E^j\).
At the generic point of \(D\), \(B|_E\) has multiplicity \(\ell\).
If \(\ell>0\), the pair \((E,B|_E/\ell)\) is log canonical there,
since the local ring is a discrete valuation ring. Inversion of
adjunction for the smooth Cartier divisor \(E\) gives
\(v(B)\le\ell(A_Y(v)-v(E))\). For \(\ell=0\) the residual section
is a unit there. The blowup formula gives \(A_Y(v)-v(E)=\gamma\),
nonnegative by log canonicity of \((Y,E)\). Thus in both cases
\[
                          v(s)\le aj+\gamma\ell .
\]

Let \(h_m=h^0(X,mH)\). For given \(m,j\) let \(W_{m,j}\) be the
space of initial forms of sections with order at least \(j\), modulo
order greater than \(j\). It is a vector space of degree-\(j\)
polynomials. Write \(J_m=\sum_j j\dim W_{m,j}\). Let \(L_m\) be
the sum, over \(j\), of the vanishing orders of a basis of
\(W_{m,j}\) adapted to order along \(D\). A basis of \(H^0(X,mH)\)
simultaneously adapted to \(v\) and to order at \(p\) induces bases
of these spaces. Their total \(D\)-order is at most \(L_m\), by the
dimension formula for the jumping numbers of a filtration.
The sum of all \(v\)-jumping numbers is therefore at most
\(aJ_m+\gamma L_m\).

Choose coordinates \(x_1,\dots,x_n\) on \(T_pX\) so the last
coordinate point lies outside \(D\). Use the least exponent in
lexicographic order, comparing the \(x_1\)-exponent first. The
initial monomial of an equation for \(D\) is then \(x_n^e\).
In each \(W_{m,j}\) take a basis simultaneously adapted to this
monomial order and to \(D\)-order. The exponent vectors \(\alpha\)
are distinct within each space. They satisfy
\(e\ell\le\alpha_n=j-\sum_{i<n}\alpha_i\) by divisibility of
polynomials and multiplicativity of leading monomials.
For each \(i<n\) lift \(x_i\) to a hyperplane section through
\(p\), possible by very ampleness. Multiplication by its \(k\)th
power carries the leading exponents in \(W_{m-k,j-k}\) injectively
to exponents with \(\alpha_i\ge k\) in \(W_{m,j}\), increasing that
coordinate by \(k\). All occurring exponents are counted by an
adapted basis, since the monomial filtration has one-dimensional
successive quotients. Summing over \(j\) and \(1\le k\le m\) gives
\[
             e L_m\le J_m-(n-1)\sum_{k=1}^m h_{m-k}.
\]
Now \(J_m/(mh_m)\to\mu\) by the definition of expected order and
\(\sum_{k=1}^m h_{m-k}/(mh_m)\to 1/N\) by the Hilbert polynomial.
This proves (AD1).

For clarity, \(N\mu\le n+d-1\) also follows by taking a general
smooth flag of hyperplane sections through \(p\), ending in a
smooth plane curve of degree \(d\) and the point \(p\) on it.
The ordinary order of a section is at most the sum of its flag
exponents, by its expansion in regular local parameters. Use a
basis adapted to both filtrations. Until reaching the curve, the
layer of hyperplane order \(i\) in section degree \(m\) restricts
to the complete system of degree \(m-i\) on the next section, by
the hypersurface sequence. On a \(j\)-dimensional rung the mean of
\(i/m\) tends to \(1/(j+1)\), leaving mean residual degree
\(j/(j+1)\). On the final curve the mean order divided by degree
of the system tends to \(d/2\) by Riemann--Roch. Iteration gives
\((n-1+d)/N\). Substituting in (AD1) yields
\[
 N(A_X(v)-rS_H(v))\ge
       (d-1)(d-2)a+\left(N-\frac{rd}{e}\right)\gamma .
\]
The first term is strictly positive.
\end{proof}

\subsection{Regular tangent forms}
\label{sec:regular-angular-divisor}
There is a stronger restriction in every degree at a point for which
the homogeneous tangent forms \(q_2,\ldots,q_d\) form a regular
sequence in the polynomial ring on \(T_pX\).

\begin{proposition}\label{prop:regular-angular-divisor}
At such a point, every divisor in Theorem~\ref{thm:angular-divisor-bound}
with \(e\ge2\) satisfies \(A_X(v)>rS_H(v)\). More precisely,
writing \(H_{d-1}^{\mathrm{harm}}=\sum_{i=1}^{d-1}1/i\), one has
\[
          N S_H(v)\le (n+H_{d-1}^{\mathrm{harm}})a
                                      +\frac{2}{e}\gamma .
 \tag{AD2}
\]
Thus a possible minimizing divisor over a regular point whose
center on the first blowup has codimension one in \(E\) must
have a hyperplane there as center.
\end{proposition}

\begin{proof}
We give the refinement of the preceding dimension count. Suppose
\(n-1\) sections, of respective section degrees \(b_i\), have ordinary
initial forms of degrees \(h_i'\) whose restrictions \(g_i\) form a
homogeneous system of parameters in \(R=\mathbf C[x_1,\ldots,x_n]/(P)\),
where \(P\) is the irreducible equation of \(D\). The ring \(R\) is
Cohen--Macaulay and finite free over \(\mathbf C[g_1,\ldots,g_{n-1}]\).
Indeed homogeneous lifts of a basis of \(R/(g_1,\ldots,g_{n-1})\)
generate over the polynomial subring by graded reduction, and the
regular-sequence Hilbert series proves independence.

In each layer of \(W_{m,j}\) of exact \(D\)-order \(\ell\), divide
by \(P^\ell\) and restrict to \(D\). This embeds the layer as a
subspace of \(R_{j-e\ell}\). Expand it in the free basis just
described. Row reduction in a monomial order, with a fixed tie order
on the free basis, counts one pivot for each dimension of the
layer. For each pivot \(g^\alpha u\) we have
\(\sum_i h_i'\alpha_i\le j-e\ell\), since the basis element
\(u\) has nonnegative degree. Multiplication by the \(k\)th power
of the \(i\)th chosen section preserves \(\ell\), raises \(j\) by
\(kh_i'\), and takes the pivot to \(g_i^kg^\alpha u\). It follows,
as in the proof of (AD1), that
\[
 e L_m\le J_m-
    \sum_{i=1}^{n-1} h_i'\sum_{k=1}^{\lfloor m/b_i\rfloor}h_{m-kb_i}.
 \tag{AD3}\label{eq:angular-parameter-count}
\]
In particular the asymptotic subtraction after multiplication by
\(N/(mh_m)\) is \(\sum_i h_i'/b_i\).

The hypertangent sections with initial forms \(q_j\) have section
degree \(j-1\) and ordinary order \(j\), up to sign. They are
obtained by truncating the homogeneous equation at degree \(j-1\)
on the affine chart at \(p\) and homogenizing to that degree.
Indeed the first omitted term on \(X\) has initial form \(q_j\),
since the normal coordinate has order at least two. Regularity
gives \(N\mu=n+H_{d-1}^{\mathrm{harm}}\), by the graded presentation
in Section~\ref{sec:nov27-pow-regular-point}. The presentation there
uses only regularity and applies also to \(d=3,4\).

For \(i=3,\ldots,d\), take a general degree-\(i\) combination
of \(q_2,\ldots,q_i\), with homogeneous polynomial coefficients
of degrees \(i-2,\ldots,0\). Each combination lifts to a section
of degree \(i-1\) by multiplying the hypertangent sections by
products of tangent hyperplanes. On the projective hypersurface
\(D\subset\mathbf P(T_pX)\), the
base set of this system is \(D\cap V(q_2,\ldots,q_i)\).
Since the latter \(i-1\) equations form a regular sequence before
restriction, its codimension in \(D\) is at least \(i-2\),
with an empty base set allowed. Inductively the \(i\)th general
choice therefore cuts every component left by the previous
\(i-3\) choices properly. Complete these \(d-2\) choices by
\(n-d+1\) general linear forms. This gives the required system
of parameters on the affine cone. The total subtraction in (AD3)
is
\[
      (n-d+1)+\sum_{i=3}^{d}\frac{i}{i-1}
                    =n-2+H_{d-1}^{\mathrm{harm}} .
\]
The two units remaining from \(N\mu\) prove (AD2).
For \(e\ge2\), \(2r/e\le r<N\). Also
\(r(n+H_{d-1}^{\mathrm{harm}})<nN\), since \(r\le n-1\)
and \(H_{d-1}^{\mathrm{harm}}<d-1\). Both coefficients in the
resulting discrepancy comparison are strictly positive.
\end{proof}

There is a useful refinement for the linear case left by (AD2).
Suppose \(n\ge d\) and that \(D\) contains no irreducible component
of \(V(q_2,\ldots,q_d)\subset E\). Then
\[
      N S_H(v)\le(n+H_{d-1}^{\mathrm{harm}})a+\gamma/e .
\]
Indeed adjoining \(P\) to the full regular sequence still has the
expected height, so \(q_2,\ldots,q_d\) restrict to a regular
sequence on the cone over \(D\). Complete these parameters by
\(n-d\) general linear forms and apply (AD3). The subtraction
becomes \(n-1+H_{d-1}^{\mathrm{harm}}\). Thus a hyperplane which
is the center of a minimizing divisor over a regular point with
\(\delta(X,H)\le r\) must contain a component of that complete
intersection. If the complete intersection is integral, it spans
\(E\) because its homogeneous ideal is generated in degrees at
least two. In this case every such minimizing divisor has center
of codimension at least two in \(E\): the case of center \(E\)
itself is the ordinary order at \(p\), already strictly positive.

\subsection{Tangent sections of multiplicity one below the degree}
\label{sec:contact-blowup}
Suppose \(n\ge3\), \(3\le d\le n+1\), and the tangent hyperplane
section of \(X=X_d^n\) at \(p\) has multiplicity exactly \(d-1\).
Write \(E\) for the exceptional divisor of the ordinary blowup at \(p\),
\(N=n+1\), \(r=n+2-d\), and \(\beta_X(v)=A_X(v)-rS_H(v)\).
The calculation in Chapter~\ref{ch:fivefold} extends to these points
in every degree.

\begin{proposition}\label{prop:contact-blowup}
At such a point,
\[
 B:=N S_H(E)=n+d-2+\frac1{d-1},\qquad
 \beta_X(E)=
 \frac{(d-2)\bigl(n+(d-1)(d-2)+1\bigr)}{(d-1)(n+1)}>0 .
\]
\end{proposition}
\begin{proof}
Take coordinates \([T:u_1:\cdots:u_n:z]\) with \(p=[1:0:\cdots:0]\)
and tangent hyperplane \(z=0\). The equation has the form
\[
 zB_{d-1}(T,u,z)+Tf(u)+g(u)=0 ,
 \qquad \deg f=d-1,\quad \deg g=d,
\]
where \(f\ne0\) and the coefficient of \(T^{d-1}\) in \(B_{d-1}\)
is one. The singular locus of the tangent section is finite. Indeed,
it lies in a fiber of the Gauss map of \(X\), whose hyperplane
pullback \((d-1)H\) is ample. If \(f\) and \(g\) had a common
factor, its zero set and the residual factor of \(Tf+g\) would
meet in dimension at least \(n-2\) in the tangent \(\mathbf P^n\),
and the section would be singular there. Thus \(f,g\) are coprime.

Write \(B_{d-1}=TB'+B_0(u,z)\) and set \(w=zB'+f\), a section
of \((d-1)H\). The identity \(Tw=-g-zB_0\) gives the initial
form \(-g\) for \(w\) on \(T=1\), since \(z\) has order \(d-1\)
and \(zB_0\) has order at least \(2d-2>d\). In particular \(w\)
has order \(d\). With the two gradings given by section degree
and order at \(p\), the associated graded section ring is
\[
 \mathbf C[T,u_1,\ldots,u_n,Z,U]/(T^{d-2}Z+f,\ TU+g).
\]
Here \(Z,U\) are the initial sections of \(z,w\), of bidegrees
\((1,d-1),(d-1,d)\), respectively. To verify the presentation,
reduce the two relations modulo \(T\). They give the regular
sequence \(f,g\), so the two relations form a regular sequence
and \(T\) is a nonzerodivisor in the quotient. After inverting
\(T\), one can eliminate \(Z,U\). The remaining generators map
injectively because \(u_i/T\) form a regular system of parameters
at \(p\) and \(T\) records section degree. The map before
localization is therefore injective. Its source has section-degree
Hilbert series \((1-t^d)/(1-t)^{n+2}\), the Hilbert series of \(X\),
so it is an isomorphism.

In the logarithmic derivative of the bigraded Hilbert series with
respect to the order variable, evaluated at one, each generator
of bidegree \((b,h)\) contributes
\(ht^b/(1-t^b)\), and each relation contributes its negative.
If \(h_m=h^0(X,mH)\), the coefficient of \(t^m\) in the product
of \(t^b/(1-t^b)\) with the section Hilbert series is
\(\sum_{k=1}^{\lfloor m/b\rfloor}h_{m-kb}\). Dividing by
\(mh_m\) gives the limit \(1/(bN)\) by the Hilbert polynomial.
Hence \(NS_H(E)=n+(d-1)+d/(d-1)-1-1=B\). Substitution of
\(A_X(E)=n\) gives the asserted discrepancy difference.
\end{proof}
For the quartic fivefold this gives \(S_H(E)=11/9\) and
\(\beta_X(E)=4/3\).

\begin{corollary}\label{cor:quartic-contact-centers}
Suppose \(d=4\) and the tangent section has multiplicity three.
Let \(v\) be a divisorial valuation whose center on the blowup is
an integral hypersurface \(D=(P=0)\subset E=\mathbf P^{n-1}\).
Each of the following hypotheses implies \(\beta_X(v)>0\).
\begin{enumerate}[label=\textup{(\roman*)},nosep]
\item \(\deg D\ge3\),
\item \(\deg D=2\) and \(P\) does not divide the tangent cubic \(f\),
\item \(D\) is a hyperplane containing no component of \(V(f,g)\subset E\).
\end{enumerate}
\end{corollary}
\begin{proof}
Put \(e=\deg D\), \(a=v(E)>0\), and \(\gamma=A_X(v)-na\ge0\) as
in Theorem~\ref{thm:angular-divisor-bound}. By the basis estimate
\eqref{eq:angular-parameter-count}, if the initials of \(n-1\)
sections of respective degrees \(b_i\) and orders \(h_i\) form
a homogeneous system of parameters modulo \(P\), then
\[
 N S_H(v)\le Ba+\frac{\gamma}{e}
                  \left(B-\sum_{i=1}^{n-1}\frac{h_i}{b_i}\right).
\]
Here \(B=n+7/3\). By the preceding presentation \(f,g\) occur,
up to sign, with ratios \(h_i/b_i=3,4/3\), respectively.
If \(P\mid f\), use \(g\) modulo \(P\), which is nonzero by
coprimality. If \(P\nmid f\), use \(f\). In case \textup{(iii)}
use both \(f,g\), since \((P,f,g)\) has height three.
In each case the specified forms are a regular sequence modulo
\(P\) and general linear forms complete them to parameters.
The linear forms have ratio one. The resulting value of
\(B-\sum_i h_i/b_i\) is \(3\) when only \(g\) is used, \(4/3\)
when only \(f\) is used, and \(1\) when both are used. Denote this
value by \(C\). Then
\[
 \beta_X(v)\ge a\beta_X(E)+\gamma\left(1-\frac{rC}{Ne}\right).
\]
In each of \textup{(i)}--\textup{(iii)}, \(C\le e\).
Since \(r=n-2<N\), Proposition~\ref{prop:contact-blowup}
gives strict positivity.
\end{proof}

\subsection{Higher-codimension centers on the exceptional tangent space}
\label{sec:nov30-audit-angular-centers}
Keep the notation of Theorem~\ref{thm:angular-divisor-bound}. Let the center of a divisorial valuation $v$ on the point blowup $Y$ be any integral $Z\subset E$ of codimension $c$ in $E$, $1\le c\le n-1$. Write $a=v(E)$, $\gamma=A_X(v)-na$, and $\mu=S_H(\operatorname{ord}_p)$. For a general linear $c$-plane $\Lambda\subset E$ the intersection $\Gamma=Z\cap\Lambda$ consists of $\deg Z$ distinct smooth points of $Z$. Define its Waldschmidt constant by
\[
 \sigma(\Gamma)=\inf_{b\ge1}\frac{\min\{t:H^0(\Lambda,\mathcal I_{\Gamma}^{(b)}(t))\ne0\}}b,
\]
where the symbolic power requires vanishing to order at least $b$ at each point. Then
\begin{equation}\label{eq:nov30-audit-arbitrary-angular}
 N S_H(v)\le N\mu a+\frac{N\mu-(n-c)}{\sigma(\Gamma)}\,\gamma .
\end{equation}
In particular $N\sigma(\Gamma)\ge r(d+c-1)$ implies $A_X(v)>rS_H(v)$. For a divisor $Z$ the points lie on $\mathbf P^1$, where $\sigma(\Gamma)=\deg Z$, recovering (AD1).

Let $F$ be the ordinary initial form of a section of order $j$ at $p$ and let $b$ be its multiplicity along the generic point of $Z$ in $E$. We first check $v(s)\le aj+\gamma b$. At a general point of the smooth locus of $Z$, restriction to a general smooth curve transverse to a hypersurface germ of multiplicity $b$ has that same multiplicity. Successive inversion of adjunction from the curve up through a general smooth local flag gives the bound $\operatorname{lct}\ge1/b$ near this point. Apply inversion once more along $E$ as in (AD1). For $b=0$ the residual divisor misses the generic point. This establishes the valuation bound for all $v$ under consideration.

Choose a general linear flag down to $\Lambda$, cut by $n-c-1$ linear forms $x_1,\dots,x_{n-c-1}$. Its cuts of $Z$ are transverse along the general smooth locus, and at the final intersection all points are smooth and transverse. Choose $x_{n-c}$ restricting to a hyperplane of $\Lambda$ missing every point of $\Gamma$. Complete these forms to coordinates. For a homogeneous initial $F$, take its least lexicographic monomial in these coordinates, with exponents $\alpha_i$. For each successive flag hyperplane, divide by its maximal power and then restrict. At each step the restricted residual polynomial has multiplicity at least $b$ along every component of the cut of $Z$. Indeed the dividing form does not contain any such component of the \emph{preceding} cut, hence division preserves generic multiplicity there. At every smooth point of any such component the same lower bound on ordinary multiplicity holds, since the partial derivatives of orders below $b$ vanish on the component. Restriction preserves ordinary multiplicity as a lower bound. The general flag ensures that the generic points in each new cut, and every point in the last cut, are in this transverse smooth open set. For $c=n-1$ start directly on $\Lambda=E$. At the last step simply divide by the maximal power of $x_{n-c}$, which is a unit at every point of $\Gamma$. The resulting nonzero homogeneous polynomial has degree $j-\sum_{i=1}^{n-c}\alpha_i$ and multiplicity at least $b$ at each point. Hence
\[
 j-\sum_{i=1}^{n-c}\alpha_i\ge\sigma(\Gamma)b.
\]
Here the indicated successive powers are precisely the first $n-c$ exponents of the lexicographic initial monomial. In each $W_{m,j}$ take a basis adapted to both this monomial flag and multiplicity along $Z$. Multiplication by the $k$th power of a lift of $x_i$ produces distinct pivots with $\alpha_i\ge k$ from section degree $m-k$, as in (AD1). Thus the sum of the $b$-jumping numbers, multiplied by $\sigma(\Gamma)$, is at most $J_m-(n-c)\sum_{k=1}^m h_{m-k}$. Applying the previously proved valuative bound to a basis adapted to $v$ and point order, and using the maximality of the $b$-adapted sums on $W_{m,j}$, gives \eqref{eq:nov30-audit-arbitrary-angular}. The numerical consequence follows from $N\mu\le n+d-1$ and $nN-r(n+d-1)=(d-1)(d-2)>0$.

\paragraph{A complete intersection at a regular point.}
If $Z$ is a scheme-theoretic complete intersection of forms $P_1,\ldots,P_c$ with degrees $e_i$, let $e=\min e_i$. When $q_2,\ldots,q_d$ are a regular sequence, there is a stronger bound
\begin{equation}\label{eq:nov30-audit-regular-ci-angular}
 NS_H(v)\le (n+H_{d-1}^{\mathrm{harm}})a+
 \frac{c+H_{\min(c,d-1)}^{\mathrm{harm}}}{e}\,\gamma .
\end{equation}
Thus $Ne\ge r(c+H_{\min(c,d-1)}^{\mathrm{harm}})$ suffices to exclude such a divisor. For an integral complete intersection the transverse point set in \eqref{eq:nov30-audit-arbitrary-angular} already has $\sigma(\Gamma)=e$: its homogeneous ideal is generated by a regular sequence of degrees $e_i$ and all its powers are saturated, as follows also from the associated graded calculation below and the depth-one property of its homogeneous ring.

Write $S=\mathbf C[x_1,\dots,x_n]$, $I=(P_1,\dots,P_c)$. The associated graded ring of a regular sequence is $(S/I)[U_1,\ldots,U_c]$; the $U_i$ have homogeneous degrees $e_i$ and ideal order one. The complete-intersection hypothesis is used at this step. Since $S/I$ is a domain all ordinary powers of $I$ are primary, as seen by multiplying in this associated graded domain, so they equal the symbolic powers. For the transverse reduced finite scheme, $S/I$ is instead a reduced one-dimensional Cohen--Macaulay ring. A general linear form avoiding those points is a nonzerodivisor on $S/I$, on the displayed graded ring and hence on all $S/I^b$. Thus $I^b$ is saturated, and locally at the points it is the required power of the point ideal, proving the assertion about $\sigma(\Gamma)$.

Choose $n-c$ homogeneous parameters $g_i$ of the complete intersection ring $S/I$ as follows. If $c\le d-2$, take general degree-$j$ combinations of $q_2,\ldots,q_j$ for $c+2\le j\le d$, as in the proof of (AD2). Their base locus on $Z$ has codimension at least $j-1-c$, since the original $j-1$ forms are a regular sequence in $S$. Thus each new system cuts every component remaining from its $j-c-2$ predecessors. Complete by $n-d+1$ general linear forms. If $c\ge d-1$, take only general linear forms. All parameters lift to sections as before, a choice of degree $j$ having section degree $j-1$. The ring $S/I$ has a homogeneous free-module basis over the polynomial ring in these parameters by Cohen--Macaulayness, exactly as in (AD3).

Each layer with angular multiplicity $b$ now has initial forms inside $I^b/I^{b+1}$. Expand in its monomials $U^\lambda$ with $|\lambda|=b$ times the free basis and parameter monomials $g^\alpha$. For every pivot in section order $j$ this gives $eb+\sum h_i'\alpha_i\le j$, where $h_i'=\deg g_i$. Multiplication by a lifted parameter preserves angular multiplicity and shifts its own pivot exponent by one. The tail count of (AD3) therefore subtracts $n-c+H_{d-1}^{\mathrm{harm}}-H_{\min(c,d-1)}^{\mathrm{harm}}$ in units of $N$. Using $N\mu=n+H_{d-1}^{\mathrm{harm}}$ yields \eqref{eq:nov30-audit-regular-ci-angular}. The $a$ coefficient has strict discrepancy margin as in (AD2), proving the numerical consequence.

\subsection{Local monomial orders and a hyperplane flag}
\label{sec:nov30-hyperplane-monomial}
There is another way to obtain the sharp point bound for a class of
valuations, without a presentation of the associated graded ring.
Let \(p\) be smooth and let \(u_1,\ldots,u_n\) be regular parameters
there, obtained by dividing hyperplane sections by a hyperplane not
through \(p\). Give them positive weights \(w_1\ge\cdots\ge w_n\)
and let \(v\) be the monomial valuation in these parameters. Assume
that, in the indicated order, cutting by the first \(n-1\) sections
gives a flag of integral scheme-theoretic hyperplane sections.
The members of this flag are smooth at \(p\) by independence of the
parameters. Then
\[
 N S_H(v)\le \sum_{i=1}^n w_i+(d-1)w_n,\qquad
 \frac{A_X(v)}{S_H(v)}\ge\frac{n(n+1)}{n+d-1}>r .
 \tag{HF}
\]
Here one can take rational weights and scale them to integers.

Indeed the exponent vector \(\alpha\) for this flag is the least
lexicographic exponent in the power series expansion at \(p\).
The term of that exponent occurs, so the monomial valuation is
at most \(\sum_iw_i\alpha_i\). Take a basis adapted to the flag
and to \(v\). On each integral hypersurface rung the exact sequence
for a hyperplane section identifies the hyperplane-order layers
with complete systems on the next rung, as in the proof of (AD1).
Integrality identifies divisibility at the generic point of the
hyperplane section with divisibility in the graded section ring:
the quotient by its linear equation is the homogeneous domain
of that section. On the integral plane curve at the bottom, the
last mean order is \(d/2\) in degree-one units. This remains true
if the curve is singular away from \(p\), by Riemann--Roch on its
normalization: the lift of the space of sections misses at most
the length of the normalization quotient, uniformly in degree.
Consequently the normalized mean flag exponents are \(1/N\) for
each of the first \(n-1\) entries and \(d/N\) for the last.
Since \(A_X(v)=\sum_iw_i\ge nw_n\), these estimates prove (HF).

For example the integrality hypothesis can be verified by the
following geometric test even when the original flag has reducible
members. Let \(z=0\) be the tangent hyperplane at \(p\), with \(z\)
written in the same local trivialization. Suppose \(w_1\le2w_n\)
and the line \(\Lambda\) cut out by \(u_1,\ldots,u_{n-1},z\) in
the ambient projective space is not contained in \(X\). Then (HF)
holds. In fact one can use \(u_i+\lambda_i z\) for \(i<n\) and
retain \(u_n\). Here \(v(z)\ge2w_n\) because \(z\) has ordinary
order at least two on \(X\). This substitution preserves the
monomial valuation. When equality of weights occurs, the
weight-leading term added to a variable of largest weight depends
only on variables of strictly lower weight, since it has ordinary
degree at least two. The substitution thus induces an invertible
triangular map on the weighted graded polynomial ring.
As the constants vary, the final planes sweep an open chart of all
planes containing \(\Lambda\). A general such plane cuts \(X\)
smoothly away from its finite intersection with \(\Lambda\) by
successive Bertini, and at each of these finitely many points the
plane can be chosen transverse to \(X\). Its smooth plane curve
is integral. Every earlier restricted homogeneous equation is
then irreducible and reduced, since a nontrivial homogeneous
factorization would give a factorization on the final plane or
would make the equation there zero. This supplies the required
flag. The argument proves a bound on the indicated point
valuations even if new generators occur beyond degree one in
their section filtrations.

\section{Weights, generators and complete intersections}

\subsection{Sharp bounds for valuations induced by ambient coordinate weights}
\label{sec:nov30-coordinate-weight}
Let \(X=(F=0)\subset\mathbf P^{n+1}\) be smooth of degree
\(3\le d\le n+1\), and let \(H=\mathcal O_X(1)\).
Fix any projective coordinate system, with integral weights
\(0=a_1\le\cdots\le a_{n+2}\) which are not all zero. Write
\(m=\#\{i:a_i=0\}\), \(U=\sum a_i\), and let \(b\) be the least
weight of a monomial of \(F\) with nonzero coefficient. If
\(b>0\), then
\[
              U\ \ge\ \frac{n+2+m(d-2)}{d}\,b .
\tag{CW1}
\]
This estimate holds for every initial form. When the initial form
\(\operatorname{in}_a F\) is an irreducible polynomial, the induced
weight filtration of the homogeneous coordinate ring defines a
divisorial valuation \(v\) up to positive scaling, using a
weight-zero coordinate to trivialize \(H\). Its threshold satisfies
\[
 \frac{A_X(v)}{S_H(v)}
    \ \ge\ (n+1)\frac{n+2-d+m(d-2)}{n+1+m(d-2)}
    \ \ge\ \frac{n(n+1)}{n+d-1}
    \ >\ n+2-d .
\tag{CW2}
\]
If \(b=0\) and the initial form is irreducible, the same ratio
equals \(n+1\). In particular on a smooth quartic ninefold a
coordinate-weight valuation with prime initial equation and
weight-zero locus a line has ratio at least \(55/7\).

For all the coordinate-weight valuations with irreducible initial
equation considered here, (CW2) gives the scale-invariant estimate
\[
 \frac{A_X(v)-rS_H(v)}{A_X(v)}
 \ \ge\ 1-\frac{(n+2-d)(n+d-1)}{n(n+1)}
 \ =\ \frac{(d-1)(d-2)}{n(n+1)} .
\]
For a quartic ninefold this lower bound is \(1/15\).

To prove (CW1), suppose \(a_j<b/d\) and consider the projective
\((j-1)\)-plane spanned by the first \(j\) coordinate points.
All monomials confined to these coordinates have weight smaller
than \(b\), so the plane lies in \(X\). A partial derivative with
respect to \(x_i\) vanishes identically on this plane unless
\(a_i+(d-1)a_j\ge b\). As \(X\) is smooth, these surviving forms
of degree \(d-1\) have empty common zero set on the plane.
There must therefore be at least \(j\) surviving indices, and
\[
              a_{n+3-j}\ge b-(d-1)a_j>b/d .
\]
Pair each index of weight below \(b/d\) with its distinct
complementary high index. Each pair sums to at least
\(b-(d-2)a_j\ge2b/d\). Each of the \(m\) zero-weight pairs
contributes at least \(b\). Every unpaired index has weight at
least \(b/d\), proving (CW1).

For the valuative assertion, divide all weights by their common
divisor. On the chart \(x_1=1\), extract the toric divisor of the
weight vector in the remaining affine coordinates, allowing weights
zero in the directions of its center. Its ambient log discrepancy
is \(U\) and the equation has order \(b\). The initial hypersurface
is integral and meets the dense torus of the exceptional divisor:
an irreducible homogeneous polynomial of degree \(d>1\) has no
coordinate factor. On that torus the strict transform meets the
exceptional divisor in a generically smooth prime divisor, by
irreducibility and characteristic zero. The toric variety is
smooth at the generic point of its exceptional divisor and on the
torus orbit there. Hypersurface adjunction therefore gives
\(A_X(v)=U-b\) in the indicated normalization. The order of a
class is its induced quotient weight: its leading form restricts
to zero generically on the exceptional hypersurface exactly when
it is divisible by \(\operatorname{in}_a F\). Equivalently its
associated graded algebra is the quotient by this single initial
equation.

For ordinary degree variable \(t\) and weight variable \(q\), the
character of that quotient is
\[
       \frac{1-t^dq^b}{\prod_{i=1}^{n+2}(1-tq^{a_i})}.
\]
Its logarithmic \(q\)-derivative at \(q=1\), divided by the ordinary
Hilbert series, is \(Ut/(1-t)-bt^d/(1-t^d)\). Thus the series of
total weights has leading term \(d(U-b/d)(1-t)^{-n-2}\),
whereas the ordinary Hilbert series has leading term
\(d(1-t)^{-n-1}\). The ratio of total weight in degree \(k\)
to \(k h^0(X,kH)\) tends to \(S_H(v)=(U-b/d)/(n+1)\).
For \(b>0\) the ratio \((U-b)/(U-b/d)\) increases with \(U\).
Equation (CW1) yields the first inequality in (CW2). The right
side increases with \(m\), and the last strict difference has
numerator \((d-1)(d-2)\).

The first bound in (CW2) is attained whenever
\(1\le m\le\lfloor(n+2)/2\rfloor\). Use coordinates
\(x_1,\ldots,x_m,z_1,\ldots,z_m\) and the remaining \(y_j\),
with respective weights \(0,d,1\), on the smooth hypersurface
\[
    \sum_{i=1}^m(x_i^{d-1}z_i+z_i^d)+\sum_j y_j^d=0 .
\]
The two partial derivatives in each \((x_i,z_i)\) block have
common affine zero only at the origin. The initial equation drops
the pure \(z_i^d\) terms and is primitive and linear in \(z_1\):
the constant term is coprime to \(x_1^{d-1}\) since there is
another pair or a \(y_j\). Thus the initial equation is irreducible.
The case without \(y_j\)'s is rescaled to primitive weights, which
preserves the ratio. Here (CW1) is an equality.

This computation also applies if a divisorial filtration on the
section ring has associated graded algebra generated by the
initial classes of a basis of \(H^0(X,H)\). Choose that basis
adapted to the filtration and normalize the weights by subtracting
their minimum times the section degree. In a hyperplane
trivialization with valuation zero this is the preceding
normalization. The initial classes induce a graded surjection
from the polynomial ring in \(n+2\) variables. Its kernel contains
\(\operatorname{in}_a F\), and the ordinary Hilbert functions of
the associated graded ring and of the quotient by that polynomial
both equal the Hilbert function of \(X\). Thus the rings agree,
the polynomial is prime, and (CW2) applies. In particular a
divisorial valuation with \(A_X(v)/S_H(v)\le n+2-d\) requires
additional generators beyond degree one in its associated
graded ring.

\subsection{A uniform Donaldson--Futaki bound for linear degenerations}
\label{sec:nov30-ambient-futaki}
For every degree and dimension in the question, the same weight
inequality gives a quantitative statement about actual test
configurations. Let \(q=n+2,\ N=n+1,\ r=q-d\). Diagonalize a
one-parameter subgroup on the ambient coordinates, subtract its
least section weight, and denote the resulting weights by \(a_i\).
Use the notation \(U=\sum_i a_i\), \(b=\operatorname{wt}_{\min}F\)
of Section~\ref{sec:nov30-coordinate-weight}. The family defined
by \(t^{-b}F(t^{a_1}x_1,\ldots,t^{a_q}x_q)=0\) has central
equation the ordinary degree-\(d\) weight-initial polynomial.
It is flat, since the equation is primitive over the base line.
Give it the polarization \(\mathcal O(r)\).

\begin{proposition}
For these families, including those with reducible or nonreduced
initial equation, the Hilbert and weight invariants in the
conventions of Section~\ref{sec:kd-definitions} satisfy
\[
 \operatorname{DF}
  =\frac{(d-1)(dU-qb)}{dN}
  \ \ge\
   \frac{(d-1)(d-2)}{r\{n(n+d)+1\}}\,J^{\mathrm{NA}} .
\]
In particular this proves a uniform positive lower bound on the
ratio for every nontrivial normal test configuration induced by
an ambient one-parameter subgroup, in all the Fano dimensions.
\end{proposition}
\begin{proof}
Put \(C_k=\binom{k+q-1}{q-1}\). In section degree \(k\) of the
hyperplane polarization the dimension is \(h_k=C_k-C_{k-d}\).
The total weight of all degree-\(k\) monomials in \(q\) variables
is \(kUC_k/q\), by symmetry. Multiplication by the initial equation
is an injection from the degree-\((k-d)\) ambient polynomial
space with a character shift of weight \(b\). Its cokernel has
total weight \(w_k\) given by
\[
       \frac{w_k}{h_k}
          =\frac{kU}{q}+
             \left(\frac{dU}{q}-b\right)\frac{C_{k-d}}{h_k}.
\]
Expanding the binomials gives
\(
 C_{k-d}/h_k=(k-q(d-1)/2)/(dN)+O(k^{-1}).
\)
The linear coefficient of \(w_k/h_k\) is
\(\overline a=(U-b/d)/N\). Minus twice its constant coefficient
is \(\operatorname{DF}\), yielding the stated formula. Replacing
\(k\) by \(rk\) gives the anticanonical polarization.

Suppose some weight is positive. By (CW1) in the preceding
coordinate-weight calculation, \(b\le dU/(n+d)\), also valid
when \(b=0\). Thus \(b<d\max_i a_i\). Every power of a
maximum-weight variable survives modulo the initial polynomial:
divisibility of such a pure power would force that polynomial
to be a scalar multiple of the \(d\)th power of the same variable.
It follows that \(J^{\mathrm{NA}}=r(\max_i a_i-\overline a)\).
Using \(\max_i a_i\le U\), we obtain
\[
 \operatorname{DF}\ge\frac{(d-1)(d-2)U}{N(n+d)},\qquad
 J^{\mathrm{NA}}\le rU\left(1-\frac{n+d-1}{N(n+d)}\right).
\]
Combining the two inequalities proves the result. If all weights
are zero then both invariants vanish. The calculation uses the
coordinate ring of the indicated flat family itself. A general
test configuration can require new section-ring generators of
higher degree, as explained in Section~\ref{sec:nov30-coordinate-weight}.
\end{proof}

\subsection{Weights on nonlinear equations of an ambient complete intersection}
\label{sec:nov30-weighted-ambient-ci}
Let \(X=(F=0)\subset\mathbf P^{n+1}\) be smooth, of degree
\(3\le d\le n+1\), and write \(N=n+1\), \(r=n+2-d\).
Suppose a nonempty integral \(Z\subsetneq X\) has homogeneous
prime ideal in the ambient polynomial ring generated by a
homogeneous regular sequence \(g_1,\ldots,g_c\), of degrees
\(e_i\). Give these generators arbitrary positive integer weights
\(w_i\) and put \(m=\dim Z+1=n+2-c\).
Let \(b>0\) be the weight order of \(F\) in the associated
filtration of the ambient ring. Suppose the initial equation
generates a prime ideal in
\[
   \mathrm{gr}_w\mathbf C[x_0,\ldots,x_{n+1}]
        =\bigl(\mathbf C[x_0,\ldots,x_{n+1}]/I_Z\bigr)
                                   [u_1,\ldots,u_c].
\]
Then the induced order on \(X\) is divisorial up to scaling
and satisfies
\[
  \frac{A_X(v)}{S_H(v)}
     \ \ge\ (n+1)\frac{n+2-d+m(d-2)}{n+1+m(d-2)}
     \ \ge\ \frac{n(n+1)}{n+d-1}
     \ >\ r .
\]
Thus the sharp coordinate-weight exclusion also applies when
the transverse parameters are nonlinear forms of arbitrary degrees.

Choose graded vector-space lifts of a basis of the ambient
coordinate ring of \(Z\). The \(g^\alpha\) times these lifts
give a basis of the ambient polynomial ring, by reduction modulo
the regular sequence and its Hilbert series. In fixed degree a
polynomial has an expansion \(\sum c_\alpha g^\alpha\) in these
lifts. The least \(\sum_i\alpha_iw_i\) with \(c_\alpha\ne0\)
gives the indicated filtration. Multiplying coefficients modulo
\(I_Z\) and the regular sequence presentation show the asserted
polynomial associated graded ring over the domain of \(Z\).
In particular the initial ideal of the principal ideal \((F)\)
is generated by its initial equation. Primality gives the true
induced order after quotienting. At the generic point of \(Z\)
there is a transverse smooth chart with regular parameters
\(g_i\), where one divides by local powers of a weight-zero
hyperplane section that is generically nonzero on \(Z\).
The weighted blowup in these parameters gives the exceptional
initial hypersurface. If its equation has at least two transverse
monomials, it meets the dense torus over the generic point of
\(Z\), and is generically smooth there in characteristic zero.
Adjunction gives, with the same normalization of weights,
\[
             A_X(v)=\sum_i w_i-b .
\]
If its equation has only one transverse monomial, primality in
the full graded ring forces it to be a nonzero scalar times a
single \(u_j\): a coefficient of positive ordinary degree is a
nonunit. In this case \(e_j=d,\ b=w_j\), and the transverse
parameters on \(X\) are the other \(g_i\), giving the same
formula. A common weight scaling has no effect on the ratio.
The Hilbert series of the quotient has weight factor
\((1-t^d z^b)/\prod_i(1-t^{e_i}z^{w_i})\) times the ordinary
Hilbert series of \(Z\)'s homogeneous coordinate ring.
Its logarithmic derivative at \(z=1\) gives
\[
               N S_H(v)=\sum_i w_i/e_i-b/d.
\]

To bound \(W=\sum_i w_i\), adjoin \(m\) formal zero weights
and sort the resulting \(n+2\) numbers as
\(a_1\le\cdots\le a_{n+2}\). For each \(j\) such that
\(t=a_j<b/d\), let \(K=\{i:w_i>t\}\).
On the projective complete intersection
\(V_K=(g_i=0:i\in K)\), the polynomial \(F\) vanishes.
Indeed a monomial without a high generator has weight at most
\(dt<b\), since its number of generator factors is at most
its ordinary degree. Group the remaining monomials as
\(F=\sum_{i\in K}g_i H_i\). Modulo the ideal of \(V_K\),
the coefficient \(H_i\) can be nonzero only if
\[
                  w_i+(d-e_i)t\ge b ;
\]
terms with two high factors give zero there, and the other
factors have total ordinary degree at most \(d-e_i\).
Smoothness of \(X\) implies the nonzero restrictions
\(H_i|_{V_K}\) have no common zero: at such a zero all
partial derivatives of \(F\) would vanish.

If one surviving coefficient has degree zero, then \(e_i=d\)
and \(w_i\ge b\). Therefore
\(\sum_j w_j(1-1/e_j)\ge(1-1/d)b\), so the exact formulas
give \(A_X(v)/S_H(v)\ge N\), which suffices. Otherwise all
surviving sections have positive degree. Their number is at
least \(\dim V_K+1=n+2-|K|\ge j\), by intersection with a
positive-dimensional projective complete intersection, or simply
the projective dimension theorem. The regular-sequence hypothesis
ensures this dimension for every subset of generators. In
particular the \(j\)-th largest weight is at least
\(b-(d-1)t>b/d\).

Pair each weight below \(b/d\) with its complementary large
weight. Each pair sums to at least \(2b/d\), and the partners
of the \(m\) zero weights each have weight at least \(b\).
Each unpaired weight is at least \(b/d\). Consequently
\[
                   W\ge\frac{n+2+m(d-2)}d\,b .
\]
Finally \(\sum_i w_i/e_i\le W\), and
\((W-b)/(W-b/d)\) increases with \(W\). Combining the exact
formulas gives the displayed bound. Its last strict difference
from \(r\) is \((d-1)(d-2)/(n+d-1)\) at \(m=1\); the first
bound increases with \(m\).

\subsection{A divisorial filtration with one additional generator}
\label{sec:nov30-one-extra-generator}
\begin{theorem}\label{thm:nov30-one-extra-generator}
Let \(X=(F=0)\subset\mathbf P^{n+1}\) be smooth, with
\(3\le d\le n+1\), \(H=\mathcal O_X(1)\), \(N=n+1\)
and \(r=n+2-d\). A nontrivial divisorial valuation \(v\) with
\(A_X(v)\le rS_H(v)\) needs at least two additional generators
beyond the degree-one subspace in
\(\mathrm{gr}_v\bigoplus_{j\ge0}H^0(X,jH)\).
In particular a special test configuration whose filtration
has just one additional generator has positive Futaki invariant.
For quartics this excludes a finite double cover of the quadric
degree-one image, including a singular quadric, whenever its
graded ring is generated in degrees \(1,2\) with one new quadric
section.
\end{theorem}
\begin{proof}
Use primitive integral valuation weights on section divisors and
a homogeneous basis of \(H^0(X,H)\) adapted to them. Denote its
weights by \(w_i\), \(1\le i\le n+2\), and their sum by \(W\).
At least one \(w_i=0\), since the complete hyperplane system is
free at the generic center. We can work in a weight-zero
trivialization of \(H\). If the initial classes of degree one
already generate, the coordinate-weight theorem applies. Suppose
there is one further minimal homogeneous generator: it may be
taken homogeneous for the valuation weight as well. Indeed the
quotient by products of lower ordinary degrees inherits this
weight grading and has just one new dimension. Denote it by
\(y\), of ordinary degree \(a\ge2\), and lift it to a section
represented by a homogeneous polynomial \(\widetilde Q(x)\)
of degree \(a\). Choose this representative to maximize its
minimum ordinary coordinate weight \(h=\mathrm{wt}_w\widetilde Q\)
modulo \((F)\), possible in fixed degree. Then
\(\mathrm{wt}_v(y)=h+p\) with \(p>0\): at equal weights its
initial would be supplied by a polynomial in the degree-one
initials. Write \(Q=\operatorname{in}_w\widetilde Q\).

Suppose first \(a<d\). The kernel of the bigraded presentation
\(\mathbf C[x,y]\to\mathrm{gr}_v R\), \(R=\bigoplus H^0(jH)\),
contains \(Q\). The homogeneous Hilbert series of \(R\) is
\((1-t^d)/(1-t)^{n+2}\). Thus there is no relation in ordinary
degree below \(a\) and exactly one relation in degree \(a\).
The kernel is prime, so \(Q\) is irreducible. The quotient by
\(Q\) alone has the ordinary Hilbert series of \(\mathbf C[x]\).
Choose a weight-homogeneous kernel element \(J\) in ordinary
degree \(d\) outside \((Q)\), of weight \(b\). Then
\(Q,J\) form a regular sequence, and their quotient has exactly
the Hilbert series of \(R\), giving the full initial presentation.
Bigraded surjectivity makes the original polynomial presentation
strictly compatible with the filtrations in each section degree
(lift the initial class, subtract to raise the weight, and use
finite dimension). Thus \(J\) lifts to a genuine ordinary degree-\(d\)
relation whose terms all have weight at least \(b\) for
\(\mathrm{wt}(x_i)=w_i,\ \mathrm{wt}(y)=h+p\). Substituting
\(y=\widetilde Q\) into it yields a nonzero scalar times \(F\):
if it gave zero as a polynomial, the relation would be a multiple
of \(y-\widetilde Q\), whose initial ideal is \((Q)\).
Consequently, after rescaling, we have an actual expansion
\[
     F=\sum_{j=0}^{\lfloor d/a\rfloor}
       \widetilde Q^{\,j} H_j(x),\qquad
       \deg H_j=d-aj,\qquad
       \mathrm{wt}_w H_j+j(h+p)\ge b .                     \tag{1}
\]
We assign infinite weight to zero coefficients.

At a generic point of the initial complete intersection all
coordinates and \(y\) are nonzero, since it is integral and all
their initial classes are nonzero. On the affine chart of a
weight-zero coordinate, the toric modification for the weight
vector in \((x,y)\) extracts the divisor defined by the initial
complete intersection in its exceptional torus. This intersection
is generically smooth by integrality, the regular-sequence
presentation and characteristic zero. The two initial equations
are independent in the conormal space there, so the strict
transform meets the exceptional divisor transversely at that
generic point. Adjunction gives the ambient sum of weights minus
the two equation weights \(h,b\). Section orders agree with the
induced filtration by its domain presentation. The bigraded
Hilbert series has ordinary generator degrees \(1^{n+2},a\)
and equation degrees \(a,d\). Its logarithmic weight derivative
at weight parameter one, or the calculation in the coordinate
case, gives
\[
      A_X(v)=W+p-b,\qquad N S_H(v)=W+p/a-b/d.               \tag{2}
\]
All identities use the common normalization of the initial weights.
If the toric weight vector needs to be made primitive, the same
positive scaling divides both sides in the discrepancy and expected
order comparisons.

Here is the numerical bound obtained from smoothness of the
original equation:
\[
    W\ge (n+2)b/d-\frac{p}{\min(a,d-a)}.                   \tag{3}
\]
Put \(t=b/d\), \(\nu=h/a\), \(\nu'=(h+p)/a\) and
\(\gamma=(b-h-p)/(d-a)\), \(\gamma'=(b-h)/(d-a)\).
For \(0<u<t\), let \(L(u)\) count the coordinates of weights
at most \(u\) and \(K(u)\) those of weights at least
\(b-(d-1)u\). Let \(P_u\) be the projective plane on the low
coordinates; its dimension is \(L(u)-1\). A degree-\(e\) form
of weight above \(eu\) vanishes on \(P_u\). Its derivative
in coordinate \(i\) vanishes there if its weight exceeds
\(w_i+(e-1)u\).

If \(u\ge\nu'\), all \(H_j\) vanish on \(P_u\) by (1),
since \(b>du\) and \(j(h+p)\le aju\). The derivatives of
the \(H_j\)'s there vanish in every coordinate of weight below
\(b-(d-1)u\), by the preceding degree bound. Thus \(P_u\subset X\)
and at least \(L(u)\) high derivatives are needed to avoid a
common zero there, by smoothness.

If \(u<\nu'\), restrict to \(\widetilde Q=0\) on \(P_u\)
when \(u\ge\nu\); for \(u<\nu\) the restriction is already
zero. There \(F=H_0=0\) and
\(dF=dH_0+H_1d\widetilde Q\). Only high derivatives of
\(H_0\) can contribute. The coefficient \(H_1\) is zero on
\(P_u\) if \(u<\gamma\); if \(u\ge\gamma'\), a derivative
of \(\widetilde Q\) can contribute only in a coordinate
with weight at least
\(h-(a-1)u\ge b-(d-1)u\). In the remaining interval
\(\gamma\le u<\gamma'\) impose the further equation \(H_1=0\).
It has positive degree since \(a<d\). On this projective
intersection all derivatives except the indicated high ones
vanish. If its imposed equations number \(k\), then it has
dimension at least \(L(u)-1-k\) when the latter is nonnegative.
The remaining degree-\(d-1\) derivatives have no common zero,
hence \(K(u)\ge L(u)-k\), with the inequality automatic when
the right hand side is nonpositive.

In summary \(K(u)\ge L(u)-k(u)\), where \(k(u)\) is zero
for \(u\ge\nu'\) and otherwise the sum of the indicators of
\([\nu,\nu')\) and \([\gamma,\gamma')\).
The excess of weights above \(t\) is at least
\((d-1)\int_0^t K(u)\,du\), by the change of variable
\(b-(d-1)u\); the deficit of weights below \(t\) equals
\(\int_0^t L(u)\,du\). Since
\((d-1)K(u)-L(u)\ge-k(u)\), the deficit in (3) is at most
\(\int_0^t k(u)\,du\). If \(t\ge\nu'\), then
\(\gamma\ge\nu'\), leaving only the first interval of length
\(p/a\). If \(t\le\nu\), only the second can contribute, with
length at most \(p/(d-a)\). Otherwise put \(D=\nu'-t\),
\(0<D<p/a\). Their contributions are at most
\(p/a-D\) and \(t-\gamma=aD/(d-a)\). Their sum is bounded
by the larger of \(p/a\) and \(p/(d-a)\), proving (3).

From (2) and (3) we conclude
\[
 N(A_X(v)-rS_H(v))\ge
   p\left(N-\frac r a-\frac{d-1}{\min(a,d-a)}\right)>0 .
\]
If the minimum is \(a\), the expression in parentheses is
\(N(1-1/a)\). Otherwise use \(N\ge d\) and
\(N-r/a\ge d-1/a\), whereas \((d-1)/(d-a)\le d-1\).

For \(a\ge d\), let \(F_*=\operatorname{in}_wF\), of weight
\(b\). This is irreducible: it belongs to the prime kernel
and there is no relation in degree below \(d\). By maximality
of the chosen lift, \(Q=\operatorname{in}_w\widetilde Q\)
is not divisible by \(F_*\), since subtracting a suitable
multiple of \(F\) would strictly raise its raw weight.
Both initial equations involve only \(x\) and again form a
regular sequence with the full Hilbert series, now directly
from the original equations \(F,y-\widetilde Q\). Thus (2)
applies with the actual coordinate weight \(b\) of \(F\).
Smoothness implies \(W\ge(n+2)b/d\) by the elementary
coordinate-weight estimate (equivalently GIT semistability).
Equation (2) now gives \(N(A_X(v)-rS_H(v))\ge p(N-r/a)>0\).
\end{proof}

\subsection{Several successive generators in a complete-intersection degeneration}
\label{sec:nov30-triangular}
Let \(q=n+2\), \(N=n+1\), and \(r=q-d\).
Present a smooth hypersurface by adjoining auxiliary coordinates
\(y_1,\ldots,y_s\) of ordinary degrees \(2\le a_i<d\).
Take homogeneous equations
\[
 y_i=Q_i(x,y_1,\ldots,y_{i-1}),\qquad H(x,y)=0 ,
\]
where each \(Q_i\ne0\) has degree \(a_i\) and \(H\) has degree
\(d\). Eliminating the auxiliary coordinates gives the original
equation \(F(x)=0\). Assign nonnegative integral weights \(w_j\)
to the \(q\) ordinary coordinates, of minimum zero, and weights
\(\lambda_i\) to \(y_i\). Let \(h_i\) be the minimum monomial
weight of \(Q_i\), and \(b\) that of \(H\). Assume
\(p_i=\lambda_i-h_i>0\). In particular all the weights are
nonnegative by induction. Suppose the \(s+1\) weight-initial
equations form a regular sequence in the polynomial ring in
\(x,y\). Their quotient then has the Hilbert series of
\(\mathbf C[x]/(F)\), so these equations generate the full
initial ideal. The corresponding flat family is given its
anticanonical polarization on the generic fiber. Use the
Hilbert and weight invariants of this family, equivalently those
of a sufficiently divisible Veronese with its exponent taken
into account.

\begin{theorem}\label{thm:nov30-triangular}
Put \(W=\sum_j w_j\), \(P=\sum_i p_i\), and
\[
 C_0=\frac{(d-1)(d-2)}{n+d},\qquad B_0=N(d/2-1).
\]
For every such family,
\[
 \operatorname{DF}\ \ge\
 \frac{J^{\mathrm{NA}}}
 {N\bigl((r+2B_0)/C_0+1\bigr)} .
\]
Moreover \(N\operatorname{DF}\ge rP/2\). The estimate is independent
of the number of auxiliary coordinates and of the factorization of
the central equations.
\end{theorem}
\begin{proof}
The ordinary Hilbert series is
\(D(z)=(1-z^d)/(1-z)^q\). The series of total weights is
\[
 D(z)\left\{\frac{Wz}{1-z}
   +\sum_i\frac{p_i z^{a_i}}{1-z^{a_i}}
   -\frac{bz^d}{1-z^d}\right\}
\]
by differentiating the bigraded complete-intersection series.
Write \(C=W+\sum_i p_i/a_i-b/d\) and \(E=W+P-b\). At \(z=e^{-v}\)
the expression in braces is \(C/v-E/2+O(v)\), while
\(D(e^{-v})=d v^{-N}(1+rv/2+O(v^2))\).
All poles away from one are at most simple. Comparing the first
two Hilbert and weight coefficients therefore gives mean weight
\(Ck/N+(rC/N-E)/2+O(1/k)\) in degree \(k\). In particular
\begin{equation}
 N\operatorname{DF}=NE-rC
 =(d-1)(W-qb/d)+\sum_i(N-r/a_i)p_i .
 \label{eq:nov30-triangular-futaki}
\end{equation}

We bound the negative part of the first term by projective
smoothness. Put \(t=b/d\). For \(0<u<t\), let \(l(u)\) count the
coordinates \(x_j\) with \(w_j\le u\), and \(m(u)\) those with
\(w_j\ge b-(d-1)u\). Work on the projective coordinate space of
the first group. In polynomial expressions on this space always
substitute the graph equations recursively. For each \(i\) define
the homogeneous polynomial of degree \(d-a_i\)
\[
 D_i=H_{y_i}+\sum_{j>i}D_j (Q_j)_{y_i}
\]
on the graph, in descending order. These are the coefficients
appearing in the chain rule. Impose the following equations on
the coordinate space whenever the indicated inequalities hold:
\[
\begin{array}{ll}
 Q_i=0 & \text{if }h_i\le a_i u<\lambda_i,\\
 D_i=0 & \text{if }h_i<b-(d-a_i)u\le\lambda_i .
\end{array}
\]
Let \(k(u)\) count the imposed equations, including repetitions.
They have positive ordinary degrees. On their common zero set,
every coordinate \(y_i\) with \(\lambda_i>a_i u\) vanishes.
Indeed it either vanishes by the first imposed equation, or
\(h_i>a_i u\), in which case every monomial of \(Q_i\) contains
an earlier variable of weight greater than \(u\) times its degree.
This proves the assertion by ascending induction. Consequently
any homogeneous polynomial all of whose monomials have weight
strictly above \(u\) times its degree vanishes there.

On this same set \(D_i=0\) whenever \(h_i<b-(d-a_i)u\).
If the second equation was not imposed, then
\(\lambda_i<b-(d-a_i)u\), so \(H_{y_i}=0\) by the preceding
monomial observation. In a remaining summand with \(D_j\ne0\)
descending induction gives \(h_j\ge b-(d-a_j)u\), hence
\(h_j-\lambda_i>(a_j-a_i)u\). Thus \((Q_j)_{y_i}=0\) as well
(a negative ordinary degree means the derivative is zero).
The identical argument for an ordinary coordinate of weight
\(w_j<b-(d-1)u\) gives \(F_{x_j}=0\), now using the chain
rule \(F_{x_j}=H_{x_j}+\sum_i D_i(Q_i)_{x_j}\).
Also \(F=0\) on this set since \(b>du\). Smoothness implies that
imposing in addition the \(m(u)\) remaining ambient derivatives
gives an empty projective intersection. By the projective dimension
theorem \(m(u)+k(u)\ge l(u)\), including the case \(l(u)=0\).

Integration of the low and high weight counts gives
\[
 W-qt\ \ge\ \int_0^t\bigl((d-1)m(u)-l(u)\bigr)\,du
 \ \ge\ -\int_0^t k(u)\,du .
\]
Indeed the negative parts below \(t\) are integrated exactly, and
the positive parts above \(t\) are counted only up to \(dt\).
For the \(i\)th generator the lengths of the two indicated bands,
clipped to \(0<u<t\), have sum at most \(p_i/\min(a_i,d-a_i)\).
To see this, partition the interval \([h_i,\lambda_i]\) at
\(a_i t\). Its portion below \(a_i t\) contributes its length
divided by \(a_i\) through the first band, and its portion above
\(a_i t\) contributes at most its length divided by \(d-a_i\)
through the second. Substituting into
\eqref{eq:nov30-triangular-futaki}, the coefficient of \(p_i\)
is bounded below by
\[
 N-r/a_i-\frac{d-1}{\min(a_i,d-a_i)}
 \ \ge\ r(1-1/a_i)\ \ge\ r/2 .
\]
For \(b=0\) this conclusion follows directly from the same
identity.

For the asserted norm estimate let \(b_0\) be the ordinary
coordinate minimum weight of \(F\). Substituting in descending
order, each occurrence of \(y_i\) can lose at most \(p_i\) in
weight, and each degree-\(d\) monomial has at most \(d/a_i\)
such occurrences at that step. Thus \(b\le b_0+d\sum_i p_i/a_i\).
The ordinary coordinate bound (CW1) gives
\((d-1)(W-qb_0/d)\ge C_0 W\), also when \(b_0=0\).
Equation~\eqref{eq:nov30-triangular-futaki} now gives
\(N\operatorname{DF}\ge C_0 W-B_0P\).
Finally every weight per unit degree is at most \(W+\sum_i p_i/a_i\)
by ascending induction in the graph equations. All section
weights are nonnegative. Hence \(J^{\mathrm{NA}}\le r(W+P/2)\).
Using successively \(P\le2N\operatorname{DF}/r\) and
\(C_0W\le N\operatorname{DF}+B_0P\) proves the claimed constant.
\end{proof}
For a quartic this tests arbitrarily many quadratic and cubic
generators whenever their successive graph relations have the
stated regular initial equations. Finite generation of a
destabilizing filtration alone does not supply this presentation.

\begin{corollary}\label{cor:nov30-all-ci-hyperplane-degenerations}
Let a multiplicative integral filtration of the hyperplane section
ring of a smooth degree-$d$ Fano hypersurface, $3\le d\le n+1$,
be exhaustive and separated with finitely many jumps in each
section degree. Suppose its associated graded ring is a finitely
generated complete intersection. Its Donaldson--Futaki invariant
and non-Archimedean $J$-norm satisfy the same uniform estimate as
in Theorem~\ref{thm:nov30-triangular}. In particular no nontrivial
destabilizing special test configuration of such a hypersurface
can have a complete-intersection associated graded hyperplane ring.
There is no bound on the degrees or number of its generators.
\end{corollary}
\begin{proof}
Use the notation of Theorem~\ref{thm:nov30-triangular}. Let
$R=\mathbf C[x_0,\ldots,x_{n+1}]/(F)$ be the original section
ring. Choose a filtration-adapted basis of $R_1$, shifting all
weights by a multiple of section degree so that its least weight
is zero. Extend the initial forms of this basis to a minimal
bihomogeneous system of algebra generators of $S=\operatorname{gr}R$.
Write the additional generators as $y_i$, of section degrees
$a_i\ge2$ and filtration weights $\lambda_i$ for $1\le i\le s$.
Lift each one to a section of $R$ and then to an ordinary
homogeneous polynomial in $x$. In the polynomial ring $P$ on
$x,y$, the kernel $I$ of the induced map to $R$ is the graph ideal
of these polynomials together with $F$.

The central ideal $J=\ker(P\to S)$ equals the weight-initial
ideal of $I$. Indeed a weight-homogeneous relation of $J$
evaluates in $R$ to a section of higher filtration weight; its
initial form is a polynomial in the chosen generators. Subtract
a lift of that expression and repeat. This procedure terminates
in each section degree since its filtration has finitely many
jumps. It therefore lifts each bihomogeneous relation of $J$ to
one of $I$ with that initial form. Minimality gives
$J\subset(P_{>0})^2$, with $P_{>0}$ the homogeneous maximal
ideal for section degree. The complete-intersection hypothesis
gives a minimal bihomogeneous regular sequence generating $J$.
Indeed its localization at the homogeneous maximal ideal is
generated by its height many elements in this minimal polynomial
presentation; homogeneous Nakayama and Cohen--Macaulayness give
the assertion. There are $s+1$ equations. If $b_j$ are their
section degrees, equality of ordinary Hilbert series gives
\[
  \prod_j(1-z^{b_j})=(1-z^d)\prod_i(1-z^{a_i}).
\]
The multisets of exponents agree, by comparing the coefficient at
the least degree where their multiplicities could differ.

Lift these central equations to $E_j\in I$ as above. They
generate $I$: to express any relation in them, subtract an
expression lifting the corresponding initial relation and
iterate in its fixed section degree. The minimum monomial weight
strictly increases, in a finite set. In section degree $a\ge2$,
write the linear coefficients of the $E_j$ in the variables
$y_i$ of degree $a$ as a matrix, with those equations as rows.
It has full column rank. In fact $I$ contains the graph equation
of every such variable, and products of lower-degree equations
with polynomials of the complementary positive degree cannot
have a term linear in a degree-$a$ variable. The number of rows
is the number of these variables plus one precisely when $a=d$.

Process the rows of each such matrix in decreasing order of their
minimum monomial weights. Subtract earlier pivot rows to clear
their columns in subsequent rows, and rescale the pivots. These
operations leave each minimum weight unchanged: a row of larger
weight cannot affect the initial form, and within a tie the
invertible row operations act on linearly independent central
equations. The resulting initials still form a minimal regular
sequence. Every pivot row has the form $y_i-Q_i$; every variable
in its linear part has weight strictly above the row minimum,
since the initial form contains no linear term. Thus
$h_i=\operatorname{wt}_{\min}Q_i<\lambda_i$. The graph is
triangular: order variables by increasing section degree and,
within a degree, by reverse order of extraction of the pivots.
A pivot row can involve only lower-degree variables and later
pivot columns of the same section degree. The single pivot-free
row, of degree $d$, has the form $H(x,y)=0$. It contains no
variables of section degree at least $d$. Eliminating the graph
equations gives the original $F$ up to nonzero scalar because
the lifted rows generate $I$.

If every $a_i<d$, this is precisely the presentation of
Theorem~\ref{thm:nov30-triangular}. In general its Hilbert and
weight identity~\eqref{eq:nov30-triangular-futaki} still holds
for the full regular sequence. Variables of degrees $a_i\ge d$
occur neither in $H$ nor in graph equations of degree below $d$.
The smoothness and band-length argument therefore uses exactly
the indices with $a_i<d$. The other indices retain the entire
positive coefficient $N-r/a_i\ge r/2$ of $p_i=\lambda_i-h_i$.
This gives $N\operatorname{DF}\ge r\sum p_i/2$. The norm argument
is unchanged: in substitution into $H$ the omitted variables
occur zero times, and the bound on generator weights per unit
degree applies to every graph equation by triangular induction.
If $s=0$ use the ordinary coordinate estimate. Finally, finite
generation of the associated graded here also gives finite
generation of the corresponding Rees algebra, by lifting its
generators and adjoining the parameter, using the same finite
reduction in each section degree.
\end{proof}

\subsection{Several extra quadratic generators in a quartic filtration}
\label{sec:nov30-bin-extra-quadrics}
Let \(X\subset\mathbf P^{n+1}\) be a smooth Fano quartic,
\(N=n+1, r=n-2\). Suppose a divisorial valuation on its
hyperplane section ring has associated graded generated by the
degree-one classes and \(k\) further classes of degree two.
If the quadratic relations among the degree-one classes form a
prime complete intersection of codimension \(k\), then
\(A_X(v)>rS_H(v)\). This always holds when \(k=2\).
It also holds for \(k=3\) whenever the three coordinate
quadrics have height three.

Indeed write \(m=n+2\), choose the adapted linear coordinates
\(x_i\) with weights \(w_i\), and write \(W=\sum w_i\).
The vector space of quadratic sections is the full degree-two
polynomial space, since the original equation is quartic.
Split its raw minimum-coordinate-weight filtration and its
valuation filtration simultaneously. The latter dominates
the former. In a common adapted basis, the elements whose
weights agree give independent degree-two products in the
target graded ring. The remaining elements \(P_i(x)\) have
raw weights \(h_i\) and valuation weights \(b_i=h_i+p_i\),
with \(p_i>0\). Their number is exactly \(k\). Their raw
initials \(Q_i\) give a basis of the degree-two coordinate
kernel, and their actual initials \(y_i\) give the additional
generators.

Assume first that these \(Q_i\) define a prime complete
intersection. The ring with the \(y_i\) adjoined and with
just these relations has the Hilbert series of the polynomial
algebra in \(m\) variables. In degree four there is one more
bihomogeneous relation \(J\), of weight \(b\), outside
\((Q_i)\). It is regular by primeness and the quotient now
has the whole quartic Hilbert series. The lifts \(y_i-P_i\)
generate the kernel of substitution alone, with initial
ideal \((Q_i)\): the inclusion and equality of Hilbert
series follow from their regular sequence and the polynomial
ring after substitution. Strict lifting in each ordinary
degree, as in the one-generator argument, gives an actual
degree-four relation lifting \(J\). Its substitution is
nonzero because its initial is outside \((Q_i)\); rescale
to obtain the original quartic \(F\). Write its expansion
\[
   F=\sum_{\boldsymbol\alpha}H_{\boldsymbol\alpha}(x)
                         \prod_i P_i(x)^{\alpha_i},\qquad
   \mathrm{wt}_w H_{\boldsymbol\alpha}
                     +\sum_i \alpha_i b_i\ge b .
\]
The coefficient degree is \(4-2\sum_i\alpha_i\).
Toric adjunction at the generic smooth torus point of the
initial complete intersection and the bigraded Hilbert
derivative give
\[
       A_X(v)=W+\sum_i p_i-b,\qquad
       N S_H(v)=W+\tfrac12\sum_i p_i-b/4 .
\]

Put \(t=b/4\); at \(0<u<t\) call a pivot active if
\(b_i>2u\). In the projective coordinate space of variables
of weight at most \(u\), impose \(P_i=0\) for each active
pivot, charging one equation only if \(h_i\le2u\).
Every coefficient term with only inactive factors now has
its coefficient vanishing by the weight bound \(b>4u\);
its differentiated coefficient can only contribute in a
coordinate of weight at least \(b-3u\). For an active
pivot the coefficient of \(dP_i\) in \(dF\) is the
substitution of \(\partial_{y_i}\) of the lifted quartic.
It vanishes on this intersection whenever \(b_i<b-2u\):
all remaining factors are inactive and have weights at
most their degrees times \(u\). If \(h_i\ge b-2u\), the
differential \(dP_i\) itself has only contributions in
coordinates of weight at least \(b-3u\). Otherwise impose the
additional coefficient equation, of ordinary degree two.
With \(L(u),K(u)\) the low and high coordinate counts of the
one-generator argument, smoothness and projective dimension
give
\[
 K(u)\ge L(u)-\sum_i k_i(u).
\]
Here \(k_i(u)=0\) for an inactive pivot and otherwise is the
sum of the indicators of \([h_i/2,b_i/2)\) and
\([(b-b_i)/2,(b-h_i)/2)\). The integrated charge from each
pivot is at most \(p_i/2\): relative to \(t=b/4\) these
two intervals are reflections, and the second is counted
only below \(b_i/2\). Equivalently apply the one-generator
interval calculation with \(d=4,a=2\) to each pivot.
Integrating \(3K(u)-L(u)\) therefore yields
\[
          W\ge m b/4-\tfrac12\sum_i p_i,\qquad
          N(A_X(v)-rS_H(v))\ge\tfrac N2\sum_i p_i>0 .
\]

For \(k=2\), two independent coordinate relations are
irreducible since the coordinate kernel is prime and has no
linear forms. They define a degree-four complete intersection.
Any component containing the nondegenerate degree-one image
is nondegenerate of codimension two, hence has degree at least
three. Degree three would give a variety of minimal degree;
its ideal contains three independent quadrics, all in the
coordinate kernel, which has only two. Thus this component
has degree four and multiplicity one. The complete intersection
is Cohen--Macaulay, so it has no embedded components and is prime.

For \(k=3\) of height three let \(D=\mathbf C[x]/(Q_i)\).
Choose a linear Noether normalization \(P\) of this complete
intersection. It is finite graded free over \(P\); moreover
\(\operatorname{Hom}_P(D,P)\simeq D(3)\).
Indeed reduction by the normalizing linear forms is an
Artinian intersection of three quadrics, with a perfect
socle pairing of top degree three. Lift a homogeneous free
basis and its top functional. Multiplication into this
functional induces the claimed isomorphism by the graded
Nakayama lemma and equality of ranks.

If \(D\) had several minimal primes, pick one \(\mathfrak p\)
different from a minimal prime \(\mathfrak q\) containing
the degree-one image geometrically (so \(\mathfrak q\) is
contained, as an ideal, in the coordinate kernel). The field
trace from \(D/\mathfrak p\) to \(P\) is regular since \(P\)
is normal and the extension integral. It is homogeneous of
degree zero by equivariance, and nonzero in characteristic
zero. Composing with the quotient, the duality identifies
this functional with a nonzero \(s\in D_3\) annihilated by
\(\mathfrak p\). An element of \(\mathfrak p\setminus
\mathfrak q\) then shows \(s\in\mathfrak q\). When \(D\) has
one minimal prime but is nonreduced, it is generically
nonreduced by Cohen--Macaulayness. At its generic Artinian
local ring the same annihilated element cannot be a unit,
so again it belongs to the selected prime. Either case
gives a nonzero extra cubic in the coordinate kernel modulo
\((Q_i)\). But after adjoining the three \(y_i\), the raw
complete intersection already gives the correct Hilbert
function through degree three. This contradiction proves
primality.

\paragraph{The sole remaining configuration with three extra quadrics.}
At height two, the three quadrics generate the ideal of a
nondegenerate degree-three variety of codimension two. Indeed
take any codimension-two component of their common zero locus.
If nondegenerate its degree is at least three, and degree four
would make it the entire prime complete intersection of the
first two quadrics, contradicting the third independent
equation. A variety of minimal degree has quadratically
generated ideal with exactly three generators here. If the
component is non-linear but contained in a hyperplane, its
ideal is \((l,g_e)\), for \(e\ge2\). Some nonzero combination
of the three quadrics would be divisible by \(l\), which is
impossible in the prime coordinate kernel without linear forms.

It remains to exclude a linear component \((l_1,l_2)\).
Write \(Q_i=l_1 A_i+l_2 B_i\), and consider the four rows
\(r_0=(-l_2,l_1)\), \(r_i=(A_i,B_i)\) of linear forms.
They are independent since \(r_0\wedge r_i=-Q_i\) are
independent. On the degree-one image the first row is
generically nonzero, so all pairwise determinants vanish
there. The determinant map from the exterior square of
this four-space of rows therefore has image of dimension
at most three. Its kernel contains a nonzero decomposable
wedge, since its projectivization meets the Pl\"ucker
quadric. Two independent linear-form rows with identically
zero determinant must be products of two distinct linear
forms with the same constant vector: if a row has coprime
linear entries, the other is a constant multiple, by unique
factorization. The determinant of either of the two
independent rows with \(r_0\) is then a nonzero reducible
quadric in the coordinate kernel, which is impossible.
This proves the assertion.

\paragraph{Equal-degree generators on hypersurfaces of arbitrary degree.}
More generally, let a smooth Fano hypersurface have degree
\(d\), index \(r=n+2-d\), and exactly \(k\ge1\) additional
generators all in degree \(a\), where \(2\le a<d\). As above,
split the two filtrations on the full degree-\(a\) polynomial
space. Let \(Q_i\) be the resulting coordinate relations and
write \(b_i=h_i+p_i\) for their actual lifting weights, with
\(h_i\) their raw coordinate weights. If the \(Q_i\) form
a prime complete intersection of codimension \(k\), then
\(A_X(v)>rS_H(v)\). It is enough that they have height \(k\)
whenever \(k(a-1)<d\).

In fact after adjoining the new degree-\(a\) variables the
raw complete intersection has the polynomial Hilbert series,
as before. In the prime case the first independent relation
is the final one of degree \(d\) and some weight \(b\).
Its strict lift substitutes to a nonzero multiple of the
original equation. In the derivative estimate call a pivot
active when \(b_i>au\) and set its lift to zero on the space
of coordinates of weight at most \(u\). Apart from the
interval \([h_i/a,b_i/a)\) this equation is automatic.
Vanishing of the coefficient of its differential only costs
another equation, of positive degree \(d-a\), in the interval
\([(b-b_i)/(d-a),(b-h_i)/(d-a))\) while the pivot is active.
Inactive factors have assigned weights at most their ordinary
degrees times \(u\), so the one-generator derivative proof
applies term by term. Integrating at \(0<u<b/d\) charges
each pivot at most \(p_i/\min(a,d-a)\), as in that proof.
Toric adjunction and the bigraded Hilbert derivative yield
\[
 N(A_X(v)-rS_H(v))\ge
  \left(N-\frac{r}{a}-\frac{d-1}{\min(a,d-a)}\right)
    \sum_i p_i>0 .
\]
Finally, for a height-\(k\) intersection the Artinian linear
reduction of the coordinate ring has socle in degree
\(k(a-1)\). The trace argument above in that degree would
give an extra independent coordinate relation if the
complete intersection were not prime. Its ordinary degree
is strictly below \(d\) under the stated assumption, while
the two Hilbert functions agree in every such degree.

\section{Finite hyperplane maps on special test configurations}

\subsection{The complete section ring of a finite hyperplane map}
\label{sec:nov30-finite-tower}
\begin{theorem}\label{thm:nov30-finite-tower}
Let \(S\) be a normal Cohen--Macaulay positively graded complex
domain with \(S_0=\mathbf C\) and
\[
  \operatorname{Hilb}_S(z)=\frac{1-z^d}{(1-z)^{n+2}}.
\]
Suppose \(S\) is finite over the subring \(T\) generated by
\(S_1\). Then \(T\) is a normal hypersurface ring, and \(S\) is
a graded complete intersection. More precisely, beginning with
\(T\), one obtains \(S\) by successively adjoining integral
elements with one monic homogeneous relation at each step.
If the degree of the degree-one image is \(e\), the first
additional generator, when needed, has section degree \(e\).
After a field extension of degree \(l\), the next generator,
if needed, has degree \(el\), and so on. The final product is
\(d\).
\end{theorem}
\begin{proof}
Choose a general linear Noether normalization
\(A=\mathbf C[u_0,\ldots,u_n]\subset T\); it exists as \(T\)
is standard graded. Both \(T\) and \(S\) are finite over \(A\).
The \(u_i\) form a regular sequence on \(S\) by the
Cohen--Macaulay assumption. The finite \(A\)-module \(S\)
is graded free, for example by the graded Auslander--Buchsbaum
formula. Its basis degrees, counted with multiplicities, are
\(0,1,\ldots,d-1\) by the Hilbert series.

The ring \(T\) has dimension \(n+1\) and \(n+2\) independent
linear generators, so its unique defining equation has some
degree \(e\ge2\). Equivalently \(T=A[x]\) for a remaining
linear generator \(x\), whose defining equation is monic over
\(A\) by the choice of the normalization. Hence \(T\) is free
over \(A\) with basis degrees \(0,\ldots,e-1\).

Consider any intermediate graded subring \(A\subset B\subset S\)
which is free over \(A\) with basis degrees \(0,\ldots,k-1\).
In particular \(k\le d\) by field degrees. The two Hilbert
series imply \(B_j=S_j\) for \(j<k\).
Let \(K=\operatorname{Frac}(B)\) and \(C=S\cap K\).
The ring \(S\) is the integral closure of \(A\) in its fraction
field, since it is finite over \(A\) and normal. Thus \(C\) is
the integral closure of \(A\) in \(K\), and also the normalization
of \(B\). The field trace from \(\operatorname{Frac}(S)\) to \(K\),
divided by the field degree, restricts to an \(A\)-linear retraction
\(S\to C\). Indeed traces of integral elements are integral and
the map is the identity on \(K\). This retraction preserves
ordinary degrees, since the scalar action defining the grading
preserves the field extension and commutes with trace.
Consequently \(C\) is graded projective over \(A\), hence graded
free. This last implication also follows by lifting a homogeneous
basis modulo the irrelevant ideal of \(A\) and using projectivity
at that maximal ideal. Its rank is \(k\), the dimension of \(K\)
over \(\operatorname{Frac}(A)\).

Since \(B_j=C_j=S_j\) for \(j<k\), the vector spaces of minimal
\(A\)-module generators of \(B\) and \(C\) agree in those
degrees. These already account for all \(k\) generators of
the free module \(C\). Thus \(B=C\). In particular \(B\) is
normal, a conclusion we can use at each stage starting with
\(B=T\).

If \(k<d\), choose \(y\in S_k\setminus B_k\), possible since
the two free modules first differ there by one dimension.
The element is integral over \(B\). Normality gives
\(y\notin K\), so its field degree \(l\) over \(K\) is at
least two. Its monic minimal polynomial has coefficients in
\(B\): they are integral over \(B\) and belong to \(K\).
Uniqueness under the scalar action makes it homogeneous of
section degree \(kl\) on assigning degree \(k\) to its
variable. Monic division and minimality show that this is
the full relation over \(B\). Hence \(B[y]\) is free over
\(B\) on \(1,y,\ldots,y^{l-1}\), and is free over \(A\)
with basis degrees \(0,\ldots,kl-1\). If there is also a
commuting weight grading preserving \(B\), \(y\) can be chosen
as a weight vector, and the relation is weight homogeneous.
The new rank \(kl\) divides \(d\) by the tower law for fields.
Iteration terminates at rank \(d\), when the two Hilbert
series force equality with \(S\). Adjoining one monic equation
at a time to the hypersurface presentation of \(T\) gives a
regular sequence. This proves the theorem.
\end{proof}

\begin{corollary}\label{cor:nov30-basepoint-necessary}
For every \(3\le d\le n+1\), a nontrivial special test
configuration of a smooth \(X_d^n\) with nonpositive
Donaldson--Futaki invariant has a base point of the limiting
degree-one system. Equivalently, its graded hyperplane section
ring fails to be finite over the algebra generated by its
degree-one part. Furthermore this graded ring cannot be a
complete intersection. For smooth quartics of dimension at least
four a special configuration of this kind cannot have terminal
factorial central fiber.
\end{corollary}
\begin{proof}
In the punctured-product trivialization, close a hyperplane to an
invariant Weil divisor on the total space. Its \(r\)th multiple
agrees with the relative anticanonical class up to the integral
principal special fiber. Denote its restriction as a reflexive
class by \(L\). As in the proof of
Proposition~\ref{prop:nov30-factorial-quartic-replacement}, the
rank-one reflexive sheaves of its multiples are Cohen--Macaulay
on the klt total space, and their restrictions are the reflexive
sheaves of \(jL\) on the normal klt special fiber \(Z\).
The restriction assertion follows from depth and agreement at
codimension-one points of \(Z\). Vanishing for \(j\ge0\) and
base change therefore identify \(S=\bigoplus_{j\ge0}H^0(Z,jL)\)
with the associated graded ring of the induced hyperplane
filtration, up to the linear character in each degree.
Here one may use the Weil form of Kawamata--Viehweg vanishing
for integral \(\mathbf Q\)-Cartier divisors on a klt variety.
In particular \(S\) has the hypersurface Hilbert series. It
is normal as a full Weil section ring and finitely generated
since a multiple of \(L\) is ample Cartier. Passing to
anticanonical multiples gives the test configuration and its
invariants, with the polarization convention of the introduction.
Corollary~\ref{cor:nov30-all-ci-hyperplane-degenerations}
excludes the complete-intersection case.

If the ideal \(S_1 S\) has empty projective zero set, \(S/S_1S\)
is finite dimensional and graded Nakayama makes \(S\) finite over
\(\mathbf C[S_1]\); conversely finiteness gives this emptiness.
In that case the degree-one sections trivialize \(L\) on the
corresponding affine opens: a nonzero value of a sufficiently
divisible power gives a section with zero Weil divisor there.
Thus \(L\) is Cartier and basepoint free. Its full section ring
is Cohen--Macaulay: for \(j\ge0\) ordinary Fano vanishing applies,
and for \(j<0\) use Serre duality and vanishing for \(K_Z-jL\)
to obtain \(H^i(Z,jL)=0\) for \(0<i<n\). Here \(K_Z\sim-rL\)
is Cartier and \(Z\) is Cohen--Macaulay. The standard local
cohomology description of a section ring gives the assertion at
the vertex, with depth at least two from normality. Away from
the vertex the ring is locally a line-bundle torsor over \(Z\)
because the degree-one sections cover. The theorem would now
make \(S\) a complete intersection, a contradiction. The last
claim also uses the global generation assertion of
Proposition~\ref{prop:nov30-factorial-quartic-replacement}.
\end{proof}
The corollary concerns special test configurations of a fixed
smooth hypersurface. A K-polystable replacement of a varying
family is addressed separately in
Section~\ref{sec:nov30-special-fibers}.

% ===== Surface sections, point valuations and filtrations =====
\chapter{Surface sections, point valuations and filtrations}\label{ch:foundations}
\section{Basis divisors on the surface}
\begin{proposition}\label{foundations:split}
Let $S$ be a smooth complex projective surface of Picard number one, let $H$ be an ample integral divisor with $H^2=d$, and let $p\in S$. Suppose an integral curve $C\equiv hH$ has $k=\operatorname{mult}_p C$ and $\tau=k/h>\sqrt d$. Set $\varepsilon=d/\tau$. For every divisorial valuation $v$ whose center contains $p$,
\[
 3S_H(v)\ \le\ \varepsilon A_S(v)+\frac{\tau-\varepsilon}{k}v(C).
\]
In particular, if $\lambda=\operatorname{lct}_p(S,C)$, then
\[
 \delta_p(S,H)\ \ge\frac{3}{\varepsilon+(\tau-\varepsilon)/(k\lambda)}.
\]
The inequality for the ordinary blowup is an equality.
\end{proposition}

\begin{proof}
Intersection with $C$ shows that for an integral $B\equiv bH$ different from $C$, one has $\operatorname{mult}_p(B)\le dbh/k=\varepsilon b$. Thus $\tau$ is the pointwise maximum of the ordinary multiplicity divided by the $H$-degree coefficient, and $\varepsilon$ is the Seshadri constant. In particular $\operatorname{mult}_p D\le m\tau$ for $D\in|mH|$.

Write $V_m=H^0(S,mH)$, $N_m=\dim V_m$, and $F^jV_m$ for the sections of order at least $j$ at $p$. Put $r_j=\dim F^jV_m/F^{j+1}V_m$ for $0\le j\le\lfloor m\tau\rfloor$. For a section with divisor $D$ of order precisely $j$, let $b$ be the multiplicity of its component $C$. As $D-bC$ has no component $C$, local intersection bounds give
\[
 h(m-bh)d=(D-bC)\cdot C\ \ge\ k(j-bk).
\]
Thus a divisor $D$ with order at least $j$ contains at least
\[
 a_j=\left\lfloor\frac{(j-m\varepsilon)_+}{k-h\varepsilon}\right\rfloor
\]
copies of $C$. Here $x_+=\max(x,0)$. Divide out these copies in the whole filtered piece $F^jV_m$. Leading forms, in local parameters at $p$, give an injection
\[
 F^jV_m/F^{j+1}V_m\hookrightarrow
 H^0(E,\mathcal O_E(q_j)),\qquad q_j=j-ka_j,\quad E=\mathbf P(T_pS).
\]
Indeed each leading form in degree $j$ is divisible by the same degree-$ka_j$ tangent form of $a_j C$, and the division has kernel precisely the next filtered piece. One has $q_j\ge0$ for the entire displayed range. Uniformly in $j$,
\[
 q_j=m g(j/m)+O(1),\quad
 g(t)=\begin{cases}t&0\le t\le\varepsilon,\\
  \varepsilon(\tau-t)/(\tau-\varepsilon)&\varepsilon\le t\le\tau.
 \end{cases}
\]
Since $\int_0^\tau g(t)\,dt=\varepsilon\tau/2=d/2$ and $N_m=dm^2/2+O(m)$ by Riemann--Roch and Serre vanishing, the dimension deficits satisfy
\[
 0\le\sum_j (q_j+1-r_j)=O(m).
\]

Choose a basis simultaneous for the two filtrations $F^\bullet$ and the valuation filtration of $v$. Such a basis exists for every pair of finite flags, and the filtration jumps of either filtration on a finite-dimensional space give a finite flag. Denote the degree-$j$ basis vectors by $s_{ji}$ for $1\le i\le r_j$. Define effective rational divisors by
\[
 D_m=\frac1{mN_m}\sum_{j,i}\operatorname{div}(s_{ji})=R_m+b_m C,
 \qquad b_m=\frac1{mN_m}\sum_j r_j a_j.
\]
The moving remainder $R_m$ has $\rho_m:=\operatorname{mult}_p R_m=(mN_m)^{-1}\sum_j r_jq_j$. On the blowup $\pi:\widetilde S\to S$, write $\widetilde R_m$ for its strict transform. For any $x\in E$, the degree-$j$ leading forms of the residual basis vectors are linearly independent in $H^0(E,\mathcal O_E(q_j))$. The sum of their vanishing orders at $x$ is at most $r_jq_j-r_j(r_j-1)/2$: in a degree-$q_j$ polynomial space the subspace of orders at least $a$ has dimension $\max(0,q_j+1-a)$. It follows uniformly for all $x\in E$ that
\[
 \operatorname{mult}_x(\widetilde R_m|_E)\le s_m,
 \qquad s_m:=\frac{1}{mN_m}\sum_j\left(r_jq_j-\frac{r_j(r_j-1)}2\right).
\]
In any of these normalized sums one can replace $r_j$ by $q_j+1$ with error $O(1/m)$: the individual weights and the derivatives of the quadratic expressions are $O(m)$ and the sum of all dimension deficits is $O(m)$. Consequently
\[
 \rho_m\to\frac{2}{d}\int_0^\tau g(t)^2\,dt=\frac{2\varepsilon}{3},
 \qquad s_m\to\frac{1}{d}\int_0^\tau g(t)^2\,dt=\frac{\varepsilon}{3}.
\]
Likewise $\rho_m+kb_m=(mN_m)^{-1}\sum r_jj$ tends to $\frac2d\int_0^\tau t g(t)\,dt=(\tau+\varepsilon)/3$ by the centroid of the triangle. Hence $b_m\to(\tau-\varepsilon)/(3k)$.

Put $\eta_m=\max(\rho_m/2,s_m)$. For all sufficiently large $m$ this is positive and tends to $\varepsilon/3$. The pair $(E,\eta_m^{-1}\widetilde R_m|_E)$ is log canonical by the pointwise bound. Inversion of adjunction gives log canonicity of $(\widetilde S,E+\eta_m^{-1}\widetilde R_m)$ around $E$. The coefficient on $E$ in the log pullback of $(S,\eta_m^{-1} R_m)$ equals $\rho_m/\eta_m-1\le1$. Thus $(S,\eta_m^{-1} R_m)$ is log canonical at $p$, including the generic centers of all curves through $p$. In particular $v(R_m)\le\eta_m A_S(v)$. By simultaneous basis compatibility,
\[
 S_m(H;v)=v(D_m)\le \eta_m A_S(v)+b_m v(C).
\]
The limit of $S_m(H;v)$ is $S_H(v)$, giving the proposition. For $v=\operatorname{ord}_E$, $v(C)=k$, $A_S(v)=2$ and $3S_H(v)=\tau+\varepsilon$ by the same centroid calculation.
\end{proof}

\paragraph{The discriminant estimate for a plane curve.}

\begin{proposition}
Let $C$ be a reduced complex plane curve germ of multiplicity $k\ge2$, with $b$ local branches, Milnor number $\mu$, and local genus defect $\delta$. Set $I=2\delta+k-b=\mu+k-1$. Then
\[
 \operatorname{lct}_0(C)\ \ge\ \frac1k+\frac{k-1}{I}.
\]
\end{proposition}
\begin{proof}
Choose a general coordinate $x$ as projection parameter, and use Weierstrass preparation to write the local equation as a monic polynomial of degree $k$ in $y$. After substituting $x=z^m$, the Newton--Puiseux theorem gives distinct roots $\phi_1(z),\ldots,\phi_k(z)$, all convergent power series vanishing to order at least $m$. Write $c_{ij}=m^{-1}\operatorname{ord}_z(\phi_i-\phi_j)$ for $i\ne j$ and $c_{ii}=+\infty$. The discriminant has order $\sum_{i\ne j}c_{ij}=I$ in $x$. To see the indicated genus formula, the finite flat $\mathbb C\{x\}$-algebra $A=\mathcal O_C$ has power basis and trace discriminant equal to the polynomial discriminant. The discriminant of its normalization is $k-b$: the normalization is the product of regular branches of ramification degrees $e_\alpha$ with sum $k$, and each smooth branch has discriminant order $e_\alpha-1$ in characteristic zero. Replacing the normalization lattice by its sublattice $A$ adds twice its colength $\delta$ to the discriminant order. Milnor's relation $\mu=2\delta-b+1$ gives the other expression for $I$.

Let $v$ be any divisorial valuation centered at the origin, and extend it to a divisorial valuation $w$ along the finite map $(z,y)\mapsto(z^m,y)$, scaling $w$ so that it restricts to $v$. For clarity, the logarithmic ramification formula reads $A(v)=A(w)+(m-1)w(z)$ with this scaling. It follows either from the finite-cover discrepancy formula in characteristic zero, or by computing canonical ramification at a smooth model extracting $v$ and the normalization of that model in the extension: upstairs the coefficient of the relative ramification divisor over its prime divisor is $e-1$ for the ramification index $e$, so the ordinarily normalized upstairs discrepancy for the log pullback pair is $e A(v)$.

Choose $i$ maximizing $w(y-\phi_i(z))$, and write $h=w(y-\phi_i)/(m w(z))$. Since $z$ and $y-\phi_i(z)$ are regular parameters, the pair with both coordinate divisors of coefficient one is log canonical. Hence $A(v)\ge m w(z)+w(y-\phi_i)=v(x)(1+h)$. For all $j$, the ultrametric property and maximality give
\[
 \frac{w(y-\phi_j)}{m w(z)}=\min(h,c_{ij}).
\]
Indeed unequal values give equality directly, and if $c_{ij}=h$ a higher value for $w(y-\phi_j)$ would contradict the choice of $i$. Let $q=v(C)/v(x)$, and let $N(t)=\#\{j:c_{ij}\ge t\}$. Then $q=\int_0^hN(t)\,dt$. At each $t$ all pairs $j,l$ counted by $N(t)$ have $c_{jl}\ge t$. Integrating their ordered count gives
\[
 I\ge\int_0^h N(t)(N(t)-1)\,dt\ \ge\ (k-1)(2q-kh).
\]
The second inequality is $(k-N)(k-1-N)\ge0$, valid for integer $0\le N\le k$. Define $T=I/[k(k-1)]\ge1$. The displayed estimate and $N\le k$ imply
\[
 \frac qk\le\min\{h,(T+h)/2\}\le\frac{T(1+h)}{1+T}.
\]
The last step uses $h\le T$ in the first branch and $h\ge T$ in the second. Consequently $A(v)/v(C)\ge(1+T)/(kT)$ as required. For a valuation with curve center through the origin, a component of the reduced divisor gives the value one, which is at least the same bound since $T\ge1$ and $k\ge2$.
\end{proof}

\section{Valuations centered at a closed point}\label{foundations:sec:angular-next}

Let \(X=X_d^n\) be smooth of index \(r=n+2-d\ge3\), let \(p\in X\) be arbitrary, and write \(E=\mathbf{P}^{n-1}\) for the exceptional divisor of its blowup. For any finite-discrepancy valuation centered at \(p\), scale so that \(v(\mathfrak m_p)=1\) and put \(q=A_X(v)-n\). Set
\[
 s=S_H(E),\qquad \gamma=\frac{n-1}{n+1},\qquad
 s_*=\frac{n+d-1}{n+1}.
\]
Let \(\pi\) be the blowup. The nef class \(\pi^*H-E\) gives the exact volume \(d-z^n\) for \(0\le z\le1\), and intersecting its \((n-1)\)-st power with \(\pi^*H-zE\) gives pseudoeffective endpoint at most \(d\). The left tangent at one bounds the concave volume root from above on the remaining big interval. Hence \(\operatorname{vol}(\pi^*H-zE)\le(d-z)^n/(d-1)^{n-1}\) there. Integrating gives \(s\le s_*\). Write \(T\in|H|\) for the tangent hyperplane section at \(p\) and \(b=\operatorname{mult}_p T\). Its singular locus is finite by finiteness of the Gauss map. Since \(\dim T\ge2\) it is integral: two components of its hypersurface equation would meet in positive dimension in the dimensions used below.

\begin{lemma}\label{foundations:lemma:point-angular}
For \(2\le b\le d\), put \(e=d/b\) and
\[
 s_b=\frac{n-2+b+e}{n+1}.
\]
Then \(s\le s_b\) and \(S_H(v)\le s+(s-\gamma)q\). If \(b^2>d\), put
\[
 J_b=e-s_*+\frac{(d-1)^2}{d(n+1)}
                \left(\frac{d-e}{d-1}\right)^{n+1}.
\]
One has
\[
 S_H(v)\ \le\ s+(s-\gamma)q-
 \frac{(s-e+J_b)_+}{b-e}
 \big((b-1)q+b-n+r-2\big)_+ ,
 \qquad z_+=\max(0,z).
\]
\end{lemma}
\begin{proof}
First consider a valuation \(w\) on \(\mathbf{P}^N\) scaled to have log discrepancy one. Let \(F(z)\) be the normalized volume of sections of \(\mathcal{O}(1)\) filtered by normalized \(w\)-order at least \(z\). Its \(N\)-th root is concave on its big interval, \(F(0)=1\), and \(F(z)=0\) for \(z>1\) since \(\alpha(\mathbf{P}^N,\mathcal{O}(1))=1\). Also \(\int F\le1/(N+1)\) by K-semistability of \(\mathbf{P}^N\). Subtracting the affine function \(1-z\) from \(F(z)^{1/N}\) shows that \(F(z)-(1-z)^N\) on \([0,1]\) changes sign at most once, from positive to negative. The integral inequality therefore implies \(\int_y^1 F\le \int_y^1(1-z)^N\,dz\) for every \(y\in[0,1]\). It follows that for a subspace of asymptotic dimension proportion \(\theta\), the mean of its normalized \(w\)-weights is at most
\[
 \frac1\theta\int_0^1\min\{\theta,F(z)\}\,dz
 \le \inf_{0\le y\le1}\left(y+\frac1{\theta}\int_y^1(1-z)^N\,dz\right)
 =1-\frac{N}{N+1}\theta^{1/N}.
\]
The first expression bounds any subspace, by containment in the full filtered space.

Take the refinement \(W_{m,j}\subset H^0(E,\mathcal{O}_E(j))\) of sections of \(mH\) by ordinary order \(j\) at \(p\). Write \(g(t)^{n-1}\) for its slice volume at \(t=j/m\). By the slice formula, equivalently by Brunn--Minkowski on its Okounkov body, \(g\) is concave and
\[
 \int g^{n-1}=\frac dn,\qquad
 \frac{\int t g^{n-1}}{\int g^{n-1}}=s.
\]
All jets of degree at most \(m\) are generated by products of ambient linear forms after choosing a nonvanishing frame. Hence \(g(t)=t\) for \(0\le t\le1\). For a normalized \(w\) on \(E\), the doubly filtered slice volume at \((t,u)\) is at most \(\min(g(t)^{n-1},t^{n-1}F_w(u/t))\), by containment in the complete system and the unfiltered slice bound. One can first check the containment on rational rays and pass to the interior by continuity. Integration as in \cite[Corollary~2.22 and Lemma~2.23]{ar2003.13788} and the subspace estimate above give
\[
 S(W;w)\le s-\frac{n-1}{n}\frac{\int g^n}{\int g^{n-1}}\le s-\gamma .
\]
For the last inequality put \(M=\sup g\ge1\), and let \(L(x)\) be the length of the interval where \(g\ge Mx\), \(0<x<1\). This length is concave, so \(L(x)/(1-x)\) is nondecreasing. Tilting the density proportional to \(x^{n-2}(1-x)\) by a nondecreasing function increases its mean. The layer-cake formulas for the two moments of \(g\) give \(\int g^n/\int g^{n-1}\ge nM/(n+1)\).

Take a basis simultaneously adapted to \(v\) and \(\operatorname{ord}_E\). The restriction of the average of its strict transforms to \(E\) is a refined basis divisor. By the uniform estimate in \cite[proof of Lemma~2.21]{ar2003.13788}, its log canonical threshold is asymptotically at least \(1/(s-\gamma)\). Applying inversion of adjunction with coefficient one on \(E\), and separating the ordinary order of the full basis divisor, gives \(S_H(v)\le s+(s-\gamma)q\).

To bound \(s\), adapt a basis simultaneously to the ordinary order and the \(T\)-order. The normalized mean \(T\)-coefficient is \(1/(n+1)\) in the limit. Removing \(i\) copies of \(T\) gives a basis on restriction in the corresponding degree \((m-i)H_T\), by projective normality. The ordinary orders of the lifts after removing \(T^i\) are bounded above by the maximal-ideal orders of these restrictions. On the blowup of \(T\) at \(p\), with exceptional Cartier divisor \(F\), the filtration by \(F\)-order dominates the latter filtration. The mean \(F\)-order is at most \((n-2+d/b)/n\): \(H_T-F\) is nef, the point has multiplicity \(b\), and \(\operatorname{vol}(H_T-zF)=d-b z^{n-1}\) on \(0\le z\le1\). Intersecting with \((H_T-F)^{n-2}\) gives pseudoeffective endpoint at most \(d/b\). Integrating the tangent upper bound for the volume root at one gives the asserted mean just as for the ordinary point on \(X\). When \(b=d\), the interval \([0,1]\) is already the whole support. Weighting the mean on \(T\) by \(n/(n+1)\) and adding \(b/(n+1)\) gives \(s_b\).

Finally \(v(T)\le A_X(v)-(r-2)\). Cut by \(r-2\) general independent hyperplanes through \(p\) to get a smooth \(d\)-fold section \(X'\) through \(p\). Its tangent divisor is a degree-\(d\) hypersurface in \(\mathbf{P}^d\) with isolated singularities, by general slicing of \(T\). For \(c<1\), the non-klt locus of this polynomial pair in \(\mathbf{P}^d\) is finite. Its multiplier ideal twisted by \(\mathcal{O}(-1)\) has vanishing \(H^1\) by Nadel, since \(d-cd>0\). Since \(\mathcal{O}(-1)\) has no sections the multiplier quotient is zero. Thus the polynomial pair is log canonical at coefficient one. The hypersurface germ in \(X'\) is analytically the same as in its tangent \(\mathbf{P}^d\): eliminate the transverse coordinate by the implicit function theorem. Inversion of adjunction along the slicing hyperplanes gives log canonicity on \(X\) with these \(r-2\) hyperplanes and \(T\) present. Each slicing hyperplane has \(v\)-order at least one.

Suppose \(b^2>d\). Any section in ordinary grade \(j\) is divisible by
\[
 T^{a_{m,j}},\qquad
 a_{m,j}=\max(0,\lceil(j-em)/(b-e)\rceil).
\]
Indeed after removing all its tangent components, intersect the residual and \(T\) with general hyperplanes through \(p\) to get the multiplicity bound \(e\) times the residual class. Remove just these forced copies in the basis adapted to \(v,\operatorname{ord}_E\). Write \(\bar a\) for the limiting mean of \(a_{m,j}/m\) under the ordinary graded dimensions; it exists by the slice formula. The residual ordinary mean is \(s-b\bar a\). At each \(t\) the normalized degree of its initial forms on \(E\) decreases by \(b(t-e)_+/(b-e)\), and their dimension stays the same. The negative term in the complete-space bound, namely \(-[(n-1)/n]g(t)\) in normalized order units, is independent of the ambient degree of these forms. Thus the same slice integration gives the bound \(s-b\bar a-\gamma\) for the residual restricted \(S\)-invariant per unit discrepancy. For finite-level uniformity, the residual basis weights on \(E\) equal the original weights minus \(\bar a_m w(T_{\rm in})\), where \(T_{\rm in}\) is the degree-\(b\) tangent form. Here \(\bar a_m\to\bar a\) and \(w(T_{\rm in})\le b A_E(w)\). The preceding uniform comparison for the original series therefore carries over to the residual bases. Inversion of adjunction gives
\[
 S_H(v)\le s+(s-\gamma)q+\bar a[v(T)-b(1+q)].
\]
Combine this with the estimate without subtraction by taking the better one.

Under the ordinary limiting distribution the expectation of \(t\) is \(s\) and the survival function is \(\operatorname{vol}(H-tE)/d\). Consequently
\[
 (b-e)\bar a=\mathbb E[(t-e)_+]=s-e+\mathbb E[(e-t)_+].
\]
The last nonnegative expectation equals \(e-\int_0^e\operatorname{vol}(H-tE)\,dt/d\). The volume bound \(\operatorname{vol}(H-tE)\le(d-t)^n/(d-1)^{n-1}\) on \([1,d]\), together with the exact initial interval, shows it is at least \(J_b\). This proves the lemma. For \(b=2\) in the degrees used below the estimate without subtraction applies. For \(b=d\), \(e=1\) makes the displayed formulas well defined.
\end{proof}

\begin{lemma}\label{foundations:lemma:simplex-tangent}
If \(b^2>d\), the value \(J_b\) in Lemma~\ref{foundations:lemma:point-angular} can be replaced by the maximum of the given value and
\[
 e-s_b+\frac{(b-1)^2}{b(n+1)}
                    \left(\frac{b-e}{b-1}\right)^n .
\]
\end{lemma}
\begin{proof}
We retain the full distribution in the proof of \(s\le s_b\). Take a degree-\(m\) basis simultaneously adapted to ordinary and \(T\)-orders. In the \(T\)-order-\(i\) block write \(U=i/m\). Its limiting law is \(\operatorname{Beta}(1,n)\): the block with this \(i\) has dimension \(h^0(T,(m-i)H)\) with leading coefficient proportional to \((m-i)^{n-1}\). Within the block the ordinary order divided by \(m\) is at most \(bU+(1-U)Z_\mathrm{res}\), where \(Z_\mathrm{res}\) is the normalized order of the restriction. For each cutoff the number of residual basis restrictions of order at least the cutoff is bounded by the dimension of the corresponding \(F\)-filtered space on \(T\), by independence and the restriction inequality. The volume tangent bound in the proof of Lemma~\ref{foundations:lemma:point-angular} shows the residual survival probabilities are bounded by \(1-z^{n-1}/e\) on \([0,1]\), and by \((e-z)^{n-1}/[e(e-1)^{n-2}]\) on \([1,e]\). These are exactly the survival probabilities of the affine function with vertex values \(0,1,\ldots,1,e\) on the uniform \((n-1)\)-simplex: below the second lowest knot and above the second highest knot, the removed and retained regions are simplices with explicitly scaled edges. If \(e=1\), only the initial interval is required. Passing to the weighted sum over \(i\) gives stochastic domination of the ordinary limiting law by the affine function on the uniform \(n\)-simplex with vertex values
\[
 0,\ \underbrace{1,\ldots,1}_{n-2},\ e,\ b .
\]
Indeed conditioning a uniform simplex point on its last coordinate gives the same beta law and a residual independent rescaled uniform simplex. One can justify the passage using finitely many bins in \(i/m\) and the survival cutoff and then shrinking the bin widths; the displayed model survival is continuous. The mean of the dominating law is \(s_b\). For \(z\ge e\), the simplex region above \(z\) has relative volume \((b-z)^n/[b(b-e)(b-1)^{n-2}]\), by linear edge scaling from the unique remaining vertex. Integrating from \(e\) to \(b\) gives \((b-e)^n/[b(n+1)(b-1)^{n-2}]\). Under the dominating law the expectation of \((e-z)_+\) is \(e-s_b\) plus this integral. Since this test function decreases, the actual expectation is at least its model value. This proves the assertion.
\end{proof}

\paragraph{The estimates for an arbitrary prime divisor.}

The estimates extend to any prime divisor \(P\in|hH|\) through \(p\) by taking \(b=\operatorname{mult}_p(P)/h\), \(e=d/b\). The mean \emph{number} of copies of \(P\) in a compatible basis divided by its degree level is \(1/[h(n+1)]\). The class proportion \(U=hi/m\), where \(i\) is its \(P\)-order in degree \(m\), has limiting law \(\operatorname{Beta}(1,n)\) from the dimensions \(h^0(P,(m-hi)H)\). Restriction in this sequence is surjective by the vanishing of intermediate line-bundle cohomology on \(X\). The ordinary mean contribution of these copies is \(b/(n+1)\). The degree and local multiplicity of \(P\) are \(hd\) and \(hb\), respectively, so their ratio in the residual volume comparison is exactly \(e=d/b\). This proves both the arbitrary-prime bound \(s\le s_b\) and stochastic domination by the uniform-simplex law with knots \(0,1,\ldots,1,e,b\). Integrality of \(P\) suffices for the volume comparison on its blowup. In the forced-power calculation for \(b^2>d\), the number of copies in a grade-\(j\) section is at least \(\lceil (j-em)_+/[h(b-e)]\rceil\). The saving thus uses \(v(P)/h\).

One has \(v(P)/h\le A_X(v)-(r-2)\). On a smooth index-two hypersurface \(X'\) of dimension at least three every \(B\in|hH_{X'}|\) satisfies \(\operatorname{lct}(X';B)\ge1/h\). Indeed the codimension-one cycle bound of \cite[Proposition~2.1]{ar1507.00285} implies that the locus \(\operatorname{mult}_x B>h\) is finite. Elsewhere the local log canonical threshold is at least \(1/\operatorname{mult}_x B\). Thus for every \(c<1/h\) the multiplier ideal of \(cB\) can fail to be trivial only at finitely many points. Its twist by \(-H\) has vanishing \(H^1\) by Nadel since \((1-ch)H\) is ample, and \(\mathcal{O}(-H)\) has no sections. A nonzero zero-dimensional multiplier quotient is impossible. Apply this on a general smooth pointed index-two section of \(X\) to \(P|_{X'}\), which is a nonzero effective Cartier divisor, and use inversion of adjunction along the \(r-2\) cutting hyperplanes.

\section{Hyperplane flags through an arbitrary point}
Let $X=X_d^n$ be a smooth complex hypersurface of dimension $n\ge3$ and degree $d\ge4$, with $H=\mathcal O_X(1)$ and $p\in X$. Put $T=\tau_p(H)$ and $B=\max(T,\sqrt{6d})$. For every divisorial valuation centered at $p$, scaled to have $v(\mathfrak m_p)=1$ and with $A_X(v)=n+q$, one has
\[
 (n+1)S_H(v)\le B(q+3)+n-3.\tag{*}
\]

Choose a very general linear flag $X=X_n\supset X_{n-1}\supset\cdots\supset X_3=Y$ through $p$. The members are smooth by Bertini with prescribed transverse tangent spaces. By successive applications of Abban--Zhuang's lifting lemma (Lemma 4.2 in \emph{Seshadri constants and K-stability of Fano manifolds}), $\tau_p(H|_Y)\le\max(T,\sqrt d)$. More explicitly the contrapositive at a strict larger threshold propagates that threshold up all dimensions at least four. As usual the occurrence of a prime with fixed integer degree and multiplicity is tested on countably many jet incidences, so generic violations can be spread along the flags.

On $Y$ take a very general smooth $S\in|2H\otimes\mathfrak m_p|$. One has $\operatorname{Pic}(S)=\mathbb Z H|_S$ at \emph{every} prescribed $p$. Apply Di Gennaro--Franco \cite[Theorem 1.1(b)]{ar1205.4528}, with $Z=\{p\}$ and generation degree $1$ for its ideal. The degree $2$ of the new hypersurface is at least $1+1$. The vanishing cohomology has a nonzero $(2,0)$-part since $H^{2,0}(Y)=0$ and $K_S=(d-3)H|_S$ has a nonzero section. The zero-dimensional base has no curve component. Lefschetz gives $H^2(Y,\mathbb Z)=\mathbb ZH$ and $H^1(S,\mathcal O_S)=0$.

Suppose a very general $S$ carries an integral $C\in|hH_S|$ of multiplicity $k$ at $p$, for fixed integers $h\ge1,k$. For a family dominating the parameter space of $S$, one has either $k/h\le\tau_p(H_Y)$ or
\[
 k(k+1)\le2d h(h+2).
\]
Indeed, if these integral curves do not cover $Y$, their closed union is an integral surface $P$. The linear system of quadrics through $p$ restricts to an irreducible general curve on $P$: in the $2H$ embedding the projection from $p$ on $P$ is generically finite, as the embedded $Y$ contains no lines. Bertini therefore gives $C=P|_S$ on a general $S$. The equality is scheme-theoretic by general reducedness off $p$; $P$ is Cartier on the smooth threefold. Its class is $hH$, and restricting to a general quadric through $p$ preserves its multiplicity there. If the curves cover $Y$, their strict transforms cover the blowup $Y'$ of $p$. Choose a general two-dimensional slice of the parameter space retaining dominance: at a general smooth pair, two general base directions span the normal space modulo the tangent to the fiber. Resolve the resulting generically finite evaluation. On a general normalized fiber $F$ its differential determinant injects $K_{Y'}|_F$ into $K_F$, since the normal bundle of a smooth fiber is trivial. Thus $2g(F)-2\ge2dh(d-5)+2k$. On the blowup $S'$ of the surface the strict transform of $C$ has arithmetic genus given by $2p_a(C')-2=2dh(h+d-3)-k(k-1)$. The inequality $p_a(C')\ge g(F)$ gives the claimed bound.

Counting possible $h,k$ and taking a high-ratio prime component now proves $\tau_p(H_S)\le B$, since $k/h<\sqrt{2d(1+2/h)}\le\sqrt{6d}$ in the covering case. The rank-one surface inequality gives $\delta_p(S,H_S)\ge3/B$. Applying adjunction to $S\in|2H_Y|$ yields $\delta_p(Y,H_Y)\ge4/B$. In fact the asymptotic $S$-coefficient of a $H$-basis divisor is $1/8$, and the residual restricted series is almost complete in the class $3H_S/4$, with no fixed component. To verify this without a normalization convention, the $i$-th restricted space in degree $m$ is the full $H^0(S,(m-2i)H_S)$ by the hypersurface exact sequences. The conditional means of $i/m$ and $(m-2i)/m$ are therefore $1/8$ and $3/4$. The two terms in adjunction are $8$ and at least $4/B$.

Finally retain the linear members as boundaries. On $X_j$, write $\Delta_j=(1-1/B)\sum_{l=4}^j L_l|_{X_j}$, where $L_l$ are ambient hyperplanes, indexed so $L_j$ cuts $X_{j-1}$ on $X_j$. Then $\delta_p(X_j,\Delta_j;H)\ge(j+1)/B$ by induction from $\Delta_3=0$. To see this choose each basis divisor compatibly with $X_{j-1}$. Its asymptotic $X_{j-1}$ coefficient is $1/(j+1)$, its residual series on $X_{j-1}$ has class $jH/(j+1)$ and is almost complete without fixed part. At the scaled coefficient $(j+1)/B$, including $\Delta_j$ gives total coefficient one on $X_{j-1}$, and the remaining restricted inequality is exactly the induction hypothesis. This follows from adjunction for a pair with effective boundary \cite[Theorem 3.3 and Corollary 3.4]{ar2003.13788}, or by taking coefficients just below their limits in inversion of adjunction. All other boundary hyperplane restrictions are simple normal crossings for a general flag. Therefore $(n+1)S_H(v)\le B A_X(v)-(B-1)\sum_{l=4}^n v(L_l)$. Each $v(L_l)\ge1$, proving (*).

\section{An ordinary-order estimate for a subspace of sections}

Write \(N=n-1\) and normalize \(v(\mathfrak m_p)=1\), \(A_X(v)=n+q\). For a prime divisor \(P\in |hH|\), write \(b=\operatorname{mult}_p(P)/h\).

\paragraph{The angular estimate on a subspace.}

Write $N=n-1$. If $M$ distinct degree-$j$ monomials in $N+1$ variables are chosen, the mean exponent of any one coordinate is at most
\[
 j-\frac N{N+1}(N!M)^{1/N}+O_N(1).
\]
Indeed among all the monomials the ones with deficiency at most $l$ in this coordinate number $\binom{l+N}{N}$ as long as $l\le j$. Take the largest $l$ for which this number is at most $M$. Among the monomials of deficiency at most $l$, the average deficiency is $Nl/(N+1)$, by the simplex symmetry or the binomial identity. Each subsequent monomial has higher deficiency. Thus every such set has mean deficiency at least $Nl/(N+1)$, with $l>(N!M)^{1/N}-N-1$.

Let $R_m=H^0(X,mH)$, $D_m=\dim R_m$, and let $U_m\subset R_m$ be \emph{any} subspace. Write $M_j$ for the dimension of its $j$-th graded piece for $\operatorname{ord}_p$. Suppose that every vector in this piece has a representative divisor divisible by $P^{a_j}$ in the following strong sense: every $f\in U_m$ of order $j$ is divisible by $P^{a_j}$. Define
\begin{align*}
 s_{U,m}&=\frac{1}{mD_m}\sum_j jM_j,&
 a_{U,m}&=\frac{1}{mD_m}\sum_j h a_j M_j,\\
 K_{U,m}&=\frac{N}{(N+1)mD_m}\sum_j M_j(N!M_j)^{1/N}.&&
\end{align*}
The largest sum of the $v$-weights of a basis of $U_m$, divided by $mD_m$, is at most
\[
 (1+q)s_{U,m}-qK_{U,m}
 -a_{U,m}\bigl[b(1+q)-v(P)/h\bigr]+O_n(q/m).
 \tag{1}
\]
To see this, fix linear local parameters on $X$ at $p$. Refine the ordinary order flag by the order of least monomials in a fixed lexicographic order among each homogeneous degree. This gives a single finite flag with one-dimensional leaves. Choose a basis of $U_m$ compatible both with this flag and with $v$. Its elements of order $j$ have distinct initial monomials. Divide each such basis vector by the fixed power $P^{a_j}$. The residual homogeneous forms on the exceptional $E=\mathbf P^N$ have degree $j-hba_j$ and still have pairwise distinct initial monomials within each $j$: division subtracts the same exponent vector in that degree. A one-parameter diagonal torus degeneration simultaneously specializes all these forms to their initial monomials, since a finite set of lexicographic comparisons admits an integral system of weights realizing it. The preceding monomial bound controls \emph{every} coordinate coefficient of this specialized residual basis divisor by $s_{U,m}-b a_{U,m}-K_{U,m}+O_n(1/m)$, with normalization $mD_m$. The specialized coordinate divisors have simple normal crossings, so semicontinuity of the log canonical threshold gives the corresponding bound on the original restriction to $E$. This can be done at once on projective $E$ by properness. Inversion of adjunction then bounds the $v$-value of the averaged residual strict transform by $q$ times the same quantity. Its ordinary $E$-coefficient is $s_{U,m}-b a_{U,m}$, and the separated copies of $P$ contribute $a_{U,m}v(P)/h$. This proves (1) without uniformity theorems for a non-graded family $U_m$.

\paragraph{Gluing on the intersection of two pointed quadrics.}
Let \(Y\subset\mathbf P^4\) be a smooth complex hypersurface threefold
of degree \(d\ge4\), fix any \(p\in Y\), and write \(S_t=Y\cap Q_t\)
for a very general quadric \(Q_t\) through \(p\). Write \(H\) for the
hyperplane class. Then
\[
 \tau_p(H|_{S_t})\ \le\
 \max\left\{\tau_p(H|_Y),\sqrt{2d+8\sqrt{2d}}\right\}.
\]
In particular, on a very general pointed threefold section of a smooth
degree-\(d\) hypersurface \(X^n\), one may use the signed quadric flag
with
\[
 B=\max\left\{\tau_p(H_X),\sqrt{2d+8\sqrt{2d}}\right\}.
\]

The pointed Noether--Lefschetz theorem gives \(\operatorname{Pic}(S_t)
=\mathbf ZH\). Fix a dominating family of integral \(C_t\in|hH|_{S_t}|\)
with multiplicity \(k\) at \(p\), with
\(k/h>\max\{\tau_p(H|_Y),\sqrt{2d}\}\). As in the pointed quadric flag,
the family covers \(Y\) and
\[
                 k(k+1)\le2dh(h+2).
\]
The family may be taken over a dominant generically finite cover of the
pointed quadric system: two distinct integral curves with ratios
above \(\sqrt{2d}\) on the same smooth \(S_t\) would have local
intersection exceeding their total intersection. We claim that
\[
                       k\le4h^2\qquad(h<d).
\]
Every section cutting out \(C_t\) lifts uniquely up to scalar to a
section of \(\mathcal O_{Q_t}(h)\). Indeed \(Q_t\) is arithmetically
Cohen--Macaulay, so \(H^1(Q_t,\mathcal O_{Q_t}(h-d))=0\), and
\(h-d<0\) gives uniqueness. Its divisor \(G_t\subset Q_t\) is
integral. On the smooth quadric threefold each prime divisor has class
\(aH\) for \(a>0\); each such component therefore restricts to a
nonzero effective divisor on \(S_t\) (no component can contain
\(S_t\), as \(h<d\)). The primeness of \(C_t\) excludes a multiple
component or more than one component.

For two general independent parameters \(s,t\), the surface
\(V=Q_s\cap Q_t\) is smooth and integral, and
\(A_s=G_s|_V\), \(A_t=G_t|_V\) are integral curves. Indeed, for
fixed general \(G_s\), the quadrics through \(p\) induce on it the
hyperplane system through a point of its quadratic Veronese embedding.
This embedding contains no line, since the degree of every curve in it
is even. The projection from the marked point therefore has
two-dimensional image, so Bertini gives irreducibility and generic
reducedness away from the base point. The intersection is a pure
Cohen--Macaulay curve, hence reduced. The assertion holds for both
curves by symmetry and independence. Similarly \(D=Y\cap V\) is
smooth at \(p\) and integral by successive Bertini and transversality
at \(p\). It cannot be either \(A_t\) or \(A_s\), since their degrees
are \(4h<4d\). On \(S_t\), restriction to the smooth curve \(D\)
of the equation of \(C_t\) has order at least \(k\) at \(p\).
Consequently \(i_p(A_t,D)\ge k\), and the same holds for \(s\). If
\(A_t,A_s\) are distinct, their local intersection algebra surjects
onto its quotient by the equation of the smooth \(D\), of length
\(\min\{i_p(A_t,D),i_p(A_s,D)\}\ge k\). Their total intersection
on \(V\) equals \(4h^2\). Thus \(k>4h^2\) would imply that the two
sections on \(V\) agree up to a nonzero scalar for general pairs.

For completeness, such proportionality glues the general sections.
Fix \(N=\lfloor h/2\rfloor+1\) general quadrics \(Q_i\) and their
sections. On every triple intersection, their restricted sections are
nonzero: for the first fixed pair the nonempty integral curve \(A_i\)
cannot contain the general complete-intersection curve obtained by
cutting \(V\) by the third pointed quadric. Indeed the latter curves
move and cover a dense open of \(V\). Apply this argument also to a
further general parameter. Hence, after scalar normalizations using
the first anchor, the sections agree on every pair of anchors, and a
further general section can be normalized to agree with every anchor.
An ambient form \(F\) of degree \(h\) interpolates the \(N\) anchor
sections. Start by lifting the first one to \(\mathbf P^4\). Once
the first \(j-1\) have been matched, the difference on \(Q_j\)
vanishes along each of the distinct integral Cartier divisors
\(Q_i|_{Q_j}\) for \(i<j\). Divide on the smooth \(Q_j\) by their
product, lift the resulting section of degree \(h-2(j-1)\)
to \(\mathbf P^4\), and add its product with the first \(j-1\)
quadric equations to the existing interpolant. The restriction
maps to \(Q_j\) are surjective by the quadric exact sequence. On
a further general quadric \(Q_t\), the difference of the normalized
section and \(F|_{Q_t}\) vanishes on \(N\) distinct such divisors,
so is divisible by a section of degree \(2N>h\) and must be zero.
This puts the covering family into \(Y\cap(F=0)\), impossible since
\(F\ne0\) and \(h<d\). This proves the claim.

For a ratio larger than \(\sqrt{2d}\), the claim forces
\(h>\sqrt{2d}/4\) when \(h<d\), and the same is automatic for
\(h\ge d\). The differential inequality now gives
\[
    (k/h)^2 < 2d(1+2/h)<2d+8\sqrt{2d}.
\]
Passing to all primes on a very general surface gives the asserted
threshold. The last statement uses pointed hyperplane lifting in
dimensions at least four to bound
\(\tau_p(H|_Y)\le\max\{\tau_p(H_X),\sqrt d\}\).

\paragraph{Angular entropy after separating a prime divisor.}

Suppose \(P\in |hH|\) has multiplicity \(hb\) at \(p\), with \(b>\sqrt d\), and put \(e=d/b\), \(N=n-1\), \(q=A_X(v)-n\).

The threshold at $p$ is precisely $b$. In fact, every other integral
prime of ratio $b'$ has $bb'\le d$, by intersecting the two divisors
and general hyperplanes through $p$. The local multiplicity of a
proper intersection is at least the product of the two multiplicities.
Since $e<b$, taking prime components gives $T=b$.
The ordinary expected order satisfies the previously established bound
$s\le(n-2+b+e)/(n+1)$, and $v(P)/h\le d+q$.

Set $K=\lceil(n+1)b\rceil$. In degree $m$, choose a basis compatible
with the $v$-filtration and with the full filtration by powers of $P$.
A vector of $P$-order $\nu$ has a class fraction $w=h\nu/m$.
Group the basis vectors by the integer $i=\lfloor Kw\rfloor$ and
let $W_i$ be their spans. Since the initial basis was $v$-adapted,
the sum of the maximal $v$-weights of these spaces is the full
maximal sum. In $W_i$ every element is divisible by
$P^{a_i}$ where $a_i=\lceil mi/(Kh)\rceil$. The mean of $ha_i/m$,
weighted by the dimensions of $W_i$, differs from the full mean of
$w$ by a number in $[-1/K,0]$. The full mean tends to $1/(n+1)$,
as follows by summing $h^0(X,(m-hj)H)$ over the positive powers $j$.

Write $M_{ij}$ for the ordinary grade-$j$ rank of $W_i$ and
$p_{ij}=M_{ij}/D_m$, $D_m=h^0(X,mH)$. Apply the finite angular
estimate separately on each $W_i$, with constant forced power $a_i$.
There are at most $K+1$ spaces, fixed independently of $m$.
Within each, the necessary compatible flags of ordinary monomials
and $v$ can be chosen afresh. The union of these bases is a basis of
the full space, so the sum of their ordinary weights is at most the
full ordinary weight sum. The total angular moment subtracted with
coefficient $q$ is
\[
  D^{\rm ang}_m=\frac{N}{n}\frac{(N!D_m)^{1/N}}m
                         \sum_{i,j}p_{ij}^{1+1/N}.
\]
For a discrete probability vector let
$\mathcal H(p)=-\sum p_\alpha\log p_\alpha$, with zero terms omitted.
Jensen gives $\sum p_\alpha^{1+1/N}\ge\exp(-\mathcal H(p)/N)$.
Since there are at most $bm+1$ ordinary grades,
$\mathcal H(i,j)\le\mathcal H(i)+\log(bm+1)$.
The limiting upper bound on the mean of $i$ is
$K/(n+1)\le b+1/(n+1)$. A nonnegative integer random variable of
mean $\mu$ has entropy at most
$(1+\mu)\log(1+\mu)-\mu\log\mu\le 1+\log(1+\mu)$.
This follows from nonnegativity of relative entropy with respect to
the geometric distribution of mean $\mu$ and
$\mu\log(1+1/\mu)\le1$; the zero-mean case is immediate.
Consequently $\limsup\exp(\mathcal H(i))<3(b+2)$.
Using $N!D_m\sim(d/n)m^n$, we get
\[
  \liminf D^{\rm ang}_m\ge D^{\rm ang}
  :=\frac{N}{n}\left(\frac{e}{3n(b+2)}\right)^{1/N}.
\]

Put $G=b(1+q)-d-q$. Whether $G$ has positive or negative sign,
the total forced-class term in the finite angular estimate is at most
$-G/(n+1)+G_+/K+o(1)$ using $v(P)/h\le d+q$.
For positive $G$ use the lower bound for the mean class, and for
negative $G$ its upper bound. The resulting estimate is
\[
 S_H(v)\le\frac{n-2+e+d}{n+1}
       +q\left(\frac{n-1+e}{n+1}-D^{\rm ang}\right)
       +\frac{G_+}{K}.
\tag{1}
\]

\paragraph{A trial density for the angular contribution.}

Put \(K\ge1\), \(0<c<1\), and divide the \(P\)-power filtration into intervals of length \(1/K\) in the normalized divisor class.

The limiting bin dimensions are
\[
   w_i=(1-i/K)^n-(1-(i+1)/K)^n,\quad 0\le i<K.
\]
Their mean forced class is $A=\sum(i/K)w_i$; the negligible final
endpoint of class $1$ contributes zero in the limit.

After removing this fixed binwise power of $P$, the residual ordinary
grade divided by $m$ is $y\ge0$. The limiting graded dimensions
in bin $i$ define a measure $\rho_i(y)dy$ of mass $w_i$ such that
\[
 0\le\rho_i(y)\le(n/d)y^N,\qquad
   \sum_i\int y\rho_i(y)dy\le
   \frac{n-2+b+d/b}{n+1}-bA=:\mu_b .
\]
For rigor take any common subsequential weak limit; support is
bounded, and the jet rank in residual degree $j$ is at most
$\binom{j+N}{N}$ after dividing by the actual forced copies.
The ordinary mean bound follows from the ordinary-adapted bases
of the bins forming a basis of the full series. Jensen on finitely
many intervals and then refinement of intervals bounds the limiting
angular moment below by
\[
  \frac{N}{n}(d/n)^{1/N}\sum_i\int\rho_i^{1+1/N}.
\]

Here is a trial density estimate for this expression. Put
$a_i=(dw_i)^{1/n}$, $a=\sum_iw_i a_i$, with the chosen $0<c<1$. Let
$g_i$ be the tent rising with slope one from zero to $t_i=c a_i$,
then falling linearly to zero at $u_i=c^{-N}a_i$.
For $\rho_i^0=(n/d)g_i^N$ direct integration gives mass $w_i$,
aggregate mean $a(c^{-N}+Nc)/(n+1)$ and aggregate angular
moment $\Gamma=Nac/(n+1)$. The slopes of the descending sides
are common, so write their affine extensions as $L_i(y)=A_i-Cy$
with a common $C>0$. If the aggregate trial mean is at least the
actual mean, the angular moment of the trial is a lower bound.
Indeed, convexity of $x^{1+1/N}$ bounds the difference of the
actual and trial moments below by a positive multiple of
$\sum\int g_i(\rho_i-\rho_i^0)$. This is at least
$\sum\int L_i(\rho_i-\rho_i^0)$: on the rising side $g_i\le L_i$
and $\rho_i\le\rho_i^0$ by the jet cap, while beyond $u_i$ one has
$g_i=0\ge L_i$ and $\rho_i^0=0$.
The final expression is nonnegative by the mass and mean comparisons.

% ===== Positive centers and the range \(n\ge2r-1\) =====
\chapter{\texorpdfstring{Positive centers and the range \(n\ge2r-1\)}{Positive centers and the uniform index range}}\label{ch:global}
\section{The dimension range}
Let $X=X_d^n\subset\mathbf P^{n+1}_{\mathbf C}$ be a smooth hypersurface,
let $H=\mathcal O_X(1)$ and let $r=n+2-d$, so $-K_X=rH$.
We first treat the range
\[
 \boxed{\quad d\ge4,\quad 1\le r,\quad n\ge2r-1
      \quad\Longrightarrow\quad \delta(X,-K_X)>1 .\quad}
\]
This chapter proves the range up to \(n=2d-3\) at fixed degree. Chapters~\ref{ch:boundary} and~\ref{ch:sextic} extend it.

\paragraph{Notation.}
For a divisorial valuation $v$, $A_X(v)$ is its log discrepancy and
$S_H(v)$ its expected order in the complete section ring of $H$. A
$v$-compatible basis in degree $m$ defines a basis divisor normalized by
$m h^0(X,mH)$ of value $S_{m,H}(v)$. These expectations converge to
$S_H(v)$ and satisfy the uniform upper approximation of Blum--Jonsson
\cite{ar1706.04548}. If $\delta(X,H)\le r$, some divisorial valuation computes it
\cite{ar2102.09405}. For a klt boundary $B$ of class $\lambda H$ we abbreviate
$\inf_v(A_X(v)-v(B))/S_H(v)$ by $\delta(X,B;H)$.

\section{Rational normal jets and low-degree curves}

Let $X=X_d^n\subset\mathbf P^{n+1}_{\mathbf C}$ be smooth, $n\ge4$, and $H=\mathcal O_X(1)$. For a reduced integral curve $B\subset X$ put $e=H\cdot B$ and $c=n-1$. Define
\[
 \theta(e)=1-\frac{n-1}{n}
          \left(\frac{d}{n e(d-1)}\right)^{1/(n-1)}.
\]
We prove the localized estimate
\begin{equation}\label{global:eq:jet-curve}
 \delta_B(X,H)\ \ge\ \min\{n-1,\,1/\theta(e)\}.
\end{equation}
Here the center of a test valuation may contain $B$. The bound depends only on the degree of the curve. It also applies to a curve singular at some closed points.

In $R_m=H^0(X,mH)$ take the filtration by ordinary order at the generic point of $B$. Every nonzero section has order $j\le m$ by the divisor case of Suzuki's multiplicity theorem for a smooth hypersurface \cite[Proposition 2.1]{ar1507.00285}. Shrink around the generic point so that $B$ is smooth. Write $K=\mathbf C(B)$. The degree-$j$ initial part of a section is a homogeneous polynomial of degree $j$ on the $c$-dimensional normal space over $K$, with values in the fiber of $mH$ over $K$.

We choose coordinates on that normal space in which each individual monomial coefficient, for sections in a fixed grade $j$, ranges in a complex vector space of dimension at most
\[
                         L_m=e(d-1)m+1.
\]
Write $F$ for the homogeneous equation of $X$ and $F_i$ for its coordinate derivatives. On the affine cone use the homogeneous derivations $F_i\partial_k-F_k\partial_i$ of degree $d-2$. They kill $F$ and span the tangent space to its smooth affine cone wherever its gradient is nonzero. Their images in the projective normal space to $B$ consequently span that space. Select $c$ constant linear combinations whose images give a basis at its generic point. On a chart $x_0\ne0$ meeting $B$ these define regular directions after homogeneous trivialization and passage modulo the Euler direction. For a degree-$m$ section of order at least $j$, the values on $B$ of the $j$-fold iterated derivatives give its degree-$j$ symmetric normal tensor, up to fixed factorials and trivializations. This follows from the Leibniz rule on products of $j$ local generators of the ideal. Derivatives falling on a coefficient, and lower-order corrections arising in dehomogenization, vanish modulo the ideal. Each value is the restriction of a section of degree $m+j(d-2)$. On the normalization of $B$, a line bundle of degree $b\ge0$ has at most $b+1$ independent sections. Thus any selected coefficient has the asserted capacity, since $j\le m$.

Let $v$ be a divisorial valuation whose center contains $B$. Refine the ordinary-order flag by a fixed monomial flag in these coordinates on each normal grade. Choose a complex basis of $R_m$ compatible with the refined flag and with $v$ (simultaneous bases for two filtrations). Its initial monomials in the normal grades have total degrees at most $m$. Each occurs at most $L_m$ times, since the corresponding leading coefficients are complex linearly independent in the indicated coefficient space. Put $D_m=\dim R_m$.

Fix one of the $c$ normal coordinates. Among all monomials of total degree at most $m$, there are $\binom{l+c}{c}$ with exponent at least $m-l$ in this coordinate, for $0\le l\le m$. In the full group of these monomials the mean deficiency from $m$ is $cl/(c+1)$, by the elementary simplex count (use a slack variable to represent unused degree). Fill copies of the smallest deficiencies first, using at most $L_m$ copies of each monomial. The average exponent in the selected coordinate is consequently at most
\[
               m-\frac{c}{c+1}
                      \left(\frac{c!\,D_m}{L_m}\right)^{1/c}
                     +O_c(1).
\]
Indeed take the largest $l$ for which $L_m\binom{l+c}{c}\le D_m$, with a bounded endpoint adjustment. All additional copies have deficiency at least $l$. The bound on $l$ from the next binomial coefficient gives the stated error. Since $D_m=dm^n/n!+O(m^{n-1})$, every normalized coordinate coefficient has limsup at most $\theta(e)$.

To translate this count to discrepancies, blow up the smooth generic locus of $B$, with exceptional divisor $E$. The generic fiber of $E\to B$ is $\mathbf P^{c-1}_K$. Let $D$ be the basis divisor just chosen and $D'$ its strict transform. Its restriction to that fiber is the product divisor of the normal initial forms, with normalization $1/(mD_m)$. One diagonal one-parameter subgroup in the chosen fixed coordinates specializes all these forms to their leading monomials: there are only finitely many monomial comparisons at each $m$. The resulting divisor is supported on coordinate hyperplanes and is simple normal crossings. For each fixed $\theta'>\theta(e)$ its coefficients are at most $\theta'$ for all sufficiently large $m$. Semicontinuity of log canonical thresholds, applied over an algebraic closure of $K$ on the proper projective fiber, gives log canonicity of the original restricted divisor with coefficient $1/\theta'$. This gives the same conclusion along the generic fiber over $\mathbf C$ by spreading a log resolution over an open part of $B$. All monomial counts used the fixed derivative coordinates before this specialization.

Inversion of adjunction gives log canonicity of $E+D'/\theta'$ around the generic fiber. The order $\rho_m$ of $D$ along $B$ is at most one. On pulling back $tD$, the coefficient of $E$ in the log pullback is $t\rho_m-(c-1)\le1$ if $t\le c$. Thus $(X,tD)$ is log canonical at the generic point of $B$ for every fixed $t<\min(c,1/\theta(e))$ and large $m$. This also tests valuations centered on any larger subvariety containing $B$, by properness of the blowup and openness around the generic point of $B$. Since the basis was $v$-compatible, $v(D)=S_m(H;v)\to S_H(v)$. This proves \eqref{global:eq:jet-curve}.

\paragraph{Sharper generic jets for a rational curve.}
Let $X=X_d^n\subset\mathbf P^{n+1}$ be smooth, $n\ge3$, and let
$f:\mathbf P^1\longrightarrow C\subset X$ normalize an integral
rational curve of degree $e$. Put $c=n-1$ and let
$N^*=\ker(f^*\Omega_X\to\Omega_{\mathbf P^1})$. The image of this
map is an invertible ideal times $\Omega_{\mathbf P^1}$, so $N^*$ is
a subbundle of rank $c$, including above singular points of $C$.
Write $N^*\otimes\mathcal O(e)=\bigoplus_{i=1}^c\mathcal O(a_i)$.
Then
\[
                   \sum_i(a_i)_+\le (d-1)e.
\]
Indeed the positive summands form a subbundle $P$ of $f^*\Omega_X(1)$.
Its preimage in $f^*\Omega_{\mathbf P^{n+1}}(1)$ is an extension of
$P$ by $\mathcal O(-(d-1)e)$, by the conormal sequence. The Euler
sequence embeds that preimage into a trivial bundle on $\mathbf P^1$,
so its degree is at most zero by taking its top exterior power.

If $s\in H^0(X,mH)$ vanishes to order at least $j$ at the generic
point of $C$, its ordinary $(j-1)$-jet vanishes along the whole of
$f(\mathbf P^1)$, by density and local freeness of the bundle of
principal parts on smooth $X$. The exact sequence of principal
parts therefore identifies the pullback of its $j$-jet with a global
section of $\operatorname{Sym}^j f^*\Omega_X\otimes\mathcal O(em)$.
At a general smooth point of $C$ its symbol belongs to the subbundle
$\operatorname{Sym}^j N^*\otimes\mathcal O(em)$, by the local
equations of the smooth curve in $X$. This inclusion holds over
all of $\mathbf P^1$, because the quotient is locally free. The
resulting map from the order-$j$ graded piece of $H^0(X,mH)$ into
its generic fiber is injective. In the splitting frame the coefficient
of the monomial $\beta=(\beta_1,\ldots,\beta_c)$, $|\beta|=j$,
thus belongs to a fixed complex vector space of dimension at most
\[
       h^0(\mathbf P^1,\mathcal O(e(m-j)+a\cdot\beta))
       \ \le\ e(m-j)+\sum_i(a_i)_+\beta_i+1 .
\]
Here $j\le m$ by Suzuki's multiplicity bound for each fixed positive-
dimensional center. All monomial slots are taken in the same splitting
frame independently of any valuation.

Fix a divisorial valuation $v$ on $X$ whose center contains $C$.
In $R_m=H^0(X,mH)$ refine the ordinary-order filtration along $C$
by a diagonal weight on each order-$j$ quotient in the splitting
frame. Choose the weight to separate all exponent vectors of the
same degree $j\le m$. Take a complex basis of $R_m$ simultaneously
adapted to this refined filtration and to $v$. The order-$j$ symbols
of its members give a basis compatible with the diagonal weight.
Those with a fixed leading slot $\beta$ have linearly independent
coefficients in the indicated coefficient space. Thus their number
is bounded by its dimension. On the generic exceptional projective
space of the blowup of $C$, the diagonal action degenerates the
product of these symbols to a sum of coordinate hyperplanes, counted
with the normalized total exponents of their coordinates. Lower
semicontinuity of the log canonical threshold gives the same lower
bound for the product of the original symbols. One can check this
over $\overline{\mathbf C(C)}$ and spread the conclusion over a
dense open of $C$, as in the preceding proof.

Let $x_i=\beta_i/m$, $x_0=1-\sum_{i=1}^cx_i$, and take uniform
probability on the $c$-simplex. The limiting number of lattice slots
divided by $m^c$ is $1/c!$. The total number of basis sections divided
by $m^{c+1}$ tends to $d/n!$. The limiting capacity of the tail
$x_1\ge u$, divided by the total dimension, is at most
\[
   \frac{ne}{d}(1-u)^{n-1}
          \left((d-1)u+\frac dn(1-u)\right).
\]
Indeed the tail simplex has fraction $(1-u)^c$, and its barycenter
has coordinate $u+(1-u)/n$ on the first vertex and $(1-u)/n$
on all other vertices including vertex zero. Apply this to the
capacity $e x_0+\sum (a_i)_+x_i$, using its coefficient bounds
$\sum(a_i)_+\le(d-1)e$ and $(a_1)_+\le(d-1)e$. The additive unit
in the finite-level slot bound contributes $O(1/m)$ after summation.
Integrating the capacity above an arbitrary $t\in[0,1]$, and using
the total capacity bound $1$ below $t$, bounds the normalized mean
exponent by
\[
  \theta_{n,d,e}(t)=t+\frac ed(1-t)^n
       \left((d-1)t+\frac{2d-1}{n+1}(1-t)\right).
\]
The same bound holds for each coordinate. Put
$\Theta(e)=\inf_{0\le t\le1}\theta_{n,d,e}(t)$.
For every fixed $c'<\min\{n-1,1/\Theta(e)\}$ the restricted strict
transform of the chosen basis divisor is therefore log canonical
with coefficient $c'$ on the generic exceptional fiber for all
sufficiently large $m$. Its mean ordinary order along $C$ is at most
one. Inversion of adjunction on the blowup, exactly as in the preceding
proof, yields log canonicity of $c'$ times this basis divisor near
the generic point of $C$. This tests $v$, and the chosen basis computes
$S_{m,H}(v)$. Taking limits gives
$\delta_C(X,H)\ge\min\{n-1,1/\Theta(e)\}$.
The optimum is the unique $t\in[0,1)$ at which the displayed tail
capacity equals one, including $t=0$ when $e=1$. For $e=1$, the value
is $(2d-1)/(d(n+1))$, agreeing with the line estimate.

\section{Linear sections containing a rational curve}
Let $C\subset X$ be an integral rational curve of degree $e$ and projective span $\mathbf P^a$, so $a\le e$. For any $j\ge3$ with $a+1\le j\le n$, the following bound is available, and one may always take $j=n$ on the original smooth hypersurface:
\begin{equation}\label{global:eq:slicejet}
 \frac{(n+1)S_H(w)}{A_X(w)}
 \le (j+1)\max\{1/(j-1),\Theta_{j,d}(e,t)\}
 \quad(0\le t\le1,\ C\subset c_X(w)).
\end{equation}
Here $\Theta_{j,d}(e,t)$ denotes the polynomial $\Theta$ in the rational-jet calculation with dimension $j$. To verify the slicing assertion let $Y=X\cap\mathbf P^{j+1}$ be general among linear sections containing $\mathbf P^a$. Off the span Bertini applies. On the span the singularities can occur only at points where the span lies in $T_xX$ and the $j+1-a$ further independent directions belong to $T_xX$. At a fixed such point the latter conditions have codimension $j+1-a$ in the Grassmannian (or are impossible). Thus $\dim\operatorname{Sing}Y\le a-(j+1-a)\le j-3$, and $\operatorname{Sing}Y$ misses the entire curve, as $j+1-a\ge2$. The same holds for intermediate general members. The hypersurfaces are integral and normal by this singular-locus bound and the hypersurface criterion. Suzuki's Proposition~2.1 allows a singular hypersurface in $\mathbf P^{j+1}$ when $2+\dim\operatorname{Sing}Y<j+1$ and the positive-dimensional center misses its singular locus. Thus the order cap in the proof of the jet estimate still holds. The conormal sequence and the jet bundles used there are taken on the open smooth along the entire curve.
 
On each intermediate $k$-fold use a basis adapted to the next hyperplane and the valuation filtration. Its hyperplane coefficient tends to $1/(k+1)$. The graded restriction spaces are the full spaces in the residual degrees, by the hypersurface restriction sequence. Their mean normalized residual degree tends to $k/(k+1)$. Local basis adjunction at the generic point of $C$ therefore multiplies the terminal threshold by $(n+1)/(j+1)$, with truncation at $n+1$; the sections are Cartier and smooth there. This local use of adjunction follows directly by inversion of adjunction and disjointly taking bases in the full residual degrees as in \cite{ar2003.13788}. Equivalently one may use the multigraded basis adjunction theorem: the normality and Cartier conditions hold on each member, and plt is needed around the generic point of $C$. The upper bounds on the right of \eqref{global:eq:slicejet} exceed one, so the truncation is harmless.

We also use the following stronger low-degree estimates:
\begin{equation}\label{global:eq:smalldeg}
 \frac{(n+1)S_H(w)}{A_X(w)}
 \ \le\ \begin{cases}
   (2d-1)/d,& e=1,2,\ \text{or $C$ is a plane cubic},\\
   8/3,& C\ \text{is a twisted cubic},\quad d\ge4 .
 \end{cases}
\end{equation}
The plane-pencil and two-quadric proofs given next work in these degrees for any ambient dimension $n\ge5$. In the twisted cubic calculation $(4d-3)/(3d-1)<4/3$ for all $d>1$. The SNC and complete-intersection restriction arguments are unchanged. Each assertion concerns centers containing $C$.

For a line cut by a general linear flag through
it down to a smooth surface. On the resulting surface
put \(F=H-C\). The planes through \(C\) give a basepoint-free pencil
of residual curves: if the equation is \(xA+yB=0\) with
\(C=(x=y=0)\), then \(A|_C,B|_C\) have no common zero by smoothness.
One has \(F^2=0,\ C^2=2-d\) and \((H-tC)F=(d-1)(1-t)\), so the
pseudo-effective threshold of \(C\) is one. On \(0\le t\le1\) the
class \(H-tC\) is nef. Integrating its square gives
\(S_H(\operatorname{ord}_C)=(2d-1)/(3d)\). Flag adjunction gives
$\delta_C(X,H)\ge(n+1)d/(2d-1)$. Suppose $C$ spans a plane $M$
and has degree $e=2$ or $3$. Take a general linear threefold
$Y=X\cap\mathbf P^4$ through $M$, and a general degree-$e$
surface section $S$ on $Y$ through $C$. Write $(u,v)$ for the
ideal of $M$ in $\mathbf P^4$, $f$ for an integral plane equation
of $C$, $F$ for the equation of $X$, and $G$ for that of the
degree-$e$ section. Then $F|_M\in(f)$, including the case
$F|_M=0$, and $G|_M=\lambda f$ with $\lambda\ne0$ for general
$G$. Consequently $S\cap M$ is exactly the reduced Cartier
divisor $C$ scheme-theoretically, since its ideal on $M$ is
$(F|_M,G|_M)=(f)$. Hence $(u,v)\mathcal O_S=\mathcal O_S(-C)$.
Dividing the two pencil generators by any local Cartier equation for
$C$ leaves two generators with no common zero. This proves both
that the fixed multiplicity is exactly one and that $H_S-C$
is globally generated by a pencil. The surface
is smooth. Indeed outside $M$ the linear flag is generally smooth.
On its base, singularities of $Y$ arise only where the two general
normal derivatives for the additional directions vanish, on the
locus $x\in X\cap M$, $M\subset T_xX$. Their system on that locus
is generated since $X$ is smooth. Consequently such singularities
are finite, disjoint from $C$ (or empty), and $S$ avoids them as
its restriction to $M$ equals $C$. The ideal of the plane curve is
generated in degrees at most $e$; Bertini and general generation
of the rank-two conormal along its smooth part give smoothness
of $S$ there and off $C$. At the finitely many singular points of
a plane cubic there is still a nonzero linear conormal direction
on the smooth $Y$ since the Zariski tangent space of the plane
curve has dimension two. A general section is smooth there too.
Since $(H_S-C)^2=0$, $H_S^2=ed$, $H_S C=e$, integration of
$((1-t)H_S+t(H_S-C))^2$ from zero to one gives
$S_H(C)=(2d-1)/(3d)$, and higher classes are not pseudoeffective
by intersection with $H_S-C$. The point-free generic valuation
on $S$ whose center contains $C$ is just its divisor order.
Linear and hypersurface adjunction along the flag supplies
$\delta_C(X,H)\ge(n+1)d/(2d-1)$.

In the remaining case $C$ is a twisted cubic spanning a three-plane
$M$. Take a general $\mathbf P^5$ through $M$, $Y=X\cap\mathbf
P^5$, and two general quadrics $Q_1,Q_2$ through $C$. The surface
$S=Y\cap Q_1\cap Q_2$ is smooth: $\operatorname{Sing}Y$ is at
most one-dimensional and misses $C$, by the same two normal
derivatives on $X\cap M$. The two general quadrics avoid that
singular locus. Along $C$ generation of its quadratic ideal gives
two general sections of the rank-three conormal in $Y$; their
degeneracy condition has codimension two and misses the curve.
Bertini handles the complement. Intermediate flag varieties are
complete intersections, normal and smooth at the generic point of
$C$; for instance the first quadric meets the singular locus of $Y$
in at most finitely many points. The section and restriction
sequences are the full complete-intersection ones. The two quadrics
on $M$ cut $C\cup L$, where $L$ is a general secant line meeting
$C$ transversely twice. This follows from the net of quadratic
minors defining the twisted cubic, or its secant parametrization:
quadrics through a given secant form a pencil, whose base has the
displayed degree-four complete intersection. As the secants of $C$
fill $M$, the general such $L$ is not contained in $X$ unless
$M\subset X$. In the latter case the intersection on $M$ is the
reduced curve $C\cup L$ scheme-theoretically.

On the smooth surface $S$ the pencil through $M$ has fixed divisor
$C$ if $M\not\subset X$, and $C+L$ otherwise. Indeed its base is
the above intersection with $X|M$, and it is generically reduced on
each surviving curve. Removing the fixed divisor gives a nef class
in either case, since the moving pencil has no fixed curve.
Adjunction gives $H_S^2=4d$, $K_S=(d-2)H_S$, $H_S C=3$ and
$C^2=4-3d$. When $L$ survives, $C\cdot L=2$ by the reduced
transverse secant intersections. Intersecting $H_S-tC$ with the
nef residual pencil yields respectively
\[
 \tau_H(C)\le(4d-3)/(3d-1)\quad\hbox{or}\quad\tau_H(C)\le4/3 .
\]
The complete-series inequality $S_H(C)\le2\tau_H(C)/3$ on
$S$ and flag adjunction give
$\delta_C(X,H)\ge(n+1)/(3 S_H(C))\ge3(n+1)/8$.
For clarity these flag comparisons apply to every valuation whose
center contains $C$. Use Abban--Zhuang, \emph{K-stability of Fano
varieties via admissible flags}, Theorem~3.5, or successive inversion
of adjunction at the generic point of $C$ for a simultaneous
valuation- and flag-adapted basis. All flag members are smooth in
a neighborhood of the generic point of $C$ by the checks above,
and each is a prime Cartier divisor in its predecessor. Their full
section rings are the homogeneous complete-intersection rings.
Restriction onto every residual graded piece $H^0(V',(m-jk)H)$
is surjective by the arithmetically Cohen--Macaulay restriction
sequence. In particular the flag refinements have full graded
pieces. On a $k$-fold a full series in class $mH$ has mean
order along a degree-$j$ Cartier member asymptotic to
$m/[j(k+1)]$, by summing $(m-ja)^k$ for the spaces vanishing
to order at least $a$. Its mean residual degree is therefore
$km/(k+1)+o(m)$. Iteration gives asymptotic mean degree
$3m/(n+1)$ on the terminal surface, and normalized flag-member
coefficients $1/[j(n+1)]$. At the generic point of $C$ the only
nontrivial valuation on the surface is a multiple of
$\operatorname{ord}_C$. Applying its ordinary full-series mean
order in each final graded piece gives normalized coefficient
$3S_H(C)/(n+1)$; bounded-degree pieces have $o(1)$ normalized
contribution. Successive local inversion of adjunction gives the
minimum of the reciprocal coefficients, as asserted. For this
local argument smoothness of the flag near the generic point of
$C$ suffices.

\section{Positive log canonical centers and isolated points}

\paragraph{Bounded degree curves in positive-dimensional non-klt centers.}
Let $X$ be smooth projective of dimension $n$, let $H$ be very ample,
and suppose $-K_X\equiv rH$. If $0\leq B\equiv bH$ is an effective
$\mathbf Q$-divisor, let $t$ be its smallest log canonical threshold
over centers of dimension at least one. This minimum is read off a log
resolution: discard the finitely many exceptional or strict-transform
divisors whose centers are points and take the minimum of $A/v(B)$
over the remaining divisors with $v(B)>0$. Away from the discarded
points, the usual simple-normal-crossings computation of the threshold
proves the assertion. Suppose $t<\infty$ and $\eta=r-tb>0$. Every
positive-dimensional center $Z$ attaining $t$ contains an integral
rational curve $C$ of degree at most $2\dim Z/\eta$.

Indeed $(X,tB)$ is log canonical off a finite set and is log canonical
and non-klt generically along $Z$. In Fujino, \emph{Cone theorem and
Mori hyperbolicity}, J. Differential Geom. \textbf{129} (2025),
617--693, Example~4.11 identifies the lc centers of the effective
normal pair with the qlc centers of $[X,\omega=K_X+tB]$, and the
non-qlc scheme with the non-lc scheme. Its set-theoretic support $N$
is finite. Theorem~4.6(i) equips $X'=Z\cup N$ with the induced
quasi-log scheme structure having the same non-qlc scheme and
quasi-log canonical class $\omega|_{X'}$. As $-\omega\equiv\eta H$
and $N$ is finite, both $\overline{\mathrm{NE}}(X')_{\omega\ge0}$
and $\overline{\mathrm{NE}}(X')_{-\infty}$ are zero. The Mori cone of
$X'$ is nonzero since $\dim Z>0$. Parts (i) and (iii) of Fujino's
Theorem~1.6 therefore supply a rational curve $C$ in $X'$ with
$0< -\omega\cdot C\le 2\dim Z$. Its support belongs to $Z$, proving
the claim. This usage permits worse isolated singularities of $(X,tB)$.

For application to a residual basis divisor $B_m=D_m-a_m H_0$ of
the complete $mH$-series with $a_m\to 1/(n+1)$, a fixed target
$T$ with $Tn/(n+1)<r$ yields a uniform degree bound whenever
$t_m\le T$. If for every curve $C$ of degree below that bound the
asymptotic threshold $\delta_C(X,H)$ exceeds $T$ by a fixed positive
margin, such levels cannot persist. Indeed the divisor $E$ attaining
$t_m$ has center $Z\supset C$, and for every complete basis
the uniform valuation comparison gives
\[
 t_m^{-1}=\frac{v_E(B_m)}{A_X(E)}\le
 \frac{S_{m,H}(E)}{A_X(E)}\le
 (1+o(1))\frac{S_H(E)}{A_X(E)}\le
 \frac{1+o(1)}{\delta_C(X,H)}.
\]
The comparison is uniform over all valuations for the fixed complete
series. In particular the curve, center and residual basis may all
depend on the level.

\paragraph{A shorter rational curve in the point-centered step.}
Suppose one has already proved the uniform lower bound
$\delta_Z(X,H)\ge r+\xi$ for all positive-dimensional $Z\subset X$,
where $\xi>0$. Fix a valuation $v$ centered at a closed point $p$
with $rS_H(v)\ge A_X(v)$. For each sufficiently divisible level take
a basis $D_m$ simultaneously adapted to $v$ and to a general fixed
hyperplane avoiding $p$, and remove its hyperplane coefficient to get
$\Gamma_m$. Thus $v(\Gamma_m)=S_{m,H}(v)\to S_H(v)$. Let $t_m$ be
the first positive-center threshold of $\Gamma_m$, as in the previous
paragraph. The complete-series uniform estimate, $\Gamma_m\le D_m$,
and the stated positive-center gap give $t_m\ge r+\xi/2$ for large $m$.
Consequently $(X,t_m\Gamma_m)$ is not log canonical at $p$. Suppose
$-K_X-t_m\Gamma_m\equiv\eta_m H$, $\eta_m>0$. Its non-lc scheme $N$
is finite, nonempty, and contains $p$.

Choose a minimal positive-dimensional lc center $Z$ under inclusion.
Quasi-log adjunction endows $X'=Z\cup N$ with the structure described
above. Fujino's Theorem~4.6(ii), with the trivial Cartier divisor on
$X$, gives $H^1(X,\mathcal I_{X'})=0$ since $-\omega$ is ample.
Therefore constants surject onto $H^0(X',\mathcal O_{X'})$, so $X'$
is connected. In particular the whole finite support of $N$ lies in
$Z$, since any external point would constitute a separate component.
By the minimality of $Z$ and Theorem~4.6(i), the non-qklt locus of
$X'$ consists of $N$ and the proper zero-dimensional qlc strata, so
it is finite, nonempty, and meets a complete curve of $Z$. Apply
O.~Fujino and K.~Hashizume, \emph{Existence of log canonical
modifications and its applications}, Eur. J. Math. \textbf{9}
(2023), Article 13, Theorem~1.8, to $X'\to\operatorname{Spec}\mathbf C$.
It furnishes an integral rational curve $C\subset Z$ with
\[
                       0<\eta_m (H\cdot C)\le1.
\]
Thus the point-centered contradiction, after the initial positive-center
gap has been established, requires checking curve jets only up to
$e\le1/\eta_m$. The positive-center gap itself can be established at
a smaller multiplier with larger ample slack, using the unrestricted
bound $e\le2(n-1)/\eta_m$ from the ordinary quasi-log cone theorem.

The short-curve proof applies with a fixed klt boundary as follows. At the first positive threshold for a residual basis, add an arbitrarily small general pointed hyperplane to make the original point optimizer non-lc if necessary. The uniform estimate $S_{m,H}\le(1+o(1))S_H$ ensures that the difference between its discrepancy and the value of the current boundary tends to zero from above, or is already negative. Such a hyperplane is free off the point and retains the previous positive lc centers. The ample deficit decreases by $o(1)$. Thus strict bounds on the degree of a short curve persist. If $\Xi$ is the final rational boundary with ample $-(K_X+\Xi)$ and klt away from a finite set, Nadel vanishing with the trivial line bundle gives $\operatorname{length}(\mathcal O_X/\mathcal J(\Xi))\le1$. In particular $\mathfrak m_p\subset\mathcal J(\Xi)_p$. For every divisorial $w$ centered at $p$ the multiplier-ideal criterion yields $w(\Xi)<A_X(w)+w(\mathfrak m_p)$. We use the weak inequality after taking limits of basis levels. A real threshold that occurs in an auxiliary inequality may be approached by rational coefficients within its strict margins.

\section{The dimension of a minimizing center}
Let $X$ be smooth over $\mathbf C$, let $-K_X=rH$ with $H$ very ample, and let $B\sim_{\mathbf Q}\lambda H$ be an effective klt boundary. Suppose a prime divisor with valuation $w=\operatorname{ord}_E$ and center $W$ of dimension $s$ computes $\beta'=\delta(X,B;H)$. Then
\begin{equation}\label{global:eq:dim-center}
      r-\lambda-\frac{n}{n+1}\beta'\le s+1 .
\end{equation}
In particular the strict reverse inequality for $s=0$ excludes a point center without estimating its discrepancy.

Choose a hyperplane $H_0$ avoiding the generic point of $W$. Use a basis adapted to both $w$ and $H_0$ in degree $m$, and subtract its coefficient $a_m H_0$ from the normalized basis divisor, leaving $\Gamma_m$. Then $a_m\to1/(n+1)$ and $w(\Gamma_m)=S_{m,H}(w)\to S_H(w)$. For every fixed rational $t<\beta'$ close to $\beta'$, the pair $(X,B+t\Gamma_m)$ is klt for sufficiently large $m$, since $u(\Gamma_m)\le S_{m,H}(u)\le(1+o(1))S_H(u)$ uniformly in $u$. Suppose the inequality in \eqref{global:eq:dim-center} fails. Fix a small rational $\nu>0$ preserving a strict margin there, then choose $t$ so close to $\beta'$ that $(\beta'-t)S_H(w)<\nu/2$.

For each such $m$ take a general linear section $V$ of codimension $s$. The variety $V$ is smooth, meets $W$ transversely in points of its general smooth locus, and $(V,(B+t\Gamma_m)|_V)$ is klt by successive Bertini on a log resolution. On a resolution carrying $E$, the inverse images of these general hyperplanes have no exceptional fixed component. Over a general mark $p\in V\cap W$, a prime component $E_p$ of the intersection of $E$ with the pulled-back hyperplanes is a divisor over $V$, centered at $p$. Take a smooth resolution including the supports of $B,\Gamma_m$. Generic smoothness over $W$ and adjunction on the sliced resolution give $A_V(E_p)=A_X(E)$; no other exceptional divisor contains its generic point. The value there of the restricted boundary is at least $w(B)+t w(\Gamma_m)$. For $s=0$ simply take $V=X$ and $E_p=E$. A general hyperplane through $p$ on $V$, with coefficient $\nu<1$, keeps this pair klt off $p$; at $p$ it increases its $E_p$-value by at least $\nu$, making the pair non-lc for large $m$.

Write $\Xi_m$ for this boundary on $V$. Since $-K_V=(r-s)H_V$, the divisor $-H_V-(K_V+\Xi_m)$ is ample by our choices and $a_m\to1/(n+1)$. Nadel vanishing gives $H^1(V,\mathcal J(V,\Xi_m)\otimes\mathcal O_V(-H_V))=0$. Restrictions of sections of $-H_V$ must therefore surject onto the sections of its restriction to the nonempty zero-dimensional multiplier subscheme supported at $p$. The former space is zero, and the latter has positive dimension. This proves \eqref{global:eq:dim-center}.

\section{The rational image of an optimizing center}

\paragraph{Short covering curves on an optimal center.}
Let \(X\) be a smooth complex Fano hypersurface of dimension \(n\),
\(-K_X=rH\), and let \(\Delta\ge0\) be rational with
\(\Delta\sim_{\mathbf Q}\lambda H\) and \((X,\Delta)\) log Fano.
Put
\[
 \beta=\inf_w(A_X(w)-w(\Delta))/S_H(w),\qquad
 0<\beta\le r-\lambda.
\]
Suppose a divisorial minimizer has center \(Z\) of dimension \(s>0\).
There is an integral rational curve of \(H\)-degree at most
\((s+1)n/(s\beta)\) through a general point of \(Z\).

\emph{A quantitative tangent slope.}
Let \(Y\) be an arbitrary klt \(\mathbf Q\)-Fano \(n\)-fold, let
\(L=-K_Y/r\), and put \(\delta_Y=\delta(Y,-K_Y)\).
For \(0<t<\min(1,\delta_Y)\) the twisted equation
\[
       \operatorname{Ric}(\omega)=t\omega+(1-t)\omega_0,
       \qquad [\omega]=[\omega_0]=c_1(-K_Y)
\]
has a solution with bounded potentials, smooth on \(Y_{\rm reg}\),
for a smooth reference K\"ahler form \(\omega_0\). Here is one way
to deduce this fixed-parameter assertion using log Fano pairs.
It suffices to use rational \(t\). For sufficiently divisible
\(m\gg1\) take \(D\in|-mK_Y|\) general. The pair \((Y,D)\) is
lc by Bertini on a log resolution: the general divisor pulls
back to a smooth divisor transverse to the exceptional loci, with
no exceptional component. In particular \(w(D)\le A_Y(w)\) for every
divisorial \(w\). Thus
\[
 \delta\bigl(Y,(1-t)D/m\bigr)
     \ge \frac{1-(1-t)/m}{t}\,\delta_Y>1 .
\]
The Ding functional of this uniformly K-stable log Fano pair is
coercive by Li--Tian--Wang, \emph{The uniform version of
Yau--Tian--Donaldson conjecture for singular Fano varieties}, Theorems
1.1 and 3.12. After writing potentials in \(c_1(-K_Y)\)
and dividing the Ding functional in \(t\,c_1(-K_Y)\) by \(t\),
its expression, up to an additive constant, is
\[
 -E(\varphi)-\frac1t\log\int_Y e^{-t\varphi}
      \|s_D\|^{-2(1-t)/m}\,d\mu_0 .
\]
Here \(\mu_0\) is the adapted canonical measure of curvature
\(\omega_0\), and \(s_D\) has the metric on \(-mK_Y\) of curvature
\(m\omega_0\). Constant rescaling lets us assume \(\|s_D\|\le1\).
Consequently the smooth-twist functional, obtained by replacing its
factor \(\|s_D\|^{-2(1-t)/m}\) by \(1\), is also coercive.

For clarity, apply the same coercive bound at any rational
\(t<t'<\min(1,\delta_Y)\). On a normalized minimizing sequence
\(\sup\varphi_j=0\) for the smooth-twist \(t\)-functional,
coercivity bounds \(-E(\varphi_j)\). A weakly convergent subsequence
has a finite-energy limit by upper semicontinuity of \(E\).
The lower bound for the \(t'\)-functional yields
\(\sup_j\int e^{-t'\varphi_j}d\mu_0<\infty\), so the \(t\)-integrals
converge by uniform integrability after taking a subsequence.
Here the klt adapted measure pulls back to a measure absolutely
continuous with respect to smooth volume, and normalized psh weak
convergence gives a subsequence converging almost everywhere.
Thus the lower semicontinuous \(t\)-functional attains its minimum.
The variational Euler--Lagrange argument gives
\[
       \omega_\varphi^n/
       {\textstyle\int_Y\omega_0^n}
            = e^{-t\varphi}\mu_0/
                 {\textstyle\int_Y e^{-t\varphi}d\mu_0}.
\]
The finite-energy potential has zero Lelong numbers on a resolution.
The klt canonical density has some \(L^{1+\varepsilon}\) bound there.
Hölder's inequality and Skoda integrability give an \(L^{1+\varepsilon'}\)
bound for the right hand side for some \(\varepsilon'>0\), so the
potential is bounded by the klt pluripotential estimate; regularity
on the smooth locus follows from the local regularity theorem for
the smooth twist. These are the variational and regularity arguments
of Berman--Boucksom--Eyssidieux--Guedj--Zeriahi,
\emph{K\"ahler--Einstein metrics and the K\"ahler--Ricci flow on log Fano
varieties}, Proposition 1.4, Lemmas 1.5 and 3.6, Theorem 4.8 and
Appendix B. The factor \(t\) in the
exponent simply cancels the \(1/t\) when differentiating the logarithm.

Now fix a smooth resolution \(\pi:\widetilde Y\to Y\). For any rank
\(k\) subsheaf \(\mathcal F\subset T_{\widetilde Y}\), the
fixed-parameter inequality (18) in Druel--Guenancia--Păun,
\emph{A decomposition theorem for \(\mathbf Q\)-Fano K\"ahler--Einstein
varieties}, with \(\gamma=1-t\), gives
\[
 \mu(\mathcal F)\le
       \bigl[1+(1-t)(n/k-1)\bigr]\mu(T_{\widetilde Y}).
\]
The proof through (18) uses only the one given bounded twisted
solution; the later step taking \(\gamma\to0\) is not needed.
The polarization for slopes is \(\pi^*c_1(-K_Y)\); the exceptional
discrepancy terms have zero degree. Subtracting the subsheaf degree
from that of the tangent sheaf and dividing by the quotient rank
shows that every rank-\(s\) torsion-free quotient has slope at least
\(t\mu(T_{\widetilde Y})\). In terms of \(P=\pi^*L\), for every dominant
morphism \(h:\widetilde Y\to W\) to a smooth projective \(s\)-fold
this yields
\begin{equation}\label{global:eq:optimizer-covering-slope}
       -P^{n-1}h^*K_W\ \ge\ \frac{str}{n}P^n .
\end{equation}
Indeed the image of the differential is a rank-\(s\) quotient of the
tangent sheaf and is contained in \(h^*T_W\), so taking its determinant
can only decrease its degree. For the semistable endpoint \(t=1\)
one may also apply Theorem 6(i) and Remark 4 of the same paper.

\emph{The projective image and its curves.}
By Blum--Liu--Zhou, \emph{Optimal destabilization of K-unstable Fano
varieties}, Theorems 4.5 and 4.7, the chosen minimizing divisor
induces a special degeneration \((Y,\Delta_Y)\) of \((X,\Delta)\)
preserving its anticanonically normalized threshold. Put
\(\mathcal L=-(K_{\mathcal X/\mathbf A^1}+\Delta_{\mathcal X})/(r-\lambda)\)
on the total space. The closure \(\Delta_{\mathcal X}\) satisfies
\(\Delta_{\mathcal X}\sim_{\mathbf Q}\lambda\mathcal L\), since this
holds on the punctured product and the only vertical prime is the
principal central fiber. Hence \(\Delta_{\mathcal X}\) is
\(\mathbf Q\)-Cartier, its restriction is \(\Delta_Y\), and
\(-K_Y\sim_{\mathbf Q}rL\) for \(L=\mathcal L|_Y\). Thus \(Y\)
is klt \(\mathbf Q\)-Fano, and \(r\delta(Y,-K_Y)\ge\beta\)
since dropping \(\Delta_Y\) can only increase discrepancies.

For sufficiently divisible \(a\), the section ring of \(aL\) on \(Y\)
is the associated graded of the complete series of \(aH\) for the
valuation, with degreewise shifts of the linearization allowed.
Choose the linearization for which the weights are orders of
divisors of sections and hence nonnegative. The weight-zero
subalgebra in each degree is exactly the image of the restriction
map to \(Z\): a section has strictly positive valuation precisely
when its divisor contains the generic point of \(Z\). Its Proj
is \(Z\); for this identification, the homogeneous ideal is exactly
the vanishing ideal of \(Z\), and all sufficiently high degrees of
the restrictions are complete by Serre vanishing. This subalgebra
naturally includes in the associated graded as its weight-zero sum.
It gives a dominant rational map \(Y\dashrightarrow Z\), with
the very ample coordinates in \(|aH_Z|\) coming from a subsystem
of \(|aL|\) after making \(a\) large. Lift to a resolution \(W\to Z\)
and a common smooth resolution \(\pi:\widetilde Y\to Y\),
\(h:\widetilde Y\to W\). Write \(H_W\) for the pullback of \(H_Z\).
Resolving the subsystem gives \(P-h^*H_W\sim_{\mathbf Q}E\ge0\).

Choose very general complete curves cut by sufficiently divisible
very ample multiples of \(L\) on \(Y\). They avoid the fixed
codimension-two exceptional images of \(\pi\). Their strict
transforms have class proportional to \(P^{n-1}\), and a general
one maps to a nonconstant curve in \(W\) by
\eqref{global:eq:optimizer-covering-slope}. The family covers \(W\).
Apply Jovinelly--Lehmann--Riedl, \emph{Optimal bounds in
Bend-and-Break}, arXiv:2509.08065v2, Theorem 1.1, to the irreducible
image curve in the smooth variety \(W\), using the nef divisor \(H_W\).
The ratio of the two intersection numbers is unchanged on replacing
the image cycle by its reduced curve. The rational curves provided by
the theorem through general points are not contracted to \(Z\), since
these points are in the isomorphism locus of \(W\to Z\). From
\(P^{n-1}h^*H_W\le P^n\) and \eqref{global:eq:optimizer-covering-slope}
we find
\[
 H\cdot R\le \frac{(s+1)n}{s\,tr}.
\]
Taking \(t\nearrow\beta/r\) and using integrality gives
\(H\cdot R\le(s+1)n/(s\beta)\). The families of curves of
each bounded integer degree are parameterized by finitely many
components, so the limiting bound holds through a general point.

\paragraph{The canonical extension and eliminating the dimension factor.}
Use the notation of the short-covering-curve argument: the minimizing
divisor has center \(Z\) of dimension \(s>0\), threshold \(\beta\)
in \(H\)-units, and central klt \(\mathbf Q\)-Fano \(Y\) with
\(-K_Y\sim_{\mathbf Q}rL\) and \(r\delta(Y,-K_Y)\ge\beta\).
If \((s+1)\beta\ge n+1\), then \(Z\) is covered by rational curves
of \(H\)-degree at most \((n+1)/\beta\).

We give the fixed-parameter canonical extension estimate behind this
improvement. On the common smooth resolution \(\pi:\widetilde Y\to Y\)
put \(P=\pi^*L\) and let \(\mathcal V\) be the extension of
\(T_{\widetilde Y}\) by \(\mathcal O_{\widetilde Y}\) of class
\(\pi^*c_1(-K_Y)\). For each rank \(l\) torsion-free quotient
\(\mathcal Q\) of \(\mathcal V\),
\[
       c_1(\mathcal Q)\cdot P^{n-1}\ge
       \frac{ltr}{n+1} P^n,\qquad 0<t<\min(1,\delta(Y,-K_Y)).
\]
Indeed use the bounded fixed-\(t\) twisted K\"ahler--Einstein metric
from the preceding argument. In Druel--Guenancia--Păun, proof of
Theorem 11, Section 3.3, put \(\gamma=1-t\). Here is the mean-curvature
calculation, with the convention from their (11) that wedging with
\(\omega^{n-1}/\omega^n\) contracts a form by \(\Lambda_\omega/n\).
In a unitary frame, normalize the extension form to
\(\sqrt c\sum_j d\bar z_j\otimes dz_j\), where rescaling the vertical
metric permits any \(c>0\). The parallelism of the K\"ahler form kills
the off-diagonal curvature blocks, as in their (31). The quadratic
corrections in the two diagonal blocks after contraction are \(+c\)
on the vertical line and \(-c/n\) on each horizontal coordinate:
the former sums all \(n\) squares, and the latter retains one square.
Their total trace vanishes. Since the horizontal Ricci scalar in (11)
coming from \(t\omega\) is \(t/n\), taking \(c=t/(n+1)\) gives
\(c=(t-c)/n=t/(n+1)\) in both blocks. This computation on the
regularized metrics uses their (12) to split off the same scalar.
The remaining smooth-twist contribution is a
nonnegative endomorphism on the horizontal block, and its integrated
trace with their slope normalization is \(\gamma\) times the total
volume. Projection onto any subsheaf can only reduce this trace.
All other integrated error terms are controlled at this fixed twist
by the estimates of their Section 2.2 for Terms (I) and (III).
Those estimates bound contractions with any orthogonal projection:
the compression of an extension-bundle projection to the horizontal
block is a positive contraction and satisfies the same trace bounds.
The second fundamental form contributes nonpositively. In particular,
if \(\alpha=\pi^*c_1(-K_Y)\) and \(k=n+1-l\), the degree of any
rank-\(k\) subsheaf is at most
\((kt/(n+1)+1-t)\alpha^n\) in this polarization. The auxiliary
exceptional twists used in (25)--(36) have zero intersection against
\(\alpha^{n-1}\). Subtracting from
\(c_1(\mathcal V)\alpha^{n-1}=\alpha^n\) and changing polarization to
\(P=\alpha/r\) yields the quotient bound.

Let \(h:\widetilde Y\to W\) resolve the weight-zero projection to a
smooth model of \(Z\), with \(H_W\) the pullback of \(H|_Z\).
For sufficiently divisible \(a\) take the subsystem of \(|aL|\)
generated by products of \(a\) weight-zero degree-one hyperplane
coordinates. The zero-weight part of each product is its restriction
to \(Z\), so this system maps to the \(a\)-Veronese of \(Z\) by the
original projective embedding. Write the resolved fixed part as
\(aD\), so \(D\ge0\), \(P\sim_{\mathbf Q}h^*H_W+D\). For all
nonexceptional primes, \(D\) has integral coefficients. In detail,
for each projective coordinate \(s\) not identically zero on \(Z\),
the invariant section extending \(s^a\) of \(a\mathcal L\) on the
special degeneration, normalized by its valuation, restricts to a
nonzero central section. Its zero divisor is exactly \(a\) times
the horizontal closure of \(\operatorname{div}(s)\): there is no vertical
component. At the generic point of any prime divisor on \(Y\), the
two-dimensional local ring of the total space is regular, since its
quotient by the parameter is a regular discrete valuation ring. Thus
the horizontal closure is Cartier in codimension one along \(Y\),
and every central order of \(s^a\) there is divisible by \(a\).
At any such prime the order of a mixed product is the arithmetic mean
of the orders of the corresponding pure \(a\)-th powers, since its
\(a\)-th power is their product. The minimum order in the subsystem
is therefore attained on a pure power and is divisible by \(a\).
It follows that
\[
                 P^{n-1}\cdot D_{\rm red}\le P^{n-1}\cdot D
\]
also on the resolution, because its exceptional divisors have zero
degree in this intersection.

For a line bundle \(M\) let \(\operatorname{At}(M)\) be its Atiyah
extension of the tangent bundle by \(\mathcal O\). Its extension
class is \(c_1(M)\) up to the common sign convention. On
\(\widetilde Y\) the extension for \(aP\) is isomorphic to \(\mathcal V\)
by rescaling the vertical coordinate. Away from \(D\), its rational
morphism to \(h^*\operatorname{At}(aH_W)\) is induced by the
line-bundle isomorphism \(aP=h^*(aH_W)\) and the differential of
\(h\). Along the fixed divisor, in compatible local trivializations,
the horizontal block is the regular \(dh\), and only the vertical
output acquires the term \(d\log f\) contracted with the input
vector, where \(f\) is a local defining rational function for
\(aD\). Thus this generically surjective map has at most a simple
pole in a single row along each fixed prime. The image, regarded
as a torsion-free quotient of the source, consequently satisfies
\[
  c_1(\operatorname{im})\cdot P^{n-1}
       \le(-h^*K_W+D_{\rm red})\cdot P^{n-1}.
\]
This assertion can equivalently be read directly from the maximal
minors at the codimension-one discrete valuation rings.

Let \(t\nearrow\beta/r\) in the quotient inequality with \(l=s+1\)
and put \(b=(s+1)\beta/(n+1)\). Since \(b\ge1\), we obtain
\[
 -h^*K_W\cdot P^{n-1}
   \ge bP^n-D\cdot P^{n-1}
   \ge b\,h^*H_W\cdot P^{n-1}.
\]
Jovinelly--Lehmann--Riedl, \emph{Optimal bounds in Bend-and-Break},
Theorem 1.1, applied as before on \(W\) gives covering curves of
degree at most \((s+1)/b=(n+1)/\beta\).

\section{Uniform estimates for point valuations}

The principal-factor and ordinary angular estimates are developed in Chapter~\ref{ch:foundations}. The following estimates explain the constants in the numerical tests.

\paragraph{A dimension-independent point estimate.}
Let $X=X_d^n$ be a smooth hypersurface, $4\le d\le n$, $n\ge10$,
and $r=n+2-d$. Normalize the divisorial valuation at an arbitrary
closed point by $v(\mathfrak m_p)=1$, and write $A_X(v)=n+q$.
Set
\[
 \begin{gathered}
 B=\sqrt{2d}+4,\qquad
 J=\max\left\{\frac{3B+n-3}{n+1},
                  \frac{n-2+d+\sqrt{d/2}}{n+1}\right\},\\
 G=\max\{B,\sqrt{d/2}+\log(3n(d+2))\}.
 \end{gathered}
\]
Then for every $q\ge0$,
\[
                    S_H(v)\le J+\frac{G}{n}q.                  \tag{1}
\]
In particular $J<3$. When $n>3r$, all such valuations satisfy
$rS_H(v)<n+q$ for $0\le q\le Q$, where one may take
\[
                      Q=\frac{n(n-3r)}{rG}.
\]
If $rG\le n$ and $rJ<n$, inequality (1) itself suffices for every
$q$. If $rG>n$ and $rJ<n$ one can take any
$Q<n(n-rJ)/(rG-n)$. In particular, for $n\asymp r\log n$ one has
$Q\asymp \sqrt n\log n$ for the displayed finite choice, since
$d/n\longrightarrow1$.

When $T=\tau_p(H)\le B$, the pointed quadric overlap and signed
adjunction give $(n+1)S_H(v)\le B(q+3)+n-3$.
Otherwise take a prime divisor $P\in |hH|$ of ratio $b>B$ and set
$e=d/b$. The principal-bin angular argument is independent of the
relation between dimension and index: put $K=\lceil(n+1)b\rceil$,
$N=n-1$, $G_q=b(1+q)-d-q$ and
\[
                 D_{\rm ang}=\frac{N}{n}
                        \left(\frac{e}{3n(b+2)}\right)^{1/N}.
\]
Its intermediate bound, including either sign of $G_q$, is
\[
 S_H(v)\le\frac{n-2+e+d}{n+1}
  +q\left(\frac{n-1+e}{n+1}-D_{\rm ang}\right)
  +\frac{(G_q)_+}{K}.
\]
For clarity, uniqueness of the high-ratio prime gives $T=b$, so the
full ordinary grades range from $0$ to $bm$. Jensen over their
principal bins gives the stated $D_{\rm ang}$ because the entropy of
the bin variable is bounded by $1+\log(b+2)$ and the number of grades
is at most $bm+1$. The mean separated class is within $1/K$ below
$1/(n+1)$. This uses only the general bounds
$s\le(n-2+b+e)/(n+1)$ and $v(P)/h\le d+q$.
Set $A=3n(b+2)/e$. Since $A\ge1$,
$1-D_{\rm ang}\le(1+\log A)/n$. Also $(G_q)_+\le bq$ since
$b\le d$. The coefficient of $q$ is therefore at most
\[
  \frac{e-1}{n+1}+\frac{1+\log A}{n}
   \le\frac{e+\log A}{n}\le
     \frac{\sqrt{d/2}+\log(3n(d+2))}{n}.
\]
This proves (1). The intercept bound follows from $d\le n$ and
$3\sqrt{2n}+6<2n$ for $n\ge10$.

\paragraph{A linear flag at an arbitrary point.}
For a smooth hypersurface $X_d^n$, $d\ge4$, $n\ge4$, fix $p\in X$
and let $T=\tau_p(H)$. A very general flag of pointed linear sections
down to a smooth threefold $Y$ and surface $S$ satisfies
\[
 \varepsilon_p(H|_S)\ge d/B,\qquad B=\max(T,\sqrt{2d}).
\]
By Abban--Zhuang's pointed lifting lemma, $\tau_p(H|_Y)\le
\max(T,\sqrt d)$. Fix a family of integral curves on the varying
pointed hyperplane sections of $Y$, with fixed degree $a$ and
multiplicity $k\ge1$ at $p$, dominating the hyperplane parameters.
These families suffice by countability of Hilbert strata.

If the swept locus is proper, it is an integral surface $D\subset Y$.
Unless it is a cone from $p$, its general pointed hyperplane section
is integral by Bertini and equals the chosen curve. Since $D\sim hH$
we get $a=dh$, and general slicing preserves $k=\operatorname{mult}_pD$.
Thus $dk/a\le\tau_p(H|_Y)$. If projection from $p$ has image
dimension one, then $D$ is a cone of lines lying in $T_pY$. The
tangent hyperplane section is integral: the singular points of that
projective surface lie in a fiber of the finite Gauss map of $Y$,
whereas any reducible or multiple surface hypersurface in $\mathbf P^3$
has a positive-dimensional singular locus. Hence the cone is the
whole tangent section and $\tau_p(H|_Y)=d$. In that case
$B\ge d$ and the very ample bound $a/k\ge1$ suffices.

If the curves cover $Y$, their strict transforms cover the blowup
$Y'$ of $p$. A resolved two-dimensional subfamily of their
normalizations has generically finite dominant evaluation to $Y'$:
at a general point keep two parameter directions whose evaluated
derivatives complement the fiber tangent line. Restricting the
determinant of this evaluation to a general complete normalized fiber
gives
\[
                   2g-2\ge(d-5)a+2k .
\]
Here the intersection with the exceptional divisor is $k$, since the
smooth surface section through $p$ contains the curve and its
strict transform is the point blowup of that surface. On $S$, adjunction
on its point blowup and the Hodge index theorem give
$2g-2\le a^2/d+(d-4)a-k(k-1)$. Consequently
$k(k+1)\le a^2/d+a$. Furthermore an integral projective curve of
degree $a$ has normalized genus at most $(a-1)(a-2)/2$ by general
birational plane projection. The lower bound on $2g-2$ implies
$a(a-d+2)\ge2k>0$, hence $a\ge d-1$. If $a=d-1$ one has
$a-k\le d-2<a^2/d$; for $a\ge d$ the same strict inequality is
immediate. It follows that $k^2<2a^2/d$. This proves the asserted
bound on $\varepsilon_p$ for all curve classes, with no Picard-rank
hypothesis on the pointed surface.

Abban--Zhuang's surface inequality gives $\delta_p(S,H)\ge3/B$.
Retaining the $n-2$ transverse cutting hyperplanes with coefficient
$1-1/B$ in basis adjunction lifts this inequality to
\[
           (n+1)S_H(v)\le B(q+2)+n-2,
     \qquad q=A_X(v)-n,\quad v(\mathfrak m_p)=1.
\]
Indeed on each dimension-$j$ rung the distinguished section has mean
coefficient $1/(j+1)$, the residual restricted basis series has full
pieces and first Chern class $jH/(j+1)$, and adding the given boundary
puts coefficient one on the cut at multiplier $(j+1)/B$. All coefficients
used in inversion of adjunction are nonnegative, and the remaining
cuts are transverse. Each cutting hyperplane has $v$-order at least one.

\section{Tilting the basis of an optimizer}
Suppose $\beta=\delta(X,H)\le r$ and let $v$ be a divisorial minimizer, $S=S_H(v)$. The distribution of its normalized vanishing orders is the pushforward of uniform volume on the Okounkov body by the first coordinate, up to the scaling of the divisorial valuation: start the flag with the prime on a smooth model of $X$ and use the complete pulled-back section ring. This also proves weak convergence of the finite weight distributions and convergence of bounded quantile integrals \cite{ar0805.4559,ar1706.04548}. Its distribution function $F$ is log-concave by Brunn--Minkowski on the sublevel slices, and continuous on its support, since that body is full-dimensional. Write $Q(u)=F^{-1}(u)/S$ for the normalized quantile, of mean one.

A general finite linear projection $\pi:X\to\mathbf P^n$ of degree $d$ is \'etale at the generic point of the center of $v$. The restriction $v_0$ to the base field is divisorial up to scale, with $A_{\mathbf P^n}(v_0)=A_X(v)$ for this scaling. This follows by \'etale base change of a smooth model carrying the base divisor: at the generic point over the chosen center the map is \'etale and the relative canonical divisor contributes zero. Equivalently use the finite-morphism ramification formula for log discrepancies, with ramification order zero here. The polynomial subseries from the projection has asymptotic dimension fraction $1/d$ and mean at most $A_X(v)/(n+1)$ by K-semistability of $\mathbf P^n$. Its own distribution function $F_0$ is continuous and log-concave by the same Okounkov argument. If $U$ has this distribution, then $\mathbb E\log F_0(U)=-1$, so Jensen's inequality gives $F_0(\mathbb EU)\ge\exp(-1)$. The quotients by the high-order subspaces of the polynomial subseries inject into the corresponding complete-series quotients. Passing to the limiting continuous distributions gives
\[
                Q(p)\le z,\qquad p=1/(3d),\qquad z=\beta/(n+1).
\]

Here is an explicit bound on the cost of discarding a lower fraction $f$ with $p<f<1$. Put $L=\log(f/p)$, $D=\log(1/p)-1-f(L-1)$. Then
\begin{equation}\label{global:eq:quant}
 Q(f)\le z+\big[1-(1-f)z\big]L/D
 \qquad\text{when }0\le z\le1 .
\end{equation}
Indeed $D=\int_f^1\log(u/p)\,du>0$. If $Q(f)>z$, convexity of $a\mapsto Q(e^a)$ implies
$Q(u)\ge z+(Q(f)-z)\log(u/p)/L$ for $u\ge f$. Integrate and use the nonnegativity of the remaining quantiles. If $Q(f)\le z$, the inequality is immediate.

More generally, suppose the bottom $f$ fraction has normalized mean at most $z_*<1$. Choose rational $K>0$, $f\in\mathbf Q$ and set $\lambda=K(1-f)/f$. At a large level, take an adapted basis for $v$, and let $D_m^+$ be the sum of the divisors of the upper $1-f$ fraction, divided by $m$ times their number. Thus $D_m^+\sim_{\mathbf Q}H$. For any fixed $\varepsilon>0$, choosing $m$ large gives the fixed rational boundary $B=\lambda D_m^+$ and $k=\lambda+K+\varepsilon$ with
\begin{equation}\label{global:eq:tilt}
 w(B)\le k S_H(w)\quad\text{for all }w,\qquad
 \beta-k\le \beta'=\delta(X,B;H)
 \le\beta-\lambda-K(1-z_*)+\varepsilon .
\end{equation}
The first bound follows because every selected vector extends to the full basis and $S_{m,H}(w)\le(1+o(1))S_H(w)$ uniformly. The upper bound for $\beta'$ uses $v$ and convergence of its top-fraction mean; rounding at most one endpoint is harmless. If $k<\beta$, the pair is klt since $A_X(w)\ge\beta S_H(w)$. If also $\beta'<r-\lambda$, it is log Fano and has a divisorial optimizing valuation by \cite{ar2102.09405}. For every such optimizer $w$,
\begin{equation}\label{global:eq:sumcap}
   A_X(w)\le(\beta'+k)S_H(w).
\end{equation}

\paragraph{Two half-basis estimates.}
The empirical full-series weights of a divisorial valuation $v$ have
a continuous limiting distribution. Indeed one may start an
Okounkov flag on a model with the divisor of $v$; the limiting law
is the first-coordinate projection of Lebesgue measure on the convex
body. Both its lower and upper cumulative functions are log-concave.
The induced valuation $w$ on $\mathbf P^n$ from a general finite
projection has the same discrepancy. With $F_0$ its lower cumulative
function, Jensen gives
\[
 F_0(\mathbb E W_0)\ge \exp(\mathbb E\log F_0(W_0))=\mathrm e^{-1},
\]
since $F_0(W_0)$ is uniform. As $\delta(\mathbf P^n,\mathcal O(1))=n+1$,
the degree-$d$ finite projection subseries proves
$F(a)\ge p$ for $a=A_X(v)/(n+1)$ and $p=1/(3d)$.

Write $S=\mathbb E W$, let $m$ be the median, and set
$D_- =\mathbb E(m-W)_+$, $D_+=\mathbb E(W-m)_+$,
$L=\log(1/(2p))$. The relative excess of the top-half mean is
$(D_-+D_+)/S$. The geometric chord on the lower side of $m$ gives
$D_-\ge l(m-a)$ for $l=(1/2-p)/L$. On the upper side, the quantile
$F^{-1}(\exp z)$ is convex in $z$. Extrapolating its chord between
$\log p$ and $\log(1/2)$ gives $D_+\ge h_0(m-a)$, where
$h_0=(\log 2-1/2)/L$; this follows by integration of $\log(2u)$ for
$1/2\le u\le1$. If $a\ge m$, both displayed lower bounds hold
automatically. Using $S=m+D_+-D_-$, the resulting gain is
\[
 \frac{D_-+D_+}{S}\ge
 \frac{(\log2-p)(1-a/S)}{L+\log2-1+p}.
\]
Indeed the difference after multiplication by $S(1+h_0-l)$ is
$(1+2h_0)[D_--l(m-a)]+(1-2l)[D_+-h_0(m-a)]\ge0$.

For any fixed polarized variety the largest normalized mean of a
half-basis satisfies, uniformly over divisorial $w$,
\[
 S^{\rm top}_{m,1/2}(w)\le(1+\log2+o(1))S_H(w).
\]
Fix an Okounkov flag independently of $w$. The compatible bases in
Blum--Jonsson Lem.~2.9 have normalized finite weights bounded
slotwise by the concave transform $G_w$ sampled at the semigroup's
degree-$m$ lattice points. Put $g_w=G_w/T_H(w)$ so $0\le g_w\le1$.
On each compact subset of the body's interior the bounded concave
functions $g_w$ are uniformly Lipschitz. Consequently
$(g_w-z)_+$ are equicontinuous there for every $w$ and $0\le z\le1$.
The integral of each such function against the empirical normalized
lattice measure differs from its integral against uniform body
measure by $o(1)$ uniformly in $w,z$: first remove a sufficiently
thin boundary layer (the body has boundary of measure zero), then
apply compactness of the remaining equicontinuous family and weak
convergence of the lattice measures. Continuous cutoffs between
two interior layers justify the restriction of measures. For a
bounded scalar random variable $U$ the top-half mean, with fractional
endpoints if necessary, equals
$\inf_z[z+2\mathbb E(U-z)_+]$. The ceiling or floor of half the
finite dimension changes the normalized estimate by $o(T_H(w))$.
Hence the discrete top mean is bounded by the limiting top mean
plus $o(T_H(w))$, uniformly in $w$. Here $T_H(w)\le(n+1)S_H(w)$.

For the last comparison, write the limiting $W$ in distribution as
$f(E)$ for a standard exponential $E$. Its log-concave survival
implies $f$ is nonnegative, nondecreasing and concave, using the
generalized inverse of minus the log survival. Hence $f(t)/t$ is
nonincreasing. Relative to the probability measure $t e^{-t}dt$,
its covariance with $\mathbf1_{t\ge\log2}$ is nonpositive. The share
of the total mean in the top half is therefore at most
$\int_{\log2}^\infty t e^{-t}dt=(1+\log2)/2$. This proves the
uniform half-basis estimate.

\paragraph{Curve images of dimension one.}
If the center $W$ of the tilted optimizer is a curve, it is rational
by the negative canonical degree on its resolution. At a fixed
general point $x\in W_{\rm sm}$ the complete degree-$k$ subsystem
on $W$ has $ke+O(1)$ sections, $e=H\cdot W$. It contains a section
vanishing to order $ke-O(1)$ at $x$. Pulled back through the rational
map $Y\dashrightarrow W$, this gives a section in $|kL_Y|$ with
order at least $ke-O(1)$ on a fixed prime fiber component $F$
meeting the regular locus of the map, since its fixed part is
effective. Thus $T_{L_Y}(F)\ge e$ and
$\beta'\le\delta(Y,L_Y)\le(n+1)/e$ using
$S_{L_Y}(F)\ge T_{L_Y}(F)/(n+1)$. For $s=\dim W\ge2$ the
Jovinelly--Lehmann--Riedl bend-and-break bound gives a curve of
degree $(1+1/s)n/\beta'$ or less. Consequently a rational curve of
degree at most $3n/(2\beta')$ exists for every positive center.

\section{The all-index linear ranges}\paragraph{The first range.}
The preceding results imply that every smooth degree-$d$ Fano
hypersurface of index $r=n+2-d\ge3$ satisfying $n\ge6r$ is
uniformly K-stable. We may assume $6r\le n<r(r-1)$ by the quadratic
range, so $r\ge8$, $n\ge48$, and $d<r^2$.

Suppose $\beta=\delta(X,H)\le r$ and choose the divisorial optimizer.
The projection subseries and the log-concavity of the lower cumulative
distribution give the relative top-half gain at least
\[
 g_* =\frac{(1-2r/(n+1))(1-1/d)}{2\log d+1-1/d}.
\]
Indeed the lower cumulative function $F$ has $F(2A/(n+1))\ge1/(2d)$.
With median $m$, mean $S$ and top-half mean $(1+u)S$, one has
$uS=\mathbb E|W-m|$ and $m\ge(1-u)S$. On the interval between
$2A/(n+1)$ and $m$, log-concavity bounds $F$ below by its geometric
chord with endpoint values $1/(2d)$ and $1/2$. Its integral is
$(m-2A/(n+1))_+(1-1/d)/(2\log d)$, at most $uS$.
Since $A/S=\beta\le r$, rearranging gives $u\ge g_*$.

Fix any $0<g<g_*$, put $k=1/16$, $h=13/100$ and choose a sufficiently
high level top-half basis divisor $D\sim_{\mathbf Q}H$ with
$v(D)>(1+g)S_H(v)$ and $k\beta w(D)\le hA_X(w)$ for every divisorial
$w$. The latter follows from uniform finite-level comparison since
$2k<h$. Write $B=k\beta D$ and
\[
 \beta'=\delta(X,B;H)<b\beta,\qquad b=1-k(1+g)<15/16.
\]
The tilted pair has a divisorial optimizer $v'$. If its center is
positive, the underlying Fano of the tilted special degeneration yields
a rational $C$ in that center with $e\le2n/\beta'$. The threshold
bound for $C$ forces either $\beta'\ge(1-h)(n-1)$ or
$\beta'\Theta(e)\ge1-h$. Both are excluded by the jet test at $4/n$:
\[
 \beta'\Theta(2n/\beta')\le 4br/n+2(1-4/n)^n c_{n,d},\qquad
 c_{n,d}=\frac{4(d-1)}d+\frac{(2d-1)(n-4)}{d(n+1)}<6.
\]
Indeed its upper bound is $5/8+12\exp(-4)<0.87=1-h$.

Hence $c_X(v')=p$ is a point. Normalize $v'(\mathfrak m_p)=1$ and
put $\mu=A_X(v')$. Take $Q\ge0$ with $\mu\ge n+Q$. Choose
$a>n/[(1-h)(n+Q)]$ and set $t=(1+a/n)\beta'$, $T=(1+a/n)br$,
\[
 \eta=r\{1-k-b(n+a)/(n+1)\},\qquad E=\lfloor1/\eta\rfloor
 \quad (\eta>0).
\]
For a $v'$-compatible basis also adapted to a hyperplane avoiding $p$,
remove that hyperplane to obtain $\Gamma_m\sim (1-a_m)H$ with
$a_m\to1/(n+1)$ and $v'(\Gamma_m)\to S_H(v')$. The ample slack
$r-k\beta-t(1-a_m)$ has lower limit at least $\eta$ by the bounds for
$\beta$ and $\beta'$. If there is a positive non-klt center at
coefficient $t$ for infinitely many $m$, let $t_m\le t$ be the first
such threshold. The uniform complete-series comparison gives
$t_m\ge\beta'/(1+o(1))$. Thus adding an arbitrarily small general
pointed hyperplane coefficient $\epsilon_m\to0$ makes the pair
non-lc at $p$, preserving its positive lc centers off $p$, as in the
short-curve lemma. That lemma produces a rational curve of degree
at most $E$ on one of the original positive centers. For its testing
valuation $w$ one has $w(B)\le h A_X(w)$ and
$w(\Gamma_m)\le(1+o(1))S_H(w)$ uniformly. Therefore this is
impossible as soon as $E=0$, or as soon as
\[
     T<(1-h)(n-1),\qquad
     \frac{T}{n}\left[4+E(1-4/n)^n c_{n,d}\right]<1-h .
\]
Under these tests $(X,B+t\Gamma_m)$ is thus klt in a punctured
neighborhood of $p$ for large $m$. The point Nadel inequality gives
\[
   (1+1/\mu)A_X(v')\ \ge\ v'(B)+tS_H(v')
    = A_X(v')+\frac an(A_X(v')-v'(B))
    \ \ge\ (1+a(1-h)/n)A_X(v'),
\]
which contradicts the choice of $a$.

For $r\ge40$ one can take $Q=0$, $a=29/25$ and
$g=0.65/(4\log r+1)<g_*$. As $r/(4\log r+1)$ is increasing and
$\log40<3.7$, direct substitution gives
\[
  \eta= rkg - (a-1)br/(n+1)
    >\frac{0.65r}{16(4\log r+1)}-0.025>1/13.
\]
Thus $E\le12$, $T<0.943r$. The displayed jet test follows from
$4+72\exp(-4)<5.32$ and $n\ge6r$.

For $8\le r<40$ the program \path{r2-lit-six-index.py} in the archive
described in Chapter~\ref{ch:checks} supplies $Q$ and verifies the above
strict inequalities for all $15136$ remaining pairs. It uses the
universal point envelope $S_H(w)\le 2+G(A_X(w)-n)/n$ where
\[
 G\ge\max\{\sqrt{2d}+4,\ \sqrt{d/2}+\log(3n(d+2))\}.
\]
Indeed the high-branch intercept is below two since
$d+\sqrt{d/2}<n+2$ for $n<r(r-1)$, and the low-branch intercept is
below two since $3\sqrt{2n}+9<n+2$ for $n\ge48$.
Set $R=br/(1-h)>0$. At the optimizer $RS_H(v')\ge A_X(v')$, hence,
when $n>2R$ and $RG>n$, the point envelope gives
\[
      Q=\frac{n-2R}{RG/n-1}.
\]
If $RG\le n$, $n>2R$ excludes the point at once. The script uses
$a=1.01n/[(1-h)(n+Q)]$ and rounds square roots and logarithms
upward. All thresholds, powers, floors and inequalities are evaluated
with integers or rational fractions.

\paragraph{The all-index bound $n\ge4r$.}
The curve bound $3n/(2\beta')$, a variable fraction of a basis,
and the span-sliced normal estimate improve the all-index linear
range to $n\ge4r$. We retain the notation and the point-optimizer
argument of the $6r$ proof. Only $4r\le n<\min(6r,r(r-1))$
requires verification, so $n\ge24$ and $r\ge6$.

Take a fraction $\tau$ of the highest-weight sections, with $u=1-\tau$.
For $p=1/(2.72d)$, $F(A/(n+1))\ge p$ by the projection argument.
Set $m=F^{-1}(u)$, $L=\log(u/p)$, $D_\pm=\mathbb E(\pm(W-m))_+$.
Integrating the geometric chord below $m$ and the quantile chord
above $m$, as in the half-basis estimate, gives
\[
 D_-\ge l(m-A/(n+1)),\quad D_+\ge h_0(m-A/(n+1)),\quad
 l=(u-p)/L,\quad h_0=(-\log u-\tau)/L.
\]
Here $u>p$ in all the choices used. The top mean minus the mean is
$D_-+(u/\tau)D_+$. Substituting $S=m+D_+-D_-$ and using $l\le u$
shows that the relative gain is at least
\[
    g_* =\frac{-u\log(u)/\tau-p}{\log(1/p)-1+p}
             \left(1-\frac r{n+1}\right).
\]
Indeed with $c=(l+(u/\tau)h_0)/(1+h_0-l)$ the difference after
subtracting $c(S-A/(n+1))$ is the sum of the two displayed
nonnegative remainders with coefficients $1+c$ and $u/\tau-c$.
The same uniform lattice comparison and exponential survival proof
give the finite top mean upper bound $(1-\log\tau+o(1))S_H(w)$
uniformly over all divisorial $w$: the normalized upper tail mean
of the exponential is $1-\log\tau$.

Fix $H_0$ avoiding the center of the original optimizer. The complete
basis may be chosen compatible with its valuation and with $H_0$.
Remove all copies of $H_0$ from the selected top fraction. The
valuation of the resulting effective boundary $D$ is unchanged.
In fact $D\sim_{\mathbf Q}c_m H$ with $c_m<1-\alpha$ for all
sufficiently large $m$, where we may use
\[
       \alpha=\frac{1+u\log(u)/\tau}{n+1} .
\]
The proportion in the full basis divisible by $H_0^{\lceil mz\rceil}$
tends to $(1-z)^n$ for $0\le z\le1$. Thus the normalized mean
removed from a proportion $\tau$ of its members has lower limit at
least $\tau^{-1}\int_0^1[(1-z)^n-u]_+\,dz$, or
$[1-(n+1)u+nu^{1+1/n}]/[(n+1)\tau]$. This is strictly greater
than $\alpha$ since $u^{1/n}>1+(\log u)/n$. Consequently the
tilted boundary $B=k\beta D$ has class at most $k\beta(1-\alpha)H$.
For rational $g<g_*$ and $h>k(1-\log\tau)$ the established
conclusions $w(B)\le hA(w)$ and $\beta'<b\beta$,
$b=1-k(1+g)$, apply. Its underlying-central-fiber argument uses
the actual rational class $k\beta c_m H$.

For $e\le23$ use the reciprocal threshold bound $U(e)$ from a
span flag to $j=\max(3,e+1)$, namely
\[
 U(e)=\frac{j+1}{n+1}\max\left(\frac1{j-1},\ \Theta_{j,d,e}(t)\right)
\]
at any chosen rational $t$. The bound at $e$ also controls every
curve of smaller degree by taking the same $j$ and $t$.
Positive optimizers are excluded whenever
$U(e)\min(br,3n/(2e))<1-h$ for $1\le e\le23$.
For $e\ge24$ one has $n/\beta'\ge16$; the full-dimensional
trial $4/n$ bounds $\beta'\Theta(e)$ by
$4/16+9\exp(-4)<0.42$. The full codimension term is included in
the checker. For the point optimizer use $R=br/(1-h)$ and the
point envelope as before. For $\mu\ge n+Q$ and
$a>n/[(1-h)(n+Q)]$, put
\[
 T=(1+a/n)br,\qquad \eta=r\{1-k(1-\alpha)-b(n+a)/(n+1)\}.
\]
When $\eta>0$ and $E=\lfloor1/\eta\rfloor$, the point argument
requires $E=0$ or $TU(E)<1-h$.

For completeness, the finite checker slightly improves the universal
point slope without changing its proof. In the low branch
$s=S_H(v)\le (2\sqrt{2d}+n-2)/(n+1)+\sqrt{2d}\,q/n$ by the
Picard-free point flag. At a high principal ratio $b_0>\sqrt{2d}$
put $e_0=d/b_0$, $N=n-1$. The intercept is at most
$(n-2+d+\sqrt{d/2})/(n+1)$. The entropy of the principal bin obeys
$\exp \mathcal H(i)\le \mathrm e(1+b_0+1/(n+1))<
(11/4)(b_0+11/10)$ in the limit, since its mean is at most
$b_0+1/(n+1)$. Therefore the earlier bin formula bounds $n$ times
the slope by
\[
 e_0+\log A - c_0,\qquad A=\frac{11n(b_0+11/10)}{4e_0},\qquad
 c_0=\frac{\log^2(5n)}{2(n-1)}
              \left(1-\frac{\log(3n(d+2))}{3(n-1)}\right).
\]
Indeed the slope is at most
$(n+e_0)/(n+1)-(n-1)A^{-1/(n-1)}/n$. Use the cubic Taylor lower
bound for $\exp(-x)$ and $5n<A<3n(d+2)$. The functions
$d/b_0+\log(11nb_0(b_0+11/10)/(4d))$ attain their maximum at
one of the endpoints $\sqrt{2d},d$: the sign of the derivative
is that of a quadratic with a single positive root. The checker
rounds both endpoint values up and $c_0$ down. This gives
$S_H(v)\le J+Gq/n$ with the stated intercepts and slopes, so
$Q=(n-RJ)/(RG/n-1)$ when $RG>n>RJ$.

For $r\ge256$ take $\tau=3/4$, $k=1/40$, $h=0.033$, $a=1.04$,
$Q=0$ and $g=0.34/(2\log r+1)$. Since $n<r(r-1)$, $n\ge4r$,
and $-u\log(u)/\tau=\log 4/3$, one has $g<g_*$. Without the
bonus $\alpha$, the slack is still at least
\[
      \frac{0.34r}{40(2\log r+1)}-\frac{0.04(39/40)}4>1/6.
\]
The right side is increasing, and $\log256<5.55$. Hence $E\le5$.
The large-index table verifies the positive-center test for all
$e\le23$ and the point test for $e\le5$ using $d=\infty$ in the
increasing coefficient ratios $(d-1)/d$ and $(2d-1)/d$.

For $6\le r<256$ the program \path{r2-lit-four-index.py} tests every
$4r\le n<\min(6r,r(r-1))$. Logarithm bounds use the positive
$\operatorname{arctanh}$ series after extracting a power of two,
with a geometric upper bound on the remainder. Square roots are
rounded outward by the integer square root. Its witnesses evaluate
the normal-jet polynomials at exact rational arguments.

\paragraph{The all-index range $4n\ge15r$.}
Only $15r/4\le n<\min(4r,r(r-1))$ is new. In this strip $n\ge19$.
We give the two strengthened estimates used by the additional finite
checker.

First, if the center $W$ of the tilted optimizer is a surface, let
$S$ denote a smooth resolution of $W$ with big nef Cartier
polarization $H_S$. The determinant slope inequality gives a moving
curve class $C_0$ on $S$ such that $(K_S+cH_S)C_0<0$ for any fixed
$0<c<2\beta'/n$. Run the surface cone theorem. A negative extremal
ray is $K_S$-negative by nefness of $H_S$. Every birational extremal
contraction from a smooth surface therefore blows down a $(-1)$-curve
$E$ to a smooth point. The pushforward $H'$ is again big nef Cartier,
and $H_S=\pi^*H'-aE$ with $a=H_S E\ge0$. Then
$K_S+cH_S=\pi^*(K'+cH')+(1-ca)E$ and $1-ca>0$. As $C_0$ is
moving, its pushforward still has negative intersection and is moving.
After finitely many blowdowns the chosen negative contraction is a
Mori fiber space. Over a curve its general fiber is rational of
$H'$-degree $<2/c$. It avoids all contracted points and gives the
same degree on $W$.

Over a point the final surface is $\mathbf P^2$ with $H'=h\mathcal O(1)$
and $h<3/c$. If one of the blowup multiplicities is positive, at
least one proper base point of the birational morphism to $\mathbf P^2$
has positive multiplicity. Indeed the intersection of the nef $H_S$
with each strict exceptional curve is its own multiplicity minus the
sum of the multiplicities proximate to it. The strict transform of a
general line through such a proper point has positive $H_S$-degree
at most $h-1$, by mobility and bigness of $H_S$; its image is a
rational curve on $W$.

If all multiplicities vanish, $H_S$ is the pullback of $h\mathcal O(1)$.
Powers of a general line on $\mathbf P^2$ yield sections of $|kH_S|$
vanishing to order $kh$ along the corresponding moving prime on $S$.
They descend to the normalization of $W$ by the projection formula.
To descend to $W$, multiply by a fixed nonzero section of a sufficiently
high twist of the conductor ideal of its normalization. After this
fixed degree shift the corresponding zero-weight sections of $Y$ give
$T_{L_Y}(F)\ge h$ on a prime divisor $F$ on $Y$ mapping to the general
line, exactly as for curve centers. Hence $h\le(n+1)/\beta'$. Letting
$c$ tend up to $2\beta'/n$ and using integral degrees shows, including
the other center dimensions treated earlier, that some rational curve
in every positive optimizing center has degree $e$ with
\[
 \beta'\le M_e:=\max\{(n+1)/e,\ 4n/(3e),\ 3n/(2(e+1))\}.
\]
Here centers of dimension at least three use Jovinelly--Lehmann--Riedl
with factor at most $4/3$. If $U=(b+h)r$, the discrepancy ratio at
the optimizer is at most $\min(U,M_e+hr,M_e/(1-h))$, since the
tilted boundary has $w(B)\le h\beta S_H(w)$.

Second, in each principal bin of the point-envelope proof all nonzero
elements have $P$-order fraction $w$ in its interval
$[i/K,(i+1)/K]$. Indeed the original basis is adapted to $P$-order,
so its grouped vectors retain independent $P$-initial forms in each
grade. The residual divisor after removing exactly that power meets
$P$ properly. By intersection multiplicities, its ordinary order is
at most $e_0 m(1-w)$, where $P\sim h_P H$,
$b_0=\operatorname{mult}_p(P)/h_P$ and $e_0=d/b_0$. It follows
that the ordinary grades $j/m$ of this bin lie in an interval of
length at most $L=e_0+(b_0-e_0)/K$. Use the conditional entropy of
ordinary grade given bin; its upper bound is $\log(mL+2)$. The
angular moment is at least $(n-1)A^{-1/(n-1)}/n$ where
$A=nL\exp(\mathcal H(i))/d$, with limiting upper bounds understood.
Furthermore the exact difference between the full separated mean and
the binwise forced mean tends to
\[
 \frac1{n+1}-\frac1K\sum_{i=1}^{K-1}(1-i/K)^n\le\frac1{2K}.
\]
The identity uses the limiting beta law of the full $P$-order class
fraction and the inequality is the trapezoid rule for a convex function.
Thus choose $K=\lceil(n+1)b_0/2\rceil$. The entropy estimate
$\exp\mathcal H(i)\le\mathrm e(1+K/(n+1))$ then gives the uniform
upper bound
\[
 A\le A_0=2.72n\left(\frac12+\frac{1+1/(n+1)}{\sqrt{2d}}\right)
              \left(1+\frac2{n+1}\right)
\]
for $b_0>\sqrt{2d}$. The term from the positive part divided by
$2K$ contributes at most $q/(n+1)$. The high-branch slope times
$n$ is therefore bounded by
\[
 1+\frac{n(\sqrt{d/2}-1)}{n+1}+\log A_0
 -\frac{(\log A_0)^2}{2(n-1)}\left(1-\frac{\log A_0}{3(n-1)}\right).
\]
This uses the cubic Taylor bound after first replacing $A$ by the
fixed upper bound. The low-branch slope times $n$ is at most
$n\sqrt{2d}/(n+1)$ and the two intercepts are as in the $4r$ proof.

For the finite cases take $Q$ from these intercepts and slopes with
$R=U$. Write $\eta_0=r[1-k(1-\alpha)-bn/(n+1)]$. If $\eta_0>1$,
the point-center dimension bound directly excludes the optimizer.
Otherwise take a rational $\Delta>0$, with added pointed hyperplane
coefficient tending to $\eta=\eta_0-\Delta n/(n+1)>0$, such that
\[
          \Delta(n+Q)/U+\eta>1.
\]
Use $t=\beta'+\Delta$. If the positive centers at $t$ have been
excluded by the short-curve tests up to degree $\lfloor1/\eta\rfloor$,
adding the free pointed hyperplane with this limiting coefficient
keeps the pair klt off $p$ and ample with slight finite shrinkage.
The Nadel inequality at $p$ gives the contradiction displayed above
because $S_H(v')\ge(n+Q)/U$. For testing positive centers it suffices
that $(U+\Delta)S_H(w)<A_X(w)$ on every resulting curve. The
hyperplane is added only in the last Nadel step. Every endpoint in
the tests is strict.

For $r\ge512$ use $\tau=17/20$, $k=1/40$, $h=0.03$, $Q=0$ and
$a=1.04$ in the simpler point test at $(1+a/n)\beta'$. One may take
$g=0.24/(2\log r+1)$ in the quantified gain since
$-(1-\tau)\log(1-\tau)/\tau>0.334$, $n\ge15r/4$,
and $d<r^2$. Without the additional frozen hyperplane, the slack is
already bigger than
\[
 \frac{0.24r}{40(2\log r+1)}-\frac{0.04(39/40)}{15/4}>1/5
\]
(the first term is increasing, and $\log512<6.24$). Thus the point
curve has degree at most four. The large-index table uses the
dimension-free limits of the sliced normal inequalities for
$1\le e\le15$; $e\le3$ uses the plane and twisted-cubic estimates.
For $e\ge16$ the full-dimensional trial $3/n$ suffices, since
$\tfrac32(3/16+5\exp(-3))/(1-h)<1$. For $5\le r<512$ the
program \path{r2-lit-fifteen-quarters.py} verifies every remaining integer pair.

\section{A projection estimate for the bottom of a basis}
Let $\beta=\delta(X,H)\le r$. Normalize the quantile $Q(u)$ of a divisorial optimizer by its mean, so $\int_0^1Q=1$. The subseries of any general finite linear projection has dimension fraction $1/d$ and mean at most $\beta/(n+1)$ in these units. At a finite level the ordered weights of a subspace interlace those of the full space: its quotient by sections of weight at least $t$ injects into the full quotient. Thus the mean of the lowest $1/d$ is at most $\beta/(n+1)$ in the limit. Replacing the exact subseries dimension by $\lfloor h^0(mH)/d\rfloor$ costs $o(1)$ in normalized mean for this fixed valuation, whose maximal normalized weight is bounded.

Put $D(u)=u^{-1}\int_0^uQ(p)\,dp$. Since the lower distribution function is log-concave, the function $Q(e^z)$ is convex for $z<0$. Hence
\[
       D(e^z)=\int_0^\infty Q(e^{z-t})e^{-t}\,dt
\]
is convex as well. The chord joining $u=1/d$ and $u=1$, and monotonicity below $1/d$, now yield the following guaranteed relative gain for the top $1-f$ of a basis:
\begin{equation}\label{global:eq:average-gain}
 g_f=\left(1-\frac{\beta}{n+1}\right)
       \frac f{1-f}\min\left(1,\frac{\log(1/f)}{\log d}\right).
\end{equation}
This argument needs only the mean of the projection subseries. All finite constructions below use rational parameters with a strict allowance for convergence.

\section{Integral pairs at moderate index}
\paragraph{A curve in a surface center.}
In the surface contraction proof of the $15r/4$ strip, stop at the first birational contraction with $a=H_S\cdot E>0$. Every earlier contraction had $a=0$, so its pullback preserves the degree of the strict transform of $E$. Since $(K_S+cH_S)E=-1+ca<0$, the original surface has a rational curve of positive degree $<1/c$. Its image in the original center is a rational curve of degree no larger. If no such contraction intervenes, the polarization pulls back from the final smooth surface. For a Mori fibration use its general rational fiber, and for $\mathbf P^2$ use the general-line section and conductor argument already given. Letting $c$ increase to $2\beta'/n$ proves the existence of a rational curve of degree at most $(n+1)/\beta'$ in every surface optimizing center. Existence of one curve suffices since the threshold test includes every valuation whose center contains that curve.

\paragraph{Subdividing the original threshold.}
The following specifies the quotient and point tests used in the finite table. Use the notation of the slope and adjunction sections. Put $m=n/(n+1)$ in these formulas. On a rational interval $a\le\beta/r\le b$, multiply the top fraction by $\ell\beta$. Use \eqref{global:eq:average-gain} with $br$ in place of $\beta$, call its lower bound $g$, and put
\[
 \alpha=\frac{1-f\log(1/f)/(1-f)}{n+1},\quad
 h=\ell(1-\log(1-f))+\epsilon,\quad t=1-\ell(1+g)+\epsilon,
 \qquad\epsilon=10^{-7}.
\]
Only $\ell>0$ uses these allowances; for the original boundary use $h=0,t=1$. The hyperplane subtracted from the upper fraction is the one avoiding the center of the original optimizer. We may use a fixed sufficiently large level. All positive search choices have positive gain, and the strict finite upper bound $\beta'<r-\lambda$ follows by using the actual gain before the allowances; $\lambda$ denotes the actual class coefficient of the boundary. Put
\[
 B_-=ar(1-h),\quad \lambda_+=\ell br(1-\alpha),\quad
 U=(t+h)br,\quad z=r-\lambda_+-m tbr,\quad c=(1-h)^{-1}.
\]
Then $B_-\le\beta'\le tbr$, $\lambda\le\lambda_+$, and $A_X(w)/S_H(w)\le\beta'+h\beta\le U$ at a tilted optimizer $w$. We also have $\beta'+h\beta\le c\beta'$ since $\beta'\ge(1-h)\beta$. If $s$ is the center dimension, its bound below is $\beta'\ge(r-\lambda_+-s-1)/m$ and $z\le s+1$.

For a positive center let $e$ be the degree of a rational curve supplied by the slope arguments. For $s=1,2$, or $(s+1)\beta'\ge n+1$, the upper bound on $\beta'$ is $(n+1)/e$. Otherwise $\beta'<(n+1)/(s+1)$ and $\beta'\le (s+1)n/(se)$. For $s\ge9$ one may use $(n+1)/10$ and $10n/(9e)$ in the latter case. Apply the dimension inequalities before retaining any case. Let $M$ be the resulting maximum admissible upper bound on $\beta'$ among the center dimensions. The optimizer is impossible in degree $e$ if $B_->M$ or
\[
      \min(U,cM,M+hbr)\,D_{n,d}(e)<1,
\]
where $D_{n,d}(e)$ is the upper bound on $S_H/A_X$ from the span sections: use the line/conic and twisted cubic values up to degree three, then $j=\min(n,e+1)$. For $e\ge41$, $\beta'e\le4n/3$ in all cases ($n\ge16$). The full-dimensional trial $4/n$ bounds the product by $(4/3)c(4/41+6/54)$. The radial term is tested separately by $U<n-1$.

For a point the dimension inequality suffices when $z>1$. Otherwise write the point envelope $S_H(w)\le J+Gq/n$ with $w(\mathfrak m_p)=1$, and take $Q=0$ or the improved $Q=(n-UJ)/(UG/n-1)$ when $n>UJ$, $UG>n$; if $n>UJ$, $UG\le n$ the point is excluded immediately. Take $\Delta=(1-z+\epsilon)U/(n+Q-mU)$, $\eta=z-m\Delta$. Require $\eta>0$ and $(U+\Delta)D_{n,d}(e)<1$ for every integer $1\le e\le\lfloor1/\eta\rfloor$. At the first positive threshold for the residual bases the short-curve lemma then excludes a positive center. The Nadel inequality at the point gives $1\ge\Delta S_H(w)+\eta\ge 1+\epsilon$. If the actual slack exceeds one one may shrink the added hyperplane coefficient to a number less than one still giving the strict contradiction; for the displayed choice $\eta\le1$ whenever $z\le1$. Endpoints can be perturbed by rational amounts because all decisive comparisons are strict.

For $32\le r<50$ and $2r\le n<15r/4$ we allow a separate witness on each dyadic interval in $\beta/r$. The archive described in Chapter~\ref{ch:checks} contains the witnesses in \path{middle-witnesses.json}. The program \path{middle.py} verifies the inequalities and checks that the intervals cover $[0,1]$ for every such pair. Floating arithmetic proposes parameters; all checks that can certify a pair use rational bounds. Different values of $n$ at a fixed $r$ are checked individually. The low-index arguments and the proof of the quadratic range appear in Chapters~\ref{ch:fivefold}, \ref{ch:foundations} and the final sections of this chapter.

\section{Three-dimensional centers and the seventeen-fifths range}

\paragraph{A terminal Fano threefold pseudoindex lemma.}
Let $V$ be a complex projective terminal threefold of Picard number one with
$-K_V$ ample. If every integral rational curve has $-K_V$-degree $>3$, then
$V\simeq\mathbf P^3$. We include the detail needed to pass from an integral
length assumption to the strictly weaker rational length hypothesis here.

Use the canonical index-one covering stack $\mathcal V$. It is proper, lci,
agrees with $V$ off isolated points, and $K_{\mathcal V}$ pulls back $K_V$.
For every nonidentity inertia element its virtual tangent age is strictly
less than two. At a quotient $\tfrac1m(a,-a,1)$ it is $1+\{k/m\}$.
At an index-one-cover hypersurface the virtual tangent character is
the four ambient tangent characters minus the equation character.
Mori's normal forms, listed for instance by J.~A.~Chen,
\emph{Explicit resolution of three dimensional terminal singularities},
arXiv:1310.6445, Sections~3--6, give the following list up to permutation.
Here $w$ is the ambient weight vector, $b$ the equation weight, and
$1\le k<m$.
\[
\begin{array}{c|c|c|c}
 &m&w&b\\ \hline
 cA/m&m&(a,-a,1,0)&0\\
 cAx/2,\ cD/2,\ cE/2&2&(1,1,1,0)&0\\
 cAx/4&4&(1,3,1,2)&2\\
 cD/3&3&(0,2,1,1)&0
\end{array}
\]
In each case $\sum_i\{kw_i/m\}-\{kb/m\}$ is strictly less than two,
also for a nonprimitive element. Gorenstein points have trivial inertia.

Proposition~3.2 of Chen--Tseng,
\emph{Note on characterization of projective spaces}
(\href{https://arxiv.org/abs/math/0509649}{arXiv version}), gives a twisted rational
stable map of anticanonical degree $\ell\le4$ with at most one twisted
point. Its coarse map is birational to its image. Indeed a multiple cover
would have degree greater than six, as would any degeneration with two
nonconstant components. In particular its coarse family is unsplit.
Fix one ordinary point $0$ with smooth image. If a twisted point
$\infty$ is present, its image is fixed under deformation because the
nontrivial inertia is supported over isolated coarse points. The lci
deformation estimate and orbifold Riemann--Roch in Chen--Tseng,
\emph{Cone theorem via Deligne--Mumford stacks}, arXiv:math/0505043,
Lemma~4.6 and its proof, give
\[
 \dim_{[f]}\operatorname{Mor}(\mathcal C,\mathcal V;f|_0)
     \ge \ell-\operatorname{age}_{\infty}(T_{\mathcal V}^{\mathrm{vir}})>1.
\]
Here the tangent-obstruction Euler characteristic is computed from the
dual of the ambient conormal presentation on the source; the numerical
term subtracted is precisely the virtual age in the table. There is just
one dimension of reparametrizations fixing $0$ and $\infty$. The map
to coarse stable maps is quasi-finite, so this gives a moving coarse
rational curve with two distinct fixed images. Bend-and-break contradicts
unsplitness, as in Proposition~3.5 and Lemma~3.7 of the projective-space
paper. Thus the curve has no twisted point. Its anticanonical degree is
integral and hence is four.

For precision, the family argument of Chen--Tseng, loc.~cit.,
Corollaries~3.8 and~3.10, Lemma~3.12 and Proposition~4.2, and Propositions~4.1,
4.3 and~4.4 now applies to this degree-four class. Here are its inputs:
the complete family of these maps is unsplit and has no multiple covers,
because any rational component has degree greater than three. At the
untwisted map the morphism space fixing one ordinary point has dimension
at least four. The subfamily additionally meeting any fixed singular
point has dimension at most two, including parametrizations, by
two-point bend-and-break. Thus the general member avoids the singular
locus. The resulting proper normalized family of rational curves through
a general point has dimension two and a general curve is very free.
The cited propositions identify its tangent map with the exceptional
$\mathbf P^2$ of the blowup at the marked smooth point. Its universal
family is $\mathbf P(\mathcal O\oplus\mathcal O(-1))$ over $\mathbf P^2$,
and its evaluation contracts the negative section and factors through a
finite surjection $g:\mathbf P^3\to V$. In those cited steps the length
hypothesis serves only to prohibit splittings and multiple covers of
this class; both have just been excluded for all degree-four rational
stable maps. The lines through the marked point map birationally to
the degree-four family. Hence $g^*(-K_V)=4\mathcal O(1)$ and the
ramification formula makes $g$ quasi-\'etale. Removing the inverse image of
$V_{\rm sing}$ leaves a simply connected open subset of $\mathbf P^3$
(the removed set is finite). By purity over $V_{\rm reg}$ this is the universal
finite \'etale cover there. Its deck group acts on $\mathbf P^3$ by
normalization, and $V=\mathbf P^3/G$. If $G$ were nontrivial, the projective
line spanned by two eigenvectors of any nonscalar element with distinct
eigenvalues would map with degree at least two. Its rational image would have
anticanonical degree at most two, contrary to hypothesis. This proves
the lemma.

\paragraph{Small-dimensional rational images.}
For center dimension $s=2$ the first contraction of positive $H$-degree
in the surface MMP, if it is birational, contracts a $(-1)$-curve of
degree $<1/c$. All preceding contractions are $H$-trivial and pull
back the successive Cartier polarizations. If there is no such step,
a Mori fiber over a curve has degree $<2/c$, and over a point the last
surface is $\mathbf P^2$, with the high powers of a line giving
$h\le(n+1)/\beta'$, as in the preceding surface argument.
Thus one may always use $e\le(n+1)/\beta'$ when $s=2$.

For $s=3$ take a relative $K$-minimal model $V\to W^\nu$ of the smooth
resolution over the normalization of the optimizing center. It exists by
the threefold MMP, since the morphism is birational, and is
$\mathbf Q$-factorial terminal. Its big nef Cartier divisor $H_V$ is
pulled back from $W^\nu$. For $0<c<3\beta'/n$, a general moving curve
from the slope argument still has negative intersection with $K_V+cH_V$.
Indeed on a common resolution the pullback of the initial canonical
divisor is the pullback of $K_V$ plus an effective divisor and the
moving curve avoids its support generically. Choose a
$(K_V+cH_V)$-negative extremal ray. It is $K_V$-negative and
$H_V$-positive, since $K_V$ is nef over $W^\nu$. If the ray contraction
has positive-dimensional base, a positive-dimensional component of a fiber
has dimension at most two. Jovinelly--Lehmann--Riedl, Theorem~1.2,
therefore supplies a rational curve with $-K_V$-degree at most three;
its image on $W$ has positive $H$-degree $<3/c$. If the contraction
has target a point, $V$ is terminal Fano with Picard number one.
The pseudoindex lemma provides the same conclusion unless $V=\mathbf P^3$.
In that case $H_V=h\mathcal O(1)$, $V\to W^\nu$ is finite and birational
and hence is an isomorphism. Powers of a hyperplane on $V$, multiplied
by a fixed conductor section to descend to $W$, yield $T\ge h$ over a
general hyperplane and $h\le(n+1)/\beta'$, as for the surface case.
Letting $c\nearrow3\beta'/n$ proves $e\le(n+1)/\beta'$ also for
three-dimensional centers.

\paragraph{The range $5n\ge17r$.}
For a positive optimizing center of dimension $s\le3$, or when
$(s+1)\beta'\ge n+1$, write the degree constraint as
$\beta'\le(n+1)/e$. In any remaining dimension $s\ge4$ one has
\[
 \beta'\le\min\{(n+1)/(s+1),(s+1)n/(se)\}.
\]
Here we invoke the strong canonical extension inequality for the first
alternative, and the ordinary quotient slope and optimal bend-and-break
for the second. In particular, in all dimensions the right side is at most
$5n/(4e)$ when $n\ge4$. Let $M$ be the largest permissible right side
over the alternatives. The untilted small case $\beta<(r-1)(n+1)/n$
has no point optimizer by the center-dimension lemma. Its positive
optimizers are also excluded by the following uniform comparisons with
$U=\beta\le r$, $h=0$.

For $r\ge200$, take omitted fraction $f=1/4$ and boundary scale
$k=1/10$. In the new strip $17r/5\le n<15r/4$, the projected bottom-mean
inequality gives
\[
  g\ge\frac{(1/3)\log4\,(1-r/(n+1))}{\log d}
      >\frac{0.325}{\log(4r)}.
\]
Use $h=0.129$ in $w(B)\le h\beta S_H(w)$, $b=1-k(1+g)$, and
$U=(b+h)r\le1.029r$. The fixed finite-level boundary allows
strict arbitrarily small margins in these budgets. Even without
subtracting a free hyperplane from $B$ the point dimension slack is
\[
 r(1-k)-r b\frac{n}{n+1}
 =\frac{0.9r}{n+1}+0.1r\frac{n}{n+1}g>1 .
\]
The last term exceeds $0.96$ at $r=200$ and the lower bound
$0.0325r(680/681)/\log(4r)$ increases thereafter. The first term is
greater than $0.23$. Thus every optimizer is positive-dimensional.
For $e\le3$ multiply $U/(n+1)$ by $2$ or $8/3$ respectively.
For $4\le e\le6$ use
\[
 \frac{A_X(w)}{(n+1)S_H(w)}
 \le\frac1{1-h}\max\{1/e,\min(1/5,5/(4e))\}
\]
at a hypothetical optimizing valuation. The sliced jet factors in
these degrees, at the fixed trial parameters of the $4r$ table, give
products $<1$. The function \texttt{large} in \path{r2-lit-seventeen-fifths.py} evaluates these products with rational arithmetic.
For $e\ge7$ the full-dimensional trial $3/n$ gives the upper
product $\tfrac54(3/7+1/4)/(1-h)<1$, using $5\exp(-3)<1/4$.
The independent ordinary term is harmless since $U<n-1$. These
comparisons also exclude every untilted positive optimizer for the
small-$\beta$ branch in all dimensions in the new strip.

For $5\le r<200$ the same program verifies the finite range
$17r/5\le n<\min(15r/4,r(r-1))$. In the remaining high-$\beta$ branch
it uses $\beta'\ge(1-h)(r-1)(n+1)/n$. A weak alternative of
dimension $s$ can be discarded if this lower bound exceeds its cap,
if the dimension slack exceeds $s+1$, or if that cap violates the
center-dimension inequality. For each remaining degree the budget
is $\min(U,M/(1-h),M+hr)$. For the pointed stage the checker uses the
previously proved sharp slope envelope and the Nadel inequality with
the auxiliary pointed hyperplane. In the finite beta-law estimates
one may replace an occurrence of $n(1-e^{-z/n})$ by its upper bound
$z-z^2/(2n)+z^3/(6n^2)$: this applies both to the uniform upper mean on a
top fraction and, with $z=\log(1/f)$, to the removed upper hyperplane mean.
For the first bound, the inverse upper quantile is concave in
$u=1-(1-p)^{1/n}$ by concavity of the survival volume root, vanishes
at the origin, and its ratio to $u$ is decreasing. Consequently its
top-fraction share is bounded by that of $u$, whose conditional top
mean divided by its full mean is $1+n(1-(1-f)^{1/n})$. The same
calculation for the hyperplane filter, whose survival law is exactly
$(1-t)^n$, computes the bottom retained hyperplane coefficient. The
uniform finite-level approximation is the one in the top-half lemma.
Every successful finite witness is checked with rational arithmetic.

\paragraph{The quartic refinement at $(n,r)=(34,10)$.}
Here $d=26$. For an integral rational quartic of span dimension at
most three one can slice to dimension $j=4$ in the ordinary jet
formula. Otherwise the curve is a rational normal quartic in its
span $\mathbf P^4$; slice to a hypersurface of dimension $j=5$ in
$\mathbf P^6$. The twisted conormal of the curve in $\mathbf P^6$ is
$\mathcal O_{\mathbf P^1}(-2)^{\oplus3}\oplus\mathcal O_{\mathbf P^1}^{\oplus2}$
by the classical normal-bundle formula for a rational normal curve and
the split tangent sequence of a linear subspace. The twisted conormal
of the curve in the hypersurface is its vector-bundle quotient. Thus
each summand has degree $a_i\ge-2$, and
$\sum_i a_i=(d-1)e-(2e-2)$ with $e=4$. The previous signed-slot formula
has exact leading capacity $[e x_0+\sum_i a_ix_i]_+$.
For a trial $t>0$, if $a_1 t\ge2(1-t)$, this linear form is nonnegative
throughout each tail $x_1\ge u\ge t$. Its integral over that tail is
the previous upper bound minus at least
$(2e-2)(1-u)^j/d$ after normalization. If $a_1t<2(1-t)$, use positive
parts as before. In the positive-parts tail bound the loss in the
first-vertex term also gives this improvement whenever
\[
       (d-1)e\,t-2(1-t)\ge(2e-2)(1-t)/j.
\]
The trial in the checker satisfies this inequality. Integrating in $u$
gives a reduction $(2e-2)(1-t)^{j+1}/d$ in the normalized angular
factor $(j+1)\Theta_{j,d}(e,t)$; the independent ordinary term is
kept. The checker takes the larger of this result and the
dimension-four estimate for the smaller span. This closes the boundary
case using retained fraction $0.85$ and $k=0.20$.

\section{Twice the index}

\paragraph{A uniform cutoff of $50$.}
Use the surface-center bound $e\le(n+1)/\beta'$ and the
averaged-quantile gain. With $R=50$, $r\ge R$, and $x=n/r\in[a,b_1]$,
the following rows prove $\delta(X,H)>r$ throughout $2\le x\le4$.
\[
\begin{array}{c|c|c@{\qquad}c|c|c}
{}[a,b_1]&f&L&[a,b_1]&f&L\\\hline
{}[2,65/32]&8/r&23/250 &[65/32,33/16]&8/r&23/250\\
{}[33/16,17/8]&9/r&9/100 &[17/8,9/4]&11/r&43/500\\
{}[9/4,5/2]&14/r&41/500 &[5/2,11/4]&13/r&21/250\\
{}[11/4,3]&5/r&59/500 &[3,7/2]&13/20&11/250\\
{}[7/2,29/8]&20/r&23/500 &[29/8,15/4]&10/r&7/125\\
{}[15/4,4]&4/r&27/250 &&&
\end{array}
\]
Here $L\beta$ multiplies the upper-fraction basis and the finite-level
upper tail allowance in $h$ is $10^{-6}/r$. To explain the interval
checks without summing over dimensions, write
\[
 b_f=\frac{f\log(1/f)}{1-f},\quad
 g=\left(1-\frac r{n+1}\right)\frac f{1-f}
               \min\left(1,\frac{\log(1/f)}{\log d}\right),
 \quad h=L[1-\log(1-f)]+10^{-6}/r.
\]
Use rational gains arbitrarily close to the indicated lower bound.
The tilt and class subtraction give $\beta'\le[1-L(1+g)]\beta$ and
$B\sim_{\mathbf Q}\lambda H$ with
$\lambda\le L\beta[1-(1-b_f)/(n+1)]$. Consequently the bounds on
$\beta'+h\beta$ and on $r-\lambda-n\beta'/(n+1)$ are respectively
\[
 U=r+Lr[-\log(1-f)-g]+10^{-6},\qquad
 \eta_0=\frac1{x+1/r}+L\left[\frac n{n+1}rg-
                           \frac{b_f}{x+1/r}\right].
\]
The former is an upper bound, the latter a lower bound.

For $f=F/r$ in the table we have $f\ge1/d$. As $r$ increases,
$b_f$ decreases and $-r\log(1-F/r)$ decreases. Indeed
$\log t/(t-1)$ decreases for $t>1$ and $-\log(1-t)/t$ increases
for $0<t<1$. The expression
\[
 (1-1/a)\frac{F\log(r/F)}{(1-F/r)\log[(b_1-1+2/R)r]}
\]
is a lower bound for $rg$, increasing on $r\ge R$ in every row.
To check the last assertion put $u=\log(r/F)$ and
$k=\log[F(b_1-1+2/R)]$. The sign of the derivative of
$u/[(1-e^{-u})(u+k)]$ is that of $k(e^u-1-u)-u^2$.
The checker verifies the nonnegative sign at $R$ with outward rational
logs. Since $(e^u-1-u)/u^2$ increases, the sign persists. For the one
constant-$f$ row use the lower bound
\[
 (1-1/a) r\min\{f/(1-f),\ b_f/\log[(b_1-1+2/R)r]\},
\]
which is likewise increasing. In that row one may discard $g$ in the
upper bound for $U/n$. These monotonicities give rational bounds
$u_0\ge U/n$, $z\le\eta_0$ and $c\ge(1-h)^{-1}$ on each whole row.

For degrees one and two check $U<(n+1)d/(2d-1)$. For variable $f$
use the lower bound $(a-2)R/2+1/2+a/[4(b_1-1)+6/R]$ for the excess of
this mark over $r$. For higher degrees use $j=e+1$ and the tabulated
span-jet trial from the exact checker up to $e=40$, with $8/3$ for
$e=3$. For positive center dimension $s\le2$ the cap for
$(\beta'+h\beta)/(n+1)$ is $\min(u_0,c/e)$. The same cap applies
in any dimension when $(s+1)\beta'\ge n+1$. For $3\le s\le8$,
the other case can occur only when
$z\le s+1$ and $1-L-(s+1)/R<b_1/(s+1)$, by the center inequality;
then use $\min(u_0,c/(s+1),c(s+1)/(se))$. The remaining dimensions
have corresponding weak cap $\min(u_0,c/10,10c/(9e))$. Each cap
times its reciprocal-threshold coefficient is less than one. For
$e\ge41$ it suffices to check $(3/2)c(4/41+6/54)<1$ and $U<n-1$,
by the full-dimensional trial $4/n$.

A point optimizer is excluded immediately if $z>1$. Otherwise take
$\Delta=(1-z+10^{-5})u_0/(1-u_0)$: the normalized valuation has
$S_H(w)\ge n/U$ and hence
$\Delta S_H(w)+r-\lambda-n(\beta'+\Delta)/(n+1)>1$.
The checker verifies $z-\Delta>0$ and the curve tests for
$U+\Delta$ through degree $\lfloor1/(z-\Delta)\rfloor$.
The short-curve and Nadel argument therefore excludes the point.
The deterministic command \texttt{python r2-flag-index-boxes.py}
checks every row with rational arithmetic.

\section[The adjacent row at large index]{The adjacent row $n=2r-1$ at large index}
For every integer $r\ge10^{10}$, every smooth degree-$(r+1)$
hypersurface in $\mathbf P^{2r}$ is uniformly K-stable.
Write $n=2r-1$, $d=r+1$, $H=\mathcal O_X(1)$, and suppose
$\beta=\delta(X,H)\le r$. The theorem of Liu--Xu--Zhuang gives a
divisorial optimizer. We use the center-dimension, short-curve,
canonical-extension and span-sliced estimates established above.

\paragraph{The case $\beta\le r-1$.}
Use the original pair with zero boundary. A minimizing center of
dimension $s$ satisfies $r-n\beta/(n+1)\le s+1$. Since
\[
 r-\frac n{n+1}\beta\ \ge\ 1+\frac{r-1}{2r}>1,
\]
the optimizer has positive center. If $s=1$, the curve-image bound
gives a rational curve of degree $e\le(n+1)/\beta$. If
$2\le s\le8$, the dimension bound gives
$\beta\ge(n+1)(r-9)/n$, and consequently $(s+1)\beta>n+1$.
The canonical-extension bound gives the same degree estimate in
these dimensions. For $s\ge9$ the ordinary quotient slope and
Jovinelly--Lehmann--Riedl give $e\le10n/(9\beta)$.
Thus in all cases $e\beta/(n+1)\le10/9$.

Such a curve cannot lie in an optimizing center. For $e=1,2$ the
curve threshold is at least
\[
   \frac{(n+1)d}{2d-1}=r+\frac{r}{2r+1}>r-1.
\]
For $e=3$ the upper bound on $\beta S_H(w)/A_X(w)$ is at most
$(10/9)(8/9)<1$ by the plane and twisted-cubic estimates.
For $e=4$ take $j=5$ and $t=454/1000$ in the span-sliced normal
estimate. Its numerator is at most
\[
  6\max\{1/4,\ t+4(1-t)^5[t+(1-t)/3]\}<3.465.
\]
The corresponding product is at most $(10/9)(3.465/4)<1$.
For $e\ge5$, the full-dimensional trial $t=3/n$ gives
$n\Theta(e,t)\le3+5e\exp(-3)$. Its contribution is bounded above by
\[
    (1+1/n)(10/9)(3/5+5\exp(-3))<1.
\]
The other full-dimensional
term satisfies $\beta/(n-1)<1$. All possibilities contradict
the optimizer equality.

\paragraph{A boundary from the complete projection subseries.}
Now suppose $r-1<\beta\le r$. Put $y=\sqrt r$, $f=1/d$,
and choose rationals $l,\Delta$ satisfying
\[
 \frac13-\frac{0.02}y-\frac1{r^2}<l\le\frac13-\frac{0.02}y,
 \qquad
 \frac13-\frac{0.04}y-\frac1{r^2}<\Delta\le\frac13-\frac{0.04}y.
\]
Set $\lambda=lr$. On a general finite linear projection
$X\to\mathbf P^n$ that is \'etale at the generic center of the
optimizer $v$, the restricted valuation has discrepancy $A_X(v)$.
The polynomial subseries has limiting dimension fraction $f$ and
normalized mean at most $A_X(v)/(n+1)\le S_H(v)/2$. The full bottom
$f$ fraction has mean no larger. Indeed in each degree the sum of
the smallest $\dim U$ full weights is at most the sum for an adapted
basis of any subspace $U$ of that dimension, by the filtration
quotient rank inequalities. Replacing this number by the desired
rounded fraction changes $o(\dim H^0(X,mH))$ slots, whose normalized
weights are uniformly bounded for the fixed $v$.

In a basis compatible with $v$ and a hyperplane $H_0$ avoiding its
generic center, let $D^+$ be the normalized upper fraction $1-f$.
Remove its $H_0$ copies to obtain $D\sim_{\mathbf Q}c H$ with
\[
 c<1-\alpha,\qquad
 \alpha=\frac{1-\log(r+1)/r}{2r},
\]
by the free-hyperplane computation. Its $v$-value is unchanged.
The limiting mean of the retained fraction is at least
$(1+1/(2r)) S_H(v)$ since $f/(1-f)=1/r$. For the fixed rational
boundary $B=\lambda D$ at a sufficiently large level, put
$k=\lambda+l+1/r$ and $\beta'=\delta(X,B;H)$. Uniform complete-
series comparison on a basis extending the retained vectors gives
$w(B)\le kS_H(w)$ for all $w$, since $\lambda/(1-f)=\lambda+l$.
Testing $v$ and using $\beta>k$ gives
\[
  \beta-k\le\beta'\le\beta-\lambda-l/2+1/r,
  \qquad U:=\beta'+k\le r+l/2+2/r<r+0.17.
\]
Thus $(X,B)$ is klt, and it has a divisorial optimizer $w$ since
$\beta'<r-\lambda c$. In particular $A_X(w)\le US_H(w)$.
Also $\beta'>2r/3+0.02y-1-l-1/r>2r/3$. A positive center of
dimension one therefore yields $e\le(n+1)/\beta'<3$, and any
center of dimension $s\ge2$ satisfies $(s+1)\beta'>n+1$, giving
the same bound by the canonical-extension argument. The line and
conic threshold $r+r/(2r+1)>U$ rules out a positive center.

\paragraph{The possible point center.}
Normalize $w(\mathfrak m_p)=1$ and write $\mu=A_X(w)=n+q$.
The universal point envelope $S_H(w)\le J+Gq/n$ has
$J\le1.501$ and $G\le1.42y$ here. In the high branch the
intercept is $(3r-2+\sqrt{(r+1)/2})/(2r)$ or smaller, and the
low intercept is $1+O(1/y)$, explicitly at most
$(3(\sqrt{2d}+4)+n-3)/(n+1)$. For the slopes use
$\sqrt{2d}+4<1.42y$ and
$\sqrt{d/2}+\log(3n(d+2))<1.42y$ for $r\ge10^{10}$; the latter
follows from $\log(3n(d+2))<\log(12r^2)$ and monotonicity of
$\log(12r^2)/\sqrt r$. Since $A_X(w)\le U S_H(w)$, the case
$UG\le n$ is impossible by $UJ<n$. Otherwise
\[
 q\ge\frac{n-UJ}{UG/n-1}>0.6y,
 \qquad \frac{n+q}{U}>2+0.59/y.
\]
For the first inequality the numerator is at least $0.498r$ and
the denominator at most $0.711y$; the second follows from
$U<r+0.17$ and $0.01y>1.35$.

For bases adapted to $w$ and to a hyperplane avoiding $p$, let
$\Gamma_m$ be the usual residual, of class $(1-a_m)H$ with
$a_m\to1/(n+1)$ and $w(\Gamma_m)\to S_H(w)$. The limiting ample
slack at $\beta'$ satisfies
\[
 r-\lambda c-\frac n{n+1}\beta'
 \ \ge\ \frac12+\frac l2-
     \frac{1+l/4+(l/2)\log(r+1)}r
 \ >\ \eta_0:=\frac12+\frac l2-0.001/y.
\]
For the last bound use that
$(1.084+\log(r+1)/6)/\sqrt r$ decreases for $r\ge10^{10}$ and
is already below $0.001$ there. At $t=\beta'+\Delta$ this leaves
the limiting lower slack
\[
   \eta=\eta_0-\frac n{n+1}\Delta>1/3+0.028/y.
\]
If $(X,B+t\Gamma_m)$ has a positive non-klt center for infinitely
many $m$, take the first positive threshold $t_m\le t$. The
Blum--Jonsson uniform comparison gives $t_m\ge\beta'/(1+o(1))$.
The short-curve argument (including the arbitrarily small general
pointed hyperplane if needed at $p$) gives a rational curve of degree
at most $\lfloor1/\eta\rfloor\le2$ in a positive center attaining
$t_m$. For its place $u$, the coefficient comparison uses
$u(B)+t_m u(\Gamma_m)\le(k+t+o(1))S_H(u)$. This is strictly
smaller than $A_X(u)$, since the curve threshold is at least
$r+r/(2r+1)$, whereas
\[
      k+t=U+\Delta\le r+1/2-0.05/y+2/r
                   <r+1/2-1/(4r+2).
\]
Thus $B+t\Gamma_m$ is klt in a punctured neighborhood of $p$ for
all large $m$. Add a general hyperplane through $p$ with rational
coefficient tending to $\eta$ from below; the strict ample slack
and $\eta<1$ preserve the needed ampleness and punctured klt property.
Abban--Zhuang's point Nadel inequality then implies
\[
  1\ge\Delta S_H(w)+\eta
    \ge \eta_0+\Delta\left(\frac{n+q}{U}-\frac n{n+1}\right)>1.
\]
For clarity, the strict lower bound follows already from
$\eta_0>2/3-0.012/y$, $\Delta>1/3-0.041/y$ and
$(n+q)/U-n/(n+1)>1+0.59/y$: their combination exceeds
$1+0.14/y$ for $y\ge10^5$. This contradiction finishes the row.

\subsection{Lowering the adjacent-row cutoff}
Here we improve the cutoff in the adjacent-row argument. On $n=2r-1$ the preceding untilted proof for $\beta\le r-1$ works already for $r\ge65$. Indeed $3(r-9)>2r-1$; the remaining four degree comparisons hold at $n\ge129$. Suppose $r-1<\beta\le r$ and take $\epsilon=r^{-4}$ and
\[
 \ell=(\beta-2r/3-3\epsilon)/(r+1),\quad
 \lambda=r\ell,\quad u=\frac n{n+1},\quad k=\lambda+\ell+\epsilon .
\]
The original optimizer is dreamy, hence $\beta$ is rational. Take the same upper $1-1/d$ fraction and free hyperplane subtraction as in the adjacent-row proof. Then
\[
  \beta'>2r/3,\qquad \beta'\le\beta-r\ell-\ell/2+\epsilon,\qquad
  \beta'+k\le\beta+\ell/2+2\epsilon .
\]
In the following formulas evaluate $\ell$ at $\beta=r$, calling this value $l$. The right endpoint maximizes the discrepancy budget, now denoted $U=r+l/2+2\epsilon$, and minimizes the following lower bound on the ample slack:
\[
 z=\tfrac12+l g-u\epsilon,\qquad
 g=\tfrac12\big(u-\log(r+1)/r\big).
\]
Indeed before evaluation this lower bound is $r-u\beta+\ell g-u\epsilon$, a decreasing affine function of $\beta$. A positive center yields $e\le2$ by the canonical extension for $s\ge2$ and the curve bound for $s=1$. The line/conic threshold $r+r/(2r+1)$ exceeds $U$.

For a point, seek $\Delta>0$ satisfying
\[
   \Delta<\frac r{2r+1}-(U-r),\qquad z-u\Delta>1/3,\qquad
   z+\Delta\left(\frac{n+Q}{U}-u\right)>1 .
\]
The first two inequalities exclude positive centers at $\beta'+\Delta$ with the line and conic estimate. The third is precisely the contradiction to Nadel at the point optimizer. For $65\le r<1024$ the program \path{r2-root-adjacent.py} evaluates both point-envelope constants, chooses the smaller of the two upper bounds on $\Delta$ minus $\epsilon$, and verifies the last inequality strictly with rational arithmetic.

For $r\ge1024$, put $y=\sqrt r$ and $E=0.26+\log(r+1)/6$. Here $U<r+0.17$, $z\ge 2/3-E/r$. Indeed $1/3-l=(1/3+3\epsilon)/(r+1)$ and $1/2-g=[1/2+\log(r+1)]/(2r)$. The simpler universal envelope gives $J\le1.512$ and $G\le1.55y$. For the low intercept use $(3(\sqrt{2(r+1)}+4)+n-3)/(n+1)$; the high intercept is $(3r-2+\sqrt{(r+1)/2})/(2r)$. The slope follows from $\sqrt{2(r+1)}+4\le1.55y$ and $\sqrt{(r+1)/2}+\log(12r^2)\le1.55y$; the ratios of the logarithms to $y$ decrease here. Thus $n-UJ\ge0.486r$ and $UG/n-1<0.777y$. Unless the envelope already excludes the point, $Q>0.62y$, and $(n+Q)/U>2+0.57/y$, using $0.05y>1.35$. Take $\Delta=1/3-(E+0.1)/r$, with a sufficiently close rational perturbation if needed. The first two tests are strict since $r/(2r+1)>1/2-1/(4r)$, $l<1/3$, and $E+0.1>0.36$. For the third it suffices that
\[
 0.19y>2E+0.1+0.57(E+0.1)/y .
\]
This holds at $y=32$ using $\log(1025)<7$ and persists since $E/y$ decreases. In particular the adjacent row holds for every $r\ge65$.

\subsection{The adjacent row at the remaining indices}
Write $n=2r-1$, so $d=r+1$. For $\beta\le r-1$ the center-dimension inequality rules out a point. Every positive optimizer with $s\le2$, or with $(s+1)\beta\ge n+1$, has a rational curve of degree $e\le(n+1)/\beta$. Otherwise $s\ge3$ and $\beta/(n+1)<1/4$. The sliced bounds for $e=3,4,5,6$ are respectively at most $8/3$, $3.465$, $3.721$, and $3.939$, using rational trial values from the uniform interval checker. Each factor is less than both $4$ and $e$. These choices use $j=e+1\le n$ for $r\ge4$. In degree one or two the mark is $r+r/(2r+1)>r$. In every case $e\beta\le4n/3$ and for $e\ge7$ the full-dimensional trial $3/n$ gives at most $(4/3)(3/7+5\exp(-3))<1$. The radial term is strictly below one since $\beta<n-1$. These comparisons give the untilted branch at every $r\ge4$.

In the tilted branch the preceding formulas for $U,z,\Delta$ apply unchanged. For $4\le r\le64$ the finite angular estimate below bounds $q$ from below at a hypothetical optimizer with $A_X/S_H\le U$. It suffices to prove
\[
                  Q> U\left(\frac{1-z}{\Delta}+\frac n{n+1}\right)-n .
\]
The following calculation derives the low- and high-multiplicity bounds that determine $Q$ and verifies the strict comparison for each of these indices. Together with the large-index estimate this establishes the adjacent row for every $r\ge4$. For $r=3$ it is the quartic fivefold of Chapter~\ref{ch:fivefold}.

\paragraph{A finite angular estimate on the adjacent row.}
In the effective adjacent-row argument the point test holds throughout
$4\le r\le64$, with $n=2r-1$ and $d=r+1$. Write $u=n/(n+1)$, and
use $\epsilon,l,U,g,z,\Delta$ of that argument, with
\[
 \Delta=\min\{r/(2r+1)-(U-r),(z-1/3)/u\}-\epsilon.
\]
Here $g$ and $z$ may be replaced by lower rational bounds, as in the
checker. Put $Q_{\rm req}=U((1-z)/\Delta+u)-n$ and round it up
strictly to the next multiple of $1/100$ for $4\le r<12$ or of
$1/4$ for $12\le r\le64$. We show that any point-centered divisor
with $A_X(v)/S_H(v)\le U$ and $v(\mathfrak m_p)=1$ has
$q=A_X(v)-n>Q$. Then $z+\Delta((n+q)/U-u)>1$ is the required
contradiction. The proof of the finite test follows.

Choose a rational $B$ with $B^2>d$ and
$(B^2/d-1)(2d-5)+B-d\ge0$. If $\tau_p(H)\le B$, the signed pointed
linear flag gives
\[
                  (n+1)S_H(v)\le B(q+2)+n-2 .
\]
For the improved surface constant, consider an integral curve of
degree $a$ and multiplicity $k$ on a very general pointed
linear surface in a smooth threefold section $Y$. The proper
swept-locus case is the Picard-free point flag: either the curve
comes from a prime on $Y$ with $dk/a\le\tau_p(H|_Y)\le B$, or the
tangent hyperplane section is a cone, forcing $B\ge d$.
In the covering case the canonical determinant on the point blowup
of $Y$, adjunction on the blown-up surface, and Hodge index give
\[
  2g-2\ge(d-5)a+2k,\qquad k(k+1)\le a^2/d+a .
\]
If the curves are planar, $k\ge2$ forces their planes to be the
tangent planes of the corresponding surfaces at $p$. All then lie
in the fixed hyperplane $T_pY$, contrary to coverage. Thus $k=1$
and the plane genus estimate gives $a\ge d-1$, which implies
$dk/a<B$. If nonplanar, Castelnuovo's bound
$g\le\lfloor(a-2)^2/4\rfloor$ gives $a>2d-6$, hence $a\ge 2d-5$.
The inequality on $B$ now implies
$(Ba/d)(Ba/d+1)\ge a^2/d+a$, so $dk/a\le B$.
This proves $\varepsilon_p(H|_S)\ge d/B$ and gives the displayed
signed bound as in the Picard-free point flag. The two genus
inequalities apply already for $d\ge4$.

If $\tau_p(H)>B$, there is a prime $P\in|hH|$ of ratio
$b=\operatorname{mult}_p(P)/h>B$; any other prime has ratio at most
$d/b$ by proper intersection, so $\tau_p(H)=b\le d$.
Use the principal-bin angular estimate with $K$ bins in the
class fraction $h\operatorname{ord}_P/m$. Write
\[
 w_i=(1-i/K)^n-(1-(i+1)/K)^n,\quad
 A=\sum_{i=0}^{K-1}(i/K)w_i,\quad R=\sum_{i=0}^{K-1}w_i(dw_i)^{1/n}.
\]
For $b$ in a closed interval $[b_-,b_+]\subseteq[B,d]$, put
\[
 M=\max_{b_-,b_+}\left\{\frac{n-2+b+d/b}{n+1}-bA\right\}.
\]
Convexity makes this an upper bound for the ordinary mean after the
fixed powers have been removed. If $0<c<1$ and
$R(c^{1-n}+(n-1)c)\ge(n+1)M$, the principal-bin trial tents give
\[
 S_H(v)\le(1+q)M+(d+q)A-q\frac{n-1}{n+1}Rc .
\]
Indeed the residual grade density in bin $i$ has mass $w_i$ and is
at most $ny^{n-1}/d$. The trial density is $n/d$ times the
$(n-1)$st power of a tent rising with slope one to
$c(dw_i)^{1/n}$ and then descending to zero at
$c^{1-n}(dw_i)^{1/n}$. It has the same mass; the aggregate trial
mean is $R(c^{1-n}+(n-1)c)/(n+1)$. Convexity of
$x^{1+1/(n-1)}$ and the common slope of the descending sides show
that its angular moment $(n-1)Rc/(n+1)$ is a lower bound for the
actual angular moment: replace each tent in the first-variation
integral by the affine continuation of its descending side. The
jet cap on the rising side and nonnegativity beyond the support
give the required sign. Integrating those affine continuations
uses just the bin masses and the upper bound for the actual
aggregate mean. Finally the removed principal class has mean $A$
and $v(P)/h\le d+q$, giving the displayed estimate.

The program \path{r2-local-adjacent-exact.py} in the archive described in Chapter~\ref{ch:checks}
checks both affine comparisons at $q=0,Q$. It rounds each radial
root down to $10^{-5}$, chooses $B$ just below the low-branch
endpoint, and takes $K=3n$ on $[B,d]$ whenever possible. This
single interval works for every $12\le r\le64$. For smaller $r$
the code can subdivide by intervals of width $1/8$ and
independently choose $K,c$ on each interval. The calculation uses
rational arithmetic throughout, including an upper bound for
$\log(d)$ from the positive atanh series. The resulting point
conclusion holds for all 61 indices, including $(r,n,d)=(4,7,5)$
and $(5,9,6)$.

\section{Finite completion of the dimension range}

\paragraph{The remaining strip down to $n=2r-1$.}
On the linked smooth surface associated to a twisted cubic, put
$H^2=4d$, $H C=3$ and $C^2=4-3d$. We sharpen the previous estimate
of the expected vanishing by integrating its volume. If the residual
secant line $L$ is absent, the class $H-C$ is nef and hence $H-tC$
is nef for $0\le t\le1$. When $L$ is present,
$HL=1$, $CL=2$, $L^2=-d$ by adjunction, and $H-C-L$ is nef.
For $0\le t\le1$ the Zariski positive part is
\[
    P_t=H-tC-\frac{(2t-1)_+}{d}L.
\]
Indeed it intersects $L$ nonnegatively, with equality whenever the
subtracted multiple is positive. On other curves it is nonnegative
because it is a nonnegative combination of $H$, $H-C-L$ and $L$.
Thus the volume on this first interval is
\[
  4d-6t+(4-3d)t^2+\frac{(2t-1)_+^2}{d}
\]
when $L$ is present, with the last term omitted otherwise. Its integral
is at most $3d-5/3+1/(6d)$.

Write $V=d-2+1/d$ and $D=3d-1-2/d$. At $t=1$ the volume is $V$ and
$P_1C=D$ in the secant case. In the other case the volume is smaller
and the intersection with $C$ is larger. The square root of volume is
concave on the big cone. Its left derivative at $1$ is
$-P_1 C/\sqrt{P_1^2}$ by the displayed formula. The tangent line
gives a tail integral at most $V^2/(3D)$. Including the
surface-to-ambient factor three, this proves
\[
 \Xi_3(d)=\frac{3}{4d}\left(3d-\frac53+\frac{1}{6d}
                  +\frac{(d-2+1/d)^2}{3(3d-1-2/d)}\right)<\frac73.
\]
For plane cubics retain $2-1/d$, which is smaller than this expression
for $d\ge4$. The surface and intermediate flags are the same ones as
in the previous twisted-cubic calculation.

Apply the preceding tilted test with this $\Xi_3(d)$,
the strong degree bound in center dimensions at most three, and the
canonical extension bound in the applicable remaining dimensions. If
$\beta<(r-1)(n+1)/n$, the point optimizer is excluded by the dimension
bound. The untilted positive optimizer is excluded in the present strip,
as follows. For degree one or two,
$\beta(2-1/d)/(n+1)<2(r-1)/n\le1$ for $n\ge2r-1$.
In degree three the positive cap $M/(n+1)\le1/3$ gives a product
less than $7/9$. The degree-four to degree-six caps are exactly
those of the preceding large-index table with $h=0$, and the
full-dimensional degree-at-least-seven estimate also applies.
If $j=e+1>n$ in the degree-four to degree-six slice bound for $n=5$,
the untilted case $\beta<(r-1)(n+1)/n$ and the full-dimensional
$j=n$ estimate can instead be checked directly with the same rational
trials as the finite certificate.

For the high-$\beta$ branch, use the previously proved fixed-bin
point estimates or, when the last column of \path{r2-lit-final-strip-witnesses.txt} equals $1$,
the intervalwise principal-bin estimate from the angular proof for
the adjacent row. That proof works here for all $d\ge4$ and $n\ge5$.
In its pointed surface estimate the nonplanar Castelnuovo bound gives
$(d-5)a+2k\le a(a-4)/2$, so $a>2d-6$ also at $d=4$.
For the planar $k=1$ case the plane genus bound gives
$a(a-d+2)\ge2$, hence $a\ge d-1$. The remainder of the pointed
flag and bin estimates use $n\ge4$. The resulting lower bound $Q$ can
be combined with the fixed-bin bound by taking their maximum.

For $3\le r\le31$, $2r-1\le n<\min(17r/5,r(r-1))$, the program
\path{prior-range/checks/r2-lit-final-strip.py} verifies the 716 records of this table with exact arithmetic, including
$(r,n)=(3,5),(4,7),(4,8),(4,9),(5,9)$. It also verifies that the
complement comprises exactly the eight adjacent pairs $n=2r-1$
with $r\in\{9,17,18,24,25,26,27,28\}$. These eight are covered by
the separate exact adjacent-row argument. For example,
$(r,n)=(4,8)$ uses omitted fraction $1/4$ and tilt $3/50$.
Together with the results for $n\ge17r/5$, $n\ge r(r-1)$ and the
already established $r\ge32$ strip, this proves the claimed
full range $n\ge2r-1$.

\section{A range determined by two rational inequalities}
Write $r=n+2-d$, and assume $4\le d\le n-1$ and $r\ge3$. For $e\ge1$ set
\[
 \Theta(e,t)=t+\frac ed(1-t)^n
   \left((d-1)t+\frac{2d-1}{n+1}(1-t)\right),\qquad
 \Theta(e)=\min_{0\le t\le1}\Theta(e,t).
\]
The minima only enter through upper bounds. Any two rational trial values in $[0,1]$ therefore give an exact sufficient test.
Put
\begin{align}
  e_0&=\left\lfloor\frac{2(n^2-1)}{r}\right\rfloor,&
  \eta&=\frac{r}{nd+1},&
  e_1&=\lfloor1/\eta\rfloor,&
  T&=\frac{n+1}{n}(r-\eta).
  \label{global:eq:constants}
\end{align}
Then every smooth $X_d^n$ is uniformly K-stable whenever
\begin{equation}\label{global:eq:test}
  2r(n+d-1)<(n+1)^2,\qquad
  r\Theta(e_0)<1,\qquad T\Theta(e_1)<1 .
\end{equation}
In this range $r<T<n-1$, since $r\le n-2$. In particular the codimension term in the normal-jet estimate exceeds both test coefficients.

\paragraph{Positive-dimensional centers.}
For every integral rational curve $C$ of $H$-degree $e$, the normal-jet calculation gives $\delta_C(X,H)\ge\min(n-1,1/\Theta(e))$, where valuations with center containing $C$ are included. The bound decreases with $e$. Take a fixed $T_0>r$ sufficiently close to $r$ that $T_0 n/(n+1)<r$ and the curve bound for $e\le e_0$ exceeds $T_0$. We may also require
\[
        \frac{2(n-1)}{r-T_0n/(n+1)}<e_0+1 .
\]
This is possible by strictness and the definition of the floor, including when $2(n^2-1)/r$ is integral.

Let $H_0$ be a general hyperplane. From any basis divisor $D_m$ of $|mH|$, normalized by $m h^0(mH)$ and compatible with order along $H_0$, subtract its $H_0$-coefficient $a_m$ to obtain $\Gamma_m=D_m-a_m H_0$. One has $a_m\to1/(n+1)$ by summing $h^0((m-i)H)$ over $i\ge1$. We claim that for all sufficiently large $m$, $(X,T_0\Gamma_m)$ is klt at every positive-dimensional generic point. Otherwise the first threshold $t_m$ over positive-dimensional centers is at most $T_0$. The quasi-log cone theorem gives a rational curve $C_m$ inside a center attaining $t_m$, with
\[
     H\cdot C_m\le 2(n-1)/[r-t_m(1-a_m)]<e_0+1 .
\]
The uniform valuation estimate for the complete series, $S_m(H;w)\le(1+o(1))S_H(w)$, holds uniformly in $w$ \cite[Corollary 3.6]{ar1706.04548}. Thus $w(\Gamma_m)\le S_m(H;w)$ and the strict curve bound exclude this threshold for large $m$. This also proves $\delta_Z(X,H)\ge T_0$ for every $Z$ of positive dimension: for a given $w$ choose $H_0$ avoiding its generic center and a basis compatible with both $w$ and $H_0$. Then $w(\Gamma_m)=S_m(H;w)$.

\paragraph{Closed points.}
Assume that $\delta(X,-K_X)\le1$. A divisorial valuation $v$ computing it exists by \cite[Theorem 1.2]{ar2102.09405}. The preceding step forces its center to be a point $p$. Normalize $v(\mathfrak m_p)=1$, and write $A_X(v)=n+q$, $q\ge0$. Fix $H_0$ avoiding $p$ and use the same compatible bases and residual divisors. Thus $rS_H(v)\ge A_X(v)$ and $v(\Gamma_m)\to S_H(v)$.

If the first positive-center threshold $t_m$ is at most $T$ for infinitely many $m$, for those $m$ it is at least $T_0/(1+o(1))$ by the uniform comparison and the positive-center bound. In particular $(X,t_m\Gamma_m)$ is non-lc at $p$ for all large such $m$. Its non-lc scheme is finite. The short rational curve theorem above gives a rational curve in a positive-dimensional center attaining $t_m$, with degree at most
\[
                1/[r-t_m(1-a_m)]\le 1/(\eta+o(1))<e_1+1.
\]
The strict test $T\Theta(e_1)<1$ and the same complete-series comparison give a contradiction. Thus $(X,T\Gamma_m)$ is klt off a finite set for large $m$.

For $q>1$ use an auxiliary general $H'_m\in|H\otimes\mathfrak m_p|$. Choose $0\le\lambda_m<1$ tending to $\eta$ and satisfying $\lambda_m<r-T(1-a_m)$ so that $-(K_X+T\Gamma_m+\lambda_m H'_m)$ is ample. Bertini for klt pairs preserves klt on a punctured neighborhood of $p$, since the hyperplane system is basepoint free there. Abban--Zhuang's Nadel vanishing lemma \cite[Lemma 4.24]{ar2003.13788} with $\mu=n+q$ gives in the limit
\[
        A_X(v)\ \ge\ \frac{\mu}{\mu+1}
                           \big(T S_H(v)+\eta\big).
\]
Consequently $(T/r-1)(n+q)+\eta\le1$. Direct substitution in \eqref{global:eq:constants} gives $(T/r-1)(n+1)+\eta=1$ and $T>r$, excluding $q>1$.

For $0\le q\le1$ there is an elementary estimate. On the blowup of $p$ the hyperplanes through $p$ generate $H-E$. For $0\le z\le1$ the volume of $H-zE$ is $d-z^n$. Concavity of its $n$th root bounds the volume beyond one by $(d-z)^n/(d-1)^{n-1}$ up to $z=d$, and the nef intersection with $(H-E)^{n-1}$ excludes any larger big value. Integration gives
\[
       s=S_H(\mathrm{ord}_p)\le \frac{n+d-1}{n+1}.
\]
For any divisor germ of ordinary order $j$ at $p$, its strict transform meets $E=\mathbf P^{n-1}$ in a degree-$j$ divisor. Its log canonical threshold there is at least $1/j$. Inversion of adjunction with boundary $E$ therefore gives $v(f)\le j+jq$. The bound includes a unit directly. Applying it to an adapted basis yields $S_H(v)\le(1+q)s$. At $q=1$ the first inequality of \eqref{global:eq:test} gives $r(1+q)s<n+q$. At $q=0$ the numerator of the margin is $(d-1)(d-2)>0$. Affinity in $q$ proves strictness throughout. This excludes the last possible computing valuation and proves the criterion.

\paragraph{Closing the quadratic range independently.}
We use the preceding two-stage test to cover the unbounded region needed by the linear-range proof. Put $x=n/r$. If $x\ge4$ and $e^x\ge 2e^2rx^3$ (here $e$ is Euler's number), both degree caps in \eqref{global:eq:constants} are at most $2rx^2$. Taking $t=(x-2)/n$ gives $(1-t)^n\le\exp(2-x)$ and
\[
 n\!\left[(1-1/d)t+(2-1/d)(1-t)/(n+1)\right]\le x .
\]
Thus $n\Theta\le x-2+2\exp(2)rx^3\exp(-x)\le x-1$ in both tests. Also $T<r(1+1/n)$, so $T\Theta<1$. The remaining inequality $2r(n+d-1)<(n+1)^2$ holds when $n\ge4r$.

For $r\ge15$ it suffices that $n\ge r(r-1)$. Indeed $x\ge r-1$, $e^x/x^3$ increases for $x>3$, and $e^{r-1}/[r(r-1)^3]$ increases for $r\ge15$. At $r=15$ the required inequality is $e^{12}\ge 30\cdot14^3$, which follows already from $e>2.7$. For $3\le r\le14$ the same exponential test holds whenever $x\ge16$, by checking $x=16,r=14$. Between $r(r-1)$ and $16r-1$ the 724 remaining pairs at these indices pass the finite-tilt test in the preceding section. The program \path{r2-root-independent-check.py} verifies the inequalities in \path{r2-root-independent.txt} and checks that its records cover exactly these pairs. Thus the quadratic portion of the proof can be recovered entirely from the curve and point arguments used here.

\section{Coverage and the next dimension}
For $r\ge3$, combine the following overlapping proofs. The bound $n\ge r(r-1)$ follows independently from the preceding section.
The six-index, four-index and fifteen-quarters arguments cover
$4n\ge15r$. The terminal threefold and linked-surface calculations
extend this to $5n\ge17r$. For $3\le r\le31$ the finite strip proof and
the adjacent-row estimates cover every remaining integer $n\ge2r-1$.
For $32\le r<50$, the finite interval covers give $2r\le n<15r/4$.
For $r\ge50$, the uniform interval proof covers $2r\le n\le4r$.
The adjacent-row calculation covers $n=2r-1$ at every such index.
Finally the published cases $r=1,2$ complete the stated dimension range.
Chapter~\ref{ch:checks} gives the Python commands for verifying both
the interval covers and the pointwise inequalities.

\paragraph{An alternative check using surface centers.}
As an additional check on the use of the terminal threefold argument, consider
$3\le r<50$ and $2r\le n<\min(15r/4,r(r-1))$. In 2,130 of the 2,147 pairs,
the surface-center bound suffices: allow the ordinary weak alternative
already at $s=3$, retaining the canonical-extension bound whenever
$(s+1)\beta'\ge n+1$. On the untilted branch, degrees $3$ to $6$ use
$\Xi_e<\min(4,e)$ and higher degrees use $e\beta\le4n/3$ as in the adjacent-row
proof. For the tilted branch recompute $M$ with $s=3$ allowed among
the weak dimensions, and use $4n/3$ also in the high-degree cutoff. The file
\path{surface-cover.json} gives the witnesses for the 2,130 pairs.
The program \path{r2-root-surface-check.py} verifies the witnesses and
lists the 17 remaining pairs. The theorem above
includes those 17 through its primary proof.

% ===== Local thresholds along lines and conics =====
\chapter{Local thresholds along lines and conics}\label{ch:models}
\section{Local threshold estimates along lines}\label{models:sec:local-line-models}
Write $H=\mathcal O_X(1)$, $N=n+1$ and
$\delta_C=\inf_{C\subseteq c_X(v)}A_X(v)/S_H(v)$.
We work over the geometric generic point of $C$, with the original
complex dimensions as multiplicities in all section spaces. In the
local computations a normalized mean is $N$ times the asymptotic
mean order per section degree, and a normalized threshold is
$\delta_C/N$. A normalized denominator is the upper bound for the
normalized mean divided by discrepancy. Finite-level divisors
are tested below any limiting threshold and the coefficient is then
allowed to increase to the endpoint.

For a line \(L\), the \(n\) normal polar forms on \(L\) span a
basepoint-free subspace of dimension \(k\ge2\) in binary forms of
degree \(d-1\). In the quotient
\(\mathcal O^{n}/\mathcal O(1-d)=\bigoplus_{i=1}^{n-1}\mathcal O(a_i)\)
defined by this basepoint-free vector, all \(a_i\ge0\), their sum is
\(d-1\), and exactly \(k-1\) are positive. The last count follows by
dualizing and taking global sections. For ordinary order \(j\) the
split monomial of exponent \(\gamma\), \(|\gamma|=j\), has coefficient
space of dimension \(m-j+1+a\cdot\gamma\) for \(m-j\ge d-2\).
This is the symmetric-power presentation of the quotient bundle:
the kernel is a sum of copies of \(\mathcal O(m-j-d+1)\).
The associated graded is cut out by the primitive linear polar
equation over the binary base. Its ordinary Hilbert series and
generic localization agree with the full initial ring. The remaining
bounded strip of residual degrees has \(O(m^{n-1})\) sections, each
of angular degree at most \(m\). The simplex averages with density
\(x_0+\sum a_ix_i\), \(x_0=(m-j)/m\), give mean coordinate
coefficients \((1+a_i/d)/(n+1)\). The radial expected order is
\((n-1/d)/(n+1)\), with discrepancy \(n-1\). Thus for \(k\ge3\)
and \(n\ge6,d\ge4\),
\[
  \delta_L(X,H)\ge (n+1)\frac d{2d-2}>
  (n+1)\frac{d+1}{2d}.
\]
For \(k=2\), the common linear space of the two polar directions
has dimension \(n-1\). It cannot lie in \(X\): their two positive-degree
normal derivatives would have a common zero there. The weighted
calculation below therefore always applies.

\paragraph{Several transverse directions in the line calculation.}
Write $N=n+1$, and express expected orders in units of $1/N$ in this
paragraph, so the threshold bounds below are in units of $N$.
Along a line on a smooth degree-$d$ hypersurface let $k$ be the
dimension of its vector space of binary normal polars, and let
$l=n-k$. Choose polar variables $z_1,\ldots,z_k$ and zero-polar
variables $x_1,\ldots,x_l$. If the linear space $z=0$ is not
contained in $X$, let $h$ be the first degree in $x$ of the equation
there, and $G(u,v,x)$ that nonzero homogeneous term. Thus $2\le h\le d$.
Give the $z$'s weight $h$ and the $x$'s weight one. The initial
equation is $\sum f_i(u,v)z_i+G$, primitive and linear over the
binary polynomial ring. Localizing at the base and comparing ordinary
Hilbert functions identifies the whole weighted initial ring, just as
for the quartic. The discrepancy is $l+(k-1)h$, the weighted mean
$l+(k-1)h+(1-1/d)h$. Write $\alpha=1-1/d$.

Suppose a constant torus degeneration on the light space, possibly
after restricting $G$ to a linear subspace, has nonzero initial term
$a(u,v)M$, with $M=\prod_{i=1}^s x_i^{b_i}$ and $\sum b_i=h$.
Choose the first two $f_i$ as a basepoint-free pencil, possible since
the base is a curve. Split all other $z$'s and the remaining light
directions in the simultaneous degeneration. Concretely one first
specializes the eliminated coordinates to zero, retaining $f_1z_1+
f_2z_2+G|_{\rm subspace}$ as the initial relation; then gives both
retained $z$'s the weight of the chosen initial monomial of $G$.
Every surviving relation is primitive. This degeneration is geometric
on the independent generic angular coordinates, eliminating $z_2$
there; primitivity and the ordinary Hilbert series give precisely all
initial pieces. Let $b=\max b_i$. Then
\[
 \frac{\delta_L(X,H)}{N}\ \ge\
 \min\left\{\frac1{\alpha+1/h},\frac{1+b}{1+(1+\alpha)b},
 \frac{l+(k-1)h}{l+(k-1)h+\alpha h}\right\}.                    \tag{L1}\label{models:eq-mixed-monomial}
\]
The last term is at least the middle term: $l\ge s$, $h\le sb$ and
$k\ge2$. In particular taking $s=1$ always gives the uniform
bound $(d+1)/(2d)$ (the finite-point term exceeds this for $d\ge3$).

Here are details for the mixed monomial. Put $Q=M$. The ring in the
light variables is free over $\mathbb C[Q]$ with basis the monomials
not divisible by $M$, by unique monomial division. The subtorus fixing
$M$ splits the angular pieces into a residual monomial times a binary
heavy system of some degree $j$. Indeed two exponent remainders of
the same character differ by an integral multiple of the primitive
direction of $(b_i)$. Equality of total angular degrees modulo $h$
makes that difference an integral multiple of $(b_i)$ itself. In the
heavy binary system $f_1z_1+f_2z_2+aQ=0$, with $Q$ of section degree
$h$. Degenerating to $f_1z_1+aQ$ for $f_1$ coprime to $a$ gives
limiting mean coordinate exponents $1$ and $\alpha+1/h$ on the coarse
binary line, by the binary calculation below. This follows also by
using degrees $1,1,h$ for $z_1,z_2,Q$ and subtracting the degree-$d$
relation in the logarithmic Hilbert derivative. The remainder
monomials carry binary degree zero. Thus the coefficient of any finite
factor $z_1-\lambda M$ of the binary polynomials is at most
$\max(1,\alpha+1/h)=\alpha+1/h$, uniformly (including $\lambda=0$).
Each coordinate $x_i=0$ has coefficient at most one, using divisibility
in the original split ring: modulo $x_i$ the relation is the same
primitive $f_1z_1+f_2z_2$. This includes the roots at infinity.

A toric log resolution of the monomial ideal $(z_1,M)$ separates
all these finite roots: on each torus orbit where the two monomials
have the same order their primitive ratio is a torus character, and
its nonzero distinct fibers are disjoint and transverse to the toric
strata. The exponent of $z_1$ in this character is one, ensuring
primitivity. Refine the fan to a smooth fan, or use smooth quotient
charts on the weighted projective stack. For nonnegative toric orders
$U_i$ on the lights and $W$ on $z_1$, the mean on the wall
$W=\sum b_iU_i$ is exactly $\sum U_i+(1+\alpha)W$ by homogeneity
of the relation and its Hilbert derivative. For $W=0$ it is at most
$\sum U_i$, and for all $U_i=0$ it is at most $(\alpha+1/h)W$.
It is linear on either side of the wall because every nonzero finite
root has order $\min(W,\sum b_iU_i)$, and coordinate components
give linear orders. Divide by the discrepancy $\sum U_i+W$; on the
wall the maximum occurs for $b_i=b$. Free coordinate factors each
have mean one and add independent simple normal crossings divisors.
This proves (L1), including near the common angular axis before its
weighted blowup. In each finite piece use a basis adapted to the tested valuation and the diagonal order. Proper semicontinuity passes the lower threshold from its initial basis to this basis.

\paragraph{Two disjoint transverse monomials.}\label{models:par:local-disjoint}
Retain a basepoint-free polar pencil. Suppose a weight and torus
specialization yields $fz_1+gz_2+aM(x)+bP(y)$, where $M,P$ are
nonconstant monomials on disjoint variable sets, both of extraction
weight that of the heavy coordinates $z_i$. Let their ordinary
degrees be $p,q$, and their largest exponents $m_1,m_2$. All light
variables here have section degree one. Assume the gcds of the
exponents of $M$ and of $P$ are relatively prime to one another.
Then the normalized line-center threshold is at least
\[
 \min\left\{1,\frac1{1+1/p-1/d},\frac1{1+1/q-1/d},
           \frac{1+1/m_1+1/m_2}{2-1/d+1/m_1+1/m_2}\right\}.       \tag{DP}\label{models:eq-disjoint}
\]
Independent free variables cause no loss. For the primary radial
ratio, normalize the heavy weight to one; the sum of the positive
weights on the support of $M$ is at least $1/m_1$ and similarly
for $P$. Any other common weight adds the same term to the
discrepancy and mean, giving this same radial comparison.

Write $U=M,V=P$. Each coordinate polynomial ring is free over its
single monomial subring, by unique monomial division with remainders
not divisible by that monomial. Split the characters of the torus
preserving $M$ and $P$. In each finite angular-weight/character
block the two remainder indices are fixed: if the exponent gcds
are $g_1,g_2$, respective exponent differences have directions
$r\,\mathrm{exp}(M)/g_1$ and $s\,\mathrm{exp}(P)/g_2$. Equality
of angular weights modulo one forces $r/g_1+s/g_2$ integral;
coprimality implies both are integral. Thus the remaining factors
are homogeneous ternary systems in $U,V,z_1$ over the geometric
generic field, with $z_2$ eliminated by the displayed relation.
The normalized mean heavy degree is $2+1/p+1/q-1/d$. A variable
$x_i$ of $M$-exponent $b_i$ has mean remainder exponent $1-b_i/p$:
the character for the remainders is the polynomial ring character
times $(1-t^p)$, with the monomial marked by its exponent on $x_i$.
The corresponding expression for $V$ has $c_j/q$ in place of $b_i/p$.

Fix a threshold coefficient $\rho$ strictly below (DP), using units
of $n+1$ as usual. The product of a split basis consists of the
coordinate remainder powers just computed and of ternary polynomials
pulled back from $\mathbf A^3_{U,V,z}$. It suffices to prove the
klt threshold upstairs by integration. The direct image under
$U=\prod x_i^{b_i}$ of the polydisc density
$\prod |x_i|^{-2\rho(1-b_i/p)} dV_x$ is bounded near $0$ by a
constant times $|U|^{-2(1-\epsilon_U+\nu)}dV_U$ for each $\nu>0$,
where
\[
         \epsilon_U=\min_i\frac{1-\rho(1-b_i/p)}{b_i}
                   =\frac{1-\rho}{m_1}+\frac\rho p .
\]
Indeed in negative logarithmic radii the measure is a convolution
of exponential densities with those rates (up to the common factor
two); repeated minimal rates introduce at worst a fixed polynomial
in the logarithmic radius. The angular variables have constant Haar
pushforward. This also covers any small polydisc by enlargement.
Use $\epsilon_V=(1-\rho)/m_2+\rho/q$ for the second disjoint block.
This argument may be made over a general complex base point at each
fixed level after spreading out; thresholds on the geometric generic
fiber agree with those at general parameter. The $o(1)$ errors in
the remainder means are absorbed in arbitrarily small $\nu$.

It remains to check on $\mathbf A^3$ the ternary cone series with
boundary $(1-\epsilon_U)(U=0)+(1-\epsilon_V)(V=0)$; small additions
to these coefficients are permitted in the strict range. The
origin blowup tests $1+\epsilon_U+\epsilon_V$ against
$\rho(2+1/p+1/q-1/d)$, which gives exactly the last bound in (DP).
Its directions are the ternary projective plane. Away from the two
axes there, a coarse torus specialization to $fz_1+aU$ (choose the
polar $f$ coprime to $a$) yields independent coordinates with means
$1+1/p-1/d,1,1/q$. Primitivity gives the full initial spaces.
Along $U=0$ use its actual divisor first: the quotient by $U$ has
primitive equation $fz_1+gz_2+bV$, so exact divisibility gives its
normalized coefficient $1/p<\epsilon_U/\rho$. The residual system
on this line is binary with finite coefficients at most
$\max(1,1+1/q-1/d)$. At its point $V=0$ the exact $V$-coefficient
is $1/q<\epsilon_V/\rho$. Interchange the coordinates for the
other axis. The primitive quotients follow from coprimality of the
polar pencil; for the finite-point torus choose the retained polar
coprime to the coefficient of the monomial in question. The same
finite-level semicontinuity and basis-adaptation argument as for
the single binary proof applies. Blowup of the origin and polar
coordinates therefore give integrability of the pulled-back bases
in the common affine cone. In particular it holds on each weighted
angular stack chart, including points of all coordinate axes.

\paragraph{Two coprime pure powers.}
In the disjoint-monomial estimate \eqref{models:eq-disjoint}, take
$M=x^p$ and $P=y^q$ with $1\le p\le q\le d$ and $\gcd(p,q)=1$.
The largest exponents are $m_1=p$ and $m_2=q$, and the exponent
gcds are coprime as required. Since
$1+1/p-1/d\ge 1+1/q-1/d\ge1$, the estimate reduces to
\[
 \frac{\delta_L(X,H)}{n+1}\ge
 \min\left\{\frac{1}{1+1/p-1/d},\
       \frac{1+1/p+1/q}{2-1/d+1/p+1/q}\right\}.
 \tag{L2}\label{models:eq-pure-powers}
\]
In particular this applies to the consecutive exponents $p=d-1$ and $q=d$.

\paragraph{A fractional line gain at maximal pencil contact.}\label{models:par:local-maximal-contact}
Let $L$ be a line on a smooth degree-$d$ hypersurface $X\subset
\mathbf P^{n+1}$, $d\ge3$, and assume its first normal polars span
a basepoint-free pencil. In coordinates $u,v$ on $L$, choose heavy
normal coordinates $z_1,z_2$ realizing the pencil and $k=n-2$
remaining light coordinates $x$. Assume maximal contact, namely the
nonzero restriction to $z_1=z_2=0$ has generic order $d$ along $L$.
It therefore equals $G(x)$ for an ordinary homogeneous polynomial of
degree $d$ independent of $u,v$. The projective hypersurface
$Y=(G=0)\subset\mathbf P^{k-1}$ is smooth. In fact a singular
direction would span with $L$ a projective plane along which every
light and base derivative vanishes; on that plane the two remaining
degree-$(d-1)$ derivatives must have a common zero, contradicting
smoothness of $X$. We assume $k\ge3$ in what follows.

Let $c>0$ satisfy $c\le k$, $c\le k\delta(Y,\mathcal O(1))/(k-1)$,
and put $b=1-1/d$. Then the threshold over valuations whose centers
contain $L$ is at least $n+1$ times
\[
 \min\left\{\frac{k+d}{k+2d-1},\ \frac{k}{k+d-1},
             \left(1+b\frac{1-c/k}{1-c/(dk)}\right)^{-1}\right\}.
\]
In particular if $k/d\to2$ and the lower-dimensional smooth $Y$
is K-semistable, one can take $c=k(k-d)/(k-1)$, and the displayed
lower bound is $(2/3-o(1))(n+1)$.

Write the weight-$d$ initial of the defining equation for light weight
one and heavy weight $d$ as $fz_1+gz_2+G(x)$, with coprime binary
forms $f,g$ of degree $d-1$. The weighted
extraction over $\eta_L$ has discrepancy $k+d$ and expected order
$(k+2d-1)/(n+1)$. Indeed the primitive displayed equation produces
an integral graded algebra, embeds after localization at the base,
and has exactly the ordinary Hilbert series of the smooth hypersurface.
On the weighted angular fiber $\mathbf P_K(1^k,d)$, $K=\mathbf C(L)$,
scale the independent heavy coordinate $z_1$ toward zero. The resulting
equation is $gz_2+G=0$, primitive and prime. It again gives all complex
initial spaces by localization and Hilbert series, with the extraction
order and the free $z_1$ exponent retained. The normalized mean
exponent of $z_1$ is $1/(n+1)$. After removing it the remaining
polynomials are homogeneous in $x$; their normalized mean degree is
$(k+d-1)/(n+1)$, by the order-series computation. On the affine root
chart at $x=0$ the ordinary blowup has discrepancy $k$ and precisely
this mean as order. Away from this vertex the homogeneous polynomials
pull back from $\mathbf P^{k-1}_K$. On the blowup of the root chart
the same is true, including along its exceptional divisor, which is
transverse to the pullbacks. The quotient chart is quasi-\'etale, and
the free divisor $z_1=0$ away from the vertex is a transverse
coordinate. Thus it remains to estimate those series on
$\mathbf P^{k-1}_K$.

Fix the residual section degree $M$ after removing the $z_1$ exponent,
and a fixed angular $x$-degree $D$. On a base chart with parameter
$b_1=v/u$ and frame $u$, the space is
\[
 \sum_{0\le j\le D/d} U_j\, G(x)^j\mathbf C[x]_{D-dj},\qquad
 U_j=g(b_1)^{-j}\mathbf C[b_1]_{\le M-D+(d-1)j},
\]
where a negative upper degree means zero. These coefficient spaces
are nested increasing, because multiplication by $g$ uses at most
$d-1$ degrees. The spaces in $x$ are nested decreasing. Choose an
ordinary base mark $b_0$ with $g(b_0)\ne0$. The vanishing indices of
$U_j$ at $b_0$ are exactly $0,\ldots,\dim U_j-1$. For each index $r$
let $j_r$ be the least $j$ with $r<\dim U_j$. In $U_{j_r}$ choose an
element of exact order $r$. These chosen elements give a basis of
each $U_j$ by their distinct indices; tensoring each with the nested
$x$-space shows that the whole associated graded for base order has
one block $G^{j_r}\mathbf C[x]_{D-dj_r}$ per $r$. Its nonzero scalar
coefficient can be discarded for testing the geometric fiber.

This base-order degeneration is legitimate even for geometric
valuations defined using nonconstant functions of $b_1$. Apply
$b_1\mapsto b_0+t(b_1-b_0)$ to every coefficient of a basis
simultaneously adapted to the given valuation and base order. For
each $t\in\mathbf C^*$ this substitution induces a field automorphism
of $K$ and preserves geometric log canonical thresholds. Dividing
each section by its least $t$-power yields a nonzero limiting
polynomial with scalar factor $(b_1-b_0)^r$. Lower semicontinuity on
the proper projective fiber gives the asserted comparison. As in the
quartic calculation one first fixes the new test valuation on the
residual space; replacing an incoming basis by the adapted one can
only increase its value sum.

In the resulting blocks the limiting mean $j_r$ is
$(1-1/d)/(n+1)$ and the mean complete polynomial degree after
subtracting $dj_r$ is $k/(n+1)$. For a direct count, multiplicity at
index $j$ is $\dim U_j-\dim U_{j-1}$, using zero for the predecessor
at $j=0$. This counts the normal monomials of the ordinary-degree
ring $\mathbf C[u,v,x,z_1,z_2]/(u^{d-1}z_2)$ with fixed $z_1$ exponent,
$x$ degree and $z_2$ exponent $j$: at $j=0$ all base monomials are
allowed and at $j>0$ the exponent of $u$ is less than $d-1$. Marking
the $z_2$ exponent in this monomial Hilbert series subtracts $1/d$
from the free-generator contribution of one; marking the total
$x$-degree gives the free contribution $k$. Both are divided by
$n+1$ for normalized means. Thus for every angular valuation $w$ on
$\mathbf P^{k-1}$ the asymptotic mean is bounded by
\[
 \frac{k S_{\mathcal O(1)}(w)+(1-1/d)w(G)}{n+1}.
\]
The complete polynomial blocks may be estimated after extending
scalars; their bases still consist of independent scalar polynomials.
Bounded-degree blocks contribute $o(1)$, and for the remaining ones
use the uniform finite-level comparison to their complete-series
expected orders.

Finally write $A=A_{\mathbf P^{k-1}}(w)$. Projective-space
K-semistability gives $kS(w)\le A$. For centers contained in the
smooth hypersurface $Y$, signed basis adjunction also gives
\[
             cS(w)+(1-c/(dk))w(G)\le A.
\]
Indeed the mean power of $G$ is $1/(dk)$, and the residual restriction
is the full series on $Y$, with mean normalized degree $(k-1)/k$,
using the hypersurface restriction sequence. The stated choice of $c$
therefore permits inversion of adjunction. For centers outside $Y$
one has $w(G)=0$. In the first case maximize
$kS(w)+b w(G)$ subject to the two inequalities; its slope in $S$
is nonnegative because $c\le k$ and $b=1-1/d$. This gives at most
$[1+b(1-c/k)/(1-c/(dk))]A$. Together with the two radial
discrepancies and inversion over $\eta_L$ this proves the estimate.

To express this maximal-contact bound in the form (L3), let $t=k$,
let $K_0\ge1$, and suppose $\delta(Y,H)\ge(t-1)/K_0$. Put $c=t/K_0$
and $\alpha=1-1/d$. As $(t+d)/(t+2d-1)\ge t/(t+d-1)$, the bound is
\[
 \frac{\delta_L}{N}\ge\min\!\left\{\frac{t}{t+d-1},
     \frac{1}{1+\alpha(K_0-1)/(K_0-1+\alpha)}\right\}. \tag{L3}
\]
In particular the smooth sextic transversal in $\mathbf P^8$ has
$t=9$; the range for its sevenfold section gives $\delta(Y,H)\ge3$
and $K_0=8/3$, so that both entries are $9/14$.

\paragraph{Penultimate contact and a nonvanishing top form.}\label{models:par:local-penultimate}
In the penultimate polar-pencil case let \(e\) be a common zero of the
two light forms of degree \(d-1\). Assume a line through \(e\)
has finite contact order \(b\ge2\) and the top degree form has
nonzero value at \(e\).
With $y$ in direction $e$
and transverse $x$ on this line, set $c'=d-1-b$. Cut to these two
light directions and the two general heavy variables by a localized hyperplane flag through the line, retaining a basepoint-free pair of polars so each cut is smooth along the line. Give $y,x,z_i$ respective weights
$b,b+1,db$. The primitive initial relation is
\[
       f_1z_1+f_2z_2+a x^b y^{c'}+a_0y^d,
                       \qquad a\ne0,\quad a_0=C_d(e)\ne0 .       \tag{L5}
\]
This follows by minimizing weight at the smallest $x$ exponent $b$
in the penultimate term. Again the integral initial ring injects
after localization at the base and has the full Hilbert series.
The radial ratio in normalized units is
\[
 \frac{1+(2+1/b)/d}{2-1/d+(2+1/b)/d}\ge\frac{d+2}{2d+1}.
\]

On the angular weighted plane put $U=x^b$. Its $b$-th Veronese
has weights $(1,b+1,d)$ on $(y,U,z_1)$: in fixed degree modulo
$b$ the $x$ exponent has a unique remainder. Stack-theoretically
this replacement records the $b$-th-root boundary along $U=0$;
off that boundary it is an isomorphism, also in quotient charts.
The remainder $x^s$, $0\le s<b$, is bounded and contributes no
asymptotic order there. On this complement further degenerate
(L5) to $f_1z_1+f_2z_2+a U y^{c'}$ by giving $U,z_1,z_2$ lower
weight than $y^d$. This diagonal action is geometric on the entire
new weighted plane and primitivity again gives all initial spaces.
For a lower bound on the complement one may use the global threshold
of this last model without the root boundary along $U=0$; the test
point is taken downstairs where the discarded root structure is
absent. The mixed-monomial calculation of (L1) applies there with
section degree $b$ for $U$, and one for $y,z_i$. In more detail,
its finite binary coefficient is $1-1/d+1/(b+c')$ because $Uy^{c'}$
has section degree $b+c'=d-1$. The two coordinate means on the
weighted plane are at most $1/b$ and one by exact divisibility.
On the Newton wall $W=U_0+c' V_0$ for orders on $z_1,U,y$ its mean
is $U_0/b+V_0+(2-1/d)W$, whereas discrepancy is $U_0+V_0+W$.
The maximal denominator along the two extreme walls is consequently
\[
       \max\left\{\frac{1/b+2-1/d}{2},\quad
               \frac{1+(2-1/d)c'}{1+c'}\right\}.
\]
The first entry is at most $5/4$, the second at most $[1+(2-1/d)(d-3)]/(d-2)$
for the exponent range in that lemma. Off the wall one again interpolates the linear
orders of the binomial divisors and coordinates. For $c'=0$ the
unused $y$ is free, giving its coordinate bound one.

At $U=0$ use the original stack-chart $x$ filtration before this
last torus action. The primitive quotient is $f_1z_1+f_2z_2+a_0y^d$,
so the exact upstairs $x$-coefficient is one in normalized units.
The restricted angular systems on the projective binary line are
precisely the earlier two-variable pencil calculation with $q=d$:
at finite points its coefficient is at most $1-1/d+1/d=1$,
and at the remaining point exact $y$-divisibility also gives one.
The common affine binary axis is empty on this projective boundary.
This argument uses smooth stack charts to account for the orbifold
orders. In particular on the boundary the $y,z_1$ weights share a
common factor $b$ before normalization; its generic stabilizer acts
only normally, and the induced effective weighted line has weights
$(1,d)$ with the same finite and infinity orders as above. Filtered
inversion of adjunction at $x=0$ and adjunction from the first
extraction give the claimed threshold for this local model.

\section{Local threshold estimates along conics}\label{models:sec:local-conic-models}

\subsection{Conics whose plane cuts a proper divisor}
Let $X_d^n\subset\mathbf P^{n+1}$ be smooth in characteristic zero,
$d\ge3$, let $C$ be a smooth plane conic on $X$, and suppose its plane
$\Pi$ is not contained in $X$. With $H=\mathcal O_X(1)$, one has
\[
                  \delta_C(X,H)\ge (n+1)\frac{d+2}{2d}.
\]
Here the local threshold ranges over all valuations whose center contains
$C$. Write $X|_\Pi=kC+R$ with $C\not\subset\operatorname{Supp}R$ and
$1\le k\le d/2$. Throughout, a general linear flag through $\Pi$ from
dimension $n$ to dimension $b$ scales complete residual thresholds by
$(n+1)/(b+1)$, up to the additional lower bound $n+1$ contributed by
the hyperplanes. This follows from the complete restrictions of every
graded piece, the mean hyperplane coefficient $1/(j+1)$ in dimension
$j$, and localized basis inversion of adjunction at the generic point
of $C$.

\paragraph{Multiplicity one.} Take a general three-space through the
plane, with section surface $S$. At the smooth points of the plane
divisor the section is automatically transverse. At the finitely many
points of $C\cap R$ a general transversal polar is nonzero, by
smoothness of $X$, so $S$ is smooth along $C$. On any nonreduced
component of $R$ generality ensures that the transversal polar does
not vanish identically, again by smoothness. General Bertini off the
plane therefore makes $S$ a hypersurface with at most isolated
singularities. It is integral and normal: a reducible hypersurface
surface is singular along an intersection curve, and a hypersurface
is Cohen--Macaulay. Use on a resolution of $S$ the nef pullback $H$
and the unchanged conic. By adjunction $H^2=d$, $HC=2$, and
$C^2=6-2d$. The effective divisor $E=H-C$ coming from the specified
plane section consists of the residual curves and any exceptional
components in the pullback. On each residual curve $L$ mapping to a
curve we have $C\cdot L=2\deg L$ by plane B\'ezout: their intersection
is supported where $S$ is smooth and their scheme intersection is
given by their two equations inside the plane. Exceptional components
are disjoint from $C$ and orthogonal to $H$.

Consequently $2H-C=H+E$ is nef, big with square $2d-2$, and orthogonal
to each component of $E$. Thus $H-tC$ is nef for $t\in[0,1/2]$ by
convexity. Its positive part for $1/2\le t\le1$ is $(1-t)(2H-C)$,
with orthogonal negative part $(2t-1)E$ by the negative definiteness
of a collection of curves orthogonal to a big nef class. For $t>1$
intersection against $2H-C$ excludes pseudoeffectivity. Thus
\[
 S_H(C)=\frac1d\left[\int_0^{1/2}(d-4t+(6-2d)t^2)\,dt+
     (2d-2)\int_{1/2}^1(1-t)^2\,dt\right]
 =\frac{3d-2}{6d}.
\]
On $S$ the only nontrivial valuation with center containing $C$ is its
generic divisor order up to scaling. Hence $\delta_C(X,H)\ge
(n+1)2d/(3d-2)$, which is at least the stated bound since
$4d^2-(d+2)(3d-2)=(d-2)^2$.

\paragraph{Multiplicity at least two.} Here $d\ge4$. For the Fano range
$d\le n+1$ we have $n\ge3$. Take a general four-space through $\Pi$
and let $Y$ be its section threefold. With homogeneous coordinates $u$
on the plane and $x,y$ transverse, its equation is
\[
      g(u)^k c_0(u)+x A_0(u)+y B_0(u)+O((x,y)^2),
\]
where $g$ defines $C$. The normal polars of $X$ are basepoint-free
along $C$ since the plane differential vanishes there. Thus two general
combinations $A,B$ have no common zero on $C$, so $Y$ is smooth
there. Taking the combinations successively, we may also assume $A$
is nonzero at all zeros of $c=c_0|_C$. Off the plane use Bertini;
on the other components of the plane divisor the two combinations
are simultaneously nonzero wherever needed after a general choice,
or at least give only isolated singularities. In particular the flag
is integral and normal and smooth near $C$.

At the generic point of $C$ give $g$ weight one and the free normal
parameter weight $k$. The weighted blowup has log discrepancy $k+1$.
Its bidegree initial ring of absolute sections is
\[
 R_C[x,y,w]/(Ax+By+cw^k),\qquad R_C=\mathbf C[u_0,u_1,u_2]/(g),
\]
where $w$ is the symbol of $g$, of section degree two, and $x,y$
have section degree one. The presentation follows by adjoining the
symbol $w$, taking the two displayed initial relations $g(u)=0$ and
$Ax+By+cw^k=0$, and comparing Hilbert functions. The latter quotient
is a domain: over the punctured conic cone its linear coefficients
have no common zero and give an affine-space fibration, and the
fiber over the vertex has dimension at most three in a pure
four-dimensional Cohen--Macaulay ring. Thus the localization at the
generic base point, where one heavy normal is eliminated, computes
the divisorial order there. Ordinary Hilbert series agrees with that
of $Y$ since the added degree-two generator and relation cancel.
Hilbert-series differentiation gives the mean exceptional order
\[
          S_Y(E)=\tfrac14\bigl(\tfrac12+2k-k/d\bigr).
\]

The generic exceptional fiber is the weighted line $\mathbf P(1,k)$.
Its coarse binary coordinates are $x,y$, since $c\ne0$ in $K=\mathbf
C(C)$, and its stack point at infinity is $w=0$, with stabilizer $k$.
Set $Q=w^k$. The graded ring is free with basis $1,w,\ldots,w^{k-1}$
over the ring with $Q,x,y$ and the relation $Ax+By+cQ=0$. In each
stack-weight slice $kj+s$, $0\le s<k$, the remaining factor is a
binary system of degree $j$ on the coarse line, with complex base
coefficient spaces. For finite points the $w^s$ factor gives no
vanishing; the normalized mean $j/m$ equals $S_Y(E)/k+o(1)$.

On this coarse binary line, diagonal scaling of $y$ against $x$
specializes the heavy ring to $Ax+cQ=0$. Its affine-cone presentation
again injects at the generic point because $A,c$ have no common zero
on the punctured conic cone and the vertex fiber has lower dimension.
It has the same multigraded Hilbert function. On the generic binary
line $Q$ specializes to $-Ax/c$, so the initial sections are products
of powers of $x$ and $y$ with nonzero coefficients in $K$. Since $y$
splits freely in the limit, its mean exponent is $1/4$ on $Y$, and
the mean complementary exponent of $x$ is $S_Y(E)/k-1/4$. These are
upper bounds for the normalized basis-divisor coefficients on the
original coarse line in the sense that their maximum bounds the
maximum order at any geometric point: choose a complex basis for
each slice simultaneously adapted to this point and the diagonal
filtration, and apply semicontinuity of the log canonical threshold
to the limiting coordinate divisor. The finite-dimensional diagonal
coefficients are independent of the chosen adapted bases, giving
uniformity in the geometric point. At infinity use the upstairs
order of $w$ on the stack. The quotient by $w$ has initial relation
$Ax+By=0$, still prime by absence of base points, and so actual
vanishing along $w=0$ on the generic fiber is equivalent to
divisibility by $w$ in the graded ring. For the product of all
slices the order-$i$ subspace has dimension $h^0(Y,(m-2i)H)$.
Its normalized mean is therefore $1/8+o(1)$. This upstairs order
already accounts for the different $1-1/k$ on the coarse line.

The radial and angular estimates, after scaling to dimension $n$,
give the following upper bound on $(n+1)/\delta_C(X,H)$:
\[
  \max\left\{\frac{1/2+k(2-1/d)}{k+1},\ 1,\
                1-1/d+1/(2k),\ 1/2\right\}.
\]
Indeed for any smaller threshold than the reciprocals the angular
basis divisors are lc along the entire geometric generic fiber;
weighted adjunction and properness test all centers containing $C$.
The first entry increases in $k$ and equals $2d/(d+2)$ at $k=d/2$.
The third entry is at most $5/4-1/d<2d/(d+2)$ for $d\ge4$; the
remaining entries are smaller as well. This proves the claim.

\subsection{Conics on a contained plane}
Let $C\subset\Pi\subset X_d^n\subset\mathbf P^{n+1}$ be a smooth
conic on a contained plane, with $X$ smooth and $d\ge3$. Let $e$ be
the dimension of the span of the $n-1$ normal polars of the plane,
restricted to $C$. They form a basepoint-free system in
$H^0(\mathbf P^1,\mathcal O(2d-2))$, so $e\ge2$.
Put $E=\operatorname{coker}(\mathcal O(2-2d)\to\mathcal O^{n-1})$
on $\mathbf P^1$, where the map has these entries. If
$E=\bigoplus\mathcal O(a_i)$ then $a_i\ge0$, $\sum a_i=2d-2$, and
the number of positive summands is exactly $e-1$ by dualizing and
taking global sections. The ordinary extraction gives
\[
 \delta_C(X,H)\ge (n+1)\min\left\{
 \frac{n-1}{n-1/2-1/d},\,\frac{2d}{2d+\max a_i}\right\}.
\]
Indeed the conic ideal has $n-1$ linear generators and its plane
quadratic $g$. Its graded ring is
$R_C[x_1,\ldots,x_{n-1},w]/(\sum f_i x_i)$, where $w$ represents
$g$ and $R_C$ is the homogeneous conic ring. This follows from the
complete-intersection conic ideal in the polynomial ring, primitivity
over the punctured conic cone, the dimension of the vertex fiber,
and equality of ordinary Hilbert series. In section degree $m$, fix
the power $s$ of $w$ and joint $x$ degree $k$; put $p=m-2s-k$.
For $p\ge d-1$ the symmetric-power quotient sequence of $E$ has
vanishing $H^1$ in its kernel after twisting by $\mathcal O(2p)$,
since $2p-(2d-2)\ge0$.
Its split monomial indexed by $\beta$, $|\beta|=k$, then has exactly
$2p+1+a\cdot\beta$ independent coefficients. The finitely many
boundary values of $p$ have $O(m^{n-1})$ total dimension. Each residual
angular polynomial has degree at most $m$, so their normalized effect
on a valuation of discrepancy $A$ is $O(A/m)$. Integrating on the
simplex of variables $(\beta,2s,p)/m$ with density $a\cdot\beta+2p$
gives $\bar\beta_i=(1+a_i/(2d))/(n+1)$ and $\bar s=1/(2(n+1))$.
Diagonal splitting over the geometric generic fiber bounds all
angular valuations by the largest of these coefficients per unit log
discrepancy. The ordinary extraction has discrepancy $n-1$ and
mean $(n-1/2-1/d)/(n+1)$. Local inversion of adjunction at its entire
generic fiber gives the display for all centers containing $C$.
In particular if $e\ge4$, or $e=3$ and the smaller positive $a_i$ is
at least two, then
\[
                \delta_C(X,H)\ge(n+1)d/(2d-2).
\]
For the radial comparison one can use $n\ge4$ (a smooth hypersurface
of degree at least two containing a plane has at least three transverse
polars) and $d\ge3$.

\paragraph{The sharper pencil estimate.} Suppose $e=2$, and change
normal coordinates to have two heavy variables $z_1,z_2$ spanning
the nonzero polars on $C$ and $n-3$ remaining variables $x_i$ whose
polars there vanish. Let $h$ be the order of $F|_{z_1=z_2=0}$ in
the complete-intersection ideal $(g,x_i)$ of $C$. We have
\[
                      2\le h\le1+\lfloor(d-1)/2\rfloor .
\]
In fact the restrictions of the two heavy polars to the entire plane
are positive-degree forms and have a common zero there. Smoothness
of $X$ forces some light polar to be nonzero as a plane form.
It vanishes on $C$ but is divisible by $g$ to order at most
$\lfloor(d-1)/2\rfloor$ by its degree $d-1$. Its product with its
single light normal variable has order at most the displayed bound.
There is no cancellation with terms of other normal monomials.

As in the weighted conic calculation, give the $x_i$ and $g$ weight
one and the $z_i$ weight $h$. The leading equation is
$f_1z_1+f_2z_2+G(x_i,w)$ over $R_C$, with $G\ne0$ homogeneous of
weight $h$. By replacing $w$ with the symbol of $g+\sum L_i(u)x_i$
for suitable linear forms and changing the $x_i$ by a constant matrix,
we may assume $G(x_1,0,\ldots,0)=a(u)x_1^h$, $a\ne0$. Over the
generic base choose a constant direction of the $x_i$ for which the
residual polynomial in $w$ is nonzero, then choose its shear to miss
the finitely many roots. The complete-intersection presentation,
codimension of the vertex fiber and general choice of $f_1$ avoiding
the zeros of $a$ prove the primitivity and Hilbert-function assertions
in both diagonal initial rings as in the proper-plane calculation.

The $n-4$ remaining $x_i$ and $w$ split off with respective mean
exponents $1/(n+1)$ and $1/(2(n+1))$. The three variables of the
binary block, $x_1,z_1,z_2$, contribute total normalized weighted
degree $B/(n+1)$ with $B=1+h(2-1/d)$. Its preceding bidegree
calculation on $\mathbf P(1,h)$ gives angular coefficients at most
$\max(1,1-1/d+1/h)/(n+1)$ at finite points, and $1/(n+1)$ upstairs
at $x_1=0$, since here $x_1$ has section degree one. Where both
binary coordinates vanish, their weighted blowup, of discrepancy
$h+1$, separates the binary factors, giving the additional normalized
bound $B/(h+1)$. Its two quotient charts give transverse strict
transforms, also in the independent split coordinates.
The original extraction has discrepancy $n-2+h$ and normalized mean
$n-5/2+h(2-1/d)$. Its quotient by that discrepancy is no bigger
than $\max(1,B/(h+1))$, by adding $n-3$ extra unit radial discrepancies
to the binary ones and at most $n-3$ corresponding normalized order
($n-4$ linear and one degree-two generator). Angular adjunction and
properness therefore give
\[
 \frac{n+1}{\delta_C(X,H)}\le
   \max\{1,1-1/d+1/h, [1+h(2-1/d)]/(h+1)\}.
\]
For $d\ge4$ the bound on $h$ implies $h\le d-2$. The last entry
is then at most $2-2/d$; the finite-point entry also satisfies this
inequality whenever $d\ge3$ and $h\ge2$. Hence the target
$(n+1)d/(2d-2)$ holds for every contained-plane pencil with $d\ge4$.

\paragraph{A three-dimensional conic polar system.}
Let $C\subset\Pi\subset X_d^n$ be a smooth conic on a contained
plane of a smooth hypersurface, with $n\ge6$ and $d\ge3$. If the
polar restrictions on $C$ span a net, then
\[
                 \delta_C(X,H)\ge(n+1)(d+1)/(2d).
\]
This estimate includes every center containing $C$. Choose three
heavy normal variables $z_1,z_2,z_3$ so their polars span the net,
and write $x_i$ for the remaining $n-4$ light linear normals. The
space $W=(z_1=z_2=z_3=0)$ has dimension $n-2$ and cannot be
contained in $X$: three positive-degree normal polars on that
space have a common zero. The nonzero polynomial $F|_W$ has order
$2\le h\le d$ in $(g,x_i)$ for the conic equation $g$, by the
vanishing of its first jets on $C$ and the ordinary degree of the
terms. Give $z_i$ weight $h$ and light variables $x_i,g$ weight one.
The conic-ideal associated graded of the ambient polynomial ring is
$R_C[x_i,w,z_1,z_2,z_3]$ where $R_C=\mathbf C[u_0,u_1,u_2]/(g)$
and $w$ represents $g$, with section degree two. The initial equation
of $X$ is $f_1z_1+f_2z_2+f_3z_3+G(x_i,w)$, $G\ne0$. The quotient
is integral since the $f_i$ have no common zero on the punctured
conic cone, giving an affine-space fibration there; the smaller
vertex fiber cannot contribute a component in a Cohen--Macaulay
hypersurface. Hence the ring injects into its generic localization,
which is the weighted order in the regular local parameters after
one heavy variable is eliminated. Its Hilbert functions agree with
those of the original section ring (the relation is one initial
nonzerodivisor). Thus it computes every weighted-order slice.

Choose the basis of the polar net so $f_1,f_2$ span a pencil without
base point on $C$. In fact a basepoint-free net maps $C$ to a plane
curve and two general lines meet off that image. Replace $w$ by
the symbol of $g+\sum L_i(u)x_i$, for linear $L_i$, and make a
constant change of $x_i$ so $G(x_1,0,\ldots,0)=a(u)x_1^h$ with
$a\ne0$. Over $\mathbf C(C)$ choose a constant $x$-direction where
the polynomial in $w$ is nonzero; varying the value of the linear
shear on that direction avoids its finitely many roots. Degenerate
all complementary light variables, $w$, and the third heavy variable
$z_3$ to weight-separated free coordinates, leaving
\[
                       f_1z_1+f_2z_2+a x_1^h=0.
\]
On the independent generic coordinates one performs this by scaling
the split variables after eliminating $z_2$; primitivity from the
pencil and the same Hilbert functions identify all limiting spaces.
The free $z_3$ has section degree one, so its normalized mean
exponent is $1/(n+1)$ despite its heavy weight $h$. Each other free
linear coordinate has the same mean; $w$ has mean $1/(2(n+1))$.

The remaining binary block has mean weighted degree $B/(n+1)$,
$B=1+h(2-1/d)$, by logarithmic Hilbert-series differentiation.
On its core $\mathbf P(1,h)$ set $Q=x_1^h$, split off the bounded
parity $x_1^s$, $0\le s<h$, and use $z_1,z_2$ as coarse binary
coordinates. Choose $f_1$ general in the pencil so it has no common
zero with $a$ on $C$. Diagonal scaling gives the relation
$f_1z_1+aQ=0$, still primitive by the punctured-cone argument,
and frees $z_2$. The two normalized mean coordinate exponents are
$1$ and $B/h-1=1-1/d+1/h$, with means multiplied by $n+1$ here
and below. Simultaneous adaptation to the point-order and diagonal
filtrations and semicontinuity bound the order at every finite
binary point by their maximum. At infinity $x_1=0$ the quotient
has equation $f_1z_1+f_2z_2$, again prime. Thus upstairs order at
infinity equals divisibility by $x_1$ in the section ring; this
has mean one, since $x_1$ has section degree one. This upstairs
calculation includes the coarse different.

On the binary axis blow up its two normal directions with weights
$(1,h)$; the log discrepancy is $1+h$, its mean order is $B$,
and the angular factors separate into distinct binary points,
transverse to the exceptional divisor and to the independent free
coordinate divisors on the usual covering charts. This holds also
when only a heavy free coordinate is nonzero: there the local chart
is a cyclic quotient of affine space, and the same calculation
applies to its covering coordinates; these give the projective stack
and its induced different. The initial extraction over $C$ adds
the weights of all the free parameters to both the binary discrepancy
and its mean order, except that $w$ contributes only half its weight
to the latter. Therefore its ratio of normalized mean to discrepancy
is no larger than $\max(1,B/(1+h))$. The only constants to test are
thus $1$, $1-1/d+1/h$ and $(1+h(2-1/d))/(1+h)$. The last increases
in $h$ and equals $2d/(d+1)$ at $h=d$, while the others are at most
$2d/(d+1)$ for $d\ge3,h\ge2$. Angular semicontinuity for bases
simultaneously adapted to the first specialization, followed by
weighted inversion of adjunction near the entire generic fiber and
properness, gives the asserted local threshold.

\paragraph{The mixed leading form in the extreme conic pencil.}
Let $C\subset\Pi\subset X_d^n$, with $d=4$ or $5$, and assume the
normal plane polars on $C$ form a pencil. Use the two heavy normals
$z_1,z_2$, the light linear normals, and the light conic symbol $w$
of section cost two as in the conic calculation. If $d=4$, the first
light order is $h=2$; if $d=5$, it is $2$ or $3$. These bounds follow
because the two heavy polars have a common zero on the plane, where
some light polar is nonzero. For $(d,h)=(4,2),(5,3)$ at least one
nonzero light plane polar is $g^{h-1} a$, where $g$ is the conic
equation, and $a$ has degree at most one. Indeed all such polars
are divisible by $g^{h-1}$ by the assumed order, and the remaining
degree is $d-1-2(h-1)$. In particular one leading light monomial
is $a x w^{h-1}$ for a suitable linear normal $x$.

Keep $x,w$ and split off all complementary light coordinates by a
torus. Now scale $x,z_1,z_2$ to auxiliary weight one and $w$ to zero.
The pure term $w^h$ was zero since the plane is contained in $X$.
Thus the exact new relation is
\[
              f_1 z_1+f_2 z_2+a x w^{h-1}=0.
\]
It is primitive since $f_1,f_2$ have no common zero on $C$; the
same generic localization and ordinary Hilbert-series calculation
as for the pencil gives all initial spaces. In independent geometric
coordinates $x,w,z_1$ the specialization is a diagonal action after
eliminating $z_2$ over $\mathbf C(C)$. Scale $x$ oppositely to
$w^{h-1}$ to split characters in this new ring. Each character
block factors as $x^i w^j$ times a system of ordinary binary
polynomials in $z_1,z_2,Q=xw^{h-1}$ with the displayed linear
relation. Here $\min(i,\lfloor j/(h-1)\rfloor)=0$, by taking out
the largest power of the product. This factorization is unique:
the ring adjoining $x,w$ with $xw^{h-1}=Q$ is free on the indicated
monomials over the binary ring, by monic polynomial division in the
mixed monomial. All coefficients below are normalized by $n+1$.

At finite binary points choose $f_1$ avoiding the zeros of $a$.
After the usual diagonal binary degeneration $z_2$ splits freely,
with normalized exponent one. The total mean binary degree $k$
satisfies $(h-1)\bar k\le1/2+(h-1)(2-1/d)$: for the grading
assigning one to $w$, zero to $x$, and $h-1$ to both heavies,
each binary polynomial has degree $(h-1)k$ plus its nonnegative
outside power of $w$. Thus all finite binary coefficients are at
most $\max(1,1-1/d+1/[2(h-1)])$. At $x=0$ and $w=0$ the true
upstairs coefficients are at most one and $1/2$, respectively,
since the quotient by either variable has primitive relation
$f_1z_1+f_2z_2=0$. Hence divisibility in the section ring computes
the corresponding orders on the generic fiber.

The binary factors are $z_1-\lambda xw^{h-1}$ in the finite
chart. Where $w$ is a unit, blow up $x=z_1=0$ with weights $(1,1)$.
Its normalized mean is $1+(2-1/d)$ by the auxiliary grading above
which instead puts one on $x,z_1,z_2$ and zero on $w$. The discrepancy
is two. Where $x$ is a unit use the $(1,h-1)$ blowup of $w=z_1=0$;
its mean is $1/2+(h-1)(2-1/d)$ and its discrepancy is $h$. In
each case the distinct finite factors separate and their transforms,
coordinate divisors and exceptional divisor are transverse on the
covering charts. Where $x=w=0$ and a heavy coordinate is nonzero,
only the coordinate divisors meet. Where all three vanish, the
$(1,1,h)$ extraction has normalized mean $3/2+h(2-1/d)$ and
discrepancy $2+h$, and its angular weighted plane is the one just
analyzed. Free coordinates contribute normalized coefficients at
most one. The same bound persists at the original extraction and
before specialization by weighted adjunction and semicontinuity.

For $(d,h)=(4,2)$ the largest ratio is $11/8$, from the first
transverse blowup; for $(d,h)=(5,3)$ it is $7/5$ at the same blowup.
The ordinary binary bound for $(d,h)=(5,2)$ is $23/15$. Thus all
quartic contained-plane conic pencils satisfy
$\delta_C\ge(n+1)8/11$, and all quintic ones satisfy
$\delta_C\ge(n+1)15/23$. In particular the classical single-binary-
axis estimate is not sharp in the two extreme pencil cases.

The same reasoning gives a uniform improvement in every degree
$d\ge3$. Put $k=\lfloor(d+1)/2\rfloor$, $J=2-1/d$. Smoothness
ensures $h\le k$ since some nonzero light plane polar of degree
$d-1$ is divisible by $g^{h-1}$. When $h=k$, its quotient has
degree zero or one and so is nonzero on $C$. The mixed calculation
applies with no change. All its ratios are at most
\[
                P_d=\frac{1+(k-1)J}{k}.
\]
Indeed the first transverse ratio is $(1+J)/2\le P_d$ for $k\ge2$;
the second is $(1/2+(k-1)J)/k<P_d$; and the full ratio
$(3/2+kJ)/(k+2)\le P_d$ follows from $k(J-1/2)\ge2(J-1)$.
The finite-point ratio $1-1/d+1/(2(k-1))$ is at most $(1+J)/2$.
For $h<k$ the ordinary binary-axis ratio $(1+hJ)/(1+h)$ is at
most $P_d$ by monotonicity; the ordinary finite-point coefficient
$1-1/d+1/h$ is also at most $(1+J)/2$ since $h\ge2$. Thus
$\delta_C(X,H)\ge(n+1)/P_d$ for every contained-plane pencil.

\section{A moving point of a polar net}\label{models:sec:local-net-product}
Let $C\subset\Pi\subset X_d^n$ be a smooth conic on a contained
plane, with $n\ge5$, $d\ge4$. Suppose the plane polars restricted to
$C\simeq\mathbf P^1$ form a net with splitting $(1,2d-3)$ after
discarding the trivial factors. Then
\[
                \delta_C(X,H)\ge(n+1)d/(2d-2).
\]
Write $[s:t]$ on $C$. In suitable constant normal coordinates its three
nonzero polars are $(t\ell,-s\ell,b_0)$ with $\deg\ell=2d-3$,
$\deg b_0=2d-2$ and no common zero of $\ell,b_0$, since the linear
syzygy has the form $(s,t,0)$. Write $z_1,z_2,z_3$ for their normal
variables. The remaining normal variables and the conic equation are
light. The $(n-2)$-plane $z_1=z_2=z_3=0$ is not contained in $X$:
a contained $k$-plane on a smooth hypersurface of degree $d\ge2$
requires $n+1-k$ positive-degree polars with no common zero, so
$2k\le n$, whereas $2(n-2)>n$ here. The first nonzero form $G$ on
this plane has conic-ideal order $2\le h\le d$. Give the heavy
variables order $h$ and the light ones order one. As in the pencil
proof, after a shear of the degree-two conic symbol and constant
change of the light linear variables, one such variable $x$ yields
$G(x,0,\ldots,0)=a x^h\ne0$. There is a light linear variable
because $n\ge5$. Split off the remaining light variables by diagonal
degeneration. The primitive initial equation, its injection over the
generic base and its ordinary Hilbert series give the exact limiting
spaces, as in the pencil case. The remaining core has equation
\[
      a x^h+\ell(tz_1-sz_2)+b_0z_3=0
\]
over the conic cone; $a$ is a nonzero section of $\mathcal O_C(d-h)$.
The exceptional core is $\mathbf P_K(1,h,h)$ over $K=\mathbf C(C)$.
For $h\le d-2$ split off a general heavy variable leaving a
basepoint-free pencil, whose existence follows from a generic
projection of the net. The binary-axis estimate
$[1+h(2-1/d)]/(h+1)\le 2-2/d$ and the finite binary estimates prove
the result. We treat $h=d-1,d$. The following calculation also works
for $d=h=3$.

The coarse core is $\mathbf P^2_K$ with independent $z_1,z_2,z_3$.
Its boundary is $(1-1/h)L$ for the line $L=(\ell T+b_0z_3=0)$,
$T=tz_1-sz_2$, since $Q=x^h=-(\ell T+b_0z_3)/a$. On the cone
the algebra is free with basis $1,x,\ldots,x^{h-1}$ over the heavy
ring with $z_1,z_2,z_3,Q$ and the displayed linear relation. Their
section costs are $1,1,1,h$, respectively, and they all have coarse
degree one. The bounded $x$ parity factor contributes $O(1/m)$ in
normalized order upstairs. In this proof all coefficient bounds are
reported after multiplying by $n+1$.

\paragraph{Comparison with complete product systems.}
In a block of coarse degree $k$ and remaining section degree $p+k$,
the largest possible exponent of $Q$ is at most $\lfloor p/(h-1)\rfloor$.
Multiplication by $a$ to that power clears all denominators, and is a
unit on the geometric generic fiber. Hence the block embeds in the
full product of degree $k$ in the three $z$'s and of degree $2rp+O(1)$
on $C$, where
\[
 (r,R)=(1,1)\ (h=d),\qquad
 (r,R)=\left(\frac{d-1}{d-2},\frac{(d-1)^2}{d(d-2)}\right)
       \ (h=d-1).
\]
Indeed $\deg a=2(d-h)$. Enlarging to $\lfloor 2rp\rfloor$ if necessary
only requires multiplying by another base frame over $K$. In the sum
of all blocks the ratio of the enlarged to original dimensions tends
to $R=rh/d$. Compared with the full ring on the $z$'s over the conic
cone, the original heavy ring has leading Hilbert coefficient factor
$d/h$ (extra generator of degree $h$, one relation of degree $d$),
whereas giving base spaces degree $2rp$ multiplies that coefficient
by $r$. All free and bounded-power factors appear in both sums and
cancel. For any valuation an original basis can be extended to a
basis of its corresponding enlarged block, so its expected order is
at most $R$ times the enlarged expectation, when each is normalized
by its own dimension. Nonnegative orders here use a local regular
frame on the generic projective fiber; denominators depending only
on $C$ are units. Keeping the free-degree and bounded-power factors
in the block weights gives the usual free-variable means $1/(n+1)$.

\paragraph{The line and its distinguished point.}
Off $L$ the enlarged system consists of full plane polynomials,
whose normalized mean order is at most $A$ after the $n+1$ scaling
(the mean coarse degree is $3$ in this scaling). Along $L$ the
enlarged system has normalized mean divisor order
\[
                    1/D,\qquad D=1+(d-1)/r.
\]
Indeed divisibility by the primitive family equation $\ell T+b_0z_3$
costs one plane degree and $2d-2$ base degrees, thus section cost $D$.
Boundary rounding changes the limit by zero. The discrepancy is $1/h$.
The restricted blocks on $L$ after peeling off its powers come from
the rank-two split quotient $\mathcal O(1)\oplus\mathcal O(2d-3)$.
Away from the finitely many boundary values in base degree this is the
complete symmetric-power quotient tensored with the full base system,
by the symmetric-power exact sequence as before. The high-summand
coordinate generically equals $z_3$ times a unit, so it vanishes only
at $P=[s:t:0]$ on $L$. At any geometric point of $L\setminus\{P\}$
dehomogenize this high coordinate to one. The order of any polynomial
in the low coordinate is at most its degree in that coordinate. A
basis adapted simultaneously to this degree filtration and the point
order gives upper bound equal to the normalized mean low exponent,
\[
                         1+1/[2(r+d-1)].
\]
For this last computation the simplex density in low exponent $i$,
high exponent $j$, base cost $p$ and peeled line cost $Dq$ is
$i+(2d-3)j+2rp$; other free variables have density coefficient zero.
It has total vertex coefficient $2(r+d-1)$, giving the claimed
low mean by elementary simplex moments. The boundary strips in base
degree have asymptotically zero mass and bounded-degree sections give
a uniform $o(1)A$ correction. Localized adjunction along $L$ now
tests every center there except $P$. The two normalized ratios,
including the comparison factor, are $Rh/D=R$ and
$R(1+1/[2(r+d-1)])$, both at most $2-2/d$. For $h=d$ use $d\ge3$;
for $h=d-1$ use $d\ge4$, $R\le9/8$, and $r+d-1\ge4$.

\paragraph{Blowing up the moving point.}
At $P$ use its ordinary blowup on the generic coarse plane. Its
discrepancy for the pair is $1+1/h$. Global normal generators for
the graph of $P$ on $C\times\mathbf P^2$ are $z_3,T$ of bidegrees
$(0,1),(1,1)$. Filter first by their total order $j$ at $P$. In
each block of base degree $2rp$ and plane degree $k$, the normal
monomial $z_3^aT^b$, $a+b=j$, has independent image in the complete
coefficient space of degree
\[
      D_b=2rp-b+k-j,\qquad 0\le j\le k,\quad 0\le b\le2rp.
\]
Indeed first split by the $z_3$ exponent; the remaining two-variable
binary expression is divisible by $T^b$ precisely when it vanishes
to order $b$ on the graph $[s:t]$. After division its restriction
to this graph is the multiplication of complete binary systems of
degrees $2rp-b$ and $k-j$, which is surjective. Thus all the stated
coefficients vary independently in each order-$j$ piece. In the
following integrals bounded changes caused by $\lfloor2rp\rfloor$
are immaterial. Put $c=k-j$ and $e=2rp-b$. The four nonnegative
parameters $a,b,c,e$ have section cost
\[
           a+(1+u)b+c+ue,\qquad u=1/(2r),
\]
with dimension density $c+e$; the other free indices are unchanged.
The normalized mean $a$ is $1$ and the mean $b$ is $1/(1+u)$.
This gives radial ratio $R(1+1/(1+u))/(1+1/h)$. For $h=d$ it
is $5d/[3(d+1)]\le2(d-1)/d$ for $d\ge3$. If $h=d-1$, the same
inequality is equivalent to $d^3-4d^2-d+6\ge0$; setting $d=4+j$
gives $j^3+8j^2+15j+2$. On the exceptional line all points other
than the direction of $L$ have discrepancy one, and admit normalized
order upper bound $R\max(1,1/(1+u))=R$ by the diagonal action on
the two independent jet monomials and threshold semicontinuity.

The direction of $L$ has discrepancy $1/h$ and is defined by the
primitive $b_0z_3+\ell T$. A jet polynomial divisible by its $q$th
power has remaining coefficients of degrees $D_b-(2d-2)q$. Their
$T$ exponents must satisfy $b\le\min(j,2rp)-q$ by the leading
nonzero coefficient in $T$. These bounds follow by taking sections
of the corresponding line bundle on the ruled exceptional over $C$:
the weighted base cost of $T$ exceeds that of $z_3$ by one, the
relative linear equation has total base cost $2d-2$, and its
coefficients have no common zero. Thus the dimension is bounded
above by the sum of $(D_b-(2d-2)q+1)_+$ over these exponents.

For clarity the resulting normalized upper integral has a short
closed form. For the simplex with the displayed section costs one
may replace the section-cost-one uniform simplex by independent
exponentials of rates equal to the costs before normalizing: the
integrand for the dimension is homogeneous of degree one and that
for the summed vanishing is homogeneous of degree two. The extra
radial factor in their quotient is exactly $n+1$, already included
in our normalization. Equivalently this is the ratio of Laplace
transforms of the two homogeneous counting integrals. The extra
free coordinates cancel. Thus take $a,c$ exponential of rate one,
$e$ exponential of rate $u$, and integrate the remaining $q$
between zero and $\min(a,e)$. After shifting $a,e$ by $q$ the
unweighted density contributes $\exp(-(1+u)q)$, while the dimension
weight is $(c+e-(2d-3)q)_+$. Since
\[
 \mathbb E[(c+e-t)_+]=[e^{-ut}/u-u e^{-t}]/(1-u),\qquad
 \mathbb E[c+e]=1+1/u,
\]
the normalized mean is at most
\[
 J=\frac1{1-u^2}\left[\frac1{1+(2d-2)u}
                          -\frac{u^2}{2d-2+u}\right].
\]
For both choices of $h,r$ we have $1+(2d-2)u=h$. Therefore the
ratio is $RhJ<R/(1-u^2)$. For $h=d$, this is at most $4/3$ before
the strict correction. For $h=d-1$, $u\le1/2$ and $R\le9/8$ give
$R/(1-u^2)\le3/2$. These are at most $2-2/d$ in the respective
ranges. Adjunction on the exceptional line has now covered $P$.

\paragraph{The remaining coordinate axis.} Over the locus where
the ternary directions all vanish but a free light coordinate is
nonzero, blow up the transverse parameters with weights $(1,h,h)$.
The discrepancy is $1+2h$ and the normalized mean core order is
$1+h(3-1/d)$. Their quotient is at most $2-2/d$ since the difference
after clearing $d$ is $(d-2)+h(d-3)\ge0$. On this exceptional
weighted plane the preceding angular calculation applies with the
free indices carried along, including the different on $x=0$;
the weighted charts are products with the free coordinate factors,
whose normalized coefficients are at most one. The complete first
extraction at $C$ adds the free radial discrepancies and their
no larger normalized order, so it satisfies the same estimate.
The specialized bounds are uniform on the geometric generic fiber
and apply before the first diagonal specialization by semicontinuity.
Weighted adjunction and properness downstairs prove the claim.

\paragraph{A polar net on a line.}\label{models:par:local-line-net}
In rank three with splitting $(1,d-2)$ the polars may be written
$(t\ell,-s\ell,b_0)$ in line parameters $[s:t]$, with degrees
$d-2,d-2,d-1$ on $\ell,\ell,b_0$ and $\ell,b_0$ coprime. Apply the
three-heavy-variable product comparison and moving-point calculation
just established, with base polarization degree
one for a line in place of two for a conic. The normalized constants
are given here to specify every change. Put
\[
 r=(d-1)/(h-1),\quad R=rh/d,\quad u=1/r,\quad D=1+(d-1)/r=h.
\]
The pure light symbol $Q=x^h$ has coefficient of base degree $d-h$.
After clearing its denominators the enlarged $(k,p)$ blocks thus
have base degree $rp+O(1)$ (instead of $2rp$) and the same total
asymptotic dimension ratio $R$. Their line $\ell T+b_0z_3=0$ costs
$(d-1)$ base degrees. Off the marked point, the normalized low-summand
mean is $1+1/(r+d-1)$, since the three nonzero simplex weights in the
restricted line computation are $1,d-2,r$. The line discrepancy test
has denominator $Rh/D=R$. At the moving point $[s:t:0]$, an order-$j$
normal monomial $z_3^aT^b$ has coefficient-space degree
$rp-b+k-j$. Thus the two nonconstant mean jet exponents, radial
denominator and final exceptional-direction denominator are exactly
\[
 1,\quad 1/(1+u),\qquad
 R\frac{1+1/(1+u)}{1+1/h},\qquad RhJ,
\]
where
\[
 J=\frac{1}{1-u^2}\left(\frac1{1+(d-1)u}-\frac{u^2}{d-1+u}\right)
 \quad(u<1),\qquad dJ=1+\frac{d-2}{2d}\quad(u=1).
\]
Indeed the residual degree after a further direction order $q$ loses
$(d-1)q$ and after shifting the two exponential variables the
positive-part threshold is $(d-2)q$. For $u=1$ the two remaining
rate-one exponentials have positive-part mean $(2+t)e^{-t}$, proving
the formula by direct integration against $e^{-2q}$ at $t=(d-2)q$.
Lastly the normalized ternary-axis denominator is
$[1+h(3-1/d)]/(1+2h)$, and the initial radial denominator is no
larger, since additional free coordinates have unit cost.

For the sextic elevenfold these formulas give the following
normalized upper denominators. The four columns are respectively
$R(1+1/(r+d-1))$, $R(1+1/(1+u))/(1+1/h)$, $RhJ$, and the cost of the
ternary axis. The first column bounds $R$ as well.
\[
\begin{array}{c|rrrr}
h&\text{off the marked point}&\text{marked radial term}&
       \text{direction of the polar line}&\text{axis}\\ \hline
3&17/12&45/28&1025/756&19/14\\
4&23/18&13/9&325/252&37/27\\
5&29/24&875/648&8125/6264&91/66\\
6&7/6&9/7&4/3&18/13
\end{array}
\]
The largest entry is $45/28$. For $h=2$ use the binary estimate (L1)
with $d=6$. Thus every degree-six line net whose splitting is $(1,4)$
has $\delta_L(X,H)\ge12(28/45)=112/15$. For a quintic ninefold the
same expressions at $h=3,4,5$ give marked radial terms $3/2,704/525,5/4$
and each of the four denominators is at most $3/2$. At $h=2$ use (L1).

\section{Conic thresholds in degrees four, five and six}
\label{models:sec:local-small-conic-table}
For a smooth conic $C$ on a smooth quartic eightfold, quintic ninefold,
or sextic elevenfold, the local estimates give respectively
\[
       \delta_C(X,H)\ \ge\
       \frac{369}{61},\qquad\frac{32950}{5083},\qquad\frac{4176}{553}.
       \tag{SCC}
\]
We prove these estimates simultaneously, with the normalized denominator
$b$ equal to $61/41,5083/3295,553/348$ in the three degrees.
If the conic plane is a proper section, Section~\ref{models:sec:local-conic-models}
gives denominator $2d/(d+2)\le b$. For a contained plane and a polar
pencil the mixed extreme-pencil proof in that section gives
\[
 P_d=\frac{1+(k-1)(2-1/d)}{k},\qquad k=\lfloor(d+1)/2\rfloor ;
\]
this is at most $b$ in each degree. The ordinary split estimate suffices
as soon as the largest positive splitting degree is at most $2d-5$,
because $1+(2d-5)/(2d)\le b$. The positive degrees sum to $2d-2$.
Thus the remaining types are nets with positive summands $s,2d-2-s$,
$s=1,2$, and four-dimensional systems with summands $1,1,2d-4$.
For the latter a linear syzygy allows two independent polar coordinates
$t\ell,-s_0\ell$ in line parameters $[s_0:t]$ on $C$; their space
can be completed to a net by choosing a third polar that avoids the
finitely many roots of $\ell$. This is a basepoint-free net with
linear syzygy. Give all active normals weight $h$, where $h$ is the
initial light order on their common kernel; this kernel cannot lie
in $X$ since it has dimension at least $n-3$ and the at most four
transverse polars of positive degree would have a common zero. In
the four-polar case keep the fourth heavy coordinate as an independent
free variable by the same diagonal specialization as in the weak-net
calculation. The two-light and three-heavy calculations therefore use
one of these nets at some $2\le h\le d$.
The remaining core on its coarse plane, including the boundary, is
\[
 ax^h+L(z)=0,\qquad L=f_1z_1+f_2z_2+f_3z_3,\qquad
          (\mathbf P^2,(1-1/h)(L=0)).
\]
Here $\deg a=2(d-h)$ on the conic in its binary coordinates.
The complete axis is controlled by $[1+h(3-1/d)]/(1+2h)$,
which is smaller than $b$ in all the entries of the table below.

\paragraph{Two simultaneous bounds on the actual plane.}
Without the boundary $L$, all valuations of the coarse plane have
normalized mean at most $c_0 A$, $c_0=1+1/h-1/d$. Indeed choose
$f_1$ avoiding the zeros of $a$ and diagonally specialize both
$z_2,z_3$ out of the relation $aQ+L$. Primitivity of $aQ+f_1z_1$
and equality of the Hilbert functions give all initial spaces. The
two free normalized exponents equal one, and the remaining exponent
equals $1+1/h-1/d$ by differentiation in heavy degree. Threshold
semicontinuity after simultaneous diagonal and valuation adaptation
gives the assertion on the whole coarse plane.

The actual mean order along $L$ is $1/h$ up to $o(1)$: setting $x=0$
leaves the primitive relation $L=0$, so upstairs divisibility is by
the degree-one symbol $x$. Fix a plane valuation and choose a basis
adapted simultaneously to it and to the actual $L$-order. After the
$L$-component is removed, let $\Gamma$ be the remaining normalized
basis divisor. At any marked geometric point $P\in L$ its
intersection with $L$ is at most
\[
                     I=1+(2d-2-s)/(2d)+o(1).
\]
In fact the quotient after peeling $x$ is the complete symmetric
power of the binary split $(s,2d-2-s)$ for all sufficiently large
base degrees, and the other degrees give asymptotically negligible
mass. Split monomials have coefficient dimension
$2p+1+s\alpha+(2d-2-s)\beta$; the simplex density yields the two
means $1+s/(2d)$ and $1+(2d-2-s)/(2d)$. On the binary line this
bounds the order of every basis at every point by diagonal
semicontinuity. If $s=1$, away from the high coordinate zero
$P=[s_0:t:0]$ take the high coordinate as a unit: the order of
every residual binary polynomial is bounded by its degree in the
low coordinate. There the smaller $1+1/(2d)$ suffices and ordinary
adjunction along $L$ gives the desired bound.

\paragraph{Multiplicity of the residual at the distinguished point.}
Suppose $s=1$. The following calculation applies at each $h$; if $a$ is constant, the positive-pole sector has zero mass.
On $C\times\mathbf P^2$ the graph of $P$ has the two normal
generators $z_3$, $T=tz_1-s_0z_2$, of base/plane bidegrees
$(0,1)$ and $(1,1)$. The normalized mean total multiplicity at
$P$ in the actual plane system is at most
\[
            \frac hd\frac53+\left(1-\frac hd\right)(2+1/h).
\]
To see this, filter each heavy-degree-$k$, section-degree-$(p+k)$
block by pole order $i$ in $Q=L/a$. The regular level consists
of the full product with base degree $2p$ and plane degree $k$.
Its fraction of the total dimensions tends to $h/d$ by the heavy
Hilbert series (one degree-$h$ generator and one degree-$d$ relation).
Its normalized mean point multiplicity with its own dimension is
$1+2/3=5/3$. Indeed split the $z_3$ powers, then the $T$ powers;
they cost one and $3/2$ section degrees. Equivalently the $z_3^jT^l$
jet coefficients have degree $2p-l+k-j-l$ with $2p\ge l,k\ge j+l$,
by surjectivity of multiplication of the binary sections on the graph.
These costs give the asserted two means. All bounded-power factors and
free variables are included in the common normalized grading.

At each pole level $i\ge1$, leading residue on $a=0$ injects the
quotient into $L^i$ times polynomials of plane degree $k-i$ over
the zero scheme of $a$, of base degree $2p-2(h-1)i$. It is the
complete residue space once that base degree is at least $2(d-h)-1$.
Indeed the top numerator is $L^i F_i$, and $L$ is a relative linear
coordinate at every base point of $a=0$, by basepoint freeness. A
zero residue implies $F_i$ divisible by $a$, reducing the pole by
substituting $aQ=L$. The finitely many boundary base degrees have
asymptotically negligible dimension. If the original section vanishes
to order $j$ along the generic graph, its regularized numerator
vanishes along its entire $j$th infinitesimal neighborhood: the
graph is a relative smooth section and all its thickened structure
sheaves are flat over the base. Reduction at $a=0$ preserves this
condition, including at nonreduced base roots. The residual $F_i$
therefore has graph-point order at least $(j-i)_+$ on that finite
scheme. One can choose relative coordinates at each root with one
normal coordinate equal to $L$, whose linear part is nonzero in
every fiber. On the complete plane polynomial system of degree $k-i$
the maximal mean multiplicity at a point is $2(k-i)/3$, by counting
monomials in two affine normal coordinates. The same count holds
over an Artinian base multiplied by its length. This argument bounds
each order-$j$ space by its regular portion and pole quotients,
so it bounds the maximal mean, regardless of basis compatibility.
In the positive pole sector the top power $i$ costs $h$, and the
three remaining plane variables each cost one. The bounded-length
base residue in high base degrees supplies one free cost-one index.
The elementary simplex means, retaining the common free and parity
factors, give normalized $\bar i=1/h$ and $\overline{k-i}=3$
within this sector. Its resulting normalized bound is $2+1/h$;
weighting by the two fractions proves the display. Subtracting the
actual $1/h$ copies of $L$, we obtain the residual multiplicity
bound
\[
     M=\min\{I,\,2-h/(3d)-1/d\}\qquad(s=1).
\]
The additional estimate by $I$ follows since $L$ is smooth and local
intersection with $L$ bounds point multiplicity. For $s=2$ this
argument gives $M=I$ at every point of $L$, which is sufficient.

\paragraph{A two-coordinate inequality on the smooth plane.}
For a valuation centered exactly at $P$, normalize $v(\mathfrak
m_P)=1$, let $A$ be its discrepancy without boundary, and set
$u=v(L)\ge1$, so $A\ge u+1$. After dropping asymptotically vanishing
errors, the normalized order is bounded above by both
\[
            c_0 A,\qquad u/h+(A-1)\min(M,I/u).
\]
For the second, choose local coordinates $x=L,y$ with $y$ general,
$v(y)=1$. After the Puiseux cover $y=t_0^e$ each curve of $\Gamma$
splits formally into factors $x-\phi(t_0)$ times a unit, the total
number of factors equaling its multiplicity. The ramification
formula and the smooth coordinates $t_0,x-\phi$ give
$v(x-\phi)\le A-1$, since $v(t_0)=1/e$. If $v(\phi)\ne u$ its
order is $\min(u,v(\phi))$, and if $v(\phi)=u$ use $A-1$.
Since $u\le A-1$, each factor is also at most
$(A-1)v(\phi)/u$; their $v(\phi)$ sum to intersection with $L$.
This proves the formula for every effective rational $\Gamma$
without component $L$. Equivalently one may use the smooth-surface
valuation intersection inequality. All competing divisors have
already been separated for the very same adapted basis.

Put $c=b(1-1/h)$, $r=(b-c_0)/c$ and $\rho=c+1/h$.
The first bound proves the desired estimate $bA-cu$ when $A\ge u/r$.
On $u+1\le A\le u/r$ the second requires
$(m-b)A+\rho u-m\le0$, $m=\min(M,I/u)$. It is enough to test
the two endpoints in $A$. At $A=u+1$ the worst occurs at $u=I/M$,
with slack $b+(b-\rho)I/M-I$. At $A=u/r$ the expression increases
up to $u=I/M$ in the case $m>b$; beyond that point its derivative
is at most $I+\rho-b/r<0$. The slack at $I/M$ is
$M-[(M-b)/r+\rho]I/M$. Where $m\le b$ only the first endpoint
is needed.

Consequently, with
\[
 I=2-\frac{2+s}{2d},\quad c_0=1+1/h-1/d,\quad
 \rho=b-(b-1)/h,\quad r_e=\frac{b-c_0}{b(1-1/h)},
\]
the endpoint test used below is
\[
\begin{gathered}
 b+(b-\rho)I/M\ge I,\qquad r_e>0,\\
 M\le b\quad\text{or}\quad
 \left\{M\ge\big[(M-b)/r_e+\rho\big]I/M,\quad b/r_e\ge\rho+I\right\}.
\end{gathered} \tag{E}
\]
These inequalities test all $A,u$ in the two-coordinate estimate;
in the preceding derivation $r$ means $r_e$.
The generic valuation of $L$ itself has ratio one since its
coefficient is $1/h$ and its discrepancy is $1/h$.

Two other tests can be used when they give a smaller upper bound.
The binary specialization of Section~\ref{models:sec:local-conic-models}
needs
\[
  \max\left\{1+1/h-1/d,\frac{1+h(2-1/d)}{1+h}\right\}\le b . \tag{B}
\]
For $s=1$ the product comparison of Section~\ref{models:sec:local-net-product}
gives the following sufficient test:
\[
 \max\!\left\{R\left(1+\frac1{2(t+d-1)}\right),
       R\frac{1+1/(1+u)}{1+1/h},\frac{R}{1-u^2}\right\}\le b,\quad
 t=\frac{d-1}{h-1},\quad R=\frac{th}{d},\quad u=\frac1{2t}.
 \tag{P}
\]
Indeed the same computation works at any $2\le h\le d$:
in a block of base cost $p$ the maximal exponent of $Q=x^h$
is $\lfloor p/(h-1)\rfloor$; multiplying by this power of $a$
embeds the block in one of base degree at most
$2p+2(d-h)\lfloor p/(h-1)\rfloor\le2tp$.
Thus the dimension ratio is $R$, and $1+(2d-2)u=h$ gives
the asserted bound for the special direction on the blowup.
The first entry also bounds $R$. The other centers and free axes were
treated in Section~\ref{models:sec:local-net-product}.

Here are rational certificates for the finite cases. In each entry the
letters in order correspond to $h=2,\dots,d$: substitute $d,h,s,b$
in the displayed test with $M$ as given above. The axis, pencil,
ordinary and proper-plane terms have already been accounted for.
\[
\begin{array}{c|c|cc}
d&b&s=1&s=2\\\hline
4&61/41& E,E,E & E,E,E\\
5&5083/3295& B,E,P,P & B,E,E,E\\
6&553/348& B,E,P,P,P & B,E,E,E,E
\end{array}
\]
For example, the entry determining the endpoint in each row is
$E$ for $(s,h)=(1,2),(1,3),(1,3)$ respectively. There $b$ is
obtained by solving the second endpoint equality in (E):
\[
       b=\frac{c_0(M^2-I/h)}
            {M^2-I/h+I(1-1/h)(c_0-M)}.
\]
The checking consists solely of the comparisons (B), (E), (P) at the
listed rational arguments; the same comparisons with exact arithmetic
appear in the conic certificate. Equalities here give lower bounds by
using an arbitrarily smaller test coefficient and taking the limit.
Applying weighted adjunction above the entire geometric generic fiber
gives (SCC). In particular $32950/5083>58/9$ since the cross-products
differ by $1736$.

% ===== The quintic ninefold =====
\chapter{The quintic ninefold}\label{ch:quintic}
\section{The range and the new case}
Let \(X=X_d^n\subset\mathbf{P}^{n+1}_{\mathbb C}\) be smooth and let \(H=\mathcal{O}_X(1)\).
Adjunction gives \(-K_X=rH\), where \(r=n+2-d\). Write \(S=S_H\) for
expected vanishing. For a closed subvariety \(C\), use the convention
\[
 \delta_C(X,H)=\inf_{C\subseteq c_X(v)}A_X(v)/S(v).
\]
Thus a lower bound at a curve also tests divisors whose center strictly
contains that curve.

\begin{theorem}\label{quintic:main}
Using the uniform range, the slope estimates of Chapter~\ref{ch:global} and the local adjunction lemmas of Chapter~\ref{ch:models}, every smooth
quintic ninefold in \(\mathbf{P}^{10}\) is uniformly
K-stable. Together with the preceding chapters this gives the following sufficient range
for \(d\ge4\):
\[
 d-1\le n\le 2d-2,\qquad\text{or}\qquad(d,n)=(4,7),(5,9).
\]
\end{theorem}
Here \(n=2d-2\) is \(n=2r-2\). Each of the two further cases has
\(n=2r-3\). Chapter~\ref{ch:boundary} gives the proof for the
quartic sevenfold, including its own short-curve threshold \(16/3\).
Chapter~\ref{ch:sextic} adds
\((d,n)=(6,11)\) by a new calculation for lines of contact five.
The results of Chapters~\ref{ch:eight} and~\ref{ch:nov30-replacements}
cover all quartics of dimension at least eight. For quintics,
Chapters~\ref{ch:secondrow} and~\ref{ch:beyond} extend the present
ninefold result to dimensions ten and eleven, respectively. The extensions beyond the published cases in Section~\ref{sec:established-cases}
are research arguments of this paper.

\section{Point valuations and the lowest one fifth of the series}
\begin{lemma}\label{quintic:bottom}
On a smooth polarized \(n\)-fold with \(H^n=d\), let \(w\) be centered
at a closed point \(p\) and normalized by \(w(\mathfrak m_p)=1\).
For fixed \(0<f<1\), the average of the lowest fraction \(f\) of its
jumping numbers on \(H^0(X,mH)\), divided by \(m\), is at least
\(\frac n{n+1}(df)^{1/n}-o(1)\). The error is uniform in \(p,w\).
\end{lemma}
\begin{proof}
Every germ satisfies \(w(g)\ge\operatorname{ord}_p(g)\). The codimension of the
subspace of sections of \(w\)-order at least \(j\) is at most
\(\ell(\mathcal{O}_{X,p}/\mathfrak m_p^j)=\binom{j+n-1}{n}\), by the jet map.
The complete system has dimension \(d m^n/n!+O(m^{n-1})\).
Consequently the increasing degree-normalized quantile at index
\(x\) is bounded below by \((dx)^{1/n}-o(1)\), uniformly on compact
subintervals of \((0,1)\). Integrate from \(0\) to \(f\), using
nonnegativity at zero. Both dimension bounds are independent of \(p,w\).
\end{proof}

\section{Reduction to a numerical inequality}
We now work on a smooth quintic ninefold. Put \(N=10\), \(u=9/10\).
Sections~\ref{quintic:lines} and~\ref{quintic:conics} establish the common mark
\[
                 M=58/9 \ \le\ \delta_C(X,H)
       \qquad (C\text{ a line or an integral conic}).               \tag{1}
\]
We recall the global deductions explicitly, as the point estimate and
the universal comparison have different constants.

Suppose \(\beta=\delta(X,H)\le6\). The divisorial minimizer theorem
\cite{ar2102.09405} applies to the log Fano polarization after scaling by six.
Take a closed point on the center of an original minimizer \(v\).
A general hyperplane section through that point is a smooth quintic
eightfold with index five, in the range of Chapter~\ref{ch:boundary}. Local hyperplane
basis adjunction \cite[Lemma 4.6]{ar2003.13788} gives
\(\beta\ge\min(10,(10/9)\delta(X\cap H,H))>50/9\).
Fix \(f=1/5,\ a=151/400,\ \varepsilon=10^{-8}\). For \(t\in(0,1)\)
let \([t^{1/9}]_-\) be the downward rounding to a multiple of \(10^{-12}\).
Define
\begin{align*}
 k_0&=9(1-[(1-f)^{1/9}]_-),&
 \gamma&=(1-6/10)f/(1-f)=1/10,\\
 \alpha&=\frac{1-10f+9f[f^{1/9}]_-}{10(1-f)}-\varepsilon,&
 \lambda&=\frac{a(\beta-10/3-3\varepsilon)}{1+k_0}.
\end{align*}
Use a basis adapted both to \(v\) and to the order of a hyperplane
avoiding the generic center of \(v\). Retain the largest \(1-f\)
fraction for \(v\), remove each member's copies of that hyperplane,
average in degree-normalized units over the retained members and
multiply by \(\lambda\). Denote the resulting effective boundary by
\(D\sim_{\mathbb Q}\lambda cH\). The original dreamy minimizer has
rational \(\beta\), so the selected constants give a rational boundary.
At a sufficiently high basis level the following bounds hold:
\begin{align}
 c&\le 1-\alpha,\qquad
 \nu(D)\le kS(\nu),\qquad k=\lambda(1+k_0)+\varepsilon,\tag{A}\label{quintic:all}\\
 \beta-k&\le\beta':=\delta(X,D;H)
          \le\beta-\lambda(1+\gamma)+\varepsilon .\tag{B}\label{quintic:two}
\end{align}
For clarity, the tail comparison and its uniform finite-dimensional
version are proved in Chapter~\ref{ch:global}, in the top-half lemma and the
subsequent beta-law calculation. Concavity of the \(n\)-th root of
the survival volume bounds the upper \(1-f\) conditional mean by
\((1+n(1-(1-f)^{1/n}))S(\nu)\). For the discarded hyperplane factors
their full survival function is \((1-t)^9\), so their average in any
retained fraction is at least \(\int_0^1((1-t)^9-f)_+dt/(1-f)\).
At \(v\), the subseries from a general finite linear projection has
asymptotic fraction \(f\) and mean at most \(A_X(v)/10\), since
the projection is \'etale at the generic center and
\(\delta(\mathbf{P}^9,\mathcal{O}(1))=10\). Interlacing with that subseries bounds
the mean of the lowest fifth from above. The complementary mean at
\(v\) gains at least \(\gamma S(v)\). These three estimates give
\eqref{quintic:all}--\eqref{quintic:two}, with room for the indicated finite-level errors.

For an expression depending on \(\beta\), subscript \(6\) means
evaluation at \(\beta=6\). Set
\begin{align*}
 U &=6+\lambda_6(k_0-\gamma)+2\varepsilon,\\
 z &=6/10+\lambda_6\{u(1+\gamma)-1+\alpha\}-u\varepsilon .
\end{align*}
The underlying comparison at a \(\beta'\)-optimizer is \(A_X/S\le U\),
and \(6-\lambda c-u\beta'\ge z\). Indeed the corresponding affine
functions \(\beta+\lambda(k_0-\gamma)+2\varepsilon\) and
\(6-u\beta+\lambda\{u(1+\gamma)-1+\alpha\}-u\varepsilon\)
are respectively increasing and decreasing. Furthermore
\(\beta-k>9/2\) and \(\beta'<6-\lambda c\) on \(50/9<\beta\le6\);
these endpoint comparisons are included in the certificate. In particular
\((X,D)\) is log Fano and again has a divisorial minimizing valuation \(w\).

Its center must be a point. Indeed for a positive-dimensional
minimizing center of dimension \(s\ge2\), the ordinary tangent-slope
argument of Chapter~\ref{ch:global} supplies an integral rational curve of
degree at most \(9(s+1)/(s\beta')<3\). For \(s=1\) the same slope
gives rationality, and the center has degree at most \(10/\beta'<3\):
the weight-zero system gives a projection of the special fiber onto
the curve; at a general point of the curve, the asymptotic maximum
vanishing of its hyperplane series is its degree. At a prime above
that point use \(S\ge T/10\) and discrepancy at most one for the
log pair. These deductions, including the preservation of the
threshold on the special fiber, are given in Chapters~\ref{ch:global} and~\ref{ch:boundary}. Equation (1) contradicts \(A_X(w)/S(w)\le U<M\).

Normalize \(w(\mathfrak m_p)=1\) at its point center and write
\(A_X(w)=9+q\), \(q\ge0\). Applying Lemma~\ref{quintic:bottom} to the full
\(w\)-jump sequence gives \(u-o(1)\) for the mean of its bottom fifth.
The sum of orders of any fixed number of independent sections is at
most the sum of that many largest jumps. This comparison applies to
the actual \(v\)-selected sections, with no \(w\)-adaptation assumption.
Removing effective hyperplanes can only decrease their order. Consequently
\[
 w(D)\le(\lambda/(1-f)+\varepsilon)S(w)-\lambda uf/(1-f).
\]
Here uniform \(S_m\le(1+o(1))S\) and the uniform point estimate
absorb the errors into \(\varepsilon S(w)\), since \(S(w)\) is
uniformly bounded below by Lemma~\ref{quintic:bottom}. Combining with
\eqref{quintic:two} gives
\[
 S(w)\ge (A_X(w)+P)/U_p,\qquad
 U_p=6+\lambda_6(f/(1-f)-\gamma)+2\varepsilon,\qquad
 P=\lambda_6 uf/(1-f).                                            \tag{2}
\]
In passing to the endpoint one must move numerator and denominator
together. At \(\beta\), replace \(6,\lambda_6\) by \(\beta,\lambda\)
in these expressions. The quotient \((A+P(\beta))/U_p(\beta)\)
decreases: \(P'(\beta)U_p(\beta)<9\le A\) and \(U'_p(\beta)\ge1\).

Set
\[
 \Delta=\min\{M-U,(z-1/3)/u\}-\varepsilon,\quad
 \eta=z-u\Delta,\quad Q=53/20.
\]
We claim the necessary multiplier-ideal inequality
\[
                         1\ \ge\ \Delta S(w)+\eta .                \tag{3}
\]
In a new large degree, use a basis adapted to \(w\) and a hyperplane
missing \(p\). Remove its hyperplane factors to get
\(\Gamma_m\sim_\mathbb Q(u+o(1))H\), with \(w(\Gamma_m)\to S(w)\)
and \(\nu(\Gamma_m)\le(1+o(1))S(\nu)\) uniformly. At
\(t=\beta'+\Delta\), a first positive-dimensional non-klt center
at \(t_m\le t\) would have \(t_m\ge\beta'/(1+o(1))\). If necessary,
a general hyperplane through \(p\) with coefficient tending to zero
makes the pair non-lc at \(p\) while preserving those centers
generically. By the choice of the first threshold, the non-lc locus
is finite. The ample remainder has coefficient at least
\(\eta-o(1)>1/3\). Quasi-log adjunction and vanishing on the union
of a minimal positive lc center and the non-lc scheme make that
union connected with finite nonempty non-qklt locus. The length-one
theorem \cite[Theorem 1.8]{ar2103.01417} produces a rational curve of degree
at most two on this center. Its lc place before the hyperplane
perturbation would have \(A_X/S\le k+t+o(1)\le U+\Delta+o(1)<M\).
This contradicts (1). Full details of this quasi-log argument are
in Chapter~\ref{ch:boundary}. The multiplier ideal at \(t\) therefore has finite
cosupport. Add a general pointed hyperplane of coefficient approaching
\(\eta<1\) from below. The augmented ideal still has finite cosupport;
Nadel vanishing with ample remainder makes constants surject onto the
quotient. Its length is at most one, so it contains \(\mathfrak m_p\)
locally. Applying the valuative multiplier-ideal inequality gives
\(1\ge tS(w)+\eta-A_X(w)+w(D)\) in the limit, proving (3).

\section{The discrepancy at the optimizing point}
We prove \(q>Q\). Let \(b=\operatorname{mult}_p(T)\) for the projective tangent
hyperplane section \(T\). When \(b<5\), general pointed linear cuts
give a smooth quintic surface whose tangent plane cuts an integral
quintic \(C\) of multiplicity \(b\). Indeed the tangent section has
isolated singularities by finiteness of the Gauss map; projection
from \(p\) is generically finite when \(b<5\). Pointed Bertini then
gives the integral curve. The surface estimate
\(3S_{H|_S}\le B A_S\) at \(p\) lifts along seven smooth pointed cuts
with the retained hyperplane coefficient, giving
\[
                       10S(w)\le B(q+2)+7.                        \tag{4}
\]
This is the signed flag in Chapter~\ref{ch:boundary}. For \(b=2\), B\'ezout with
the tangent curve gives Seshadri constant at least two, so the
surface estimate in \cite{ar2101.09246} gives \(B=5/2\).
For \(b=3,4\), if \((S,cC)\) is lc near \(p\) one can take
\[
                        B=5/b+(1-5/b^2)/c .
\]
Indeed a degree-\(m\) section of order \(j\) must contain at least
\(\lfloor(bj-5m)_+/(b^2-5)\rfloor\) copies of \(C\). After dividing
by these copies the binary leading forms have the tent-shaped
degree bound derived in Chapter~\ref{ch:boundary}, section on the two small indices.
Its area is \(5/2\); integrating the binary orders gives the displayed
coefficient. For \(b=3\), the integral quintic has arithmetic genus
six. The plane-curve discriminant bound in Chapter~\ref{ch:boundary} gives
\(c\ge1/3+2/(2\cdot6+3-1)=10/21\), whence \(B\le13/5\).
For \(b=4\), if the initial quartic has no factor of multiplicity
greater than two, a general direction line has all roots of
multiplicity at most two. Blow up the surface point. The strict
transform with coefficient \(1/2\) restricts to the exceptional
line with coefficients at most one, so inversion of adjunction
gives \(c=1/2\). This yields \(B=21/8\).

For the remaining \(b=4\) case an initial linear factor has
multiplicity \(s=3\) or \(4\). Write coordinates
\([T_0:y_1:\ldots:y_8:L:z]\), with tangent hyperplane \(z=0\).
There the equation is \(T_0L^s Q(y,L)+g_5(y,L)\), with
\(Q(y,0)\ne0\). Put \(G=g_5(y,0)\). Its zero locus in \(\mathbf{P}^7\)
has only isolated singularities: a positive-dimensional singular
locus meets \(\partial_Lg_5(y,0)=0\), producing a line of singular
points as \(T_0\) varies in the tangent section. The same argument
excludes \(G=0\) identically. Thus \((\mathbf{P}^7,G)\) is lc. Indeed,
scale one coordinate to degenerate it to the cone over a general
smooth hyperplane section. The vertex has blowup discrepancy \(7\)
and multiplicity \(5\), and the resolved cone has smooth strict
transform transverse to the exceptional divisor.

Give \(y_i,L\) weights \(s,s+1\) in the pointed smooth chart,
so \(z\) has weight \(5s\). The weighted divisor \(E\) has
\[
 A_X(E)=9s+1,\qquad S(E)=(13s+1)/10.
\]
The section-ring relation \(T_0^4 Z+T_0L^sQ(y,0)+G=0\) is
primitive. Localization at \(T_0\) and equality of Hilbert series
give the full associated graded. We spell out the angular constants
for the weighted extraction of Chapter~\ref{ch:boundary}. Scaling \(L\) out
leaves the prime relation \(T_0^4 Z+G=0\). In \(L\)-exponent \(k\),
total \(y\)-degree \(D_0\), the spaces are
\(L^k G^j\mathbb C[y]_{D_0-5j}\) with
\(j=\max(0,\lceil(D_0+k-m)/4\rceil)\). The section-normalized
means of \(k,D_0,j,D_0-5j\) are \(1/10,12/10,2/25,8/10\).
For example use standard monomials \(T_0^aL^k Z^h P_l(y)\) with
\(P_l\) in a complement to \((G)\); then \(j=(h-a/4)_++O(1)\).
On the homogeneous simplex the two slots \(h,a\) give mean
\((4/5)/10\) for this positive part. Since \(G\) is lc and
\(\delta(\mathbf{P}^7,\mathcal{O}(1))=8\), the mean base order per unit discrepancy
is at most \(2/25+(8/10)/8=9/50\). The independent \(L\) exponent
has mean \(1/10\). In the \(L\ne0\) chart the radial inequality
\((50/9)(12/10)<8\) holds. Weighted basis adjunction on the
uniformizing charts, as in Chapter~\ref{ch:boundary}, therefore gives
\[
 A_X(w)-(50/9)S(w)
 \ge (A_X(E)-(50/9)S(E))\,w(E)
 \ge \frac{16s+4}{9(s+1)}\ge\frac{13}{9}.
\]
We used \(w(E)\ge1/(s+1)\), since the largest local coordinate
weight is \(s+1\). Together with \(A_X(w)\le U S(w)\) this gives
\(A_X(w)\ge(13/9)U/(U-50/9)>9+Q\).

If \(b=5\), the tangent section is a cone with smooth quintic
sevenfold base; a singular point downstairs would give an entire
singular line upstairs. The base is in the range of Chapter~\ref{ch:global}.
The generalized Eckardt formula \cite[Corollary 4.19]{ar2003.13788} gives
\(\delta_p(X,H)=90/13>U\). This excludes \(b=5\).

In all the cases with (4), use \(B=21/8\). The rational certificate
checks
\[
        U_p\{(21/8)(q+2)+7\}<10(9+q+P)
                       \qquad(q=0,\ 53/20).
\]
Both sides are affine, so (2) and (4) exclude the entire closed interval
\(0\le q\le Q\). Finally,
\[
 z+\Delta\big((9+Q+P)/U_p-u\big)>1+10^{-4}.
\]
Equations (2) and (3) now contradict each other. Thus (1) gives
\(\delta(X,H)>6\), equivalently \(\delta(X,-K_X)>1\). This proves
Theorem~\ref{quintic:main} once the curve estimates below are supplied.

\paragraph{Size of the strict inequalities.}
The accompanying verifier uses rational arithmetic, including integer
enclosures of the ninth roots. For orientation it gives
\(U=6.09931390\), \(U_p=6.12372994\), \(P=0.18559489\),
\(z=0.64384697\), \(\Delta=0.34501514\) to the displayed precision.
The margin in the affine surface comparison at \(Q\) exceeds \(0.74\);
the final Nadel margin exceeds \(0.00015\).
\section{The line threshold on a smooth quintic ninefold}\label{quintic:lines}
Every line $L$ on a smooth complex quintic ninefold satisfies
\[
                  \delta_L(X,H)\ge 10\frac{29}{45}
                              =\frac{58}{9}>\frac{32}{5}.          \tag{Q9}
\]
All constants below are in units of $N=n+1=10$. Any polar rank at
least four has ordinary lower bound $5/7$, and the rank-three
splitting $(2,2)$ has the same bound. At any polar rank contact two
has lower bound $15/23$ from (L1) with one light direction. It
suffices to treat the pencil and the exceptional rank-three splitting
$(1,3)$ in contacts three, four and five.

\paragraph{The pencil at contact five.}
Its transverse leading form is a smooth quintic $F\subset\mathbf
P^6$ by finite Gauss. Through every point lies a line in $F$, since
four forms of degrees $2,3,4,5$ have a common point in the $\mathbf
P^4$ of tangent directions. The projective line scheme has a covering
component. Its general member on the smooth fivefold is free, by
generic smoothness of evaluation in characteristic zero: evaluation
surjects on the normal fiber, hence the normal bundle has no
negative summands. The rank-four $N^*_{B/F}(1)$ of that line has
degree four and all its splitting degrees are at most one, so all
four equal one. Its five normal quartic polars are therefore linearly
independent (dualize the conormal sequence to count constant syzygies)
and give the complete binary quartic system. Choose a transverse
direction to $B$ giving binary polar $s^2t^2$. On its plane with
$B$, the initial expression of $F$ under normal-order degeneration
is $x s^2 t^2$. This is an instance of (L1) with $b=2$, yielding
$15/23$.

\paragraph{The pencil at contact four.}
On the seven light variables write $G_4=uP+vQ$ for quartics $P,Q$.
Write $C$ for the top quintic on the light $\mathbf P^6$. One may
assume $C$ smooth: the complete intersection $X|_{z_1=z_2=0}$ has
singular locus of dimension at most one, so a general codimension-two
section disjoint from $L$ is smooth. Express that section as
$u=L_1(x),v=L_2(x)$ and shift the base coordinates. This preserves
the leading normal form and the weighted computation. We shall choose
$e\in V(P,Q)\setminus C$ and a line through $e$ in the common tangent
space on which $P,Q$ do not both vanish identically. Its contact is
then $2\le b\le4$. The weighted mixed/pure model (L5) has $c'=4-b$
and its actual radial ratio is
\[
     \frac{1+(2+1/b)/5}{9/5+(2+1/b)/5}\ge 29/45 .
\]
On the complement of the root divisor its coarse mixed wall for
$c'\le2$ has threshold at least $15/23$, and the other wall has
denominator at most $5/4$. The coordinate boundary reduces to the
finite binary calculation as in (L5).

Here is the choice of $e$. The smooth quintic $C$ cannot contain
any component of the codimension-two complete intersection of two
coprime forms of degrees below five. Indeed it has Picard group
$\mathbb ZH$, so any such reduced component would be $(C,A)$
scheme-theoretically, and both lower-degree forms would be divisible
by $A$. The same fact holds on any fixed hyperplane $H'\simeq
\mathbf P^5$ for the quintic $C|_{H'}$: by finite Gauss it has only
isolated singularities, is integral and normal, and its dimension
four complete-intersection local rings are locally factorial by
Grothendieck parafactoriality at the isolated singularities.
The projective complete-intersection Picard theorem gives
$\operatorname{Pic}(C|_{H'})=\mathbb ZH$. In both cases section
restriction from projective space is surjective. A shared divisor
of lower degree cannot be contained in the corresponding integral
quintic either.

If $P,Q$ are coprime, their base is pure of dimension four. A nonlinear
component gives the desired general $e$: its common tangent has
dimension at least four, and if the whole tangent were in the base,
it would be a linear component through $e$. The same conclusion
holds at a general point on any component where the two derivatives
have rank less than two. In the remaining case all components are
generically reduced four-planes. The complete intersection has degree
16 and is unmixed, so there are at least two. Some two meet in a
three-plane: otherwise intersection with a general $\mathbf P^3$
would give a disconnected union of disjoint lines, although a complete
intersection curve of two positive-degree surfaces in $\mathbf P^3$
is connected. This three-plane is not contained in a smooth quintic
fivefold (three positive-degree transverse derivatives must have a
common zero). Choose $e$ in it off $C$. Its common tangent contains
the span of the two neighboring four-planes and so cannot be in
the four-dimensional base.

If a nonlinear prime divides $P,Q$, take a general smooth point on
it off $C$. The common tangent is either the whole space or its
tangent hyperplane. The latter varies nontrivially in characteristic
zero for a nonlinear hypersurface and cannot generally be one of the
finitely many hyperplane components of $V(P,Q)$. If a common linear
factor $l$ occurs, write $P=lP_1,Q=lQ_1$. On $l=0$ choose a common
zero of the cubics $P_1,Q_1$ outside $C$ by the hyperplane version
just proved (allowing zero polynomials and shared components). Both
original differentials vanish there, and the whole light space is
not their base. This also handles proportional $P,Q$.

\paragraph{The pencil at contact three.}
Write $G_3$ using its three quadratic-in-$u,v$ coefficient cubics.
A common zero not a common cone vertex gives a mixed leading cubic
of exponents $(1,2)$, hence threshold $15/23$. Otherwise the common
linear vertex $V$ has codimension $s\le3$ and equals their entire
common zero set. The forms depend only on $s$ active variables.
If the degree-four term does not vanish on $V$, retain a vertex
direction $y$ where it is nonzero and an active $x$. Their initial
powers $x^3,y^4$ give normalized radial ratio $95/143>29/45$ and
both finite ratios are larger; use the disjoint pure-power argument
\eqref{models:eq-pure-powers} with degree-$3,4$ symbols and ambient $d=5$. The two exponents
are coprime, and each kept polar can be chosen coprime to the other
coefficient where required. In general the radial expression for
coprime powers $p,q$ is $(1+1/p+1/q)/(2-1/d+1/p+1/q)$.

Suppose the degree-four term vanishes on $V$ and again take the top
quintic $C$ smooth by shifting $u,v$; this preserves that vanishing
since all the degree-three coefficients vanish identically on $V$.
The restriction of $C$ to $V$ is nonzero, since otherwise the join
with $L$ is a linear subspace of dimension at least five inside $X$.
We choose $y\in V\setminus C$ and a light direction $x$ outside $V$
annihilating the derivatives of the two degree-four coefficients
at $y$. If $s=3$ this is immediate by dimensions: the common kernel
in the seven-dimensional vector space has dimension at least five,
whereas the vector vertex has dimension four. For $s=2$ the two
derivatives in the two active directions form a $2\times2$ matrix
of cubics on $V=\mathbf P^4$. Its determinant is either identically
zero or a sextic, whose zero locus cannot be contained in $C|_V$.
Indeed the latter is an integral quintic by finite Gauss (singular
locus at most one), and six is not a multiple of five. At such a
zero the common kernel strictly contains the vertex. For $s=1$ the
two active derivatives are cubics on the hyperplane $V=\mathbf P^5$;
choose a common zero off $C$ by the preceding factorial hyperplane
argument. Now $x$ has a nonzero cubic coefficient by the common-zero
description. On the two retained directions $x,y$, every quartic
coefficient has $x$ exponent at least two, whereas the remaining
degree-five term contains nonzero $y^5$. Weights $1/3,1/5,1$ on
$x,y,z_i$ give the pure powers $x^3,y^5$: all intermediate quartic
terms have weight at least $2/3+2/5>1$. The radial ratio is
$23/35>29/45$ and the finite bounds are larger by the same coprime
two-power calculation.

\paragraph{The near-pencil net.}
In rank three with splitting $(1,d-2)$ the polars may be written
$(t\ell,-s\ell,b_0)$ in line parameters $[s:t]$, with degrees
$d-2,d-2,d-1$ on $\ell,\ell,b_0$ and $\ell,b_0$ coprime. Apply the
three-heavy-variable product comparison and moving-point calculation
of Chapter~\ref{ch:models}, with base polarization degree
one for a line in place of two for a conic. The normalized constants
are given here to specify every change. Put
\[
 r=(d-1)/(h-1),\quad R=rh/d,\quad u=1/r,\quad D=1+(d-1)/r=h.
\]
The pure light symbol $Q=x^h$ has coefficient of base degree $d-h$.
After clearing its denominators the enlarged $(k,p)$ blocks thus
have base degree $rp+O(1)$ (instead of $2rp$) and the same total
asymptotic dimension ratio $R$. Their line $\ell T+b_0z_3=0$ costs
$(d-1)$ base degrees. Off the marked point, the normalized low-summand
mean is $1+1/(r+d-1)$, since the three nonzero simplex weights in the
restricted line computation are $1,d-2,r$. The line discrepancy test
has denominator $Rh/D=R$. At the moving point $[s:t:0]$, an order-$j$
normal monomial $z_3^aT^b$ has coefficient-space degree
$rp-b+k-j$. Thus the two nonconstant mean jet exponents, radial
denominator and final exceptional-direction denominator are exactly
\[
 1,\quad 1/(1+u),\qquad
 R\frac{1+1/(1+u)}{1+1/h},\qquad RhJ,
\]
where
\[
 J=\frac{1}{1-u^2}\left(\frac1{1+(d-1)u}-\frac{u^2}{d-1+u}\right)
 \quad(u<1),\qquad dJ=1+\frac{d-2}{2d}\quad(u=1).
\]
Indeed the residual degree after a further direction order $q$ loses
$(d-1)q$ and after shifting the two exponential variables the
positive-part threshold is $(d-2)q$. For $u=1$ the two remaining
rate-one exponentials have positive-part mean $(2+t)e^{-t}$, proving
the formula by direct integration against $e^{-2q}$ at $t=(d-2)q$.
Lastly the normalized ternary-axis denominator is
$[1+h(3-1/d)]/(1+2h)$, and the initial radial denominator is no
larger, since additional free coordinates have unit cost.
For $d=5$, in rising $h=3,4,5$ the line-off-point denominators are
$7/5,\,19/15,\,6/5$ (the middle can be bounded by $19/15$), the
marked radial denominators are $3/2,\,704/525,\,5/4$, and the
special direction and ternary-axis denominators are all at most
$3/2$. The other points have denominator no larger than $R$.
Thus the net gives at least $2/3$ in these contacts. All these
estimates test the full geometric generic fiber of the extraction,
hence cover every ambient valuation with center containing the line.

\section{The conic threshold}\label{quintic:conics}
For \(d=5,n=9\), the conic calculation in Chapter~\ref{ch:models}
gives \(\delta_D(X,H)\ge 10(3295/5083)>58/9\). The numerator in
the cross-multiplied strict comparison is
\(9\cdot32950-58\cdot5083=1736\). Together with the line estimate
this gives the common mark used in the global argument.

% ===== The sextic elevenfold =====
\chapter{The sextic elevenfold}\label{ch:sextic}
\section{Statement and proof structure}
\begin{theorem}\label{sextic:main}
Let $X=X^ {11}_6\subset\mathbf{P}^{12}_{\mathbb C}$ be smooth and put
$H=\mathcal O_X(1)$. Then $\delta(X,H)>7$. In particular $X$ is
uniformly K-stable.
\end{theorem}
Adjunction gives $-K_X=7H$. For a curve $D$, we write
$\delta_D(X,H)=\inf_{D\subset c_X(v)} A_X(v)/S_H(v)$,
allowing equality in the containment. The positive-center and multiplier-ideal estimates are in Chapter~\ref{ch:global}, the curve estimates in Chapter~\ref{ch:models}, and the bottom-fraction jet lemma in Chapter~\ref{ch:quintic}. The global argument here is in hyperplane units. In the local calculations
the normalized line mark is $\delta_L(X,H)/12$.

\begin{center}
\begin{tabular}{ll}
\toprule
Part of the argument & Bound in hyperplane units\\
\midrule
Sufficient mark on each line and conic & $7433/1000$\\
Every integral conic & $4176/553$\\
Contact five, coprime light quintics & $12(42/67)$\\
Contact five, order $3$--$5$ off the transverse sextic & $12(41/66)$\\
Contact five, order $3$--$5$ transverse to that sextic & $12(204/329)$\\
\bottomrule
\end{tabular}
\end{center}
\section{The global short-curve criterion}

\begin{theorem}\label{sextic:criterion}
Let $X\subset\mathbf{P}^{12}_{\mathbb C}$ be a smooth sextic elevenfold.
If every integral line and conic $C$ satisfies
$\delta_C(X,H)\ge 7433/1000$, then $\delta(X,H)>7$.
Here $H=\mathcal O_X(1)$, and the valuations in the local threshold
have centers containing the given curve.
\end{theorem}
Use the two-short-curve Nadel argument from Chapter~\ref{ch:boundary} with the point jet improvement. A general smooth pointed
hyperplane section is a sextic tenfold of index six, uniformly
K-stable by the preceding row. Hyperplane basis adjunction gives
$\beta=\delta(X,H)>(12/11)6=72/11$. Suppose $\beta\le7$.
First consider $\beta\le33/5$. An original positive-dimensional
minimizing center contains a line or conic: in dimension $s\ge2$
the covering degree is at most $11(s+1)/(s\beta)<3$; in dimension
one it is at most $12/\beta<3$. Either contradicts the curve mark.
At a point optimizer, normalized by $w(\mathfrak m_p)=1$, one has
$S_H(w)=A_X(w)/\beta\ge11/\beta$. In the untilted short-curve
Nadel argument take $\Delta=1/4$ and
$\eta=7-(11/12)(\beta+\Delta)$. Then
$\beta+\Delta\le137/20<7433/1000$ and $1/3<\eta<1$.
The necessary inequality $1\ge\eta+\Delta S_H(w)$ contradicts
$\eta+11\Delta/\beta\ge91/80$, evaluated at $\beta=33/5$.
We may now assume $33/5<\beta\le7$. The program \path{r7-lit-sextic-finite.py} in the archive described in Chapter~\ref{ch:checks} uses
$f=1/6$, share $167/400$ and $\epsilon=10^{-8}$. The tilted
threshold satisfies $\beta'\ge\beta-(167/400)(\beta-4)-\epsilon
\ge11029/2000-\epsilon>11/2$. Thus a positive-dimensional
tilted optimizer again contains a line or conic by the same degree
formulas. The upper discrepancy ratio in the short-curve test
is strictly below the mark by the comparison in the program's \texttt{budgets} function, so the
optimizer is centered at a point. If $q=A_X(w)-11$
at the tilted optimizing point with $w(\mathfrak m_p)=1$, the Nadel
inequality is contradicted once $q>3.2306$. On $0\le q\le3.2306$ our
pointwise estimate implies $A_X(w)\le U S_H(w)$ with $U<7.027715$.

Take a general smooth pointed linear section $S$ of dimension two.
The tangent section on $X$ has multiplicity $2\le b\le6$ and isolated
singularities by finiteness of the Gauss map. For $b<6$ the section
is integral and is not a cone at $p$, and the general pointed
restriction cuts an integral plane sextic of the same multiplicity.
The scalar surface estimate $3S_{H|_S}\le B A_S$ implies
$12 S_H(w)\le B(q+2)+9$ by signed hyperplane adjunction.
At both endpoints $q=0,3.2306$ the upper estimate with $B=26/9$
contradicts $11+q\le U S_H(w)$.

\paragraph{Multiplicity two.}
We can use $B=26/9$. Consider an irreducible curve of degree $a$,
point multiplicity $k$ and $6k/a>B$ on the general pointed surface $S$.
It moves in the family of pointed surface sections of the general
pointed sextic threefold $Y$. The moving-curve genus inequalities proved
in Chapter~\ref{ch:global} apply when the curves sweep out $Y$:
$a\ge5$, $6k(k+1)\le a^2+6a$, and $2g-2\ge a+2k$.
B\'ezout with the integral tangent curve gives $a\ge2k$ unless that
curve is the tangent section itself; the latter has ratio two.
For any counterexample in the swept case the integer inequalities leave
only $(a,k)=(6,3)$: substituting $k>13a/27$ gives
$a<756/95<8$.
The genus lower bound at this pair and adjunction on the point blowup
give $D^2\ge6$, while Hodge gives $D^2\le a^2/6=6$.
Thus $D$ is numerically $H$. On a smooth hypersurface surface in
$\mathbf{P}^3$ numerical equivalence to $H$ implies linear equivalence:
the surface is simply connected, its Picard variety is zero, and
$H^2(S,\mathbb Z)$ is torsion free. Only the tangent hyperplane
section has multiplicity at least two, and its multiplicity here is two.

In the proper swept-locus case Chapter~\ref{ch:global} shows that the
curve $D$ is a general pointed section of an integral divisor on $Y$,
so $D\sim hH$ for an integer $h$. Its cone exception would give
tangent multiplicity six. Write $t=k/h$. Since $t>B$,
$h\ge2$ and $B<t\le3$. The positive and negative parts of
$H-sE$ on the blown-up surface are determined by $D$ once $s>6/t$.
Indeed every other integral curve obeys the nef bound $6/t$ by
intersection with $D$. Subtracting the forced copies of $D$ therefore
gives exactly the same angular tent computation as for a tangent
curve, rescaled by $h$. In particular
\[
 3S_{H|_S}(v)\le (6/t) A_S(v)+(1-6/t^2)v(D)/h .
\]
For completeness let $I$ be the discriminant exponent of the
multiplicity-$k$ curve in a general projection parameter. The local
delta-invariant and genus give $I\le6h^2+12h+k+1$.
In a splitting cover choose the closest Puiseux root for $v$, writing
$s$ for its contact divided by $v(x)$ and $e=v(D)/v(x)$.
The cluster counting function $m(u)$ on $[0,s]$ has integral $e$,
$1\le m(u)\le k$ and $\int m(u)(m(u)-1)du\le I$.
Writing $m=e/s$, Cauchy--Schwarz gives $s m(m-1)\le I$. Thus
\[
 \frac{e}{1+s}\le\frac{I m}{I+m(m-1)}
                 \le\frac{I}{2\sqrt I-1}\le 2h .
\]
For the middle inequality use $m^2+I\ge2m\sqrt I$.
For the last one, $1\le I\le6h^2+15h+1\le55h^2/4$, since $k\le3h$
and $h\ge2$. Consequently $2h-\sqrt{4h^2-2h}<1\le\sqrt I$ and
$\sqrt I\le (\sqrt{55}/2)h<(2+\sqrt3)h\le2h+\sqrt{4h^2-2h}$.
These two bounds give $I\le2h(2\sqrt I-1)$.
Since $A_S(v)\ge(1+s)v(x)$
by the finite-cover discrepancy argument, this proves
$v(D)/(h A_S(v))\le2$. Hence
$B'=6/t+2(1-6/t^2)\le8/3<26/9$.
If no curve violates $B$, the usual Seshadri surface bound supplies
$3S_{H|_S}\le B A_S$.

\paragraph{Multiplicity three.}
The plane sextic has arithmetic genus ten; its discriminant exponent
is at most $22$. The three-root estimate from Chapter~\ref{ch:quintic} gives $\operatorname{lct}_p(C)\ge1/3+2/22=14/33$ (or use the cluster proof
with $m(m-1)\ge2(2m-3)$).
The surface tent bound gives $B=2+(1/3)(33/14)=39/14$.

\paragraph{Multiplicity four.}
Write the tangent equation $T^2F_4+TF_5+F_6$.
If the projective quartic $F_4$ has no component of multiplicity
three or four, a general binary direction has root multiplicities at
most two. In the Puiseux count $e\le\min(4s,2s+2)$, so that
$\operatorname{lct}_p(C)\ge1/2$. The scalar bound is then $11/4$.
Otherwise $F_4=L^jQ$ for a linear form and $j=3$ or $4$.
If $L\nmid F_5$, the plane strict transform on the ordinary point
blowup is smooth at this multiple direction, with smooth contact
$j$ with the exceptional line. Its local threshold on $S$ is at
least $\min(2/4,(2+1/j)/5)\ge9/20$, by resolving the two smooth
branches in $j$ successive point blowups. All other directions are
simple or absent. This gives $B\le 3/2+(5/8)(20/9)=26/9$.

Suppose $L\mid F_5$. Write $j'=\operatorname{ord}_L(F_5)$ with
$1/j'=0$ if $F_5=0$. Weight the ten transverse variables $y$ by $a$
and $L$ by $a+h$, choosing positive coprime integers with
$h/a=\max(2/j,1/j')\le1$. The weights of all terms of the tangent
equation are at least $6a$, and the weight-$6a$ symbol contains
$G(y)=F_6(y,0)$. This projective sextic in $\mathbf{P}^9$ has at most
isolated singularities. Indeed on any positive-dimensional component
of its singular locus the remaining normal derivative is
$T(\partial_LF_5)+\partial_LF_6$. If the first coefficient does not
vanish identically, it has a root in $T$ at general points. Otherwise
the second coefficient has a zero and the whole corresponding
$T$-line is singular in the tangent section. Both violate finite
Gauss. This also excludes $G=0$. Degeneration to the cone over a
general smooth hyperplane section proves $(\mathbf{P}^9,G)$ lc.
The weighted symbol $T^5Z+$ (weight-$6a$ tangent terms) is prime
since it contains $G$. Scaling $L$ out of the angular series leaves
the same primitive relation $T^5Z+G$. As in the preceding weighted
quintic calculation, normal monomials $T^i L^k Z^j P_l$, with $P_l$
in a graded complement to $(G)$, give normalized mean $L$-order
$1/12$, mean fixed $G$-order $5/72$ and residual complete $y$-degree
$10/12$. On $\mathbf{P}^9$ the angular coefficient is at most
$5/72+(10/12)/10=11/72$ times discrepancy. At the $L$-axis the
ordinary blowup has discrepancy ten and mean radial degree
$(10/12)+6(5/72)=15/12$, below $10/(72/11)$.
Extraction adjunction gives with $c=72/11$
\[
 A_X(w)-cS_H(w)\ge(11a+h-c(16a+h)/12)/(a+h).
\]
For $A_X(w)\le U S_H(w)$ this implies $A_X(w)>14.231$ throughout
$0<h/a\le1$.

\paragraph{Multiplicities five and six.}
For $b=5$ the tangent section is the unique maximal multiplicity
prime by B\'ezout. The general principal-bin bound with fixed
coefficient upper $v(T_p)\le6+q$ excludes $0\le q\le3.2306$.
The function \texttt{bin\_check} in \path{r7-lit-sextic-finite.py} verifies the specified rational nonuniform partition with margin
greater than $0.013$ at the upper endpoint.
For $b=6$ the tangent is a cone over a smooth sextic ninefold.
The latter has index five and is K-stable by $n\ge2r-1$.
The Eckardt local formula gives $\delta_p(X,H)=132/16>U$.
This exhausts the possible multiplicities.

\section{Conics and lines before the contact-five case}
For an integral conic over $\mathbb C$ the plane curve is smooth.
The degree-six conic estimates in Chapter~\ref{ch:models} (proper plane section
and contained plane, including the polar net cases) give
$\delta_D(X,H)\ge4176/553$. For $d=6$, \path{r6-flag-conic-general-check.py} verifies the
upper bound $553/348$ for the normalized denominator, giving
$12/(553/348)=4176/553$. We consider lines in the rest of the proof.

The normal polars along a line $L$ are the restrictions of the
normal partial derivatives of the defining sextic. They are binary
quintics without a common zero on $L$. The dimension of their span
is the polar rank. In constant normal coordinates, the normals
corresponding to independent polars are called heavy coordinates,
and those in the kernel are called light coordinates. The least
nonzero light degree on the common heavy-zero space is the contact
order. For a pencil this space cannot be contained in $X$, by the
linear-space bound for smooth hypersurfaces.

For rank at least four, or rank three with positive conormal
splitting $(2,3)$, the ordinary estimate gives
$\delta_L/12\ge6/9$. For the other net the positive splitting is
$(1,4)$. Use the moving polar calculation in Chapter~\ref{ch:quintic}
with $d=6$. For $h=3,4,5,6$ substitution gives a maximal
denominator at most $45/28$ for its four bounds. The contact-two
single-power estimate (L1) is $9/14$. Thus every net gives
$\delta_L\ge112/15$. For a pencil at order six the transverse sextic
in $\mathbf{P}^8$ is smooth: a singular point gives a singular plane on
the two-heavy section, contrary to its singular-locus dimension
bound one by finite Gauss. This sextic sevenfold has index three. Its dimension satisfies
\(7\ge 2\cdot3-1\), so Chapter~\ref{ch:global} gives uniform K-stability. In (L3) take nine light variables and
use $\delta(Y,H)\ge3$; the lower ratio is $9/14$.
The next calculation treats the two intermediate orders three and four.
\paragraph{A Newton edge at a common vertex.}
Let a pencil of normal quintics along a line of a sextic elevenfold
be fixed. Retain two heavy variables $z_1,z_2$ with coprime
polars $f_1,f_2$, and light coordinates $x,y$ along a two-plane,
by the actual hyperplane flag through the line. Write $G$ for
the light restriction. Suppose its lowest total degree is $p$ and
its degree-$p$ part consists of $a x^p$ for $a\ne0$ a binary
base coefficient. Assume a higher term involving $y$ occurs.
The binary Newton polygon then has a first lower edge starting
at $(p,0)$ and ending at $(i,j)$ with $i<p$, provided a higher
term of $x$ exponent $<p$ occurs. Write $(p-i)/j=l/m$ in
lowest terms; $0<l<m$. Put $k=j/m$. Give $x,y,z_i$ weights
$m,l,pm$. On the weighted exceptional, the initial relation is
\[
      f_1z_1+f_2z_2+\sum_{u=0}^{k}a_u x^{p-lu}(y^m)^u=0,
                    \qquad a_0a_k\ne0 .
\]
Primitivity plus the ordinary Hilbert series computes all
initial pieces, also in all subsequent specializations.
In fixed weighted degree the exponent of $y$ has a unique
remainder modulo $m$, since $(m,l)=1$. Thus $V=y^m$ gives coarse
weighted plane $\mathbf P(1,l,p)$ on $x,V,z_1$, with root
boundary $V=0$ and section costs $1,m,1$. The bounded remainder
has asymptotically zero normalized effect.

Off $V=0$, drop its root boundary and diagonally specialize to
$a_k x^i V^k$ (give $x$ weight one, $V$ zero and both $z_i$
weight $i$). In the mixed-monomial calculation the finite binary
coefficient is at most $\max(1,5/6+1/(i+mk))$. By exact
divisibility the two coordinate coefficients are at most $1,1/m$.
On the mixed Newton wall the normalized mean is
$u+v/m+(11/6)w$ when $w=iu+kv$, with discrepancy $u+v+w$.
Its two maximal denominators are therefore
\[
       \frac{1+(11/6)i}{1+i},\qquad
       \frac{1/m+(11/6)k}{1+k}.
\]
Use just the nontrivial wall if $i=0$. The common affine axis
and the coordinate strata are included in the mixed-monomial
calculation, so its bound tests the whole coarse model and
semicontinuity applies. On $V=0$ take the actual $y$ filtration
upstairs. The quotient relation $f_1z_1+f_2z_2+a_0x^p$ is
primitive, giving coefficient one for $y$; the residual weighted
projective binary line has finite coefficient at most
$\max(1,5/6+1/p)$ and infinity coefficient one by exact
$x$-divisibility. No common affine axis lies on this line.
Adjunction includes the original radial threshold ratio, in units
of $12$, equal to
\[
       \frac{1+(1+l/m)/p}{11/6+(1+l/m)/p}.
\]

For contact three, all pencil cases give
$\delta_L(X,H)/12\ge17/27>149/240$. The four cubic coefficient
forms of the light third-order term have a common zero. If some
common zero is not a common cone vertex, a general line through it
has finite contact one or two, giving the ordinary mixed-monomial
bound with largest exponent at most two. Otherwise their zero set
is their common linear cone vertex $W$ of projective dimension at
least four. Choose general $y\in W$ and general active $x$. At
least one higher restricted monomial has $x$ exponent $<3$:
otherwise the quartic, quintic and sextic light coefficients all
vanish to order at least three along $W$ (use the Taylor expansion
in arbitrary $x$ at general $y$). Then the whole join of $W$ and
the original line would be singular on the two-heavy section,
contradicting finiteness of the Gauss map, which bounds its
singular locus by one. The desired detectable monomial and
$a_0\ne0$ are simultaneous open conditions. Apply the edge
calculation. Here $i\le2$, $k\le3$, $m\ge2$ and
$l/m\ge1/4$ (use $j\le6-i$). The wall denominators are at
most $14/9$ and $3/2$. The radial ratio is at least $17/27$.

The same argument for contact four works whenever the first
edge after $x^4$ has far endpoint $i\le2$. Here $l/m\ge1/2$,
$k\le4$, and the radial ratio is at least
$(1+3/8)/(11/6+3/8)=33/53>149/240$.
The second wall is at most $47/30$, and the first at most $14/9$.

\paragraph{Tangent directions in contact four.}
Let $C$ be the smooth top sextic on the light $\mathbf P^8$,
obtained by the general base-coordinate shift. The three coefficient
quartics $q_1,q_2,q_3$ have a common zero off $C$. Indeed, if their
tuple were a nowhere-zero section of $\mathcal O(4)^3$ there, then
$64H^3=0$ in $\mathrm{CH}^3(\mathbf P^8\setminus C)$. Localization
would give a cycle of codimension two on the smooth sevenfold $C$
having degree 64. Weak Lefschetz gives $H^4(C,\mathbf Z)=\mathbf ZH^2$,
so its degree must be divisible by six.

In fact at some such $e$ the common projective tangent linear
space $T$, annihilated by the $dq_i(e)$, contains a line through
$e$ not contained in $B=V(q_1,q_2,q_3)$. Choose a maximal-dimensional
component $W$ meeting $\mathbf P^8\setminus C$. Unless the general
tangent space $T$ is contained in $B$, the claim follows. If it
is contained, it lies in $W$ at a general point of $W$ and contains
the tangent space to $W$; equality of dimensions implies $W$ is
a linear space $\mathbf P^{8-c}$, $1\le c\le3$. The derivatives
transverse to $W$ give a $c\times3$ matrix of cubics on $W$.
This matrix drops rank somewhere outside $C$. Indeed choose a
general $A=\mathbf P^{8-2c}$ in $W$; the singular locus of
$C\cap W$ has dimension at most $c-1$ since the finite Gauss
map of $C$ sends it to the $\mathbf P^{c-1}$ of hyperplanes
containing $W$. Then $D=C\cap A$ is a smooth sextic hypersurface
of dimension $7-2c$. The matrix cannot induce a surjection
$\mathcal O^3\to\mathcal O(3)^c$ on $A\setminus D$, as its kernel
of rank $3-c$ would have $c_{4-c}=0$ whereas the inverse total
Chern class $(1+3H)^{-c}$ gives $-27,27,-9$ times $H^{4-c}$
for $c=1,2,3$ respectively. None is divisible by six; localization
and weak Lefschetz again give divisibility for cycles of the
required codimension $3-c$ on $D$, since $2(3-c)<\dim D$.
At a rank-drop point off $C$, $T$ has dimension greater than
$\dim W$, so it cannot lie in $B$ by maximality of $W$.

Retain a light line in $T$ with nonzero quartic restriction,
with coordinates $y$ at $e$ and $x$ transverse to it. Let $b$ be
the least $x$ exponent in the restricted quartics, so $2\le b\le4$.
For $b=2$ the mixed quartic $x^2y^2$ gives $\delta_L/12\ge9/14$.
Otherwise let $p=5$ if at least one quintic coefficient is nonzero
at $e$, and $p=6$ if both vanish. The pure $y^p$ coefficient is
nonzero. Take the Newton edge abutting $(0,p)$. To find its far
endpoint maximize $(p-j)/i$ over all other exponent pairs $(i,j)$
with $i>0$, breaking equality by larger $i$. Among any fixed degree
below $p$ the maximum occurs at the least $i$. For $p=5$ the far
endpoints are $(b,4-b)$. For $p=6$ writing $s$ for the least
$x$ exponent in the fifth degree, one compares $1+2/b$ with
$1+1/s$. For $b=3,4$ this gives $(b,4-b)$ or $(1,4)$, with the
latter precisely when $s=1$. When $b=s=2$ either the original
mixed quartic already suffices or one takes the larger quartic
endpoint under the tie rule.

All endpoints $(b,4-b)$ with $b=2,3,4$ give the target
$\delta_L/12\ge149/240$. Set $g=\gcd(b,p-4+b)$,
$m=b/g$, $l=(p-4+b)/g$ so the weights on $x,y,z$ are
$(l,m,mp)$. After taking $U=x^m$, off the root boundary
$U=0$ specialize the edge to its extreme term $U^g y^{4-b}$.
The wall denominators are at most $(1/m+(11/6)g)/(1+g)$
and $(1+(11/6)(4-b))/(5-b)$. In these ranges $g\le2$,
$4-b\le2$ and the bounds are at most $14/9$. On the boundary
the actual $x$ quotient is pure $y^p$. When $m=1$ (only $p=6,b=2$)
the same argument needs no root. The full radial ratio is at least
$33/53$, since $l/(mp)+1/p\ge1/4+1/6$. Thus the only case left
by this argument is $p=6$, $b\ge3$, $s=1$, when the exposed
light monomials are $xy^4$ and $y^6$. If any choice of allowed
$(e,x)$ kills the two quintic differentials the other endpoints
give the asserted threshold. Here is the remaining calculation.

Suppose $b=3$ or $4$ and the quintic restriction has least
$x$ exponent one. Use the opposite edge, between $x^b y^{4-b}$
and $xy^4$. Their coefficients $a,c$ have base degrees two and
one, respectively. For $b=4$ give $x,y,z$ weights $(4,3,16)$.
The boundary $y=0$ has pure residual $x^4$; modulo the bounded
$y$ remainder the coarse variable $V=y^4$ has section cost four.
Off the boundary the diagonal limit $xV$ has both wall exponents
one. The radial ratio is $23/(7+16\cdot11/6)=69/109$.

For $b=3$ give $x,y,z$ weights $(3,2,11)$. The initial light
form is $a x^3y+cxy^4$, and the full radial ratio is $96/151$.
On the angular open $xy\ne0$, use ordinary projective coordinates
\[
             U=x^2,\quad V=y^3,\quad W_i=z_i/(xy).
\]
The point in $\mathbf P^2$ is $[U:V:W_2]$, with
$f_1W_1+f_2W_2+aU+cV=0$. This identifies the two open
sets $xy\ne0$ and $UV\ne0$: the two invariant ratios
$x^2/y^3,z_2/(xy^4)$ are a basis of the character lattice modulo
the original $(3,2,11)$ weights. In angular weight $q$, after
substitution $z_i=xyW_i$ the common remainders of the $x,y$
exponents modulo $2,3$ are $q\bmod2,2q\bmod3$. Removing
them gives complete polynomials in the coarse coordinates.
Take a feasible diagonal action with $V$ the unique lowest
coordinate, choosing $f_1$ coprime to $c$, and separate $U,W_2$
as well. Its special relation is $f_1z_1+cxy^4$ in the original
ring. It is primitive and has the original ordinary Hilbert
function. The resulting coordinate means are
\[
   \bar U=(3-1/6)/2=17/12,\quad
   \bar V=(1+4+1-4/6)/3=16/9,\quad \bar W_2=1 .
\]
The diagonal can equivalently be pulled to the original weights
after clearing denominators: assign half the $U$ weight to $x$,
a third of the $V$ weight to $y$, and their sum plus the $W_2$
weight to $z_2$. The monomial Hilbert derivative subtracts
the relation weight divided by six. These computations retain
the original complex coefficient dimensions.

On $UV\ne0$ a fixed power of $V$ can be removed from every
weight-$q$, degree-$M$ system without affecting any threshold
there. In fact every representative monomial contains at least
$(q-3M)_+/8$ heavy variables: the maximal other-symbol weight
per section degree is three, and each heavy has weight eleven.
Each heavy contributes a factor $y$ before passage to the
coarse coordinates, giving common $V$ exponent at least
$(q-3M)_+/24-O(1)$. Remove this power before applying a global
diagonal degeneration on the ordinary coarse plane to the
resulting honest polynomials. The actual hyperplane flag retained
two light and two heavy variables, so the new hypersurface
dimension is four and its normalization factor is five.
The $5$-normalized expected $q/M$ is $151/6$, hence the
normalized mean removed power is at least
$(151/6-15)/24=61/144$. The remaining $V$ coefficient is at
most $16/9-61/144=65/48$; all three coarse coefficients
are smaller than $240/149$. Semicontinuity gives the desired
bound at every light-torus center.

For completeness, when $x\ne0,y=0$ use the projective
diagonal with $y,z$ of positive weight one and $x$ of
weight zero. Its primitive special relation has light term
$a x^3y$. At $y=z=0$ on the covering chart $x=1$ the
residual factors are $z-\lambda y$ and the relevant mixed
wall denominator is $17/12$. A point with $y=0,z\ne0$ has
inverse orbit tending to the heavy vertex as $t\to0$, and
there the finite factors $z-\lambda x^3y$ are units. The
coordinate factors have coefficients at most one and meet
transversely on its quotient covering chart. Properness of
the projective family or local semicontinuity along the
limiting point gives the bound. The heavy vertex itself is
fixed. For $x=0,y\ne0$ take the diagonal with $x,z$ of weight
one, $y$ of weight zero, specializing instead to $cxy^4$;
the identical chart calculation has exponent one in $x$.
This uses inverse orbits only on the indicated boundary.
Thus every sextic elevenfold pencil with contact four has
normalized threshold strictly greater than $149/240$.
\section{The contact-five case for a sextic elevenfold}
Let \(L\) be a line on a smooth sextic elevenfold whose five-degree
normal polars span a pencil. Put \(u,v\) on \(L\). After choosing
two heavy normal coordinates spanning the pencil, consider the
restriction of the sextic to their common zero space. Contact five
means this restriction has the form
\[
                 uP(x)+vQ(x)+C(x),\quad x=(x_1,\ldots,x_9),
\]
where \(P,Q\) are quintic forms, at least one is nonzero, and \(C\)
has degree six. After a linear shift of \(u,v\), one can take the
transverse sextic \(C\subset\mathbf P^8\) smooth. Indeed, the
singular locus of the section by the two heavy coordinates has
dimension at most one by the finite Gauss map. A general
complementary \(\mathbf P^8\) misses that locus and is transverse
on the smooth locus by Bertini.

The mixed/pure model (L5) of Chapter~\ref{ch:models} applies when some
$e\in V(P,Q)\setminus C$ has finite binary contact $b=3,4$ or $5$.
With $j=5-b$ its second wall denominator is at most
$[1+(11/6)j]/(1+j)\le14/9$ and its first wall is smaller.
The radial ratio is minimized at $b=5$, giving $41/66$.
We next allow the mark $e$ to lie on $C$. This is used in both
shared-factor arguments.
\paragraph{Third contact at a point of the top sextic.}
Let $e\in C$ and choose a light line whose contact with the penultimate
quintic pair at $e$ has finite smaller order $b\in\{3,4,5\}$.
Suppose that line is transverse to the smooth sextic $C$ at $e$.
Then $\delta_L(X,H)/12\ge 204/329>7433/12000$. In particular
the existence of any point where both quintics have multiplicity
at least three suffices, whether or not it lies on $C$: take a
general line there, and use (L5) if the point is outside $C$.

Use binary lights $x,y$ with $e=[0:1]$, and retain both heavy
polar variables $z_1,z_2$ by the actual hyperplane flag. Put
$a=b-1$, $h=5-b$ and $w=6b-5$. Give $x,y,z_i$ weights
$b,a,w$. Among the quintic terms the minimum occurs at $x^b y^h$;
among the sextic terms it occurs at $xy^5$ by transversality.
Both weights are $w$. Thus the initial relation is
\[
      fz_1+gz_2+c x^b y^h+kxy^5=0,
\]
where $f,g$ are the two basepoint-free quintic polars, $c$ is a
nonzero binary linear form on the original line and $k\ne0$ is
constant. Primitivity and localization at the base compute all its
initial pieces. The radial ratio in units of the new dimension plus
one (five) is $(2b-1+w)/(2b-1+11w/6)$. Its three values are
$108/173,156/251,204/329$, respectively.

On the angular open $xy\ne0$, set $U=x^a$, $V=y^b$ and
$W_i=z_i/(xy^h)$. These have common weight $ab$, with relation
$fW_1+gW_2+cU+kV=0$. This identifies the light torus chart
with $UV\ne0$ in an ordinary projective plane with coordinates
$[U:V:W_2]$: $U/V$ is primitive on the light torus, and
$W_2/V=z_2/(xy^5)$ has heavy exponent one. In weight $J$ the
remainders of the total $x,y$ exponents modulo $a,b$ after putting
$z_i=xy^hW_i$ are fixed. Removing these bounded common remainders
leaves ordinary homogeneous polynomials. In each fixed section degree
choose integral diagonal weights on the coarse plane with $V$ its
unique lowest coordinate and avoiding the finitely many hyperplanes
where two distinct equal-degree exponent vectors have equal weight.
This makes all individual residual initial sections monomials;
their coefficients over the generic base still carry their original
complex multiplicities. The original special relation is
$fz_1+kxy^5$, primitive and prime since $k$ is a nonzero constant.
The resulting coordinate means are
\[
     \bar U=(3-1/6)/a,\qquad
     \bar V=(1+2h+b-5/6)/b,\qquad \bar W_2=1.
\]
For example $z_1=xy^h$ times a term with initial $V$, and the second
heavy contributes one $xy^h$; the relation subtracts $1/6$ of
the exponents in $xy^5$. The diagonal lifts to original symbols by
giving $x,y$ respectively $1/a,1/b$ times the assigned $U,V$
weights, giving $z_1,z_2$ the weights of $xy^h V$, $xy^h W_2$.
The ordinary Hilbert function and the complex base multiplicities
are preserved. All displayed means are in the dimension-plus-one
normalization.

In degree $m$, weight $J$, every ordinary monomial has at least
$(J-bm)_+/(w-b)$ heavy factors. Each heavy contributes $y^h$.
Before the coarse diagonal specialization remove the common power
$\max(0,\lfloor h(J-bm)_+/[b(w-b)]\rfloor)$ of $V$, a unit on the
light torus. The bounded light remainders introduce at most a
constant error per section. The normalized mean heavy lower bound
is at least $T_b=(2b-1+11w/6-5b)/(w-b)$, by convexity and the mean
of $J/m$. Thus the remaining $V$ coefficient is at most
$(1+2h+b-5/6-hT_b)/b$, namely $22/15,53/45,31/30$ for the three
values of $b$. The $U$ coefficient is at most $17/12$. These and
the free coefficient one all lie below $329/204$. Upon removal of
the units the sections are global homogeneous coarse polynomials.
The coordinates of their diagonal limit have simple normal crossings
on the full projective plane, so threshold semicontinuity gives the
claimed estimate on the torus with no limiting-center restriction.

At $x=0,y\ne0$ give $x$ and the heavies weight one and $y$ weight
zero. The special light term is $kxy^5$. Points with $z=0$ remain
fixed; at such points the mixed wall of (L1) on the covering chart
$y=1$ has exponent one in $x$. Points with $z\ne0$ on this boundary
have inverse orbit limiting to the heavy vertex; on that chart the
finite factors are units and the remaining coordinate divisors have
coefficients at most one. The heavy vertex itself is dealt with in
the same chart. At $y=0,x\ne0$ use weights one on $y$, $h$ on the
heavies and zero on $x$. The special term is $cx^b y^h$ and the
mixed wall in the chart $x=1$ has exponent $h\le2$, with denominator
at most $14/9$. For $h>0$ the same fixed and heavy-vertex limits
apply. For $h=0$ both $x$ and the heavy coordinate have weight zero,
so the whole boundary is fixed. Its residual binary systems for
$fz_1+gz_2+c x^5$ have finite coefficient at most
$\max(1,5/6+1/5)=31/30$ by the pure case of (L1), while $y$
is free. These arguments restrict adjacent finite factors by their
actual exponent on the indicated unit charts. Adjunction from the
angular stack yields the stated bound.
\paragraph{A coprime pencil of quintic light forms on a sextic elevenfold.}
Suppose the two quintic light forms $P,Q$ have no common factor in
$\mathbf P^8$. A general light line misses their common zeros, and
their restrictions to it give a basepoint-free binary pencil.
Cut by the linear flag retaining this line of directions, the two
heavy coordinates and the original line. It is smooth along the
original line because the two heavy polars $f,g$ form a basepoint-free
binary quintic pencil there. The resulting hypersurface has dimension
four, with two independent light normals $x,y$ and one independent
heavy normal. The quartic-sixfold filtered flag argument lifts every
bound in units of dimension plus one, up to a cutoff of one. Use
units of five on this fourfold. Give $x,y$ weight one, $z_1,z_2$
weight five and $u,v$ weight zero. The initial equation is
\[
       f(u,v)z_1+g(u,v)z_2+G=0,\qquad G=uP(x,y)+vQ(x,y).
\]
This is primitive and embeds over the generic base. The weighted
radial discrepancy-to-mean ratio is $7/(2+10-5/6)=42/67$.
The angular surface is $\mathbf P(1,1,5)$ over $K=\mathbf C(u/v)$.
We show its normalized threshold at least $12/19$.

The heavy vertex $x=y=0$ on its smooth uniformizing chart can be
blown up in the two light variables. The initial ring has equation
$fz_1+gz_2=0$, primitive and prime, with both lights free. Their
combined mean is two, equal to the radial discrepancy, and its
angular binary systems are full with two normalized coordinate
coefficients one. The usual two-filtration basis comparison and
inversion of adjunction give normalized threshold one there.
Likewise the common heavy zero $z_1=z_2=0$ is a separable degree-five
closed point over $K$, since $G$ is irreducible over $K$. In its
normal filtration the ring has primitive equation $G=0$ and both
heavy coordinates are free. Their combined mean is two, and the
independent binary systems on the blowup give the same bound.

\paragraph{The divisor of $G$ away from these two loci.}
In any $(m,J)$ section degree and angular-weight block, write each
section over $K$ as a polynomial in $z=z_1$ after eliminating $z_2$.
Use the largest $z$-degree $j$ as an increasing filtration. The
leading-in-$z$ ring has prime equation $fz_1+gz_2$ and all the
initial pieces by Hilbert series and injection. Thus the scaled mean
of $j$ is $2-1/6=11/6$ and that of $J$ is $67/6$. The coefficient
of $z^j$ has binary light degree $J-5j$. A section divisible by $G^i$
over the generic light ring therefore satisfies $5i\le J-5j$.
Use a simultaneous complex basis for these two filtrations, showing
the scaled mean of $i$ is at most $2/5$. This controls each geometric
prime above $G$ since the roots are separable and all conjugate.

Before eliminating $z_2$, every monomial on the block uses at least
$j_0=\max(0,\lceil(J-m)/4\rceil)$ heavy variables, by its ordinary
cost. Over $K$, substituting $z_2=-(fz+G)/g$ shows that the
coefficient of $z^a$ has a factor $G^{(j_0-a)_+}$. Consequently
the restriction of a section of $G$-order $i$, after removing
$G^i$, has $z$-degrees between $\max(0,j_0-i)$ and its actual top
degree $j$. At a finite nonzero $z$ the restricted vanishing is at
most $j-j_0+i$. Indeed the mandatory common monomial is a unit
there. To control arbitrary adapted residual bases at each $i$,
use the filtration from maximum degree $j$ on the quotient at $i$.
Its weighted dimensions, summed over $i$, give exactly the original
degree mean, by the simultaneous-basis identity. Since the scaled
mean $j_0$ is at least $(67/6-5)/4=37/24$, the scaled residual cost
is at most $11/6-37/24+2/5=83/120$. Restricting to $G$ is legitimate
on the nonvertex charts since its geometric components are smooth
and disjoint there. With its actual coefficient at most $2/5$,
adjunction gives threshold at least one at these centers.

\paragraph{The complement of $G$ and the heavy vertex.}
For a test valuation choose a geometric closed point on its center
with light vector nonzero and $G\ne0$, and spread this point over
a finite cover of the line parameter. A general base point $p$
specializes it to a point with the same properties, avoiding
finitely many poles. Change $[u:v]$ so that $p=[1:0]$, normalize
the polars so $f(1,0)=1$, $g(1,0)=0$, and use $P,Q$ for the new
binary quintics. The quintic $P$ is squarefree since a general member
of a coprime binary pencil is separable. On the specialized angular
plane the solved symbol $z_1$ is $-P$ in the frame $u=1$. Change the
free coordinate $z_2$ by a constant multiple of $z_1$ so that $z_2$
vanishes at the specialized point. This is possible since $P$ is
nonzero there, and preserves the section cost of one. It preserves
$f(1,0)=1$ and $g(1,0)=0$.

Write $W=u^4z_1+P$ of section degree five and angular degree five.
It acquires base order one by the original equation. For the
filtration by $v$-order of generic angular coefficients the
base-specialized presentation is the complete intersection
\[
 u^4z_1+P=0,\qquad
 uW+v(Q+c_1u^4z_1+c_2u^4z_2)=0 .
\]
Indeed $W$ is solved to first order by the original equation, with
$c_1,c_2$ constants. On inverting $u$ this quotient is an ordinary
polynomial ring. At $u=0$ it has $P=0,vQ=0$, of codimension two
among remaining variables, since $P,Q$ are coprime. Thus there is
no component above $u=0$. The complete intersection is integral,
its symbols inject in the generic angular ring with the ordinary
and base frames retained, and its ordinary graded Hilbert series is
that of the original sextic. This gives all initial spaces in every
grading. Now filter by the exponent of $z_2$ with its lowest power
first. The system has the same presentation with $c_2=0$ by the same
argument. This second torus specialization fixes the chosen point
because $z_2=0$ there. In applying semicontinuity one first chooses
a simultaneous basis for the test valuation and base order. To
bound the resulting specialized pair near the chosen point, bound
the maximum for any valuation there by bases adapted to it and the
second torus filtration in each base-order piece. This suffices for
the original pair by relative lower semicontinuity along the chosen
multisection on the proper angular surface. With finite-level split
coefficients independent of $p$, all limiting bounds are uniform.

In the latter ring put $Q'=Q-c_1P$ and use $z_1=-P$, $W=-Q'$ in
the specialized frame. The free $z_2$ has normalized mean one.
Fix a block with $z_2$ exponent $j$, base order $i$, light angular
degree $D=J-5j$ and total section degree $m$. Its binary system
is the sum of all
\[
 P^c (Q')^t\mathbf C[x,y]_{D-5c-5t},\qquad
 0\le t\le i,\qquad 4c+t\ge i+j+D-m.
\]
Indeed the underlying monomial has $z_1^cW^t v^{i-t} z_2^j$ and
$u$ exponent $m-(i+j+D)+4c+t$. The forms $P,Q'$ are a homogeneous
regular sequence, and $\mathbf C[x,y]$ is graded free over
$\mathbf C[P,Q']$ on lifts of the quotient's homogeneous basis
in degrees at most eight. For any such basis element of degree $l$
compatible with $D$, put $k=(D-l)/5$. A monomial of $P$ exponent
$q$ and complementary $Q'$ exponent $k-q$ is present precisely when
\[
  4q+\min(i,k-q)\ge i+j+D-m.
\]
The left side increases with $q$. For the possible $l$, the integer
$k$ differs by at most one. Thus each component contains all exponents
from a lower cutoff onwards, and these cutoffs differ by at most one,
allowing cutoff $k+1$ for an empty component. With $s$ the least
actual $P$ exponent in a nonzero block, the block contains
$P^{s+1}\mathbf C[x,y]_{D-5(s+1)}$ and lies in
$P^s\mathbf C[x,y]_{D-5s}$. In particular after its fixed part it
has codimension at most five in the complete binary moving system.

For computing the *mean* fixed exponent, use a monomial basis in
the complete-intersection presentation. Equivalently use the ideal
$z_1$. Modulo it the equations become $P=0,uW+vQ'=0$. This quotient
injects after inverting $u$ since the fiber at $u=0$ has codimension
one inside it, by the same coprimality argument. Thus divisibility
of a binary section by $P^a$ in this evaluated, fixed block agrees
with its algebraic divisibility by $z_1^a$ in the graded ring. In
particular the mean of the minimum block exponents $s$ is at most
the scaled mean of actual $z_1$ exponents of an adapted basis, which
is one. The scaled mean of $D$ is $67/6-5=37/6$. At any point of
the light binary projective line, the scaled mean vanishing is at
most $\overline{s}+(\overline D-5\overline{s})/2+o(1)$ if it is
a root of $P$, and at most $(\overline D-5\overline{s})/2+o(1)$
otherwise. Here the uniform $O(1)$ error in unnormalized binary
means comes from the codimension-five comparison. The amount of
*additional total* $P$-divisibility in a nearly complete degree-$q$
system must also be accounted for in determining $\overline{s}$.
For a complete binary degree-$q$ system it is $q/10+O(1)$, by
summing dimensions in the $P$-power filtration. The bounded
codimension changes this by $O(1)$. Hence
\[
  1=\overline{s}+(\overline D-5\overline{s})/10,
  \qquad \overline{s}=23/30.
\]
The resulting light-point bound is $\overline{s}+(37/6-5\overline{s})/2
=29/15$, which does not give the sextic elevenfold target at a root
of $P$. At the selected specialized point, however, $P$ is a unit.
There the bound is $(37/6-5\overline{s})/2=7/6$. The independent free
coordinate $z_2$ has coefficient one; the heavy vertex is absent
at that point. This completes the complement using the local bound
$7/6$ and the actual-$G$ calculation at its roots. In particular the
normalized threshold in the claimed coprime case is at least $42/67$.
\paragraph{A common factor of degree two, three, or four.}
Suppose that the homogeneous quintic pair has an irreducible common
factor $R$ of degree $2\le s\le4$ and has no common linear divisor.
Write $P=RA$, $Q=RB$, where each nonzero cofactor has positive degree.
Take any irreducible component $Z$ of $V(R,A,B)$. It is nonempty
of dimension at least five by projective dimension; a zero cofactor
does not affect the argument. If $Z\subset\operatorname{Sing} R$,
any $e\in Z$ has ambient multiplicity at least three for $(P,Q)$.
Take a general line there, using (L5) off $C$ and the
transverse-on-$C$ estimate on $C$.

Otherwise take a general $e\in Z$. The tangent hyperplane $T_eR$
is contained in neither $R$ nor $V(A,B)$, since $R$ is irreducible
of degree greater than one and the two cofactors have no common
linear divisor. A general line through $e$ in $T_eR$ has finite
contact between three and five with the quintic pair. When
$e\notin C$, (L5) applies. If $Z\subset C$, the two tangent
hyperplanes at general $e$ are distinct. Otherwise the nowhere
vanishing gradient $\nabla C$ on $Z$ would express $\nabla R$ as
its multiple by a global section of $\mathcal O_Z(s-6)$. Indeed
the proportionality equations hold everywhere on the reduced $Z$
and the local ratios patch wherever a partial of $C$ is nonzero.
This negatively twisted section vanishes since $Z$ is integral and
positive dimensional, contradicting generic smoothness of $R$ along
$Z$. Thus the line in $T_eR$ can also be chosen transverse to $C$,
and the transverse-on-$C$ estimate applies. The case of an additional
linear divisor is part of the common-hyperplane argument.\paragraph{A common irreducible factor of degree at least three.}
Let $P,Q$ be the two light quintics in $\mathbb P^8_{\mathbb C}$
at penultimate contact on a sextic elevenfold, and let $C$ be the
smooth transverse sextic. Suppose that $P,Q$ have a common irreducible
factor $R$ of degree $3\le s\le5$. We show there are
$e\in V(P,Q)\setminus C$ and a line through $e$ on which the smaller
vanishing order of $P,Q$ is finite and belongs to $[3,5]$.
This includes a zero member of the pair provided the other is nonzero.

At a general smooth $e\in V(R)\setminus C$ at least one nonzero
cofactor is a unit, since $2s>5$. Thus $V(P,Q)$ locally equals
$V(R)$. If any direction of second-order projective contact with
$V(R)$ there (intersection multiplicity at least three) is realized
by a line not contained in $V(R)$, that line has the required
property. Otherwise at a general point all such lines are contained
in $V(R)$. In the notation of J.~M.~Landsberg and C.~Robles,
\emph{Fubini's theorem in codimension two}, Theorem~1.1,
\cite{fubini}, this implies
$\mathcal C_{2,e}=\mathcal C_{3,e}$.
Their statement of the classical hypersurface Fubini theorem gives a
quadric if the Gauss image has dimension at least two. Since $s\ge3$,
the dual variety is a curve (a constant Gauss map gives a hyperplane).

Let $B$ be the normalization of that dual curve, whose projective span
is $\mathbb P^k$ in the dual projective space. The tangent lines of
the immersed smooth locus extend to a morphism
$\psi:B\longrightarrow\operatorname{Gr}(2,k+1)$ by properness.
Projective biduality in characteristic zero identifies $V(R)$
with the cone over the hypersurface dual to the curve in its span.
The degree of this hypersurface equals $\deg\psi^*\mathcal O(1)$:
take a general pencil of hyperplanes in $\mathbb P^k$, with
codimension-two axis $A$. A line $\psi(b)$ meets $A$ precisely
when some member of the pencil contains $\psi(b)$. The pencil
misses tangent hyperplanes contributed by any exceptional finite
set in $B$, and at every intersection with the dual hypersurface
it is transverse at a smooth point, with a unique point of tangency
by conormal biduality. The pencil can also avoid members with
exceptional tangency and the locus of tangent lines contained in the
axis, by generality. Thus pulling back the Pl\"ucker hyperplane of
lines meeting $A$ counts $s$ tangencies with multiplicity one.

At an ordinary point on $B$ the coordinate linear series admits a
basis $f_0,\ldots,f_k$ with consecutive orders $0,\ldots,k$ in a
local frame, by the characteristic-zero Wronskian criterion.
The Pl\"ucker coordinates $f_i f'_j-f_j f'_i$, $0\le i<j\le k$,
have orders $i+j-1$. These include all integers from zero through
$2k-2$, so at least $2k-1$ pulled-back sections are linearly
independent. On a smooth projective curve a line bundle of degree
$s$ has at most $s+1$ independent sections, giving $s\ge2k-2$.
In particular $k\le3$.

The cone vertex of $V(R)$ therefore contains a linear space of
dimension $7-k\ge4$. The smooth sextic $C\subset\mathbb P^8$
cannot contain a four-plane: its four normal derivatives there are
positive-degree forms with a common zero, which would be singular on
$C$. Take a vertex point $e$ off $C$. Since the equation of $V(R)$
depends only on the coordinates transverse to its cone vertex, it
has multiplicity $s$ at $e$. A general line through $e$ is contained
in neither the common base nor $V(R)$, so the required finite order
is between $s$ and five. Applying (L5) gives
$\delta_L(X,H)/12\ge41/66$.
\paragraph{A simple linear factor and the second jets of its cofactors.}
Write the penultimate quintics as $P=\ell A$, $Q=\ell B$, with
$H=(\ell=0)\simeq\mathbb P^7$, and let $A_0,B_0$ be the restrictions
of the two quartics to $H$. Either member of a pair may be zero,
with order taken to be infinity. Suppose there exist a point $e$ on
a line $M\subset H$ for which
\[
     \min\{\operatorname{ord}_e(A_0|_M),
             \operatorname{ord}_e(B_0|_M)\}=2.
\]
No condition on the top sextic at $e$ is needed. Retain as lights
the plane spanned by $M$ and an arbitrary direction transverse
to $H$, in addition to the two basepoint-free heavy polars. Put
$M=(z=0)$ in this plane, take $z=\ell$, and give $M$ coordinates
$[x:y]$ with $e=[0:1]$. After the ordinary light-weight-five
extraction, the first angular specialization gives weight one to
$z$ and both heavies and zero to $x,y$. Its initial relation is
\[
       f Z_1+g Z_2+z(uA_0|_M+vB_0|_M).
\]
In a second angular specialization give $x$ weight one, $y,z$
weight zero and the heavies weight two. The light term becomes
$a(u,v)z x^2y^2$, where $a$ is a nonzero linear form along the
original line. These are global diagonal actions on the successive
proper weighted angular spaces over the generic base, using the
independent coordinates $x,y,z,Z_1$ after eliminating $Z_2$.
Both initial polynomials are linear in the heavy pair with coprime
coefficients, so are primitive and prime. Localization at the base
and equality of Hilbert series gives all complex initial pieces.
The final monomial model of (L1) has a global angular estimate,
which permits specialization with any limiting center. Its exponents
are $1,2,2$, and (L1) gives $\delta_L(X,H_X)/12\ge9/14$ for the
hyperplane polarization $H_X$.

\paragraph{A common linear factor.}
Let $P=\ell A$, $Q=\ell B$ with quartics $A,B$, and put
$H=V(\ell)=\mathbf P^7$, $A_0=A|_H$, $B_0=B|_H$. We repeatedly
use two observations. A common point of ambient multiplicity
at least three for $(P,Q)$ suffices: a general line there has finite
order $3\le b\le5$ and, if the point is on $C$, is transverse to
$C$. Also, if the restricted pair of quartics has exact order two
at a point on some line in $H$, first taking the $\ell$-linear
symbol transverse to that line and then its first binary term gives
the symbol $\ell x^2y^2$. In that case (L1) gives $9/14$.

If $A_0=B_0=0$, a common zero on $H$ of the cubics
$A/\ell,B/\ell$ gives the first observation. If a linear form $z$
on $H$ divides $A_0,B_0$, choose coordinates beginning with $\ell,z$.
On $\mathbf P^6=V(\ell,z)$ the four cubics $A_0/z$, $B_0/z$,
$\partial_\ell A$, $\partial_\ell B$ have a common zero, where
both ambient quartic gradients vanish. A common nonlinear irreducible
factor of minimum multiplicity one gives exact order two along a
general simply tangent line: its projective second fundamental form
cannot vanish generically in characteristic zero, since its Gauss
map would then be constant. If the minimum multiplicity is greater
than one, it is exactly two by degree, and a general transverse line
works. These arguments allow precisely one zero restriction by
assigning to it infinite order.

Suppose now $A_0,B_0$ are coprime and take any reduced irreducible
component $Z$ of their base scheme. Its dimension is five.
At a general smooth point of the affine cone on $Z$ let
$T=\ker dA_0\cap\ker dB_0$. If either Hessian has nonzero value on
$(w,w)$ for some $w\in T$, the corresponding line has exact order
two. Otherwise both Hessians vanish on $T\times T$ by polarization.
On the dense open of constant rank $j$ the linear span $N$ of
$dA_0,dB_0$ is then constant along $Z$: for every tangent direction
$v$ to the cone on $Z$ we have $v\in T$, so each of
$d^2A_0(v,-),d^2B_0(v,-)$ annihilates $T$ and belongs to $N$.
This is the vanishing of the differential of the conormal-space
map to the Grassmannian. In characteristic zero zero differential
makes the map constant. Thus everywhere on $Z$ both restricted
gradients lie in the same fixed vector space $N$ of dimension $j$.

If $j\le1$, the ambient gradients of $A,B$ along $Z$ have at most
$2(j+1)\le4$ potentially nonzero cubic coordinates, in $N$ and
the additional $\ell$ direction. These cubics have a common zero
on the projective fivefold $Z$. There $A,B$ have multiplicity at
least two, giving the first observation.

If $j=2$, Euler's identity puts $Z$ in the fixed annihilator
$\mathbf P(N^\perp)=\mathbf P^5$, whence equality. This linear
space is not contained in $C$ (even a four-plane in the smooth
sextic of $\mathbf P^8$ is impossible by its normal derivatives).
At a general point $e\in Z\setminus C$ the ambient space
$T'=T_eV(A,B)$ is a $\mathbf P^6$ transverse to $H$, by independence
of the restrictions of the gradients to $H$. If $T'$ is not contained
in $V(A,B)$, a general line in it through $e$ transverse to $H$
gives finite contact at least three and (L5) applies. Otherwise
$W=T'$ is a linear component of the dimension-six $V(A,B)$, and
$K=W\cap H=Z$ by local smoothness at $e$. The two constant
normal directions to $W$ give a $2\times2$ matrix $M$ of cubic
normal derivatives on $W$, with nonzero determinant at $e$.
Let $D$ be a reduced irreducible component of
$V(\det M)\cap K$; its dimension is four. For $p\in D$ the
ambient tangent space $T_pV(A,B)$ has dimension at least seven,
contains $W$, and cannot lie in $V(A,B)$. At $p\notin C$ a general
line in this space transverse to $H$ works. If $D\subset C$, then
at a general $p\in D$ the same tangent space is not contained in
$T_pC$. Indeed if its generic codimension is zero this is immediate;
in codimension one the opposite inclusion would force both
$\nabla A,\nabla B$ to be multiples of $\nabla C$ along all of $D$.
The multipliers are sections of $\mathcal O_D(3-5)$ and vanish,
contradicting rank one. Thus a general line in $T_pV(A,B)$ through
this point can meet both $H$ and $C$ transversely and have finite
contact with $(P,Q)$ of order at least three. The transverse-on-$C$
estimate applies.
\section{Conclusion}
At contact five the nonzero pair of degree-five coefficients is
coprime or has a common irreducible factor of degree one through
five. The preceding local arguments cover each case, including a
zero member of the pair. Together with the remaining polar ranks,
contact orders and conics, the argument gives the common lower
bound $\delta_D(X,H)\ge2448/329$ and supplies the curve mark of Theorem~\ref{sextic:criterion}.
In the tightest displayed line comparison,
\[
   \frac{2448}{329}-\frac{7433}{1000}=\frac{2543}{329000}>0 .
\]
Theorem~\ref{sextic:criterion} yields $\delta(X,H)>7$ and hence
$\delta(X,-K_X)>1$ as claimed in Theorem~\ref{sextic:main}.

% ===== Hypersurfaces of dimension \(2d-1\) in higher degree =====
\chapter{\texorpdfstring{Hypersurfaces of dimension \(2d-1\) in higher degree}{Hypersurfaces of dimension 2d-1 in higher degree}}\label{ch:higher}
\section{Results and dependencies}
Put \(H=\mathcal O_X(1)\), \(r=n+2-d\).
This chapter develops the global criterion in dimension \(2d-1\) for \(d\ge14\) and verifies several sufficient local conditions. Chapter~\ref{ch:fullrow} supplies the estimate for arbitrary lines.

\begin{theorem}
Let \(d\ge14\), \(n=2d-1\). Write a line as \(z_1=\cdots=z_n=0\)
with homogeneous coordinates \(u,v\). The restrictions to the line
of the \(n\) derivatives \(\partial F/\partial z_i\) of an equation
of \(X\) are binary forms of degree \(d-1\). If their span has
dimension at least three on every line, then \(X\) is uniformly
K-stable. More generally the conclusion holds if, for every line of
rank two, the restriction of \(F\) to the common linear kernel of its
two independent polars has order at most \(\lfloor d/2\rfloor+1\) or exactly \(d\)
along the line. At contact $d-1$, either of the further tests in the subsection on penultimate contact also suffices.
\end{theorem}
Here the common kernel is obtained by setting the two active normal
coordinates equal to zero after a constant change of normal coordinates.
For this theorem all conics, including those on planes contained in
\(X\), are handled by the estimates below. The criterion with rank
at least three is open: the incidence of lines on \(X\) with rank
at most two is determinantal on the relative Grassmannian, and its
projection is proper. The more explicit incidence calculations at
the end give examples in subfamilies with a prescribed plane or line.
The short-curve bound in Chapter~\ref{ch:fullrow} removes the condition on lines.

\section{A sufficient threshold on the next dimension}

\begin{theorem}
Let $X$ be a smooth hypersurface of degree $d\ge14$ and dimension
$n=2d-1$ over $\mathbb C$, and write $H=\mathcal O_X(1)$,
$r=d+1$, $N=n+1=2d$. Suppose that for every line and integral conic $C$
the local threshold, taking all valuations whose centers contain $C$, satisfies
\[
             \delta_C(X,H)\ \ge\ M:=\frac{N(d+2)}{2d+1}.
\]
Then $\delta(X,H)>r$. In particular $X$ is uniformly K-stable.
\end{theorem}

We use the canonical-extension positive-center estimate, principal-bin
point estimate and two-low-degree Nadel test from Chapter~\ref{ch:boundary}. These inputs include uniform K-stability for a smooth hypersurface
of index $s$ and dimension $2s-2$. Therefore if $\beta=\delta(X,H)\le r$,
a pointed smooth hyperplane section of index $r-1$ shows
$\beta>N(r-1)/n>r-1$. Write $u=n/N$, $b=r/N$.
The target is $M=r+1/2-3/[2(2d+1)]$.

\paragraph{Finite levels.}
Retain the upper $1-f$ fraction, $f=1/d$, of the basis adapted to an
original optimizer and a free hyperplane missing its generic center.
Subtract that hyperplane and multiply the remainder by
$\lambda=a(\beta-N/3-3\epsilon)/h$, for the recorded share
$0<a\le1$ and $\epsilon=10^{-8}$. Put
\[
 h\ge 1+n(1-(1-f)^{1/n}),\qquad
 g=(1-b)/(d-1),\qquad
 \alpha\le\frac{1-Nf+nf^{1+1/n}}{N(1-f)}.
\]
The existing uniform tests give $\beta'>N/3$, and for
$\lambda_r=\lambda|_{\beta=r}$ allow the choices
\begin{align*}
 U&=r+\lambda_r(h-1-g)+2\epsilon,&
 z&=b+\lambda_r\{u(1+g)-1+\alpha\}-u\epsilon,\\
 U_p&=r+\lambda_r((d-1)^{-1}-g)+2\epsilon,&
 P&=u\lambda_r/(d-1).
\end{align*}
Here $U$ bounds the comparison against all positive centers, and $z$
is a lower bound on the initial ample slack. At a point optimizer $w$
normalized by $w(\mathfrak m_p)=1$, the additional inequality is
\[
          S_H(w)\ge (A_X(w)+P)/U_p.
\]
Indeed the lowest $1/d$ fraction of the weight list at any such $w$
has mean at least $u-o(1)$ in degree-normalized units. The number of its
weights below $mt$ is at most the length of the smooth-point order-$mt$
jet ring, and thus at most $(mt)^n/n!+O(m^{n-1})$ for bounded $t$.
Integrating the resulting quantile lower bound proves the assertion.
Every collection of independent basis vectors of the retained size has
sum of weights at most the sum of the largest jumps. Subtracting the
omitted jet contribution gives the point inequality at the actual
$\beta$. Its right-hand side decreases with $\beta$ since
$P' U_p<n\le A_X(w)$ and $U_p'\ge1$. Finite-level errors are absorbed
by the displayed allowance on $U_p$ since $A_X(w)/S_H(w)$ is uniformly
bounded here and $A_X(w)\ge n$.

Take $\Delta=\min(M-U,(z-1/3)/u)-\epsilon$.
For $15\le r\le105$ the accompanying exact calculation rounds
\[
 Q>U_p\big((1-z)/\Delta+u\big)-n-P
\]
upwards on a grid of denominator 1000, setting $Q=0$ if this upper
rounding is negative. If $A_X(w)=n+q$ with $0\le q\le Q$, the point
inequality gives $A_X(w)/S_H(w)\le U_p(n+Q)/(n+Q+P)$. The previous
principal-bin estimate at the latter parameter excludes this entire
interval. Thus $z+\Delta((n+q+P)/U_p-u)>1$.
The positive-center cap and the ample slack $z-u\Delta>1/3$ exclude
positive centers as in the two-low-degree Nadel argument. Nadel at the
remaining point contradicts this last display.

\paragraph{Larger indices.}
For $r\ge106$ use the former, coarser cutoff
$x=(r-N/3-3r^{-4})/d$, with
\[
 U=r+bx+2r^{-4},\qquad
 z=b+x\{u(1-b)-\log(d)/N\}-ur^{-4}.
\]
For $106\le r<1024$ the previous affine point envelope
$S_H(w)\le J+Gq/n$ suffices by rational arithmetic. For $r\ge1024$ set
$y=\sqrt r$, $E=0.2\log r$, and $\Delta=1/3-(1.5+E)/r$.
Since $x<1/3+1/d$, one has $U-r<1/6+0.67/d$. Also $x\ge1/3$ and
$u(1-b)=1/2-3/(4d)+1/(4d^2)$, which give $z>2/3-E/r$.
Thus $U+\Delta<M$ and $z-u\Delta>1/3$. As in the uniform range,
\[
 \begin{aligned} J&\le\max\left(\frac{2\sqrt{2d}+n-2}{N},
                 \frac{n-2+d+\sqrt{d/2}}{N}\right)<1.515,\\
 G&\le\max\big(\sqrt{2d}+4,\sqrt{d/2}+\log(12r^2)\big)<1.56y .\end{aligned}
\]
Consequently $n-UJ>0.4815r$ and $UG/n-1<0.782y$. The affine
envelope gives $q>0.615y$ unless the point is already excluded.
It follows that $(n+q)/U>2+0.5/y$, since $0.115 y>3.34+0.085/y$.
The final inequality
$z+\Delta((n+q)/U-u)>1$ follows from
\[
       y/6>2E+1.5+0.5(E+1.5)/y .
\]
At $y=32$ it holds directly; dividing by $y$ makes the right-hand
side decreasing thereafter. This is again the two-low-degree Nadel
contradiction.
\section{All conics on the next row in degree at least fourteen}
Let $X=X_d^{2d-1}$ be a smooth hypersurface over $\mathbf C$,
$d\ge14$, $N=2d$, and put $b=(2d+1)/(d+2)$. Every smooth plane
conic $C\subset X$ satisfies $\delta_C(X,H)\ge N/b$. Throughout,
mean orders are normalized by $N$ and tested over the full geometric
generic fiber. Thus valuations with centers strictly containing $C$
are included. The proper-plane conic estimate $(d+2)/(2d)$ in
$N$-units exceeds $1/b$. On a contained plane the ordinary
radial term is smaller than $b$ and angular splitting gives
$1+a_{\max}/(2d)$. This suffices as soon as the polar rank is
at least six. A conic pencil at every light order satisfies the
mixed-pencil bound $[1+(k-1)(2-1/d)]/k\le b$ with
$k=\lfloor(d+1)/2\rfloor$. We treat ranks three, four and five.

Give all active normals weight $h$, the order of $X$ in the conic
ideal after setting these normals to zero, and the light symbols
weight one. The light linear space, of codimension at most five,
cannot be contained in $X$ since $2d-5>(2d-1)/2$. Thus
$2\le h\le d$. For $h\le d/2$ the weighted binary calculation,
splitting all but two active polars, gives the bound since
$(1+h(2-1/d))/(1+h)\le b$; the finite binary coefficient is at
most $3/2-1/d$. For $h>d/2$ select any independent basepoint-free
net of the active polars; two general combinations already have no
common zero on $C$. Split off the other heavy variables. By shearing
the light conic symbol with light linear normals as in the weighted
conic presentation, retain one light linear variable $x$ where the
first polynomial restricts to $a x^h$, $a\ne0$. The coarse plane
of the resulting ternary core has boundary $(1-1/h)L$ with
$L=f_1z_1+f_2z_2+f_3z_3$. The primitive ring has remaining
generator $Q=x^h=-L/a$ of section cost $h$, and three $z_i$ of cost
one. All four have coarse degree one. The normalized mean order of
the ternary axis divided by its discrepancy is
$(1+h(3-1/d))/(1+2h)<3/2$. The initial extraction adds the free
weights to discrepancy and no more than these to mean order.

Let $s,2d-2-s$ be the low and high degrees in the net quotient
bundle; $1\le s\le d-1$. The actual $L$ coefficient is $1/h$ by
upstairs divisibility, and the global normalized coarse mean without
boundary is at most $c_0 A$ for $c_0=1+1/h-1/d$, by the primitive
diagonal limit $f_1z_1+aQ$. After simultaneous adaptation to a
test valuation and to $L$, remove its $L$ component, leaving
$\Gamma$. On $L$ the complete two-summand quotient after peeling
gives the intersection bound $I=2-(2+s)/(2d)$ at every point.
It is the larger of the two split means, by diagonal degeneration.
Thus for $s\ge4$ ordinary adjunction along $L$ already suffices.
When the split degrees differ, away from the zero $P$ of the
high-summand coordinate take that coordinate as unit. The order is
bounded by the degree in the low coordinate. For the maximum at
any given point choose a basis simultaneously adapted to its order
and to the degree filtration in that low coordinate; leading low
exponents then have the fixed multiplicities of the split blocks.
Their normalized mean is
$1+s/(2d)<b$, so only $P$ needs attention. It is a rational moving
point with degree-$s$ basepoint-free homogeneous coordinates on the
conic: the quotient of the rank-two bundle by its high subbundle
is $\mathcal O(s)$, whose three generating sections give $P$. In
particular it is nonconstant. These facts are independent of the
choice of the three polars in higher polar rank.

\paragraph{Uniform multiplicity gain at a moving point.}
For any nonconstant rational graph $P\colon\mathbf P^1\to\mathbf
P^2$ the full-product portion with base/plane degrees $(2p,k)$
has normalized mean point multiplicity at most $5/3$ after the
usual section-cost normalization $p+k$ and inclusion of common
free factors. Choose homogeneous plane coordinates so the projection
$[P_1:P_2]$ is nonconstant and $P_2\ne0$ generically. On the chart
$z_2=1$, write $t=z_3/z_2$, $y=z_1/z_2$. In the block of polynomials
with multiplicity at least $j$ at $P$, filter by the highest exponent
$l$ of $t$. Differentiating $l$ times with respect to $t$ lowers
multiplicity by at most $l$ and leaves only $l!$ times the leading
coefficient in $y$. It consequently vanishes to order $(j-l)_+$ at
the generic binary coordinate $P_1/P_2$. The primitive binary graph
equation has bidegree $(e,1)$ for some $e\ge1$ after removing common
factors. Therefore the leading coefficient of base degree $2p$
and binary degree $k-l$ belongs to a space of dimension at most
\[
 (2p-e(j-l)_++1)_+ (k-l-(j-l)_++1)_+.
\]
Summing in $l$ gives the order-$j$ dimension precisely for the
model graph $[P_1:P_2:0]$, where powers of $t$ split and each
remaining coefficient vanishes on the binary graph. It is in turn
bounded by the same formula at $e=1$, the linear graph model.
That model has normalized mean $1+2/3=5/3$: the transverse
coordinate $z_3$ has cost one, and the binary graph equation has
cost $1+1/2$; equivalently use the complete jet blocks with density
$(2p-l')+(k-l-l')$ after $z_3^lT^{l'}$. This proves the bound
including for nonembedded or covered graphs.

Use now the regular/pole comparison for the actual primitive system
$Q=-L/a$. The regular product sector has dimension fraction $h/d$.
The fraction $1-h/d$ of positive pole residues on $a=0$ has normalized
mean multiplicity at most $2+1/h$. Indeed residue level $i$ imposes
a fixed $L^i$ and otherwise has a complete plane system of degree
$k-i$ over the finite scheme of base roots. The marked graph lies
on $L$, which is a smooth relative linear coordinate even over
those roots. The per-block bound is $i+2(k-i)/3$, with normalized
sector means $1/h,3$ for $i,k-i$. The regular/pole lemma's argument
with flat infinitesimal graph neighborhoods works for every morphism
$P$, including repeated roots of $a$. Subtracting the actual $L$
gives
\[
          \operatorname{mult}_P\Gamma\le
          (h/d)(5/3)+(1-h/d)(2+1/h)-1/h
                = M:=2-h/(3d)-1/d .
\]
For $h>d/2$ and $d\ge14$ one has $M<b$, since
$1/6+1/d>3/(d+2)$.

\paragraph{Split degrees two and three.}
For a valuation at $P$ normalized by $v(\mathfrak m_P)=1$, put
$u=v(L)\ge1$ and let $A$ be its ordinary discrepancy, so $A\ge u+1$.
The smooth-surface intersection inequality bounds the normalized
order of $\Gamma$ by $(A-1)\min(M,I/u)$. Use $M'=\min(M,I)$.
To obtain $S\le b(A-(1-1/h)u)$ for the preceding adapted basis
it suffices, since $M'\le b$, to check the smallest $A=u+1$.
With the actual $L$-coefficient $1/h$ this reduces to
\[
           \min(M'u,I)\le b+\frac{b-1}{h}u .
\]
It holds for every $u\ge1$ if $(b-1)/h\ge I-b$: if the right side
is already at least $I$ at $u=1$, it remains so afterwards.
For $s\ge2$ one has
$I-b\le (d-4)/(d(d+2))<(b-1)/d$, as required. This proves
all such cases, including any encountered with equal split degrees.

\paragraph{Split degree one.}
For $d\ge19$ one may apply the established product comparison
for the moving linear graph $[s_0:t:0]$. Put $r=(d-1)/(h-1)$,
$R=rh/d$, $u_0=1/(2r)$. Its normalized ratios are no greater than
\[
 R,\quad R(1+1/(2(r+d-1))),\quad
 R(1+1/(1+u_0))/(1+1/h),\quad R/(1-u_0^2).
\]
Here $r\le9/4$, $u_0\ge2/9$, $h^2\ge d$ implies
$R/(1+1/h)\le1$, and $R\le9/8$ since $h\ge9$. Thus the third
entry is at most $20/11<b$, the fourth at most $3/2$, and the
others at most $(9/8)(1+1/(2d))<3/2$. The same four direct
fractions already work for $14\le d\le18$ with $h>d/2$, as the
short certificate checks (or one can use the preceding $M<b$
endpoint where sufficient). This proves the bound on the coarse
plane. All radial, free, and branch discrepancies have been included;
local weighted adjunction over the entire generic fiber and
properness downstairs give the claimed bound for all conics.
\section{Line polar systems on the next row}
Let $d\ge14$, $n=2d-1$, $N=2d$ and $b=(2d+1)/(d+2)$. Every line
on a smooth $X_d^n$ with at least three independent normal polars
satisfies $\delta_B(X,H)\ge N/b$. The same holds for polar pencils
whose light-order contact is at most $d/2$ or is equal to $d$.

For polar rank at least four the ordinary split quotient has at
least three positive summands of total degree $d-1$. Hence
$a_{\max}\le d-3$ and its angular order after multiplying by $N$
is at most $1+a_{\max}/d\le2-3/d<b$. The radial ratio is smaller.
The single-binary estimate covers all $h\le d/2$ irrespective of
rank. Assume rank three and $h>d/2$. The common linear space of
the three polar equations is too large to lie in $X$, so $h\le d$.
Keep the three full heavy variables and one light direction with
nonzero pure symbol $a x^h$. The normalized ternary-axis numerator
is $1+h(3-1/d)$ against discrepancy $1+2h$. All free variables
cost the same normalized mean as their added weight. The coarse
plane $\mathbf P^2$ has boundary $(1-1/h)L$ with $L=f\cdot z$;
the extra generator $Q=x^h=-L/a$ has section cost $h$ and plane
degree one, and $\deg a=d-h$ on the line. The $z_i$ have cost
and plane degree one. The bounded $x$ remainder does not affect
the asymptotic coefficients. If the two positive split degrees
are $s,d-1-s$, the larger intersection bound after peeling $L$
is $I=2-(1+s)/d$ and its own ratio is one. The smaller bound
$1+s/d$ suffices away from the high-coordinate zero, by the
binary coordinate degree. Off $L$ the usual coarse degeneration
gives $1+1/h-1/d$. If $s\ge2$ these estimates already suffice
by adjunction.

For $s=1$ write $P$ for the high-coordinate zero. The regular/pole
argument of the conic proof now uses base polarization one. The
regular product fraction $h/d$ consists of blocks of bidegree
$(p,k)$ on the line and the plane, with total section cost $p+k$.
The point $P$ is the linear graph $[s_0:t:0]$, so its two normal
coordinates $z_3$ and $tz_1-s_0z_2$ have costs one and two.
Its normalized mean point multiplicity in these full products
is therefore $1+1/2=3/2$. Explicitly the separated $z_3^j T^l$
jets have full base coefficient space of degree $p-l+k-j-l$ by
multiplication on the diagonal binary graph. The same subdivision
and simplex calculation as for the conic gives the two means.
Among the pole residues on the finite base scheme $a=0$, the
fixed level $i$ costs $h$ and the remaining complete ternary
system has order per block at most $2(k-i)/3$. With sector means
$1/h,3$ this gives at most $2+1/h$, unchanged from the conic.
If $a$ is constant there is no pole sector. At repeated roots use
the relative graph jets over their Artin rings exactly as above.
After removing the actual $L$ coefficient $1/h$, the residual
multiplicity is at most
\[
          M=(h/d)(3/2)+(1-h/d)(2+1/h)-1/h
                 =2-h/(2d)-1/d<7/4-1/d<b .
\]
For a plane valuation normalized at $P$, $v(L)=u\ge1$ and
$A\ge u+1$. The residual bound $(A-1)\min(M,I/u)$ and $M<b$
reduce the needed inequality to
$\min(Mu,I)\le b+(b-1)u/h$. This holds for every $u\ge1$
since $I-b=(d-4)/(d(d+2))<(b-1)/d$. Thus all rank-three lines
satisfy the threshold, by weighted adjunction on the entire
geometric generic fiber and properness.

For a pencil of maximal contact $d$ the transverse light hypersurface
$Y$ is smooth by finite Gauss and has dimension $2d-5$ in
$\mathbf P^{2d-4}$. Its Fano index is $d-3$, and the previously
established uniform range gives $\delta(Y,H)>d-3$. In the
fractional maximal-contact line estimate (L3) take $k=2d-3$ and
$c=k(d-3)/(k-1)$. Its three normalized threshold entries simplify
to $3/4,(2d-3)/(3d-4),(2d-3)/(3d-4)$. Each exceeds
$(d+2)/(2d+1)$ for $d\ge14$, the last comparison having numerator
$d^2-6d+5=(d-1)(d-5)>0$. This proves the maximal-contact claim.

Combined with the all-conic and conditional next-row theorems,
at this stage, the remaining cases in an arbitrary smooth hypersurface of
this dimension are polar pencils on lines with strict intermediate
light order $d/2<h<d$. In particular, requiring rank at least
three on every line gives a Zariski-open sufficient condition
with no hypotheses on conics or planes. The incidence calculation
for prescribed planes remains valid, a fortiori, at forbidden
rank at most two.
\subsection{One further contact order for every polar pencil}
Let $d\ge14$, $n=2d-1$, and let $L$ be a line whose normal polars
span a pencil. The short-curve target
\[
 \frac{\delta_L(X,H)}{n+1}\ \ge\ \frac{d+2}{2d+1}
\]
holds for every pencil of contact $h=\lfloor d/2\rfloor+1$.
Together with the existing low-contact estimate this gives it for
all $h\le\lfloor d/2\rfloor+1$.

Write $W=(z_1=z_2=0)$ for the common polar kernel. Its light
coordinates span $\mathbf P^{2d-4}$. The degree-$h$ light form
$G_h$ has $q=d-h+1$ scalar coefficients (some may be zero) in the
original binary line variables. Their common zero locus $B$ has
dimension at least $2d-4-q=d+h-5$. At a point $e\in B$ that is
not a common cone vertex, a general line through $e$ has finite
order $a$ in the coefficient tuple with $1\le a<h$. Indeed, a
homogeneous degree-$h$ form has multiplicity $h$ at $e$ exactly
when it is independent of the $e$-coordinate. Taking the first
binary term on that line gives $x^a y^{h-a}$; both exponents are
at most $h-1\le d/2$. In (L1) the wall comparison
\[
 \frac{1+(2-1/d)t}{1+t}\le\frac{2d+1}{d+2}
       \quad\Longleftrightarrow\quad t\le d/2
\]
and the finite coefficient $1-1/d+1/h$ give the assertion.

Suppose henceforth that all points of $B$ are common cone vertices.
Then $B$ is the common linear vertex: it is the projectivization
of the common kernel of the linear maps which take a vector to
the directional derivative of each coefficient form. There exist
general $e\in B$ and a general light vector $a$ such that the
restriction to the two light coordinates $x,y$ in the directions
$a,e$ has a nonzero $x^h$ initial term and some term $x^iy^j$ of
higher total degree with $i<h$. Indeed if no such term existed,
every higher scalar coefficient would have multiplicity at least
$h$ along $B$. This follows by setting $y=1$ and varying $a$ in
its Taylor expansions at general $e$. The same already holds for
the degree-$h$ coefficients. Hence the entire join of $B$ and $L$
would be singular on $X\cap W$: every term and each of its first
derivatives would vanish there. The dimension of that singular
locus is at most one, since the finite Gauss map of smooth $X$
sends it to the projective line of hyperplanes containing $W$.
This contradicts the lower bound on $\dim B$. The two required
nonvanishings on $(a,e)$ are open and hold simultaneously.

Retain the two selected lights and the two heavies by the actual
linear flag through $L$. Take the first lower Newton edge from
$(h,0)$ with far endpoint $(i,j)$, $i<h$. Write
$(h-i)/j=l/m$ in lowest terms and $j=mk$. Then $0<l<m$, since
all other terms have total degree greater than $h$.
Give $x,y,z_1,z_2$ weights $m,l,hm,hm$. The radial ratio is
\[
 \frac{1+(1+l/m)/h}{2-1/d+(1+l/m)/h}.
\]
It is at least the target precisely if $(1+l/m)/h\ge2/d$.
Since $j\le d-i$, we have
\[
 \frac lm\ge\frac{h-i}{d-i}\ge\frac1{d-h+1}
                 \ge\frac{2h}d-1
\]
for $h=\lfloor d/2\rfloor+1$, including equality in the last
comparison for even $d$.

For completeness the angular estimate uses the same Newton-edge
calculation as for the sextic contact-three pencil. In each fixed
weight, the exponent of $y$ has a fixed remainder modulo $m$.
Setting $V=y^m$ gives coarse weights $(1,l,h)$ on $x,V,z$;
the root divisor is $V=0$. Away from it, remove the bounded
remainder and specialize to the far-end monomial $x^iV^k$ by
giving $x$ weight one, $V$ zero, and both heavies weight $i$.
The equation retaining both basepoint-free polars is primitive.
The section costs of $x,V,z$ are $1,m,1$. On the mixed toric wall
$w=is+kt$ for orders $w,s,t$ in $z,x,V$, the normalized mean is
$s+t/m+(2-1/d)w$. Thus its two extremal denominators are
\[
 \frac{1+(2-1/d)i}{1+i},\qquad
 \frac{1/m+(2-1/d)k}{1+k}.
\]
Each is at most $(2d+1)/(d+2)$, since $i\le h-1\le d/2$ and
$k=j/m\le d/2$. At $i=0$ the first coordinate is free.
The individual coordinate coefficients are at most one. The finite
binary coefficient is at most $\max(1,1-1/d+1/(i+mk))$ by the
mixed-monomial calculation of (L1), the monomial having section
cost $i+mk=i+j>h$. Diagonal specialization can use the global
coarse threshold on this complement. At the root divisor, take the
actual $y$-order upstairs. The primitive quotient has pure light
term $x^h$, giving $y$ coefficient one and residual finite binary
coefficient at most $\max(1,1-1/d+1/h)$. Its infinity coordinate
has coefficient at most one by exact divisibility. These entries
all lie below $(2d+1)/(d+2)$ for $d\ge14$. Adjunction from the
weighted extraction and the retained linear flag proves the claim.

Consequently the explicit incidence open of the previous next-row
criterion can use forbidden contact $e=\lfloor d/2\rfloor+2$ in
place of $\lfloor d/2\rfloor+1$. In this open all lines and conics
reach the target. The excluded incidence has codimension at least
\[
 \sum_{j=0}^{e-1}(d-j+1)\binom{2d-4+j}{j}-(8d-8)
\]
in the parameter space of equations, with the same flag dimension
and independent jet count as in the preceding computation.

There is also a conditional gain for $h=\lfloor d/2\rfloor+2$.
Suppose the degree-$h$ light coefficients cut out a subscheme with
a nonlinear generically reduced irreducible component $Z$.
At a general smooth point of $Z$ the differentials of these
coefficients span its conormal space, as the local scheme is then
smooth. The second fundamental form of $Z$ is nonzero at a general
point; otherwise its Gauss map has zero differential and is constant
in characteristic zero, making $Z$ a linear space. A tangent vector
therefore has nonzero quadratic coefficient in at least one of
the equations, by the conormal description of the second fundamental
form and polarization. Its tangent line has finite order exactly two
for the coefficient tuple. The resulting mixed symbol has exponents
$2,h-2$, each at most $d/2$, proving the same line mark by (L1).
\subsection{Two further tests for pencils of penultimate contact}
Let $d\ge14$, $n=2d-1$, $N=n+1$, and let $L\subset X_d^n$ be a line
with a pencil of normal polars and first light order $h=d-1$. Write
the order-$h$ term on the common polar kernel as $uP+vQ$, and its
degree-$d$ term as $C_d$. Either of the following conditions suffices
for $\delta_L(X,H)/N\ge(d+2)/(2d+1)$:
\begin{enumerate}
\item The degree-$h$ forms $P,Q$ have no common factor.
\item Some $e\in V(P,Q)\setminus V(C_d)$ admits a projective light
line on which the smaller of the two vanishing orders at $e$ is
finite and satisfies $b\ge\lceil(d-2)/2\rceil$.
\end{enumerate}

For the first condition restrict to a general line of light directions.
It misses $V(P,Q)$, so its binary $P,Q$ define a degree-$h$ morphism
to $\mathbf P^1$. Retain these two light variables $x,y$ and the two
heavy coordinates by localized flag adjunction in Chapter~\ref{ch:sixfold}. The retained pencil of polars is basepoint-free, so the cuts
are smooth along $L$. The resulting fourfold has initial relation
\[
 fz_1+gz_2+G=0,\qquad G=uP(x,y)+vQ(x,y),
\]
with weights $1,1,h,h$ on $x,y,z_1,z_2$. Here $f,g$ are binary
polars of degree $d-1$ in $u,v$. In the angular calculation below a
bar denotes the asymptotic mean multiplied by five and divided by
section degree. The angular weight $J$ has $\bar J=2+h(2-1/d)$ by
the Hilbert derivative, and the discrepancy of the extraction is
$2+h$. The angular surface is $\mathbf P(1,1,h)$ over
$K=\mathbf C(u/v)$. We give the changes in the coprime binary proof
for the sextic elevenfold.

At the heavy vertex, taking initial light order leaves the prime
relation $fz_1+gz_2=0$ with two free light generators. At the common
heavy zero the initial relation is $G=0$ with both heavy generators
free. The latter locus is a separable degree-$h$ point over $K$ by
the morphism just constructed. In each normal blowup the two free
exponents have mean one each, the discrepancy is two, and torus
splitting on its exceptional line bounds every directional
coefficient by one. On the remaining part of $G=0$, let $j$ be the
top $z_1$ degree after eliminating $z_2$, and $i$ the order of $G$.
Then $hi\le J-hj$, $\bar j=2-1/d$, and $\bar i\le2/h$. On removing
$G^i$ the smallest possible $z_1$ exponent is at least $j_0-i$,
where $j_0=\max(0,\lceil(J-m)/(h-1)\rceil)$ in section degree $m$.
At a nonzero finite $z_1$ the restriction therefore has maximal
normalized order at most
\[
 \bar j-\frac{\bar J-5}{h-1}+\frac2h
     =\frac{1+1/d}{h-1}+\frac2h<1.
\]
These averages use the induced top-degree filtration on each
$G$-order quotient. Conjugacy gives the same estimate at each root.

Off $G$ and the vertex, spread a test point over a finite base cover
and choose a general specialized point $[u:v]=[1:0]$ with
$P(x,y)\ne0$, $f(1,0)=1$, $g(1,0)=0$. Shear $z_2$ by a constant
multiple of $z_1$ to make $z_2=0$ at this point. For base order use
the generator $W=u^{h-1}z_1+P$; its section degree is $h$. The
two initial relations are
\[
 u^{h-1}z_1+P=0,\qquad
 uW+v(Q+c_1u^{h-1}z_1+c_2u^{h-1}z_2)=0.
\]
The lowest-$z_2$ specialization removes the $c_2$ term and fixes the
marked point. Both presentations are complete intersections. Their
$u=0$ fibers have $P=vQ=0$ of codimension two among the remaining
variables; after inverting $u$ the solved presentations are polynomial.
Thus the rings inject into their generic presentations, and their
ordinary Hilbert series agree with that of the fourfold. The two
successive basis comparisons and semicontinuity along the marked
multisection apply as in the sextic proof.

Put $Q'=Q-c_1P$. In the second ring fix $z_2$ exponent $j$, base
order $i$, and residual light degree $D=J-hj$. A block is the sum
of $P^c(Q')^t\mathbf C[x,y]_{D-h(c+t)}$ over
\[
      0\le t\le i,\qquad (h-1)c+t\ge i+j+D-m.
\]
The binary ring is free over $\mathbf C[P,Q']$ with homogeneous
basis degrees $0\le l\le2h-2$. For a compatible $l$ put $k=(D-l)/h$.
The condition on the total $P$ exponent $a$ becomes
$(h-1)a+\min(i,k-a)\ge i+j+D-m$. The possible $k$ differ by at most
one, and so do the corresponding lower cutoffs for $a$. If $s$ is
the smallest occurring exponent, division by $P^s$ leaves a system
containing $P\mathbf C[x,y]_{D-h(s+1)}$ of codimension at most $h$
in the complete degree-$(D-hs)$ system. Moreover the ideal $(z_1)$
is saturated with respect to $u$: modulo $z_1$ the equations are
$P=0$ and $uW+vQ'=0$, with no component above $u=0$, by coprimality.
It agrees blockwise with divisibility by $P$, and iteration gives the
same statement for all powers. The normalized mean actual $P$ exponent
is therefore one by the cost-one Hilbert shift for $z_1$. In a
complete binary system of degree $a$ the mean exponent of $P$ is
$a/(2h)+O(1)$ by summing the dimensions divisible by its powers. A
codimension at most $h$ changes this by $O_h(1)$. Thus
\[
  1=\bar s+\frac{\bar D-h\bar s}{2h},\qquad
  \bar D=\bar J-h,
\]
since the free $z_2$ exponent has mean one. At the chosen point $P$
is a unit, so the residual binary point order is at most
$(\bar D-h\bar s)/2=\bar D-h=1+1/d$. Its fiber direction meets the
free divisor $z_2=0$ transversely on this light chart. Consequently
the angular threshold is at least $d/(d+1)$. By extraction and flag
adjunction the line threshold divided by $N$ is at least
\[
 \min\left\{\frac{d}{d+1},\frac{d(d+1)}{2d^2-d+1}\right\}
       >\frac{d+2}{2d+1}.
\]
The cross-product difference in the radial comparison is $2d-2$.

For the second condition use the mixed and pure initial model (L5) in Chapter~\ref{ch:models}, with $j=d-1-b$. Its radial ratio is
$(1+(2+1/b)/d)/(2-1/d+(2+1/b)/d)>(d+2)/(2d+1)$.
The root divisor has normalized coefficient at most one, and the
finite binary coefficient on its complement is at most
$1-1/d+1/(d-1)$. The two remaining extremal denominators are
\[
 \frac{1/b+2-1/d}{2},\qquad \frac{1+(2-1/d)j}{1+j}.
\]
The first is at most $5/4$, and the second is at most
$(2d+1)/(d+2)$ precisely when $j\le d/2$. This is the stated
condition on $b$. The root divisor is tested before the last torus
specialization in (L5), so this computation requires no assumption
on the limiting position of a point outside it.
\section{Examples with prescribed planes and lines}
The all-conic lemma removes the conic rank conditions
from this open and proves the same stability conclusion under
the line rank condition alone. The incidence argument below records
the stronger genericity assertion for hypersurfaces through a plane.
Let $d\ge14$, $n=2d-1$, $X=X_d^n\subset\mathbf P^{2d}$ smooth.
Write $H=\mathcal O_X(1)$ and $M=2d(d+2)/(2d+1)$. The global criterion in this chapter proves K-stability once every
line and smooth conic has local threshold at least $M$, testing
all valuations whose center contains the curve.

The largest simple projective incidence open supplied by the present
jet test is as follows. Set $e=\lfloor d/2\rfloor+2$. Consider flags
$B\subset W\subset\mathbf P^{2d}$ consisting of a line and a
codimension-two linear space. Exclude forms satisfying
$F|_W\in H^0(\mathcal I_{B/W}^{e}(d))$ on some such flag.
This is a closed projection of a projective incidence. Its complement
among smooth forms consists exactly of the hypersurfaces with every
pencil line of light order less than $e$. Indeed the order condition
already implies the transverse differential on $W$ vanishes along
$B$. Smoothness forces at least two independent normal polars on the
line, so their rank is exactly two and $W$ is their common kernel.
The flag dimension is $8d-8$. For fixed flag the linear conditions
have codimension
\[
     \sum_{j=0}^{e-1}(d-j+1)\binom{2d-4+j}{j}.
\]
This counts the independent monomials on $W$ by order $j$ along
$B$, with $2d-3$ light normal variables. Thus the codimension of
the excluded locus in the space of hypersurface equations is at
least the displayed number minus $8d-8$. By the all-conic and
line results every smooth member of this explicit open is uniformly
K-stable. Its definition does not exclude planes or ordinary pencil
lines. The closed exceptional set can still contain stable members,
as the fixed-jet calculations below illustrate.

Define $\mathcal U_d$ by the following two rank conditions:
\begin{enumerate}
\item On every contained line, the normal polars restricted to
the line span a vector space of dimension at least four.
\item For every contained plane and every conic divisor on that
plane (including reducible and double conics), the normal polars
restricted to the conic span a vector space of dimension at least
six.
\end{enumerate}
Then $\mathcal U_d$ is Zariski open and every one of its hypersurfaces
is uniformly K-stable. This open is dense even in the linear family
of hypersurfaces containing any previously specified plane.

Indeed the split degrees of the line-normal quotient sum to $d-1$
and the number of positive summands equals the polar span dimension
minus one. If the span has dimension at least four, the largest
degree is at most $d-3$, and the ordinary extraction's normalized
angular order (here normalized by $2d$) is at most $2-3/d$.
On a proper conic plane section the conic threshold is at least
$2d(d+2)/(2d)$ by the proper-plane lemma. For a contained plane,
the positive split degrees of its active normal polars along a
smooth conic add to $2d-2$ and their number is one less than the
polar dimension. Six polars therefore give largest degree at most
$2d-6$, and normalized angular order at most $2-3/d$. The ordinary
radial orders are smaller, since the relevant $n-1$ discrepancy
divides $n-1/d$ for lines and $n-1/2-1/d$ for conics. The inequality
$2-3/d < (2d+1)/(d+2)$ proves the claim. The rank constraints are
determinantal on the relative projective Grassmannian of lines
and on the relative projective bundle of quadratic divisors over
the Grassmannian of contained planes. In the latter the restriction
to a conic of any plane form of degree $d-1$ is a vector bundle
of rank $2d-1$, by the hypersurface restriction sequence including
nonreduced conics. Thus each bad locus projects to a closed set.

\paragraph{Incidence with a prescribed plane.}
Fix $\Pi_0\simeq\mathbf P^2$ and work in the linear system
$\mathbf P H^0(\mathcal I_{\Pi_0}(d))$. A general member is smooth:
the normal polars on $\Pi_0$ are $n-1$ unrestricted degree-$(d-1)$
forms; for a general choice they have no common zero, and away from
the plane use Bertini. For any other plane $\Pi$, let $t$ be the
dimension of the vector-space intersection of the two planes'
affine cones, $0\le t\le2$. Vanishing on $\Pi$ imposes
$\binom{d+2}{2}-h^0(\mathbf P^{t-1},\mathcal O(d))$ independent
conditions, using coordinate subspaces to obtain every polynomial
vanishing on the intersection. This is at least $d(d+1)/2$. The
whole extra plane family has dimension $3(n-1)=6d-6$, strictly
smaller for $d\ge14$. Hence a general member has no second plane.

For each conic on $\Pi_0$ its $n-1=2d-2$ independent normal
polars are arbitrary in the $(2d-1)$-dimensional conic system,
so rank at most five has codimension $(2d-7)(2d-6)>5$, the
dimension of the space of conics. For a line $L\subset\Pi_0$
the same $n-1$ polars are arbitrary in $H^0(L,\mathcal O(d-1))$;
rank at most three has codimension $(n-4)(d-3)>2$. For a line
disjoint from $\Pi_0$, requiring it to lie on $X$ costs $d+1$,
and all $n$ normal polars there are freely prescribable thereafter.
For a line meeting $\Pi_0$ at one point the containment costs $d$;
at least $n-2$ of its normal polars can then be chosen freely in
$H^0(L,\mathcal O(d-1))$, using normals to the span of the two
linear spaces. These are direct monomial interpolation facts in
simultaneous linear coordinates: multiply arbitrary degree-$(d-1)$
line forms by the specified linear normal coordinates. For disjoint
spaces the line variables themselves vanish on $\Pi_0$. In the
intersecting case the prescribed normal coordinates vanish on
$\Pi_0$. The respective codimensions of the bad flags thus are
at least $d+1+(n-3)(d-3)$ and $d+(n-5)(d-3)$, both greater than
the dimension $2n$ of all lines. This proves the asserted density.
It also shows that imposing the absence of planes is unnecessary.

\paragraph{A weaker pointwise geometric hypothesis.}
Membership in the above open is sufficient, and the proof extends
to each hypersurface satisfying this less restrictive incidence
test. For a low-rank line (rank at most three) take the linear
space $W$ obtained by setting to zero every active polar normal.
For a low-rank conic on a contained plane (rank at most five)
define $W$ analogously. These spaces are intrinsic: modulo the
original linear span, their direction spaces are the respective
kernels of the polarized differential; in the conic case the
plane tangential polars vanish identically. Require that the initial
conic- or line-ideal order of the equation of $X|_W$ is at most
$\lfloor d/2\rfloor$ for every such flagged curve, except that
polar pencils of conics need no assumption. The light-space
restriction is nonzero for these low ranks: a contained
$k$-plane requires $2k\le n$ by smoothness of positive-degree
normal polars. The weighted-pencil calculation provides normalized
angular order at most the maximum of $1$, $1-1/d+1/h$, and
$[1+h(2-1/d)]/(1+h)$ at light order $h$, after splitting off
all but a general basepoint-free pencil. In the conic case shear
the quadratic generator toward a light linear variable as in the
contained-conic calculation. The first two constants and ordinary
radial tests are smaller than $(2d+1)/(d+2)$, and the third
passes precisely at $h\le d/2$. For a conic pencil the stronger
mixed estimate uses $k=\lfloor(d+1)/2\rfloor$ and gives instead
$[1+(k-1)(2-1/d)]/k$, always no larger than this target. This
proves the claimed extension without asserting that the pointwise
test itself is closed under specialization of rank strata.

\paragraph{A prescribed polar pencil.}
The new line estimate also gives examples outside the rank-three
open. Fix a line $B$ and prescribe any basepoint-free polar pencil
as its full linear normal jet. Work in the degree-$d$
affine linear system of this prescribed jet, translated by
$H^0(\mathcal I_B^2(d))$. A general member of that system is
smooth since the prescribed polars have no common zero on $B$
and the square-ideal linear system is generated away from $B$.
Its first light order is two, since the square-ideal variations
include light quadratic monomials multiplied by arbitrary base
forms of degree $d-2$. A general member has no second polar-pencil
line. Indeed on a line disjoint from $B$ the variations give
arbitrary values and first normal jets in simultaneous coordinates,
since the other line's coordinates vanish on $B$ and $d\ge3$.
For a line $B'$ meeting $B$ at $p$, use the $n-1$ normal linear
equations to the plane $\langle B,B'\rangle$, multiplied by
arbitrary degree-$(d-1)$ forms on $B'$ vanishing at $p$ that
extend to forms vanishing along $B$. These products are in
$\mathcal I_B^2(d)$ and vanish identically on $B'$ itself. Thus
under the condition of containing $B'$ its normal polars in these
$n-1$ directions vary independently through a translated copy of
$H^0(B',\mathcal O(d-1)(-p))^{n-1}$. Project the polar coefficients
onto that $(d-1)$-dimensional vector space along any complement;
the projected matrix is unrestricted and rank at most two if the
original polar rank is at most two. The codimension of this bad
condition is at least $(n-3)(d-3)>2n$, the dimension of all lines.
The disjoint case has still larger codimension. Incidence then
excludes every additional bad line on a dense open of the fixed-jet
system. The fixed line uses light order two, every other line has
rank at least three, and all conics satisfy the unconditional bound.
Thus a nonempty open of each such prescribed-pencil system satisfies
the next-row stability criterion.

One can prescribe higher jets as well. Fix any normal jet through
order $j\le d-1$ along $B$ which is smooth on $B$, and use the
translated system $F_0+H^0(\mathcal I_B^{j+1}(d))$. It contains
all degree-$d$ polynomials in normal coordinates to $B$, hence is
generated away from $B$. The general member is smooth. For every
line $B'\ne B$ choose $p'\in B'\setminus B$, a normal linear
coordinate $T$ on $B$ nonzero at $p'$, and independent linear forms
$\ell_i$ vanishing on the full span of $B$ and $B'$. The variations
$\ell_i\ell_jT^{d-2}$ preserve the entire first jet on $B'$, so on
each fiber of that jet they vary the indicated quadratic coefficients
independently at $p'$. If $B'$ is contained and of polar rank two,
its light linear subspace $W'$ is already fixed by that first jet.
The restrictions of the $\ell_i$ to $W'$ span at least $n-4$ normal
directions of $B'$: the common span of two lines has dimension at
most three and intersecting with the codimension-two $W'$ loses
at most two further linear forms. At $p'$ the condition of contact
at least three along $B'$ on $W'$ consequently has codimension at
least $\binom{n-3}{2}$. This exceeds $2n+1$, the dimension of all
flagged pairs $(B',p')$. A constructible incidence bound, stratifying
by the normal jet or its polar kernel if needed, therefore excludes
every further pencil line of order at least three on a dense open of
the fixed-jet system. This works even for $j=d-1$, since the displayed
variations have full normal degree $d$ at $B$.

For example prescribe on $B$ a basepoint-free polar pencil and
light order $h$ for any $d/2<h<d$. One may prescribe its leading
light form to be a nonzero binary coefficient in the base line
times $X^{\lfloor h/2\rfloor}Y^{\lceil h/2\rceil}$. Both exponents
are at most $d/2$. The mixed-monomial bound (L1) gives $N/b$ at
$B$. On a dense open of this family all other pencil lines have
order two, the remaining lines are covered by the net result,
and every conic satisfies $N/b$. This proves stability even on
nonempty open subsets of these arbitrarily prescribed high-contact
polar-pencil strata whenever their fixed light jet has such a
balanced leading monomial.

\chapter{\texorpdfstring{The full boundary \(n=2d-1\)}{The full boundary n = 2d - 1}}\label{ch:fullrow}
\section{The short-curve threshold}
The point argument in Chapter~\ref{ch:higher} reduced the higher degrees
to a local estimate along lines and conics. We prove a common bound that
also gives enough room in the finite point comparisons for degrees seven
through thirteen.

\subsection{A local threshold for polar pencils on lines}
Let \(X_d^n\) be smooth of degree \(d\ge4\), with \(n>4\). For every
line \(L\) on which the first normal polars span a pencil we have
\[
             \frac{\delta_L(X,H_X)}{n+1}\ \ge\ \frac{d+2}{2d+1},
             \qquad H_X=\mathcal O_X(1).
\]
We use the
linear and weighted flag adjunction and the normalization of local
thresholds used in the preceding chapters; all asymptotic inequalities
are first tested strictly below the displayed threshold.

\paragraph{A smooth plane germ.}
Let \(B\) be a smooth (possibly only locally specified) curve germ
in \(\mathbf P^2\) at \(p\), with contact \(b\ge2\) with its tangent
line. For all valuations whose centers contain \(p\), interpreting
the germ near that point, one has the signed bound
\[
       3S_{\mathcal O(1)}(v)+v(B)
            \ \le\ \left(1+\frac{b}{b+1}\right)A(v).
       \tag{*}
\]
Indeed work in affine linear coordinates \(x,y\) at \(p\) with
\(y=0\) tangent to \(B\), and blow up with weights \(1,b\).
The exception has discrepancy \(1+b\); its orders on the scaled
complete-series basis divisor and \(B\) tend respectively to
\(1+b\) and \(b\). For a basis compatible with the extraction,
the initial series in weighted degree \(J\) consists of the monomials
\(x^{J-bj}y^j\), with
\[
 j_{\min}=\max(0,\lceil (J-M)/(b-1)\rceil)\le j
        \le \lfloor J/b\rfloor=j_{\max}
\]
in ordinary section degree \(M\). On the exceptional chart with
\(s=y/x^b\) this is a fixed monomial times a complete polynomial
system of degree \(j_{\max}-j_{\min}\). The normalized mean vanishing
at a finite nonzero \(s\) is at most \(1/b\). To see the constant,
write \(a,j,u\) for the exponents of \(x,y,1\), adding to \(M\);
then the block width is
\(\lfloor a/b\rfloor+\min(j,\lfloor u/(b-1)\rfloor)\).
In the uniform simplex \(3\mathbf E[a/M]=1\) asymptotically and
\(3\mathbf E[\min(j,u/(b-1))/M]\to 1/b\), the latter by integrating
the simplex cut \(j>s,\ u>(b-1)s\). The sum of vanishing indices
at a finite nonzero point of a polynomial system of width \(w\)
is at most \(w(w+1)/2\); its mean is at most \(w/2\).
At the two coordinate points the normalized multiplicities are
at most one, using the full monomial lists and their ordinary
coordinate exponents (on the stack chart at infinity).
The strict transform of \(B\) restricts to a single simple point
at finite nonzero \(s\), since its weighted initial form is
a nonzero multiple of \(y-cx^b\) with \(c\ne0\).
The exceptional angular denominator is at most \(1+1/b\).
It is at most \(1+b/(1+b)\) since \(b\ge2\).
Weighted inversion of adjunction and simultaneous choice of
valuation/extraction-adapted bases prove \((*)\), by first testing
strictly below the limiting threshold.

In particular take a smooth affine hypersurface \(W\) in the standard
affine chart of \(\mathbf P^3\) at a point \(p\). Suppose that
the smallest nonzero order beyond the linear term in its local
graph over its tangent plane is \(b\), with \(2\le b\le d/2\).
A general affine plane through \(p\) cuts a smooth curve of contact
\(b\) with the tangent line. Indeed its tangent direction can be
chosen generally in the tangent plane to \(W\), where that first
nonzero form of degree \(b\) is nonzero. Complete hyperplane flag
adjunction, retaining \(W\) as additional boundary, now gives
\[
      4S_{\mathbf P^3,\mathcal O(1)}(v)+v(W)
               \le \left(1+\frac{d}{d+2}\right)A_{\mathbf P^3}(v)
                    \qquad (p\in c(v)).
       \tag{**}
\]
For the flag the normalized hyperplane coefficient is one and the
residual complete \(\mathbf P^2\)-series has normalized mean degree
three. The test coefficient in \((**)\) is less than one, permitting
raising the hyperplane coefficient to one in inversion of adjunction.
Outside \(W\) the ordinary projective-space bound suffices.

\paragraph{Two light directions.}
Write \(u,v\) for homogeneous coordinates on \(L\), and choose two
heavy linear normal variables realizing the basepoint-free pencil of
normal polars of degree \(d-1\). The common polar kernel \(W\) cannot
lie in \(X\), because a contained linear space has dimension at most
\(n/2\): its positive-degree normal derivatives must have no common
zero. Let \(G\) be the first light term of \(F|_W\), of light degree
\(2\le h\le d\) and base degree \(d-h\). We shall restrict to two
constant light directions \(x,y\), keeping the two heavy variables.
In ordinary projective space this is a linear flag through \(L\);
every cut remains smooth along \(L\) since the retained polars
have no common zero. Each restricted hypersurface is integral. In a homogeneous factorization
its normal order one places one factor in normal order zero. The
binary initial of that factor divides both coprime first polars, so
it is constant. Homogeneity then forces the factor itself to be
constant. Each next linear equation is consequently regular, and
the hypersurface exact sequence gives the full restricted graded
pieces. Inversion of adjunction takes place in the smooth open
along the line. The linear flag estimate reduces the assertion
to the resulting fourfold, in units of five.

If the geometric generic polynomial \(G_\eta\) has no component
of multiplicity greater than \(d/2\), take a general line in the
projective light space. It meets the reduced components transversely
at distinct smooth points and preserves their multiplicities.
The line can be defined over \(\mathbf C\), since complex points
of its Grassmannian are Zariski dense after extension to
\(\overline{\mathbf C(L)}\). Take light weights \(m=l=1\) and
heavy weight \(H=h\); write \(\Phi\) for this homogeneous binary
restriction of \(G\). The geometric generic root multiplicities
of \(\Phi\) are at most \(d/2\). This holds outside finitely many
base points by projective closedness of the repeated-root condition.

Otherwise the unique component of multiplicity \(s>d/2\) must be
a hyperplane, since \(h\le d\). It is defined over \(\mathbf C(L)\)
by uniqueness and Galois invariance. Take its primitive equation
\(\ell(u,v,x)\) over the binary polynomial ring, with bidegree
\((e,1)\). Gauss's lemma gives \(\ell^s\mid G\), whence
\(se\le d-h<s\); the hyperplane is constant. Use its equation
as the coordinate \(x\). Choose a general complementary vector
with a nonzero degree-\(h\) value, giving a pure term \(x^h\).
There is a light vector in \((x=0)\) giving a pure term \(y^t\)
of degree \(t\le d\). Indeed otherwise the join of this kernel with
\(L\) would be a contained linear space of dimension \(n-2>n/2\).
Here and below pure coefficients are binary forms nonzero at the
generic point of the base.

In the restricted Newton polygon at weight \((1,1)\), every
minimal term has \(x\)-exponent at least \(s>d/2\). As the
\(x\)-weight divided by the \(y\)-weight increases, eventually
every minimum has \(x\)-exponent zero. Some lower Newton edge thus
has endpoint \(x\)-exponents \(i_0\le d/2<i_1\), with primitive
positive integer normal \((m,l)\), \(m>l\). Denote its weight by
\(H\) and its face polynomial by \(\Phi\). The other endpoint
coordinates include the minimum \(y\)-exponent \(j_1\le d-i_1<d/2\).
These endpoint coefficients are nonzero outside a finite base set.

In both cases give heavy variables the weight \(H\).
The normalized extraction denominator is
\[
             \frac{m+l+(2-1/d)H}{m+l+H}\le\frac{2d+1}{d+2}.
\]
Indeed the required \(m+l\ge2H/d\) follows from \(h\le d\) for
equal weights and from \(H\le mh,\ H\le lt\) for the crossing edge.
The weighted initial equation is \(fz+gz_2+\Phi=0\) for a choice
\(f,g\) of pencil generators. Its primitivity and its ordinary
Hilbert series identify the full weighted initial ring, as in the
preceding pencil estimates.

\paragraph{Ordinary polynomial series from the pole blocks.}
Choose \(g\) with \(d-1\) simple roots \(p_a\) outside the
finite excluded set just described. Over the angular stack
\(\mathbf P(m,l,H)\) retain \(z\) and eliminate \(z_2\) using
\((fz+\Phi)/g\). At section degree \(M\) and angular weight \(J\),
the term of \(z_2\)-exponent \(j\) and remaining monomial
\(x^A y^B z^C\) has \(mA+lB+HC=J-Hj\) and binary base degree
\(M-j-A-B-C\). Filter by maximum pole order at the roots and refine each
layer by successive vanishing of its leading coefficients at \(p_1,\ldots,p_{d-1}\).
The subspace with maximum pole order at most \(i\) is spanned
by the terms with \(j\le i\). Indeed in a representative with
highest \(j>i\), the corresponding coefficient polynomial vanishes
at every \(p_a\), since \(f(p_a)z+\Phi(p_a,x,y)\) is nonzero in
the polynomial domain. Dividing the coefficient by \(g\) and
multiplying by \(fz+\Phi\) replaces it at level \(j-1\), preserving
ordinary and angular degrees. Thus the
positive block at \(p_a\) and order \(j>0\) is \(W_a^j\) times
the full monomial subspace
\[
 \left\langle x^A y^B z^C:\ mA+lB+HC=J-Hj,\quad
                 A+B+C\le M-j-a+1\right\rangle ,
 \]
Here \(W_a=z+\Phi(p_a,x,y)/f(p_a)\).
The regular group at base order \(e\) has the same subspace with
\(j=0,\ A+B+C\le M-e\) and no \(W_a\). These assertions use only
interpolation of the binary base coefficients at distinct points.

For a fixed test valuation on the geometric generic angular space
choose complex bases in each angular degree simultaneously adapted
to it and to this refined pole flag. There are two filtrations here on a finite-dimensional
complex vector space, so a basis splitting both exists. For each root group
regard its chosen sections as a family of weighted homogeneous polynomials
over an open part of the original base line. At that group's root multiply
each section by a power of the local parameter so that it has regular
coefficients and nonzero leading residue. After shrinking around the root,
the product defines a relative effective Cartier divisor on the constant
proper weighted projective stack. For any rational coefficient strictly
below the threshold of the special-fiber product, openness of the klt
locus and properness give klt on the entire geometric generic fiber.
The scaling factors are units on this fiber. The same argument works at
the ordinary point used for the regular group. Each group estimate thus
bounds the order at the originally fixed test valuation. Adding the
inequalities does not require a common special fiber.
Although each block is homogeneous in the angular weights,
for fixed positive \(j,a\) their union over all \(J\), after
removing \(W_a^j\), is a basis of \emph{all ordinary polynomials}
of degree at most \(M-j-a+1\) on \(\mathbf A^3_{x,y,z}\).
Likewise for fixed regular \(e\) the union is a basis of all
ordinary polynomials of degree at most \(M-e\). Homogenize them
on ordinary \(\mathbf P^3\), where they constitute a basis of
the corresponding complete series. This does not change the
divisors on the affine chart. Homogeneous angular products can
be tested on the punctured affine cone: it maps smoothly and
surjectively to the weighted projective stack.
Consequently ordinary projective-space valuative bounds on this
punctured chart imply the requisite angular bounds, including the
different on coarse quotient charts.

For completeness,
\(\sum_{j\ge1}\binom{M-j-a+4}{3}=M^4/24+O(M^3)\) per positive
group, with the analogous formula for the regular group. Each
of the \(d\) groups has mass \(1/d\) against the total leading
dimension \(dM^4/24\). Conditional on the group the mean complete
degree, divided by \(M\), tends to \(4/5\); the mean positive
fixed exponent \(j/M\) tends to \(1/5\). This follows by summing
the binomials with the indicated extra factors, or by the uniform
four-simplex. Thus in units of five, the regular group contributes
at most \((1/d)\sup_v 4S_{\mathbf P^3}(v)/A(v)\le1/d\).
Each positive group contributes at most \((1/d)\) times the
supremum of \((4S_{\mathbf P^3}(v)+v(W_a))/A(v)\) on the
punctured affine chart. One uses the uniform finite-level complete-series bound of \cite[Corollary~3.6]{ar1706.04548} here and discards a vanishing contribution from
bounded ordinary degrees.

\paragraph{Tangent contacts on the graph.}
Fix \(p_a\), absorb \(1/f(p_a)\) in its constant face polynomial,
and write \(W=(z+\Phi(x,y)=0)\). A graph point in the punctured
cone has \((x_0,y_0)\ne(0,0)\). Weighted Euler differentiation
gives
\[
 \operatorname{Hess}(\Phi)\binom{mx}{ly}
     =\binom{(H-m)\Phi_x}{(H-l)\Phi_y}.
\]
Here \(H>m,l\): for equal weights \(h\ge2\), and on a crossing
edge \(H\ge mi_1>m\). Wherever \(d\Phi\ne0\) its Hessian is
therefore nonzero and the tangent contact of the graph is two.
At a critical point weighted Euler gives \(\Phi=0\). There the
contact equals the multiplicity of the weighted homogeneous binary
polynomial. Under equal weights, its distinct linear factors meet
only at the origin and all have multiplicity at most \(d/2\).
For a crossing face the two coordinate factors have multiplicities
\(i_0\le d/2\) and \(j_1<d/2\). On the torus the remaining factors
have the form \(x^l-\lambda y^m\), \(\lambda\ne0\). They are smooth
and disjoint there, and each multiplicity is at most \(d/m\le d/2\)
by the ordinary degree bound. Thus every graph point under
consideration has tangent contact \(2\le b\le d/2\).

Estimate \((**)\) applies at all these points. The positive-group
denominator is at most \(1+d/(d+2)\); adding the regular group
and the \(d-1\) positive groups gives
\[
          \frac1d+\frac{d-1}{d}
                         \left(1+\frac{d}{d+2}\right)
               = \frac{2d+1}{d+2}.
\]
The limiting products used for each group are angular-homogeneous,
so evaluation on the punctured affine cone proves their threshold
bounds on the proper weighted projective stack. Using the root-group
specializations, then weighted and linear flag adjunction, proves
the stated line threshold.

\subsection{Lines in odd degree}
Let \(d\ge7\) be odd and \(n=2d-1\). Write \(N=2d\), \(r=d+1\),
and use the short-curve target
\[
 M=r+\frac{d-1}{2d-1},\qquad
 D=\frac{N}{M}=1+\frac{1-1/d}{1+1/(d-1)+1/d}.
\]
For polar rank at least four the ordinary angular bound \(2-3/d\)
and radial bound \((n-1/d)/(n-1)<3/2\) are smaller than \(D\).
For rank two or three, let \(W\) be the projective linear space
of all common light directions together with the line. Its codimension
is the polar rank, hence \(\dim W\ge n-2>n/2\). Take the first nonzero
light term of degree \(h\) on \(W\), retaining a general
basepoint-free pencil of active polars. A rank-three net admits
such a pencil: its associated basepoint-free map to the projective
plane has image of dimension one, so one may project from a point
outside the image. Cut away the remaining one heavy direction
in an actual linear flag. Every flag cut is smooth along the line
by basepoint-freeness. The ordinary hypersurface ring and section
restrictions compute its graded systems even if a prescribed cut
is singular elsewhere: the normal initial \(fz_1+gz_2\) is
primitive and irreducible, precluding a homogeneous factorization
since the order-zero factor would divide both polars in binary
degree. The same applies on earlier cuts retaining additional polars.

If \(h<d\), use the two-light Newton proof above. In its repeated
hyperplane case the additional hyperplane in \(W\) has dimension
at least \(n-3>n/2\), so the second pure direction exists also in
polar rank three. Equal weights give \((m+l)/H=2/h\), whereas
on the crossing face \(H\le mh,lt\) for \(h\le d-1,t\le d\).
In either case the radial denominator is at most \(D\).
Since \(d\) is odd, each graph point has finite tangent contact
at most \(k=(d-1)/2\). Its angular denominator after averaging
the pole groups is at most \(1+(1-1/d)k/(k+1)<D\).

If \(h=d\), the equation on \(W\) is a homogeneous polynomial
\(G\) of the light coordinates alone. Its projective singular
locus in the space of light directions has dimension at most zero.
Indeed the join of any singular component with the original
line is contained in the singular locus of \(X\cap W\). At those
points the Gauss map of smooth \(X\) takes values in the space of
hyperplanes containing \(W\), of dimension at most two.
The Gauss map is finite since its pullback of \(\mathcal O(1)\)
is the ample line bundle \(\mathcal O_X(d-1)\). The join has
dimension two plus the dimension of the component. A general
projective plane of light directions avoids the finite singular
locus and cuts a smooth plane curve by Bertini.
Retain its three light coordinates, together with the same two heavy
directions, and take weights \((1,1,1,d,d)\). The previous pole proof
works with three light coordinates. The punctured affine chart for its
angular stack is now \(\mathbf A^4\setminus\{0\}\), and each pole
group in ordinary degree \(q\) gives all polynomials on \(\mathbf A^4\)
of degree at most \(q-j-a+1\). The simplex integration gives regular
contribution \(1/d\) and positive-group normalized masses \(1/d\),
with weighted complete degree and graph exponent respectively five
and one. At a punctured graph point \(dG\ne0\): simultaneous
vanishing of its derivatives would give a singular projective point
by Euler. The Hessian identity gives tangent contact two. A general
flag of affine planes gives
\(5S_{\mathbf P^4}+v(z+G/f)\le(1+2/3)A\) by the plane-germ
lemma. The total angular denominator is at most
\(1+(1-1/d)2/3\le D\), and the radial denominator is
\(1+(1-1/d)/(1+3/d)<D\).

\subsection{The sharpened short-curve mark in even degree}
Let \(d=2k\ge6\), \(n=2d-1\), \(N=2d\). Put
\[
 A_d=d+1+\frac{d-1}{2d-1},\qquad
 D_d=\frac{N}{A_d}=1+\frac{1-1/d}{1+1/(d-1)+1/d}.
\]
For every line \(L\) on a smooth \(X_d^n\), we have
\(\delta_L(X,H)\ge A_d\). 

For lines of polar rank at least four use the ordinary split and
radial estimates \(2-3/d\) and \((n-1/d)/(n-1)\), both smaller
than \(D_d\). For polar rank two or three let \(W\) be the common
kernel in the ambient projective space; its codimension is the polar
rank. As in the two-light proof keep two general heavy directions
whose polars are basepoint-free, using ordinary linear flags.
Write \(h\) for the first nonzero light order on \(W\).

There is one preliminary restriction if \(h=k\) and the first light
term \(G\) has a component of multiplicity \(k\) over the geometric
generic point of \(L\). Its entire degree is \(k\), so this component
is linear and unique. Its primitive equation \(\ell\) is defined
over the binary base and has base degree at most one, since
\(\ell^k\mid G\) and \(\deg_{\rm base}G=d-k=k\).
In base coordinates \(u,v\) it is therefore either one constant
linear form in the lights or \(u\ell_0+v\ell_1\), with constant
linear forms \(\ell_i\). Restrict the light space to their common
kernel, retaining all of \(L\) and the two heavy directions for
the subsequent flag. The new common light space \(W'\) has ambient
codimension \(c\le5\), so it cannot be contained in \(X\) by the
linear-space bound. Its first light order \(h'\) is strictly
greater than \(k\). In other cases put \(W'=W\), \(h'=h\),
\(c=\operatorname{codim}W\). In every case \(2c<n\).
Each linear flag member in these restrictions is smooth along \(L\),
by basepoint freeness of the retained polars. Its equation is
irreducible: in a factorization one factor would have normal order
zero on \(L\), since the total normal order is one. The restriction
of this positive-degree factor to \(L\) would divide all the retained
polars. The hypersurface restriction sequence gives the complete
graded pieces on each integral flag member, so local inversion of
adjunction along \(L\) applies also to these prescribed sections.

Suppose first that \(h'<d\). If the geometric generic form has
all component multiplicities at most \(k-1\), take a general
two-dimensional light direction space, exactly as in the two-light
proof. The radial comparison uses \(2/h'\ge1/(d-1)+1/d\);
the graph tangent contacts are at most \(k-1\).
In the remaining case \(h'>k\). Indeed the preliminary restriction
has removed precisely the problematic case at \(h'=k\).
The component of multiplicity at least \(k\) is unique and linear
as \(h'<2k\). It is a constant hyperplane since its primitive
base degree \(e\) satisfies \(ke\le d-h'<k\). In two general
light directions put \(x=0\) for this hyperplane, with the other
direction furnishing a pure \(y^t\)-term for \(t\le d\).
Such a direction exists: the corresponding hyperplane of \(W'\)
has dimension \(n-c>n/2\) and hence is not contained in \(X\).
The complementary direction gives a pure \(x^{h'}\) term. In
the lower Newton polygon take the crossing edge with
\[
         i_{\min}<k\le i_{\max}
\]
for the \(x\)-exponents at its ends. Its primitive weights \(m,l\)
satisfy \(m>l\), since minimal total-degree monomials have
\(x\)-exponents at least \(k\). Its weighted degree \(H\) obeys
\(H\le mh',lt\), so the radial estimate holds as
\[
              (m+l)/H\ge 1/h'+1/t\ge1/(d-1)+1/d .
\]
On this face both coordinate factors have multiplicity less than
\(k\). For the \(y\)-factor, at the end \(i_{\max}\) one has
\(j_{\min}=H/l-(m/l)i_{\max}\le d-(m/l)k<k\).
If a torus factor \(x^l-\lambda y^m\) had multiplicity \(k\),
its ordinary degree \(km\) would be at most \(2k\), so
\((m,l)=(2,1)\). The entire polynomial would be a nonzero scalar
times \((x-\lambda y^2)^k\), whose monomial \(x^k\) contradicts
the least total light order \(h'>k\). These bounds persist at the
general chosen polar roots by excluding the finitely many zeros
of the relevant coefficient and discriminant data.
The Hessian identity and the two-light graph proof now bound all
tangent contacts by \(k-1\); this includes noncritical contacts of
order two since \(k\ge3\). Thus the normalized angular denominator
is at most \(1+(1-1/d)(k-1)/k<D_d\).

Finally let \(h'=d\). The equation on \(W'\) is a degree-\(d\)
homogeneous form \(G\) solely in the light variables. The joint
span of \(L\) with any projective singular component of \(G\) lies
in the singular locus of the section \(X\cap W'\). The Gauss
map of \(X\) is finite, and takes this join into the projective
\((c-1)\)-space of hyperplanes containing \(W'\). The singular
locus of \(G\) in \(\mathbf P^{n-1-c}\) thus has dimension at
most \(c-3\). As \(2c<n\), a general plane of light directions
avoids it and meets \(G=0\) in a smooth plane curve. Retain these
three lights and the same heavy pencil. The three-light pole
calculation for \(h=d\) in the odd-degree argument gives radial
and angular denominators at most
\[
       1+\frac{1-1/d}{1+3/d},\qquad 1+(1-1/d)\frac23,
\]
both smaller than \(D_d\). This proves the line assertion.

\subsection{A uniform short-curve estimate}
\label{sec:nov27-uniform-short-mark}
Let \(d\ge6\), \(n\ge2d-1\), and write
\[
      D_d=\frac{d(2d-1)}{d^2+d-1}, \qquad H=\mathcal O_X(1).
\]
For every line and every integral conic \(C\) on a smooth complex
degree-\(d\) hypersurface \(X^n\),
\[
        \delta_C(X,H)\ \ge\ \frac{n+1}{D_d}.
\]
In particular the bound on the row \(n=2d-1\) is
\(d+1+(d-1)/(2d-1)\). For \(n>2d-1\) take a general linear
section of dimension \(2d-1\) through the linear span of the curve.
All intermediate sections are smooth off this span by Bertini. At
any point of the span lying on \(X\), either transversality is
automatic or all the remaining transverse directions must lie in
a fixed hyperplane to give a singular section. The codimension of
the latter condition is at least \(2d-2>2\), which suffices since
the span has dimension at most two. Complete hyperplane flags
carry the normalized bound from \(2d-1\) back up to \(n\).

We therefore assume \(n=2d-1\). For lines the two-light argument
and its three-light case at \(h=d\) in the odd-degree proof apply
to every odd \(d\ge7\): their only numerical requirements are
\(2\le(d-1)/2\) and enough light dimensions for the indicated
cuts, both satisfied here. For even \(d\ge6\) use the removal of
the degree-\(d/2\) repeated linear component proved above.

For conics the proper-plane denominator \(2d/(d+2)\) is smaller
than \(D_d\). The mixed pencil denominator is
\(1+(1-1/d)(j-1)/j\), \(j=\lfloor(d+1)/2\rfloor\); it is at most
\(D_d=1+(1-1/d)/(1+1/(d-1)+1/d)\) since
\(j-1\le (d-1)/2\le d(d-1)/(2d-1)\).
At polar rank at least six use \(2-3/d<D_d\), and the ordinary
radial denominator is smaller than \(3/2<D_d\). At ranks three
through five use the ternary reduction, with split degrees
\(s,2d-2-s\) and light order \(2\le h\le d\), as in the earlier
conic proof. The free axes, ternary axis, coarse off-branch estimate
and unexceptional points on the branch line are bounded by \(3/2\).
Equal split degrees also give branch intersection \(I<3/2\). For
\(h<d/2\), the binary test (B) suffices since
\[
       h\le(d-1)/2\le d(d-1)/(2d-1).
\]

For \(6\le d\le14\) use the tests (B), (E), (P) displayed and proved in
Section~\ref{models:sec:local-small-conic-table}, with the binary
construction from Section~\ref{models:sec:local-conic-models} and
the product comparison from Section~\ref{models:sec:local-net-product}.
Take \(b=D_d\),
\[
       I=2-(2+s)/(2d),\qquad M=\min\{I,2-h/(3d)-1/d\}.
\]
The argument of Chapter~\ref{ch:higher} under \emph{Uniform multiplicity gain at a moving point}, applied to any positive split degree, justifies this \(M\). In the following table each letter
corresponds to one consecutive value \(h=2,\ldots,d\); in the
last column the tests for \(s=2\) work for all \(s\ge2\),
since the upper bound on residual order decreases with \(I,M\).
\[
\begin{array}{c|c|cc}
d&D_d&s=1&s=2\\\hline
6 &66/41 &\mathrm{BEPPP}&\mathrm{BEEEE}\\
7 &91/55 &\mathrm{BBPPPP}&\mathrm{BBEEEE}\\
8 &120/71&\mathrm{BBEPPPP}&\mathrm{BBEEEEE}\\
9 &153/89&\mathrm{BBBPPPPP}&\mathrm{BBBEEEEE}\\
10&190/109&\mathrm{BBBPPPPPP}&\mathrm{BBBEEEEEE}\\
11&231/131&\mathrm{BBBBPPPPPP}&\mathrm{BBBBEEEEEE}\\
12&276/155&\mathrm{BBBBPPPPPPP}&\mathrm{BBBBEEEEEEE}\\
13&325/181&\mathrm{BBBBBPPPPPPP}&\mathrm{BBBBBEEEEEEE}\\
14&378/209&\mathrm{BBBBBPPPPPPPP}&\mathrm{BBBBBEEEEEEEE}
\end{array}
\]

It remains to consider \(d\ge15,\ h\ge d/2\). If \(s\ge2\),
the multiplicity bound satisfies
\[
     M\le 11/6-1/d\le D_d .
\]
The latter inequality reduces to
\(d^3-11d^2+17d-6\ge0\). Moreover,
\[
 I-D_d\le\frac{d^2-4d+2}{d(d^2+d-1)}
       <\frac{(d-1)^2}{d(d^2+d-1)}=\frac{D_d-1}{d}.
\]
Thus for \(u=v(L)\ge1\), normalized by \(v(\mathfrak m_P)=1\),
one has \(\min(Mu,I)\le D_d+(D_d-1)u/h\).
As \(M\le D_d\), the smooth-surface intersection inequality at
the least log discrepancy \(A=u+1\) gives the required bound
at every \(A\ge u+1\).

For \(s=1\) use the product test (P) with
\(t=(d-1)/(h-1), R=th/d, w=1/(2t)\). Since
\(h\ge d/2\), one has \(1\le t\le 2+2/(d-2)\le9/4\) and
\(R\le(d-1)/(d-2)\le9/8\). Also \(R/(1+1/h)\le1\) since
\(h^2\ge d\). Consequently the three product bounds are at most
\[
  (9/8)(1+1/(2d)),\qquad 20/11,\qquad 3/2
\]
respectively. All are smaller than \(D_d\) for \(d\ge15\):
for the middle bound the difference has numerator
\(2d^2-31d+20>0\) with denominator \(11(d^2+d-1)\).
This finishes the conic estimate.

\section{The global threshold on the boundary}
\begin{theorem}\label{thm:nov27-first-row}
Let \(d\ge7\). Every smooth complex hypersurface of degree \(d\)
and dimension \(2d-1\) has \(\delta(X,H)>d+1\). In particular it
is uniformly K-stable.
\end{theorem}
\begin{proof}
Put \(n=2d-1,\ N=2d,\ r=d+1,\ u=n/N,\ b=r/N\).
The preceding short-curve estimate gives
\[
 M=r+\frac{d-1}{2d-1}\le\delta_C(X,H)
 \qquad(C\text{ a line or an integral conic}).
\]
For \(d\ge14\), this exceeds \(N(d+2)/(2d+1)\), so the global
criterion of Chapter~\ref{ch:higher} applies. We give the additional
point comparison for \(7\le d\le13\).

Suppose \(\beta=\delta(X,H)\le r\). A general smooth pointed
hyperplane section belongs to the range \(n'=2d-2\) of
Chapter~\ref{ch:boundary}. Complete hyperplane adjunction gives
\(\beta> N(r-1)/n\). In particular \(\beta>N/3\).
Use the upper-fraction construction of Chapters~\ref{ch:quintic}
and~\ref{ch:higher}. Put \(f=1/d,\ \epsilon=10^{-8}\) and write
\([x]_- = 10^{-8}\lfloor 10^8 x\rfloor\). Choose the share \(a\)
from the table below and set
\begin{align*}
 h&=1+n(1-[(1-f)^{1/n}]_-),&
 g&=(1-b)/(d-1),\\
 \alpha&=(1-Nf+nf[f^{1/n}]_-)/(N(1-f)),&
 \lambda_\beta&=a(\beta-N/3-3\epsilon)/h .
\end{align*}
In a sufficiently high degree, the construction gives a boundary
\(D\sim_{\mathbf Q}\lambda_\beta cH\) with the estimates
\[
 c\le1-\alpha+o(1),\quad v(D)\le
 (\lambda_\beta h+\epsilon)S_H(v),\quad
 \beta-\lambda_\beta h-\epsilon\le\beta'
 \le\beta-\lambda_\beta(1+g)+\epsilon ,
\]
where \(\beta'=\delta(X,D;H)\). In the bound for \(c\) the
integral in the hyperplane-tail calculation strictly exceeds this
downward-rounded \(\alpha\), so one may take \(c\le1-\alpha\)
for all sufficiently large levels. The certificate checks at both
endpoints of \([N(r-1)/n,r]\) that \(\beta'>N/3\) follows from the
lower estimate and that the upper estimate is strictly less than
\(r-\lambda_\beta(1-\alpha)\). These are affine comparisons in
\(\beta\). The resulting log Fano pair has a divisorial minimizer.

Write \(\lambda=\lambda_r\) and put
\begin{align*}
 U&=r+\lambda(h-1-g)+2\epsilon,&
 z&=b+\lambda\{u(1+g)-1+\alpha\}-u\epsilon,\\
 U_p&=r+\lambda(1/(d-1)-g)+2\epsilon,&
 P&=u\lambda/(d-1),\\
 \Delta&=\min\{M-U,(z-1/3)/u\}-\epsilon,&
 \eta&=z-u\Delta .
\end{align*}
The affine derivatives for the estimates in Chapter~\ref{ch:quintic}
have the same signs here. In particular \(U\) bounds \(A_X/S_H\)
at the tilted optimizer, and \(z\) bounds
\(r-\lambda_\beta c-u\beta'\) from below. The certificate checks
\(U+\Delta<M\) and \(1/3<\eta<1\). The positive-center slope
estimate of Chapters~\ref{ch:global} and~\ref{ch:boundary} with
\(\beta'>N/3\) would place a curve of degree at most two in any
positive-dimensional minimizing center. Its local threshold would
be at most \(U<M\). Hence the minimizer \(w\) is centered at a
point. Normalize it by \(w(\mathfrak m_p)=1\) and set
\(A_X(w)=n+q\), \(q\ge0\).

The bottom-fraction jet bound of Chapter~\ref{ch:quintic} gives
\(S_H(w)\ge(n+q+P)/U_p\). In passing to \(\beta=r\) the
certificate checks \(P'U_p<n\) and \(U_p'\ge1\), with derivatives
taken before endpoint substitution. The finite values used are
\[
\begin{array}{c|rrrrrrr}
d&7&8&9&10&11&12&13\\ \hline
100a&55&60&65&68&71&74&76\\
1000Q&3111&3244&3384&3464&3550&3640&3699
\end{array}
\]
Specifically \(Q\) is the next multiple of \(1/1000\) above
\(U_p((1-z)/\Delta+u)-n-P\). If \(q\le Q\), the ratio
\(A_X(w)/S_H(w)\) is at most \(U_p(n+Q)/(n+Q+P)\).
For this rational parameter the signed surface flag and
principal-factor interval test of Chapters~\ref{ch:foundations}
and~\ref{ch:global} exclude \(0\le q\le Q\) in each column.
The exact substitutions, interval endpoints, bin counts, moment
tests and outward root enclosures are specified by
\path{nov27-first-row.py} and its accompanying
\texttt{r2-local-adjacent-exact.py} in the archive described in Chapter~\ref{ch:checks}.
Thus \(q>Q\), giving \(\eta+\Delta S_H(w)>1\).

On the other hand, in the second basis construction the first
positive-dimensional non-klt center would contain a line or conic
by quasi-log adjunction and the length-one theorem, since the
ample remainder exceeds \(1/3\). Its lc place would give
\(A_X/S_H\le U+\Delta+o(1)<M\). The remaining multiplier ideal
has finite cosupport, and the pointed hyperplane perturbation
and Nadel vanishing in Chapter~\ref{ch:quintic} give
\(\eta+\Delta S_H(w)\le1\). This contradiction proves the result.
\end{proof}

\chapter{\texorpdfstring{The boundary \(n=2d\)}{The boundary n = 2d}}\label{ch:secondrow}
\section{Lines and conics}
In dimension \(n=2d\), write \(N=2d+1\) and \(b=(2d+2)/(d+3)\).
We first establish the short-curve bound \(N/b\). The final section
combines this bound with the pointed hyperplane result of
Chapter~\ref{ch:fullrow}.

\subsection{A mixed estimate for every conic polar pencil}
\label{sec:nov27-mixed-conic-pencil}
Let \(C\) be a smooth conic on a plane \(\Pi\) contained in a smooth
degree-\(d\) hypersurface, \(d\ge4\), and set \(J=2-1/d\). Suppose that
some linear normal direction to the plane has a nonzero polar on
\(\Pi\) whose restriction to \(C\) vanishes. Then
\[
 \frac{n+1}{\delta_C(X,H)}\le
 \max\left\{\frac{1+J}{2},\frac{2d^2-2d+1}{d(d+1)}\right\}
 <\frac{2d+2}{d+3}.
\]
In particular this holds for every conic at which the plane polars
span a pencil on \(C\). Indeed two positive-degree polar forms on
\(\Pi\) have a common zero. Smoothness forces the full system of
plane polars to have no common zero, so the kernel of restriction
of their span to \(C\) is nonzero in the pencil case.

Let \(g=0\) be the conic equation, and take two heavy normals realizing
a basepoint-free pencil on \(C\). Such a pencil exists at every
rank by general projection of the one-dimensional polar image.
Let \(x\) be a linear normal direction in the stated kernel, with
plane polar \(g^p a_0\), where \(1\le p\le\lfloor(d-1)/2\rfloor\) and
\(g\nmid a_0\). In an actual linear flag through \(\Pi\), retain only
these two heavy directions and \(x\). All flag sections are smooth
along \(C\). They are integral: if the restricted equation factored,
its conic-ideal order one and its two basepoint-free heavy polars
would imply that a positive-degree factor restricted to \(C\) was
nonzero and divided both of them, an impossibility. Their standard
hypersurface restriction sequences give complete residual series,
and local flag adjunction carries the normalized result to \(X\).
The final hypersurface in this flag has dimension four.

Use the conic-ideal symbol \(w\) of \(g\), of section degree two, and
choose an integer \(m>p\). Give \(w,x,z_1,z_2\) weights
\(1,m,m+p,m+p\), respectively. The plane is contained in \(X\), so
there is no pure term in \(w\); terms with at least two linear lights
have weight at least \(2m>m+p\). The exact initial equation is therefore
\[
          f_1z_1+f_2z_2+a xw^p=0, \qquad a=a_0|_C\ne0 .
\]
Here \(f_1,f_2\) are basepoint-free in \(H^0(C,\mathcal O_C(d-1))\).
In the complete-intersection presentation over the conic ring this
equation cuts out an integral hypersurface: the punctured conic cone
gives an affine-space bundle, and its vertex fiber has lower dimension
in a Cohen--Macaulay ring. After generic localization the relation
eliminates one heavy variable. Its ordinary Hilbert series agrees
with the original one, and the added quadratic conic symbol cancels
its degree-two relation. Thus the initial ring computes the actual
weighted order over \(C\), with discrepancy \(2m+p+1\) and normalized
mean \(m+1/2+(m+p)J\).

The angular ring has the same mixed monomial presentation and
ordinary section costs as in the paragraph on the extreme conic
pencil in Section~\ref{models:sec:local-conic-models}. The auxiliary
grading separating the characters of \(x,w\) refines the new grading
as well, since \(z_i\) and \(xw^p\) all have weight \(m+p\). Thus the
same factorization into powers of \(x,w\) and binary factors in
\(z_1,z_2,Q=xw^p\) applies. For completeness, the finite binary
point coefficient is at most \(\max(1,1-1/d+1/(2p))\), and the
coefficients of \(x=0,w=0\) are at most \(1,1/2\) by actual
divisibility. On the punctured three-dimensional covering cone the
transverse axis with \(w\) a unit has radial denominator
\(R_x=(1+J)/2\); the one with \(x\) a unit has denominator
\(R_w=(1/2+pJ)/(1+p)\). These are the respective weighted blowups
of \(z_1=x=0\) and \(z_1=w=0\) with weights \((1,1)\) and \((p,1)\).
Their strict transforms separate the finite binomials and are
transverse on the smooth covering charts. At \(x=w=0,z_1\ne0\),
the two coordinate divisors meet transversely. The punctured cone
maps smoothly onto the weighted projective stack of the first
extraction, so these bounds apply with all coarse differents.
The finite point coefficient is at most \(R_x\), while the full
extraction ratio equals
\[
 \frac{m+1/2+(m+p)J}{2m+p+1}=
 \frac{2mR_x+(p+1)R_w}{2m+p+1}.
\]
Finally \(R_w\) increases with \(p\), with the stated upper value
at \(p=(d-1)/2\). The differences of \((2d+2)/(d+3)\) from
\(R_x\) and this last upper value have positive numerators
\((d-1)(d-3)\) and \(7d-3\) with respective denominators
\(2d(d+3)\) and \(d(d+1)(d+3)\). Weighted adjunction over the
whole geometric generic fiber proves the assertion.

\subsection{All conics when \(n=2d\) and \(d\ge13\)}
\label{sec:nov27-all-conics-2d}
Let \(d\ge13\), \(n=2d\), \(X_d^n\) smooth, \(H=\mathcal O_X(1)\), and
\(b=(2d+2)/(d+3)\). Every integral conic \(C\subset X\) satisfies
\[
              \delta_C(X,H)\ge (n+1)/b.
\]
The same normalized estimate holds for \(n>2d\) by taking a general
linear flag through the span of the conic down to dimension \(2d\),
as in Section~\ref{sec:nov27-uniform-short-mark}.
A proper plane gives denominator \(2d/(d+2)<b\); a contained-plane
polar pencil is handled by Section~\ref{sec:nov27-mixed-conic-pencil}.
In rank at least eight there are at least seven positive splitting
degrees of sum \(2d-2\); the largest is at most \(2d-8\), giving
ordinary angular denominator \(2-4/d<b\). The ordinary radial ratio
is smaller than \(3/2\). For ranks three to seven repeat the earlier
ternary conic reduction. The common light space of dimension at
least \(2d-6>n/2\) ensures a finite first contact \(2\le h\le d\).
The ternary axis, free directions, off-branch denominator,
low-summand direction and equal-split cases all have denominator
at most \(3/2<b\). Let \(s\le d-1\) be the low positive split degree
and \(P\) the high-coordinate zero on the branch line \(L\).
On the coarse plane the actual normalized branch coefficient is
\(1/h\), and the residual \(\Gamma\) after removing it has intersection
on \(L\) at most \(I=2-(2+s)/(2d)\). For \(s\ge6\),
\(I\le2-4/d<b\), which suffices by adjunction. We explain the
additional comparison for the remaining split degrees.

\paragraph{Multiplicity and the first infinitely near point.}
The moving graph \(P\colon C=\mathbf P^1\to\mathbf P^2\) has
degree \(s\), in the sense that \(P^*\mathcal O(1)=\mathcal O(s)\).
Project from a constant point outside its image; a general projection
to a projective line is nonconstant and has binary degree exactly
\(s\), including for a covered image. In the moving-point leading
coefficient argument of the preceding conic chapter the graph
equation in this binary coordinate accordingly has bidegree
\((s,1)\). Filtering by the remaining plane variable bounds the
product-sector dimensions by those for the model on a fixed line
with this binary graph. Its transverse variable costs one and its
binary graph equation has section cost \(1+s/2\). Splitting out
powers of these generators is exact: the remaining coefficient
after dividing by the graph equation is a binary polynomial of
the remaining base and plane degrees, modulo its multiples. Thus
the two normalized mean exponents are \(1\) and \(2/(s+2)\)
by the same independent cost shifts as for \(s=1\). The positive
pole fraction has the earlier multiplicity bound \(2+1/h\). In
particular, for the same basis adapted to a test valuation and
to actual \(L\)-order, write \(m=\mathrm{mult}_P\Gamma\); then
\[
                 m\le M:=\min\{I,2-hs/((s+2)d)-1/d\}.
\]
All order and multiplicity bounds refer to normalized asymptotic
basis divisors, or to finite-level bounds with their vanishing errors.

Blow up \(P\) on the smooth geometric generic plane. Denote by
\(P_1\) the point of the exceptional line marked by the strict
transform of \(L\), and by \(m_1\) the multiplicity at \(P_1\) of the
strict transform of \(\Gamma\). Set
\[
     t=\frac{d-1}{h-1},\quad R=\frac{th}{d},\quad
     K_s=\frac{s}{2d-2},\quad J=\frac{R(4t+3s)}{2d-2}.
\]
We claim
\[
          m_1\le J-K_s(m+1/h)-1/h.                     \tag{T1}
\]
Use the regularized product comparison from
Section~\ref{models:sec:local-net-product}: each block of base
cost \(p\) and plane degree \(k\) embeds in the full product of
bidegree \((\lfloor2tp\rfloor,k)\), by clearing the maximal power
of the coefficient of \(x^h\). This only multiplies by a scalar on
the geometric generic plane. The ratio of total enlarged to actual
dimensions tends to \(R\), by the ordinary leading Hilbert
coefficients; this comparison holds for every split degree.

Consider one section of bidegree \((B,k)\) in that full product,
where \(B=\lfloor2tp\rfloor\), and let \(j\) be its multiplicity at
the generic point of the graph. Denote
\(\mathcal K=\ker(\mathcal O_C^3\to\mathcal O_C(s))\), where the
map evaluates plane hyperplanes on the graph. The first nonzero
jet is a nonzero section of
\[
       \operatorname{Sym}^j\mathcal K\otimes
                       \mathcal O_C(B+s(k-j)).
\]
Indeed the associated graded of the graph ideal is the symmetric
algebra of its relative conormal
\(\mathcal K\otimes\mathcal O_C(-s)\), while the original twist
restricts to \(\mathcal O_C(B+sk)\). Generic divisibility by the
graph ideal extends over the base since its powers define flat
infinitesimal neighborhoods. The primitive linear equation of
\(L\) gives a subbundle \(\mathcal O_C(-(2d-2))\subset\mathcal K\):
the polar net has no base point and its equation vanishes on
the graph. Let \(q\) be the vanishing of the first jet at the
direction \(P_1\). Its homogeneous jet polynomial is divisible by
the \(q\)-th power of this subline, also over the base since the
subline is locally a direct summand. The quotient is a section of
\[
 \operatorname{Sym}^{j-q}\mathcal K\otimes
               \mathcal O_C(B+s(k-j)-(2d-2)q).
\]
The symmetric power embeds into a trivial bundle. Nonvanishing
forces \((2d-2)q\le2tp+sk-sj\). This is an individual bound
and requires no simultaneous adaptation to two more flags.

For the enlarged full product the normalized mean base cost is
two and the normalized mean plane degree is three. These are
the simplex means for the two binary and three plane coordinates,
or the logarithmic derivatives of their graded Hilbert series;
rescaling base degree by \(t\) multiplies the leading coefficient
without changing the cost means. In original blocks the sums of
the nonnegative quantities \(2tp+sk\) are bounded by their sums
in full blocks, giving the upper normalized bound \(R(4t+3s)\).
The actual normalized mean \(j\) is \(m+1/h\). On the first blowup
the residual strict transform has intersection at \(P_1\) bounded by
the sum of the jet orders \(q\) minus the actual branch coefficient
\(1/h\). Its multiplicity there is at most that intersection. This
proves (T1).

\paragraph{A sufficient numerical test.} Set
\[
  S=(1-K_s)M+J-(1+K_s)/h.
\]
In addition to (B), (E), (P) of
Section~\ref{models:sec:local-small-conic-table}, the following
inequalities suffice:
\[
 M\le b,\quad J\le b,\quad S\le b+(b-1)/h,\quad
                        I\le b+2(b-1)/h.                  \tag{T}
\]
Normalize the valuation by \(v(\mathfrak m_P)=1\), put
\(u=v(L)\ge1\), and write \(A\) for ordinary log discrepancy. For
\(u=1\) use the earlier \((A-1)M\) residual bound. For \(u\ge2\),
the earlier bound \((A-1)\min(M,I/u)\), of slope at most \(b\),
suffices at \(A=u+1\), because \(I\le b+(b-1)u/h\). In the interval
\(1<u<2\) the center on the first blowup is \(P_1\), with minimal
parameter value \(u-1\) and ordinary discrepancy \(A-1\) there.
The smooth-surface multiplicity inequality applied on this blowup
gives \(v(\Gamma)\le m+m_1(A-u)\). By (T1) its slope is at most
\(J\le b\), and at \(A=u+1\) its value is at most
\(m+m_1\le S\le b+(b-1)u/h\). Thus in all cases including the actual
\(u/h\) branch contribution the normalized order is at most
\(b[A-(1-1/h)u]\), as required.

\paragraph{Large degrees and finite arithmetic.} Suppose \(d\ge100\).
For \(h\le d/3\) use (B), since the binary radial denominator at
\(h=d/3\) equals \(b\). Otherwise \(t\le3(d-1)/(d-3)<3.1\),
\(R\le(d-1)/(d-3)<1.03\). Also \(h^2\ge d\). For \(s=1\), writing
\(u_0=1/(2t)>1/7\) gives \(R/(1+1/h)\le1\), so the three terms of
(P) are at most \(1.03(1+1/(2d)),15/8,1.03(4/3)\), each smaller
than \(b\). For \(s\ge2\), \(M\le11/6\). For \(s\ge4\) one has
\(I-b\le(d-9)/(d(d+3))<(b-1)/d\), and the old multiplicity and
intersection test at \(A=u+1\) suffices for every \(u\ge1\). For
\(s=2,3\) we have
\[
 J<1.03(4\cdot3.1+9)/198<0.112,\qquad
 M+J<11/6+0.112<202/103\le b,
\]
and \(I-b\le(2d-6)/(d(d+3))<2(b-1)/d\). This proves (T).

For the remaining finite degrees the following exact-arithmetic
pseudocode gives the certificate. Loop over \(13\le d\le99\),
\(1\le s\le6\), \(2\le h\le d\); in that order of priority test (B),
then (P) if \(s=1\), then (E) with the degree-\(s\) value of \(M\),
then (T). The numbers passing each first applicable test are
\[
       \begin{array}{c|rrrr}
       \text{test}&B&P&E&T\\\hline
       \text{number of triples}&9048&3273&14394&1995
       \end{array}
\]
The counts sum to \(28{,}710\). See \path{nov27-lit-next-conics.py}.
This exhausts
the distinguished-point tests and hence all conics by weighted
adjunction over the full geometric generic fiber.

\subsection{The smaller degrees in the conic estimate}
\label{sec:nov27-small-conic-2d}
For every \(d\ge6,\ n\ge 2d\), every integral conic on a smooth
\(X_d^n\) satisfies \(\delta_C(X,H)\ge(n+1)(d+3)/(2d+2)\).
It suffices to treat \(6\le d\le12,\ n=2d\). Use the reductions and
notation of Section~\ref{sec:nov27-all-conics-2d}. For \(d=6\),
polar rank at least six already gives ordinary angular denominator
at most \(2-6/(2d)=3/2\); the remaining common kernel has dimension
at least \(n+1-5>n/2\). For \(d\ge7\) the preceding dimension argument
applies directly.

Write \(D=2d-2,\ t=(d-1)/(h-1),\ R=th/d,\ u=1/(2t)\), and let
\(b=(2d+2)/(d+3)\). For the previously untreated cases the split
degree \(s\) is one or two. On blowing up the distinguished point
\(P\), the normalized actual radial coefficient is at most
\(m_0=M+1/h\), where
\[
 M=\min\{2-(2+s)/(2d),\,2-hs/((s+2)d)-1/d\}.
\]
The discrepancy there is \(1+1/h\). We may use this actual radial
bound together with estimates on the exceptional line made using the
enlarged product. Indeed, in each block multiplication by a base
scalar is an injection into the enlarged block preserving order at
\(P\), hence an injection of the associated graded spaces over
\(\mathbf C\) for the order at the moving point. Their coefficients
lie in the same jet spaces on the geometric generic exceptional
line. Thus every filtered jet space embeds as well. Summing these
dimension bounds gives at most \(R\) times the normalized enlarged
mean for any angular point, regardless of which radial bound is
used. Away from this blowup the existing actual estimates on and
off the branch suffice. The distinguished direction \(P_1\)
on the exceptional line has discrepancy \(1/h\), and every other
direction has discrepancy one. Anywhere on the exceptional line
its angular order is bounded by \(m_0\), since the restriction of
the strict transform of an individual nonzero section of point
multiplicity \(j\) is a homogeneous jet of degree \(j\).

For \(s=1\), the product calculation already gives normalized
angular coefficients at most \(R\) off \(P_1\) and at most \(RJ\)
at \(P_1\), where
\[
 J=\frac{1}{1-u^2}\left(\frac1{1+Du}-\frac{u^2}{D+u}\right).
\]
Consequently the three inequalities \(m_0\le b(1+1/h),\
R\le b,\ RhJ\le b\) suffice.

For \(s=2\), we use \(m_0\le b\) away from \(P_1\). The moving
map \(C\to\mathbf P^2\) of degree two either covers a line
or parametrizes a smooth plane conic. For a map to a line of
degree \(s\), the graph ideal is generated by a constant plane
linear form \(z\) and a binary graph equation \(T\) of bidegree
\((s,1)\). The equation of the branch is \(b_0z+\ell T\) with
\(\deg b_0=D,\ \deg\ell=D-s\), and the coefficients have no
common zero. The enlarged product calculation gives
\[
 J_s=\frac{1}{1-(su)^2}
       \left(\frac1{1+Du}-\frac{(su)^2}{D/s+su}\right)
\]
as an upper bound before multiplication by \(R\); for \(su=1\)
its continuous value is
\((1+(D-s)/(2(D+s)))/(1+Du)\).
Here are the changes from the linear graph calculation. With \(a,l\)
the jet powers of \(z,T\), \(c=k-a-l\) and \(e=B-sl\) the remaining
plane and base degrees, the section cost is
\(a+(1+su)l+c+ue\) and the dimension density is \(e+sc\).
The restrictions of the remaining binary forms to the graph fill
this coefficient space for \(e\ge s-1\), by the hypersurface
sequence of the graph of bidegree \((s,1)\); smaller base degrees
are a strip of lower asymptotic dimension. Divisibility by the
\(q\)-th power of the marked equation requires \(a\ge q,\ e\ge sq\)
in the coefficient bound, with \(Dq\) base degrees used. After
shifting \(a,e\) accordingly, the remaining density is
\((e+sc-(D-s)q)_+\) with exponential weight
\(\exp(-(1+su)q)\). Integration of the two independent exponential
variables \(e,sc\) of rates \(u,1/s\), divided by the original
mean of \(e+sc\), yields the displayed formula. This is the
ordinary homogeneous simplex calculation in Laplace coordinates,
including its common normalized section-degree factor.

If the graph is a smooth conic, its degree-two parametrization is
the full binary quadratic system. Hence
\(\ker(\mathcal O_C^3\to\mathcal O_C(2))=\mathcal O_C(-1)^2\),
by the two linear syzygies of the quadratic Veronese. For
bidegree \((B,k)\), the order-\(j\) leading jet is therefore in
\[
 H^0(\mathbf P^1\times\mathbf P^1,\,
                       \mathcal O(B+2k-3j,j)).
\]
The first factor parametrizes the base conic, the second the
constant projective line of the two linear syzygies. In these
coordinates the branch direction is a primitive equation of
bidegree \((D-1,1)\), since the primitives of the syzygies cost
one base degree. Therefore the dimension whose order at the marked
direction is at least \(q\ge1\) is bounded by
\[
 (j-q+1)_+\,(B+2k-3j-(D-1)q+1)_+,\qquad j\le k .
\]
This uses only containment in the full jet space. Base scalars are
units for the generic angular valuation.

Put \(v=1/2,\ g(z)=1/((1+z)^2(1+Dz))\). The normalized enlarged
angular bound furnished by these inequalities is
\[
 V=\frac{u}{v-u}\left(\frac vu g(u)-\frac uv g(v)\right)
 \quad (u\ne v),\qquad
 V=g(u)\left(2+u\left(\frac2{1+u}+\frac D{1+Du}\right)\right)
 \quad (u=v).
\]
Indeed set \(a=j-q,\ c=k-j,\ B=2tp\). The full dimension of
the enlarged product in the exponential integral is
\(\int e^{-p-k}2tp\,k^2/2\,dp\,dk=2t\).
The sum of the filtered angular dimensions is bounded by
\[
 \int_{\mathbf R_+^4} e^{-p-c-a-q}\,
                    a(2tp+2c-a-Dq)_+\,dp\,dc\,da\,dq.
\]
The other free indices and bounded factors cancel. The original
homogeneous simplex expectation in the established normalization
equals this ratio of Laplace integrals. For independent unit
exponentials \(p,c\), the random variables \(2tp,2c\) have rates
\(u,v\). Their sum has expected excess above \(w\ge0\) equal to
\(((v/u)e^{-uw}-(u/v)e^{-vw})/(v-u)\) for \(u\ne v\).
Integrate at \(w=a+Dq\) with weight \(ae^{-a-q}\); division by
\(2t=1/u\) gives \(V\). The equal-rate case follows by continuity.
Thus it suffices to check \(Rh\max(J_2,V)\le b\) when \(s=2\)
and \(m_0\le b\). One may also handle an unmarked direction when
\(m_0>b\) by the following two bounds. For a graph lying on a
line, the enriched jets in the normal factors \(z,T\) have
independent coefficient spaces. Diagonal specialization at a point
of the exceptional line bounds the mean by the larger of the two
coordinate means \(1,1/(1+su)\), before multiplication by \(R\),
as in the linear graph calculation. For a Veronese graph, an order
\(q\) at any given direction in a homogeneous binary polynomial of
degree \(j\) costs at least \(q\) independent coefficient positions,
even for a point over an algebraic extension of the function field:
after a constant linear change avoiding infinity, projection onto
the top \(j-q+1\) coefficients is injective on the filtered space,
since a nonzero polynomial of degree less than \(q\) cannot vanish
to order \(q\). Thus the dimension at fixed \(B,k,j,q\) is at most
\((j-q+1)_+(B+2k-3j+1)_+\). The same integral as above bounds the
mean by \(R V_0\), with \(V_0\) obtained from \(V\) by putting
\(D=1\) in its formula. In particular
\(m_0\le b(1+1/h),\ R\max(1,V_0)\le b,\ 
Rh\max(J_2,V)\le b\) is sufficient in any degree.

Over \(6\le d\le12,\ 1\le s\le d-1,\ 2\le h\le d\) there are
476 triples. The previous tests pass 459. The remaining eleven
with \(s=1\) and six with \(s=2\) pass the inequalities just
proved, by direct rational substitution. The complete check is
\path{nov27-lit-small-conics.py}. Weighted adjunction gives
the assertion for every conic.

\subsection{Lines of polar rank at least three in dimension \(2d\)}
\label{sec:nov27-next-lines}
For \(d=4\) and every \(d\ge6\), \(n=2d\), every line \(B\subset X_d^n\)
on a smooth complex hypersurface whose first normal polars have
rank at least three satisfies
\[
      \delta_B(X,H)\ge(n+1)/b,\qquad b=\frac{2d+2}{d+3}.
\]
At rank at least five the largest of the positive split degrees
of sum \(d-1\) is at most \(d-4\), so the ordinary angular ratio
is at most \(2-4/d<b\), and the radial ratio \((n-1/d)/(n-1)<b\).
At rank three or four, the common polar kernel has dimension
at least \(n-3>n/2\), so the first light order satisfies
\(2\le h\le d\). Keep a basepoint-free triple of active polars,
splitting off any further heavy and complementary light directions
as in the conic-net argument. Write \(s,d-1-s\),
\(1\le s\le(d-1)/2\), for the split degrees of this triple.
The ternary axis ratio is at most
\((1+h(3-1/d))/(1+2h)\le b\): it increases with \(h\), and the
inequality at \(h=d\) reduces to \((d-1)(d-2)>0\). The branch
line ratio is one, off-branch ratio is at most \(1+1/h-1/d\le b\),
and the low-summand direction has denominator
\(1+s/d\le1+(d-1)/(2d)\le b\). At the distinguished point the
residual intersection and multiplicity satisfy
\[
      I=2-\frac{1+s}{d},\qquad
      M=\min\left(I,2-\frac{hs}{(s+1)d}-\frac1d\right).
\]
For the multiplicity use the regular/pole comparison of the line-net
calculation with the graph projection taken outside its image, as in
Section~\ref{sec:nov27-all-conics-2d}. The binary graph equation
in the full-product regular sector now has section cost \(1+s\),
since the base polarization has degree one. The positive-pole sector
and the actual normalized branch coefficient \(1/h\) are unchanged.
For equal split degrees ordinary adjunction with \(I=1+s/d\)
suffices, as it also does whenever \(s\ge3\) since then \(I<b\).

For \(s=1,2\) use the binary test (B) and surface endpoint test (E)
of Section~\ref{models:sec:local-small-conic-table} with these
line-specific \(I,M\). For \(s=1\) the exact line-graph product
test in paragraph~\ref{models:par:local-line-net} is
\[
\max\left\{R\left(1+\frac1{t+d-1}\right),
        R\frac{1+1/(1+u)}{1+1/h},\,RhJ\right\}\le b,
\quad t=\frac{d-1}{h-1},\quad R=\frac{th}{d},\quad u=1/t,
\]
where
\[
 J=\frac{1}{1-u^2}\left(\frac1{1+(d-1)u}
                 -\frac{u^2}{d-1+u}\right)\ (u<1),\qquad
 J=\frac1d\left(1+\frac{d-2}{2d}\right)\ (u=1).
\]
For \(d=4\) and \(6\le d\le99\), successive use of (B), the
line product test, and (E), in this order, gives respectively
\(3040,3319,3326\) successful tests among all \(9685\) triples
with \(s=1,2\), \(s\le(d-1)/2\), \(2\le h\le d\). They are
the direct rational substitutions in \path{nov27-lit-next-lines.py}.

For \(d\ge100\) the binary test suffices when \(h\le d/3\).
Otherwise for \(s\ge2\) one has \(M\le 2-2/9<b\) and
\(I-b\le(d-9)/(d(d+3))<(b-1)/d\), so the smooth-surface
intersection bound suffices at every \(u=v(L)\ge1\) as before.
For \(s=1\) use the first infinitely near point from the preceding
conic proof. With base polarization degree one its jet has maximum
coefficient denominator \(d-1\) in place of \(2d-2\), and the
regularized base degree is \(tp\), not \(2tp\). Thus its coarse
unsubtracted normalized bound is \(J_1=R(2t+3)/(d-1)\).
Here \(t<3.1, R<1.03\), so \(J_1<0.096\), whereas
\(M\le11/6\), giving \(M+J_1<202/103\le b\). Finally
\(I-b\le(2d-6)/(d(d+3))<2(b-1)/d\). These inequalities imply
the sufficient test (T) on the line model with room to spare,
finishing the claim by localized weighted adjunction.

\subsection{A three-light estimate for a power of an irreducible quadric}
Suppose that along a polar-pencil line of a degree-$d$ smooth
hypersurface the first light form has degree $h=2s\le d$ and is a
scalar times an $s$th power of an irreducible quadric over the geometric
generic point. For $d\ge9$ and $s\ge2$, we claim
\[
             \frac{\delta_L(X,H)}{n+1}\ge\frac{d+3}{2d+2},
\]
provided that three general constant light directions can be retained.
A general projective plane of light directions cuts a smooth conic
in the quadric since its rank is at least three. The same holds at
the general pole roots, by openness and uniqueness of its reduced form.
Use three light variables of weight one and two active heavy variables
of weight $h$. The normalized radial denominator is
$1+(1-1/d)/(1+3/h)\le 1+(d-1)/(d+3)$. The pole groups now give
ordinary complete polynomials on the punctured affine space of three
light variables $x$ and one heavy variable $z$. The regular group
has normalized denominator one. For a positive group it suffices to
estimate $5S_{\mathbf P^4}(v)+v(W)$ with $W=(z+cq(x)^s=0)$,
where $c\ne0$ and $q$ is a smooth ternary quadratic form. Off $W$
the ordinary projective-space inequality suffices. At graph points with
$q\ne0$ the Euler identity gives contact two, so that the previous
smooth-germ inequality gives $5/3\le1+d/(d+3)$.

Consider a graph point with $q=z=0$. Choose a linear form $L$ on
the light variables with $L(p)=1$, and cut by the ordinary affine
hyperplane $L(x)=1$. The normalizing complete-degree factor is four,
and the hyperplane coefficient is one. On the resulting ordinary
$\mathbf P^3$ the curve $C=(Z=Q=0)$ is a smooth conic: $Z(z)=z$
and $Q$ is the homogenization in that plane of $q|_{L=1}$, which
recovers the original ternary form. Extend $Q$ to be independent of
the $Z$-coordinate. Denote the graph section in the affine chart by
$W$ as well. We prove at this point that
\[
        4S_{\mathbf P^3}(v)+\lambda v(W)\le2A(v),
                \qquad \lambda=1+\frac{3}{2s}.               \tag{1}
\]
Take the weighted blowup of $C$ with orders $s$ on $Z$ and one
on $Q$. It has log discrepancy $s+1$, and $W$ has order $s$.
In ordinary degree $m$, the associated graded is the direct sum
of $Z^iQ^jH^0(C,\mathcal O_C(t))$, $t=m-i-2j\ge0$,
at weight $si+j$. Indeed expand in $Z$, then take powers of $Q$
on the plane. Each restriction to the conic is surjective by the
plane hypersurface exact sequence. Independence in each weight follows
in its normal bundle, whose two independent normal parameters have
weights $s$ and one. Since $\mathcal O_C(t)=\mathcal O_{\mathbf P^1}(2t)$,
index a monomial by $0\le e\le2t$, the order of a conic section at
$p$. The barycentric variables $i,2j,e/2,t-e/2$ give an asymptotically
uniform simplex with sum $m$. In units of four the normalized means
of $i,j,e$ are therefore $1,1/2,2$. The radial numerator in (1)
is $s+1/2+\lambda s=2(s+1)$ as required.

On the exceptional weighted projective bundle, take the fiber at
$p$ as an ordinary flag. Its normalized coefficient is two by the
same simplex count. The strict transform of $W$ does not contain this
fiber; its restriction is a single simple noncoordinate point at
$w=Z/Q^s=-c'$ with $c'\ne0$, in local trivializations on $C$.
The fiber-graded piece has the monomial interval of $Z^iQ^j$ with
\[
              si+j=J,\qquad i+2j\le m-e/2.
\]
Indeed an order-$e$ coefficient in $H^0(C,\mathcal O_C(t))$ is
arbitrary when $2t\ge e$, independently in each term. Writing
$a=t-e/2$, the width of this interval as $i$ varies is
\[
       \lfloor j/s\rfloor +\min\{i,\lfloor a/(2s-1)\rfloor\}.
\]
At a finite nonzero angular point the mean vanishing is at most half
this width, since a consecutive interval is a monomial times a complete
polynomial system. The scaled simplex mean of the first summand is
$1/(2s)$ and that of the second is also $1/(2s)$: its tail at $t'$
cuts the two simplex coordinates $i>t',\ a>(2s-1)t'$, hence cuts
total length $2st'$. The normalized coefficient of any such point
in the complete-series product is at most $1/(2s)$. Including $W$
raises it to at most $\lambda+1/(2s)=1+2/s\le2$. At the coordinate
points the mean vanishing is at most one and $1/2$, from the same
full monomial lists in the stack charts; $W$ misses them. The
primary and fiber flags, with simultaneous adapted bases and inversion
of adjunction applied strictly below coefficient $1/2$, prove (1).

Interpolating (1) with $4S\le A$ gives
\[
          4S(v)+v(W)\le (1+1/\lambda)A(v)
             =\left(1+\frac{2s}{2s+3}\right)A(v)
             \le\left(1+\frac d{d+3}\right)A(v).
\]
The statement on $L=1$ lifts by the ordinary hyperplane flag. Adding
the regular one-$d$th group and $d-1$ positive groups of equal mass
gives the angular denominator $1+(d-1)/(d+3)$, as asserted.

\paragraph{A repeated quadric with a residual factor.}
The three-light inequality also holds if the geometric generic
first form is $q^s R$, where $s>d/3$, $q$ is an irreducible
quadric of rank at least five, and $q\nmid R$. Write $t=\deg R$,
so that $2s+t=h\le d$ and $t<d/3$. The large component is unique,
and its primitive base degree is zero as in the preceding cases.
The quadric is factorial: it is smooth if its dimension is three
and its rank five, and otherwise its singular locus has codimension
at least four; its class group is generated by a hyperplane.
Equivalently one may use the usual calculation of the class group
of a projective cone over a smooth quadric of dimension at least three.
Consequently each prime component of $(R|_q)$ has multiplicity at
most $t$, since each has positive integral hyperplane class and
the divisor has class $tH$. On a general plane of light directions
the conic is smooth, and its intersection with $R$ has order at
most $t$ at each point. Indeed the plane meets the components
transversely in their smooth loci, avoiding their intersections,
the singular loci and the exceptional locus of the original quadric.
The same generic property holds at the finitely many chosen pole
roots, by taking generic roots after fixing the general plane.

Keep three light variables as in the pure-quadric proof. At a graph
point with $qR\ne0$ the contact is two by Euler. At $q\ne0$,
$R=0$ the contact is either two or is bounded by $t<d/3$, since
multiplicity is bounded by the degree. A graph point with $q=0$
and $R\ne0$ is covered by the pure-quadric proof: a local unit
multiplying $q^s$ changes only the value of the nonzero graph
direction on the exceptional fiber at the prescribed point.
All estimates were local over that point and the generic order
along the conic still equals $s$.

At the remaining graph points first cut by $L=1$ as above, and
choose affine linear plane coordinates $x,y$ centered at the
point, with $y=0$ tangent to the conic. For weights $1,2$ its
initial form is a nonzero scalar times $Q_2=y-a x^2$, $a\ne0$.
Write $B$ for the weighted order of $R$ there. Parametrizing the
conic by $x$, its $y$-coordinate has order two; hence
$B\le\operatorname{ord}_p(R|_q)\le t$. Give the heavy coordinate
$z$ weight $H=2s+B$. For the complete series on ordinary
$\mathbf P^3$ the scaled radial mean is its discrepancy $1+2+H$
(sum the exponents of the three affine coordinates in the uniform
degree simplex); the graph adds order $H$. Thus the signed
radial denominator is $1+H/(H+3)\le1+d/(d+3)$.
The signed angular bound is checked on the punctured affine
$\mathbf A^3$. The complete degree-$M$ ordinary monomials split
by weighted degrees and the union of the full initial blocks
is precisely a basis of all polynomials of degree at most $M$.
Thus on that punctured chart the ordinary $4S_{\mathbf P^3}$
bound can again be used with the new graph
$z+Q_2^sR_B(x,y)$, where $R_B$ is the weighted initial form.
At a noncritical graph point its tangent contact is two, by
weighted Euler with $H>2$. At a critical point $(x,y)\ne(0,0)$
one obtains the multiplicity of the binary weighted polynomial.
Its coordinate factors have multiplicities at most $B\le t<d/3$.
Each other factor has the form $y-bx^2$, $b\ne0$, and any two
such factors meet only at the origin. If its multiplicity is $e$
then $2e\le H\le d$. At every point of this parabola in the
punctured chart, the other factors are units. Its projective
closure is a smooth plane conic, so the preceding conic-flag
estimate $4S+v(W)\le(1+2e/(2e+3))A$ applies at that point
when $e\ge2$; if $e=1$ the gradient is nonzero. For the coordinate
factors use the ordinary smooth-graph contact lemma. All angular
ratios are at most $1+d/(d+3)$. Applying the second weighted
adjunction and the first hyperplane adjunction, then averaging the
original pole groups, proves the claim.

\subsection{Three-light tests and the linear loci that can obstruct them}
Fix a rank-two polar line on a smooth degree-$d$ hypersurface, $d\ge9$,
and write $h$ for the first nonzero light degree and $G$ for its
geometric generic form on the projective light space $P$. Take
\[
 T_d=\frac{d+3}{2d+2},\qquad c=3/d .
\]
If $h\le d/3$ the one-light threshold directly gives
$\delta_L(X,H)/(n+1)\ge T_d$. If $h\le 2d/3$ and every
component of $G$ has multiplicity at most $d/3$, use two general
light directions. Their binary form has all root multiplicities at
most $d/3$ and $2/h\ge3/d$, so the ordinary-degree pole and
plane-germ proof gives the same result.

If $h\le d$ and $(P,cG)$ is log canonical, keep three general
light directions. Their plane section is log canonical with the
same coefficient by repeated general-section adjunction. Equivalently
choose the general constant plane over the original complex field
and the general distinct pole roots: geometric generic choices hold
in a nonempty Zariski open. Take two heavy weights $h$ and three
light weights one. The radial denominator is
$1+(1-1/d)/(1+3/h)\le1/T_d$. For the angular argument use the
weighted complete-series decomposition for $P'=\mathbf P(1,1,1,h)$
and $Z=(z=0), W=(z+G/f=0)$ at a pole root. 
Here \(k=3\). In the ordinary section degree \(m\), with angular weight
\(J\), apply the pole and residue filtration of Chapter~\ref{ch:fullrow}.
The base degree of a term with \(j\) eliminated and \(\ell\) retained
heavy factors is \(B_{j+\ell}=m-J+(h-1)(j+\ell)\). Vanishing at the
previous \(a-1\) points leaves a nonzero evaluation precisely when
\(B_{j+\ell}\ge a-1\).
Consequently the block at \(p_a\) is
\[
          W_a^j z^K\, \mathbf C[x,z]_{J-h(j+K)}^{\rm wt},
 \quad W_a=z+G(p_a,x)/f(p_a),\quad
 K=\max\left(0,\left\lceil
       \frac{a-1-m+J}{h-1}\right\rceil-j\right).
\]
Here the superscript denotes all weighted homogeneous polynomials of
the indicated degree. At \(j=0\), filter the regular sections by base
order \(e\) at one ordinary point; the analogous expression omits
\(W_a^j\) and uses \(e\) in place of \(a-1\).

\paragraph{The complete angular series.}
The \(d\) groups have equal asymptotic mass. Explicitly count normal
monomials avoiding \(u^{d-1}z_2\). The regular group has \(j=0\),
free variables \(u,v,x,z\), and \(u\)-exponent \(e\) for base
vanishing order. Positive group \(a\) has \(u\)-exponent \(a-1\),
free variables \(v,x,z,z_2\), and \(z_2\)-exponent \(j\ge1\).
The \(v\)-exponent in each positive group is \(B_{j+\ell}-(a-1)\);
in the regular group it is \(B_\ell-e\). These monomials biject
with the bases of the residue blocks including their inner complete
systems. Conditional on any group, in units of \(k+3=6\) (the cut has dimension \(k+2=5\)) the
normalized means of each free exponent equal one, by the standard
simplex integral. In every group the common fixed exponent of \(z\)
is \(K=\max(0,\ell-\lfloor v/(h-1)\rfloor)\), reading \(v\) as
its exponent. Its conditional mean is \(1-1/h\) in these units:
in a uniform simplex the tail volume for
\(\ell>t+v/(h-1)\) is \((1-1/h)(1-t)^D\), where \(D\) is the
simplex dimension. This follows by the substitution
\(\ell=t+v/(h-1)+\ell'\), then \(v'=hv/(h-1)\).
The normalized conditional mean complete weighted degree is therefore
\(k+h(1-(1-1/h))=k+1\). The positive groups additionally have
mean \(W_a\)-exponent one.

Write \(S_P=S_{\mathcal O_P(1)}\) and \(Z=(z=0)\). Hence the regular
group has normalized expected-order upper bound
\(\frac1d\sup_v D_0(v)/A_P(v)\), and group \(a\) has upper bound
\(\frac1d\sup_v(D_0(v)+v(W_a))/A_P(v)\), where
\[
                D_0(v)=(k+1)S_P(v)+(1-1/h)v(Z).
\]
For the complete residual spaces use their uniform finite-level
expected-order estimate; bounded weighted degrees make \(o(1)\)
contribution. Equivalently one may pass to levels divisible by
the Cartier index, handle the finitely many twists by multiplication
with a fixed section and use the same asymptotic estimate.

The pair \((P,(1-1/h)Z)\) is log K-polystable: its crepant cyclic
cover is the ordinary projective space under
\([x_1:\cdots:x_k:t]\mapsto[x_1:\cdots:x_k:t^h]\),
and its log anticanonical class is \((k+1)\mathcal O_P(1)\).
Thus \(D_0(v)\le A_P(v)\); the change of variable \(z\mapsto W_a\)
gives this inequality also with \(W_a\) in place of \(Z\).
Write \(w=v(W),z_0=v(Z)\). Combining the two cyclic-cover inequalities
gives \(D_0+w\le A+\min\{w,w/h+(1-1/h)z_0\}\).
Inversion of adjunction from $Z,W$ gives
$\max(v(Z),v(W))+c\min(v(Z),v(W))\le A$. These inequalities imply
$D_0+v(W)\le[1+1/(1+c)]A$ since $h\ge1+c$. Averaging the one
regular and $d-1$ positive groups gives precisely $1/T_d$.

The following sharp structural reduction requires no assumption
on smoothness of the degree-$h$ form. Suppose every component
of $G$ has multiplicity at most $d/3$ and $(P,cG)$ fails to be
log canonical. Choose a rational $c'<c$ sufficiently close to $c$
that the non-klt locus remains nonempty. Its codimension is at least
two. Take a general plane cutting a codimension-two component, if any.
On the plane the non-klt locus is finite, and Nadel vanishing for
its multiplier ideal $\mathcal J$ implies that constants surject onto
$H^0(\mathcal O/\mathcal J)$ since $3-c'h>0$. Hence the locus
is a single reduced point in the multiplier-ideal scheme. General
plane section commutes with the multiplier ideal, so the component
in $P$ is unique of codimension two and degree one, thus linear.
There must be such a component: if all components have codimension
at least three, take a general three-plane through a general point
of a chosen component. The non-klt scheme is zero-dimensional
on this slice, with a point at the prescribed point by repeated
inversion of adjunction (or generic transverse slicing at that
component). But Nadel vanishing gives
$H^1(\mathbf P^3,\mathcal J(-1))=0$ since $3-c'h>0$, contradictory
to the surviving sections of the nonempty finite scheme. On a slice
away from the prescribed point Bertini precludes new divisorial or
positive-dimensional non-klt loci. 

If a component of $G$ instead has multiplicity $s>d/3$, it is a
hyperplane or quadric. In the latter case it is unique and its
primitive base degree $e$ obeys $se\le d-h<s$, thus the quadric
is constant. A unique multiple hyperplane has $se\le d-h<2s$ and
hence is either constant or of base degree one. Its coefficient
forms span at most two constant linear normals; cutting their
common kernel raises $h$. If two high hyperplanes occur and are
conjugate, their product is rational; a primitive product has
base degree $e$ with $se\le d-h<s$, so both hyperplanes are
constant over the algebraically closed constant field. If they
are individually rational, each is constant by the same bound.
Any retained cut whose light space has ambient codimension less
than $n/2$ again has a finite first light order, by the bound on
contained linear spaces. Thus moving divisorial components have an
explicit and inexpensive restriction; this alone does not bound
the total number of successive such cuts.

In the high range $h>2d/3$, every horizontal irreducible component
of the non-klt locus in the binary base and the projective light
space is constant. To see this consider the divisor of bidegree
$(d-h,h)$ on $\mathbf P^1\times P$ and take $c'<3/d$ sufficiently
close to $3/d$. Write $Z$ for the multiplier-ideal scheme. For
$j$ sufficiently large, Nadel with $\mathcal O(-1,j)$ and
$1-c'(d-h)>0$ gives
$H^1(\mathcal J(c'G)(-1,j))=0$, hence
$H^0(\mathcal O_Z(-1,j))=0$. For the projection $\pi$ to $P$,
Serre's theorem now gives $\pi_*\mathcal O_Z(-1,0)=0$.
Any irreducible component $T$ is an associated subvariety of
$\mathcal O_Z$. It supports an inclusion of a rank-one coherent
sheaf on $T$ into $\mathcal O_Z$ (for instance multiply a generic
associated section by enough local equations of its poles and use an
ample twist). If $T\to\pi(T)$ were generically finite, this sheaf
after tensoring by $\mathcal O(-1,0)$ would have nonzero pushforward
at the generic point of $\pi(T)$, contradiction. Hence every
component has fibers of dimension one generically, which forces
$T=\mathbf P^1\times\pi(T)$. The general-fiber multiplier ideal is
the restriction of the total ideal by restriction to a general base
fiber. In particular the unique offending codimension-two linear
subspace in the high range is independent of the base point.

\subsection{All polar pencils on the row \(n=2d\)}
\label{sec:nov27-all-pencils-2d}
Let \(d\ge9\), \(n=2d\), and let \(B\) be a line of polar rank two
on a smooth complex \(X_d^{2d}\). Then
\[
      \frac{\delta_B(X,H)}{n+1}\ \ge\ \frac{d+3}{2d+2}.
\]
We use the three-light log-canonical test, the structural description
of its failures at high contact, and the repeated-quadric estimates
proved above. The new points are the numerical control of successive
linear restrictions and a two-variable argument at intermediate
contact. Write \(b=(2d+2)/(d+3)\). Keep two heavy normal directions
along the line giving a basepoint-free polar pencil. Initially their
common light linear space \(V\subset\mathbf P^{n+1}\) has codimension
\(c=2\). At later steps \(V\) is a fixed linear subspace of that
kernel still containing \(B\), and \(c\) denotes its codimension in
the original ambient projective space. The retained hypersurface
is the restriction to the span of \(V\) and the heavy directions.
All cuts are smooth along \(B\), and the basepoint-free two-polars
argument proves they are integral. Their complete graded pieces
restrict by the hypersurface exact sequence, as required for
localized hyperplane flags. Whenever \(c<d\) the space \(V\) has
dimension greater than \(n/2=d\) and hence is not contained in
\(X\). Write \(j\) for its first light order and \(G_j\) for its
first nonzero homogeneous light form at the generic point of \(B\).

\paragraph{A bounded sequence of higher-contact restrictions.}
Suppose the current data \(j,c\) satisfy
\[
                         c+2(d-j)<d.                       \tag{H}
\]
At order \(j=d\), \(G_d\) is a homogeneous polynomial in the light
variables alone. The join of \(B\) with the projective singular
locus of \(G_d\) maps under the finite Gauss map of \(X\) to the
hyperplanes containing \(V\). Its dimension implies
\(\dim\operatorname{Sing}(G_d)\le c-3\) in \(\mathbf P^{n-1-c}\).
As \(2c<n\), a general plane avoids the singular locus and cuts a
smooth curve. Retain its three light directions. The three-light
radial denominator is \(1+(1-1/d)/(1+3/d)=b\), and the signed
graph bound at the punctured angular points has tangent contact two,
giving averaged denominator \(1+(1-1/d)2/3\le b\).

At \(j<d\), the following tests either finish the proof or give
a fixed linear subspace of \(V\) of codimension at most two,
containing \(B\), on which \(G_j\) vanishes identically.
If \(j\le d/3\), the single-binary calculation suffices. If
\(j\le2d/3\) and every component has multiplicity at most \(d/3\),
the two-light equal-weight calculation suffices. If
\(j>2d/3\) and every component has multiplicity at most \(d/3\),
apply the three-light log-canonical test with coefficient \(3/d\).
If it fails, the Nadel result above gives a codimension-two linear
locus, constant over the binary base, along which \(G_j\) vanishes.

Otherwise a component has multiplicity \(s>d/3\). It is a hyperplane
or a quadric; a quadric with this multiplicity is unique and constant
by the primitive base-degree bound \(s e\le d-j<s\). A quadric of
rank at least five is handled by the repeated-quadric and residual
estimates. A constant quadric of rank three or four contains a
constant linear subspace of codimension two in the light space:
put its form into hyperbolic coordinates over \(\mathbf C\) and
use a maximal totally isotropic space including its radical.
If a repeated hyperplane is rational over the binary base, it
has primitive base-degree \(e\le1\), since \(j>d/3\) gives
\(d-j<2s\). Its coefficients span at most two constant linear
forms. Their common kernel kills \(G_j\). If a repeated hyperplane
has a nontrivial Galois orbit, the orbit has exactly two members.
Their primitive quadratic product is constant since \(j\ge2s\)
gives \(d-j<s\); therefore these hyperplanes are themselves
constant. This covers every outcome in the tests.

In each unresolved case restrict to the indicated subspace. The
new ambient codimension is at most \(c+2\), and the new light
order is at least \(j+1\). Inequality (H) persists and ensures at
each step that the restricted polynomial is nonzero. After at most
\(d-j\) strict increases one reaches the terminal calculation or
an earlier successful test. With initial \(c=2\), this already
treats every initial order \(h>(d+2)/2\).

\paragraph{Intermediate contact and a crossing with controlled exponents.}
We may now assume \(d/3<h\le(d+2)/2<2d/3\) initially, with a
component of multiplicity \(s>d/3\). It must be a unique rational
hyperplane, since \(h<2d/3\). If it has primitive base-degree one,
write its equation \(u\ell_0+v\ell_1\) with constant linear forms
in the light variables. Cut by one general constant linear combination
of \(\ell_0,\ell_1\). The moving equation restricts to a nonzero base
factor times a constant hyperplane. The initial restricted form is
still nonzero of light degree \(h\): only finitely many members of the
chosen pencil can divide its nonzero residual factor at the generic
point. For a constant hyperplane no preliminary cut is necessary.
Write \(x=0\) for the resulting hyperplane in \(V\) and \(V_0\) for
its join with \(B\) inside \(V\). It has ambient codimension at most
four, so is not contained in \(X\). The actual first light order on
\(V_0\) is some \(h<t\le d\).

Choose a general direction in \(V_0\) giving a pure term \(y^t\)
and a complementary direction giving a pure \(x^h\). In these
two light variables every least total-degree monomial has
\(x\)-exponent at least \(s>d/3\). As the weight of \(x\) relative
to \(y\) increases, a lower Newton edge is met whose endpoint
\(x\)-exponents satisfy \(i_0\le d/3<i_1\). Let its primitive
positive weights be \((m,l)\), \(m>l\), and its weight be \(H\).
The pure monomials show \(H\le mh,lt\), whence
\[
       \frac{m+l}{H}\ge\frac1h+\frac1t.
\]
If \(1/h+1/t\ge3/d\), this gives radial denominator at most \(b\).
Writing \(r=m/l\), the minimum \(y\)-exponent on that face is
\[
    j_1=H/l-ri_1\le\min\{r(h-i_1),\,t-ri_1\}
             \le t(1-i_1/h)\le d/3 .
\]
For the middle inequality split at \(r=t/h\); for the last use
\(i_1>d/3\) and \(1/h+1/t\ge3/d\). Thus both coordinate factors
of the face form \(\Phi\) have multiplicity at most \(d/3\).
The endpoint coefficients remain nonzero at general pole roots.

On the punctured angular covering cone, each positive pole group
has graph \(z+\Phi(x,y)=0\). At every noncritical graph point
weighted Euler gives nonzero Hessian, since \(H\ge m i_1>m,l\).
At a critical point the tangent contact is the multiplicity of
the unique vanishing factor on the punctured binary chart. At the
axes the ordinary signed-contact bound therefore gives at most
\(1+d/(d+3)\). On the torus each factor is
\(x^l-\lambda y^m\), \(\lambda\ne0\), of multiplicity \(e\le d/m\).
For \(m\ge3\) the same contact bound suffices. Otherwise \((m,l)=(2,1)\)
and the factor is a smooth affine parabola of degree two, with
\(2e\le d\). It defines with \(z=0\) a smooth conic in ordinary
\(\mathbf P^3\), by homogenizing \(x-\lambda y^2\). Near the selected
point the other factors are units. For \(e\ge2\), apply the conic
flag estimate for \(z+U(x,y)(x-\lambda y^2)^e\) proved in the
repeated-quadric calculation: its weighted initial along that conic
is \(z+U|_C\,Q^e\), with nonzero coefficient at the marked point;
the proof on its marked exceptional fiber is unchanged for a unit
\(U\). The signed bound is \(1+2e/(2e+3)\le1+d/(d+3)\). If \(e=1\)
the point is noncritical. After averaging one regular group with
denominator one and \(d-1\) positive groups with this signed bound,
the angular denominator is at most \(1+(d-1)/(d+3)=b\).

Finally suppose \(1/h+1/t<3/d\). Then
\[
        t>\frac{dh}{3h-d}\ge\frac{d(d+2)}{d+6}>d-4,
\]
since the middle function decreases with \(h\). Restrict the light
space to \(V_0\). Its new order \(t\ge d-3\) and ambient codimension
\(c\le4\) give \(c+2(d-t)\le10<d\) for \(d\ge11\). In degrees
\(9,10\) the same displayed lower bound gives \(t\ge d-2\), hence
\(c+2(d-t)\le8<d\). The preceding bounded sequence therefore
applies. This proves the pencil estimate in every case.

\section{\texorpdfstring{The point estimate in dimension \(2d\)}{The point estimate in dimension 2d}}\label{sec:nov-second-row}
Suppose \(n=2d\), \(d\ge20\), \(N=2d+1\) and \(r=d+2\).
If every line and integral conic \(C\subset X\) satisfies
\[
 \delta_C(X,H)\ \ge\ M:=\frac{(2d+1)(d+3)}{2d+2}
                =r+\frac12-\frac1{d+1},
\]
then \(X\) is uniformly K-stable.
Here and below the short-curve threshold tests all valuations whose
centers contain \(C\). This criterion assumes the all-member theorem
on the preceding row \(n=2d-1\).

Suppose \(\beta=\delta(X,H)\le r\). The smooth pointed hyperplane
estimate and the preceding row give \(\beta>N(r-1)/n\).
Consequently all the finite-fraction tilted-basis comparisons and
the positive-center Nadel contradiction in the preceding chapter
apply verbatim with \(u=n/N,\ b=r/N\), since \(\beta>N/3\).
The formulas there for \(h,g,\alpha,\lambda,U,U_p,P,z,\Delta,Q\)
use \(n,N,r,d\) explicitly and do not otherwise fix their relation.
Use \(f=1/d\), \(\epsilon=10^{-8}\) and downward rounding on a grid of denominator \(10^9\) for the roots in \(h,\alpha\). For \(25\le d\le 162\) use
the share \(a=1-7/d\). For \(d=20,\ldots,24\) take \(1000a=
649,664,678,690,702\), respectively. Round \(Q\) upwards
to the next multiple of \(10^{-4}\), or to zero if this is larger.
For \(d<25\) take \(L=1000\), otherwise \(L=100\), and put
\[
 B=\frac{\left\lfloor L\min_{q=0,Q}
     \frac{N(n+q+P)/U_p-(n-2)}{q+2}\right\rfloor-1}{L}.
\]
The checker verifies \(B^2>d\) and
\((B^2/d-1)(2d-5)+B-d>0\), as well as the signed-flag comparisons
for \(0\le q\le Q\) when \(\tau_p(H)\le B\).
For larger \(\tau_p(H)=t\), the principal-bin endpoint and bin count
\((t_{\rm end},K)\) for each successive interval are as follows.
\[
\begin{array}{c|c|c}
d&\hbox{starting value}&(t_{\rm end},K)\\ \hline
20&13/2&(105/8,80),\ (20,160)\\
21&13/2&(111/8,84),\ (21,168)\\
22&13/2&(29/2,88),\ (22,176)\\
23&13/2&(121/8,92),\ (23,184)\\
24&13/2&(63/4,96),\ (24,192)\\
25&B&(65/4,100),\ (25,200)\\
26&B&(189/8,156),\ (26,208)\\
27&B&(49/2,162),\ (27,216)\\
28,\ldots,162&B&(d,8d)
\end{array}
\]
In the bin check one may omit terms of mass below \(10^{-16}\)
from the nonnegative radial lower sum and round the retained roots
down to \(10^{-6}\). The tent parameter \(c\) is chosen on the grid
of width \(10^{-5}\) and its moment constraint checked exactly.

For \(20\le d\le24\) and \(B<t<13/2\), let the exceptional high-ratio
prime be \(D\in|jH|\), of multiplicity \(k=jt\). Use the
discriminant surface estimate of the preceding row. For \(j\ge12\)
the continuous test takes
\[
 E=d+\frac{d(d-4)+13/2}{12}+\frac1{144},\qquad
 L_*=\frac{E}{2\sqrt E-1/12}.
\]
The maximum of \(d/t+L_*(1-d/t^2)\) on \([B,13/2]\) occurs at
\(2L_*\) clamped to that interval and is less than \(B\).
For \(j\le11\) test the finite list \(jB<k<13j/2\) using
\(I=dj(j+d-4)+k+1\). For \(d=20\), \(j=1,2\), use the additional
upper bound
\[
 I_0\ \le\ (k-1)\left[1+
       \left(\frac{j(d+j-2)^n}{k-1}\right)^{1/(n-1)}\right]
\]
on the discriminant of a general plane section. Here is a proof.
Write \(D=X\cap V(G)\), \(\deg G=j\), using homogeneous projective
normality. For \(j=1\), \(V(G)\) is smooth. For \(j=2\), if \(V(G)\)
were singular at \(p\), its affine equation at \(p\) would be
homogeneous quadratic. Since \(k>2\), that quadratic vanishes on
the tangent space to \(X\). It would have the tangent linear equation
as a factor, contradicting primeness. Thus \(V(G)\) is smooth at \(p\).
On the reduced singular locus of \(D\) the proportionality
\(\nabla G=\sigma\nabla F\) defines a regular section of
\(\mathcal O(j-d)\), since the coordinates of \(\nabla F\) have
no common zero. On every positive dimensional integral component
this section vanishes, so that such a component lies in
\(\operatorname{Sing}V(G)\). Hence \(D\) has an isolated hypersurface
singularity at \(p\) in the smooth germ \(V(G)\).

Use \(n\) independent ambient affine coordinates \(x_i\) as
analytic coordinates on \(V(G)\); let \(z\) be the remaining affine
coordinate with \(G_z(p)\ne0\). Generic linear flags in the \(x_i\)
give legitimate pointed linear sections; the \(x_i\) themselves
may first be taken general. The local Milnor number \(\mu_n\) is the
intersection multiplicity of \(G=0\) and the \(n\) polynomials
\(G_zF_{x_i}-F_zG_{x_i}\) of degree at most \(d+j-2\). The
Jacobian critical scheme is isolated at the point. B\'ezout bounds
\(\mu_n\le j(d+j-2)^n\).
Teissier's log-convexity of the sectional Milnor numbers gives
\(\mu_2/\mu_1\le(\mu_n/\mu_1)^{1/(n-1)}\), where
\(\mu_1=k-1\). See B.~Teissier, \emph{Cycles \'evanescents,
sections planes et conditions de Whitney}, Ast\'erisque
\textbf{7--8} (1973), 285--362. Finally the reduced plane-curve
discriminant exponent is \(I_0=\mu_2+k-1\), as in the earlier
surface proof. This proves the displayed additional bound.
All the resulting point comparisons at \(q=0,Q\) are strict.

For \(163\le d\le1023\) the coarse comparison suffices. Put
\(\epsilon=r^{-4}\), \(x=(r-N/3-3\epsilon)/d\),
\[
 U=r+bx+2\epsilon,\quad
 z=b+x\{u(1-b)-\log(d)/N\}-u\epsilon,\quad
 \Delta=\min\{M-U,(z-1/3)/u\}-\epsilon.
\]
The elementary affine point envelope \(S_H\le J+Gq/n\)
excludes \(q\le Q=\max(0,U((1-z)/\Delta+u)-n)\). The arithmetic
for \(20\le d\le1023\) is contained in
\path{nov27-bin-second-row.py}.

For \(d\ge1024\), write \(y=\sqrt d\) and \(E=\log(d)/6\). One has
\(1/3\le x<1/3+5/(3d)\), \(b=1/2+3/(2(2d+1))\),
\(u(1-b)>1/2-1/d\). It follows that
\[
 U<r+1/6+1.09/d,\quad z>2/3-E/d .
\]
Take \(\Delta=1/3-(1+E)/d\). Since \(E>1.09\), this gives
\(U+\Delta<M\), \(z-u\Delta>1/3\). The point envelope on this row
has \(J<1.515,\ G<1.56y\); these follow directly by substituting
\(n=2d\) in the two maxima of the previous large-index proof, or
from \(G\le\max\{\sqrt{2d}+4,\sqrt{d/2}+\log(12(d+2)^2)\}\).
Thus \(n-UJ>0.48d\), \(UG/n-1<0.782y\). Unless the point is already
excluded, \(q>0.613y\), so
\((n+q)/U>2+0.46/y\). Indeed \(0.153y>4.336+0.998/y\) at \(y=32\)
and thereafter. The final point contradiction
\(z+\Delta((n+q)/U-u)>1\) follows from
\[
   (0.46/3)y>2E+1+0.46(1+E)/y .
\]
After division by \(y\) each term on the right decreases for
\(y\ge32\), and the endpoint is strict. This proves the point criterion.

\subsection{The degree-nineteen endpoint in dimension thirty-eight}
\label{sec:nov27-cubic-degree19}
The point calculation in dimension \(n=2d\) with the local mark
\(M=(2d+1)(d+3)/(2d+2)\) also succeeds for \(d=19\).
Use \(a=633/1000\), \(B=5.404\), \(T=13/2\) and the two
principal-bin data \((25/2,76),(19,152)\). In the finite surface
test take the sectional Milnor bound of the previous subsection
for \(j\le2\). For \(3\le j\le11\) its genus bound already gives
every comparison except \((j,k)=(3,19)\). For this pair it
suffices to prove that the discriminant exponent \(I\) of a general
pointed plane curve section of the prime \(D\in|3H|\) is at most
\(884\). We prove this directly. As before, represent \(D\) as
\(X\cap(G=0)\). The cubic \(G\) is irreducible, since every
nonconstant homogeneous factor defines an effective divisor on \(X\).
If \(G\) is smooth at \(p\), the preceding isolated Milnor argument
with \(j=3\) suffices.

Suppose \(G\) is singular at \(p\). Write \(z=0\) for the tangent
hyperplane of \(X\). Multiplicity three at \(p\) would make its
affine equation homogeneous cubic, divisible by \(z\) because
\(k>3\). At multiplicity two its quadratic term likewise must be
\(z\ell\). If \(\ell\) is proportional to \(z\), the cubic term
also vanishes on \(z=0\), since \(z|_X\) has order at least two.
This again makes \(G\) reducible. Thus \(z,\ell\) are independent.

We claim that the singular locus of \(G\) through \(p\) contains a
linear subspace \(W\) of ambient codimension two. If it had local
codimension at least three, a general pointed affine three-space
would cut a cubic surface with isolated singularity at \(p\).
Indeed Bertini gives smoothness away from the original singular
locus and from the base point at each cut. On this three-space take
local analytic coordinates \((u,v,t)\) with \(u=0\) defining \(X\).
The restriction of \(G\) to \(u=0\) has multiplicity \(k\), and
its quadratic part transverse to this restriction has two independent
factors. We may thus write \(G=uv+R(v,t)\), \(R\in(v,t)^k\),
up to a remainder in \((u,v,t)^k\): more explicitly, first write
\(G=uV+R_0\) with \(R_0\) of order at least \(k\), then replace
one of the coordinates by \(V\), which has independent linear part.
In these coordinates the two derivatives in \(u,v\) eliminate
\(u,v\) as power series of order at least \(k-1\) in \(t\). The
remaining derivative has order at least \(k-1\), so the Milnor
number is at least \(k-1=18\). Three affine partials of the cubic
have degrees at most two, and B\'ezout gives the upper bound eight
on their isolated local intersection. This contradiction furnishes
a codimension-two component \(W\) through \(p\). A general
projective plane cuts an integral cubic curve, by Bertini, which
has at most one singular point by its genus. Transversality to
\(W\) thus gives \(\deg W=1\). Since all first derivatives vanish
along this linear space, \(G\in I_W^2\). The kernel of the Hessian
at \(p\) contains the tangent space of \(W\). Their dimensions agree,
so in the affine chart at \(p\), \(W=(z=\ell=0)\).

Take a general linear projection of that chart to \(n\) coordinates
whose kernel is transverse to \(z=0\), and complete its coordinates
as \((w,y_1,\ldots,y_{n-1})\), with \(I_W=(z,w)\). Here \(w\)
is the nonzero linear combination of the two generators vanishing
on the chosen projection direction. Write \(z=tw\). The proper
transform has polynomial equation
\[
  G_*(t,w,y)=G(tw,w,y)/w^2=0,\qquad \deg G_*\le4 .
\]
At the exceptional point corresponding to the tangent \(z=0\), this
polynomial has a simple root \(t=0\), since the two original tangent
hyperplanes are distinct. The general choice of projection makes the
other root finite and simple as well. Both roots lift analytically
to \(t_1(w,y),t_2(w,y)\). The two smooth branches
\(H_i=(z=t_i(w,y)w)\) of the cubic meet transversely along \(W\).
This exhausts the local equation since it has order two and the
product of the two distinct smooth factors divides it. The first
branch is tangent to \(X\), and the second is transverse. The
multiplicities of their divisors on smooth \(X\) are respectively
\(k-1\) and one. Since \(X\) and \(H_1\) are smooth graphs over
\((w,y)\), the same difference of graph functions, up to a unit,
defines their intersection on either graph, also on sections by
linear coordinates in \((w,y)\).

The singular locus of the divisor \(X|_{H_1}\) is contained near
\(p\) in \(W\). Indeed any positive-dimensional component not
contained in \(W\) would give a positive-dimensional component of
the projective singular locus of \(D\) whose general point is smooth
on \(G\). The proportional-gradient section of \(\mathcal O(3-d)\)
in the previous subsection excludes it. Moreover \(W\) has tangent
space contained in the tangent space of \(H_1\) at each point, so
any such singularity on \(W\) has \(T_xX\supset W\) projectively.
The Gauss map of \(X\) is finite, as its hyperplane pullback is
\((d-1)H\). Hyperplanes containing \(W\) form a pencil. The local
singular locus therefore has dimension at most one. A general
hyperplane through \(p\) in the projection coordinates cuts an
isolated singularity on \(H_1\), by Bertini off the base point
and proper intersection with any singular curve through it.

After this cut, the ambient incidence chart has \(n\) coordinates,
and the smooth degree-at-most-four hypersurface \(G_*=0\) has
dimension \(s=n-1=37\). Its local coordinates are the \(s\)
remaining independent projection coordinates. The pulled-back
polynomial \(F(tw,w,y)\) has total degree at most \(2d\).
Clearing the nonzero \(t\)-partial of \(G_*\) in the Jacobian of
its restriction gives \(s\) polynomials of degree at most
\(D_*=4+2d-2=40\). The isolated Milnor number satisfies
\(\mu_s\le4D_*^s\) by B\'ezout. Its one-dimensional sectional
Milnor number is \(k-2\); Teissier's inequality gives for the
general plane section of this tangent branch
\[
 I_1=\mu_2+k-2\le(k-2)\left[1+
       \left(\frac{4D_*^s}{k-2}\right)^{1/(s-1)}\right].
\]
The branch from \(H_2\) is smooth on the general section. The mutual
intersection multiplicity of the two branches is the intersection
multiplicity of \(X\) with the corresponding line in \(W\), which
is at most \(d\). Indeed \(W\) cannot lie on \(X\) since its
dimension exceeds \(n/2\), and the general line through \(p\) in
\(W\) does not lie on \(X\). The selected ambient three-spaces are
general pointed three-spaces when the projection center and the
flag in its quotient are allowed to vary. The plane-curve product
discriminant formula consequently gives
\[
            I\le I_1+2d<779<884 .
\]
This proves the required surface estimate in the last cubic case.
The remaining rational tests are given in
\path{nov27-bin-degree19.py}. Consequently the line and conic
mark in the preceding dimension-\(2d\) theorem implies uniform
K-stability here as well.

\subsection{Improved integer estimates in small degrees}\label{sec:nov27-integer-second}
Let \(n=2d\). The following table gives a local threshold \(M=(2d+1)/b\)
for every line and every integral conic and proves uniform K-stability
of every smooth hypersurface in the indicated degrees. Put
\[
              b=1+(1-1/d)\frac{T}{T+1}.
\]
The last three columns give parameters for the point criterion.
\[
\begin{array}{c|c|c|r|r|l}
d&T&b&1000a&1000B&(t_{\rm end},K)\\ \hline
10&28/9&311/185&575&4230&(6,30),(10,50)\\
11&7/2&169/99&560&4277&(6,22),(11,55)\\
13&4&113/65&750&6141&(13,39)\\
14&9/2&271/154&641&5134&(14,70)\\
16&5&57/32&794&7365&(16,48)\\
17&11/2&397/221&698&5974&(17,85)\\
19&6&241/133&825&8569&(19,57)
\end{array}
\]
We first refine the rank-two line calculation in
Section~\ref{sec:nov27-all-pencils-2d}. Use \(1/T\) as the coefficient
in the log-canonical test. At any nonterminal light order \(j\le d-1\)
one has \(j\le3T\). A single light suffices for \(j\le T\). For
\(j\le2T\) without components of multiplicity greater than \(T\),
use two general lights. The radial reciprocal-weight ratios are at
least \(1/T\) and in the latter case the binary root multiplicities
are at most \(T\). For \(j>2T\) in the absence of a high multiplicity
component, use the three-light log-canonical test. Its radial ratio
satisfies \(3/j\ge1/T\). If log canonicity fails, Nadel with
\(c'<1/T\) sufficiently close to \(1/T\) gives the same unique
codimension-two linear locus since \(c'j<3\). Furthermore \(j>2T\)
implies \(d-j\le T\) in each row, so \(1-c'(d-j)>0\) and the
base-independence proof applies.

For every integer \(s>T\), \(s\le j<d\), the table has
\(3s>j\) and \(2s>d-j\). If \(j\ge2s\), it also has \(s>d-j\).
Thus the same primitive-degree argument handles every repeated
hyperplane and every constant low-rank quadric. For a quadric of
rank at least five its residual degree \(u=j-2s<T\). In the marked
conic calculation the new radial graph weight is at most
\(2s+u=j\le3T\). Each noncoordinate parabolic factor has
multiplicity \(e\) with \(2e\le j\le3T\), so its signed denominator
\(1+2e/(2e+3)\) is at most \(1+T/(T+1)\). The coordinate factors
have multiplicities at most \(u<T\), and the noncritical contact two
is at most \(T\). This proves the required estimates for the
repeated-quadric calculations with the sharper constant.

Keep the same potential condition \(c+2(d-j)<d\) for ambient
codimension \(c\). It guarantees termination of each sequence of
restrictions costing at most two per unit increase in \(j\).
At terminal order \(j=d\), retain \emph{four} general lights. Indeed
the previous Gauss-map bound on \(\operatorname{Sing}G_d\) gives it
codimension at least \(n+2-2c\ge4\) in projective light space, as
\(c\le d-1\). A general three-plane cuts a smooth surface. The radial
ratio satisfies \(4/d\ge1/T\). At a graph point away from its zero
locus Euler's identity gives tangent contact two, and on its zero
locus smoothness gives nonvanishing gradient; Euler again gives
contact two. The signed graph estimate therefore applies at every
point on the punctured cover.

For intermediate initial orders where \(2+2(d-h)\ge d\), a component
of multiplicity \(s>T\) is necessarily a unique rational hyperplane,
since then \(h<2s\) in the table. Make it constant by the same initial
cut, if needed. On its kernel the next order is \(t\), \(h<t\le d\),
and \(c\le4\). If \(t=d\) use the terminal computation directly.
If \(1/h+1/t<1/T\), the finite integer checks give \(4+2(d-t)<d\)
and hence return us to the terminating sequence, with one exception,
\((d,h,t)=(10,6,7)\), treated below. In the remaining
case choose the crossing edge with endpoint exponents
\(i_0\le T<i_1\) as before. Its smallest second exponent satisfies
\(j_1\le t(1-i_1/h)\le T\) and its radial reciprocal-weight ratio
is at least \(1/h+1/t\ge1/T\). Each torus factor
\(x^\ell-\lambda y^m\) of multiplicity \(e\) on the face has
\(me\le t\): its largest \(y\)-exponent is at most the face weight
divided by \(\ell\), at most \(t\) by the pure term \(y^t\).
Here \(t\le d-1\le3T\), giving \(e\le T\) when \(m\ge3\) and
the requisite \(2e\le3T\) when \((m,\ell)=(2,1)\). The conic-flag
bound applies in that last case. Thus the positive pole groups have
denominator at most \(1+T/(T+1)\) and the averaged denominator is
at most \(b\).

In the exceptional degree-ten case a constant repeated hyperplane
costs just one cut, giving \(c=3\) and \(c+2(10-7)=9<10\).
For a moving repeated hyperplane the base degree of the sextic is
four. Its equation is therefore
\[
        G_6=(uA+vB)^4R_2 ,
\]
with independent constant linear forms \(A,B\), binary base
coordinates \(u,v\), and a constant quadric \(R_2\).
Choose two complementary directions for \(A,B\) and one general
light direction \(y\) in their common kernel. The latter has a
nonzero seventh-order coefficient by the assumption \(t=7\).
Let \(k=4+\operatorname{ord}_{(A,B)}R_2\).
The selected three-space preserves this order: choose \(y\) at
which one of the lowest nonzero transverse coefficients does not
vanish. With weights \((k+1,k+1,k)\) on \(A,B,y\) the face is
\[
      \Phi=y^{6-k}\big(H_k(A,B)+\gamma y^{k+1}\big),
      \qquad \gamma\ne0,\qquad 4\le k\le6 .
\]
Indeed every degree-six term of least weight has transverse degree
\(k\), only the pure seventh-degree term has weight \(7k\), and
higher degrees have larger weight. The binary factor \(H_k\) is
\((uA+vB)^4\) times a constant form of degree \(k-4\). Thus every
root over the geometric generic base has multiplicity \(e\le4\).
The same holds at a general pole root. The radial inequality is
\(7k\le (28/9)(3k+2)\).

At a nonzero point with \(y\ne0\) the zero hypersurface of \(\Phi\)
is smooth. On \(y=0\) at a root of \(H_k\), its transverse equation
is analytically \(y^{6-k}(x^e+\gamma y^{k+1})\), after absorbing
the unit in a local root coordinate \(x\). The compact edge is
nondegenerate on the torus. Its primitive weighted ray, together
with the coordinate divisors and the strict transform, computes the
local threshold as the minimum of
\[
        1,\qquad \frac1{6-k},\qquad \frac{k+1+e}{7e},
\]
where a zero denominator is omitted. This follows directly by
subdividing the two cones adjacent to the edge normal; on each
cone the order and discrepancy are linear and the strict transforms
meet the strata transversely. Each term is at least \(9/28\).
Away from the roots the coordinate multiplicity is at most two.
The ordinary three-variable signed graph bound now applies with
\(T=28/9\). To recall its numerical requirements, slice transversely
to the weighted orbit at a critical point. The resulting ordinary
plane polynomial has degree at most seven and threshold at least
\(1/T\). At an equal-weight tangent face its binary multiplicities
are at most \(T\), unless a unique direction has multiplicity
greater than \(T\), in which case take the crossing edge in this
direction. The monomial valuation and the threshold bound give
\(H\le T(m+\ell)\); hence both axis multiplicities on that edge
are at most \(T\). Every torus factor \(x^\ell-\lambda y^m\) has
\(me\le7\le3T\), so for \(m\ge3\) its multiplicity is at most
\(T\), and for \(m=2\) the conic-flag signed estimate applies.
These checks give signed denominator \(1+T/(T+1)\) by weighted
adjunction on the plane graph and ordinary hyperplane adjunction.
Weighted Euler handles a noncritical graph point by tangent contact
two. Averaging the pole groups gives the asserted degree-ten bound.

For the remaining short curves the proper-plane, high-rank,
ordinary-axis and conic polar-pencil bounds are respectively bounded
by \(2d/(d+2),\,2-4/d,\,3/2,\,(2d^2-2d+1)/(d(d+1))\),
whose maximum is at most the listed \(b\). Apply the previously
proved ternary line and conic tests with this \(b\) for all
\(2\le h\le d\), with \(1\le s\le2\) for lines and \(1\le s\le6\)
for conics. Greater \(s\) already has residual intersection at most
\(2-4/d\). Include the additional conic tests at \(s=1,2\).
All finite comparisons pass. Finally apply the point criterion with
the new threshold \(M\). Its formulas hold for any chosen
short-curve threshold satisfying the displayed numerical comparisons.
Use \(f=1/d\), \(\epsilon=10^{-8}\), and the tabulated \(a,B\).
One has \(B^2>d\) and \((B^2/d-1)(2d-5)+B-d>0\) throughout.
The two point endpoints for \(\tau_p(H)\le B\) and the indicated
principal-bin comparisons from \(B\) to \(d\) all hold strictly.
These rational substitutions, and every integer inequality in the
pencil reduction, are reproduced by
\path{nov27-pow-integer-second.py}. The preceding-row hyperplane
estimate and the positive-center argument now imply
\(\delta(X,H)>r=d+2\).

\subsection{Degrees divisible by three}\label{sec:nov27-divisible-second}
Suppose \(d=3k\), \(k\ge3\), \(n=2d\), and put
\(T=k-1/6\). Every polar-pencil line has normalized denominator
at most \(b=1+(1-1/d)T/(T+1)\).

At terminal contact use the four-light argument above whenever
the ambient codimension \(c<d\). Indeed \(4/d\ge1/T\).
At a nonterminal order \(j\ge2k\) one has \(j\le d-1<3T\).
A log-canonical light form at coefficient \(1/T\) therefore
passes the three-light test. For \(j>2k\), its base degree
\(d-j\le k-1<T\) ensures that any codimension-two Nadel locus is
constant, by the product Nadel argument. If a component instead
has multiplicity \(s>T\), then \(s\ge k>d-j\); it is a
constant hyperplane or a constant quadric, with conjugate
hyperplanes handled by their primitive product. The repeated-quadric
estimates apply when the rank is at least five since
\(j\le3T\) and \(j-2s<T\); a low-rank quadric admits a
constant codimension-two linear cut. Thus at \(j>2k\) every
unsuccessful test can be followed by at most two cuts raising \(j\).

At \(j=2k\) a repeated irreducible quadric has multiplicity exactly
\(k\) and is the entire form up to scalar. The pure-quadric argument
also works when it varies in the binary base: a general constant
three-dimensional light subspace cuts its geometric generic
quadric in a smooth conic since its rank is at least three,
and the same holds at general pole roots. Here \(3/(2k)\ge1/T\)
and the signed increment \(2k/(2k+3)\le T/(T+1)\).
If instead a hyperplane component of multiplicity at least \(k\)
occurs, choose one at each of two distinct general base points.
Their intersection is a constant codimension-at-most-two space
on which the full light form vanishes. In fact a fixed point there
gives a base polynomial with two roots of multiplicities at least
\(k\), whereas its degree is \(k\). For a conjugate component
take a local branch after an étale base change at each general
base point. Divisibility by the power of its local equation gives
the same root orders, as the local base map is unramified.
In the remaining case every component has multiplicity at most
\(k-1\). The plane-curve gap proved below shows that again
a failed test can be followed by a constant codimension-two cut.
Consequently \(c+2(d-j)<d\) suffices for termination in every
case with \(j\ge2k\). Initially \(c=2\), so this condition holds
for \(j\ge2k\) since \(2+2k<3k\).

For an initial \(h\le T\) one light suffices. For
\(k\le h\le2k-1\) two lights suffice in the absence of a component
with multiplicity greater than \(T\), since \(h\le2T\).
Otherwise there is a unique rational hyperplane of multiplicity
\(s\ge k\). Its primitive base degree \(e\) is at most one
for \(h>k\) and at most two for \(h=k\), by \(se\le d-h\).
Make this hyperplane fixed by imposing at most \(e\) general
linear equations among its coefficient forms: the subspaces on
which the coefficient vector has a prescribed projective direction
cover a dense open subset, so a general such subspace preserves
the nonvanishing of the order-\(h\) form. One further cut by that
hyperplane gives an order \(h<t\le d\), since the resulting
ambient codimension is at most \(5<d\).
At \(t=d\) use the terminal computation directly.
If \(t<d\) and \(1/h+1/t\ge1/T\), use the two-variable crossing
as above; now each torus factor satisfies \(me\le t<3T\).
Otherwise \(t\ge2k+1\), since
\[
       \frac1{2k-1}+\frac1{2k}>\frac1{k-1/6}.
\]
Moreover \(h>k\), as \(1/k+1/(3k)>1/T\). Thus the new codimension
is at most four, and \(c+2(d-t)\le4+2(k-1)<3k\).
The high-order termination argument applies. This exhausts all
polar-pencil lines.

\paragraph{A gap for plane curves of degree \(2k\).}
Let \(k\ge 2\) and let \(D\) be an effective plane divisor of degree
\(2k\), each irreducible component appearing with multiplicity at most
\(k-1\). If \(\operatorname{mult}_pD\le k\), then
\[
            \operatorname{lct}_p(D)\ge\frac1{k-1/6}.
\]
For multiplicity less than \(k\) use the elementary inverse-multiplicity
bound. Suppose the multiplicity equals \(k\). For a general linear
projection \(x\) centered at \(p\), write the distinct convergent
Puiseux roots through \(p\) as \(y_i(x)\), of weights \(w_i\) equal
to their component multiplicities. Thus \(\sum_iw_i=k\),
\(w_i\le k-1\) and at least two distinct roots occur. Put
\[
             t_0=\min_{i\ne j}\operatorname{ord}_x(y_i-y_j).
\]
Write the global components as \(D_\alpha\) of degrees \(d_\alpha\)
and multiplicities \(m_\alpha\). The sum of the contacts over ordered
pairs of distinct roots of the same component through \(p\) is at
most \(d_\alpha(d_\alpha-1)\): this is the local discriminant order
of its Weierstrass polynomial, bounded by intersection with a general
polar under the chosen generic projection. The sum over roots of
two different components, in one ordering, is their local
intersection multiplicity, at most \(d_\alpha d_\beta\). Therefore
\[
\begin{aligned}
t_0\left(k^2-\sum_i w_i^2\right)
 &\le \sum_{i\ne j}w_iw_j\operatorname{ord}_x(y_i-y_j)\\
 &\le \sum_\alpha m_\alpha^2d_\alpha(d_\alpha-1)
       +2\sum_{\alpha<\beta}m_\alpha m_\beta d_\alpha d_\beta
 \le (2k)(2k-1).
\end{aligned}
\]
The smallest possible positive \(k^2-\sum_iw_i^2\) is \(2(k-1)\),
attained at the partition \(k=(k-1)+1\). In particular
\(t_0\le k(2k-1)/(k-1)\le6k-1\).

For a divisorial valuation centered at \(p\), normalize \(v(x)=1\).
Extend it to the finite local cover \(x=u^N\) splitting the Puiseux
roots, keeping this normalization, and set
\(Y=\max_i v(y-y_i(u^N))\). The log discrepancy on the plane is at
least \(1+Y\). Indeed on the smooth cover the pair of coordinate
divisors \(u=0,\ y-y_i(u^N)=0\) has simple normal crossings;
the ramification contributes \((N-1)v(u)\), giving the lower bound
\(Nv(u)+v(y-y_i(u^N))\) for every \(i\). One of the root distances
is at most \(t_0\), by the ultrametric inequality on a pair realizing
that minimum. Hence
\[
       \frac{v(D)}{A(v)}
       \le \frac{(k-1)Y+\min(Y,t_0)}{1+Y}
       \le \max\left\{k-1,\frac{k t_0}{1+t_0}\right\}
       \le k-\frac16.
\]
Divisors already on the plane have coefficient at most \(k-1\).
This proves the gap.

Apply it to the moving light form of bidegree \((k,2k)\), in the
absence of a component with multiplicity at least \(k\).
If the three-light test fails at \(T=k-1/6\), the Nadel plane
argument gives a unique codimension-two linear locus \(W_\eta\)
over the geometric generic binary base. Its transverse plane section
at a general point fails log canonicity with coefficient \(1/T\);
restriction of the multiplier ideal to a general plane gives this
claim. Global component multiplicities on that plane are at most
\(k-1\) by Bertini and generic transversality, so the gap proves
that the generic order of the light form along \(W_\eta\) is at least
\(k+1\). Uniqueness descends the linear space to the function field.
After shrinking the binary base it spreads to a relative linear
subspace \(W\), with the light equation in \(\mathcal I_W^{k+1}\).
Indeed the relative normal jets of orders less than \(k+1\) vanish
generically and the successive jet bundles on \(W\) are locally
free. Fix a general base value \(u_0\) and any point \(w\in W_{u_0}\).
The polynomial \(u\mapsto G(u,w)\) has vanishing order at least
\(k+1\) at \(u_0\) and degree at most \(k\), so it vanishes identically.
Thus \(W_{u_0}\) is a constant codimension-two restriction on which
the entire generic form vanishes.

\paragraph{Uniform K-stability in degrees \(9,12,15,18\).}
The preceding line estimates complete the following point calculations
at \(n=2d\), with \(M=(2d+1)/b\).
\[
\begin{array}{c|c|c|r|r|c|l}
 d&T&b&1000a&1000B&t_{\rm start}&(t_{\rm end},K)\\ \hline
 9&17/6&343/207&482&3695&37/8&(6,27),(9,45)\\
 12&23/6&601/348&590&4564&B&(35/4,36),(12,60)\\
 15&29/6&133/75&662&5416&B&(15,75)\\
 18&35/6&1333/738&713&6253&B&(18,90)
\end{array}
\]
The ordinary high-rank, proper-plane and conic polar-pencil bounds
are at most \(b\) in every row. The finite ternary line and conic
tests with these \(b\) pass over their full integer ranges as in
Section~\ref{sec:nov27-integer-second}, including the additional
tests for split degrees one and two. In each row the signed point
tests for \(\tau_p(H)\le B\) and the principal bins from
\(t_{\rm start}\) through \(d\) are strict. In degree nine, between
\(B\) and \(37/8\), apply the discriminant comparison from
Section~\ref{sec:nov-second-row} with the upper endpoint \(37/8\)
in place of \(13/2\), and \(H=12\). For \(j=1,2\) include the
sectional Milnor upper bound proved there and round the positive
\((n-1)\)-st root upwards to a multiple of \(10^{-6}\).
The comparisons for the finitely many \(j\le11\) and the single
continuous estimate for \(j\ge12\) are strict. Their rational
certificate is \path{nov27-pow-integer-second.py}. This proves
\(\delta(X,H)>d+2\) in all four degrees.

\subsection{The octic in dimension sixteen}\label{sec:nov27-octic-lines}
Let \(X=X_8^{16}\), \(T=5/2\), and \(b=1+(7/8)T/(T+1)=13/8\).
We prove \(\delta_C(X,H)\ge 17/b=136/13\) for every line or integral
conic \(C\).
For lines other than polar pencils, and for conics, apply the
comparisons at the end of Section~\ref{sec:nov27-integer-second}.
The four separate bounds there are
\(8/5,\,3/2,\,3/2,\,113/72\), all at most \(13/8\).
The finite ternary comparisons are recorded in \path{nov27-pow-octic.py}.

For polar pencils use the same notation as in
Section~\ref{sec:nov27-all-pencils-2d}. The common light space has
initial ambient codimension \(c=2\). Whenever \(c<8\) its dimension
exceeds eight, so it has a first light form \(G_h\), of degree \(h\)
in the lights and base degree \(8-h\). At the terminal order \(h=8\)
the codimension of its singular locus is at least \(18-2c\ge4\)
in projective light space. Retaining four general lights gives a
smooth surface and radial ratio \(4/8>1/T\); the noncritical signed
graph calculation applies. At \(h=2\) one light suffices. At
\(3\le h\le5\) two lights suffice provided every component has
multiplicity at most two.

At \(h=6,7\), when every component has multiplicity at most two,
use the three-light test if \((\mathbf P_{\rm light},G_h/T)\)
is log canonical. Otherwise the bidegree Nadel argument gives a
constant codimension-two linear center, since \(h<3T\) and
\(8-h<T\). A component of multiplicity at least three has degree
at most two. A repeated hyperplane at these orders is constant,
as follows by applying the primitive base-degree bound to its
Galois orbit. A repeated quadric is likewise constant. In ranks
three and four it contains a fixed codimension-two linear subspace;
in greater rank use the repeated-quadric estimate of
Section~\ref{sec:nov27-integer-second}: its residual form has degree
at most one, \(h\le3T\), and \(T\ge2\). Thus either the estimate
holds or a restriction of codimension at most two raises \(h\).
If \(c+2(8-h)<8\), these restrictions terminate with \(c<8\).
This treats every initial \(h\ge6\).

At \(h=3,4,5\) a component of multiplicity \(s\ge3\) is a unique
hyperplane. Its primitive base degree is at most one. At \(h=3,4\)
make it constant by at most one preliminary cut and let \(t\) be
the next order on its kernel, where \(c\le4\). If \(t=8\), apply
the terminal test. Otherwise whenever \(1/h+1/t\ge1/T\) take the
crossing edge of Section~\ref{sec:nov27-integer-second}. Its torus
factors of multiplicity \(e\) satisfy \(me\le t\le7\le3T\).
Among \(3\le h\le4\), \(h<t\le7\), the crossing fails only at
\((h,t)=(4,7)\), where \(4+2(8-7)<8\).
At \(h=5\) a constant high hyperplane costs one cut, and then
\(c+2(8-t)\le3+2(8-6)<8\). For a moving high hyperplane, the
entire base degree three is used by its cube. Write
\[
              G_5=(uA+vB)^3R_2,\qquad K=(A=B=0),
\]
where \(A,B\) are independent constant linear light forms and \(R_2\)
is a constant quadric. The first order \(t\) on \(K\) exists,
since \(c=4\). If \(t\ge7\), use \(4+2(8-t)<8\). We treat \(t=6\)
below. When \(G_6|_K\) has a component of multiplicity at least three,
a constant hyperplane costs one cut, giving \(c=5,t'\ge7\),
and a repeated irreducible quadric accounts for the entire sextic.
Its rank is at least three and its pure-quadric three-light calculation
applies: \(3/6\ge1/T\), \(6/(6+3)\le T/(T+1)\), and \(2\le T\).
We may therefore suppose that every geometric generic component
of this sextic has multiplicity at most two.

If \(R_2|_K=0\), let \(j=3+\operatorname{ord}_{(A,B)}R_2\in\{4,5\}\).
Take two complementary \(A,B\) directions and one general \(K\)-light
\(y\). The weights \((j+1,j+1,j)\) have face
\[
                y^{5-j}\bigl(H_j(A,B)+\gamma y^{j+1}\bigr),
                \qquad \gamma\ne0 .
\]
Here \(H_j\) is the moving cubic power times a fixed binary form
of degree \(j-3\); at a general pole each root has multiplicity
\(e\le3\). The general directions preserve both the indicated
transverse order and the nonzero pure sextic. The radial check is
\(6j\le T(3j+2)\). As in the exceptional degree-ten calculation,
the punctured threshold at a root on \(y=0\) is the minimum of
\(1,1/(5-j),(j+1+e)/(6e)\), with zero denominators omitted.
This is at least \(4/9>1/T\). Elsewhere the zero divisor is smooth
off \(y=0\) and has multiplicity at most one at the other points
of \(y=0\). The signed plane-graph argument after transverse
restriction applies because the total degree is at most \(6<3T\).

If \(R_2|_K\ne0\), retain two general \(K\)-lights \(y_1,y_2\)
as well as \(A,B\). The binary quadric \(r=R_2|_{\langle y_1,y_2\rangle}\)
is nonzero. By generic transversality the restricted sextic \(g\)
has root multiplicities at most two. These assertions hold at
general pole roots by openness. With weights \((4,4,3,3)\) the
least weight is 18. The quintic contributes
\((uA+vB)^3r(y_1,y_2)\), the sextic contributes \(g(y_1,y_2)\),
and higher forms have weight at least 21. The radial ratio is
\(14/18>1/T\). After a linear change in the two equal-weight
directions at a pole, the graph takes the form
\[
     Z+f(x,y_1,y_2)=0,\qquad f=x^3r(y_1,y_2)+g(y_1,y_2),
\]
and is independent of the fourth light \(v\).

Apply ordinary hyperplane adjunction at the marked point to
\(v=v_0\). The ordinary signed graph term becomes
\(5S_{\mathbf P^4}+\operatorname{ord}(Z+f)\).
First suppose \((x,y_1,y_2)\ne0\). At any noncritical point
weighted Euler gives a nonzero Hessian, so tangent contact is two.
At a critical point \(f=0\) and the local threshold of \(f\) is
at least \(1/2\). Indeed, for \(y_1=y_2=0\), \(x\ne0\), the
multiplicity equals two, which suffices for this threshold. For
\((y_1,y_2)\ne0\) fix one nonzero \(y\)-coordinate on a section
transverse to the weighted orbit and write \(w\) for the binary
projective parameter. Up to units the section is
\[
      x^3 w^\alpha A(w)+w^\beta B(w),\qquad
      0\le\alpha,\beta\le2,\qquad A(0)B(0)\ne0.
\]
For \(\beta<\alpha\) the equation is a coordinate power times a
unit. For \(\alpha=\beta\) the second factor, when zero, is smooth
and transverse to \(w=0\), by its nonzero \(x\)-derivative.
For \(\alpha<\beta\), \(x=0\), absorb the units to get
\(w^\alpha(x^3+w^{\beta-\alpha})\).
Its threshold is the minimum of \(1,1/\alpha\) and
\((3+\beta-\alpha)/(3\beta)\), hence is at least \(1/2\).
The case \(x\ne0\) here gives a coordinate power. Weighted scaling
exhibits the punctured germ as a local product of its transverse
section with a smooth parameter, up to an invertible factor.
The ordinary affine plane transverse to the orbit therefore inherits
this threshold. Its degree is at most six. Apply the signed
plane-graph estimate as above, using \(6\le3T\) and \(T\ge2\).

Finally the marked point may project to the origin of
\((x,y_1,y_2,Z)\). Extract that origin with primitive weights
\((4,3,3,18)\). The complete-series expected order in its
\(5S_{\mathbf P^4}\) normalization is 28, and the graph order is 18.
The radial signed ratio is \(1+18/28\le1+T/(T+1)\). The graph is
itself weighted homogeneous. In the angular calculation the
ordinary monomial basis splits by weighted degree and its aggregate
is the complete ordinary polynomial series of degree at most \(m\),
as in the residual-quadric calculation. On the punctured affine
cover every point of this graph projects to a nonzero triple
\((x,y_1,y_2)\). Its signed comparison is the one just proved.
Weighted inversion of adjunction at the origin now gives the same
signed bound. Including the free light by ordinary adjunction
and averaging the positive pole groups of total share \(7/8\)
give the asserted denominator \(b\).

\paragraph{An octic point comparison.}
For \(d=8,\ n=16,\ N=17\), suppose that
\(\delta_C(X,H)\ge136/13\) for each line and each integral conic.
Then the proof of the dimension-\(2d\) criterion applies with this
mark. In the tilt use \(f=1/8,\ a=54/125,\ \epsilon=10^{-8}\).
The downward-rounded signed surface bound is \(B=17/5\).
Above \(B\) write \(t=k/j\) for the order-to-degree ratio of the
unique exceptional prime. At \(t=8\) the single-value principal bin
with \(K=70\) applies in every section degree \(j\). For \(B<t<8\)
use the discriminant test of Section~\ref{sec:nov-second-row}. Put \(H=12\) and
\(E=8+40/H+1/H^2\) for \(j\ge H\). For \(j<H\) use the genus
upper bound \(I=8j(j+4)+k+1\), taking its minimum with the
sectional Milnor bound when \(j\le2\). The discrete test omits
\(k(k-1)>I\), since distinct Puiseux roots have pairwise
contacts at least one. All pairs then pass the surface comparison
except possibly \(j=1\), \(k=5,6,7\). The respective
single-value principal-bin checks have \(K=29,40,51\).
The exact arithmetic is in \path{nov27-bin-octic.py}.
This concludes the point and Nadel comparisons under the stated
short-curve hypothesis.

\subsection{The septic in dimension fourteen}\label{sec:nov27-septic-lines}
Let \(X=X_7^{14}\) be smooth and set
\[
 T=\frac{21}{10},\qquad b=1+\frac67\,\frac{T}{1+T}=\frac{49}{31}.
\]
We prove \(\delta_C(X,H)\ge 15/b=465/49\) for every line and
every integral conic \(C\). For the lines outside the polar-pencil
case and for the conics, use the short-curve comparisons of
Section~\ref{sec:nov27-integer-second}. The four separate bounds are
\(14/9,10/7,3/2,85/56\), all smaller than \(b\). The finite ternary
comparisons are reproduced in \path{nov27-q8-septic.py}.

For a line whose normal polars span a pencil, write \(c\) for the
ambient codimension of the current constant light subspace and \(G_h\)
for its first nonzero light form, as in
Section~\ref{sec:nov27-all-pencils-2d}. Initially \(c=2\). Whenever
\(c<7\) the next light order exists. At terminal order \(h=7\),
the Gauss-map estimate bounds the codimension of the singular locus
in projective light space below by \(16-2c\ge4\).
Four general lights therefore cut out a smooth surface in
\(\mathbf P^3\). Their radial ratio is \(4/7>1/T\), and the
noncritical signed graph estimate applies.

At \(h=5,6\), in the absence of a component of multiplicity at least
three, the three-light test suffices if
\((\mathbf P_{\rm light},G_h/T)\) is log canonical, since \(h\le3T\).
In the remaining case the bidegree Nadel argument gives a constant
codimension-two linear center, using \(h<3T\) and \(7-h<T\).
If a component of multiplicity at least three occurs, its degree is at
most two. A hyperplane component is constant by the primitive
base-degree bound applied to its Galois orbit, since \(7-h\le2\).
An irreducible quadric component forces \(h=6\) and contributes the
whole form as a cube. Its rank is at least three, so three general
lights cut a smooth conic. The pure-quadric calculation applies:
\(3/6\ge1/T\), \(6/(6+3)\le T/(1+T)\), and \(2\le T\).
Thus each outstanding high order admits a cut of codimension at most
two raising its order. The potential \(c+2(7-h)\) never increases
under these cuts. If it is smaller than seven, all cuts stay in \(c<7\)
through termination. This handles every initial \(h\ge5\).

A single light suffices for \(h=2\). For \(h=3,4\) two lights
suffice when all component multiplicities are at most two. Otherwise
there is a unique hyperplane of multiplicity \(s\ge3\). Its primitive
base degree is at most one since \(se\le7-h\). Impose one
preliminary constant linear condition if needed to make the
hyperplane constant, retaining the nonzero form \(G_h\). On cutting
its kernel the next order \(t\) exists, since \(c\le4\). Use the
terminal test if \(t=7\). If \(h=3\) and \(t\le6\), then
\(1/h+1/t\ge1/2>1/T\), so the crossing-edge estimate of
Section~\ref{sec:nov27-integer-second} applies. Its torus factors have
\(me\le t\le6\le3T\). For \(h=4,t=6\), the potential is at most
\(4+2(7-6)<7\). It remains to treat \(h=4,t=5\).

Write \(K\) for the constant kernel of the repeated quartic
hyperplane when it is fixed. In the moving case one has
\[
 G_4=(uA+vB)^3R,\qquad K=(A=B=0),
\]
where \(A,B\) are independent constant linear forms and \(R\) is a
constant linear form, since the cube uses the entire base degree.
The quintic \(G_5|_K\) is nonzero. When the quartic hyperplane is
fixed and \(G_5|_K\) is a fifth power at the geometric generic
base point up to scalar, its root is a constant hyperplane by the
primitive bound \(7-5<5\). Cutting this hyperplane gives \(c=4\)
and an order at least six, where the potential test applies.

After the preliminary cut retain one complementary coordinate \(x\)
and two general lights \(y_1,y_2\) on \(K\). The transverse
\(x\)-order of the quartic is \(j=3\) or \(4\). The complement and
the binary section may be chosen generally to preserve this order.
Give the three lights weights \((j+1,j,j)\). The face of weight
\(5j\) is
\[
              f=x^jr(y_1,y_2)+g(y_1,y_2).
\]
Here \(r\) is nonzero and linear when \(j=3\), and is a nonzero
constant when \(j=4\). The form \(g\) is the binary restriction of
\(G_5|_K\). The radial ratio \((3j+1)/(5j)\) exceeds \(1/T\).
The following computations take place at a general pole root,
where the generic nonvanishing and transversality conditions hold.

Weighted Euler gives tangent contact two at a noncritical graph
point. At \(y_1=y_2=0,\ x\ne0\), either \(f\ne0\) when \(j=4\),
or the derivative in the \(r\)-direction is nonzero when \(j=3\).
For \((y_1,y_2)\ne0\), fix a nonzero \(y\)-coordinate to cut out
a plane transverse to the weighted orbit. In a binary projective
coordinate \(w\), the resulting germ has equation
\[
          x^j w^\alpha U(w)+w^\beta V(w),
          \qquad U(0)V(0)\ne0,\quad 0\le\alpha\le1,\quad
          0\le\beta\le5 .
\]
If \(\beta<\alpha\) this is a coordinate power times a unit. For
\(\beta=\alpha\) any zero of the residual factor is simple in \(x\)
and transverse to \(w=0\). If \(\beta>\alpha\) and \(x\ne0\), again
only a coordinate factor can vanish. At \(x=0,\ \beta>\alpha\) one
absorbs the units to obtain \(w^\alpha(x^j+w^{\beta-\alpha})\).
Its threshold by the compact Newton edge is the minimum of
\(1,1/\alpha,(j+\beta-\alpha)/(j\beta)\), omitting a zero
denominator. This is at least \(1/2\) except possibly for
\(j=3,\alpha=1,\beta=5\) and \(j=4,\beta=5\).
A general binary section has a root of multiplicity five only if
\(G_5|_K\) is a hyperplane to the fifth power up to scalar:
a general projective line meets the reduced components transversely
and avoids their intersections. For \(j=3\), sharing that root
with the restriction of the residual linear form on a general line
requires the residual hyperplane itself to be the root. The preceding
codimension reduction covers fifth powers for a fixed quartic hyperplane.

For a moving quartic hyperplane with \(R|_K=0\), retain two
complementary \(A,B\) directions and one general \(K\)-light \(y\).
The weights \((5,5,4)\) have face \(H_4(A,B)+\gamma y^5\),
\(\gamma\ne0\), and radial ratio \(14/20>1/T\).
The binary quartic is the cube of a general pencil member times
a fixed linear form. The two roots are distinct at a general pole.
At a critical point off the origin one has \(y=0\) and meets a
root of \(H_4\) of multiplicity \(e\le3\). The local transverse
threshold is \(\min(1,1/e+1/5)\ge8/15>1/T\).

It remains to treat the moving case with \(r=R|_K\ne0\) and with
\(G_5|_K\) proportional to \(r^5\). Make the preliminary cut so that
\(uA+vB\) becomes a nonzero base multiple of \(x\). Retain \(R\)
itself as another coordinate, denoted again by \(r\), and take a
general light \(y\) in \(K\cap(r=0)\). Its first pure order is six
or seven because this full kernel has ambient codimension five.
Initially assign \((x,r,y)\) weights \((4,3,3)\). The minimum face
consists precisely of \(x^3r\) and \(r^5\), with nonzero coefficients.
Decrease the weight \(a\) of \(y\) to the first wall; the pure power
of \(y\) ensures \(a\ge15/7\). Here are all possible new monomials
at weight 15.
\[
\begin{array}{c|l@{\qquad}c|l}
 a&\text{new monomials}&a&\text{new monomials}\\ \hline
 11/4&xy^4&7/3&x^2y^3\\
 8/3&xry^3&9/4&r^2y^4\\
 5/2&xr^2y^2,\ y^6&11/5&xy^5\\
 12/5&ry^5&15/7&y^7
\end{array}
\]
Indeed for a new term \(x^i r^j y^k\) one has \(k>0\),
\(5\le i+j+k\le7\), \(a=(15-4i-3j)/k\in[15/7,3)\);
at ordinary degree five the only term independent of \(x\) is \(r^5\).
The radial ratio \((7+a)/15>1/T\) suffices. At \(y=0\) the punctured
zero divisor is smooth by the \(x,r\) derivatives of \(x^3r+cr^5\)
with \(c\ne0\). At \(y\ne0\) normalize \(y=1\) by weighted scaling.
The transverse plane polynomial is
\[
 F=x^3r+cr^5+M(x,r),\quad
 M\in\{Ax,Axr,Axr^2+B,Ar,Ax^2,Ar^2,B\},
\]
where \(M\ne0\).
At \(r\ne0\), multiplicity at least three requires
\(F_x=F_{xx}=0\). For \(M=Ax^2\) these equations are incompatible.
Otherwise they force \(x=0\), and any nonzero root of
\(cr^5+M(0,r)\) is simple. At \(r=0,\ x\ne0\), the cases
\(M=Ax,Ax^2,B\) give no zero. For \(M=Axr^2+B\) with \(B=0\),
and for \(M=Ar^2\), one has \(F_r=x^3\ne0\). For \(M=Axr\) or
\(Ar\), if \(F_r=0\) then \(F_{xr}=-2A\) or \(3x^2\), respectively,
and is nonzero. This proves multiplicity at most two at all such
zeros. At the origin a nonzero constant or a term of degree at most
two suffices except when \(M=Axr^2,\ A\ne0\).
In this case the Newton polygon has vertices \((3,1),(1,2),(0,5)\);
the two compact faces are nondegenerate on the torus. Their normals
\((1,2),(3,1)\) give ratios \(3/5,4/5\), and the coordinate
factor \(r\) has coefficient one. Hence the local threshold is \(3/5\).
The multiplicity bound at the other zeros gives threshold at least
\(1/2\) there.

All the transverse planes at critical points in the order-\((4,5)\)
argument have threshold at least \(1/2>1/T\) and degree at most six.
In particular the sole possible septic term \(y^7\) on the last face
becomes constant on \(y=1\). The signed plane-graph calculation of
Section~\ref{sec:nov27-integer-second} applies: the threshold controls
the radial order and coordinate factors, the torus factors satisfy
\(me\le6\le3T\), and contact two suffices at a noncritical graph
point. Weighted and ordinary adjunction give the positive pole
denominator \(1+T/(1+T)\). Averaging the regular and positive groups,
of shares \(1/7\) and \(6/7\), respectively, gives \(b=49/31\).

\paragraph{The point estimate for the septic fourteenfold.}
Let \(d=7\), \(n=14\), \(r=9\), and suppose that each line and each
integral conic \(C\) in \(X\) satisfies \(\delta_C(X,H)\ge465/49\).
In the finite-fraction argument of Section~\ref{sec:nov-second-row}
take \(f=1/7\), \(a=447/1000\) and \(\epsilon=10^{-8}\).
The tilt gives \(Q=7317/2500\) and \(B=409/125\) after outward
rounding. The two signed point comparisons at \(q=0,Q\) are strict
when \(\tau_p(H)\le B\); the surface cutoff applies since
\(B^2>7\) and \(9(B^2/7-1)+B-7>0\).

For \(\tau_p(H)=t>B\) let the exceptional prime be \(D\in|jH|\)
of multiplicity \(k=jt\). At \(t=7\), use the principal-bin comparison
directly with \(50\) bins. When \(B<t<7\) the general pointed curve
section of the prime is integral by projection from \(p\); a cone
with vertex \(p\) would have multiplicity \(7j\). For \(1\le j<12\)
enumerate \(Bj<k<7j\) and put \(I=7j(j+3)+k+1\). For \(j=1,2\)
replace \(I\) by its minimum with the upper sectional Milnor bound
of Section~\ref{sec:nov-second-row}. The local discriminant exponent
is at least \(k(k-1)\), since the \(k\) Puiseux roots have pairwise
contact at least one in a general projection. Discard the pairs with
\(I<k(k-1)\). With
\[
  m=\min(k,\sqrt I),\qquad
  L=\frac{I m}{j(I+m(m-1))},
\]
one has \(7/t+(1-7/t^2)L<B\) except possibly at \((j,k)=(1,5),(1,6)\).
The principal-bin comparison there uses \(29\) and \(40\) bins
respectively. For \(j\ge12\) use
\[
 E=7+\frac{28}{12}+\frac1{144},\qquad
 L_*=\frac{E}{2\sqrt E-1/12}.
\]
The maximum of \(7/t+(1-7/t^2)L_*\) on \([B,7]\) is strictly below
\(B\), as follows by evaluating at the nearest point of this segment
to \(2L_*\). The rational root enclosures, signed inequalities and
principal-bin comparisons are recorded in \path{nov27-bin-septic.py}.
The pointed hyperplane section has dimension thirteen and satisfies
Theorem~\ref{thm:nov27-first-row}. The positive-center and point
contradictions of Section~\ref{sec:nov-second-row} now give
\(\delta(X,H)>9\) under the stated short-curve bound.

\subsection{Signed graph comparisons below ordinary tangent contact two}
\label{sec:nov27-parabolic-flag}
We use the normalized signed numerator
\((m+1)S_{\mathcal O_{\mathbf P^m}(1)}(v)+v(W)\)
for a smooth affine graph \(W\). Suppose first that at the marked
point the Hessian of the graph over its tangent hyperplane has
rank at least two. Then its signed denominator is at most \(3/2\).

Indeed, by successive ordinary affine-linear flags reduce to
\(\mathbf P^3\) so that the restriction of the quadratic form to
the two surviving light directions is nondegenerate. Subtract the
tangent affine-linear equation in the graph coordinate. The
weighted extraction at the origin with weights \((1,1,2)\) has
discrepancy four, normalized complete-series order four, and graph
order two. The initial graph is \(z+xy\), after a linear change of
the two lights. Its ordinary projective closure \(Q\subset\mathbf P^3\)
is a smooth quadric surface. Along \(Q\) the normalized coefficient
of the complete degree series is \(1/2\), so including \(Q\) gives
exactly \(3/2\). Its restrictions are the complete systems of
\(\mathcal O_Q(t)\) for \(t=m-2j\). Their mean \(t/m\) is \(3/4\)
because the degree-\(t\) dimension has quadratic leading term.
On \(Q\simeq\mathbf P^1\times\mathbf P^1\), the polarization is
\(\mathcal O(1,1)\) and \(\delta(Q,\mathcal O(1,1))=2\); this is
also the standard complete product-series calculation. Thus the
normalized residual series \(3S_{Q,\mathcal O(1,1)}\) has
denominator at most \(3/2\). Ordinary divisor adjunction gives
\(4S_{\mathbf P^3}+v(Q)\le(3/2)A\) over \(Q\); off \(Q\) ordinary
projective semistability gives the result. The ordinary complete
sets of monomials split and reassemble by their weighted degrees
as in Section~\ref{sec:nov27-all-pencils-2d}. Consequently the
punctured angular graph of the original point extraction has this
same signed bound. This and the radial calculation prove the assertion.

A second improvement applies to graph germs in three affine coordinates.
Consider \(W=(z+x^2+y^k=0)\), where \(k\ge3\) is an integer.
Let
\[
 \lambda_k=\min\{3/2+1/k,\ 2-1/k\}.
\]
At all points of its ordinary affine chart,
\[
             4S_{\mathbf P^3}(v)+v(W)
                \le (1+1/\lambda_k)A(v).
\]
Off the graph projective semistability suffices. When \(y\ne0\)
the graph Hessian is nondegenerate, as just treated. On
\(y=0\) its graph points lie on the smooth projective plane conic
\[
             C=(Y=0,\ Q=0),\qquad Q=ZX_0+X^2 .
\]
Extract with weights one on \(Y\) and \(k\) on \(Q\). Here the
restriction argument is made over the smooth affine locus of \(C\)
where \(X_0\ne0\). The graph has order \(k\) and its marked fiber
direction has both coordinates nonzero. In the notation of the
earlier conic-flag proof, the initial pieces are
\(Y^iQ^j H^0(C,\mathcal O_C(t))\), \(t=m-i-2j\); the ordinary
complete-series mean values in units of four are \(1,1/2,2\)
for \(i,j,e\) respectively, where \(e\) is the vanishing order
of a coefficient section at the marked point. The radial comparison
for the stronger signed inequality \(4S+\lambda_k v(W)\le2A\) is
\[
                1+k/2+\lambda_k k\le 2(1+k).
\]
The ordinary fiber over the marked point has coefficient two and
the graph misses its generic point. On this weighted projective
fiber the monomials satisfy \(i+kj=J\),
\(i+2j\le m-e/2\). Put \(a=t-e/2\). At a finite nonzero point
the average order in the consecutive monomial interval in \(j\) is
bounded by half its width
\[
               \lfloor i/k\rfloor+
                     \min\{j,\lfloor a/(k-2)\rfloor\}.
\]
In the limiting simplex with coordinates \(i,2j,e/2,a\) the two
summands have normalized mean \(1/k\) each: the second tail cuts
length \(2s+(k-2)s=ks\) at level \(s\). Thus the normalized
complete-series order at this point is at most \(1/k\). The graph
gives one simple mark; the total coefficient there is
\(\lambda_k+1/k\le2\). Each coordinate point has normalized
coefficient at most one, read in the corresponding stack chart,
and is disjoint from the graph mark. Weighted and ordinary fiber
adjunction therefore prove \(4S+\lambda_kv(W)\le2A\).
Interpolating with \(4S\le A\) gives the asserted inequality.
For \(k=3,4,5,6\) its maximum denominator is \(8/5\).

In particular on the punctured affine chart of the graph
\(z+x^h+y^k=0\), \(2\le h,k\le6\), the positive-pole
denominator is at most \(8/5\). Here punctured means that the
light pair \((x,y)\) does not vanish. If both are nonzero, the
Hessian has rank two. If just one is nonzero, say \(x=x_0\ne0\),
subtract the tangent affine-linear term and translate \(x\) by
\(x_0\). In the resulting \(X=x-x_0\) the light part is
\(aX^2+O(X^3)+y^k\) with \(a\ne0\). For \(k=2\) use rank two.
For \(k\ge3\) extract with weights \(k,2,2k\) on \(X,y,z\)
after division by their gcd. The normalized radial ratio is
\(1+2k/(3k+2)\le8/5\) here, and its punctured initial graph
is the preceding \(z+X^2+y^k\) after rescaling. The other axis
follows with the exponents interchanged. If additionally
\(1/h+1/k\ge2/3\), the same \(8/5\) bound holds also at their
common origin using the weights \(k,h,hk\).

\subsection{Noncritical graphs of weighted forms of degree at most six}
\label{sec:nov27-six-weighted-noncritical}
Let \(P\) be a polynomial of ordinary degree at most six on an affine
space of light coordinates. Suppose its Hessian at \(p\) has rank one,
and let \(L\) be a linear equation of its Hessian kernel centered at
\(p\). If \(P-T_pP\) is not divisible by \(L^2\), where \(T_pP\) is
the affine tangent polynomial, the signed graph denominator at \(p\)
is at most \(13/8\).

Indeed cut to a plane with coordinates \(X,Y\) through \(p\), taking
\(X=L\). A general direction on \(L=0\) and a general complementary
direction preserve at least one nonzero pure-\(Y\) or \(X\)-linear
coefficient of \(P-T_pP\). Write \(k\) for the first pure-\(Y\) degree
and \(j\) for the first \(Y\)-degree of an \(X\)-linear term, using
infinity for an absent term. Then \(k\ge3,\ j\ge2\), and finite orders
satisfy \(k\le6,\ j\le5\). The quadratic form is \(aX^2,\ a\ne0\).
With the weight of \(Y\) one, increase the weight of \(X\) to the first
wall \(\min(k/2,j)\), and take the weight of the graph coordinate equal
to twice that value. At \(k/2<j\) the signed bound is at most \(8/5\)
by the disjoint-power calculation. Otherwise the first form, after
rescaling, is
\[
           F=X^2+BXY^j+CY^{2j},\qquad B\ne0,\quad 2\le j\le5 .
\]
Here \(C=0\) if \(2j>6\). The radial denominator is
\(1+2j/(3j+1)\le13/8\).

At a punctured point with \(Y=0,\ X\ne0\), the Hessian has rank two
if \(j=2\); for \(j\ge3\), after subtracting the affine tangent its
first pure term on the Hessian kernel is \(BXY^j\) with the original
\(X\) evaluated at the point. The weights \((j,2,2j)\) leave a
quadratic-plus-\(j\)-th-power graph and have radial denominator at
most \(8/5\). For \(Y\ne0\), rank is at least one by \(F_{XX}=2\).
Writing \(X=qY^j\), its rank is one exactly if
\[
         2B(j-1)q=B^2j-4(2j-1)C.
\]
At such a point \(e=(-BjY^{j-1}/2,1)\) spans the kernel and
\[
       D_e^3F=\frac{j^2(2j-1)}{2}(4C-B^2)Y^{2j-3}.
\]
When the cubic derivative is nonzero the quadratic-plus-cubic flag
applies. In the remaining case \(F=(X+BY^j/2)^2\) and rank one
forces \(X+BY^j/2=0\). At \(Y\ne0\) its reduced equation is smooth
and has nonzero second derivative in the curve tangent direction.
In affine linear tangent coordinates \(U,V\), the weights \((2,1)\)
therefore give its initial equation \(U+bV^2,\ b\ne0\).
The corresponding weighted graph has first light part
\((U+bV^2)^2\). At a point of the parabola use the conic flag of
Section~\ref{sec:nov27-plane-curve-flag}; at its complement the
Hessian is nondegenerate, by the degree-four determinant computation.
Its radial denominator is \(1+4/7\). The Hessian estimate and the
quadratic-plus-power estimates used here are in
Section~\ref{sec:nov27-parabolic-flag}. Weighted adjunction and
ordinary linear restriction prove the claim.

Suppose now that \(P\) is weighted homogeneous of weight \(H\),
strictly greater than each light weight, with at least two independent
light directions. Unless \(P\) is a scalar times a power of a linear
form, at \emph{every} noncritical point its graph has signed
denominator at most \(13/8\). Weighted Euler differentiation shows
the Hessian has positive rank at a noncritical point. In rank one
it remains to exclude a hyperplane \(A\) on which all first
derivatives of \(P\) are constant, which would result from the
divisibility in the preceding calculation. Weighted scaling by
\(t\in\mathbf C^*\) carries \(A\) to another hyperplane with constant
gradient there, transforming the \(i\)-th gradient value by
\(t^{H-w_i}\). If \(A\) is fixed these constants vanish, since
\(H>w_i\). If generic translates are nonparallel to \(A\), each
meets \(A\); at an intersection their gradient constants agree,
so again all constants vanish. Both alternatives contradict
noncriticality. Otherwise the translates are parallel and their
common linear normal is a weight eigenvector. If they vary, their
union is dense in the affine space. In coordinates \((u,v)\) with
these hyperplanes \(u=\mathrm{const.}\), constancy of all derivatives
along each such hyperplane implies
\(P=B(u)+\sum_i a_i(u)v_i\). Constancy of \(\partial P/\partial u\)
then forces the \(a_i\) constant. Weighted homogeneity and \(H>w_i\)
remove those tangent linear terms and the affine terms, and \(B(u)\)
is a single power since \(u\) can be chosen as a linear weight
eigenvector through the origin. This is exactly the excluded case.

\subsection{Quadratic and cubic initial forms in a sextic polar pencil}
\label{sec:nov27-sextic-quadratic}
Use the signed graph inequality and the two estimates for a quadratic
tangent form and for \(z+x^2+y^k\) in
Section~\ref{sec:nov27-parabolic-flag}. All the following graphs are
on the ordinary affine chart, and all changes before a weighted
extraction are affine linear. The signed denominator is the least
constant \(D\) in \((m+1)S_{\mathcal O_{\mathbf P^m}(1)}(v)+v(W)
\le D A(v)\).

\paragraph{A noncritical homogeneous form.}
Let \(P\) be a homogeneous polynomial of degree \(2\le h\le6\)
which is not proportional to a power of one linear form. At any
point \(p\) of its light space with \(dP_p\ne0\), the graph of \(P\)
has signed denominator at most \(8/5\). A Hessian of rank at least
two gives \(3/2\). Otherwise the Euler identity
\(\operatorname{Hess}_p(P)(p,-)=(h-1)dP_p\) forces rank one;
its kernel is \(\ker dP_p\), and \(P(p)\ne0\). Indeed for a rank-one
symmetric form \(b\), \(b(p,-)\ne0\) implies \(b(p,p)\ne0\).
Choose an affine plane by the two light directions \(p\) and a
general \(v\in\ker dP_p\). Its polynomial is
\[
  P((1+x)p+yv)=\sum_{i=0}^h c_i(1+x)^{h-i}y^i .
\]
One has \(c_0\ne0,\ c_1=c_2=0\). Some \(c_k\ne0\) for \(k\ge3\):
otherwise the polynomial on \(\langle p\rangle\oplus\ker dP_p\)
would be a pure \(h\)-th power, by varying \(v\). Subtract the
affine-linear tangent value in the graph coordinate and take the
weights \(k,2,2k\) on \(x,y,z\), where \(k\) is the first such index.
The initial graph is \(z+a x^2+c_k y^k\), \(a\ne0\). The radial
denominator \(1+2k/(3k+2)\) is at most \(8/5\), and the punctured
graph has the same bound by Section~\ref{sec:nov27-parabolic-flag}.
The same conclusion holds in the original dimension by the general
linear flag through the indicated plane (or successive restriction
to fixed affine-linear hyperplanes with this intersection).

\paragraph{A weighted quartic.}
For constants \(B,C\) and \(B\ne0\), put
\[
                 P=x^2+Bxy^2+Cy^4 .
\]
At every affine graph point the signed denominator is at most \(8/5\).
On the punctured chart its Hessian determinant is
\(4Bx+(24C-4B^2)y^2\). At its zeros \(y\ne0\), since the origin is
omitted. The Hessian has rank one, with kernel vector \(e=(-By,1)\);
the third derivative in that direction is
\((24C-6B^2)y\). When this is nonzero, subtract the tangent value,
write the quadratic part as \(aX^2\), and extract with weights
\((3,2,6)\) on \(X,Y,z\), taking \(Y\) in the kernel direction.
The initial graph has light part \(aX^2+bY^3\) with \(ab\ne0\).
Its radial denominator \(1+6/11\) and its angular denominator
are at most \(8/5\). In the remaining case \(C=B^2/4\), the
rank-one locus is \(Q=x+By^2/2=0\) and \(P=Q^2\). The ordinary
projective closure of \(Q=0\) in the light plane is a smooth conic;
the conic flag of Section~\ref{sec:nov27-plane-curve-flag} with
\(e=s=2\) gives \(1+4/7\) there. At the origin use weights
\((2,1,4)\); they preserve \(P\), their radial denominator is
\(1+4/7\), and the punctured chart was just checked.

\paragraph{Applications to a sextic line pencil.}
Put \(d=6,\ n=12\) and keep the line pencil notation of
Section~\ref{sec:nov27-all-pencils-2d}. For \(T\ge3/2\) the desired
positive-pole denominator \(1+T/(1+T)\) is at least \(8/5\).
At terminal light order six, four general lights cut the smooth
projective surface in \(\mathbf P^3\) whenever the ambient codimension
of the constant light kernel is at most five. This follows from
the lower bound \(14-2c\ge4\) for the singular locus codimension.
The homogeneous sextic is not a pure power, and smoothness of its
projective zero set and Euler's identity give nonzero gradient at
every punctured point. The first calculation applies. The radial
ratio is \(4/6\ge1/T\).

At initial light order two, a quadratic of rank at least two remains
nondegenerate on two general light directions. Its Hessian gives
the required graph bound everywhere. Otherwise the quadratic is a
square up to scalar, with rational root since it is unique. The
primitive base degree of the root is at most two, since the residual
base degree is four. Make the root constant on a subspace, retaining
the quadratic, by at most two preliminary linear cuts as in
Section~\ref{sec:nov27-divisible-second}; then its kernel has
ambient codimension at most five. Let \(t\ge3\) be its first
nonzero light order. If \(t=6\), the preceding terminal test applies.
If \(t\le5\), keep one general direction \(y\) in the kernel and a
complementary direction \(x\). Increase the weight of \(x\) from one,
with the weight of \(y\) fixed, until \(x^2\) first meets another
monomial. The pure \(y^t\) ensures the weight of \(x\) at this wall
is at most \(t/2\). A mixed term on the wall is \(xy^j,\ j\ge2\);
then \(2j\le t\), forcing \(j=2\). The first face therefore is
\(x^2+C y^t\) with \(C\ne0\), or the preceding weighted quartic
with \(B\ne0\), allowing \(C=0\) in the latter. The radial ratio
is at least \(1/2+1/t\ge7/10\), and the two signed estimates
give \(8/5\).

Suppose next \(T\ge5/3\), and the initial cubic at the geometric
generic base point has a double hyperplane and a distinct residual
hyperplane. The double root is rational and its primitive base
degree is at most one. After making it constant, write \(x=0\) for
its equation. On its kernel the residual linear equation is nonzero
at a general pole unless the preliminary cut made the restricted
cubic a pure cube. In the following calculation assume that the two
roots stay distinct there. The case without a preliminary cut
always satisfies this assumption. Let the next order on \(x=0\)
be \(t\), so \(4\le t\le6\). At \(t=6\) use the terminal flag.
At \(t=4,5\), choose a single general kernel light \(y\) and any
complement \(x\), in which the cubic is \(x^2(ay+bx)\), with
\(a\ne0\). Give \(x,y\) the weights \(t-1,2\) and the graph
coordinate weight \(2t\). The face has light part, after rescaling,
\[
          F=x^2y+C y^4\quad(t=4),\qquad
          F=x^2y+Bxy^3+Cy^5\quad(t=5),\qquad C\ne0 .
\]
Indeed any mixed monomial \(xy^j\) of higher ordinary order has
\(j\ge3\). The pure \(x^3\) is strictly above the displayed face.
The radial ratio \((t+1)/(2t)\ge3/5\) suffices.

At \(y=0,\ x\ne0\) the Hessian determinant is \(-4x^2\).
At \(t=4\) the rank-one locus with \(y\ne0\) has \(x^2=6Cy^3\).
The third derivative on the kernel vector \((-x/y,1)\) is
\(6x^2/y^2+24Cy=60Cy\ne0\); thus the quadratic-plus-cubic flag
applies. For \(t=5\) write \(x=q y^2\) at a point with \(y\ne0\).
The rank-one condition and the third derivative in its kernel
direction \(e=(-(q+3B/2)y,1)\) are respectively
\[
       4q^2=40C-9B^2,\qquad
       D_e^3F=6q(2q+B)y^2 .
\]
If \(q=-B/2\), then \(C=B^2/4\ne0\) and \(F=y(x+By^2/2)^2\).
Apply the smooth conic flag with unit \(y\). If \(q=0\), then
\(C=9B^2/40\) forces \(B\ne0\); here
\(D_e^4F=-9B^2y\ne0\). After subtracting the tangent equation,
extract with weights \((2,1,4)\) in coordinates where the Hessian
is \(aX^2\) and the kernel direction is \(Y\). The face has light
part \(aX^2+bXY^2+cY^4\) with \(ac\ne0\), and was treated above,
also when \(b=0\) by the disjoint-power estimate. When
\(q(2q+B)\ne0\) use the quadratic-plus-cubic estimate directly.
This proves the desired bound in all the stated cases.

\paragraph{All cubic initial forms.}
Suppose \(d=6,n=12\) and \(T\ge12/7\).
The positive graph target is then at least \(31/19>13/8\).
For a reduced cubic keep two general lights. They give a homogeneous
binary cubic with simple roots, whose punctured graph is noncritical
and satisfies Section~\ref{sec:nov27-six-weighted-noncritical}.
The radial ratio is \(2/3\). If the cubic has a repeated factor, its
unique repeated hyperplane is rational over the function field of the
base and the primitive degree of its equation is at most one, by its
multiplicity and the complementary base degree three. Use at most one
preliminary linear cut to make it constant, retaining a nonzero cubic.
The case of a still distinct simple hyperplane factor was covered in
Section~\ref{sec:nov27-sextic-quadratic}. Otherwise its equation is
\(x^3\) up to a scalar, and the constant kernel \(K=(x=0)\) has
ambient codimension at most four. Write \(t\) for the first order on
\(K\). At \(t=6\) use the smooth terminal test. If \(t=4\), use one
general \(K\)-light \(y\); the weights \(4,3\) leave the disjoint powers
\(x^3+c y^4\), with radial ratio \(1/3+1/4=7/12\) and the punctured
signed estimate of Section~\ref{sec:nov27-parabolic-flag}.

Let \(t=5\), and denote by \(B_3\) the coefficient of the \(x\)-linear
part of the quartic form on the fixed light section. It is a cubic on
\(K\) with base degree two. If \(B_3=0\), keep two general \(K\)-lights
\(y_1,y_2\) and use weights \(5,3,3\). The resulting face is
\(x^3+g_5(y_1,y_2)\); the normalized radial ratio \(11/15\) suffices.
The only punctured critical points have \(x=0\) at a multiple root of
\(g_5\). In the affine transverse plane obtained by setting a nonzero
\(y\)-coordinate equal to its marked value, the binomial light part
after the next extraction is \(x^3+a w^\beta\), \(a\ne0\), where
\(\beta\) is the root multiplicity. Its radial ratio is at least
\(1/3+1/4\) for \(\beta\le4\), and its punctured graph has the
disjoint-power signed bound. If a general binary restriction contains
a root of multiplicity five, the full quintic on \(K\) is a fifth
power of a hyperplane up to scalar, since a general projective line
meets all reduced components transversely, away from their pairwise
intersections. The complementary base degree one forces that hyperplane
to be constant. Its kernel has ambient codimension at most five, and
there the terminal test applies.

If \(B_3\ne0\), keep again two general \(K\)-lights and give \(x,y_1,y_2\)
weights \(3,2,2\). The weight-nine face is \(x^3+xB_3(y_1,y_2)\)
after a scalar normalization. The radial ratio \(7/9\) suffices.
The punctured critical locus has \(x=0,\ B_3=0\). Indeed Euler at a
critical point forces the face to vanish, and the equations
\(x^3+xB_3=3x^2+B_3=0\) give the assertion. On a transverse
plane at a simple root of \(B_3\) the Hessian has rank two. At a
double root the ordinary initial light part on this plane is
\(x^3+b xw^2\), \(b\ne0\), a reduced binary cubic already treated.
At the noncritical angular points in these calculations apply
Section~\ref{sec:nov27-six-weighted-noncritical}. None of the stated
first forms is a pure linear power.

A triple root on a general binary section forces \(B_3\) to be the
cube of a hyperplane on \(K\) up to scalar. Its base degree is at
most two, so the hyperplane is constant. Write \(L\) for its equation.
If the quintic vanishes on \(K\cap(L=0)\), this kernel has ambient
codimension at most five and immediately gives the terminal test.
Otherwise retain \(x,L\) and a general light \(u\) in that intersection
with nonzero pure quintic coefficient. The weights \(15,10,9\)
give the face
\[
                  x^3+a xL^3+b u^5,\qquad ab\ne0 .
\]
Indeed the only \(x\)-linear quartic terms come from \(xL^3\);
quartic terms with two \(x\)'s have weight at least \(48\);
and \(u^5\) is the only possible term of ordinary degree at least
five of weight at most \(45\). Its radial ratio \(34/45\) suffices.
The derivatives force \(u=0\), then \(xL^2=0\), then \(x=L=0\) at
any critical point. The punctured graph bound therefore follows
from Section~\ref{sec:nov27-six-weighted-noncritical}.

\subsection{Quartic and quintic initial forms in a sextic line pencil}
\label{sec:nov27-sextic-high}
Let \(d=6,\ n=12,\ T=12/7\). Write \(c\) for the ambient
codimension already used by the constant light subspace. We use the
positive-pole signed denominator \(1+T/(1+T)=31/19\). The terminal
test of Section~\ref{sec:nov27-sextic-quadratic} applies at \(h=6\)
and \(c\le5\). The arguments below also work for all larger \(T\).
At ordinary or weighted noncritical points we use the bound \(13/8\)
of Section~\ref{sec:nov27-six-weighted-noncritical} for the
non-pure weighted forms occurring below. The graph bounds are tested
with general pole roots after choosing the constant linear subspaces.

\paragraph{Quintic order.}
If the geometric generic quintic is reduced and the pair at
coefficient \(7/12\) is log canonical, use three lights and
\(3/5\ge7/12\) in the ordinary three-light argument. If it fails log
canonicity, the bidegree Nadel argument gives a constant linear locus
of codimension two: indeed \(5/T<3\) and the base degree \(1<T\).
Cutting this subspace gives the terminal order and costs two. In
particular every reduced quintic at \(c\le3\) is treated.

For a nonreduced quintic, each repeated geometric factor is constant.
Indeed the product over its Galois orbit has a primitive equation of
base degree \(e\) with \(2e\le1\), hence factors over the constant
algebraically closed field. If it is a hyperplane cut it; this also
treats \(c=4\). Otherwise the quintic is \(q^2R\), with \(q\) an
irreducible constant quadric and \(R\) a base-linear family of
hyperplanes up to scalar. If \(R\) is fixed, cut it. Assume it
moves. At a general base point \(q|_{R=0}\) has rank at least two.
This is automatic for \(\operatorname{rank}q\ge4\). In rank three,
a hyperplane not containing the radical retains rank three.
The exceptional hyperplanes that contain the radical and have
restricted rank one correspond to the dual conic of the nonsingular
ternary quadric. This conic contains no projective line, whereas
the image of the nonconstant base-linear family is a projective
line. Thus three general lights give a smooth plane conic whose
intersection with the selected residual line is transverse. For
example, first take two general directions in the residual
hyperplane on which the quadratic restricts nondegenerately; a
general third direction gives a nondegenerate ternary restriction.

The homogeneous radial ratio is \(3/5\). On the punctured chart
at noncritical points use the noncritical bound. At a zero of \(q\)
where \(R\) is a unit, use an ordinary affine plane chart of the
projective light space, transverse to the radial direction, and
the smooth conic flag with \(s=e=2\). At a common zero of \(q,R\),
choose dehomogenized affine-linear coordinates \(x,y\), with \(x=R\)
and \(y=0\) tangent to \(q\). The weights \((1,2)\) leave the conic
equation \(y-a x^2,\ a\ne0\), after rescaling. The first graph has
light part \(x(y-a x^2)^2\); its radial ratio is \(3/5\).
On its punctured chart the only critical points are on the smooth
parabola away from the axis, where the remaining factor is a unit.
The conic-flag denominator is \(1+4/7\). All other points are
noncritical, so the noncritical estimate and weighted adjunction
finish this case. This proves the sufficient quintic assertions
also when the base coefficient of \(R\) is specialized at any of
the finitely many general pole roots.

\paragraph{A reduced quartic.}
When the geometric generic quartic \(G\) is log canonical at
coefficient \(7/12\), use three lights and \(3/4\ge7/12\).
Otherwise the Nadel plane-slice argument gives a unique
codimension-two linear non-lc center \(W\) in projective light
space over the geometric generic point. It descends to the base
function field by uniqueness. We initially work in one geometric
fiber; dependence on the base is addressed below.

Take a general two-plane meeting \(W\) transversely at a general
point. Its quartic is reduced by Bertini, has the same multiplicity
\(m\) there, and fails the coefficient \(7/12\) at that point by
generic restriction of the multiplier ideal (choose a coefficient
slightly below \(7/12\) first). Necessarily \(m\ge3\). For \(m=2\),
under a general projection let \(s\) be the contact of its two
distinct Puiseux roots. The local discriminant exponent \(2s\) is
at most \(4\cdot3\) by intersecting with a general cubic polar.
For a valuation normalized by \(v(x)=1\) on the projection
coordinate, write \(Y\) for its larger order of these two branches
on a split Puiseux cover. Finite ramification and the two regular
parameters at either split branch give \(A(v)\ge1+Y\); their
sum of orders is at most \(Y+\min(Y,s)\). The resulting normalized
quotient is at most \(\max(1,2s/(1+s))\le12/7\). Smooth points
satisfy the same assertion.

If \(m=4\), the quartic depends on just the two normals to \(W\).
Over the full base \(\mathbf P^1\) let \(E\) be the rank-two bundle
of the linear equations of these subspaces, obtained by extending
their Grassmannian map. It is a subbundle of a trivial bundle and
splits as \(\mathcal O(-a)\oplus\mathcal O(-b)\), \(0\le a\le b\).
The bidegree equation is a section of \(\operatorname{Sym}^4 E(2)\):
membership holds generically and the quotient of the corresponding
symmetric bundles is locally free. If \(b\ge1\), the coefficients
with exponent at least three on the second generator vanish, since
their line bundles have negative degree. Every term would be
divisible by the square of the first generator over the function
field, contrary to reducedness. Hence \(a=b=0\), giving two fixed
normals \(x_1,x_2\). Their common kernel now has \(c=4\). At the
next order six use the terminal test; at the next order five retain
one general kernel light \(y\). The weights \((5,5,4)\) have face
\(H_4(x_1,x_2)+\gamma y^5,\ \gamma\ne0\). The binary quartic has
four distinct roots. The face has no punctured critical point, and
its radial ratio \(14/20\) suffices.

Let \(m=3\). In transverse affine coordinates the quartic has
equation \(H_3(x,y)+Q_4(x,y)\) with homogeneous summands. Three
distinct tangent directions give threshold \(2/3\) by one blowup.
For a double and a simple direction write \(H_3=x^2y\). At the
simple direction the strict transform is smooth and transverse.
At the double direction put \(y=u,x=uv\); the total equation is
\(u^3(v^2+uB(v))\). If \(B(0)\ne0\), the compact edge in these
coordinates has normal \((2,1)\). The discrepancy relative to the
original plane is \(2\cdot2+1=5\), the order is eight, and the
nondegenerate edge together with the exceptional divisor and the
strict transforms gives threshold at least \(5/8\). If \(B(0)=0\),
reducedness forces \(\operatorname{ord}_0 B=1\), since otherwise
\(x^2\) divides the original equation. The germ on the chart is
\(u^3v(v+u\,\mathrm{unit})\), with three distinct tangent
directions; blowing up their intersection gives ratio \(3/5\),
and resolves the crossings. All these thresholds exceed \(7/12\).

Consequently at the general marked point the tangent cubic is a
triple line, say \(x^3\). A nonzero \(y^4\) coefficient in \(Q_4\)
gives the nondegenerate compact edge \((3,0),(0,4)\) and threshold
\(7/12\). Hence this coefficient vanishes, and reducedness gives a
nonzero \(xy^3\) coefficient. To express this fiberwise on the whole
light space choose homogeneous coordinates \((U,V;Y)\) where
\(U,V\) define \(W\). The terms of transverse order three have
coefficients linear in \(Y\). As \(Y\) varies, these binary cubic
forms all lie in the affine cone of cubes, since this holds at a
general marked point. The twisted cubic contains no projective
line, hence their span is one-dimensional. Thus
\[
     G=L(U,V)^3R(Y)+Q_4(U,V)
\]
for a linear form \(R\). The missing fourth coefficient along
\(L=0\) gives \(L\mid Q_4\); it is unaffected by choosing oblique
normal directions because the change of homogeneous complement
adds a multiple of \(L^3\) to \(Q_4\). We obtain \(G=L C_3\)
where the residual cubic has a transverse cusp at a general point
of \(W\): after a choice \(x=L\), its transverse plane equation
has initial compact face \(a x^2+b y^3\), \(ab\ne0\).
The cubic is geometrically irreducible, since if it were reducible
a factor of degree at most two would have to carry that reduced
irreducible cusp on a general transverse plane. Thus \(L\) is the
unique linear component and is rational over the base function field.

In fact \(L\) is constant. Otherwise its primitive base degree is
at least one, and the complementary family \(C_3\) has base degree
at most one. If constant, its singular locus has codimension at
least two and has finitely many codimension-two components, so
\(W\) is constant. If varying, it is a pencil of irreducible
cubics without a fixed divisor. Bertini gives smoothness of its
general member outside its fixed codimension-two base locus.
Again the possible linear components \(W\) of that dimension
form a finite set, so \(W\) is fixed. In fixed normals at a
general constant marked point the nonzero transverse tangent cubic
of \(G\) is a family of pure cubes of base degree at most two.
Its primitive root has base degree \(e\) with \(3e\le2\), showing
\(L\) is constant after all. Its kernel costs just one cut from
\(c=2\). The resulting quintic at \(c=3\) or the terminal sextic
was already treated.

\paragraph{A nonreduced quartic.}
A constant hyperplane factor gives the same one-cut argument.
A repeated hyperplane of multiplicity at least three is unique and
constant, by the base degree bound. A geometric irreducible
repeated quadric gives its square; three general lights cut a smooth
conic. The homogeneous radial ratio is \(3/4\), the noncritical
bound applies on its complement, and the conic flag with a unit
applies at its zeroes on a dehomogenized transverse plane.

A nonconstant rational repeated hyperplane therefore has multiplicity
exactly two, and
\[
       G=(uA+vB)^2 Q ,
\]
where \(A,B\) are fixed independent linear forms and \(Q\) is a
constant quadric. Put \(K=(A=B=0)\), of ambient codimension four.
Keep one general complementary direction \(x\) and the subspace
\(K\); thus the restricted moving linear factor is a scalar multiple
of \(x\). The quartic has transverse \(x\)-order \(m\) for some
\(2\le m\le4\). Its leading term is \(x^m q(Y)\) with \(q\)
a constant nonzero form of degree \(4-m\) up to a base scalar.
If the first order on \(K\) is six, the terminal test applies.
For order five, put \(g=G_5|_K\). If \(g\) is nonreduced, the quintic
argument above works at \(c=4\); assume \(g\) reduced. Two general
\(K\)-lights \(y_1,y_2\) give a squarefree binary restriction of \(g\)
at a general pole, and preserve the general restriction of \(q\).
Weights \(m+1,m,m\) leave precisely the face
\[
      F=x^m q(y_1,y_2)+g(y_1,y_2).
\]
Its normalized radial ratio \((3m+1)/(5m)\) is at least \(13/20\).
At \(y\ne0\), a critical point cannot have \(x=0\) by squarefreeness.
For \(x\ne0\) criticality gives \(q=g=0\) by Euler and \(F_x\),
so \(dq\ne0\) by squarefreeness of \(g\) and \(d_yF=0\). Then
\(F_{xx}=0\) and \(F_{xy}=m x^{m-1}dq\ne0\), giving Hessian rank at
least two. At \(y=0,x\ne0\) we either have nonzero gradient
(\(m=3,4\)), or \(m=2\) with light Hessian of the same rank as \(q\).
If the binary restriction of \(q\) has rank two this suffices. In
the remaining case the original \(q\) has rank one and equals a
scalar times a fixed \(\ell^2\). When \(\ell\) divides \(g\) on \(K\),
one further fixed cut gives the terminal form at \(c=5\). Otherwise
on a general binary plane \(g|_{\ell=0}\) is a nonzero fifth power
of its single remaining coordinate \(Y\) up to scalar. At \(x\ne0\)
fix \(x\) to its marked value and extract in \(\ell,Y\) with
weights \(5,2\); the face is the disjoint sum of a square and
a fifth power. Its radial ratio is \(7/10\) and its punctured graph
has the bound in Section~\ref{sec:nov27-parabolic-flag}. This
finishes the rational moving case.

Finally suppose the repeated hyperplanes are conjugate. Their
product has primitive base degree at most one and its square is
the quartic, up to a base scalar. In the nonconstant case that
product defines a pencil of reducible quadrics with no fixed factor:
a fixed factor would be defined over the constants and could not
belong to a nontrivial conjugate pair. The normalization of the
parameter curve of the individual hyperplanes is irreducible.
For two general members of this curve with distinct base points,
their hyperplanes meet in a codimension-two linear component of the
pencil base. Only finitely many components are possible; since the
space of general pairs is irreducible, this intersection is a fixed
subspace \(K\). Every moving hyperplane contains \(K\). Its ambient
codimension is four, and retaining one general complement \(x\)
therefore restricts the quartic to a scalar times \(x^4\).
The preceding order-five or order-six proof applies with \(m=4\).
This exhausts the quartic cases.

In particular, in combination with the quadratic, cubic and terminal
tests, a rank-two normal polar pencil along a line on any smooth
sextic twelvefold satisfies the short-curve mark
\((n+1)/(1+(5/6)T/(1+T))=247/29\).

\subsection{Sextic twelvefold lines of higher polar rank}
\label{sec:nov27-sextic-line-nets}
Set \(n=12,\ d=6,\ T=12/7\). For every line in a smooth sextic
twelvefold, the line denominator is at most \(29/19\).
Here a positive angular graph bound of \(31/19\) suffices. We retain
the complete polar kernel and two heavy directions giving a
basepoint-free pencil, an actual constant linear flag smooth along
the whole line. The dimension of the resulting section changes,
but the normalized denominator under this flag and the regular and
positive pole proportions are as in
Section~\ref{sec:nov27-all-pencils-2d}. For rank at least four the
ordinary syzygy estimate suffices: its positive integers have sum
five with at least three summands, giving maximum at most three
and denominator at most \(3/2\). Rank two was treated above. At rank
three the ternary estimates leave only initial orders two and three.
Here the original light kernel has ambient codimension three.

\paragraph{The last possible kernel.}
All kernels below have ambient codimension \(c\le6\). Their ambient
linear section has dimension \(13-c\ge7\), so cannot be contained
in the smooth sextic twelvefold. If its first order is six,
the terminal sextic on its lights has singular set of codimension
at least \(14-2c\ge2\) by the tangent-hyperplane argument. Thus it
is reduced; for \(c\le5\) it still admits the preceding smooth
projective three-space section.

\paragraph{The quadratic case.}
Rank two or more for the initial quadratic immediately gives the
Hessian bound. In rank one the root has base degree at most two and
can be made constant, preserving the quadratic, by at most two
cuts. Its kernel has \(c\le6\), hence admits a first nonzero term of
some order \(3\le t\le6\). On the binary section with one general
kernel light \(y\) and a complement \(x\), give \(y\) weight one
and increase the weight of \(x\) to the first Newton wall. As
before, either the first face has light part \(x^2+C y^t\),
\(C\ne0\), or it has light part
\[
       x^2+Bxy^j+C y^{2j},\qquad j=2,3,\quad B\ne0 .
\]
The radial ratio is at least \(1/2+1/6\). At noncritical punctured
points each face is covered by
Section~\ref{sec:nov27-six-weighted-noncritical}. For the binomial
there are no other punctured points. For the mixed form the only
possible punctured criticals occur when \(C=B^2/4,\ y\ne0\), along
\(x+(B/2)y^j=0\), where the squared factor is a smooth plane curve
of tangent contact two, since \(y\ne0\). In affine tangent
coordinates its weights \((2,1)\) therefore give a squared parabola,
as in Section~\ref{sec:nov27-six-weighted-noncritical}. This gives
the degree-two signed bound \(1+4/7\).

\paragraph{The cubic case.}
In the proof for a cubic initial form, all arguments made before a
second constant kernel work at present rank: its first kernel \(K\)
has \(c\le5\). At a second kernel \(M\subset K\), one has \(c\le6\).
We detail the two cases where the preceding proof called for a
smooth terminal three-space section on \(M\). Write \(x\) for the
fixed complement of \(K\) and \(L\) for the fixed complement of
\(M\) in \(K\). Keep in addition one general \(M\)-light \(y\).
The terminal sextic on \(M\) is reduced, so its value on general
\(y\) is nonzero.

First take \(t=5,\ B_3=0,\ g_5=L^5\) up to scalar. The light weights
\(5,3,5/2\) on \(x,L,y\) give exactly the face
\[
 x^3+aL^5+A x^2 y^2+Bxy^4+C y^6,\qquad aC\ne0 .
\]
Indeed every quartic term has \(x\)-order at least two, and each
quintic term either has \(x\)-order at least one or is \(aL^5\).
The least weight is fifteen. The ratio \( (5+3+5/2)/15=7/10\)
suffices. At \(y=0\) all punctured graph points are noncritical.
Set \(y=1\) at the other points by weighted homogeneity.
The monic cubic \(P(X)=X^3+AX^2+BX+C\) has no triple root for
general \(y\) at a general pole. Indeed its constant coefficient
is the terminal reduced sextic evaluated at \(y\), while \(A\)
is the value at \(y\) of a quadratic form. An identity
\(27C=A^3\) on \(M\), necessary for a triple root, contradicts
reducedness. At a critical point one has \(L=0\) and a double root.
After taking the ordinary affine section \(y=1\) and translating
\(x\), its next weighted face has light part \(X^2+b L^5,\ b\ne0\).
Its radial ratio \(1/2+1/5\) and its punctured estimate suffice.

It remains to take \(t=5,\ B_3=L^3\) up to scalar, \(g_5|_M=0\).
Give \(x,L,y\) the respective weights \(3,2,w\) and decrease \(w\)
from two, until \(x^3+a xL^3\) acquires another term of weight
nine. The possibilities are \(w=7/4,5/3,3/2\). The quartic terms
other than \(a xL^3\) have \(x\)-order at least two, the quintic on
\(K\) is divisible by \(L\), and the pure sextic coefficient is
nonzero. At the first two walls a single new term \(bL^k y^{5-k}\)
occurs, \(k=1,2\), respectively. At \(y=1,L\ne0\), a critical
zero is impossible: apply the Euler operator of weights
\(3,2\) to \(x,L\) and subtract nine times the polynomial,
obtaining \((2k-9)bL^k\). At \(L=0\) the only critical
candidate is \(x=0\), smooth when \(k=1\) and with binomial
initial part \(x^3+bL^2\) when \(k=2\).

At the last wall, put \(y=1\) to obtain
\[
      P(x)+(ax+b)L^3,\qquad P(x)=x^3+Ax^2+Bx+C,\quad C\ne0 .
\]
Again \(P\) has no triple root by the reduced sextic argument.
At a critical zero with \(L\ne0\) one has \(ax+b=0\);
the Hessian has nonzero mixed entry \(3aL^2\) and zero \(LL\) entry,
thus has rank two. At \(L=0\) any critical \(x_0\) is a double
root of \(P\). If \(ax_0+b\ne0\), the next light face is
\(X^2+cL^3\) up to nonzero scalar factors. Otherwise the weights
\(3,1\) give \(X^2+cXL^3,\ c\ne0\). For the last face all
punctured graph criticals are absent: its two derivatives are
\(2X+cL^3,\ 3cXL^2\). Its radial ratio is \(2/3\).
For every wall, at \(y=0\) any punctured point is already
noncritical for \(x^3+axL^3\). The noncritical graph estimate
applies throughout, and the first radial ratio
\((3+2+w)/9\ge13/18\) suffices.

\subsection{Conics on the sextic twelvefold}
\label{sec:nov27-sextic-twelve-conics}
Let \(X_6^{12}\subset\mathbf P^{13}\) be smooth, \(C\subset X\)
an integral conic, and \(b=29/19\). Then \(\delta_C(X,H)\ge13/b=247/29\).
The proper-plane bound is \(3/2\), and a nonzero plane polar
divisible by the conic equation is covered by the mixed pencil
estimate, since both \(17/12\) and \(61/42\) are smaller than \(b\).
In particular this covers polar pencils. An ordinary positive
splitting degree at most six gives angular denominator at most
\(3/2\), and the ordinary radial bound is smaller. At rank at
least four, if the largest degree is at least seven, each partition
of ten into at least three positive parts with this property contains
a one. The corresponding linear syzygy supplies two independent
polars \(t\ell,-s_0\ell\), in base coordinates \([s_0:t]\); a
third polar avoiding the roots of \(\ell\) gives a basepoint-free
net with low split degree one. Retain these three heavy coordinates
and the original common kernel by an actual linear flag through
the conic plane. Every cut is smooth along \(C\), and the light
order \(h\) is preserved. The common kernel has codimension at
most five in all remaining cases, and hence dimension greater
than six. Thus \(2\le h\le6\) as in the preceding conic estimates.

For the ternary axes the ratios \([1+h(3-1/6)]/(1+2h)\)
are smaller than \(b\). Write \(s\) for the smaller positive
split degree on a net. Tests (E),(P) of the conic chapter and
the two marked-point tests of Section~\ref{sec:nov27-small-conic-2d}
give the following table for \(h=2,3,4,5,6\) at \(b=29/19\).
The labels \(1,2\) denote the single- and double-degree graph
tests in that section. Entries \(*,\dagger\) are treated below.
\[
\begin{array}{c|ccccc}
s&2&3&4&5&6\\\hline
1&*&1&1&P&P\\
2&*&2&E&2&\dagger\\
3,4,5&E&E&E&E&E
\end{array}
\]
The off-branch, low-coordinate and free estimates in those tests
are all smaller than \(b\). For \(s=1,2\) there remains the
distinguished point \(P\) on the generic branch line \(L\).

\paragraph{Retaining the moving point in the pole specialization.}
Suppose \(h=2\). Write \(S\) for the normalized actual plane mean
and \(A\) for ordinary plane discrepancy. For every geometric
generic test valuation centered at \(P\) we claim
\[
       S+C_s v(\mathfrak m_P)\le (4/3)A,\qquad
              C_1=4/9,\qquad C_2=1/3.                  \tag{SC1}
\]
Use the actual blocks of the conic chapter, with \(Q=x^2=-L/a\),
\(\deg a=8\). The regular sector is the complete product of base
degree \(2p\) and plane degree \(k\). It has mass \(1/3\), with
conditional normalized mean plane degree three. The positive pole
sector has mass \(2/3\). Refine the increasing pole filtration at
level \(i\ge1\) by the successive coefficients of the leading
residue over the finite divisor \((a=0)\), ordering its roots
and then their jets. If the root \(p_0\) has multiplicity \(e\),
these jets have indices \(j=0,\ldots,e-1\). A representative
whose leading residue vanishes through jet \(j-1\) has leading
term at this root, up to a base scalar,
\[
             L(p_0)^i F,\qquad F\in H^0(\mathbf P^2,\mathcal O(k-i)).
\]
Indeed the preceding pole levels have pole order at most
\((i-1)e\), strictly below \(ie-j\), and \(L(p_0)\) is a nonzero
linear form. The quotient gives the whole displayed coefficient
space by binary interpolation whenever the remaining base degree
is at least seven. This also allows prescribed vanishing at the
earlier roots and their jets. The bounded boundary strip has
vanishing relative dimension; its contribution to normalized
orders is \(o(1)A\), since any degree-\(k\) plane form has order
at most \(kA\) and \(k\) is bounded by the section degree.
For every fixed root and jet the conditional normalized means of
\(i\) and \(k-i\) are \(1/2\) and three. Indeed each positive
root-jet block has \(\binom{k-i+2}{2}\) coefficients and the
same range of \(i,k,p\) away from that strip, so the
positive-sector simplex computation applies separately.

Choose bases simultaneously adapted to the tested valuation and
to this refined flag of complex vector spaces, in each graded
block. For a root \(p_0\), substitute \(b_0\mapsto
p_0+\tau(b_0-p_0)\) in a base coordinate \(b_0\) containing
the root in its affine chart. For \(\tau\ne0\) this is a base
field automorphism. Rescaling each basis element by its leading
scalar in \(\tau\) specializes the elements in that root-jet group
to bases of the displayed complete systems times \(L(p_0)^i\).
Carry the ideal of the moving section \(P\) in this same family.
It specializes to the ideal of \(P(p_0)\in L(p_0)\), since \(P\)
is a morphism on the base and \(L\) is relatively basepoint-free.
The normalized complete contribution to the limit is at most
\((1+o(1))A\), by the uniform complete-series basis thresholds
of \(\mathbf P^2\); the fixed line coefficient tends to \(1/2\).
The pair with line \(L(p_0)\) of coefficient one and ideal
\(\mathfrak m_{P(p_0)}\) of exponent one is log canonical, as
seen by blowing up that point. Lower semicontinuity of thresholds
for the product with the marked ideal, on the proper plane, now
gives the conditional estimate
\[
          S_{\rm pos}+\tfrac12v(\mathfrak m_P)\le\tfrac32 A .
\]
At finite level test strictly below the limiting coefficient.
The separate specializations compare the same original valuation
with each conditional bound, so their inequalities can be added.

For the regular product conditioning, a further base-order flag
at an ordinary point specializes each residue to the complete
plane series. Thus \(S_{\rm reg}\le A\). If \(s=1\), the graph
on the plane is a line parametrization and one has the refinement
\[
             S_{\rm reg}+\tfrac13v(\mathfrak m_P)\le A.
\]
Let \(Z=(z_3=0)\) be the constant line containing the image of
\(P\). Choose a basis adapted to the tested valuation and to
order along \(Z\). At order \(a\) on \(Z\), a regular block of
bidegree \((2p,k)\) restricts, after dividing by \(z_3^a\), to
the full product of binary systems of bidegree \((2p,k-a)\).
Its subspace whose order at \(P\) is at least \(j\) consists
exactly of the multiples of the primitive graph equation
\(T^j=(tz_1-s_0z_2)^j\). Its dimension is
\((2p-j+1)_+(k-a-j+1)_+\).
The normalized mean \(a\) is one. The normalized restricted
order is at most \(2/3\): in these dimensions the shift for
divisibility by \(T^j\) has cost \(j+j/2\) in \(k-a,p\),
respectively. Equivalently, integrate the two tail dimensions
of the product over the section-degree simplex. A shift of unit
cost gives normalized mean one and cost \(3/2\) gives \(2/3\);
rounding \(j/2\) affects only the lower-order term.
On the geometric generic \(Z\) the marked ideal restricts to
\(\mathcal O_Z(-P)\). Its exponent \(1/3\) brings the total
coefficient at \(P\) to at most one. Inversion of adjunction on
the smooth line \(Z\), testing slightly below the limiting
coefficients at finite degree, proves the refinement. Averaging
the regular and positive sectors proves (SC1).

\paragraph{The numerical plane comparison.}
Normalize \(v(\mathfrak m_P)=1\) and put \(u=v(L)\ge1\), so
\(A\ge u+1\). The actual residual intersection and multiplicity
estimates give
\[
 S\le u/2+(A-1)m,\qquad m=\min(M,I/u),\qquad
 (M,I)=
 \begin{cases}(31/18,7/4)&s=1,\\(5/3,5/3)&s=2.\end{cases}
\]
We require \(S\le b(A-u/2)\). Estimate (SC1) suffices when
\(A\ge A_*(u)=(bu/2-C_s)/(b-4/3)\). At \(A=u+1\)
the residual estimate has slack at least
\(B=b+(b-1)/2-I>0\), since \(um\le I\).
It continues to suffice for all \(A\ge u+1\) when \(m\le b\).
Otherwise \(u\le I/b\); moving from \(u+1\) to \(A_*(u)\)
can lose at most
\[
       D=(M-b)\max\{0,A_*(I/b)-I/b-1\},
\]
because \(A_*(u)-u-1\) is increasing. The exact values are
\[
\begin{array}{c|ccc}
s&B&D&B-D\\\hline
1&3/76&43349/2618352&60007/2618352\\
2&7/57&3820/54549&2879/54549
\end{array}
\]
which proves the two starred cases.

\paragraph{Terminal rank-three nets.}
The only remaining entry has \(s=2\), polar rank precisely three,
and \(h=6\). Cut by the quadric cone over the conic as in
Section~\ref{sec:nov27-conic-quadric-flag}, retaining the whole
common light kernel and two basepoint-free heavy directions. The
order on that kernel after imposing the conic equation is at least
six. Its vertex has dimension \(12-3-2=7>6\) and cannot lie in
\(X\). Thus its first form is the nonzero base-independent sextic
\(G_6\). The tangent Gauss-map argument of that section with
\(c=3\) gives singular codimension at least \(12-2c-1=5\) in
the projective light space. Four general constant lights therefore
cut a smooth projective sextic surface. Give the lights weight
one and the two heavy directions weight six. The radial denominator
is \(1+(1-1/6)/(1+4/6)=3/2\). At every punctured point of a
positive pole graph its light polynomial has nonzero gradient,
by Euler's identity and the smooth zero surface; it is not a pure
power. Section~\ref{sec:nov27-six-weighted-noncritical} bounds its
signed denominator by \(13/8<31/19\). The quadric-flag proportions
\(1/6,5/6\) yield an angular denominator at most \(b\).
The linear and quadric flags and weighted adjunction now give the
claimed threshold for all conics.

\paragraph{The short-curve criterion for the sextic twelvefold.}
Let \(X\subset\mathbf P^{13}\) be a smooth sextic of dimension
\(12\). If every line and every integral conic \(C\) on \(X\)
satisfies \(\delta_C(X,H)\ge1599/188\), then \(\delta(X,H)>8\).
Indeed the pointed hyperplane section is a smooth sextic elevenfold,
and Theorem~\ref{sextic:main} and hyperplane flag adjunction give
\(\beta=\delta(X,H)>(13/12)7\).
Apply the point and positive-center argument of
Section~\ref{sec:nov-second-row} with \(N=13,\ r=8,\ f=1/6\),
share \(a=21/50\) and \(\epsilon=10^{-8}\). The outward-rounded
constants are \(Q=4951/2000,\ B=1501/500\). The surface cutoff and
both signed comparisons at \(q=0,Q\) hold up to \(t=B\).

For \(B<t<6\) the exceptional prime of class \(jH\) has multiplicity
\(k=jt\), and its general pointed plane curve section is integral
by the projection argument. For \(1\le j<20\) use
\(I=6j(j+2)+k+1\), improved for \(j\le2\) by the sectional
Milnor bound, and discard \(k(k-1)>I\). For
\(m=\min(k,\sqrt I)\) the comparison
\[
  \frac{6}{t}+
   \left(1-\frac{6}{t^2}\right)
   \frac{I m}{j(I+m(m-1))}<B
\]
holds except at \((j,k)=(1,4),(1,5)\), where respectively \(16,25\)
principal bins suffice. For \(j\ge20\) put
\(E=6+18/20+1/400\) and \(L_*=E/(2\sqrt E-1/20)\).
The maximum of \(6/t+(1-6/t^2)L_*\) on \([B,6]\), obtained
by clamping \(2L_*\) to this interval, is less than \(B\).
Finally \(t=6\) passes the principal test with \(35\) bins.
The exact comparisons are given in \path{nov27-bin-sextic.py}.

\section{The quintic tenfold}\label{sec:nov27-quintic-ten}
Let \(X_5^{10}\subset\mathbf P^{11}\) be smooth. We prove that the
local threshold at every line and integral conic is at least \(385/51\)
in hyperplane units. The point argument below explains how this
bound implies \(\delta(X,H)>7\).

\paragraph{A short-curve criterion for a quintic tenfold.}
Let \(X\subset\mathbf P^{11}\) be a smooth quintic of dimension ten.
If every line and integral conic \(C\subset X\) satisfies
\(\delta_C(X,H)\ge385/51\), then \(\delta(X,H)>7\).
Indeed \(N=11\), \(r=7\), and a general smooth pointed hyperplane
section satisfies Theorem~\ref{quintic:main}. Apply the criterion of
Section~\ref{sec:nov-second-row} with \(f=1/5,\ a=429/1000\) and
\(\epsilon=10^{-8}\). It gives \(B=1373/500,\ Q=18231/10000\).
The surface cutoff and both signed endpoint comparisons below \(B\)
are strict.

When \(B<t=\tau_p(H)<5\), the exceptional prime of degree class
\(jH\) and multiplicity \(k=jt\) has an integral general pointed
curve section. For \(1\le j<25\), use \(I=5j(j+1)+k+1\), taking
its minimum with the sectional Milnor bound when \(j\le2\), and
discard \(k(k-1)>I\). The comparison of
\(5/t+(1-5/t^2)I m/[j(I+m(m-1))]\) with \(B\), where
\(m=\min(k,\sqrt I)\), is strict except at \((j,k)=(1,4)\);
there fifteen principal bins suffice. For \(j\ge25\) use
\(E=5+10/25+1/625\), \(L_*=E/(2\sqrt E-1/25)\).
Clamping \(2L_*\) to \([B,5]\) gives
\(\max(5/t+(1-5/t^2)L_*)<B\).

At \(t=5\) the degree of the exceptional prime equals its multiplicity,
so it is a projective cone with vertex \(p\). Every generating line
lies in \(T_pX\); hence its support is the tangent hyperplane section.
That section is integral since its singular locus is finite by the finite
Gauss map and its dimension is nine. The cone base is smooth by the
same isolation and is a quintic eightfold, uniformly K-stable in the
preceding range. The generalized Eckardt formula
\cite[Corollary~4.19]{ar2003.13788} gives
\(\delta_p(X,H)=110/14>385/51>U\). This excludes the endpoint.
The finite comparisons are reproduced by \path{nov27-bin-quintic.py}.

\paragraph{The Hessian on a general four-dimensional space of lights.}
Let \(G\) be homogeneous of degree \(e\ge2\) on a vector space
\(V\), with \(\mathbf P(V)=\mathbf P^N\), \(N\ge4\). Suppose the
projective singular locus \(B=\operatorname{Sing}(G=0)\) has
codimension at least five, with the empty set allowed. For a
general four-dimensional linear subspace \(U\subset V\), the
Hessian of \(G|_U\) has rank at least two at every nonzero vector
of \(U\). In particular its ordinary affine graph has signed
denominator at most \(3/2\) at every punctured graph point by the
rank-two quadratic estimate.

Write \(\sigma=\dim B\), taking \(\sigma=-1\) if \(B\) is empty.
The rational gradient map on \(\mathbf P^N\setminus B\) has
fiber dimensions at most \(\sigma+1\). Indeed choose a polar
coordinate nonzero at the value of the map. Its degree-\((e-1)\)
equation does not vanish along the fiber, so its intersection with
the projective closure of a positive-dimensional fiber component
lies in \(B\) and has pure codimension one on that closure.
The stratum of points where the full affine Hessian has rank \(r\)
therefore has dimension at most \(r+\sigma\). In fact Euler
differentiation sends the radial vector to \((e-1)dG\), so the
projective differential has rank \(r-1\). On the smooth open of
each stratum component the dimension of its image is at most
\(r-1\), by separability in characteristic zero. The fiber bound
then gives the assertion.

A general projective three-space avoids \(B\) and cuts the smooth
locus of \((G=0)\) transversely. If \(G(p)=0\) on this section,
its restricted gradient is nonzero and
\(\operatorname{Hess}G(p)(p,p)=0\) by Euler. A rank-one symmetric
form with isotropic vector \(p\) has \(p\) in its radical, which
contradicts \(\operatorname{Hess}G(p)(p,-)=(e-1)dG_p\) on \(U\).

Suppose \(G(p)\ne0\), so the radial line is nonisotropic and its
Hessian orthogonal complement is \(T=\ker dG_p\). Write \(k=r-1\)
for the rank on \(T\). Among four-spaces \(U\) through \(p\),
their intersections \(Y=U\cap T\) parametrize all three-spaces in
\(T\). The restricted Hessian has rank one exactly when \(Y\) is
totally isotropic. For \(k=0\), the \(p\)-locus has dimension at
most \(\sigma+1<N-3\) and is avoided. For \(1\le k\le2\), a
totally isotropic \(Y\) is contained in one of at most two fixed
hyperplanes of \(T\). This imposes at least three conditions on
\(Y\), and the \(p\)-locus has dimension at most
\(\sigma+3<N\). For \(k\ge3\), total isotropy imposes at least
four conditions. On ordered independent generating triples the
three equations for the vectors to be null cut out an open subset
of the third product of an irreducible quadric cone, of codimension
three. The pairing between the first two vectors is not identically
zero on this product, yielding a fourth condition. Passing to the
Grassmannian preserves this codimension. Finally, containing \(p\)
already imposes \(N-3\) conditions on \(U\). In each stratum the
incidence of the bad pairs \((p,U)\) thus has dimension smaller
than the Grassmannian, proving the claim.

\paragraph{The sharp ordinary graph bound for a square and a cube.}
On the affine chart of \(\mathbf P^3\), put \(W=(z+x^2+y^3=0)\).
Then \(4S_{\mathcal O(1)}(v)+v(W)\le(17/11)A(v)\) at every
point in this chart. On \(W\) with \(y\ne0\) the Hessian has rank
two. The quadratic flag therefore gives \(3/2\). The translation
\(x\mapsto x+a\) with the corresponding affine-linear change of
\(z\) preserves the graph, so it remains to work at the origin.

The projective graph closure \(D=(T^2z+Tx^2+y^3=0)\) is a normal
cubic surface (its only singular point is at infinity). At the
origin \(\delta_0(D,\mathcal O_D(1))\ge5/3\). Here is the direct
calculation. In affine coordinates the sections of degree \(m\)
are spanned by \(x^a y^b Q^c\), \(a+b+c\le m\), \(Q=x^2+y^3\).
Indeed all surface sections are restrictions by its hypersurface
sequence. Extract with weights \(3,2\) on \(x,y\). At weight
\(j\) set \(q=\lceil(j-3m)_+/3\rceil\), \(E=j-6q\). The
weighted piece is exactly \(Q^q\) times the full binary weighted
degree-\(E\) space, interpreted as zero when necessary. The
inequality \(j=3a+2b+6c\le3m+3c\) forces the divisor. Conversely
write \(y^{3i+r}=y^r(Q-x^2)^i\), \(0\le r<3\), to express any
binary monomial of weight \(E\) in section degree at most
\(\lceil E/3\rceil\le m-q\). The extraction thus has mean weight
\(3m+o(m)\) and discrepancy five. To see all means at once, put
\(t=j/m\); the weight dimension is \(E/6+O(1)\). Hence the density
in \(t\in[0,6]\) is proportional to \(w(t)=\min(t,6-t)\), giving
\(\overline{E}/m=2+o(1)\), \(\overline q/m=1/6+o(1)\).
At the two coordinate points of the exceptional weighted line,
the mean vanishing of a complete weight-\(E\) binary basis in
the smooth stack coordinates is at most \(E/6+O(1),E/4+O(1)\),
by averaging their respective exponents. Neither coordinate point
lies on \(Q\). Their asymptotic coefficients are at most \(1/2\).
At a noncoordinate point its complete basis has mean at most
\(E/12+O(1)\), since the consecutive monomials give an ordinary
complete system of length \(E/6+O(1)\) in a torus coordinate.
The factor \(Q\) adds at most the mean \(1/6\) there, as its mark
is simple. Thus all angular coefficients are at most \(1/2\).
Weighted adjunction gives the stated threshold at the origin.

Refining the complete system of \(\mathbf P^3\) by the cubic
surface, its coefficient in units of four is \(1/3\) and its
residual series has mean degree three in those units, by the
hypersurface restriction sequence. The surface threshold and
divisor adjunction consequently give
\[
 4S(v)+(9/5-1/3)v(W)\le(9/5)A(v).
\]
For this local application the origin is smooth on \(D\).
Interpolating with \(4S\le A\) yields the coefficient \(17/11\).
It is sharp at the point with weights \((3,2,6)\), where the
complete normalized mean is eleven and the graph order is six.

\paragraph{The sharp bound for a quartic weighted face.}
Take \(Q=x^2+Bxy^2+Cy^4\), with \((B,C)\ne(0,0)\), on the
ordinary affine chart of \(\mathbf P^3\), and set \(W=(z+Q=0)\).
Then
\[
             4S_{\mathcal O(1)}(v)+v(W)\le \frac{11}{7}\,A(v).
\]
At the origin the weights \((2,1,4)\) give equality. Let \(D\)
be the projective graph closure, \(H=\mathcal O_D(1)\), and
\(e=\deg D\), equal to four if \(C\ne0\), three otherwise.
Near the origin \(D\) is smooth; only its section spaces and this
smooth neighborhood are used. Its hypersurface sequence identifies
the restricted degree-\(m\) space with the span of
\(x^a y^b Q^c\), \(a+b+c\le m\). On the surface extract
with weights \((2,1)\), so the discrepancy is three.

If \(C\ne0\), set \(q=\lceil(j-2m)_+/2\rceil,\ E=j-4q\)
at weight \(j\). Each piece is contained in \(Q^q S_E\), with
\(S_E\) the complete binary weighted system, by
\(j\le 2m+2c\). It contains \(Q^{q+2}S_{E-8}\).
Indeed, division by the monic equation expressing \(y^4\) in
terms of \(x,Q,y\) gives a homogeneous basis \(1,y,y^2,y^3\)
over \(\mathbf C[x,Q]\). Its degree-\(E'\) monomials have ordinary
section degree at most \(E'/2+3/2\). Multiplication by \(Q^{q+2}\)
therefore suffices. For \(E'<0\) take the zero space.
The containing and contained spaces differ by a bounded dimension
uniformly in \(m,j\). At any angular point the vanishing sum of a
subspace is bounded by that of its containing space, and the
difference of the total dimensions of all weights is \(O(m)\).
Thus the leading means are computed by the containing systems.

With \(t=j/m\) the density is proportional to \(\min(t,4-t)\) on
\([0,4]\); the normalized means of \(j,E,q\) are \(2,4/3,1/6\).
The stack coordinate orders of a complete binary basis at \(x=0\)
and \(y=0\) have respective means \(E/4+O(1)\) and \(E/2+O(1)\).
The factor \(Q\) misses both coordinate points. At a torus point
the residual series in the coordinate \(y^2/x\) has mean at most
\(E/4+O(1)\), and a root of \(Q\) has multiplicity at most two.
Consequently every angular coefficient is at most \(2/3\),
including for a double root. The radial coefficient is \(2\) with
discrepancy three, so local weighted adjunction gives
\(\delta_0(D,H)\ge3/2\).

When \(C=0,\ B\ne0\), put \(j=2u+\epsilon,\ 0\le\epsilon<2\)
and \(w=y^2/x\). Standard representatives modulo the relation
\(Bxy^2=Q-x^2\) have either \(a=0\) or \(b\le1\).
This reduction never increases ordinary section degree.
Dividing by \(y^\epsilon x^u\), the weight piece is
\[
       (1+Bw)^q\mathbf C[w]_{\le L-q},\qquad
 q=\max(0,u+\epsilon-m),\quad
 L=u-\max(0,\lceil(2u+\epsilon-m)/3\rceil).
\]
Here only nonempty weight pieces are considered. To see equality,
the terms with \(a\ge1,\ b=\epsilon\) give
\((1+Bw)^c\) with \(q\le c\le\lfloor(u-1)/2\rfloor\).
The terms with \(a=0\) give
\(w^{u-2c}(1+Bw)^c\) with \(u-L\le c\le\lfloor u/2\rfloor\).
Their degrees \(u-c\) supply the remaining successive degrees up
to \(L\); every term has the stated common factor. If the latter
interval is empty, the weight piece is empty or has bounded
dimension near the extreme weight, which has no effect on the
limit.

Put \(t=u/m\), \(l=\min(t,(1+t)/3)\), \(q_0=(t-1)_+\).
The weight density is proportional to \(l-q_0\) on \([0,2]\).
Integration on \([0,1/2],[1/2,1],[1,2]\) gives
\[
 \overline{t}=17/18,\quad\overline l=17/27,\quad
 \overline{q_0}=4/27,\quad\overline{l-q_0}=13/27 .
\]
The mean radial order is \(2\overline t=17/9\).
The angular coefficients at \(y=0,x=0\) are respectively
\(\overline{l-q_0}=13/27\) and
\(\overline t-(\overline l+\overline{q_0})/2=5/9\), by the
two stack coordinate orders on the monomial basis. At any
torus point the coefficient is at most
\((\overline l+\overline{q_0})/2=7/18\). All are at most
\(17/27\), proving \(\delta_0(D,H)\ge27/17\).

In both cases the degree-\(m\) pieces after refining
\(\mathbf P^3\) along \(D\) are the complete restrictions in
degree \(m-ei\), by the hypersurface sequence. In units of four
the mean of \(i/m\) is \(1/e\) and the residual degree has mean
three. Take \(g=2\) for \(e=4\), \(g=17/9\) for \(e=3\).
Divisor adjunction at the smooth origin gives
\(4S+(g-1/e)v(W)\le gA\).
Interpolating with \(4S\le A\) gives \(11/7\) in both cases.

For completeness the bound at the origin implies the bound
throughout the affine chart. The diagonal action with weights
\((2,1,4)\) preserves \(W\) and contracts the chart to the origin.
At degree \(m\), choose a complete basis compatible with the
given divisorial order and the weight filtration simultaneously.
After scalar rescaling its initial sections form a complete
basis. Its basis divisor plus the fixed multiple of \(W\) is the
flat limit of the translates of the original divisor with boundary.
The finite-degree estimates above have errors tending uniformly to
zero after normalization (the exceptional angular errors are
\(O(1)\) per section, and dimension errors are \(O(m)\) among
\(O(m^2)\) sections). For each fixed multiplicative loss in the
claimed threshold, all sufficiently large \(m\) thus give klt at
the origin for the limiting basis divisor with boundary. Openness
of klt along the contracting orbit gives the corresponding inequality
at the original point. Then let \(m\) tend to infinity and the loss
tend to zero.

\paragraph{A square plus a quintic monomial.}
Let \(Q=x^2+u^3v^2\) and \(W=(z+Q=0)\) in the affine chart of
\(\mathbf P^4\). Then \(5S_{\mathcal O(1)}(E)+E(W)\le(14/9)A(E)\)
at the origin. Write \(D\) for the degree-five projective graph
closure. We prove \(\delta_0(D,\mathcal O_D(1))\ge2\), using only
the smooth affine chart of \(D\). Its complete degree-\(m\)
sections have basis
\[
  x^a u^b v^c Q^q,\qquad a+b+c+q\le m,\qquad b<3\ \text{or}\ c<2.
\]
Indeed \(u^3v^2=Q-x^2\) is a monic reduction that lowers ordinary
degree; distinct normal forms are independent by polynomial
division, and restriction is complete by the hypersurface sequence.

Use the two-dimensional torus for which \(u^3v^2\) and \(x^2\)
have the same character and put \(w=u^3v^2/x^2\). Within a fixed
character, the normal exponents \(b,c\) are unique: all candidate
pairs differ by an integer multiple of \((3,2)\), and precisely
the minimal nonnegative pair satisfies the normal-form condition.
Thus \(B=a+2q\) is fixed and the character piece is
\[
  x^B u^b v^c(1+w)^{q_0}\mathbf C[w]_{\le l},\quad
  q_0=\max(0,B+b+c-m),\quad l=\lfloor B/2\rfloor-q_0.
\]
Empty pieces are omitted. In the basis counting measure the
limiting normalized means, dividing exponents by \(m\), are
\[
   (\bar a,\bar b,\bar c,\bar q)=(1/4,1/10,3/20,1/4),
   \qquad \overline{l/2}=1/8 .
\]
The \(c=0,1\) layers have two fifths of the mass with \(a,b,q\)
free, while the \(b=0,1,2\) layers have three fifths with \(a,c,q\)
free; their overlap has lower order. Each free exponent averages
\(1/4\) by the three-simplex. The last identity follows because
\(\lfloor B/2\rfloor-q=\lfloor a/2\rfloor\) and \(q\) ranges over
the consecutive interval from \(q_0\) to \(q_0+l\).

For a divisorial order on the affine chart put
\((\alpha,\beta,\gamma)=(E(x),E(u),E(v))\) and
\(\omega=3\beta+2\gamma-2\alpha\). The discrepancy is at least
\(\alpha+\beta+\gamma\) by the coordinate normal crossings.
If \(\omega>0\), the complete \(w\)-polynomials in each character
have maximal mean at \(w=0\) equal to \(l\omega/2\).
Since \(1+w\) is a unit, the resulting normalized mean is at most
\[
 \tfrac34\alpha+\tfrac1{10}\beta+\tfrac3{20}\gamma+
             \tfrac18\omega
   =\tfrac12\alpha+\tfrac{19}{40}\beta+\tfrac25\gamma .
\]
If \(\omega<0\), the pole order is the degree in \(w\), so the
mean is at most
\[
 \tfrac14\alpha+\tfrac{17}{20}\beta+\tfrac{13}{20}\gamma
 \le\tfrac12(\alpha+\beta+\gamma).
\]
The inequality uses \(2\alpha\ge3\beta+2\gamma\).
When \(\omega=0\), the monomial part gives
\(3\alpha/4+\beta/10+3\gamma/20\), also at most half the sum.
If \(w\) specializes to a constant \(t\ne0\), put
\(\lambda=E(w-t)>0\). The additional mean is at most
\(\bar q\lambda=\lambda/4\) at \(t=-1\), and at any other
constant at most \(\overline{l/2}\lambda=\lambda/8\).
Here \(A(E)\ge\alpha+\beta+\gamma+\lambda\). Indeed replacing
\(d\log v\) by \(d\log(w-t)\) in
\(d\log x\wedge d\log u\wedge d\log v\) introduces the factor
\((w-t)/(2w)\). A wedge of logarithmic differentials has at most
a simple pole along any prime on a smooth model, giving precisely
the stated inequality. Nonconstant residue of \(w\) contributes
no further order. These bounds have uniform \(o(1)A(E)\) errors,
by the finite simplex sums and the bounded integer remainders.

To apply the character computation to an arbitrary valuation
centered over the origin, at level \(m\) choose a positive integral
cocharacter separating the finitely many characters of that level,
and a basis simultaneously adapted to its filtration and the given
valuation. Its torus initial basis splits into the character pieces
above. The limiting basis divisor is log canonical locally at the
origin at coefficient \(2-o(1)\), by the estimates for arbitrary
orders just given. Lower semicontinuity at the torus-fixed origin
transfers every strictly smaller coefficient to the original basis
divisor. Thus \(\delta_0(D,\mathcal O_D(1))\ge2\).
On restricting the complete ambient series to \(D\), its normalized
fixed coefficient is \(1/5\) and its residual degree coefficient is
four. Divisor adjunction therefore yields
\(5S+(2-1/5)E(W)\le2A\); interpolation with \(5S\le A\)
proves \(14/9\).

For the whole punctured graph the uniform bound \(11/7\) suffices.
Translate \(x\) to zero by an affine-linear change of \(z\).
At \(u=v=0\) the preceding calculation applies. At \(uv\ne0\)
the Hessian has rank at least two, since the binary monomial has a
nonzero Hessian and \(x^2\) is a disjoint nondegenerate block.
When \(u\ne0,v=0\) the full Hessian again has rank two.
When \(u=0,v\ne0\), fix \(v\) to its nonzero value and apply the
square-plus-cube bound \(17/11\) on the section. Ordinary linear
adjunction recovers the required estimate upstairs.

\paragraph{A smooth quintic surface in four light variables.}
For a nonsingular homogeneous quintic \(P\) in four variables,
every punctured point of the graph \(z+P=0\) has signed denominator
at most \(11/7\). At a rank-two or higher Hessian point use the
quadratic flag. Rank one at a nonzero point forces \(P(p)\ne0\)
and \(dP_p\ne0\), by smoothness and Euler's identity.
Write the lights as \((1+x)p+Y\), with \(Y\in\ker dP_p\).
The expansion is
\[
 c(1+x)^5+(1+x)^2P_3(Y)+(1+x)P_4(Y)+P_5(Y),\quad c\ne0 .
\]
If \(P_3\) or \(P_4\) is nonzero, a general one-dimensional
\(Y\)-space yields the square-plus-cube or square-plus-quartic
face after subtracting the affine tangent, with signed bound
\(17/11\) or \(11/7\). Otherwise \(P_5\) defines a smooth
plane quintic, since the original hypersurface is smooth.
A general tangent line to this plane curve has contact exactly
two in characteristic zero (equivalently, the finite Gauss map of
the curve is generically unramified). Use its two-dimensional
vector space for \(Y\). Choose binary coordinates \(u,v\) with
the double root at \(v=0\), so the first binary monomial is
\(b u^3v^2\), \(b\ne0\). Give \(u,v\) weights \(1,1+\epsilon\)
and \(x\) weight \((5+2\epsilon)/2\), for small positive rational
\(\epsilon\). After subtracting the affine tangent, the first
graph has precisely the preceding light part \(a x^2+b u^3v^2\),
\(ab\ne0\). The sum of the light weights divided by the graph
weight is \(1/2+(2+\epsilon)/(5+2\epsilon)>3/4\). Its radial
and punctured estimates prove the claim by weighted adjunction.

\subsection{Quintic tenfold conics at the short-curve threshold}
\label{sec:nov27-quintic-conic-strong}
Every integral conic \(C\) on a smooth quintic tenfold satisfies
\(\delta_C(X,H)\ge385/51\). Put \(b=51/35\). The proper plane
denominator \(10/7\), both mixed pencil terms \(7/5,41/30\),
and every ordinary angular denominator with largest split degree
at most four are less than \(b\). In rank at least four, any
partition of eight with at least three summands and largest
degree at least five contains a one (the possibilities with
three parts are \(6,1,1\) and \(5,2,1\)). Thus reduce to a
basepoint-free net with a linear syzygy, retaining the original
common kernel as in Section~\ref{sec:nov27-sextic-twelve-conics}.
The previously displayed tests in this paper, including the
single- and double-degree graph tests, suffice for \(h=2,3,4,5\)
in every net except \((s,h)=(1,2),(2,2),(2,5)\). The ordinary
radial and ternary-axis ratios are also smaller than \(b\).

For \(h=2\) the marked pole comparison in
Section~\ref{sec:nov27-sextic-twelve-conics} applies with regular
and positive masses \(2/5,3/5\), and with \(\deg a=6\).
The positive conditional bound is
\(S_{\rm pos}+v(\mathfrak m_P)/2\le3A/2\).
For the regular product at both split degrees we have
\[
   S_{\rm reg}+\frac{s}{s+2}v(\mathfrak m_P)\le A.
\]
Here is the additional argument for \(s=2\). The radial
coefficient on blowing up \(P\) is at most \(3/2\) by the
projected graph multiplicity estimate. For a double cover of a
line the two jet coordinates have exponent means \(1,1/2\),
so all angular coefficients are at most one. For a Veronese
graph the unmarked angular estimate of
Section~\ref{sec:nov27-small-conic-2d}, with \(t=1,\ u=v=1/2\)
and \(D=1\) in \(V_0\), is \(8/9\). The marked ideal adds
only \(1/2\) to the radial coefficient. This proves the bound
by adjunction. It gives the combined estimates
\[
   S\le(13/10)A-C_s,\qquad C_1=13/30,\quad C_2=1/2,
\]
in the normalization \(v(\mathfrak m_P)=1\). Write \(u=v(L)\)
and \(A\ge u+1\). The residual bound remains
\(S\le u/2+(A-1)\min(M,I/u)\), with
\((M,I)=(5/3,17/10),(8/5,8/5)\) for \(s=1,2\), respectively.

For \(s=2\) the marked estimate suffices for \(1\le u\le I/b\):
the required lower bound on \(A\) is
\((bu/2-C_2)/(b-13/10)\), whose difference from \(u+1\)
is increasing and equals \(-106/561\) at \(u=I/b\).
For \(u\ge I/b\) the residual slope is at most \(b\) and the
endpoint slack at \(A=u+1\) is at least \(b+(b-1)/2-I=3/35\).

Suppose \(s=1\). At \(u=1\) the first marked estimate suffices,
with slack \(2/105\) at \(A=2\). If \(u>1\), the center on
the first blowup is the branch direction \(P_1\).
The regular radial coefficient is at most \(5/3\); its strict
transform has intersection with the exceptional line at \(P_1\)
at most
\[
     J=\frac1{1-1/4}\left(\frac1{1+8/2}
                    -\frac{1/4}{8+1/2}\right)=58/255
\]
by the marked product computation with \(t=1\). Inversion of
adjunction for the exceptional line gives
\(S_{\rm reg}\le5/3+J(A-2)\). Indeed the log discrepancy on
the blowup is \(A-1\) and the exceptional divisor has value one.
Use a basis of this regular subgroup adapted to both tested
and radial orders, independently of the residue specializations
in the positive sector. Averaging now gives \(S\le KA+R\),
where \(K=2527/2550,\ R=157/850\).
For \(1<u\le17/16\) this implies \(S\le b(A-u/2)\);
the worst endpoint \(A=u+1,\ u=17/16\) has slack \(9/3400\).
For \(u\ge17/16\) one has \(b+(b-1)u/2\ge I\), so the residual
bound suffices at \(A=u+1\). When \(I/u\le b\) it also suffices
thereafter. Otherwise it suffices to test it up to
\(A_*(u)=(R+bu/2)/(b-K)\), beyond which the bound \(KA+R\)
applies. At that endpoint the excess over \(b(A-u/2)\) is
\[
 -\frac{43636}{41605}u+\frac{173043}{83210}-\frac{272}{265u}.
\]
This decreases for \(u\ge1\), since the ratio of the absolute
last and first coefficients is \(10676/10909<1\); at \(17/16\)
it equals \(-27/33284\). This proves both order-two cases.

For \(h=5,\ s=2\) the polar rank is three. Let \(V\) be its full
linear polar kernel and \(\Pi\) the conic plane. Write the
restriction of its defining polynomial to \(V\) by homogeneous
light degree \(i\) in coordinates normal to \(\Pi\). Each
coefficient has plane degree \(5-i\), and since the ordinary
conic-ideal order is five, it is divisible by the conic equation
to power \(5-i\) when \(i<5\). Hence it vanishes, and
\(F|_V=G_5(x)\) is base independent. The vertex light projective
space has dimension \(11-3-3=5\). The quintic \(G_5\) is nonzero
since \(V\) has dimension eight and cannot lie in \(X\). It
is smooth: if it had a projective singular point, its join
with \(\Pi\) would be a three-space on which the tangent
hyperplanes of \(X\) all belong to the fixed dual plane of
hyperplanes through \(V\), contrary to finiteness of the Gauss map.
Now take the quadric conic flag, a basepoint-free pair of the
original heavy polars and four general light directions. The
general Hessian-rank estimate for forms with singular codimension
at least five implies Hessian rank at least two on every punctured
graph point of this smooth terminal form. The positive signed
graph denominator is therefore \(3/2\), below \(11/7\).
The radial denominator \(1+(1-1/5)/(1+4/5)=13/9\)
is less than \(b\), and the pole proportions \(1/5,4/5\)
prove the angular target. This completes the conic estimate.

\subsection{Quintic tenfold lines of polar rank at least three and quadratic pencils}
\label{sec:nov27-quintic-line-quadratic}
Set \(d=5,n=10,b=51/35\). The positive graph denominator
for the two-heavy pole flag is \(11/7\), and its required light
weight sum divided by graph weight is \(3/4\).

\paragraph{Ordinary and marked nets.}
At polar rank at least four the positive split degrees sum to
four and there are at least three of them, so their maximum is
at most two. The ordinary angular ratio is at most \(7/5\).
At rank three and balanced splitting \(2,2\), adjunction along
the branch line gives the same angular bound at its marked points;
the off-branch bound \(1+1/h-1/5\) and ternary axis bound
\((1+h(3-1/5))/(1+2h)\) are also smaller than \(b\).
For the splitting \(1,3\) and orders \(h=4,5\), the three terms
of the line product test in Section~\ref{sec:nov27-next-lines}
are \(19/15,704/525,512/399\) and \(6/5,5/4,13/10\), respectively.
For \(h=2,3\) let \(P\) be the moving point on \(L\) and use the
marked pole specialization of
Section~\ref{sec:nov27-sextic-twelve-conics}, now on the line.
The regular product has base degree \(p\), plane degree \(k\),
and mass \(h/5\). On the constant line through the moving graph
the primitive binary graph equation costs one plane degree
and one base degree. Its normalized mean order on that line
is at most \(1/2\), by peeling that equation in the full binary
product just as in the conic flag. Thus the conditional means satisfy
\[
        S_{\rm reg}+\tfrac12v(\mathfrak m_P)\le A,\qquad
        S_{\rm pos}+\tfrac1h v(\mathfrak m_P)\le(1+1/h)A .
\]
For the second estimate the leading pole residues are taken root
by root and jet by jet on the divisor of the base coefficient of
degree \(5-h\). They give complete plane systems of conditional
normalized degree three, plus \(1/h\) copies of the specializing
line through the specializing point. The pair consisting of that
line and the marked point ideal, both with coefficient one,
is log canonical on the plane by blowing up the point.
The argument of the cited section applies with base degrees \(p\).

Normalize \(v(\mathfrak m_P)=1\), set \(u=v(L)\ge1\), and use
\(A\ge u+1\). Averaging gives \(S\le K A-C\), with
\((K,C)=(13/10,1/2),(17/15,13/30)\) at \(h=2,3\), respectively.
The residual bound is \(S\le u/h+(A-1)\min(M,I/u)\) with
\(I=8/5\) and \(M\le I\). We need \(S\le b(A-(1-1/h)u)\).
At \(A=u+1\) the residual bound has slack at least
\(b+(b-1)/h-I>0\). Thus for \(u\ge I/b\) it suffices at all
\(A\ge u+1\), by its slope. For \(1\le u\le I/b\) the marked
bound suffices as soon as
\[
             A\ge A_*(u)=\frac{b(1-1/h)u-C}{b-K}.
\]
The affine function \(A_*(u)-u-1\) is increasing for each
\(h\), with values at \(I/b\) equal to \(-106/561,-29/204\)
respectively. This covers the remaining centers of the net.

\paragraph{Quadratic initial term in a polar pencil.}
Keep the full original light kernel of ambient codimension two.
A quadratic of rank at least two remains nondegenerate on two
general light directions; its signed graph denominator is \(3/2\)
and its radial light ratio is \(2/2=1\).
A rank-one quadratic is a scalar times the square of a rational
linear form over the function field of the line. A primitive
equation has base degree at most one, since its square divides
a form with base degree three. As in
Section~\ref{sec:nov27-divisible-second}, at most one preliminary
constant linear cut makes that equation constant while preserving
the quadratic. The hyperplane kernel \(V_0\) of the resulting
equation \(x=0\) has ambient codimension \(c\le4\).
Its dimension is at least seven, so \(V_0\) is not contained
in the smooth quintic tenfold. Its first light order is
\(t\in\{3,4,5\}\).

At \(t=3,4\) keep a general kernel light \(y\) and a complement
\(x\). For \(t=3\) the weights \(3,2,6\) on \(x,y,z\) give
the initial light part \(x^2+C y^3\), \(C\ne0\), up to rescaling.
At \(t=4\) the weights \(2,1,4\) give
\(x^2+Bxy^2+Cy^4\), \(C\ne0\). The normalized radial ratios
are \(5/6\) and \(3/4\). The sharp cubic and quartic graph
bounds proved above are \(17/11\) and \(11/7\).
At \(t=5\) restrict directly to \(V_0\) together with the two
basepoint-free heavy directions. On \(V_0\) the terminal degree-five
form \(G_5\) is base independent. The finite Gauss map bounds the
projective singular codimension of \(G_5\) by at least
\(12-2c\ge4\), so four general lights cut a smooth quintic in
\(\mathbf P^3\). The signed bound \(11/7\) at all punctured
graph points of a smooth projective quintic surface is proved
above. The normalized radial ratio \(4/5>3/4\) suffices.
In each case the regular and positive two-heavy pole masses
\(1/5,4/5\) give denominator at most \(b\).

\subsection{The higher orders of a quintic tenfold line pencil}
\label{sec:nov27-quintic-line-high}
Let \(C\) be a line whose restricted polars span a pencil on a
smooth quintic tenfold. Use the full constant polar kernel, of
ambient codimension two, and put \(T=4/3\). If its first order
\(h\) is at least three, then \(\delta_C(X,H)\ge385/51\).
Indeed it suffices, at each flag, that the sum of light weights
divided by the heavy weight is at least \(3/4\), and that every
positive pole graph has signed denominator at most \(11/7\):
the normalized line denominator is then at most
\(1/5+(4/5)(11/7)=51/35\).

For a constant kernel of ambient codimension \(c\le5\), the
next nonzero order exists since its projective dimension is
\(11-c>5\). At terminal order five the Gauss-map argument gives
singular codimension at least \(12-2c\) on its projective light
space. For \(c\le4\) take four general lights, which cut a smooth
quintic surface; the preceding smooth-quintic graph calculation
and \(4/5\ge3/4\) suffice. At \(c=5\) the terminal quintic is
still reduced.

We record the elementary graph tests used below. A Hessian of
rank at least two gives \(3/2\). In rank one, if the restriction
of the graph polynomial minus its affine tangent to the Hessian
kernel has a nonzero term of order at most four, take a general
direction there and one complementary direction. The first face
is a square plus a cube, or \(X^2+BXY^2+CY^4\) with \(C\ne0\),
up to rescaling. Their signed bounds are \(17/11,11/7\).
The second bound also holds if \(C=0,B\ne0\).
For a nonpure homogeneous cubic or quartic at a noncritical
point, the rank-one condition identifies the Hessian kernel with
the gradient kernel, and the graph is covered by these tests:
on the plane of the radial vector and a general kernel direction
its expansion is \(\sum c_i(1+X)^{h-i}Y^i\), with
\(c_0\ne0,c_1=c_2=0\) and some \(c_i\ne0,3\le i\le h\);
otherwise the form is a pure power. Affine sections mentioned in
these tests lift by complete hyperplane flag adjunction.

\paragraph{Quartic order in the original kernel.}
The quartic \(G\) has base degree one and \(c=2\).
A constant hyperplane factor costs one cut and gives the terminal
case. A nonconstant rational hyperplane factor costs at most two,
by imposing the two linear coefficients of its base-linear
equation. The same conclusion holds for a nontrivial orbit of
hyperplane factors. Its product is either constant or defines a
pencil with no fixed factor. In the latter case consider the
irreducible curve parametrizing the factors. Two general points
over distinct base points give independent hyperplanes whose
codimension-two intersection lies in the pencil base. It must be
one of finitely many irreducible components of this dimension.
Irreducibility of the general pairs forces this intersection to
be constant, so every moving hyperplane contains it. Cutting there
raises the order to five at \(c=4\).

When \(G\) is irreducible over the geometric generic base, its
general plane section is irreducible by Bertini irreducibility,
and has arithmetic genus three. Its multiplicity-three locus
contains at most one point, since an integral plane triple point
contributes at least three to the genus defect. If a general plane
contains a point of multiplicity at least three, the full
multiplicity-three locus has a unique linear component \(W\) of
codimension two. Indeed a general plane meets each codimension-two
component transversely at its degree-many distinct general points,
with the same multiplicity; components of higher codimension are
avoided. This equality of multiplicities follows by generic
restriction of initial normal forms on smooth equimultiple
strata. The component is constant: for a varying irreducible
pencil, its general member is smooth away from the fixed
codimension-two base; if the family is constant the assertion
is immediate. Cut \(W\) to reach the terminal case.
Otherwise the general plane quartic is integral with all
multiplicities at most two. If the full form is a product of two
distinct irreducible quadrics, a general plane likewise gives
two smooth conics and has all multiplicities at most two.
Use its three lights, with radial ratio \(3/4\). Off the critical
locus apply the homogeneous graph test. At a singular point pass
to the affine light plane by fixing a nonzero coordinate. In
multiplicity two either the Hessian has rank two or its kernel
tangent line is not a component of the quartic, giving pure order
at most four there. Finally a double irreducible quadric yields
a smooth conic squared on this section; at the zeros the conic
flag gives \(11/7\). This covers every quartic order at \(c=2\).

\paragraph{A reduced cubic.}
Three general lights give a reduced plane cubic. Unless it
has a triple point, its local multiplicity-two germs in an affine
chart have rank two, or a double tangent. In the latter case an
irreducible cubic has a pure cubic term along the tangent line.
If that line is a component, the curve is a line tangent to a
smooth conic: the multiplicity-two germ then has first face
\(X^2+bXY^2\), \(b\ne0\). Three reduced line components with no
triple point have only nodes. The homogeneous noncritical test
covers the rest of the punctured graph.

The full triple-point locus of a cubic is a linear subspace:
its vectors form the kernel of contraction into the symmetric
cubic polarization. If a general plane has a triple point, the
full cubic is therefore a binary form in two independent normals
(a single normal is excluded by reducedness). The general plane
detects the full multiplicities as above. Write its normal
subbundle over the base as
\(E=\mathcal O(-a)\oplus\mathcal O(-b)\), \(0\le a\le b\).
The cubic is a section of \(\operatorname{Sym}^3 E(2)\).
Thus \(a=0\); for \(b\ge1\) the cube of the second generator
is absent and the first, constant normal divides the cubic.
Moreover reducedness implies that the term with the first normal
to exponent one is present, so \(2b\le2\). In this case use
the following direct alternative, which also covers \(b=0\):
make both normals constant by at most one preliminary cut if
\(b=1\), keeping a nonzero reduced binary cubic whenever the
two restricted normal directions stay distinct. This is possible
by choosing the restriction of the varying second normal to
have a nonzero component independent of the fixed first normal.
Indeed the span of the coefficient vectors of its base-linear
equation has dimension at most two; impose one generic member of
this span, chosen so the remaining member is independent of the
fixed first normal. If their span contains that fixed normal,
both original normals already lie in a constant two-space, so
their Grassmannian point is constant and \(b=0\) in the intrinsic
normal bundle. The common kernel consequently has \(c\le5\).

Retain both complementary directions \(X_1,X_2\) and one general
kernel light \(Y\). The restricted base cubic \(G_3(X)\) has
three distinct roots at each general pole. If the next order is
four, weights \((4,4,3)\) give \(G_3(X)+aY^4\), \(a\ne0\).
At kernel order five, weights \((5,5,3)\) give \(G_3(X)+aY^5\)
unless a quartic term linear in \(X\) is present, in which case
weights \((3,3,2)\) give \(G_3(X)+\ell(X)Y^3\), \(\ell\ne0\).
At \(X\ne0\) fix \(Y\); the remaining binary polynomial is
the squarefree cubic plus an affine term, with the homogeneous
noncritical cubic bound since its gradient before that affine
addition vanishes only at zero. At \(X=0,Y\ne0\), in the
disjoint faces the radial \(Y\)-direction gives a nonzero square
and a general line in \(X\) gives the cube. In the last face the
mixed derivatives with \(Y\) give Hessian rank two. The radial
ratios of the three faces are \(11/12,13/15,8/9\).

\paragraph{A cube or a fixed double hyperplane.}
For a pure cubic its root \(x\) is fixed by the base degree
bound \(3e\le2\). For a double root consider first a fixed root
\(x\), the cubic then being \(x^2L\). The kernel \(K=(x=0)\)
has \(c=3\); if its next order is five use its terminal test,
so write \(g\) for its nonzero quartic. Any fixed hyperplane
factor of \(g\) yields a terminal cut at \(c=4\). In particular
we may suppose no geometric factor of \(g\) has multiplicity
at least three: it would be a hyperplane, and its Galois orbit
would have primitive base degree zero by the base degree one.

For the cubic \(x^3\) use two general \(K\)-lights and weights
\((4,3,3)\). The face is \(x^3+g(Y_1,Y_2)\), with all binary
root multiplicities at most two. At \(x\ne0,Y\ne0\) its Hessian
has rank two, since \(g\) has a nonzero Hessian at such \(Y\)
(Euler and the multiplicity bound). At \(Y=0,x\ne0\) take the
radial quadratic and a general quartic kernel direction.
At \(x=0,Y\ne0\), fix \(x\) at a noncritical point of \(g\);
at a double root use \(x^3\) and the transverse square on an
affine slice. The radial ratio is \(10/12\).

For \(x^2L\) with \(\ell=L|_K\ne0\), use weights \((3,2,2)\);
the face is \(x^2\ell(Y)+g(Y)\). At \(\ell\ne0\) the Hessian
has a nonzero \(xx\) entry. If it has rank one its normal
has a nonzero \(x\)-coefficient. The degree-four term \(g(Y)\)
cannot vanish on that kernel, which projects onto all \(Y\),
so the rank-one quartic test applies. At \(\ell=0,x\ne0\) the
mixed derivatives and zero \(xx\) entry give rank two. At
\(x=0,Y\ne0,\ell=0\), a general binary restriction can encounter
a multiple root of \(g\) only if \(\ell\) is a repeated factor
of the full quartic; this is fixed by its base degree and was
excluded. Otherwise the noncritical homogeneous quartic estimate
on \(x=0\) suffices. The radial ratio is \(7/8\).

\paragraph{A moving double hyperplane.}
The cubic has the form \(zL^2\) with \(z\) fixed, since the
moving root uses the whole base degree two. On \(K=(z=0)\)
the terminal case is already settled. Let \(g\) be the quartic
there and put \(\ell=L|_K\), a nonzero linear form at a general
pole. A fixed hyperplane factor of \(g\) again suffices. Giving
\(z\) weight two and the kernel lights weight one gives
\(z\ell(Y)^2+g(Y)\). Use two general kernel lights when
\(\ell\nmid g\), and three in the remaining case. The radial
ratio is at least one. If \(\ell\ne0\), the zero \(zz\) entry
and nonzero mixed derivatives give Hessian rank two.

At \(Y\ne0,\ell(Y)=0\), suppose first \(g(Y)\ne0\) or
\(d g_Y\ne0\). If the Hessian has rank one, its kernel contains
the \(z\)-direction but is not equal to \(\ker d\ell\).
Indeed Euler on \(Y\) gives a nonzero radial quadratic if
\(g(Y)\ne0\), and its contraction gives \(3d g_Y\ne0\) in
the second case; in both cases \(d\ell(Y)=0\).
Thus a general vector in the Hessian kernel has a nonzero
\(\ell\)-component; varying its \(z\)-component gives a nonzero
pure kernel cubic from the summand \(z\ell^2\). At \(Y=0,z\ne0\)
when \(\ell\nmid g\), fix \(z\) to obtain the square and a
nonzero pure quartic on its kernel. These arguments handle the
whole two-kernel-light case, since generically \(g\) is nonzero
on the one projective zero of \(\ell\) in the chosen binary space.

If \(\ell\mid g\), we may assume \(\ell\) is moving, since a
fixed factor is already settled. Its base degree is one and
\(g=\ell R\) for a fixed cubic \(R\). If \(R\) is reducible,
it has a fixed hyperplane factor. Assume it irreducible.
Its restriction to the general hyperplane of the pencil of
\(\ell\)'s is not a triple hyperplane. Indeed the general pencil
member cuts the smooth locus of \((R=0)\) transversely outside
the fixed pencil base by Bertini, and contains no divisorial
component of its singular locus outside that base. If the general
section were a triple hyperplane it would therefore be supported
in that linear codimension-two base. On each pencil hyperplane
the cubic would be a scalar multiple of the cube of a pencil
coordinate; varying the hyperplane shows \(R\) depends only on
the two pencil coordinates, contrary to irreducibility. Thus two
general directions on \(\ell=0\) give a binary cubic \(r\) with
root multiplicities at most two. Three general fixed \(K\)-lights
induce precisely such a section at general pole roots.

At \(Y=0,z\ne0\) set \(z\) to its nonzero value and give
\(X=\ell\) weight three and its other two lights weight one.
The face is \(aX^2+Xr(U,V)\), \(a\ne0\), with radial ratio
\(5/6\). Its \(XX\) derivative is nonzero. At \(dr\ne0\), in
rank one its Hessian kernel is \(2a\,dX+dr=0\); the highest
term \(Xr(U,V)\) restricts nontrivially there, so the rank-one
quartic test applies. At \((U,V)=0,X\ne0\), the pure kernel
cubic is \(Xr\) with \(X\) evaluated at the point. At a nonzero
multiple root of \(r\), fixing a nonzero coordinate gives a
double zero transverse to the radial direction. For \(X\ne0\)
this gives rank two, while for \(X=0\) the first transverse
face is \(aX^2+b Xw^2\), \(b\ne0\).

It remains on the original three-light face to consider
\(Y\ne0,\ell(Y)=0,R(Y)=0\) with \(dg_Y=0\). Dehomogenize
by fixing a nonzero \(Y\)-coordinate and put \(X=\ell\);
write \(w\) for the other transverse coordinate. The local light
part is \(X^2z+X R(X,w)\). At a simple root of \(R(0,w)\)
its mixed Hessian has rank two. At a double root, a nonzero
quadratic term gives on fixing \(z\) a face \(aX^2+bXw^2\).
When the quadratic coefficient vanishes, retain the variation
of \(z\): the ordinary initial cubic, after a linear change, is
\(X^2 Z+bXw^2\), \(b\ne0\). It is a line times a smooth conic
tangent to that line in the projective light plane. Its radial
ratio is one, and the reduced-plane-cubic argument above gives
the required punctured bound. This completes all cubic cases,
and hence all higher pencil orders.

\section{\texorpdfstring{Conclusion in dimension \(2d\)}{Conclusion in dimension 2d}}
\begin{theorem}\label{thm:nov27-second-row}
For every \(d\ge5\), every smooth complex hypersurface \(X_d^{2d}\subset
\mathbf P^{2d+1}\) satisfies \(\delta(X,H)>d+2\) and is uniformly K-stable.
\end{theorem}
\begin{proof}
For \(d\ge20\) the bounds in Sections~\ref{sec:nov27-all-conics-2d},
\ref{sec:nov27-next-lines} and~\ref{sec:nov27-all-pencils-2d} imply
\(\delta_C(X,H)\ge (2d+1)(d+3)/(2d+2)\) for every line or integral conic.
For \(d\ge20\) Section~\ref{sec:nov-second-row} applies with this bound.
For \(9\le d\le19\), Sections~\ref{sec:nov27-integer-second}
and~\ref{sec:nov27-divisible-second} verify the point criterion using
their stronger estimates on curves.
For \(d=8\) Section~\ref{sec:nov27-octic-lines} supplies the
common curve bound \(136/13\) and excludes the remaining point centers.
For \(d=7\), Section~\ref{sec:nov27-septic-lines} gives the common
curve bound \(465/49\). Its point comparison uses a general septic
thirteenfold as a hyperplane section and gives \(\delta(X,H)>9\).
For \(d=6\), Sections~\ref{sec:nov27-sextic-line-nets}
and~\ref{sec:nov27-sextic-twelve-conics} give
\(\delta_C(X,H)\ge247/29\) on every line and conic. The point
calculation requires only \(1599/188\). Their difference
is \(65/5452>0\), so \(\delta(X,H)>8\) in this case.
For \(d=5\), the conic estimate and the line estimates for higher
polar rank and for pencils in Section~\ref{sec:nov27-quintic-ten}
give the common bound \(385/51\). The strict endpoint comparisons
of that section yield \(\delta(X,H)>7\).
Smoothness along a line forces the first polars to have rank at least
two, since a single nonconstant binary form has a zero. Thus the point
criterion applies to every smooth member. The equality
\(-K_X=(d+2)H\) gives the claimed anticanonical threshold.
\end{proof}

\chapter{Sufficient conditions in larger dimensions}\label{ch:beyond}
\section{Induction using pointed hyperplane sections}
For dimensions above the two preceding boundaries, the local
estimates impose additional conditions on the hypersurface.
In this chapter \(H=\mathcal O_X(1)\).

\subsection{A short-curve criterion on a strip of growing width}
Suppose \(d\ge 10^4\) and \(n=2d+k\), where
\(0\le k\le\sqrt d/8\). Write \(N=n+1\), \(r=n+2-d\), and put
\[
               M_{d,k}=\frac{(n+1)(d+k+3)}{2d+k+2}
                       =r+\frac{d-1}{n+2}.
\]
Assume a general smooth hyperplane section through each closed
point of \(X\) is K-stable. If every line and every integral
conic \(C\subset X\) has \(\delta_C(X,H)\ge M_{d,k}\), then
\(X\) is uniformly K-stable. In particular, if the local bound
holds for every smooth member successively in these dimensions,
the hyperplane hypothesis can be propagated by induction starting
with the theorem in dimension \(2d-1\).

Suppose \(\beta=\delta(X,H)\le r\). Take a divisorial optimizer.
Hyperplane flag adjunction through a point of its center gives
\(\beta>N(r-1)/n>N/3\). Apply the upper \(1-1/d\) projection
tilt and free-hyperplane subtraction of
Section~\ref{boundary:numbers}, with
\(u=n/N,\ b=r/N,\ \epsilon=r^{-4}\). The center and point estimates
there are valid in arbitrary dimension once the initial strict lower
bound on \(\beta\) has been established. Their required numerical inputs are
\[
 \begin{gathered}
 x=(r-N/3-3\epsilon)/d,\qquad
 U=r+bx+2\epsilon,\qquad
 z=b+x\{u(1-b)-\log(d)/N\}-u\epsilon,\\
 U+\Delta<M_{d,k},\qquad z-u\Delta>1/3,\qquad
 z+\Delta\big((n+q)/U-u\big)>1 ,
 \end{gathered}
\]
where the last inequality need only hold at any hypothetical point
optimizer with \(A_X=n+q\le U S_H\) and \(v(\mathfrak m_p)=1\).
We verify them uniformly. Put \(y=\sqrt d\); thus \(y\ge100\).
Then \(b=1/2+(k+3)/(2N)\) and
\(1/3<x<1/3+(2k+5)/(3d)<0.336\). Hence
\[
                U<r+1/6+0.43k/d+1.1/d .
\]
Indeed the contribution of \(b-1/2\) is at most
\(0.084(k+3)/d\), and the additional \(2\epsilon\) fits in the stated
strict bound. Next consider \(b+xu(1-b)\), increasing with \(b\)
and \(x\) in the indicated range. Its minimum there is at least
\(1/2+u/6\). Consequently
\[
 z>2/3-\{1/12+0.168\log d\}/d-\epsilon
     >2/3-0.017/y .
\]
For the last comparison divide the numerator of the error by \(y\);
\(\log(d)/y=2\log(y)/y\) decreases for \(y\ge100\).

Now take
\[
      \Delta=\frac13-\frac{0.69k+2.2}{d}-\frac{0.017}{y}
               >\frac13-\frac{0.126}{y}.
\]
In fact \((d-1)/(n+2)=1/2-(k+4)/(2(n+2))
\ge 1/2-k/(4d)-1/d\), so the estimates on \(U,z\) imply the
first two strict comparisons. It remains to check the point.
The universal point envelope is \(S_H\le J+Gq/n\), with
\[
 \begin{aligned}
 J&=\max\left\{\frac{3(\sqrt{2d}+4)+n-3}{N},
               \frac{n-2+d+\sqrt{d/2}}{N}\right\},\\
 G&=\max\{\sqrt{2d}+4,\sqrt{d/2}+\log(3n(d+2))\}.
 \end{aligned}
\]
In this strip \(J<1.504\), \(G<1.46y\) and
\(U<d+k+2.17\le1.002717d\). For the log bound use
\(3n(d+2)<7d^2\) and the decreasing function
\(\log(7y^4)/y\). In particular \(n-UJ>0.491d\),
\(UG/n-1<0.733y\). Unless the point is already excluded,
\(A_X\le US_H\) gives \(q>0.669y\). Since \(k\le y/8\),
\[
 n+q-2U>0.669y-k-4.34>0.5y,\qquad
                         (n+q)/U>2+0.49/y.
\]
The final strict point comparison now follows from
\[
 z+\Delta\big((n+q)/U-u\big)>
 (2/3-0.017/y)+(1/3-0.126/y)(1+0.49/y).
\]
The right side exceeds one since $0.49/3-0.017-0.126>0.02$
and $0.126(0.49)/y<0.001$ for $y\ge100$. This excludes the
point as well and completes the criterion.

\subsection{A wider strip from the two point estimates}\label{beyond:wider-strip}
The short-curve criterion in the preceding subsection holds for
\[
 d\ge10^4,\qquad n=2d+k,\qquad 0\le k\le \tfrac12\sqrt d .
\]
More generally, for every fixed \(c<2\sqrt2/3\) it holds for
\(0\le k\le c\sqrt d\) and \(d\) sufficiently large in terms of \(c\).
It suffices to assume \(\delta_C(X,H)\ge M_{d,k}\) on every line and
integral conic of \(X\).

\begin{proof}
First we give the initial lower bound used in the tilt.
For any \(n\ge m=2d-1\), a general smooth pointed linear section
of dimension \(m\) has threshold \(>d+1\) in \(H\)-units by the
all-member theorem in that dimension. At every intermediate dimension
\(j\), hyperplane flag adjunction at the chosen point \(p\) gives
\[
 \frac{\delta_p(X_j,H)}{j+1}\ge
 \min\left\{1,\frac{\delta_p(X_{j-1},H)}{j}\right\}.
\]
Equivalently one may induct on \(j\) using the global threshold on a
smooth pointed section for the second term. To keep strictness when
taking the global threshold, apply this flag at a point of the center
of a divisor computing any putative \(\beta=\delta(X,H)\le r\).
The result is
\[
                  \beta>\frac{N(d+1)}{2d}>\frac N2 .
\]
Thus all intermediate sections can be used without assumptions on
their stability. In the tilt put
\(x_\beta=(\beta-N/3-3\epsilon)/d\) and
\(\lambda=(d-1)x_\beta\). We have \(x_\beta>N/(6d)-3\epsilon/d\).
For \(d\ge10^4,\ k\le\sqrt d\), this gives \(x_\beta>0.33\) and
\(1-b>0.49\). In particular \(\beta'>N/3\), and the upper estimate
\(\beta'\le\beta-\lambda-x_\beta(1-b)+\epsilon\)
lies strictly below \(r-\lambda c\) for the residual class factor
\(c\le1\) of Section~\ref{boundary:numbers}. The endpoint substitutions
in \(U,z\) use only \(\beta\le r\) and \(0<u(1-b)-\log(d)/N<ud\).
This supplies the same tilt and center comparisons as before.

Retain \(N,r,b,u,\epsilon,x,U,z\) from the preceding subsection and
put \(y=\sqrt d\). First assume \(y\ge100\) and \(k\le y/2\). We have
\[
 \tfrac13+\tfrac{2k}{3d}\le x<
 \tfrac13+\tfrac{2k+5}{3d}<0.337,\qquad
 b=\tfrac12+\tfrac{k+3}{2N}.
\]
Consequently \(U<r+1/6+0.43k/d+1.1/d<r+0.17\).
Indeed \(x(b-1/2)\le0.08425(k+3)/d\), and
\(x/2\le1/6+k/(3d)+5/(6d)\). In particular,
\(U/d<1.00523\) and \(U/n<0.503\). Furthermore
\[
                 z>\tfrac23+0.49k/d-0.017/y .
\]
To see this, write \(b+xu(1-b)=1/2+xu/2+(b-1/2)(1-xu)\).
The first two terms are at least \(2/3-1/(6N)+uk/(3d)\);
the last term exceeds \(0.663k/(4.011d)\) when \(k>0\).
Use \(u>0.9999\) and
\(1/(6N)+0.337\log(d)/N+\epsilon<0.017/y\).
The latter follows by evaluating at \(y=100\) after multiplication
by \(y\), since \(\log(y)/y\), \(1/y\) and \(y^{-7}\) decrease there.

Take \(\Delta=1/3-0.69k/d-0.022/y>0\).
The identity \(M_{d,k}-r=1/2-(k+4)/(2(n+2))\)
gives \(M_{d,k}-U>1/3-0.68k/d-2.1/d>\Delta\).
The estimate for \(z\) gives \(z-u\Delta>1/3\).
It remains to prove the last point comparison. In the
dimension-independent point estimate, the bounds
\(S_H(v)\le J_i+G_iq/n\) apply with the following \emph{paired}
intercepts and slopes:
\[
 \begin{array}{c|cc}
 &J_i&G_i\\ \hline
 \tau_p(H)\le\sqrt{2d}+4&
 (n-3+3(\sqrt{2d}+4))/N&\sqrt{2d}+4\\
 \tau_p(H)>\sqrt{2d}+4&
 (n-2+d+\sqrt{d/2})/N&\sqrt{d/2}+\log(3n(d+2)).
 \end{array}
\]
Indeed the first pair is the pointed quadric bound and the second
pair follows from the principal-bin bound in its proof.
In the first case \(J_i<1.022\) and \(G_i<1.455y\);
in the second \(J_i<1.504\) and \(G_i<0.912y\).
For the latter slope use \(3n(d+2)<7d^2\) and the decreasing
function \(\log(7y^4)/y\) on \(y\ge100\).
The hypothetical inequality \(n+q\le US_H(v)\) would give
\[
                  q(UG_i/n-1)\ge n-UJ_i .
\]
In the first case the right side exceeds \(0.972d\) and the
coefficient on the left is less than \(0.732y\).
In the second case these bounds become \(0.488d\) and \(0.459y\).
Thus \(q>1.06y\) in both cases. It follows that
\(n+q-2U>1.06y-k-4.34>0.516y\), whence
\((n+q)/U-u>1+0.51/y\). Since
\(\Delta\ge 1/3-0.367/y\), the final comparison is
\[
 z+\Delta\big((n+q)/U-u\big)>
 1+\frac{0.51/3-0.017-0.022}{y}
       -\frac{0.20k}{d}-\frac{0.367\cdot0.51}{y^2}>1 .
\]
This proves the explicit range.

For the asymptotic assertion all errors below are uniform for
\(0\le k\le y\). Expanding the displayed formulas for \(x,b,U,z\)
and \(M_{d,k}\), using \(k^2/d^2\le1/d\), gives
\[
 \begin{aligned}
 U-r&=\tfrac16+\tfrac{5k}{12d}+O(d^{-1}),&
 M_{d,k}-r&=\tfrac12-\tfrac{k}{4d}+O(d^{-1}),\\
 z&=\tfrac23+\tfrac{k}{2d}+O((\log d)/d).
 \end{aligned}
\]
In particular
\(\Delta=1/3-2k/(3d)-C(\log d)/d\) is positive and gives
\(U+\Delta<M_{d,k}\) and \(z-u\Delta>1/3\) for an absolute, sufficiently
large \(C\). In the first point case
\(n-UJ_i=d+O(y)\) and \(UG_i/n-1=y/\sqrt2+O(1)\).
In the second these quantities are \(d/2+O(y)\) and
\(y/(2\sqrt2)+O(\log d)\), respectively. Division in either case
gives \(q\ge\sqrt2\,y-O(\log d)\). Substituting this lower bound
for \(q\), expanding \(U=d+k+O(1)\), yields
\[
 \frac{n+q}{U}-u\ge 1+\frac{\sqrt2\,y-k}{d}
                                  -O((\log d)/d).
\]
Consequently
\[
 z+\Delta\big((n+q)/U-u\big)
 \ge 1+\frac{\sqrt2\,y/3-k/2}{d}-O((\log d)/d).
\]
The right side exceeds one for all sufficiently large \(d\) when
\(k\le c y\) with \(c<2\sqrt2/3\), as required.
\end{proof}

\paragraph{The first higher dimension.}
For every \(d\ge40\), the following implication in dimension
\(n=2d+1\) holds: if every line and every integral conic satisfies \(\delta_C(X,H)\ge M=(2d+2)(d+4)/(2d+3)\), then \(X\) is uniformly K-stable.
Since \(n-1=2d\), the smooth pointed hyperplane section belongs
to the range already proved. For \(d\ge10^4\) this is the first
positive width in the preceding strip. Here are the finite
comparisons for the remaining degrees. For \(230\le d\le9999\)
use the coarse point estimate and \(\epsilon=r^{-4}\) from the
preceding proof, taking \(\Delta\) to be the smaller of \(M-U\)
and \((z-1/3)/u\) minus \(\epsilon\). For \(40\le d\le229\) take
the finite-fraction tilt with \(a=1-10/d\), \(f=1/d\) and
\(\epsilon=10^{-8}\). Use the point and center definitions of
the dimension-\(2d\) calculation with the new \(n,N,r,M\), and
round \(Q\) upwards on the grid \(10^{-4}\). Define \(B\) by the
same downward rounding on the grid \(10^{-3}\). For
\(\tau_p(H)\ge B\) a single principal-bin comparison on \([B,d]\)
with \(5n\) bins suffices except for \(d=40\), when the two
intervals \([B,38],[38,40]\) use \(4n,5n\) bins, respectively.
All inequalities, including the two affine endpoint assumptions
of the tilt, are checked in \path{nov27-bin-third-row.py}.
The local threshold premise still applies to each line and conic.

\paragraph{Degrees thirty-three through thirty-nine.}
The same point criterion in dimension \(n=2d+1\) holds for
\(33\le d\le39\) at the same short-curve mark
\(M=(2d+2)(d+4)/(2d+3)\). In the preceding finite-fraction
calculation use
\[
\begin{array}{c|rrrrrrr}
 d&33&34&35&36&37&38&39\\ \hline
 1000a&702&710&717&724&731&737&743\\
 1000B&7383&7559&7738&7915&8090&8268&8444
\end{array}
\]
with the same rounding rules for \(B,Q\). The cutoff and the two
signed comparisons below \(B\) are strict. For \(B<t<17/2\) use
the unique prime of class \(jH\) and multiplicity \(k=jt\).
Here \(t<d\), so its general pointed curve section is integral.
Use the discriminant bounds of Section~\ref{sec:nov-second-row}.
For \(1\le j<20\) take \(I=dj(j+d-4)+k+1\), taking its minimum
with the sectional Milnor bound when \(j\le2\). For \(j=3,\ k>9\)
also use the cubic argument of Section~\ref{sec:nov27-cubic-degree19}.
In the smooth case its bound is
\((k-1)(1+[3(d+1)^n/(k-1)]^{1/(n-1)})\); in the singular
case it is
\((k-2)(1+[4(2d+2)^{n-1}/(k-2)]^{1/(n-2)})+2d\).
Indeed \(k-1>8\) forces the same codimension-two singular locus
and the incidence-chart Jacobian degrees are at most \(2d+2\).
Take the maximum of the two bounds before intersecting with the
genus bound. Discard
\(k(k-1)>I\). The test
\[
 \frac d t+\left(1-\frac d{t^2}\right)
 \frac{Im}{j(I+m(m-1))}<B,\qquad m=\min(k,\sqrt I),
\]
holds except for \((j,k)=(4,33)\) at \(d=33\), where the
single-value principal-bin test with \(3n\) bins applies.
For \(j\ge20\) put
\[
 E=d+\frac{d(d-4)+17/2}{20}+\frac1{400},\qquad
 L_* =\frac{E}{2\sqrt E-1/20}.
\]
Evaluation at \(2L_*\) clamped to \([B,17/2]\) gives
\(d/t+(1-d/t^2)L_*<B\) throughout. Finally, the intervals
\([17/2,d/2]\) and \([d/2,d]\) pass the principal-bin test with
\(3n\) and \(10n\) bins respectively. These comparisons and all
higher-rank tests in the same seven degrees are reproduced by
\path{nov27-bin-third-small.py}.

\subsection{The codimension of the high-contact exceptions}
\label{sec:nov27-high-contact-strip}
Let \(d\ge10^4\), \(0\le k\le\sqrt d/2\), \(n=2d+k\). Set
\[
 t=\lfloor d/(k+3)\rfloor,\qquad
 c_1=k+4,\quad c_2=2k+7,\qquad m_i=n+1-c_i.
\]
In the projective parameter space of degree-\(d\) hypersurfaces,
every smooth member outside a closed subset of codimension at least
\[
 \min\{D_1-c_1(m_1+1)-2(m_1-1),\quad
           D_2-c_2(m_2+1)-(3m_2-1)\}
\]
is uniformly K-stable, where
\[
 \begin{aligned}
 D_1&=\sum_{a=0}^{t}\binom{m_1-2+a}{a}(d-a+1),\\
 D_2&=\sum_{a+j\le t}\binom{m_2-3+a}{a}
                             \bigl(2(d-a-2j)+1\bigr)
                      \quad(a,j\ge0).
 \end{aligned}
\]
In particular the codimension is at least
\(3^{\lfloor\sqrt d\rfloor}-(2k+10)(n+2)\) on the whole strip.
The same statement holds for each fixed \(0<c<2\sqrt2/3\)
when \(0\le k\le c\sqrt d\) and \(d\) is sufficiently large in
terms of \(c\), using Section~\ref{beyond:wider-strip} as the point
input. In both ranges \(k+3\le\sqrt d\) and \(m_2-2\ge d\).

The two exceptional sets are the closures of the loci where
there is a line, respectively a smooth conic, \(C\) in a linear
\(W=\mathbf P^{m_i}\) such that the restricted equation of \(X\)
belongs to \(I_{C/W}^{t+1}\). We verify their sufficiency for
all the pointed hyperplane sections in the preceding strip criterion.
Fix a smooth linear section \(Y\) of ambient codimension \(j\),
where \(0\le j\le k\), and put \(k'=k-j\). It has dimension
\(2d+k'\). For a line of polar rank at least \(k'+5\), or a
conic of polar rank at least \(2k'+8\) whose plane restriction has
multiplicity at least two, use the preceding ordinary-rank bounds.
A conic of multiplicity one on its plane is similarly covered.
In every outstanding case the common polar kernel \(V\) through
the span of \(C\) has ambient codimension at most \(j+k'+4=c_1\)
for a line and at most \(j+2k'+7\le c_2\) for a conic.
In the conic case the assumption on plane multiplicity includes
the identically zero plane equation. Shrinking \(V\) if necessary
gives a \(W\) of the appropriate dimension containing \(C\).
Outside the indicated exceptional loci the first curve-ideal
order \(h\) of the equation on \(V\) is consequently at most
\(t\), since restriction preserves membership in ideal powers.
It is at least two by the definition of the common polar kernel.

Retain a basepoint-free pencil in the active polar system and take
the linear section spanned by \(V\) and its two chosen normal
directions. It and the successive linear cuts are integral and
smooth along \(C\): the retained binary sections have no common
zero there, and every extra projective hypersurface component would
meet the curve and give a singularity. Thus the same complete
hyperplane flags and localized adjunction as in Chapter~\ref{ch:models}
apply. For a line use the one-light calculation with this order
\(h\). For a conic the light variables are the linear normals in
\(V\) and the symbol \(w\) of its plane quadratic equation. The
leading relation is \(f_1z_1+f_2z_2+G_h(x,w)\) in the ordinary
curve-ideal graded polynomial ring over the conic ring. This
description applies also when the plane restriction is nonzero,
since all higher terms of the plane equation appear within the
curve-ideal filtration. On the punctured conic cone the relation
defines an affine bundle by basepoint freeness, and its vertex
fiber is too small to contribute a component. Hilbert series and
generic localization therefore give exactly the weighted slices.
A homogeneous shear of the degree-two symbol by a base linear
form times a light linear coordinate, followed by a general
linear light direction, makes the pure coefficient of \(x^h\)
nonzero: over the generic base its value tests the nonzero
polynomial \(G_h\) on a Zariski dense set of light directions.
There exist light linear directions because \(m_i>2\). Splitting
the remaining variables now gives the binary core of
Chapter~\ref{ch:models} irrespective of the original polar rank.
For both curves the normalized denominators are bounded by
\[
 \max\left\{1,\ 1+1/h-1/d,\ 
             1+(1-1/d)/(1+1/h)\right\}
       \ \le\ 1+\frac{d-1}{d+k'+3}.
\]
Indeed \(h\le t\le d/(k'+3)\), and the second term is below
\(3/2\) while the right hand side exceeds \(3/2\).
This verifies the required short-curve estimate on every such
smooth section. Starting at the established all-member row
\(n=2d-1\), induction with the pointed strip criterion proves the
assertion for \(X\).

For the dimension count fix \(C\subset W\). The restriction of the
ambient degree-\(d\) equation to \(W\) is surjective. The
associated graded ring for a line ideal has \(m_1-1\) independent
linear normal variables over its binary coordinate ring. Keeping
monomials of normal degree at most \(t\) gives the \(D_1\)
independent equations. For the conic ideal there are \(m_2-2\)
linear generators and one quadratic generator forming a regular
sequence. Its associated graded is polynomial in their symbols
over the conic coordinate ring. Order \(a+j\) with \(a\)
linear symbols and \(j\) quadratic symbols gives the dimension
in \(D_2\), since \(d-a-2j\ge d-2t\ge0\). These ideals have
the same powers in the homogeneous coordinate ring and at the
generic point in the sense needed here: each successive graded
quotient is torsion-free over the curve coordinate ring, so a
nonzero initial symbol remains nonzero over its function field.
Choices of \(W\) have dimension \(c_i(m_i+1)\), and choices
of \(C\) in \(W\) have dimensions \(2(m_1-1)\) and \(3m_2-1\).
The contact conditions for each flag are linear of the calculated
constant ranks. The dimension of the closure of their incidence
images therefore gives the claim. For the final lower bound use
the single summand with \(a=t,j=0\). In each case the number
\(s\) of linear normal generators is at least \(d\), and
\(\binom{s+t-1}{t}\ge(s/t)^t\ge3^t\).
Finally \(k+3\le\sqrt d\) and the two flag parameter dimensions
are at most \((2k+10)(n+2)\).

\subsection{The higher polar ranks on successive dimensions}
\label{sec:nov27-net-strip}
Let \(n=2d+k\), \(k\ge1\), \(K=k+3\), and \(d\ge32K^3\).
Write \(N=n+1\), \(b=1+(d-1)/(d+K)\). Every line and every
integral conic \(C\) of normal-polar rank at least three on a smooth
\(X_d^n\) satisfies
\[
                \delta_C(X,H)\ge N/b=N(d+k+3)/(n+2).
\]
For a conic of simple multiplicity on its plane the estimate holds
without the rank hypothesis. In the other conic cases the normal
polars mean the plane-normal polars restricted to the conic.

Set \(e=\deg C\in\{1,2\}\). A conic of simple plane multiplicity
has denominator below \(3/2<b\) by the surface flag. For the other
cases, as in Section~\ref{sec:nov27-high-contact-strip}, the
common polar kernel contains the entire span of the curve, the
plane equation of the conic supplies an additional light symbol
of section degree two, and proper higher-multiplicity plane
restrictions give the same type of ternary model as contained
planes. By the ordinary estimates we may assume that the rank
is at most \(K+1\) for a line or \(2K+1\) for a conic
(retaining one additional rank would do no harm). Indeed the
ordinary angular bound is \(2-(q-1)/d\) for line polar rank
\(q\), and \(2-q/(2d)\) for conic polar rank \(q\);
for a line \(q\ge K+2\) works, for a conic \(q\ge2K+2\)
works. Both radial denominators are below \(3/2\). The
common kernel consequently has dimension greater than \(n/2\).
Its first curve-ideal order \(h\) is finite. Keep a general
basepoint-free net among the polars and split off the extra heavy
variables. The model of Sections~\ref{sec:nov27-all-conics-2d}
and~\ref{sec:nov27-next-lines} has low split degree \(s\ge1\),
\(s\le e(d-1)/2\). The single-binary bound suffices for
\(h\le d/K\). Suppose \(h>d/K\). Every term away from the
distinguished branch point has denominator at most \(3/2<b\).
The residual intersection and multiplicity bounds \(I,M\) are
\[
 I=2-(e+s)/(ed),\qquad
 M=\min\{I,2-hs/((s+e)d)-1/d\}.
\]
When \(s\ge eK\), ordinary adjunction suffices since \(I<b\).
Otherwise use the first infinitely near point on the strict
transform of the branch. With \(t=(d-1)/(h-1)\) and \(R=ht/d\)
the upper bound for the residual multiplicity \(m_1\) there,
after dropping its nonpositive terms, is
\[
      m_1\le J=R(2et+3s)/(e(d-1)).
\]
These are the same degree-one or degree-two base jets visualized on
the geometric generic plane in the preceding local calculation.
We record how the bound controls directions beyond the first
interval. Normalize the valuation by \(v(\mathfrak m)=1\).
Write \(u=v(L)\ge1\), \(A\) for the ordinary plane discrepancy,
and \(m\le M\) for the multiplicity of the residual \(\Gamma\).
For \(u>1\) the center on the first blowup is the branch point
\(P_1\). There the maximal ideal has value \(\min(1,u-1)\),
since the regular parameters may be chosen with values \(1,u-1\).
Applying the ordinary smooth-surface multiplicity inequality there
to the strict transform of the residual gives
\[
 v(\Gamma)\le m+m_1(A-1-\min(1,u-1)).
\]
At the minimum possible discrepancy \(A=u+1\) the right-hand
side is at most \(M+\max(1,u-1)J\). At all higher discrepancies
its slope is at most \(J\). For \(u=1\) the original bound is
\((A-1)M\). The intersection comparison on the original plane
also gives \((A-1)\min(M,I/u)\) for every \(u\).

The following inequalities are uniform:
\[
 t<1.01K,\quad R<1.01,\quad J<6K/d,\quad M\le2-1/(3K).
\]
Indeed \(t\le K(d-1)/(d-K)\), \(R\le(d-1)/(d-K)\),
\(s<eK\), and \(s/(s+e)\ge1/3\). For \(1\le u\le K+2\) we have
\[
 M+(K+1)J\le 2-1/(3K)+6K(K+1)/d
            <2-(K+1)/(d+K)=b.
\]
For the strict inequality it suffices that
\(d>18K^2(K+1)+3K(K+1)\), a consequence of \(d\ge32K^3\)
and \(K\ge4\). For \(u\ge K+2\) one has
\(I\le b+(b-1)u/h\) even with the crude bounds \(I\le2,h\le d\),
since \((K+2)(d-1)\ge(K+1)d\). The intersection comparison
therefore suffices at \(A=u+1\), and its slope is at most
\(M<b\). In both ranges we obtain at every \(A\ge u+1\)
\[
             v(\Gamma)+\frac{u}{h}
                \le b\,[A-(1-1/h)u].
\]
Weighted flag adjunction completes the local calculation.

\paragraph{Higher polar ranks on a square-root strip.}
The same local conclusion holds for \(d\ge10^4\) and
\(1\le k\le3\sqrt d/10\). It also holds for each fixed
\(c<1/3\) on \(1\le k\le c\sqrt d\) for sufficiently large \(d\).
Set \(x=h/d,\ z=s/e\). The binary case covers \(x\le1/K\),
and ordinary adjunction covers \(z\ge K\). The remaining polar ranks
have common kernels of dimension greater than \(n/2\), since
their ranks are at most \(2K+1\). For \(x>1/K,\ 0<z<K\), set
\[
 U=\frac{(I-b)h}{b-1}<\frac{(K-z)x}{b-1}.
\]
For \(u\ge\max(1,U)\), the original intersection estimate
suffices at \(A=u+1\). For \(1<u<U\), the first blowup gives
\[
     v(\Gamma)\le M+\max(1,u-1)J\le M+UJ .
\]
With \(d\ge10^4\), write \(y=\sqrt d\). Then \(K\le0.33y\),
\(K/d\le0.0033\), and
\[
 R,\ t/K\le (d-1)/(d-K)<1.0034,\qquad
 1/(b-1)<1.0035,\qquad d/(d-1)<1.0002 .
\]
Consequently \(J<1.01(2K+3z)/d<b\) and
\(UJ\le1.012x(K-z)(2K+3z)/d\) when \(U>1\).
Moreover \(2-b-1/d<K/d<xK^2/d\). Thus \(M+UJ\le b\)
in that case follows from
\[
       \frac{K^2+1.012(K-z)(2K+3z)}d\le\frac z{1+z}.
\]
If \(z\ge1\), use \((K-z)(2K+3z)\le25K^2/12\), giving
an upper bound \((1+1.012\cdot25/12)0.33^2<1/2\).
The only remaining value is \(z=1/2\). Its left-hand side is at
most \(3.024\cdot0.33^2+0.506\cdot0.0033<1/3\).
The same inequalities imply \(M<b\) for every \(U\),
so the plane bounds extend to all discrepancies and the original
multiplicity comparison works at \(u=1\).
For the asymptotic statement the factor \(1.012\) becomes
\(1+O(d^{-1/2})\) uniformly. The two estimates become
\((1+25/12)c^2+o(1)<1/2\) and \(3c^2+o(1)<1/3\).

As a consequence, use \(m=n-1\) and \(t=\lfloor d/K\rfloor\).
On these two strips one can take the minimum of the two
high-contact codimension bounds of
Section~\ref{sec:nov27-third-rank-two-contact} with these
\(m,t\). Indeed the point criterion of
Section~\ref{beyond:wider-strip} requires the short curves only
on \(X\). All possible exceptions are now rank two, their common
polar kernel is an \(m\)-plane, and the one-light calculation
handles contact \(h\le t\). The codimension calculation then
has exactly the same flag dimensions and graded dimensions
as there.

\paragraph{All degrees at least forty on the next dimension.}
When \(k=1\) the same conclusion holds for every \(d\ge40\).
We detail rational substitutes for the last coarse estimates.
Here \(K=4\). For \(d\ge d_0=200\), partition
\(x=h/d\in[1/4,1]\) into intervals \([L,H]=[j,j+1]/32\),
\(8\le j\le31\). Put \(b_0=1+(d_0-1)/(d_0+4)\),
\[
 t_0=(d_0-1)/(d_0L-1),\quad R_0=t_0L,\quad
 J_0=R_0(2e t_0+3s)/(e(d_0-1)),\quad M_0=2-Ls/(s+e).
\]
These upper-bound \(t,R,J,M\) respectively, by monotonicity
in \(x,d\). Put \(U_0=H(4-s/e)/(b_0-1)\).
Since \(I-b\le(4-s/e)/d\), the intersection test applies
at \(u\ge U_0\). Each of the 240 elementary substitutions
\(e=1,2\), \(1\le s<4e\), \(8\le j\le31\) satisfies
\(M_0\le b_0,\ J_0\le b_0\) and either \(U_0\le1\) or
\(M_0+\max(1,U_0-1)J_0\le b_0\). Thus the first-blowup test
covers the entire complementary interval.

For \(40\le d\le199\), use the previous binary, product, endpoint
and conic tests first. For the remaining entries use the signed
version of the first jet,
\[
  m_1\le J-K_s(m+1/h)-1/h,\qquad K_s=s/(e(d-1)).
\]
Write \(a=(b-1)/h\). The original bound \((A-1)\min(M,I/u)\)
suffices at \(A=u+1\) when \(u\le b/(M-a)\) or
\(u\ge (I-b)/a\). In the remaining interval
\([u_0,u_1]=[\max(1,b/(M-a)),(I-b)/a]\) test
\[
 (1-vK_s)_+M+v(J-(1+K_s)/h)\le b+au,
             \qquad v=\max(1,u-1).
\]
Its left side bounds \(m+v m_1\) by the signed jet inequality.
On either side of \(u=2\) the difference from the right hand side
is convex, so the endpoints and \(u=2\), if present, suffice.
Slope domination follows from \(M,J\le b\) at every entry to which
this test is applied. For \(1\le s<4e\), \(2\le h\le d\), there
are 189,600 triples. The number first passing the binary test,
preceding conic tests, endpoint test and signed first-jet test are,
respectively, 45,600, 91,031, 24,647 and 28,322. Both the 240
upper-box inequalities and all the finite tests are implemented in
\path{nov27-lit-third-row-nets.py} with exact rational arithmetic.

\subsection{A short-curve criterion in degrees nine through thirty-two}
\label{sec:nov27-low-third-row}
Let \(9\le d\le32,\ n=2d+1,\ r=d+3\), and let \(X_d^n\) be smooth.
Define \(M=r+c_d\) by
\[
\begin{array}{c|cccccc}
d&9&10,11&12&13,14&15,16&17,\ldots,32\\ \hline
c_d&11/20&3/5&11/20&21/40&1/2&19/40 .
\end{array}
\]
Suppose \(\delta_C(X,H)\ge M\) for every line of normal-polar rank
two and for every integral conic of plane-normal polar rank two
whose plane is contained in \(X\) or intersects \(X\) with
multiplicity at least two along the conic. Then \(X\) is uniformly
K-stable. In particular the conclusion holds when there are no
such lines or conics.

Write \(e=\deg C\) and \(b=(n+1)/M>3/2\). For a simple proper
conic use the plane flag. The ordinary angular bounds are
\(2-(q-1)/d\) for a line of polar rank \(q\) and
\(2-q/(2d)\) for a conic of plane-normal polar rank \(q\).
At the remaining ranks the common kernel has dimension greater
than \(n/2\), giving finite contact \(2\le h\le d\). For rank
at least three retain a basepoint-free net. Use the local
estimates of Sections~\ref{sec:nov27-next-lines},
\ref{sec:nov27-small-conic-2d} and~\ref{sec:nov27-net-strip}
with \(1\le s\le e(d-1)/2\) and the present \(b\). Those tests
only use the dimensions for the existence of the common kernel.
The remaining cases have \((e,s)=(1,1),(2,3)\).

\paragraph{A jet estimate at the direction of the branch.}
Use the point \(P\) of the moving graph of degree \(s\) on the
angular plane and the direction \(P_1\) of the branch \(L\) on the
first blowup. Let \(D=e(d-1)\), \(t=(d-1)/(h-1)\), \(R=th/d\),
\(u=1/(et)\), and \(v=1/s\). The upper bound on the normalized
total multiplicity at \(P\) is \(M_0+1/h\), where
\[
 I=2-\frac{e+s}{ed},\qquad
 M_0=\min\!\left(I,\,2-\frac{hs}{(e+s)d}-\frac1d\right).
\]
Each actual filtered block embeds after multiplication by a base
scalar in the product with bidegree \((B,k)=(\lfloor etp\rfloor,k)\).
Multiplication by this scalar preserves point orders over the
geometric generic base. It therefore induces injections on all
point jets and their subspaces with a specified directional order.
The asymptotic ratio of enlarged to actual dimensions is \(R\).
When the graph maps onto a line, the two graph generators and
the residue calculation of Section~\ref{sec:nov27-small-conic-2d}
give an upper bound \(RJ_s\) for its normalized mean directional
order at \(P_1\), and \(R\) at any other direction, with
\[
 J_s=\frac{1/(1+Du)-(su)^2/(D/s+su)}{1-(su)^2}.
\]
At \(su=1\) the value is \((1+(D-s)/(2(D+s)))/(1+Du)\).
Indeed the two generator powers have costs \(1,1+su\) in the
regularized product, their residual coefficient density is
\(B'+sc'\) with \(B',c'\) the free base and plane degrees,
and the marked power costs \(D\) base degrees. These are the
same product integrals with the indicated \(s,e\).
This treats every graph with \(s=1\).

For \(s=3\) and image spanning the plane, its three coordinate
sections generate \(\mathcal O_{\mathbf P^1}(3)\) without a
constant relation. The kernel of the evaluation splits as
\(\mathcal O(-1)\oplus\mathcal O(-2)\). The branch equation
gives a subline \(\mathcal O(-D)\) with locally free quotient,
since the net is basepoint-free. The order-\(j\) leading jet
belongs to the \(j\)-th symmetric power of the kernel tensored by
\(\mathcal O(B+s(k-j))\). A jet vanishing at the marked direction
to order \(q\) is divisible by the \(q\)-th power of the subline.
Writing \(x+y=j-q\), its dimension is therefore bounded by
\[
 \sum_{x+y=j-q}\bigl(B+s(k-j)-x-2y-Dq+1\bigr)_+ .
\]
At another direction one can replace \(Dq\) by \(q\). Indeed, in
the splitting basis over the generic base, if the marked root is
finite, the highest \(j-q+1\) coefficients determine such a jet:
a polynomial of degree less than \(q\) vanishing there to order
\(q\) is zero. The retained monomials have kernel degrees obtained
by shifting \(x+y=j-q\) by \(q\) units of either exponent; each
shift costs at least \(q\) base degrees. At infinity use the
opposite ordering. This argument also works at geometric points.

For \(T=D,1\) put \(g_T(z)=1/((1+z)(1+2z)(1+Tz))\) and define
\[
 V_T=\frac{u}{v-u}
 \left\{\frac{v}{u}g_T(u)-\frac{u}{v}g_T(v)\right\}.
\]
Take the continuous value \(2g_T(u)-ug_T'(u)\) at \(u=v\).
Then \(RV_D\), respectively \(RV_1\), bound the two normalized
average angular orders. To see this write \(c=k-j\) and use the
homogeneous simplex integral in Laplace coordinates. The full
product dimension integral is
\(\int e^{-p-k}etp\,k^2/2\,dp\,dk=1/u\).
The numerator has weight \(e^{-p-c-x-y-q}\) and integrand
\((etp+sc-x-2y-Tq)_+\) on the nonnegative orthant.
Integrating first in \(p,c\) gives the expected excess of a sum
of independent exponentials of rates \(u,v\) above \(x+2y+Tq\).
The remaining three integrations give exactly \(V_T/u\).

Let \(J\) be \(RJ_s\) for a line image and \(RV_D\) otherwise.
The direct blowup adjunction test uses
\(M_0+1/h\le b(1+1/h)\), \(hJ\le b\), and the unmarked
bound \(\min(M_0+1/h,R)\le b\) for a line image or
\(\min(M_0+1/h,RV_1)\le b\) otherwise. We give a second test
when the marked estimate alone falls short. Normalize the plane
valuation by \(v'(\mathfrak m_P)=1\) and put \(w=v'(L)>1\).
The first exceptional \(E\) has value one, and its marked point
is the new center. If an individual section has order \(j\) at
\(P\) and its first nonzero jet has order \(q\) at the direction
\(P_1\), its strict transform has intersection number \(q\)
with \(E\) there. For \(q>0\), inversion of adjunction on the
smooth surface at \(E\) gives
\[
 v'(\mathrm{div}\,f)\le j+q(A-2),
 \qquad A=A_{\mathbf P^2}(v').
\]
Indeed \(A_{\mathrm{Bl}_P\mathbf P^2}(v')=A-1\) and the
coefficient-one exceptional divisor accounts for one more unit.
For \(q=0\) the strict transform avoids the new center.
Take a basis simultaneously adapted to \(v'\) and point order.
The leading jets are independent in each order, so their mean
directional order is at most \(J\) by the dimension estimate.
The total normalized mean is at most \(M_0+1/h+(A-2)J\).
This bound and the earlier intersection bound for the same
expected vanishing may be computed with different adapted bases.
When \(M_0,J\le b\) it suffices to test at \(A=w+1\).
With \(a=(b-1)/h\), the original intersection estimate works
at \(w\le b/(M_0-a)\) or \(w\ge(I-b)/a\), including \(w=1\)
by its original-point version. At both endpoints of the intervening
interval above \(w=1\) we test
\[
        M_0+1/h+(w-1)J\le b(1+w/h).
\]
The two sides are affine in \(w\). These tests verify all the
remaining higher ranks and hence supply the required bound on
each short curve. The exact substitutions appear in
\path{nov27-lit-low-third-row.py}.

The pointed hyperplane sections belong to the dimension-\(2d\)
theorem. For \(10\le d\le32\) use the preceding finite-fraction
tilt and point criterion with the displayed curve mark. In
ascending order the tilt shares in thousandths are
\[
\begin{gathered}
689,727,651,623,647,610,629,587,602,615,628,640,\\
650,660,670,678,687,694,701,708,715,721,726 .
\end{gathered}
\]
For \(B<t<17/2\) the previous surface and sectional Milnor
bounds apply, including the singular cubic bound at \(k>9\);
each remaining individual pair \(j<20\) passes the single-value
principal-bin test with \(3n\) bins. For \(j\ge20\) use the
same continuous envelope. The intervals with endpoints
\(17/2,d/2,d\) pass their previous comparisons (one may drop
the first interval if \(d/2<17/2\)). The above certificate
checks all the point comparisons in these degrees. The degree-nine
point calculation is given separately in the nonic criterion.

\paragraph{A nonic nineteenfold.}
Let \(X\subset\mathbf P^{20}\) be a smooth hypersurface of degree nine.
If \(\delta_C(X,H)\ge251/20\) for every line and integral conic
\(C\subset X\), then \(\delta(X,H)>12\).
A smooth pointed hyperplane section belongs to the dimension-\(2d\)
range proved above. Use the same point criterion with \(n=19,\ N=20,\)
\(r=12,\ a=544/1000,\ f=1/9,\ \epsilon=10^{-8}\). At the
stated mark its rounded comparison constants are
\(B=1061/250,\ Q=23701/10000\); the surface cutoff and both
signed endpoint inequalities below \(B\) are strict.
For \(B<t=k/j<9\) use the pointed curve-section discriminant bound
\(I=9j(j+5)+k+1\), taking its minimum with the sectional Milnor
bound for \(j\le2\). As above discard \(k(k-1)>I\). For
\(1\le j<25\), with \(m=\min(k,\sqrt I)\), the surface test
\[
 \frac9t+\left(1-\frac9{t^2}\right)
 \frac{Im}{j(I+m(m-1))}<B
\]
holds except for \((j,k)=(1,7),(1,8)\), where the single-value
principal-bin test passes with \(56,69\) bins, respectively.
For \(j\ge25\) use \(E=9+54/25+1/625\) and
\(L_*=E/(2\sqrt E-1/25)\). Evaluating at \(2L_*\) clamped to
\([B,9]\) gives \(\max(9/t+(1-9/t^2)L_*)<B\).
If \(t=9\), the exceptional divisor on \(X\) computing \(t\) is
a cone of degree equal to multiplicity, lying in \(T_pX\) since
all its generating lines pass through \(p\). The tangent section
is integral by finiteness of the Gauss map, so it equals this cone.
Its base is a smooth nonic seventeenfold, uniformly K-stable by the
preceding dimension-\((2d-1)\) range. Consequently the generalized
Eckardt formula \cite[Corollary~4.19]{ar2003.13788} gives
\(\delta_p(X,H)=380/27>251/20\), as required. All numerical
inequalities are verified by \path{nov27-bin-nonic.py}.

\subsection{The possible exceptional locus in dimension \(2d+1\)}
\label{sec:nov27-third-rank-two-contact}
Let \(d\ge9,\ n=2d+1,\ m=n-1=2d\), and \(N=n+1\). For \(d\le32\)
let \(M\) have the value in Section~\ref{sec:nov27-low-third-row};
for \(d\ge33\) put \(M=N(d+4)/(2d+3)\). Define
\[
 b=N/M,\qquad t=\left\lfloor\frac{d(b-1)}{2d-1-db}\right\rfloor .
\]
Every smooth member outside a closed set of codimension at least
\[
\min\left\{
  \sum_{a=0}^{t}\binom{m-2+a}{a}(d-a+1)-4m,\
  \sum_{a+j\le t}\binom{m-3+a}{a}\bigl(2(d-a-2j)+1\bigr)-5m-1
  \right\}
\]
is uniformly K-stable, where \(a,j\) are nonnegative. For example,
the codimension is at least \(1315\) for degree nine in dimension
nineteen, and at least \(1831\) for degree ten in dimension twenty-one.
For every \(d\ge33\), the formula uses \(t=\lfloor d/4\rfloor\).

Indeed the third-row criterion, including the finite low-degree
tests, reduces the problem to curves of polar rank two. For a
line or a conic of multiple plane contact with this rank, its
common polar kernel \(W\) is a projective \(m\)-plane. Write
\(h\) for the order of the restricted equation along the curve.
We have \(2\le h\le d\), since a smooth hypersurface cannot contain
this \(W\) of dimension greater than \(n/2\). The one-light
calculation of Section~\ref{sec:nov27-high-contact-strip} gives
normalized denominator at most
\[
 \max\{1,\ 1+1/h-1/d,\ 1+(d-1)h/[d(h+1)]\}.
\]
Here \(3/2<b<2-1/d\), by the displayed marks. The maximum is
at most \(b\) when \(h\le t\). Hence any remaining exceptional
hypersurface admits a flag \(C\subset W\) with its restricted
equation in \(I_{C/W}^{t+1}\). For \(d\ge33\) the expression
inside the floor is \(d/4\); for \(9\le d\le32\), substitution
gives \(2\le t\le d/2\). The conic ideal in \(W\) is a complete
intersection generated by \(m-2\) linear equations and a quadric.
Counting the graded pieces of orders \(0,\ldots,t\) for the line
and conic ideals as in Section~\ref{sec:nov27-high-contact-strip}
gives the two sums above. The \(m\)-planes vary in dimension
\(2(m+1)\); their lines and smooth conics vary in dimensions
\(2(m-1)\) and \(3m-1\), respectively. Subtracting the flag
dimensions and taking the closures of the resulting incidence
images proves the bound.

\paragraph{Point estimates in degrees six, seven, and eight.}
Let \(d\in\{6,7,8\}\), \(n=2d+1\), and let \(X_d^n\) be smooth.
Assume that a general smooth pointed hyperplane section has
\(\delta(Y,H)>d+2\). Suppose that every line and integral conic
\(C\subset X\) has \(\delta_C(X,H)\ge M\), with \(M\) as follows.
Then \(\delta(X,H)>d+3\).
\[
\begin{array}{c|ccc}
d&6&7&8\\ \hline
M&48/5&263/25&11501/1000\\
1000a&468&388&402\\
B&323/100&81/25&69/20\\
Q&16099/10000&13547/5000&2012/625
\end{array}
\]
Indeed, use the comparison of Section~\ref{sec:nov-second-row}
with the displayed tilt share \(a\), \(f=1/d\), and
\(\epsilon=10^{-8}\). The surface cutoff and both signed comparisons
up to \(B\) are strict. For \(B<t=k/j<d\) take
\(I=dj(j+d-4)+k+1\), intersecting with the sectional Milnor bound
when \(j\le2\), and discard the impossible values \(k(k-1)>I\).
For \(j<25\), apply the same surface comparison as in degree nine.
Its complete list of exceptions consists of the pairs \((1,k)\)
in the following table. The entry attached to each \(k\) gives the
number of principal bins which proves the signed comparison there.
\[
 \begin{array}{c|lll}
  d & \multicolumn{3}{c}{(k,\text{number of bins})}\\ \hline
  6 &(5,26)\\
  7 &(5,29)&(6,39)\\
  8 &(5,30)&(6,42)&(7,54)
 \end{array}
\]
For \(j\ge25\) use
\(E=d+(d(d-4)+d)/25+1/625\) and \(L_*=E/(2\sqrt E-1/25)\).
The maximum of \(d/t+(1-d/t^2)L_*\), computed by clamping
\(2L_*\) to \([B,d]\), is strictly smaller than \(B\).
At \(t=d\) the tangent section is a cone exactly as in the
degree-nine argument. Its smooth degree-\(d\) base has dimension
\(2d-1\), which is in the preceding all-member range. The
generalized Eckardt formula therefore gives
\(\delta_p(X,H)=n(n+1)/(n+d-1)>M\) in each case. The arithmetic
is certified by \path{nov27-bin-third-six-eight.py}.

\subsection{Explicit polar criteria in dimensions thirteen, fifteen and seventeen}
\label{sec:nov27-third-low-polar}
Let \(d=6,7,8\) and \(n=2d+1\). In the smooth parameter space
there is a closed subset of codimension at least \(25,43,65\),
respectively, outside which every hypersurface is uniformly
K-stable. Let \(M\) be the mark of the preceding point
calculation, \(b=(n+1)/M\), and \(A_0=\lfloor2d(b-1)\rfloor\).
The required geometric criteria are the following.
\begin{enumerate}
\item Every line has normal-polar rank at least three.
\item For each smooth conic of plane multiplicity at least two
(including a contained plane), every positive splitting degree of
its plane-normal polar kernel is at most \(A_0\).
\end{enumerate}
For a line of rank at least five the ordinary angular denominator
is at most \(2-4/d\le b\). For a conic satisfying the second
hypothesis it is at most \(1+A_0/(2d)\le b\).
The ordinary radial denominators and the simple plane conic
denominator \(3/2-1/d\) are also smaller than \(b\).

For a line of rank three or four keep a general basepoint-free
net of polars and the whole common kernel of dimension at least
\(n-3>n/2\). Its first light order satisfies \(2\le h\le d\).
Write \(s\) for the smaller of the two positive splitting degrees,
so \(1\le s\le(d-1)/2\).
The line tests of Section~\ref{sec:nov27-next-lines} suffice for
\(s\ge2\). For \(s=1,\ h\ge3\) use either the line product test or
the marked graph test of Section~\ref{sec:nov27-low-third-row}.
The exact comparisons, including the axis and the off-branch
and low-summand estimates, appear in
\path{nov27-lit-third-six-eight.py}.

It remains to treat \(s=1,\ h=2\) at the distinguished point \(P\)
of the angular plane. Normalize \(v(\mathfrak m_P)=1\), let
\(w=v(L)\ge1\) for the branch line, and let \(A\ge w+1\) be the
ordinary plane discrepancy. Put \(p=2/d,\ I=2-2/d\) and \(D=d-1\).
The actual normalized plane mean \(S\) has regular mass \(p\)
and positive-pole mass \(1-p\). As in the marked pole specialization
of Section~\ref{sec:nov27-sextic-twelve-conics}, the binary
coefficient \(a\) now has degree \(d-2\): the relation is
\(x^2=-L/a\) over the generic line. Filter the positive-pole
blocks by the roots and successive jets of \(a\), specializing the
moving point together with the basis. Each residue consists of
the complete plane series times \(L(p_0)^i\) at its root \(p_0\),
with conditional normalized means three for the free degree
and \(1/2\) for \(i\). Thus \(S_{\rm pos}+1/2\le 3A/2\), since
the marked point specializes onto that line. The same calculation
applies with binary base polarization of degree one: a root
jet costs one binary degree, and each root jet has the same
conditional simplex as in that specialization.

The degree-one graph parametrizes a fixed plane line. In the
regular product at bidegree \((j,k)\), restriction after removing
\(l\) copies of this fixed line gives the full product with
bidegree \((j,k-l)\). Order \(q\) at the moving point requires
divisibility by the graph equation of bidegree \((1,1)\) to
power \(q\), with dimension \((j-q+1)_+(k-l-q+1)_+\).
The mean fixed-line coefficient is one. Shifting both indices by
\(q\) costs \(2q\) in the section simplex, giving normalized mean
restricted order at most \(1/2\). Adding the marked ideal to
exponent \(1/2\) and applying inversion along that fixed line
gives \(S_{\rm reg}+1/2\le A\).
The regular mean point order is at most \(3/2\) by graph
projection. If \(w>1\), the valuation centers on the branch
direction after the first blowup. The regular jet at that
direction has normalized order at most
\[
           J=\frac1d\left(1+\frac{d-2}{2d}\right)
\]
by the line graph formula with \(t=1\) for the regular product.
Inversion along the exceptional line gives
\(S_{\rm reg}\le3/2+(A-2)J\).
The regular and pole estimates bound the same tested valuation,
using bases adapted to that valuation and to the flag used on
each sector.

We must prove \(S\le b(A-w/2)\). Put
\[
 K_0=1+(1-p)/2,\quad K=3(1-p)/2+pJ,\quad
 R=p(3/2-2J)-(1-p)/2 .
\]
The preceding inequalities give \(S\le K_0A-1/2\) and, at
\(w>1\), \(S\le KA+R\). Here \(K_0,K<b\) and
\(2K_0-1/2\le3b/2\), settling \(w=1\). The original
residual intersection estimate is
\[
       S\le w/2+(A-1)I/w ,
\]
since the residual multiplicity and intersection bounds coincide
at \(I\). At \(A=w+1\) the slack is at least
\(B=b+(b-1)/2-I\). If \(w\ge I/b\) the residual slope
is at most \(b\). Otherwise write
\(T(w)=(R+bw/2)/(b-K)-w-1\), an affine function.
At distance at least \(T(w)\) from \(A=w+1\) the marked bound
suffices. Up to that distance the residual estimate can lose
at most \(E=(I-b)\max(0,T(1),T(I/b))\). Direct substitution gives
\[
\begin{array}{c|ccc}
 d &6&7&8\\ \hline
 B&1/48&247/3682&4491/46004\\
 E&215/13944&0&0
\end{array}
\]
and \(B>E\) throughout. This proves the local line bound.
The pointed hyperplane theorem in dimension \(2d\) and the
preceding point calculation in dimension \(2d+1\) now give
uniform K-stability under the stated two criteria.

Containing a fixed line costs \(d+1\) conditions. Its \(n\)
normal polars are unrestricted binary forms of degree \(d-1\).
Rank at most two costs an additional \((n-2)(d-2)\) conditions.
Subtracting the \(2n\) parameters of the line gives the respective
bounds \(25,43,65\).

For the conic put \(D'=2d-2,\ t=A_0-1\). The fixed plane equation
must be divisible by the square of the conic, at cost \(4d-2\).
Independently the \(n-1\) polars \(F_i\in S_{D'}\) vary freely.
Failure of the splitting bound is equivalent to failure of
surjectivity of \(\bigoplus_i S_t\to S_{D'+t}\),
\((G_i)\mapsto\sum F_iG_i\): twist the evaluation sequence by
\(\mathcal O(t)\) and take cohomology. A nonzero functional
\(\lambda\in S_{D'+t}^*\) must then annihilate every \(F_iS_t\).
At any such basepoint-free tuple, the pairing
\((P,B)\mapsto\lambda(PB)\) on \(S_{D'}\times S_t\) has rank
\(t+1\). Otherwise some nonzero \(B\in S_t\) is in its kernel.
Basepoint freeness gives a linear combination \(F\) of the polars
coprime to \(B\). But \(B S_{D'}+F S_t=S_{D'+t}\) by the
dimension formula, since the intersection is generated by \(BF\).
The functional would be zero. For fixed \([\lambda]\) of the
required rank the polar equations therefore cost \((n-1)A_0\),
and the projective functional varies in dimension \(D'+t\).
Finally the smooth plane conic varies in dimension \(3n+2\).
The resulting lower bound on codimension is
\[
      4d-2+(n-1)A_0-(D'+t)-(3n+2)
             = 2d-1-3n+(n-2)A_0 .
\]
Here \(A_0=5,7,9\), giving \(27,59,99\), respectively. The closures
of these incidence images have the same dimension bounds, proving
the asserted codimensions.

\subsection{A flag along a smooth plane curve of arbitrary degree}
\label{sec:nov27-plane-curve-flag}
Let \(C\subset\mathbf P^3\) be a smooth plane curve cut out by a linear
form \(Z\) and a form \(Q\) of degree \(e\ge2\) in the coordinates of
its plane. Let \(p\in C\), and choose an ordinary affine chart about
\(p\) with coordinates in which \(Z,Q\) denote their dehomogenizations.
For an integer \(s\ge e\) and a unit \(U\) at \(p\), consider the
divisor germ \(W=(Z+UQ^s=0)\). Every divisorial valuation whose
center contains \(p\), with \(W\) evaluated as a germ, satisfies
\[
        4S_{\mathcal O_{\mathbf P^3}(1)}(v)+v(W)
        \le \left(1+\frac{es}{es+e+1}\right)A(v).
\]
In particular a smooth cubic component of multiplicity \(s\ge3\) on
an ordinary plane section gives the factor \(1+3s/(3s+4)\) at a
point where the remaining factor is a unit.

\begin{proof}
Put \(\lambda=e-1+(e-1/e)/s\). We prove the signed inequality with
numerator \(4S+\lambda v(W)\) and denominator \(eA\). Apply the
weighted blowup along \(C\) with weights \(s,1\) on \(Z,Q\). Its
log discrepancy is \(s+1\) and \(W\) has order \(s\). The ordinary
section ring in degree \(m\) has associated graded pieces indexed
by monomials \(Z^iQ^j\) with coefficients in
\(H^0(C,\mathcal O_C(t))\), \(t=m-i-ej\ge0\). Indeed expand first
in powers of the independent variable \(Z\), then in the plane use
successive quotients by powers of \(Q\). All coefficients lift by
the complete restriction sequence for a plane hypersurface, and the
distinct pairs are independent in the two normal coordinates. They
retain these properties when collected at equal weight \(si+j\).

Further filter each piece by vanishing order \(a\) on \(C\) at \(p\).
If \(g\) denotes the genus of \(C\), all integers \(0\le a\le et-2g\)
are realized, by Riemann--Roch applied after subtracting \(ap\) and
\((a+1)p\). Every realized order is at most \(et\). Thus omission
from the list \(0,\ldots,et\) occurs only within a bounded distance
of the upper endpoint, also for the finitely many small \(t\). The
number of omitted triples \((i,j,a)\) is \(O_e(m^2)\) against a
basis size of order \(m^3\), and their orders are \(O_{e,s}(m)\).
The barycentric coordinates \(i,ej,a/e,t-a/e\) have sum \(m\).
Their asymptotic means, after multiplication by \(4/m\), are one
by simplex integration. The means of \(i,j,a\) in these units
are consequently \(1,1/e,e\), respectively. The signed primary
coefficient is \(s+1/e+\lambda s=e(s+1)\). On the exceptional
bundle the next flag member is the fiber above \(p\), of normalized
coefficient \(e\). The strict transform of \(W\) misses both coordinate
points and cuts this fiber simply at a noncoordinate point, by the
nonvanishing of \(U(p)\).

For fixed \(m,si+j,a\), the full possible exponents \(i\) form a
consecutive interval under \(i,j\ge0,\ i+ej\le m-a/e\). Its width,
viewed from any of its indices, is
\[
 \lfloor j/s\rfloor+
 \min\{i,\lfloor (t-a/e)/(es-1)\rfloor\}.
\]
The two summands have scaled means \(1/(es)\) and \(1/(es)\).
For the second, the simplex conditions are \(i>v\) and
\(t-a/e>(es-1)v\), with combined cost \(esv\).
At a noncoordinate point of the fiber any independent system
of \(r\) polynomials in an interval of width \(w\) has total
order at most \(rw-r(r-1)/2\). Completing the missing monomials
changes the global normalized bound by \(O_{e,s}(1/m)\), since there
are \(O_e(m^2)\) omissions, each of cost \(O_{e,s}(m)\). Thus the
scaled mean order at such a point is at most half the mean width,
namely \(1/(es)\). At the marked point inclusion of \(W\) gives
\(\lambda+1/(es)=e-1+e/s\le e\). At the coordinate points, the
scaled bounds \(1,1/e\) follow from the means of the two exponents
on the uniformizing charts. Successive local inversion of adjunction,
testing strictly below \(1/e\) at finite level as in the preceding
conic-flag proof, proves \(4S+\lambda v(W)\le eA\). Combine it
with \(4S\le A\) with coefficients \(1/\lambda\) and \(1-1/\lambda\)
to obtain the claim.
\end{proof}

\subsection{A signed graph bound from a plane threshold}
\label{sec:nov27-binary-signed-general}
Let \(G(x,y)\) be a polynomial of ordinary degree at most \(d\) at a
point \(p\) of the affine plane with \(G(p)=dG(p)=0\). Suppose
\(T\ge\max\{2,\sqrt d\}\) and \(\operatorname{lct}_p(G)\ge1/T\).
For the graph \(W=(z+G=0)\) in the ordinary affine chart of
\(\mathbf P^3\), one has, at the point over \(p\),
\[
       4S_{\mathcal O(1)}(v)+v(W)
             \le\left(1+\frac{T}{T+1}\right)A(v).
\]
In particular this applies to \(T=d/4\) for \(d\ge16\).

First the smooth plane-curve flag of
Section~\ref{sec:nov27-plane-curve-flag} works unchanged for an
integral plane curve of degree \(m\) smooth at the specified point,
with no smoothness condition elsewhere. Indeed powers of its prime
equation and the plane equation still give the same graded coefficient
spaces, since the plane hypersurface sequence gives surjectivity of
all homogeneous restriction maps. For the orders
at the smooth Cartier mark use Riemann--Roch and Serre duality on the
integral Gorenstein plane curve, with its arithmetic genus. A line
bundle of negative degree has no global sections, so all orders up to
\(mt-2p_a\) occur in section degree \(t\), as in that proof. The
extraction and the adjunction on the marked fiber are used only on
the smooth open of the curve around that point. Thus at a smooth
point of the curve germ \(Q\), if \(e\ge m\), the graph
\(z+UQ^e\) with a local unit \(U\) has signed denominator at most
\(1+me/(me+m+1)\).

Translate \(p\) to the origin and use the two-variable Newton
crossing argument. Its multiplicity is at most \(2T\) by the ordinary
point blowup. If each root of its homogeneous tangent form has
multiplicity at most \(T\), use weights \(1,1\). Otherwise put the
unique higher-multiplicity linear factor equal to \(x=0\). A term
with \(x\)-exponent at most \(T\) exists by the divisor test for
the assumed threshold. Increasing the relative weight of \(x\)
therefore meets a compact lower edge of weights \(m>l\ge1\),
coprime, with its two endpoint \(x\)-exponents \(i_0\le T<i_1\).
The log canonical condition applied to this weight gives face weight
\(H\le T(m+l)\), so the smaller endpoint \(y\)-exponent
\((H-mi_1)/l<T\). Give \(z\) weight \(H\). The normalized radial
mean of the complete series is \(m+l+H\), its discrepancy, and the
graph adds \(H\); hence the radial signed denominator is as asserted.
Also \(H>m,l\), by ordinary multiplicity two or by \(i_1>T\).

On the punctured angular cone noncritical graph points have contact
two by weighted Euler, and angular critical points on the two axes
have contact at most \(T\) by the endpoint tests. The same is true
for the equal-weight roots. At a weighted torus point the unique
vanishing irreducible factor of the face is
\(Q=x^l-\lambda y^m\). Let \(e\) be its multiplicity. If \(e\le T\)
the smooth-contact bound suffices. Otherwise the ordinary degree gives
\(me\le d\), thus \(m<d/T\le T<e\). Its degree-\(m\) projective
closure is an integral plane curve, smooth at this nonzero torus
point. Moreover \(ml e\le H\le T(m+l)\) by the weight of the factor,
so \(me\le T(1+m/l)\le T(m+1)\). The plane-curve flag with the
other factors as a local unit now gives
\(1+me/(me+m+1)\le1+T/(T+1)\), including singularities of the
projective closure at infinity. As in the preceding weighted plane
argument, grouping the complete ordinary polynomial basis by its
weights preserves that complete ordinary system. The signed bounds
on the punctured smooth covering charts and the radial inequality
therefore give the assertion by weighted adjunction.

\subsection{Disjoint monomial faces without a coprimality condition}
\label{sec:nov27-disjoint-arbitrary}
Let \(M_1,\ldots,M_a\) be nonconstant monomials in mutually disjoint
nonempty sets of variables, each of ordinary degree at least two. Let
\(b_i\) be the largest exponent appearing in \(M_i\), set
\(c=\sum_i1/b_i\), and let \(k\) count all active variables. For nonzero
constant coefficients \(u_i\), consider the graph
\(W=(z+\sum_i u_iM_i=0)\) in the affine chart of ordinary
\(\mathbf P^{k+1}\). Throughout this affine chart every divisorial
valuation whose center contains the tested point satisfies
\[
       (k+2)S_{\mathcal O(1)}(v)+v(W)
       \ \le\ D A(v),\qquad
       D=\max\left\{\frac53,1+\frac1{1+c}\right\}.
\]

We induct on the number of active variables. Give the variables in
block \(i\) a common positive weight \(m_i=H/\deg M_i\), taking \(H\)
a common multiple, and give \(z\) weight \(H\). In particular all
\(m_i<H\). Off the graph use \((k+2)S\le A\). At a noncritical graph
point weighted Euler's identity applied to each first derivative gives
a nonzero Hessian entry; after subtracting the tangent affine linear
function the graph has contact two. Its smooth-graph bound is \(5/3\).
At a critical point each monomial vanishes and has order at least two
there: its partial derivatives cannot cancel those of another block.
If any active coordinate is nonzero, fix every such coordinate to its
current value by successive affine-linear hyperplanes, taking their
projective closures. The restricted graph still has all \(a\) blocks,
each of degree at least two, and its block maxima have only decreased.
The induction hypothesis bounds it with denominator at most \(D\).
Complete hyperplane flag adjunction lifts this signed estimate: in
dimension-plus-one units each cut has mean coefficient one, the
residual complete coefficient drops by one, and the graph does not
contain the cutting hyperplane. Since \(D>1\) all radial comparisons
for these linear cuts hold.

At the origin use the weighted blowup with the preceding weights,
dividing by their common gcd if needed. Its log discrepancy is
\(A_0=H+\sum_{i,j}m_i\), where \(j\) ranges over the variables in
block \(i\). The complete ordinary series has normalized mean \(A_0\)
there, and the graph has order \(H\). Since
\(H\le b_i\sum_j m_i\) in each block, the signed radial denominator
is at most \(1+1/(1+c)\). On the angular punctured cone the graph
equation is unchanged by weighted homogeneity. The union of the
independent fixed-weight bases in a given ordinary degree, after
ordinary homogenization, is a full basis of the ordinary polynomial
system. It is therefore enough to use the estimates already proved
at every nonzero graph point above, on the smooth punctured cover.
For this reduction test strictly below \(1/D\) at finite levels,
using the uniform complete-series upper bounds and then weighted
adjunction. If all coordinates vanish, the point of the graph is
the origin; no such point occurs on this angular punctured chart.
Thus each application of induction at an angular critical point
strictly lowers the number of active variables. In the base case
there is no critical punctured point with all surviving coordinates
zero. This proves the estimate. Additional unused affine coordinates
may always be fixed by further ordinary flags.

Suppose at a rank-two polar line two heavy variables have the weight
of a light face \(\sum_i u_iM_i\) of the preceding kind, with
coefficients nonzero at general polar roots. If \(d\ge8\) and
\(\sum_i1/b_i\ge4/d\), the pole-block comparison gives
\[
             \frac{\delta_L(X,H_X)}{n+1}\ge\frac{d+4}{2d+3}.
\]
In the primary radial computation the sum of light weights divided
by their common face weight is at least \(\sum_i1/b_i\), by the
same inequalities, even if the primary weights chosen on the line
are unequal within a block. Each positive residue has the signed
ordinary graph bound \(D\le 1+d/(d+4)\), since \(d\ge8\).
One regular and \(d-1\) positive groups have equal weights, so their
mean denominator is at most \(1+(d-1)/(d+4)=(2d+3)/(d+4)\).
The direct base specialization and complete ordinary series arguments
of Chapter~\ref{ch:fullrow} apply. This conclusion does not require
coprimality among the blocks' exponent gcds.

For example let \(q,r\) be constant quadrics such that every nonzero
linear combination has rank at least fifteen. There is a four-space
of light vectors on which \(q=x_1x_2\) and \(r=x_3x_4\) after
rescaling. Indeed on any vector space with pencil rank at least
seven take an \(r\)-isotropic vector \(v\) with \(q(v)\ne0\).
On \(v^{\perp_q}\cap v^{\perp_r}\) every nonzero combination still
has rank at least three. Hence the restricted \(r\) is irreducible
and does not divide the restricted \(q\); choose a second vector
on its zero set where \(q\ne0\). Their span has identically zero
\(r\) and nondegenerate \(q\). Pass to the joint orthogonal
complement of this two-space for both forms; its codimension is at
most four, so the restricted pencil still has rank at least seven.
Repeat with the two quadrics interchanged. The two planes are
disjoint and mutually orthogonal for both forms. A hyperbolic
basis in each gives the assertion. In particular, if the first
light term is \(a q^s\), the first subsequent term modulo \((q^s)\)
is \(b r^t\) of ordinary degree \(2t>2s\), and the terms of higher
degree are arbitrary, then give \(x_1,x_2\) weight \(t\) and
\(x_3,x_4\) weight \(s\). The precise face is the sum of the two
disjoint powers. Its mark follows from \(s,t\le d/2\).

\subsection{A repeated quadric and the next light coefficient}
\label{sec:nov27-quadric-next-form}
Put \(T=d/4\) and suppose \(d\ge16\). Consider a rank-two polar
line with a constant quadratic form \(q\) of rank at least seven on
the retained light space. Suppose the first light form of order
\(h\) is \(q^sR\), \(q\nmid R\), with \(s>3d/8\), and the form at
order \(t=h+1\le d\) has nonzero restriction \(P\) to the quadric
\(Q=(q=0)\). Over the geometric generic base take the divisors of
\(R|_Q\) and \(P\) on \(Q\). For each prime component of \(P=0\)
write \(a\) for its multiplicity in \(P\), \(b\ge0\) for its
multiplicity in \(R|_Q\). If
\[
          b\ge a\qquad\hbox{or}\qquad sa\le T(s+a-b)
                                                        \tag{Q}
\]
at every component, then
\[
           \delta_L(X,H_X)/(n+1)\ge(d+4)/(2d+3).
\]
Every component for which (Q) fails is the section of \(Q\) by a
constant hyperplane, with \(a>T\). The argument holds equally for
\(t>h+1\) if each intervening term at order \(j\) vanishes to
order strictly greater than \(s(t-j)/(t-h)\) along \(Q\).

Indeed choose a general two-dimensional totally isotropic vector space
in the nonvertex part of \(Q\), with coordinates \(z,w\). On its
projectivization, \(R,P\) restrict to binary forms such that each
root of \(P\) has the pair \((a,b)\) at its simple transverse
intersection with a unique component. This is general transversality
on the smooth quadric orbit: its linear group acts transitively on
the smooth locus, and a general translate of a fixed nonvertex line
avoids codimension-two intersections and singular loci, meeting each
divisor transversely on its smooth locus. The same choice works at
the general simple roots of the eliminated polar. Inside the
orthogonal complement of this two-space select a nondegenerate
hyperbolic pair independent of it, with coordinates \(x,y\).
This is possible since the quotient by the two-space and the
radical has rank at least \(\operatorname{rk}q-4\ge3\).
On these four lights \(q=xy\) after rescaling, without terms in
\(z,w\). Give \(z,w\) weight one and \(x,y\) weight
\(m=1+(t-h)/(2s)>1\), clearing denominators. The face weight is
\(t\) in these relative units. Its precise form, up to nonzero
scalars and writing the two binary restrictions again as \(R,P\), is
\[
              \Phi=(xy)^s R(z,w)+P(z,w),\qquad
              \deg R=B=h-2s<T,\quad \deg P=t.
\]
This follows by keeping only \(z,w\) in the residual factors. An
intermediate term divisible by \(q^j\) gains at least
\(j(t-h)/s\) over its ordinary degree, giving the stated extended
hypothesis. All higher terms have weight greater than \(t\).
The primary sum of light weights divided by the face weight is
\((2m+2)/t\ge4/d\), giving the required radial comparison.

We prove the signed ordinary graph estimate with denominator
\(D=1+T/(T+1)\) at every punctured point of the graph
\(v+\Phi=0\) in its affine chart (\(v\) denotes the graph coordinate
in this paragraph). Off the graph one uses projective semistability,
and at noncritical graph points weighted Euler gives contact two,
since each surviving weight is strictly below the face weight.
Thus consider a critical point. First suppose \((z,w)\ne(0,0)\).
Fix a nonzero linear combination of \(z,w\) by an ordinary affine
hyperplane flag. In a coordinate \(\tau\) at the resulting binary
slope the polynomial is \((xy)^s R_0(\tau)+P_0(\tau)\).
If \(x_0y_0\ne0\), criticality and weighted Euler imply that both
coefficients vanish there. Fix also \(y=y_0\). If the vanishing
orders \(a,b\) differ, the multiplicity is \(\min(a,b)\le B<T\).
If they agree, it is \(\tau^b\) times a local unit or times a germ
smooth and transverse to \(\tau\), since its derivative in \(x\)
is nonzero at a canceled nonzero root. Its plane threshold is then
at least \(\min(1,1/b)\ge1/T\), so the binary signed lemma applies.

If \(xy=0\), then \(P_0(0)=0\). If \(b\ge a\), the local expression
is \(\tau^a\) times a unit, with \(a\le B<T\); its tangent contact
suffices. Assume \(b<a\). If one of \(x,y\) is nonzero, fix it by
an affine hyperplane. In the two remaining variables the polynomial
is analytically \(\tau^b(X^s+\tau^{a-b})\) times a unit, by taking
an analytic root of a unit to scale \(X\). Its threshold is at least
\(1/T\) by (Q) and \(b<T\). Indeed its Newton polygon gives the
three bounds \(1,1/b,(s+a-b)/(sa)\), with \(1/b\) omitted for
\(b=0\); the binomial on its compact edge has reduced smooth torus
zeros in characteristic zero. A toric resolution thus gives the
stated plane threshold. Apply the binary signed lemma in the original
ordinary coordinates, where its polynomial degree is at most \(d\).

When \(x=y=0\), use at this point the weighted blowup in the three
light coordinates with equal positive weights on \(x,y\) and weight
one on \(\tau\), the former equal to \((a-b)/(2s)\). The graph
coordinate has weight \(a\). A critical point here has \(a\ge2\).
The signed radial ratio is at most \(D\) precisely by (Q), since the
face weight divided by the sum of light weights is
\(as/(s+a-b)\). The new graph polynomial is
\(\tau^b(c_1(xy)^s+c_2\tau^{a-b})\) with nonzero constants.
At a punctured point with \(\tau\ne0\) its zero set has \(xy\ne0\)
and a nonzero derivative. At \(\tau=0\) at least one of \(x,y\)
is nonzero; fix it and use the same binary threshold or the contact
\(b\) when both are nonzero and the point is critical. This proves
the signed estimate there and hence at the origin of this extraction.

It remains to consider \(z=w=0\), with one of \(x,y\) nonzero.
Fix that coordinate. If both are nonzero the initial binary
homogeneous polynomial is a nonzero constant times \(R\), of
degree \(B<T\), giving tangent contact at most \(B\) at a critical
point (when \(B=0\) there is no critical zero). If the surviving
other coordinate vanishes, write it \(y\) and use weights
\(((t-B)/s,1,1)\) on \((y,z,w)\). The graph polynomial
\(c y^s R(z,w)+P(z,w)\) is weighted homogeneous. Its signed
radial test follows from
\[
        \frac{st}{t-B+2s}=\frac{st}{t-h+4s}\le t/4\le T.
\]
At punctured critical zeroes either \(z,w\) are not both zero,
when the same fixed binary chart above applies, or \(y\ne0\) and
the tangent contact is at most \(B\). This exhausts all punctured
points. At each secondary extraction the normalized complete mean
is the sum of the positive coordinate weights, and the angular
pieces reassemble to an ordinary full polynomial basis. All the
linear restrictions lift by complete projective flag adjunction.
Applying the independent pole-block estimates to the initial pencil
now gives the displayed line threshold.

Finally the high-rank quadric cone is factorial: localizing at the
first coordinate in a hyperbolic pair gives a Laurent polynomial
ring, and that coordinate defines a prime divisor since the remaining
quadratic form has rank at least three, so Nagata's lemma applies.
If the component has class \(jH_Q\), \(j\ge2\), then
\(a\le t/2\le d/2\), \(b\le B/2\le d/2-s\). These imply (Q), since
\((s-d/4)(a-d/2)\le0\) rearranges to the desired inequality
with the largest possible \(b=d/2-s\). If \(a\le T\) and \(b<a\),
(Q) holds immediately. A failing component is therefore a hyperplane
section with \(a>T\). The whole Galois orbit of any such component
consists of high-multiplicity factors; its primitive product over
the binary base has base degree zero, since \(a\) times this degree
is at most \(d-t<T\). Equivalently take the primitive factors
in the factorial ring of the quadric with the two polynomial base
variables adjoined. A base-independent product splits over the
algebraically closed constant field, proving constancy of each
hyperplane. This concludes the reduction.

\subsection{A quadric flag for conic polar pencils}
\label{sec:nov27-conic-quadric-flag}
Let \(n=2d+1,\ d\ge33\), and suppose that a conic
\(C\subset X_d^n\) has polar rank two. The conic may lie in
a contained plane or have arbitrary multiplicity at least two
in its proper plane section. Choose plane coordinates
\(u_0,u_1,u_2\) and a complementary space of normal coordinates.
Let \(Q\) be the quadric cone cut out by the conic equation using
only these plane coordinates. Let \(V\) be the common polar
kernel of ambient codimension two through the plane of \(C\).
On \(Q\cap V\), write \(u=(s^2,st,t^2)\) in Veronese
coordinates. The polynomial
\[
                 F(s^2,st,t^2,x)\big|_V
\]
has a least nonzero degree \(h\) in the remaining linear normal
variables \(x\). Its degree-\(h\) form \(G_h(s,t,x)\) has
binary degree \(2(d-h)\). The normalized denominator
of \(\delta_C(X,H)\) is at most \((2d+3)/(d+4)\) in any of
the following cases:
\begin{enumerate}
\item \(h\le d/4\).
\item \(h\le d/2\) and every component of the geometric generic
      form has multiplicity at most \(d/4\).
\item \(h\le 3d/4\) and a general projective plane section of
      \((G_h=0)\) is log canonical with coefficient \(4/d\).
\item A general projective three-space cuts \((G_h=0)\) smoothly.
      In particular this holds whenever \(h=d\).
\end{enumerate}
These tests remain valid after further constant linear restrictions
of the light space provided its ambient codimension \(c\le d-2\),
retaining the two active normals. All tests refer to the restricted
polynomial modulo the conic equation. In particular the value of a
pure power of that equation in the original plane section imposes
no upper restriction here.

\begin{proof}
Over the conic one has independent equations \(dQ\) and \(dF\):
the former is a nonzero plane covector, while the latter has
basepoint-free components in the normal directions to the plane.
Cut by \(Q\) and take further linear sections through \(C\)
retaining two active normals. All these flags are smooth along
\(C\). Their coordinate ring is the even-degree Veronese part of
the weighted hypersurface in \(\mathbf P(1,1,2,\ldots,2)\) cut
out by \(F(s^2,st,t^2,x,z_1,z_2)\), with \(H=\mathcal O(2)\).
Any positive weighted degree factor either contains \(C\) or
meets \(C\), since its restriction to the parametrizing binary
line has positive degree. The full hypersurface contains \(C\),
so reducibility or a repeated factor would contradict smoothness
along \(C\). These sections are integral complete intersections.
The complete restricted graded pieces are the hypersurface
coordinate spaces, by projective normality of complete intersections.
For the first quadratic flag the normalized fixed coefficient is
\(1/2\), because the divisor has degree two. The complete-series
residual mean degree on it is \(n/(n+1)\), as follows from the
successive complete restrictions of \((m-2i)H\). Thus its
normalization factor is the usual dimension factor. Further
linear flags and weighted adjunction apply at the generic point
of \(C\).

The vertex given by \(u=0\) inside the restricted kernel has
dimension \(n-c-2>n/2\) for the indicated \(c\). It cannot be
contained in the smooth \(X\), so the substitution on the kernel
is nonzero and \(2\le h\le d\). Retain \(t_0=1,2,3,4\) general
constant lights in the respective cases. The rule that gives them
weight one and gives both heavy variables weight \(h\) has
initial relation
\[
           f(s,t)z_1+g(s,t)z_2+G_h(s,t,x)=0,
              \qquad \deg f=\deg g=2d-2.
\]
On the punctured conic cone the first two terms make the
hypersurface an affine bundle, and its vertex fiber has smaller
dimension. Hence the complete-intersection initial ring is a
domain, agrees in Hilbert functions, and embeds over the generic
point of \(C\), where it is the stated monomial filtration after
eliminating a heavy variable. Its radial denominator is
\[
         1+(1-1/d)/(1+t_0/h).
\]
The factor \(1/d\) comes from differentiating the initial relation
of ordinary \(H\)-degree \(d\). The required bound follows when
\(h\le t_0d/4\).

Here are the changes to the angular pole proof for the doubled
degrees of the binary coefficient ring. Choose the second polar
\(g\) with \(2d-2\) distinct general roots avoiding all excluded
values. For an ordinary angular monomial of total degree \(v\)
and pole index \(j\), the coefficient base degree at section level
\(M\) is \(2(M-j-v)\). At the \(a\)-th root, previous roots
cost \(a-1\) binary degrees, giving the bound
\(v\le M-j-\lceil(a-1)/2\rceil\) on the complete residual
polynomials. Each root contributes one positive group with fixed
graph factor to order \(j\). In the regular group the base
vanishing order \(i\) gives \(v\le M-\lceil i/2\rceil\).
There are asymptotically two regular copies (even and odd \(i\))
and \(2d-2\) positive copies, each of identical simplex mass.
The same residue divisibility argument holds in the even
Veronese ring, since \(g\) has even binary degree. Consequently
the total regular fraction is \(1/d\) and the positive fraction
\((d-1)/d\), with exactly the original conditional mean complete
and fixed coefficients. The ordinary-space angular bounds used
for one light, two lights, and for three lights with the
log-canonical coefficient \(4/d\) therefore apply verbatim.
The noncritical tangent contact two is admissible because
\(2\le d/4\). For four lights with a smooth projective zero
hypersurface, every point of the punctured graph has a nonzero
gradient in the light variables, by projective smoothness on the
zero set and Euler's identity on other levels. Euler
differentiation gives tangent contact two there. This proves all
four sufficient conditions with the claimed angular denominator.

Finally let \(h=d\), so \(G_d\) is constant in the base.
Let \(S\) be its projective singular locus inside the light
\(\mathbf P^{n-c-2}\). Fix a point \(p\in C\) and join it with
\(S\) inside the light kernel. At each point of the join off the
vertex, the restriction of \(F\) to \(V\cap Q\) and its
differential both vanish, since this restriction equals
\(G_d(x)\). Every tangent hyperplane of \(X\) there thus
belongs to the projectivized span of the \(c\) constant conormals
and \(dQ(p)\), of dimension \(c\); the direction of \(dQ\)
is fixed along the join off the vertex. Finiteness of the Gauss
map of the original smooth hypersurface gives \(\dim S+1\le c\).
It follows that \(S\) has codimension at least
\(n-2c-1\ge4\). A general projective three-space avoids it
and meets its complement transversely, proving the terminal claim.
\end{proof}

\begin{corollary}\label{cor:nov27-third-rank}
Let \(d\ge33\), \(n=2d+1\). A smooth \(X_d^n\) is uniformly K-stable
if each of the following conditions holds.
\begin{enumerate}
\item Each line has normal polar rank at least three, or satisfies a
rank-two criterion in Sections~\ref{sec:nov27-disjoint-arbitrary}
and~\ref{sec:nov27-quadric-next-form} (either criterion suffices).
\item Each integral conic occurs with multiplicity one in its plane
section, has plane-normal polar rank at least three, or satisfies one
of the tests in Section~\ref{sec:nov27-conic-quadric-flag}.
\end{enumerate}
\end{corollary}
\begin{proof}
The higher-rank calculations and the indicated rank-two tests supply \(\delta_C(X,H)\ge
(2d+2)(d+4)/(2d+3)\) at every line and conic in the hypotheses.
The preceding point criterion applies since the general smooth
pointed hyperplane section belongs to Theorem~\ref{thm:nov27-second-row}.
\end{proof}

\section{A line criterion up to quadratic dimension}
\label{sec:nov27-line-only-range}
Let \(N=n+1\), \(r=n+2-d\), and \(H=\mathcal O_X(1)\). For a line
\(C\subset X\), its polar rank is the dimension of the span of the
\(n\) first derivatives in the normal directions, restricted to \(C\).
It is an integer between two and \(d\).

\begin{theorem}\label{thm:nov27-line-only-range}
Suppose \(d\ge10\) and \(4d\le n\). Put
\[
 L_*=r+2,\qquad
 \ell=1+\left\lceil d\left(2-\frac{N}{L_*}\right)\right\rceil .
\]
Every smooth \(X_d^n\) whose lines all have polar rank at least
\(\ell\) is uniformly K-stable. When \(n\le d^2/3-4\), one has
\(\ell\le d-2\) and the complement of this rank condition in the
smooth parameter space is closed, of codimension at least
\[
 (n-\ell+1)(d-\ell+1)+d+1-2n\ \ge\ n-2d+10.
\]
For the larger range \(4d\le n\le(d^2-d-8)/2\), a sufficient
condition is that the Fano scheme of lines on \(X\) be smooth of
its expected dimension. Its failure locus also lies in an explicit
proper closed subset of the smooth parameter space.
\end{theorem}

\begin{proof}
Assume \(\beta=\delta(X,H)\le r\) and choose a divisorial minimizer.
Through a closed point of its center take successive general smooth
hyperplane sections to a hypersurface \(Y\) of dimension \(2d-2\).
Theorem~\ref{boundary:main} gives \(\delta(Y,H|_Y)>d\).
The complete hyperplane flag estimate applied at the chosen point
then gives \(\beta>N d/(2d-1)>N/2\). In particular the gap above
\(N/2\) exceeds \(N/(4d-2)\), sufficient for the small errors below.

We recall how the free hyperplane tilt of
Section~\ref{boundary:numbers} specializes when one needs only lines.
Set \(u=n/N\), \(b=r/N\), \(\epsilon=r^{-4}\), and put
\[
 x=1,\quad
 U=r+bx+2\epsilon,\quad
 g=u(1-b)-\frac{\log d}{N},\quad z=b+xg-u\epsilon.
\]
In constructing the tilted boundary use
\(x_\beta=\min\{1,(\beta-N/2-3\epsilon)/d\}>0\),
\(\lambda=(d-1)x_\beta\), retaining the upper fraction \(1-1/d\)
of a basis and erasing a free hyperplane through none of the generic
center. The finite projection and free hyperplane estimates of that
section apply to these quantities without a relation between \(n,d\).
In particular the boundary has class \(\lambda c H\) and valuation
order at most \((\lambda+x_\beta+\epsilon)S_H\), and its tilted
threshold \(\beta'\) satisfies
\[
 N/2<\beta-(\lambda+x_\beta+\epsilon)\le\beta'
 \le\beta-\lambda-x_\beta(1-b)+\epsilon<r-\lambda c.
\]
For the last strictness, \(c<1\) and
\(x_\beta\ge\min\{1,(N/(4d-2)-3\epsilon)/d\}\), so
\(x_\beta(1-b)>\epsilon\) here. Thus its divisorial optimizer \(w\)
has \(A_X(w)\le U S_H(w)\). The residual ample coefficient has
lower bound \(r-\lambda c-u\beta'\ge z\): before substituting the
worst value \(\beta=r\), this bound is \(r-u\beta+x_\beta g-u\epsilon\),
which decreases both before and after the cap at one since
\(0<g<ud\). At \(\beta=r\), \(n\ge4d\) gives \(x_r=1\).

Every positive-dimensional optimizer contains a rational curve
of \(H\)-degree at most \(N/\beta'<2\) by the minimizing-center
estimates of Chapter~\ref{ch:global}; for a center of dimension
\(s\ge2\) the hypothesis \((s+1)\beta'\ge N\) holds, and the
curve-center case is included separately there. Thus that curve is
a line. At a point center normalize \(w(\mathfrak m_p)=1\), giving
\(A_X(w)\ge n\) by the ordinary blowup. Let \(L\) be a common
lower bound for \(\delta_C(X,H)\) on all lines. The following three
comparisons therefore suffice:
\[
        U+\Delta<L,\qquad z-u\Delta>1/2,\qquad
        z+\Delta(n/U-u)>1.                       \tag{*}
\]
Indeed positive centers of \(w\) are already excluded by \(U<L\).
At a point apply the same first positive-center and multiplier-ideal
argument of Section~\ref{boundary:numbers}. The ample coefficient
in the perturbed argument is now strictly above \(1/2\); hence the
Fujino--Hashizume length-one curve has integral \(H\)-degree less
than two. Its lc place would give \(A_X/S_H\le U+\Delta+o(1)<L\).
For finite cosupport the remaining Nadel estimate reads
\(1\ge z-u\Delta+\Delta S_H(w)\) in the limit (take \(z-u\Delta<1\)).
The lower bound \(S_H(w)\ge n/U\) contradicts the final comparison.

To verify these inequalities note that \(u>0.975\) and
\(z<b+u(1-b)<1\). Directly expanding gives
\[
 2N\left(z-\tfrac12-\frac{U}{2N}\right)
 =(d-1)(2ux-1)-x(2\log d+b)-(2Nu+2)\epsilon>0.
\]
For the sign use \(x=1\), \(2u-1>0.95\),
\(0.95(d-1)>2\log d+2\) for \(d\ge10\), and
\((2N+2)\epsilon<1\) since \(r> N/2\).
At \(\Delta_0=(z-1/2)/u>0\), the last expression in (*) equals
\(1/2+(z-1/2)n/(uU)>1\). Moreover
\(\Delta_0<1/(2u)<0.52\) and \(U<r+1+2\epsilon\), so
\(U+\Delta_0<L_*\). A positive \(\Delta\) just below
\(\Delta_0\) gives (*) with \(L=L_*\) and \(1/2<z-u\Delta<1\).

It remains to check this mark on lines. The ordinary extraction of
Chapter~\ref{ch:models} has radial denominator \((n-1/d)/(n-1)\).
If the polar rank is \(q\), the normal quotient has \(q-1\) positive
splitting degrees summing to \(d-1\), so its angular denominator is
at most \(2-(q-1)/d\). The radial denominator is at most \(N/L_*\):
indeed \((n-1/d)/(n-1)<1+1/(n-1)\), while
\(N/(r+2)=1+(d-3)/(r+2)>1+1/(n-1)\) for \(d\ge10\).
The angular denominator is at most \(N/L_*\) whenever \(q\ge\ell\).
Thus \(\delta_C(X,H)\ge L_*\) for all lines and uniform K-stability
follows from (*). Furthermore
\[
       n\le d^2/3-4\quad\Longleftrightarrow\quad
       \frac{3N}{d+3}\le d-3
       \quad\Longleftrightarrow\quad N/(1+3/d)\ge L_*.
\]
This yields \(\ell\le d-2\) on that interval.

For a fixed line containment imposes \(d+1\) independent equations.
Its \(n\) first normal polars are independently prescribed elements
of the \(d\)-dimensional binary space. Rank at most \(\ell-1\)
therefore has additional codimension \((n-\ell+1)(d-\ell+1)\).
The Grassmannian of lines is projective of dimension \(2n\), giving
the closed image and the stated codimension. At \(\ell=d-2\) the
display is \(n-2d+10\), and decreasing \(\ell\) increases it.

For the last assertion the normal sequence of a line \(C\) gives
\(N_{C/X}\simeq\bigoplus_i\mathcal O_C(1-a_i)\).
Thus \(H^1(N_{C/X})=0\) exactly when every \(a_i\le2\).
Equivalently its tangent dimension \(h^0(N_{C/X})\) has the
expected value \(2n-d-1\), by Riemann--Roch. The Fano scheme is
the zero scheme of a section of a rank-\(d+1\) vector bundle on
the \(2n\)-dimensional Grassmannian, hence each component has at
least this dimension. This establishes the equivalence with smoothness
of the expected dimension. Under this condition every angular
denominator is at most \(1+2/d\), which is at most \(N/L_*\)
exactly when \(2N/(d+2)\le d-3\), equivalently
\(n\le(d^2-d-8)/2\). The additional range follows from (*).
For completeness, on a fixed line the exceptional polars are the
tuples \(f_i\) for which multiplication by the two binary variables
fails to span \(H^0(C,\mathcal O_C(d))\). This is the cokernel on
sections in the normal sequence. Equivalently some nonzero functional
\(t\) on this space annihilates all \(uf_i,vf_i\).
Its two contractions on degree \(d-1\) are independent: in the
dependent case, \(t\) annihilates all multiples of a linear form,
so is evaluation at its unique zero, contradicting basepoint-freeness
of the polars. For fixed projective \(t\) the conditions on the
tuple therefore have codimension \(2n\); \([t]\) varies in dimension
\(d\). Adding the \(d+1\) containment conditions and subtracting
\(2n\) for the choice of line proves codimension at least one.
The rank conditions and projectivity of the incidence give closedness.
\end{proof}

\subsection{A larger range from the scheme of lines}
\label{sec:nov27-optimized-lines}
Let \(d\ge10\), put \(D=d-1,\ N=n+1,\ r=N-D,\ u=n/N,\ b=r/N\),
and suppose \(n\ge 2d-2\). The following numerical test uses only
the local thresholds on lines. Set
\[
 x_0=\frac{D}{2uD-2\log d-b},\qquad
 L_0=(r+b x_0)\left(1+\frac1{2n}\right).
\]
Assume the denominator of \(x_0\) is positive and
\(x_0<(r-N/2)/d\). If every line has \(\delta_C(X,H)\ge L\) for
a common constant \(L>L_0\), then \(X\) is uniformly K-stable.

In particular, if the scheme of lines is smooth of its expected
dimension, uniform K-stability holds throughout the explicit range
\[
 4d\le n\le \frac{d^2-6}{2}-\frac{(d+2)\log d}{2(d-2)}.
                                                        \tag{4}
\]
For example, in degree ten this includes dimensions forty through
forty-five. The upper endpoint from Theorem~\ref{thm:nov27-line-only-range}
in that degree was forty-one. A further explicit range, allowing
dimensions below \(4d\), is obtained by putting
\[
  c_d=\frac{d-1}{2(d-1)-2\log d-2}.
\]
Under the same hypothesis on the scheme of lines one may take
\[
  2d-3+2dc_d\le n\le\frac{d+2}{2}(d-3/2-c_d)-1 .
                                                        \tag{5}
\]
The numerical test itself gives every \(32\le n\le46\) in degree ten.
The logarithm is natural.

\begin{proof}
Use the point and positive-center argument for the scheme of lines
in Theorem~\ref{thm:nov27-line-only-range}. The pointed section of
dimension \(2d-2\) gives \(\beta>Nd/(2d-1)>N/2\). Fix a rational cap
\(x\) strictly between \(x_0\) and \((r-N/2)/d\) and use the actual
share \(x_\beta=\min\{x,(\beta-N/2-3\epsilon)/d\}\) in that proof.
The initial gap is uniform for this fixed \(n,d\). Thus for sufficiently
small \(\epsilon>0\), \(x_\beta>0\), \(x_\beta(1-b)>\epsilon\) and
\(x_r=x\). The ample-remainder estimate still takes its minimum
at \(\beta=r\), as the derivative on the nonconstant branch is
\(-u+[u(1-b)-\log(d)/N]/d<0\). The three strict tests (*) of the
line theorem therefore apply. At \(\epsilon=0\) their quantities are
\[
 U=r+bx,\qquad z=b+x\{u(1-b)-\log(d)/N\},\qquad
 2N(z-1/2)-U=x(2uD-2\log d-b)-D .
\]
At \(x=x_0\) the right hand side is zero. At
\(\Delta_0=(z-1/2)/u\), the first expression \(U+\Delta_0\) there
equals \(L_0\), and the point expression is
\(1/2+(z-1/2)N/U=1\). Take \(x>x_0\) sufficiently close.
The point expression then strictly exceeds one and \(U+\Delta_0<L\).
Decreasing \(\Delta_0\) slightly makes the ample remainder strictly
greater than \(1/2\) while preserving the other two strict bounds.
All comparisons persist at small positive \(\epsilon\) and with
rational parameters by continuity. This gives the numerical test.

For the explicit interval use instead \(x=1\). When \(n\ge4d\)
this cap is permissible and \(2N(z-1/2)>U\), by the verification
in the line theorem. At zero error one has \(U=r+b<n\) and
\(1-z=(1-b+\log d)/N>0\). The value of \(\Delta\) at equality in
the point test is
\[
 \Delta_1=\frac{1-z}{n/U-u}
          =\frac{U(1-b+\log d)}{n(D-b)}.
\]
The inequality just obtained gives \(0<\Delta_1<\Delta_0\).
Also \(D-b>d-2\), so
\[
 U+\Delta_1<r+b+\frac{1-b+\log d}{d-2}
             <r+1+\frac{\log d}{d-2}.
\]
The ordinary line estimate for a smooth scheme of the expected
dimension is at least \(Nd/(d+2)\). Its angular bound follows
from \(a_i\le2\) in the normal bundle calculation of the line
theorem. Its radial denominator is \((n-1/d)/(n-1)<1+2/d\).
Condition (4) is exactly \(Nd/(d+2)\ge r+1+\log(d)/(d-2)\).
Thus one can increase \(\Delta_1\) slightly while keeping all three
tests strict. Small positive errors again preserve the inequalities.

Finally \(0<c_d<1\) for \(d\ge10\), as \(d-1>2\log d+2\).
For \(n\ge2d-2\) one has \(2uD-b=2D-1-D/N>2D-2\),
so \(x_0<c_d\). The lower inequality in (5) gives
\(c_d\le(r-N/2)/d\). Also \(r+bx_0<n\), whence
\(L_0<r+c_d+1/2\). The upper inequality in (5) is precisely
\(Nd/(d+2)\ge r+c_d+1/2\), proving the claim. At \(d=10\),
\(32\le n\le46\), use \(2uD-b=17-9/N\ge17-9/33\) and
\(\log 10<7/3\) to get \(x_0<3/4\le(r-N/2)/10\).
Here \(b\le38/47<5/6\), so \(L_0<r+9/8\), whereas
\(10N/12-r=9-N/6\ge7/6\). This proves the entire stated interval.
\end{proof}

\subsection{A polar-rank calculation for every smooth conic}
\label{sec:nov27-all-conic-normal-rank}
Let \(C\) be an integral conic on a smooth degree-\(d\) hypersurface
\(X^n\), \(n\ge4\), and let \(\rho\) be the dimension of the span
of all ambient first derivatives of its equation restricted to \(C\).
Put \(e=\min\{n-1,\rho\}\). Then
\[
             \delta_C(X,H)\ \ge\ \frac{n+1}{2-e/(2d)}.
\]
Indeed \(C\) is a smooth plane conic, and all its restricted polars
have no common zero by smoothness of \(X\). A general \((n-1)\)-dimensional
space of directions complementary to the three plane directions maps
onto an \(e\)-dimensional space of polars and has no common zero on
\(C\). The first assertion follows from ordinary linear algebra and
openness of being a complement. For the second, the condition that
all members vanish at a fixed point on \(C\) has codimension \(n-1\)
in the Grassmannian, which exceeds the dimension of \(C\).

Choose such normal directions, write \(u\) for three plane coordinates,
\(z_1,\ldots,z_{n-1}\) for the complementary normals, and \(g(u)=0\)
for the conic equation. The plane restriction of the hypersurface
equation belongs to \((g)\) in every case. Adjoin a generator \(w=g\)
of section degree two. Assign primary weights zero to \(u\), one to
each \(z_i\), and two to \(w\). The two initial equations are
\[
          g(u)=0,\qquad \sum_{i=1}^{n-1}f_i(u)z_i=0,
\]
where the chosen \(f_i\) restrict to the basepoint-free normal polars.
Here a simple factor of \(g\) in the plane restriction contributes
weight two, so it also disappears from the weight-one relation.
The two equations form a complete intersection with the ordinary
Hilbert series of \(X\). Its scheme over the punctured affine conic
cone is an affine bundle, and its fiber over the vertex has dimension
at most \(n\), less than the complete-intersection dimension \(n+1\).
Hence the scheme is integral, since it is generically reduced and has
no embedded components. Localizing at the generic point of the conic
eliminates one normal coordinate. This identifies the indicated weight
filtration with the monomial divisorial extraction of weights
\((1^{n-2},2)\) in normal parameters of \(X\) there; equality of the
ordinary Hilbert series identifies its full initial spaces.

The discrepancy is \(n\) and its expected weight in units of \(n+1\)
is \(n-1/d\), by the weighted Hilbert series. The angular computation
is the one for the conic on a contained plane in Chapter~\ref{ch:models}.
We spell out why the change of primary weights preserves its bound.
Within each weight split out the free \(w\)-powers, and then diagonalize
the split quotient bundle
\(\mathcal E=\mathcal O_C^{n-1}/\mathcal O_{\mathbf P^1}(2-2d)
=\bigoplus\mathcal O_{\mathbf P^1}(a_i)\) determined by the chosen
polars. The spaces in section degree \(m\), \(w\)-power \(s\), and
split-normal multi-index \(\gamma\) have coefficient dimension
\(2p+1+\sum_i a_i\gamma_i\) for \(p=m-2s-|\gamma|\ge d-1\).
This is the same coefficient calculation as in that chapter. The
remaining bounded strip contributes \(O(m^{n-1})\) sections of
angular ordinary degree at most \(m\). Tests on the smooth charts
of the punctured cone give them vanishing \(O(mA)\) at a valuation
of discrepancy \(A\), using the degree bound for affine polynomials.
Summing over the primary weights imposes only \(2s+|\gamma|=J\) in
one block, so all simplex means coincide with the previous conic
calculation. In units of \(n+1\) the mean exponents are \(1/2\) on
\(w\) and \(1+a_i/(2d)\) on the split normal coordinates. The
diagonal degenerations are homogeneous for the present weights and
yield these coordinate monomials on the whole weighted angular
fiber \(\mathbf P(1^{n-2},2)\). As in the preceding angular
calculations one uses a common complex basis for the tested filtration
and the total flag refining the weight by the coordinate exponents.
Checking the resulting coordinate divisor on the smooth punctured
cover also accounts for the order-two stack chart. Weighted
adjunction and the strict radial comparison give a denominator at
most \(1+\max_i a_i/(2d)\). Exactly \(e-1\) summands of the quotient
are positive, with total \(2d-2\), so \(\max_i a_i\le 2d-e\).
This proves the assertion.

\subsection{An explicit open set in dimensions between \(2d\) and \(4d\)}
\label{sec:nov27-middle-range}
Suppose \(d\ge10^4\) and \(n=2d+k\) with \(0\le k\le2d\).
Let \(N=n+1\), \(r=n+2-d\), and set
\[
      L=r+4,\qquad \ell=1+\left\lceil d(2-N/L)\right\rceil,
      \qquad e=\left\lceil 2d(2-N/L)\right\rceil.
\]
If every line in \(X\) has normal-polar rank at least \(\ell\) and
the ambient polars have rank at least \(e\) on every integral conic,
then \(X\) is uniformly K-stable. The complement of the rank conditions
is contained in a closed subset of the smooth parameter space of
codimension at least
\[
 \min\big\{(n-\ell+1)(d-\ell+1)+d+1-2n,\quad
         (n-e)(2d-e)+2d-3n-1\big\}.
\]

Indeed \(2-N/L=(k+11)/(d+k+6)\). By the ordinary line estimate and
Section~\ref{sec:nov27-all-conic-normal-rank}, the indicated ranks give
\(\delta_C(X,H)\ge L\) on every integral curve of degree one or two.
For the line radial check, \((n-1/d)/(n-1)<1+1/d<N/L\) on the stated
range. One also has \(e\le n-1\) since \(d\ge10^4\) and
\((k+11)/(d+k+6)\le(2d+11)/(3d+6)<0.668\).

Suppose \(\beta=\delta(X,H)\le r\). As in the proof of
Theorem~\ref{thm:nov27-line-only-range}, a smooth pointed flag to
dimension \(2d-2\) through a point on the center of a divisorial
minimizer gives \(\beta>N d/(2d-1)>N/2\). Thus in the tilt with
cutoff \(N/3\) there is always a strict initial gap. Use
\[
 \begin{gathered}
 u=n/N,\qquad b=r/N,\qquad\epsilon=r^{-4},\qquad
 x=(r-N/3-3\epsilon)/d,\\
 U=r+bx+2\epsilon,\qquad z=b+x\{u(1-b)-\log(d)/N\}-u\epsilon .
 \end{gathered}
\]
These are the formulas of Section~\ref{boundary:numbers}. In the
construction one uses \(x_\beta=(\beta-N/3-3\epsilon)/d>0\);
\(x_\beta(1-b)>\epsilon\) by the preceding lower bound for \(\beta\).
The resulting tilted optimizer satisfies \(A_X\le U S_H\) and the
ample coefficient satisfies \(r-\lambda c-u\beta'\ge z\), since the
expression before substituting \(\beta=r\) is decreasing in \(\beta\).
Hence the three sufficient tests are \(U+\Delta<L\),
\(z-u\Delta>1/3\), and \(z+\Delta((n+q)/U-u)>1\) at a putative
point optimizer with \(A_X=n+q\), as proved there. No pointed
hyperplane hypothesis in the preceding dimension is needed once
the displayed strict initial bound for \(\beta\) is available.

For \(0\le k\le\sqrt d/8\) the point-envelope calculation on the
growing-width strip above proves all three tests, even with the
smaller common mark \(r+(d-1)/(n+2)<L\). For larger \(k\) the
third test already holds with \(q=0\). In fact
\[
\begin{split}
 3N\left(z-\frac13-\frac{2U}{3N}\right)
   &=(d-1)(3ux-1)-x(3\log d+2b)-(3Nu+4)\epsilon,\\
 3x&=1+(2k+5)/d-9\epsilon/d.
\end{split}
\]
Since \(N\ge2d+1\), the first term is
\((d-1)u(2k+5-9\epsilon)/d-(d-1)/N>1.99k+4\), and
\((3N+4)\epsilon<0.01\). For \(k\le d/2\) one has \(3x<2.01\);
for \(\sqrt d/8<k\) the bound \(1.99k>2.01\log d\) follows since
\(\log(d)/\sqrt d\) decreases from its value at \(10^4\).
The right side of the first identity is therefore positive, using
\(2b/3<1\). For \(k>d/2\) one has \(3x<5.01\), and
\(1.99k>0.995d>5.01(\log d+1)\), which proves the same strict
comparison. Take \(\Delta\) just below \((z-1/3)/u\). At that
endpoint the last test with \(q=0\) is
\(1/3+(z-1/3)n/(uU)>1\), precisely the comparison just proved.
The endpoint is positive and smaller than \(1.34\): from
\(x<2N/(3d)\) we have \(z<5/3\) and \(u>0.9999\).
Also \(x<1.67\) and \(U<r+1.68\). The first test thus holds with
strict room, and the small decrease makes the second strict as
well. This proves the individual-member statement.

The line incidence count is as in Theorem~\ref{thm:nov27-line-only-range}.
A fixed smooth conic \(C\) imposes \(2d+1\) independent conditions
on the hypersurface equation. Independently its \(n-1\) plane normal
polars in a fixed complement are arbitrary degree-\(2d-2\) binary
forms, because restriction of degree-\(d-1\) plane forms to \(C\)
is surjective. If all ambient polars have rank at most \(e-1\), this
normal tuple has rank at most \(e-1\), of determinantal codimension
\((n-e)(2d-e)\). The space of smooth conics has dimension
\(3(n-1)+5=3n+2\). Subtract and take the closure of the incidence
image to obtain the codimension asserted. It is positive, since
\(\ell\le0.668d+2\) and \(e\le1.336d+1\) whereas \(n\ge2d\).

\section{A local threshold from a regular system of tangent forms}
\label{sec:nov27-pow-regular-point}
Let \(X_d^n\) be smooth, \(d\ge5\), \(r=n+2-d\), \(k=d-2\),
\(N=n+1\). At \(p\in X\) write \(q_2,\ldots,q_d\) for the homogeneous
terms of the tangent hyperplane section on \(V=T_pX\). Assume that
they form a regular sequence. Suppose there are vector spaces
\(L\subset U\subset V\), of dimensions \(k,k+1\), such that
\[
 q_2|_U=0,\qquad U\cap\operatorname{rad}(q_2)\subset L,\qquad
 \mathbf P(L)\cap V(q_3,\ldots,q_d)=\varnothing.             \tag{R}
\]
In particular (R) holds whenever \(\operatorname{rank}(q_2)\ge2d-2\).
Indeed the smooth locus of its quadric is homogeneous under the
orthogonal group of the quotient together with the radical shears.
There exists an isotropic \((k+1)\)-space disjoint from the radical.
A general translate of its isotropic \(k\)-subspace misses the
common zero set there. Its projective dimension is \(k-1\) while
the latter set has codimension \(k\) on the quadric by regularity.

More generally put \(R=\operatorname{rank}(q_2)\), \(K=\operatorname{rad}(q_2)\),
\(s=\max\{0,d-1-\lfloor R/2\rfloor\}\) and
\(B=V(q_2,\ldots,q_d)\subset\mathbf P(V)\). For \(R\ge2\), (R)
holds if \(\dim K\ge s\) and an \(s\)-dimensional subspace \(J\subset K\)
misses \(B\) projectively. For \(s>0\) the latter condition holds
precisely when \(\dim(B\cap\mathbf P(K))+s<\dim K\), with an empty
intersection satisfying the dimension test by convention.
To prove the assertion project \(B\) away from \(\mathbf P(J)\).
Its image has dimension at most \(n-k-2\) on the quotient quadric
in \(\mathbf P(V/J)\), hence has codimension at least \(k-s\).
For \(s<k\) choose an isotropic \((k-s)\)-space disjoint from
\(K/J\) and extend it to an isotropic \((k-s+1)\)-space still
disjoint from \(K/J\). Both exist since \(2(k-s+1)\le R\).
On the smooth locus of the quotient quadric the orthogonal group
and radical shears act transitively. By incidence dimension,
a general translate of the projectivized \((k-s)\)-space misses
the image of \(B\), since that plane has dimension \(k-s-1\).
Its inverse image and that of its extension give \(L,U\).
If \(s=k\), take \(L=J\) and extend by one nonradical isotropic
vector. This proof includes \(R=2\), for which \(s=k\).

Let \(A=(A_3,\ldots,A_d)\) be uniform on the standard probability
simplex of dimension \(k-1\), and put
\[
 P(t)=\mathbb P\left(\sum_{i=3}^d A_i/i\ge1-1/t\right),\quad t>0.
\]
Set
\[
 D_d=\max\left\{\frac32,\frac{2(k+1)}{d}\int_0^{3/2}
 t^k P(t)\big[k-(k-1)P(t)^{1/(k-1)}\big]\,dt\right\}.
\]
Writing \(H_{d-1}^{\rm harm}=\sum_{j=1}^{d-1}1/j\), the local
threshold, including divisors whose centers contain \(p\), satisfies
\[
 \delta_p(X,H)\ge\min\left\{
 \frac{nN}{n+H_{d-1}^{\rm harm}},\,\frac{N}{D_d}\right\}.          \tag{1}
\]
The first entry always exceeds \(r\) since
\(n(d-1)>rH_{d-1}^{\rm harm}\).
If (R) and regularity hold everywhere and \(rD_d<N\), then
\(X\) is uniformly K-stable. The following upper endpoints of
dimension satisfy this strict inequality; the same holds at every
smaller Fano dimension in which the geometric hypotheses are met:
\[
\begin{array}{c|rrrrrrrr}
 d&5&6&7&8&9&10&11&12\\ \hline
 n&10&13&16&18&21&24&26&29\\
 d&13&14&15&16&17&18&19&20\\ \hline
 n&31&33&36&38&41&43&45&48
\end{array}
\]
These are memberwise sufficient conditions.

\paragraph{Degeneration and the residual series.}
For the ordinary blowup of \(p\), the residual series in section
degree \(m\) and vanishing order \(j\) is
\([(q_2,\ldots,q_d)^{(j-m)_+}]_j\).
Indeed write the full equation as \(T^{d-1}z+
\sum_{i=2}^dT^{d-i}f_i(y,z)\), with \(f_i(y,0)=q_i(y)\).
The successive hypertangent symbols \(Z_{i-1}\) have section
degree \(i-1\), order \(i\), and initial relations
\(TZ_{i-1}+q_i=0\). Regularity of the \(q_i\) makes \(T\) a
nonzerodivisor on this presentation. Inverting \(T\) identifies the
symbols in the ordinary graded ring of the parameters \(y/T\);
the ordinary Hilbert series \((1-t^d)/(1-t)^{n+2}\) then proves
equality in each degree. Substituting \(Z_{i-1}=-q_i/T\) gives
the residual ideal power: products with at least \(j-m\) factors
and arbitrary remaining tangent monomials correspond exactly to
nonnegative exponents of \(T\) in the original section degree.
The regular presentation also gives ordinary mean
\((n+H_{d-1}^{\rm harm})/N\), by the degree simplex of its symbols.

Choose coordinates \(a_1,\ldots,a_k,s\) on \(U\), with \(L=(s=0)\).
By (R) a complementary light direction \(z'\) pairs to \(s\) under
\(q_2\). That direction has weight one, the other complementary
directions have weight two, and the coordinates on \(U\) have weight
zero. Initializing for least weight and then for least \(s\)-degree
gives \(sz',f_3(a),\ldots,f_d(a)\), where \(f_i=q_i|_L\).
These final forms are a regular sequence since \(f_i\) have no
common projective zero. The intermediate sequence is regular as
well, by the same filtered argument from its regular initials.
In each representation of a power of the ideal, products of the
initial generators belong to its limit, and the Hilbert functions
of powers of either homogeneous regular sequence agree (use
the corresponding free polynomial presentation). This proves that
the indicated limiting residual ideal is exact. Finite-dimensional
simultaneous adaptation to the valuation and each successive weight
and semicontinuity for the products on projective space transfer
their threshold bounds back to the original residual series.

There are \(\ell=n-k-2\) other tangent coordinates.
Split the characters in all variables outside \(a\). A block of
ordinary core degree \(D\), with exponents \(i,j'\) on \(s,z'\)
and multi-index \(\alpha\) on the \(\ell\) other variables,
is the fixed monomial times
\[
   [(f_3,\ldots,f_d)^e]_D,\quad
   e=(D+u-m)_+,\qquad u=|\alpha|+\max(i,j').
\]
Indeed the usable powers of \(sz'\) are at most
\(\min(i,j')\). The normalized mean exponents, in units \(N/m\),
of \(s,z'\) are both \(3/2\), and those of the \(\ell\) others are
one. To see this, the monomial basis over \(\mathbf C[sz']\) consists
of one pure power of \(s\) or \(z'\); these sectors have equal
asymptotic mass. The \(sz'\) symbol of section cost one contributes
one exponent to each, and the free pure power contributes \(1/2\)
to each. The remaining \(f_i\) have section costs \(i-1\).
The bounded shifts of a free basis of \(\mathbf C[a]\) over
\(\mathbf C[f_3,\ldots,f_d]\) do not affect these
means.

\paragraph{The complete-system upper-tail bound on the core.}
In degree \(D\) the relative dimension of the core block in all
homogeneous polynomials in \(a\) is asymptotically \(P(t)\),
\(t=D/(m-u)\). The full core ring is finite free over the forms
\(f_i\). In each free summand, monomials
\(\prod f_i^{b_i}\) have \(\sum i b_i=D+O(1)\) and the ideal power
requires \(\sum b_i\ge e\). The proportions \(ib_i/D\) tend to
uniform measure on the simplex, giving the asserted probability.
All bounded free-module shifts, lattice faces and blocks of bounded
ordinary degree have negligible mass in the total dimension.
Equivalently the count is the Riemann sum for a bounded polytope
with finitely many simplex cuts. One has \(0\le u\le m\) and
\(0\le t\le3/2\) up to these shifts.

On ordinary \(\mathbf P^{k-1}\), for a test valuation normalized
to discrepancy one, the largest mean vanishing of any fraction
\(P\) of the degree-\(D\) series is at most
\(D(1-(k-1)P^{1/(k-1)}/k)+o(D)\). This is the complete-subseries
upper-tail estimate proved above by concavity and
\(\delta(\mathbf P^{k-1},\mathcal O(1))=k\).
For completeness, the normalized survival volume \(F(x)\) on
\([0,1]\) has \(F^{1/(k-1)}\) concave on its positive support
and \(\int F\le1/k\); hence each integrated upper tail is bounded
by that of \((1-x)^{k-1}\), using the single sign change.
Minimizing \(z+P^{-1}\int_z^1(1-x)^{k-1}dx\) proves the formula.
Here is a finite-degree estimate uniform in both the valuation and the
chosen subspace. Put \(p'=k-1\), \(h_D=\binom{D+p'}{p'}\), and let
\(F_D(x)\) be the fraction of the complete degree-\(D\) system with
order at least \(Dx\). For \(D>0\), at continuity points one has
\[
 F_D(x)\le \frac{(D+p')^{p'}}{p'!\,h_D}\,
                F\left(\frac{Dx}{D+p'}\right).
\]
To see this, choose affine coordinates from a projective chart meeting
the generic center of the valuation. A basis of the filtered subspace
can be chosen with distinct leading monomials of exponents \(\gamma\)
in a fixed monomial order. Write its polynomials as \(s_\gamma\).
In degree \(q(D+p')\), homogenize \(s_\gamma^q x^\nu\) for
\(0\le\nu_i<q\). Each has order at least \(qDx\), and their leading
exponents \(q\gamma+\nu\) are pairwise distinct. Dividing the count
\(h_D F_D(x)q^{p'}\) by the complete dimension and taking
\(q\to\infty\) proves the bound. The support is in \([0,1]\) by
\(\alpha(\mathbf P^{p'},\mathcal O(1))=1\).
The prefactor is \(1+O_{p'}(D^{-1})\). Integrating from any
\(z\in[0,1]\) and using the integrated tail comparison above gives
\[
 \int_z^1 F_D(x)\,dx\le\int_z^1(1-x)^{p'}\,dx+O_{p'}(D^{-1}),
\]
uniformly in \(z\). For a subspace of dimension \(P h_D\), the sum
of orders in any basis, divided by \(Dh_D\), is at most
\(Pz+\int_z^1F_D(x)\,dx\). Taking \(z=1-P^{1/p'}\) gives the
claimed upper-tail contribution with an additive error \(O_{p'}(h_D)\)
in the unnormalized sum. Over all nonempty blocks \(\sum h_D=O_n(m^n)\),
since the noncore exponents and \(D\) are \(O_d(m)\). The total error
after division by \(m h^0(X,mH)\) thus tends to zero, uniformly even
for block fractions tending to zero.

Here is the normalized mean of this bound over the original basis.
Scale \(m\) to one. The multiplicity of the noncore indices at \(u\)
is proportional to \(u^{n-k-1}\), since exactly \(2c+1\) pairs
have \(\max(i,j')=c\), with \(\ell\) remaining free indices.
Consequently the joint mass in \(u,t\) is proportional to
\[
 u^{n-k-1}(1-u)^k t^{k-1}P(t)\,du\,dt,\qquad 0<u<1.
\]
In this law \(u,t\) are independent and \(N\mathbb E(1-u)=k+1\).
Moreover
\[
  \int t^{k-1}P(t)\,dt=\frac1k\mathbb E\big(1-\sum A_i/i\big)^{-k}
  =\frac1k\prod_{i=3}^d(1-1/i)^{-1}=\frac d{2k}.
\]
For the expectation identity integrate \(\exp(-\sum (1-1/i)x_i)\)
over the positive orthant, using \(x_i=\rho A_i\).
Multiplying the preceding upper-tail bound by \(N\) and averaging
now gives precisely the integral in \(D_d\). It also dominates
\(N\mathbb E(D)/k\), since \(k-(k-1)P^{1/(k-1)}\ge1\).

Blow up the coordinate vertex \(a=0\) in \(\mathbf P^{n-1}\).
Its radial coefficient \(N\mathbb E(D)\) is at most \(kD_d\)
and its discrepancy is \(k\). The homogeneous core divisors
pull back from \(\mathbf P^{k-1}\) on the resolved projection.
The remaining coordinate hyperplanes and the exceptional divisor
are relatively normal crossings for this locally trivial smooth
projection. On its local product charts the threshold therefore
follows from the preceding bound on the base and the individual
coordinate coefficients at most \(D_d\); this can also be checked
on a product log resolution. The cone vertex has positive codimension
in the setting considered here. Thus the normalized angular
denominator is at most \(D_d\). Combined with the ordinary radial
bound this gives (1) by local basis adjunction.

\paragraph{Exact numerical comparisons.}
For \(v_i=1/i\), the simplex-slice formula gives
\[
 P(t)=\sum_{i=3}^d
 \frac{(v_i-1+1/t)_+^{k-1}}{\prod_{j\ne i}(v_i-v_j)}.
\]
This follows by successive integration over the barycentric simplex
(or the divided-difference formula for its volume). It is one on
\(t\le1\). Partition \([1,3/2]\) into 64 equal intervals and
subdivide at the knots \(i/(i-1)\). On each interval compute the
exact \(t^kdt\)-weighted mean \(\bar P\) by integrating
\(t(1-(1-v_i)t)_+^{k-1}\). The function
\(x\mapsto x(k-(k-1)x^{1/(k-1)})\) is concave on \([0,1]\);
Jensen gives a rational upper bound using its value at \(\bar P\)
with the root rounded downward on the grid \(10^{-8}\).
These outward comparisons are given by
\path{nov27-pow-regular-point.py}.

\subsection{Codimension of the exceptional tangent systems}
\label{sec:nov27-regular-jets-codimension}
Let \(\mathcal P\) be the projective space of equations, and take
\(d\ge4,\ n\ge2d-1\). Outside a closed subset of the smooth locus
of codimension at least \(C(d,n)\) below, the tangent forms at
\emph{every} point are regular and satisfy the isotropic flag
condition (R). Thus the preceding numerical local estimate proves
uniform K-stability on this complement whenever \(d\ge5\) and
\(rD_d<n+1\). For the quartic eightfold use the strict
point-centered flag estimate and the exclusion of positive-dimensional
minimizing centers in Chapter~\ref{ch:eight}.

For each \(2\le R\le2d-3\) put
\[
 \kappa=n-R,\quad s=d-1-\lfloor R/2\rfloor,\quad
 t=\lfloor R/2\rfloor,\quad b=\kappa-s.
\]
Define
\[
 M_R=\binom{\kappa+1}{2}+
 \max\left\{\min_{3\le i\le s+2}\binom{\kappa+2}{i},
             \sum_{i=3}^{t+2}\binom{b+i}{i}\right\}.
\]
Here \(\kappa\ge s\). One may take
\[
 C(d,n)=\min\left\{\binom n2,\ \min_{2\le i\le d}\binom{n+1}{i},
                 \ \min_{2\le R\le2d-3}M_R\right\}-n .
\]
For the following degree and dimension pairs this method gives
\[
 \begin{array}{c|rrrrrrr}
 (d,n)&(4,8)&(5,10)&(6,12)&(6,13)&(7,15)&(7,16)&(8,18)\\ \hline
 C(d,n)&20&35&54&65&90&104&135 .
 \end{array}
\]
For instance the locus of smooth quartic eightfolds that fail uniform
K-stability is contained in a closed subset of codimension at least
twenty in the space of smooth quartics. For \(d\ge5\) this
codimension argument uses the tangent-form bound directly.

To prove the formula, fix the marked point and its tangent
hyperplane. Once the nonzero normal derivative is normalized,
the tangent forms of all degrees \(2,\ldots,d\) can be prescribed
independently. Whenever the first regularity failure occurs at \(q_i\),
that form must vanish identically on some component of the common
zero scheme of the preceding forms in \(\mathbf P^{n-1}\), of
dimension \(n-i+1\). For a fixed integral projective variety of
dimension \(a\), vanishing of a form of degree \(i\) imposes at
least \(\binom{a+i}{i}\) independent conditions, by a finite
linear projection onto \(\mathbf P^a\). For fixed preceding forms
there are only finitely many components. Thus failure of regularity
has codimension at least \(\min_{2\le i\le d}\binom{n+1}{i}\).
Ranks zero or one impose at least \(\binom n2\) independent
conditions in codimension, by the symmetric determinantal formula.

For rank \(R\ge2d-2\), regularity suffices. At each remaining
rank \(R\ge2\) the symmetric determinantal stratum has codimension
\(\binom{\kappa+1}{2}\). Its radical \(K\) is fixed after prescribing
\(q_2\). The restrictions of the higher forms to \(K\) remain
arbitrary and independent. The lower-rank criterion above
applies whenever the common zero scheme of these restrictions has
codimension at least \(s\) in \(\mathbf P(K)\), with the empty
scheme allowed. If this fails, the first \(s\) restrictions fail
regularity. Counting at the first failure, with degrees now
starting at three in a projective space of dimension \(\kappa-1\),
gives the first term in the maximum.

For the other term, on this bad locus follow a chain of components
from \(\mathbf P(K)\) down to a component of the final intersection
of dimension at least \(b=\kappa-s\). At each step the new equation
either contains the chosen preceding component or cuts its dimension
by one on the chosen descendant. There can be at most \(s-1\)
dimension drops. At least \(d-2-(s-1)=t\) of the equations therefore
vanish identically on chosen components of dimension at least \(b\).
At each step the possible preceding component depends only on the
earlier forms and has finitely many choices. The fiber dimension
bound for the independent next equation subtracts at least
\(\binom{b+i}{i}\) at each vanishing step of degree \(i\). One can
make these counts on constructible strata for the finitely many
relative component choices, or induct on the successive choices
using their finite generic fibers. Summing the smallest \(t\)
possible contributions gives the stated lower bound. Finally,
the universal smooth marked-hypersurface family has relative
dimension \(n\) over the parameter space. Projecting the bad
incidence and taking closure loses at most \(n\) in codimension,
as claimed. The integer substitutions are recorded in
\path{nov27-pow-regular-codim.py}.

\subsection{A smooth ternary quadratic block and the fixed cubic factor}
\label{sec:nov27-orthogonal-cubic}
Continue with \(k=d-2,\ N=n+1\), now allowing \(d=4\).
At a point \(p\) assume there are disjoint subspaces \(L,B\subset T_pX\)
with \(\dim L=k,\ \dim B=3\) and
\[
 q_2|_L=0,\quad B\subset L^\perp,\quad q_2|_B
 \text{ nondegenerate},\qquad
 \mathbf P(L)\cap V(q_3,\ldots,q_d)=\varnothing .
\]
Write \(f_i=q_i|_L\). Suppose \(v(f_3)\le c A_{\mathbf P(L)}(v)\)
for all divisorial valuations, using local frames, with \(c\le3\).
One can always take \(c=3\). One can take \(c=1\) when
\((\mathbf P(L),(f_3=0))\) is log canonical. For \(d=4\) the latter
condition means that the binary cubic has three distinct roots.
Write \(D_d\) for the constant of Section~\ref{sec:nov27-pow-regular-point},
using the same integral at \(d=4\). Put
\[
 J_d=\frac1{18}\prod_{i=5}^d\frac{i-1}{3i-9},\qquad
 K_{d,c}=\max\left\{\frac{112}{81},\
  \frac1k\sum_{i=3}^d\frac{i}{i-1},\ D_d-(3-c)J_d\right\}.
\]
Then
\[
 \delta_p(X,H)\ge \min\left\{
 \frac{nN}{n+H_{d-1}^{\mathrm{harm}}},\
 \frac{N}{K_{d,c}}\right\}.
\]
This local statement includes centers containing \(p\).
Consequently if the hypotheses hold everywhere with \(c=1\), uniform
K-stability follows up to dimensions \(9,12,14,16,19\) in degrees
\(4,5,6,7,8\), respectively. For a chosen \(L\), an appropriate
\(B\) exists if
\(\operatorname{rank}q_2-2(\dim L-\dim(L\cap\operatorname{rad}q_2))\ge3\),
by Witt decomposition. These are tests on the individual equation.

\paragraph{The two projective factors.}
Give the coordinates on \(L,B\) weights zero and one, respectively,
and the remaining coordinates weight three. The initial polar forms
are \(Q(x)=q_2|_B,\ f_3(a),\ldots,f_d(a)\), a regular sequence.
The original polar forms are consequently regular, and their powers
have the exact initial powers by the Hilbert-function argument of
the preceding section. Split the free coordinates \(y\) and the
two ordinary degrees in \(x,a\) by their torus characters. In a block
of \(x\)-degree \(E\), \(a\)-degree \(D\), free exponent \(\alpha\)
and section degree \(m\), set \(e=(E+D+|\alpha|-m)_+\),
\(I=(f_3,\ldots,f_d)\). On \(\mathbf P(B)\times\mathbf P(L)\)
the block after removing the free monomial is
\[
 V=\sum_{0\le j\le E/2}
 Q^j\mathbf C[x]_{E-2j}\otimes [I^{(e-j)_+}]_D.                  \tag{2}
\]
Resolve the two linear projections by the projective bundle over
these two factors associated to their two tautological lines and
the trivial bundle of free coordinates. It maps birationally to
the original projective space by taking the vector sum. Locally
on the base, write each \(x\) or \(a\) as its corresponding fiber
coordinate times its angular coordinates. The coordinate fiber
hyperplanes contributed by the two groups have log discrepancies
\(3,k\) over the original projective space, as seen at the respective
generic linear vertices. The free coordinate hyperplanes have
discrepancy one, and all are relatively normal crossings.
In units \(N/m\) their mean multiplicities are, respectively,
\(4,\sum_{i=3}^d i/(i-1),1\), by the count below. The remaining
divisors pull back from the two projective factors. A product
resolution therefore reduces the estimate to (2).

\paragraph{Separating the two estimates in a coupled block.}
On each full space \(\mathbf C[a]_D\) choose a complete flag of
subspaces \(G\) with one-dimensional quotients, as in the proof of
Z. Zhuang, \emph{Product theorem for K-stability}, Adv. Math.
371 (2020), 107250, Theorem~3.1. If the tested prime dominates the
base \(\mathbf P(L)\), use successive evaluations at general points.
Otherwise its restriction to the base is a positive multiple of
a divisorial valuation \(w\). On a model realizing \(w\), first
take its orders, and inside each order take successive evaluations
of its leading coefficients at general points of the base prime.
For each quotient slot \(\nu\) let \(h_\nu\) be the greatest \(h\)
such that \([I^h]_D\) uses this slot. Simultaneous adaptation to
\(G\) and the nested powers gives the \(x\)-space in that slot as
\[
 U_\nu=Q^{j_\nu}\mathbf C[x]_{E-2j_\nu},\qquad
                  j_\nu=(e-h_\nu)_+ .
\]
This is the largest of the nested \(x\)-spaces occurring in its
sum in (2). Empty slots are omitted.

Refine the tested valuation filtration to one-dimensional quotients
on each \(V\) and choose general simultaneous lifts for this and
the base filtration. After dividing out their base orders the
group of \(\dim U_\nu\) lifts specializes to a basis of \(U_\nu\)
over its defining evaluation point. Thus on a general fiber the
threshold of this group's divisor is bounded below by the minimum
over all bases of \(U_\nu\), by lower semicontinuity in the smooth
family \(\mathbf P(B)\) over the base prime (or over the base).
Only this threshold bound is transported to the general fiber.
One may choose a common general fiber for the finitely many
groups, and add their bounds with the dimension weights. The
coefficient of the base prime itself is the weighted sum of its
orders on the slots. If that coefficient divided by discrepancy and
the general fiber threshold denominators are both at most \(K\),
inversion of adjunction gives the bound \(K\) for the tested
valuation: work on the base model, raise the coefficient of the
base prime to one at a slightly smaller multiplier, and apply
inversion along the prime at its generic fiber. This is the
argument in equations (3.4)--(3.8) of Zhuang's proof, with the
lower semicontinuity comparison at each different slot. When the
tested prime dominates the entire base only the general fiber
step is necessary.

\paragraph{The finite-dimensional counts.}
Resolve a nonempty \(x\)-space by its \(Q\)-orders \(j\ge j_\nu\),
where its graded dimension is \(h^0(Q,\mathcal O_Q(l))=2l+1\),
\(l=E-2j\). On the base, the slots of height at least \(e-j\)
have precisely the orders of a basis for
\([I^{(e-j)_+}]_D\) adapted to the base valuation. Thus the
total base coefficient, obtained by summing the orders of the
slots with the indicated \(Q\)-graded multiplicities, is exactly
the corresponding sum of the base-valuation coefficients on
these ideal powers.
Here their exponent is \((D+u-m)_+\), \(u=|\alpha|+j+l\).
The multiplicity density at \(u\) is proportional to
\(u^{n-k-1}\), since there are \(n-k-3\) free exponents,
one \(j\), and one \(l\) of multiplicity \(2l+1\).
Consequently the normalized joint distribution of \(u/m\) and
\(D/(m-u)\) agrees with that in the proof of the preceding
point estimate. The same normalized angular upper-tail integral
therefore controls the base coefficient.

There is a further gain. For \(h=(D+u-m)_+\), every element of
\([I^h]_D\) is divisible by \(f_3^i\) with \(i=(4h-D)_+\),
since each other generator has degree at least four. Write
\(D'=D-3i\) and let \(P\) be the relative dimension in the degree
\(D\) complete system. The residual dimension fraction in degree
\(D'\) is \(P(D/D')^{k-1}\) asymptotically. Adding the bound \(ci\)
for the fixed factor, the complete-system upper-tail expression
for the residual is exactly the earlier expression in degree
\(D\) minus \((3-c)i\). Indeed its nonlinear term is unchanged
because \(D'[P(D/D')^{k-1}]^{1/(k-1)}=DP^{1/(k-1)}\).
The uniform finite-degree estimate of the preceding section applies
also before taking limits, with the binomial dimensions of the two
complete systems; the contribution of degrees \(D'=o(m)\) is
controlled directly by their degree bound.

For explicit means use the free-module monomials
\(\prod_{i=3}^d f_i^{b_i}\) with the bounded homogeneous basis shifts
of \(\mathbf C[a]\) over \(\mathbf C[f_3,\ldots,f_d]\).
Their exact ideal-power heights are \(\sum b_i\) and
\(D=\sum i b_i+O(1)\). Each slot of that height occurs with the
same multiplicity when summed over all slots, regardless of the
chosen base filtration. Adding the \(Q\)-order \(j\) and \(l=E-2j\)
with multiplicity \(2l+1\), the simplex coordinates in section
degree are
\[
 j,\quad l_1,l_2,\quad \alpha,\quad
 ((i-1)b_i)_{3\le i\le d},\quad x_0 ,
 \qquad l=l_1+l_2 ,
\]
up to the bounded shifts. The factor \(2l+1\) has the same leading
density as two copies of the possible splits of \(l\).
Here \(x_0\) is the remaining nonnegative section degree.
These \(N\) barycentric coordinates have normalized mean one.
In particular this gives the fiber hyperplane multiplicities
stated above. The compulsory \(f_3\)-order, up to a bounded error,
is \((b_3-4x_0-\sum_{i=5}^d(i-4)b_i)_+\).
Uniform simplex averages of a degree-one homogeneous expression
multiplied by \(N\) equal its expectation on independent standard
exponential coordinates, by conditioning on their sum. Integrating
the positive tail of the \(b_3\) coordinate of mean \(1/2\) gives
\(N\,\mathbb E(i/m)=J_d\) in the limit. This proves the asserted
base denominator \(D_d-(3-c)J_d\).

\paragraph{The conic on the other factor.}
On \(\mathbf P(B)=\mathbf P^2\), write \(S=S_{\mathcal O(1)}\).
The smooth conic \(Q\) and \(\mathbf P^2\) are K-polystable.
Zhang--Zhou \cite[Theorem~1.7]{ar2003.06839} give K-semistability
of \((\mathbf P^2,tQ)\) for \(0\le t\le3/4\), equivalently
\((3-2t)S(v)+t v(Q)\le A(v)\). Also \(v(Q)\le A(v)\).
For all \(0\le2q\le E\) these inequalities imply
\[
 (E-2q)S(v)+qv(Q)
 \le \max\{E/3,\ (E-2q)/6+q\}\,A(v).                              \tag{3}
\]
If \(q\le E/4\), use \(t=3q/E\); otherwise use \(t=3/4\)
and add the remaining nonnegative multiple of \(v(Q)\le A(v)\).
This gives the required individual bound for every basis of
\(U_\nu\) on the defining evaluation fiber, up to a uniform
\(o(m)\) times discrepancy by the complete-series approximation.
Bounded residual degrees can instead use their ordinary degree
bound.

In the preceding simplex, \(j_\nu=(j-x_0)_++O(1)\), because
\(e-h_\nu=E+D+|\alpha|-m-\sum b_i=j-x_0+O(1)\) whenever positive.
The mean of the first branch of (3) in units \(N/m\) is \(4/3\).
The additional term in taking the maximum is
\[
 \tfrac13(2j_\nu-E/2)_+
 =\tfrac13\big(j-2x_0-(l_1+l_2)/2\big)_+ +O(1).
\]
Its normalized limiting mean is
\(1/[3(1+2)(1+1/2)^2]=4/81\). Hence the general fiber
denominator is at most \(112/81\). This proves the angular bound
and the claimed local threshold by ordinary point adjunction.
The outward rational numerical evaluations are given in
\path{nov27-pow-orthogonal-cubic.py}.

\paragraph{A parameter-space consequence.}
For the following pairs the smooth members failing uniform K-stability
are contained in a closed subset of the indicated codimension or
higher, inside the parameter space of smooth equations:
\[
 \begin{array}{c|rrrrrr}
 (d,n)&(4,9)&(5,11)&(5,12)&(6,14)&(7,16)&(8,19)\\ \hline
 \text{codimension}&19&29&43&64&80&117 .
 \end{array}
\]
In particular the preceding tests hold on nonempty open sets
in these dimensions. To see this fix \(p\) and its tangent
hyperplane as before. Put \(Q_0=(q_2=0)\), \(R=\operatorname{rank}q_2\),
\(K=\operatorname{rad}q_2\), \(\kappa=n-R\), and \(Z=V(q_2,\ldots,q_d)\).
Ranks at most two have codimension \(\binom{n-1}{2}\). Nonregularity
costs at least \(\min_{2\le i\le d}\binom{n+1}{i}\).
For \(R\ge3\) we further assume that the singular locus of the cubic
section of \(Q_0^{\rm sm}\) has codimension at least \(k\) there.
The locus of cubic forms failing this condition, for fixed \(q_2\),
has codimension at least \((n-k)^2\). Indeed, if the singular
locus has dimension at least \(a=n-k-1\), it contains \(a+1\)
independent projective points. At these points the cubic jets are
independent: choose linear forms \(\ell_i\) dual to the points,
and vary cubics \(\ell_i^2 u_i\) with \(u_i\) arbitrary linear.
The restriction jets on the smooth quadric are therefore independent
as well. Prescribing a singular point of the cubic section at a
fixed point imposes \(n-1\) equations. Varying the ordered independent
points costs at most \((a+1)(n-2)\) parameters, and each such
cubic with a singular \(a\)-dimensional component has an
\((a+1)a\)-dimensional family of choices. This gives \((a+1)^2\).

Take \(s=\max\{0,k-\lfloor(R-3)/2\rfloor\}\) and \(h=k-s\).
The dimensions in the table satisfy \(\kappa\ge s\). An isotropic
\(k\)-plane with radical part \(J\) of dimension \(s\) leaves
at least three nondegenerate orthogonal directions. A general
such \(L\) misses \(Z\) away from the radical and cuts a smooth
cubic on its own nonradical locus. Indeed the smooth quadric locus
is homogeneous under orthogonal transformations and radical shears;
translate an initial \(L\) and apply transversality to the smooth
part of the cubic and dimension to its singular part and to \(Z\).
If \(s=0\), this constructs the required \(L\).

For \(s>0\) impose also the following three conditions on
\(\mathbf P(K)\), allowing empty bad loci. First require
\(\operatorname{codim}(Z\cap\mathbf P(K))\ge s\), which lets a
general \(\mathbf P(J)\) miss \(Z\). Next require the locus where
the full differential \(dq_3\) in \(T_pX\) vanishes to have
codimension at least \(s\) in \(\mathbf P(K)\). Finally, if \(h<s\),
require the singular locus of \(q_3|_K\) to have codimension at
least \(s-h\) in \(\mathbf P(K)\). Here if the cubic restricts
identically to zero, its entire space is counted as singular.
We verify smoothness of \(q_3|_L\) at the remaining points of
\(\mathbf P(J)\). A general \(J\) meets the smooth part of the
nonzero restricted cubic transversely. It misses the full-gradient
zero locus. Among its points where the restricted cubic has zero
differential, the dimension is at most \(h-1\): for \(s\le h\)
this follows from \(\dim\mathbf P(J)=s-1\); for \(s>h\) it
follows from the last imposed condition. At such a point \(y\)
the nonzero covector \(dq_3(y)\) factors through \(T_pX/K\).
The image of \(L\) in that quotient is a general isotropic
\(h\)-plane \(P\), with \(2h+3\le R\). For a given nonzero
covector, the locus of such \(P\) in its kernel has codimension
at least \(h\) in the isotropic Grassmannian. For a nonisotropic
normal vector this is the dimension difference between the
isotropic Grassmannians in dimensions \(R\) and \(R-1\).
For an isotropic normal vector \(v\), the locus with \(v\notin P\)
projects to the isotropic Grassmannian in dimension \(R-2\)
with \(h\) lift parameters, again giving codimension \(h\).
The locus with \(v\in P\) is the Grassmannian of isotropic
\((h-1)\)-planes in dimension \(R-2\) and has codimension
\(R-h-1>h\). Thus a general \(P\) gives a nonzero derivative on
\(L\) at all the remaining points, by incidence dimension.
For \(h=0\) there are no such points. The conditions used on and
off the radical are simultaneous open conditions in the irreducible
orbit family (the inequality \(2h+3\le R\) avoids maximal even
orthogonal Grassmannians).

For clarity, here are the bounds used to count failures of these
three conditions after fixing the rank-\(R\) quadratic. Write
\(b_0=\kappa-s\), \(t=d-1-s\), and put
\[
 A_0=\max\left\{\min_{3\le i\le s+2}\binom{\kappa+2}{i},
                    \sum_{i=3}^{t+2}\binom{b_0+i}{i}\right\}.
\]
The component-chain count of
Section~\ref{sec:nov27-regular-jets-codimension} gives \(A_0\) for
the first failure. For the second, use \(b_0+1\) independent
points in the full-gradient zero locus. The same independent
ambient cubic jets now impose \(n\) equations per fixed point.
The points move in \(\mathbf P(K)\) and their choices on the bad
component have dimension at least \((b_0+1)b_0\). This gives
\[
 A_1=(b_0+1)(R+b_0+1).
\]
When \(s=1\), take instead
\(A_1=\binom{\kappa+2}{3}+R\binom{\kappa+1}{2}\):
the full gradient vanishes identically on \(K\) exactly when
the cubic contains no term with fewer than two nonradical
variables. If \(h<s\), the third failure gives
\(A_2=(\kappa-s+h+1)^2\), by the independent jet count applied
on \(\mathbf P(K)\); otherwise omit \(A_2\).
To the minimum of these applicable \(A_i\) add
\(\binom{\kappa+1}{2}\), the codimension of the prescribed
quadratic rank. Minimize these numbers for \(3\le R<2k+3\)
together with the preceding bounds for small rank, nonregularity,
and the singular locus on \(Q_0^{\rm sm}\).
Subtract \(n\) to forget the marked point and take closure.
These finite substitutions produce the table, as reproduced by
\path{nov27-pow-orthogonal-codim.py}.

\section{Point centers and descent in dimension \texorpdfstring{\(2d+1\)}{2d+1}}
\label{sec:nov27-third-descent}
Let \(d\ge5,\ n=2d+1,\ N=2d+2,\ r=d+3\). For every smooth
complex degree-\(d\) hypersurface \(X\) in this dimension,
\[
             \delta_C(X,H)>r
\]
for every line and every integral conic \(C\). Every divisorial
minimizer of \(\delta(X,-K_X)\le1\) therefore has a closed point as
its center. Consequently \(X\) is uniformly K-stable under either
of the following hypotheses: a finite subgroup of \(\operatorname{Aut}(X)\)
has no common fixed point on \(X\), or \(X\) has a model without rational
points over a subfield of \(\mathbf C\).

\begin{proof}
Write \(b_2=1+(d-1)/(d+3)\), so \(N/b_2=r\).
A normalized short-curve bound on a smooth hyperplane section through
the span of \(C\) lifts to \(X\) by the complete hyperplane flag
of Section~\ref{sec:nov27-uniform-short-mark}. For the quintic,
this lifts \(385/51\) in dimension ten to
\((12/11)(385/51)=420/51>8\) in dimension eleven. For \(6\le d\le19\),
the sharpened bounds in Chapter~\ref{ch:secondrow} have denominator
strictly smaller than \(b_2\). For the sextic this is
\(29/19<14/9\), for the septic it is \(49/31<8/5\), and for the
octic it is \(13/8<18/11\). In every remaining entry of the two
small-degree tables one has \(T<d/3\).

Let \(d\ge20\). The third-row net bounds in
Section~\ref{sec:nov27-low-third-row} and
Corollary~\ref{cor:nov27-third-rank} give a threshold strictly
greater than \(r\) at higher-rank lines and at higher-rank conics
in contained planes. Here one uses just the local estimates in the
proof of the latter corollary. Proper conic planes have denominator
at most \(2d/(d+2)<b_2\), by Chapter~\ref{ch:models}, and polar
pencils in contained planes have denominator strictly smaller than
\(b_2\) by Section~\ref{sec:nov27-mixed-conic-pencil}.

It remains to make the bound in Section~\ref{sec:nov27-all-pencils-2d}
strict for polar-pencil lines. One may work in dimension \(2d\)
by taking a general hyperplane through the line. If rank two drops
under restriction, its surviving one-dimensional polar system has
a base point, contradicting smoothness along the line. When \(3\mid d\),
Section~\ref{sec:nov27-divisible-second} already gives the strict
bound with \(T=d/3-1/6\) for all \(d/3\ge3\).

If \(3\nmid d\), choose \(T<d/3\) sufficiently close to \(d/3\).
We verify the argument of Section~\ref{sec:nov27-integer-second}
with target \(1+(1-1/d)T/(T+1)<b_2\). The integer cutoffs at
\(T,2T\) are unchanged, because neither \(d/3\) nor \(2d/3\)
is an integer. For nonterminal orders \(j\le d-1\) one has
\(j<3T\). When \(j>2T\) one has \(d-j<T\), as needed for the
constant non-lc locus. For integers \(s>T,\ s\le j\), one has
\(3s>j,\ 2s>d-j\); when \(2s\le j\) one also has \(s>d-j\).
These inequalities give precisely the previously proved repeated
hyperplane and quadric reductions. At orders where
\(2+2(d-j)\ge d\), one has \(j\le(d+2)/2<2T\); hence any
high-multiplicity component there is a unique hyperplane.
At terminal order use four lights as in
Section~\ref{sec:nov27-integer-second}, since \(4T>d\).

In the intermediate crossing calculation take \(h\le(d+2)/2\)
and \(t>h\). Whenever \(1/h+1/t>3/d\), the bound \(1/h+1/t\ge1/T\)
persists upon choosing \(T\) sufficiently close to \(d/3\), since
there are finitely many such pairs. If \(t=d\), use the terminal
case directly. If \(1/h+1/t\le3/d\) with \(h>d/3\), then
\[
      t\ge\frac{dh}{3h-d}\ge\frac{d(d+2)}{d+6}>d-4 .
\]
After the at most two cuts producing this order, its ambient
codimension is \(c\le4\); hence \(c+2(d-t)\le10<d\), and the
high-order iteration applies. If \(h\le d/3\) the one-light
test already works. In the remaining crossings the graph exponents
use \(t\le d-1<3T\) and \(1/h+1/t\ge1/T\), as in the cited
integer argument. This proves the strict pencil bound.

For the global assertion put \(\beta=\delta(X,H)\le r\). A general
pointed hyperplane section belongs to Theorem~\ref{thm:nov27-second-row}.
Hyperplane adjunction gives \(\beta>N(d+2)/(2d+1)>N/2\).
A positive-dimensional minimizing center contains a rational curve of
degree at most two by the covering-curve inequalities of
Chapter~\ref{ch:global}: for a curve center use \(N/\beta<2\),
and in dimension \(s\ge2\) use \(n(s+1)/(s\beta)<3\).
Such a curve contradicts \(\delta_C(X,H)>r\). Thus all the
divisorial minimizers in question have point centers, including
in the case \(\beta=r\). The uniqueness and equivariant
polystable-degeneration arguments in the proof of
Theorem~\ref{thm:descent-stability} now apply verbatim.
\end{proof}

For example every smooth, strictly positive real form in even
degree \(d\ge6\) satisfies the descent hypothesis in this dimension.
There are also examples over \(\mathbf Q\) in every degree in the
theorem. Lift a norm form of
\(\mathbf F_{\ell^d}/\mathbf F_\ell\) at a prime \(\ell\) to a
homogeneous integral form \(P\) on a \(d\)-tuple; on a triple
write \(P_3\) for its restriction to the first three coordinates.
The \(2d+3\) ambient coordinates split into two \(d\)-tuples
\(x,y\) and a triple \(z\). Every smooth integral equation
congruent to \(P(x)+\ell P(y)+\ell^2P_3(z)\) modulo \(\ell^3\)
has no \(\mathbf Q_\ell\)-point: a primitive zero would force
\(x,y,z\) in succession to vanish modulo \(\ell\).
The congruence class is a translate of a full lattice of
coefficients and is Zariski dense, so smooth examples exist.

\subsection{Descent and the scheme of lines in higher dimensions}
\label{sec:nov27-descent-line-scheme}
Let \(d\ge4,\ n\ge 2d-2\), \(N=n+1\), \(r=n+2-d\), and let
\(X=X_d^n\) be smooth. Assume either that a finite automorphism
group has no common fixed point on \(X\), or that a model of \(X\)
over a subfield of \(\mathbf C\) has no rational point. If
\[
                \delta_C(X,H)>r\quad\text{for every line }C\subset X,
\]
then \(X\) is uniformly K-stable. In particular, if the scheme of
lines of \(X\) is smooth of its expected dimension, this holds for
\[
                     2d-2\le n<\frac{d^2+d-4}{2}.
\]
For instance, under the stated descent and line-scheme hypotheses,
in degree ten this includes dimension fifty-two.

\begin{proof}
At \(n=2d-2\) use Theorem~\ref{boundary:main}. Otherwise put
\(\beta=\delta(X,H)\), and suppose \(\beta\le r\).
For any closed point, apply pointed complete hyperplane adjunction
to a smooth linear section of dimension \(2d-2\). Its \(H\)-threshold
exceeds \(d\) by the same theorem. Equivalently one may choose the
point on the center of a global divisorial minimizer. Thus
\(\beta>Nd/(2d-1)>N/2\).
Let a divisorial minimizer have positive-dimensional center of
dimension \(s\). When \(s=1\), the curve-image estimate of
Chapter~\ref{ch:global} gives a rational curve there of degree at
most \(N/\beta\). When \(s\ge2\), apply the canonical-extension
estimate of that chapter: its hypothesis \((s+1)\beta\ge N\) holds,
and it gives the same bound. Its use of the slope parameter
\(t<\beta/r\) applies also at \(\beta=r\). This degree is less
than two. The center therefore contains a line \(C\), contradicting
\(\beta\ge\delta_C(X,H)>r\). Consequently every such divisorial
minimizer has a point center, and the proof of
Theorem~\ref{thm:descent-stability} gives uniform K-stability.

For the explicit interval, write the kernel of the first normal
polars of a line as \(\bigoplus_i\mathcal O_{\mathbf P^1}(-a_i)\);
then \(a_i\ge0\) and
\(N_{C/X}=\bigoplus_i\mathcal O_{\mathbf P^1}(1-a_i)\).
Smoothness of the scheme of lines of its expected dimension implies
\(H^1(C,N_{C/X})=0\) by its tangent-space description and
Riemann--Roch. Hence every \(a_i\le2\). The ordinary line
calculation of Chapter~\ref{ch:models} has angular denominator
at most \(1+2/d\) and radial denominator
\((n-1/d)/(n-1)\le 1+2/d\). This gives the common bound
\(\delta_C(X,H)\ge Nd/(d+2)\). Substituting \(r=n+2-d\),
the inequality \(Nd/(d+2)>r\) is exactly
\(2n<d^2+d-4\).
\end{proof}

There are Zariski-dense examples over \(\mathbf Q\) throughout
this interval. Fix a prime \(\ell\), partition the \(n+2\) variables
into \(m=\lceil(n+2)/d\rceil\le d\) blocks of size at most \(d\),
and lift the degree-\(d\) finite-field norm form of
\(\mathbf F_{\ell^d}/\mathbf F_\ell\) integrally. Let \(P_j\) denote
its restriction to the coordinate subspace of block \(j\).
Impose the coefficient congruence
\[
             F\equiv\sum_{j=0}^{m-1}\ell^jP_j\pmod{\ell^m}.
\]
If the first \(j\) blocks of a primitive \(\ell\)-adic solution
were divisible by \(\ell\), their terms would be zero modulo
\(\ell^{j+1}\), since \(j+1\le m\le d\). Reducing modulo this
power shows the next block is also zero modulo \(\ell\), contradicting
primitivity after all blocks. The congruence specifies a full lattice
translate, hence is Zariski dense. It meets the open locus where
the hypersurface and its expected-dimensional scheme of lines are
smooth. For this last assertion the incidence of a line on a degree-\(d\)
hypersurface is a projective bundle over the Grassmannian, since
restriction to a fixed line is surjective. Its projection to the
coefficient space is proper. Generic smoothness gives the assertion
if it is dominant (with relative dimension \(2n-d-1\)), and otherwise
a general hypersurface has no lines and the line criterion is vacuous.

\paragraph{Descent for quartic ninefolds.}
The point-center and descent conclusions of
Section~\ref{sec:nov27-third-descent} now hold also for
\(d=4,n=9\). Every line or integral conic \(C\) on a smooth quartic
ninefold satisfies \(\delta_C(X,H)\ge220/31>7\). In fact a general
hyperplane containing the span of \(C\) cuts a smooth quartic
eightfold: off the base it is Bertini, and the set of tangent
hyperplanes at points of the base in the span has dimension at
most two, strictly less than the hyperplanes containing that span.
The normalized complete hyperplane flag lifts the eightfold curve
mark by \(10/9\). If \(\beta=\delta(X,H)\le7\), a general smooth
pointed hyperplane section and the eightfold uniform stability result
give \(\beta>(10/9)6=20/3\) at a divisorial minimizer. For a curve
center the covering-curve bound is \(10/\beta<2\); for a center of
dimension \(s\ge2\) it is \(9(s+1)/(s\beta)<3\). A line or
integral conic in the center contradicts its local mark \(220/31\).
Consequently every such divisor has a point center. The same descent
argument gives uniform K-stability for quartic ninefolds with a finite
symmetry group having no common fixed point, or with a model without
rational points over a subfield of \(\mathbf C\).

\section{Interpolating the two basis divisors at the original minimizing point}
\label{sec:nov28-direct-point}
For each pair in the following table, a smooth \(X_d^n\) is uniformly
K-stable if each line and each integral conic satisfies
\(\delta_C(X,H)\ge M\).
\[
\begin{array}{c|rrrrrrr}
 (d,n)&(4,8)&(4,9)&(5,11)&(6,13)&(7,15)&(8,17)&(9,19)\\ \hline
 M&32/5&299/40&43/5&191/20&2103/200&23/2&1249/100\\
 10^4 t&62583&71030&70469&78872&87348&94752&102193\\
 10^4\lambda&1104&2898&12720&14078&15434&17870&20320
\end{array}
\]
\begin{proof}
Let \(N=n+1,\ u=n/N,\ r=n+2-d,\ b=r/N\) and suppose
\(\beta=\delta(X,H)\le r\). The all-member result in dimension
\(2d-1\) and iterated pointed hyperplane adjunction at a point in
the center of a divisor computing \(\beta\) give
\(\beta>N(d+1)/(2d)>N/2\). The curve-center and covering-curve
bounds then give an integral curve of degree at most two in any
positive-dimensional optimizing center: for curves the degree is
at most \(N/\beta<2\), and in dimensions \(s\ge2\) it is at most
\(n(s+1)/(s\beta)<3\). The curve hypothesis therefore puts the
original optimizer \(v\) at a closed point \(p\).

Write \(A_X(v)=n+q\) with \(v(\mathfrak m_p)=1\), and set
\(S=S_H(v)\). We first explain how a lower bound \(S\ge S_0\) will
be used. Put \(f=1/d\) and take
\[
 h=1+n(1-(1-f)^{1/n}),\quad g=(1-b)/(d-1),\quad
 \alpha=\frac{1-Nf+nf^{1+1/n}}{N(1-f)}.
\]
Use a basis compatible with \(v\) and a hyperplane avoiding \(p\),
and remove that hyperplane from the full basis divisor \(D_m\)
and its upper \(1-f\) basis divisor \(D_m^+\). Denote their
residuals by \(R_m,R_m^+\). The usual projection and hyperplane
estimates give their respective asymptotic \(H\)-class coefficients
\(u\) and at most \(1-\alpha\), while their \(v\)-orders are
\(S+o(1)\) and at least \((1+g)S-o(1)\). For all test valuations
\(w\), their orders are at most \((1+o(1))S_H(w)\) and
\((h+o(1))S_H(w)\), respectively, by the uniform partial-basis
estimate in Chapter~\ref{ch:global}.
For the indicated rational \(t,\lambda>0\), the following strict
comparisons suffice:
\[
    t+h\lambda<M,\qquad 1/3<\eta:=r-ut-(1-\alpha)\lambda<1,
\quad a:=t+(1+g)\lambda-r>0,\qquad aS_0+\eta>1.       \tag{1}
\]

Fix a small positive rational \(e\), take a sufficiently high
basis degree, and interpolate linearly from \((\beta-e)R_m\) to
\(\Omega_m=tR_m+\lambda R_m^+\), perturbing the initial
coefficient rationally if necessary. The initial boundary is klt
by uniform comparison. Suppose the endpoint has a positive-dimensional
non-klt center. On a log resolution of the two fixed divisors,
take the first parameter at which there is a positive-dimensional
non-klt locus. By linearity of discrepancies its positive centers
are generically lc and its non-lc scheme has finite support. The
value at \(v\) throughout the segment is at least
\(A_X(v)-eS-o(1)\) since it has this value at the initial endpoint
and a larger value at the final endpoint. A general hyperplane
through \(p\) with coefficient tending to zero with \(e\) therefore
makes the interpolation boundary non-lc at \(p\), preserving its
positive centers at their generic points. The ample remainder is
larger than \(H/3\), by convexity, (1), and
\(r-u\beta\ge r/N>1/3\). Quasi-log adjunction and the length-one
theorem, applied to a minimal positive lc center and the finite
non-lc scheme as in Section~\ref{boundary:numbers}, put a line or
conic in that center. At its lc place the ratio \(A_X/S_H\) is at
most \(\max(\beta,t+h\lambda)+o(1)<M\). This contradicts the
curve hypothesis. Thus the endpoint is klt off a finite set.

Now add a general hyperplane through \(p\) with coefficient tending
to \(\eta\) from below, keeping the remainder ample. Bertini
preserves the finite cosupport of the multiplier ideal. Nadel
vanishing on \(X\) makes its quotient have length at most one.
Applying the multiplier-ideal inequality at \(v\) and passing to
the limits gives \(aS+\eta\le1\), contrary to (1).

For the two quartic columns one may take \(S_0=7/4,49/31\),
respectively. At quadratic contact and at reduced cubic contact
the signed-surface flag gives
\(NS\le(13/6)(q+2)+n-2\) in both dimensions. Combining with
\(n+q=\beta S\le rS\) gives the claimed values. The generalized
Eckardt point has ratio \(nN/(n+3)>r\) by the all-member result on
its smooth base. At a nonreduced cubic tangent form, the weighted
linear-factor extraction of the quartic eightfold also applies
when \(n=9\). In its notation, the ordinary weighted discrepancy
and expected order become \(ns+1,((n+3)s+1)/N\), leaving the
positive discrepancy difference \((6s+3)/N\) when multiplied by
\(r=n-2\). The angular means on \(\mathbf P^{n-2}\) are
\(1/N,3/(4N),(n-1)/N\) for the linear coordinate, fixed quartic,
and residual free degree. At the axis the radial discrepancy is
\(n-1\), while \(r(n+2)/N<n-1\). On \(\mathbf P^{n-2}\) put
\(c=3r/(4N)\). A general hyperplane section of the quartic
divisor \(G\) is smooth of dimension \(n-4\). For \(n=8,9\) this gives a quartic
fourfold or fivefold, respectively. The fourfold case follows from \cite[Theorem~4.21]{ar2003.13788}
and the fivefold case from Chapter~\ref{ch:fivefold}. The log interpolation bound for its ambient
\(\mathbf P^{n-3}\) applies since
\(c<1-((n-2)/4-1)/(n-3)\). The projective cone with coefficient
\(4c/(n-2)\) at infinity then has log anticanonical coefficient
\((n-1)(1-4c/(n-2))=r(n-1)/N\).
Dropping the hyperplane after openness gives the required angular
comparison at coefficient \(r\). Thus the extraction also excludes
this point when \(n=9\).

For the remaining columns put \(S_0=(n+Q)/r\), using
\[
\begin{array}{c|rrrrr}
 d&5&6&7&8&9\\ \hline
 10^3Q&1551&2167&2794&3285&3667\\
 10^4B&27671&30220&32269&34424&36664
\end{array}
\]
Apply the point estimates to \(A_X(v)\le r S\) at both \(q=0,Q\).
Use the signed surface estimate for \(\tau_p(H)\le B\) and the
discriminant estimate for \(B<\tau_p(H)<d\), with cutoff \(j=25\)
for a prime of class \(jH\) and the sectional Milnor improvement
for \(j\le2\), as in Section~\ref{sec:nov-second-row}. Only the
following finite \((k,K)\) trials require the principal-bin bound,
all for \(j=1\). Here \(k\) is the multiplicity and \(K\) the
bin count.
\[
\begin{array}{c|l}
5&(4,16)\\6&(5,26)\\7&(5,28),(6,34)\\
8&(5,30),(6,32),(7,38)\\9&(5,32),(6,30),(7,36),(8,44)
\end{array}
\]
The extremal tangent contact uses the generalized Eckardt formula
as before. Rational certificates for the entire table, including
(1) with outward root bounds, are in \path{nov28-bin-direct.py}.
\end{proof}

\paragraph{A line on the cubic tangent form.}
Consider a smooth quartic eightfold with cubic tangent contact at a
point \(p\). Write \(f,g\) for the cubic and quartic forms in the
tangent hyperplane. Suppose \(F=(f=0)\subset\mathbf P^7\) contains
a line \(\ell\) with a normal quadratic polar having two distinct
roots on \(\ell\), neither a zero of \(g|_\ell\). Then the following
signed upper bound holds for \(v(\mathfrak m_p)=1,\ A_X(v)=8+q\):
\[
 9S_H(v)\le \frac{31}{3}
   +\left(\frac{11}{6}+\frac{(e-2)_+}{12}\right)q,
 \qquad e=\max_{u\in\ell}\operatorname{mult}_u(g|_\ell)\le4.
\]
In particular \(A_X(v)\le6S_H(v)\) implies \(S_H(v)\ge17/9\),
improving to \(S_H(v)\ge2\) when \(e\le3\).

For the proof use the ordinary point presentation with symbols
\(T,Z,W\) of section degrees \(1,1,3\) and ordinary orders \(0,3,4\)
and the eight linear tangent symbols of order one. Its two relations
are \(T^2Z+f,TW+g\), as in Section~\ref{sec:q8-tangent-tests}.
First degenerate the tangent forms to their restriction to the plane
spanned by \(\ell\) and the chosen normal direction. Taking the
least power of the coordinate \(z\) transverse to \(\ell\) next
gives \(f_0=zxy,g_0=g(x,y)\), with the roots of the polar as
coordinate points. The two initials form a regular sequence by the
root assumption. Thus the graded angular dimensions are preserved
under both degenerations by the complete-intersection presentation.
Specialization of bases for a test valuation and each weight
filtration gives the local threshold comparison used in the angular
tests of that section.

Split the monomials in \(z\) and the five complementary coordinates
and retain the binary block. In section degree \(m\) and ordinary
order \(j\) put \(l=j-m\). If the absorbed exponent of \(z\) is
\(s\) and the remaining binary degree is \(D\), the block is
\[
 \sum_{i\le s,\ 2i+h\ge l}
 (xy)^i g^h\,\mathbf C[x,y]_{D-2i-4h}\qquad (i,h\ge0).
\]
It is divisible by \((xy)^I g^J\) with
\(I=\max(0,\lceil(4l-D)/6\rceil)\) and \(J=(l-2s)_+\).
After removing these powers the difference from the complete binary
dimension is bounded by an absolute constant. Indeed the binary
ring is finite free over \(\mathbf C[xy,g]\) with a homogeneous
remainder basis of degrees \(b\le4\), by the regular sequence of
degrees two and four. A basis element \((xy)^k g^t R_b\) in degree
\(D\) belongs to the block precisely when \(2\min(k,s)+t\ge l\).
This is \(t\ge J\) and \(k\ge(4l-D+b)/6\), with nonnegativity
understood. On the degree line \(2k+4t=D-b\) only a bounded
number are lost compared with the stated two fixed powers.
The mean binary order at any marked point is therefore at most
\(D/2+(e-2)_+J+O(1)\). Here \(xy\) is squarefree; removing
its fixed copies reduces the free half-degree by exactly their
multiplicity at either root.

We record all means in units scaled by nine and divided by \(m\).
In the normal basis \(T^a Z^i W^h R\) with \(R\) complementary
to \((zxy,g)\), one has \(l=2i+h-a\) and \(s\ge i\), so
\(J\le(h-a)_+\). The ordinary degree simplex gives
\(9\overline{J/m}\le\mathbb E[(E_1/3-E_0)_+]=1/12\), for
independent rate-one exponential variables \(E_0,E_1\). The
residual ternary curve ring \(\mathbf C[x,y,z]/(zxy,g)\)
has degree twelve. Its \(z\)-axis has multiplicity eight, since
inverting \(z\) gives the binary intersection \((xy,g)\) of
length eight. Its remaining degree four lies in \(z=0\).
Thus its single free simplex index contributes mean \(2/3\) to
\(z\) and \(1/3\) to the binary degree. This also follows by
logarithmic differentiation in the two character weights in the
complete-intersection Hilbert series. The additional \(f_0^i\)
contributes one to \(z\) and two to the binary degree, while
\(g^h\) contributes \(4/3\) to the latter. The total scaled means
are consequently \(5/3\) for \(z\), \(11/3\) for \(D\), and
one for each other free coordinate. Along the binary vertex the
scaled radial ratio is \((11/3)/2\); on its exceptional binary
line the marked order bound just proved applies. Ordinary
adjunction therefore gives the claimed angular slope. Finally
\(9S_H(\operatorname{ord}_p)=31/3\) from the two-symbol ring.

When every possible optimizing point of cubic contact either
satisfies this line criterion or is already excluded by the other
tangent tests, hypothesis \(\delta_C(X,H)\ge198/31\) on all lines
and integral conics suffices. To check the last step, the quadratic
contact case gives \(S_H(v)\ge2\) at an optimizer of ratio at
most six; thus the uniform lower bound here is \(17/9\).
Use the untilted isolated-point Nadel comparison with \(t=127/20\)
and \(\Delta=t-6=7/20\). The curve budget is strict, the ample
slack is \(6-8t/9=16/45>1/3\), and
\(\Delta(17/9)+16/45=61/60>1\). Positive-dimensional original
optimizers are excluded by the original line and conic estimates.

\paragraph{Quartic ninefolds and separate hypotheses on lines.}
For a smooth quartic ninefold, the common short-curve mark
\(M=1491/200\) suffices. Alternatively it suffices that every line
on the ninefold have \(\delta_C(X,H)\ge954/125\).
For a smooth quartic eightfold, the analogous sufficient
line-only mark is \(164/25\).
If, on the ninefold, every irreducible cubic tangent form at points
of cubic contact is smooth, the sufficient common mark is any
\(M>125/17\), and the sufficient line-only mark is \(751/100\).

Indeed, increasing \(n=8\) to \(9\) adds a single independent
linear variable to the angular presentations. Thus the scaled
radial intercept in (2) is now \(34/3\), the angular slope is still
\(17/8\), and the ordinary normalization denominator is ten.
Every hypothetical optimizing valuation at reduced cubic contact
then has \(S\ge187/117\). The argument producing a transverse free
line and the two binary forms applies verbatim, with one additional
ambient coordinate. At quadratic contact the signed surface flag
gives \(10S\le2(9+q)-7\), whence \(S\ge7/4\). Nonreduced
cubics and generalized Eckardt points are excluded by the 
dimension-nine comparisons of Section~\ref{sec:nov28-direct-point}.

For the statement with smooth irreducible cubic tangent form, choose
a general free covering line on the cubic \(F=(f=0)\) itself. Lines
cover this smooth sevenfold by the same direction count, and the
quartic \(g|_F\) is reduced because its singular points give singular
points of the tangent hyperplane section and are finite. Smooth
evaluation from the free lines gives four transverse intersection
points. Choose a normal quadratic polar with two simple roots
avoiding all four, which is possible by basepoint freeness and
Bertini in characteristic zero. The angular calculation using a line
on the cubic then gives \(10S\le34/3+(11/6)q\); hence
\(S\ge31/17\). Reduced reducible cubics have a linear factor;
the same weighted linear-factor extraction excludes a point of
this type in dimension nine. In detail, take any linear factor of
multiplicity \(s\ge1\), and the other factor is nonzero on its
hyperplane. The initial relation and angular means in the proof
of Section~\ref{sec:nov28-direct-point} are unchanged for \(s=1\),
with the same positive radial discrepancy gap. Thus one may use
the smaller \(S_0=7/4\) uniformly under the smoothness hypothesis.
The untilted bound has \(\Delta>6/17\), chosen arbitrarily close
to \(6/17\). Then \(\eta=7/10-(9/10)\Delta>1/3\) and
\(\eta+(7/4)\Delta>1\), giving the stated common mark.

For the other assertions use the interpolation tested by (1).
For a line-only criterion its second inequality becomes \(\eta>1/2\).
The first positive-dimensional non-klt center in that argument
then contains a line. The initial optimizing center is itself a
point: \(\beta>N/2\) by the pointed flag from dimension seven,
so the canonical-extension bound in Chapter~\ref{ch:global}
puts a covering curve of degree \(\le N/\beta<2\) in every
positive-dimensional minimizing center. For a curve center the
same conclusion follows from the degree bound there. At the initial
interpolation boundary the ample coefficient is at least \(r/N>1/2\).
Here are all required strict comparisons. In the last column use
the smaller value of \(S_0\) specified for the corresponding case,
with \(\eta\) and \(a\) as in (1).
\[
\begin{array}{c|c|c|r|r|c}
n&M&\text{required }\eta&10^4t&10^4\lambda&S_0\\ \hline
9&1491/200&>1/3&71903&2062&187/117\\
9&954/125&>1/2&53964&17422&187/117\\
9&751/100&>1/2&59393&12240&7/4\\
8&164/25&>1/2&44757&16250&16/9
\end{array}
\]
In direct substitution it suffices to use
\(h<1.28258,\ \alpha>0.063969\) for \(n=8\), and
\(h<1.28315,\ \alpha>0.05716\) for \(n=9\).

In fact, for \(8\le n\le11\) a smooth cubic tangent form \(F\)
excludes an optimizing point of ratio \(\le r=n-2\) at cubic
contact. Repeat adjunction along \(F\subset\mathbf P^{n-1}\)
in Section~\ref{sec:q8-tangent-tests}. Before division by
\(N=n+1\), the mean fixed coefficient of \(F\) is one, the mean
of the restricted quartic is at most \(1/12\), and the moving
degree on \(F\) is \(n-1\). The two-symbol calculation is
unchanged upon adding independent linear coordinates, with unit
scaled mean for each. The isolated-singularity argument in that
section gives \(w(g|_F)\le A_F(w)\): the ordinary vertex blowup
of the cubic cone has discrepancy \(n-4\ge4\), and a smooth
general hyperplane section meets the quartic transversely on the
resolution. The smooth cubic \(F\) has \(\delta(F,H)>n-3\) by
the cubic theorem. Therefore the scaled angular comparison on
\(F\) holds since
\[
       \frac{r}{N}\left(\frac1{12}+\frac{n-1}{n-3}\right)\le1,
       \qquad (n-6)(n-11)\le0.
\]
Off \(F\), the scaled bound \(4/3\) is at most \(N/r\) in the
same dimensions, and the \(F\)-coefficient itself is \(r/N<1\).
Finally the ordinary radial discrepancy exceeds the product of
the index and the mean radial order, as
\(nN-r(n+7/3)=(2n+14)/3>0\).
This supplies a strict discrepancy margin for every valuation
centered at that point.

\paragraph{A normality alternative for cubic contact on the quartic ninefold.}
For an irreducible cubic tangent form \(F=(f=0)\subset\mathbf P^8\),
the scaled angular slope is at most two in each of the following
cases:
\begin{enumerate}
\item \(\dim\operatorname{Sing}F\le4\), including the empty locus.
\item The cubic has at most four essential variables, meaning it is
the pullback of a cubic polynomial under a linear map of rank at
most four.
\item \(F\) is nonnormal.
\end{enumerate}
Consequently if every irreducible cubic contact satisfies one of
these alternatives, then a common line and conic mark strictly
larger than \(170/23\) suffices for uniform K-stability.

For (1), a general \(\mathbf P^3\) cuts a smooth cubic surface on
\(F\), and its restriction of \(g\) is nonzero on that surface
since \(f,g\) are relatively prime. Its quartic section has degree
twelve. At least one of the twenty-seven lines is therefore not
contained in \((g=0)\). All points on the cubic surface are smooth
points of \(F\), so the ambient normal quadratic polars of that
line are basepoint-free. Choose a squarefree polar avoiding the
finite set of zeros of \(g\) on the line. The point calculation
using a line on the cubic gives the assertion because the maximum
binary quartic multiplicity is at most four.

For (2), let \(m\le4\) be the number of essential variables and
let \(R\) be their common vector kernel. The restriction of \(g\)
to \(R\) is nonzero. Otherwise all radical derivatives of \(g\)
would vanish on \(\mathbf P(R)\), where \(df=0\) as well. The
remaining \(m\) cubic derivatives have a common zero there since
\(\dim\mathbf P(R)=8-m\ge m\). This would give a whole singular
\(T\)-line of the tangent hyperplane section. Take a general
two-dimensional vector subspace of \(R\) with nonzero binary
restriction \(G(y)\) of \(g\), and a general two-dimensional
complementary section of the essential space with squarefree
binary restriction \(F(x)\) of \(f\). First specialize to this
four-space and then give the essential variables positive weight.
The resulting regular sequence is \(F(x),G(y)\) in disjoint
binary blocks.

Here is the angular computation for these disjoint forms. In a
fixed multidegree with binary degrees \(D_1,D_2\), the series is
the sum of the spaces
\[
 F^i\mathbf C[x]_{D_1-3i}\otimes
 G^h\mathbf C[y]_{D_2-4h},\qquad 2i+h\ge j-m',
\]
where \(m'\) is the section degree. In the ordinary simplex the
two residual binary rings modulo their respective forms contribute
one free degree each. Thus the scaled means are \(4\) for \(D_1\)
and \(7/3\) for \(D_2\). Their two radial quotients, after blowing
up the binary vertices, are at most two. On the angular
\(\mathbf P^1\times\mathbf P^1\), first fix a fiber over the
first factor. Its coefficient is at most the mean of \(D_1/2\).
Indeed adapt the second binary space to its \(G\)-order; every
coefficient space on the first factor is then a fixed power of
the squarefree \(F\) times the complete residual binary space.
Conversely the first binary space has a basis simultaneously
compatible with \(F\)-order and the marked first-factor order;
all latter jumps are distinct. For each such basis vector of
\(F\)-order \(i\), the residual coefficient space on the marked
fiber is exactly \(G^J\mathbf C[y]_{D_2-4J}\) with
\(J=(j-m'-2i)_+\), if nonzero. Its mean at every marked point
is at most \(D_2/2+2J\). The mean of \(J\) over all the blocks
is at most the mean of their actual \(G\)-orders, which is \(1/3\)
in scaled units. For this last identity use normal monomials
\(T^a Z^i W^h R_1R_2\), with \(R_1,R_2\) binary remainders
modulo \(F,G\). The expansion is free over the subalgebras
generated by the individual forms; it is adapted to \(G\)-order,
whose exponent is \(h\) with scaled mean \(1/3\). Thus the tagged
point on the angular fiber has scaled mean at most
\(7/6+2/3=11/6\). Adjunction along the fiber and the independent
coordinate factors gives the claimed slope two.

For (3), a component \(S\) of the singular locus has dimension
six. The secant variety of \(S\) is in \(F\), because a cubic
vanishing doubly at two points vanishes on their joining line.
It follows that \(S\) is a linear six-plane. In detail, if its
span is a seven-plane and it is non-linear, its secants fill the
span, forcing \(F\) to contain a hyperplane. If its span is all
of \(\mathbf P^8\) its secants have dimension at least eight:
otherwise the joins of \(S\) with a general point of \(S\) all
have dimension seven and coincide with the secant variety.
That variety is a cone from every general point of \(S\), so its
linear vertex contains the span. These are contradictions.
Write \(S=(x=y=0)\). Double vanishing now gives
\(f=f_3(x,y)+\sum z_i Q_i(x,y)\). If the normal quadratics
span at most two dimensions then there are at most four essential
variables, already handled. In the remaining case the quadratics
span \(x^2,xy,y^2\). Shears of the three corresponding linear
coordinates absorb \(f_3\). The remaining cubic is consequently
\(z_0x^2+z_1xy+z_2y^2\), independent of the other coordinates.
There is a covering family of lines on its smooth locus: take
\(z_i\) linear in \(x,y\) subject to the identity
\(z_0x^2+z_1xy+z_2y^2=0\), and all remaining coordinates
arbitrary linear in \(x,y\). The two linear syzygies
\((-y,x,0),(0,-y,x)\) span the kernel at every point of
\(\mathbf P^1\). These lines therefore cover the complement of
\(S\) and are disjoint from \(S\), the entire singular locus.
Some such line is not in \((g=0)\), and the forward-polar
argument of (1) applies.

Finally \(10S\le34/3+2q\) at these cubic points gives
\(S\ge5/3\) at a putative optimizer of ratio at most seven.
The quadratic contact bound \(7/4\) is stronger, and the
linear-factor extraction and Eckardt test handle the remaining
tangent multiplicities. Choose \(t\) slightly larger than
\(170/23\) below the given curve mark, and take \(\lambda=0\) in
the interpolation. At equality \(t-7=9/23\), the ample remainder
is \(8/23>1/3\), and \((t-7)(5/3)+8/23=1\). Increasing \(t\)
gives strictness, proving the assertion.

\paragraph{All cubic tangent forms in dimension nine.}
For every smooth quartic ninefold, the bound on lines and conics
\(\delta_C(X,H)\ge M\) for any \(M>125/17\) implies uniform
K-stability. The only geometric hypothesis here concerns the curves
of the ninefold. We prove that in all cases an original minimizing
point of ratio \(\beta\le7\) must satisfy \(S_H(v)\ge7/4\).

First refine the forward-polar argument using a line in \(F\) with
a squarefree binary quadratic polar \(Q\) coprime to \(g|_\ell\).
It gives scaled slope \(11/6\) for every multiplicity pattern of
the latter quartic. Indeed at a root of the binary quartic, \(Q\)
is nonzero, so the mean marked order of the binary block in that
calculation is at most \(D/2-I+2J+O(1)\). At its other points
the mean is already at most \(D/2+O(1)\). The scaled mean of
\(J\) is at most \(1/12\). The scaled mean of \(I\) is at least
\(2/5\): in the \(z\)-axis part, of mass \(2/3\), the residual
binary degree is bounded. There \(l=2i+h-a\), \(D=2i+4h+O(1)\),
and \(I\ge(i-2a/3)_+-O(1)\). The two unit simplex indices have
\(\mathbb E[(E_1-2E_0/3)_+]=3/5\) for independent rate-one
exponentials, proving the claim. Hence every marked point and
the binary vertex have scaled slope at most \(11/6\).

Next, in the case of at most four essential variables, retain a
general three-space of the essential variables and a binary
space in the radical. Relative to the preceding disjoint-block
argument the first form is now an integral plane cubic, by Bertini
and irreducibility. Its threshold is at least \(5/6\) because its
only possible singularity is a node or an ordinary cusp; the latter
has weighted order six and discrepancy five. The second form is
a nonzero binary quartic as before. Normalize means in the same
units. The ternary cubic block now has total degree five coming
from its two residual simplex slots and the three degrees of
\(f^i\). The binary quartic block still has degree \(7/3\).
Refine first by a marked fiber of the second factor. To bound its
coefficient adapt the ternary factor to its \(f\)-order; the
resulting binary spaces are fixed powers of the quartic times
complete binary spaces. As before the coefficient is at most
\(7/6+2/3=11/6\). With the binary space adapted both to the fiber
and to its quartic order, restriction to a given marked order
instead leaves \(f^I\) times the complete ternary space of degree
\(D_1-3I\). Its mean at any plane valuation \(w\) is at most
\(\{D_1/3+I/5+o(m')\}A(w)\): use the preceding threshold
and \(\delta(\mathbf P^2,\mathcal O(1))=3\). The inequality is
uniform by the complete-series comparison, treating bounded
residual degrees with vanishing normalized contribution. The
mean of \(I\) is at most one by the normal monomial basis, adapted
to \(f\)-order just as for the two binary blocks. This gives a
scaled slope at most \(5/3+1/5=28/15\). The two radial ratios
\(5/3,7/6\) are smaller. The same bounds consequently apply after
ordinary adjunction to the binary and ternary vertices.

We show that these alternatives cover every irreducible \(F\).
The preceding cubic-surface and nonnormal arguments leave only
normal \(F\subset\mathbf P^8\) with a five-dimensional component
\(Y\) of the singular locus. We may take \(G=(g=0)\) smooth by
the linear change of \(T\) already proved. No such component is
contained in \(G\). To check this, a general three-plane section
of \(F\) is an integral cubic surface with isolated singularities.
It has at most four singular points: three collinear ones would
force every quadratic partial derivative to vanish along their
line. Four coplanar ones with no three collinear force the plane
cubic restriction to vanish identically (use its factorization
in the reducible cases and the plane genus formula otherwise).
This contradicts integrality. Four singular points can therefore
be taken as the coordinate vertices. The cubic is then a linear
combination of the four squarefree cubic monomials. A fifth
singular point must have four nonzero coordinates, where division
by their product makes the equation \(\sum a_i/x_i=0\). Its
differential does not vanish. Slicing \(Y\) transversely gives
\(\deg Y\le4\). If \(Y\subset G\), weak Lefschetz on the smooth
sevenfold \(G\) gives \([Y]\in H^4(G,\mathbf Z)=\mathbf Z H^2\),
so its degree is divisible by four and hence equals four.
If the span of \(Y\) is a six-plane, every quadratic partial
vanishing on this quartic hypersurface vanishes on that six-plane,
contradicting normality. A seven-plane span is excluded as in the
secant argument above: the secants fill that span. In the full
eight-plane span, \(Y\) is a five-dimensional variety of minimal
degree. The classical classification of varieties of minimal degree
(used also in Chapter~\ref{ch:cubics}) gives a rational normal
scroll, possibly conical, with a four-plane ruling, or a cone over
the Veronese surface with nonempty vertex. A smooth quartic in
\(\mathbf P^8\) contains no four-plane, since its four normal cubics
would have a common zero there. Nor can it contain this nondegenerate
cone: the Zariski tangent space at any vertex point is the whole
ambient space. The remaining span of dimension five is a linear
space of degree one. This proves the assertion.

Let \(m\ge5\) be the number of essential variables. At a general
\(s\in Y\) the Hessian rank is at least two. Here are the details
for this rank assertion. Contraction of the cubic polarization
with a vector gives a linear Hessian map with kernel exactly the
inessential vector space. On the essential \(m\)-space it is
injective. If the generic Hessian on \(Y\) has rank zero, then
\(Y\) belongs to the radical, forcing \(m\le3\). If its rank is
one, its projective Hessian image is contained in the quadratic
Veronese image of a variety \(B\) of linear forms, where
\(\dim B\ge5-(9-m)=m-4\) by linear projection. If \(m\ge7\),
the quadratic coordinate functions on any variety of dimension
at least three include ten independent monomials, by linear
Noether normalization. This exceeds \(m\le9\), the dimension of
the linear image of the Hessian map. If \(m=6\), a variety of
dimension at least two has at least six such monomials, with
equality only if it is a plane. Indeed a general projection to
a plane supplies the six, and an additional independent linear
coordinate supplies more: equivalently the standard minimal
Hilbert-function bound for a nondegenerate \(k\)-fold in \(\mathbf
P^a\) is \(H(2)\ge\binom{k+2}{2}+(a-k)(k+1)\), obtained by
successive general hyperplane sections down to spanning points.
Thus the Hessian map would have its whole image in quadratic
forms on three variables, forcing the original cubic to use at
most three by polarization. For \(m=5\), the same Hilbert bound
shows that \(B\) spans at most a plane. In a plane, a nonlinear
curve has at least five independent quadratic coordinates. At
equality the entire image would again be supported on three
variables, impossible. Therefore \(B\) spans a line, whose
squares span all of \(\operatorname{Sym}^2\langle x,y\rangle\).
Take the three preimages \(a_0,a_1,a_2\) of \(x^2,xy,y^2\).
Symmetry of the cubic tensor, comparing contractions on the
pairs \((a_0,a_2),(a_0,a_1),(a_2,a_1)\), gives
\(x(a_i)=y(a_i)=0\) for all three. For example the first comparison
gives \(x(a_2)x=y(a_0)y\) and the other two give the remaining
four vanishings. Thus these preimages span the common kernel of
\(x,y\) in the essential five-space. The corresponding derivative
quadratics are precisely \(x^2,xy,y^2\), so the cubic has the
form \(z_0x^2+z_1xy+z_2y^2+f_3(x,y)\). Its singular locus contains
the codimension-two linear subspace \(x=y=0\), contradicting
normality. This proves the rank assertion.

Choose \(s\) general with this rank and \(g(s)\ne0\). The Hessian
annihilates the affine tangent to the cone over \(Y\), so its rank
is two or three. Write \(f=Tq_s(y)+c_s(y)\) from this double point.
Line directions on \(F\) through \(s\) are \((q_s=c_s=0)\) in
\(\mathbf P^7\). There is a five-dimensional component of directions
outside the quadric radical: for rank three the radical has dimension
four; for rank two write \(q_s=xy\). If both cubic restrictions to
the two hyperplanes have support only on their common hyperplane,
then \(c_s=Ax^3+By^3+xyL\), so after a linear shear \(m\le3\).
A restriction identically zero gives a linear factor of \(f\).
These alternatives are excluded. Since every five-dimensional
singular component has nonzero restriction of \(g\), the projection
from \(s\) of \(\operatorname{Sing}F\cap G\) has dimension at most
four. A general line of the indicated component therefore avoids
that intersection. Along this line the normal quadratic polars
all vanish at \(s\), and some have a simple zero there since the
tangent direction is outside the radical of \(q_s\). A general
such quadratic is consequently squarefree, and it avoids all zeros
of \(g\) on the line because at each such zero the line is in the
smooth locus of \(F\). We have obtained the required forward-polar
line.

The worst scaled slope of the cubic cases is therefore at most
\(28/15\). With radial intercept \(34/3\), the inequality
\(9+q\le7S\) then gives \(S\ge41/23>7/4\) at irreducible cubic
contact. The other contacts have the bound or exclusion already
proved; quadratic contact has \(S\ge7/4\). The untilted comparison
at \(M>125/17\) completes the proof.

At the exact common mark \(M=125/17\) one obtains K-semistability.
Indeed if \(\beta<7\), then at the optimizing point the additional
coefficient in the direct interpolation is \(t-\beta\), rather
than the weaker lower bound \(t-7\). At \(t=M\) the limiting
Nadel expression \((t-\beta)S+7-(9/10)t\) is therefore strictly
larger than one; the ample remainder is \(13/34>1/3\).
Both strict inequalities persist with \(t<M\) arbitrarily close
to \(M\), where the curve comparison is strict as well.

\paragraph{The possible quadratic contact at the exact mark.}
At the common mark \(125/17\) any divisor computing seven in
hyperplane units must have quadratic tangent contact. It must
satisfy \(S_H(v)=7/4\), \(A_X(v)=49/4\) after normalization
\(v(\mathfrak m_p)=1\): the signed surface inequality gives the
lower bound on \(S\), and otherwise the untilted Nadel comparison
is strict for a curve budget just below \(125/17\). At cubic
contact the preceding lower bound is strictly stronger. The
following quadratic contacts are also incompatible with equality:
\begin{enumerate}
\item the quadratic tangent form has rank one or at least six;
\item its three tangent forms constitute a regular sequence and
admit the isotropic plane and line of Section~\ref{sec:q8-tangent-tests};
\item the quadratic divides the cubic and the two-symbol line test
of that section applies.
\end{enumerate}
For (2), the scaled intercept increases by one from the eightfold
to \(65/6\), while the angular slope stays at most \(3/2\).
For (3) the new intercept is eleven with the same scaled slope.
Either upper bound for \(10S\) is strictly smaller than \(35/2\)
at \(q=A_X(v)-9=13/4\). At rank at least six the regular and
divisible cases and the required flags follow by the same
dimension count on the quadric as in that section: a general
isotropic line avoids each codimension-two base in the smooth
locus, and the nondegenerate quotient has a three-dimensional
isotropic vector space. If instead \(q_2=L^2\) and \(L\nmid q_3\),
the cubic and quartic are coprime modulo \(L\): otherwise
the linear space \(L=0\) of the tangent hyperplane would cut a
reducible quartic \(T\bar q_3+\bar q_4\). That restriction to
\(\mathbf P^8\) is integral, since the two normal derivatives of
the original smooth hypersurface and the two factors would have
a common zero. A general line on \(L=0\) therefore avoids their
common zero set. Degenerate the tangent triple to \(L^2,f,g\)
with those coprime binary restrictions of degrees three and four.
The free \(s,z\) sectors of the computation for \(sz,f,g\) in
Section~\ref{sec:q8-tangent-tests} are now replaced by a single
quadratic symbol of scaled \(L\)-mean two and bounded parity
remainder. The binary means and their bound \(3/2\) are unchanged,
since \(I=(i-4a)_++O(1)\) holds with the same symbol degrees.
This gives \(10S\le65/6+2q<35/2\) at the proposed equality.
Finally, if \(q_2=L^2\) and \(L\mid q_3\), give \(L\) weight two
and the other eight tangent parameters weight one. The common
linear-factor calculation has \(A(E)=10\), \(S_H(E)=13/10\),
with the fixed quartic, moving series, and opposite vertex
controlled as in the dimension-nine extraction of
Section~\ref{sec:nov28-direct-point}. Its positive discrepancy
margin already excludes \(A_X(v)\le7S_H(v)\).

\paragraph{A mixed quadratic contact in the equality case.}
In the notation of the preceding endpoint test, suppose
\(q_2=xy\) and \(q_3\) vanishes on \(R=(x=y=0)\). Write
\[
 q_3\equiv xA+yB\pmod{(x,y)^2},\qquad G=q_4|_R.
\]
Assume there is a line in \(\mathbf P(R)\) on which \(A,B\) are
nonzero proportional binary quadratics, with common form \(U\)
coprime to the nonzero restriction of \(G\). Then at this point
\[
                   10S_H(v)\le65/6+2q.
\]
It follows that every optimizer of ratio at most seven has
\(S_H(v)\ge43/24>7/4\).

In particular the line exists when \(A,B\) are nonzero and
proportional, or when they are coprime and their pencil contains
a member of rank at least four. To justify this, the restriction of
\(X\) to the seven-plane spanned by \(p\) and \(\mathbf P(R)\) is
the cone \(G=0\). It is integral because two factors there and
the three normal derivatives would have a common zero on this
seven-plane. Its singular scheme has support of dimension at most
two: a three-dimensional component meets all three normal cubics,
contradicting smoothness. Thus the quartic in \(\mathbf P(R)=
\mathbf P^6\) has singular locus of dimension at most one. It is
locally factorial by parafactoriality for a complete intersection
regular in codimension three, and its coordinate ring is factorial
by Grothendieck--Lefschetz and projective normality. A common
irreducible factor of \(A,B\) modulo \(G\) therefore lifts to a
common homogeneous factor, because its degree is at most two and
all graded pieces in degrees below four are ambient polynomials.
Consequently in the coprime case \((A,B,G)\) has codimension three.
On a general rank-at-least-four pencil member its base with \(G\)
thus has codimension two. Its smooth locus contains a family of
lines disjoint from its vertex and covering that locus: lift lines
from a smooth quotient quadric, and use all linear shears in the
radical. A general translate avoids the base by transitivity.
The nonzero proportional case uses a general line in \(\mathbf P(R)\).

This also covers a pair with a common linear factor and every
coprime pencil with a nonendpoint member of rank at most two.
In the first case write \(A=LA_1,B=LB_1\) with independent
linear forms \(A_1,B_1\). On a general hyperplane
\(V=(A_1+tB_1=0)\) the restrictions are proportional. The
fixed four-plane \(A_1=B_1=0\) does not lie on the factorial
quartic, since a divisor of degree one there would contradict
its Picard group. For general \(t\) the four-plane
\(V\cap(L=0)\) does not lie on the quartic either: if it
varies, its union fills \(L=0\), and if constant it equals
the same fixed four-plane. Thus a general line in \(V\) works.
In the coprime rank-at-most-two case use a hyperplane component
\(V\) of the pencil member. It is not a component of \(A\);
the proper base \(A=B=G=0\) has dimension three, hence a
general line in this five-plane avoids it.

Finally suppose the general pencil member has rank three and
its vertex varies. The vertex of a general member is a
three-plane contained in \(A=B=0\). Indeed for each vector
\(v\) of its kernel extend \(v\) locally to a section
\(v(t)\) of the kernel bundle. Differentiating
\((A+tB)v(t)=0\) and pairing with \(v(t)\) gives \(B(v,v)=0\).
The union of these varying vertices has dimension at least four,
so a general vertex does not lie in the dimension-three base
\(A=B=G=0\). Its intersection with that base has dimension
at most two. On the quadric cone use the family of lines in its
four-plane rulings through general points of the vertex. All
lines in a four-plane ruling meet its vertex hyperplane.
A general such line avoids the base off the vertex, since the
family on the smooth locus is homogeneous and the base has
codimension two. It also misses the indicated intersection
within the vertex. Because the base \(A=B=0\) has codimension
one on the cone and the lines cover it, \(A\) restricts nontrivially
to a general line. This gives the required coprime binary pair.

To prove the angular estimate, cut to the two heavy and two chosen
light directions and then take the least total heavy degree. The
three initial forms are \(F=xy,\ LU,\ G\), where \(L\) is a linear
form with nonzero coefficients of both \(x,y\), and \(U,G\) are
coprime binary forms of degrees two and four in the lights. This
is a regular sequence. The original ordinary tangent presentation
and the dimension comparison under the two specializations
therefore apply. Split the two total binary degrees \(D_H,D_L\)
and all independent outside characters. For the slice of ordinary
order \(j\) and section degree \(m\), put \(e=(j-m)_+\). Its
space consists of the part of the indicated bidegree in
\((F,LU,G)^e\).

Both binary rings are finite free over their two indicated forms,
namely \(F,L\) and \(U,G\). Use homogeneous remainder bases, of
degrees \(\rho\le1\) and \(b\le4\). For fixed light powers \(c,h\)
with \(2c+4h=D_L-b\), the heavy space is spanned by the monomials
\(F^aL^k R_\rho\) with \(2a+k=D_H-\rho\) and
\(a+\min(k,c)+h\ge e\). Thus, apart from bounded endpoint shifts,
it is a complete binary space with fixed factors
\[
 F^{a_0}L^{k_0},\qquad a_0=(e-h-c)_+,\qquad
                          k_0=(2e-2h-D_H)_+ .
\]
Indeed the exact constraints are \(a\ge a_0\) and
\(k\ge2e-2h-D_H+\rho\). Since \(F\) is squarefree, its maximum
marked heavy order is governed by the normalized mean of
\(D_H/2+k_0/2+O(1)\). Here and below the order bounds concern the
mean of a compatible basis, by the complete binary order filtration.

For clarity, the light calculation after fixing heavy marked order
is also complete up to negligible boundary terms. If the marked
heavy point is off \(FL=0\), the allowed heavy orders \(p\) range
from zero to \(D_H-2a_0-k_0\). At a root of \(F\) they range
from \(a_0\) to \(D_H-a_0-k_0\); at the root of \(L\) they range
from \(k_0\) to \(D_H-2a_0\). These statements permit bounded
endpoint shifts. With \(c=(D_L-b-4h)/2\), both \(a_0,k_0\) are
maxima of affine functions of \(h\). Thus for fixed \(p\) each
set of permitted \(h\)'s is an interval cut out by finitely many
linear inequalities. The light series consists of consecutive
monomials \(U^cG^h R_b\) with bounded shifts, giving the full
residual light binary series after taking out fixed powers
\(U^{c_{\min}}G^{h_{\min}}\) at the two interval endpoints.
One can replace all these bounded shifts by the same outward
constant in the linear inequalities: in the \((p,h)\)-polygon
the added lattice points lie in fixed-width strips about finitely
many lines of length \(O(m)\), so their total number is \(O(m)\)
per binary bidegree, including a polygon that collapses to a
segment. At fixed section degree there are \(O(m^{n-2})\) groups:
the independent
coordinates are \(j,D_H,D_L\) and the \(n-4\) outer exponents,
subject to their one total angular degree relation. The resulting
dimension error is \(O(m^{n-1})\), negligible against \(h^0(mH)\).
Every such vector has order \(O(m)A\) uniformly on the product
charts. Therefore the enlarged spaces may be used to estimate all
normalized means. At a root of \(U\) the mean is at most
\(D_L/2+c_{\min}+O(1)\), and at a root of \(G\) it is at most
\(D_L/2+2h_{\min}+O(1)\). These root sets are disjoint.

Scale the section-normalized means by ten. In the normal monomial
representation \(T^t Z_1^i Z_2^k Z_3^h R\), the three symbols
have degrees \(1,2,3\) and means \(1,1/2,1/3\). The remainder in
the four binary variables modulo \((F,LU,G)\) has degree 24.
The heavy \(x\)- and \(y\)-axes have multiplicity eight each by the
transverse complete intersection \((U,G)\); the light roots of
\(G\) carry the remaining degree eight. Equivalently, use the free
remainder bases over \(F,L,U,G\), set \(F=G=0\) and impose
\(LU=0\): the free remaining coordinate is \(L\) in heavy degree
one or \(U\) in light degree two, in proportions two to one.
Hence the means of \(D_H,D_L\) are \(19/6,8/3\). The scaled mean
of \(k_0\) is at most that of \((k-2t)_+\) since
\(e\le i+k+h-t\) whenever \(e>0\), and \(D_H\ge2i+k\).
This upper mean is \(\mathbb E[(E_1/2-2E_0)_+]=1/10\).
Summing over a fixed heavy order only permutes the light-column
masses, so the means of \(h_{\min},c_{\min}\) are bounded by the
means of the actual monomial powers \(h,c\). They are \(1/3\)
and \(1/2+1/6=2/3\), the latter including the free \(U\)-summand
in the remainder. Thus the heavy marked bound is at most
\(19/12+1/20<2\), and the light marked bound at most
\(4/3+2/3=2\). The two radial ratios are half the respective total
degrees. Adjunction along the chosen heavy fiber on the angular
product and on the two binary vertices proves the claimed slope.

The same conclusion holds when \(A=B=0\) and the least pure orders
on \(y=0\) and \(x=0\) coincide. Here pure order means the
exponent of the one surviving heavy coordinate. A common linear
factor with \(q_2\) is already handled by its own extraction.
Modulo \(xy\), suppose first both pure orders are two. Then
\(q_3\equiv x^2 A_1+y^2 B_1+ax^3+by^3\), where \(A_1,B_1\)
are nonzero linear forms in the radical variables. Choose a
hyperplane on which they are nonzero and proportional. On it
choose a general binary section where the common linear factor
\(U\) is coprime to \(G\). In the independent case the base
\(A_1=B_1=0\) is a four-plane and cannot belong to \(G\), since
on the factorial quartic fivefold it would be a Weil divisor of
degree one. In the proportional case irreducibility of \(G\)
suffices. Subtracting a linear multiple of \(q_2\) from \(q_3\)
preserves all the ideals \((q_2,q_3,q_4)^e\). The initial triple
is consequently \(F=xy,\ P U,\ G\), where
\(P=ax^2+by^2\) has two distinct roots disjoint from \(F=0\).
The bigraded calculation applies with degrees \(\deg P=2,\deg U=1\).
For fixed light column \(c,h\), the admissible heavy exponents
satisfy \(a+\min(k,c)+h\ge e\) and \(2(a+k)=D_H-\rho\).
Thus feasibility requires \((D_H-\rho)/2+h\ge e\), followed
only by the lower bound \(a\ge(e-c-h)_+\). The heavy space has
the fixed power of \(F\) times complete forms with the bounded
residue convention above. More formally, for the binary free
remainder basis over \(F,P\), each remainder has degree at most
two. The feasibility cutoff in \(h\) can be rounded outwards by
one to include all remainder parities; its single exceptional
column has \(O(m)\) entries. The same interval and boundary-strip
argument gives complete residual light spaces in each heavy
order. The heavy- and light-axis proportions are now one to two,
as the remaining free exponents have degrees two and one.
The respective binary means are \(10/3,5/2\). The heavy marked
bound is at most \(5/3\), while the light marked bound is at most
\(5/4+\max\{(1/2+2/3)/2,2/3\}=23/12\): the fixed linear
factor has mean exponent at most \(1/2+2/3\), and the fixed
quartic has mean exponent at most \(1/3\).

If both pure orders are three, subtraction of a linear multiple
of \(xy\) makes the cubic a squarefree binary form
\(P=ax^3+by^3\) with \(ab\ne0\). Keep any two radical
directions on which the quartic remains nonzero. The separated
triple is now \(F,P,G\). With \(2a+3k=D_H-\rho\), the condition
in each fixed light column is \(a+k+h\ge e\), an interval in
the heavy exponents. Both heavy forms are squarefree of degree
at least two, so taking out either endpoint fixed power does not
increase the half-degree bound. The residual light blocks are
again intervals in the quartic order. The residual ring gives
the full free simplex slot to the lights. The two degree means
are \(7/2,7/3\), and the light quartic exponent has mean at most
\(1/3\). Thus the respective marked bounds are \(7/4,11/6\).
All these estimates are strictly below two.

\paragraph{The quartic tenfold and two higher dimensions.}
If every line and integral conic on a smooth quartic tenfold has
\(\delta_C(X,H)\ge M\) for some \(M>319/38\), then it is uniformly
K-stable. At \(M=319/38\) it is K-semistable. The corresponding
criteria in dimensions eleven and twelve hold at the following
marks, provided the smooth base of the tangent cone at each
generalized Eckardt point is K-semistable. This condition is
vacuous if the hypersurface has no generalized Eckardt point.
\[
\begin{array}{c|ccc}
n&10&11&12\\ \hline
M_*&319/38&66/7&5785/553\\
S_0&8/5&3/2&121/85
\end{array}
\]
The strict common mark \(M>M_*\) gives uniform K-stability and
the exact mark gives K-semistability.

To see this, the universal cubic-contact comparison just proved
extends to these dimensions, reading
\[
       (n+1)S_H(v)\le n+7/3+(28/15)q
\]
at an irreducible cubic tangent form. The angular computations
acquire only free linear coordinates. We spell out the changes
in the geometric classification. The cubic is now in
\(\mathbf P^{n-1}\). When its singular locus has dimension at most
\(n-5\), a general three-plane again cuts a smooth cubic surface.
If the cubic is nonnormal, the singular component of dimension
\(n-3\) is linear by the same secant argument and the cubic has
at most five essential variables, with the same normal forms.
In the remaining normal case a singular component \(Y\) of
dimension \(n-4\) has degree at most four by the transverse
three-plane calculation. It cannot belong to the smooth quartic
\(G\), since it has codimension two there. As in dimension nine,
its hypersurface span case would make all quadratic polars vanish
on that span; its hyperplane span would make its secants fill
the hyperplane. The nondegenerate degree-four variety of minimal
degree is again excluded: its scroll ruling has projective
dimension \(n-5>(n-2)/2\), or the variety is a Veronese cone
with full ambient tangent at a vertex. For \(m\) essential
variables the dimension of the rank-one projective Hessian image
is at least \((n-4)-(n-m)=m-4\). If \(m\ge7\), the quadratic
Noether-normalization bound \(\binom{m-2}{2}>m\) applies, and
the other values of \(m\) are exactly the small cases already
treated. Through a general \(s\in Y\setminus G\) the directions
on the cubic have dimension \(n-4\), whereas the projection of
\(\operatorname{Sing}F\cap G\) has dimension at most \(n-5\).
The rank-two quadric argument supplies a component outside its
radical, or at rank three the radical has smaller dimension. This
completes the same forward-polar construction in all three columns.

For a cubic with a linear factor, the weighted extraction compared
with \(r=n-2\) has positive radial gap \((6s+3)/(n+1)\), and the
angular interpolation uses a smooth quartic of dimension \(n-4\),
which belongs to the all-member range since \(n\le12\).
For full tangent contact use the generalized Eckardt formula.
When \(n=10\) its base is a smooth quartic eightfold, which is
uniformly K-stable by the preceding theorem. The point threshold
\(n(n+1)/(n+3)\) is strictly greater than \(n-2\) in every
column. The quadratic contact surface bound gives
\(S\ge(n-2)/(n-5)\) at the original hypothetical optimizer.
The cubic bound gives \(S\ge(13n-35)/(13n-71)\). The smaller
of these two numbers is exactly the displayed \(S_0\).

Finally, \(M_*\) is the solution of
\( (M_*-r)S_0+r-nM_*/(n+1)=1\). The ample coefficient
\(r-nM_*/(n+1)\) equals \(7/19,\,5/14,\,190/553\), respectively,
all strictly larger than \(1/3\). Original positive centers are
excluded by hyperplane adjunction from dimension seven and the
short-curve bound as before. Apply the untilted isolated-point
Nadel comparison just above \(M_*\) when the curve mark is
strictly larger. If \(\beta<r\), the actual increase \(t-\beta\)
gives the required strictness just below \(M_*\) as well.

\subsection{A signed estimate for disjoint monomials}
\label{sec:nov28-disjoint-general}
Let \(M_i=\prod_{j=1}^{k_i}x_{ij}^{a_{ij}}\) be monomials in
disjoint nonempty blocks, with all \(a_{ij}>0\),
\(h_i=\sum_j a_{ij}\ge2\), \(b_i=\max_j a_{ij}\),
\(d=\max_i h_i\), and \(m_0=\sum_i k_i\). On the ordinary affine
chart in \(\mathbf P^{m_0+1}\) let \(W=(z+\sum_i M_i=0)\).
Nonzero constant coefficients on the monomials may be absorbed
into the coordinates. At every affine test center,
\[
 (m_0+2)S_{\mathcal O(1)}(v)+v(W)\le D A(v),\qquad
 D=1+\max\left\{
        \frac{1}{1+\sum_i1/b_i},\ \frac1d,\
        \frac{d}{2(d+1)}\right\}.
\tag{DG}
\]

\paragraph{Character reduction and discrepancies.}
Set \(H=\operatorname{lcm}(h_i)\). Use weight \(H/h_i\) on
each variable in block \(i\) and weight \(H\) on \(z\).
Its discrepancy is \(A_0=H(1+\sum_i k_i/h_i)\), the same as
the mean weighted order of the ordinary complete system in units
of \(m_0+2\). Since \(k_i/h_i\ge1/b_i\), the primary radial
denominator \(1+H/A_0\) is at most \(D\).

Let \(G\) be the full, possibly disconnected, diagonalizable
group for which \(z,M_1,\ldots,M_s\) have the same character.
To test a strict angular failure one may split the pieces by
the characters of \(G\), including its finite characters.
Indeed the angular system in section degree \(m\), in the
normalized units above, together with the fixed graph as a
one-dimensional additional system, is a \(G\)-invariant boundary
of linear systems in the sense of \cite[Definitions~4.2--4.3]{ar2004.09413}.
On the coarse weighted exceptional space with its generic
stabilizer boundary in codimension one the log anticanonical
class is \(A_0\) times the positive angular class. The class of
the boundary of systems is \(A_0+H\) in these units. This is
exact at finite level since every ordinary coordinate exponent
in the uniform degree-\(m\) simplex, including its homogenizing
coordinate, has mean \(m/(m_0+2)\). Thus an angular threshold
strictly below \(1/D\) leaves ample remainder. The cited
Theorem~4.4 gives a \(G\)-invariant divisorial minimizer, whose
valuation filtration splits into characters. In particular it is
enough to check every character-adapted basis with uniform errors
tending to zero. This argument on the coarse klt pair is
equivalent to the smooth weighted stack calculation by the
codimension-one stabilizer different.

Write \(Y_0=z\), \(Y_i=M_i\). Each character piece is a fixed
remainder monomial \(\prod x_{ij}^{r_{ij}}\) times an ordinary
homogeneous polynomial space in the \(Y_i\), with fixed degree
\(l\), monomial costs \(h_i\) (put \(h_0=b_0=1\)) and cost
bound \(m-\sum r_{ij}\). Indeed if an original exponent vector
is \(n_{ij}\), set
\(q_i=\min_j\lfloor n_{ij}/a_{ij}\rfloor\) and
\(r_{ij}=n_{ij}-a_{ij}q_i\). The complete character group is
the ambient character lattice modulo the \emph{integral} relations
\(M_i/z\). Thus two monomials share a character precisely when
their remainder vectors agree and \(q_0+\sum_iq_i\) agree,
where \(q_0\) was the original \(z\) exponent. This uses all
components of \(G\), so does not demand a saturated relation
lattice. The original radial degree is fixed by that character.
The aggregate means of the \(q_i\) in normalized units are
\(1/h_i+O(1/m)\), by integrating the simplex where every
exponent in block \(i\) is at least \(a_{ij}t\). Those of the
\(r_{ij}\) are \(1-a_{ij}/h_i+O(1/m)\).

We give the discrepancy comparison on the affine cones since
the monomial fibers can have positive dimension. Write
\(\alpha_{ij}=v(x_{ij})\). If the restriction \(v'\) to the
field of the \(Y_i\)'s is nontrivial, it is a multiple of a
divisorial valuation: the residue transcendence degree drops by
at most the relative transcendence degree. One has
\[
 A(v)\ge A(v')+\sum_{ij}\alpha_{ij}-\sum_{i\ge1}v'(Y_i).
\]
In the ambient volume form replace one differential in each
block by \(dY_i\), using logarithmic differentials for all
the other variables in that block. Relative to their wedge
and the ordinary downstairs volume form the factor is, up to
a constant, \(\prod_{ij}x_{ij}/\prod_{i\ge1}Y_i\). At a model
of \(v'\), the downstairs volume is a unit times the local
parameter to its log discrepancy, times a wedge of logarithmic
differentials. Such a wedge after pullback has at most a simple
pole at the upstairs divisor. This proves the comparison. If
the restriction is trivial then all \(\alpha_{ij}=0\) and no
polynomial or graph contributes any order.

Let \(S\) be the normalized basis order of just the polynomial
pieces downstairs. To prove (DG) it therefore suffices to show
\[
       S(v')+v'(W)+\sum_{i\ge1}c_i v'(Y_i)\le D A(v'),
 \qquad c_i=D-\frac{D-1}{b_i}-\frac1{h_i}.
\tag{DG1}
\]
Indeed in combining the discrepancy inequality with the
remainder exponents the cost for \(\alpha_{ij}\) is
\(D(a_{ij}-1)+1-a_{ij}/h_i\le a_{ij}c_i\), since
\(a_{ij}\le b_i\). Errors are uniform by
\(\sum_{ij}\alpha_{ij}\le A(v)\). Put \(c_0=0\), as is also
given by the formula at \(h_0=b_0=1\).

\paragraph{The ordinary quotient cone.}
For (DG1) the graph is the linear form \(W=\sum_iY_i\) and
its polynomial spaces have homogeneous degree \(l\) and cost
\(\sum h_iq_i\) bounded. All \(c_i\) are nonnegative. The
radial blowup of the origin passes precisely when
\(1\le(D-1)\sum_{i\ge0}1/b_i\), which holds by the first
term in (DG). On the exceptional ordinary projective space
we verify (DG1) by the following finite induction.

The active indices retain their costs and coefficient \(c_i\)
unless their coordinate divisor has been dropped because it
missed the point of interest. All the polynomial pieces still
have a fixed homogeneous degree and an upper cost bound, both
possibly different in every piece. Their aggregate exponent
means remain \(1/h_i+O(1/m)\). At a point outside \(W\) make a
diagonal specialization after dropping the graph. The largest
coordinate coefficient is at most
\(c_i+1/h_i=D-(D-1)/b_i\le D\). At a point of \(W\) having a
nonzero coordinate of maximal cost \(d\), discard its missing
coordinate boundary and replace that coordinate by \(W\).
This triangular substitution preserves every piece in either
direction since all other costs are at most \(d\). Its new
coefficient, graph included, is \(1+1/d\le D\), and all
other coordinates pass the same toric test.

Otherwise every index of cost \(d\) vanishes at the chosen point
on \(W\). At least two coordinates are nonzero there. Use one
of least cost among them as pivot \(Y_s\). Subtract its suitable
multiples from the other nonzero coordinates, after dropping
their missing boundary components. Since the pivot is no more
expensive, both coordinate changes preserve all pieces. The
resulting affine parameters vanish at the point, and the graph
equation is their sum with zero constant term. Refine by its
ordinary point blowup. At any secondary order the discarded
pivot exponent is fixed to the old degree minus that order.
Thus the remaining monomials again have fixed total degree and
the former cost bound reduced by the pivot cost times its fixed
exponent. All leading monomials are independent, and aggregating
the secondary orders simply partitions the same list. The
induction reproduces on the exceptional projective space with
one fewer index.

In this step let \(R\) consist of the remaining indices with
retained boundary and \(J\) of the rest. There is an index
\(j\in J\) of cost at least two: a second nonzero coordinate
is unbranched, and cannot have cost one since then the pivot
also has cost one, whereas there is only one such index.
There is an index \(i\in R\) of cost \(d\): every such
coordinate vanishes at each nonterminating step and stays
unmodified. The secondary radial comparison is
\[
 1\le (D-1)\sum_{i\in R}1/b_i+
                      \sum_{j\in J}(D-1/h_j).
\]
Its right side is at least \((D-1)/d+D-1/2\), which is at
least one by the last term in (DG). When fewer than three
indices survive no further nonterminating step is possible.
All divisor orders and coefficients transfer by successive
ordinary flag adjunction. At every separation of monomials one
may start with a basis adapted to the tested and specializing
orders. The finite induction preserves the uniform errors and
proves (DG1). The original weighted extraction yields (DG) at
the origin. Its positive contracting weighted action preserves
the graph and transfers each strictly smaller finite-level
coefficient from the origin of the initial basis to its whole
affine chart by openness, proving the assertion everywhere.

For example the bounds for \(xy+u^3,xy+u^4\) are \(10/7,13/9\).
For \(x^2+u^3+v^3\) the bound is \(19/13\); the primary
valuation attains it. Thus a quartic-ninefold rank-one quadratic
pencil with nonzero reduced cubic on the common kernel of its
root coefficients satisfies \(\delta_L(X,H)\ge52/7\): select two
general kernel lights, whose binary cubic is squarefree and so
linearly equivalent to \(u^3+v^3\) by projective three-
transitivity. Weights \((3,2,2)\) on the three lights have
radial ratio \(7/6\), and the pole denominator is
\(1/4+(3/4)(19/13)=35/26\).

\subsection{A criterion on the projective hypersurface of a polar form}
\label{sec:new-translated-branch}
Let \(G\) be a homogeneous light form of degree \(h\ge2\) on \(k\)
retained variables in the two-polar calculation of
Chapter~\ref{ch:fullrow}. Write \(H\) for the ordinary hyperplane
class on \(\mathbf P^{k-1}\). Suppose that \(1+1/h\le D<2\) and
that every divisorial valuation on this projective space satisfies
\[
          kS_H(v)+(2-D)v(G)\le D A(v).
\]
Then the positive angular group has normalized denominator at most
\(D\). A sufficient condition is K-semistability of
\((\mathbf P^{k-1},cG)\), where \(c=(2-D)/D\), together with
\(k(D-1)\ge h(2-D)\).

To prove this, let \(P=\mathbf P(1^k,h)\), put \(Z=(z=0)\),
\(W=(z+G=0)\), and write
\(D_0=(k+1)S_{\mathcal O_P(1)}+(1-1/h)v(Z)\).
The pole decomposition in Chapter~\ref{ch:fullrow} requires
\(D_0+v(W)\le D A_P(v)\). Away from \(Z\cap W\) the two cyclic-cover
bounds for the sections \(Z,W\) give the result: if \(v(W)=0\)
use \(D_0\le A_P\), and if \(v(Z)=0\) the translated bound is
\(D_0+v(W)\le A_P+v(W)/h\le (1+1/h)A_P\).
Both sections avoid the weighted vertex.

Now put \(w=v(W)\), \(z_0=v(Z)\), and \(\lambda=2-D\). In
the flag of \(Z\), the fixed coefficient of the complete weighted
series in the normalization of \(D_0\) is \(1/h\), and the residual
is \(k\) times the complete hyperplane series on \(Z\).
These coefficients follow directly by summing monomials of
degree \(m-hi\) after extraction of \(z^i\). Thus adding
\((D-1)z_0\) to \(D_0\) makes the coefficient on \(Z\) equal
to \(D\). Filtered inversion of adjunction and the assumed
projective inequality give
\[
       D_0+\lambda v(G)\le D A_P(v)-(D-1)z_0 .
\]
Translation \(z\mapsto z+G\) preserves the complete weighted series
and \(G\). Using both inequalities bounds \(D_0+w-D A_P(v)\) by
\[
 -\lambda v(G)+
 \min\{w-(D-1)z_0,\ (1-D+1/h)w+(1-1/h)z_0\}.
\]
For \(w\le z_0\) the first term in the minimum is at most
\(\lambda w\). For \(w\ge z_0\) the second is at most
\(\lambda z_0\), since \(1-D+1/h\le0\).
The order of \(G=W-Z\) is at least \(\min(w,z_0)\), as required.
K-semistability of the indicated pair gives the sufficient condition
by \(D(k-ch)\ge k\).

For a quartic ninefold this treats every terminal pencil quartic after
constant linear cuts of total ambient codimension \(c_0\le3\).
The terminal form is independent of the point on the line. The finite
Gauss estimate gives it singular locus of codimension at least
\(11-2c_0\ge5\) in the projective light space, as in
Section~\ref{sec:nov28-q8-cusp-mark}. There are at least five
light variables, so five general ones give a smooth quartic threefold
in \(\mathbf P^4\). All such quartics are K-stable by the index-one
case. The interpolation in the proof of
\cite[Theorem~5.1]{ar2001.10557} makes
\((\mathbf P^4,(5/13)G)\) K-semistable, since
\(5/13<1-(5/4-1)/4\).
Take \(D=13/9\) in the criterion. The regular and positive masses
give \(1/4+(3/4)(13/9)=4/3\), and the radial denominator is
\(1+(3/4)/(1+5/4)=4/3\). It follows that
\(\delta_L(X,H)\ge15/2\).
At order three the same calculation with \(k=4,D=10/7\)
applies whenever four general lights cut a K-semistable cubic
surface \(G\). The needed coefficient \(2/5\) is below
\(1-(4/3-1)/3\). Both radial and combined angular denominators
are \(37/28\), giving \(\delta_L(X,H)\ge280/37\).

\paragraph{A sharper bound for quartic polar nets on lines.}
On a smooth complex quartic ninefold, every line of polar rank
at least three satisfies \(\delta_L(X,H)\ge37/5\).
Set \(b=50/37\). At rank four the ordinary denominator is at
most \(5/4\); at rank three the splitting is \(1,2\), and
the full common kernel has first contact \(2\le h\le4\).
The ternary radial denominator is at most \(4/3\), the branch
and low-summand ratios are at most \(5/4\), and the off-branch
bound is \(1+1/h-1/4\le5/4\). For \(h=4\) the ordinary full
product test of Section~\ref{sec:nov27-next-lines} has entries
\(5/4,6/5,5/4\).

It remains to bound the plane order \(S\) at the moving branch
point \(P\) in the net. Normalize \(v(\mathfrak m_P)=1\) and
put \(u=v(L')\ge1\) for the branch line \(L'\), so \(A\ge1+u\).
For \(h=2,3\) the conditional normalized orders of the pole and
regular sectors, with masses \(1-h/4,h/4\), satisfy
\[
 S_{\rm pos}\le (1+1/h)A-1/h,\qquad S_{\rm reg}\le A-1/2
\]
by the line version of Section~\ref{sec:nov27-quintic-line-quadratic}.
At \(u=1\) these give \(S\le(5/4)A-1/2\) for \(h=2\) and
\(S\le(13/12)A-11/24\) for \(h=3\), which suffice for
\(S\le b[A-(1-1/h)u]\) by testing \(A=2\) and the slopes.
When \(u>1\), the center on the first blowup of \(P\) is the
branch direction. The regular group is already a full product of
base and plane degree, so on this group apply the single marked
line graph calculation with \(t=1\). It gives radial coefficient
\(1+1/2=3/2\) and angular coefficient at the branch direction
\[
          J=\frac14\left(1+\frac{3-1}{2\cdot4}\right)=5/16.
\]
Indeed the high polar form has base degree three, and the
moving binary equation base degree one, so this is exactly the
\(u=1\) endpoint of the expression for \(J\) in
Section~\ref{sec:nov27-next-lines} (the variable there is the
inverse of the product enlarging factor). Adjunction on the
exceptional line therefore gives
\(S_{\rm reg}\le3/2+(5/16)(A-2)\). Averaging now gives
\[
 S\le(29/32)A+3/16\quad(h=2),\qquad
 S\le(109/192)A+55/96\quad(h=3).
\]
The separate branch intersection comparison is
\(S\le u/h+(A-1)I/u\), \(I=3/2\).

For \(h=2\) the first bound passes at \(A=1+u\) up through
\(u_0=I/b=111/100\): its slack is affine decreasing and is
positive at \(u_0\), since
\[
 (50/37)(1+111/200)>(29/32)(211/100)+3/16 .
\]
Its slope in \(A\) is below \(b\). For \(u\ge u_0\) the
intersection slope is at most \(b\), and its minimal discrepancy
comparison reduces to \(I\le b+(b-1)u/2\), valid even at
\(u=1\). For \(h=3\) use the first bound up through \(u=4/3\);
its endpoint slack is positive since
\((50/37)(1+4/9)>(109/192)(7/3)+55/96\).
For \(u\ge4/3\) both \(I/u\le b\) and
\(I\le b+(b-1)u/3\). This covers each branch center and proves
the asserted line mark.

\subsection{A sufficient short-curve bound for quintic elevenfolds}
\label{sec:nov28-quintic-eleven-point}
Let \(X\subset\mathbf P^{12}\) be a smooth quintic elevenfold and
\(H=\mathcal O_X(1)\). If every line and every integral conic on
\(X\) satisfies \(\delta_C(X,H)\ge857/100\), then \(X\) is uniformly
K-stable. In particular the rational mark \(60/7\) suffices.

Assume that \(\beta=\delta(X,H)\le8\). The smooth pointed
hyperplane section in dimension ten and the covering curve bound
give a line or a conic in any positive-dimensional optimizing
center, as in Section~\ref{sec:nov28-direct-point}. Thus an
optimizing divisor is centered at a closed point \(p\). Normalize
\(v(\mathfrak m_p)=1\), put \(A_X(v)=11+q\), and write
\(S=S_H(v)\). We claim \(q\ge2\), hence \(S\ge13/8\).

\paragraph{The surface alternatives.}
Take \(B=131/50\). The signed pointed surface bound
\[
                 12S\le B(q+2)+9
\]
applies if \(\tau_p(H)\le B\). In the covering case of its proof
the only additional issue is the inequality for the integral
degree \(a\) and multiplicity \(k_0\) of a nonplanar covering curve.
Here \(a\ge 2d-5=5\) and \(k_0(k_0+1)\le a^2/5+a\).
For \(a=5,6\) this gives respectively \(k_0\le2,3\), and
\(5k_0/a\le5/2<B\). For \(a\ge7\) one has
\((B^2/5-1)a+B-5>0\), which gives \(5k_0/a\le B\) directly.
The planar and proper swept-locus cases are unchanged.
The displayed bound contradicts \(11+q\le8S\) throughout
\(0\le q\le2\), by checking its two endpoints.

If \(t=\tau_p(H)>B\), let \(j\) and \(k\) be the degree coefficient
and multiplicity of the high-ratio prime. The general-section
discriminant argument of Section~\ref{sec:nov-second-row} gives
\[
 t=k/j,\qquad k(k-1)\le I\le5j(j+1)+k+1 .
\]
Thus \((k-1)^2\le5j(j+1)+2\). For \(j\ge6\) this gives
\(k/j\le1/6+\sqrt{5+5/6+2/36}<B\). Direct checking at
\(j=1,\ldots,5\) gives the only possibilities below \(t=5\):
\[
               (1,3),\ (1,4),\ (2,6),\ (3,8),\ (4,11).
\]
At \((1,3)\), the sectional Milnor estimate gives
\(I\le2[1+(4^{11}/2)^{1/10}]<11\), hence \(I\le10\).
At the three entries \(j\ge2\), use respectively \(I\le37,69,112\).
The discriminant contribution per unit degree is at most
\(L=I/[j(2\sqrt I-1)]\), by maximizing the usual expression in
the allowed multiplicity. In the four cases just discussed one
can take \(L\le2,17/10,3/2,7/5\), respectively.
In each case \(5/t+L(1-5/t^2)<B\). The signed surface argument
therefore gives the same contradiction. At \(t=5\) the point is
generalized Eckardt. The base is a smooth quintic ninefold, so Theorem~\ref{quintic:main} gives K-stability.
The point threshold \(132/15\) exceeds eight.

\paragraph{Quartic tangent contact.}
At the remaining entry \((j,k)=(1,4)\) the tangent hyperplane
section has equation \(Tf(y)+g(y)=0\), with eleven tangent
variables, quartic \(f\) and quintic \(g\). Suppose first that
no component of \(f\) has multiplicity greater than two.
A general replacement \(g+Lf\) cuts out a smooth quintic
\(G\subset\mathbf P^{10}\), by the finite Gauss map and the base
point argument of Section~\ref{sec:nov28-reverse-polar}.
Every irreducible factor of \(f\) restricts to an integral divisor
on \(G\): \(\operatorname{Pic}(G)=\mathbf ZH\), all sections of
the required degrees lift, and the ideal of \(G\) has no part in
degree at most four. Distinct factors cut distinct divisors by the
same lifting argument. Lines cover \(G\), since their direction
equations at a smooth point have degrees two, three, four and five
on \(\mathbf P^8\). A general line of a covering component is
free, and its smooth universal evaluation gives transverse
intersections with the support of these divisors as in the cited
section. Choose a quartic normal polar \(P\) on the line avoiding
all marked roots of \(f_0=f|_{\ell}\). Projection to the plane of
the line and normal, followed by taking the least \(z\)-power as
there, gives angular initials \(f_0,zP\).

For clarity the two ordinary hypertangent symbols have degrees
one and four and orders four and five. Their relations
\(T^3 Z+f,\ TW+g\) are a regular sequence modulo \(T\);
localization and ordinary Hilbert series prove the presentation
as in the cubic-contact argument. The order-\(j_0\) piece on the
ordinary exceptional space is
\[
 \sum_{3i+h\ge j_0-m} f^ig^h
                    \mathbf C[y]_{j_0-4i-5h}.
\]
The regular sequence property of the two initial forms preserves
all of these dimensions. In binary degree \(D\) and transverse
\(z\)-exponent \(s_0\), put \(l=j_0-m\). Over the finite free
algebra \(\mathbf C[f_0,P]\), a binary monomial
\(f_0^aP^cR_b,\ D=4(a+c)+b\), belongs to the block precisely when
\(3a+\min(c,s_0)\ge l\). Equivalently
\[
      a\ge(l-s_0)/3,\qquad 2a\ge l-(D-b)/4.
\]
The remainder degrees \(b\) are bounded by six. Hence there is a
fixed \(f_0^I\) with
\(I=\max(0,\lceil(l-s_0)/3\rceil,\lceil(l-D/4)/2\rceil)\);
after it is removed the binary system has bounded codimension in
degree \(D-4I\), since only a bounded number of endpoint monomials
for each remainder can be missing. Every root has multiplicity
at most two, and the normalized marked mean in this binary block
is consequently at most \(2I+(D-4I)/2+O(1)=D/2+O(1)\).

In units of twelve the two symbol exponent means are one and
\(1/4\). The residual ternary complete intersection
\((f_0,zP)\) has degree twenty, sixteen on the \(z\)-axis and
four on \(z=0\). Indeed its zero \(z\)-degree part is
\(\mathbf C[x,y]/(f_0)\), of eventual Hilbert function four,
and every positive \(z\)-degree part has binary length sixteen.
The scaled binary and transverse means are thus
\(4+1+1/5=26/5\) and \(1/4+4/5=21/20\).
The other eight independent coordinates have mean one. On the
blowup of the binary vertex both its radial bound and the
marked binary orders are therefore at most \(13/5\); the coordinate
divisors are transverse to this resolution. The ordinary radial
order has scaled mean \(9+4+5/4=57/4\). Filtered inversion of
adjunction now yields
\[
              12S\le57/4+(13/5)q .
\]
Together with \(11+q\le8S\) this gives \(q\ge45/22>2\).

Finally suppose a linear factor \(L\) occurs in \(f\) to
multiplicity \(s=3\) or four. Give the ten other tangent variables
weight \(s\), and \(L\) weight \(s+1\). The normal tangent coordinate
has weight \(5s\). Put \(G=g|_{L=0}\), now a quintic in ten
variables. It is nonzero and has only isolated projective
singularities. Otherwise a complete singular curve meets
\((\partial_Lg)|_{L=0}=0\), and the entire \(T\)-line over that
point is singular on \(Tf+g\), contrary to finite Gauss. If
\(G=0\), the tangent section has the factor \(L\), which contradicts
its integrality. Specializing \(L\) to zero in the initial
relation gives \(T^4Z+G=0\). Primitivity and the Hilbert series
justify this specialization as in the weighted linear-factor
extraction of Section~\ref{sec:q8-tangent-tests}.
The extraction has \(A=11s+1,\ 12S_H=(15s+1)\); hence its
discrepancy at coefficient eight is positive. On the angular
space the free \(L\)-coefficient is \(1/12\), the quintic
coefficient \(1/15\), and the mobile \(H\)-degree \(10/12\).
Indeed in a normal basis \(T^aL^bZ^hR_e\) the forced exponent
of \(G\) in its fixed weight and \(L\)-character is
\((h-a/4)_++O(1)\), with mean in units of twelve
\(\mathbf E[(E_1-E_0/4)_+]=4/5\) for independent unit exponential
variables. Removing those powers leaves the stated mobile
degree. The opposite \(L\)-axis has radial discrepancy ten,
exceeding \(8(14/12)\).

On \(\mathbf P^9\) it therefore suffices to show
\((8/15)v(G)+(20/3)S_H(v)\le A(v)\).
A general hyperplane cuts \(G\) in a smooth quintic sevenfold
\(R\subset\mathbf P^8\), K-stable by the preceding range.
The log interpolation inequality gives K-semistability of
\((\mathbf P^8,(8/15)R)\), since \(8/15<1-(9/5-1)/8\).
Degenerating to its projective cone with coefficient \(8/27\)
at infinity gives, by the cone criterion and openness, a
K-semistable pair with the original \(G\), of anticanonical
coefficient \(10(1-8/27)=190/27>20/3\).
Dropping the divisor at infinity proves the required inequality.
The weighted extraction therefore excludes this optimizing point.
This finishes the proof of \(S\ge13/8\).

\paragraph{Interpolation.}
Use Section~\ref{sec:nov28-direct-point} with \(f=1/5\),
\(t=71/10,\ \lambda=6/5\). Here \(1+g=13/12\); the upper testing
factor and the lower removed hyperplane coefficient satisfy
\(h<49/40,\ \alpha>1/20\). These follow by the integer comparisons
\(5\cdot431^{11}<4\cdot440^{11}\) and \(5\cdot47^{11}<55^{11}\).
Thus \(t+h\lambda<857/100<60/7\), while
\(a=t+(1+g)\lambda-8=2/5\), and
\[
 \eta=8-(11/12)t-(1-\alpha)\lambda>211/600>7/20.
\]
The upper bound \(\eta<1\) follows, for example, from
\(\alpha\le1/12\). Consequently \(a(13/8)+\eta>1\), the required
contradiction to the Nadel inequality.

\subsection{All higher-rank line polars on a quintic elevenfold}
\label{sec:nov28-quintic-eleven-nets}
Every line \(C\) of polar rank at least three on a smooth quintic
elevenfold satisfies \(\delta_C(X,H)\ge60/7\). The normalized
denominator to prove is \(b=7/5\). If the largest split degree is
at most two, the ordinary denominator \(1+2/5\) and radial
denominator \((11-1/5)/10\) suffice. The positive degrees sum
to four. Thus the only case to consider is a net of splitting
\((1,3)\).

Use the marked net of Section~\ref{sec:nov27-quintic-line-quadratic}
and its notation \(h,P,L\). The first light order is between two
and five since the full polar kernel of ambient codimension three
cannot lie on \(X\). The ternary-axis denominator is
\((1+14h/5)/(1+2h)<7/5\), and the other unmarked elementary
terms are at most \(7/5\). The line product test of
Section~\ref{sec:nov27-next-lines} handles \(h=4,5\):
its three terms are respectively
\[
 (19/15,704/525,512/399),\qquad(6/5,5/4,13/10).
\]

For \(h=3\) set \(t=2,\ R=6/5,\ u_0=1/2\) in that product test.
The comparison off \(P\) is
\(R(1+1/(t+4))=7/5\). On blowing up \(P\) the unmarked
direction comparison is \(R=6/5\); the marked angular comparison
is \(RhJ=4/3\) since \(J=10/27\). Use the actual series for the
radial comparison on this blowup. The marked regular and positive
group bounds of Section~\ref{sec:nov27-quintic-line-quadratic}
give
\[
 S_{\mathrm{actual}}(\operatorname{ord}_P)
    \le(17/15)\,2-13/30=11/6 < (7/5)(1+1/3).
\]
For precision, after clearing scalar base denominators an original
block is a subspace of the enlarged block of the product test.
Filter both by the order at \(P\) on the geometric generic plane.
This order is independent of the scalar base frame, hence their
associated graded blocks inject in every order. Any adapted
original basis yields independent leading jets in the corresponding
enlarged block. Their nonnegative angular orders at a given point
of the exceptional line can be bounded by extending them to a full
basis there and using the established product computations with
the enlargement factor \(R\). This lets one test the original
radial order and the enlarged angular orders separately. Blowup
adjunction now gives the asserted denominator at \(h=3\).

At \(h=2\) normalize \(v(\mathfrak m_P)=1\) on the coarse plane
and write \(u=v(L)\), \(A=u+1+x\), with \(u\ge1,\ x\ge0\).
The required inequality is \(S\le(7/5)(A-u/2)\). If \(u=1\)
the marked estimate \(S\le(13/10)A-1/2\) suffices since \(A\ge2\).
For \(u>1\) the center after the first blowup is the direction of
\(L\). In the regular product alone its radial coefficient is
at most \(3/2\) and the angular order at this direction is at most
\(J=(1/5)(1+3/10)=13/50\), by the regular-product line formula
with base expansion factor one. Taking a simultaneous radial-
and valuation-adapted basis, inversion along the exceptional line
gives \(S_{\rm reg}\le3/2+(13/50)(A-2)\). Indeed its discrepancy
for the reduced exceptional boundary is \(A-2\) under the
normalization above. The positive bound is
\(S_{\rm pos}\le(3/2)A-1/2\). Averaging with masses \(2/5,3/5\)
yields \(S\le(251/250)A+23/250\). Its slack against the desired
bound is
\[
               \bigl(99x-76(u-1)\bigr)/250 .
\]
The independent residual intersection bound,
\[
            S\le u/2+(A-1)\min(8/5,8/(5u)),
\]
has slack at least \(((u-1)-x)/5\), by using the second term
for its value at \(x=0\) and the first for its slope in \(x\).
This covers \(0\le x\le u-1\), and the previous bound covers
\(x\ge u-1\). All comparisons are in the normalized units of
the net extraction; weighted adjunction proves the claim.

\subsection{Removing part of a marked hyperplane at an optimal point}
\label{sec:nov28-two-hyperplanes}
Write \(N=n+1,\ -K_X=rH\). Suppose a divisorial minimizer \(v\)
of \(\beta=\delta(X,H)\le r\) is centered at a closed point \(p\).
Normalize \(v(\mathfrak m_p)=1\), and suppose \(\beta>N/2\).
Take a hyperplane \(T\) missing \(p\) and a general hyperplane
\(\ell\) through \(p\), with \(v(T)=0,\ v(\ell)=1\).
Then, for each rational \(0<\theta<\beta/(2n)\), there are
asymptotically \(v\)-optimal basis divisors in sufficiently
divisible degrees from which one can remove the usual coefficient
\(1/N+o(1)\) of \(T\) and also \(\theta/N\) of \(\ell\).
Their upper fractions for \(v\) retain the ordinary hyperplane
estimate for removal of \(T\).

\paragraph{The fixed and moving components on the optimal degeneration.}
Let \(Y\) be the klt central fiber of the special degeneration
induced by \(v\). Its polarization \(L=-K_Y/r\) satisfies
\(\delta(Y,L)=\beta\) by Blum--Liu--Zhou. The hyperplane class
has an integral Weil extension: take the closure of \(T\) on the
total space over the punctured line. A sufficiently divisible
multiple differs from the integral Cartier test polarization by
a divisor supported on the central fiber, hence by an integral
multiple of that fiber. It is therefore \(\mathbf Q\)-Cartier.
The total space is regular at the codimension-one points of \(Y\),
since \(Y\) is a normal Cartier divisor there. Restriction gives
an integral Weil class on \(Y\) representing \(L\). In a
trivialization on this locus, divide the two original sections by
their respective orders in the base parameter. Their nonzero initial
sections \(T_*,\ell_*\) thus have integral Weil zero divisors.
A common integral base-parameter twist only shifts both weights.
The weights of these two sections differ, so they are independent.

Each of these zero divisors is reduced. In fact if a member of
\(|L|\) contains a prime divisor \(D\) with multiplicity \(a\), its
powers in divisible degrees give \(S_L(D)\ge a/N\). Equivalently,
for \(0\le x\le a\) the pullback section yields the volume lower
bound \(\operatorname{vol}(L-xD)\ge(1-x/a)^nL^n\), computed on
a resolution if needed. Since \(A_Y(D)=1\), \(a\ge2\) contradicts
\(\beta>N/2\).

Resolve the pencil \([T_*:\ell_*]\) by \(\pi:\widetilde Y\to Y\)
and denote its morphism by \(h:\widetilde Y\to\mathbf P^1\).
Write \(P=\pi^*L,\ M=h^*\mathcal O_{\mathbf P^1}(1)\), and let
\(D'\) be the sum of the components of \((T_*=0)\) on which
\(\ell_*\) does not vanish identically. Since \((T_*=0)\) is
reduced, the movable member on \(\widetilde Y\) differs from the
strict transform of \(D'\) by exceptional divisors. Hence
\[
       F=\frac{L^{n-1}D'}{L^n}=\frac{P^{n-1}M}{P^n}
              \ \ge\ \frac{\beta}{2n}.
\]
For the inequality take a general prime fiber \(E\) of the Stein
factorization of \(h\). It maps onto a prime divisor on \(Y\),
has \(A_Y(E)=1\), and is nef. Set \(x=P^{n-1}E/P^n\), so
\(0<x\le F\). The algebraic Morse inequality gives
\(\operatorname{vol}(P-tE)/P^n\ge 1-nxt\).
Integrating between \(0\) and \(1/(nx)\) gives
\(1/\beta\ge S_L(E)\ge1/(2nx)\), as required.

\paragraph{Counting the simultaneous factors.}
Choose \(k\) sufficiently divisible that the \(k\)-th Veronese
of the graded degeneration is the section ring
\(A=\bigoplus_{b\ge0}H^0(Y,bkL)\), generated in degree one.
It is the associated graded in these section degrees. Put
\(U=T_*^k,\ V=\ell_*^k,\ B=A/(U),\ C=B/(0:_B V^\infty)\).
The Hilbert multiplicity ratio \(e(C)/e(B)\) is \(F\).
Indeed the saturated annihilator removes exactly the top-dimensional
components on which \(V\) vanishes. At every other component its
localization is zero. The associativity formula for Hilbert
multiplicity counts the generic lengths, and the Cartier divisor
of \(U\) has the same multiplicity \(k\) at every prime of its
reduced support.

In total degree \(b\) use the successive \(U\)-order layers
\(U^i B_{b-i}\). Within \(B_{b-i}\) filter by
\(V^j B_{b-i-j}\). These spaces are homogeneous for the valuation
weight. Choose a homogeneous basis of all successive quotients and
represent each selected vector by \(U^iV^j a\) with \(a\)
homogeneous for that weight. Lift \(a\) to an original section
with that initial form and multiply by the actual \(T^{ki}\ell^{kj}\).
The resulting sections give a basis with maximal \(v\)-order,
because their initial classes are a basis of \(A_b\). Its tagged
orders \(ki,kj\) are lower bounds for actual divisibility.
The mean tagged \(T\)-coefficient in units of \(kb\) tends to
\(1/N\) by the successive hypersurface layers. Its actual
coefficient has the same limit since no basis exceeds the full
\(T\)-order mean. Hence the distribution of tagged \(T\)-orders
gives the usual hyperplane upper-fraction bound, equivalently by
the uniform distribution tails from the exact powers of \(U\).
For \(j\) beyond a fixed integer, the annihilator of \(V^j\) in
\(B\) equals its saturated annihilator. Thus
\(\dim(V^jB_{b-i-j})=\dim C_{b-i-j}\). Summing the positive
tagged orders over all \(i\) and \(j\ge1\) gives
\[
 \frac{1}{b\dim A_b}\sum_{i,j\ge0,\ j\ge1}\dim(V^jB_{b-i-j})
                \ \longrightarrow\ \frac{e(C)}{N e(A)}
                =\frac{F}{N}.
\]
Here \(e(A)=e(B)\) since \(U\) is regular of degree one. The
finitely many exceptional \(j\)'s contribute \(O(1/b)\)
after normalization, and the remaining sum is the double
Hilbert summation. This proves the extraction assertion.

\paragraph{The line criterion.}
In the direct interpolation keep the initial boundary
\((\beta-e)R_m\) with only \(T\) removed. From its endpoint
\(tR_m+\lambda R_m^+\) subtract a coefficient \(t\theta/N\)
of \(\ell\) from the full basis term. This can be done for all
large \(m\) by the count just proved. It reduces the endpoint
class and its \(v\)-order by the same \(t\theta/N\), since
\(v(\ell)=1\). The ample remainder changes from \(\eta\) to
\(\eta'=\eta+t\theta/N\), and the limiting Nadel expression
is still \(aS+\eta\). The upper bound at every testing valuation
only decreases. If \(1/2<\eta'<1\) and \(aS+\eta>1\), the
endpoint has \(v\)-order greater than \(A_X(v)\) and the
tiny-point perturbation applies throughout the segment. The
initial ample coefficient is greater than \(1/2\). The first
positive center therefore has a rational curve of \(H\)-degree
less than two by the length-one theorem. The original optimizer,
if positive-dimensional, already contains a line by the
canonical-extension bound \(H\cdot C\le N/\beta<2\).
Thus the direct interpolation requires only the line bound.

The line obtained from a point optimizer passes through \(p\).
Indeed let \(Z_0\) be the union of the finite non-qlc scheme of
the perturbed pair at the first positive center and all its
zero-dimensional qlc strata. Quasi-log adjunction allows this
union, and quasi-log vanishing with the ample remainder gives
\(H^1(X,\mathcal I_{Z_0})=0\). The constants therefore surject
onto \(H^0(Z_0,\mathcal O_{Z_0})\). Since \(p\) is non-lc,
the finite scheme has support exactly \(\{p\}\). For a minimal
positive-dimensional qlc center \(W\), apply the same vanishing
to \(W\) union the non-qlc scheme. That union is connected, so
\(p\in W\). Its induced quasi-log structure has non-qklt locus
supported at \(p\), by minimality and the assertion about
zero-dimensional strata. Take a curve in \(W\) through \(p\)
and apply the length-one theorem of Fujino--Hashizume to the
contraction to a point. The closure of the rational curve from
their Theorem~1.8 meets this non-qklt locus. Its degree bound
therefore supplies the required line through \(p\).

In dimension nine for quartics take \(t=125/17,\lambda=0\),
\(S_0=7/4,\theta=1/4\). One has
\[
       \eta=13/34,\qquad \eta'=77/136,\qquad
       (t-7)(7/4)+\eta=1 .
\]
Increasing \(t\) slightly gives uniform K-stability under any
common line bound \(M>125/17\). At \(M=125/17\), assuming
\(\beta<7\) makes \((t-\beta)S+\eta>1\) already at equality,
so take \(t<M\) sufficiently near \(M\) to prove K-semistability.
For quartic eightfolds the limiting value \(t=51/8,\lambda=0\)
gives \(\eta'=49/96\) at \(\theta=1/4\); the small positive
perturbation of \(t\) itself gives strictness with \(\lambda=0\).
For quintic elevenfolds use \(t=71/10,\lambda=6/5\)
and \(\theta=3/11=N/(4n)<\beta/(2n)\). This first gives
\[
       \eta'>211/600+71/440>1/2,
\]
and \(\eta'<1\) from \(\alpha\le1/12\). The sufficient line-only
mark is \(857/100\), in particular \(60/7\).

\paragraph{A sharper elementary moving-degree estimate.}
Take a general prime fiber \(E\) in the Stein factorization of
the pencil on \(\widetilde Y\), and put \(x=P^{n-1}E/P^n\le F\).
Its image on \(Y\) is a prime divisor and \(A_Y(E)=1\). The
fiber divisor \(E\) is nef, \(P\) is nef, and \(P-E\) is
pseudoeffective since it is the sum of a fixed effective
divisor and the other fibers of a general pencil member.
The algebraic Morse inequality for two nef divisors gives
\(\operatorname{vol}(P-tE)/P^n\ge1-nxt\). One can also see
this directly by bounding the cost of vanishing on \(j\)
distinct fibers: on each the number of restrictions is at
most \(m^{n-1}P^{n-1}E/(n-1)!+O(m^{n-2})\), uniformly on
general fibers. For the volume argument the numerical classes
of the fibers are equal.
For \(0\le a<1,\ nxa\le1\), concavity of the volume to the
power \(1/n\) between \(a\) and \(1\) gives
\[
 \frac1\beta\ \ge\ S_L(E)\ \ge\
 a-\frac{nxa^2}{2}
       +\frac{(1-nxa)(1-a)}{n+1}.
\]
In dimension eleven, \(x\le3/8\) and \(a=1/6\) would give
\(1/\beta\ge151/1152\). Thus \(\beta\ge31/4>1152/151\)
implies \(F>3/8\).

It follows that for quintic elevenfolds a common line bound
strictly larger than \(144/17\) gives uniform K-stability,
and the exact bound \(144/17\) gives K-semistability. Indeed
take \(\lambda=0,\ S_0=13/8,\ t_0=144/17\).
For \(\beta\ge31/4\), at \(t_0\)
\[
       \eta=4/17,\qquad (t_0-8)S_0+\eta=1,\qquad
       \eta+\frac{(3/8)t_0}{12}=\frac12 .
\]
Choose \(\theta\) a little above \(3/8\) and below the true
fraction \(F\); the tagged-factor count permits every such
\(\theta\). Increasing \(t\) slightly makes the Nadel expression
strict while preserving \(\eta'>1/2\). If \(\beta<8\), replacing
\(8\) by \(\beta\) already gives strictness, so \(t\) can be
decreased slightly below \(t_0\). For \(6<\beta\le31/4\), take
\(t=t_0-1/8\), \(\theta=1/4<\beta/(2n)\). Direct substitution
gives \(1/2<8-(11/12-1/48)t<1\) and
\((t-31/4)S_0+8-(11/12)t>1\). This handles both conclusions.

\paragraph{Four further degree rows.}
For \(d=6,7,8,9\), put \(n=2d+1\), \(N=n+1\), \(r=d+3\).
If every line \(C\) on a smooth \(X_d^n\) satisfies \(\delta_C(X,H)\ge M\)
with the entry in the first row below, then \(X\) is uniformly K-stable.
\[
\begin{array}{c|rrrr}
d&6&7&8&9\\ \hline
1000 M&9510&10491&11485&12484\\
10^4t&86870&93143&98684&104007\\
10^4\lambda&6966&10198&14264&18639\\
1000\theta&371&357&346&338\\
1000Q&2167&2794&3285&3667\\
\lfloor10^8(1/d)^{1/n}\rfloor&
87124883&87833488&88486508&89079282\\
\lfloor10^8(1-1/d)^{1/n}\rfloor&
98607315&98977591&99217598&99382006
\end{array}
\]
Use the previously proved point estimates \(S_0=(n+Q)/r\).
Write \(b=r/N,\ f=1/d\), and use the \(u,g,\alpha,h,\eta\)
of Section~\ref{sec:nov28-direct-point}. If \(\beta\ge r-1/8\),
then \(F>\theta\) for the indicated value. In fact if \(x\le\theta\),
the prime-fiber inequality evaluated at \(a=(1-\theta)/((n-1)\theta)\)
would give
\[
 \frac{1}{\beta}\ge
 \frac{1}{N}\left(1+
   \frac{n(1-\theta)^2}{2(n-1)\theta}\right)>
   \frac{1}{r-1/8}.
\]
Here \(\theta>1/n\) ensures \(a<1,\ n\theta a<1\). With the
parameters in the table, all four columns satisfy
\[
 t+h\lambda<M,\qquad t+(1+g)\lambda>r,\qquad
 (t+(1+g)\lambda-r)S_0+\eta>1+3/10^4,
\]
\[
  1/2+2/10^4<\eta+t\theta/N<51/100 .
\]
The last two integer rows give rational bounds for these comparisons
by raising consecutive integers to power \(n\) and substituting in
the definitions of \(h\) and \(\alpha\).
If \(\beta\le r-1/8\), use \(t'=t-1/8\) and
\(\theta'=1/4<\beta/(2n)\), with the same \(\lambda\).
Then \(t'-\beta\ge t-r\), the unmarked ample coefficient increases
by \(u/8\), and the new marked ample coefficient satisfies
\[
  53/100<\eta+u/8+(t-1/8)/(4N)<58/100.
\]
The line criterion applies in both intervals.

\section{Polar pencils on quintic elevenfolds}
This section studies rank-two lines in dimension eleven. Here the
\emph{light space} is the common kernel of the normal polars along
the line. The first nonzero term on it has normal degree two,
three, four, or five. Section~\ref{sec:nov28-two-hyperplanes}
reduces the global problem to a uniform estimate on these lines.
The degree of that first nonzero term organizes the local tests
below.

\subsection{Signed bounds for diagonal graphs}
\label{sec:nov28-diagonal-graphs}
Let \(e_1,\ldots,e_k\ge2\), \(d=\max e_i\) and
\(W=(z+\sum_{i=1}^k x_i^{e_i}=0)\) in the affine chart of
\(\mathbf P^{k+1}\). Write
\[
 D=1+\max\left\{\frac1{1+\sum_i1/e_i},\quad\frac1d,\quad
                      \frac{d}{2(d+1)}\right\}.
\]
Every affine center satisfies
\[
       (k+2) S_{\mathcal O(1)}(v)+v(W)\le D A(v).
\tag{Q9.2}
\]
The primary radial valuation attains the first of these three
excesses. In particular the bound is sharp for \((2,3,3)\) and
\((2,3,4)\), giving \(19/13\) and \(37/25\), respectively.

Put \(h=\operatorname{lcm}(e_i)\). Give each \(x_i\) weight
\(h/e_i\) and give \(z\) weight \(h\). The ambient complete
series has radial mean \(A_0=h(1+\sum_i1/e_i)\) in units of
\(k+2\), equal to the discrepancy, while the graph contributes
\(h\). This gives the first bound. On the angular stack take
the branch-coordinate power map to ordinary \(\mathbf P^k\)
with coordinates \(Y_0=z, Y_i=x_i^{e_i}\). Its boundary is
\(\sum_{i=1}^k(1-1/e_i)(Y_i=0)\). The graph becomes the
hyperplane \(W=(\sum_{i=0}^kY_i=0)\). We specify briefly why
one may split the root residues on the original angular piece.
Add the weight-\(h\) graph as a one-dimensional system with
coefficient one to the angular boundary of linear systems. The
degree of its total class is exactly \(A_0+h\) in units of the
positive angular class, since in the full degree-\(m\) simplex
each ordinary coordinate has mean \(m/(k+2)\). The log
anticanonical degree on the coarse weighted space with its
stabilizer boundary equals \(A_0\). As \(1/D\le A_0/(A_0+h)\),
any strict failure of the desired threshold leaves an ample
remainder. The finite-invariant minimizer argument proved above,
using \cite[Theorem~4.4]{ar2004.09413}, now gives a minimizer
invariant under the coordinate roots of unity. Its filtration
splits into those characters. The bounded root monomial can be
removed at \(O(1/m)A\) cost. Each full remaining piece is the
pullback of the following monomial space in ordinary homogeneous
degree \(H\):
\[
   Y^{\alpha},\qquad \sum_{i=0}^k\alpha_i=H,
       \qquad \sum_{i=0}^k e_i\alpha_i\le m',\qquad e_0=1,
\]
where \(m'\) is the section degree minus the bounded residue cost.
The finite crepant formula gives exactly the displayed branch
boundary. Across all pieces the normalized exponent of \(Y_i\)
has mean \(1/e_i+O(1/m)\). We can therefore work on \(\mathbf P^k\)
with its ordinary discrepancies and this branch boundary.

We verify the angular inequality by successive ordinary point
blowups and linear changes of coordinates. The induction retains
a projective coordinate set \(I\) with its original costs
\(e_i\), whose monomial pieces have the same displayed form
with a possible new fixed degree and cost bound for each piece.
There is still the hyperplane given by the sum of all surviving
coordinates. Some surviving root coordinates carry their original
boundary, while the others are unbranched at the current point;
we may drop any divisor missing that point. All retained exponent
means remain \(1/e_i\) up to the same bounded errors. At a point
off \(W\) the toric estimate is immediate after dropping \(W\):
the coefficient of any coordinate is at most
\(1/e_i+D(1-1/e_i)\le D\). This applies to every basis by a
separating diagonal specialization. At a point on \(W\) with a
nonzero coordinate of largest cost \(d\), drop its missing root
boundary and replace that coordinate by \(W\). This triangular
change and its inverse preserve the spaces since every cost is
at most \(d\). Its coefficient including the graph is
\(1+1/d\le D\); the other coordinates obey the same toric
test. This argument also applies if another step has previously
dropped additional branch components.

It remains that all cost-\(d\) coordinates vanish at the point
on \(W\). At least two other coordinates are nonzero since their
sum vanishes. Choose a nonzero pivot \(Y_s\) of least cost among
the nonzero coordinates. For every other nonzero coordinate
replace it by \(Y_i-c_iY_s\), using its value divided by the
pivot value at the point. The cost of the pivot is at most its
cost, so the spaces are unchanged in either direction. Their
missing root divisors can be dropped first, while each coordinate
through the point with a retained boundary is untouched. On the
pivot chart these are affine parameters all vanishing at the
point. The graph equation there is the sum of these parameters
with zero constant term. Refine by the ordinary blowup at the
point. At secondary order \(l\), the old pivot exponent is fixed
to \(H-l\); all remaining monomials have ordinary degree \(l\)
and the same cost bound reduced by \(e_s(H-l)\). Distinct
monomials give all initials. Thus on the new exceptional ordinary
projective space the stated induction reproduces in one fewer
dimension. Aggregating all orders \(l\) has partitioned exactly
the same monomials, so the other exponent means are unchanged.

For its radial test write \(R\subset I\setminus\{s\}\) for the
indices carrying the branch through the chosen point, and \(J\)
for the other surviving indices. There is \(j\in J\) of cost
at least two. Indeed besides the pivot another nonzero coordinate
loses its boundary; if it had cost one, minimality would give two
original cost-one coordinates, whereas only \(Y_0\) has cost one.
At least one cost-\(d\) index still lies in \(R\): no earlier
point made such a coordinate nonzero, since that would terminate
the induction, and no earlier coordinate change affected it.
The radial ratio at this point is
\[
 \frac{1+\sum_{i\in R\cup J}1/e_i}{\sum_{i\in R}1/e_i+|J|}
 =1+\frac{1-\sum_{i\in J}(1-1/e_i)}{\sum_{i\in R}1/e_i+|J|}
 \le 1+\frac{1/e_j}{1+1/d}\le1+\frac{d}{2(d+1)} .
\]
There cannot be another nonterminating step once fewer than
three coordinates survive. The angular proof therefore terminates
after finitely many ordinary extractions. The estimates apply to
arbitrary bases using successive compatibility with the tested
and extraction orders; all errors are uniform since only the
original bounded root exponents have been discarded. Weighted
adjunction gives (Q9.2) at the origin. It holds throughout the
affine chart by its contracting weighted action: at each finite
degree an adapted basis has a limiting complete basis, the graph
is preserved, and openness of klt along the orbit transfers every
strictly smaller coefficient at the fixed origin to the test
point.

For a quartic ninefold a rank-one quadratic pencil whose cubic
on the common kernel of the coefficients of its root is reduced
and nonzero consequently has line threshold at least \(52/7\).
Take a squarefree binary section of that cubic. Any squarefree
binary cubic is linearly equivalent over the pole to \(u^3+v^3\),
by the three-transitivity of the projective linear group on the
line. With the extra quadratic light, weights \((3,2,2)\) give
the diagonal face \(x^2+u^3+v^3\), and the pole calculation
gives \(10/(1/4+(3/4)(19/13))=52/7\). The radial light ratio
is \(7/6\), giving the same bound.

\subsection{Quartic eightfolds with a rank-one quadratic along a line}
\label{sec:nov28-q8-mixed-cubic}
For the quartic eightfold the exceptional codimension can be increased
to at least \(87\). The additional estimate uses the cubic on the
common kernel of the coefficients of the quadratic root.

\paragraph{The signed graph of a square and a mixed cubic.}
Let \(W=(z+x^2+u^2v=0)\) in the affine chart of \(\mathbf P^4\).
Then at every affine point one has
\[
            5S_{\mathcal O(1)}(E)+E(W)\le\tfrac32 A(E).
\tag{Q8.1}
\]
Off the graph one uses the ordinary complete-series bound.
At points with \((u,v)\ne(0,0)\) the nonzero Hessian of \(u^2v\)
and the disjoint square give Hessian rank at least two. Translation
in \(x\), with an affine-linear correction to \(z\), preserves the
graph. It remains to treat the origin. Write \(D\) for the cubic
projective graph closure and \(Q=x^2+u^2v\). We claim that its
local threshold in hyperplane units is at least \(12/5\).

The cubic hypersurface restriction sequence and division by the
leading monomial \(u^2v\) give the following basis on its affine
chart in section degree \(m\):
\[
 x^a u^b v^c Q^q,\qquad a+b+c+q\le m,\qquad b<2\quad\text{or}\quad c=0.
\]
Use the two-torus for which \(x^2\) and \(u^2v\) have the same
character and put \(w=u^2v/x^2\). The character relation
\((-2,2,1)\) is primitive. Each character fixes its normal pair
\((b,c)\) and \(B=a+2q\). Its section space is
\[
  x^B u^b v^c(1+w)^{q_0}\mathbf C[w]_{\le l},\quad
  q_0=\max(0,B+b+c-m),\quad l=\lfloor B/2\rfloor-q_0.
\]
After normalization by \(m\), the mean vector of \((a,b,c,q,l/2)\)
converges to \((1/4,1/12,1/6,1/4,1/8)\).
Indeed the \(c=0\) layer has one third of the mass, the two
\(b=0,1\) layers have two thirds, and each free exponent in the
three-simplex has mean \(1/4\). For the last entry the character
average of \(l/2\) equals that of \(\lfloor B/2\rfloor-q=\lfloor a/2\rfloor\).

Set \((\alpha,\beta,\gamma)=(E(x),E(u),E(v))\) and
\(\omega=2\beta+\gamma-2\alpha\). When \(\omega>0\), the maximal
mean at \(w=0\) of the complete polynomial block is \(l\omega/2\).
The limiting total is at most
\[
   \tfrac12\alpha+\tfrac13\beta+\tfrac7{24}\gamma
   \ \le\ \tfrac5{12}(\alpha+\beta+\gamma),
\]
using \(2\alpha\le2\beta+\gamma\). When \(\omega<0\), using
the block's ordinary degrees at infinity gives
\[
   \tfrac14\alpha+\tfrac7{12}\beta+\tfrac5{12}\gamma
   \ \le\ \tfrac5{12}(\alpha+\beta+\gamma),
\]
since \(\alpha\ge\beta\). At equality \(\omega=0\) the same linear
bound covers the monomial factors. If \(w\) has constant residue
\(t\ne0\), the extra mean is at most \(E(w-t)/4\): the coefficient
is the mean \(q\) at \(t=-1\), and the mean \(l/2\) elsewhere.
Replacing \(d\log v\) by \(d\log(w-t)\) in the coordinate
logarithmic volume form proves
\(A(E)\ge\alpha+\beta+\gamma+E(w-t)\). For nonconstant residue
there is no extra mean. These estimates have uniform
\(o(1)A(E)\) errors by the simplex sums. As in the square-plus-quintic
calculation, a basis adapted to the tested order and a separating
cocharacter has initial basis spanning the character blocks.
Lower semicontinuity at the fixed origin transfers the resulting
local threshold \(12/5\) from the initial basis divisors.
Refinement by the cubic graph gives fixed coefficient \(1/3\)
and mean residual degree four in units of five. Divisor adjunction
thus gives \(5S+(5/3-1/3)E(W)\le(5/3)A\). Interpolate with
\(5S\le A\), with weights \(3/4,1/4\), to prove (Q8.1).

In particular an ordinary affine graph with rank-one Hessian at
a point satisfies the signed denominator \(3/2\) if the cubic
term on the Hessian kernel is nonzero and is not a cube of a
linear form up to scale. Indeed cut by linear flags to one
quadratic direction \(x\) and two general kernel directions.
After subtracting the affine tangent, the quadratic is a square.
The resulting binary cubic on the kernel has a simple projective
root: a nonpure cubic has a multiplicity-one component and a
general line meets that component transversely. Put this root at
\(v=0\) in coordinates \((u,v)\). Give these coordinates weights
\(1,1+\epsilon\), \(x\) weight \((3+\epsilon)/2\) and the graph
coordinate weight \(3+\epsilon\), with small positive rational
\(\epsilon\). The initial light polynomial is \(a x^2+b u^2v\)
with \(ab\ne0\); terms of ordinary order at least four have larger
weight. The sum of the three light weights exceeds the graph
weight. The radial signed denominator is therefore at most
\(3/2\), and (Q8.1) bounds the punctured graph. Weighted and
ordinary flag adjunction give the claim.

\paragraph{The two cubic coefficients along the line.}
First suppose the quadratic light form is nonzero of generic rank
one and its primitive linear root has base degree \(e\le1\).
For \(e=1\), take one constant linear cut through the line making
the root constant up to a base scalar while preserving its square,
as in the quintic-tenfold pencil proof. Concretely annihilate one
general combination of its two independent coefficient forms.
Keep a quadratic direction \(x\) and the common kernel \(U\) of
all coefficients of the primitive root. Its dimension is \(4\)
for \(e=1\) and \(5\) for \(e=0\). The cubic on \(U\) is
\(G_t=t_0G_0+t_1G_1\). If the geometric generic member is
nonzero and nonpure, retain two general constant directions in
\(U\) on which its restriction is nonpure. Give \(x\) weight
three, these remaining lights weight two, and the two heavy
coordinates weight six. At every general pole the light part is
\(a x^2+G(u,v)\), where \(a\ne0\) and the binary cubic is
nonpure. At every \((u,v)\ne(0,0)\) the Hessian of \(G\)
is nonzero: a vector in the radical of a binary cubic tensor
would make it a cube in the complementary variable. Thus the
disjoint square gives rank at least two. At \((u,v)=(0,0)\)
the preceding ordinary rank-one calculation applies after
translation in \(x\). The initial light weight sum divided by six
is \(7/6\). The two-heavy angular proportions \(1/4,3/4\)
and the positive signed bound \(3/2\) thus give denominator
at most \(11/8<2250/1601\), including the radial test.

Consequently an unresolved nonzero quadratic forces the entire
cubic pencil on \(U\) to lie in the cone of cubes. As its
projectivization (the third Veronese variety) contains no line,
the ordered pair \((G_0,G_1)\) has the form
\((aL^3,bL^3)\), including the zero case. The two normal-degree-three
coefficient forms restrict freely and independently to this fixed
\(U\) once the quadratic and polars are fixed. For \(e=1\) this
costs at least \(2\binom63-(3+2)=35\) conditions, in addition
to the codimension \(52\) already obtained after projection.
For \(e=0\) the corresponding figures are
\(2\binom73-(4+2)=64\) and \(56\), respectively.

Finally suppose the restricted quadratic is identically zero.
Its incidence after projecting the line has codimension at least
\(5+12+63-16=64\). On all six lights its cubic pencil has
arbitrary independent cubic coefficients. Let \(V\) denote this
six-dimensional light space. Except on a closed subset of
codimension at least \(16\) in \(\operatorname{Sym}^3V^*\), every
nonzero point of \(V\) has Hessian rank at least two, or rank
one with a nonpure nonzero cubic restriction to its Hessian
kernel. Indeed fix a point \([p]\in\mathbf P(V)\). The Hessian
evaluation is surjective onto the symmetric forms, so rank zero
costs \(21\) conditions. Rank exactly one costs \(15\). After
fixing the rank-one value with kernel \(K\), the free remaining
cubic coefficients are precisely the cubics on \(V/\langle p\rangle\).
Choose four directions in \(K\) projecting independently there
(possible whether or not \(p\in K\)). Their cubic restriction varies
freely in \(20\) dimensions and the cone of pure cubes has dimension
four. Hence rank one with pure or zero restriction costs at least
\(15+16\). Varying \([p]\) costs five dimensions. Taking closures
of the two resulting constructible sets gives the claimed bound.
For any cubic outside this closed exceptional set retain the six
lights: the radial ratio is two and each punctured graph has
signed denominator at most \(3/2\) by the two local graph tests.
Thus failure of the required line mark forces both cubic coefficients
into this closed subset, at a cost of at least \(32\) more conditions.
The three quadratic strata therefore have projected exceptional
codimension at least \(87,120,96\), respectively. The finite union
of their closures contains all smooth quartic eightfolds for which
uniform K-stability could fail.

\paragraph{A square plus two cubic variables.}
On the affine chart of \(\mathbf P^4\), let
\(W=(Z+x^2+g(y,u)=0)\), where \(g\) is a binary cubic which is
not the cube of a linear form times a scalar. At every point of \(W\)
and for every divisorial order \(v\) centered there,
\[
                 5S_{\mathcal O(1)}(v)+v(W)\le \tfrac32 A(v).
\]
First take \(g=y^2u\), \(Q=x^2+y^2u\), and let \(D\) be the cubic
projective graph closure. Its affine chart is smooth. We claim
\(\delta_0(D,\mathcal O_D(1))\ge12/5\). Its complete restrictions in
degree \(m\) have the basis \(x^a y^b u^c Q^q\) with
\(a+b+c+q\le m\) and \(b<2\) or \(c<1\), by monic reduction of
\(y^2u=Q-x^2\). In the limiting normalized basis measure,
\[
 (\bar a,\bar b,\bar c,\bar q)=(1/4,1/12,1/6,1/4).
\]
Indeed the single \(c=0\) layer has mass one third with \(a,b,q\)
free, and the two \(b=0,1\) layers have mass two thirds with
\(a,c,q\) free. Their overlap is of lower order.

Use the two-dimensional torus preserving the character identity
\(x^2=y^2u\), and write \(w=y^2u/x^2\). Each character has a unique
normal pair \((b,c)\): all candidates differ by multiples of
\((2,1)\). Writing \(P=a+2q\), its piece is
\[
 x^P y^b u^c(1+w)^{q_0}\mathbf C[w]_{\le l},\qquad
 q_0=\max(0,P+b+c-m),\quad l=\lfloor P/2\rfloor-q_0 .
\]
The mean of \(l/2\) after division by \(m\) is \(1/8\), since
the values of \(q\) range over a consecutive interval and
\(\lfloor P/2\rfloor-q=\lfloor a/2\rfloor\).
Write \(\alpha=v(x),\beta=v(y),\gamma=v(u)\),
\(\omega=2\beta+\gamma-2\alpha\). For \(\omega>0\) the maximum
complete-polynomial mean is \(l\omega/2\). Consequently the
normalized maximal mean in the character-split system is at most
\[
 \tfrac34\alpha+\tfrac1{12}\beta+\tfrac16\gamma+\tfrac18\omega
       =\tfrac12\alpha+\tfrac13\beta+\tfrac7{24}\gamma
       \le\tfrac5{12}(\alpha+\beta+\gamma).
\]
The inequality uses \(\alpha\le\beta+\gamma/2\).
For \(\omega<0\) the order in a polynomial is determined by its
degree, so one obtains
\[
 \tfrac14\alpha+\tfrac7{12}\beta+\tfrac5{12}\gamma
       \le\tfrac5{12}(\alpha+\beta+\gamma)
\]
from \(\alpha\ge\beta+\gamma/2\). For \(\omega=0\) either expression
gives the monomial part. If \(w\) has nonzero constant residue \(t\),
there is an additional mean at most \(\lambda/4\), with
\(\lambda=v(w-t)>0\): for \(t=-1\) use the mean \(\bar q=1/4\);
otherwise use \(\overline{l/2}=1/8\). In that event
\(A(v)\ge\alpha+\beta+\gamma+\lambda\) by replacing \(d\log u\)
by \(d\log(w-t)\) in the coordinate logarithmic volume form.
With a nonconstant residue there is no additional order.
The ordinary log discrepancy bound is always
\(A(v)\ge\alpha+\beta+\gamma\).
The errors are \(o(1)A(v)\) uniformly by the finite simplex sums
and bounded integral remainders. An arbitrary valuation is handled
by splitting the basis by a generic torus cocharacter and its
order simultaneously and applying lower semicontinuity to the
limiting basis divisor, as in the square-plus-quintic calculation.
This proves the claimed local threshold.

The mean power of \(D\) in the complete ambient series, multiplied
by five, is \(1/3\); the restricted moving degree multiplied by
five is four, by the hypersurface sequence. Signed divisor
adjunction with \(g_0=4/(12/5)=5/3\) yields
\(5S+(4/3)v(W)\le(5/3)A\). Interpolation with \(5S\le A\)
gives the asserted bound at the origin. For a general nonpure
binary cubic, put a simple root on the axis \(u=0\). Its first
binary term in increasing \(u\)-degree is \(c y^2u\), \(c\ne0\).
Choose weights \(1,1+\epsilon,(3+\epsilon)/2,3+\epsilon\)
on \(y,u,x,Z\). The initial graph is
\(Z+x^2+cy^2u\); a projective diagonal torus acts on the whole
ambient series, and simultaneous basis adaptation and lower
semicontinuity at the fixed origin transfer the signed estimate.
On the plane \(y=u=0\) translation of \(x\), followed by an
affine-linear shear in \(Z\), gives the origin case. At all other
points the binary Hessian is nonzero: contraction with a nonzero
point could vanish only if the cubic were a pure cube. Together
with the \(x^2\) block the Hessian has rank at least two, and the
quadratic flag gives the same bound. The corresponding ambient
higher-dimensional estimate follows by a linear flag whenever
a square and a nonpure binary cubic are retained.

\paragraph{A square, a cube, and a fourth power.}
On the affine chart of \(\mathbf P^4\), let \(W=(Z+x^2+y^3+u^4=0)\).
At the origin we claim
\[
  5S_{\mathcal O(1)}(v)+v(W)\le \tfrac{37}{25} A(v).
\]
At every other point of \(W\) the uniform bound \(3A(v)/2\) holds.
Write \(D\) for the projective graph closure, \(Q=x^2+y^3+u^4\).
Its affine chart is smooth and its homogeneous degree-\(m\)
restrictions have basis \(x^a y^b u^cQ^k\) for \(c<4\) and
\(a+b+c+k\le m\), by division by the monic \(u^4=Q-x^2-y^3\).
The extraction over the origin on \(D\) has weights \(6,4,3\)
and discrepancy 13. The normalized means in this basis of
\(a,b,k,T=m-a-b-c-k\) are all \(1/4\), with bounded \(c\);
the mean extraction weight is consequently
\((6+4+12)/4=11/2\). The radial quotient is \(26/11\).

We establish an angular threshold at least \(8/3\) on its
exceptional weighted projective space. At fixed extraction weight
\(J\), the residues \(s=b\bmod 3\) and \(h=c\bmod 2\) are fixed:
reduce \(J=6a+4b+3c+12k\) modulo three and two. Factor \(y^s u^h\);
its bounded coefficients give a uniform \(O(1/m)A\) error after
normalization. Use the three pure-power variables
\[
   A=x,\qquad B=u^2,\qquad C=y^3,\qquad Q=A^2+B^2+C.
\]
Their weights after dividing by six are \(1,1,2\).
The natural finite map of smooth weighted projective stacks from
\(\mathbf P(6,4,3)\) to \(\mathbf P(1,1,2)\) has ramification
along the root coordinates \(B,C\), with respective indices two
and three. Equivalently on the smooth stack charts this is the
coordinate power map in \(u,y\). Hence testing on the target with
boundary \(B/2+2C/3\) proves the angular bound upstairs by the
finite crepant formula. Each entire \(J\)-piece after removing the
two fixed residues is the pullback of one space on this target;
no splitting by a residual finite-group character is used.
Hereafter the letters \(B,C\) in a divisor expression denote
their coordinate divisors.

On the target, put \(m'=m-s-h\). The basis at fixed weighted
degree consists of monomials
\[
     B^e A^a C^j Q^k,\qquad e=0,1,\quad
     a+2e+3j+k\le m',
\]
whose target weighted degree is \(a+e+2j+2k\).
Taken over the original pieces, the normalized means of
\(a,j,k,T=m'-a-2e-3j-k\) are \(1/4,1/12,1/4,1/4\); the
four quantities \(a,3j,k,T\) are asymptotically uniform on the
three-simplex after dividing by \(m\). The bounded residues do
not affect these statements.

First work off \(B\), where only the \(2C/3\) boundary is needed.
Degenerate \(C\) by diagonal scaling towards zero, preserving
its coordinate divisor. This replaces \(Q\) by \(Q_0=A^2+B^2\).
The displayed monomials remain independent: for fixed \(j,e\)
different \(k\)'s have distinct highest \(B\)-powers \(2k+e\)
in their binary factors. These therefore give exactly the initial
spaces. Split by the \(C\) exponent also in the limiting system.
For an original test valuation one uses bases adapted to its
order and the scaling filtration, then lower semicontinuity
applied to each weighted-degree piece together in the basis
divisor. This is legitimate on the whole target with boundary
\(2C/3\); its angular bound below holds also at the vertex.

At fixed \(j\) put \(D'=a+e+2k\) for the residual binary degree.
When nonempty, the binary space is contained in
\(Q_0^\kappa\mathbf C[A,B]_{D'-2\kappa}\), with
\(\kappa=(D'+3j-m')_+\), and contains
\(Q_0^{\kappa+1}\mathbf C[A,B]_{D'-2\kappa-2}\).
Indeed its monomials in the basis free over \(\mathbf C[A,Q_0]\)
are those with \(e=0,1\) and \(k\ge D'+3j+e-m'\).
Thus the difference from the containing binary space has bounded
dimension. The mean of \(\kappa\) is \(1/8\), since
\(\kappa=(k-T-e)_+\) and the mean of the positive difference
of two uniform simplex coordinates \(k/m,T/m\) is \(1/8\).
The mean of \(D'\) is \(3/4\). At either of the two distinct roots
of \(A^2+B^2\), the maximal normalized mean order is at most
\(\mathbb E(D')/2=3/8\), and at any other point it is at most
\(\mathbb E(D'-2\kappa)/2=1/4\).
The coefficient of \(C\) in the limiting basis divisor is
\(\bar j=1/12\). The binary hyperplanes meet \(C\) transversely.
Simultaneously resolving their common vertex on the smooth
stack chart gives exceptional discrepancy two and normalized
coefficient \(\bar D'=3/4\); the vertex misses \(C\).
It follows that coefficient \(8/3\) works with \(2C/3\).
In the bounded binary-codimension estimates the aggregate
vanishing-sum error is \(O(m^3)\) among \(O(m^3)\) original
sections of degree \(m\), and so vanishes after normalization
by \(m h^0(mH)\).

At centers on \(B\) use a second basis before degeneration.
Replace \(C^j\) by \((B^2)^j\) while keeping \(A,Q,B^e\).
This leaves every piece unchanged: the substitution
\(C=Q-A^2-B^2\) and its inverse preserve weighted degree and
the cost upper bound when both \(C\) and \(B^2\) have cost three,
\(A\) has cost one and \(Q\) has cost one. Only one bounded
factor \(B^e\) occurs. Refine by the divisor \(B\). At \(B=0\)
the residual polynomials \(A^a(A^2+C)^k\) at fixed
\(a+2k=D''\) have consecutive \(k_0\le k\le\lfloor D''/2\rfloor\);
they are linearly independent, so \(2j+e\) are the exact fixed
orders along \(B\). Their mean is \(1/6\).
On the stack \(B\simeq\mathbf P(1,2)\) the root at \(A=0\)
has normalized coefficient \(\bar a=1/4\) in its uniformizing
parameter \(A\). At \(A^2+C=0\) the coefficient is
\(\bar k=1/4\). At any other point the coefficient is at most
\(\overline{(\lfloor D''/2\rfloor-k_0)/2}=1/8\), because
\(\lfloor D''/2\rfloor-k=\lfloor a/2\rfloor\), and the consecutive
interval basis becomes a complete series in the local coordinate
\((A^2+C)/A^2\). In particular the point \(C=0\) has coefficient
at most \(1/8\). Now
\[
   \tfrac12+\tfrac83\cdot\tfrac16\le1,\qquad
   \tfrac23+\tfrac83\cdot\tfrac18=1 .
\]
All other residual coefficients are at most \(1/4\) with zero
remaining boundary. Divisor adjunction along \(B\) proves the
claim on this chart. All the finite-dimensional simplex and
binary interval estimates give uniform \(o(1)A\) errors.
Thus the extraction gives
\(\delta_0(D,\mathcal O_D(1))\ge\min(26/11,8/3)=26/11\).

The complete ambient series has normalized fixed \(D\)-coefficient
\(1/4\) and residual degree four in units of five. Signed
adjunction with \(g=4/(26/11)=22/13\) yields
\(5S+(g-1/4)v(W)\le gA\). Interpolating with \(5S\le A\)
gives \(1+(g-1)/(g-1/4)=37/25\), as claimed.
When \((y,u)\ne0\) the Hessian of \(x^2+y^3+u^4\) has rank at
least two, giving the bound \(3/2\). At \(y=u=0\) a translation
of \(x\) and an affine-linear shear of \(Z\) reduce to the origin.

\paragraph{The quartic face with a mixed quadratic coefficient.}
The preceding square--cube--quartic estimate extends on the whole
affine chart to the graph of
\[
     x^2+y^3+b x u^2+c u^4,\qquad c\ne0.
\]
At the origin its signed denominator is still at most \(37/25\),
and on the whole chart at most \(3/2\). In the proof use the same
variables \(A=x,B=u^2,C=y^3,Q=P_2(A,B)+C\), now with
\(P_2=A^2+bAB+cB^2\). All weight-piece and simplex formulas, the
\(C\)-adic binary quotient and the \(B\)-adic computation remain
valid: the leading \(B^2\)-coefficient of \(P_2\) is nonzero and
the mixed term \(AB\) costs three in the two triangular changes
between \(C\) and \(B^2\). On \(B=0\) the polynomial is again
\(Q=A^2+C\). For completeness, away from \(B\cup C\) discard the
root boundary locally: the whole angular system is monomial in
the independent weighted coordinates \(A,B,Q\), whose mean
exponents are \(1/4,1/6,1/4\). A separating diagonal degeneration
in these coordinates gives the angular bound \(8/3\) there.
At centers on \(C\) away from \(B\), refine by its actual
order before specializing: the monomials on the remaining binary
line stay independent by the nonzero leading \(B^2\)-term just
as above. Inversion along \(C\) uses its coefficient
\(2/3+(8/3)(1/12)\le1\), and gives the previous residual bound
except at a double root of \(P_2\).

At this remaining point set \(B=1\), \(R=A/B-r\), \(Z=Q/B^2\),
where \(r\ne0\) is the double root. The branch is
\((2/3)(Z-R^2=0)\). At fixed angular degree \(J'\) and fixed
\(Q\)-power \(k\), the independent monomials in \(A,B,Q\)
have \(B\)-exponents \(0\le h\le h_{\max}\), where
\[
 h_{\max}=\min(J'-2k,\;2(m'-J'+k)).
\]
Indeed if \(h=2j+e\), \(e=0,1\), its section cost is \(3j+2e\)
and its weighted degree is \(h\). The binary block therefore
has an adapted basis of orders \(i=0,\ldots,h_{\max}\) at the
nonzero slope \(R=0\), times a common unit. Give \(R,Z\) weights
one and two. The initials are distinct monomials \(R^i Z^k\),
and the limiting mean of \(i\) equals that of \(h\), namely
\(1/6\), while that of \(k\) is \(1/4\). The radial basis order
is \(2/3\), against log discrepancy \(3-2(2/3)=5/3\).
On the new exceptional line \(\mathbf P(1,2)\) the two coordinate
orders have means \(1/6,1/4\). At fixed \(J',i+2k\) the possible
\(k\)'s form an interval by the displayed formula. Thus at any
finite nonzero slope the residual order is at most half its
width; in the aggregate this is at most \(\overline{i}/2=1/12\),
since \(\lfloor i/2\rfloor\) bounds the distance to the upper
endpoint. The strict transform of \(Z=R^2\) contributes \(2/3\)
at one such slope. Hence coefficient \(5/2\) works at the double
root, including all angular centers. This exceeds the primary
radial threshold \(26/11\), so the original graph-closure and
signed estimates at the origin follow.

At any other point with \(y\ne0\) its Hessian has rank at least
two by the independent \(xx,yy\) entries. At \(y=0\) the \(xx\)
entry is nonzero. If its Hessian has rank one, subtract the affine
tangent and complete the quadratic square by an affine-linear
change of \(x\), retaining the local variables \(u,y\). On the
Hessian kernel the cubic is \(y^3+\gamma u^3\). When \(\gamma\ne0\)
the square-plus-mixed-cubic estimate applies. Otherwise the
remaining local polynomial is precisely the displayed mixed
model at the origin: translation and this linear change preserve
the pure fourth coefficient \(c\), the coefficient of \(xu^2\),
and the disjoint \(y^3\). This proves the punctured estimate.

\paragraph{The coupled face \(x^2+y^3+xu^3\).}
For \(W=(Z+x^2+y^3+xu^3=0)\subset\mathbf A^4\subset\mathbf P^4\)
the signed denominator is at most \(3/2\) at every point.
The graph closure \(D\) has degree four. In its complete series
the normal monomial basis \(x^a y^b u^c Q^k\), \(Q=x^2+y^3+xu^3\),
has \(a=0\) or \(c<3\) and ordinary cost \(a+b+c+k\le m\).
Take weights \(3,2,1\), so \(Q\) has weight six. Standard simplex
counting gives mean weighted order \(21/8\), and the radial
threshold is \(6/(21/8)=16/7\).

Use the angular finite-symmetry reduction of the preceding paragraph
for the group scaling \(u,y\) by independent cube roots of unity.
In each character the bounded \(u,y\) residues can be removed.
The quotient coordinates \(A=x,B=u^3,C=y^3\) have weights
\(1,1,2\), and the branch boundary is \(2B/3+2C/3\). This follows
from the coordinate power map; the hyperplane generator pulls
back to three times the weighted generator. Write
\(Q=A^2+AB+C\). After bounded residue shifts the complete spaces
are spanned by
\[
 A^i B^j C^h Q^k,\qquad ij=0,\qquad i+3j+3h+k\le m' .
\]
In the \(j=0\) sector the simplex has weight \(3/4\), and in the
\(i=0\) sector weight \(1/4\). The normalized means
\((\bar i,\bar j,\bar h,\bar k)\) are
\((3/16,1/48,1/12,1/4)\). To justify this after character splitting,
equivalently start with the unrestricted normal basis above. Its three
\(c<3\) layers have \(a,b,k\) free, its \(a=0\) layer has \(b,c,k\)
free, and the bounded residue classes partition each layer uniformly.

At centers on \(C\) away from \(B\), refine by \(C\). The residue
at its exact order \(h\) consists on the binary line of the
restrictions of \(A^i B^j(A(A+B))^k\).
Writing \(D'=i+j+2k\), \(m''=m'-3h\), and \(\gamma=B/A\)
yields \(\gamma^j(1+\gamma)^k\). For \(j=0\), \(k\) runs over
\((D'-m'')_+\le k\le\lfloor D'/2\rfloor\). For \(i=0\), the
remaining \(k\)'s have \(j=D'-2k\) and
\(k\ge\lceil(3D'-m'')/5\rceil\); this is at least
\(D'-m''\) since \(D'\le2m''\). Their degrees \(D'-k\)
extend the first sequence by consecutive integers. Thus the
space is the complete binary interval after a fixed power
\((1+\gamma)^{q_0}\), \(q_0=(D'-m'')_+\).
Writing \(T=m''-i-3j-k\) gives \(q_0=(k-2j-T)_+\).
Its normalized means in the two sectors are \(1/8,3/40\):
in the second use
\(\mathbb E[(k-T-(2/3)(3j))_+]=1/[4(1+1)(1+2/3)]\)
for the normalized simplex. The total is \(9/80\).
The maximum mean at any finite \(\gamma\ne-1\) is consequently
\(\overline{j+k}-\bar q_0=19/120\). At \(-1\) it is
\(\overline{j+k}=13/48\). At \(A=0\) it is
\(\overline{i+k}=7/16\). These bounds prove the \(16/7\) estimate
on this part of \(C\), including the fixed coefficient
\(2/3+(16/7)(1/12)<1\).

To handle \(B\), replace \(C^h\) in the basis by \((AB)^h\).
The change \(C\leftrightarrow Q-A^2-AB\) preserves costs in
both directions when \(AB\) is assigned cost three. Thus in
the independent weighted coordinates \(A,B,Q\) the pieces
have the monomial basis \(A^{a'}B^{b'}Q^k\) with cost
\(a'+3b'-\min(a',b')+k\). The means of \(a'=i+h,b'=j+h\)
are \(13/48,5/48\). Refining by \(B\) on the stack, the
residual space at fixed \(D''=a'+2k\) and fixed \(b'\)
has consecutive \(k\)'s: the cost is a convex piecewise-linear
function of \(k\) after \(a'=D''-2k\) is substituted.
Use the affine coordinate \(q=Q/A^2\) on \(B\).
At \(B\cap C\) one has \(q=1\); the residual mean is at most
the average half-width of this interval, itself at most
\(\bar a'/2=13/96\). Indeed its upper endpoint is at most
\(\lfloor D''/2\rfloor\), and the distance to it averages half
the width. At \(A=0\) use \(A\) as uniformizing coordinate and
the bound \(\bar a'=13/48\), at \(Q=0\) use \(\bar k=1/4\).
At remaining nonzero finite values use the same \(13/96\).
Now \(2/3+(16/7)(5/48)<1\) and \(13/96<7/48\).
Divisor adjunction yields \(16/7\) all along \(B\).
Off both branch divisors the monomial coordinates \(A,B,Q\)
give the estimate directly by torus specialization with no boundary;
each of the three coordinate coefficients is below \(7/16\).
These finite interval and simplex estimates have uniform
\(o(1)A\) errors. The finite-invariance argument and the radial
check prove \(\delta_0(D,H)\ge16/7\). As for the other quartic
graphs, signed adjunction with \(g=7/4\) gives \(3/2\).

At a nonzero point of Hessian rank one, necessarily \(y=0\)
and \(u(4x-3u^3)=0\). If \(u=0,x\ne0\), the cubic term on
the Hessian kernel is \(dy^3+x\,du^3\). Otherwise
\(x=3u^3/4\), \(u\ne0\); the kernel has
\(dx=-3u^2du/2\) and its cubic term is
\(dy^3-(15/4)u^3du^3\). These are both nonpure binary cubics.
The square-plus-nonpure-binary-cubic flag gives \(3/2\),
and the rank-two quadratic test handles all other points.

\paragraph{The coupled fourth-order face \(x^2+y^3+yu^3\).}
For \(W=(Z+x^2+y^3+yu^3=0)\) in the affine chart of \(\mathbf P^4\),
every point satisfies \(5S_{\mathcal O(1)}(v)+v(W)\le(3/2)A(v)\).
It suffices at the origin to give the projective graph closure \(D\)
the bound \(\delta_0(D,H)\ge16/7\). The graph has degree four, so
signed divisor adjunction with \(g=4/(16/7)=7/4\) gives
\(5S+(g-1/4)v(W)\le gA\); interpolation with \(5S\le A\)
gives the required \(3/2\). Away from the plane \(y=u=0\) the
Hessian of the binary part is nonzero (its \(yy,yu\) entries are
\(6y,3u^2\)), hence the full Hessian has rank at least two.
Translation of \(x\) and an affine-linear shear of \(Z\) handle
the other points.

Put \(Q=x^2+y^3+yu^3\). In affine coordinates the complete
graded pieces of \(D\) have the normal monomial basis
\(x^a y^b u^c Q^k\), \(a+b+c+k\le m\), \(b=0\) or \(c<3\).
Division by the monic relation for \(yu^3\) does not increase
ordinary degree and proves the assertion. Take weights
\((9,6,4)\) on \(x,y,u\), giving weight 18 to \(Q\). The
discrepancy is nineteen and the mean is \(65/8\):
the limiting means of \(a,b,c,k\) are respectively
\(1/4,3/16,1/16,1/4\). Indeed the \(c=0,1,2\) layers have mass
\(3/4\) with \(a,b,k\) free, and the \(b=0\) layer has mass
\(1/4\) with \(a,c,k\) free. Thus the radial threshold is
\(152/65>16/7\).

At a fixed weighted order \(J\), the residues \(a\bmod2,c\bmod3\)
are fixed by reducing modulo two and three. Factor these bounded
powers and use \(B=y,C=u^3,A=x^2\), whose weights after division
by six are \(1,2,3\), with \(Q=A+B^3+BC\). The crepant boundary
on the target stack \(\mathbf P(1,2,3)\) of this coordinate root
map is \(A/2+2C/3\). After bounded residue and cost shifts, the
angular basis is
\[
      A^i B^b C^j Q^k,\qquad
      b=0\text{ or } j=0,\qquad 2i+b+3j+k\le m'.
\]
In the sector \(b=0\) we here include only monomials outside the
other sector; its boundary overlap has lower dimension. Means of
\((i,b,j,k)\) are \((1/8,3/16,1/48,1/4)\). As before removing
the bounded residues costs \(o(1)A\) after normalization.

For centers outside \(A\), omit its boundary and specialize the
pencil to \(Q_0=A+BC\), by giving \(B,C,A\) auxiliary weights
\(1,0,1\). The \(B^3\) term has larger weight. The normal forms
remain independent by monic reduction of \(BC=Q_0-A\), so this
is the complete initial system in each degree. The boundary \(C\)
is preserved, and the following check works on the entire target
with \(2C/3\) alone, allowing lower semicontinuity at any moving
test center. After lifting to the smooth root charts, dropping the
boundary \(A/2\) gives the contribution \(+\alpha\) to the allowed
discrepancy, where \(\alpha=v(x)\). In such a chart one coordinate
is a unit. Write \(\beta=v(y),\gamma=v(u)\). The discrepancy with
this extra allowance is at least \(2\alpha+\beta+\gamma\), by
coordinate normal crossings. Apply the binomial torus calculation
to \(Q_0=x^2+yu^3\), using the invariant \(w=yu^3/x^2\) and
\(\omega=\beta+3\gamma-2\alpha\). We give the details to include
arbitrary, possibly nonconstant test orders on each fiber.

The normal pair of exponents \(b,c\) is unique in each character,
differing from any other candidate by a multiple of \((1,3)\).
The common remaining character exponent is \(P=a+2k\).
In degree \(m\) the character space is a monomial factor times
\((1+w)^{k_0}\mathbf C[w]_{\le l}\), where
\(k_0=(P+b+c-m)_+\), \(l=\lfloor P/2\rfloor-k_0\), with empty
pieces omitted. Its extraction grade is constant on a character
since \(yu^3\) and \(x^2\) have the same weight 18.
The normalized mean of \(l/2\) is \(1/8\), using
\(\lfloor P/2\rfloor-k=\lfloor a/2\rfloor\).
For \(\omega>0\) the maximal normalized mean is bounded by
\[
   \tfrac34\alpha+\tfrac3{16}\beta+\tfrac1{16}\gamma
              +\tfrac18\omega
      =\tfrac12\alpha+\tfrac5{16}\beta+\tfrac7{16}\gamma
      \le \tfrac7{16}(2\alpha+\beta+\gamma).
\]
For \(\omega<0\) the pole orders are determined by the degrees and
the bound is
\[
   \tfrac14\alpha+\tfrac7{16}\beta+\tfrac{13}{16}\gamma
      \le \tfrac7{16}(2\alpha+\beta+\gamma),
\]
since \(2\alpha\ge\beta+3\gamma\). At \(\omega=0\) the common
monomial part has the same bound. A nonzero finite constant residue
\(t\) of \(w\) adds at most \(v(w-t)/4\) to the mean and at least
\(v(w-t)\) to the discrepancy, as in the square-plus-binary-cubic
calculation. In a uniformizing chart one may fix one nonzero
coordinate and replace a logarithmic differential of either
remaining coordinate by \(d\log(w-t)\); all coefficients in the
monomial defining \(w\) are nonzero. This gives the stated
discrepancy inequality, including on the angular stack.
Generic torus weights separate the finitely many characters at
each degree, preserve \(Q_0\) up to character, and preserve the
coordinate boundaries. A basis simultaneously adapted to the
test and the weights followed by lower semicontinuity proves the
same bound for arbitrary valuations, as in the previous binomial
calculation. The specialization with only boundary \(2C/3\) is
taken downstairs, where that boundary is effective.

For centers on \(A\), refine the original series by \(A=0\).
Its fixed normalized coefficient is \(1/8\). Indeed the residual
polynomials \(B^b C^j(B^3+BC)^k\), \(b=0\) or \(j=0\), found
with a given exponent \(i\), remain independent by their distinct
degrees in \(\gamma=C/B^2\), as verified just below. The residual
weighted line is \(\mathbf P(1,2)_{B,C}\). Put
\(J'=b+2j+3k\), \(m''=m'-2i\); divide the sections by \(B^{J'}\).
For \(j=0\) their polynomials are \((1+\gamma)^k\) for
\[
 \max(0,\lceil(J'-m'')/2\rceil)\le k\le\lfloor J'/3\rfloor .
\]
For \(b=0\) they are \(\gamma^j(1+\gamma)^k\), with
\(j=(J'-3k)/2\) an integer and \(3j+k\le m''\).
As \(k\) decreases by two in the second sector, the leading degree
\((J'-k)/2\) increases by one. It starts at the end or one above
the end of the first sector (discard any duplicated \(j=b=0\)
term). The smallest \(k\) in the second sector is at least that
in the first, since its extra lower condition is
\(k\ge(3J'-2m'')/7\) and \(J'\le3m''\). Thus the combined space
is \((1+\gamma)^{q_0}\mathbf C[\gamma]_{\le L-q_0}\), where
\(q_0=\max(0,\lceil(J'-m'')/2\rceil)\) and \(L\) is its greatest
degree. Bounded endpoint changes do not affect limiting means.
The normalized average leading degree is
\(\overline{j+k}=13/48\). Writing
\(T=m'-(2i+b+3j+k)\), one has
\(q_0=(k-j/2-T/2)_++O(1)\), with mean
\[
 \tfrac34\cdot\tfrac16+\tfrac14\cdot\tfrac17=\tfrac9{56}.
\]
For the first sector use the simplex mean of \((k-T/2)_+\).
For the second use that \(3j\) and \(T\), together with \(k,2i\),
are uniform simplex coordinates. In exponential coordinates
the mean of \((k-aT-b(3j))_+\), divided by \(m\), is
\(1/[4(1+a)(1+b)]\). This proves the two values.

At any finite point other than \(\gamma=-1\), including \(C=0\),
the maximal normalized residual order is therefore at most
\(13/48-9/56=37/336\). At \(\gamma=-1\) it is \(13/48\).
At the stack point \(B=0\), leading degrees give order
\(J'-2(j+k)\) in a uniformizing parameter and mean
\(\overline{b+k}=7/16\). Consequently the coefficient \(16/7\)
is allowed everywhere on the residual line: at \(C=0\) the
remaining discrepancy is \(1/3\), and
\((16/7)(37/336)<1/3\). Along \(A\) itself
\(\tfrac12+(16/7)(1/8)<1\).
Divisor adjunction proves the angular claim. All bounded shifts
and simplex remainder terms give \(o(1)A\) after normalization.
Together with the radial check this proves the desired origin
threshold and the signed denominator.

\paragraph{The face \(x^2+y^3+xu^2\).}
On the affine graph \(W=(Z+x^2+y^3+xu^2=0)\) in
\(\mathbf P^4\), the signed denominator is at most \(3/2\) at
every point. Write \(Q=x^2+y^3+xu^2\), and let \(D\) be the cubic
graph closure. At the origin take weights \(6,4,3\) on \(x,y,u\).
The complete graph series has the basis \(x^a y^b u^c Q^k\),
\(0\le b<3,\ a+b+c+k\le m\), by monic division by \(y^3\).
Its normalized means of \(a,c,k\) are \(1/4\) each, hence
the extraction gives discrepancy thirteen and mean \(21/4\).
Its radial quotient \(52/21\) exceeds \(12/5\).

At fixed weight, remove the fixed bounded residues in \(y\bmod3\)
and \(u\bmod2\) as in the disjoint fourth-power calculation.
On \(\mathbf P(1,1,2)\) use \(A=x,B=u^2,C=y^3\), with boundary
\(B/2+2C/3\) and \(Q=A^2+AB+C\). Since \(C\) itself is a coordinate
of weight two, \((A,B,Q)\) are also weighted projective coordinates.
The pieces are monomial in these coordinates with costs \(1,2,1\).
Off \(C\), retain the boundary \(B/2\) and use the diagonal torus
in the coordinates \(A,B,Q\). The normalized binary exponents of
\(A,B\) are \(1/4,1/8\), and the \(Q\)-exponent is \(1/4\).
A basis adapted to a test order and a general diagonal weight
degenerates to the monomial basis. The three coordinate stack
divisors have simple normal crossings. Thus the angular threshold
with this reduced boundary is at least four on the entire stack.

At centers on \(C\) refine by \(C=0\). Introduce \(C\) at cost three
and use \(AB=Q-A^2-C\) to eliminate the product \(AB\).
This preserves weighted degree and each ordinary cost bound, in
both directions. Thus an equivalent full basis consists of
\(C^h A^a B^j Q^k\), with \(a=0\) or \(j=0\) and
\(3h+a+2j+k\le m'\), where \(m'\) includes bounded residue shifts.
The fixed order \(h\) is exact because, as seen below, at a fixed
weighted degree the residual \(A^aB^j(A(A+B))^k\) are independent.
The limiting simplex masses of the \(j=0\) and \(a=0\) sectors
are \(2/3\) and \(1/3\), respectively, since their free variables
have cost products three and six. Hence
\[
   (\bar h,\bar a,\bar j,\bar k)=(1/12,1/6,1/24,1/4).
\]

On \(C=\mathbf P^1_{A,B}\) set \(D'=a+j+2k\) and
\(\gamma=B/A\). After removing \(A^{D'}\) the residual polynomials
are \(\gamma^j(1+\gamma)^k\). Put \(m''=m'-3h\).
The sector \(j=0\) gives the consecutive powers with
\(q_0=(D'-m'')_+\le k\le\lfloor D'/2\rfloor\).
The other sector satisfies \(j=D'-2k\) and
\(k\ge\max(0,\lceil(2D'-m'')/3\rceil)\ge q_0\);
here \(D'\le2m''\) whenever the space is nonempty.
As \(k\) decreases, the degrees \(D'-k\) extend the initial interval
by consecutive integers. This proves independence and identifies
the combined space as \((1+\gamma)^{q_0}\mathbf C[\gamma]_{\le L-q_0}\),
where \(L\) is the greatest degree. With
\(T=m''-a-2j-k\), one has \(q_0=(k-j-T)_+\).
In the two sectors its respective normalized means are \(1/8\)
and \(1/12\): in the latter use \(2j,k,T\) as three of the
simplex coordinates, and in the former \(j=0\). Thus
\(\bar q_0=1/9\). The maximum mean at finite \(\gamma\ne-1\)
is at most \(\overline{j+k}-\bar q_0=13/72\);
at \(-1\) it is \(\overline{j+k}=7/24\). At \(A=0\)
the mean is \(\overline{D'-j-k}=5/12\).
The \(B/2\) boundary meets \(C\) only at \(\gamma=0\),
where \((12/5)(13/72)<1/2\); the fixed coefficient on \(C\)
satisfies \(2/3+(12/5)(1/12)<1\). Divisor adjunction along \(C\)
gives angular threshold at least \(12/5\). The finite residue
errors and weighted-simplex remainders are \(o(1)A\) as above.
Together with the radial bound this proves \(\delta_0(D,H)\ge12/5\).

The cubic hypersurface adjunction argument with \(g=4/(12/5)=5/3\)
gives \(5S+(4/3)v(W)\le(5/3)A\). Interpolation gives \(3/2\).
Away from the origin the only points where the Hessian of the
light part has rank one satisfy \(y=0,\ x=u^2,\ u\ne0\).
Its kernel there is \(dx=-u\,du\), with \(dy\) free.
The restriction of the cubic term after translation to this
kernel is therefore \(dy^3-u\,du^3\), a nonpure binary cubic.
Use the square-plus-nonpure-binary-cubic estimate after the usual
weighted flag of weights three on a complementary direction
and two on the kernel. This flag has normalized light-weight
sum \(7/6>1\), so its radial signed quotient is less than \(3/2\).
Points of Hessian rank at least two are covered by the quadratic
flag, completing the proof.

\paragraph{All pencils with nonzero quadratic light form.}
Let a line on a smooth quartic eightfold have two-dimensional
polar span. If its first light form has order two, its line
threshold is at least \(1601/250\). Rank at least two and the
nonpure cubic restriction after a rank-one quadratic were treated
above. For a rank-one form retain its quadratic direction \(x\)
and the common kernel \(U\) of its fixed linear coefficients,
after one constant cut when the primitive root moves. The
dimension of \(U\) is five for a fixed root and four for a
moving root. The quadratic on this section is \(a(t)x^2\),
with \(a\) nonzero at a general pole. A preliminary cut in the
moving case can make the retained \(x\)-direction any general
combination of the two active constant directions.

Suppose first the cubic restriction to \(U\) is nonzero but pure
over the geometric generic base. Its pencil is \(b(t)y^3\) for
a fixed linear form \(y\) on \(U\), since the third Veronese
contains no line. Set \(K=(y=0)\subset U\). If its quartic
restriction is nonzero, keep one general direction \(u\) of
\(K\). The weights \((6,4,3)\) on \(x,y,u\) give face
\(a x^2+b y^3+c x u^2+e u^4\), \(abe\ne0\). The mixed
square--cube--quartic estimate gives positive signed denominator
\(3/2\); the light ratio is \(13/12\). A zero restriction on
\(K\) cannot occur for a fixed root: the join of the line with
\(\mathbf P(K)\) would be a five-space contained in the smooth
quartic eightfold. For a moving root that join is a contained
\(M=\mathbf P^4\). The five normal derivatives of its equation
on \(M\) are linearly independent. Otherwise at most four
cubics on \(M\) would cut out their common nonempty zero set,
contradicting smoothness of \(X\). In particular the derivative
in the \(y\)-direction is nonzero on \(M\). It consists solely
of a normal-degree-four term \(y R_3(K)\): all lower terms on
\(U\) were the pure cubic, or zero. If a degree-three term
\(x R_2(K)\) survives for some retained active direction, take
that direction and a general \(u\) where it is nonzero. The same
weights now give \(a x^2+b y^3+c x u^2\), \(abc\ne0\), with
the signed denominator \(3/2\) of the cubic mixed graph. Otherwise
use a general \(u\) where \(R_3\ne0\). Weights \((9,6,4)\)
give \(a x^2+b y^3+c y u^3\), \(abc\ne0\): the cubic term
\(x y u\) has weight nineteen and all quartic \(x\)-terms have
weight at least twenty-one. The signed bound is again \(3/2\),
and the light ratio is \(19/18\). This handles every pure
nonzero cubic restriction.

It remains that the entire cubic vanishes on \(U\). The terminal
quartic \(C\) there is nonzero and reduced: the join of the line
with \(\mathbf P(U)\) has ambient codimension \(c\le4\) and
dimension greater than four, and the terminal Gauss bound gives
projective singular codimension at least \(10-2c\ge2\).
Write the cubic part linear in \(x\) as \(xB_t(Y)\) on \(U\).
Keep two general directions \(Y=(u,v)\) in \(U\). Their binary
quartic \(C\) is squarefree, by generic transversality to the
reduced hypersurface. Weights \((2,1,1)\), with heavy weight
four, give light part
\[
             a x^2+xB(Y)+C(Y),\qquad a\ne0,
\]
at general poles. Its radial ratio is one. At a pole, let
\(D=C-B^2/(4a)\). When the full \(D\) is nonzero and nonpure,
the initial general constant binary restriction also makes it
nonpure at general poles. At \(Y=0,x\ne0\) the Hessian has
rank at least two if \(B\ne0\); for \(B=0\) use the square-plus-
binary-quartic estimate. For \(Y=p\ne0\), subtract the affine
tangent. Write \(s=x+B(Y)/(2a)\) and use as the new \(x\)
the affine-linear part of \(s-s(p)\) at that point, so the
square is completed only in the quadratic Taylor term. If the
Hessian has rank one, its kernel now projects isomorphically
onto the \(Y\)-plane, and its cubic restriction is exactly the
cubic Taylor coefficient of \(D\) at \(p\), a nonzero constant
multiple of \(\partial_pD\). Indeed its remaining terms in the
expansion of \(a s^2+D\) are a multiple of \(xB(Y-p)\) in
degree three and the original \(C(Y-p)\) in degree four. A
binary quartic has zero derivative in a nonzero direction only
when it is a fourth power up to scale, or zero. Thus this cubic
is nonzero. If nonpure use the ordinary cubic test. If pure,
say it is \(b y^3\), then a nonzero pure quartic on its kernel
or a surviving cubic \(x u^2\) is handled by the first two
mixed graphs above. In the remaining case the squarefree
\(C(y,u)\) vanishes to order exactly one in \(y\), so its
\(y u^3\) coefficient is nonzero, and the \((9,6,4)\) graph
applies. The three secondary extractions all have light ratio
at least one. Rank at least two was already sufficient.

If instead the binary \(D\) is zero, \(C=B^2/(4a)\) contradicts
squarefreeness. If \(D\) is a nonzero pure fourth power, over the
pole field the whole light polynomial is, after rescaling,
\((x+B'(Y))^2+L(Y)^4\). Squarefreeness of \(C=B'^2+L^4\)
implies \(B'|_{L=0}\ne0\). In binary coordinates \(y=L\) one
can therefore complete a square in the other light, preserving
\(y\), and rescale it to obtain exactly
\((x+u^2+b y^2)^2+y^4\). The completed-square graph estimate
gives signed denominator at most \(3/2\), uniformly in \(b\).
All the coordinate choices for these secondary ordinary graph
computations are made over the respective pole. The first
linear cuts and the first weight vector are constant on the
line. Estimates at the finitely many chosen general poles suffice.
The radial inequalities and the regular and positive masses
\(1/4,3/4\) now give denominator at most \(11/8\), proving the
required mark for every line of quadratic initial order.

\paragraph{The square of a smooth affine quadric graph.}
In the affine chart of \(\mathbf P^4\), the graph
\(W=(Z+(x+uv)^2=0)\) satisfies \(5S+v(W)\le(3/2)A(v)\).
Away from the graph the projective-space bound suffices. Where
\(x+uv\ne0\) on it the quadratic term of the centered affine
polynomial has rank at least two: after using the affine-linear
part of \(x+uv\) as coordinate, its quadratic part is a square
plus \(2(x+uv)\,du\,dv\).
It remains to test points of the smooth affine quadric
\(C=(Z=x+uv=0)\). Its closure in \(\mathbf P^4\) is the smooth
complete intersection \(Z=q=0\) of degrees \(1,2\), with
\(q=xT+uv\) and \(T\) the homogenizing coordinate.
Blow up \(C\) with normal weights \(2,1\) for \(Z,q\).
The associated normal monomials \(Z^iq^j\) have independent
coefficients in
\(H^0(C,\mathcal O_C(t,t))\), \(t=M-i-2j\), by the regular
sequence, the quadric restriction sequence and projective
normality. Here \(C=\mathbf P^1\times\mathbf P^1\).
The normalized means in units of \(5/M\) are
\(\bar i=1,\ \bar j=1/2,\ \bar t=3\). The weighted radial
discrepancy is three; its signed order is
\(2+1/2+2=9/2\).

On the exceptional weighted line bundle over the affine
quadric, take a ruling fiber through the point being tested,
then its marked point on that fiber. Multiplicities along these
ordinary smooth flags may be computed by bases of the two
binary factors of \(\mathcal O_C(t,t)\): their orders \(a,b\)
range independently from zero to \(t\). Each therefore has
normalized mean \(3/2\). The angular graph \(ZT^3+q^2=0\)
is dominant over both base flags and adds no fixed base order.
At these flags one can work in the regular frames where \(T\ne0\).
The degrees of the coefficient spaces and the indicated ranges
also compute any divisorial valuation centered generically over
a ruling or over the whole quadric by omitting the later flags.

At fixed base orders \(a,b\) all normal monomials of cost
\(i+2j\le M-\max(a,b)\) reassemble in their angular weights.
At fixed weight \(2i+j\), the exponents \(i\) form a consecutive
interval. At any noncoordinate point of the weighted line, the
mean vanishing of its complete interval is half its width.
This equals the average distance from \(i\) to its upper
endpoint and hence is at most the average of \(j/2\), namely
\(1/4\) after summing the blocks in the same normalization.
At the single graph mark one adds one to this coefficient.
At the coordinate marks use the generic stack parameters and
the original means \(1,1/2\) respectively. Thus all are at
most \(3/2\). Successive filtered inversion of adjunction on
the two ordinary base flags and the initial weighted extraction
proves the assertion.

The graph \(Z+x^2+xuv=0\) also has signed denominator at most
\(3/2\). Its Hessian has rank at least two away from \(x=u=v=0\).
At the origin specialize by a diagonal weight preserving the
origin to \(Z+xuv=0\), giving \(x\) weight greater than the
sum of the \(u,v\) weights and \(Z\) the weight of \(xuv\).
The cubic graph closure in the latter case is toric on its
smooth affine chart. In degree \(M\) reduce the monomials
\(x^a u^b v^c Q^j\), \(a+b+c+j\le M,\ Q=xuv\), using
\(\min(a,b,c)=0\). Each omitted-coordinate layer has one third
of the mass. Each of \(a,b,c\) therefore has normalized
mean \(1/6\), whereas \(j\) has mean \(1/4\). The resulting
three coordinate exponents have means \(5/12\). Diagonal
specialization to monomial bases gives local graph closure
threshold at least \(12/5\). Ordinary divisor adjunction and
interpolation with projective space give the signed denominator
\(3/2\), as in the other cubic graphs.

\paragraph{The square--cube--quintic graph.}
In the affine chart of \(\mathbf P^4\), the graph
\(W=(Z+x^2+y^3+u^5=0)\) satisfies
\[
             5S_{\mathcal O(1)}(v)+v(W)\le(3/2)A(v).
\]
We prove that its degree-five projective closure has
\(\delta_0(W,\mathcal O_W(1))\ge20/9\). On its smooth affine chart
write \(Q=x^2+y^3+u^5\). Degree-\(m\) restrictions have the
independent normal basis \(x^i y^j u^e Q^\ell\), with
\(i+j+e+\ell\le m,\ 0\le e<5\): reduction by the monic \(u^5\)
relation lowers ordinary degree. Restriction is complete by the
hypersurface sequence. The weights \((15,10,6)\) on \(x,y,u\)
give discrepancy 31. The three unbounded exponents have mean
\(m/4+O(1)\), so the radial mean divided by \(m\) is \(55/4\).
Its threshold \(124/55\) exceeds \(20/9\).

The coarse weighted exceptional space is the plane with
coordinates \(A=x^2,B=y^3,C=u^5\), and \(Q=A+B+C\).
Its different is \(\frac12(A=0)+\frac23(B=0)+\frac45(C=0)\)
by the weighted quotient charts. A weighted degree determines
the residues of \(i\) modulo two, \(j\) modulo three, and \(e\)
modulo five, since the three exponents are pairwise coprime.
After these bounded common factors are removed, the degree-\(J\)
parts have independent monomials \(A^aB^bQ^c\) with \(a+b+c=J\)
and costs \(2a+3b+c\le m-O(1)\). Their normalized mean
exponents are \(1/8,1/12,1/4\).
Away from \(C=0\) the two remaining branch lines are coordinate
in \(A,B,Q\). A diagonal specialization preserving these two
lines and the character spaces gives normalized coefficients at
most these monomial means. After multiplication by \(20/9\) and
addition of the boundary, all three coordinate coefficients are
at most one. Equivalently apply the coordinate strata calculation
to the plane with just these two branch lines, whose log comparison
holds globally.

Along \(C=0\), replacing \(B\) by \(C=Q-A-B\) is a triangular
change of the monomial filtration in each degree, since \(B\)
has the largest cost. Hence the mean fixed order along \(C\)
is \(1/12\) and \(\frac45+(20/9)(1/12)<1\).
After imposing \(b=\operatorname{ord}_C\) in degree \(J\), the
binary exponents of \(A\) form the entire interval
\[
       0\le a\le\min(J-b,m'-J-2b),
\]
where \(m'=m-O(1)\) includes the fixed residues. At a finite
nonzero slope in \([A:Q]\) the maximum mean binary vanishing is
half the interval width, whose normalized aggregate is \(1/8\).
At \(A=0\) this same mean applies, and at \(Q=0\) the normalized
mean is \(1/4\). In particular the restricted coefficients are
at most one at both further branch marks: at \(B=0\) the slope
is \(A=Q\), giving \(2/3+(20/9)/8<1\), and at \(A=0\) use
\(1/2+(20/9)/8<1\). Inversion along \(C\) and weighted
adjunction give the claimed local threshold. The normalized
finite-degree errors vanish by the simplex estimates in each
of the finitely many residue classes.

In the ambient complete series the fixed \(W\)-coefficient in
units of five is \(1/5\), and the moving class on \(W\) is four.
Divisor adjunction therefore gives
\(5S+(9/5-1/5)v(W)\le(9/5)A\) at the origin.
Interpolating with \(5S\le A\) gives the asserted signed
denominator. If \(y\ne0\) or \(u\ne0\), the affine Hessian has
rank at least two and the quadratic flag suffices. The remaining
\(x\)-axis reduces to the origin by translation and an affine-linear
change of \(Z\).

\subsection{High-order polar pencils on a quintic elevenfold}
\label{sec:nov28-quintic-eleven-high}
Let \(C\) be a line of polar rank two on a smooth \(X_5^{11}\).
The pencil criterion for \(\delta_C(X,H)\ge60/7\) asks for radial
light-weight sum at least the graph weight and signed positive
denominator at most \(3/2\). Their combination with the regular
and positive masses \(1/5,4/5\) gives denominator \(7/5\).

If the initial form of order \(h\) is log canonical on the
geometric generic projective light space and at least \(h\)
lights are available, the criterion holds. Take \(h\) general
unit lights, preserving log canonicity, and use the translated
section argument on \(\mathbf P(1^h,h)\) with coefficient one.
In detail, \(z_0=v(Z)\), \(w=v(W)\), \(A=A_P(v)\) satisfy
\(z_0+w\le A\) by inversion from the two sections. For
\(D=(h+1)S_{\mathcal O(1)}+(1-1/h)z_0\) the two cyclic-cover
bounds give
\[
    D+w\le A+\min\{w,w/h+(1-1/h)z_0\}\le 3A/2 .
\]
For the last inequality use the first branch when \(w\le z_0\)
and the second when \(w\ge z_0\), since \(h\ge2\).
At the weighted vertex the two section orders vanish.

In particular the terminal order \(h=5\) succeeds after any
constant light restrictions of total ambient codimension \(c\le4\).
The finite Gauss morphism gives singular codimension at least
\(13-2c\ge5\) in the projective light space. Five general lights
therefore give a smooth quintic on \(\mathbf P^4\).

\paragraph{Quartic squares.}
If the order-four form is the square of a quadric of rank at
least four, take four lights on which it is nondegenerate.
On \(P=\mathbf P(1^4,4)\) with \(Z=(z=0)\),
\(W=(z+q^2=0)\), let \(D=5S_{\mathcal O(1)}+3v(Z)/4\).
The flag \(Z,Q=(q=0)\) gives the signed inequality
\[
              D+v(Q)\le 3A/2-v(Z)/2 .
\]
Indeed in the normalization of \(D\) the fixed coefficient
of \(Z\) is one, the subsequent quadric coefficient is \(1/2\),
and the residual on the smooth quadric surface is three times
the complete hyperplane series. Its hyperplane stability threshold
is two. Moving \(v(Z)/2\) to the left and applying flag inversion
proves the inequality. At a center on only one member the
corresponding shorter flag suffices; at the vertex the two
graph orders vanish. The change \(z\mapsto z+q^2\) transfers the
inequality to \(W\). Therefore
\[
 D+w\le 3A/2-v(Q)+\min\{w-z_0/2,\,-w/4+3z_0/4\}\le 3A/2 ,
\]
because the minimum is at most \(\min(w,z_0)/2\) and
\(2v(Q)\ge\min(w,z_0)\). The radial ratio is one.
A square quadric of rank at most three in the original kernel
is constant by the base-degree bound. Its maximal isotropic
subspace costs at most two hyperplanes, so the terminal case applies.

Any linear factor of the order-four form in the original kernel
also yields the terminal estimate. A rational factor has base
degree at most one, hence its coefficients have a common kernel
of codimension at most two. For a nontrivial orbit the same
constant-kernel conclusion follows by the moving-hyperplane
incidence in the quintic-tenfold proof.

\paragraph{A GIT-semistable four-light section at order four.}
The same estimate holds whenever the geometric generic quartic
restricts on some four-dimensional light space to a GIT-semistable
quartic surface \(G\subset\mathbf P^3\). This condition is open
in the Grassmannian, so one may take a general constant four-space
and the same open condition holds at general pole roots. Ascher,
DeVleming and Liu, \emph{K-stability and birational models of moduli
of quartic K3 surfaces}, Theorems~1.1(1) and 5.14,
\href{https://arxiv.org/abs/2108.06848}{arXiv:2108.06848}, identify
the quartic GIT stack with the K-moduli stack for rational boundary
coefficients \(0<c<1/3\). The valuative inequalities for
\((\mathbf P^3,cG)\) tend at \(c=1/3\) to
\[
            4 S_{\mathbf P^3,H}(v)+v(G)/2\le3A(v)/2 .
\]
The \(Z\)-flag in \(\mathbf P(1^4,4)\) now gives
\(D+v(G)/2\le3A/2-z_0/2\). The translated flag and
\(v(G)\ge\min(w,z_0)\) give the same \(3/2\) bound as for
the square quadric.

\paragraph{A fixed cut for a family of cuspidal cubics.}
Suppose that at order three the geometric generic cubic is
irreducible and fails log canonicity. Such a cubic is a cone over
a cuspidal plane cubic by the essential-variable classification.
Its three-space of essential normals defines a subbundle \(E\)
of the constant dual light bundle on the binary base, with
splitting \(\mathcal O(-a)\oplus\mathcal O(-b)\oplus\mathcal O(-c)\),
\(0\le a\le b\le c\). The cubic is a section of
\(\operatorname{Sym}^3E(2)\) because its base degree is two.
If \(b>0\), every allowed monomial contains the first generator,
contradicting irreducibility. Hence \(a=b=0\). When \(c>0\) the
third generator cannot occur cubed, so the two constant first
generators have a common kernel on which the cubic vanishes.

When \(c=0\), write the family in three constant normal coordinates
as \(u^2 A+uvD+v^2 E_0\), with \(A,D,E_0\) ternary cubics. If
these have no common projective zero, they define a morphism
\(\phi\colon\mathbf P^2\to\mathbf P^2\) pulling back
\(\mathcal O(1)\) to \(\mathcal O(3)\). The unique singular point of
the general plane cubic gives a rational map from the binary base,
extending to \(p\colon\mathbf P^1\to\mathbf P^2\). Differentiating
the identity of vanishing at \(p(t)\) gives zero base derivative
there since the spatial gradient vanishes. Thus the binary quadratic
obtained by evaluation at \(p(t)\) has a double root \(t\), and is
nonzero by basepoint freeness. The composite \(\phi\circ p\) is
therefore the dual Veronese conic, whose pullback has degree two.
This contradicts \(3\deg p^*\mathcal O(1)=2\).
A common ternary zero gives a constant linear cut of codimension two
in the original lights. This argument for constant three-dimensional
essential space also applies to a reduced tangent line and conic,
which have a unique singular point. When the resulting light order
is five in the original kernel, the terminal test applies.

\paragraph{A binary cubic and a second disjoint binary form.}
Let \(f(X_1,X_2)\) be squarefree of degree three, and let
\(g(Y_1,Y_2)\) have degree four or five with all root multiplicities
at most two. The ordinary graph \(Z+f+g=0\) has signed denominator
at most \(3/2\) at every affine point. At nonzero vectors both
binary Hessians are nonzero by Euler's identities and the
multiplicity bounds. When both vectors are nonzero their disjoint
blocks give rank at least two. When \(X=0,Y\ne0\), retain a
direction of nonzero Hessian for \(g\) and both cubic directions.
After subtracting the affine tangent the first weighted face is
a square plus the full nonpure binary cubic.

At \(Y=0,X\ne0\) consider a rank-one Hessian of \(f\).
At a projective zero, Euler would put the radial vector in its
radical, contradicting smoothness of the binary zero. Thus the
radial direction is nonisotropic. In radial and gradient-kernel
coordinates \(f=aR^3+bS^3\), with \(ab\ne0\), since the mixed
terms vanish by the first and second derivatives. Retain the
variation of \(R\), the variable \(S\), and a direction where
\(g\ne0\). The first weighted graph, after removing the affine
tangent, is a square plus a cube plus a fourth or fifth power.
Its light ratio is at least \(1/2+1/3+1/5>1\), and the graph
calculations for the exponents \((2,3,4)\) and \((2,3,5)\)
give \(3/2\). At the origin use weights \(1/3\) on both
\(X\)-variables and \(1/\deg g\) on both \(Y\)-variables, with
graph weight one. Their sum exceeds one. The punctured estimates
complete weighted adjunction.

In particular suppose a geometrically reduced non-lc pencil
cubic has a constant essential normal space of dimension two or
three. Keep two general essential directions and let \(K\)
be the common kernel of the full essential space. The restricted
binary cubic is squarefree. If the next order on \(K\) is four
and the quartic has no component of multiplicity at least three,
two general kernel lights give the disjoint face just treated.
Suppose instead that the next order on \(K\) is five and both
first derivatives of the quartic in the retained normal directions
vanish on \(K\). The ambient codimension of \(K\) is at most five,
so two general kernel lights give a squarefree quintic and again
the disjoint test applies.

\paragraph{Mixed squares and binary quartics with double roots.}
Let \(C(Y)\) be a nonzero binary quartic, let \(B(Y)\) be binary
quadratic and \(a\ne0\). Assume \(B\) is nonzero at all projective
roots of \(C\) of multiplicity at least three.
The ordinary graph of \(a x^2+xB(Y)+C(Y)\) has signed denominator
at most \(3/2\) on the entire affine chart.

Use the proof of the squarefree case in the quartic-eightfold
calculation. We give the additional arguments for the double
roots. Let \(D=C-B^2/(4a)\). At \(Y=0,x\ne0\) the Hessian has
rank at least two when \(B\ne0\), and otherwise the
square-plus-binary-quartic lemma applies. At \(Y=p\ne0\) with
rank-one Hessian, translate the point and make the affine-linear
change given by the differential of \(x+B/(2a)\). The cubic on
the Hessian kernel is proportional to \(\partial_pD\).
If \(D\) is nonpure, this derivative is nonzero. A nonpure
derivative gives the square-plus-nonpure-binary-cubic estimate.
For a derivative \(b y^3\) the mixed square--cube tests apply
if \(C|_{y=0}\ne0\), \(B|_{y=0}\ne0\), or the coefficient of
\(y u^3\) in \(C\) is nonzero, where \(u\) completes \(y\).

Assume these three coefficients vanish. The binary quartic
\(D\) vanishes on \(y=0\). If \(y(p)\ne0\), integrating
\(\partial_pD=b y^3\) in a frame consisting of \(y\) and a
linear coordinate annihilating \(p\) gives \(D=b' y^4\), since
the remaining fourth power must vanish on \(y=0\). Hence
\(y(p)=0\) in the nonpure case, and
\[
       D=y^3(cu+e y),\qquad B=y(b_1u+b_2y),\qquad b_1\ne0 .
\]
The last assertion follows because otherwise \(C=D+B^2/(4a)\)
would be divisible by \(y^3\) at a root of \(B\). The Hessian of \(D\) at \(p\)
vanishes. In coordinates adapted to the differential of
\(s=x+B/(2a)\), rank one therefore forces \(s(p)=0\), since
the remaining Hessian is \(s(p)\operatorname{Hess}(B)\) and
\(\operatorname{Hess}(B)\ne0\). The centered polynomial is
\[
         a(X+\alpha y U+\beta y^2)^2+y^3 V(y,U),
                 \qquad \alpha\ne0,
\]
with \(U=u-u(p)\). Assign weights \(1,1,0,2\) to \(X,y,U,Z\)
respectively. This is a projective diagonal degeneration fixing
the affine origin and preserving the ordinary complete series.
Its initial graph is \(Z+a(X+\alpha y U)^2=0\), covered by
the squared smooth quadric calculation. Semicontinuity transfers
the signed bound at the fixed origin.

If \(D=0\), then \(B\) has rank two by the multiplicity
hypothesis on \(C\), and the squared smooth quadric calculation
applies directly. If \(D\) is a nonzero fourth power, write
the light equation up to scalars as \((x+B')^2+y^4\).
When \(B'|_{y=0}\ne0\) use the completed-quartic estimate of the
quartic-eightfold proof. Otherwise \(B'=y(bu+ey)\) with \(b\ne0\)
by the same multiplicity hypothesis. At \(y\ne0\) its quadratic
part in the centered normal frame \(s=x+B'\) is
\[
       (ds)^2+s\,\operatorname{Hess}(B')
                      + 6y^2(dy)^2
\]
up to the harmless choices of quadratic-form convention.
Its rank is at least two, since the mixed \(u,y\) coefficient
of \(\operatorname{Hess}(B')\) is nonzero, so the last two
summands cannot cancel. At \(y=0\), rank one forces \(s=0\),
and the same degeneration with weights \(1,1,0,2\) applies.
These arguments cover every affine point; the origin itself
is treated by the primary weights \(2,1,1,4\) and the punctured
comparisons just proved, since their radial ratio is one.

Consequently the rank-one order-two quintic pencil test succeeds
when the common kernel of the quadratic root has first order four
and every geometric component of that quartic has multiplicity
at most two. Indeed two general kernel lights have precisely that
binary property at general pole roots. The usual weights \(2,1,1\)
on the light variables give the displayed face, even when the
initial quadratic root required one preliminary constant cut.

\paragraph{A square plus a smooth ternary quintic.}
If \(P(Y)\) defines a smooth quintic in \(\mathbf P^2\), the
ordinary affine graph of \(x^2+P(Y)\) has signed denominator
at most \(3/2\) throughout. For \(Y\ne0\) the Hessian of \(P\)
is nonzero by smoothness and the differentiated Euler identity.
Together with the square this gives rank at least two.
For \(Y=0\) translate \(x\) with an affine-linear shear of the
graph coordinate. At the origin use weights \(5,2,2,2,10\)
on \(x,Y,Z\). The light ratio is \(11/10>1\). The preceding
ordinary estimates cover punctured directions with \(Y\ne0\)
by reassembly of the ordinary complete bases in all angular
weights. The remaining graph direction has \(x\ne0\).

As in the square-plus-binary-quartic calculation, take the
graph flag on the exceptional stack. At level \(M\) its
\(k\)-th power costs \(5k\) ordinary degrees and \(10k\)
weight. The resulting fixed coefficient is \(1/5\) and the
residual complete graph degree is five, in the ambient units
of six. Thus the fixed coefficient including the graph mark is
\((2/3)(1+1/5)<1\), and the moving threshold required there
is \(10/3\).

Let \(R_l\) be a complement to \(P\mathbf C[Y]_{l-5}\) in
\(\mathbf C[Y]_l\), of dimension \(5l+O(1)\). Ordinary
division by the degree-five equation gives spanning independent
blocks \(x^a(x^2+P)^q R_l,\ a+q+l\le m\) of the graph series.
Their limiting normalized means are \(\bar a=\bar q=1/5\),
\(\bar l=2/5\) by the simplex count weighted by this linear
dimension function. On the chart \(x=1\) at fixed weight
\(J=5a+2l+10q\), the allowed \(l\) satisfy \(2l\equiv J\bmod5\).
The complete range in \(q\) is
\[
 \max(0,(J+3l)/5-m)\le q\le\lfloor(J-2l)/10\rfloor .
\]
Eliminating these consecutive powers of the unit \(1+P\) gives
independent initial blocks \(P^sR_l\) for every \(s\) from zero
to the interval width. At fixed order \(D=l+5s\), the largest
endpoint becomes \(s+\lfloor(J-2D)/10\rfloor\), and the
smallest is \(\max(0,(J+3D)/5-m-3s)\). Hence the allowed
\(s\) form the full interval from a lower endpoint \(s_0\) to
\(\lfloor D/5\rfloor\). Their direct sum is
\(P^{s_0}\mathbf C[Y]_{D-5s_0}\). The mean of \(s\) equals
the mean of \(q_{\max}-q=\lfloor a/2\rfloor\), namely \(1/10\).
Consequently the secondary ordinary blowup has mean order
\(\bar D=2/5+5/10=9/10\) against discrepancy three.
On its exceptional \(\mathbf P^2\), the complete moving sections
give at most \((D-5s_0)A/3\). Smoothness of \(P\) gives
\(v(P)\le A\le 5A/3\), so the total on this plane is bounded
by \(DA/3\). After taking means it is at most \(3A/10\), as
required. All normalized counts have uniform \(o(1)\) errors
as in the binary case. This proves the assertion.

\paragraph{The remaining order-two pencils.}
If the order-four form on the quadratic kernel has a component
of multiplicity at least three, this is a unique hyperplane \(L\)
and its primitive equation is constant by the available binary
base degree one. Keep the quadratic direction \(x\) as above
and let \(x B\) be its normal-degree-three term linear in \(x\).
When \(B|_L\ne0\), two general kernel directions meet the
triple or quadruple root at a nonzero value of \(B\). The
preceding mixed-square test applies. Otherwise impose the single
constant cut \(L\) on the kernel while retaining the direction
\(x\). The new kernel has ambient codimension \(c\le5\);
its dimension is at least seven, so it is not contained in \(X\).
Its order is five, and its quintic has projective singular
codimension at least \(13-2c\ge3\). Three general kernel
directions give a smooth ternary quintic \(P\). Since \(B|_L=0\),
weights \(5,2,2,2\) on \(x,Y_1,Y_2,Y_3\) and heavy weight ten
give exactly the face \(a x^2+P(Y)\), with \(a\ne0\) at general
poles. Its radial and angular comparisons have just been proved.
Together with the nonpure cubic test and the mixed
square--cube and \(E_8\) graph tests when the restriction to
the quadratic kernel starts at order three, and the terminal
test at \(c\le4\) when that order is five, this proves the
mark \(60/7\) for every order-two quintic elevenfold pencil.

\paragraph{Two non-log-canonical quartic strata.}
Continue with an order-four pencil on a quintic elevenfold in
the original light kernel. The four-light Nadel classification
applies with coefficient one: use \(c'<1\) sufficiently close to
one, so \(4-4c'>0,\ 1-c'>0\) on the light and binary factors.
All horizontal components of the non-lc locus are constant. The
following two strata always give the \(60/7\) mark.

First suppose there is a nondegenerate codimension-two non-lc
component of degree three. It is a variety of minimal degree,
with a general three-space section a smooth twisted cubic \(T\)
in \(\mathbf P^3\). It is a codimension-two component of the
singular locus, so the restricted quartic vanishes doubly along
\(T\). Every such quartic is \(J(Q_0,Q_1,Q_2)\), where \(J\)
is a quadratic and \(Q_i\) are the three standard minors cutting
out \(T\). Indeed the blowup of \(T\) maps by \(|2H-E|\) to
\(\mathbf P^2\). Its general fiber is the strict transform of
a secant line: the minors at its general point record the binary
quadratic with roots its two points on the rational normal
cubic, equivalently the kernel of its two-by-three Hankel
matrix. This morphism has connected fibers by Stein factorization.
Thus \(H^0(4H-2E)\) is the pullback of the plane quadratics.
At \(p=[s^3:s^2t:st^2:t^3]\), the three linear normal terms
of the minors have image the line
\(t^2Q_0-stQ_1+s^2Q_2=0\), by the two Hankel syzygies and
smoothness of \(T\). These lines trace a nonsingular dual conic.

Since the quartic is non-lc at a general smooth point of \(T\),
the restriction of \(J\) to the corresponding line has rank
at most one. Rank two would give a transversal node and
log canonicity there. A rank-two conic has this rank drop only
along lines through its singular point and the two components,
whose dual locus does not contain the displayed dual conic.
A rank-one \(J\) gives a double quadric. In that case the
full generic quartic is itself a double quadric by general
linear section and Bertini, and the square-quadric test applies.
If \(J\) has rank three it is the unique conic with this dual.
Thus the quartic is the discriminant of binary cubics, preserved
by the natural \(\mathrm{SL}_2\) acting irreducibly on the
four coordinates. It is \(\mathrm{SL}_4\)-semistable: a destabilizing
one-parameter subgroup would have an instability parabolic
containing its \(\mathrm{SL}_2\) stabilizer by Kempf's theorem,
impossible by irreducibility. The quartic GIT comparison applies.

Next suppose there is no codimension-two non-lc component and the
unique codimension-three component \(K\) has generic multiplicity
two. It is linear by the four-light Nadel classification. A
general projective three-space cuts it in the unique non-lc
point \(p\) on the resulting quartic surface. To justify uniqueness,
possible further components of higher codimension are avoided,
and general linear section commutes with the multiplier ideal.
If this quartic is unstable, a flag giving one of the three
maximal negative Hilbert--Mumford supports has flag point \(p\).
Indeed the four weights in each support are strictly decreasing;
in an affine chart at its flag point the positive weights
\(w_0-w_i\), \(1\le i\le3\), have sum \(4w_0\), and each allowed
monomial has order greater than \(4w_0\). This point is non-lc.
The support \((8,-1,-3,-4)\) has multiplicity at least three
and is excluded. For the supports \((4,3,-2,-5)\) and
\((6,2,-1,-7)\) the quadratic at \(p\) is a multiple of
\(x_3^2\). On its tangent plane \(x_3=0\), every allowed monomial
is divisible by \(x_2^3\). In fact if \(a,b\) are the exponents
of \(x_0,x_1\) there, the two inequalities are respectively
\(6a+5b<8\) and \(7a+3b<4\).

Use three fixed normal coordinates \(X\) to \(K\) and coordinate
vectors \(Y\) along it. The quadratic normal term \(J(X,Y)\)
has bidegree \((2,2)\) and binary base degree one. Over a
general point of \(K\) it has rank one; rank two would give
a log canonical transverse germ. Over the polynomial unique
factorization domain of the base and longitudinal variables,
its primitive square root is linear in \(X\); writing
\(J=qL^2\) and comparing base degrees shows that \(L\) is
constant in the binary base. Its longitudinal degree is zero
or one. In the latter case it is bilinear in \(X,Y\), and \(q\)
has longitudinal degree zero. If this bilinear form has rank
at least two, restrict to two normal and two longitudinal
directions on which its rank is two, then scale the normal
directions to see the square of a rank-four quadric as the
limit of that quartic section. This proves GIT semistability.
Otherwise \(J=q'(Y)L_0(X)^2\) with \(L_0\) a constant linear
normal. The case where \(L\) already had longitudinal degree
zero has the same description.

The tangent plane to the doubly tangent quartic at \(p\) in
the three-space is now its intersection with the fixed hyperplane
\(H=(L_0=0)\). As the three-space varies generally, this plane
varies generally in \(H\). The Hilbert--Mumford support calculation
therefore gives the zero form or a triple line in a general
plane section of the full \(G|_H\). In the latter case \(G|_H\)
has a triple hyperplane component, by Bertini: a general plane
cuts each geometrically integral component of degree at least
two in an integral plane curve, and different components have
distinct general plane sections. The triple hyperplane is unique
and constant because the binary base degree is one. Thus at
most two constant restrictions reach the terminal order in
ambient codimension at most four.

\paragraph{A rank-three quadric component in a hyperplane.}
Let a codimension-two non-lc component be the fixed integral
quadric \(B=(H=q=0)\) of rank three on \(H\). Choose three fixed
coordinates \(X\) on which \(q=q(X)\) is nondegenerate. Since the
quartic is in the symbolic square of this complete-intersection
ideal, it is in its ordinary square. This also follows by first
restricting to \(H\), where \(q^2\) divides the restriction, then
dividing the remainder by \(H\) and restricting again.
Write its base-linear family as
\[
          G_t=a_t q^2+H q M_t+H^2R_t .
\]
If \(a_t=0\) identically, the fixed hyperplane suffices. At each
general point of \(B\), the transverse binary quadratic has rank
one for a general \(t\), since rank two gives a transversal node.
A linear space of binary squares has dimension at most one, so
\(M_t|_B/a_t\) is independent of \(t\). The component spans \(H\),
and \(M_t\) is linear. Thus a fixed linear \(M\) gives
\(M_t\equiv 2a_t M\bmod H\). Completing the square leaves
\(G_t=a_t Q^2+H^2R'_t\) where \(Q=q+HM\) is constant and
\(R'_t|_B=0\). Hence for linear coefficients in \(t\),
\[
          G_t=a_t Q^2+k_t H^2 Q+H^3 L'_t .
\]
If \(M|_H\) uses a coordinate from the radical of \(q|_H\), we
may take \(Q=q(X)+H u\) in fixed coordinates. Restrict to
\(H,u\) and two \(X\)-coordinates on which \(q\) is nondegenerate.
Give \(H,u\) opposite weights \(1,-1\), fixing the others.
The quartic tends to the square of the resulting rank-four
quadric, so this four-space is GIT-semistable. Otherwise a
constant shear of \(X\) gives
\[
   G_t=a_t q(X)^2+k_t H^2 q(X)+H^3 L_t+c_t H^4,                \tag{*}
\]
where \(L_t\) involves no \(H\). If \(L_t\) uses a direction outside
\(X\), there is a four-space over the geometric generic field
on which \(L_t\) becomes a chosen multiple of \(H\) with
nonzero total fourth-power coefficient. The resulting form
\(a q^2+kH^2q+cH^4\), \(ac\ne0\), is semistable: its stabilizer
contains \(\mathrm{SO}(q)\) acting irreducibly on the three
\(X\)-coordinates, and its centralizer scales this space and
\(H\). The weights on the two nonzero end terms have opposite
signs. The centralizer criterion for a reductive stabilizing
subgroup follows also from the existence of a centralizing
optimal one-parameter subgroup in Kempf's instability theorem.
By openness a general constant four-space then suffices.

We may now take all four essential coordinates \(X,H\) constant.
If \(L_t=0\) identically, the preceding semistability applies
when \(c_t\ne0\), and \(q\) is a fixed factor otherwise.
Suppose the general \(L_t\ne0\) and every four-light section
is unstable. At \(H=1\) the polynomial is \(a_t(q+\lambda_t)^2
+L_t+\text{constant}\). At \(q+\lambda_t=0\) its gradient is
nonzero, and elsewhere its Hessian has rank at least two since
it is a nondegenerate ternary form plus a rank-one form.
Thus any non-lc flag point is on \(H=q=0\), where the quadratic
tangent is a double plane. In a hyperbolic frame at such a point
write \(q=xy+z^2\), with point on the \(x\)-axis and tangent
plane \(y=0\). From the two applicable degree-four negative
Hilbert--Mumford supports, its tangent plane section must contain
a triple line. If \(L_t\) has a nonzero \(x\)-coefficient,
its coefficient in (*) on that plane is a multiple of \(H^3\).
Any triple line must then be independent of \(x\) and divide
\(H\), inconsistent with the \(z^4\) term. Hence \(L_t(p)=0\).
Write its restriction to the tangent plane as \(\beta_t z\).
The quartic in this plane is binary:
\[
       a_t z^4+k_t z^2H^2+\beta_t zH^3+c_t H^4.
\]
A triple root and the vanishing cubic coefficient give
\(k_t^2+12a_t c_t=0\) and \(27a_t\beta_t^2+8k_t^3=0\),
by expansion of \(a_t(z-uH)^3(z+3uH)\).
In the chosen hyperbolic frame \(\beta_t^2\) is the dual
quadratic norm of \(L_t\), since \(L_t(p)=0\).
These invariant identities hold identically in the binary base.
As \(a_t,k_t,c_t\) are linear forms and \(a_t\ne0\), we have
\(k_t,c_t\) proportional to \(a_t\). The second identity gives
\(\|L_t\|_{q^{-1}}^2\) proportional to \(a_t^2\).

Write \(L_t=a_t L_0+b_t L_1\) in independent base coordinates
\(a_t,b_t\). If \(L_1\ne0\), the norm identity makes
\(p_*=q^{-1}L_1\) isotropic and \(L_0(p_*)=0\). Its fixed
tangent plane is \(L_1=0\); on it (*) becomes \(a_t\) times
a fixed binary quartic as in the preceding hyperbolic frame.
A linear factor of that binary quartic gives a common line for
the entire pencil. For \(L_1=0\) the family is already a fixed
quartic up to scalar. Choose any point of
\((q=0)\cap(L_0=0)\) in the \(X\)-plane to get the same binary
tangent-plane restriction and a line. Lifting this fixed line
from four coordinates gives a constant codimension-two cut and
the terminal computation.

\paragraph{The transverse cubic and quartic along a linear center.}
Return to the unique linear non-lc center \(K\) of codimension
three and suppose its multiplicity is at least three.
Let \(X\) denote its three normal coordinates, \(Y\) its
longitudinal coordinates. If the multiplicity is four, the form
depends on \(X\) alone, and the two ternary quartic coefficients
have a common point by the projective dimension theorem. Its lift
gives a constant cut of codimension two. In multiplicity three
the form is \(C_{t,Y}(X)+P_{4,t}(X)\), with \(C\) a ternary cubic
linear separately in the binary base \(t\) and in \(Y\).
The generic \(C_{t,Y}\) is plane-GIT-unstable if all general
quartic three-space sections are unstable. Indeed take a
three-space with this single longitudinal direction and the three
normals and let its quartic degenerate to \(T C_{t,Y}(X)\).
If the cubic were semistable this product would be semistable:
its preserving block torus acts scalarly on the ternary space
and on \(T\) with determinant one, and the centralizer test
reduces precisely to ternary semistability.

A nonreduced cubic has a unique repeated line, constant along
every general parameter line in the binary base or the
longitudinal projective space by the bound one on the degree
of the cubic coefficients along that parameter line. Thus this
line is constant simultaneously in both parameters. For reduced
unstable ternary cubics the possibilities are a cusp, a tangent
line and conic, or three distinct concurrent lines. On any general
line in either parameter factor their unique singular point is
constant. For a pencil of generically integral cuspidal cubics,
Bertini puts the general singular point in the finite pencil
base. A constant pencil gives the same conclusion directly. For
a line times a smooth conic the base-degree sum is at most one,
so one factor is fixed. On a fixed line, the rank-one binary
quadratic restrictions of a conic pencil have fixed double
root since the Veronese conic contains no line. The smooth
dual conic similarly permits no pencil of lines tangent to a
fixed conic. For the three concurrent lines, their essential
normal subbundle on the parameter line splits as
\(\mathcal O(-a)\oplus\mathcal O(-b)\), \(0\le a\le b\).
A squarefree cubic in its third symmetric power twisted by
\(\mathcal O(1)\) forces \(a=b=0\): for \(b>0\) each
allowed monomial is divisible twice by the first generator.
Hence that pencil's common point is fixed as well.
These linewise conclusions give constancy on the parameter product
by varying general points successively in its two factors.

In the multiplicity-two plane-cubic cases the unique tangent
line at this common point is also fixed, since its quadratic
term is a family of squares linear along each parameter line.
Restrict to this tangent line in the normal projective plane.
The resulting binary cubic is a multiple of \(y^3\), where
\(y=0\) is its marked point. For the concurrent-lines case use
any constant line through the point. For nonreduced cubics use
the fixed repeated line already found, where the restriction is
zero. After this one constant cut the full quartic has the form
\[
             y^3 L_t(Y)+P_t(y,z).
\]
If \(P_t(0,z)=0\), the additional cut \(y=0\) gives terminal
order in codimension four. Otherwise its \(z^4\) coefficient
is nonzero at general poles. On \(K\) the pure quintic is
base independent and has projective singular codimension at
least \(13-2\cdot5=3\). Three general \(Y\)-directions give
a smooth ternary quintic \(Q\).

When \(L_t\ne0\), weights \(1/4,4/15,1/5\) respectively on
\(z,y,Y\), with graph weight one, give face
\[
             z^4+y^3L(Y)+Q(Y)
\]
after rescaling coefficients at the pole. Their sum is \(67/60>1\).
For this ordinary graph, wherever \(z\ne0\) its Hessian has rank
at least two unless \(y=Y=0\): the \(y,Y\) mixed derivatives are
nonzero at \(y\ne0\), and the Hessian of \(Q\) is nonzero at
\(Y\ne0\), by smoothness and homogeneity. At the exceptional
axis fix \(y=0\) and use the radial \(z\)-square with the
smooth ternary quintic.

At \(z=y=0,Y=p\ne0\), a rank-one Hessian of \(Q\) occurs only
off its zero set, as Euler and smoothness show. If \(L(p)\ne0\),
retain its radial square, the cubic \(y^3\), and \(z^4\) using
the \((2,3,4)\) graph test. If \(L(p)=0\), choose
\(e\in\ker dQ_p\) with \(L(e)\ne0\): the radial vector pairs
nontrivially with \(dQ_p\) and trivially with \(L\). On the affine
linear section \(Y=(1+R)p+b y e\), a nonzero kernel cubic from
\(Q\) gives the same square--cube--quartic test. If it vanishes,
the pure \(y^4\) correction is \(bL(e)+O(b^4)\) and is nonzero
for a general \(b\). The face is a square plus a squarefree
binary quartic in \(z,y\). Here the mixed quadratic and cubic
terms \(R y,Ry^2\) from \(Q\) vanish by its rank-one Hessian and
homogeneity.

If \(y\ne0\) and \(L(Y)=0\), the \(yy\) entry vanishes and the
mixed \(y,Y\) row is nonzero. At \(yL(Y)\ne0\), use just the
radial \(Y\)-direction with variable \(u=L(Y)\). On the
\(y,u\) plane the light part is \(y^3u+a u^5\) up to scalars.
If its Hessian has rank one, its determinant gives
\(a u^4=3y^3/40\). Along its kernel, normalized by \(du=1\),
we have \(dy=-y/(2u)\); the cubic Taylor coefficient is
\(11y^3/(8u^2)\ne0\). The square, this cube, and the disjoint
fourth power of \(z\) give the required signed estimate.
At the origin the displayed weights and all punctured estimates
apply. This completes the case \(L_t\ne0\).

For \(L_t=0\) choose weights \(1/4\) on \(y,z\), \(1/5\) on the
three \(Y\)'s. A root of \(P_t\) of multiplicity at least three
is fixed by its base degree one and gives the additional cut
to terminal order. Otherwise each of the two blocks of
\(P_t(y,z)+Q(Y)\) has a nonzero Hessian at any nonzero argument.
Where both arguments are nonzero their disjoint Hessians give
rank two. On either axis retain a square in the nonzero block
and use the square--binary-quartic or square--ternary-quintic
test on the other one. The light-weight sum \(1/2+3/5>1\)
covers their origin. Together with the other strata this proves
the \(60/7\) bound for every initial pencil quartic in the
original light kernel of a quintic elevenfold.

\paragraph{A cube and a smooth plane quintic.}
On the ordinary affine chart of \(\mathbf P^5\), the graph of
\(F=x^3+P(Y)\), with \(P\) a homogeneous smooth ternary quintic,
has signed denominator at most \(44/29\). Give \(x\) weight five,
the three \(Y\) coordinates weight three and the graph coordinate
weight fifteen. The radial discrepancy is \(29\) and the graph has
order fifteen, giving \(1+15/29=44/29\). We bound each punctured
angular graph point by \(3/2\).

At a point with \(xY\ne0\), the disjoint Hessian blocks both have
positive rank by Euler's identity and smoothness. At \(x\ne0,Y=0\)
the translation of \(x\) and subtraction of the affine tangent have
first face a square plus \(P(Y)\), using weight five on the translated
coordinate and two on each \(Y\). That graph has denominator at most
\(3/2\) by the preceding square calculation, and the secondary radial
denominator \(1+10/21\) also suffices.

At \(x=0,Y\ne0\) only the case of a rank-one Hessian requires
attention. Then \(P(Y)\ne0\): at a smooth zero its radial vector is
isotropic and pairs nontrivially with a gradient direction. Take a
radial coordinate \(R\) and coordinates \(e\) along the gradient
kernel. By Euler and the rank-one condition the first nonzero
homogeneous derivative on the kernel has degree \(s=3,4\), or five.
Indeed the homogeneous quintic in these coordinates has the form
\[
 aR^5+\sum_{j=3}^5 R^{5-j}P_j(e)
\]
with \(a\ne0\), and at least one \(P_j\) is nonzero by smoothness.
Retain a general kernel line and let \(s\) be the smallest such
\(j\). At the given point subtract the affine tangent;
the weights proportional to \(1/2,1/3,1/s\) on the translated
radial coordinate, \(x\) and this line give a square plus a cube
plus an \(s\)-th power. Their light weights sum to at least one,
so the secondary radial bound is at most \(3/2\).
The signed bound for this graph is at most \(3/2\) by
the disjoint-monomial estimate (DG), since
\(1/2+1/3+1/s\ge1\) for \(3\le s\le5\).
Reassembly in the ordinary complete series at the punctured charts
and weighted adjunction prove the assertion.

On a quintic elevenfold, a polar pencil admitting this first face
therefore has line threshold at least \(348/41\).
Indeed the outer radial denominator and the angular denominator are
both at most
\[
 1+\frac{4/5}{1+14/15}
    =\frac15+\frac45\,\frac{44}{29}
    =\frac{41}{29}.
\]
Consequently \(12/(41/29)=348/41>144/17\). This estimate is useful
after a linear restriction has left only three independent terminal
quintic directions.

\paragraph{The mixed quartic correction to a reduced binary cubic.}
Suppose the weighted light face at a general pole is
\[
  F(x,y,U)=a x^3+b y^3+xP(U)+yQ(U),
       \qquad U=(u,v),\quad ab\ne0 ,
\]
where \(P,Q\) are binary cubics, at least one nonzero, and
they have no common triple root (the zero form vanishes everywhere
to every order).
Then the ordinary graph of \(F\) has signed denominator at most
\(3/2\). In particular the line estimate on a quintic elevenfold
is at least \(60/7\). The light and graph weights are
\((3,3,2,2;9)\), giving light ratio \(10/9\). It is enough to
verify the ordinary signed estimate at every punctured light vector
by the punctured-graph criterion. In the following Hessians and higher
Taylor terms are taken after an affine graph shear subtracts the tangent.

If \(xy\ne0\), the Hessian on the first two coordinates already has
rank two. Suppose \(x\ne0,y=0\), the other alternative being
symmetric. We need only treat a rank-one Hessian. When \(U=0\)
and \(P\ne0\), its kernel cubic restricts to \(b y^3+xP(U)\),
which is nonpure. This gives the square-plus-nonpure-binary-cubic
bound by taking two general kernel lights. If \(P=0\), retain one
direction where \(Q\ne0\). Weights \((9,6,4)\) about the point
give the coupled fourth-order face consisting of a square, \(by^3\)
and a nonzero \(yu^3\), whose signed bound is \(3/2\).
When \(U=U_0\ne0\), the row of the Hessian
in the \(y\)-direction has zero diagonal and therefore must vanish.
In particular \(dQ(U_0)=0\). Put \(p=P(U_0)\). For \(p\ne0\)
restrict \(U\) to the radial line \((1+t)U_0\). The rank-one
condition on \((x,t)\) is \(p=4a x^2\), and the kernel contains
the direction \((dx,dt)=(-2x,1)\). The cubic Taylor term on
this direction equals
\(-8ax^3-6px+px=-7px\). The kernel also contains the
\(y\)-direction with nonzero pure cubic; no term is quadratic
in \(y\) and linear in the first direction. Thus their binary cubic
is nonpure. For \(p=0\), the rank-one condition on the Hessian blocks
gives \(x\operatorname{Hess} P=dP\otimes dP/(6ax)\). Contracting
with \(U_0\) gives \(dP=0\) by Euler, then \(\operatorname{Hess}P=0\).
If \(P\ne0\), then \(U_0\) is a triple root. The Hessian kernel
contains all \(U\)-variables and \(y\), and its cubic on a transverse
\(U\)-direction and \(y\) has both a nonzero pure cube of that
direction and \(by^3\), with no term quadratic in \(y\). This
gives the same estimate. If \(P=0\), then \(Q\) has a double
root here since its derivative vanishes and a common triple root is
excluded. Its second derivative in a transverse direction contributes
a nonzero \(y t^2\) to the kernel cubic, again giving a nonpure form.

Finally take \(x=y=0, U=U_0\ne0\). Any nonzero row of the mixed
derivative matrix \((dP,dQ)\) forces rank at least two. Otherwise
both binary cubics have a multiple root here, and at least one has
a double root by hypothesis. Retain a transverse direction \(t\)
at this root. The lowest Taylor term is
\[
                    ax^3+by^3+L(x,y)t^2 ,
\]
with \(L\ne0\). Its plane cubic is smooth if \(L\nmid ax^3+by^3\):
at \(t=0\) use squarefreeness, and at \(t\ne0\) the singularity
equations would force \(L=ax^3+by^3=0\), while at \([0:0:1]\)
the differential of \(L\) is nonzero. If \(L\) divides the binary
cubic, the curve is a line and a smooth conic meeting at two distinct
points since the complementary binary quadratic is nonzero at the
root of \(L\). The homogeneous three-light log-canonical graph test
applies. This finishes the punctured verification.

There is also a test for the common triple-root exception. Let the
two-dimensional normal frame before diagonalization be constant,
and let its constant radical space \(K\), including the line, have
ambient codimension at most five. Suppose its first order is five
and the nonzero mixed quartic has a common triple hyperplane on
\(K\) in its two derivative cubics. It has the form \(L(x,y)V^3\)
at the geometric generic pole, with \(V\) constant on the base.
Indeed the common hyperplane is uniquely determined and descends;
in a constant normal frame its primitive cube divides a nonzero
polynomial of binary base degree one. The constant cut \(V=0\)
on \(K\) has ambient codimension at most six and dimension at least
six, so it is not contained in the smooth quintic elevenfold.
Its order is five. Retain a general light \(T\) there on which the
quintic is nonzero, and a complementary light measured by \(V\).
With weights \((1/3,1/3,2/9,1/5;1)\) on \(x,y,V,T,Z\),
the new face is
\[
                    ax^3+by^3+L(x,y)V^3+cT^5,\qquad c\ne0 .
\]
In particular a quartic term with two or more normal coordinates
has weight at least \(2/3+2/5>1\). The radial light ratio exceeds
one. Where \(T\ne0\), its Hessian supplies an independent
square unless the Hessian of the first three variables vanishes.
That vanishing requires \(x=y=V=0\) since \(L\ne0\), and the
kernel cubic then contains the nonpure \(ax^3+by^3\). Where
\(T=0,V\ne0,x=y=0\), the mixed block has rank two. At
\(T=0,x\ne0,y=0\) or its symmetric counterpart, the proof
above applies verbatim on the single \(V\)-direction when
\(P(V)\ne0\) and \(V\ne0\), via the radial kernel calculation.
If the coefficient of \(xV^3\) vanishes and \(V\ne0\), the
\(y,V\) mixed entry gives rank two. At \(V=0\), the nonzero
\(xV^3\) supplies a cube independent of \(y^3\); in its absence
one uses the coupled square--cube--\(yV^3\) calculation as above.
If \(xy\ne0\), the first two derivatives already have rank two.

Thus all nonzero mixed quartic corrections are covered for a
squarefree binary cubic on a constant normal two-space with
\(K\) as above. Indeed the absence of a common triple hyperplane
gives the first test on two general kernel lights by Bertini, since
any triple root on the general line comes from the triple
hyperplane of a cubic. The two tests, together with the disjoint
binary cubic and squarefree binary quintic test in the absence of
a mixed term, cover all possibilities at kernel order five:
the Gauss bound gives singular codimension at least three on
the projective light space in that disjoint case.

\paragraph{Retaining a triple quartic hyperplane beside a binary cubic.}
Suppose the first cubic on a fixed two-space and a fixed light
kernel \(K\) is squarefree binary, in suitable pole coordinates
\(f=ax^3+by^3\), \(ab\ne0\). Suppose \(K\), including the
original line, has ambient codimension at most five. The \(60/7\)
line mark holds. The kernel order five and the order-four case in
which the kernel quartic has no triple component are covered by the
mixed and disjoint binary cubic calculations. In the remaining case
the quartic on \(K\) is \(V^3R\), where the unique triple or
quadruple hyperplane \(V\) is constant by the binary base degree
one. On \(K'=(V=0)\) the terminal quintic is nonzero: its ambient
codimension is at most six, hence its dimension is greater than
half the hypersurface dimension. Retain a general direction \(Y\)
there, with nonzero quintic coefficient, and a complementary
direction in \(K\) measured by \(V\).

If the quartic derivative in the binary-normal directions is nonzero
on \(K'\), its restriction along \(Y\) is \(L(x,y)Y^3\),
\(L\ne0\). Set \(w(x)=w(y)=1/3,w(Y)=2/9\). Take
\(w(V)=7/27\) when \(R|_{K'}\ne0\), or \(1/4\) when
\(R|_{K'}=0\). The face is respectively
\[
              f+L Y^3+eV^3Y,\qquad f+L Y^3+eV^4 ,
                   \qquad e\ne0 .
\]
The light sums exceed one. On \(V=0\) retain the other coordinates
at any nonzero punctured point. The mixed binary cubic test applies
also with a single \(Y\) at these nonzero points: when only
\(x\ne0\) or only \(y\ne0\) its missing second cube is supplied
by the coupled square--cube--\(yY^3\) test; for nonzero \(Y\)
the radial-kernel calculation applies verbatim. If \(V\ne0\) in
the fourth-power face its square is independent of the other Hessian;
the latter vanishes only at \(x=y=Y=0\), where its kernel contains
\(f\). In the \(V^3Y\) face at \(VY\ne0\) a nonzero \(x\) or
\(y\) gives a radial square independent of the \(VV\) entry
since \(F_{xV}=F_{yV}=0\). At \(x=y=0,Y\ne0\) a mixed
entry from \(LY^3\) and the zero normal block give rank two.
At \(Y=0,V\ne0\) the \(Y,V\) block has nonzero determinant.

Suppose now the quartic derivative restricts to zero on \(K'\).
Write \(L(x,y)VY^2\) for its next possible term on the retained
directions. Take \(w(Y)=1/5\). When \(R|_{K'}\ne0\), take
\(w(V)=4/15\) to obtain
\[
             f+L(x,y)VY^2+eV^3Y+cY^5,\qquad ec\ne0 .
\]
When \(R|_{K'}=0\) and \(L\ne0\) use the same weights and face
with \(e=0\). When both vanish, weight \(1/4\) for \(V\)
instead yields the diagonal form \(f+e'V^4+cY^5\), which has
the desired estimate. The indicated faces contain all weight-one
terms; terms with at least two binary normals or a quintic with
any coordinate of weight greater than \(1/5\) have larger weight.

For the coupled face with \(L=0,e\ne0\), the Hessian of
\(eV^3Y+cY^5\) is nonzero for every nonzero \((V,Y)\).
If both blocks are nonzero, the full Hessian has rank at least two.
If just the second block is nonzero, retain one square there and
the nonpure cubic \(f\). If just the first block is nonzero, rank
one in that block gives a square and a further independent cube;
on a general \((Y,V)\) line the first pure term is quartic,
so the square--cube--quartic test applies.

For \(L\ne0\) it remains to consider rank at most one. With
\(Y\ne0\), both \(x,y\) cannot be nonzero, while if both are
zero a mixed normal entry from \(LVY^2\) with the \(V\)-direction
gives rank two. Assume by symmetry \(x\ne0,y=0\). The row of
the Hessian in \(y\) must vanish, so its mixed \(y,V\) entry
forces the \(y\)-coefficient of \(L\) to vanish. The \(by^3\)
is now completely independent. The Hessian on \((x,V,Y)\)
has nonzero \(xx\) entry, so its kernel projects onto \(dY\).
Along any kernel line with \(dY\ne0\), the Taylor remainder after
subtracting its affine tangent has order \(k=3,4,5\), since its
degree-five term is \(c(dY)^5\). Retain this direction \(s\) and
a complementary square variable \(u\), as well as the independent
\(y\). For \(k=3\) one gets a square and two cubes. For \(k=4\)
the face is \(u^2+by^3+Bu s^2+Cs^4\), \(C\ne0\).
For \(k=5\) the same weights give \(u^2+by^3+Bu s^2\) if
\(B\ne0\); otherwise \(u^2+by^3+Cs^5\) is the face.
All these signed denominators are at most \(3/2\).

At \(Y=0,V\ne0\), if \(e\ne0\) the \((Y,V)\) block
has nonzero off-diagonal entry and zero \(VV\) entry. If \(e=0\)
and \(x=y=0\), the cubic on \(x,y,Y\) is
\(f+V L Y^2\), the smooth or transverse line-conic cubic of the
mixed-cubic lemma. At \(x\ne0,y=0\) the Hessian rank is two unless
\(L(x,0)=0\). In that last case \(L\) is a nonzero multiple of
\(y\), so the kernel cubic contains the nonpure \(by^3+c'yY^2\).
Finally, at \(V=Y=0\) with just \(x\ne0\) the restriction
\(V=0\) has a radial square in \(x\), the independent cube
in \(y\), and the fifth power in \(Y\). This proves the estimate
at all punctured points; the light sum exceeds one throughout.

In particular every reduced non-lc cubic with constant essential
normal space is covered: its dimension is two or three, so the
radical has ambient codimension at most five, and two general
essential directions cut a squarefree binary cubic. A moving
binary essential space with splitting \((0,1)\) can be made
constant by one preliminary hyperplane as in the quintic
tenfold proof; its resulting radical still has codimension
at most five.

\paragraph{A single cubic direction with quartic corrections.}
On the ordinary affine chart take the graph of \(x^3+xP(Y)\),
where \(P\) is a ternary cubic defining a reduced plane curve with
only ordinary nodes as singularities. Its signed denominator is
at most \(3/2\). The weights \((3,2,2,2;9)\) have light ratio
one. At \(x=0, Y=Y_0\ne0\) with \(dP(Y_0)\ne0\) the mixed
Hessian has rank two. At a projective node fix a nonzero radial
coordinate of \(Y\) and retain its two transverse parameters.
The initial cubic on three lights is \(x(x^2+uv)\) up to
rescaling, a transverse line and conic. The homogeneous
three-light log-canonical test applies.

At \(x\ne0,Y=0\) the square in the first variable and nonpure
kernel cubic \(xP(Y)\) give the bound. Suppose now \(x\ne0\),
\(Y_0\ne0\) and the Hessian has rank one. If \(p=P(Y_0)=0\),
contracting \(x\operatorname{Hess}P=dP\otimes dP/(6x)\) with
the radial vector gives \(dP=0\), then \(\operatorname{Hess}P=0\),
contrary to nodality. Hence \(p\ne0\). Choosing radial and
gradient-kernel coordinates gives \(P=pR^3+B(S,T)\):
the \(R^2(S,T)\) and \(R(S,T)^2\) terms vanish by the
differential and the rank condition. The binary \(B\) is nonzero
since \(P\) is reduced. The radial calculation of the mixed
quartic lemma gives a nonzero cubic in the Hessian kernel of
the \((x,R)\)-plane, while a general direction in \((S,T)\)
gives the independent cubic \(xB\). There are no mixed terms
between these two directions in the kernel cubic. This proves
all punctured estimates and therefore the assertion.

There is also a disjoint face \(x^3+G_4(Y)\) of signed
denominator at most \(3/2\) if the ternary quartic \(G_4\)
defines a reduced projective curve all of whose points have
multiplicity at most two, or the square of a smooth conic.
The light weights divided by the graph weight are
\(1/3,1/4,1/4,1/4\), with sum greater than one.
When both variables are nonzero each block has nonzero Hessian:
vanishing of all second derivatives of a homogeneous quartic
would give a triple point. At \(x\ne0,Y=0\) retain two general
\(Y\)-directions. Their binary quartic has no root of multiplicity
greater than two, so the square-plus-binary-quartic calculation
applies after translating \(x\) and taking its radial square.

It remains to consider \(x=0,Y=Y_0\ne0\) with rank-one
quartic Hessian. If \(G_4(Y_0)\ne0\), the radial direction
supplies a square. Choose a general kernel direction. On its
span with the radial direction, after making the gradient-kernel
change the quartic has the form \(aR^4+bRs^3+cs^4\), with
\((b,c)\ne(0,0)\) since the full quartic is not a pure power.
The first face after translation at \(R=1\) is a square plus
the independent \(x^3\) and a third or fourth power, as required.
At a smooth projective zero the Hessian has rank at least two
by contraction with the isotropic radial vector. At a double
point with rank-one Hessian fix the nonzero radial coordinate
and use transverse parameters \(u,v\) in which the quadratic
part is \(u^2\). A nonzero pure cubic on \(u=0\) gives the
square plus two cubes. Otherwise the next candidates in the
transverse binary quartic are \(buv^2+cv^4\).
If at least one coefficient is nonzero the weighted face with
the independent \(x^3\) is the mixed fourth-power test (including
the case \(c=0\)). If both vanish, the quartic's restriction to
the tangent line vanishes identically. In the reduced case the line
is a simple component, so the \(uv^3\)-coefficient must be
nonzero (otherwise every term is divisible by \(u^2\)).
The remaining face is \(u^2+x^3+c'uv^3\), covered by the coupled
fourth-order estimate with the mixed term in the square variable.
In the double-conic case the local conic equation has a nonzero
linear normal term and nonzero pure quadratic on its tangent;
its square therefore already falls in the \(c\ne0\) case.
This completes the punctured check.

\subsection{Reduced cubic polar pencils on a quintic elevenfold}
\label{sec:nov28-quintic-eleven-moving-cubic}
Every line of polar rank two and initial light order three with
geometrically reduced cubic satisfies \(\delta_C(X,H)>144/17\).
The homogeneous lc test gives \(60/7\) when the cubic is lc.
The non-lc classification leaves a squarefree binary cubic, a
ternary cuspidal cubic, or a ternary line tangent to a smooth conic,
in each case with an arbitrary radical. The binary cases, including
three concurrent lines, and all constant ternary cases were treated
in the triple-hyperplane lift. We give the remaining moving cases.
All constants in a graph at a general pole are taken over its
algebraic closure.

\paragraph{The two normal forms.}
For an irreducible ternary cusp the essential subbundle has splitting
\((0,0,c)\), \(1\le c\le2\), using the binary base degree two
as in the preceding essential-variable argument. Namely all three
positive degrees are impossible, a positive second degree would
force the first generator to divide the cubic, and the third occurs.
Write the two fixed generators as \(x,y\) and the remaining
degree-\(c\) generator as \(z\). If \(z\) occurs linearly,
its coefficient is the tangent cone at \([0:0:1]\), hence a square.
Its primitive root is fixed since its base degree is at most
\(2-c\le1\). Calling the root \(s\), and the other fixed normal
\(t\), the cubic is \(s^2 W+B(s,t)\) and
\(B(0,t)\ne0\) by irreducibility.

When \(z\) occurs quadratically, necessarily \(c=1\) and
\[
                 F=L z^2+Q z+B
\]
has \(L\) a fixed binary linear form, \(Q\) of binary base degree
one and \(B\) of binary base degree two. The point \([0:0:1]\)
is smooth. Projection from that point sends the cusp to a triple
root of \(\Delta=Q^2-4LB\). Indeed the line through those two
points is transverse to the tangent at the cusp, since the cusp
tangent already meets the cubic there with multiplicity at least
three. This line also differs from the tangent at the smooth point
by degree. Consequently the local quadratic equation over the
base-direction line at the image of the cusp is monic; completing
its square identifies an ordinary cusp with a discriminant zero of
order three (the local equation is a square plus a univariate
series, whose order is three for an ordinary cusp). This unique
root \(x\) is fixed since \(\Delta\) has base degree two.
It is independent of \(y=L\). Absorb the terms of \(Q\) divisible
by \(y\) by a shift of \(z\) by a base-linear binary linear form.
Then \(Q=kx^2\). Divisibility of
\(k^2x^4-4yB\) by \(x^3\) gives \(B=bx^3\), \(b\ne0\).
The shifted \(z\) is still a base-linear moving normal, and one
constant preliminary cut makes it a nonzero scalar times a fixed
normal. To see this write it as \(uU+vV\) for the two base
coordinates; if \(U,V\) are dependent no cut is needed. If the
coefficient span meets the fixed two-space of \(x,y\), the full
essential three-space was already constant. Otherwise impose a
general linear combination of \(U,V\), keeping the fixed two-space.
We have reached \(s^2W+ks t^2+bt^3\) with \(s,t\) fixed
and \(W\) independent on their common kernel.

For a tangent line and smooth conic, first suppose the unique line
factor \(s\) is constant. The conic restricted to \(s=0\) is a
square; the primitive root has base degree at most one. Extending
its equation to the full space writes the cubic as
\(a s t^2+s^2 W\) with \(t\) fixed or base-linear.
When needed, one constant restriction makes \(t\) fixed and
independent of \(s\). The complementary \(W\) remains independent
on the resulting kernel: otherwise the original essential
three-space lies in the constant span of \(s\) and the two
coefficients of \(t\), and was already constant. If the coefficient
span of \(t\) already contains \(s\), choose a cut that keeps
\(s,t\) independent by shifting the chosen hyperplane by a
complementary light if necessary; more directly, modulo \(s\)
the primitive square root then has degree zero and its fixed
representative suffices in the displayed normal form.

For a moving line factor, the residual conic has base degree at
most one. If it is constant, so is the essential three-space:
the tangent line in its essential plane pulls back to a hyperplane
containing the full radical. Otherwise its essential subbundle
\(E\) has splitting \((0,0,1)\) in the moving case, since a
rank-three quadratic in \(\operatorname{Sym}^2 E(1)\) has
nonzero determinant in \((\det E)^2(3)\). The cubic has that same
essential subbundle. If its third generator occurs linearly then
\([0:0:1]\) is the unique double point, and the tangent cone
there is the square of a fixed normal as in the cusp argument.
The tangent line is this fixed normal, a case already handled.
Otherwise use the quadratic presentation above. The smooth point
\([0:0:1]\) cannot lie on the line component since any line
through it which divides the cubic divides the fixed \(L\).
Projection from this point gives \(\Delta\) a quadruple root:
the line through the tangency point and the smooth projection
point is distinct from both curve tangents by degree on the conic.
The two smooth branches at the tangency point therefore give two
local graphs in this projection with difference of order two,
whose square is the discriminant up to a unit. Its unique linear
root \(x\) is constant and independent of \(y=L\). The same
shift gives \(Q=kx^2\), now \(B=0\) because the discriminant is
a scalar times \(x^4\). Thus \(F=z^2 y+kzx^2\), \(k\ne0\).
The previous one-cut argument makes \(z\) a fixed normal.

In all cases the two fixed normals \(s,t\) have common kernel
\(K\) of ambient codimension at most five, and the retained
cubic is one of
\[
 s^2 W+B(s,t),\quad B(0,t)=bt^3,\ b\ne0;
 \qquad s^2W+a st^2+\text{terms divisible by }s^2 ,
       \quad a\ne0 ,
\]
with \(W|_K\ne0\). Terms divisible by \(s^2\) in the second
formula can be absorbed into \(W\). In the first formula all
binary terms except \(bt^3\) have greater weight in the
extractions below. The first order on \(K\) is four or five.
At order five its quintic has projective singular codimension at
least \(13-2\operatorname{codim}K\ge3\) by the finite Gauss map.
It is therefore geometrically irreducible and reduced.

\paragraph{The cuspidal faces.}
At kernel order four, keep a general direction \(Y\) in \(K\)
and use weights \((3/8,1/3,1/4)\) on \(s,t,Y\), with unit
graph weight. The face up to nonzero rescaling is
\(s^2Y+t^3+eY^4\). Its light sum is \(23/24\).
At \(Y\ne0\) the \(s,s\) Hessian entry is nonzero; for
\(t=0\) a rank-one \((s,Y)\)-block requires \(s^2=6eY^3\).
The kernel direction \((ds,dY)=(-s/Y,1)\) has cubic coefficient
\(s^2/Y^2+4eY=10eY\), giving with \(t^3\) a nonpure kernel
cubic. At \(Y=0,s\ne0\) the mixed Hessian has rank two. With
just \(t\ne0\) use its radial square and the nonpure \(s^2Y\).

At kernel order five, if the quartic derivative in \(t\) is
nonzero there, use \((7/18,1/3,2/9)\). The face is
\(s^2Y+t^3+e tY^3\), \(e\ne0\), of light sum \(17/18\).
When \(Y\ne0\) its \(ss\) entry is nonzero; at \(t\ne0\)
the \(tt\) entry is independent, and at \(t=0\) the \(tY\)
entry is nonzero while \(st=tt=0\). When \(Y=0,s\ne0\) use
the \(s,Y\) mixed entry with \(ss=0\). With just \(t\ne0\)
the kernel cubic \(s^2Y+etY^3\) is nonpure.
If this derivative is zero, weights \((2/5,1/3,1/5)\) give
\(s^2Y+t^3+b'sY^3+cY^5\), \(c\ne0\), of light sum
\(14/15\). At \(Y\ne0\), a rank-one Hessian has a square
in \(s\) and requires \(t=0\); a kernel line in \((s,Y)\)
has nonzero \(dY\), so its remainder after subtracting the
affine tangent has order between three and five by \(cY^5\).
The independent cube in \(t\) and this remainder use precisely
the mixed E6, coupled-square, and E8 alternatives detailed in the
triple-hyperplane lift. At \(Y=0,s\ne0\) the mixed Hessian
gives rank two, and with just \(t\ne0\) use \(s^2Y\) as before.

\paragraph{The tangent-conic faces.}
Write \(x=s,y=t\), absorbing nonzero scalars. Kernel order four
with weights \((3/8,5/16,1/4)\) gives
\(x^2Y+xy^2+cY^4\), \(c\ne0\), of light sum \(15/16\).
Its Hessian has rank at least two at every nonzero point: at
\(x\ne0\) use \(F_{yy}=F_{xY}=2x\), \(F_{yY}=0\);
at \(x=0,y\ne0\) use \(F_{xy}\ne0,F_{yy}=0\);
when just \(Y\ne0\) the \(xx,YY\) entries are nonzero.
At order five, a nonzero quartic derivative in \(y\) gives
\(x^2Y+xy^2+c yY^3\) with weights
\((5/13,4/13,3/13)\) of sum \(12/13\). For \(Y=0\)
use the same mixed entries. For \(Y\ne0\) a rank-one Hessian
requires \(x=y^2/Y\) and simultaneously
\(3cY^4=2y^3\), \(3cY^4=y^3\) from its \(yY\) and
\(YY\) entries relative to its nonzero \(xx\) entry; this
contradicts \(c\ne0\).

If that quartic derivative vanishes, keep two general kernel lights
\(U\). With weights \((2/5,3/10,1/5,1/5)\) the complete face is
\[
       F=L(U)x^2+xy^2+xP(U)+y^2B(U)+G(U),
\]
where \(L=W|_U\) is a nonzero linear form, \(B,P,G\) have
degrees two, three and five, and \(G\) is the terminal quintic
restricted to these lights. The light sum is \(11/10\). The
binary choice can ensure that \(G\ne0\) at the root of \(L\)
and that \(T=P-2LB, R=G-BP+B^2L\) have no common root of
multiplicities at least two and three, respectively, with zero
forms interpreted as vanishing to all orders. For the first
assertion use irreducibility of the full terminal quintic.
For the second, a common prime factor \(D\) of those multiplicities
in the full light space would give \(G=R+BT+B^2L\) singular along
\(D=B=0\). If \(B\) is zero or divisible by \(D\),
\(D^2\) divides \(G\); otherwise the intersection is nonempty
of codimension two, contradicting the singular codimension at
least three. General line restriction detects the multiplicities
of prime divisors and avoids their pairwise intersections by
Bertini, proving the binary assertion.

At \(U=0\), a nonzero \(x\) gives \(F_{xx}=0\),
\(F_{xU}=2x\,dL\ne0\), and a nonzero \(y\) with \(x=0\)
gives nonzero \(F_{xy}\) and zero \(F_{xx}\).
Consider \(U\ne0\) with Hessian rank at most one. If \(y\ne0\),
\(F_{xy}=2y\ne0\) shows it has rank one with both \(x,y\)
coefficients of its normal nonzero. The kernel line in the
\((x,y)\)-plane has both coordinates nonzero; the cubic Taylor
term there is \(dx(dy)^2\ne0\). If the full kernel cubic is
nonpure we are done. If it is a pure cube \(w^3\), the simultaneous
kernel of the Hessian normal and \(w\) projects onto all the
\(U\)-directions, since these two functionals are independent
on the \((x,y)\)-plane. Choose a direction \(s'\) there with
\(G(dU)\ne0\). Its pure Taylor remainder has order four or
five since its cubic vanishes and its term of degree five is
\(G(dU)\). Retain a square variable, the cube direction, and
\(s'\). At order four use the square--cube--mixed-fourth face.
At order five an earlier square-variable times \((s')^2\) gives
the coupled square face; in its absence a nonzero \(w(s')^3\)
gives the coupled cube face with weights \((9,6,4;18)\).
If neither occurs use the square--cube--fifth face. These are
all possible earlier mixed terms by those weights, and each
has signed denominator at most \(3/2\).

If \(y=0,L\ne0\), rank one requires \(x=-B\). Its kernel
contains the \(y\)-direction, and the only terms in its cubic
involving \(dy\) are \((dy)^2(dx+dB)\). When this is nonzero
the kernel cubic is nonpure, having no pure \((dy)^3\).
Otherwise the Hessian normal is proportional to \(dx+dB\).
At this point set \(h=x+B\) and write exactly
\(F=Lh^2+h(y^2+T)+R\). The \(x,U\) mixed terms imply
\(dT=0\) in the \(U\)-directions after removing
\(2L(dx+dB)^2\). Hence \(T=0\) by Euler, and the residual
\(U,U\) Hessian is precisely \(\operatorname{Hess}R=0\).
This contradicts the chosen binary-root condition.

Finally suppose \(y=L=0\). If \(F_{yy}\ne0\) the \(y\)-row
has no mixed entries, so rank one puts all \((x,U)\)-directions
in the kernel. Its cubic has the nonzero \((dx)^2 dL\) but
no \((dx)^3\), hence is nonpure. If \(F_{yy}=0\), both the
\(x,y\)-rows have zero diagonals and must vanish for rank at
most one. At rank one the kernel cubic on these two directions
is the nonpure \(dx(dy)^2\). At rank zero the \(x,U\) row
gives \(2x\,dL+dP=0\), hence \(P(U)=0\) by Euler. The
second derivative in the radial \(U\)-direction is now
\(6xP(U)+20G(U)=0\), contrary to \(G\ne0\) at the root of
\(L\). This completes the punctured check.

Each punctured signed bound is at most \(3/2\). The smallest
primary light-to-graph ratio in these nonconstant cases is
\(12/13\), strictly greater than \(23/25\). The ordinary
primary signed denominator is therefore at most
\(1+13/25=38/25<73/48\). Including the regular and positive
line masses gives the uniform normalized denominator
\(1/5+(4/5)(38/25)=177/125<17/12\), hence in all reduced
cubic cases the common line mark \(500/59>144/17\).

\paragraph{Pure-cube initial lines on a quintic elevenfold.}
If the cubic initial on the full two-polar kernel is a pure cube,
the line satisfies \(\delta_C(X,H)\ge348/41>144/17\). Its root
\(x\) is constant since the binary base degree is two. On
\(x=0\) the ambient codimension is three. Its next order five
gives the terminal test directly. Otherwise write \(G\) for its
quartic, of base degree one.

When \(G\) is nonreduced, a repeated hyperplane is fixed: an
orbit of repeated hyperplanes of size at least two fills the
degree as the square of a primitive quadratic, then this square
root is constant in the binary base. Cutting a fixed factor
gives terminal order in codimension four. The other case is a
constant irreducible quadric squared. At rank at least four use
the four-light square-quadric test. At rank three keep its three
directions, giving \(x^3+q^2\) with light ratio
\(1/3+3/4>1\). At \(x\ne0\) and nonzero quadric vector the
two Hessian blocks both have positive rank. Indeed at a nonzero
isotropic vector the quadric has nonzero derivative; at a
nonisotropic vector its square has nonzero radial derivative.
On the \(x\)-axis its translated radial square and a nondegenerate
binary restriction of \(q\) give the square-plus-binary-quartic
test. At \(x=0,q\ne0\), the rank of the Hessian of \(q^2\)
is at least two since it is a nonsingular ternary quadratic
form times a scalar plus a rank-one form. On \(q=0\) at a
nonzero isotropic vector, fixing a nonzero hyperbolic coordinate
by an affine flag writes \(q\) as \(A+B^2\). The coupled
square--cube--quartic calculation applies to \((A+B^2)^2+x^3\).

Now let \(G\) be reduced. Three general lights give a reduced
plane quartic at general poles, and \(x^3+G(Y)\) is the first
face for light weights \(1/3,1/4,1/4,1/4\). At \(xY\ne0\),
both Hessian blocks are nonzero unless \(Y\) is a point of
multiplicity at least three, treated separately below. At \(Y=0\)
and \(x\ne0\), the translated square is disjoint from a general
squarefree binary quartic. At \(x=0\) a rank-one Hessian of
the quartic off its zero set gives a radial square and a nonzero
term of degree three or four in the gradient kernel. Indeed
in the radial frame the homogeneous form there is
\(aT^4+T P_3+P_4\), \(a\ne0\); both residual forms cannot vanish
for a reduced quartic. Use the square, \(x^3\), and the first such
kernel term, giving a square and a nonpure cubic or the
\((2,3,4)\) face. On the smooth zero set the Hessian already has
rank at least two by Euler. At a double point with rank-one
Hessian, fix a nonzero homogeneous coordinate and use \(y=0\)
for the double tangent. If the quartic stays nonzero along its
tangent line, the remaining pure order is three or four. The
rank-one cubic or mixed square--cube--quartic tests with \(x^3\)
apply. If that line is a component, write the local quartic as
\(y(ay+\cdots)\), with \(a\ne0\). Its complementary factor is
not divisible by \(y\); its restriction to \(y=0\) begins at
order two or three. The face is thus \(ay^2+x^3+b y z^k\),
\(b\ne0,\ k=2\) or three. The coupled square--cube computations
prove the \(3/2\) estimate for both faces.

Suppose a general plane quartic has a triple point. A reduced
quartic cannot have two such points: their joining line is a
component by B\'ezout, and removing it leaves a cubic with two
double points on the same line, forcing that line again.
General transverse linear restriction to smooth equimultiple
strata therefore shows that the full triple locus has a unique
linear codimension-two component \(K\). This component is
constant. The general transverse plane divisor at multiplicity
at least three is non-lc (its point blowup gives threshold at
most \(2/3\)); it follows by transverse log restriction at the
generic stratum that \(K\) is a codimension-two component of a
multiplier scheme with coefficient strictly less than one.
There are no divisorial components since \(G\) is reduced.
The base-line Nadel argument for the quartic pencil makes all
its horizontal components constant. In particular \(K\) has
ambient codimension five. Its terminal quintic \(Q\) has
projective singular codimension at least three, hence three
general \(K\)-lights give a smooth ternary quintic.

For multiplicity exactly three write the cubic tangent on the
transverse light plane in two variables \(y,z\). At \(x\ne0\)
above the triple point its nonpure tangent gives the square
and a nonpure binary cubic. At \(x=0\), if all three tangent
lines are distinct, \(x^3+C_3(y,z)\) defines a smooth plane
cubic, and its homogeneous cubic test suffices. For a double
and a simple tangent take \(C_3=y^2z\). If the quartic term
on \(y=0\) is nonzero, the first face is
\[
          x^3+y^2z+z^4, \qquad
          (w_x,w_y,w_z)=(1/3,3/8,1/4).
\]
Its light ratio is \(23/24\). On the \(x\)-axis use its radial
square and the nonpure \(y^2z\). Elsewhere the binary Hessian is
nonzero, giving rank two if \(x\ne0\). At \(x=0\) and rank one
we have \(y^2=6z^3\), \(z\ne0\); on the kernel direction
\(dz=1,\ dy=-y/z\) the binary cubic coefficient is \(10z\ne0\).
Together with \(x^3\) this gives the square-plus-nonpure-cubic
test. If instead the quartic vanishes on \(y=0\), the coefficient
of \(yz^3\) is nonzero because the curve is reduced. Its first
face is
\[
          x^3+y^2z+yz^3,\qquad
          (w_x,w_y,w_z)=(1/3,2/5,1/5).
\]
Its ratio is \(14/15\). The same axis comparison applies.
Rank one of the binary Hessian at a nonzero point requires
\(4y^2+9z^4=0\). A kernel vector with \(dz=1\) has
\(dy=-(2y+3z^2)/(2z)\ne0\). Therefore the quartic term on
the kernel is nonzero. If its cubic term is nonzero it joins
\(x^3\) as a nonpure cubic; otherwise use the coupled
square--cube--quartic test. Thus all punctured chart comparisons
are at most \(3/2\), and the radial denominators for these two
inner weights are at most \(1+1/(1+14/15)=44/29\).

If the transverse tangent at a general point of \(K\) is a pure
cube, its unique normal \(y\) is constant in the base and along
\(K\), because its cube divides an equation linear separately in
these two sets of coordinates. Dropping the original \(x\),
the full quartic in the original \(x=0\) is
\[
                    y^3 L_t(U)+P_t(y,z).
\]
If \(P_t(0,z)=0\), cutting \(y=0\) gives terminal order at
codimension four. Otherwise the nonzero \(z^4\) and the smooth
ternary quintic on \(K\) give exactly the transverse cubic
calculation of the preceding quartic strata proof, with or
without \(L_t\).

Finally suppose the generic multiplicity along \(K\) is four.
Then \(G=C_t(y,z)\) in two fixed normal coordinates, since
multiplicity four along \(K\) gives membership in the fourth
power of its ideal. The binary quartic \(C_t\) is squarefree.
Discard the original \(x\). Give \(y,z\) weight \(1/4\) and
three general directions on \(K\) weight \(1/5\). The face is
\(C_t(y,z)+Q(U)\) for a smooth ternary quintic \(Q\): all
mixed quintic terms have larger weight. This is the disjoint
binary-quartic and ternary-quintic test in the quartic strata
proof, of radial ratio \(1/2+3/5>1\) and signed denominator
at most \(3/2\). The radial and pole proportions in all cases
give denominator at most \(1/5+(4/5)(44/29)=41/29\), proving
the claimed threshold.

\paragraph{A forward-polar improvement at a quintic point.}
At a quartic tangent contact on a smooth quintic elevenfold write
\(Tf+g\) for the tangent hyperplane section. Suppose \(f=0\)
contains a line \(\ell\) having a squarefree normal cubic polar
\(P\), with \(P\) and \(g|_\ell\) coprime. Then a point optimizer
normalized as above and of ratio at most eight satisfies \(S\ge7/4\).
Indeed use the constant three-space spanned by the line and its
polar direction. As in the forward-polar flag the two angular
initials become \(zP(x,y),G(x,y)\), where \(G=g|_\ell\).
They are a regular sequence. In fixed binary degree \(D\), fixed
\(z\)-exponent \(s\), and section-order excess \(l\), the binary
block is
\[
\sum_{i\le s,\ 3i+h\ge l} P^iG^h\,
         \mathbf C[x,y]_{D-3i-5h}.
\]
It contains the full series after removing \(P^I G^J\) up to
bounded codimension, where \(I=\max(0,\lceil(5l-D)/12\rceil)\)
and \(J=(l-3s)_+\). To verify this, use the free homogeneous
remainder basis over \(\mathbf C[P,G]\). For each bounded
remainder degree \(b\) the monomials \(P^iG^h\) of degree \(D-b\)
are allowed exactly when \(h\ge J\) and \(12i\ge5l-D+b\).
The bounded difference in the lower endpoint costs only bounded
codimension. At a root of \(G\) the normalized marked mean before
averaging is at most \(D/2-3I/2+5J/2+O(1)\). At other points
the mean is at most \(D/2+O(1)\) by squarefreeness of \(P\).
Here and below the normalized aggregate means are in units of
\(12/m\), with section degree \(m\).

The normal basis \(T^a Z^i W^h R\) has \(l=3i+h-a\), and
\(s\ge i\). Hence \(\overline{J}\le1/20\), using the simplex
mean of \((E_1/4-E_0)_+\) for independent standard exponentials.
The residual ternary ring has degree twenty, with multiplicity
fifteen on its \(z\)-axis and five off it. Its positive-\(z\)
monomials have bounded binary remainder degree since \(P,G\)
form a binary regular sequence. In this axis mass one has
\(
 5l-D=12i-5a+O(1)
\);
hence \(\overline I\ge(3/4)\mathbf E[(E_1-5E_0/12)_+]=9/17\).
All other sectors contribute nonnegatively. Total mean binary
degree is \(3+5/4+1/4=9/2\); mean \(z\)-degree is \(1+3/4=7/4\).
The binary vertex and all its marked directions therefore have
normalized denominator at most \(9/4\), as
\(-3(9/17)/2+5(1/20)/2<0\).
The complementary coordinates each have mean one.
The ordinary point radial mean is \(57/4\), so adjunction gives
\[
 12S_H(v)\le57/4+(9/4)q,\qquad A(v)=11+q.
\]
When \(A(v)\le8S_H(v)\) this gives \(q\ge3\), as claimed.
The line hypothesis cannot be imposed on every quartic tangent
form. For \(f=Q^2\) every first normal polar on a line in \(f=0\)
vanishes identically. The reverse polar and repeated-factor
calculations above still give the bound \(S\ge13/8\).

\paragraph{Coupling a nowhere-zero initial cubic to its optimizing point.}
Suppose every line of a smooth quintic elevenfold has
\(\delta_C(X,H)\ge1380/163\), and that all lines except a class
\(\mathcal E\) have a common bound strictly larger than \(144/17\).
Suppose each line in \(\mathcal E\) has a rank-two polar pencil,
with initial light cubic nonzero as a form on its constant kernel
at every point of the line. Then the elevenfold is uniformly
K-stable. This nowhere-zero hypothesis holds in particular for
a cubic \(L^2M_t\) with \(L\) fixed and \(M_t\) a primitive
base-quadratic moving hyperplane: the coordinate quadratics of
\(M_t\) have no common zero on \(\mathbf P^1\), and the binary
degree has been exhausted, leaving a nonzero constant prefactor.

Indeed if \(\beta\le8\), positive-dimensional original optimizing
centers contain lines and contradict the common bound above eight.
Fix an optimizing point \(p\). By the proof of the line criterion
at \(144/17\), using the first comparison just above that value,
there must be a line through \(p\) of mark below the common bound
for lines outside \(\mathcal E\). To see the pointed assertion,
let \(Z\) be the minimal positive center and \(N\) the finite
non-qlc locus containing \(p\). Vanishing for the finite union
of \(N\) and all zero-dimensional qlc centers contained in \(Z\)
makes constants surject onto its global sections, so its support
is \(\{p\}\). The curve supplied by the length-one theorem meets
this support and hence passes through \(p\).

The tangent hyperplane at \(p\) contains the constant light
kernel of this line. On its affine slice through \(p\) the cubic
does not vanish identically, by hypothesis, so the multiplicity
of the tangent hyperplane section is two or three. In the first
case the signed surface bound has \(B=5/2\); in the second one
may take \(B=13/5\). For the latter bound a general pointed
quintic surface has integral plane tangent curve of multiplicity
three and arithmetic genus six. The plane-curve discriminant
bound gives local threshold at least \(10/21\); the resulting
coefficient is at most \(5/3+(1-5/9)(21/10)=13/5\).
The nine pointed flags give
\[
  12 S\le(13/5)(q+2)+9,\quad 11+q=\beta S\le8S,
  \qquad q\ge23/11,\quad S\ge18/11 .
\]
Rerun the untilted interpolation at this same point with
\[
  t_0=\frac{1+8(18/11-1)}{18/11-11/12}
      =\frac{804}{95}<\frac{1380}{163},\qquad
       (t_0-8)(18/11)+8-11t_0/12=1 .
\]
For \(\beta\ge31/4\), choose \(\theta=3/8<F\) using the
prime-fiber estimate. Then the marked ample coefficient at
\(t_0\) is \(8-(11/12-1/32)t_0>1/2\), since \(t_0<144/17\).
Increase \(t\) slightly, keeping it strictly below \(1380/163\).
For \(\beta\le31/4\), use \(t=t_0-1/8,\ \theta=1/4\).
Here \(t-\beta>t_0-8\) and the unmarked ample coefficient also
increases. The marked ample coefficient
\(8-(11/12-1/48)(t_0-1/8)\) lies strictly between \(1/2\) and
one. Thus the first positive-center test contradicts the bound
\(1380/163\) on every line, and the isolated-point Nadel
comparison gives the final contradiction.

\paragraph{A nonzero line cubic and tangent quintic contact.}
In the preceding two-pass criterion one may replace \(1380/163\)
by any common mark strictly larger than \(312/37\) if each of
the exceptional lines has a nonzero initial light cubic
at every point. This hypothesis holds for \(L_{[s:t]}^2 M\)
when \(L_{[s:t]}=sa+tb\) and the three linear forms \(a,b,M\)
are independent, and also if \(M\) belongs to the span of
the independent \(a,b\). In particular \(76/9>312/37\) is sufficient.

Only the cubic-contact estimate at the point supplied by the
first pass changes. A nonpure restriction makes the full cubic
tangent form nonpure. Its restriction to two general tangent
directions has a simple root and is either squarefree or of
type \((2,1)\). Indeed a general projective line meets the
reduced components of its divisor transversely at smooth
points away from their intersections. A general pointed
surface section therefore has integral plane tangent quintic
with such a cubic tangent at \(p\). In the squarefree case
the plane germ has threshold \(2/3\) by one ordinary blowup.

In type \((2,1)\) the completed plane germ splits as a smooth
branch and a reduced double germ transverse to the first
branch's tangent. This is the coprime-factor lifting of the
tangent cones. Use the smooth branch itself as coordinate
\(x=0\), and take \(y=0\) tangent to the other germ. By
Weierstrass preparation in \(y\), completion of the square
by a translation of \(y\), and multiplication by a unit,
the full germ becomes \(x(y^2+x^m)\), with \(m\ge3\).
To normalize the coefficient of \(x^m\) one may absorb a
unit root into \(x\), preserving the smooth branch. The
delta invariant is
\[
  \delta_p=\left\lfloor m/2\right\rfloor+2 .
\]
Here the double germ has the displayed first summand (from
its two branches for even \(m\), or its parametrization for
odd \(m\)), and its intersection multiplicity with \(x=0\)
is two. Integrality of the projective plane quintic bounds
\(\delta_p\le p_a=6\), so \(m\le9\). The local log canonical
threshold of \(x(y^2+x^m)\) is \((m+2)/(2m+2)\ge11/20\).
Indeed the weights \(2,m\) give that ratio. The Newton edge
has squarefree roots on the dense torus. The positive toric
divisors on a regular subdivision have ratio
\((a+b)/(a+\min(2b,ma))\), minimized at \(2b=ma\);
the coordinate component has threshold one.

Thus the signed surface coefficient in cubic contact is at most
\[
       \frac53+\left(1-\frac59\right)\frac{20}{11}
          =\frac{245}{99}<\frac52 .
\]
The other terms in the surface flag are at most \(5/2\), as
in the quadratic-contact case. Hence \(12S\le(5/2)(q+2)+9\)
and \(11+q\le8S\) give \(q\ge5/2,\ S\ge27/16\).
The untilted interpolation at the same point now has equality
at \(t_0=312/37\), since
\[
                 (t_0-8)(27/16)+8-11t_0/12=1 .
\]
The large-\(\beta\) marked comparison still exceeds \(1/2\)
as \(t_0<144/17\). For \(\beta\le31/4\), \(t=t_0-1/8\)
and \(\theta=1/4\) give
\(1/2<8-(11/12-1/48)(t_0-1/8)<1\). The remaining inequalities
are strict when \(t\) is slightly increased in the first
case, exactly as in the preceding two-pass proof.

For completeness, if the tangent form at the point is a pure
cube, the surface tangent cone is \(y^3\) in affine linear
coordinates. Integrality of the plane quintic and B\'ezout give
intersection order \(I=4\) or \(5\) with \(y=0\).
When \(I=4\) the Newton boundary is the single nondegenerate edge
with initial \(y^3+ax^4\), \(a\ne0\). When \(I=5\) the only
monomial strictly below the segment from \((0,3)\) to \((5,0)\)
is \(x^3y\): every term other than \(y^3\) has ordinary degree
at least four. In its absence this segment has nondegenerate
initial \(y^3+ax^5\). In its presence the two initials
\(y^3+bx^3y\) and \(bx^3y+ax^5\) have simple torus zeros and
primitive inward normals \((2,3)\) and \((1,2)\).
The corresponding toric discrepancy ratios are \(5/9\) and
\(3/5\). A regular toric subdivision resolves these plane
germs, since all face polynomials have simple zeros on the
torus and their edge endpoints have nonzero coefficients.
The minimum of the ratios in the three cases is \(8/15\),
attained at weights \((3,5)\) in the single quintic-edge case.
It gives the surface coefficient \(5/3+(4/9)(15/8)=5/2\).
Thus precisely the same improved comparison holds at this last
tangent type.

\paragraph{Restoring a factor of a nonreduced cubic.}
Suppose the first light form of a two-polar line on a smooth
quintic elevenfold is \(L^2 M\), with distinct factors over the
geometric generic pole field. Their multiplicities distinguish
them over the rational function field of the binary base.
On choosing primitive representatives the base degrees satisfy
\(2\deg L+\deg M\le2\), so at least one factor is constant.
The following tests retain a direction transverse to such a
factor. All the asserted sections and first weights can be
chosen over the constant field whenever their stated conditions
hold on a general constant section. Secondary coordinates and
Taylor expansions are taken separately at each general pole.

First let the double factor be constant, with coordinate \(x\).
Suppose that on a constant light subspace in \(x=0\) the first
nonzero form is \(G_h(Y)\), \(h=4\) or \(5\).
Let \(l=M|_Y\ne0\) and retain a transverse \(x\)-direction.
For \(h=5\) assume the first normal derivative of the quartic
vanishes on this section. The weights
\(((h-1)/(2h),1/h,\ldots,1/h;1)\) have face
\[
                      x^2l(Y)+G_h(Y).
\]
Take at least three \(Y\)-directions. Assume that \(G_h\) and
its restriction \(g=G_h|_{l=0}\) have nonzero affine Hessian
at every nonzero vector of their respective spaces. Equivalently
their projective hypersurfaces have no point of multiplicity at
least three, by the homogeneous Euler identities.
For \(h=4\) the condition on \(g\) can be weakened to \(g\ne0\).
More generally at \(h=4\) one may admit projective triple points
of \(G_4\) away from \(l=0\) provided their tangent cubics are
not pure cubes.
Then the signed denominator of this ordinary graph is at most
\(3/2\), and the line mark is at least \(60/7\).

To see the assertion, write \(F=x^2l+G_h\) and subtract the affine
tangent at each graph point. At \(x\ne0,l=0\), the \(xx\) entry
of the quadratic vanishes and its \(x,Y\) row is \(2x\,dl\ne0\),
giving rank at least two. At \(x=0,l\ne0\) the \(xx\) entry and
the nonzero Hessian of \(G_h\) likewise give rank at least two.
At an admitted triple point choose one general kernel plane
on which its cubic stays nonpure. With weights \(1/2,1/3,1/3\)
the \(x^2l\) term gives a separate square and all its corrections
in the kernel directions have weight greater than one. This is
the square--nonpure-cubic test.
At \(x=l=0,Y=p\ne0\) the Hessian of \(g\) is nonzero. If the
whole Hessian has rank one, its kernel contains the \(x\)-axis
and a direction \(e\) in the \(Y\)-space with \(dl(e)\ne0\).
Indeed its unique supporting covector is nonzero on \(\ker dl\).
On this kernel plane the cubic is \(a x^2s+b s^3\), \(a\ne0\),
so the square-plus-nonpure-cubic test applies.

For \(h=4\), suppose at the latter point the restriction of the
Hessian to \(l=0\) vanishes. The whole Hessian remains nonzero,
so only rank one with unique supporting covector \(dl\) remains.
The restriction \(g\) has multiplicity \(j=3\) or \(4\) at \(p\).
Retain a general direction \(u\) detecting this multiplicity in
\(\ker dl\), the transverse coordinate \(s=l\) and \(x\).
Weights \(1/2,1/4,1/j\) on \(s,x,u\) give
\[
  s^2+s x^2+c u^3\quad (j=3),\qquad
  s^2+s x^2+b s u^2+c u^4\quad (j=4),\qquad c\ne0,
\]
after rescaling. A nonzero cubic with one \(s\) and two \(u\)
has weight greater than one in the first case. We check their
punctured graphs. For \(j=3\) and \(u\ne0\) the \(s\)- and
\(u\)-entries give rank two. At \(u=0\) the binary
\(s,x\) block has rank at least two unless \(s=x^2\).
At the latter punctured point \(x\ne0\). Its kernel direction
is \(ds=-x\,dx\), with cubic coefficient \(-x(dx)^3\).
The independent \(u^3\) gives a nonpure kernel cubic.

For \(j=4\) the Schur complement of the \(ss\) entry forces,
at rank one, \(s=x^2,\ bxu=0\) and
\(2bs+(12c-2b^2)u^2=0\). If \(b=0\) then \(u=0\), and
\(x\ne0\) gives the preceding transverse square and kernel
cube together with the independent \(cu^4\), a \((2,3,4)\)
face. If \(b\ne0\) these equalities force \(x=s=0,\ u\ne0\)
and \(b^2=6c\). Restrict the kernel to \(dx,du\). There
\(ds=-bu\,du\) and the cubic contains
\(-bu\,du\,(dx)^2\), with no \((dx)^3\) term, so is nonpure.
All punctured points and the secondary radial light sums,
which are at least one, give signed denominator \(3/2\).

Suppose \(x\ne0,l(p)\ne0\) and the Hessian has rank one. Normalize
\(l(p)=1\), and write a variation of \(Y\) as \(rp+U\),
\(l(U)=0\), with \(r=1\) at the given point. The Schur complement
of the \(xx\) entry gives
\[
        \operatorname{Hess}G_h(p)=2x^2\,dl\otimes dl .
\]
Euler's identity gives \(h(h-1)G_h(p)=2x^2\) and
\(dG_h(p)|_{\ker dl}=0\). Consequently
\[
  F=x_{\rm var}^2 r+c r^h+
           \sum_{j=3}^{h} r^{h-j}G_j(U),\qquad c=G_h(p)\ne0.
\]
Here the last coefficient is \(g(U)\ne0\). In the radial plane
the Hessian kernel is \((dx_{\rm var},dr)=(-x,1)\). Its cubic
coefficient is \(x^2+\binom h3c=(h+1)x^2/3\ne0\).
The Hessian supplies a transverse square. If \(G_3\ne0\), a
general direction \(U\) gives a second cube, with no mixed
term on the two kernel directions, and the kernel cubic is
nonpure. Otherwise the first nonzero \(G_j\) has \(j=4\) or
\(5\). A general \(U\)-direction and the weights \(1/2,1/3,1/j\)
in the square, radial kernel, and \(U\)-directions give the
disjoint \((2,3,j)\) face. The factors \(r^{h-j}\) are units
with constant leading coefficient. Both signed estimates are
\(3/2\). The first radial light sum is at least
\((h+5)/(2h)\ge1\). The punctured-graph and line adjunction
therefore give the assertion.

There is a second test when the simple factor is constant, with
coordinate \(x\). Put \(l=L|_Y\ne0\) on a constant subspace of
\(x=0\) whose first order is \(h=4\) or \(5\). The weights
\(((h-2)/h,1/h,\ldots,1/h;1)\) give the face
\[
                       x l(Y)^2+G_h(Y).
\]
No restriction on the quartic first derivative is needed when
\(h=5\), since its weight is then \(6/5\). Suppose the Hessian
of \(g=G_h|_{l=0}\) is nonzero at all its nonzero vectors, and
retain two directions in \(\ker l\) for \(h=4\), or three for
\(h=5\). In the latter case also assume the ternary quintic
\(g\) has log canonical threshold at least \(3/5\). This holds
in particular for a smooth plane quintic. The signed denominator
is again at most \(3/2\).
Indeed \(l\ne0\) gives a nonzero \(x,Y\) row with zero \(xx\)
entry. At \(l=0,Y=p\ne0\) a rank-one Hessian has kernel
containing \(\partial_x\) and a direction \(e\) with \(dl(e)\ne0\),
by the Hessian hypothesis on \(\ker l\). Its cubic on those
two directions is \(a t s^2+b s^3\), \(a\ne0\), and the
square-plus-nonpure-cubic estimate applies. At \(Y=0,x\ne0\)
fix \(x\) and take weights \(1/2\) on \(l\) and \(1/h\) on
the \(\ker l\) directions. The resulting face is the square
plus \(g\), covered by the square--binary-quartic or
square--ternary-quintic calculation. In the ternary proof the
terminal estimate \(v(g)\le5A/3\) is precisely the indicated
log canonical threshold condition, and the punctured nonzero
\(Y\)-directions use its Hessian hypothesis. Its radial sum is
at least one. The original sum of
light weights is also at least one. These estimates give line
mark \(60/7\).

If the quartic vanishes entirely on the chosen constant factor
hyperplane, the terminal test already applies there in ambient
codimension three. If its quartic \(G\) has a constant hyperplane
factor, one further cut reaches terminal ambient codimension four.
If instead \(G\) has a rational moving hyperplane factor \(P\),
its primitive coefficients have base degree one. Write
\(P=t_0a+t_1b\) with independent constant \(a,b\); base-degree
additivity gives \(G=P Q\) for a fixed cubic \(Q\). Put \(K=(a=b=0)\)
in the enlarged light space. If \(Q|_K\ne0\), take three general
directions \(U\) on \(K\) and one fixed complementary direction
\(z\). The weights \(1/5\) on \(U\) and \(2/5\) on \(z\)
give face
\[
                       z B_3(U)+C_5(U).
\]
Here \(B_3\ne0\) is a multiple of \(Q|_U\) at general poles,
and \(C_5\) is a smooth ternary quintic by the terminal Gauss
bound in ambient codimension five. The original cubic vanishes
on the entire chosen hyperplane, so no lower terms survive.

If the ternary cubic \(B_3\) is log canonical, this face has
signed denominator \(3/2\). At \(U=0,z\ne0\) keep \(z\) fixed.
The first homogeneous form is \(B_3\); three ordinary lights
and the homogeneous log-canonical cubic test apply. At a point
\(U=p\ne0\), the only possible low-rank Hessian has \(dB_3(p)=0\),
since a nonzero \(z,U\) row and zero \(zz\) entry give rank two.
Euler then gives \(B_3(p)=0\). Such a point is necessarily a
double point with nonzero Hessian \(J\), as a plane cubic with a
triple point fails log canonicity there. We have \(Jp=0\), and
\((zJ+\operatorname{Hess}C_5(p))p=4dC_5(p)\ne0\). Therefore the
\(U\)-block of the Hessian has positive rank. In rank one its
unique supporting covector \(\theta\) satisfies \(\theta(p)\ne0\).
It follows that \(J|_{\ker\theta}\ne0\). On the kernel of the
full Hessian the cubic contains \(sJ(e,e)t^2/2\) for \(s\) the
\(z\)-variation and a suitable kernel vector \(e\) with
parameter \(t\). Its binary restriction is nonpure, giving the
square--cubic estimate. The outer radial light ratio is one.

\paragraph{Three independent coefficients in the simple factor.}
Let the cubic be \(x^2M\) with \(x\) constant. Suppose the
coefficient linear forms of the primitive binary expression
for \(M|_{x=0}\) span a three-space. Then every such quintic
elevenfold line has threshold at least \(60/7\).

Here \(M|_{x=0}\) is primitive of binary degree two. Its
primitive polynomial cannot divide a nonzero quartic on \(x=0\),
whose binary degree is at most one. If the quartic vanishes on
that hyperplane, the terminal test already applies at ambient
codimension three. Otherwise let \(G\) be this quartic. A fixed
repeated hyperplane of \(G\) gives the terminal test after one
more cut. The other nonreduced possibility is a constant
irreducible quadric squared, by the base degree one factor
argument in the pure-cube case. A general plane cuts it in
a smooth conic squared, which has no triple point.

For reduced \(G\), if its general plane section has a triple
point the triple locus has a fixed linear codimension-two
component \(W\), unique in this codimension, as proved in
the pure-cube calculation. If the multiplicity there is
three and the normal cubic is nonpure, a general transverse
plane sees this same tangent type. At its marked point \(p\)
one has \(l(p)\ne0\) for \(l=M|_{x=0}\) at a general pole:
otherwise all three independent fixed coefficient forms would
vanish on the codimension-two \(W\). The restored-factor graph
therefore has signed denominator \(3/2\), both at this
point and elsewhere. In applying that graph criterion the
quartic on the general moving line \(l=0\) in the plane is
nonzero: the full restriction modulo \(l\) is nonzero by the
binary degree comparison, and a general constant plane
preserves that property at general poles. The same graph test
handles every general plane quartic without a triple point.

If instead the transverse cubic is pure, discard the original
\(x\): the remaining quartic has the same form
\(z^3L(U)+P_4(z,y)\) and the same smooth ternary terminal
quintic in \(W\) as in the pure-cube calculation. Its
transverse cubic and quartic graph gives \(3/2\). Finally,
multiplicity four along \(W\) makes the reduced quartic
a squarefree binary quartic \(C(y,z)\) in its fixed normals.
Retain two such normal coordinates and three general
\(W\)-directions with a smooth ternary terminal quintic \(P\),
available because \(W\) has ambient codimension five. With
weights \(1/4,1/4,1/5,1/5,1/5\) in \(x=0\) the face is
\(C+P\). Each block has positive Hessian rank at any of its
nonzero vectors. When both are nonzero this gives rank at least
two. When just one is nonzero, retain one direction there with
a nonzero quadratic Taylor coefficient, and the full other
block. The square--binary-quartic or square--ternary-quintic
test proves the comparison. The first radial light sum \(11/10\)
is sufficient, as in all the preceding cases.

\paragraph{Every fixed double factor in the initial cubic.}
Let \(x^2M\) be the initial cubic of a two-polar line on a
smooth quintic elevenfold, with \(x\) constant and \(M\) distinct
from \(x\). Its line threshold is at least \(1860/219=620/73\).
More precisely the only new denominator is the inner
\(1+16/31=47/31\), and after the \(1/5,4/5\) pole proportions
the mark is
\[
       \frac{12}{\,1/5+(4/5)(47/31)\,}
              =\frac{620}{73}>\frac{144}{17}.
\]
The outer light ratio is \(3/8+3/4=9/8\).

Put \(l=M|_{x=0}\), dropping any binary scalar from this
restriction when discussing its primitive equation. Use the
same quartic \(G\), reduced plane sections and fixed
multiplicity-three locus \(W\) as in the three-independent-
coefficient case. That proof already treats the nonreduced
quartics, a pure transverse cubic, and generic multiplicity
four, except for the need for \(G|_{l=0}\ne0\), addressed below.
Let a general plane meet \(W\) at \(p\) and have a nonpure
binary tangent cubic \(f\) there. If \(l(p)\ne0\) the preceding
test applies. Otherwise take an affine transverse plane through
\(p\); its two parameters \(y,z\) identify \(l\) with a
nonzero linear form. At \(x=0\) the first cubic on \(x,y,z\)
is
\[
                          x^2l(y,z)+f(y,z).
\]
It defines a log canonical plane cubic unless the repeated
root of \(f\) is the root of \(l\). To check this, on \(x=0\)
a multiple root not on \(l=0\) gives a node by the independent
\(x^2\) and root-transverse square. When \(x\ne0\), the
\(x\)-derivative at a singular point forces \(l=0\). The point
\([1:0:0]\) is smooth by the differential of \(l\). A singularity
elsewhere on this line requires \(l\) to divide \(f\).
If \(l\) is a simple factor the equation is \(l(x^2+f/l)\);
its two distinct intersection points on \(l=0\) are transverse
and the binary quadric may itself split into distinct lines
after adding \(x^2\). Its singular point in this last case is
a node away from \(l=0\). At \(x\ne0\) above the original
triple point of \(G\), with \(l=0\), the mixed Hessian row
of \(x^2l\) already gives rank two.

The exceptional situation therefore has \(l=y,\ f=y^2z\)
after rescaling. If the quartic restriction in \(x=0,y=0\)
is nonzero on this transverse section, it supplies \(cz^4\)
with \(c\ne0\). Weights \(5/16,3/8,1/4\) on \(x,y,z\)
give the second face \(x^2y+y^2z+cz^4\). Its Hessian has
rank at least two at every nonzero vector. When \(y\ne0\),
the minor in rows \(x,z\), columns \(x,y\), is \(4y^2\).
When \(y=0,x\ne0\) the \(x,y\)-minor equals \(-4x^2\);
on the remaining \(z\)-axis the two nonzero entries \(2z,12cz^2\)
suffice. The radial light ratio is \(15/16\), giving precisely
the claimed inner denominator.

If the repeated-tangent exception occurs at a general pole
and general \(p\in W\), then \(l\) vanishes on \(W\).
The complete order-three polynomial in its two normal
coordinates, linear in the longitudinal variables, is
divisible by \(l^2\). Indeed this divisibility holds upon
evaluation at general longitudinal points and hence holds
identically. The remaining fourth-order term depends only
on the two normals. If its restriction to \(l=0\) is zero,
then \(l\mid G\) on the full fixed factor hyperplane.
If it is nonzero, a transverse plane sees exactly the
\(cz^4\) just used.

It remains to treat every case \(l\mid G\). A fixed \(l\)
gives a terminal cut in ambient codimension four. For a
nonconstant primitive \(l\), base degrees force \(\deg_t l=1\)
and \(G=lQ\) up to scalar with \(Q\) a fixed cubic.
Any linear factor of \(Q\) again provides the terminal cut,
so take \(Q\) integral. On a general constant plane \(Q\)
is a fixed integral plane cubic by Bertini, and \(l\)
runs in a basepoint-free pencil of linear forms as polynomials
in the binary parameter, defining a pencil of projective
lines on that plane with a single fixed base point.
A general member cuts \(Q\) in simple points away from
the fixed base point: projection from that point has
separable generic fibers on the normalization of \(Q\)
and avoids its finitely many singular or ramified values.
At a smooth fixed base point it also cuts transversely,
since just one pencil line is tangent. Consequently any
multiple intersection there is a singular point of \(Q\).
Such a point is double and a general line of the pencil
avoids its finitely many tangent directions.

Apply the restored graph \(x^2l+lQ\). For \(x=0,l\ne0\)
and \(G\) of multiplicity at most two use its nonzero
Hessian as before. A triple point with \(l\ne0\) would give
one of \(Q\) itself; the integral plane cubic cannot have
one. For \(x=l=0\), a simple root of \(Q|_{l=0}\), or a
point off \(Q\), gives Hessian rank at least two: the row
of \(d l\otimes dQ+dQ\otimes dl\) has rank two whenever
\(dQ|_{\ker dl}\ne0\), and off \(Q\) this last differential
is nonzero by Euler. At a singular base point of \(Q\), the
cubic tangent of \(G\) is \(lq_2\) with \(q_2|_{l=0}\ne0\)
for a general pencil line. It was among the log canonical
ternary cubics treated above. For \(x\ne0,l=0\) the nonzero
mixed row of \(x^2l\) suffices. Finally for \(x l\ne0\)
the radial-kernel proof of the first restored-factor lemma
still applies: in its rank-one normal form some \(G_j\) with
\(j\ge3\) must be nonzero even though \(G|_{l=0}=0\),
because \(G=lQ\) with \(Q\) integral cannot be \(c l^4\).
This supplies every missing punctured estimate.

\paragraph{A square plus a binary quintic.}
Let \(g(B)\) be a binary quintic whose projective roots have
multiplicity at most two. Its ordinary affine graph
\(W=(Z+x^2+g(B)=0)\) in \(\mathbf P^4\) satisfies
\[
       5S_{\mathcal O(1)}+v(W)\le(29/19)A .
\]
At nonzero \(B\), Euler's identities and the multiplicity bound
show that its Hessian is nonzero. Together with the \(x\)-square
this gives rank at least two. Translation on the \(x\)-axis and
an affine-linear graph shear reduce the other points to the origin.

On the degree-five graph closure the extraction of weights \(5,2,2\)
has discrepancy nine. In the complete graph restriction the normal
monomials are \(x^a(x^2+g)^jR_l\), \(a+j+l\le m\), with \(R_l\)
a homogeneous complement to \(g\,\mathbf C[B]_{l-5}\).
These latter spaces have eventual dimension five. The normalized
means of \(a,j,l\) are \(1/4\) each, giving radial mean \(17/4\).
In the uniformizing chart \(x=1\) the secondary ordinary blowup
at \(B=0\) has the same binary interval calculation as the ternary
quintic proof. More explicitly, fix \(J=5a+10j+2l\). For fixed
\(l\), the upper endpoint for \(j\) is
\(\lfloor(J-2l)/10\rfloor\) and the lower endpoint is
\(\max(0,(J+3l)/5-m)\). Eliminating consecutive powers of
\(1+g\) gives \(g^sR_l\) over the full range from zero to the
width. At \(D=l+5s\) these pieces reassemble as
\(g^{s_0}\mathbf C[B]_{D-5s_0}\) for the least permitted \(s_0\),
by the same endpoint substitution. The mean of \(s\) equals that
of \(\lfloor a/2\rfloor/m\) in the normalized count, hence
\(\bar D=1/4+5/8=7/8\). The radial ratio on this chart is \(7/16\).
On its exceptional binary line each marked mean is at most
\(2s_0+(D-5s_0)/2\le D/2\), with the same aggregate \(7/16\),
smaller than \(17/36\). On the ambient weighted extraction of weights (5,2,2,10), the ordinary signed radial ratio is (1+10/19). The punctured charts with (B
e0) are covered by the Hessian comparison through reassembly of the ordinary basis. In the remaining angular chart (x=1,B=0), divisor adjunction at the graph gives the moving threshold (36/17), while the fixed graph coefficient is (1/5). Since the binary residual calculation has aggregate coefficient (7/16<17/36), the moving graph contribution is below (17/9=4/(36/17)). Interpolation with the complete-system bound yields the signed inequality, because ((17/9-1)/(17/9-1/5)=10/19). Weighted adjunction finishes.

\paragraph{A repeated binary cubic and its quartic corrections.}
Suppose the cubic has been made \(x^2y\) in a constant normal
two-space, with constant complementary kernel \(K\) of ambient
codimension at most five, including the original line. A moving
double factor of base degree one has this form after one
preliminary constant cut. Indeed the double factor uses at most
two constant normals, the simple factor is fixed, and a general
cut equates those two normals up to scale while retaining the
simple factor independently. If the three normals are dependent,
choose an original constant frame of size two; the cubic is
of the required shape in pole-dependent coordinates on that
constant two-space. Its two transverse coordinates receive the
same first weight in the estimates below, so these latter linear
changes may be made at the poles.

If the quartic on \(K\) has no triple hyperplane component, two
general directions there give a nonzero binary quartic \(C\)
with no triple root, and weights \(1/3\) on \(x,y\), \(1/4\)
on these lights give \(x^2y+C\). If the quartic vanishes on
\(K\) together with both first normal derivatives, take three
general directions there and \(1/5\) weights instead. The
terminal Gauss estimate gives \(C\) a smooth ternary quintic.
Both graphs have signed denominator \(3/2\). The Hessian of
each of their separate blocks is nonzero away from its origin.
When just the cubic block is nonzero its Hessian has rank at
least two for \(x\ne0\), and for \(y\ne0,x=0\) fixing \(y\)
gives the square--binary-quartic or square--ternary-quintic test.
When just \(C\) has a nonzero argument, retain a direction having
nonzero quadratic coefficient there and the two independent
cubic variables. The square--nonpure-cubic test applies.
The radial light sums exceed one.

Suppose instead the quartic vanishes on \(K\), its first normal
derivative in the \(y\)-direction has nonzero nonpure cubic
restriction, and hence on two general \(K\)-directions gives a
binary cubic \(Q\) without a triple root. Allow any binary cubic
\(P\) for the other derivative, including zero. Weights \(1/3\)
on \(x,y\) and \(2/9\) on these two kernel directions give face
\[
                     x^2y+xP(U)+yQ(U).
\]
It has signed denominator \(3/2\). At \(x\ne0\) the \(x,y\)
Hessian minor has determinant \(-4x^2\). At \(x=0,y\ne0,U=0\)
fix \(y\) and use the first face consisting of its square and
the nonpure \(yQ(U)\). Elsewhere on \(x=0,U=p\ne0\), a nonzero
\(dQ(p)\) gives Hessian rank two since its \(y\)-row has zero
diagonal. If \(dQ(p)=0\) the root of \(Q\) is exactly double,
so \(J=\operatorname{Hess}Q(p)\ne0\). At \(y\ne0\) the rank-one
condition forces \(J=dP\otimes dP/(2y^2)\) by the Schur complement
of the \(xx\) entry. Choose \(\xi\) with \(dP(\xi)\ne0\) and
let \(e\) be its kernel lift with \(dx=-dP(\xi)/(2y)\).
The restriction of the cubic Taylor term to \(s\,\partial_y+t e\)
has coefficient
\[
        (dP(\xi)/(2y))^2+J(\xi,\xi)/2
             =(dP(\xi))^2/(2y^2)\ne0
\]
on \(st^2\), and has no term divisible by \(s^2\).
This gives a nonpure kernel cubic.

At \(x=y=0,U=p\ne0\), a nonzero row \(dP(p)\) already gives
rank at least two. When both gradients vanish, fix the radial
coordinate and let \(u\) be transverse to \(p\). The lowest
cubic on \(x,y,u\) is
\[
                    x^2y+(a x+b y)u^2,\qquad b\ne0.
\]
For \(a=0\) it is three lines forming a triangle. For \(a\ne0\)
it is integral, being linear in \(y\) with relatively prime
coefficients \(x^2+b u^2\) and \(axu^2\). Its only singularity
is the ordinary node at \([0:1:0]\): at \(u\ne0\), the three
partial equations give \(x^2=-bu^2\), \(by=-ax\), and then
\(3au^2=0\), impossible. The ternary log-canonical cubic
test completes this last case. The first weight sum is \(10/9\),
so all the line estimates in this paragraph give \(60/7\).

The same mixed test applies if \(Q\) is a nonzero pure cube,
provided \(P\) does not vanish identically on the hyperplane
supporting \(Q\). A general binary restriction then gives
\(Q=v^3\), \(P(0,t)=\gamma t^3\) on its root direction,
\(\gamma\ne0\). At \(U=p\ne0\) on this root \(dP(p)\ne0\)
and \(\operatorname{Hess}Q(p)=0\), so the Schur complement
at \(x=0,y\ne0\) has positive rank. At \(x=y=0\) the nonzero
\(x,U\) mixed row suffices. The only other new point is
\(x=0,y\ne0,U=0\). Fix \(y\) there, translate the point, and
use weights \(1/2,1/3,1/6\) on \(x,v,t\). They give
\(a x^2+b v^3+\gamma x t^3\), with all three
coefficients nonzero, precisely the coupled square--cube test
with the quartic term in the square variable. Its radial light
sum is one and its signed denominator is \(3/2\).

\paragraph{A fixed binary double factor with terminal radical.}
Suppose the retained cubic at the poles is \(x^2y\) up to nonzero
scalar, with both coordinates constant up to scalar along the base.
Let their common kernel \(K\), including the original line, have
ambient codimension at most five, and let its first order be five.
Then the common line threshold is at least \(76/9>312/37\).
For a moving double root whose two coefficients are independent
modulo the constant simple root, one general constant preliminary
cut achieves these coordinates and \(K\) has ambient codimension
five. If all three coefficient forms instead span a two-space,
\(K\) has codimension four and its terminal test applies directly.
The terminal quintic on \(K\) has singular codimension at least
three, so is integral. Write \(xP+yQ\) for the quartic terms
linear in the two normals, with cubics \(P,Q\) on \(K\).
The existing binary-double tests settle the vanishing of both,
every nonpure \(Q\), and the pure \(Q\) whose root does not
divide \(P\). All cubic factors discussed below are in the
geometric generic light space.

Suppose \(Q=H^3\) up to scalar and \(H\mid P\).
The root \(H\) is constant by the binary base degree one.
Choose a general direction \(u\) in \(H=0\) and retain a
complementary direction measured by \(H\); the quintic on \(u\)
is nonzero by integrality. Let \(s\) be the multiplicity of
\(H\) in \(P\), assigning \(s=3\) if \(P=0\).
Take \(w_u=1/5\),
\[
   w_x=\frac{3(s+2)}{5(2s+3)},\qquad
   w_y=1-2w_x,\qquad w_H=2w_x/3 .
\]
Except for the lower term specified below, the face is
\[
         x^2y+yH^3+a xH^s u^{3-s}+c u^5,\qquad c\ne0,
\]
where only \(a\) may vanish when \(P=0\).
Its light sum exceeds one. At \(x\ne0\) use nonzero \(F_{xy}\)
and zero \(F_{yy}\); at \(x=0,H\ne0\) use nonzero \(F_{yH}\).
At \(x=H=0,u\ne0\), a nonzero \(y\) gives independent
\(xx,uu\) entries, while at \(y=0\) use the square in \(u\)
and the nonpure kernel cubic \(x^2y\) unless the additional
mixed terms already give rank two. Indeed for \(s=1\) the
\(xH\)-entry is nonzero; otherwise the first nonzero additional
mixed term has order at least three and vanishes upon setting
\(H=0\). On the remaining \(y\)-axis the face after fixing
\(y\) is a square, \(H^3\), and \(u^5\); the mixed monomial
has larger weight for \((w_x,w_H,w_u)=(1/2,1/3,1/5)\).
At \(s=1\) a term \(b y^2u^2\) has lower weight and
\(y^2Hu\) has weight one. If either occurs, use weights \((12,10,8,7;34)\)
on \(x,y,H,u,Z\), giving
\[
                 x^2y+yH^3+a xHu^2+b y^2u^2,\qquad a\ne0 .
\]
Now \(y^2Hu\) and \(u^5\) both have weight \(35/34\).
For \(b\ne0\), at \(u\ne0\) rank at most one relative to the nonzero \(yy\)
entry implies
\[
 x^2=byu^2,\quad ab u^4=3xH^2,\quad
 y=3H^3/(4bu^2),\quad aHu^2=2xy ,
\]
the latter three being inconsistent since \(xH\ne0\).
For \(u=0\), a nonzero \(x\) or \(H\) gives a nonzero
\(y\)-row with zero diagonal; on the \(y\)-axis the \(xx\)
and \(uu\) entries are independent.
If \(b=0\), the zero \(yy\)-entry forces the whole \(y\)-row to
vanish at rank at most one, hence \(x=H=0\). There the \(xH\)-entry
is \(au^2\) and \(HH=0\), leaving only the \(y\)-axis.
Fixing its nonzero coordinate gives, after rescaling,
\(x^2+H^3+A xHu^2\) with \(A\ne0\). Its weights are
\((1/2,1/3,1/12)\), of sum \(11/12\). At \(u\ne0\) rank one
relative to \(xx=2\) would require
\(6H=A^2u^4/2,\ 2Axu=A^2Hu^3,\ 2AxH=2A^2H^2u^2\),
which are inconsistent. At \(u=0,H\ne0\) the \(xx,HH\) entries
give rank two; at \(u=H=0,x\ne0\) the kernel cubic
\((dH)^3+A x\,dH\,(du)^2\) is squarefree on its two variables.

It remains that \(Q=0,P\ne0\). On \(x=0\) the quartic is
divisible by \(y^2\). If it is zero, this hyperplane has ambient
codimension at most four (in this application there is at most one
preliminary cut) and the terminal test suffices. Otherwise let
\(s=2,3,4\) be its multiplicity along \(y=0\) and \(B\)
its first coefficient of degree \(4-s\) on \(K\).
Take two general kernel lights \(U\). Up to nonzero scalars the
three faces and their weights \((w_x,w_y,w_U)\) are
\[
\begin{array}{c|c|c}
s& F & (w_x,w_y,w_U)\\ \hline
2&x^2y+xP(U)+y^2 B(U)&(5,4,3)/14\\
3&x^2y+xP(U)+y^3 B(U)&(7,5,4)/19\\
4&x^2y+xP(U)+b y^4&(3/8,1/4,5/24).
\end{array}
\]
For forms gaining a common high root the modified choices with
one extra terminal direction are given below. General binary
restrictions preserve divisorial multiplicities and avoid
pairwise intersections, so a triple root of \(P\) means that
the full \(P\) is \(aH^3\) for a constant \(H\). A common
double quadratic root at that triple point then means that the
full quadratic is \(bH^2\).

For \(s=2\), suppose first \(B(U)\ne0\) and the Hessian has
rank at most one. Its nonzero \(yy\)-entry implies
\(x^2=yB\) and \(dP=2xy\,dB/B\), hence \(3P=4xy\).
The radial \(U,U\) entry would be both
\(6xP+2y^2 B=10y^2B\) and \((4yB)^2/(2B)=8y^2 B\).
Thus \(x=y=0\), \(dP=0\). A double root of \(P\) gives a
nonpure kernel cubic \(c x v^2\) with \(v\) transverse to
the binary root. At a triple root the centered first weighted
face is \(y^2+x^2 y+a x v^3\), with weights \((1/2,1/4,1/4)\)
on \(y,x,v\) and \(a\ne0\). Its \(yy\)-entry is nonzero;
at \(v\ne0\) the \(xv\) and zero \(yv\) entries exclude rank
one. At \(v=0\) a rank-one point has \(y=x^2\); the kernel
direction \((dy,dx)=(-x,1)\) gives a nonzero cube at \(x\ne0\),
independent of the cube from \(ax v^3\).

At \(B(U)=0,U\ne0\), rank at most one forces \(x=0\),
\(y\,dB=0\) from the \(y\)-row. When \(y\ne0\) the root
of \(B\) is double. Its nonzero Hessian and the \(xx\)-entry
force in rank one \(y^2\operatorname{Hess}B=dP\otimes dP/(2y)\).
Take a kernel vector over a direction \(e\) with \(dP(e)\ne0\),
so \(dx=-dP(e)/(2y)\). On its span with \(\partial_y\) the
cubic has nonzero coefficient
\((dx)^2+y\operatorname{Hess}B(e,e)=3(dP(e))^2/(4y^2)\)
of the term linear in the \(y\)-variation and quadratic in
the other direction, with no term quadratic in the \(y\)-variation.
When \(y=0\), any nonzero \(dP\) gives rank two. Otherwise fix
the nonzero radial coordinate, write \(v\) for the transverse
direction, and put \(k=\operatorname{mult}P\), \(m=\operatorname{mult}B\).
The four possibilities are as follows.
\begin{itemize}
\item For \((k,m)=(2,1)\) the ordinary cubic is
\(x^2y+a xv^2+b y^2v\), smooth in \(\mathbf P^2\).
Its three exponent vectors are independent, ruling out a torus
singularity; at coordinate points one partial is nonzero.
\item For \((3,1)\) the weighted face is
\(x^2y+a xv^3+b y^2v\), of light sum \(12/13\), the
punctured-rank-two tangent-conic face proved above.
\item For \((2,2)\) the face is
\(x^2y+a xv^2+b y^2v^2\), with light weights
\((2/5,1/5,3/10)\) and sum \(9/10\). Its Hessian has rank
at least two on the punctured space. Indeed \(v\ne0\)
would at rank one give \(x^2=byv^2\), \(av^2=2xy\), and
\(ax=3by^2\) from the \(xx,xv,vv\) Schur equations, contrary
to \(a\ne0\). On \(v=0\) use the \(xy\) entry if \(x\ne0\)
or the \(xx,vv\) entries on the \(y\)-axis.
\item For \((3,2)\) the full \(P,B\) are \(aH^3,bH^2\).
Retain a new general \(u\in(H=0)\) with quintic nonzero,
of weight \(1/5\), keeping weight \(3/14\) on \(H\).
The face becomes \(x^2y+a xH^3+b y^2H^2+c u^5\).
On the punctured \((x,y,H)\)-space its first block has Hessian
rank two by the preceding radial Schur computation at \(H\ne0\)
and by the \(xy\) or \(xx,HH\) entries at \(H=0\).
On the complementary \(u\)-axis use a square and \(x^2y\).
\end{itemize}
At \(U=0\) in the general \(s=2\) face the \(xy\) mixed
derivative gives rank two if \(x\ne0\), and a nonzero
\(y\) gives a square and the nonzero independent binary
quadratic \(y^2B\).

For \(s=3\) write \(B=L\) linear. At \(y\ne0\) a rank-one
Hessian relative to \(xx=2y\) would give
\(x\ne0\), \(x^2=3y^2L\),
\(dP=3y^3dL/x\), hence \(P=xy/3\).
The radial \(U,U\) entry would be simultaneously
\(6xP=2x^2y\) and \((3P)^2/(2y)=x^2y/2\), impossible.
At \(y=0,x\ne0\) use the \(xy\) entry. At
\(x=y=0,U\ne0\) only a multiple root of \(P\) needs testing.
A double root with \(L(U)\ne0\) gives cubic
\(x^2y+c y^3+a xv^2\), log canonical by the squarefree
binary cubic plus a linear form times a square test.
If \(L(U)=0\), use \(x^2y+a xv^2+b y^3v\) of light sum
\(12/13\), another permutation of the punctured-rank-two
tangent-conic face.
For a triple root \(P=aH^3\) use a new \(u\in(H=0)\)
with nonzero quintic and weight \(1/5\). If \(L\) does not
vanish identically there retain
\(x^2y+a xH^3+b y^3u+c u^5\) with weights
\(11/30,4/15,19/90,1/5\). At \(y u\ne0\) the independent
\(xx,uu\) entries give rank two. At \(y\ne0,u=0\) the
\(yu\) mixed entry and zero \(uu\) do. At \(y=0\) use the
\(xy\) entry when \(x\ne0\), the \(xH\) entry if \(x=y=0,H\ne0\),
or else the radial square on the \(u\)-axis and \(x^2y\).
If \(L\) vanishes identically on \(H=0\), keep the original
\((7,5,4)/19\) weights on \(x,y,H\), giving instead
\(x^2y+a xH^3+b y^3H+c u^5\). The one-\(H\) Schur
computation already made proves the first block has rank at
least two at each punctured point. On the \(u\)-axis add
the radial square to the nonpure cubic.

For \(s=4\) the \(y\ne0\) rank-one condition gives
\(x^2=6by^3\), then \(dP=0\) and
\(\operatorname{Hess}P=0\) by the mixed and \(U\)-entries.
On the kernel of the \(x,y\) block the direction
\((dx,dy)=(-x/y,1)\) has cubic coefficient \(10by\ne0\).
Since the nonzero homogeneous cubic \(P\) has zero first and
second derivatives at this \(U\), a general \(U\)-variation
supplies an independent nonzero cube. At \(y=0,x\ne0\) there
is a rank-two mixed block. At \(x=y=0,U\ne0\) a simple or
nonroot is handled by \(dP\ne0\); a double root gives
\(x^2y+a xv^2+b y^4\), the punctured-rank-two face of
light sum \(15/16\). A triple root has \(P=aH^3\); add
the previous \(u\)-direction to obtain the disjoint fifth
power beside \(x^2y+a xH^3+b y^4\). All punctured points
of the latter block were tested above, and the \(u\)-axis
again gives a square and the nonpure cubic.

All supplementary directions above exist with a nonzero terminal
quintic because the quintic on \(K\) is integral and has no
linear factor. For the modified weights with the extra
\(u\)-direction, terms containing a pair of normals including
\(x\) have weight greater than one, as do all omitted higher
\(y\)-terms. Thus the listed initial equations are complete.
The smallest required light-to-graph ratio is \(9/10\), giving
signed denominator at most \(29/19\); averaging the \(1/5,4/5\)
groups gives \(12/(1/5+(4/5)(29/19))=76/9\).

\paragraph{Every moving double root in the cubic pencil.}
Suppose the cubic is \(L_t^2M\), with \(L_t\) a primitive
moving normal of binary base degree one and \(M\) fixed.
Every such line on a smooth quintic elevenfold has threshold
at least \(76/9\). Write \(K\) for the common constant kernel
of the two coefficients of \(L_t\) and \(M\).
It has ambient codimension four or five. A general normal
two-space restricts the cubic to \(x^2y\) at the general pole,
with \(x=L_t,y=M\) after rescaling. In the independent case
one general linear relation between the two coefficients of
\(L_t\) makes \(x,y\) fixed normals on that section and leaves
\(K\) unchanged. In the dependent case all tests below retain
equal weights on these two normals so changes at the separate
poles suffice. The only test requiring distinct first weights
is the terminal-radical lemma, which in this dependent case
already follows from the original terminal test at codimension four.

Write \(G_4\) for the light quartic. If \(G_4|_K=0\), apply
the preceding terminal-radical calculation. If \(G_4|_K\)
has no triple hyperplane component, apply the disjoint
binary-double and binary-quartic calculation. Suppose therefore
that \(G_4|_K=V^3R\) with \(V\) constant, by its binary
base degree one. Let \(K'=(V=0)\) inside \(K\). This constant
space has original ambient codimension at most six, and so
cannot be contained in the smooth elevenfold. Its order is
five, with nonzero base-independent terminal quintic \(C\).

If some quartic derivative in the original normal directions
to \(K\) within the light space is nonzero on \(K'\), choose
the normal two-space so at least one persists on it. In the
independent case its kernel in the span of the two coefficient
directions of \(L_t\) is a freely chosen line, while the
\(M\)-direction is also available. Thus a general such choice
detects any nonzero derivative. Take one general direction
\(u\) in \(K'\) and one complementary direction \(V\) in
\(K\). Give \(x,y\) weight \(1/3\), \(u\) weight \(2/9\).
Give \(V\) weight \(7/27\) when \(R|_{K'}\ne0\),
and \(1/4\) otherwise. The complete faces are respectively
\[
   x^2y+(ax+by)u^3+e V^3u,\qquad
   x^2y+(ax+by)u^3+e V^4,\qquad e\ne0,\quad(a,b)\ne(0,0).
\]
Both light sums exceed one. Write \(J\) for the last monomial.
At \(x\ne0\) the \(xy\)-entry and zero \(yy\) give rank
two. At \(x=0,u\ne0\), if \(b\ne0\) the \(yu\)-entry is
nonzero with \(yy=0\). If \(b=0,y=0\), use the \(xu\)
entry with \(xx=0\). If \(b=0,y\ne0,V=0\), use the same
\(xu\) entry, this time with \(uu=0\).
For \(b=0,y\ne0,uV\ne0\), the \(xx,VV\) entries are
independent. At \(x=u=0,V\ne0\) the \(V^3u\) face has
rank two from its \(Vu\) entry and zero \(VV\). The fourth
power instead supplies a square; if \(y\ne0\) there is a
second square from \(x\), and otherwise the kernel cubic
contains \(x^2y\).

On the remaining \(y\)-axis fix its nonzero value by a linear
flag. For \(b\ne0\) the retained square and cube with \(J\)
give the square--cube--quartic or the coupled \(uV^3\) test;
the term \(axu^3\) has larger weight. For \(b=0\) the
secondary graph is, up to units and constant rescaling,
\[
                    p^2+p u^3+eV^4
                  \quad\hbox{or}\quad p^2+p u^3+eV^3u .
\]
Use weights \(1/2,1/6,1/4\) or \(1/2,1/6,5/18\)
on \(p,u,V\); the sums are \(11/12\) and \(17/18\).
For the first face, when \(V\ne0\) its square adds to
the \(pp\)-entry. At \(V=0\) a rank-one punctured
\((p,u)\) block at \(u\ne0\) requires \(p=3u^3/4\);
on its kernel \(dp=-(3/2)u^2du\) the cubic coefficient
is \(p-(9/2)u^3\ne0\). With the independent \(V^4\)
this is the square--cube--quartic face. For \(u=0,p\ne0\)
use the cube with coefficient \(p\) and the same fourth power.
In the second face \(V\ne0\) always gives rank at least two,
since the \(pV\)-entry is zero while \(VV\ne0\) if \(u\ne0\),
and \(Vu\ne0,VV=0\) if \(u=0\). At \(V=0,u\ne0\)
the same kernel cube is independent of the \(V^3\) term.
At \(V=u=0,p\ne0\) its square, the \(pu^3\) cube and
\(eV^3u\) give the known coupled test. These estimates cover
all punctured points in this first derivative case.

Suppose now all the original quartic normal derivatives vanish
on \(K'\). If \(C\) has no triple hyperplane component, choose
two general directions in \(K'\). Their binary quintic has
no root of multiplicity greater than two. With both cubic
normals of weight \(1/3\) and kernel lights of weight \(1/5\)
the face is \(x^2y+C(u,v)\). Its signed denominator is at
most \(29/19\): off the axes use the two independent Hessian
blocks; with only the quintic present use a square and the
nonpure cubic; on the \(y\)-axis fix \(y\) and use the
square-plus-binary-quintic bound. The \(x\ne0\) block
already has rank two. The light sum exceeds one.

If \(C\) has a triple hyperplane \(T\), it is fixed since
\(C\) is independent of the binary base. In ambient codimension
at most five for \(K'\) the terminal Gauss bound makes \(C\)
integral, so only codimension six occurs here, with three independent
original light normals to \(K\). Put \(W=(T=0)\) in
\(K'\), an actual projective five-plane in the hypersurface.
It has four light coordinates \(U\) in addition to the two
binary coordinates on the original line. On \(W\) the cubic
and all its partials vanish. The quartic also has zero partials:
along \(K\) this follows from \(V^3R\), and in the three
remaining light directions it is the hypothesis of this case.
The partial derivative in \(T\) of the terminal quintic
vanishes by its triple factor. Thus the four remaining
nonheavy derivatives, in the three original normals to \(K\)
and in \(V\), restrict to quartics in \(U\) alone. They
have no common projective zero in \(\mathbf P^3_U\).
Such a zero would make every nonheavy derivative vanish on
its joining plane with the original line. The other two
heavy derivatives are quartics on that \(\mathbf P^2\)
and have a common zero there, contradicting smoothness of \(X\).
Consequently their general linear combination \(G(U)\)
is a smooth quartic surface by Bertini.

Retain one general normal direction in this four-space of
normals, called \(z\), and the four \(U\)-directions on
\(W\). Its cubic coefficient \(a\) is nonzero at general
poles, since the direction in the three original normals can
be general. Give \(z\) weight \(1/3\) and \(U\) weight
\(1/6\), of total light weight one. The complete face is
\[
                   F=a z^3+z^2 B(U)+zG(U),
\]
with \(B\) an arbitrary quadratic: the missing zero and first
normal quartic coefficients vanish on \(W\), the zero normal
quintic coefficient vanishes, and all later terms have weight
greater than one.

At \(z=0,U\ne0\) the nonzero \(z,U\) row \(dG\) and zero
\(U,U\) block give Hessian rank at least two. At \(z\ne0\)
the Hessian cannot vanish identically. For \(U=0\) use the
\(zz\) entry. Otherwise the \(z,U\) row and contraction of
the \(U,U\)-block with the radial vector would give
\(2z\,dB+dG=0\) and \(z^2dB+3z\,dG=0\), contrary to
\(dG\ne0\). Consider a rank-one Hessian with kernel \(E\)
in the five tangent variables. A nonpure cubic on \(E\)
gives the usual square-plus-nonpure-cubic bound. If that cubic
is a nonzero pure cube, its zero hyperplane inside \(E\)
is a three-space projecting injectively onto the four \(U\)
variables: a vertical \(z\)-direction in \(E\) has nonzero
cubic \(a(dz)^3\). Its projected projective plane cannot
lie in the smooth quartic surface. Choose a general direction
there with \(G(dU)\ne0\). If \(dz\ne0\), the top Taylor
coefficient is \(dz\,G(dU)\) of degree five; otherwise its
degree-four coefficient is \(zG(dU)\ne0\). Thus the first
pure power on this direction has order four or five. A square
variable and the pure cube give the mixed E6, coupled-square,
coupled-cube, or E8 alternatives as in the moving reduced-cubic
proof, all with denominator at most \(3/2\).

Finally, if the cubic on \(E\) is zero, \(E\cap(dz=0)\)
projects onto a subspace of dimension at least three in \(U\).
Choose a projective plane in it. Its section of \(G\) is
reduced: a repeated curve would map to this one tangent
hyperplane by the finite Gauss morphism of the smooth quartic
surface. Two general directions in the three-space therefore
give a squarefree binary quartic. After translation their
quartic is the nonzero scalar \(zG(dU)\). With a transverse
square variable \(v\) the resulting face is
\(v^2+vB'_2(dU)+zG(dU)\), covered by the mixed
square--binary-quartic estimate. This finishes the punctured test.
The terminal-radical case had smallest allowed light ratio
\(9/10\); all new ratios are at least this large. The line
denominator is thus at most \(1/5+(4/5)(29/19)=27/19\),
proving the stated mark.

\subsection{The quintic elevenfold}
\begin{theorem}\label{thm:nov28-quintic-eleven}
Every smooth complex quintic elevenfold \(X\subset\mathbf P^{12}\)
is uniformly K-stable.
\end{theorem}
\begin{proof}
Here \(-K_X=8H\). The following are common lower bounds for
\(\delta_C(X,H)\) for the indicated lines \(C\). For a polar pencil,
\(h\) is the first nonzero normal degree on its constant kernel.
\[
\begin{array}{lc}
\text{line type}&\text{lower bound}\\ \hline
\text{polar rank at least three}&60/7\\
\text{polar rank two},\ h=2,4,5&60/7\\
h=3,\ \text{geometrically reduced cubic}&500/59\\
h=3,\ \text{pure cube}&348/41\\
h=3,\ L^2M,\ L\text{ fixed},\ M\ne L\text{ up to scalar}&620/73\\
h=3,\ L^2M,\ L\text{ moving}&76/9
\end{array}
\]
The rank cannot be one since the quartic normal polars have no
common zero on \(C\). The rank at least three case is
Section~\ref{sec:nov28-quintic-eleven-nets}. For a pencil,
the order cannot exceed five: its full kernel would otherwise
give a linear space of codimension two contained in \(X\),
contrary to smoothness. The high-order, square and quartic
calculations above give the second row. For \(h=3\) every nonreduced
cubic over the geometric generic field is either a pure cube or
\(L^2M\) with distinct factors. The different multiplicities make
the latter factors descend to the base function field. The reduced
case is Section~\ref{sec:nov28-quintic-eleven-moving-cubic}; the
three remaining rows are the pure-cube, fixed-double-root and
moving-double-root calculations above.

Set \(\mathcal E=\{C:h=3,\ L^2M,\ L\ {\rm moving}\}\).
The first five bounds have minimum \(500/59\), which exceeds
\(144/17\) by \(4/1003\). Every line has threshold at least
\(76/9>8\). If \(\delta(X,H)\le8\), a positive-dimensional
minimizing center is excluded by the line bound and the
center estimate of Section~\ref{sec:nov27-third-descent}.
At a minimizing point \(p\), the length-one argument of
Section~\ref{sec:nov28-two-hyperplanes} forces a line in
\(\mathcal E\) through \(p\), by taking the first interpolation
coefficient strictly between \(144/17\) and the common bound
on the other lines and sufficiently close to \(144/17\).
On such a line the primitive moving double root has binary
degree one and the remaining factor is fixed: the coefficient
degree two of the cubic is already exhausted. Thus this cubic
is nonzero as a form at every point, including \(p\).
The tangent hyperplane section at \(p\) has multiplicity two
or three. The plane-quintic threshold computation above,
including the pure tangent, gives \(S_H(v)\ge27/16\) there,
with \(v(\mathfrak m_p)=1\). The second comparison at this
same point consequently requires only a common line mark
strictly greater than \(312/37\). The margin of the weakest
row is \(76/9-312/37=4/333>0\). The two-pass pointed argument
therefore contradicts \(\delta(X,H)\le8\) and proves
\(\delta(X,-K_X)>1\).
\end{proof}

% ===== Explicit families in larger dimensions =====
\chapter{Explicit families in larger dimensions}\label{ch:families}
\begin{theorem}\label{families:main}
Fix \(m\in\{4,5,6,7,8\}\) and let \(X\) be a smooth quartic of dimension \(n\ge m\). If, in some projective coordinates, the equation of \(X\) is even separately in \(n-m\) coordinates, then \(X\) is K-stable.
\end{theorem}

More explicitly, set \(k=n-m\), use \(m+2\) coordinates \(x_0,\ldots,x_{m+1}\) and \(k\) coordinates \(y_1,\ldots,y_k\). The theorem covers every smooth equation
\begin{equation}\label{families:family}
        f_4(x)+\sum_{i=1}^k q_i(x)y_i^2+
          \sum_{1\le i\le j\le k} a_{ij}y_i^2y_j^2=0 ,
\end{equation}
where the \(q_i\) are quadrics and the \(a_{ij}\) are constants. Taking \(m=8\) uses Theorem~\ref{thm:nov28-quartic-eightfold}. The choice \(m=4\) uses the published fourfold result.

\section{Two double covers}
\begin{lemma}\label{families:one}
Let \(n\ge4\) and let \(X=(y^4+q(x)y^2+f(x)=0)\subset\mathbf{P}^{n+1}\) be smooth, with \(q,f\) of degrees two and four. Both \(C_0=(f=0)\) and \(C_1=(q^2-4f=0)\) in \(\mathbf{P}^n_x\) are smooth quartics. If both are K-polystable, \(X\) is K-stable.
\end{lemma}
\begin{proof}
At a point with \(y=0\) the \(y\)-derivative vanishes, so smoothness gives \(df\ne0\) on \(C_0\). If \(x\in C_1\) and \(q(x)=0\), then \(d(q^2-4f)=-4df\ne0\). Otherwise take \(y^2=-q(x)/2\) on an affine chart. The \(y\)-derivative again vanishes and the remaining gradient equals \(d(f-q^2/4)\), so this gradient is nonzero. This also excludes an identically zero discriminant.

Put \(w=y^2\). In \(\mathbf{P}(1^{n+1},2)\) the hypersurface
\[
       Y=(w^2+qw+f=0)
\]
misses the weighted vertex. Completing the square identifies \(\pi:Y\to\mathbf{P}^n\) as the smooth double cover branched along \(C_1\). Its K-polystability follows from the cyclic-cover theorem of Liu and Zhu, since \(4\le n\) and \(C_1\) is K-polystable \cite[Thm.~5.1]{ar2001.10557}. With \(A=\pi^*\mathcal{O}(1)\), the ramification formula gives \(-K_Y=(n-1)A\). The smooth prime divisor \(D=(w=0)\) is isomorphic to \(C_0\) and is in \(|2A|\).

Use Zhang and Zhou's theorem for a smooth K-polystable divisor \(D\sim_{\mathbf{Q}}\lambda(-K_Y)\) on a K-polystable Fano manifold \cite[Thm.~1.7]{ar2003.06839}. Here \(\lambda=2/(n-1)<1\). Their upper endpoint for log K-polystability is
\[
    1-\frac{\lambda^{-1}-1}{n}=\frac{n+3}{2n}>\frac12 .
\]
Thus \((Y,D/2)\) is K-polystable. The map \(X\to Y\) is the Galois double cover defined by \(y^2=w\), branched along \(D\), with \(K_X=\pi_{X/Y}^*(K_Y+D/2)\). Liu--Zhu \cite[Thm.~1.2(2)]{ar2001.10557} give K-polystability upstairs. A smooth Fano quartic has finite automorphism group \cite{MM}, hence \(X\) is K-stable.
\end{proof}
\begin{proof}[Proof of Theorem~\ref{families:main}]
At each step the coefficient of \(y_k^4\) is nonzero. Otherwise the \(y_k\)-coordinate point is singular, because evenness excludes the terms \(y_k^3x_j\) and \(y_k^3y_i\). Normalize the leading coefficient to one and apply Lemma~\ref{families:one}. Both new quartic equations in dimension \(n-1\) remain even in the other \(k-1\) coordinates, and the lemma guarantees their smoothness. The induction ends in dimension \(m\). For \(m=4\) apply \cite[Theorem~4.21]{ar2003.13788}; for \(m=5,6,7\) use respectively Chapters~\ref{ch:fivefold}, \ref{ch:sixfold} and~\ref{ch:boundary}. For \(m=8\) use Theorem~\ref{thm:nov28-quartic-eightfold}.
\end{proof}

\begin{corollary}\label{families:binary}
If \(F(x)=0\) in \(\mathbf{P}^{m-1}\) is a smooth K-polystable quartic, \(m\ge5\), and \(g(y,z)\) is a binary quartic with four distinct roots, then \(F(x)+g(y,z)=0\) in \(\mathbf{P}^{m+1}\) is smooth and K-stable.
\end{corollary}
\begin{proof}
Send the four roots to \(0,\infty,1,t\). The involution \(u\mapsto t/u\) exchanges the roots in pairs, and neither fixed point is a root. In eigenvector coordinates this involution acts by a sign change. The invariant equation has nonzero pure fourth powers at both ends, so after rescaling it is \(y^4+\eta y^2z^2+z^4\), with \(\eta^2\ne4\). In Lemma~\ref{families:one}, the two equations become \(F+z^4\) and \(-4F+(\eta^2-4)z^4\). Each is a cyclic fourth-root cover of \(\mathbf{P}^{m-1}\) branched over \(F=0\), hence K-polystable by \cite[Thm.~5.1]{ar2001.10557}. The separated variables and the smoothness of each projective equation give smoothness of the sum by its partial derivatives. Lemma~\ref{families:one} finishes the proof.
\end{proof}

\paragraph{Why a general projection does not give this induction.}
Choose an exterior point and project \(X\to\mathbf{P}^n\) from it. After translating the fiber variable the equation becomes \(y^4+q_2y^2+q_3y+q_4=0\). A tower of two quadratic function-field extensions would give a nontrivial deck involution of the top step fixing \(\mathbf{C}(\mathbf{P}^n)\). It extends to \(X\), the normalization of \(\mathbf{P}^n\) in \(\mathbf{C}(X)\), and preserves \(H\), the pullback of the base hyperplane. It therefore acts linearly in the embedding by \(|H|\). With the natural linearization it fixes the \(x_i\) and sends \(y\) to \(-y+\ell(x)\). Since \(y\) has trace zero in the depressed equation, \(\ell=0\). Invariance now forces \(q_3=0\). Conversely \(q_3=0\) gives the tower of Lemma~\ref{families:one}.

For a fixed reflection the space of even quartics loses \(\binom{n+3}{3}+(n+1)\) independent coefficients, those of \(y q_3(x)\) and \(y^3\ell(x)\). The reflection itself varies in \(2(n+1)\) dimensions (its center and its fixed hyperplane). Hence even after varying coordinates the reflection locus has codimension at least \(\binom{n+3}{3}-(n+1)>0\). Smoothness cannot supply the missing symmetry.

\paragraph{Pure powers in any of the degrees.}
If \(Y=(f=0)\subset\mathbf P^{m+1}\) is a smooth K-polystable Fano
hypersurface in the range of the introduction, adjoining any number of
pure powers of the same degree \(d\) gives a K-stable hypersurface.
At a single step the new equation is \(f+y^d=0\); its derivatives
give smoothness, and projection defines the cyclic degree-\(d\) cover
of \(\mathbf P^{m+1}\) branched over \(Y\). Apply
\cite[Theorem~5.1]{ar2001.10557} and repeat. Finite automorphisms give K-stability.

\section{Sums of forms in disjoint groups of variables}
\label{sec:cyclic-block-hypersurfaces}
The next result gives a family in every degree, including prime degrees.
Its proof uses published results on special degenerations and local
stability thresholds. The constituent hypersurfaces may be of general
type.

\begin{theorem}\label{thm:cyclic-block-hypersurfaces}
Let $d\ge3$ and $d\le n+1$. Suppose a smooth hypersurface of degree
$d$ in $\mathbf P^{n+1}$ has an equation
\[
                     \sum_{i=1}^s f_i(x^{(i)})=0,
       \qquad s\le d,\qquad 1\le a_i:=\dim x^{(i)}\le d-1,
       \qquad \sum_i a_i=n+2,
\]
where the variable sets are disjoint. Then it is uniformly K-stable.
In particular the theorem gives explicit families through dimension
$d(d-1)-2$. For $d=5$, it includes every smooth sum of at most five
quintics in at most four variables each, through dimension $18$.
\end{theorem}

\begin{lemma}\label{lem:cyclic-delta-fixed}
Let $X$ be a smooth complex Fano variety with finite automorphism group,
and let $g$ be an automorphism of finite order. If
$\delta_p(X,-K_X)>1$ at every fixed point of $g$, then $X$ is uniformly
K-stable. Here the local threshold uses centers containing $p$.
\end{lemma}
\begin{proof}
By Datar--Sz\'ekelyhidi, Theorem~1 \cite{blockDS}, it suffices to test
special degenerations equivariant under the compact cyclic group
$\langle g\rangle$. For any nontrivial such degeneration
$\mathcal X\to\mathbf A^1$ its central fiber $V$ is a klt
$\mathbf Q$-Fano variety. The restriction of $\operatorname{ord}_V$ to
$\mathbf C(X)$ is a positive multiple of a divisorial valuation $v$,
and the sign of the Futaki invariant is the sign of
$A_X(v)-S_{-K_X}(v)$
\cite[Definition~2.25 and Lemma~2.26]{ar2102.09405}. A product
configuration is trivial here since $\operatorname{Aut}(X)$ is finite.

The equivariant rational map $\mathcal X\dashrightarrow X$ from the
punctured trivialization is defined at the generic point of $V$ by
properness of $X$. Its restriction to $V$ has image closure
$Z=c_X(v)$ by the valuative definition of the center. Resolve the
graph of $V\dashrightarrow Z$ equivariantly to obtain a smooth
projective $\widetilde V$ with a morphism to $Z$. Klt singularities
are rational and $-K_V$ is ample, so Kawamata--Viehweg vanishing gives
$H^j(\widetilde V,\mathcal O_{\widetilde V})=0$ for $j>0$.
The trace of $g$ on $H^0(\widetilde V,\mathcal O_{\widetilde V})$
is one. Its holomorphic Lefschetz number is therefore one, and the
holomorphic Lefschetz fixed point theorem gives a fixed point on
$\widetilde V$ \cite{blockAB}. Its image is a fixed point $p\in Z$.
The local hypothesis gives $A_X(v)>S_{-K_X}(v)$. Thus every
nontrivial equivariant special degeneration has positive Futaki
invariant. Datar--Sz\'ekelyhidi give the K\"ahler--Einstein metric,
and finiteness of automorphisms gives uniform K-stability by the
equivalences in the introduction.
\end{proof}

\begin{proof}[Proof of Theorem~\ref{thm:cyclic-block-hypersurfaces}]
We induct on the number of variables. Smoothness implies that each
block has no nonzero common zero of its first partials: by Euler's
identity this would give a singular point of $X$ after setting every
other block to zero. Every nonempty subcollection of the blocks has
the same property. For a singleton, its coefficient is nonzero.

Let $g$ multiply the different blocks by distinct $d$th roots of
unity. Each fixed point of this projective transformation is supported
on just one block. A singleton gives no fixed point on $X$. Suppose
that $p$ is supported on a block of $a\ge2$ variables, and write $b$
for the number of remaining variables. The Fano inequality and
$a\le d-1$ imply $b\ge2$. Write the equation as $f(x)+G(y)=0$.
The projective degree-$d$ hypersurface $(G=0)\subset\mathbf P^{b-1}$
is smooth and still has the stated form. Whenever it is Fano, it is
K-semistable by induction. The induction only invokes a smaller
Fano hypersurface with $d\le b-1$. If no such hypersurface occurs,
no induction hypothesis is needed.

Choose a general projective line through $p$ in the active
$\mathbf P^{a-1}$. Its restriction $f_2$ has simple roots by Bertini:
the line is transverse at $p$ and meets the smooth hypersurface
transversely elsewhere. For $a=2$ use the whole line. Cutting only
the active variables gives a smooth section
\[
            X'=(f_2(u,v)+G(y)=0)\subset\mathbf P^{b+1}
\]
of dimension $b$. It lies at the end of a flag of smooth hyperplane
sections of $X$: choose a general nested flag through the chosen
line in $\mathbf P^{a-1}$. The restrictions of $f$ are smooth off
the line by Bertini and are smooth on it since its roots are simple.
The disjoint-variable gradient calculation then gives smoothness of
each hypersurface in the flag.

In $X'$ the tangent hyperplane section at $p$ is the projective cone
over $(G=0)\subset\mathbf P^{b-1}$. Indeed the binary tangent
direction is the root ray of $f_2$, and every first derivative of
$G$ vanishes at the affine origin. Thus $p$ is a generalized Eckardt
point. Abban--Zhuang, Theorem~4.18 \cite{ar2003.13788}, gives
\[
                 \delta_p(X',H)=\frac{b(b+1)}{b+d-1}.
\]
Their theorem requires precisely the preceding K-semistability
assumption on the base when $d\le b-1$.

For clarity, on a smooth $j$-fold hypersurface $T$ with smooth
hyperplane section $D$ through $p$, the hyperplane flag inequality is
\[
 \frac{\delta_p(T,H)}{j+1}\ \ge\
      \min\left\{1,\frac{\delta_p(D,H|_D)}j\right\}.
\]
This follows from adjunction for the refined series
\cite[Theorem~3.3]{ar2003.13788}: the expected coefficient of $D$ is
$1/(j+1)$, and the residual complete hyperplane series has expected
order $jS_{H|_D}/(j+1)$. Indeed, in section degree $m$ and fixed
hyperplane order $i$, the quotient is the complete space of degree
$m-i$ on $D$ by the hypersurface restriction sequence. Summing its
dimensions gives the stated mean residual coefficient.
Iterating up the flag yields, with $r=n+2-d$,
\[
     \delta_p(X,H)\ \ge\ \frac{(n+1)b}{b+d-1},\qquad
     \frac{(n+1)b}{b+d-1}-r
                   =\frac{(d-1)(d-a)}{b+d-1}>0.
\]
Lemma~\ref{lem:cyclic-delta-fixed} applies to every fixed point and
proves the theorem.
\end{proof}

\paragraph{Binary and singleton blocks.}
Blocks may be grouped before the theorem is applied. In particular,
suppose the equation is a sum of $t$ arbitrary squarefree binary
degree-$d$ forms and $u$ nonzero pure $d$th powers on distinct
coordinates. If
\[
        2t+u=n+2\le d(d-1),\qquad
        t\le d\lfloor(d-1)/2\rfloor,
\]
then it is uniformly K-stable. Allocate the binary blocks to $d$
groups with at most $\lfloor(d-1)/2\rfloor$ pairs each and fill the
remaining places with singletons. Every group has at most $d-1$
coordinates. Smoothness follows from the squarefree hypothesis.

\subsection{Klein hypersurfaces in unbounded dimension}
\label{sec:klein-all-dimensions}
\begin{theorem}\label{thm:klein-all-dimensions}
For every $d\ge5$ and $n\ge d-1$, the Klein hypersurface
\[
       X=\left(\sum_{i=0}^{m-1}x_i^{d-1}x_{i+1}=0\right)
           \subset\mathbf P^{m-1},\qquad m=n+2,
\]
with indices read modulo $m$, is uniformly K-stable.
\end{theorem}
\begin{proof}
Put $h=d-1$. For a loop in $k\ge2$ variables the partial
derivative in $x_i$ is $x_{i-1}^h+h x_i^{h-1}x_{i+1}$.
A vanishing coordinate at a common critical point forces its
predecessor to vanish, hence every coordinate is zero. If all
coordinates are nonzero, multiplication of the equations would give
$1=(-h)^k$. Thus every projective loop hypersurface used below is
smooth, including a binary loop with its reduced set of zeros.

Take $L=|(-h)^m-1|$, a primitive $L$th root of unity $\zeta$, and
the automorphism $g$ with diagonal entries $\zeta^{(-h)^i}$.
All loop monomials are invariant. The entries are pairwise distinct:
for $0\le i<j<m$ the nonzero difference of their exponents has
absolute value at most $h^{m-1}+h^{m-2}<L$. Thus the fixed points
are the coordinate vertices. The cyclic permutation lets us work at
$p=[1:0:\cdots:0]$.

In $X$ take the two linear cuts
\[
              x_{m-1}=u x_1,\qquad x_{m-2}=x_2+v x_1
\]
with general constants $u,v$, and denote the resulting section by $T$.
Here $m\ge d+1\ge6$, so all displayed coordinates apart from the
intended substitutions are distinct. On $x_1\ne0$ the first pencil
is basepoint free and hence cuts a smooth open variety, and the
second pencil is basepoint free there as well. Bertini gives smoothness
of $T$ on this open set. On $x_1=0$ its equation is the cone with
vertex $p$ over the smooth degree-$d$ loop on variables
$x_2,\ldots,x_{m-3}$, whose last summand is $x_{m-3}^h x_2$.
Thus $T$ is smooth there away from $p$, and its derivative in $x_1$
at $p$ is nonzero from the term $x_0^h x_1$. In particular $T$ is
smooth throughout. A general hyperplane between its linear span and
the full ambient space gives a smooth intermediate section: away
from $T$ this is Bertini, and on $T$ the differentials of its two
cutting equations in $X$ are independent.

We now induct on the ambient dimension. The dimension of $T$ is
$b=n-2$. Its tangent hyperplane section at $p$ is the cone just
described, and its base, when Fano, is a Klein hypersurface of
dimension $n-4$ to which induction applies. When $d>b-1$ the base
need not be Fano and Abban--Zhuang's theorem requires no stability
hypothesis on it. The generalized Eckardt formula and the two complete
hyperplane flags used in the proof of Theorem~\ref{thm:cyclic-block-hypersurfaces}
therefore give
\[
   \delta_p(X,H)\ \ge\frac{(n+1)(n-2)}{n+d-3},\qquad
   \frac{(n+1)(n-2)}{n+d-3}-(n+2-d)
          =\frac{(d-1)(d-4)}{n+d-3}>0.
\]
The same inequality holds at each fixed point. Apply
Lemma~\ref{lem:cyclic-delta-fixed}.
\end{proof}

\subsection{The quartic Klein hypersurface}
\label{sec:klein-quartic}
\begin{theorem}\label{thm:klein-quartic}
Every smooth quartic Klein hypersurface is uniformly K-stable in all
Fano dimensions $n\ge3$.
\end{theorem}
\begin{proof}
Use the published index-one and index-two results for $n=3,4$. For
$n\ge5$ put $N=n+1$, $k=n-1$ and $r=n-2=k-1$. At the coordinate
vertex $p=[1:0:\cdots:0]$ cut out its predecessor coordinate
$x_{n+1}=0$. The resulting hypersurface $T$ of dimension $k$ has
equation
\[
   U^3z+z^3y_1+G(y)=0,\qquad
   G(y)=\sum_{i=1}^{k-1}y_i^3 y_{i+1},
\]
where $U=x_0$, $z=x_1$ and $y_i=x_{i+1}$. An open chain
$\sum_{i=1}^{l-1}z_i^3z_{i+1}$ has as its entire projective singular
locus the last vertex: the last partial forces $z_{l-1}=0$, and the
other equations then force vanishing successively back to $z_1$.
In the dimensions considered here these hypersurfaces are integral and
normal, since distinct projective components meet in codimension one
and hypersurfaces are Cohen--Macaulay. In particular $T$ is smooth
near $p$. We will use hyperplane flag adjunction locally at this point.

\emph{A strict inequality on the space of tangent directions.}
Set $P=\mathbf P^{k-1}$, $H_P=\mathcal O_P(1)$, and
$c=3r/(4N)$. For every divisorial valuation $w$ on $P$ we claim
that there is a common positive constant $\epsilon$ such that
\begin{equation}\label{eq:klein-base-strict}
 A_P(w)-c\,w(G)\ \ge\left(\frac{rk}{N}+\epsilon\right)S_{H_P}(w).
\end{equation}
First, for rational $c'$ in a neighborhood of $c$, one has
\begin{equation}\label{eq:klein-base-weak}
 A_P(w)-c' w(G)\ \ge\ L(c')S_{H_P}(w),\qquad
 L(c')=k\left(1-\frac{4c'}{k-1}\right).
\end{equation}
A general hyperplane cuts $G$ smoothly. If $k=4,5$, degeneration
to the cone over that section gives $\operatorname{lct}(P;G)\ge
(k-1)/4$. On the cone the vertex blowup resolves the divisor and
has discrepancy $k-1$ and multiplicity four. Semicontinuity gives the
claim for $G$. Therefore $A_P-c'w(G)\ge(1-4c'/(k-1))A_P$, and use
$A_P\ge kS_{H_P}$.

For $k\ge6$ the particular hyperplane $y_k=y_1$ cuts out a smooth
quartic Klein hypersurface in dimension $n-4$, K-stable by induction.
Openness gives a smooth K-stable general section $G|_B$ on
$B\simeq\mathbf P^{k-2}$. The log interpolation theorem
\cite[Theorem~1.7]{ar2003.06839} gives K-semistability of
$(B,c'G|_B)$ for $c'$ in a neighborhood of $c$, since
\[
   c=\frac{3(k-1)}{4(k+2)}<\frac{3(k-1)}{4(k-2)}
       =1-\frac{(k-1)/4-1}{k-2}.
\]
Scale a coordinate complementary to $B$ to degenerate $G$ to the
cone over $G|_B$. The log projective cone criterion
\cite[Lemma~2.8]{blockZhou}
gives K-semistability with coefficient $4c'/(k-1)$ on the section
at infinity $B$. Openness gives K-semistability of
$(P,c'G+(4c'/(k-1))B)$. Its log anticanonical class is $L(c')H_P$,
so dropping the last effective boundary gives
\eqref{eq:klein-base-weak}.

Now $L(c)=rk/N<k-4c$. If \eqref{eq:klein-base-strict} failed, the
log Fano pair $(P,cG)$ would have stability threshold
$L(c)/(k-4c)<1$, attained by a prime divisor
\cite[Theorem~1.2]{ar2102.09405}. Interpolate linearly between
$A_P\ge kS_{H_P}$ at zero and \eqref{eq:klein-base-weak} at
$c'>c$. Equality at $c$ forces $A_P=kS_{H_P}$ for this divisor.
Every divisorial minimizer at threshold one for projective space
induces a special test configuration of Futaki invariant zero
\cite[Theorem~1.1(1)]{ar1907.05399}. Projective space is
K-polystable, so the configuration is a product. Its valuation is
given, up to a positive factor, by a one-parameter subgroup in
linear coordinates, with nonnegative weights $w_1,\ldots,w_k$,
minimum zero, and $A_P=\sum w_i$ in that normalization.

Eliminate a coordinate of minimum weight in a general hyperplane
section. Substituting a linear combination of the remaining coordinates
cannot decrease the minimum weight of a monomial of $G$. If the
remaining weights are not all equal, the smooth quartic section is
GIT stable, so Hilbert--Mumford gives
$w(G)\le w(G|_B)<4\sum_iw_i/(k-1)$. If the remaining weights are
all equal, this is a multiple of ordinary point order. The multiplicity
is less than four everywhere: the partial derivatives of the open
chain $G$ are linearly independent, by their distinct cubes
$y_1^3,\ldots,y_{k-1}^3$ and the remaining term $3y_1^2y_2$.
Multiplicity four at a projective point would make $G$ independent
of its coordinate. Thus $w(G)<4A_P/(k-1)$ in both cases, which
contradicts equality in \eqref{eq:klein-base-weak}. This proves
\eqref{eq:klein-base-strict}.

\emph{The point estimate.}
Blow up $p$ on $T$, with exceptional divisor $E=P$. The associated
graded section ring for ordinary point order has presentation
$\mathbf C[U,y_1,\ldots,y_k,Z]/(U^3 Z+G)$, with weights $0,1,4$
on $U,y_i,Z$. Indeed the equation is primitive and linear in $Z$,
and after localizing at $U$ it gives the initial relation in the
regular parameters $y_i/U$ at $p$. The domain before localization
has the full quartic Hilbert series, which gives the presentation in
every section degree. In degree $m$ and order $j$, the refined system
on $E$ is exactly $G^b\mathbf C[y]_{j-4b}$ for
$b=\max(0,\lceil(j-m)/3\rceil)$, whenever nonzero. This follows
by substituting for $Z$ in a monomial. Its averaged fixed coefficient,
moving degree and full point order, each divided by $m$, tend to
\[
   \overline b=\frac{3}{4(k+1)},\qquad
   \overline{j-4b}=\frac{k}{k+1},\qquad
   S_H(E)=\frac{k+3}{k+1}.
\]
These are also the complete-series computation in the proof of
\cite[Theorem~4.18]{ar2003.13788}. Directly, use basis monomials
$U^a Z^hR_l$ with $R_l$ in a graded vector-space complement to $(G)$.
The Hilbert series of its quotient is
$(1+t+t^2+t^3)/(1-t)^{k-1}$. Thus the two indices $a,h$ and $k-1$
free indices have the uniform simplex distribution after scaling,
up to bounded shifts. One has $j=4h+l$ and
$b=(h-a/3)_++O(1)$. Its normalized mean is $3/[4(k+1)]$, since for
two independent mean-one exponentials the mean positive part of
$E_1-E_2/3$ is $3/4$, and the total of $k+1$ such exponentials is
independent of their normalized vector. This proves the formulas.

Put $t=r(k+1)/N$. The radial comparison $k>t S_H(E)$ follows from
$k>r$. On $E$, the refined expected order at $w$ is
$\overline b\,w(G)+(k/(k+1))S_{H_P}(w)$. Equation
\eqref{eq:klein-base-strict} makes its threshold strictly greater
than $t$. Indeed $(P,cG)$ is klt, so $w(G)\le C A_P(w)$ with a fixed
$C$; in addition $S_{H_P}(w)\le A_P(w)/k$. Equivalently for a
uniform margin in the whole positive expression, use the klt bound
$A_P-cw(G)\ge\eta A_P$ for some $\eta>0$ together with
\eqref{eq:klein-base-strict} and take a convex combination with
coefficient sufficiently close to one on the latter inequality.
Adjunction for the refined linear series
\cite[Theorem~3.3 and Corollary~3.4]{ar2003.13788} gives
$\delta_p(T,H)>t$. This comparison is local near $p$, where $T$
is smooth. The complete hyperplane flag from $T$ back to $X$ has
all residual sections by the hypersurface sequence and gives
$\delta_p(X,H)\ge\min\{N,(N/(k+1))\delta_p(T,H)\}>r$.
Every coordinate vertex satisfies the same bound. Apply
Lemma~\ref{lem:cyclic-delta-fixed} to the diagonal automorphism in
Theorem~\ref{thm:klein-all-dimensions}.
\end{proof}

\subsection{Chain hypersurfaces}
\label{sec:chain-hypersurfaces}
\begin{corollary}\label{cor:chain-hypersurfaces}
For all $d\ge5$ and $n\ge d-1$, putting $m=n+2$, the hypersurface
\[
    x_0^{d-1}x_1+x_1^{d-1}x_2+\cdots+x_{m-2}^{d-1}x_{m-1}
      +x_{m-1}^d=0
\]
is uniformly K-stable.
\end{corollary}
\begin{proof}
Write $h=d-1$. This closed-ended chain is smooth in any number of
variables: if the first nonzero coordinate at a critical point has
index $j<m-1$, the $j$th partial forces $x_{j+1}=0$ and the next
partial then gives $x_j^h=0$. If only the last coordinate is nonzero
its pure power has nonzero derivative. The automorphism with entries
$\zeta^{(-h)^j}$, for $\zeta$ a primitive $d h^{m-1}$st root, has
distinct eigencharacters and fixes the equation. Its fixed points on
the hypersurface are precisely the coordinate points $p_i$ with
$0\le i\le m-2$. We describe a smooth linear section of codimension
at most two on which $p_i$ is generalized Eckardt and its tangent
hyperplane section is a cone over a smaller closed-ended chain.

At $i=0$ this already holds: set $x_1=0$ in the tangent section.
At $i=1$ take the single cut $x_0=u x_2$ with $u$ general. For
$2\le i\le m-3$ take the cuts
\[
       x_{i-1}=u x_{i+1},\qquad
       x_{i-2}=x_{i+2}+v x_{i+1}.
\]
The tangent equation on the section is $x_{i+1}=0$. The surviving
left chain ends at $x_{i-2}$, which has been identified with the
start $x_{i+2}$ of the right chain. Thus the residual polynomial is
a single closed-ended chain. This description includes $i=2$, where
the first variable of the original chain has no surviving term of
its own. Finally at $i=m-2$ take $x_{i-1}=u x_{i+1}$ and
$x_{i-2}=x_{i-3}+v x_{i+1}$: the last surviving term on the left
becomes $x_{i-3}^d$. Here $m\ge6$, so the indicated coordinates
are available.

For general parameters the cuts are smooth off $x_{i+1}=0$ by
successive Bertini, since the relevant pencils have no base point
on that open set. On $x_{i+1}=0$ the cone is smooth off $p_i$, and
the derivative in $x_{i+1}$ at $p_i$ is nonzero. A two-cut section
therefore fits in a smooth hyperplane flag as in
Theorem~\ref{thm:klein-all-dimensions}. If its dimension is $b$,
the generalized Eckardt formula and induction on dimension give
$\delta_{p_i}(X,H)\ge(n+1)b/(b+d-1)$. Indeed the residual
closed-ended chain, if Fano, has fewer variables. Since
$b\ge n-2>n+2-d$ for $d\ge5$, this estimate exceeds the index.
Now apply Lemma~\ref{lem:cyclic-delta-fixed}.
\end{proof}

\subsection{Pullbacks by quotients of binary forms}
\label{sec:binary-quotient-covers}
Let $e\ge2$. Partition the $n+2$ homogeneous coordinates into
singletons and pairs. On each singleton use $z\mapsto z^e$. On each pair
$(u,v)$ choose a finite subgroup $G\subset\mathrm{PGL}_2(\mathbf C)$ of
order $e$ and use two relatively prime binary forms $(P,Q)$ of degree $e$
representing the quotient map $\mathbf P^1\to\mathbf P^1/G\simeq
\mathbf P^1$. These substitutions together define a finite endomorphism
$\Phi$ of $\mathbf P^{n+1}$. Each pair may have its own group and its own
source and target coordinates.

\begin{proposition}\label{prop:binary-quotient-covers}
Let $Y=(F=0)$ have degree $m\ge1$, and let $X=(F\circ\Phi=0)$ be
smooth of Fano degree $3\le d=em\le n+1$. Then $Y$ is smooth. Set
$H=\mathcal O_Y(1)$ and $r=n+2-em$. If $\delta(Y,H)>r$, then $X$
is uniformly K-stable. In particular this holds whenever $Y$ is
K-semistable. One may take any hyperplane or smooth quadric for $Y$
for which the pullback is smooth.
\end{proposition}

\begin{proof}
The two binary forms have no common projective zero. The substitutions
give a basepoint-free degree-$e$ morphism, which is finite since it pulls
back the hyperplane class to $e$ times that class. The morphism is a
quotient by a finite group. Indeed, each element of a binary group admits
a matrix lift on which $(P,Q)$ is unchanged: an arbitrary matrix lift
multiplies both forms by the same scalar, since it fixes their ratio and
they are coprime; rescale the matrix to remove this scalar. There are
exactly $e$ such lifts of each group element. The resulting groups,
together with the $e$th roots of unity on the singletons, act transitively
on each general affine fiber block by block. After projectivizing one
divides by their common scalar subgroup. The same action preserves $X$.

Let $\Delta$ be the branch divisor with coefficients $1-1/h$, where $h$
is the ramification index. In the affine cone of the target a singleton
contributes the coordinate hyperplane with $h=e$. A binary block
contributes one line times the other affine factors for each branch
value of the quotient of $\mathbf P^1$, with the same ramification
index as that quotient. In the complement of the binary origin the
radial power map is \'etale, so this follows directly in a projective
line chart. Riemann--Hurwitz gives, in each binary block,
\[
     \sum_{b\ {\rm branch}}(1-1/h_b)=2(1-1/e).
\]
Consequently $\Delta/(1-1/e)$ is log canonical. In each binary plane
all line coefficients are at most one, and blowing up their common
origin gives an exceptional coefficient equal to one. Taking the
product of these resolutions with the singleton factors gives simple
normal crossings. The assertion descends from the punctured affine
cone to projective space by its local product with $\mathbf C^*$.

It also holds on $Y$. In a singleton block stratify the target by its
origin and complement. In a binary block take its origin, the branch
lines minus the origin, and their complement. At a lift of a point of
any product stratum, the image of $d\Phi$ on the affine cones is exactly
the tangent space to the stratum. This follows from the radial and
directional derivatives in a nonzero binary block, and the zero block
has zero differential. Smoothness upstairs makes $dF$ nonzero on this
image. Hence $Y$ is smooth and transverse to every projective stratum
it meets. Locally one may solve the equation of the cone over $Y$ for
a coordinate tangent to the stratum, while the normal coordinates
carry precisely the singleton and binary arrangements just resolved.
Thus $(Y,\Delta|_Y/(1-1/e))$ is log canonical. Transversality also shows
that $X\to Y$ is generically \'etale and the restrictions of the branch
hyperplanes have their indicated generic ramification indices. Its
general fibers are orbits, so the restricted map is finite Galois after
factoring out a possible kernel.

On the ambient projective space $\Delta\sim_{\mathbf Q}
(n+2)(1-1/e)\mathcal O(1)$ by the singleton and binary sums. Writing
$\Delta_Y=\Delta|_Y$, adjunction and ramification give
\[
      -(K_Y+\Delta_Y)\sim_{\mathbf Q}(r/e)H,\qquad
      K_X=(\Phi|_X)^*(K_Y+\Delta_Y).
\]
Log canonicity yields $A_{Y,\Delta_Y}(v)\ge A_Y(v)/e$ for every
divisorial valuation. Therefore the stability threshold of the log Fano
pair is at least $\delta(Y,H)/r$. Under the strict hypothesis this pair
is uniformly K-stable. Liu--Zhu's finite crepant Galois cover theorem
\cite[Theorem~1.2(2)]{ar2001.10557} gives K-polystability of $X$,
and finiteness of $\operatorname{Aut}(X)$ gives uniform K-stability.
Finally K-semistability of $Y$ implies
$\delta(Y,H)\ge n+2-m>n+2-em$.
\end{proof}

For a concrete instance take $e=60$ and the icosahedral group in each
binary block. Let $g_i$ be a binary form with simple roots consisting
of a single free orbit of this group; independent projective changes
of coordinates and independent orbits are allowed in the different
blocks. Such a form is a nonzero scalar multiple of a fiber equation
$a_iP_i+b_iQ_i$. Taking $Y$ to be a hyperplane proves that
\[
       \sum_i g_i(u_i,v_i)+\sum_j c_j z_j^{60}=0,\qquad c_j\ne0,
\]
is uniformly K-stable in every dimension $n\ge59$ with the indicated
number of variables. The gradient of a squarefree binary form vanishes
only at its affine origin, so this equation is smooth. The same
construction applies in degrees $12$ and $24$ using the tetrahedral and
octahedral groups. More general smooth pullbacks arise by taking $Y$ to
be a quadric and allowing its cross terms between different blocks.

\section{Larger blocks and attached chains}\label{sec:nov30-larger-families}

\subsection{Arbitrary blocks of up to \(d\) variables in prime degree}
\label{sec:nov30-prime-blocks}
\begin{theorem}\label{thm:nov30-prime-blocks}
Let \(d\ge5\) be prime and \(d-1\le n\).
Every smooth hypersurface whose equation is a sum of at most
\(d\) forms in disjoint blocks of at most \(d\) variables is
uniformly K-stable. In particular this permits \(d\) arbitrary
smooth forms, each in exactly \(d\) variables, giving
\(n=d^2-2\).
\end{theorem}
\begin{proof}
Use the diagonal symmetry with distinct \(d\)-th-root scalars
on the blocks and induction on the total number of variables
as in Theorem~\ref{thm:cyclic-block-hypersurfaces}. The gradient
test makes each summand smooth. A fixed point \(p\) is supported
on one block. If this block has fewer than \(d\) variables, the
binary cut there gives the same strict generalized Eckardt bound;
the complement, when Fano, is covered by the current induction.
Start the induction using the established index-one and index-two
results. It suffices to consider an active block of exactly \(d\)
variables and \(n\ge d+1\). Let \(m=n+2-d\ge3\),
\(k=m+1\), \(N=n+1=m+d-1\); the target in hyperplane
units is \(m\). Write \(G(Y)\) for the sum of the other blocks
in the \(m\) complementary variables. It is smooth and, if
Fano, uniformly K-stable by induction.

Restrict the active block to a general projective plane through
\(p\). The resulting degree-\(d\) plane curve is smooth by
Bertini off \(p\) and transversality at \(p\). The full
section \(T\) is a smooth hypersurface of dimension \(k=m+1\),
and intermediate general hyperplanes give a smooth flag.
Take coordinates \(U,x,z\) on the plane in which \(p=[1:0:0]\)
and its tangent line is \(z=0\). The intersection multiplicity
\(e\) with this line lies between two and \(d\), since a
smooth plane curve of degree \(d\) has no line component.
For \(e=d\) the tangent hyperplane on \(T\) is a cone over
\(x^d+G(Y)=0\), after rescaling. This base is smooth and,
if Fano, covered by induction by replacing the active block
with a singleton. The generalized Eckardt bound lifted along
the flag is \(Nk/(k+d-1)>m\).

Suppose \(2\le e<d\). We will use just \(\gcd(e,d)=1\).
After scaling coordinates give \(U,x,Y,z\) weights
\(0,d,e,de\), respectively. The initial relation for \(T\) is
\[
                 U^{d-1}z+U^{d-e}x^e+G(Y)=0.                \tag{1}
\]
Every term in the plane curve other than the displayed \(z\)
and \(x^e\) terms has strictly larger weight. The relation is
primitive and linear in \(z\), hence prime: its coefficient
is a power of \(U\) whereas \(U\nmid G\).
After setting \(U=1\), the local parameters are \(x,Y\)
and their initial expression for \(z\) is \(-x^e-G\).
Thus (1) presents the section initials, since localization at
\(U\) verifies injectivity in each fixed section degree and
the hypersurface Hilbert series gives equality of dimensions.

For \(S=S_{\mathcal O(1)}\) on \(\mathbf P^{m-1}_Y\) put
\[
      c_1=\frac{(d-1)m}{dN},\qquad
      c_2=\frac{(d-e)m}{edN}.
\]
For \(i=1,2\) there is \(\epsilon_i>0\) such that
\[
       A-c_i v(G)\ge(m-dc_i+\epsilon_i)S(v).              \tag{2}
\]
For \(m>d\), Zhang--Zhou, Theorem~1.7, applies to the smooth
K-stable hypersurface \(G=0\) and projective space. Both
coefficients lie in the open interval
\((0,(d-1)m/[d(m-1)])\), giving K-polystability; the projective
stabilizer of the smooth hypersurface is finite. When \(m=d\)
use small-angle log stability of a smooth anticanonical
hypersurface of projective space by its GIT stability and
Ascher--DeVleming--Liu, Theorem~1.4. Interpolate the strict
estimate at a sufficiently small positive coefficient with
\(A-v(G)\ge0\) at coefficient one. For \(m<d\) interpolate
\(A\ge mS\) with \(A-v(G)\ge0\); the resulting
\(m(1-c_i)\) strictly exceeds \(m-dc_i\).

Test the local threshold on \(T\) at
\(t=(k+1)m/N\). The radial weighted discrepancy is \(em+d\).
Choose a homogeneous complement \(R_j\) of
\(G\mathbf C[Y]_{j-d}\) in \(\mathbf C[Y]_j\). Polynomial
division in (1) gives the normal blocks
\(U^s x^a z^bR_j\), \(s+a+b+j=M\), in section degree \(M\).
Indeed the polynomial algebra in \(Y\) is free as a graded
module over \(\mathbf C[G]\) on the chosen complements;
eliminating \(G\) using (1) gives this presentation.
The Hilbert density of \(j\) is that of \(m-1\) simplex
coordinates. In units of \(M/(k+1)\), the mean values of
\(s,a,b,j\) are \(1,1,1,m-1\). The mean radial value at
this normalization is \(d+de+e(m-1)\); the strict comparison is
\[
     N(em+d)-m\bigl(d+de+e(m-1)\bigr)=d(d-1)>0.
\]

The angular weighted stack has weights \(d,e,\ldots,e\).
Since \((d,e)=1\), the coarse model is
\(P=\mathbf P(d,1^m)\) in coordinates \(A_0=x^e,Y\),
with boundary \((1-1/e)(A_0=0)\). Put \(Q=A_0+G\).
At fixed original weight the residue of \(a\) modulo \(e\)
is fixed and bounded. After removing it, in angular degree
\(K\) the independent blocks are
\[
       A_0^iQ^bR_j,\qquad di+db+j=K,\qquad ei+b+j\le M',
\]
where \(M'=M+O(1)\). The unused cost is
\(s=M'-ei-b-j\). The means of \(i,b,j,s\) in the same
units are \(1/e,1,m-1,1\).

Along \(A_0=0\) refine by its order \(i\).
At fixed \(i\), the leading systems in \(Y\) are
\(G^bR_{D-db}\), \(D=K-di\).
These are independent; the allowed \(b\)'s fill the interval
from \(b_0=\max(0,\lceil(ei+D-M')/(d-1)\rceil)\)
to \(\lfloor D/d\rfloor\). Thus the residual is exactly
\(G^{b_0}\mathbf C[Y]_{D-db_0}\).
Here \(b_0=(b-s/(d-1))_++O(1)\) has normalized mean
\((d-1)/d\), since independent unit exponentials satisfy
\(\mathbf E(E_1-E_0/(d-1))_+=(d-1)/d\), and their sum is
independent of the normalized simplex coordinates.
At multiplier \(t\) the coefficient on \(A_0=0\), including
the stack boundary, is \(1-1/e+m/(eN)<1\). The residual
coefficients are \(c_1\) on \(G\) and \(m^2/N=m-dc_1\)
on the complete moving hyperplane series. Thus (2) gives a
strict comparison by flag adjunction.

Off \(A_0=0\) omit its boundary. In coordinates \(Q,Y\)
scale \(Q\) while fixing \(Y\). This preserves each space,
since \(A_0\) changes to \(A_0+(\lambda-1)Q\), and the cost
of \(Q\) is one whereas that of \(A_0\) is \(e\).
Equivalently replacing \(A_0^i\) by \(G^i\) gives the
same span triangularly, in both directions. After splitting
by \(Q^b\), the residual at fixed \(b\) is exactly
\[
       G^{i_0}\mathbf C[Y]_{D'-di_0},\quad
       D'=K-db,\qquad
       i_0=\max(0,\lceil(D'+b-M')/(d-e)\rceil).
\]
This follows by summing \(G^iR_{D'-di}\) over all
\(i_0\le i\le\lfloor D'/d\rfloor\).
Consequently \(i_0=(i-s/(d-e))_++O(1)\), with normalized
mean \((d-e)/(ed)\), by the same exponential integral.
The means of \(b,D'\) are \(1,m-1+d/e\).
On the product charts \(Y\ne0\), the coefficient of \(Q=0\)
is \(m/N<1\); the residual \(G\)-coefficient is \(c_2\)
and the moving coefficient is again \(m^2/N\), since
\(m-1+d/e-d(d-e)/(ed)=m\).
This is strictly controlled by (2), as \(c_2<c_1\).
At the weighted vertex take the finite cyclic chart \(Q=1\)
and blow up the \(Y\)-origin. Its discrepancy \(m\) exceeds
\((m/N)(m-1+d/e)\); the same projective residuals give (2).
To transfer the estimate to arbitrary angular bases, choose a
basis adapted to both the tested valuation and the \(Q\)-power
filtration in each weight. Degeneration by the indicated scaling
of \(Q\) and lower semicontinuity of log canonical thresholds
transfer the common estimate for the split bases. Here the
boundary has been omitted globally in the scaling argument
and contributes zero at the tested center off \(A_0=0\).

The simplex coefficients converge, and the complete residual
hyperplane series have the uniform finite-level basis bounds.
The strict margins in (2), together with a positive klt
discrepancy margin for each of its log pairs, allow a small
common increase of \(t\); this also preserves all the strict
radial and fiber inequalities. Weighted adjunction gives
\(\delta_p(T,H)>t\). The complete hyperplane flag lifts this
to \(\delta_p(X,H)>m\). The cyclic fixed-point criterion
concludes the induction.
\end{proof}

The same argument allows a block of \(d\) variables in composite
degree wherever at every point of its hypersurface there is
a smooth plane section whose tangent contact is \(d\) or is
coprime to \(d\). All blocks of fewer than \(d\) variables
continue to require only smoothness.

\subsection{Full-size blocks in every degree}
\label{sec:nov30-all-degree-blocks}
\begin{theorem}\label{thm:nov30-all-degree-blocks}
Let \(4\le d\le n+1\). If a smooth hypersurface in
\(\mathbf P^{n+1}\) has an equation that is a sum of at most
\(d\) forms in disjoint blocks of at most \(d\) variables each,
it is uniformly K-stable. The degree need not be prime. This
permits \(n=d^2-2\), including four arbitrary smooth quaternary
quartics in dimension fourteen.
\end{theorem}
\begin{proof}
Keep the induction, finite diagonal automorphism, and notation of
Theorem~\ref{thm:nov30-prime-blocks}. At a point \(p\) supported
on a block of \(d\) variables, let \(S\subset\mathbf P^{d-1}\)
be that smooth summand and \(G\) the smooth complementary sum in
\(m=n+2-d\) variables. The induction begins with index one and
two, so \(m\ge3\). Put \(N=n+1=m+d-1\). The finite Gauss map
of \(S\), whose polarization pulls back to \((d-1)H\), implies
that its tangent hyperplane section \(D\) at \(p\) has at most
isolated singular points.

In the projective tangent directions \(\mathbf P^{d-3}\) impose
the homogeneous Taylor coefficients \(q_2=\cdots=q_{d-2}=0\)
of \(D\). These are \(d-3\) forms of positive degrees on
\(\mathbf P^{d-3}\), so there is a common direction, even if
some coefficients vanish identically. Let \(L\) be its line
through \(p\). If \(L\not\subset S\), the contact order
is \(d-1\) or \(d\). Extend \(L\) to a general plane.
It cuts \(S\) smoothly outside \(L\) by Bertini; at the
finitely many points of \(S\cap L\) choose the extra direction
transverse whenever necessary. The plane-curve proof of the
preceding theorem applies, since \(d-1\) is coprime to \(d\).

Suppose \(L\subset S\). At the general point of \(L\) the
tangent divisor \(D\) is smooth. Choose a general tangent
two-plane \(\Pi\) containing \(L\), with coordinates
\(U,x_1,x_2\) such that \(p=[1:0:0]\) and \(L=(x_1=0)\).
The induced curve on \(\Pi\) contains \(L\) with multiplicity
exactly one: a general additional direction has a nonzero
polar along \(L\). Thus its polynomial is
\(x_1 H(U,x_2)+O(x_1^2)\) with \(H\ne0\). Let \(v\) be the
smallest exponent of \(x_2\) occurring in \(H\).
Since \(p\) is a tangent point, \(1\le v\le d-1\).
The curve on \(\Pi\) has at most finitely many singular points:
along \(L\) its transverse derivative is \(H\), and away from
\(L\) use Bertini for the general linear cut through \(L\) on
\(D\) and the isolated singular locus of \(D\). When \(d=4\)
the two-plane is the tangent plane itself. Extend \(\Pi\) by
a further direction \(z\) transverse to \(S\) at all these
finitely many points. A general such three-plane cuts \(S\)
smoothly, by Bertini off \(\Pi\); for \(d=4\) use the whole
ambient three-space. This gives a smooth full hypersurface
section \(T\) of dimension \(k=m+2\), with a smooth intermediate
flag by transversality along the final cut.

Write \(h=1+v\le d\). Choose positive integer weights
\(a,b,c\) with \(a+vb=dc\), \(2a>dc\), and primitive total
weight vector, for example \(b=1,\ c=2v+1,\ a=dc-v\).
Give \(U,x_1,x_2,Y,z\) weights \(0,a,b,c,dc\).
All terms containing \(x_1^2\) have weight greater than \(dc\);
the term of least weight among those linear in \(x_1\) is the
one with exponent \(v\) in \(x_2\). After rescaling, the
initial relation is therefore
\[
       U^{d-1}z+U^{d-h}x_1x_2^v+G(Y)=0.                    \tag{3}
\]
The equation is prime and computes the section initials by
the localization and Hilbert series argument for (1) in the
preceding subsection. With \(R_j\) the same graded complements,
the normal blocks are
\(U^s x_1^\alpha x_2^\beta z^{b'}R_j\), with
\(s+\alpha+\beta+b'+j=M\). In units of \(M/(k+1)\) the
means of \(s,\alpha,\beta,b',j\) are \(1,1,1,1,m-1\).
The radial discrepancy is \(a+b+cm\). At the required
multiplier \(t=(k+1)m/N\) the strict radial inequality follows
from
\[
    N(a+b+cm)-m\bigl(a+b+dc+(m-1)c\bigr)=(d-1)(a+b)>0 .
\]

On the angular weighted stack a torus acts with characters
\(v,-1,0\) on \(x_1,x_2,Y\), preserving all angular spaces
and fixing \(A_0=x_1x_2^v\) and \(Q=A_0+G\).
In a torus character \(v\alpha-\beta\) write uniquely
\[
  (\alpha,\beta)=(\gamma_1,\gamma_2)+i(1,v),\qquad
          i=\min(\alpha,\lfloor\beta/v\rfloor).
\]
The residual exponents \(\gamma_1,\gamma_2\) depend just on
the character. In each original angular weight at that character,
after factoring \(x_1^{\gamma_1}x_2^{\gamma_2}\) the system is
the pullback of the independent blocks on
\(P=\mathbf P(d,1^m)\) in \(A_0,Y\) given by
\[
       A_0^i Q^{b'} R_j,\qquad
       di+db'+j=K,\qquad hi+b'+j\le M',\qquad
       M'=M-\gamma_1-\gamma_2 .
\]
To bound all angular bases, choose in each original weight a
basis simultaneously adapted to the tested valuation and the
torus-character filtration. Degenerate by the torus. The limiting
bases have the displayed factorization; a common strict lc bound
for all such split bases transfers by lower semicontinuity.
Here the angular stack is tested on its smooth uniformizing charts.

All means below are taken over all weights and characters, using
the same normalized units. For independent unit exponentials
\(E_1,E_2\), their minimum \(\min(E_1,E_2/v)\) is exponential
of rate \(h\). It is independent of the independent variables
for \(s,b'\). Dividing the exponential vector by its sum gives
the normal simplex, so the means of
\(i,\gamma_1,\gamma_2\) are \(1/h,v/h,1/h\). In particular
the unused cost \(s=M'-hi-b'-j\) still has normalized mean one.
The calculation for (2) in the preceding subsection now applies
to the residual projective systems on \(P\), with \(e=h\).
More explicitly, at \(A_0=0\), extracting the \(i\)-th
order gives the fixed \(G\)-power
\((b'-s/(d-1))_++O(1)\) of mean \((d-1)/d\), and moving
hyperplane degree of mean \(m\). Away from \(A_0=0\), split
further by scaling \(Q\). The fixed \(G\)-power is
\((i-s/(d-h))_++O(1)\) of mean \((d-h)/(hd)\) if \(h<d\).
If \(h=d\), replacing \(A_0\) by \(Q-G\) leaves the same
cost, all \(G\)-powers down to zero occur at fixed \(Q\)-power,
and the fixed exponent is zero. In both cases the moving degree
again has mean \(m\). The coefficient on \(Q=0\) is \(m/N\).
Thus (2) gives the strict residual estimate in both places;
at zero \(G\)-coefficient use \(A\ge mS>(m^2/N)S\).
The weighted vertex of \(P\) is controlled by the same
\(Q=1\) chart as above, since \(m-1+d/h<N\).
It follows that the residual projective estimates allow an
additional fixed coefficient \(B\ge0\) on \(A_0=0\)
whenever \(B+m/(hN)<1\); off this divisor its coefficient
is omitted in the \(Q\)-scaling argument.

We check the passage from \(P\) to all angular charts of (3).
On a chart \(Y_j=1\), if both \(x_1,x_2\ne0\), the map
to \((A_0,Y)\) is smooth and the fixed monomials are units.
At \(x_1=0\) refine first by that divisor. Its total mean
order is \(\gamma_1+i=\alpha\) in each normal block, hence
its asymptotic coefficient is \(m/N<1\). The leading
residual along it contains the other coordinate to order
\(\gamma_2+vi=\beta\) of the same normalized mean one,
times exactly the \(G\)-power and moving linear system
computed on \(A_0=0\). The latter are pulled back from the
\(Y\) coordinates and independent of \(x_2\); this is a product
chart. The same reasoning, refining first by \(x_2=0\), gives
the bound there, with \(\gamma_2+vi\) as first order.
These flag estimates apply to arbitrary split bases, since in
each character and original weight we can maximize the order of
any additional valuation using a basis adapted also to the
coordinate order.

For a chart with \(x_1\ne0\), pass to its uniformizing cover
\(x_1=1\). Then \(A_0=x_2^v\) and the systems descend to
the affine cone in \((A_0,Y)\), with fixed factor
\(x_2^{\gamma_2}\). By the finite-cover discrepancy formula
the equivalent boundary on the cone along \(A_0=0\) is
\[
           B=1-1/v+\frac{m}{hvN}.
\]
Including the mean \(A_0\)-order of the systems gives
\(B+m/(hN)=1-1/v+m/(vN)<1\).
Blow up the cone origin with weights \(d,1,\ldots,1\).
Its projective residual is just the systems on \(P\) already
controlled. After accounting for the cover ramification, its
radial discrepancy is \(m+d/v\). The fixed factor and the
remaining homogeneous systems have radial degree of normalized
mean \(m-1+d+d/v\): the terms are \(d\gamma_2/v\),
\(di\), \(db'\) and \(j\). The required difference is
\[
           N(m+d/v)-m(m-1+d+d/v)=d(d-1)/v>0 .
\]
On the chart \(x_2=1\) the map is simply \(A_0=x_1\);
the fixed factor is \(x_1^{\gamma_1}\) with coefficient
\(B=vm/(hN)\), so \(B+m/(hN)=m/N<1\).
The radial difference on the cone is \(d(d-1)>0\),
by the same formula with the cover degree one.
The weighted blowup calculation and the homogeneous projective
test also give the estimate on the punctured cone, by its
local product presentation on the uniformizing charts.

The lc estimates throughout have strict margins by (2), the
strict divisor and radial comparisons, and the klt margin
for the two fixed pairs on \(\mathbf P^{m-1}\), exactly as
for (1). Thus they persist at a slightly larger \(t\).
Weighted adjunction and the complete hyperplane flag give
\(\delta_p(X,H)>m\). Applying the cyclic fixed-point
criterion at every fixed point concludes the induction.
\end{proof}

\subsection{A further variable in each arbitrary block}
\label{sec:nov30-blocks-degree-plus-one}
\begin{theorem}\label{thm:nov30-blocks-degree-plus-one}
Let \(4\le d\le n+1\). Every smooth degree-\(d\) hypersurface
whose equation is a sum of at most \(d\) forms in disjoint blocks
of at most \(d+1\) variables is uniformly K-stable. This includes
\(d\) independently chosen smooth forms with \(d+1\) variables
each, of total hypersurface dimension \(d(d+1)-2\).
For quartics this endpoint is dimension eighteen.
\end{theorem}
\begin{proof}
Induct on the total number of variables, with the same cyclic
symmetry and the published index-one and index-two starting cases.
The arguments of the previous two subsections at a block of size
at most \(d\) require only the induction hypothesis on its smooth
complement when Fano, so it remains to treat a point \(p\) on
a block of size \(d+1\). Let \(m\) be the number of other
variables, \(r=m+1\), and \(N=n+1=m+d\). We have \(m\ge2\)
in the induction step. Their sum \(G(Y)\) is smooth and when
Fano K-stable by induction.

On the projective tangent directions \(\mathbf P^{d-2}\) at
\(p\) the locus \(C=(q_2=\cdots=q_{d-2}=0)\) has dimension
at least one, by the projective dimension theorem. If
\(q_{d-1}\) is nonzero somewhere on \(C\), this gives a
line of contact exactly \(d-1\) and the ternary graph (1)
with \(h=d-1\), coprime to \(d\). Otherwise, the degree-\(d\)
form \(q_d\) has a zero on \(C\) and the corresponding line
is contained in the block hypersurface. The same argument with
the finite Gauss map yields the primitive monomial graph (3)
with \(1\le v=h-1\le d-1\). We use the same smooth linear
sections and flags as there, of respective dimensions \(k=m+1\)
and \(k=m+2\).

We sharpen the angular estimate by retaining the contribution
of the distinguished coordinate divisor in adjunction. Write
\(P=\mathbf P(d,1^m)\) with coordinates \(A_0,Y\) and write
\(V(w)=w(A_0=0)\). Let \(S(w)\) denote the expected order on
the split residual projective systems
\(A_0^iQ^{b'}R_j\) of the preceding proofs, in units of
\(M/(k+1)\), summed with their original multiplicities; it
excludes all monomial prefactors \(x_i^{\gamma_i}\).
Set
\[
      \theta_1=\frac{m}{m+d-1},\qquad
      \theta_2=\frac{m}{m+d/h-1},\qquad
      \theta=\frac{m+1}{m+d}.
\]
We claim that every divisorial \(w\) on \(P\) satisfies
\[
 \begin{split}
    A_P(w)&\ge\theta_1 S(w)+(1-\theta_1/h)V(w),\\
    A_P(w)&\ge\theta_2 S(w).
 \end{split}
 \tag{4}
\]
Here and below these comparisons of basis systems may equivalently
be made with strictly smaller coefficients at finite level, and
passed to the limit by uniform convergence on the complete
hyperplane spaces.

For the first inequality at a center on \(A_0=0\), refine by that
divisor. Its mean order is \(1/h\). The full residual at multiplier
\(\theta_1\) has a fixed \(G\)-coefficient
\(c_1=\theta_1(d-1)/d\) and moving degree \(\theta_1m=m-dc_1\).
It is lc by (2), whose proof applies also when \(m=2\) by
the smooth binary branch and \(m<d\). Inversion of adjunction
with total coefficient \emph{one} on \(A_0=0\) leaves precisely
the extra \((1-\theta_1/h)V\) on the right: this is the finite-level
proof of Abban--Zhuang, inequality (3.6), before forgetting that
term. The second inequality suffices at centers off \(A_0=0\)
for the first one since \(\theta_2>\theta_1\).

For the second inequality, omit all boundary and use the
\(Q\)-scaling calculation on all of \(P\). At multiplier
\(\theta_2\), the coefficient of \(Q=0\) is \(\theta_2\le1\);
the base has fixed \(G\)-coefficient
\(c_2=\theta_2(d-h)/(hd)\) and moving degree
\(\theta_2m=m-dc_2\). The same projective pair comparison
applies: for \(m>d\), this \(c_2\) is below
\((d-1)m/[d(m-1)]\), and in the other cases use the preceding
interpolation. At \(h=d\) one has \(c_2=0\) and the exact
projective space inequality \(A\ge mS_{\mathcal O(1)}\).
The weighted vertex comparison is
\(m\ge\theta_2(m-1+d/h)\), an equality by definition.
Use multipliers strictly below \(\theta_2\), which make the
radial and fiber tests strict, and take their limit. This
proves (4) everywhere.

Combining the two inequalities in (4) leaves available at
multiplier \(\theta\) the boundary coefficient
\[
       B_*=(1-\theta_1/h)
              \frac{\theta_2-\theta}{\theta_2-\theta_1}
           \quad\hbox{on }A_0=0 .
\]
The inequalities \(\theta_1<\theta<\theta_2\) hold for
\(h\ge2,m\ge2\). The exact excess over \(1-1/h\) is
\[
   B_*-(1-1/h)=
    \frac{(d-h)(d-1)\,[m(h-1)-h]}
         {h m d(h-1)(m+d)} .                                \tag{5}
\]
For the ternary graph, \(h=d-1\ge3\), so (5) is positive.
Its entire stack boundary on \(P\) is \((1-1/h)(A_0=0)\);
the neighboring multipliers still give a bound for that boundary.
The radial numerator of the original weighted blowup is now
\[
      N(hm+d)-r\bigl(d+dh+h(m-1)\bigr)
                      =(d-1)(d-h)>0 .
\]

For the primitive graph put \(\nu_1=1,\nu_2=v=h-1\).
The mean remaining monomial exponents are
\(\overline\gamma_i=1-\nu_i/h\). At multiplier \(\theta\),
including the fixed factor and finite ramification on a
chart with a root of order \(\nu_i\) requires on \(A_0=0\)
only the boundary
\[
 B_i=1-\frac{1-\theta\overline\gamma_i}{\nu_i}
     =1-\frac\theta h-\frac{1-\theta}{\nu_i}
     <1-\frac1h,\qquad i=1,2.
\]
The strict inequality follows from \(\nu_i<h\) and
\(\theta<1\). Formula (5) is nonnegative even for \(h=d\)
or \(h=m=2\). Thus \(\max B_i<B_*\), and this remains true
when \(\theta\) is slightly increased.

For completeness the same maximum of the two coefficients
suffices at points where \(Y\ne0\) and both monomial
coordinates vanish. On the uniformizing chart \(Y_j=1\),
the map is \(A_0=x_1x_2^v\) with the remaining \(Y\)
unchanged. Write \(B=\max B_i\). The lc pair on the smooth
base chart with boundary \(B[A_0=0]\) and the residual
basis divisor remains lc on its product with
\((\mathbf P^1,[0]+[\infty])\). On the common function
field use \(x_2\) as the product coordinate. The identity
\[
           d\log x_1\wedge d\log x_2
                   =d\log A_0\wedge d\log x_2
\]
identifies the corresponding logarithmic canonical divisors
in the two coordinate systems. More explicitly the pullback
of \(dA_0\wedge dx_2\) is
\(x_2^v dx_1\wedge dx_2\). The product pair therefore pulls
back on this chart with boundary coefficients
\(B\) on \(x_1=0\) and \(1-v(1-B)\) on \(x_2=0\)
and the pulled-back residual divisor. Each of these coefficients
is at least \(\theta\overline\gamma_i\) by definition.
This proves the lc bound also at their intersection, for the
split bases; use a strict multiplier margin to absorb uniform
finite-level errors in the averaged prefactors.

At \(x_1=1\) and \(x_2=1\) use the weighted cone charts
from (3). For remaining cover degree \(j=v\) or \(j=1\)
the new radial difference is
\[
         N(m+d/j)-r(m-1+d+d/j)=(d-1)(d/j-1)>0 .
\]
The angular projective test uses exactly \(B_i\) above. The
original primitive radial comparison is
\[
      N(a+b+cm)-r\bigl(a+b+dc+(m-1)c\bigr)
                       =(d-1)(a+b-c)>0,
\]
for the weights \(a=dc-v,b=1,c=2v+1\) used in (3).
The proper torus degeneration transfers these strict estimates
to all angular bases as in the previous subsection. Weighted
adjunction thus gives \(\delta_p(T,H)>(k+1)r/N\). Lifting
along the smooth complete hyperplane flag proves the strict
test \(\delta_p(X,H)>r\) for every fixed point and concludes
by Lemma~\ref{lem:cyclic-delta-fixed}.
\end{proof}

\subsection{Sums with several loop and chain blocks}
\label{sec:nov30-mixed-loop-blocks}
Put \(h=d-1\). A loop of length \(l\) means
\(\sum_{i\in\mathbf Z/l}z_i^h z_{i+1}\); a closed-ended chain
of length \(l\) means \(\sum_{i=0}^{l-2}z_i^h z_{i+1}+z_{l-1}^d\).
A one-variable loop is the pure power \(z_0^d\).
In this subsection a smooth form in a set of variables means one
whose first derivatives have no common nonzero zero.

\begin{theorem}\label{thm:nov30-mixed-loop-blocks}
Let \(4\le d\le n+1\). Suppose a degree-\(d\) equation is a sum
of smooth forms in disjoint sets of variables, of the following types.
\begin{enumerate}
\item Arbitrary forms in at most \(d-1\) variables.
\item Loops of length \(l\ge5\). For \(d=4\) assume \(4\nmid l\).
\item When \(d\ge5\), closed-ended chains of length at least five.
\end{enumerate}
Let \(c\) be the number of chains in the third class, with \(c\le d\).
Assign to an arbitrary form its number of variables, and to a loop
the integer in \(\{1,2,3,4\}\) congruent to \(l\) modulo four.
If these assigned integers can be partitioned into \(d-c\) groups
with each sum at most \(d-1\), the hypersurface is uniformly K-stable.
Empty groups are permitted. In every \(d\ge5\) this includes sums of
at most \(d\) arbitrary long loops or closed-ended chains, with
arbitrary smooth small forms allowed in the indicated groups.
For each \(d\ge4\) it also includes a loop of length \(4j+1\)
together with an arbitrary smooth form in at most \(d-1\) further
variables, in unbounded dimension as \(j\) increases.
\end{theorem}

\begin{proof}
Smoothness of the sum and of the partial sums follows from the
independence of the gradient equations in disjoint variables.
Use induction on the total number of variables. Index one and
index two are covered by the published theorems in
Section~\ref{sec:established-cases}. We will use the fixed-point
criterion, generalized Eckardt formula, and complete hyperplane
flag inequality proved or cited in Section~\ref{sec:cyclic-block-hypersurfaces}.

First construct a single diagonal symmetry with separate spectra.
For a loop of length \(l\ge5\) put \(N_l=|(-h)^l-1|\). Given a
primitive \(N_l\)-th root \(\xi\), the diagonal entries
\(\xi,\xi^{-h},\ldots,\xi^{(-h)^{l-1}}\) preserve its equation.
They all have order \(N_l\). They are distinct: for \(0<k<l\)
one has \(0<|(-h)^k-1|<N_l\). The primitive roots consequently
split into \(\varphi(N_l)/l\) disjoint orbits of size \(l\).
Different lengths have different \(N_l>d\), all coprime to \(h\).
There are enough distinct orbits even for the largest possible
number \(d(d-1)\) of loops allowed by the partition.

Here is the elementary estimate used in this count.
For all \(h\ge3,\ l\ge5\),
\[
             \varphi(|(-h)^l-1|)\ \ge\ l h(h+1).
\]
Euler's product gives \(\varphi(N)\ge N^{2/3}/2\): each prime
power \(p^a\) contributes \(p^{a/3}(1-1/p)\) to the ratio, which
is at least one for \(p\ge5\), and the product of the smallest
possible factors for two and three exceeds one half. The ratio
\((h^l-1)^{2/3}/[2l h(h+1)]\) increases in each argument on the
stated domain. In \(h\), its logarithmic derivative exceeds
\((2l/3-2)/h>0\); in \(l\) it exceeds
\(2\log(h)/3-1/l>0\).
The ratio is at least one at \((h,l)=(3,7),(4,6),(7,5)\),
as is seen by squaring \(h^l-1\) and cubing \(2l h(h+1)\).
The remaining cases have \(N_l=244,728,1025,3126,7777\) with
totients \(120,288,800,1040,6000\), respectively.

Assign a distinct character \(\tau\in\mu_d\) to each chain and
each of the \(d-c\) groups. Within a group take the sum of its
arbitrary small forms as a single block and give its variables
the group's character. Loops of length at least five use the
separate primitive spectra just constructed. On a chain give
the coordinate at distance \(j\) from its terminal variable the
eigenvalue \(\tau\eta_j\), where \(\eta_0=1\),
\(\eta_j^{-h}=\eta_{j-1}\) and \(\eta_j\) has exact order \(h^j\).
Such a sequence is obtained by taking compatible primitive roots.
The monomials are fixed because \(-h\equiv1\bmod d\). Since
\((h,d)=1\), the \(h\)-primary part records the distance \(j\),
and the part of order dividing \(d\) records the chain.
Thus every fixed point of the resulting projective automorphism
is supported on a single small group, or is one coordinate vertex
of a large loop or a chain. The terminal chain vertices do not
belong to the hypersurface.

For a point on a small group, write \(a\) for the total number of
its variables and \(b=n+2-a\). The binary-section argument of
Theorem~\ref{thm:cyclic-block-hypersurfaces} gives
\[
      \delta_p(X,H)\ge\frac{(n+1)b}{b+d-1}>n+2-d
\]
since \(a\le d-1\). The complement of the entire group is smooth;
whenever it is Fano it satisfies the induction hypothesis because
deleting the small group does not increase any assigned sum.
Here and in the vertex calculations every appeal to induction
in the generalized Eckardt formula is needed only for a Fano base.

At a large-loop vertex for \(d\ge5\), use the two cuts in the
proof of Theorem~\ref{thm:klein-all-dimensions} entirely within
this block. The base of the generalized Eckardt section replaces
that loop by the loop on \(l-4\) variables. At \(l=5\) it is the
one-variable pure power: at the vertex \(z_0\), take
\(z_4=u z_1,\ z_3=z_2+v z_1\) and then set \(z_1=0\).
At \(l-4\le4\) the shortened loop is now an arbitrary small
form taking precisely the reserved number of places in its group.
For longer loops the assigned integer is unchanged.
The smoothness argument for the cuts is unchanged: away from
the tangent hyperplane use the general parameters and Bertini,
and on it the cone over the smooth residual sum is smooth off
the vertex, where the normal derivative is nonzero.

At a chain vertex apply the corresponding cuts in
Corollary~\ref{cor:chain-hypersurfaces}. None involves any
other block. At the first or second vertex the residual chain
has length \(l-2\) or \(l-3\); at the remaining vertices its
length is \(l-4\). For \(l=5\) at the penultimate vertex the
explicit cuts \(z_2=u z_4,\ z_1=z_0+v z_4\) leave \(z_0^d\).
If a residual chain has at most four variables, place it alone
in its previously reserved chain group, increasing the number
of available small groups by one as \(c\) decreases by one.
Otherwise it is still one of the \(c\) chains. Thus all these
bases satisfy the induction hypotheses. The generalized Eckardt
and linear flags at every tested vertex give
\[
       \delta_p(X,H)\ \ge\ \frac{(n+1)b}{b+d-1},
                 \qquad b\ge n-2>n+2-d \quad (d\ge5).
\]
The right side exceeds \(n+2-d\) since its difference has numerator
\((d-1)(b-(n+2-d))\).

For a quartic loop vertex use instead the strict angular comparison
in the proof of Theorem~\ref{thm:klein-quartic}. We verify its
inputs for the present sum. Take the cut making the predecessor
coordinate zero. With \(U\) the vertex coordinate, \(z\) its
successor and \(k=n-1\), the section has local form
\[
           U^3z+z^3 y_1+G(y)=0 .
\]
The form \(G\) in \(k\) coordinates is a sum of one open chain
on \(l-3\ge2\) variables and all the other smooth blocks. Its
only projective singularity is the last vertex of the open chain,
by successive differentiation; all other blocks have independent
gradient equations. The complete linear section itself is
normal and smooth at the marked point by the same argument
applied to its open chain of length \(l-1\).
The partial derivatives of \(G\) are linearly independent.
For the open chain its last \(l-4\) partials have distinct pure
cube terms and the first is nonzero; for any other block
linear dependence would make that smooth form independent of
one of its variables. Distinct blocks use disjoint variables.
Finally identify the last and first open-chain variables in a
hyperplane. This gives the smooth sum with the tested loop
replaced by the loop of length \(l-4\), a pure power when those
two open-chain variables were the only ones. Its assigned group
remains valid: by \(4\nmid l\), a newly small loop has at most
three variables.

These facts give exactly the two inputs for the strict comparison
of Theorem~\ref{thm:klein-quartic}: a general hyperplane of \(G\)
is smooth and has a K-stable section whenever Fano by induction
and openness, and \(G\) is not a cone. For \(k=4,5\) only a
smooth general section is necessary. The same hypersurface
presentation \(U^3 Z+G\) computes the angular and radial
coefficients, since \(G\) is nonzero. The cited strict estimate
therefore gives \(\delta_p(X,H)>n-2\) here as well. The fixed-point
criterion concludes the proof in every degree.
\end{proof}

\subsection{A chain ending in an arbitrary form}
\label{sec:nov30-attached-chain}
\begin{theorem}\label{thm:nov30-attached-chain}
Let \(d\ge5,\ 1\le a\le d-2,\ l\ge1,\ n=a+l-2\ge d-1\).
Let \(f\) be any degree-\(d\) form in \(a\) variables \(x\) and
let \(L(x)\) be a nonzero linear form. Every smooth hypersurface
\[
 f(x)+L(x)y_1^{d-1}
              +\sum_{j=2}^l y_{j-1}y_j^{d-1}=0
         \quad\subset\mathbf P^{a+l-1}
\]
is uniformly K-stable. The length \(l\) is unrestricted.
\end{theorem}

\begin{proof}
Write \(h=d-1\) and \(r=n+2-d\). When \(a=1\), smoothness makes
the core coefficient nonzero, and rescaling gives the closed-ended
chain of Corollary~\ref{cor:chain-hypersurfaces}. Assume \(a\ge2\).
Smoothness implies that \(f\) and its restriction to \(L=0\)
have no nonzero critical vector. For \(f\) set all the \(y\)'s
to zero. If \(L(x)=0,\ x\ne0\) and \(df_x\) is proportional to
\(dL\), choosing \(y_1^h\) to cancel that proportionality and
all \(y_i=0\) for \(i\ge2\) gives a singular point. Conversely
these two smoothness conditions suffice at every length: at a
critical point a nonzero outermost \(y_j\), \(j\ge2\), forces
\(y_{j-1}=0\) from the derivative in \(y_j\), then
\(y_j^h=0\) from the derivative in \(y_{j-1}\). The remaining
equations give one of the two core critical points just excluded
or, at \(x=0\), force \(y_1=0\).

Use induction on dimension and the published index-one and
index-two results. In the inductive step take \(n\ge d+1\);
in particular \(l\ge5\). The automorphism that fixes all core
coordinates and gives \(y_j\) a compatible primitive \(h^j\)-th
root \(\eta_j\), with \(\eta_j^{-h}=\eta_{j-1}\) and \(\eta_0=1\),
preserves the equation. Its fixed points are the chain vertices
and the points of \((f=0)\) in the core. We prove the strict
hyperplane threshold \(>r\) at each such point.

At the \(j\)-th chain vertex the tangent hyperplane coordinate
is \(T=y_{j-1}\) for \(j\ge2\), and \(T=L\) for \(j=1\). At
the outermost vertex \(j=l\), no cut is needed: on \(T=0\)
the remaining hypersurface is a cone over the same attached
chain with two fewer variables. At \(j=l-1\) make the single
cut \(y_l=uT\); the shortened chain now has length \(l-3\).
For \(3\le j\le l-2\) take
\[
      y_{j+1}=uT,\qquad y_{j+2}=y_{j-2}+vT .
\]
At \(j=2\) use the same cuts with \(y_{j-2}\) replaced by \(L\).
On \(T=0\) the two surviving pieces then meet in a single
attached chain of length \(l-4\). These descriptions include
a vanishing piece at the end of a tail.

Lastly \(j=1\) gives \(T=L\). Use
\[
              y_2=uL,\qquad y_3=y_4+vL .
\]
The remaining base is the disjoint sum of \(f|_{L=0}\) and a
closed-ended chain of length \(l-3\), with terminal pure power
\(y_4^d\). All these bases are smooth by the two core conditions
just proved, and when Fano are K-stable by induction or by
Theorem~\ref{thm:nov30-mixed-loop-blocks}. For each use of the
generalized Eckardt formula we have a smooth linear section:
off \(T=0\) the displayed pencils are successively basepoint
free; on \(T=0\) the cone is smooth away from its vertex and the
derivative in \(T\) is nonzero at the vertex. A general flag to
the smooth final complete intersection is also smooth. The
generalized Eckardt formula and hyperplane flags give
\[
       \delta_p(X,H)\ge\frac{(n+1)b}{b+d-1}>r,\qquad b\ge n-2,
\]
because \(d\ge5\). Here \(b\) is the dimension of the tested
linear section.

Now take \(p\in(f=0)\) inside the core. Cut the core to a
general projective line through \(p\) and set \(y_1=0\).
The binary form on the line is squarefree, in particular it has
a simple root at \(p\). Choose coordinates \(U,z\) there with
that root at \(z=0\). The resulting section is
\[
    T=\left(f_2(U,z)+G(y_2,\ldots,y_l)=0\right),\qquad
    G=\sum_{j=3}^l y_{j-1}y_j^h .
\]
It is smooth near \(p\), normal everywhere, and has dimension
\(k=l-1\). Indeed the open chain has just its inner endpoint
as projective singular locus, by successive differentiation.
The section \(z=0\) is a cone with vertex \(p\). If \(a<d-2\),
make one further cut \(y_2=y_3\). Its cone base is now a smooth
closed-ended chain, K-stable when Fano. There are \(a\) cuts
in total, leaving dimension \(n-a>r\). The same generalized
Eckardt and hyperplane argument proves strictness.

When \(a=d-2\), use the section \(T\) without that final cut,
so \(k=r+1\). We give the strict comparison at this boundary.
Let \(P=\mathbf P^{k-1}\), \(H_P=\mathcal O(1)\), \(N=n+1\), and
set \(c=(d-1)r/(dN)\). For rational \(c'\) near \(c\) one has
\[
    A_P(w)-c'w(G)\ \ge\ k\left(1-\frac{dc'}{k-1}\right)
                                     S_{H_P}(w).
                                                        \tag{1}
\]
For \(k-1\le d\) degenerate \(G\) to the projective cone over
a smooth hyperplane section, giving
\(\operatorname{lct}(P;G)\ge(k-1)/d\). Use \(A_P\ge kS_{H_P}\).
When \(k-1>d\), setting \(y_2=y_3\) in the open chain gives a
smooth closed-ended chain of dimension \(k-3\), K-stable by the
chain theorem. The smooth general hyperplane section is K-stable
by openness. Zhang--Zhou log interpolation and the projective
cone construction, as in Theorem~\ref{thm:klein-quartic}, give
the bound after dropping the effective hyperplane at infinity.
Their range contains a neighborhood of \(c\), because
\[
    c=\frac{(d-1)(k-1)}{dN}<
          \frac{(d-1)(k-1)}{d(k-2)}
            =1-\frac{(k-1)/d-1}{k-2}.
\]

The coefficient on the right in (1) at \(c'=c\) is \(rk/N\),
strictly less than \(k-dc\). In fact (1) is strict uniformly
at \(c\), in the sense that \(rk/N\) can be increased slightly.
If equality occurred in the infimum, the log Fano pair
\((P,cG)\) would have threshold less than one, attained by a
divisorial valuation. Interpolating (1) above \(c\) with the
projective space inequality at zero forces \(A_P=kS_{H_P}\)
for that valuation. The product configuration argument in
Theorem~\ref{thm:klein-quartic} then makes it a linear
one-parameter-subgroup valuation. Its coordinate weights have
minimum zero and satisfy \(A_P=\sum w_i\). Eliminating a
minimum-weight coordinate on a general hyperplane cannot lower
the minimum weight of \(G\). The smooth hyperplane section is
GIT stable, hence \(w(G)<dA_P/(k-1)\) unless all the
remaining weights are equal. The latter valuation is ordinary
point order. The \(k\) partial derivatives of the open chain
are independent, by their distinct pure powers and the nonzero
remaining endpoint derivative. Thus its multiplicity is less
than \(d\) at every point, and the same strict inequality holds.
This contradicts equality in (1).

At \(p\in T\) the ordinary-order associated graded ring has
equation \(U^{d-1}Z+G=0\), after scaling \(Z\), with weights
\(0,d,1,\ldots,1\) on \(U,Z,y\). The equation is prime and
localizing at \(U\) gives the regular local parameters and this
initial relation; the ordinary Hilbert series gives the full
graded presentation. As in the quartic open-chain calculation
or the proof of Abban--Zhuang, Theorem~4.18, the normalized
mean degree of the radial order is \((k+d-1)/(k+1)\); on its
exceptional \(P\) the coefficients of the fixed \(G\) and
complete moving hyperplane series are respectively
\((d-1)/(d(k+1))\) and \(k/(k+1)\). Explicitly the fixed
exponent at level \(m\), order \(j\), is
\(\max(0,\lceil(j-m)/(d-1)\rceil)\). In the normal monomial
count \(U^aZ^bR_e\) its mean is that of \((b-a/(d-1))_+\)
up to a bounded error, where the \(a,b\) simplex coordinates
have mean \(m/(k+1)\). This gives the displayed coefficient.

The radial discrepancy \(k\) exceeds
\(r(k+d-1)/N\). The strict form of (1) gives the corresponding
strict angular comparison at multiplier \(r(k+1)/N\), using
the klt margin exactly as in Theorem~\ref{thm:klein-quartic}.
The local flag theorem gives
\(\delta_p(T,H)>r(k+1)/N\). The complete hyperplane flags lift
this to \(\delta_p(X,H)>r\). They can be taken normal and smooth
at \(p\): the final section \(T\) has at most an isolated
singularity, its cutting differentials are independent on its
smooth locus, and Bertini applies outside the final linear
space. All section restrictions are complete by the hypersurface
sequence. Lemma~\ref{lem:cyclic-delta-fixed} now proves the theorem.
\end{proof}

\paragraph{Indecomposable examples.}
Already \(a=2\) gives varying examples that do not split in
disjoint groups of variables after any linear change.
For a form \(F\) on the coordinate space, consider constant
matrices \(A\) such that \(\operatorname{Hess}(F)A\) is symmetric.
Any nontrivial disjoint decomposition supplies a nonscalar matrix,
by taking the projection to one summand. General smooth binary
\(f\) have only scalar matrices in this system. For instance
\(f=U^d+V^d+tU^{d-1}V\) at general \(t\ne0\) has two
rank-one Hessian coefficients, at \(U^{d-3}V\) and \(V^{d-2}\),
forcing \(A\) diagonal. Its mixed Hessian coefficient at \(U^{d-2}\)
then equates the diagonal entries. Having scalar system is open.

Take \(L=U\) and such an \(f\). In the attached chain, the
coefficient matrix at \(y_{j-1}y_j^{d-3}\), taking \(y_0=U\),
is a nonzero rank-one square of \(dy_j\). Thus each \(dy_j\)
is an eigenvector of the dual \(A\). The off-diagonal Hessian
matrix at \(y_j^{d-2}\) then equates the successive eigenvalues
for \(j\ge2\). The core Hessian coefficients force the core
component to be scalar and its rows to have no tail terms:
a vector annihilated by all \(\operatorname{Hess}f\) would be
an inessential variable, impossible for a smooth binary form.
The coefficient at \(y_1^{d-2}\) finally equates the core and
tail scalars. Hence \(A\) itself is scalar. These examples vary
in moduli: restricting the infinitesimal projective-orbit
directions \(\sum \ell_i\partial_iF\) to \(y=0\) gives at
most the four-dimensional binary \(\mathfrak{gl}_2\)-orbit
space, whereas the binary degree-\(d\) variations have dimension
\(d+1\). Smooth choices exist, since \(f|_{U=0}=V^d\)
and a general such binary form is squarefree.

\subsection{Quartic chains with a binary end form}
\label{sec:nov30-quartic-binary-end}
\begin{theorem}\label{thm:nov30-quartic-binary-end}
The conclusion of Theorem~\ref{thm:nov30-attached-chain} holds
also for \(d=4,\ a=1,2\). In particular, in every Fano dimension
each closed-ended quartic chain and each smooth quartic
\[
    f(U,V)+L(U,V)y_1^3+\sum_{j=2}^l y_{j-1}y_j^3=0
\]
is uniformly K-stable; \(f\) can be any squarefree binary
quartic whose zero set avoids the point \(L=0\).
\end{theorem}
\begin{proof}
Use the two core transversality conditions, diagonal automorphism
and induction on the total dimension from
Theorem~\ref{thm:nov30-attached-chain}. The one-variable and
two-variable cores are included in this same induction. In
dimensions \(3,4\) the published index-one and index-two
theorems suffice. Take \(n\ge5\), write \(k=n-1,\ r=k-1,\ N=k+2\),
and let \(p\) be one of the fixed points.

For a two-variable core point the section \(y_1=0\) gives
\(f(U,V)\) plus an open chain in \(k\) variables, exactly the
\(k=r+1\) strict angular calculation in
Theorem~\ref{thm:nov30-attached-chain}. Its smooth hyperplane
section is a closed-ended quartic chain of smaller dimension,
to which the present induction applies whenever Fano. The
one-variable core has no point on the hypersurface.

At the outermost and second outermost vertices the earlier
zero-cut and one-cut generalized Eckardt estimates are strict,
since the dimensions of the sections are at least \(n-1>r\).
Suppose the tested vertex is \(y_j\) with \(j\le l-3\). Set its
predecessor \(y_{j+1}=0\). The tangent form in the resulting
one-cut section is the sum of an open chain of length
\(l-j-1\ge2\) and the smooth core with its remaining shorter
tail. For \(j=1\) the latter is just the core on \(L=0\),
which for \(a=1\) has no variables. The partial derivatives
are independent and the singular locus is isolated. There is
a smooth hyperplane section of this tangent form in the inductive
family. For \(j\ge3\) identify \(y_{j+2}\) with \(y_{j-2}\),
and at \(j=2\) identify \(y_{j+2}\) with \(L\), reattaching
the two remaining tails. For \(j=1,a=2\), identify \(y_3\)
with the one surviving core variable; its pure fourth power
closes the remote open chain. For \(j=1,a=1\) identify \(y_3\)
with \(y_4\) instead. All these hyperplanes reduce the number
of variables, and their sections are K-stable when Fano by
induction. The one-cut hypersurface is normal and smooth at \(p\):
it is the disjoint sum of a shorter smooth attached chain and
one open chain with isolated projective singularity. The strict
open-chain comparison in Theorem~\ref{thm:klein-quartic}, with
the verified general-section and noncone hypotheses, gives
\(\delta_p(X,H)>r\).

\paragraph{The third vertex from the outer end.}
Take \(j=l-2\) and set \(y_l=0\). The resulting smooth
hypersurface \(T\) has dimension \(k\). At \(p\), use \(U=y_j\)
for its nonzero coordinate, \(x=y_{j+1}\) for the new outer
end, and \(z=y_{j-1}\) for the tangent hyperplane. The remaining
\(m=k-1\) coordinates are \(Y\). Give \(U,x,Y,z\) weights
\(0,4,3,12\), respectively. The initial equation is
\[
                       U^3z+U x^3+G(Y)=0,                    \tag{1}
\]
where \(G\) is the smooth shorter attached quartic remaining
after removing the outer four tail variables from \(X\).
Whenever the projective hypersurface \(G=0\) is Fano it is
K-stable by dimension induction. The omitted monomial
\(z^3\) times its successor has weight greater than twelve.
The primitive linear relation (1) is prime, and after
localizing at \(U\) it computes the weighted filtration in
the smooth parameters \(x,Y\) at \(p\). The ordinary Hilbert
series proves the presentation of the global section initials.

On \(\mathbf P^{m-1}_Y\) write \(S=S_{\mathcal O(1)}\).
For each of the two coefficients
\[
         c_1=\frac{3r}{4N},\qquad c_2=\frac{r}{12N}
\]
there is a uniform strict inequality
\[
       A_Y-c_i v(G)\ \ge\ (m-4c_i+\epsilon_i)S,
                         \qquad \epsilon_i>0,                \tag{2}
\]
where \(A_Y\) here denotes discrepancy on the projective space
of the \(Y\) directions. If \(m>4\), the log pair is K-polystable
by Zhang--Zhou, Theorem~1.7: \(G\) is K-stable and
\(c_i<3m/[4(m-1)]\). Its projective stabilizer is finite, so
it is uniformly K-stable. If \(m=3\), the pair \((\mathbf P^2,G)\)
is lc and \((1-c_i)3>3-4c_i\). If \(m=4\), use the quartic
surface GIT comparison of Ascher--DeVleming--Liu at any small
rational positive \(b<\min(c_i,1/3)\). A smooth quartic surface
is GIT-stable, and the pair at \(b\) is uniformly K-stable.
Interpolate its strict inequality with \(A_Y-v(G)\ge0\) at
coefficient one to obtain (2). These cases exhaust \(k\ge4\).

We show \(\delta_p(T,H)>(k+1)r/N\). Let \(t=(k+1)r/N\).
The radial discrepancy of the weighted blowup is \(3k+1\).
In degree \(M\) use a graded complement \(R_e\) to
\(G\mathbf C[Y]_{e-4}\) in \(\mathbf C[Y]_e\). Ordinary
division in (1), with \(G\) as the eliminated leading form,
gives independent normal blocks
\[
       U^s x^\alpha z^b R_e,\qquad s+\alpha+b+e=M .
\]
For example one may take the complements from a monomial order
eliminating a nonzero pure \(Y_1^4\); the leading term then
comes from \(G\) and the normal basis has the asserted counts.
The \(e\)-index has the Hilbert density of \(m-1=k-2\)
simplex coordinates. In units of \(M/(k+1)\), the limiting
means of \(s,\alpha,b,e\) are \(1,1,1,k-2\). The radial mean
is therefore \((3k+10)/(k+1)\). The strict radial test is
\((3k+1)(k+2)>(3k+10)(k-1)\), with difference twelve.

The angular stack has weights \(4,3,\ldots,3\) on \(x,Y\).
Put \(A=x^3,\ Q=A+G\). Its coarse model is
\(P=\mathbf P(4,1^m)\) in \(A,Y\), with stabilizer boundary
\((2/3)(A=0)\). At each fixed original weight the residue
\(\alpha\bmod3\) is constant; its order after division by \(M\)
has vanishing contribution. After removing this bounded residue
and dividing weights by three, the angular system at weight
\(K\) is spanned by the independent blocks
\[
          A^i Q^b R_e,\qquad 4i+4b+e=K,\qquad
          3i+b+e\le M',
\]
where \(M'\) differs from \(M\) by the bounded residue. In
the same units the means of \(i,b,e\) are \(1/3,1,k-2\);
the unused cost \(s=M'-3i-b-e\) has mean one.

At \((A=0)\) refine these blocks by their \(A\)-orders. At
fixed \(i\), their leading residues are \(G^bR_{D-4b}\)
with \(D=K-4i\). This is a direct sum, and the allowed \(b\)
form the full interval from \(b_0=\max(0,\lceil(3i+D-M')/3\rceil)\)
to \(\lfloor D/4\rfloor\). The residual is consequently
\(G^{b_0}\mathbf C[Y]_{D-4b_0}\).
Here \(b_0=(b-s/3)_++O(1)\). From the uniform simplex in
\(3i,b,s\) and the \(e\) coordinates its mean is \(3/4\) in
the stated units. Indeed for independent unit exponentials
\(\mathbf E(E_1-E_0/3)_+=3/4\), and their sum is independent
of the normalized simplex coordinates. Thus the fixed
\(A\)-coefficient after testing at \(t\) is
\(2/3+r/(3N)<1\); the residual \(G\)-coefficient is \(c_1\),
and the residual hyperplane degree coefficient is
\[
              \frac rN(4+k-2-3)=\frac{r m}{N}=m-4c_1 .
\]
Inequality (2) and flag adjunction give a strict margin here.

Away from \((A=0)\) omit its boundary and give \(Q\) its own
one-parameter scaling, with \(Y\) fixed. This preserves every
angular space: in \(Q,Y\) coordinates one has \(A=Q-G\), and
replacing \(A\) by \(A+(\lambda-1)Q\) preserves the costs
\(3,1\) of \(A,Q\), in both directions. Equivalently, replacing
\(A^i\) by \(G^i\) gives the same full span triangularly at
fixed \(K\). After splitting by the \(Q\)-powers, the exact
residual at \(Q^b\) is
\[
       G^{i_0}\mathbf C[Y]_{D'-4i_0},\qquad
       D'=K-4b,\quad i_0=\max(0,D'+b-M').
\]
Indeed the allowed indices \(0\le i\le\lfloor D'/4\rfloor\)
satisfy \(3i+b+D'-4i\le M'\). Thus \(i_0=(i-s)_+\).
The normalized means of \(b,i_0,D'\) are
\(1,1/12,k-2+4/3\): the middle is
\(\mathbf E(E_1/3-E_0)_+=1/12\).
On a chart with \(Y\ne0\), the coefficient of \(Q=0\) at
multiplier \(t\) is \(r/N<1\), and the base coefficients
are \(c_2\) on \(G\) and \(rm/N\) on the moving hyperplane
series. Since \(c_2<c_1\), (2) yields the strict test.
At the weighted vertex take the cyclic chart \(Q=1\) and
blow up the \(Y\)-origin. Its discrepancy is \(m=k-1\)
against \(r(k-2+4/3)/N<m\); the same projective residual
systems give (2). For arbitrary initial angular bases, choose
them simultaneously adapted to the given valuation and the
\(Q\)-power filtration. Degenerate by this projective scaling
of \(Q\). The preceding chart computations bound every limiting
basis with the \(A\)-boundary omitted globally; semicontinuity
therefore transfers the bound at centers off \(A=0\).

All simplex comparisons have the usual uniform limiting errors
for these complete hyperplane systems. Since (2) has positive
margins and \(A_Y-c_i v(G)\) also bounds a positive multiple
of \(A_Y\) by klt, convex interpolation preserves the angular
inequalities after a small common increase in \(t\). Weighted
adjunction and the one original hyperplane flag now prove
\(\delta_p(X,H)>r\). The cyclic fixed-point criterion concludes.
\end{proof}

\subsection{Several chains attached to a smooth core}
\label{sec:nov30-several-attached}
\begin{theorem}\label{thm:nov30-several-attached}
Fix \(d\ge5\) and set \(h=d-1\). Suppose the equation of a smooth
degree-\(d\) hypersurface in \(\mathbf P^{n+1}\), \(n\ge d-1\),
is a sum of forms on disjoint variable sets. Allow a closed-ended
chain of any positive length, including a pure power. Also allow
components of the form
\[
 f(x)+\sum_{s=1}^{t}\left[L_s(x)y_{s,1}^{h}
       +\sum_{j=2}^{l_s} y_{s,j-1}y_{s,j}^{h}\right],
 \qquad x=(x_1,\dots,x_a),\quad a\ge2,
\]
where \(l_s\ge1\), the \(L_s\) are linear, and
\[
                   a+\sum_s\min(2,l_s)\le d-1 .
\]
Count each closed-ended chain once and each displayed component
\(1+t\) times. If the total count is at most \(d\), then the
hypersurface is uniformly K-stable. In particular, for \(d\ge7\)
the statement includes two chains of arbitrary lengths attached
in independent directions to a smooth binary core, whenever the
resulting equation is smooth.
\end{theorem}
\begin{proof}
Smoothness of a homogeneous form here means that its affine
gradient has no nonzero zero. Each disjoint summand has this
property. In an attached component the core form is smooth by
setting all chain variables equal to zero. The forms \(L_s\) are
independent: a dependence \(\sum c_s L_s=0\) would give a nonzero
critical vector at \(x=0\), with \(y_{s,1}^h=c_s\) and all other
coordinates zero. Deleting any whole tail or shortening any tail
preserves smoothness, since a critical vector of the shortened
form extends by zeros. There is also the following restriction
rule. Delete tail \(s\) and restrict all core variables to
\(L_s=0\). A critical vector of that restricted component would
give a differential of the component with the tail deleted that
is a multiple of \(dL_s\). Choosing \(y_{s,1}^h\) to cancel the
multiple and setting the rest of that tail to zero would give a
nonzero critical vector of the initial component. The restricted
component is consequently smooth. If its core has one variable,
independence gives at most one attached tail and smoothness makes
its core coefficient nonzero. Rescaling makes this component a
closed-ended chain.

Give the components distinct scalars \(\tau\in\mu_d\). Inside an
attached component use \(\tau\) on the core and \(\tau\eta_{s,j}\) on
its tail variables. Choose distinct nonidentity
\(\eta_{s,1}\in\mu_h\), then choose
\(\eta_{s,j}^{-h}=\eta_{s,j-1}\). The numerical hypothesis gives
\(t\le d-3<h\). For any such character at depth \(j\), its
\(h^j\)-th power is one and its \(h^{j-1}\)-th power is nontrivial.
Thus different depths have distinct characters, and different
tails at equal depth are separated by their first characters.
For a closed-ended chain regard its pure-power variable as a
one-variable core and use the same construction on its single
tail. Since \((d,h)=1\), components with different scalars have
disjoint spectra. This finite-order diagonal automorphism preserves
the equation. Each of its fixed points on the hypersurface lies
in a single core or is a single tail vertex.

Induct on the number of variables, starting with the established
index-one and index-two cases. Write \(r=n+2-d\). At each
fixed point \(p\) we find a smooth linear section of dimension
\(b>r\) on which the tangent hyperplane section is a cone with
vertex \(p\) over a smooth hypersurface of the same class. When
Fano the base is K-stable by induction. The generalized Eckardt
formula and complete hyperplane flags then give
\[
              \delta_p(X,H)\ge\frac{(n+1)b}{b+d-1}>r .
\]
Indeed the difference has numerator \((d-1)(b-r)\).

First suppose \(p\) is on an attached core. Restrict that core
to a general two-dimensional vector subspace through \(p\).
The restricted binary form has simple roots by Bertini and
transversality at \(p\). In each attached tail set \(y_{s,1}=0\).
For a length-two tail also set \(y_{s,2}=0\); for a tail of
length at least three impose \(y_{s,2}=y_{s,3}\). In the latter
case the surviving summand is a closed-ended chain of length
\(l_s-2\), with pure power \(y_{s,3}^d\). The final equation
is a squarefree binary form plus a smooth sum in the remaining
variables. Its tangent through the simple root exposes the cone
over that sum. The number of cuts is
\(a-2+\sum_s\min(2,l_s)\le d-3\), so \(b>r\). Removing
the core component of count \(1+t\) has added at most \(t\)
closed-ended chains, so the induction applies. A one-variable
pure-power core contributes no point.

At a vertex of a separate closed-ended chain of length less than
\(d\), one may instead regard the whole chain as a single smooth
core without tails, with its one unit of count. At lengths
\(\ge d\) the one- and two-cut constructions in
Corollary~\ref{cor:chain-hypersurfaces} leave a single shorter
closed-ended chain or remove the component. For length five at
the penultimate vertex, in the ordering
\(z_0^h z_1+\cdots+z_3^hz_4+z_4^d\), the cuts
\(z_2=u z_4,\ z_1=z_0+v z_4\) leave \(z_0^d\) on the base.

It remains to check a vertex \(y_j\) of one attached tail,
suppressing its index \(s\), and writing \(l=l_s,L=L_s\).
For \(j\ge2\) the tangent coordinate is \(T=y_{j-1}\).
At \(j=l\) the tangent cone already removes the last two variables.
At \(j=l-1\) impose \(y_l=uT\), leaving a tail of length \(l-3\).
At \(2\le j\le l-2\) impose
\[
        y_{j+1}=uT,\qquad y_{j+2}=y_{j-2}+vT ,
        \qquad y_0=L .
\]
On \(T=0\) the surviving tail is attached to the original core
and has length \(l-4\). A length zero means that the tail
disappears. At \(j=1\) the tangent coordinate is \(T=L\).
Impose no cut if \(l=1\). If \(l\ge2\) use \(y_2=uT\).
For \(l=3\) also use \(y_3=vT\), and for \(l\ge4\) also use
\(y_3=y_4+vT\). On \(T=0\) the old core is restricted to
\(L=0\) with this tail deleted. For \(l\ge4\) a disjoint
closed-ended chain of length \(l-3\) remains. The restriction
rule proves smoothness. Its core has lost a variable and a
tail, and a single closed-ended chain costs at most the one
unit released by deleting the tail. Treating a one-variable
residual core as a closed-ended chain only decreases the count.
Every other summand stays the same.

These vertex constructions cost at most two cuts, so
\(b\ge n-2>r\). Away from \(T=0\) their pencils are successively
free. On \(T=0\) the cone over the smooth base is smooth off
\(p\), and the derivative in \(T\) is nonzero at \(p\).
The general indicated section is therefore smooth. A general
hyperplane flag through each smooth final section is smooth by
Bertini outside it and by independence of its cutting
differentials along it. The restriction series are complete by
the hypersurface exact sequence. Finally
Lemma~\ref{lem:cyclic-delta-fixed} gives the conclusion.
\end{proof}

\begin{corollary}\label{cor:nov30-two-attached}
In Theorem~\ref{thm:nov30-several-attached}, a component with
exactly two attached tails may alternatively satisfy
\[
         a=2\qquad\hbox{or}\qquad a+3\le d-1 ,
\]
regardless of the lengths of both tails. In particular a smooth
binary quintic with two independently attached chains gives a
uniformly K-stable hypersurface in every Fano dimension in which
its equation occurs.
\end{corollary}
\begin{proof}
Enlarge the class used in the induction of
Theorem~\ref{thm:nov30-several-attached} by the indicated
components. Restriction of the core at a first tail vertex
reduces \(a\) and deletes a tail. Its remaining single tail
satisfies the original bound if \(a\ge3\), since
\((a-1)+2\le d-1\). For \(a=2\) the residual is a closed-ended
chain. The same bound \(a+2\le d-1\) holds when a whole tail
disappears at a later vertex and the core does not shrink.
All earlier arguments therefore apply except at a core point.

Write \(A_0=L_1(p),B_0=L_2(p)\), in a fixed nonzero representative
of \(p\). In a binary core both are nonzero, since smoothness of
the restricted form on each \(L_i=0\) means that this point of
\(\mathbf P^1\) is not a root. For general \(a\) first cut the core
to a general projective line through \(p\), using \(a-2\) cuts.
Choose line coordinates \(U,z\) with \(p=[1:0]\). The binary
restriction has a simple root at \(p\) and is squarefree. Its
other roots avoid \(L_i=0\): for \(a\ge3\) the smooth core
hypersurface has no linear component and meets each of these
hyperplanes in codimension two in projective space; projection
from \(p\) excludes the lines through another point of that
intersection. The binary intermediate section thus is smooth
off \(z=0\). Indeed at a critical vector each tail has all
coordinates beyond its first zero, by successive differentiation.
Every nonzero first coordinate forces its attaching form to
vanish, and Euler's identity forces the binary core coordinate
to be a root. The case where all first coordinates vanish is
excluded by squarefreeness.

If \(A_0B_0=0\), close a zero-coefficient tail at the cost of a
single cut. For length one set its first variable equal to \(uz\).
For a longer tail set its first variable equal to its second
plus \(uz\). On \(z=0\) the first attaching term vanishes, so
in the second case the remaining chain ends in a pure power.
On a nonzero-coefficient tail use the original one or two cuts,
allowing a general multiple of \(z\) on the right of each to
give a free pencil away from \(z=0\). Thus at most three cuts
in addition to the \(a-2\) core cuts suffice, and this case
can only occur under the hypothesis \(a+3\le d-1\).

Suppose \(A_0B_0\ne0\). Set \(v=y_{1,1}\) and choose
\(\lambda^h=-A_0/B_0\). Make the first cut
\[
                          y_{2,1}=\lambda v+u z .
\]
On \(z=0\) the two attaching terms cancel identically.
If both tails have a second variable, write
\(A=y_{1,2},\ B=y_{2,2}\). Make the second cut
\(v=M(A,B)+u'z\) for a general binary linear form \(M\).
On the tangent hyperplane the equation, after omitting the
vertex coordinate \(U\), has a new binary core
\[
                         M(A,B)(A^h+\lambda B^h).
\]
The subsequent variables form two tails attached along \(A,B\),
of lengths \(l_1-2,l_2-2\); discard tails of length zero.
The binary core is squarefree and is nonzero on each coordinate
axis away from the origin. This new component is smooth.
At a critical vector only its two first tail coordinates could
be nonzero beyond the core. If exactly one is nonzero, the
core would have a root on its attaching axis. If both are
nonzero the core vector must be zero, and the independence of
\(dA,dB\) then contradicts the core gradient equation.
The case of no first tail coordinate follows by squarefreeness.
If only one second variable \(A\) exists, instead put
\(v=A+u'z\): its nonzero contribution to the tangent equation
is a pure power in \(A\), starting a shorter closed-ended chain.
If neither exists put \(v=u'z\), so the original component drops
out of the base altogether.

In the nonvanishing case there were exactly two extra cuts.
For \(a=2,d\ge5\) and for \(a+3\le d-1\) the section has
dimension \(b>r\). The new binary component consumes at most the
original three units of count. Whenever tails were closed or
discarded, at most two units remain. The intermediate pencils
are successively free off \(z=0\), and on this hyperplane the
section is the cone over the smooth sum just described. At \(p\)
the derivative in \(z\) is nonzero. Smoothness of the final
section and the application of generalized Eckardt and
hyperplane flags follow as before. This finishes the enlarged
induction.
\end{proof}

\subsection{A binary quartic with two attached chains}
\begin{corollary}\label{cor:nov30-quartic-two-tails}
Let \(f\) be a binary quartic and \(L_1,L_2\) binary linear forms.
For \(l_1,l_2\ge1\), every smooth Fano quartic
\[
 f(U,V)+\sum_{s=1}^2\left[
 L_s(U,V)y_{s,1}^{3}
       +\sum_{j=2}^{l_s}y_{s,j-1}y_{s,j}^{3}\right]=0
\]
is uniformly K-stable, without restriction on the two lengths.
\end{corollary}
\begin{proof}
Induct on dimension together with the zero- and one-tail cases
of Theorem~\ref{thm:nov30-quartic-binary-end} and the
closed-ended chain case of that theorem, using the published
index-one and index-two results in dimensions three and four.
The two linear forms are independent by the critical-vector
test at the origin of the core, and the squarefree quartic has
no root on either attaching axis by restriction as in
Theorem~\ref{thm:nov30-several-attached}. Put \(n=l_1+l_2\),
\(k=n-1\) and \(r=n-2\). Give the two inner tail coordinates
the two different nonidentity characters in \(\mu_3\) and
extend by inverse cubing out the tails as in that theorem.
Only the tail vertices and the four binary roots are fixed.

For clarity, the one-hyperplane estimate established in
Theorem~\ref{thm:klein-quartic} uses the following inputs.
The section is normal and smooth at the mark; its point-order
initial equation is \(U^3Z+G(Y)\) in weights \(0,4,1\),
where \(G\) is a degree-four form in \(k\) variables.
The form \(G\) has independent first derivatives and admits a
smooth hyperplane restriction which is K-stable whenever Fano.
The proof there then gives \(\delta_p(X,H)>r\): the hyperplane
restriction yields the weak projective log inequality, and
the derivative condition excludes equality for ordinary point
order. The associated graded presentation uses primeness of
\(U^3Z+G\), localization at \(U\), and the quartic Hilbert
series. It is unaffected by terms of larger point weight.

Let \(p\) be a binary root. Write \(p=[1:0]\) in binary
coordinates \(U,z\). Neither \(L_s(p)\) vanishes. Choose
\(\lambda^3=-L_1(p)/L_2(p)\), put \(v=y_{1,1}\), and cut by
\(y_{2,1}=\lambda v+u z\) with general \(u\). On \(z=0\)
the terms in \(Uv^3\) cancel. The tangent form independent
of \(U\) is
\[
            G=v(A^3+\lambda B^3)+A y_{1,3}^3
                +B y_{2,3}^3+\cdots ,
\]
where \(A=y_{1,2},B=y_{2,2}\); omit any absent variables and
their terms. At least one of \(A,B\) exists, since \(n\ge5\).
The derivative in \(v\) has a nonzero pure-cube term. Each of
the other derivatives has its own mixed monomial consisting
of the preceding variable times the square of the differentiated
variable. Thus the derivatives are independent. A critical
point has all coordinates beyond \(A,B\) equal to zero, by
the last nonzero variable test on each tail. Its remaining
equations are \(A^3+\lambda B^3=v A^2=v B^2=0\), with the
same convention. There are finitely many projective solutions.
If both second variables exist, the hyperplane \(v=M(A,B)\)
cuts out the smooth shorter binary component with core
\(M(A,B)(A^3+\lambda B^3)\) from
Corollary~\ref{cor:nov30-two-attached}. If just one exists,
setting \(v\) equal to that variable produces a smooth
closed-ended chain. Both are K-stable whenever Fano by induction.
The one-cut hypersurface is smooth outside \(z=0\) by Bertini,
and its singular locus on \(z=0\) is contained in the cone
over the finite singular set of \(G\). It is smooth at \(p\)
by the simple root. Since \(k\ge4\), it is normal, using
Serre's criterion and the fact that different projective
hypersurface components would intersect in codimension one.
All the one-cut inputs thus hold at every core point.

At a tail vertex \(y_j\) suppress its tail subscript and
write \(l\) for this tail length and \(L\) for its attaching
form. The tangent coordinate is \(T=y_{j-1}\) with \(y_0=L\).
For \(j=l\) the tangent is a cone over a smooth shorter component,
without a cut. For \(j=l-1\) put \(y_l=uT\) to obtain a
smooth shorter component on the base. If \(j=1\) its core
is restricted to \(L=0\), so the other tail plus its one
surviving core variable is a closed-ended chain. The generalized
Eckardt calculation in these cases uses dimension \(b\ge n-1>r\).

For \(j\le l-3\) set \(y_{j+1}=0\). The one-cut hypersurface
is the disjoint sum of a shorter smooth component and an open
chain, whose sole projective singular point is its inner end.
It is normal and smooth at \(p\). On its tangent, the residual
form is the smooth component shortened by two further variables
and an open chain in \(y_{j+2},\ldots,y_l\) of length at
least two. Both summands have independent derivatives.
For \(j\ge3\) impose \(y_{j+2}=y_{j-2}\) to join the pieces;
for \(j=2\) impose \(y_4=L\). These yield a smooth shorter
component of the same type. At \(j=1\) the smooth residual
component is the previously described closed-ended chain on
the other tail and its one surviving core variable. Identify
\(y_3\) with the outermost variable of that closed-ended chain.
The combined chain is again closed-ended: the first surviving
term on the remote open chain uses the previously outermost
variable as coefficient of \(y_4^3\). This supplies a smooth
hyperplane restriction in the Fano induction, and the same
strict one-hyperplane estimate applies.

Finally take \(j=l-2\), set \(y_l=0\), and put
\(U=y_j, x=y_{j+1}, z=T\). The weighted initial equation is
\[
                    U^3 Z+U x^3+G(Y)=0
\]
in weights \(0,12,4,3\). At \(j\ge2\) the derivative
in \(z\) comes from \(zU^3\) and the residual \(G\) is the
smooth component with this tail shortened to length \(j-2\).
At \(j=1\) take \(z=L\); restricting the core and the other
attaching form to \(z=0\) gives the smooth closed-ended chain.
In either case every omitted monomial has larger weight.
The shortened section itself is smooth by the truncation rule,
and the Fano residual is K-stable by induction. This is
precisely the weighted one-cut calculation in the proof of
Theorem~\ref{thm:nov30-quartic-binary-end}. It gives
\(\delta_p(X,H)>r\). The cyclic fixed-point criterion
now proves the corollary.
\end{proof}

\subsection{Iterated cyclic covers in composite degree}
Let \(3\le d\le n\). Choose integers \(e_1,\ldots,e_k\ge2\)
with product \(d\), put \(a_0=1\), \(a_j=e_1\cdots e_j\),
and choose homogeneous forms \(q_j(x)\) of degree \(a_j\) in
\(x=(x_0,\ldots,x_n)\). Define
\[
 T_0=y,\qquad T_j=T_{j-1}^{e_j}+q_j(x),\qquad
 X=(T_k=0)\subset\mathbf P^{n+1}.
\]
For each \(j\) let \(D_j\subset\mathbf P^n\) be the degree-\(d\)
hypersurface given by setting \(T_{j-1}=0\) in the remaining
composition \(T_k\). In particular \(D_k=(q_k=0)\).
If \(X\) is smooth and all \(D_j\) are K-polystable, then \(X\)
is uniformly K-stable. Moreover, smoothness of \(X\) implies
smoothness of every \(D_j\).

For the proof introduce \(V_j\) by taking \(T_j\) as an independent
coordinate of weight \(a_j\) in the remaining composition, in
\(\mathbf P(1^{n+1},a_j)\). Thus \(V_0=X\) and \(V_k=\mathbf P^n\)
(its weight-\(d\) coordinate is zero). In every affine base chart
the map \(V_{j-1}\to V_j\) has equation \(T_{j-1}^{e_j}=T_j-q_j\).
The weighted vertex lies on none of the \(V_j\), since their equations
are monic in the weight-\(a_j\) coordinate. All maps are finite and
surjective. A singular point of an intermediate \(V_j\) would lift to
a singular point of \(X\), by the chain rule in a base chart. At a
point of \(D_j\) use \(T_{j-1}=0\) and lift further to \(X\). There
the derivative of \(T_{j-1}^{e_j}\) vanishes, so the differential of
the equation of \(X\) is precisely that of the degree-\(d\) equation
of \(D_j\). This proves the asserted smoothness (and also rules out
an identically zero branch equation).

Write \(A\) for pullback of \(\mathcal O_{\mathbf P^n}(1)\).
Adjunction gives \(-K_{V_j}=(n+1-d+a_j)A\); the branch divisor
\((T_j=q_j)\subset V_j\) has class \(a_jA\) and maps isomorphically
to \(D_j\). Its proportionality coefficient
\(\lambda_j=a_j/(n+1-d+a_j)\) lies strictly between zero and one.
Starting at \(V_k\), apply Zhang--Zhou, Theorem~1.7
\cite{ar2003.06839}, to \(V_j,D_j\). The pair with branch coefficient
\(1-1/e_j\) is K-polystable, since the difference from the theorem's
upper endpoint is
\[
  \left(1-\frac{\lambda_j^{-1}-1}{n}\right)
   -\left(1-\frac1{e_j}\right)
    =\frac{(a_{j-1}-1)n+d-1}{a_j n}>0.
\]
The crepant cyclic-cover equivalence of Liu--Zhu, Theorem~1.2(2)
\cite{ar2001.10557}, yields K-polystability of \(V_{j-1}\).
Descending to \(j=1\) gives the assertion by finiteness of
\(\operatorname{Aut}(X)\).

This supplies a range in arbitrarily large dimension for every
composite degree. Suppose all smooth degree-\(d\) hypersurfaces of
dimension \(m\ge d-1\) are K-polystable. In \(m+2\) base variables
\(x\) choose a degree-\(d\) form \(H\). For each of \(n-m\) additional
variables \(y_i\) choose any of the above compositions
\(\Psi_i(y_i,x)\) of total degree \(d\); their factorizations of \(d\)
can be different. Then every smooth hypersurface
\[
              H(x)+\sum_{i=1}^{n-m}\Psi_i(y_i,x)=0
\]
is uniformly K-stable. Indeed one may absorb \(H\) and all other
summands into the last shift \(q_k\) of the last variable. Each branch
hypersurface in dimension \(n-1\) has the same other summands and a
modified degree-\(d\) form in \(x\), and is smooth as proved above.
Induction on \(n-m\) proves the claim. One may in particular take
\(m=d\) by the published all-members index-two theorem.

In any composite degree the displayed family contains smooth members
with no pure-power coordinate after any linear change. Take the
summands \((y_i^{a_i}+p_i(x))^{b_i}\), \(a_i b_i=d\) and \(a_i,b_i\ge2\),
with \(p_i\) nonzero, and choose \(h=H+\sum p_i^{b_i}\) general.
For a pure-power coordinate there must be a nonzero direction
\(v=(v_x,v_y)\) in which the derivative of the equation is the
\((d-1)\)st power of a linear form up to scalar. Restrict to \(y=0\).
For \(v_x\ne0\) the derivative \(\partial_{v_x}h\) would belong to
the cone of powers (including zero). For fixed \([v_x]\) the derivative
on all degree-\(d\) forms is surjective; the codimension of this cone
is \(\binom{m+d}{d-1}-(m+2)>m+1\). The incidence over \([v_x]\)
therefore excludes this for general \(h\). For \(v_x=0\) the putative
linear form must depend only on \(y\), whereas the derivative has
the nonzero mixed terms \(d\,v_{y_i}p_i^{b_i-1}y_i^{a_i-1}\) for
the indicated nonzero entries of \(v_y\). General smooth examples
exist by openness from \(p_i=0\) and a smooth \(h\).

\paragraph{An explicit range in every even degree.}
Take $d=2a$, $a\ge2$. Suppose every smooth degree-$d$ hypersurface in
dimension $m\ge d-1$ is K-polystable. For $n\ge m$, assume that in some
coordinates each of $y_1,\ldots,y_{n-m}$ occurs with exponents
divisible by $a$. Every smooth $X_d^n$ of this form is uniformly
K-stable. In particular $m=d$ is available from the published
index-two theorem; one may also take the largest \(m\) in Theorem~\ref{thm:range} for the given degree.

To see this, take the last coordinate $y$. Its coefficient at
$y^{2a}$ is nonzero: otherwise the $y$-coordinate point would lie
on $X$ and all partial derivatives would vanish there, since
the divisibility excludes monomials $y^{2a-1}x_j$ (including
the other distinguished variables among the $x_j$). Thus write
the equation as $y^{2a}+p(x)y^a+h(x)$ after scaling. Completing
the square makes $q=p/2$ and $f=h-p^2/4$ in the preceding result.
Both $h=0$ and $h-p^2/4=0$ are smooth degree-$d$ hypersurfaces
in dimension $n-1$. Their equations retain divisibility by $a$ in
the remaining distinguished coordinates. Induct on $n-m$. At an
induction step $n\ge m+1\ge d$, as required.

\section{A destabilizing cone and its smooth fibers}
Let \(n\ge d\ge3\) and let \(B_d^{n-1}\subset\mathbf P^n\) be smooth
of degree \(d\), with equation \(f(x)=0\). Its projective cone
\(Y=(f(x)=0)\subset\mathbf P^{n+1}_{[x:z]}\) is a klt Fano hypersurface.
At its vertex the blowup resolves the hypersurface, with smooth exceptional
divisor \(E\simeq B\). Hypersurface adjunction gives
\(A_Y(E)=n+1-d\), which is positive. Write \(H=\mathcal O_Y(1)\).
A basis of \(H^0(Y,mH)\) splits into the spaces
\(z^{m-j}H^0(B,jH|_B)\) of exact vertex order \(j\).
Since \(h^0(B,jH|_B)=d j^{n-1}/(n-1)!+O((1+j)^{n-2})\), their
degree-normalized average tends to \(S_H(E)=n/(n+1)\). Thus
\[
 A_Y(E)-S_{-K_Y}(E)
  =n+1-d-(n+2-d)\frac{n}{n+1}=-\frac{d-1}{n+1}<0.
\]
The valuative criterion makes this cone K-unstable. When \(B\) is
K-semistable, the projective-cone formula of Zhang and Zhou
\cite[Corollary~1.6(1)]{ar2003.06839} shows that \(E\) computes
the threshold, so that
\[
       \delta(Y)=\frac{(n+1)(n+1-d)}{n(n+2-d)}.
\]
In the family
\(f(x)+t z^d=0\) every fiber for \(t\ne0\) is smooth. When \(f\) is the
Fermat polynomial, these fibers are themselves projectively equivalent
to the Fermat hypersurface and are K-stable by
\cite[Corollary~1.4]{ar2004.09413}. The vertex computation concerns
the singular central fiber.

\paragraph{A quantitative obstruction on smooth examples.}
For every fixed \(d\ge3\), the smooth Fermat examples also give a
dimension-dependent ceiling on any lower bound for all members. Take
\(p=[1:\zeta:0:\cdots:0]\) with \(\zeta^d=-1\). In homogeneous
coordinates \(T,z,y_1,\ldots,y_n\) choose \(T=1\) at \(p\) in the
usual affine chart and \(z=0\) as its tangent hyperplane. On this
hyperplane the equation is precisely \(\sum y_i^d=0\): the first
two coordinates are proportional there and their powers cancel.
The implicit function theorem makes \(y_i\) parameters at \(p\);
the function \(z\) has order \(d\), with leading form a nonzero
multiple of \(\sum y_i^d\). Assign orders \(0,d,1\) to \(T,z,y_i\)
respectively. In the ring of symbols for the point order and
ordinary section degree the initial relation is
\(aT^{d-1}z+\sum y_i^d\), \(a\ne0\). Its quotient is a domain
(the coefficients in the primitive linear relation are coprime)
and after inverting \(T\) the parameters give injectivity. Its
ordinary Hilbert series equals that of \(X\), so this presentation
computes the full order filtration in every section degree.
Its bigraded Hilbert series is therefore
\[
   \frac{1-t^d q^d}{(1-t)(1-tq^d)(1-tq)^n}.
\]
For the ordinary blowup \(E\) we obtain
\(S_H(E)=(n+d-1)/(n+1)\). Indeed logarithmically differentiating
in \(q\) at \(q=1\) multiplies the ordinary Hilbert series by
\((n+d)t/(1-t)-dt^d/(1-t^d)\), whose leading pole is
\((n+d-1)/(1-t)\). Comparing degree-\(m\) coefficients and
dividing by \(m\) gives the claimed mean. As \(A_X(E)=n\), this proves
\[
  1<\delta(X,-K_X)\le\frac{n(n+1)}{(n+d-1)(n+2-d)}
        =1+\frac{(d-1)(d-2)}{(n+d-1)(n+2-d)} .
\]
The strict lower bound uses K-stability of the Fermat hypersurface
and finiteness of its automorphism group. Thus a lower estimate
valid for every smooth hypersurface at fixed \(d\) can exceed one
by at most \(O_d(n^{-2})\) as \(n\) increases. This is also the
ordinary point valuation in the generalized Eckardt calculation
of Abban--Zhuang \cite[Theorem~4.18]{ar2003.13788}.

\chapter{Possible counterexamples in the parameter space}
\label{ch:new-counterexamples}
Write \(\mathcal P_{n,d}=\mathbf P H^0(\mathbf P^{n+1},\mathcal O(d))\)
and \(\mathcal U_{n,d}\subset\mathcal P_{n,d}\) for its smooth locus.
All codimensions in this chapter refer to equations up to scalar.
Let \(\mathcal K_{n,d}\subset\mathcal P_{n,d}\) consist of the
equations defining a uniformly K-stable klt Fano hypersurface.
For \(3\le d\le n+1\) this is open. Indeed the klt locus in the
flat projective hypersurface family is open, the relative canonical
class is Cartier, and K-stability is open in such Fano families
\cite[Theorem~4.5]{ar1907.02408}. Uniform K-stability is equivalent
to K-stability here \cite[Theorem~1.6]{ar2102.09405}. Thus the
smooth counterexample locus itself is relatively closed in
\(\mathcal U_{n,d}\).

\begin{theorem}\label{thm:new-catalan}
For every \(d\ge4\) and \(n\ge d-1\),
\[
   \operatorname{codim}_{\mathcal P_{n,d}}
       (\mathcal P_{n,d}\setminus\mathcal K_{n,d})
       \ \ge\ \operatorname{Cat}_{d-1}
       :=\frac1d\binom{2d-2}{d-1}.
\]
For \(d\ge5,6,13,31,110\), respectively, the right hand side can
be multiplied by \(2,5/2,3,10/3,7/2\), using the largest applicable
factor. For small degrees the following full-space bounds hold
uniformly in every \(n\ge d-1\):
\[
\begin{array}{c|rrrrrr}
d&4&5&6&7&8&9\\ \hline
\text{codimension at least}&9&34&106&346&1169&4016\\[2pt]
d&10&11&12&13&14&15\\ \hline
\text{codimension at least}&13998&49348&175650&630419&2278942&8290402
\end{array}
\]
The same lower bounds hold inside the smooth parameter space, with
empty complement allowed. In particular a
general projective linear system of dimension strictly smaller
than this bound consists entirely of uniformly K-stable
hypersurfaces, including its singular fibers.
\end{theorem}
The proof of the Catalan estimate begins below. The improved
factors follow in Theorem~\ref{thm:dec1-all-degree-catalan}.
Proposition~\ref{prop:dec1-index-one-neighborhood} supplies the
index-one step for quartics and quintics. For smooth quartic
ninefolds, Chapter~\ref{ch:nov30-replacements} proves that the
exceptional locus is empty. The direct local method also gives
intermediate codimension bounds. Imposing the point and line
conditions at the same marked point gives \(36\)
(Section~\ref{sec:nov28-quartic-ninefold-joint-sharp}).
Theorem~\ref{thm:dec1-q-final-ninefold} improves this to \(148\)
by examining the quartic term along the same line.
The mixed cubic and quartic derivatives in
Theorem~\ref{thm:dec1-q-mixed-normal-181} raise the smooth-locus
bound to \(181\). For comparison, in the full parameter space
the quartic ninefold bounds are \(9\le c_{9,4}\le211\).
The upper bound comes from equations with a single double point
in Section~\ref{sec:dec1-unstable-flags}.
For quintic elevenfolds, two intermediate incidence counts give
\(635\) and \(1067\) before the last line strata are analyzed.
The full argument of Theorem~\ref{thm:nov28-quintic-eleven}
excludes all smooth exceptions in that dimension.
The ninefold estimate uses the local stability arguments for
quartics. The dimension-independent estimates follow from
singularities of the branch equation and cyclic covers.
A component of the complement
of the smooth K-stable locus, if it exists, must satisfy the same
inequality. The discriminant of singular equations has codimension
one, and its general point represents a uniformly K-stable singular
hypersurface by the theorem.
There is a version through a prescribed member. Fix any uniformly
K-stable equation \([F]\) and let \(c\) be a full-space lower
bound above. Projection from \([F]\) maps the closed exceptional
locus to a closed set of dimension at most
\(\dim\mathcal P_{n,d}-c\). For \(1\le k<c\), a general
\((k-1)\)-plane of directions misses its image by the projective
dimension count. The corresponding \(k\)-dimensional linear
system through \([F]\) has uniformly K-stable fibers throughout.
The proof of this theorem uses neither the special-fiber slope
estimates nor the weighted short-curve estimates of the preceding
chapters. The next two arguments give its initial degree estimate
and its propagation in dimension.

\section{Codimension from the tangent forms of a divisor}
\label{sec:new-tangent-catalan}
For \(s\ge2\) and \(e\ge2\), call a degree-\(e\) form on
\(\mathbf P^s\) good if \((\mathbf P^s,V(f))\) is purely log
terminal (plt). In particular its divisor is normal and klt.
Conversely a normal klt hypersurface divisor is plt in the smooth
ambient space by inversion of adjunction. These equivalences and
openness of klt singularities in a flat Gorenstein family are
used also at a marked point. Include the zero form in the closed
cone of forms that fail the condition. Write \(K(s,e)\) for a
lower bound on its affine codimension.

\begin{lemma}\label{lem:new-tangent-recursion}
One may use \(K(2,e)=1\) for all \(e\ge2\). For \(s\ge3\) one
may use any bounds obtained recursively by
\[
    K(s,e)=1+\sum_{h=2}^{\min(e,s-1)} K(s-1,h).
                                                        \tag{9}
\]
In particular \(K(s,e)\ge\operatorname{Cat}_{s-1}\) for
\(e\ge s-1\).
\end{lemma}

\begin{proof}
A smooth plane curve suffices for the initial value. Fix a point
\(p\) of \(\mathbf P^s\). In an affine chart centered at \(p\)
write the equation as a sum of homogeneous polynomials \(q_h\)
of degrees \(0,\ldots,e\) in \(s\) variables. These coefficients
are independent. For each \(2\le h\le\min(e,s-1)\) choose a
vector space \(V_h\) of homogeneous forms of degree \(h\), of
dimension \(K(s-1,h)\), in which every nonzero member is good
on \(\mathbf P^{s-1}\). Such a space exists since a general
projective space of dimension \(K(s-1,h)-1\) avoids the closed
bad locus.

Consider the direct sum of the constants, all linear forms, and
the spaces \(V_h\). Every nonzero member defines a pair plt at
\(p\). This is immediate if the constant or linear term is the
first nonzero term. Otherwise let \(q_h\in V_h\) be its first
nonzero term. The homogeneous cone \(q_h=0\) is klt at its vertex:
its blowup is the total space of \(\mathcal O(-1)\) on the klt
projective hypersurface, and the relative canonical coefficient
at the exceptional divisor is \(s-1-h\ge0\). The cone is normal,
by the hypersurface criterion and normality of its punctured cone.
In the family \(t^{-h}F(tx)\) the fiber at zero is precisely this
cone. Near the marked section this is a flat family of hypersurfaces.
Openness of klt singularities applies to the central fiber, and
all nonzero fibers are isomorphic at the mark by rescaling.
This proves the assertion about the direct sum.

The forms whose pair fails to be plt at the fixed point form a
closed cone: on nonzero projective forms apply the marked openness
just used together with the open set where the point misses the
divisor. The projectivized bad cone misses the projectivization of
the direct sum. The projective dimension theorem bounds its affine
codimension below by \(1+s+\sum_h K(s-1,h)\). Allowing the point to
move costs at most \(s\), which proves (9).

To solve the basic recursion put \(a=s-1\) and
\(j=\min(e-1,a)\). It gives the ballot number
\[
              \frac{a-j+1}{a+1}\binom{a+j}{j}.
\]
Indeed its value at \(j=0\) is one, adjacent differences in a
row give the preceding row until \(j=a-1\), and its values at
\(a-1,a\) agree. These are exactly (9) and its truncation.
The last two entries both equal \(\operatorname{Cat}_{s-1}\).
\end{proof}

\section{Propagation by cyclic covers}
\label{sec:new-propagation}
In degree \(d=s+1\) a good form \(f\) is GIT stable. To check the
strict Hilbert--Mumford inequalities choose coordinates diagonalizing
any nontrivial one-parameter subgroup. Subtract its least weight
from all weights, obtaining \(w_i\ge0\), at least one zero, and let
\(m\) be the smallest weight of a monomial of \(f\). The corresponding
monomial order at the generic point of the weight-zero coordinate
space has log discrepancy \(\sum w_i\) on \(\mathbf P^s\) and
value \(m\) on \(D=(f=0)\). One may divide the weights by their
greatest common divisor. The order formula follows on an affine
chart of a weight-zero variable by using the remaining weight-zero
ratios as independent coordinates of the center. The plt hypothesis
gives \(\sum w_i>m\): its only divisorial equality is
\(\operatorname{ord}_D\), whereas this order has linear center
and \(D\) is integral of degree \(d>1\). Since \(d=s+1\), this
strict inequality is precisely the Hilbert--Mumford test with
weights normalized to sum zero.

\subsection{Log canonicity and the Hilbert--Mumford estimate}
\label{sec:nov28-lc-git}
In the affine space of homogeneous degree-\(e\) forms on \(\mathbf P^s\),
write \(b(s,e)\) for the codimension of the cone consisting of zero and
forms with non-log-canonical divisor. The bad cone is closed: log canonical thresholds in this algebraic
family are lower semicontinuous and take finitely many values, by
stratification and log resolution in characteristic zero.
For \(s=1\), \(b(1,e)\ge1\).
For \(s\ge2\), the marked-point construction in the preceding recursion
also gives
\[
 b(s,e)\ge 1+\sum_{h=2}^{\min(e,s)}b(s-1,h).
\]
Indeed the blowup of the vertex for a homogeneous initial equation of
degree \(h\) in \(s\) variables has relative pair coefficient
\(h-s+1\) on the exceptional divisor. Its strict transform is pulled
back from the projective equation on \(\mathbf P^{s-1}\) in the total
space of \(\mathcal O(-1)\). Log canonicity there, including the zero
section with coefficient one, follows locally by taking a product
with the affine line. Since \(h\le s\), this proves log canonicity
at the vertex. Semicontinuity of log canonical thresholds under
scaling transfers the assertion to the marked germ with this initial
form. The constant and linear coefficients together have dimension
\(s+1\), and moving the marked point costs \(s\) dimensions. This
proves the displayed recurrence with the projectively disjoint linear
test spaces for each homogeneous order. It gives \(b(s,s+1)\ge
\operatorname{Cat}_s\), using the same Catalan triangle shifted by
one row. For \(e=2\) the exact value is \(\binom{s+1}{2}\), since a
nonzero quadric gives a log canonical divisor precisely in ranks at
least two.

Here is a separate estimate for the GIT locus in the anticanonical
degree. Write \(C_j=\operatorname{Cat}_j\). The non-GIT-stable
forms of degree \(d\) on \(\mathbf P^{d-1}\) have codimension at least
\[
  g_d=\min_{1\le j<d}\{C_d-C_j C_{d-j}-j(d-j)\}.
\]
Indeed order the weights of a nontrivial one-parameter subgroup
\(w_1\ge\cdots\ge w_d\), with sum zero. Write \(\alpha_i\) for
the exponents of a monomial of degree \(d\) and set
\(u_j=\sum_{i=1}^j(\alpha_i-1)\). Summation by parts gives
\[
       \sum_i w_i\alpha_i=\sum_{j<d}(w_j-w_{j+1})u_j.
\]
The number of nonnegative compositions with all \(u_j\ge0\) is
\(C_d\): draw \(\alpha_i\) vertical steps followed by one horizontal
step for each index \(i\), giving precisely the lattice paths above
the diagonal in the \(d\)-square. Within these compositions the
weight is strictly positive unless the path meets the diagonal at
every boundary of its constant-weight blocks. If their sizes are
\(m_1,\ldots,m_l\), the number of these exceptional paths is
\(\prod_i C_{m_i}\). For this sign pattern, the nonstable forms
vanish on all positive-weight coefficients. The flags are
parametrized by a variety of dimension \(\sum_{i<j}m_i m_j\).
Only finitely many monomial sign patterns occur.

It remains to maximize \(\prod_i C_{m_i}+\sum_{i<j}m_im_j\) over
partitions with at least two blocks. Merging adjacent blocks of
sizes \(a,b\) cannot decrease this expression. Indeed
\(C_{a+b}-C_aC_b\ge ab\): in the centered exponent notation the
paths with a single \(+1\) in the first block and a single \(-1\)
in the second give \(ab\) paths missing the return at their common
boundary. The other Catalan factors are at least one and the
flag dimension loses exactly \(ab\). It therefore suffices to
consider the two-block partitions in the displayed formula.

In particular \(g_d\ge C_{d-1}\) for \(d\ge4\).
The increasing Catalan ratios give \(C_jC_{d-j}\le C_{d-1}\),
and \(j(d-j)\le d^2/4\). For \(d\ge5\) use
\(C_d-2C_{d-1}=2(d-2)C_{d-1}/(d+1)\ge d^2/4\),
which holds at five and persists by the consecutive Catalan ratio.
For \(d=4\) evaluate the two distinct partitions to get \(g_4=6\).
Consequently the complement of the locus of GIT stable forms with
log canonical divisor has codimension at least \(C_{d-1}\) for every
\(d\ge4\). The suspension and affine intersection argument gives
this bound in every Fano dimension \(n\ge d-1\). The stronger
input is the minimum of \(g_d\) and the recursively obtained lower
bound for \(b(d-1,d)\) using the exact quadratic entry above.

A stronger version treats the two-block case by its exact
monomial count. For a single break at \(j\) the weight is positive
exactly when \(k=\alpha_1+\cdots+\alpha_j>j\). Thus set
\[
 H(d,j)=\sum_{k=j+1}^d
   \binom{k+j-1}{j-1}\binom{2d-k-j-1}{d-j-1}-j(d-j).
\]
For three or more blocks the product of their Catalan numbers is
at most \(C_{d-2}\), by merging down to three blocks and using
the increasing successive Catalan ratios. The maximum for three
positive sizes is \(C_{d-2}C_1^2\). The flag dimension is at
most \(\binom d2\). Combining the two cases gives the lower bound
\[
 g_d'\ =\min\left\{C_d-C_{d-2}-\binom d2,\quad
                      \min_{1\le j<d}H(d,j)\right\}.
\]
One may take the larger of \(g_d\) and \(g_d'\) in every degree.

\paragraph{The exact cubic log-canonicity count.}
For every \(s\ge2\), the cubic entry in the global recurrence is
\[
                 b(s,3)=\binom{s+3}{3}-3s-2.
\]
A cubic admitting a reduced plane section with at most ordinary
nodes gives a log canonical divisor. Indeed, after a coordinate
change, scale all variables complementary to the plane to zero.
The limiting equation defines its projective cone, which is log
canonical. Transverse to the vertex this is the affine cone in
three variables over the log canonical plane divisor. Its blowup
has exceptional coefficient one and strict transform pulled back
from that divisor. Away from the vertex the claim follows by smooth
pullback. Semicontinuity transfers the bound to the original form.

For \(s\ge3\) a general plane section of a normal cubic is smooth
by Bertini and the codimension of its singular locus. If an integral
cubic is nonnormal, a general plane section is still integral, and
its arithmetic genus one allows at most one singular point. A
codimension-two component of the ambient singular locus has degree
one, since a general plane meets that component in its degree many
distinct singular points. Writing \(K=(x=y=0)\) for this linear
space and \(u\) for coordinates on it, the cubic is
\[
 A(u)x^2+B(u)xy+C(u)y^2+g_3(x,y)
\]
with \(A,B,C\) linear. A general plane section has an ordinary node
at its unique singular point if the binary quadratic there has rank
two. Otherwise the linear span of the triples \((A(u),B(u),C(u))\)
consists of binary quadrics of rank at most one. The projective
conic of squares contains no line, so its underlying vector space
has dimension at most one. Thus the cubic depends on at most three
linear variables. The analogous assertion for an integral cubic
in the starting plane itself is immediate.

For a product \(LQ\) with \(Q\) an irreducible quadric, a general
plane cuts a smooth conic and a transverse line provided
\(\operatorname{rank}(Q|_{L=0})\ge2\). Otherwise \(Q=LA+B^2\)
for linear forms \(A,B\), allowing \(B=0\). Again \(LQ\) depends
on at most three variables. A product of independent linear forms
has simple normal crossings; the other products of three linear
forms already depend on at most three variables. Hence every
non-log-canonical cubic is defined by a ternary form, and that
ternary form must itself be non-log-canonical by the cone argument.

The non-log-canonical ternary cubics have affine codimension two.
Indeed the discriminant is an irreducible divisor whose general point
is nodal, by the incidence of singular points and the freely variable
quadratic term. The cuspidal locus has codimension two: prescribing
a singular point costs three conditions and rank one of its tangent
quadratic costs one more. Varying the point adds two parameters,
and the nonzero cubic coefficient in the kernel and integrality
are open conditions realized by a cuspidal plane cubic. Therefore
this exceptional ternary space has dimension eight. Varying the
three variables costs \(3(s-2)\) parameters on the Grassmannian,
giving the lower bound claimed for \(b(s,3)\). For the converse,
a cuspidal ternary cubic has no constant derivative equal to zero,
so the three-space of its essential variables is recoverable as
the annihilator of its constant zero derivatives in the full space.
Away from its cone vertex it retains the transverse non-log-canonical
cusp by smooth pullback. Consequently its cones have the full
dimension \(8+3(s-2)\).

\paragraph{Codimension under linear sections.}
\label{par:new-lc-linear-sections}
Put \(T(s,e)=\binom{s+e}{e}\). For \(s\ge t\ge e-1\),
\[
     b(s,e)\ge\left\lceil \frac{T(s,e)}{T(t,e)}b(t,e)\right\rceil .
\]
Indeed, if a nonzero equation on one \(t\)-plane is log canonical,
scale the complementary variables to zero. The central equation is
a cone over its restriction. At the vertex the transverse affine
cone has exceptional coefficient \(e-t\le1\), and its strict
transform is pulled back from the log canonical pair on the base.
It is log canonical by the line-bundle calculation. Away from the
vertex this follows by smooth pullback. Thus non-log-canonical forms
restrict to a non-log-canonical form or to zero on every \(t\)-plane.

Fix a smooth point \(F\) of an irreducible component \(Z\) of the bad
cone in \(\operatorname{Sym}^e(V^*)\), where \(\dim V=s+1\), and write
\(N\subset\operatorname{Sym}^e V\) for its conormal vector space.
For any \((t+1)\)-dimensional \(U\subset V\) the restriction map on
\(Z\) has image dimension at most \(T(t,e)-b(t,e)\). At a smooth
point the differential rank is at most that image dimension: a
nonvanishing differential minor gives the same bound from below
on the dimension of the image by the complex rank theorem. Dualizing
that differential gives
\[
      \dim(N\cap\operatorname{Sym}^e U)\ge b(t,e).
\]
Equip \(V\) with a Hermitian form and \(\operatorname{Sym}^e V\) with
the induced form. The trace of the product of the orthogonal projections
to \(N\) and \(\operatorname{Sym}^e U\) is at least the dimension of
the intersection, as follows by extending an orthonormal basis of
the intersection to one on \(\operatorname{Sym}^e U\). Averaging the
second projection over the unitary group gives
\((T(t,e)/T(s,e))\mathrm{id}\): the standard symmetric power is
irreducible, and its trace determines the scalar by Schur's lemma.
Taking the trace against the first projection proves the claim.

The same cone argument also controls the dimension of the bad
locus inside one hypersurface. Write
\(t_a=\lceil e/a\rceil-1\). For \(s\ge t_a\), any nonzero degree-\(e\)
form on \(\mathbf P^s\) whose pair with coefficient \(1/a\)
fails to be lc has a non-lc locus of dimension at least \(s-t_a\).
Indeed, otherwise a general \(t_a\)-plane avoids this locus. Its
restricted pair is lc by Bertini on the lc locus: intersect a log
resolution with general hyperplanes, transverse to every stratum
of its normal-crossing boundary, and restrict the discrepancy
formula. The restricted equation is nonzero. Degenerating the
original form to the cone over this restriction gives an lc central
fiber, since \(e/a\le t_a+1\). This contradiction proves the claim.
Consequently one can add \(\max\{0,s+1-\lceil e/a\rceil\}\) to
the marked-point lower bounds for \(b_a(s,e)\) after subtracting
the \(s\) parameters for the point. To justify the improvement,
take the whole non-lc incidence before projecting to the parameter
space; every fiber over a nonzero image point has the just-proved
minimum dimension. Apply this correction once to the best bound
for the entire marked incidence. The zero equation itself has
full ambient codimension.

\paragraph{Stratifying the quadratic term.}
For \(a=1,2\), let \(b_a(s,e)\) denote the affine codimension of
zero together with the degree-\(e\) forms on \(\mathbf P^s\) whose
log canonical threshold is strictly less than \(1/a\). Thus
\(b_1(s,e)=b(s,e)\). Write \(T(s,e)=\binom{s+e}{e}\).
The following estimate holds for \(s\ge3\) and \(e\ge3\):
\[
\begin{split}
b_1(s,e)\ \ge\ 1+\min\left\{
  T(s-1,2)+\sum_{h=3}^{\min(e,s)} b_1(s-1,h),\quad
  \binom{s}{2}+\sum_{h=3}^{\min(e,2s-2)}b_2(s-2,h)\right\}.
\end{split}
\]
Fix a marked point and write the affine equation as
\(F_0+\cdots+F_e\) in its Taylor coefficients there. A non-lc germ
requires \(F_0=F_1=0\). If \(F_2\) has rank at least two, the
quadratic cone is lc, hence so is the original germ by scaling and
semicontinuity. If \(F_2=0\), its vanishing costs \(T(s-1,2)\)
conditions. For each \(3\le h\le \min(e,s)\), take a homogeneous
linear space of dimension \(b_1(s-1,h)\) meeting the bad cone only
at zero. General projective linear sections supply these spaces.
In their direct sum every nonzero vector has an lc initial equation
by the cone calculation. The bad set of higher terms is invariant
under \(F_h\mapsto t^{h-2}F_h\), so each component contains the
origin. The affine dimension theorem, applied to that direct sum,
gives the first bound within the marked incidence.

If \(F_2\) has rank one, the stratum of possible nonzero quadrics
has codimension \(\binom{s}{2}\). Fix one and write it as \(x^2\).
The higher terms vary freely and their bad locus again contracts
to zero. Use linear test spaces in the \(s-1\) complementary
variables \(y\), this time with dimension \(b_2(s-2,h)\) for
each \(3\le h\le\min(e,2s-2)\). In their direct sum any nonzero
\(g(y)\) has an initial form of degree \(h\) with threshold at least
\(1/2\) projectively. The transverse affine cone also has threshold
at least \(1/2\), since \(h/2\le s-1\). By semicontinuity the
same holds for \(g\) at the origin. The Thom--Sebastiani formula
for log canonical thresholds gives
\[
 \operatorname{lct}_0(x^2+g(y))
 =\min\{1,\tfrac12+\operatorname{lct}_0(g)\}=1 .
\]
For example, the formula follows from additivity of minimal
exponents in disjoint variables (M.~Saito,
\emph{On microlocal \(b\)-function}, Bull.\ Soc.\ Math.\ France
122 (1994), 163--184, Theorem~0.8) and the identity
\(\operatorname{lct}=\min\{1,\widetilde\alpha\}\).
Applying the dimension theorem in this fiber proves the second
bound. The constant and linear terms cost \(s+1\) conditions,
and moving the marked point subtracts \(s\).

The auxiliary estimates require only the same argument with
coefficient \(1/2\). For \(e\ge3\) and \(s\ge2\) they give
\[
 b_2(s,e)\ge 1+T(s-1,2)
      +\sum_{h=3}^{\min(e,2s)}b_2(s-1,h).
\]
Indeed an equation of multiplicity at most two has threshold at
least \(1/2\), by the multiplicity bound for thresholds. The
constant, linear and quadratic terms must all vanish. In higher
degrees the cone calculation applies provided \(h/2\le s\).
The starting entries are \(b_2(1,e)=2\) for \(e\ge3\) and
\(b_2(s,3)=T(s,3)-(s+1)\) for all \(s\ge1\). For the latter
identity, a nonzero pure cube has threshold \(1/3\). Every cubic
that is not a pure cube has a binary linear restriction that is not
a pure cube, by taking a general line transverse to the reduced
support and avoiding intersections between distinct components.
This binary form has projective threshold at least \(1/2\), as does
its affine cone because \(3/2\le2\). Degenerating the original cubic
to that restriction proves the assertion. The cone of pure cubes
has dimension \(s+1\).

One additional entry used below is
\(b_2(s,4)=T(s,4)-(2s+1)\) for \(s\ge2\). A quartic
with some line section whose binary quartic has no root of
multiplicity at least three has threshold at least \(1/2\), by
the following explicit degeneration. Choose coordinates in which
the line has parameters \(u,v\). Scale all the other variables by
\(t\). For \(t\ne0\) the resulting equation is projectively
equivalent to the original quartic; at \(t=0\) it is the product
\(\prod_i L_i(u,v)^{m_i}\) of binary linear factors of
multiplicity \(m_i\le2\). Away from their common vertex, all
boundary coefficients \(m_i/2\) are at most one. Transversely to
the vertex blow up the origin of \(\mathbf A^2\). The exceptional
coefficient for the pulled-back pair is \(\sum_i m_i/2-1=1\).
The strict transforms meet that exceptional curve at distinct
points, transversely, so the boundary has simple normal crossings
with all coefficients at most one. This proves log canonicity of
the central fiber, and semicontinuity proves the threshold bound
for the general fiber. For a general line the root multiplicities are the
multiplicities of the irreducible components of the plane curve,
by Bertini and transversality. Thus the nonzero exceptions are
exactly the forms \(L^3M\) with \(L,M\) linear. Their affine cone
has dimension \(2s+1\), giving in particular \(b_2(2,4)=15-5=10\).

\paragraph{The exact quartic log-canonicity seed.}
The cone of non-log-canonical quartic surfaces in \(\mathbf P^3\),
including zero, has codimension \(b(3,4)=9\). For the lower bound
we need only refine the rank-one case of the quadratic stratification.
Fix \(F_2=x^2\), so \(F=x^2+F_3+F_4\) in local variables \(x,y,z\).
By the analytic splitting lemma \(F\) is equivalent to \(x^2+g(y,z)\).
In particular \(g(y,z)=F(\phi(y,z),y,z)\), where \(F_x(\phi,y,z)=0\)
and \(\phi\) has order at least two. Its cubic term is
\(g_3=F_3(0,y,z)\). Non-log-canonicity requires
\(\operatorname{lct}_0(g)<1/2\).
If \(g_3\ne0\), it must be a cube, since every other binary cubic
has affine threshold at least \(1/2\). This restriction costs two
conditions, and we can fix \(g_3=y^3\). In succession the coefficients
\[
          [z^4]g_4,\qquad [yz^3]g_4,\qquad [z^5]g_5
\]
must vanish. Otherwise the respective weighted initial polynomials
are \(y^3+az^4,\ y^3+byz^3,\ y^3+cz^5\) with displayed coefficient
nonzero. In the middle case use weights \((3,2)\) on \((y,z)\).
These three polynomials have respectively thresholds
\(7/12,\ 5/9,\ 8/15\), from the weighted blowup of the origin and
smoothness off it. In each case weighted degeneration gives
\(\operatorname{lct}_0(g)>1/2\).

The first two equations are independent because their pure quartic
coefficients in \(F_4(0,y,z)\) are freely variable with coefficient
one. To see that the third costs a further condition, put
\[
 A=[xz^2]F_3,\qquad B=[x^2z]F_3,\qquad E=[xz^3]F_4 .
\]
Writing \(\phi(0,z)=-Az^2/2+O(z^3)\) in \(F\) gives
\([z^5]g_5=A^2B/4-AE/2\).
Here the terms involving the coefficient of \(z^3\) in \(\phi\)
cancel between \(x^2\) and \(Axz^2\).
The last equation therefore forces \(A=0\), or determines \(E\)
when \(A\ne0\), independently of the first two equations. Thus the
rank-one case costs at least \(1+3+2+2+1=9\), including the net
cost one for marking a singular point. If \(g_3=0\), its vanishing
costs four conditions, and the quartic \(g_4\) is freely variable
through \(F_4(0,y,z)\). If it has no root of multiplicity three
or more, its threshold at the origin is at least \(1/2\). The
exceptional binary quartics, including zero, have codimension two,
which gives the asserted lower bound in this case as well.

For \(F_2=0\), the six quadratic coefficients must vanish and
the ternary cubic tangent must belong to the non-lc cone, of
codimension two. This costs \(1+6+2=9\). To attain equality, take
a marked triple point whose tangent cubic is a cuspidal irreducible
plane cubic and general quartic terms. The blowup of the point has
exceptional divisor with pair coefficient one, whose adjunction
pair contains that cuspidal cubic with coefficient one. Adjunction
therefore prevents log canonicity. The affine family at a fixed
point has codimension \(4+6+2=12\). Its general member is smooth
off that point by Bertini, since the freely variable quartic terms
in the transverse variables are base-point free there. Moving the
point yields the asserted codimension nine.

\paragraph{The quintic seed.}
For quintic threefold equations in \(\mathbf P^4\), the non-lc cone
has codimension at least \(31\), that is, \(b(4,5)\ge31\).
The rank-one quadratic stratum already gives
\[
  1+6+b_2(2,3)+b_2(2,4)+b_2(2,5)\ \ge\ 1+6+7+10+8=32.
\]
Suppose now \(F_2=0\), at net cost \(1+10=11\) including the marked
singular point. Write \(C=F_3\) in four local variables. The cone
test gives log canonicity if the cubic surface defined by \(C\) is
lc. The classification in the cubic calculation leaves the six
possibilities below. The middle column is their affine codimension
inside the 20-dimensional space of cubic forms. The last column is
a lower bound for the number of additional conditions on \(F_4,F_5\)
in order for the germ to be non-lc.
\[
\begin{array}{l|r|r}
 C & \operatorname{codim}&\text{additional conditions}\\ \hline
 \text{ternary cusp} &9&11\\
 \text{ternary tangent line and smooth conic}&10&10\\
 \text{squarefree binary cubic}&12&8\\
 \text{binary cubic with a double and a simple factor}&13&7\\
 \text{pure cube}&16&4\\
 0&20&0
\end{array}
\]
For ternary cubics there are three parameters for the essential
variable subspace. Cuspidal forms in three variables have affine
dimension eight; a tangent line and conic have dimension seven.
For binary cubics the variable subspace contributes four
parameters, and the squarefree and double-root loci have dimensions
four and three. Pure cubes form a four-dimensional cone.

To check the last column, fix \(C\). The following are linear
spaces of possible quartic and quintic corrections. Each of their
nonzero elements gives an lc germ when added to \(C\). The bad cone
for these corrections is invariant under
\((F_4,F_5)\mapsto(tF_4,t^2F_5)\), so the affine dimension theorem
gives their dimensions as lower bounds.

Both ternary cases admit coordinates \(x,y,z,w\) in which
\(C=x^2z+y^3\) or \(C=x^2z+xy^2\), with unused variable \(w\).
Use the independent corrections
\[
     g_4(z,w)\in w\operatorname{Sym}^3(z,w),\qquad
     y\,g_3(z,w)\in y w\operatorname{Sym}^2(z,w),\qquad
     g_5(z,w)\in w^2\operatorname{Sym}^3(z,w).
\]
Their dimensions are \(4,3,4\). For the tangent case restrict
the quintics to the three monomials with \(2\le\deg_w\le4\).
First suppose \(g_4\ne0\). Give
\(z,w\) weight \(1/4\) and \(x\) weight \(3/8\). Give \(y\) weight
\(1/3\) for the cusp and \(5/16\) for the tangent case. The weighted
initial form is \(C+g_4\). If \(g_4=0\) and \(yg_3\ne0\), give
\(z,w\) weight \(2/9\) for the cusp and \(3/13\) for the tangent
case. Set the \(x,y\) weights to \((1-q)/2,1/3\) or
\((1-q)/2,(1+q)/4\), respectively, with \(q\) the common weight of
\(z,w\). The initial form is \(C+yg_3\). If just \(g_5\ne0\), take
\(q=1/5\) and the respective \(x,y\) weights from the same formula.
In every case the retained terms have weight one and the discarded
ones have weight strictly greater than one. The retained correction
depends on \(w\). For \(g_k(z,w)\) let \(M\ge1\) be the largest
\(w\)-exponent present, and for \(yg_3\) use the \(w\)-exponent in
\(g_3\). To retain only the extremal monomial and the two terms of
\(C\), choose any
rational \(\rho_z>0\), solve \(2\rho_x+\rho_z=0\) and either
\(3\rho_y=0\) or \(\rho_x+2\rho_y=0\). Prescribe \(\rho_w\) by
making the extremal monomial have weight zero. Then
\(\rho_w<\rho_z\), so all other correction terms have positive
weight. This gives an actual polynomial degeneration even if some
of the variable weights are negative.

Here and below the retained monomial gives a nondegenerate Newton
polynomial. Indeed its exponent has a positive power of \(w\), and
the two distinct exponent vectors of \(C\) are independent and have
zero \(w\)-coordinate. Every face polynomial therefore has linearly
independent exponent vectors, which precludes a torus critical
point. The Newton polyhedron formula for the log canonical threshold
(or the corresponding toric log resolution) implies the following
test: if exponent vectors \(a_i\) admit real coefficients
\(\lambda_i\ge0\) with \(\sum\lambda_i=1\) and
\(\sum\lambda_i a_i\le(1,\ldots,1)\) coordinatewise, the polynomial
is lc. These inequalities put the all-ones vector in the Newton
polyhedron and bound every toric order by the toric log discrepancy.
The coefficients \((\lambda_1,\lambda_2,\lambda_3)\) for the two
monomials of \(C\) and the extremal correction are, for \(g_4,g_5\),
\((1/2,1/3,1/6)\) for the cusp and \((1/4,1/2,1/4)\) for the
tangent case. For \(yg_3\) use \((1/2,1/4,1/4)\) in the cusp
and \((1/3,1/3,1/3)\) in the tangent case. The stated bounds
on \(M\) give the required coordinatewise inequalities.

For the squarefree binary case normalize \(C=x^3+y^3\), with
unused \(z,w\). Use a three-dimensional general linear space of
binary quartics on \(z,w\) with no nonzero member having a
quadruple root. Use also \(x\) times a general two-dimensional
space of binary cubics with no nonzero cube, and a three-dimensional
general space of binary quintics with no nonzero root of order four.
Projectivized dimension counts give those spaces. Take their direct
sum in the higher terms, of dimension eight. Give \(x,y\) weight
\(1/3\). If the quartic is present, give \(z,w\) weight \(1/4\).
The weighted initial is \(C+g_4\), and
\(\operatorname{lct}_0(g_4)\ge1/3\) by the binary cone calculation.
If just the quintic and \(xg_3\) remain and \(xg_3\ne0\), take
weights \(2/9\) on \(z,w\). The initial form is \(C+xg_3\).
The three binary cubic factors have multiplicities at most two,
so take a root of multiplicity \(m\in\{1,2\}\), move it to
\(z=0\), and degenerate to \(z^m w^{3-m}\), keeping \(x,y\)
fixed. This follows by giving \(z,w\) the weights \(3-m,-m\).
The three exponent vectors of \(x^3+y^3+xz^mw^{3-m}\) are
independent; their convex combination with coefficients
\((2/9,1/3,4/9)\) is coordinatewise at most one. This gives lc
by the Newton test. Finally \(C+g_5\) is lc by Thom--Sebastiani and the binary
quintic threshold bound \(1/3\).

For the double-root binary type write \(C=x^2y\). For a correction
\(g(y,z,w)\), the morphism \(X=xy\) to the affine space with
coordinates \(X,y,z,w\) pulls the equation \(X^2+yg\) back to
\(y(x^2y+g)\); its Jacobian has divisor \(y=0\). The two ambient
pairs with coefficient-one equations are therefore crepant. Thus
it suffices for \(X^2+yg\) to be lc, which follows from
\(\operatorname{lct}_0(yg)\ge1/2\) by Thom--Sebastiani.

There is a five-dimensional linear space of ternary quartics
\(g(y,z,w)\) for which every nonzero member satisfies this last
condition. Indeed the affine cone degree is five in three
variables, so it suffices to check the plane pair
\((\mathbf P^2,\tfrac12[(y=0)+(g=0)])\).
Away from the line \(y=0\) the only possible non-lc quartics at
coefficient one half are those with a triple line, a locus of
codimension ten by the binary-section proof of \(b_2(2,4)=10\)
above. At a fixed
point on \(y=0\) choose parameters \(y,t\). If \(g\) has multiplicity
zero or one, the multiplicity bound gives the result. Multiplicity
at least three costs six conditions. For multiplicity two, \(yg_2\)
must be a cube, since every other binary cubic has threshold at
least one half; hence \(g_2\) is proportional to \(y^2\), at a cost
of two conditions besides the three constant and linear equations.
The coefficient of \(t^3\) in \(g_3\) must also vanish. Otherwise
the weighted initial form is \(y(ay^2+bt^3)\) for weights
\((3,2)\), of threshold \(5/9>1/2\). We have at least six
conditions at a fixed point and at least five after moving the
point along the line. A general five-dimensional linear subspace
therefore has the stated property.

For quintic corrections use a two-dimensional space
\(g_5=yh_4(z,w)\) in which each nonzero binary quartic \(h_4\)
has no triple root. Such a space exists since triple roots have
codimension two. The product \(yg_5=y^2h_4\) has threshold at
least \(1/2\) by the binary cone calculation. In a mixture with
a nonzero quartic, give \(y,z,w\) weight \(1/4\) and \(x\) weight
\(3/8\). The quintic drops from the initial form. This proves the
claimed seven extra conditions.

Finally take \(C=x^3\). Use a four-dimensional linear space of
ternary quartics \(g(y,z,w)\), each nonzero member having threshold
at least \(2/3\). Then \(x^3+g\) is lc. To see such a space exists,
the affine cone bound \(4(2/3)\le3\) reduces the matter to the
plane quartic. At a fixed point, multiplicity at least three costs
six conditions. In multiplicity two, use the splitting lemma:
a rank-two tangent gives lc; a rank-one tangent, costing one
further condition after the three constant and linear conditions,
admits coordinates in which the germ is \(a^2+h(b)\). The
coefficients \([b^3]h\) and \([b^4]h\) must both vanish if the
threshold is less than \(2/3\), since otherwise it is at least
\(1/2+1/4\). Each condition has its respective pure coefficient
in the original cubic or quartic jet as an independent variable.
Moving the point subtracts two, leaving codimension at least four.

The minimum of the two columns in the sense of their rowwise sum
is twenty. Restoring the eleven conditions for the marked point
and vanishing Hessian proves the quintic seed.

\paragraph{Small degrees by ordering the weight gaps.}
The codimension of the complement of the GIT stable locus of
degree-\(d\) forms in \(d\) variables is exactly
\[
             \binom{2d-3}{d-2}-(d-1) \qquad (4\le d\le10).
\]
To obtain the lower bound, let \(x_i=w_i-w_{i+1}\), \(1\le i<d\),
for ordered weights with zero sum. Let \(S=\{i:x_i>0\}\).
The distinct positive values of the gaps, in decreasing order,
determine a strict chain of nonempty sets
\(A_1\subsetneq\cdots\subsetneq A_\ell=S\) so that
\[
 (x_1,\ldots,x_{d-1})=\sum_{j=1}^\ell a_j{\bf1}_{A_j},
                      \qquad a_j>0 .
\]
For a degree-\(d\) exponent \(\alpha\), put
\(u_i(\alpha)=\sum_{h=1}^i(\alpha_h-1)\) and
\(p_A(\alpha)=\sum_{i\in A}u_i(\alpha)\).
Count \(\alpha\) such that \(p_{A_j}(\alpha)\ge0\) for all \(j\),
with strict inequality for at least one \(j\). Every such coefficient
has strictly positive weight by summation by parts. Subtract the
dimension \(f(S)=\sum_{j<k}m_jm_k\) of the flag for the weight
blocks, where the block lengths \(m_j\) are successive differences
in \(\{0\}\cup S\cup\{d\}\). The result is a lower bound for the
codimension of that stratum of nonstable forms.

The integer minima over all one-break and multiple-break chains
are as follows. The second column is the exact global minimum.
\[
\begin{array}{r|rr}
 d& |S|=1 & |S|\ge2 \\ \hline
 4&7&8\\
 5&31&37\\
 6&121&145\\
 7&456&535\\
 8&1709&1942\\
 9&6427&7053\\
10&24301&25780
\end{array}
\]
For reproducibility, initialize \(D=\{\alpha:|\alpha|=d\}\),
\(E=\varnothing\) and \(A=\varnothing\). For every nonempty subset
\(I\subseteq\{1,\ldots,d-1\}\setminus A\), extend to \(B=A\cup I\)
and make the nested call
\[
 (A,D,E)\longmapsto
 \bigl(B,\ D\cap\{p_B\ge0\},\ E\cup\{p_B>0\}\bigr).
\]
At each nonempty node record \(|D\cap E|-f(A)\). This enumeration
visits each chain exactly once. The script
\path{nov28-lit-git-gap-orders.py} implements this recurrence using
integer sets represented by bit vectors.

For the opposite inequality, take only \(x_1>0\). The
\(\binom{2d-3}{d-2}\) positive coefficients are precisely those with
\(\alpha_1\ge2\). Their vanishing is equivalent to multiplicity
at least \(d-1\) at \([1:0:\cdots:0]\).
The corresponding linear system is base-point free away from
that point, since it contains the pure \(d\)-th powers of all
coordinates other than the first. By Bertini its general member is
smooth away from that point and has multiplicity \(d-1\) there.
The marked incidence therefore maps generically one-to-one onto its
image. Moving the point contributes \(d-1\) parameters. This proves
equality with the first column.

\paragraph{Numerical form of the codimension estimate.}
Combining the log-canonicity estimates with the Hilbert--Mumford
calculation gives the following lower bounds on the codimension
of the complement of the uniformly K-stable locus in the full
parameter space, valid in every Fano dimension \(n\ge d-1\).
\[
\begin{array}{r|rrr}
 d& b_1(d-1,d)\text{ bound}&\text{GIT codimension}&
                 \text{K-stability bound} \\ \hline
 4&9&7&7\\
 5&31&31&31\\
 6&94&121&94\\
 7&305&456&305\\
 8&1029&1709&1029\\
 9&3535&6427&3535\\
10&12326&24301&12326
\end{array}
\]
The file \path{nov28-lit-global-codim.py} implements the integer
recursion for the first numerical column. The GIT column is checked
independently of the singularity estimates by the weight-gap
enumeration. These codimension estimates apply to families
of arbitrary smooth hypersurfaces through intersection with
the open smooth parameter locus. More precisely, any general
projective linear family of dimension strictly below the stated
codimension consists entirely of uniformly K-stable fibers,
including its singular fibers.

For \(11\le d\le15\) the successive combined lower bounds,
using the two-block formula for the GIT estimate, are
\(43477,\ 154859,\ 556175,\ 2011877,\ 7323428\).

\subsection{K-stability under adding a pure power}
\label{sec:nov28-suspension}
All varieties in this subsection are over \(\mathbf C\). Let
\(D=(f=0)\subset \mathbf P^{q-1}\), with \(f\) homogeneous of
degree \(d\ge3\), and form the hypersurface
\[
 X=(y^d+f(x)=0)\subset \mathbf P^{q}.
\]
If \(q=d\), assume \((\mathbf P^{q-1},D)\) is log canonical and
\([f]\) is GIT stable for the standard \(\mathrm{SL}_q\)-action.
If \(q>d\), assume \(D\) is a K-stable \(\mathbf Q\)-Fano variety.
Then \(X\) is a uniformly K-stable \(\mathbf Q\)-Fano variety.

For \(q=d\), take a sufficiently small rational \(c_0>0\).
Ascher--DeVleming--Liu, Theorem~1.4
\cite{ar1909.04576}, proves the K-stability of
\((\mathbf P^{q-1},c_0D)\). Its stability threshold is strictly
greater than one by Liu--Xu--Zhuang, Theorem~1.6
\cite{ar2102.09405}. At \(c=1-1/d>c_0\) put
\(\lambda=(1-c)/(1-c_0)\). For every divisorial valuation,
\[
 A_{cD}=\lambda A_{c_0D}+(1-\lambda)A_D\ge\lambda A_{c_0D},
 \qquad S_{-(K+cD)}=\lambda S_{-(K+c_0D)}.
\]
The first inequality uses log canonicity and the equality of the
expected orders uses \(D\sim-K_{\mathbf P^{q-1}}\). Thus the
pair at coefficient \(1-1/d\) is uniformly K-stable. Projection
from the \(y\)-coordinate point makes \(X\) its crepant cyclic
Galois cover. Normality follows from the cyclic-cover equation
and generic reducedness of \(D\): the cover is regular in
codimension one and Cohen--Macaulay. Log canonicity in the
hypothesis implies that \(D\) is reduced, and the pair with
coefficient \(1-1/d\) is klt. The cover is klt by the ramification
formula. Its anticanonical bundle is \(\mathcal O_X(1)\) by
adjunction. Liu--Zhu, Theorem~1.2(2) \cite{ar2001.10557}, gives
K-polystability of \(X\). For \(q>d\) the entire K-polystability
assertion follows from Liu--Zhu, Theorem~5.1, applied with
ambient dimension \(q-1\) and covering degree \(d\).

We verify finiteness of automorphisms in both cases, including
singular covers. In the first case the projective stabilizer of
\([f]\) is finite by GIT stability. In the second it is finite by
K-stability of \(D\). Its action on \(D\) is faithful, since an
integral hypersurface of degree at least two spans its ambient
projective space and cannot be contained in the union of the
proper eigenspaces of a nonscalar automorphism. In particular
there is no nonzero constant-coefficient derivative of \(f\)
that vanishes identically, since this would make \(f\) independent
of a coordinate and give a positive-dimensional stabilizer.
An infinitesimal projective automorphism preserving \(X\) can be
represented by \(\dot x=Ax+by,\ \dot y=\ell(x)+ey\).
The mixed terms of its derivative of the equation are
\(y\,\partial_b f\) and \(dy^{d-1}\ell(x)\). Their distinct
bidegrees force \(b=0\) and \(\ell=0\). The remaining equation
requires \(\partial_{Ax}f\) to be scalar times \(f\), so \(A\) is
a scalar matrix. Comparing the pure \(y^d\) coefficient gives
the same scalar \(e\), proving that the projective Lie stabilizer
is zero. Since \(H^1(X,\mathcal O_X)=0\) for a klt Fano variety,
\(\operatorname{Aut}^0(X)\) preserves \(\mathcal O_X(1)\). This
line bundle embeds \(X\) by its complete linear system, by the
hypersurface sequence. Thus \(\operatorname{Aut}^0(X)\) is
trivial. The automorphism group is of finite type since it preserves
the ample anticanonical bundle. It follows that \(X\) is K-stable
and uniformly K-stable.

\subsection{Intersecting a component with the cyclic cover family}
\label{sec:new-affine-suspension}
Let \(B\) be an irreducible component of the closed complement of
\(\mathcal K_{n,d}\). It is invariant under
\(\operatorname{PGL}_{n+2}\), since this connected group preserves
the complement and has only finitely many components to permute
in that complement. In coordinates \(x_1,\ldots,x_{n+1},y\), the
affine chart in which the coefficient of \(y^d\) is one meets \(B\):
one can move a point outside any fixed zero set to the \(y\)-point.
Let \(c\) be the codimension of \(B\) and write \(B'\) for its
intersection with this chart. Every coefficient other than the
coefficient of \(y^d\) has strictly positive weight for
\((x,y)\mapsto(tx,y)\). Letting \(t\) tend to zero shows that
\([y^d]\in B'\).

Inside the chart take the affine linear space
\[
                    W=\{y^d+f(x)\}.
\]
It contains \([y^d]\), so the affine intersection \(B'\cap W\)
is nonempty. Each component through this point has dimension
at least \(\dim W-c\), by the dimension theorem in the smooth
ambient chart. If a closed cone of forms \(f\) of codimension \(k\)
contains all the failures of uniform K-stability in \(W\), it
follows that \(c\ge k\). The nonemptiness at the monomial is
essential to this dimension argument.

For \(n=d-1\), take as exceptional cone the union of the
non-log-canonical and non-GIT-stable forms of degree \(d\) on
\(\mathbf P^{d-1}\), including zero.
Sections~\ref{sec:nov28-lc-git} and~\ref{sec:nov28-suspension}
give uniform K-stability of every suspension outside this cone
and give \(c\ge\operatorname{Cat}_{d-1}\).
For every subsequent \(n\),
one can take the affine cone over
\(\mathcal P_{n-1,d}\setminus\mathcal K_{n-1,d}\), again including
zero. The Fano branch case of Section~\ref{sec:nov28-suspension}
then gives the same lower bound at the next dimension.
This proves the basic Catalan bound in Theorem~\ref{thm:new-catalan}.
For the stated table one also uses the ambient GIT estimate in
Proposition~\ref{prop:dec1-index-one-neighborhood} and the refined
branch estimates. The uniform factors are proved in
Theorem~\ref{thm:dec1-all-degree-catalan}.

More generally the lower bound for the complement in the full
projective coefficient space in any one Fano dimension propagates
to all higher dimensions, up to the dimension of the affine space
of branch forms if the original complement is empty. In the
applications here the complement in the full space contains the
multiple hyperplanes. Taking its intersection with the open
smooth locus can only increase its codimension.

These codimension bounds also hold in the coarse moduli space of
smooth hypersurfaces with their hyperplane polarizations. Indeed
the smooth equations are GIT stable and have finite projective
stabilizers. Their geometric quotient has fibers of dimension
\((n+2)^2-1\), and the image of the invariant relatively closed
exceptional locus is closed in this quotient. Subtracting the
common orbit dimension from both dimensions gives the assertion.

Singular equations provide an upper bound on the possible improvement
of a bound in the full parameter space. Cones over smooth
degree-\(d\) hypersurfaces in \(\mathbf P^n\) form a locus of
codimension \(\binom{n+d}{d-1}-(n+1)\). At a fixed vertex the
polynomial uses just \(n+1\) of the \(n+2\) variables, and the
vertex varies in dimension \(n+1\). It is recovered from the
equation since it is the only singular point. For \(n\ge d\)
the vertex divisor has \(\beta=-(d-1)/(n+1)\) by
Chapter~\ref{ch:families}. At \(n=d-1\) the same divisor
has log discrepancy zero. Thus these cones belong to the
full exceptional locus in every dimension of the theorem.

\subsection{A stronger quartic bound using K-semistable branches}
\label{sec:new-semistable-branch}
For every \(n\ge4\), the closed complement of the uniformly K-stable
quartic locus in \(\mathcal P_{n,4}\) has codimension at least
\(\mathbf{9}\). In particular this is a lower bound in the smooth
parameter space in each dimension beyond the quartic eightfold.

Here are two refinements of the cyclic-cover argument. First take the
base case \(q=d\) of Section~\ref{sec:nov28-suspension}. If the
anticanonical divisor \(D\) is log canonical, then
\((\mathbf P^{q-1},cD)\) is K-semistable for each \(0<c<1\).
Indeed \(A(v)\ge qS_H(v)\) by K-semistability of projective space,
and \(A(v)\ge v(D)\) by log canonicity. Their convex combination with
weights \(1-c,c\) gives the full valuative inequality.
The finite-cover theorem gives K-semistability of the cyclic cover.
Replacing the open uniformly K-stable locus by the open K-semistable
locus in Section~\ref{sec:new-affine-suspension} thus bounds the
codimension of the complement of K-semistability in the full space
\(\mathcal P_{d-1,d}\) below by \(b(d-1,d)\). Liu--Zhu's
Theorem~5.1 propagates this same semistability bound to every larger
Fano dimension. In particular the exact quartic log-canonicity seed
gives codimension at least nine for failures of K-semistability in
\(\mathcal P_{3,4}\).

Second, when \(q>d\), a K-semistable branch \(D\) with a GIT stable
equation already has a uniformly K-stable suspension. Write
\(N=q-1\) and \(\rho=q/d-1\). The cone degeneration in the proof of
\cite[Theorem~5.1]{ar2001.10557} gives K-semistability of
\((\mathbf P^N,c_1D)\) for \(c_1=1-\rho/N\).
For sufficiently small rational \(c_0>0\) the pair at \(c_0\) is
uniformly K-stable by \cite[Theorem~1.4]{ar1909.04576} and
\cite[Theorem~1.6]{ar2102.09405}. Now \(c_0<1-1/d<c_1\), since
\(c_1-(1-1/d)=(d-1)/(dN)>0\).
Write \(1-1/d=\theta c_0+(1-\theta)c_1\), with \(0<\theta<1\).
If \(\delta_0>1\) is the threshold at \(c_0\), the affine
dependence on \(c\) of
\(A_{cD}(v)-(q-cd)S_H(v)\) gives a lower bound
\(\theta(\delta_0-1)(q-c_0d)S_H(v)\) at \(1-1/d\).
This is uniformly positive relative to \(S_{-(K+(1-1/d)D)}(v)\).
The finite-cover theorem and the Lie-algebra calculation in
Section~\ref{sec:nov28-suspension} prove the claim.

It remains to bound the non-GIT-stable quartics in five variables.
The weight-gap count above applies with prefix sums
\(p_j=5\sum_{i=1}^j\alpha_i-4j\) for \(1\le j\le4\), with the
same partial-flag dimension. For transparency, the minimum lower
bound at each number of positive gaps is
\[
\begin{array}{c|rrrr}
\text{number of gaps}&1&2&3&4\\ \hline
\text{number of ordered partitions of gap supports}&4&18&52&75\\
\text{minimum of coefficient count minus flag dimension}&29&17&15&14 .
\end{array}
\]
The counts enumerate every weak ordering, including ties. Explicitly,
for each nonempty final support \(S\subseteq\{1,2,3,4\}\), form all
chains \(\varnothing\subsetneq A_1\subsetneq\cdots\subsetneq A_l=S\).
Among the 70 sorted quadruples \(1\le i_1\le\cdots\le i_4\le5\),
set \(\alpha_i=\#\{h:i_h=i\}\). Count a quadruple when
\(\sum_{j\in A_k}p_j\ge0\) for every \(k\) and at least one of these
sums is positive. Subtract \(\sum_{i<j}m_im_j\), using the block
lengths defined by \(S\), and minimize at each \(|S|\).
This gives the displayed integers, all at least nine.

Apply the affine intersection argument to \(\mathcal P_{4,4}\).
On its suspension family the exceptional equations are confined
to the union of the non-K-semistable branch cone (codimension at
least nine) and the non-GIT-stable branch cone (codimension at
least fourteen). Thus the full complement of uniform K-stability
has codimension at least nine there. Ordinary propagation of
uniform K-stability proves the assertion for \(n\ge4\).

\paragraph{A codimension bound starting with cubic fourfolds.}
\label{par:nov28-cubic-fourfold-codimension}
For every \(n\ge4\), the complement of the uniformly K-stable
locus in the projective space of cubic equations has codimension
at least seven. Here the open locus includes all the singular
\(\mathbf Q\)-Fano cubics that are uniformly K-stable.
For \(n=4\), Liu's Theorem~1.1 \cite{ar2007.14320} identifies
K-stability with GIT stability for cubic fourfolds. By Laza,
Table~1, Proposition~2.6 and Lemma~2.2
\cite{ar0704.3256}, the nonstable cubics are contained in the
\(\mathrm{SL}_6\)-saturations of eight weight subspaces. Fix the
coordinates of \(\mathbf P^5\). In the notation of that table,
a one-parameter subgroup with weights \(w_i\) imposes the
vanishing of all coefficients \(x_i x_j x_k\) of positive weight.
The remaining space is preserved by the corresponding parabolic
subgroup. Its saturation consequently has codimension at least
\(a-b\), where
\[
 a=\#\{0\le i\le j\le k\le5:w_i+w_j+w_k>0\},\qquad
 b=\#\{0\le i<j\le5:w_i\ne w_j\}.
\]
Indeed \(b\) is the dimension of the partial flag variety for
the multiplicities of the distinct weights, and each flag costs
\(a\) linear equations. The counts for Laza's weights are
\[
\begin{array}{c|r|c|c|c}
 & (w_0,w_1,w_2,w_3,w_4,w_5)&a&b&a-b\\\hline
 S_1&(2,2,-1,-1,-1,-1)&16&8&8\\
 S_2&(1,1,1,1,-2,-2)&20&8&12\\
 S_3&(2,1,0,0,-1,-2)&23&14&9\\
 S_4&(4,1,1,-2,-2,-2)&19&11&8\\
 S_5&(2,2,2,-1,-1,-4)&22&11&11\\
 S_6&(2,0,0,0,-1,-1)&18&11&7\\
 S_7&(1,0,0,0,0,-1)&16&9&7\\
 S_8&(1,1,0,0,0,-2)&25&11&14
\end{array}
\]
The bound seven is attained in the full coefficient space in
dimension four. For a fixed point, the six conditions of being
singular and the six additional conditions that the five-by-five
quadratic term have rank at most two are independent in the
indicated incidence count. Removing the five parameters of the
point gives seven when its general fiber is finite. To verify this,
fix the rank-two quadratic term. In coordinates with the point
\([1:0:\cdots:0]\), the remaining homogeneous cubic
terms form the complete space of cubic forms in \(x_1,\ldots,x_5\).
This variable space is basepoint-free away from the fixed point,
so a general full equation with this quadratic term is smooth away
from it by Bertini. Laza's case \(S_6\) shows these cubics fail
GIT stability.
The preceding affine intersection and cyclic-cover argument
propagates the minimum of seven to every larger dimension.

\section{Counting the normal-bundle splitting types of lines and conics}
\label{sec:new-splitting-codimension}
For \(d\ge10\) and \(n\ge4d\) there is a range with a codimension
bound that grows even when the permitted dimension is quadratic in
the degree. Let \(\mathcal U_{n,d}\) denote the open locus of
smooth equations in the projective space of all degree-\(d\) forms.

\begin{theorem}\label{thm:new-splitting-range}
Fix an integer \(2\le a<d\). In the range
\[
             d\ge10,\qquad 4d\le n\le d-4+\frac{d(d-3)}a ,
\]
every failure of uniform K-stability in \(\mathcal U_{n,d}\)
lies in a closed subset of codimension at least
\[
                      (a-2)(n-1)+1 .
\]
More precisely, such a hypersurface must contain a line \(C\) for
which \(N_{C/X}\) has a direct summand of degree at most \(-a\).
For example, for degree twenty and dimension eighty one may take
\(a=5\), giving codimension at least \(238\).
\end{theorem}

\begin{proof}
We first compute a splitting stratum on the projective line.
Write \(D>0\) and let \(v\) binary degree-\(D\) forms be
basepoint-free. They give a surjection
\(\mathcal O^{v}\twoheadrightarrow\mathcal O(D)\) with
kernel \(K=\bigoplus_{i=1}^{v-1}\mathcal O(-b_i)\), where
\(b_i\ge0\) and \(\sum b_i=D\).
In the open affine space of surjective tuples the codimension
of a nonempty stratum with this prescribed unordered list is
\[
    h^1(\operatorname{End}K)=\sum_{i,j}(b_i-b_j-1)_+ .
                                                        \tag{6}
\]
Indeed, injections \(K\hookrightarrow\mathcal O^v\) with
locally free quotient form an open subset of a vector space of
dimension \(v(v-1+D)\). The quotient has degree \(D\).
Choosing its isomorphism to \(\mathcal O(D)\) adds one parameter,
and changing the identification of the kernel gives precisely a
fiber of dimension \(h^0(\operatorname{End}K)\). Subtract from the
dimension \(v(D+1)\) of all tuples and apply Riemann--Roch on
\(\mathbf P^1\). If \(b_i=m=\max b_i\), the summands for this
one \(i\) give the lower bound
\[
 h^1(\operatorname{End}K)
 \ge (v-2)(m-1)-(D-m)=(v-1)(m-1)-(D-1).
                                                        \tag{7}
\]

For a fixed line on \(X\), take its \(v=n\) first normal polars of
degree \(D=d-1\). Their coefficients are arbitrary after imposing
the \(d+1\) independent equations for containment of the line.
The derivative normal sequence identifies \(N_{C/X}=K(1)\).
If every \(b_i\le a\), the ordinary line extraction of
Chapter~\ref{ch:models} gives normalized angular denominator at
most \(1+a/d\); its radial denominator
\((n-1/d)/(n-1)\) is also at most \(1+a/d\) in this range.
Thus \(\delta_C(X,H)\ge (n+1)d/(d+a)\).
The upper bound on \(n\) is precisely the comparison of this
threshold with \(r+2\). The common line-threshold criterion in
the proof of Theorem~\ref{thm:nov27-line-only-range} then gives
uniform K-stability if every line passes.

At a line that fails, take \(m\ge a+1\) in (7). Adding \(d+1\)
and subtracting \(2n=\dim\operatorname{Gr}(2,n+2)\) gives the
asserted codimension. There are finitely many possible splitting
types. Alternatively the bad splitting condition is the closed
condition \(H^1(C,K(a-1))\ne0\) in the proper incidence of lines
over the smooth parameter space. This proves closedness after
projection.
\end{proof}

For any of the improved line thresholds \(L\) in
Section~\ref{sec:nov27-optimized-lines}, the same count applies
with any integer \(a\ge2\) for which both ordinary denominators are
at most \((n+1)/L\). The count (6) also applies when conics enter
the point argument. In detail, let \(L\) be a sufficient common
threshold for all lines and integral conics on \(X\), let \(a,b\)
be nonnegative integers, and assume the following three bounds:
\[
  1+\frac ad\le\frac{n+1}{L},\qquad
  1+\frac b{2d}\le\frac{n+1}{L},\qquad
  \frac{n-1/d}{n-1}\le\frac{n+1}{L}.
\]
Then possible failures of uniform K-stability are contained in
a closed subset of \(\mathcal U_{n,d}\) of codimension at least
\[
       \min\big\{(a-2)(n-1)+1,\ (b-3)(n-2)-4\big\}.
                                                        \tag{8}
\]
Only positive values give an exceptional-set estimate. To see the
second term, fix a smooth conic and a vector space complementary
to its plane coordinates. Containment imposes \(2d+1\)
independent conditions. Independently the \(v=n-1\) plane-normal
polars restrict to arbitrary binary forms of degree \(D=2d-2\)
by surjectivity of restriction from the plane. For a smooth
hypersurface, a general choice of complementary directions makes
this tuple basepoint-free, as shown in
Section~\ref{sec:nov27-all-conic-normal-rank}. That calculation
gives angular denominator at most \(1+\max b_i/(2d)\), and the
radial denominator \((n-1/d)/n<1\). For any conic where the local
threshold is below \(L\), every such general complement must
therefore have \(\max b_i\ge b+1\).
Parameterize the complementary spaces as well as the conics.
The offending pairs have a dense open family of bad complements
of the full complement dimension, so this auxiliary dimension
cancels in the incidence estimate. Applying (7) at \(m=b+1\),
adding \(2d+1\) and subtracting \(3n+2\) for the moving conic
gives the second term in (8); taking the Zariski closure preserves
the dimension bound.

\section{Divisors obtained by blowing up a complete intersection}
\label{sec:new-complete-intersection}
Let \(N=n+1,\ r=n+2-d\) and let \(H=\mathcal O_X(1)\). For an
individual divisorial valuation \(v\), write
\(\beta_X(v)=A_X(v)-r S_H(v)\). A smooth hypersurface with
\(\delta(X,-K_X)\le1\) must have a prime divisor with
\(\beta_X(v)\le0\), by the divisorial attainment theorem
\cite[Theorem~1.2]{ar2102.09405}.
The following computations exclude particular divisors in all degrees
and dimensions under consideration.

\begin{proposition}\label{prop:new-ci-beta}
Let \(Z\subset X\) be smooth, integral and positive dimensional.
If \(Z\) is a scheme-theoretic complete intersection in \(X\) of
hypersurface sections of degrees \(a_1,\ldots,a_k\), then the
ordinary order \(E=\operatorname{ord}_Z\) satisfies
\[
 A_X(E)=k,\qquad N S_H(E)=\sum_{i=1}^k\frac1{a_i}.
\]
If \(Z\) is a scheme-theoretic complete intersection in the ambient
\(\mathbf P^{n+1}\) of degrees \(a_1,\ldots,a_c\), then
\[
 A_X(E)=c-1,\qquad N S_H(E)=\sum_{i=1}^c\frac1{a_i}-\frac1d .
\]
Here \(Z\subsetneq X\), and order along a prime divisor means its
ordinary divisorial order. In both cases \(\beta_X(E)>0\).
\end{proposition}

\begin{proof}
Write \(P=\mathbf C[x_0,\ldots,x_{n+1}]\), \(F\) for the defining
equation of \(X\) and \(R=P/(F)\). In the first case take the
homogeneous regular sequence in \(R\) defining \(Z\), with ideal
\(I\). Its associated graded ring is
\((R/I)[u_1,\ldots,u_k]\) with bidegrees
\(\deg u_i=(a_i,1)\), ordinary section degree followed by normal order.
Indeed the degree-\(j\) conormal piece of a regular sequence is the
\(j\)-th symmetric power of \(I/I^2\), with the indicated free
generators. The homogeneous quotient \(R/I\) is the section
coordinate ring of the integral complete intersection. The regular
sequence and positive dimension imply saturation, so this quotient
is a domain. Thus each nonzero initial section remains nonzero upon
localization at the generic point of \(Z\). The filtration computes
the ordinary order there, hence the order of the exceptional prime.

In the ambient case let \(g_1,\ldots,g_c\) cut out \(Z\), put
\(I=(g_1,\ldots,g_c)\subset P\) and write \(F=\sum_i h_i g_i\)
homogeneously. The coefficients \(b_i=h_i|_Z\) have no common zero
on \(Z\), since at such a zero \(dF=\sum h_i\,dg_i\) would vanish.
In \(\mathrm{gr}_I P=(P/I)[u_1,\ldots,u_c]\) the initial equation
is consequently \(L=\sum_i b_i u_i\), of bidegree \((d,1)\).
As this associated graded ring is a domain, the initial ideal of
the principal ideal \((F)\) is generated by \(L\). Its quotient is
also a domain. Away from the vertex in \(\operatorname{Spec}(P/I)\),
the equation exhibits a vector bundle of rank \(c-1\) over an
integral base. Above the vertex the fiber has dimension at most
\(c\), smaller than \(\dim Z+c\), the dimension of this quotient.
The quotient is a complete intersection, so it is unmixed and has
no embedded components. Its dense reduced open then proves the
claim. The same localization argument now identifies its graded
pieces with the true order filtration on the sections of \(X\).

In both cases \(R_m=H^0(X,mH)\) for \(m\ge0\). To compute the mean
order, logarithmically differentiate the bigraded Hilbert series in
its order variable at value one. Relative to the ordinary Hilbert
series of \(R\), the resulting factor is
\(\sum_i t^{a_i}/(1-t^{a_i})\), with the additional term
\(-t^d/(1-t^d)\) in the ambient case. At \(t=1\) the pole
coefficients are exactly the two asserted values of \(N S_H(E)\):
the ordinary Hilbert series has leading term \(d/(1-t)^{n+1}\),
and taking the order-weighted coefficient and dividing by
\(m\dim R_m\) contributes the factor \(1/(n+1)\).
Other roots of unity contribute only bounded terms to the
coefficients of this differentiated series, since the ordinary
Hilbert series is regular there. Blowing up a smooth center
gives discrepancy equal to its codimension in \(X\).

In the first case \(\sum 1/a_i\le k\) and \(r<N\).
In the ambient case, if a \(b_j\) is a nonzero constant, then
\(a_j=d\), and the remaining terms sum to at most \(c-1\).
Otherwise the nonzero \(b_i\) are all of positive degree. They
have no common zero on \(Z\), so their number is at least
\(\dim Z+1=n+2-c\) by projective dimension. Hence \(2c\ge n+2\).
Using \(a_i\ge1\) we obtain
\[
 N\beta_X(E)\ \ge\ N(c-1)-r(c-1/d)
  =(d-1)\left(c-\frac{n+2}{d}\right)>0
\]
because \(d\ge3\). This proves the proposition.
\end{proof}

The computation on \(X\) also excludes every weighted monomial
order in the equations of such a complete intersection. In the
first case of Proposition~\ref{prop:new-ci-beta}, give its defining
sections \(g_i\) arbitrary positive relatively prime integer
weights \(w_i\). At the generic point of \(Z\) use the monomial
valuation in these regular parameters, trivial on the residue
field. Then
\[
 (n+1)\beta_X(v)
       =\sum_i \left(n+1-\frac r{a_i}\right)w_i
       \ \ge\ (d-1)\sum_i w_i>0 .
\]
Indeed choose homogeneous vector-space lifts in \(R\) of a basis
of \(R/I\). The monomials in the \(g_i\) times these lifts form
a basis of \(R\): successive reductions modulo \(I\) span in
every section degree and the regular-sequence Hilbert series
gives equality of dimensions. At the generic point of \(Z\)
the weight-initial coefficients in this expansion cannot cancel
across distinct monomials, by the definition of the monomial
valuation in the regular local ring. For each individual monomial
the coefficient has nonzero residue unless its linear combination
of basis lifts was zero. Thus the weighted Hilbert series
calculation above gives \(N S_H(v)=\sum_i w_i/a_i\).
The toric discrepancy in a transverse smooth chart is
\(A_X(v)=\sum_i w_i\). This proves the formula, including
successive monomial blowups with arbitrarily large weights.

For instance a positive-dimensional linear space of ambient
codimension \(c\) on a smooth hypersurface gives equality in the
last displayed lower bound. The ordinary blowup of a closed point
also has positive \(\beta_X\) without any regularity assumption on
its tangent forms. Indeed the point-volume estimate in
Chapter~\ref{ch:foundations} gives \(N S_H(\operatorname{ord}_p)
\le n+d-1\), so
\[
  N\beta_X(\operatorname{ord}_p)
    \ge n(n+1)-(n+2-d)(n+d-1)=(d-1)(d-2)>0 .
\]
On a smooth Fermat hypersurface the point in the preceding cone
calculation attains equality. In particular even the sharp ordinary
point estimate leaves open the possibility of a different divisor
above the same point.

\section{Testing a candidate by weights on the ambient coordinates}
\label{sec:new-ambient-weight}
Let \(X=(F=0)\subset\mathbf P^{n+1}\) be smooth of degree
\(3\le d\le n+1\). Fix any projective coordinates and integers
\(w_0,\ldots,w_{n+1}\ge0\), of minimum zero and greatest common
divisor one. Assume some weight is positive. Let \(m\) be the smallest
weight of a monomial of \(F\), and \(F_*\) the sum of its terms of
weight \(m\). When \(F_*\) is irreducible, the filtration induced
on the section ring by these weights defines a divisorial valuation
\(v\) on \(X\). For this valuation the formula is
\[
 (n+1)\beta_X(v)
    =\frac{d-1}{d}\left(d\sum_i w_i-(n+2)m\right)
    \ \ge\ \frac{d-1}{d}\gcd(d,n+2)>0 .
\]
Thus an instability witness for a smooth equation cannot be obtained
from a choice of projective coordinates and weights with integral
initial hypersurface.

\begin{proof}
The initial ideal of the homogeneous principal ideal \((F)\) is
\((F_*)\), so the associated graded ring is a domain. For a section
the filtration order is the largest order of a polynomial representative
modulo \(F\). The domain property makes this order additive on
products. On a chart of a weight-zero coordinate it agrees with the
divisorial order given by restricting the ambient monomial extraction
to the strict transform of \(X\). To see this at the generic point,
make the toric modification inserting the ray of the weight vector
in the affine coordinate cone. On the dense torus of its exceptional
divisor the weighted leading equation is \(F_*=0\) up to a unit.
The Laurent hypersurface there is integral and generically smooth
by irreducibility and characteristic zero. It is nonempty since the
homogeneous irreducible polynomial of degree \(d>1\) is not a
coordinate monomial. The affine modification is smooth at points of
that torus, and the strict transform meets the exceptional divisor
transversely at the generic point of this hypersurface. This also
covers a coordinate ray using the coordinate divisor on the original
chart. A section with nonzero initial class remains nonzero at the
generic point after removing its exceptional factor. Adjunction
therefore gives \(A_X(v)=\sum_i w_i-m\).

The bigraded Hilbert series of the section ring is
\[
       \frac{1-t^d z^m}{\prod_i(1-tz^{w_i})}.
\]
Taking its logarithmic derivative in \(z\) at one as in the preceding
section gives \((n+1)S_H(v)=\sum_i w_i-m/d\).
The formula for \(\beta_X\) follows by using \(-K_X=(n+2-d)H\).

Here is a direct verification of its strict sign. Put \(q=n+2\).
The resultant of the \(q\) first partial derivatives of \(F\) is
nonzero by smoothness. For a parameter \(s\ne0\), consider
\(F_s=s^{-m}F(s^{w_0}x_0,\ldots,s^{w_{q-1}}x_{q-1})\).
Its coefficients are regular at zero. The resultant is homogeneous
of degree \((d-1)^{q-1}\) in each partial and scales by
\(\det(A)^{(d-1)^q}\) under a change of variables \(A\) on all
the partials. It follows that its value for \(F_s\) is its value
for \(F\) times
\[
       s^{(d-1)^{q-1}(d\sum_iw_i-qm)} .
\]
The exponent is nonnegative. Equality would make \(F_0=F_*\) smooth.
Its partials would then be a homogeneous regular sequence.
All their first syzygies are generated by the Koszul relations, whose
coefficients have degree \(d-1\ge2\). Ordinary and weighted Euler
identities, however, give the nonzero linear syzygy
\(\sum_i(w_i-m/d)x_i\partial_iF_*=0\), since the weights are
not all equal. Consequently the exponent is positive, and divisibility
of its parenthesized integer by \(\gcd(d,q)\) gives the claimed gap.
\end{proof}

\paragraph{A singular family at equality in every Fano dimension.}
In the full parameter space put
\[
 U(n,d)=\sum_{b=0}^{\lfloor(d-1)/2\rfloor}
       \binom{n+d-2b-1}{d-2b-1}-(2n+1).
\]
There is a locus of this codimension whose general member is klt,
has a unique singular point, and admits a divisor with
\(A_X(E)=S_{-K_X}(E)=r\). Thus \(U(n,d)\) is an upper bound
on the codimension of the complement in the full space. The family
is disjoint from the smooth parameter space.

Fix a point \(p=[1:0:\cdots:0]\) and a hyperplane \(y=0\)
through it. On the affine chart at \(p\), give \(y\) weight two
and the other \(n\) coordinates \(z_i\) weight one. Require every
term of weight below \(d\) to vanish. This is an intrinsic
condition on the flag: the weight ideals are generated by products
of the maximal ideal of \(p\), counted with weight one, and the
hyperplane equation, counted with weight two. All monomials
\(y^b z^\gamma\) of weight \(2b+|\gamma|<d\) have ordinary
degree at most \(d\), so their coefficients are independent.
They give exactly the displayed sum before subtracting the dimension
\(2n+1\) of the flag variety.

This linear system is basepoint free away from \(p\), since it
contains the \(d\)-th powers of every coordinate vanishing at
\(p\). A general member is smooth off \(p\) by Bertini. If
\(d=2s\), its first ordinary tangent form is a nonzero multiple of
\(y^s\). If \(d=2s+1\), this form is \(y^s L(y,z)\) with a
general linear form \(L\). In either case its factorization gives
at most finitely many choices for the marked hyperplane, including
the case of a cubic and two distinct tangent factors. The incidence
therefore has generically finite projection, giving its asserted
codimension.

The weight-\(d\) initial polynomial has smooth zero locus away from
the origin on the affine cone for a general choice of coefficients.
For even \(d=2s\) this holds for \(y^s+\sum z_i^d\), and for odd
\(d=2s+1\) it holds for \(y^s z_1+\sum z_i^d\), as is seen by
differentiation. The condition is open by properness on the weighted
projective space. The initial polynomial is irreducible, since two
positive-degree factors in at least three variables would meet away
from the origin. In the weighted blowup, the strict transform is
smooth and transverse to the exceptional divisor on its finite
uniformizing charts. The resulting pair with its reduced exceptional
divisor is plt near that divisor, by descent through the quotient
charts (or the toric adjunction formula). Its discrepancy over \(X\)
is \(A_X(E)=(n+2)-d=r>0\). Hence \(X\) is klt.
The Hilbert-series computation above, which requires only an integral
initial form and transversality at its generic point, now gives
\(S_H(E)=((n+2)-d/d)/(n+1)=1\). Thus \(\beta_X(E)=0\) and
the valuative criterion excludes uniform K-stability.
For cubic fourfolds the count gives \(U(4,3)=7\), agreeing with the
sharp full-space computation above. For quartic eightfolds it gives
\(U(8,4)=157\), while Theorem~\ref{thm:nov28-quartic-eightfold}
places every smooth member in the uniformly K-stable locus.

\subsection{A joint point and line bound for quartic ninefolds}
\label{sec:nov28-quartic-ninefold-joint}
In the parameter space of smooth quartic ninefolds, all failures
of uniform K-stability lie in a closed set of codimension at
least \(31\). By Section~\ref{sec:nov27-third-descent}, a
minimizer of the global threshold at most seven is centered at
a point \(p\). Fix the endpoint coefficient strictly between
\(125/17\) and \(37/5\), sufficiently close to \(125/17\).
The preceding interpolation forces a line through \(p\) with
local threshold strictly below \(37/5\). Indeed, if all lines
through this point had threshold at least \(37/5\), the fixed
positive gap absorbs the finite-level errors and excludes
every line produced by the length-one argument.
The preceding polar-net estimate forces that line to have
polar rank two. Throughout, incidence dimensions refer to the
space of equations, and closures do not increase those dimensions.

\paragraph{A pencil of singular symmetric matrices.}
For symmetric matrices of size \(m\), pairs \((A,B)\) for which
\(A+tB\) is generically singular have codimension at least
\(m+1\). Choose a nonzero polynomial kernel vector
\(u(t)=u_0+\cdots+t^e u_e\) of minimal degree. Its coefficients
are independent. Indeed a dependence ending at \(u_j\) gives,
by reversing its coefficients, a scalar polynomial with nonzero
constant term whose product with \(u(t)\) has zero coefficient
at \(t^j\). The part below this hole, nonzero since \(u_0\ne0\),
is again a kernel vector: after multiplication by the linear
pencil its support is disjoint from that of the part above the
hole. This contradicts minimality. For fixed independent
coefficients the equations are
\[
 A u_0=0,\qquad Au_i+Bu_{i-1}=0\ (1\le i\le e),\qquad Bu_e=0 .
\]
Both forms vanish on their span: by symmetry the entries of
\(B(u_i,u_j)\) are constant along each antidiagonal, and its
initial edge vanishes by \(Au_0=0\) and its terminal edge by
\(Bu_e=0\). The recurrence then gives the same vanishing for
\(A\). Each vector of a complementary basis gives \(e+2\)
independent equations. This is a total of \(m(e+2)\) linear
conditions, while the kernel tuple varies in at most
\(m(e+1)-1\) dimensions modulo a common scalar.

We also need a variant for two distinct six-dimensional
hyperplanes of a seven-space. For a freely varying pair of
quadrics on the seven-space, requiring both restricted
pencils to be generically singular costs at least nine.
On their common five-space let \(r\) be the generic pencil rank.
If \(r=5\) the first hyperplane condition costs at least seven.
For two independent members nonsingular on the common space,
singularity of their extensions to the second hyperplane gives
two further independent equations, linear with nonzero leading
coefficient in their two freely varying complementary diagonal
entries. If \(r=4\), generic singularity on the common space
costs at least six. At two independent rank-four members, the
pairing of each complementary vector with the kernel of each
form must vanish for the respective six-dimensional extension to
be singular. These four linear conditions on the two pairs of
complementary columns are independent. In these arguments the
members can be chosen on finitely many open charts, or equivalently
one can compute the extension fibers over every fixed common pencil.
If \(r=3\), take two rank-three generators. Their individual rank
conditions cost six. Tangency to the rank-three symmetric locus
requires the second form to vanish on the two-dimensional radical
of the first, by differentiating the Schur complement. This costs
three additional conditions even within rank three: express the
second form as the pullback of a nonsingular ternary form by a
rank-three map. The images of two fixed independent vectors must
belong to one isotropic line. The isotropic line and their two
coordinates on it have dimension three, against six unrestricted
entries. The rank-three parametrization has equidimensional
orthogonal-group fibers. If \(r\le2\), the individual rank
conditions already cost twelve.

\paragraph{A point of Hessian rank at most two.}
For a fixed line, containment costs five conditions, and rank
two for its matrix of nine normal cubics costs fourteen.
Lines vary in dimension eighteen, leaving codimension at least
one. At fixed such polars their common kernel \(V\) is a
seven-space. Over the line the quadratic term restricted to
it is \(Q_t=s^2Q_0+stQ_1+t^2Q_2\) in homogeneous coordinates
with \(p=[1:0]\). For each \(p\), the other coefficients
\(Q_1,Q_2\) and the two cubic light coefficients vary freely
and independently after fixing the Hessian at \(p\). Indeed
they come from the separate normal-degree-two terms multiplied
by \(st,t^2\), and the normal-degree-three terms multiplied by
\(s,t\).

On any subspace \(U\subset V\) of dimension at least six where
\(Q_0\) vanishes, the following observation applies. If the
restricted pencil of \(Q_1,Q_2\) contains a nonsingular member
and a cubic is nonzero on their common projective zero set, take
a point \([u]\) there with nonzero cubic. The common orthogonal
to \(u\) for the two quadrics inside \(U\) has codimension at
most two. The nonsingular member restricts there with rank at
least \(\dim U-4\ge2\). Choose two such transverse quadratic
directions at generic base parameter. The two-heavy flag
with weights \(3,3,2\) on them and on \(u\) has precisely the
face of a nondegenerate binary quadratic plus a nonzero cube.
Its signed denominator \(10/7\), from the disjoint monomial
estimate for \(xy+u^3\), gives the line mark \(280/37>37/5\).
The common projective zero has dimension at least three, and
vanishing of a free cubic on it costs at least twenty: projection
from a general center dominates \(\mathbf P^3\), so the cubics
in four general linear coordinates are independent on it.
This holds for the reduced underlying set and requires no
complete-intersection assumption.

At \(p\) write the Hessian in the tangent direction of the line
and eight complementary tangent directions as
\(\left(\begin{smallmatrix}0&B\\ B^t&D\end{smallmatrix}\right)\),
where \(D\) is an arbitrary symmetric matrix even after fixing
the first polars. The row \(B\) vanishes precisely at the
ramification points of the basepoint-free normal polar map
\(\mathbf P^1\to\mathbf P^1\) of degree three. In fact its
value is the differential of the polar row modulo the row at
\(p\); the eight tangent directions constitute the kernel of
that row. If \(B\ne0\), the seven-space \(V\) equals its
kernel among those directions. Rank at most two requires
\(D|_V=0\), at cost twenty-eight. All-singularity on \(V\)
costs at least another eight by the pencil lemma, and cubic
failure costs at least twenty instead. Allowing \(p\) to move
on the line loses at most one dimension.

If \(B=0\) the point is in a finite set for the fixed polars.
The condition \(\operatorname{rank}D\le2\) costs twenty-one.
When \(D|_V\) has rank two its two distinct linear factors
give two six-dimensional hyperplanes \(U\). To fail both
available tests even when neither cubic system vanishes on its
common zero set, both restrictions of the quadratic pencil must
be generically singular. The preceding variant costs nine.
If \(D\) has rank two but \(D|_V\) has rank at most one, its
two linear factors have proportional restrictions to \(V\),
including a possible zero restriction. The span of those
factors must then contain the specified equation of \(V\).
In the Grassmannian of two-planes of eight linear forms this
costs six extra conditions. One six-dimensional \(U\) with
\(Q_0|_U=0\) still exists, and failure costs at least seven
more by the pencil lemma and cubic test. Finally rank at most
one for \(D\) itself costs twenty-eight and gives the same
single-hyperplane test. The smallest total codimension in this
paragraph is therefore \(1+21+9=31\).

\paragraph{Higher point Hessian rank.}
Consider the point conditions of
Section~\ref{sec:nov27-orthogonal-cubic}. Their estimates at
fixed marked smooth point for nonregularity or the singular
cubic on the smooth quadric locus are \(45,49\). For
\(R=3,4,5,6\) the radial-kernel failures there have marked
codimension at least \(46,31,65,31\), including the quadratic
rank codimension. In rank four the sole bound below forty is
when the cubic restricted to the five-dimensional radical has
singular locus of codimension at most one. Such a cubic is zero
or has a repeated linear factor. The forms \(l^2m\) vary in
at most nine affine dimensions inside dimension thirty-five.
Thus the corresponding marked bound increases from \(31\) to
\(15+26=41\). In rank six the only cost below forty occurs if
both \(q_3,q_4\) vanish identically on the three-dimensional
radical \(K\) of \(q_2\). This costs \(6+10+15=31\) at a
marked point. All other strata have global codimension at least
\(40-9=31\) directly.

For this rank-six exception fix \((p,T_pX,q_2)\). In local
homogeneous tangent coordinates, a direction \([v]\) of a
rank-two line through \(p\) satisfies
\[
 q_2(v)=q_3(v)=q_4(v)=0,\qquad
 \dim\langle dq_2(v),dq_3(v),dq_4(v)\rangle\le1 .
\]
Indeed the normal polar in the transverse direction to \(X\)
does not vanish at \(p\); the polars in all tangent directions
have constant term zero. Their linear, quadratic and cubic
coefficients along the line are exactly the three gradient rows,
up to nonzero numerical factors. If \(v\notin K\), vanishing and
proportionality to \(dq_2(v)\ne0\) impose eight independent
first-jet conditions each on \(q_3,q_4\): one may multiply
arbitrary forms of degree one less by a linear form vanishing
on \(K\) but nonzero at \(v\). The direction varies in dimension
seven, so this adds at least nine, giving \(31+9-9=31\) on
forgetting \(p\).

If \(v\in K\), its choices have dimension two. The remaining
gradients of \(q_3,q_4\) are independent arbitrary vectors in
\((T_pX/K)^*\) before imposing their proportionality, which
costs five. They cannot both vanish for a rank-two line. Let
\(l\) span them. The seven-space of lights is \(T\oplus W\),
where \(T=K/\langle v\rangle\) has dimension two and \(W\) is
a lift of the five-dimensional hyperplane \(\ker l\subset
T_pX/K\). On \(V\) the first quadratic coefficient is \(Q_0(w)\)
of rank four or five. Write the other two as
\[
           Q_i=P_i(w)+\sum_{j=1}^2 M_{ij}(w)t_j .
\]
The four linear entries of \(M\) vary independently: they are
the coefficients of the separate terms \(v\,t_jw\) in the cubic
and \(v^2t_jw\) in the quartic. They respect the vanishing on
\(K\) and do not affect the already fixed gradients at \(v\).
If \(\det M\) is nonzero at a general zero of \(Q_0\), solve
the two affine-linear equations \(Q_i=0\) uniquely for \(t_j\).
The tangent space to these three equations projects
isomorphically onto the tangent hyperplane to \(Q_0\) in \(W\)
at every such regular nonzero \(w\). The quadratic \(Q_0\)
restricts on that hyperplane with rank at least two, since its
rank on \(W\) is at least four. A cubic light coefficient
nonzero on this common locus now gives exactly the same
quadratic-plus-cube flag as above. Failure by cubic vanishing
costs at least twenty independently: the pure cubics in \(W\)
in the coefficient at \(p\) remain freely variable and their
restrictions to the irreducible quadric in \(\mathbf P(W)\)
have dimension \(35-5=30\).

Otherwise \(\det M\) is a scalar multiple of \(Q_0\). The
zero multiple costs at least nine among the twenty arbitrary
matrix entries: a linear space of singular \(2\)-by-\(2\)
matrices has a common kernel or common image, by the two
rulings of the rank-one quadric. For a nonzero multiple
\(Q_0\) must have rank four, and all four entries of \(M\)
must lie in its fixed four-dimensional space of essential
linear forms. This alone costs four extra conditions. The
codimension after allowing \(v,p\) to vary is consequently
at least \(31+5-2+4-9=29\) from this weaker count. In fact
the ordered four entries are a conformal change of variables
for the standard nonsingular four-dimensional determinant
quadratic. Such changes form a group of dimension seven,
against twenty free coefficients. Thus the nonzero multiple
costs thirteen, and the final bound is at least
\(31+5-2+9-9=34\). This proves the joint incidence statement.

\subsection{A stronger joint incidence estimate}
\label{sec:nov28-quartic-ninefold-joint-sharp}
The exceptional codimension in
Section~\ref{sec:nov28-quartic-ninefold-joint} is at least \(36\).
We retain its notation and give the additional incidence bounds.

\paragraph{The pencil on one six-space.}
Call a pair of quadrics on a six-space intrinsically bad if no
component of their common projective zero set has a nonempty
open subset where a general pencil member restricts to rank at
least two on the common polar hyperplanes at that point. Every
component has dimension at least three. If a pair is not
intrinsically bad, failure of the quadratic-plus-cube flag from
the previous section costs at least twenty on the free cubic.
We claim the codimension of intrinsically bad pairs is at least
\(14,12,17\) when the generic rank is \(5,4,\le3\)
respectively. A nonsingular member always suffices.

First let \(A\) have rank five and \(k\) span its radical.
If \(Bk\ne0\), the common base contains a point \(u\) with
\(B(u,k)\ne0\) unless the pencil contains a rank-one member.
Indeed \(\mathbf C[x]/(A)\) is a graded unique factorization
domain: write \(A=x_1x_2+q'\) with \(q'\) irreducible,
so \(x_1\) is prime and inverting it gives a Laurent
polynomial ring, and apply the elementary factoriality criterion
by clearing powers of this prime. If the zero support of the
nonzero homogeneous quadratic \(B\) is contained in the
hyperplane \(B(-,k)=0\), it is therefore a scalar times the
square of that linear form modulo \(A\). Here a nonzero linear
form is prime since in a positively graded factorial algebra
generated in degree one it is irreducible. At a point with
\(B(u,k)\ne0\) the common polar hyperplanes project
isomorphically onto the \(A\)-orthogonal hyperplane of \(u\)
in the five-dimensional nonsingular quotient. The latter
restriction has rank three at its nonzero isotropic vector.
A pencil containing a rank-one form costs at least
\(15-1=14\), by varying the member in \(\mathbf P^1\).

Otherwise the pencil has a common kernel line, costing at least
\(12-5=7\) before passing to the five-dimensional quotient.
There set \(C=A^{-1}B\). At a common zero where
\(A(Cu,Cu)\ne0\), the span of \(u,Cu\) has Gram rank one
and dimension two. Its orthogonal has dimension three and rank
two. The same rank bound therefore holds upstairs. If this
quadratic vanishes throughout the common base, either the
quadratic \(B\) modulo \(A\) is a square, giving the rank-one
case, or its principal ideal in the factorial quadric ring
is radical. Thus \(C^2\in\langle I,C\rangle\). Among all
fifteen self-adjoint operators on the nonsingular five-space
this costs at least seven. Two distinct eigenvalues and the
choice of nondegenerate eigenspaces give dimension at most
\(2+2\cdot3=8\). In the one-eigenvalue case the square-zero
part with isotropic image of rank \(e\le2\) has dimension
\(e(5-e)\), by choosing that subspace and its invertible symmetric
map from its dual, plus one for the eigenvalue. This proves the
rank-five assertion. A proportional pencil is included in the
same bound, or can be counted directly.

For generic rank four choose a rank-four generator with nonsingular
block \(A_0\) and radical of dimension two. This costs three;
write the second symmetric form in complementary blocks
\(P,H,H^t,J\). The Schur complement for a pencil of rank at
most four imposes \(J=0\) (three conditions) and, successively,
\[
 H^t A_0^{-1}H=0,\qquad U^tPU=0,\qquad
       U^tP A_0^{-1} P U=0,\quad U=A_0^{-1}H .
\]
If the two columns span an isotropic two-plane their constraint
costs three (the maximal isotropic planes have dimension one,
and their ordered bases have dimension four). The second displayed
condition costs another three independently on \(P\). If their
span has dimension one, the isotropic-line and two-scalar
description costs four from eight entries. For a generator
\(u\) of the span of \(U\), \(u^t P u=0\) then costs one and
\((Pu)^t A_0^{-1}(Pu)=0\) a further one: the vector
\(A_0^{-1}Pu\) varies freely in \(u^{\perp_{A_0}}\) after
the first condition and the quadratic form there is nonzero.
If the two columns vanish the cost is eight instead. This gives
at least twelve in every case. For generic rank three the first
rank condition costs six and \(J=0\) on the three-dimensional
radical costs six. The three columns of the \(3\)-by-\(3\) block
\(A_0^{-1}H\) must span an isotropic line or zero, costing at
least five (dimension at most \(1+3\)). At smaller generic rank
the two individual rank bounds already cost twenty. One may
choose a maximal rank generator on finitely many open charts:
a nonzero minor in a binary pencil has bounded degree.

If the two distinct six-dimensional hyperplanes from a rank-two
\(D|_V\) both have intrinsically bad pencils, the resulting
codimension is at least fourteen. This is already true if the
first restriction has rank five or at most three. At its
remaining rank four the first restriction costs at least twelve.
If the common five-space has generic rank at least four,
generic singularity on the other six-space costs two more,
using the two distinct pencil members in the previous section.
Otherwise the common five-space contains the two-dimensional
radical of a maximal-rank member of the first pencil: a hyperplane
missing that radical surjects onto the nondegenerate quotient.
For fixed first pencil this containment costs two on the
projective parameter space of hyperplanes. The codimension is
the same for any fixed common five-space by the transitive
linear action on that parameter space. Combining with the
twenty cubic conditions when one restriction is intrinsically
good, the rank-at-most-two ramified case now costs at least
\(1+21+14=36\). When \(D|_V\) has rank at most one, use the
single six-space bound of twelve with the six extra conditions
in the preceding section (or the separate stronger rank bound
on \(D\)). At an unramified point the seven-dimensional
singular-pencil bound already gives \(1+28+8-1=36\).

\paragraph{The rank-six point away from its radical.}
The old exceptional case requires \(q_3|_K=q_4|_K=0\) on the
three-dimensional radical. For \(v\notin K\) the previous
count up to choosing the line direction had global codimension
\(31\). In the seven lights, \(Q_0\) is the induced quadratic
on \(v^{\perp}/\langle v\rangle\), of rank four with radical
\(K\). Write the space \(K\oplus W\); now \(\dim W=4\).
The two quadratic coefficients \(Q_1,Q_2\) are unrestricted:
the terms linear in \(v\) in \(q_3\) and quadratic in \(v\)
in \(q_4\) do not constrain the values or the previously
prescribed gradients at \(v\) and vanish on \(K\).
If the quadratic family has generic rank at most five, its
values \(Q_2,Q_0+Q_1+Q_2\) each have rank at most five,
costing six on these two independent matrices.

Otherwise any common zero \(u\in\mathbf P(K)\) of
\(Q_1,Q_2\) gives common polar hyperplanes of codimension at
most two. Rank at least six gives rank at least two on them.
The cubic \(q_4(v,u^3)\), with the usual polarization
normalization, must therefore vanish on the common zero set.
This costs at least four in the coefficients restricted to
\(K\), including specializations of the quadratic pair. If
their pencil contains no smooth conic, this follows already
from the singular-pencil lemma with \(m=3\). Otherwise fix
a smooth member; the other cuts a binary quartic on it. Having
\(j\) distinct roots costs at least \(4-j\), and the freely
variable cubic restricted to the conic separates all \(j\le4\)
points. An identically zero restricted quartic costs five.

There is at least one further condition, independently of these
restricted data. Fix a nonsingular nonzero isotropic vector
\(w\in W\) of \(Q_0\). For generic remaining cross and pure
coefficients the two equations \(Q_i(w+k)=0\) in \(k\in K\)
have a solution where their gradients in \(K\) are independent.
Their affine linear terms and constants vary freely, even for
fixed quadratic leading terms on \(K\). The incidence of such
a solution is dominant: choose constants giving a solution at
zero and independent linear terms, and use the surjective
Jacobian in \(k\) for smoothness of projection there. At that
solution the common tangent maps surjectively to the tangent
hyperplane of \(Q_0\) at \(w\) in \(W\), on which the restriction
has rank two. The cubic coefficient at \(p\) is generically
nonzero there, since the pure cubics in \(W\) inside \(q_3\)
still vary freely, including its value at \(w\). Thus for
every fixed set of previously restricted coefficients, the
remaining bad fiber is not dense. The total bound is at least
\(31+4+1=36\).

\paragraph{The rank-six point in its radical.}
In the notation of the preceding section only \(\det M=0\)
still needs attention: the nonzero multiple of \(Q_0\) already
costs thirteen. The case of both common image and common
kernel has at most seven affine parameters (two projective
directions and one linear form), instead of the eleven allowed
by either single ruling, so it too suffices. In either remaining
ruling write \(l_1,l_2\) for the two independent entries after
a constant change in the row or column. Away from their common
zero, projection to \(W\) sends the common quadratic equations
smoothly onto \(Q_0(w)=F(w)=0\), with affine-space fibers:
for common image eliminate one \(t\)-linear equation, so
\(F=P_2\) after changing the two equations; for common kernel
eliminate \(t_1\) leaving \(t_2\) free, with
\(F=l_2P_1-l_1P_2\). At smooth points where the two equations
in \(w\) have independent differentials, the upstairs
quadratic rank on the common tangent is at least the rank
of \(Q_0\) on this projected tangent.

Choose a general plane \(\Pi\subset\mathbf P(W)=\mathbf P^4\)
cutting \(Q_0\) in a smooth conic and with that conic avoiding
\(l_1=l_2=0\). The latter is possible unless their common
projective plane is contained in \(Q_0\). A rank-five quadric
contains no such plane; for rank four these planes must contain
its vertex and correspond to the two rulings of its smooth
quadric surface quotient, costing five in \(\mathrm{Gr}(3,5)\).
We may thus omit this more special case. In rank four also
choose \(\Pi\) missing the vertex. The restriction of \(F\)
to the conic is an arbitrary binary quartic or sextic: for
the sextic the two coprime binary quadratics \(l_1,l_2\)
times the two full quartic spaces span every sextic. Unless
the restriction has at least two distinct roots, its cost is
already at least three (the affine cone of pure fourth or sixth
powers, including zero; an identically zero restriction can
also supply arbitrary points).

Choose two distinct roots. At each keep all values on \(\Pi\)
fixed and vary the two normal derivatives of \(F\) independently
by changing \(P_i\) by plane-normal linear forms times linear
forms on \(\Pi\). Values of the latter at the two points are
independent, and at least one \(l_i\) is nonzero in the sextic
case. For rank five the projective tangent quadratic at each
isotropic point induces a nonsingular ternary form on
\(w^{\perp}/\langle w\rangle\). Its hyperplane cut by \(dF\)
has rank two unless its dual normal is isotropic. The two free
normal derivatives give a two-dimensional subspace of this
dual ternary space, since the conic on \(\Pi\) is smooth;
such a subspace cannot be totally isotropic. For rank four
the analogous ternary form has radical the projection of the
vertex. A nonzero derivative in this radical ensures a
rank-two restriction, and can be arranged since the vertex
is off \(\Pi\). Thus at either root failure of the tangent
rank test costs at least one condition, independently at the
two roots. Alternatively if the rank passes, the cubic light
coefficient must vanish at a chosen lift. Its value can be
varied independently at both roots using the pure cubics on
\(W\) in \(q_3\). The lifts depend only on the fixed values
of the quadrics at the two roots. Consequently the required
additional cost is two in every case, raising this branch
from \(34\) to at least \(36\).

\paragraph{The rank-four point.}
The only old point bound below \(45\), apart from those just
treated, is the rank-four quadratic with radical \(K\) of
dimension five and nonreduced \(q_3|_K=l^2m\). Its cost after
forgetting the point was at least \(32\). The zero cubic costs
\(15+35-9=41\) and may be omitted. For \(v\notin K\) the
two first-jet tests off \(K\) cost sixteen and the direction
varies in dimension seven, as before. For \(v\in K\) the
equation \(q_3(v)=0\) leaves at most three dimensions for
the direction. If \(dq_3(v)\ne0\), the value and first jets
of the unrestricted quartic impose eight conditions. Otherwise
the four nonradical derivatives of the cubic at \(v\) must
vanish and the quartic must vanish there too, imposing five
independent conditions. A rank-two line then requires
\(dq_4(v)\ne0\). Its light kernel in \(K/\langle v\rangle\)
has dimension at least three. This subspace \(U\) is in the
radical of \(Q_0\), and \(Q_1|_U\) is a square up to scalar,
or zero: the zero differential of \(l^2m\) forces \(l(v)=0\)
unless \(m\) is zero, and its second term along \(v\) is
\(m(v)l^2\). Take a two-space in the kernel of this square
(any two-space if it is zero). The freely variable \(Q_2\)
cuts it in a binary quadratic. At any root the common tangent
has codimension at most two in the seven lights. If \(Q_2\)
has rank at least six on all seven, its restriction there has
rank at least two. Otherwise that rank condition itself costs
three. For a binary quadratic with distinct roots the free
coefficient \(q_4(v,u^3)\) must vanish at both; its coefficients
are independent of the terms defining \(Q_2\) and the earlier
value and gradients at \(v\). A degenerate binary quadratic
already costs one and still has at least one root, giving the
second condition. This branch therefore has codimension at least
\(32+5-3+2=36\). All other fixed-point failure bounds from the
original table are at least \(45\) before subtracting nine,
as seen directly in its individual entries. This finishes the proof.

\paragraph{The locus of possible failures in degree five and dimension eleven.}
In the parameter space of smooth quintic elevenfolds, every failure
of uniform K-stability is contained in a closed subset of codimension
at least \(635\). Indeed the line-only interpolation with threshold
\(144/17<60/7\), the higher-rank line calculation and the order-two,
order-four and terminal pencil calculations reduce a failure to a
rank-two line with initial light cubic which is geometrically
non-log-canonical on the projective light space.

For a fixed line in \(\mathbf P^{12}\), containment imposes six
independent conditions. Its eleven normal quartics span a pencil in
the five-dimensional binary quartic space. The rank-two determinantal
condition costs \(27\); fix these polars, and write \(V\) for their
nine-dimensional light kernel. Vanishing of the entire quadratic
light form on \(V\) costs \(4\binom{10}{2}=180\), its coefficients
being arbitrary binary cubics. The next light coefficient is an
arbitrary member of
\[
 H^0(\mathbf P^1,\mathcal O(2))\otimes\operatorname{Sym}^3 V^*
\]
of dimension \(495\), independently of all those conditions.

By the classification of non-log-canonical cubics above, the
space of essential linear variables of the geometric generic
cubic has dimension \(k\le3\). It descends to the function field
of the base as the annihilator of all constant directions on
which the directional derivative vanishes. Its saturated extension
is a subbundle \(E\subset V^*\otimes\mathcal O_{\mathbf P^1}\),
with splitting \(\bigoplus_{i=1}^k\mathcal O(-a_i)\),
\(0\le a_1\le\cdots\le a_k\). The cubic belongs to
\(H^0(\operatorname{Sym}^3 E(2))\): membership is generic and the
inclusion of subbundles locally splits. Every generator must be
used by some monomial, so \(a_k\le2\). The family of such
subbundles has dimension at most
\[
 B=9\sum_i(a_i+1)-\sum_{i,j}\max\{0,1+a_i-a_j\}.
\]
This subtracts the dimension of the freely acting group of
automorphisms of the source from the dimension of the space of
bundle injections. For fixed \(E\) the section dimension is
\(
 D=\sum_{i\le j\le l}\max\{0,3-a_i-a_j-a_l\}.
\)
The splittings allowed by essentiality and the corresponding
upper bounds for \(B+D\) are
\[
\begin{array}{c|rrrrrrrrr}
 (a_i)&(0)&(0,0)&(0,1)&(0,2)&(0,0,0)&(0,0,1)&(0,0,2)&
 (0,1,1)&(0,1,2)\\ \hline
 B+D&11&26&29&35&48&47&49&46&51.
\end{array}
\]
If every \(a_i\ge1\), there are no sections. The only further
splitting with \(a_1=0\) is \((0,2,2)\), for which the generic
form is \(x^2\) times a linear form, of essential rank at most
two. Thus the cubic imposes at least \(495-51=444\) additional
conditions. The lines vary in dimension \(22\), so the required
bound is \(6+27+180+444-22=635\). The finitely many incidence
families are constructible, and taking closures preserves their
dimension.

The separate reduced-cubic line calculation improves this codimension
to at least \(1067\). A remaining cubic has a repeated linear
factor and can be written \(L^2M\) over the geometric generic base.
The repeated factor is unique, hence both factors descend and have
primitive binary polynomial representatives. If \(L\) is constant,
the possible cubics have dimension at most \(8+27=35\) (the
quadratic coefficients of \(M\) include the scalar content). If
\(L\) moves, it has base degree one, \(M\) is constant and the
family has dimension at most \(18+8=26\). The pure-cube case has
constant root and is already included in the first count. Thus
the cubic imposes \(495-35=460\) conditions.

For each such cubic choose a constant linear factor and cut by its
hyperplane. The eight remaining lights have an independent
quartic of binary base degree one. At geometric generic base this
quartic must be non-log-canonical or zero, as the four-light
homogeneous log-canonical criterion after the constant cut would
give the sufficient mark \(60/7\). Evaluate at two fixed distinct
base points. Their forms are independent and lie in the closed
cone of non-log-canonical quartics, including the zero form.
Each costs at least \(208\) by the non-log-canonicity recurrence
in projective dimension seven. The choice of factor is finite
once the cubic is fixed. The total codimension is at least
\(191+460+416=1067\).

\section{Singular unstable hypersurfaces}

\paragraph{Singular klt hypersurfaces in every Fano degree.}
For every $3\le d\le n+1$ the hypersurface
\[
 X=\bigl(x_0x_1^{d-1}+x_2^d+\cdots+x_{n+1}^d=0\bigr)\subset\mathbf P^{n+1}
\]
is a K-unstable klt Fano variety with exactly one singular point, including in index one. This gives the scope of the smoothness assumption in the all-member question. The partial derivatives isolate $p=[1:0:\cdots:0]$. In its affine chart the germ is the cyclic cover of $\mathbf A^n$ given by adjoining the $(d-1)$st root of $-(x_2^d+\cdots+x_{n+1}^d)$. The reduced branch divisor is a cone over a smooth projective hypersurface. Ordinary blowup shows its log canonical threshold at the origin is $\min(1,n/d)$. The branch coefficient $(d-2)/(d-1)$ is strictly below this threshold because $n\ge d-1$. The finite-cover discrepancy formula gives klt singularities on $X$.

Give the affine variables $x_1,x_2,\ldots,x_{n+1}$ weights $d,d-1,\ldots,d-1$. The equation is homogeneous of weight $d(d-1)$. The exceptional hypersurface of the weighted blowup at $p$ is irreducible: the affine weighted equation is irreducible by Eisenstein at any prime factor of the squarefree branch form in $\mathbf C[x_2,\ldots,x_{n+1}]$. Its generic point lies in the torus of the weighted exceptional space, where the primitive weights have trivial generic stabilizer. Consequently its valuation $v$ has
\[
 A_X(v)=d+n(d-1)-d(d-1)
\]
by weighted blowup and adjunction. Its section weights are the homogeneous weights with $\operatorname{wt}(x_0)=0$: after inverting $x_0$, every nonzero weight-homogeneous element of the projective section ring restricts in the affine chart to a nonzero element of the same weight. Thus its bigraded Hilbert series, with variables recording section degree and valuation weight, is
\[
 \frac{1-z^d u^{d(d-1)}}{(1-z)(1-zu^d)(1-zu^{d-1})^n}.
\]
The logarithmic $u$-derivative at $u=1$ has leading term $[d+n(d-1)-(d-1)]/(1-z)$. Dividing its coefficient means by the section degree gives
\[
 S_H(v)=\frac{n(d-1)+1}{n+1},\qquad
 A_X(v)-(n+2-d)S_H(v)=-\frac{(d-1)(d-2)}{n+1}<0.
\]
The valuative criterion proves the assertion.

\subsection{Unstable hypersurfaces with one singular point}\label{sec:dec1-unstable-flags}
Let \(c_{n,d}\) be the codimension of the complement of the uniformly
K-stable klt locus in the full projective space of degree-\(d\)
hypersurfaces. For \(3\le d\le n+1\) put \(m=\lfloor(d-1)/2\rfloor\).
Then
\[
 c_{n,d}\le
 \binom{n+d}{d-1}-(2n+1)
     -\sum_{j=1}^{m}\binom{n+d-2j-1}{d-2j-1}.
\]
The family giving this bound has general member a klt Fano
hypersurface with precisely one singular point and a strictly
negative Donaldson--Futaki invariant for an ambient degeneration.

Fix a point \(p=[1:0:\cdots:0]\) and a hyperplane through it,
with equation \(L(x)=0\), \(x=(x_1,\ldots,x_{n+1})\).
Use all members of the linear system
\[
 F=f_d(x)+\sum_{j=1}^{m}x_0^j L(x)^{j+1}g_{d-2j-1}(x).
\]
The system depends only on the flag: replacing \(x_0\) by a nonzero
multiple of itself plus a linear form in \(x\) preserves the
successive divisibilities by \(L^{j+1}\). Its coefficient space for
fixed flag has vector dimension
\[
       \binom{n+d}{d}+\sum_{j=1}^m\binom{n+d-2j-1}{d-2j-1}.
\]
The degree-\(d\) forms \(f_d\) give a base-point-free system away
from \(p\), so a general member is smooth outside \(p\) by
Bertini. Its multiplicity at \(p\) is \(d-m\ge2\). The lowest
nonzero homogeneous term at \(p\) is \(L^{m+1}g_{d-2m-1}\).
The last factor is constant in odd degree and linear in even
degree, so that \(L\) is the unique repeated factor at a general
member. Thus the member determines its flag, and the dimensions
\((n+1)+n=2n+1\) for the point and its hyperplane are gained in
the incidence. Subtracting this dimension and that of the linear
system from the dimension of all equations gives the claimed bound.

In coordinates with \(L=x_1\), give \(x_0,x_1,x_2,\ldots,x_{n+1}\)
the weights \(-1-n\epsilon,1,\epsilon,\ldots,\epsilon\), for a
sufficiently small positive rational \(\epsilon\). Their sum is
zero and the weight of every allowed monomial is positive: this
holds directly on \(f_d\), and on the remaining terms the limiting
weight as \(\epsilon\) tends to zero is at least one. Clearing
denominators supplies a one-parameter subgroup strictly destabilizing
the equation. The ambient CM-weight formula yields a negative
Donaldson--Futaki invariant on its ample degeneration; normalization
can only decrease this invariant.

The fixed-flag system contains the klt hypersurface
\(x_0x_1^{d-1}+x_2^d+\cdots+x_{n+1}^d=0\). For completeness,
its derivatives isolate \(p\); the affine germ is the cyclic cover
of \(\mathbf A^n\) of degree \(d-1\) over the cone on the smooth
Fermat degree-\(d\) hypersurface. Its branch coefficient
\((d-2)/(d-1)\) is less than \(\min\{1,n/d\}\), the threshold
given by blowing up the origin. The finite-cover discrepancy
formula proves klt. Openness of the klt condition in the flat
Gorenstein hypersurface family proves the claim for a general
member. Adjunction gives the Fano polarization since \(d\le n+1\).

More generally choose \(1\le a\le d-2\) and an \(a\)-dimensional
subspace of linear forms in the \(x\)'s, generating a linear
ideal \(I\). Put \(m_a=\lfloor(d-1)/(a+1)\rfloor\). Then
\[
\begin{split}
c_{n,d}\le {}&\binom{n+d}{d-1}-(n+1)-a(n+1-a)\\
 &-\sum_{j=1}^{m_a}
       \sum_{\ell=aj+1}^{d-j}
       \binom{a-1+\ell}{\ell}
       \binom{n-a+d-j-\ell}{d-j-\ell}.
\end{split}
\]
Indeed take the coefficient of \(x_0^j\) in
\((I^{aj+1})_{d-j}\) for each \(j=1,\ldots,m_a\), leaving
\(f_d(x)\) unrestricted. The inner sum counts its monomials
by their exact degree \(\ell\) in \(a\) adapted coordinates.
Give those variables weight one, the other \(x\)'s weight
\(\epsilon\), and \(x_0\) weight \(-a-(n+1-a)\epsilon\).
Every term is strictly positive by the same calculation. The
previous klt example belongs to the system by taking \(x_1\in I\),
since \(d-1\ge a+1\). The general member is smooth off \(p\)
and has lowest tangent multiplicity \(d-m_a\ge2\).
The lowest nonzero tangent equation in \(\mathbf P^n\)
vanishes to order at least \(am_a+1\ge2\) along the linear
subspace defined by \(I\). Its linear system is base-point-free
off that subspace, because it contains the pure
\((d-m_a)\)-th powers of each form in \(I_1\).
Bertini gives smoothness of its divisor everywhere else.
Its reduced singular locus therefore recovers the linear subspace
and hence \(I_1\). The flag contributes
\((n+1)+a(n+1-a)\) parameters, as required.

For example minimizing these upper bounds over \(a\) gives
\(c_{3,4}\le23,\ c_{4,5}\le98,\ c_{5,6}\le381\) at index
one and \(c_{9,4}\le256\) for the full quartic ninefold space.
These upper bounds concern the full space, since the members
exhibiting failure are singular.

There is a larger family whose general member fails uniform
K-stability. Allow \(1\le a\le d-1\), put
\(m'_a=\lfloor d/(a+1)\rfloor\), and replace
\((I^{aj+1})_{d-j}\) by \((I^{aj})_{d-j}\) through
\(1\le j\le m'_a\). The resulting upper bound is
\[
\begin{split}
c_{n,d}\le {}&\binom{n+d}{d-1}-(n+1)-a(n+1-a)\\
 &-\sum_{j=1}^{m'_a}
       \sum_{\ell=aj}^{d-j}
       \binom{a-1+\ell}{\ell}
       \binom{n-a+d-j-\ell}{d-j-\ell}.
\end{split}
\]
As before, the general member is klt and smooth away from \(p\).
For \(am'_a\ge2\) the singular locus of its lowest tangent recovers
the subspace \(I\), by the same base-point-free argument.
The only remaining possibility is \(d=3,a=1\). Its tangent
is a product of two linear forms, so there are at most two choices
for \(I\), still giving a generically finite marked incidence.
Use the weights \((-a,1,\ldots,1,0,\ldots,0)\), with \(a\)
positive entries. All monomials have nonnegative weight. The
general initial equation has weight zero and involves a pure
power of a zero-weight variable, since those variables exist
and \(f_d\) is general. Thus it is neither a pure highest-weight
nor a pure lowest-weight variable power. In its hypersurface
section ring the powers of both a highest-weight and a
lowest-weight variable survive in every degree. The leading
mean weight per section degree is zero: the ambient trace
vanishes, as does the weight of the single defining equation
in the hypersurface character series. Explicitly, multiplication by a
highest-weight variable is injective since the initial equation
has a pure zero-weight monomial. For \(a/(a+1)<\theta<1\), its
\(\lceil\theta k\rceil\)-th power embeds the sections of degree
\(k-\lceil\theta k\rceil\) in the degree-\(k\) sections of weight
at least \(((a+1)\theta-a)k\). Their asymptotic proportion is
\((1-\theta)^n>0\). Hence the one-parameter
degeneration has positive non-Archimedean \(J\)-norm and zero
Donaldson--Futaki invariant. Normalization preserves the norm
and can only decrease that invariant, giving failure of uniform
K-stability. For \(a=1\) one can also read the computation off
the two weights \(+1,-1\). The upper bounds for the four examples
in the preceding paragraph improve to \(17,77,303,211\),
respectively.

For every quartic Fano dimension the choice \(a=1\) gives the
particularly simple upper bound
\[
             c_{n,4}\le\binom{n+4}{3}-\binom{n+2}{2}-2n-2.
\]
The general member has exactly one singular point, a double
point of quadratic rank one, with quartic residual on its
Hessian kernel: in affine coordinates at \(p\) its equation is
\(cL^2+L g_2(x)+f_4(x)\), \(c\ne0\), and completing the
square gives a quartic leading term on \(L=0\).

\paragraph{Singular cyclic covers from independent hyperplane arrangements.}
Fix $d\ge3$, disjoint blocks of $q_i\ge2$ variables with
$\sum_i q_i\le d$, and $n\ge d-1$. In the $i$th block choose $d$
nonzero linear forms $\ell_{ij}$, with repetitions allowed. In the
remaining $n+2-\sum_iq_i$ variables $z_k$ consider
\[
 X=\left(\sum_i\prod_{j=1}^d\ell_{ij}+\sum_k z_k^d=0\right)
       \subset\mathbf P^{n+1}.
\]
Call a form outside the $\mathbf Q$-Fano locus a failure of
K-semistability. Then $X$ is K-semistable exactly when for each
block, each $1\le a<q_i$, and each $a$-dimensional vector subspace
$U$ of its dual space,
\[
                 \#\{j:\ell_{ij}\in U\}\le da/q_i.
\]
It is uniformly K-stable exactly when all these inequalities are
strict. In the product of the ordered arrangement parameter spaces
the exact codimension of failure of uniform K-stability is
\[
 \min_i\min_{1\le a<q_i}(q_i-a)(\lceil da/q_i\rceil-a).
\]
For failure of K-semistability replace $\lceil da/q_i\rceil$ by
$\lfloor da/q_i\rfloor+1$. The same codimensions apply after
independently forgetting the order in each block. For binary blocks
the tests are simply maximum root multiplicity $\le d/2$ and
$<d/2$, respectively. These are statements in the indicated split
families; the numerical codimensions refer to their own parameter spaces.

To prove sufficiency, first consider a degree-$d$ homogeneous sum of
nonzero forms in disjoint blocks of sizes $q_i$, now including any
number of one-variable blocks, with $\sum_iq_i=d$. Its projective
divisor in $\mathbf P^{d-1}$ is log canonical precisely when
\[
        \operatorname{lct}(\mathbf P^{q_i-1};(f_i=0))\ge q_i/d
                       \qquad\text{for each }i.
\]
For a singleton the projective divisor is empty. The threshold of
a homogeneous block at the affine origin is the minimum of the
projective threshold and $q_i/d$, by blowing up the origin and
pulling back a projective log resolution. At a singular point of
the projective sum, the homogeneous Euler identities give that all
active blocks are critical points of their respective zero sets.
In the chart of one nonzero block the germs in the different blocks
remain independent. For any other nonzero block its extra radial
coordinate multiplies the dehomogenized germ by a unit, preserving
its local threshold. The Thom--Sebastiani formula gives the minimum
of one and the sum of the block thresholds. This proves sufficiency.
If a block has smaller projective threshold, take only that block
nonzero at a point attaining it. The other block thresholds at the
origin are at most $q_j/d$, so the sum is strictly below one. For
products of linear forms the displayed thresholds are equivalent to
the weak configuration inequalities: blowing up the intersections
of the arrangement in increasing order gives the ratio of their
codimension to the number of hyperplanes through them, counting
multiplicity. See M.~Musta\c{t}\u{a}, \emph{Multiplier ideals of
hyperplane arrangements}, Trans. Amer. Math. Soc. 358 (2006),
5015--5023, with the same log-resolution calculation for repeated
hyperplanes.

The strict inequalities make each block GIT stable for its special
linear group. Indeed the minimum weight of a product is the sum
of the minima for its factors. Decomposing the weights of any
traceless one-parameter subgroup into the gaps of the associated
flag gives the usual subspace criterion for an equally weighted
point configuration. A disjoint sum of stable forms is stable:
use Kempf--Ness to choose a balanced representative in each block,
then rescale the block coordinates by positive scalars of total
determinant one so the squared norms of the blocks divided by
their dimensions coincide. The diagonal moment-map components now
vanish, and the off-diagonal components vanish by multidegree since
$d\ge3$. Thus the sum is polystable. An infinitesimal off-diagonal
stabilizer would give a zero constant-coefficient derivative of
one of the blocks, by separating bidegrees $(d-1,1)$. This is
impossible for a stable form, as a missing variable gives a
destabilizing one-parameter subgroup. On each block a diagonal
infinitesimal stabilizer is scalar by stability, and comparing
the summands makes the scalars equal. This proves GIT stability
of the projective sum, including singleton blocks.

Fill the blocks to $d$ variables with pure powers. For the weak
inequalities, the resulting divisor is anticanonical and lc and
its first cyclic suspension is K-semistable by
Section~\ref{sec:new-semistable-branch}. For the strict
inequalities its first suspension is uniformly K-stable by
Section~\ref{sec:nov28-suspension}. The same sections propagate
the two assertions through all further pure powers.

For necessity, suppose $k$ of the factors lie in $U$ and that
$kq_i-da\ge0$. Choose coordinates adapted to $U$ and assign
weights $q_i-a$ on $U$ and $-a$ on its complement in this block,
with weights zero outside the block. The minimum weight of the
block polynomial is $kq_i-da$, where $k$ counts all factors in $U$.
The minimum weight of the full equation is zero. If $kq_i-da>0$,
perturb by a small positive scalar weight on all the other blocks
with the compensating negative scalar weight on this block to
remain traceless. All monomials of the full form then have positive
weight. After clearing denominators, the ambient hypersurface
weight formula gives negative Donaldson--Futaki invariant. If
$kq_i-da=0$, the unperturbed nontrivial subgroup gives zero Futaki
invariant and positive non-Archimedean $J$-norm. Indeed the initial
equation contains a pure power outside the chosen block and is not
a pure power of any variable in that block, so the arbitrarily high
powers of its highest- and lowest-weight variables both survive.
The leading mean weight is zero by the hypersurface character series
(both the ambient trace and the defining equation weight are zero).
Thus the maximum weight differs linearly in section degree from the
mean. Normalization preserves the norm and only decreases the Futaki
invariant. These configurations establish both necessities.

Finally, for any chosen $k$ marks the condition that they span at
most an $a$-space has codimension $(q_i-a)(k-a)$ when $a\le k\le d$.
Indeed the incidence over the Grassmannian has the stated dimension
and the general marked set in this incidence spans its chosen
$a$-space. Taking the finite union over choices of marks, and then
the minimum over $a$ and the independent blocks, gives the formulas.

The K-polystable members can also be distinguished among the
K-semistable ones. The criterion is that every block arrangement
is GIT polystable. Indeed, replace the forms in the blocks by their
closed affine-orbit representatives using zero-weight degenerations.
The limits still split into linear factors: the projective product
map from the product of the dual projective spaces is proper.
They are semistable configurations, hence satisfy the weak count
above. The disjoint sum giving the Calabi--Yau branch is polystable
by the same block-balancing argument, which only uses polystability
for the moment map. It is log canonical by the threshold computation.
The interpolation in Section~\ref{sec:dec1-index-one-neighborhood}
yields a K-polystable first cyclic cover. All remaining pure powers
preserve K-polystability by \cite[Theorem~5.1]{ar2001.10557}.
This proves sufficiency and supplies an embedded K-polystable
degeneration $X_0$ of any K-semistable $X$ in the family. All its
nonzero fibers are isomorphic to $X$ since the block changes are
special linear. The family is a relative hypersurface with klt
fibers and relatively ample anticanonical bundle.

If the original $X$ were K-polystable it would be isomorphic to
$X_0$, by the uniqueness of the closed point in a geometric
K-moduli fiber. But if a block orbit was nonclosed then passage to
its closed orbit strictly increases its stabilizer dimension.
For a semistable affine form the affine and projective special
linear stabilizers have equal dimension, since a nonzero homogeneous
invariant at the form forces the character scaling it to have
finite image. Its block uses all variables by semistability. The
same multidegree calculation as above shows that the dimension of
the projective Lie stabilizer of the full sum equals the sum of the
dimensions of the projective Lie stabilizers of the blocks. Indeed
the block diagonal matrices preserve the line of each summand,
and equality of their scalar actions leaves only the common scalar
in the scalar matrices modulo which one passes to the projective
Lie algebra. Singleton blocks contribute zero. These dimensions
equal the dimensions of the abstract automorphism groups: the
connected group preserves the hyperplane line bundle since
$\operatorname{Pic}^0=0$, and the complete hyperplane system embeds
each hypersurface. Thus, if any orbit was nonclosed, $X_0$ has
strictly larger automorphism dimension, which excludes the claimed
isomorphism.

More generally replace each product by an arbitrary nonzero
degree-$d$ form $f_i$ in its $q_i$ variables. If
$\operatorname{lct}(\mathbf P^{q_i-1};(f_i=0))\ge q_i/d$
for every block, the same Fermat-completed hypersurface is
K-semistable. It is K-polystable exactly when all the block forms
are GIT polystable, and uniformly K-stable exactly when all are
GIT stable. Indeed the disjoint-sum calculation gives an lc
Calabi--Yau branch, which already proves semistability. Each block
is GIT semistable: the pair with boundary $(q_i/d)(f_i=0)$ is
lc log Calabi--Yau. For $q_i\ge3$, interpolating projective-space
K-semistability with this endpoint gives K-semistability for any
smaller positive coefficient, including the small GIT coefficient.
For $q_i=2$ the condition is exactly the binary-form criterion
that every root multiplicity be at most $d/2$. In a zero-weight
special-linear degeneration of that block the limiting boundary
remains lc. The crepant argument in
Section~\ref{sec:dec1-index-one-neighborhood} applies to the
rational $d$-fold pluricanonical form
$\Omega^d/f_i^{q_i}$ with the $d$-fold logarithmic differential
in the parameter. Thus the preceding closed-orbit and stabilizer
arguments apply unchanged to these blocks. At a closed orbit,
finiteness of its projective stabilizer is equivalent to GIT
stability, which gives the last assertion.

\section{Line thresholds and the exceptional locus for quartic ninefolds}\label{sec:nov30-ninefold-results}

\paragraph{Factoriality of the tangent section and the intermediate ranks.}
Let \(X^9_4\subset\mathbf P^{10}\) be smooth, and at a point \(p\)
of quadratic tangent contact write the tangent hyperplane section as
\[
             D=(T^2 q_2+T q_3+q_4=0)\subset\mathbf P^9 .
\]
Every point of Hessian rank four, and every point of rank at least
three at which \(q_2\mid q_3\), is excluded as a center of a
divisor computing \(\delta(X,H)\le7\). At rank five the same
exclusion holds unless both \(q_3,q_4\) vanish on the radical
of \(q_2\). In the remaining regular rank-five case and in
all regular rank-three cases every such computing divisor,
normalized by \(v(\mathfrak m_p)=1\), satisfies \(S_H(v)\ge43/24\).

The Gauss morphism of the smooth quartic is finite, so \(D\) has
isolated singularities. In particular the quartic is integral.
Its local rings are factorial by parafactoriality of isolated
complete-intersection singularities of dimension at least four,
since outside the closed singular point all Weil divisors are
Cartier. By the Grothendieck--Lefschetz hyperplane theorem
\(\operatorname{Pic}(D)=\mathbf Z[H]\). If \(q_2\) is irreducible,
the scheme \(D\cap(q_2=0)\) is integral: a nontrivial effective
decomposition of this Cartier divisor of class \(2H\) would have
two components counted with multiplicities, each of class \(H\).
Their sections lift to linear forms on projective space. Their
product would equal \(q_2\) up to scalar on \(D\), hence also in
the ambient polynomial ring since \(D\) has degree four. The
divisor is Cohen--Macaulay, so generic reducedness suffices.

Consequently for every rank at least three, either \(q_2\mid q_3\)
or \(q_2,q_3,q_4\) form a regular sequence. Indeed a divisorial
component \(\Gamma\) of \((q_3=q_4=0)\) on the irreducible
quadric \(Q=(q_2=0)\subset\mathbf P^8\) would give the cone
over \(\Gamma\) with vertex \(p\) as an irreducible component
of \(D\cap(q_2=0)\). If \(q_3\notin(q_2)\), that intersection
also contains points above a dense open of \(Q\) where \(q_3\ne0\).
This contradicts integrality.

In the regular case at rank four the radical cannot be contained
in \(V(q_3,q_4)\): its join with \(p\) would give a five-plane
contained in the smooth quartic ninefold. The same argument works
at smaller ranks. At rank four or five with a radical point
outside \(V(q_3,q_4)\), the isotropic plane and line construction
of Section~\ref{sec:q8-tangent-tests} gives the separated
\(sz,f_3,g_4\) model. In units of ten its bound is
\(10S_H(v)\le65/6+(3/2)q\), with \(q=A_X(v)-9\).
Thus \(A_X(v)\le7S_H(v)\) forces \(S_H(v)\ge16/3\).

In the remaining regular rank-five case take a general isotropic
projective line \(\ell\) disjoint from the vertex and from the
codimension-two base on \(Q\). Such lines exist by lifting a line
of the nonsingular three-dimensional quadric in the quotient and
using transitivity on the smooth locus. For its vector span \(S\),
the form on \(S^\perp/S\) has rank one; choose \(z\in S^\perp\)
with \(q_2(z)\ne0\). In the initial-form calculation give the two
\(S\) directions weight zero, \(z\) weight one, and all other
directions weight at least three. The resulting regular triple is
\(z^2,q_3|_\ell,q_4|_\ell\). The two binary restrictions are
coprime and the normal-basis calculation for the rank-one initial
triple gives \(10S_H(v)\le65/6+2q\). Substituting \(9+q\le7S_H(v)\)
gives the stated \(43/24\).

At rank three choose a general radical point \(r\) off the base,
as is possible by the same contained-plane obstruction. The
family of nonradical isotropic lines through this point covers
\(Q_{\rm sm}\); its stabilizer acts transitively on \(Q_{\rm sm}\)
by radical shears and the orthogonal group of the rank-three
quotient. Thus a general such line avoids the base there by
codimension two, and it meets the vertex only at \(r\).
The vector span \(S\) projects to a one-dimensional isotropic
space in the ternary quotient, so \(S^\perp/S\) also has rank
one. The same \(z^2,f_3,g_4\) specialization and estimate apply.

For the divisible case, replace \(T\) to write the tangent
equation \(T^2q_2+g=0\), and use \(g\) modulo \(q_2\) in the
two-symbol presentation. Its zero divisor on the smooth locus
of \(Q\) has at most isolated singularities by the \(T\)-fiber
and finite Gauss argument. In ranks four and five choose a
general radical point \(r\). If \(g\) vanishes on the whole
radical, then \(dg(r)\ne0\): otherwise every point of the
projective \(T\)-line over \(r\) is singular on \(D\).
Since \(Q\) spans \(\mathbf P^8\), a general nonradical
isotropic direction through \(r\) gives a simple zero of \(g\)
at \(r\). If \(g\) is nonzero at \(r\), there is no mark there.
Off \(r\) a general such line is transverse to the smooth locus
of the divisor and avoids its singular locus. Indeed the
stabilizer of \(r\) acts transitively on \(Q_{\rm sm}\) by the
orthogonal group of the quotient and the linear radical shears.
Extend the line to an isotropic plane with radical intersection
\(r\), as in the regular case. In rank at least six take a line
and then a plane disjoint from the radical instead.

At rank three take a general line lying wholly in the radical.
The restriction of \(g\) to its projectivization \(\mathbf P^5\)
is nonzero and reduced. In fact the span with \(p\) has ambient
codimension four, so its section has singular locus of dimension
at most three: a four-dimensional component would meet all four
remaining partial derivatives and contradict smoothness of \(X\).
This section is the cone over \(g|_{\mathbf P^5}\); a nonreduced
component would give a singular component of dimension at least
five, and the zero polynomial would give a contained six-plane.
Thus the general radical line cuts this quartic transversely.
Add any nonradical isotropic direction to obtain the required
plane. In every divisible case the separated squarefree binary
quartic comparison gives \(10S_H(v)\le11+(3/2)q\), hence
\(S_H(v)\ge5\) at a possible optimizer.

Finally, every line on a smooth quartic ninefold has threshold
at least \(220/31\), by the eightfold line estimate and the
normalized hyperplane flag. For \(S_H(v)\ge5\) use the
two-hyperplane interpolation at \(t=92/13<220/31\) and
\(\theta=1/4\). Its Nadel expression satisfies
\[
             (t-7)5+7-\tfrac9{10}t=\tfrac{66}{65}>1 ,
\]
and the marked ample remainder \(7-(9/10-1/40)t\) lies strictly
between \(1/2\) and one. This excludes all the asserted optimizing
points without an additional incidence hypothesis on the lines.

\subsection{A ternary comparison for every regular quadratic contact}
\label{sec:nov30-nine-ternary}
On a smooth quartic ninefold let \(p\) have quadratic contact of
rank at least two, and assume the triple \(q_2,q_3,q_4\) is regular.
For each divisorial valuation centered there, normalized by
\(v(\mathfrak m_p)=1\), put \(q=A_X(v)-9\). Then
\[
               10S_H(v)\ \le\ \frac{65}{6}+\frac{15}{8}q .
\]
In particular \(A_X(v)\le7S_H(v)\) implies \(S_H(v)\ge29/15\).

Take a general three-dimensional vector subspace in the nine
tangent variables. Its projectivization misses the codimension-three
base and its quadratic restriction \(F_2\) has rank at least two.
Kill the six complementary variables in all three forms by the
first diagonal specialization of
Section~\ref{sec:q8-tangent-tests}. Their ternary restrictions
\(F_2,F_3,F_4\) remain a regular sequence. Thus the presented
initial spaces have the same Hilbert function, and the specialized
basis-divisor estimates bound the original ones. Write \(R\) for
the ternary polynomial algebra and \(J=(F_2,F_3,F_4)\).
Splitting off the six complementary characters and the total
ternary degree \(D\), the block at section degree \(m\) and
ordinary order \(j\) is \((J^e)_D\), \(e=(j-m)_+\).
Throughout the calculation an integer ceiling in a bounded-width
endpoint changes the unnormalized orders only by \(O(1)A\).

Every element of the block has the factor \(F_2^I\) for
\(I=(3e-D)_+\), since each other generator has degree at least
three. Dividing gives exactly \((J^{e-I})_{D-2I}\): both
inclusions follow from the generators, and \(I\le e\) for a
nonempty block. At a valuation of the angular plane its fixed
factor has order at most \(I A\), since a reduced plane quadric
is log canonical. For the residual block one may extend any
valuation-adapted basis to a full basis of the ternary polynomials
of degree \(D'=D-2I\). All orders are nonnegative, and the
complete series has mean at most \((D'/3+o(m))A\) in uniform
finite-level comparisons. Thus the contribution of the residual
block to the sum of orders is at most the full space dimension
times this complete mean. The next computation counts this
enlarged contribution without a loss from any small block.

In the normal basis the ternary algebra is finite free over
\(\mathbf C[F_2,F_3,F_4]\), of rank \(24\); choose a homogeneous
complement basis in the quotient by \(J\), of bounded degrees.
Write \(a,b,c,t\) for the exponents of the degree-\(1,2,3\)
hypertangent sections and the coordinate \(T\). After dividing
by \(m\), multiplying the mean of any homogeneous degree-one
function of these four exponents by ten is equivalent in the
simplex limit to its expectation for independent exponentials
of rates \(1,2,3,1\), respectively. This follows by dividing
the ten independent unit exponentials by their sum, using the
six complementary coordinates and the three indicated costs.
In these variables, up to the bounded complement degree,
\[
     D=2a+3b+4c,\qquad e_{\rm raw}=j-m=a+b+c-t,\qquad
                      I=(a-c-3t)_+ .
\]
Consequently the scaled mean of \(I\) is
\(\mathbf E(e^{-c-3t})=3/16\).

We spell out the enlarged residual mean. At fixed
\(s=e_{\rm raw}/D\le1/2\), put
\(h=1-2(3s-1)_+\), so \(D'=Dh\). In the triangular simplex
\(2a+3b+4c=D\) let \(f(s)\) be the area fraction satisfying
\(a+b+c\ge sD\). The leading dimension of the original
block is \(f(s)D^2/2\). Indeed, for each of the \(24\)
bounded complementary basis elements the exponent lattice on
this triangle has full leading density \(D^2/48\).
Polygon areas compute these counts with uniform \(O(D+1)\)
error by bounded-width boundary strips. After dividing the
fixed factor the complementary space of all ternary sections
has dimension \(D'^2/2+O(D'+1)\).
All these enlargements, and their order errors of size at most
\(O(m)A\) per added boundary vector, contribute \(o(1)A\) in
the normalized global mean: there are \(O(m^7)\) relevant
ternary blocks at fixed \(m\), by the six complementary
exponents and \(j\) (which determines \(D\)). For the complete
series one can first discard bounded \(D'\), then use its
uniform finite-level bound, and finally let the cutoff tend
to infinity.

Under the four exponential laws the joint density in \(D,s\)
is
\[
            \frac18 D^3 f(s)e^{-D(1-s)}\,dD\,ds .
\]
Indeed the exponent is \(a+2b+3c+t=D-e_{\rm raw}\), the
Jacobian for \(t\) is \(D\,ds\), the full triangle has
coarea \(D^2/48\), and the product of the four rates is six.
Extension to the ambient ternary space thus replaces \(f(s)\)
in this measure by \(h^2\). Multiplying by its mean \(Dh/3\)
and integrating in \(D\), its total angular coefficient becomes
\[
 \int_{-\infty}^{1/2}\frac{h^3}{(1-s)^5}\,ds
     =\frac{81}{64}+\frac{27}{64}=\frac{27}{16}.
\]
The split is at \(s=1/3\). In the second interval
\(w=(1-s)^{-1}\) turns the integral into
\(27\int_{3/2}^2(2-w)^3\,dw\).
Adding the fixed quadric coefficient gives \(15/8\).
The ternary vertex has radial coefficient per unit discrepancy
\(\overline D/3=29/18<15/8\); the six other coordinate
coefficients are one. Blow up that vertex and use the angular
plane comparison and the complementary coordinate strata as
in Section~\ref{sec:q8-tangent-tests}. The original ordinary
radial intercept is \(9+1+1/2+1/3=65/6\). Filtered inversion
of adjunction proves the claimed signed inequality, and direct
substitution gives \(S_H(v)\ge29/15\).

It follows that for a smooth quartic ninefold a common bound on
the lines \(\delta_C(X,H)\ge M\), with any \(M>1490/203\), implies
uniform K-stability. Equality suffices for K-semistability.
For the quadratic contact alternatives, if a linear factor of
\(q_2\) divides \(q_3\), the common-factor extraction of the
point criterion excludes the original optimizer. For rank one
with no common factor the square and coprime binary forms give
\(S_H(v)\ge43/24\). For rank two with no common factor the
triple is regular: on each component hyperplane of \(q_2=0\),
the section \(Tq_3+q_4\) is integral by factoriality of the
tangent section, and its nonzero cubic coefficient precludes
a common factor of the restrictions. In all higher ranks this
regularity or the divisible alternative was proved by the same
factoriality argument. Divisibility is excluded by the
two-symbol line test at the actual all-member line floor.
For cubic and quartic contact use the preceding irreducible-cubic,
linear-factor and generalized-Eckardt calculations. Their
combined lower bound for all remaining optimizers is
\(S_0=41/23<43/24<29/15\).
In the two-hyperplane interpolation take \(\lambda=0,\theta=1/4\)
and \(t_0=1490/203\). Then
\[
    (t_0-7)S_0+7-9t_0/10=1,\qquad
             \frac12<7-(9/10-1/40)t_0<1 .
\]
Testing at \(t\) slightly above \(t_0\) proves the strict
assertion. At the endpoint, any \(\beta<7\) gives a strict
Nadel expression at \(t_0\), so one may test immediately below
\(t_0\), proving K-semistability.

\paragraph{A weighted estimate at a double tangent hyperplane.}
At rank-one quadratic contact on a quartic ninefold put \(q_2=x^2\)
and suppose \(x\nmid q_3\). Write \(A=A_X(v)\), \(b=v(x)\),
and \(S=S_H(v)\) with \(v(\mathfrak m_p)=1\). If \(b\ge3/2\),
then \(10S+2b\le3A-41/3\).

Take the weighted point extraction with weights three on \(x\)
and two on the other eight tangent parameters \(y\). Its
discrepancy is nineteen. With \(w\) the transverse coordinate
of the original quartic write its full equation as
\(T^3w+T^2 Q_2+TQ_3+Q_4=0\). Here the homogeneous \(Q_j\)
can involve \(w,x,y\). Set \(W=T^2w+TQ_2+Q_3\), a section of
degree three. The weighted section presentation has \(w,W\) of
weights six and eight, and relations
\[
        T^2 w+T x^2+f_3(y)=0,\qquad T W+g_4(y)=0,
\]
where \(f_3=q_3|_{x=0}\), \(g_4=q_4|_{x=0}\) are coprime.
This follows by the actual identity \(TW=-Q_4\) and then the
definition of \(W\). Modulo \(T\) the two forms are regular;
localization at \(T\) embeds the presented ring in the associated
graded ring of the smooth point, and its ordinary Hilbert series
is the quartic Hilbert series. Thus the presentation is exact.
Scale \(x\) away in the first relation, and specialize \(y\)
to a general binary subspace. The two forms remain coprime
there, so these diagonal specializations preserve the Hilbert
functions of all the residual spaces as above. The marked
coordinate divisor \(x=0\) is preserved by both operations.

The two-symbol monomials have exponents \(i,h\) on \(w,W\)
and \(t\) on \(T\). Their ordinary costs are \(1,3,1\).
The binary ring is finite free over its two forms. There are
six free \(y\)-coordinates and one free \(x\); in units of ten
each free exponent has mean one. The scaled weighted radial
mean is \(6\cdot2+3+6+8/3=71/3\).
After fixing the free characters, a binary block of degree \(D\)
is a sum of \(f_3^a g_4^c\mathbf C[y]_{D-3a-4c}\) with
\(2a+c\ge e\). Here \(e=i\,2+h-t\) on the indicated normal
monomials, since substituting the two symbols raises light
degree over section cost by two and one. Bounded binary
complement shifts have no effect on the normalized formulas.
The forced power \(I\) of \(f_3\) is
\(\max(0,\lceil(4e-D)/5\rceil)\) up to a bounded endpoint.
After its removal the complete binary remainder has bounded
codimension: in the free binary module each complementary
degree is bounded, and the admissible \(f_3\)-exponents run
from this endpoint to the full upper endpoint in steps of four.
The mean of \(D\) is \(13/3\); the mean of \(I\) equals
\(\mathbf E[(i-4t/5)_+]=5/9\) in the usual simplex units.
Thus every binary angle has scaled coefficient at most
\(13/6+(3/2)(5/9)=3\), even at a triple root. The binary vertex
has ratio \(13/6\), and every free coordinate, including \(x\),
has coefficient one.

Let \(u=v(E)\) be the weighted exceptional order. The radial
comparison and angular inversion with the two unused units of
discrepancy on the \(x\)-divisor give
\[
             10S\le (71/3)u+3(A-19u)-2(b-3u).
\]
One can see the sign directly by adding twice the strict transform
of \(x=0\) to the angular basis divisor: its resulting coefficient
is three, and the binary and other coordinate bounds remain
at most three on the product resolution. If \(b\ge3/2\),
some \(y\) has \(v(y)=1\), so \(u=\min(b/3,1/2)=1/2\).
This proves the stated weighted comparison. In particular, at an
optimizer of ratio at most seven equality in the ordinary bound
\(S\ge43/24\) requires \(b\le145/48\) (the same conclusion is
automatic below \(3/2\)).

\paragraph{The forced cubic power in the disjoint tangent model.}
In the ternary cubic and binary quartic model for irreducible
cubic contact with at most four essential variables the ternary
angular coefficient is at most \(9/5\). To check the refinement,
write \(i,h,t\) for the normal exponents of cubic, quartic,
and \(T\). Their simplex costs are \(1,3,1\), and
\(e=j-m=2i+h-t\). First refine by a marked point of the
binary quartic block. Choose a binary basis compatible with
both its marked order and the filtration by powers of the
quartic. Since the jumps for the point order in fixed binary
degree are simple, each basis vector of quartic order \(h\)
selects exactly the ternary space
\[
       f^{I}\mathbf C[x]_{D_1-3I},\qquad
                  I=\max(0,\lceil(e-h)/2\rceil).
\]
This is also seen by reversing the two filtrations in the sum
over \(2a+c\ge e\): a binary basis vector of order \(h\)
belongs to every quartic-power subspace with \(c\le h\).
The fixed multidegrees and their multiplicities are unchanged
by replacing the binary remainder basis in this fashion.
For these normal indices \(I=(i-t/2)_++O(1)\), and its scaled
mean is \(2/3\). The ternary degree mean is five; the integral
plane cubic has threshold at least \(5/6\), so the bound on
the ternary factor improves to \(5/3+(2/3)/5=9/5\).
The marked binary factor and the other forward-polar cubic
alternatives have slope at most \(11/6\). Thus at every
irreducible cubic contact
\[
                 10S\le34/3+(11/6)(A-9),\qquad
                 A\le7S\ \Longrightarrow\ S\ge31/17 .
\]

Consequently the line-only criterion in
Section~\ref{sec:nov30-nine-ternary} holds with the smaller
mark \(M>785/107\), and K-semistability holds at equality.
Indeed \(43/24<31/17<29/15\). Take
\(S_0=43/24,\ t_0=785/107,\ \theta=1/4\) there:
\((t_0-7)S_0+7-9t_0/10=1\), and
\(1/2<7-(9/10-1/40)t_0<1\). The same strict and endpoint
arguments apply. If the rank-one contacts without a common
factor are absent, the mark improves to \(1150/157\), by
the same substitution with \(S_0=31/17\).

\paragraph{Removing an actual power of the tangent hyperplane.}
\label{par:nov30-rank-one-fixed-tangent}
At a rank-one quadratic contact without a common factor, the uniform
lower bound for a possible minimizing divisor improves to
\[
                      S_H(v)\ge\frac{21929}{12225}.
\]
Consequently the quartic ninefold is uniformly K-stable whenever
every line has local hyperplane threshold at least \(M\), for any
\(M>160306/21853\). At equality it is K-semistable.

Write \(D=(w=0)\) for the tangent hyperplane section of the smooth
ninefold. We first show, with the normalization
\(v(\mathfrak m_p)=1\), that
\[
                        v(D)\le A_X(v)-19/3.                  \tag{R1}
\]
In the tangent \(\mathbf P^9\) the hypersurface \(D\) has isolated
singularities by finiteness of the Gauss map. In affine tangent
coordinates its equation is \(x^2+q_3(x,y)+q_4(x,y)\). Write
\[
 C_3=q_3(0,y)\ne0,\qquad A_2=\partial_xq_3(0,y),\qquad
 C_4=q_4(0,y)-A_2^2/4 .
\]
If \(C_3\) has a repeated factor, it is linear, say \(l\).
Then \(C_4|_{l=0}\ne0\). Otherwise the projective section
of \(D=(T^2x^2+Tq_3+q_4=0)\) on \((x=l=0)\simeq\mathbf P^7\)
would have equation \(A_2^2/4\). Its singular set has dimension
at least six, even if the restriction is the zero polynomial.
Intersecting it with the two remaining normal partials would
produce a positive-dimensional singular set of \(D\).

Cut the eight \(y\) variables to a general two-dimensional space
of directions \((u,v)\), retaining \(x\). A reduced \(C_3\)
restricts to a squarefree binary cubic by Bertini. If it is
nonreduced, the intersection of \(l=0\) with the general
projective line is a point where \(C_4\ne0\), and the
multiplicity of this root is two or three. Make the triangular
coordinate change replacing the old \(x\) by \(x-A_2/2\).
It cancels the cubic cross term and the new pure binary term
of degree four is \(C_4\). Choosing binary coordinates for
the multiple root, the weighted initials are, up to rescaling,
respectively
\[
  x^2+C_3(u,v),\qquad x^2+u^2v+v^4,\qquad x^2+u^3+v^4
\]
in weights \((3,2,2),(4,3,2),(6,4,3)\).
Indeed all remaining terms have higher weight: pure binary
terms beyond degree four have ordinary degree at least five;
mixed terms other than \(x^2\) are divisible by \(x^2(u,v,x)\)
or belong to \(x(u,v)^3\).

On each initial the hypersurface with coefficient one and
the origin ideal with coefficient \(1/3\) give an lc pair near
the origin. Here is the Newton check. Each face polynomial
is nonsingular on its torus zero set: in the first case use
squarefreeness of the cubic, choosing both pure cube coefficients
nonzero, and in the other cases the displayed faces give this
directly. A simultaneous toric subdivision for the Newton
polyhedron and the origin ideal makes their transforms normal
crossings. For coordinate weights \((a,b,c)\) put
\(s=\min(a,b,c)\) and let \(m\) be the least monomial weight.
In the first and third cases,
\[
 a+b+c\ge\max(s,m/2)+\max(s,m/3)+\max(s,m/q),
 \qquad q=3,4
\]
respectively. Subtracting \(m\) gives at least \(s/3\);
check the affine intervals with endpoints \(2s,3s,qs\) and
the final ray. In the middle case \(a\ge m/2\) and
\(2b+c\ge m\) give \(a+b+c\ge m+s/2\). These inequalities
check all toric discrepancies; the strict transform itself has
coefficient one. Weighted specialization fixes the origin and
its ideal, so lower semicontinuity, applied for instance to
the cube of the hypersurface equation times the origin ideal
at coefficient \(1/3\), gives the same lc bound on the
three-dimensional slice.

Take six independent linear hyperplanes whose intersection is
this slice. Successive inversion of adjunction gives log
canonicity near \(p\) on the tangent affine nine-space with
all six hyperplanes, the hypersurface \(D\), and the origin
ideal of exponent \(1/3\). The smooth germ of \(X\) identifies
analytically with this affine space, carrying its tangent
section to the equation of \(D\) up to a unit by the implicit
function theorem. Each hyperplane has order at least one at
the normalized valuation, proving (R1).

Use the three hypertangent generators of section degrees \(1,2,3\)
and ordinary orders \(2,3,4\), the first being the actual section
\(w\). They give the complete ordinary point-order filtration
since the tangent triple is regular. In degree \(m\), every
section of order at least \(j\) is divisible on \(X\) by
\[
                  w^{I_0(j)},\qquad I_0(j)=(2j-3m)_+ .
\]
Indeed each homogeneous normal product of order \(j'\ge j\)
satisfies this divisibility when lifted to a product of the
actual hypertangent generators: only the first generator
contributes positively to twice the order minus three times
the degree, and its contribution is one. An initial basis in
every such grade admits these lifts; descending through the
finitely many grades, their union spans the full subspace of
order at least \(j\). Thus the factor is present before
restriction to the ordinary exceptional divisor.

On that exceptional divisor the marked tangent hyperplane has
equation \(x^2=0\). In the rank-one specialization of the
regular triple to \(x^2,f_3,g_4\), the cubic and quartic use
two other variables and are coprime. After deleting the above
actual factors, the residual angular coefficient on \((x=0)\)
decreases by \(2\alpha\) in units of ten, where
\[
 \alpha=\lim_m\frac{10}{m\dim R_m}
                         \sum_j I_0(j)\dim\operatorname{gr}^j R_m
        =\frac3{1024}.
\]
To compute it take exponents \(a,b,c,t\) on the three
hypertangent generators and on \(T\), respectively. There are
six further free unit-degree variables. A complement basis
modulo the binary forms and \(x^2\) has bounded degrees.
Thus in the usual scaled simplex mean,
\[
 \alpha=\mathbf E\big[(a-c-3t-e_1-\cdots-e_6)_+\big]
       =\frac34\cdot\frac14\cdot\left(\frac12\right)^6 ,
\]
with independent exponentials having rates three for \(c\)
and one for \(a,t,e_i\). The specializations preserve \(x\);
dividing each original angular layer by its known \(x^{2I_0(j)}\)
therefore divides its limiting layer by precisely the same
factor. In the original rank-one angular calculation the
coefficient of this coordinate hyperplane was two, while all
other coordinate and binary tests were at most \(3/2\). In
particular the new common angular bound is \(2-2\alpha>3/2\).
The residual basis divisor is effective on the actual point
blowup; no adjunction for negative boundaries is needed.

For each basis section remove \(I_0(j)\) copies of the actual
tangent divisor, whose exceptional order is two. Its strict
transform has order \(v(D)-2\) for the given valuation. The
ordinary extraction and the preceding residual angular bound
give
\[
 \begin{aligned}
 10S_H(v)
 &\le\tfrac{65}{6}+(2-2\alpha)(A_X(v)-9)
                              +\alpha\big(v(D)-2\big)\\
 &\le\tfrac{65}{6}+\tfrac23\alpha+(2-\alpha)(A_X(v)-9),
 \end{aligned}
\]
where (R1) gives the second line. For \(A_X(v)\le7S_H(v)\)
this yields \((24-42\alpha)S_H(v)\ge43-58\alpha\), proving
the claimed bound. This is smaller than \(31/17\) and \(29/15\).
In the preceding line interpolation substitute
\(S_0=21929/12225,\ t_0=160306/21853\): one has
\[
  (t_0-7)S_0+7-9t_0/10=1,\qquad
  7-7t_0/8=50813/87412\in(1/2,1).
\]
Perturbing \(t\) as in the previous criterion proves the
strict and endpoint assertions.

\paragraph{Two quadratic coefficients with high generic rank.}\label{par:nov30-audit-polar-pencils}
Let $V$ be the seven-dimensional common polar kernel of a line of polar rank two on a smooth quartic ninefold. Write $Q_t=t_0^2q_0+t_0t_1q_1+t_1^2q_2$, let $s=\dim\langle q_0,q_1,q_2\rangle$, and write $\rho$ for the rank of $Q_t$ at general $t$. The line has $\delta_L(X,H)\ge15/2$ in each of the following cases:
\begin{enumerate}
\item $s=1$ and $\rho\ge4$;
\item $s=2$ and $\rho\ge6$;
\item $s=2$, $\rho=5$, and every nonzero member of the pencil of quadrics has rank at least three.
\end{enumerate}
To prove this we seek a component $W$ of the reduced common zero set in $\mathbf P(V)$ whose affine tangent space at general $u\in W$ supports $Q_t$ of rank at least two and on which some cubic coefficient or the terminal quartic is nonzero. If some cubic coefficient is nonzero we obtain the $xy+u^3$ face. Otherwise the two cubic coefficients vanish on $W$, so their differentials vanish on its tangent space. Where the quartic is nonzero we obtain the $xy+u^4$ face. For clarity, take two constant directions $x,y$ in the affine tangent on which $Q_t$ has rank two for general $t$, and take the third direction along $u$. This is possible modulo the radial vector $u$, which is orthogonal to the affine tangent for all the common quadratic equations. The quadratic coefficients of $u^2,ux,uy$ vanish for every $t$. With weights $(3,3,2)$ the degree-six face is the nondegenerate binary quadratic and the nonzero cube; with weights $(2,2,1)$ the degree-four face is the binary quadratic and the nonzero fourth power. In the latter case both components of the cubic differential vanish on $x,y$, eliminating $xu^2,yu^2$. Higher-weight terms disappear, and the other complementary directions may be split off as in the original two-polar extraction. Choose the active polar coordinates generally, so all three roots used for the pole decomposition are distinct and avoid the finitely many parameters where the binary rank or the marked coefficient drops. At each such pole a linear change of the equal-weight variables $x,y$ gives $xy$, preserving the ordinary complete polynomial spaces on the chart. Section~\ref{sec:nov28-disjoint-general} gives the signed positive-sector denominators $10/7$ and $13/9$, respectively. The regular sector has share $1/4$ and denominator one, so the two resulting denominators are $37/28$ and $4/3$. They agree with the radial denominators $1+(3/4)/(1+8/6)$ and $1+(3/4)/(1+5/4)$, respectively. Dividing ten by these denominators gives $280/37$ and $15/2$.

When $s=1$, the common quadric is integral and its tangent restriction has rank $\rho-2$ at every smooth point. When $s=2$ the quadrics are a regular sequence since a common linear factor would make every member have rank at most two. Each reduced component has codimension two, so restriction to its affine tangent space lowers rank by at most four. In cases (1) and (2) the tangent test therefore holds on every component.

In these cases the higher forms cannot vanish on all components. Otherwise $X$ contains each join of the line with a component $W$. This join has dimension at least six and degree $\deg W$. Every integral subvariety of dimension at least five in a smooth quartic ninefold has degree divisible by four, by integral weak Lefschetz applied in its even cohomological codimension. In case (1) the degree is two. In case (2), B\'ezout would force a unique degree-four component of multiplicity one. The complete intersection is then integral, since it is Cohen--Macaulay and generically reduced. On the $\mathbf P^8$ spanned by the line and the kernel the equation would belong to the ideal of the two quadrics, hence would equal $qH+q'H'$ for quadrics $H,H'$. Where all four vanish this section is singular; that locus has dimension at least four. The finite Gauss map of $X$ sends it to the projective line of hyperplanes containing $\mathbf P^8$, an impossibility. This also proves nonvanishing whenever the common zero scheme is integral.

For case (3), the singular locus of the complete intersection is contained in the union of the projective radicals of pencil members. The generic radical has projective dimension one. There are only finitely many parameters where the rank drops below five and each of their radicals has dimension at most three. Hence the fourfold intersection is generically reduced. A linear component would give a five-dimensional totally isotropic subspace for a rank-five form in dimension seven, whereas the maximal such dimension is four. A degree-two component would span a hyperplane on which both quadratic equations were proportional, giving a nonzero pencil member divisible by a linear form. These exclusions and B\'ezout imply integrality.

Choose a general rank-five member $A$ and an independent member $B$. If the tangent test fails generically, both constraints have independent gradients at general $u$, and $T=\ker A(u,-)\cap\ker B(u,-)$ has codimension two. Then $\mathrm{rad}(A)\subset T$: otherwise the image of $T$ in $V/\mathrm{rad}(A)$ has codimension at most one and its restriction has rank at least three. In particular $B(\mathrm{rad}(A),u)=0$ generically on the intersection. It spans $\mathbf P(V)$ since it is an integral complete intersection of quadrics, hence $B$ vanishes on $\mathrm{rad}(A)\times V$. Pass to $V/\mathrm{rad}(A)$ of dimension five and write $S=A^{-1}B$. There $T^\perp=\langle u,Su\rangle$. The assumed rank at most one on the three-dimensional $T$ forces $T^\perp\subset T$, whence $A(Su,Su)=0$ generically. The quadratic ideal of the intersection is precisely $\langle A,B\rangle$, so $S^2$ lies in $\langle I,S\rangle$. If $S$ has two distinct eigenvalues, it is diagonalizable with orthogonal eigenspaces and the pencil contains the form supported on its smaller eigenspace, of rank at most two. If $S$ has a single eigenvalue, the nonzero nilpotent part has square zero, hence rank at most two. Either alternative contradicts the hypothesis of case (3). This proves the tangent test and the claimed bound.

\paragraph{Three quadratic coefficients with a nonsingular member.}\label{par:nov30-audit-three-quadrics}
Suppose the three coefficients are independent, their common zero scheme is an integral complete intersection in $\mathbf P^6$, and $Q_t$ is generically nonsingular. Then the tangent test, and hence the line bound $15/2$, always holds. We prove the matrix fact that rules out the exceptional algebra in the three-coefficient test.

Let $J$ be a three-dimensional vector space of self-adjoint operators on a seven-dimensional complex vector space with a nondegenerate symmetric form. Suppose $I\in J$ and $J$ is closed under the Jordan product $(BC+CB)/2$. We claim that $J$ contains a nonzero operator of rank at most two. All powers of each element belong to $J$, so every minimal polynomial has degree at most three. If some $B$ has minimal polynomial of degree three, $J=\mathbf C[B]$. With three distinct eigenvalues, an eigenspace projection has rank at most two. With two eigenvalues, take the projection onto the simple-root generalized eigenspace if it has dimension at most two. Otherwise the complementary double-root space has dimension at most four; the nonzero nilpotent part supported there is square-zero and has rank at most two. With a single eigenvalue $\lambda$, $(B-\lambda I)^2$ is nonzero, with rank bounded by the number of three-dimensional Jordan blocks, at most two.

Now suppose every minimal polynomial has degree at most two. Subtract half the respective linear coefficients from two generators modulo $I$ to obtain $B,C$ with $B^2,C^2$ scalar. Writing $BC+CB=eI+fB+gC$ and comparing the coefficients of $B,C$ in $(B+tC)^2\in\langle I,B+tC\rangle$ for general $t$ gives $f=g=0$. Thus $(xB+yC)^2=q(x,y)I$ for a binary quadratic form. If $q$ is nondegenerate, two orthogonal normalized generators are anticommuting involutions. One exchanges the plus and minus eigenspaces of the other, forcing even total dimension. If $q$ has rank one, take an involution and a nonzero square-zero generator anticommuting with it. The latter operator exchanges the two eigenspaces. They are orthogonal and nondegenerate, so its off-diagonal blocks are adjoint and of equal rank. Its rank is even and at most three, hence at most two. If $q=0$, every member of the pencil $\langle B,C\rangle$ has rank at most three. Constant rank three would give a rank-three quotient $E$ of the trivial bundle on $\mathbf P^1$ by the kernels of the pencil. The symmetric forms induce $E\simeq E^*\otimes\mathcal O(1)$, contradicting $2\deg E=3$. This proves the claim in all cases.

For the application, let $A$ be a nonsingular member and write $AB,AC$ for complementary quadratic coefficients. At a general smooth point of the common zero set the affine tangent $T$ is the $A$-orthogonal of $\langle u,Bu,Cu\rangle$. If the tangent test fails, the rank of $A$ on this four-dimensional tangent space is at most one, so $T^\perp\subset T$. Hence the quadratic forms represented by $B^2,(BC+CB)/2,C^2$ vanish on the complete intersection. Its quadratic ideal is exactly the three-dimensional span of the coefficients. Thus $J=\langle I,B,C\rangle$ is a Jordan algebra. By the claim its pencil of quadratic equations contains a nonzero equation of rank at most two. Over $\mathbf C$ such a quadratic factors into linear forms. The homogeneous ideal is prime and contains no nonzero linear form, a contradiction. The cubic or quartic nonvanishing on the complete intersection follows also from writing its equation on $\mathbf P^8$ as $\sum_{i=1}^3 q_i H_i$ if all the higher forms vanish. Their common zero locus then has dimension at least two and is singular on that linear section, contrary to finiteness of the Gauss map as above.

The same tangent test succeeds for an integral three-quadric complete intersection when $Q_t$ has generic rank six and the three quadratic matrices have no common radical vector. Choose a general rank-six member $A$ and let $R=\ker A$. If $R$ is not contained in the affine tangent $T$ at a general point, $T$ projects to a codimension-two subspace in the six-dimensional nondegenerate quotient $V/R$; its restriction has rank at least two. A failure would therefore force $R\subset T$ at a general point $u$. In particular every defining quadric $B$ has $B(R,u)=0$ generically. Nondegeneracy of the integral complete intersection implies $B(R,V)=0$ for all three, contrary to hypothesis.

\subsection{Du Val cubic surfaces in a quartic ninefold line calculation}
\label{sec:nov30-quartic-du-val}
Consider a polar-rank-two line \(L\) in a smooth quartic ninefold.
Suppose its equation on the common light space first has order
three and its cubic form on four general constant lights defines
a Du Val cubic surface at general points of the line. Then
\[
       \delta_L(X,H)\ge\frac{125}{17}.
\]
If a general such cubic surface has no singularity of type \(E_6\),
then \(\delta_L(X,H)\ge170/23>125/17\). Both statements include
the smooth case. The constants follow from the surface thresholds
computed in E.~Denisova, \emph{\(\delta\)-invariants of Du Val del
Pezzo surfaces of degree \(3\)}, arXiv:2311.14181, Main Theorem.
In hyperplane units, equal here to anticanonical units, these
thresholds are at least \(1/3\) in all cases and at least \(9/19\)
after excluding \(E_6\).

Write \(G\subset\mathbf P^3\) for this surface, \(S\) for the
expected hyperplane order on \(\mathbf P^3\), and let \(A=A(v)\)
and \(u=v(G)\). For every divisor with center on \(G\), the complete
series and the flag along \(G\) give
\[
     4S(v)\ \le\
     \min\left\{A,\ \frac{u}{3}+B(A-u)\right\},
     \qquad B=\frac{3}{\delta(G,\mathcal O_G(1))}.
\tag{DV1}
\]
Here the first estimate is \(\delta(\mathbf P^3,\mathcal O(1))=4\).
For the second, take in degree \(k\) a basis adapted to \(v\) and
to order along \(G\). After removal of \(G^j\), the graded layer
restricts as the complete series of \(\mathcal O_G(k-3j)\), by the
hypersurface exact sequence. The mean value of \(j/k\) tends to
\(1/12\), by summing
\(\dim H^0(\mathbf P^3,\mathcal O(k-3i))\) for \(i\ge1\).
The mean residual normalized degree therefore tends to \(3/4\).
On \(G\) the residual basis divisor with its coefficient multiplied
by four has asymptotic order at most \(B\) times discrepancy, by
the finite-level convergence of basis thresholds on \(G\).
The contribution from layers of bounded residual degree vanishes
after normalization. Inversion of adjunction for the normal
Cartier divisor \(G\), which has klt singularities, gives (DV1)
by taking coefficients strictly below the limiting ones and
passing to the limit. In particular \(0\le u\le A\).

For the universal Du Val estimate one may replace \(B\) by \(9\);
for the \(E_6\)-free estimate one may use \(19/3\).
In each case put
\[
        x_0=\frac{B-1}{B-1/3},\qquad D=\frac{1+2x_0}{1+x_0}.
\]
For \(x=u/A\) the two terms of (DV1) cross at \(x_0\).
The function
\(\min\{1,B-(B-1/3)x\}+(2-D)x\) increases up to \(x_0\)
and decreases thereafter, and its value there is \(D\).
Off \(G\) use the first estimate alone. This proves on
\(\mathbf P^3\) the projective inequality
\[
             4S(v)+(2-D)v(G)\ \le\ D A(v)
\]
with \(D=37/25\), or \(D=25/17\) in the \(E_6\)-free case.
Both numbers exceed \(1+1/3\). Apply the translated-branch
criterion of Section~\ref{sec:new-translated-branch} with
four lights of weight one and heavy weight three. Its radial
denominator is \(37/28\), and its angular denominator is at
most \((1+3D)/4\). The latter is respectively \(34/25\) or
\(23/17\) and dominates the radial denominator. Dividing the
ambient normalization ten by these denominators proves the claims.

In particular the first bound applies whenever the geometric generic
cubic is normal and has at least four essential variables. A general
projective three-plane section is a normal cubic surface by Bertini.
The triple points of the original cubic form constitute its vertex,
as vanishing of its Hessian at a vector says exactly that contraction
of the cubic polarization with the vector is zero. This vertex has
codimension the number of essential variables and misses a general
three-plane. The singular locus of a normal cubic has ambient
codimension at least three. A general three-plane avoids its strata
of codimension at least four and meets its codimension-three strata,
if any, at general points. At each double point on these strata the
quadratic part restricts nontrivially on a general transverse
three-plane, since its vanishing there imposes a proper closed
condition on the three-plane directions. Bertini transversality
on the smooth locus excludes further singularities. The surface
section has only double points. The classification of normal cubic
surfaces gives Du Val singularities in this case; the non-Du Val
normal cubic surfaces are cones over smooth plane cubics. See
J.~W.~Bruce and C.~T.~C.~Wall, \emph{On the classification of
cubic surfaces}, J. London Math. Soc. (2) \textbf{19} (1979),
245--256.

\paragraph{A cubic and a transverse quartic double contact.}
At order three, suppose there are two constant light directions on
which the cubic restricts to a generically squarefree binary form,
and two further constant directions belonging to the radical of
every cubic coefficient such that the binary quartic on these
last directions has a root of multiplicity exactly two.
Then \(\delta_L(X,H)\ge52/7>37/5\). In particular the quartic
contact condition is available if its restriction to a three-space
in the common radical is a smooth plane quartic. The tangent line
at a general point of this quartic has intersection multiplicity
exactly two, since its Gauss map is nonconstant and separable.

Give the two cubic lights weights four, the two radical lights
weights three, and the two heavy variables weight twelve. The
degree-twelve light face consists of the binary cubic and the
binary quartic in disjoint variables. Indeed every cubic
coefficient annihilates the radical, while a quartic monomial
involving a cubic light has weight at least thirteen.
The radial denominator is \(1+(3/4)/(1+14/12)=35/26\).
At each of the general poles make a linear change of the two
cubic directions to write their squarefree binary cubic as
\(x^3+y^3\), absorbing its nonzero coefficients into the
coordinates. Put the double root of the quartic at \(u=0\).
It is \(a u^2v^2+b u^3v+c u^4\) with \(a\ne0\).
Give \(u,v\) additional weights \(2,1\), give \(x,y\) additional
weight two, and give the graph coordinate weight six. The
secondary graph equation is \(z+x^3+y^3+a u^2v^2=0\).
Before applying this specialization, subtract a sufficiently
large multiple of the original weights \(4,4,3,3,12\) on
\(x,y,u,v,z\) from the additional weights. The whole original
graph is homogeneous of degree twelve for the original weights,
so this alters every term by the same weight and preserves its
initial graph. The resulting additional weights \(\nu_i\) are
all negative. For each fixed ordinary degree the diagonal change
\(x_i\mapsto t^{\nu_i}x_i\) preserves the complete space of
polynomials of at most that degree. A basis simultaneously
adapted to the tested valuation and the specializing filtration
has initial forms constituting a basis of that same space.
An original affine point \(p\) has coordinates
\((t^{-\nu_i}p_i)_i\) in this family and specializes to the
affine origin. Openness of the klt
locus at a coefficient strictly below the signed bound for this
initial graph thus proves the same upper signed denominator at
the original affine center. The disjoint monomial estimate of
Section~\ref{sec:nov28-disjoint-general} gives
\[
 D=1+\max\{(1+1/3+1/3+1/2)^{-1},1/4,2/5\}=19/13 .
\]
The combined pole denominator is again \(1/4+(3/4)D=35/26\).
Linear coordinate choices after specialization to an individual
pole only concern its ordinary affine graph estimate. This proves
the asserted line mark.

\begin{corollary}\label{cor:nov30-quartic-ninefold-du-val}
Let \(X\subset\mathbf P^{10}\) be a smooth quartic ninefold. At
each line of polar rank two assume at least one of the following.
\begin{enumerate}
\item The first light form has degree three, and its restriction
to four general constant lights at a general point of the line
defines a Du Val cubic surface. This holds in particular for a
geometrically normal cubic with at least four essential variables.
\item The line satisfies one of the quadratic tangent tests of
paragraphs~\ref{par:nov30-audit-polar-pencils}
and~\ref{par:nov30-audit-three-quadrics}.
\item The line satisfies the terminal-quartic criterion of
Section~\ref{sec:new-translated-branch}, or the cubic and
quartic double-contact criterion of
Section~\ref{sec:nov30-quartic-du-val}.
\end{enumerate}
Then \(X\) is uniformly K-stable and admits a K\"ahler--Einstein
metric.
\end{corollary}
\begin{proof}
For a line of polar rank at least three the threshold is at least
\(37/5\). The first hypothesis at a polar-rank-two line gives
\(125/17\), including a surface of type \(E_6\). The second gives
\(15/2\). The two alternatives in the third give \(15/2\) and
\(52/7\), respectively. Every line therefore has local threshold
at least \(125/17\). This exceeds the sufficient line coefficient
\(785/107\) by \(30/1819\). The strict line criterion proves
uniform K-stability, and the K\"ahler--Einstein conclusion follows
from the Yau--Tian--Donaldson theorem.
\end{proof}

\subsection{Combining the regular and positive systems at a common pole}
\label{sec:nov30-joint-poles}
In the two-polar line calculation suppose the first light form has
degree \(h\) on \(k\) retained light variables, all of weight
one, and let \(d\) be the hypersurface degree. Write
\(c=(d-1)/d\). For each of the \(d-1\) general simple poles take
\(P=\mathbf P(1^k,h)\), \(Z=(z=0)\), \(W_a=(z+G_a=0)\), where
\(G_a\) is the light form there after division by its nonzero
polar coefficient. Put
\[
 D_0(v)=(k+1)S_{\mathcal O_P(1)}(v)+(1-1/h)v(Z).
\]
If \(D_0(v)+c v(W_a)\le bA_P(v)\) at every pole, for all
divisorial valuations, the angular denominator is at most \(b\).
Thus it suffices that \(b\ge1,\ b\ge c+1/h,\ b\le1+c\), and
\[
        kS_{\mathcal O_{\mathbf P^{k-1}}(1)}(v)
                +(1+c-b)v(G_a)\le b A_{\mathbf P^{k-1}}(v)
\]
for every pole. The full denominator is the maximum of the
angular one and \(1+ch/(k+h)\).

Indeed fix a test valuation at the geometric generic angular
fiber, a section degree, and an angular degree. Filter by the
largest pole order and by consecutive residue evaluations at
the ordered roots exactly as in the two-polar calculation.
Leave the regular kernel (pole order zero) unrefined for now.
For this fixed filtration the sum of valuation values in each
graded piece of a doubly adapted basis is intrinsic: it is the
sum of jumps for the induced valuation filtration on that piece.
Denote the regular sum by \(R\) and the sum for the pieces at
root \(a\) by \(P_a\), aggregating angular degrees. For each
chosen root \(a\) further refine the regular kernel by base
vanishing order at that very root. An adapted basis of the
refined filtration remains adapted to the first filtration.
In particular the numbers \(R,P_a\) are the same for all
choices of the regular refinement.

The regular kernel has no denominator from an eliminated heavy
coordinate, so its base-order residue calculation applies at
the pole itself. Normalize by the appropriate base parameter
powers the sections in the regular family and the positive
pieces at this particular root. The product of their special
fibers, counting the positive sections with multiplicity
\(d-1\), is a divisor on the same weighted projective fiber.
Its two contributions have the complete-series and fixed-part
descriptions of Chapter~\ref{ch:fullrow}: in the normalization
conditional on either group they are bounded by \(D_0\) and
\(D_0+v(W_a)\). Both groups have identical asymptotic mass
\(1/d\), including the stated complete and fixed means. Thus
the bound \(D_0+c v(W_a)\le b A_P\) gives the same asymptotic
threshold denominator \(b\) for the generic product contributing
\(R+(d-1)P_a\), by openness of klt and properness of the
weighted fiber. Repeat this comparison at each pole with the
independently chosen regular refinement. The average of these
order sums over \(a\) is precisely \(R+\sum_a P_a\); the
uniform finite-level and bounded-degree errors vanish as before.
This proves the angular assertion.

To check the sufficient projective inequality write
\(w=v(W_a), z_0=v(Z), \lambda=1+c-b\). The two cyclic-cover
comparisons and translation give the required bound at all
centers outside \(Z\cap W_a\): if \(w=0\) use \(D_0\le A_P\);
if \(z_0=0\), use
\(D_0+c w\le A_P+(c-1+1/h)w\) and \(w\le A_P\).
At a center in the intersection, flag adjunction and translation
as in Section~\ref{sec:new-translated-branch} give
\[
 D_0+cw-bA_P\le-\lambda v(G_a)+
 \min\{cw-(b-1)z_0,\ (c-b+1/h)w+(1-1/h)z_0\}.
\]
For \(w\le z_0\) the first term in the minimum is at most
\(\lambda w\); for \(w\ge z_0\) the second is at most
\(\lambda z_0\). Now use
\(v(G_a)\ge\min(w,z_0)\).

\paragraph{Quadratic light forms of rank at least three.}
On a smooth quartic ninefold any line of polar rank two whose
quadratic light form has geometric generic rank at least three
satisfies \(\delta_L(X,H)\ge100/13\).
Retain three general constant lights and choose the simple poles
so every conic \(G_a\subset\mathbf P^2\) is smooth.
At every plane valuation write \(A=A(v),u=v(G_a)\).
The complete plane series gives \(3S_H\le A\). A flag along
the conic gives \(3S_H\le u/2+2(A-u)\) when its center is
on the conic: the fixed mean in the normalization by three
is \(1/2\), the residual mean degree is two, and
\(\delta(G_a,H)=\delta(\mathbf P^1,\mathcal O(2))=1\).
The two bounds cross at \(u/A=2/3\), whence
\[
                 3S_H+(9/20)u\le (13/10)A .
\]
Use \(b=13/10,\ c=3/4,\ h=2,\ k=3\).
Here \(1+c-b=9/20\), \(c+1/h=5/4\), and the radial
denominator is also \(13/10\). Restoring the ordinary
normalization ten proves the line estimate.

\subsection{Quartic ninefolds with a deficient quadratic polar on a line}
\label{sec:nov30-quartic-rank-family}
In the parameter space of smooth quartic ninefolds, all failures
of uniform K-stability lie in a closed subset of codimension at
least \(90\).

By the point-to-line criterion, failure requires a line with
threshold strictly below \(125/17\): indeed one may choose the
common target strictly between \(785/107\) and \(125/17\).
Its polar rank is two by the previously proved \(37/5\) bound
for higher polar ranks. At a fixed line its containment costs
five conditions and the required rank of its freely varying
\(4\)-by-\(9\) matrix of first normal cubics costs another
\((4-2)(9-2)=14\). The Grassmannian of lines in \(\mathbf P^{10}\)
has dimension eighteen, so these conditions already leave
codimension at least one after forgetting the line. Smoothness
ensures that the pencil has no base point. Fix these polars and
write \(V\) for the common seven-dimensional space of lights.
The quadratic, cubic, and quartic light coefficients are
independently arbitrary, since in the full quartic polynomial
they have separate degrees in the coordinates normal to the line.
The family of quadrics is an arbitrary
\[
 Q\in H^0(\mathbf P^1,\mathcal O(2))\otimes\operatorname{Sym}^2V^*,
 \qquad \dim H^0(\mathcal O(2))\otimes\operatorname{Sym}^2V^*=84 .
\]
By Section~\ref{sec:nov30-joint-poles} its generic rank is at
most two.

For generic rank \(r=1,2\), let \(E\subset V^*\otimes\mathcal O\)
be the saturation of the image of the symmetric map
\(V\otimes\mathcal O(-2)\to V^*\otimes\mathcal O\).
Its quotient is a vector bundle, and symmetry gives
\(Q\in H^0(\operatorname{Sym}^2 E(2))\), where the twist comes
after taking the symmetric square. This follows in local frames
of the summand since its rows and columns outside that summand
vanish generically and hence identically. The induced map
\(E^*\to E(2)\) is generically nonsingular. If
\(E=\bigoplus_i\mathcal O(-a_i)\), then \(a_i\ge0\) by its
embedding in the trivial bundle. The determinant section of
\((\det E)^2\otimes\mathcal O(2r)\) gives \(\sum a_i\le r\).
In particular the coefficients of \(Q\) use at most
\(r+\sum_i a_i\le 2r\) constant linear forms on \(V\), namely
the coefficients of a split matrix for the embedding of \(E\).
Hence the constant common radical \(K\) has dimension at least
\(7-2r\).

For each type an upper dimension bound for the family of
matrices is
\[
 7\sum_i(a_i+1)-h^0(\operatorname{End} E)
        +h^0(\operatorname{Sym}^2 E\otimes\mathcal O(2)).
\]
Indeed embeddings of the split bundle form a locally closed
subset of the displayed Hom space. Isomorphisms of the source
act freely, and each orbit has dimension
\(h^0(\operatorname{End} E)\). Once its image is chosen the
tensor belongs to the last displayed vector space.
The dimensions in these three terms and the resulting upper
bounds are
\[
 \begin{array}{c|rrrr}
 (a_i)&\dim\operatorname{Hom}&\dim\operatorname{Aut}&
          h^0(\operatorname{Sym}^2 E(2))&\text{upper bound}\\\hline
 (0,0)&14&4&9&19\\
 (0,1)&21&4&6&23\\
 (0,2)&28&5&4&27\\
 (1,1)&28&4&3&27\\
 (0)&7&1&3&9\\
 (1)&14&1&1&14 .
 \end{array}
\]
At generic rank two a cubic nonzero on \(K\) yields a constant
direction \(u\in K\) with generically nonzero cubic value.
Two further constant directions give a generically nondegenerate
binary quadratic. The quadratic-plus-cube calculation then
gives \(\delta_L\ge280/37\). Both cubic coefficients must
therefore vanish on \(K\); this costs at least
\(2\binom{3+2}{3}=20\) independent conditions. The resulting
global codimension is at least \(1+84-27+20=78\).
If \(\dim K\ge4\), the same independent cubic cost is at
least forty, which already gives \(98\). Suppose then that
\(\dim K=3\), and choose a four-dimensional complement
\(W\) in \(V\). If the quartic \(T\) vanishes on \(K\), that
costs fifteen independent conditions, bringing the bound to
\(93\). Otherwise choose six points \(u_i\in\mathbf P(K)\)
outside this quartic whose evaluations on ternary quadratics
are independent. This is possible because their evaluation
vectors on a nonempty open subset span the dual of the space
of quadratics.

At each point \(u_i\) the rows
\(dC_0(u_i)|_W,dC_1(u_i)|_W\) are an arbitrary ordered
pair of four-dimensional covectors, independently for all
six points: the coefficients of the terms of type \(K^2W\)
are the full tensor \((\operatorname{Sym}^2 K^*)\otimes W^*\)
in each cubic. Dependent rows cost at least three. For
independent rows let \(B_i\subset W\) be their common kernel.
If the general quadratic restricts nondegenerately to \(B_i\),
two constant directions there together with \(u_i\) give the
quadratic-plus-quartic estimate \(\delta_L\ge15/2\).
In fact all quadratic coefficients annihilate \(u_i\), both
cubic values and both first derivatives on \(B_i\) vanish
at \(u_i\), and \(T(u_i)\ne0\), as required in the
\(xy+s^4\) test. Deficiency therefore forces all general
quadrics to restrict degenerately to \(B_i\).

At generic rank two each individual general complex member
factors as a product of two independent linear forms. The
condition on \(B_i\in\operatorname{Gr}(2,W)\) is the Schubert
divisor that \(B_i\) meet the common annihilator of those
two factors. It is irreducible, by the incidence of a line
in this annihilator and a plane containing that line. Take
two members with distinct spans of factors; such members
exist since otherwise all quadrics would use only two constant
variables, contradicting \(\dim K=3\). The corresponding
Schubert divisors are distinct, since their decomposable
Pl\"ucker forms are nonproportional. Thus their intersection
has codimension at least two. Passing to independent ordered
row pairs is smooth with equidimensional fibers, so this
condition costs at least two per \(u_i\), including the
dependent-row alternative. Independence at the six points
gives \(78+12=90\).

At rank one the quadratic family alone gives
\(1+84-14=71\), and at rank zero it gives \(1+84=85\).

At rank one a nonzero reduced cubic on the common radical
gives the \(x^2+u^3+v^3\) estimate \(\delta_L\ge52/7\).
At a fixed five-space in this radical both endpoint cubics
must therefore be nonreduced or zero. Cubics of the form
\(l^2m\) in five variables have at most \(5+5-1=9\)
parameters in the space of dimension thirty-five. Thus this
branch has codimension at least \(71+2(35-9)=123\).
When \(Q=0\) choose a fixed four-space of lights. The two
restricted cubics are independent arbitrary cubic surfaces.
Whenever the generic surface is smooth, the translated-branch
cubic calculation gives \(\delta_L\ge280/37\). More generally
a normal nonconical surface gives the Du Val mark \(125/17\)
of Section~\ref{sec:nov30-quartic-du-val}. Reducible surfaces
and cones over plane cubics each have at most thirteen parameters
in the twenty-dimensional space (linear times quadratic, or
the vertex and a ternary form). An integral nonnormal cubic
surface has a singular line: the secant lines of any irreducible
curve of singular points lie on the cubic, and a nonlinear such
curve has secants filling its plane or all of \(\mathbf P^3\).
Having a singular line costs at least
\(4+2\cdot3-\dim\operatorname{Gr}(2,4)=6\).
The three exceptional cones are closed by projective incidence.
Both independent endpoint cubics must belong to their union,
so the rank-zero branch has codimension at least \(85+12=97\).
Combining the three strata gives the
assertion, since projection and Zariski closure do not increase
the incidence dimension.

\paragraph{A pole-specific mixed quadratic and quartic flag.}
Let a polar-rank-two line in a smooth quartic ninefold have
generic quadratic light rank two, common radical \(K\) in
its seven-space of lights, and choose a constant complement
\(W\). Assume both light cubics vanish on \(K\), and the
terminal quartic \(T\) is nonzero there. Write \(Q_t,C_t\)
for the quadratic and cubic light coefficients at parameter
\(t\). For a general \(t\), suppose there exists \(u\in K\)
with \(T(u)\ne0\) such that the covector \(dC_t(u)|_W\)
either lies outside the span of the factors of \(Q_t\), or
belongs to the line of one of those factors. Then
\(\delta_L(X,H)\ge15/2\).

Give all directions of \(W\) weight two, all of \(K\) weight
one, and both heavy variables weight four. The polynomial on
the light face is precisely
\[
 Q_t(w)+dC_t(u)(w)+T(u),\qquad (w,u)\in W\oplus K.
\]
In fact the old linear polars have weight four, the quadrics
annihilate \(K\), and each cubic vanishes on \(K\).
The original radial denominator is
\(1+3/(4+2\dim W+\dim K)\le4/3\).
Choose simple poles in the open set supporting the hypothesis.
The pole filtration for arbitrary positive light weights gives
one regular group with ordinary affine denominator one and three
positive groups represented by the complete affine polynomial
systems with one fixed copy of the pole graph in the normalized
bound. This uses the unchanged ordinary degree cut on the
monomials of the residue blocks of Chapter~\ref{ch:fullrow}.

Work at an individual pole and write \(\ell=dC_t(u)|_W\).
If it lies outside the row space of \(Q_t\), its kernel maps
surjectively to the active two-dimensional quotient of this
single quadric. Retain two directions there producing \(xy\)
and retain the direction \(u\) with variable \(s\). The
resulting graph polynomial is \(z+xy+a s^4\), \(a\ne0\),
after rescaling. If \(\ell\) is a multiple of one factor,
the restricted mixed term is \(bxs^2\) after naming that
factor \(x\). Give \(x,y,s,z\) secondary weights \(3,1,1,4\)
to remove this term while keeping \(xy,z,s^4\). In the first
case take weights \(2,2,1,4\). Give every omitted direction
an arbitrarily large secondary weight in either case, so every
term involving such a direction disappears.

For clarity the comparison of the ordinary affine series here
applies to centers at every affine point. The pole graph is
homogeneous for the original positive outer weights. Subtract
a sufficiently large multiple of those outer weights from
all secondary coordinate weights: this shifts all weights on
the graph by a common constant, so gives the same initial
graph, and makes the secondary coordinate weights negative.
Under the family \(f(t^{w_i}x_i)\) an affine point \(p\) of
the original fiber now moves to \((t^{-w_i}p_i)\), tending
to the origin. The subspace of all ordinary affine
polynomials of bounded degree is preserved by this diagonal
action. In each such space choose a basis compatible both
with the tested valuation and with this secondary filtration.
Its initial terms form a basis of that same complete bounded
degree space. After normalizing the graph and each basis
member by its least parameter power, any coefficient strictly
below the signed threshold for the special product applies to
the generic product at the chosen point by openness of klt.
The coordinate choices are made after taking residues at the
individual pole and need not extend along the original line.

Thus at each of the three positive groups the signed graph
bound for \(xy+s^4\) gives denominator \(13/9\). The extra
variables are free and do not change this bound: split their
characters, each with ordinary normalized mean one, leaving
complete polynomial spaces of the residual ordinary degrees
in the four graph coordinates, with their usual simplex means.
Their independent coordinate factors are transverse on the
product log resolution. The combined denominator is at most
\(1/4+(3/4)(13/9)=4/3\). This proves the line estimate.

\paragraph{Restriction of the cubic derivative as a bundle section.}
Consequently, for every deficient line of the preceding type
and every \(u\in K\), the section
\(\ell_t(u)=dC_t(u)|_W\) lies in \(H^0(E(1))\), with \(E\)
the saturated rank-two polar image. Otherwise noninclusion
in the generic row space is open in \(t\), and one may choose
\(u\) outside the terminal quartic preserving this property.
Saturation upgrades generic membership to bundle membership.
At fixed \(Q\), before imposing the condition \(\ell_t(u)\)
ranges freely in \(H^0(\mathcal O(1))\otimes W^*\), at a
collection of \(\binom{\dim K+1}{2}\) independently evaluating
points of \(\mathbf P(K)\) outside \(T\). Hence the independent
cost in the two cubic coefficients after their vanishing on
\(K\) is at least
\[
 \binom{\dim K+1}{2}
       \big(2\dim W-h^0(E(1))\big).
\]
In the two types \((1,1),(0,2)\) with \(\dim W=4\) this
costs thirty-six. In type \((0,2)\) every section of \(E(1)\)
is in the fixed linear factor, so the second alternative of
the pole-specific test rules out this type altogether when
the quartic is nonzero on \(K\).

\subsection{The exceptional locus in the smooth quartic ninefold family}
\label{sec:nov30-ninefold-final}
\begin{theorem}\label{thm:nov30-ninefold-codimension}
In the space \(\mathcal U_{9,4}\) of smooth quartic hypersurfaces
in \(\mathbf P^{10}\), the locus of hypersurfaces failing uniform
K-stability is relatively closed and has codimension at least
\(127\), with the empty locus allowed.
\end{theorem}
\begin{proof}
Closedness follows from openness as in the beginning of this
chapter. The line criterion in Section~\ref{sec:nov30-ninefold-results}
reduces the remaining assertion to the dimension of an incidence
of lines with threshold below \(125/17\).
The estimate before using the isotropic-factor alternative for
the last terminal quartic is \(114\). We give the dimension
checks in the notation of Section~\ref{sec:nov30-quartic-rank-family}
and the preceding pole-specific mixed test.

At quadratic light rank two write \(k=\dim K\). The vanishing
of both cubic coefficients on \(K\) is necessary and contributes
\(2\binom{k+2}{3}\) conditions for each fixed quadric family.
For \(k\ge5\) these conditions alone give the overall
bound \(1+84-27+70=128\). For \(k=4\) the terminal quartic
cannot vanish identically on \(K\): the span with the original
line would then be a projective five-plane contained in \(X\).
Indeed a smooth quartic ninefold contains no such plane since
its five normal cubic derivatives on that plane have a common
zero. The three-dimensional complement \(W\) implies
\(h^0(E(1))\le3\), since type \((0,0)\), the sole type with
four sections, has constant active dimension two. The mixed
bundle restriction thus costs at least \(10(6-3)=30\),
giving \(1+84-27+40+30=128\).

It remains to consider \(k=3\). Before any additional mixed
test the overall bound is \(78\). The four independent active
forms imply that the type is \((1,1)\) or \((0,2)\). In
the former the quadratic is \(L_tM_t\) for linear pencils
with four independent coefficients: its nonsingular symmetric
form on the split rank-two bundle is a constant split binary
quadratic after the degree-two twist. In the latter
\(Q_t=L M_t\) with fixed \(L\) and with three independent
coefficients of the degree-two \(M_t\) modulo \(L\). This
follows since the square of the degree-two column has zero
coefficient and the cross coefficient is a nonzero constant.

If the terminal \(T\) is nonzero on \(K\), the bundle condition
on the derivatives costs \(6(8-2)=36\), giving \(114\).
Moreover in type \((0,2)\) both sections of \(E(1)\)
belong to the fixed factor \(L\), so the isotropic-factor
alternative excludes this case entirely. In type \((1,1)\)
we have
\[
             dC_t(u)|_W=\alpha(u)L_t+\beta(u)M_t ,
\]
where \(\alpha,\beta\) are ternary quadratic forms on \(K\):
the derivative as a function of \(u\) is a homogeneous
quadratic with values in the constant two-space \(H^0(E(1))\).
At a zero of either quadratic outside \(T\) the covector is
a multiple of one isotropic factor for every general pole.
This again gives \(15/2\). Thus neither quadratic can be
identically zero, and the support of \(\{\alpha\beta=0\}\)
is contained in \(\{T|_K=0\}\). After the cubic jets have
been fixed, the independently arbitrary ternary quartic is
divisible by the squarefree part of \(\alpha\beta\).
Write \(s\) for its degree. If \(s=4\), the quartics form a line
in a fifteen-dimensional vector space, costing fourteen.
For \(s=3\), divisibility costs twelve and there is at least one
condition on the conics, since a pair of distinct smooth conics
has squarefree product. For \(s=2\) divisibility costs nine.
In this case an irreducible common support makes the conics
proportional, costing five among ordered pairs. If the support
consists of two distinct lines, the line pair varies in dimension
four and the scalar factors in dimension two. The multiplicities
are discrete, so this costs at least six. If \(s=1\), both quadratics are
scalar multiples of the square of the same linear form. Within
the twelve parameters of the pair \((\alpha,\beta)\), this costs
at least eight: its line varies in dimension two and its two
scalars in dimension two. Divisibility of \(T|_K\) by this fixed
line costs a further five. The two ternary quadratic coefficients
are independently arbitrary after the linear bundle restriction,
since \(H^0(E(1))\) is freely spanned by \(L_t,M_t\) in this
splitting. Thus every case gives at least \(114+13=127\).

Suppose instead that \(T|_K=0\), bringing the starting
bound to \(78+15=93\). For type \((1,1)\) let \(U_1,U_2\)
be the common kernels of the two coefficient vectors of
\(L_t,M_t\), respectively. They each have dimension five,
meet in \(K\), and span the seven-space of lights. All
quadratic coefficients vanish on each of them. For both
cubic coefficients the restrictions to \(U_1,U_2\) are
independently arbitrary in the systems of cubics vanishing
on \(K\). Indeed choose a basis along \(K\) and then separate
the complements in \(U_1,U_2\); all residual cubic monomials
for the two restrictions involve disjoint sets of complementary
variables. If the general member of either restricted cubic
pencil is integral, normal, and has at least four essential
variables, its general projective three-plane section gives
the Du Val line bound \(125/17\) by
Section~\ref{sec:nov30-quartic-du-val}. Thus each of these
four independent cubic parameters must belong to the closed
exceptional cubic cone.

Here is its dimension bound in the system of cubics in
\(U_i\) through the fixed plane \(S=\mathbf P(K)\).
The system has affine dimension \(35-10=25\). In a reducible
cubic the plane \(S\) is contained in the linear or in the
quadratic factor. These alternatives have at most
\(1+15=16\) and \(4+(15-6)=13\) parameters. Cubics with at
most three essential variables have at most \(6+10=16\)
parameters even without the incidence. Every integral
nonnormal cubic in \(\mathbf P^m\) has a singular linear
\(\mathbf P^{m-2}\). Indeed normality fails in codimension
one by the hypersurface Cohen--Macaulay property. An irreducible
component \(Y\) of this dimension in the singular locus has
all its secants in the cubic. If a nonlinear \(Y\) spans a
hyperplane, its secants fill this hyperplane; if it spans the
ambient space, its secants fill that space. For the latter
assertion, if the secants have dimension at most \(\dim Y+1\),
the joins \(\operatorname{Join}(y,Y)\) for general \(y\in Y\)
already fill the secant variety. Indeed a general \(y\) is not
a vertex of the nonlinear \(Y\). The secant variety would
then have every such \(y\) as a vertex, making it linear and
equal to the span of \(Y\). This also proves the hyperplane
case. Both alternatives contradict the irreducibility of the
cubic. A linear span of dimension \(m-2\) gives the asserted
linear singular component.

In \(\mathbf P^4\) write \(D\) for this singular plane. If
\(D=S\), second-order vanishing costs twelve further conditions
in the system already containing \(S\). If \(D\cap S\) is a
point, cubics in \(I_D^2\) impose twenty-two ambient conditions;
their restrictions to \(S\) give the seven-dimensional system
of ternary cubics vanishing twice at this point. This is seen
by lifting the two generating linear forms of the restricted
ideal. After subtracting the ten given plane conditions and
allowing \(D\) to vary in dimension six, the cost is at least
\(22+7-10-6=13\). If \(D\cap S\) is a line, the corresponding
restricted dimension is three (multiples of the square of a
line), and the planes meeting \(S\) in a line vary in dimension
four. The cost is at least \(22+3-10-4=11\). Hence the
union of the exceptional cubic loci costs at least nine in
this system. These cones are closed by the multiplication,
vertex-subspace, and singular-plane projective incidences.
The four independent restrictions give \(93+4\cdot9=129\).

In type \((0,2)\) use all of the six-space \(U=\ker L\).
Each cubic restricted to it is independently arbitrary among
the \(56-10=46\) cubic forms containing the fixed plane.
A linear factor containing this plane gives at most \(2+21=23\)
parameters, and a residual quadric containing it gives at most
\(5+(21-6)=20\). An integral nonnormal cubic has a singular
three-plane in \(\mathbf P^5\), hence its whole ambient cone
has dimension at most \(56-(20+2\cdot10)+8=24\).
At most three essential variables allow only \(9+10=19\)
parameters. Failure therefore costs at least twenty-two in
each restricted cubic, giving \(93+44=137\).

At rank one, in the notation above, a pencil of cubic forms in
five variables all of whose members are zero or nonreduced has
dimension at most fourteen as an ordered pair. Indeed if its
coefficients are independent, the generic cubic has a unique
multiple linear factor over the algebraic closure of
\(\mathbf C(t)\). Uniqueness descends its projective class to
\(\mathbf C(t)\). This gives a morphism from the smooth base
\(\mathbf P^1\) to the space of linear factors after extending
across finitely many points. The residual factor defines a second
morphism. The multiplication map \((L,M)\mapsto L^2M\) pulls
back the hyperplane bundle to \(\mathcal O(2,1)\). Since the
pencil has degree one, the repeated factor is constant. A
proportional nonzero pencil already has a constant repeated
factor. The common factor varies in dimension four and the two
remaining linear forms give ten parameters. The count at rank
one therefore improves to \(71+70-14=127\). At rank zero,
choose six constant light directions; the two cubic parameters
on them are independent without any imposed plane. In their
affine space of dimension fifty-six, reducibility costs
\(56-(6+21-1)=30\); the nonnormal and inessential exceptions
cost at least \(56-24\) and \(56-19\), as just proved.
Both coefficients in a deficient pencil belong to this
closed union, giving \(85+60=145\). The smallest of all
the resulting lower bounds is \(127\). For the passage to the
parameter space take the union of the incidence conditions on
coefficients used in these counts. It is constructible: one
stratifies by the ranks, splitting types and radical dimensions,
chooses frames on Grassmannian charts, and imposes the indicated
factorizations and divisibilities by multiplication incidences
of forms of fixed degrees. All degrees and splitting types in
the counts range over finite sets. The fiber dimension estimate
therefore applies to this union. Taking its closure and
projecting from the product with the Grassmannian preserves the
codimension lower bound already computed after allowing the line
to vary. The deficient lines lie in this incidence union.
\end{proof}

\paragraph{A smooth realization of the remaining quadratic jets.}
The balanced type with terminal quartic equal to the product of
the two derivative conics occurs on smooth members.
Take homogeneous coordinates
\((s,t,w_1,\ldots,w_4,u_1,u_2,u_3,z_0,z_1)\) on
\(\mathbf P^{10}\). Let \(A,B\) be smooth plane conics in the
\(u\)-variables meeting transversely, and put
\(L=sw_1+tw_2,\ M=sw_3+tw_4\). Start with the quartic
\[
       z_0s^3+z_1t^3+(L+B)(M+A).
\]
Perturb it by general quartic forms in the \(w\)-variables,
quartic forms in the \(z\)-variables, and terms
\(\sum_i w_i h_i(u)\) with homogeneous cubics \(h_i\).
All these perturbations have strictly larger weight than four
for the weights \(0,2,1,4\) on the successive variable groups.
The resulting linear system has no base point where \(w\)
or \(z\) is nonzero, because it contains every pure fourth
power of those coordinates. Bertini and normalization of
the initial coefficient to one give a smooth general member
there. On the subspace \(w=z=0\) its equation is \(AB=0\).
When \(u\ne0\) and exactly one of \(A,B\) vanishes, an
\(u\)-derivative is nonzero. At the four crossing directions
the \(w\)-derivative vector is \((h_i(u))\); this is nonzero
at all four for a general choice, independently of \(s,t\).
The points with \(u=0\) form the given line, where the two
normal derivatives \(s^3,t^3\) have no common zero.
This proves smoothness.

Along this line the quadratic light family is \(LM\), with
four independent coefficient forms. On its common radical
the derivative rows assemble to \(A L+B M\), and the terminal
quartic is \(AB\). Thus no point of the radical with nonzero
terminal value gives an individual-factor covector. The
example verifies that smoothness alone does not guarantee
the hypothesis of this mixed-pole criterion. It does not
determine the K-stability of the constructed members.

\paragraph{The marked hyperplane inequality at its endpoint.}
Let \(X\) be a smooth quartic ninefold and \(Y\) a hyperplane
section through \(p\), smooth there. For any divisor \(v\) centered
at \(p\), the weak local bound
\(\delta_p(Y,H|_Y)\ge63/10\) suffices to conclude
\[
                  A_X(v)-7S_H(v)\ \ge\ \frac3{10}v(Y)>0 .
\tag{HE1}
\]
Indeed at a coefficient \(\lambda<7\), adapt the basis to \(v\)
and order along \(Y\). The basis coefficient on \(Y\) tends to
\(1/10\), and the complete restricted degrees have normalized
mean \(9/10\). The effective residual boundary at coefficient
\(\lambda\) is log canonical on \(Y\) near \(p\) for sufficiently
large levels, by the local threshold hypothesis and the uniform
complete-series comparison. Inversion of adjunction after increasing
the coefficient on \(Y\) to one gives
\[
    A_X(v)-\lambda S_H(v)\ \ge\ (1-\lambda/10)v(Y).
\]
Take \(\lambda\uparrow7\). Thus equality at the eightfold
threshold permits strictness at each ninefold point-centered
divisor.

In particular one can exclude \(p\) if the nine-dimensional
tangent parameter space admits a seven-dimensional subspace
\(U\) on which \(q_2,q_3\) vanish identically and whose section
\((q_4|_U=0)\subset\mathbf P(U)\) is a smooth quartic fivefold.
Let \(B\simeq\mathbf P^7\) be the span of \(p\) and
\(\mathbf P(U)\) in the projective tangent hyperplane. Its
section with \(X\) is the cone with vertex \(p\) over this
smooth fivefold. Take a general flag of two hyperplanes through
\(B\). Off \(B\) the two successive sections are smooth by
Bertini for the system of linear forms vanishing on \(B\).
At points of \(X\cap B\) other than \(p\), the differential
of the equation already restricts nontrivially to \(B\).
At \(p\) the final \(\mathbf P^8\) through \(B\) is chosen
not to lie in \(T_pX\). Thus both intermediate hypersurfaces
\(Z^7\subset Y^8\subset X\) are smooth. The point \(p\) is
generalized Eckardt on \(Z\) with tangent section \(X\cap B\).
Its quartic fivefold base is K-stable by Chapter~\ref{ch:fivefold},
so the Abban--Zhuang formula
\cite[Theorem~4.18]{ar2003.13788} gives
\(\delta_p(Z,H|_Z)=28/5\). Hyperplane adjunction from \(Z\)
gives \(\delta_p(Y,H|_Y)\ge(9/8)(28/5)=63/10\), as required.
For example the test applies when \(q_2=xy\) and
\(q_3=xA+\ell B_2\) for quadratic forms \(A,B_2\) and a
linear form \(\ell\) independent of \(x\), provided
\(q_4|_{x=\ell=0}\) is smooth. These hypotheses allow
\(q_2,q_3\) with no common factor.

More generally let \(3\le d\le n+1\) and \(r=n+2-d\ge2\).
Suppose the tangent space at a smooth point of \(X_d^n\)
contains a vector subspace of dimension \(b\ge r\) on which
all tangent forms of degrees \(2,\ldots,d-1\) vanish and the
last form cuts a smooth degree-\(d\) hypersurface, K-semistable
if it is Fano. Choose the smooth \(b\)-dimensional linear
section as above and use the generalized Eckardt formula on it.
Iterated adjunction gives
\[
       \delta_p(X,H)\ge \frac{(n+1)b}{b+d-1}.
\]
For \(b>r\) the excess over \(r\) is
\((d-1)(b-r)/(b+d-1)>0\). For \(b=r\), necessarily \(b<n\);
the penultimate threshold in the flag gives, by the same
argument as (HE1),
\[
     A_X(v)-r S_H(v)\ \ge\
          \left(1-\frac r{n+1}\right)v(Y)>0
\]
for every divisor centered at \(p\), with \(Y\) the first
hyperplane in the chosen flag.

One instance admitting the section endpoint estimate follows from
the common linear factor calculation in
Section~\ref{sec:q8-tangent-tests}. On the quartic eightfold
suppose the quadratic and cubic tangent forms have a common
linear factor with quadratic contact, or there is cubic contact
whose cubic has a linear factor. Suppose additionally that the
terminal quartic \(G\subset\mathbf P^6\) of that calculation is
smooth. Then \(\delta_p(Y,H|_Y)\ge63/10\).

The normalized fixed and moving terms of its angular series on
\(\mathbf P^6\) have coefficients \(1/12\) on \(G\) and
\(7/9\) on the hyperplane complete series. At coefficient
\(63/10\) they become \(21/40\) and \(49/10\), respectively,
with \(49/10=7-4(21/40)\). The smooth quartic fivefold \(G\) is
K-stable by Chapter~\ref{ch:fivefold}.
Zhang--Zhou's log-pair interpolation
\cite[Theorem~1.7]{ar2003.06839} gives K-semistability of
\((\mathbf P^6,(21/40)G)\), since \(21/40<1-(7/4-1)/6\).
This is exactly the required angular discrepancy bound.
The opposite vertex comparison is
\((63/10)(10/9)=7\), and the isolated coordinate has coefficient
\((63/10)/9<1\). In the quadratic-contact extraction with
linear-factor multiplicity \(s=1,2\) put \(b=4-s\); the radial
discrepancy minus the tested coefficient times expected order is
\((7+b)-(7/10)(10+b)=3b/10>0\). For cubic contact with
multiplicity \(s=1,2,3\) the difference is
\((8s+1)-(7/10)(11s+1)=3(s+1)/10>0\).
The weighted flag bounds hold for the pointed threshold by the
local adjunction theorem. Equivalently the estimates at its
exceptional divisor control all valuations whose centers contain
the marked point; add a smaller rational coefficient first and
pass to the endpoint.

There is a limitation on extending this estimate through the
ordinary eightfold tangent models. In the regular rank-one model
\(x^2,f_3(a,b),g_4(a,b)\) of the endpoint calculation, the
coordinate divisor \(x=0\) on the ordinary exceptional space
has scaled basis order exactly two. Indeed the exponent of
the symbol for \(x^2\) has normalized simplex mean one, while
the parity remainder in \(x\) is bounded. The angular threshold
of this specialized series is consequently at most \(9/2\).
A general marked hyperplane preserves a quadratic form of rank
one. This example occurs for smooth original equations: prescribe
on the tangent hyperplane \(T^2 x^2+T f_3(y)+g_4(x,y)\), with
general \(f_3,g_4\), and extend by \(z\) times a general cubic
nonzero at \(p\). The system on the tangent hyperplane is
basepoint free away from \(p\) because all \(g_4\) are available,
and the extension is smooth at \(p\) by its \(z\)-derivative.
Bertini on and off the tangent hyperplane proves smoothness
for a general extension. This bounds the strength of the
specialized angular comparison; the actual pointed threshold
of these smooth sections is not computed by this argument.

\paragraph{Orders on the original eightfold at a rank-one point.}
The angular obstruction of the specialized model admits a direct
comparison with a full family of divisorial orders on the original
eightfold. Assume \(q_2=x^2\) and \(x\nmid q_3\). In the seven
remaining tangent parameters \(y\) the forms
\(f=q_3|_{x=0},\ g=q_4|_{x=0}\) are coprime, by the integrality
of the section \(Tf+g=0\) as in the original rank-one test.
For rational \(a\ge1\), use the monomial valuation in the actual
smooth tangent parameters with \(v_a(x)=a,\ v_a(y_i)=1\).
Then \(A_Y(v_a)=7+a\), and
\[
  9 S_H(v_a)=
   \begin{cases}
        47/6+2a,&1\le a\le 3/2,\\
        28/3+a,&a\ge3/2 .
   \end{cases}
\tag{HE2}
\]
In particular the minimum of \(A_Y(v_a)/S_H(v_a)\) over this
family is \(459/65>63/10\), attained at \(a=3/2\).

Write the equation as \(T^3w+T^2Q_2+TQ_3+Q_4=0\), with transverse
coordinate \(w\) and \(Q_i|_{w=0}=q_i\). For \(1\le a<3/2\)
the three hypertangent sections \(w,\ U=Tw+Q_2,\ W=TU+Q_3\)
have section degrees \(1,2,3\) and orders \(2a,3,4\).
Their initial relations are the three of ordinary degrees \(2,3,4\),
with respective weights \(2a,3,4\). For \(a>1\) their reductions
modulo \(T\) are \(x^2,f,g\); at \(a=1\) these are
\(x^2,q_3,q_4\), also a regular sequence by coprimality modulo
\(x\). In verifying these initial relations one uses
\(TW=-Q_4\) and \(TU=W-Q_3\) before using \(Tw=U-Q_2\).
Every term involving \(w\) in the homogeneous \(Q_i\) has order
strictly larger than the least order of \(q_i\).
Regularity modulo \(T\) makes \(T\) a nonzerodivisor on the
presented complete intersection. Localizing at \(T\) realizes
it as the polynomial algebra in the actual tangent initial
parameters and \(T^{\pm1}\). Its ordinary Hilbert series is
that of a quartic eightfold, so these are all graded relations.

For \(a\ge3/2\) retain only \(w,W\) with orders \(3,4\).
The relation \(T^2w+TQ_2+Q_3=W\) has initial
\(T^2w+f=0\), with the additional term \(Tx^2\) when \(a=3/2\).
The other relation is \(TW+g=0\). Regularity modulo \(T\)
follows from the coprimality of \(f,g\), and the same localization
and ordinary Hilbert series prove this presentation.
In a weighted complete intersection the scaled expected section
order is the sum of weight divided by section degree over
generators, minus the same sum over relations. This follows
by logarithmic differentiation of the bigraded Hilbert series.
In the first presentation the sum is
\[
 7+a+2a+\tfrac32+\tfrac43
       -\left(\tfrac{2a}{2}+\tfrac33+\tfrac44\right)
     =47/6+2a .
\]
In the second it is \(7+a+3+4/3-3/3-4/4=28/3+a\).
The resulting ratio decreases to \(a=3/2\) and increases
thereafter, as asserted.

\subsection{A quartic derivative at a base point of the cubic derivatives}
\label{sec:dec1-q-final-ninefold}
\begin{theorem}\label{thm:dec1-q-final-ninefold}
The failure locus for uniform K-stability of smooth quartic ninefolds
has codimension at least $148$ in $\mathcal U_{9,4}$. If the quadratic
light family of a polar-rank-two line has the form $L_tM_t$ for two
linear pencils with four independent coefficients, its line threshold
is at least $125/17$ on every smooth member.
\end{theorem}
\begin{proof}
Use the reduction and notation of Theorem~\ref{thm:nov30-ninefold-codimension}.
A failure requires a line of polar rank two and mark below $125/17$.
Its seven-dimensional light space has quadratic, cubic and quartic
coefficients $Q_t,C_t,T$, independently variable within their specified
degrees. Write $K$ for the constant quadratic radical. In all the
weighted faces below use the pole-specific specialization from
Section~\ref{sec:nov30-quartic-rank-family}, suppressing omitted
directions by large secondary weights and subtracting a sufficiently
large multiple of the outer homogeneous weights. This justifies the
signed estimate on the entire affine chart at each general pole.
If the light weights sum to $s$ and the heavy weight is $h$, the radial
denominator is $1+3h/(4(s+h))$. A signed denominator $D$ from
Section~\ref{sec:nov28-disjoint-general}, formula (DG), gives pole
denominator $(1+3D)/4$. Divide ten by the larger denominator to obtain
the line mark. The pole coordinate is denoted by $z$.

\paragraph{Two tests in quadratic rank two.}
Both cubic coefficients vanish on $K$ by the quadratic-plus-cube
bound of Section~\ref{sec:nov30-quartic-rank-family}. First fix
$u\in K$ such that $dC_t(u)$ lies outside the row space of $Q_t$
for general $t$. In weights two on a complement $W$ to $K$, one on
$K$, and four on the heavy directions, the light face is
$Q_t(w)+dC_t(v)(w)+T(v)$. At each general pole retain the direction
$u$, two vectors in $\ker dC_t(u)$ carrying the nondegenerate binary
quadric $xy$, and a vector in $\ker Q_t$ not killed by $dC_t(u)$.
The graph becomes $z+xy+a b^2+c b^4$, with $b$ the $u$-coordinate
and $a$ the third complementary coordinate. Secondary weights
$(1,2,1,1,3)$ on $(x,y,a,b,z)$ give $z+xy+a b^2$. Formula (DG)
gives $D=7/5$ and the pole denominator $13/10$. There are at least
nine units of light weight, so the radial denominator is smaller.
The mark is at least $100/13>125/17$. Thus on a deficient line the
derivative belongs to the polar row space for all $u$ and general
$t$, even if $T|_K=0$. Saturation gives the $H^0(E(1))$ row-bundle
condition of Section~\ref{sec:nov30-quartic-rank-family}.

For the second test suppose $u\ne0$ is a common zero of the two full
cubic derivatives on $K$ and $T(u)=0$. The quartic covector
$b=dT(u)$ on the light space is nonzero. Otherwise the plane joining
the original line and $[u]$ lies on $X$ and all light derivatives
vanish there, as do derivatives along the original line coordinates.
The two remaining normal cubics on this plane have a common zero,
contradicting smoothness. Euler's identity gives $b(u)=0$. If $b$
lies outside the generic polar row space, give $u$ weight one, all
six remaining lights weight three, and the heavy variables weight
six. The face is $Q_t(w)+b(w)s^3$: the possible weight-five term
has coefficient $dC_t(u)=0$. At each pole retain two directions in
$\ker b$ carrying $xy$, a third in $\ker Q_t$ with nonzero $b$-value,
and $s$. Formula (DG) for $z+xy+a s^3$ gives $D=10/7$. The combined
denominator is $37/28$, larger than the radial bound for $s=19,h=6$.
The mark is at least $280/37$. Consequently $b$ belongs to the generic
row space whenever the line is deficient.

Write $k=\dim K$. The large-radical quadratic-pencil estimate excludes
$k\ge5$. When $k=4$, choose the constant active dual three-space
in dimension twelve. Replacing $7$ by $3$ in the split-bundle count
of Section~\ref{sec:nov30-quartic-rank-family} bounds all four types
by eleven dimensions each. Thus the old $128$-bound in the proof of
Theorem~\ref{thm:nov30-ninefold-codimension} increases by
$27-(12+11)=4$. Independently the terminal $T|_K$ must have a linear
or quadratic factor, by the final assertion of the large-radical
quadratic-pencil estimate. It is nonzero by smoothness, and this
factorization costs at least $35-\max\{4+20-1,10+10-1\}=12$ in
the quaternary quartics. Thus this stratum gives at least $144$.
For fixed polar data the row-bundle section of the two cubic
derivatives is specified by $h=h^0(E(1))\le3$ independently arbitrary
quaternary quadratics. Their common projective zero set is nonempty.
At each such zero $[u]$, the old isotropic-factor test forces $T(u)=0$,
and the second test above puts $dT(u)$ in the intersection of all
general polar row spaces. This constant intersection has dimension at
most one since otherwise every quadratic would use the same two
variables. The three normal derivatives of $T$ outside $K$ are
arbitrary quaternary cubics, independently of $T|_K$. These conditions
cost at least two per point, independently at two distinct points.
If a single point supports all common zeros, we must have $h=3$ and
the three quadrics cut out a complete intersection of length eight.
This already costs two further conditions on their coefficients.
Indeed at a fixed point their three values vanish, offsetting the
three parameters of the point. Rank at most one of the tangent
gradient matrix costs at least four. Rank two costs one. After taking
two equations as smooth parameters the residual equation has zero
gradient and its quadratic along the common tangent line must vanish
to allow length eight. Otherwise the residual intersection has length
two. This is an independent condition on its quadratic jet with the
linear jets fixed. In either case the combined additional cost is at
least four, yielding $148$.

For $k=3$ in type $(1,1)$, $Q_t=L_tM_t$ has the four independent
coefficients in the assertion. The first test gives
$dC_t(u)=\alpha(u)L_t+\beta(u)M_t$ for fixed ternary quadratics.
Take a simultaneous projective zero of $\alpha,\beta$. It can be
chosen with $T(u)=0$: this is automatic if $T|_K=0$, and otherwise
the isotropic-factor test in the proof of
Theorem~\ref{thm:nov30-ninefold-codimension} gives
$V(\alpha\beta)\subset V(T|_K)$. The second test gives a nonzero
constant $b$ in the generic span of $L_t,M_t$. Such a covector cannot
exist: the two complementary two-dimensional spans of the coefficient
pairs each contain a nonconstant pencil of lines with zero intersection.
This proves the asserted line bound in balanced type.

In type $(0,2)$ write $Q_t=L M_t$, where the three coefficients of
$M_t$ modulo $L$ are independent. The proof of Theorem~\ref{thm:nov30-ninefold-codimension} excludes $T|_K\ne0$.
For $T|_K=0$ its bound $137$ comes from restrictions of
the cubics on the six-space $\ker L$, in addition to the fixed
quadratic and $T|_K=0$ conditions. It improves first to $139$.
Indeed a cubic on that six-space containing the fixed plane
$\mathbf P(K)$ varies in dimension $46$. The exceptional reducible
and inessential cones in the previous calculation have dimension at
most $23$. Integral nonnormal cubics belong to the same incidence of cubics singular
along a three-plane of dimension $24$ before imposing the fixed
plane. This incidence is an irreducible vector bundle over the
Grassmannian and containment of the fixed plane is a proper condition:
the cube of a linear form vanishing on a chosen three-plane need
not vanish on the fixed plane. Thus the exceptional dimension in
this system is at most $23$ for each of the two independent cubics,
giving $93+2(46-23)=139$.

The first test additionally requires
$dC_t(u)=(t_0a_0(u)+t_1a_1(u))L$ since all row-bundle sections
belong to the fixed factor. For each admissible cubic pair fix any
common projective zero $[u]$ of its two ternary quadratic coefficients.
The second test forces $dT(u)\in\langle L\rangle$, since the three
independent coefficients of $M_t$ modulo $L$ have zero common row
intersection in that quotient. After $T|_K=0$, its four first normal
derivatives are freely prescribed ternary cubics; evaluation at the
fixed nonzero $u$ is surjective onto the four-dimensional conormal
space. This costs at least three further independent conditions on
the terminal coefficients, yielding $142$. The choice of $u$ does
not increase the dimension: for fixed cubic data the terminal condition
holds at \emph{every} common zero, so the fiber is contained in the
linear conditions obtained from any one of them. In fact this branch
costs at least $145$. Varying the coefficients of the two separate
monomials $L a_i(u)$ changes $a_0,a_1$ arbitrarily without affecting
either cubic restriction to $\ker L$ or the row condition. If $a_0,a_1$
both vanish identically this costs twelve coefficients. Otherwise if
they have only a single projective common zero, both are nonzero with
no common component, and B\'ezout's theorem gives intersection
multiplicity four at this zero. This exceptional conic-pair locus
has codimension at least three: if one conic is smooth, the fixed
conic has degree-two restriction series of its second conic equal
to the complete binary quartics, with a one-dimensional scalar-kernel;
being a fourth power costs three. If both conics are singular, they
already each cost one, and having all four line intersections equal
is a proper closed condition on the product of the two irreducible
singular-conic loci (away from pairs with a common component), giving
at least one further condition. In all other cases there are two
distinct projective common zeros. The two evaluations on ternary
cubics are independent, so their terminal conditions add six. Thus
either six terminal conditions hold, or at least three conic-pair
conditions and three terminal conditions hold beyond $139$.

The unconditional row restriction improves this to at least $150$
without requiring a second common zero. On $\ker L$, already in
the $46$-dimensional cubic system, it forces the two cubics to
belong to $I_S^2(3)$ for $S=\mathbf P(K)$. In this subspace the
exceptional cubic cone has dimension at most $19$. Reducible cubics
with linear factor through $S$ have dimension at most $2+(21-6)=17$,
and those with linear factor transverse to $S$ at most $5+6=11$,
since the quadratic factor then belongs to $I_S^2(2)$. Cubics in at
most three variables have dimension at most $9+10=19$. For a
singular three-plane $D$ containing $S$, its two parameters and the
sixteen cubics in $I_D^2(3)$ give dimension at most eighteen. If $D$
meets $S$ in a line, write their ideals $(a,b,c)$ and $(a,v)$.
The intersection of their squared ideals in degree three is spanned
by $a^2$ times the six coordinates and $abv,acv$. Even allowing
all eight parameters of $D$ gives at most sixteen. If they meet in
a point the two ideals have disjoint coordinate generators, with
zero intersection of their squares in degree three. Thus the two
independent constrained cubics give at least $93+2(46-19)=147$,
to which the three independent terminal derivative conditions at a
single fixed zero add three more.

\paragraph{A double factor of the radical cubic.}
For a constant quadratic row at rank one, $k=6$. The two radical
cubics are nonreduced with a fixed repeated factor by the pencil
argument of Theorem~\ref{thm:nov30-ninefold-codimension}, varying in
dimension at most $5+12=17$ inside $2\binom83=112$ parameters.
Together with the previous quadratic cost of at least $71$ this exceeds
$141$. Otherwise $Q_t=A_t^2$ uses two independent coefficient forms
on a complement $W$ to the five-space $K$. The quadratic count gives
$71$. On $K$ the cubics either both vanish, costing another
$2\binom73=70$, or restrict to $B^2 M_0,B^2 M_1$. In the former
case give $W$ weight two, $K$ weight one, and the heavy variables
weight four. At a general pole the affine graph has light polynomial
$x^2+xP(u)+aD(u)+T(u)$, where $x=A_t$ and $D$ is proportional to
$A_t\wedge dC_t|_K$. If $D$ has rank at least two, choose two vectors
$e,f\in K$ with $D(e)=0$ and $D(e,f)\ne0$. With $u=s e+v f$,
assign secondary weights $(6,3,4,5,12)$ to $(x,a,s,v,z)$.
The only terms retained are $z+x^2+c asv$, $c\ne0$. Formula (DG)
gives $D_{\rm signed}=7/5$ and line mark at least $100/13$ after
the pole and radial comparisons. Thus the associated quadratics at
two general poles must both have rank at most one. The quadratic
coefficients at the two poles are independent, and rank at most one
costs $15-5=10$ each on $K$. This branch therefore costs $161$.
In the remaining case the linear factor
$[B]$ varies in dimension four and the $M_i$ in ten, giving base
codimension $71+70-14=127$.

Put $U=\ker(B|_K)$ and $m_t=M_t|_U$. First suppose $m_t$ is not
identically zero at a general $t$. Give weights $(4,3,2)$ to $W$,
a one-dimensional complement of $U$ in $K$, and $U$ respectively,
and assign heavy weight eight. With $v$ the complementary coordinate
in $K$, the face is
\[
            A_t(w)^2+v^2m_t(u)+dC_t(u)(w)+T(u).
\]
For a general $t$, if some $e\in\ker m_t$ has $T(e)\ne0$ and
$dC_t(e)|_W$ is zero or independent of $A_t$, take a $W$-direction
in its kernel carrying $x=A_t$, and retain $e$ and one $U$-direction
with nonzero $m_t$. For coordinates $(x,v,q,s,z)$ where $s$ is the
$e$-coordinate, use secondary weights $(2,1,2,1,4)$. The quartic
and derivative terms involving $q$ have weight above four and the
coefficient of $xs^2$ is zero. The initial $z+x^2+v^2q+s^4$ gives
$D=13/9$ and combined denominator $4/3$, exceeding the radial
denominator $1+6/(8+8+3+8)$. The line mark is at least $15/2$.
Thus for each general $t$ either $T|_{\ker m_t}=0$ or the quadratic
form $F_t=A_t\wedge dC_t|_U$ vanishes on $\ker m_t$. The alternatives
are closed where $m_t\ne0$. Since the parameter line is irreducible,
one of them holds for every general $t$.

If $m_0,m_1$ are independent, the union of their pencil's hyperplanes
is dense in $U$, so the first alternative costs all $35$ coefficients
of $T|_U$. In the second, the bidegree $(2,2)$ form $F(t,u)$ is
divisible by the irreducible bilinear form $m(t,u)$. Indeed the
bilinear form has rank two and its zero set over the general $t$
is dense in its irreducible zero hypersurface. The coefficients of
$F$ are arbitrary: contracting independent $A_0,A_1$ with the two
independent derivative covectors gives a surjection onto binary
quadratics, separately for each quadratic monomial on $U$. The
ambient space thus has dimension $3\cdot10=30$, and the multiples
of the fixed bilinear form have dimension $2\cdot4=8$. This costs
$22$. All these coefficients are independent of the cubic restrictions
to $K$. Consequently this case gives at least $127+22=149$.

If the nonzero pair $m_0,m_1$ is proportional, its two-by-four
matrix has rank one, adding three conditions to the base count.
The hyperplane $\ker m_t$ is constant of dimension three. The
vanishing-$T$ alternative costs fifteen further coefficients.
The derivative alternative costs eighteen: the restriction of the
arbitrary bidegree $(2,2)$ form $F(t,u)$ to the fixed three-space
vanishes identically, which gives $3\binom42$ independent conditions.
In the vanishing-$T$ alternative restrict the same outer face at
each pole to $W$ and the fixed hyperplane $H=\ker m_t\subset U$,
discarding the complementary variables to these spaces in the pole
graph. On $H$ the resulting graph has light polynomial
$x^2+xP(u)+aD(u)$, with $D=F_t|_H$ up to a nonzero scalar. The
rank-two-conic secondary test $(6,3,4,5,12)$ above forces $D$ to
have rank at most one. This costs $6-3=3$ more at each of two
independent evaluations. Hence this alternative gives at least
$127+3+15+6=151$, and the other gives $127+3+18=148$.

Finally, if $m_0=m_1=0$, the nonzero pencil restricts to a pure
cube of $B$ on $K$, of parameter dimension at most $4+2=6$.
Give $W$ weight six, $y=B$ weight four, $U$ weight three, and the
heavy variables weight twelve. If $dC_t(u)|_W$ is independent of
$A_t$ for a general pole and some $u\in U$, retain that $u$ and
coordinates $x=A_t,a=dC_t(u)|_W$ on $W$. The light polynomial on
these directions is $x^2+c y^3+a s^2+b s^4$, $c\ne0$. Secondary
weights $(3,2,2,2,6)$ on $(x,y,a,s,z)$ remove the last term.
Formula (DG) gives $D=10/7$ and pole denominator $37/28$, above
the radial denominator for $s=28,h=12$. Hence deficiency forces
$dC_t(u)|_W\in\langle A_t\rangle$ identically on $U$. Its two-pole
evaluation costs twenty independent coefficients in addition to
$71+70-6=135$.

At rank zero use the full seven-space of lights. Its affine space of
cubics has dimension $84$. Reducible cubics have dimension at most
$7+28-1=34$, and cubics in at most three variables at most
$12+10=22$. A nonnormal integral cubic has a singular four-plane
by the earlier secant argument. Vanishing twice on such a four-plane
imposes $35+2\cdot15=65$ conditions and the four-plane varies in
dimension ten. The closed exceptional cubic union thus has dimension
at most $34$ and costs $2(84-34)=100$ in the two independent cubic
coefficients, in addition to the old $85$. The minimum is $148$.
All estimates
are on the fixed-line incidence with the earlier subtraction for the
Grassmannian of lines already included. The same finite rank and
factorization stratifications and fiber-dimension argument as in
Theorem~\ref{thm:nov30-ninefold-codimension} give the claimed
codimension after taking the closure and projecting.
\end{proof}

\paragraph{Exclusion of a four-dimensional quadratic radical.}
On every smooth quartic ninefold a polar-rank-two line whose generic
quadratic light form has rank two and whose constant radical $K$ has
dimension four satisfies $\delta_C(X,H)\ge125/17$.

Suppose otherwise, and keep the notation of
Theorem~\ref{thm:dec1-q-final-ninefold}. The terminal $T|_K$ is
nonzero, and the cubic derivative lies in the generic quadratic row
space. For varying $u\in K$ the latter gives a quadratic map to the
space $H^0(E(1))$, of dimension at most three. Let $B\subset\mathbf
P(K)$ be its base locus, the common zeros of the coefficient quadratics.
The isotropic-factor test forces $T|_B=0$, and the transverse quartic
derivative argument in the cited theorem puts the nonzero covector
$dT(u)$, for every $[u]\in B$, in the constant intersection of the
generic quadratic row spaces. This intersection has dimension at most
one. If it is zero then $B$ is empty, impossible for three quadratics
on $\mathbf P^3$. If it is a line, the restriction of $dT$ to $B$ is
a nowhere-zero section of this fixed line tensored with
$\mathcal O_B(3)$. This excludes any integral projective curve in $B$
by the positive degree of its restriction. Hence $B$ is finite and
the three coefficient quadratics must be linearly independent.

The only possible splitting with three independent sections and
three constant active variables is $E=\mathcal O\oplus\mathcal O(-1)$:
the $(0,0)$ splitting uses at most two constant variables, and the
other splittings have $h^0(E(1))\le2$. Write $L$ for the constant
generator and $M_t$ for the degree-one column of the row bundle.
The derivative is
\[
                 (t_0\alpha_0+t_1\alpha_1)L+\beta M_t,
\]
with independent quaternary quadratics $\alpha_0,\alpha_1,\beta$.
At a general pole project this covector onto the quotient by either
of the two isotropic factor lines of $Q_t$. The result, after choosing
a scalar coordinate on that quotient, is
$G_t=\lambda(t)(t_0\alpha_0+t_1\alpha_1)+\mu(t)\beta$. For at least
one factor $\lambda(t)\ne0$, since the two factor lines are distinct.
For every general complex $t$ choose such a factor. By the isotropic-
factor test $V(G_t)\subset V(T|_K)$: otherwise the pole-specific
condition would hold at a general pole and give the required mark.
More explicitly, the existence of a $u$ with $T(u)\ne0$ and a
derivative in an isotropic factor line is a constructible condition
in the base parameter (the two factors may be parametrized on a
finite cover). If it held on infinitely many base parameters it
would hold generically. Thus, outside a finite set, the containment
holds for both choices.

There are only finitely many projective quadratics with support
contained in the fixed quartic: each is a product of its finitely
many prime factors with total degree two. In contrast the ratio of
the coefficients of $\alpha_0,\alpha_1$ in $G_t$ is $t_0:t_1$.
Their independence from each other and from $\beta$ gives infinitely
many projective quadratics, a contradiction.

\paragraph{The remaining unbalanced rank-two quadratic pencil.}
\begin{proposition}\label{prop:dec1-q-all-quadratic-rank-two}
On a smooth quartic ninefold every polar-rank-two line whose quadratic
light family has generic rank two satisfies $\delta_C(X,H)\ge125/17$.
\end{proposition}
\begin{proof}
The only remaining case in Theorem~\ref{thm:dec1-q-final-ninefold}
and the preceding four-dimensional-radical argument has radical $K$
of dimension three and quadratic form $Q_t=L M_t$, of bundle type
$(0,2)$. On a deficient line the cubic and quartic vanish on $K$,
and $dC_t(u)=a_t(u)L$, with $a_t=t_0a_0+t_1a_1$ for ternary
quadratics $a_i$. Choose $x=L$, coordinates $w$ on a three-dimensional
complement of $K$ in $\ker L$, and coordinates $u$ on $K$. Denote
by $m_t$ the restriction of $M_t$ to the $w$-space. Its three binary
quadratic coefficients give a basis of the dual $w$-space. Write
$b(u)=dT(u)|_w$, a triple of ternary cubics. At a common zero of
$a_0,a_1$ the second test in the cited theorem forces the entire
light covector $dT(u)$ to be nonzero and in the constant row
intersection $\langle L\rangle$. Thus the common zeros cannot contain
a curve, since the remaining $L$-coefficient of $dT$ is a section
of $\mathcal O(3)$ and would have a zero on that curve. In particular
$a_0,a_1$ are nonzero and relatively prime.

Give $x$ weight three, $w$ weight two, $u$ weight one and the heavy
variables weight five. The light face at a pole is
\[
    xm_t(w)+x a_t(u)+P_t(u,w,w)+b(u)(w),
\]
where $P_t$ is linear in $u$ and quadratic in $w$. Suppose that at
a general pole one can choose a nonzero $e\in K$ with $a_t(e)=0$
and $b(e)$ independent of $m_t$. Keep only $u=se$ in $K$ and
choose two $w$-directions on which $y=m_t$ and $v=b(e)$ are
independent. The face graph becomes
$z+xy+sP_t(e,y,v)+vs^3$. Give $x,y,v$ secondary weight three,
$s$ weight one and $z$ weight six. This removes the middle term,
of weight seven. The resulting graph $z+xy+vs^3$ gives denominator
$37/28$ after the pole comparison already used in
Theorem~\ref{thm:dec1-q-final-ninefold}. The original light weight
sum is twelve and its radial denominator is $1+15/68<37/28$.
Thus a deficient line requires $b(e)\in\langle m_t\rangle$ for
every $e\in V(a_t)$ at a general pole.

The rational map $\mathbf P(K)\dashrightarrow\mathbf P^1$ given by
$u\mapsto[-a_1(u):a_0(u)]$ is dominant since the conics are coprime.
Substitution in the preceding necessary condition forces $b(u)$ to
be proportional over the function field to
$(a_1^2,-a_0a_1,a_0^2)$ in the coefficient basis of $m_t$. This
triple is primitive in the polynomial ring, so the proportionality
scalar is itself a polynomial by unique factorization. The entries
of $b$ have degree three, while the primitive triple has degree four.
Hence $b=0$. The join of the original line and $\mathbf P(K)$ is a
projective four-plane contained in $X$. Its three $w$-direction
derivatives vanish there, successively by the radical property of
$K$, the cubic derivative condition, and $b=0$. Its remaining normal
derivatives are those in the $x$ direction and the two heavy
directions. Three cubics on this four-plane have a common zero,
contradicting smoothness.
\end{proof}

\begin{corollary}\label{cor:dec1-q-ninefold-156}
In the parameter space of smooth quartic ninefolds the locus
failing uniform K-stability has codimension at least $156$.
\end{corollary}
\begin{proof}
Proposition~\ref{prop:dec1-q-all-quadratic-rank-two} eliminates the
quadratic-rank-two cases. Only the residual rank-one pencil counts
in Theorem~\ref{thm:dec1-q-final-ninefold} need strengthening. In
its notation $U=\ker B$ is four-dimensional, the two-dimensional
$W$ carries the moving quadratic factor $A_t$, $m_t$ is a pencil
on $U$, and $\gamma_t(u)=dC_t(u)|_W$ is a covector-valued ternary
or quaternary quadratic on the space where it is restricted. The
passing test for a general $t$ requires only some $u\in H_t=\ker m_t$
with $T(u)\ne0$ and $\gamma_t(u)$ either zero or outside the line
$\langle A_t\rangle$. Thus if the derivative alternative holds on
$H_t$, its scalar coefficient must be nonzero wherever $T(u)\ne0$.

First suppose $m_0,m_1$ are independent. The derivative alternative
in the cited theorem already costs $149$. Substitute
$[t_0:t_1]=[-m_1(u):m_0(u)]$. The two components of $A_t$ now
form a coprime pair of linear forms after a constant change of
coordinates in $W$. The two components of $\gamma_t(u)$ are cubic
and proportional to this pair, hence equal an ordinary quadratic
$P(u)$ times that pair. On $V(P)$ away from the base line of the
pencil the vanishing test forces $T=0$ whenever $t$ is general.
For the numerical estimate it suffices instead to fix one general
$t$ after fixing the cubic data: either the scalar quadratic
$a_t$ in $\gamma_t|_{H_t}=a_t A_t$ is zero or $T|_{H_t}$ contains
its squarefree part. This gives at least five new independent
quartic conditions, for a total of $154$. To get $156$, use two
distinct general parameters, hence two distinct planes
$\mathbf P(H_i)$ in $\mathbf P(U)$. If the support of the two
scalar quadratics in their respective planes is not the common line
for at least one index, choose that line first in the ordering of
the two conditions so that the second plane has this property.
After the first five conditions, quartics divisible by the defining
linear equation of the first plane are still freely available in
the kernel. On the second plane these give that common line times
all ternary cubics. Divisibility by a component other than the
common line therefore costs at least four further conditions there.
If the only support on both planes is the common line for every
general choice, then $P$ has zero locus supported on a union of
fibers of the pencil: any irreducible component dominating the base
would meet a general fiber also away from the common line. Hence
$P$ lies in $\mathrm{Sym}^2\langle m_0,m_1\rangle$, or is zero.
Among the ten freely variable coefficients of $P$ this costs seven:
vary $\gamma_t$ by $Q(u)A_t$ for any quaternary quadratic $Q$
without affecting its wedge condition. This alone gives $156$.
The alternative $T|_U=0$ already costs $162$.

If $m_0,m_1$ are proportional and nonzero then $H_t=H$ is a fixed
three-space. In the derivative alternative, linearity of both
$\gamma_t$ and $A_t$ in $t$ gives $\gamma_t|_H=a(u) A_t$ for a
fixed ternary quadratic, with all six coefficients independently
arbitrary after imposing the wedge condition. Beyond the earlier
$148$, if $a$ has rank at least two then its squarefree equation
has degree two and the quartic divisibility costs nine. If $a$ has
rank one this costs three on $a$ and five on $T|_H$. If $a=0$ the
quartic vanishes on $H$. Thus this alternative gives at least $156$.
In the other alternative $T|_H=0$ the old graph test gives rank at
most one for the whole binary-quadratic family $A_t\wedge\gamma_t|_H$.
This family is arbitrary in eighteen coefficients before imposing the
rank restriction. Its saturated rank-one row bundle is
$\mathcal O(-e)$ with $e=0$ or $1$, hence the family dimension is at
most $\max_e\{3(e+1)-1+(3-2e)\}=6$. Its cost is therefore at least
$127+3+15+12=157$.

For a nonzero pure cube on $K$ the previously established identity
$\gamma_t(u)\in\langle A_t\rangle$ holds on the whole $U$. For
each of the ten quadratic monomials the corresponding linear family
in a fixed two-space is a scalar multiple of $A_t$, costing three
independently. This gives at least $135+30=165$. The zero-pair,
constant-row, and rank-zero cases already have bounds $161$, $166$,
and $185$ respectively. These prove the assertion.
\end{proof}

\subsection{Transverse quartic derivatives on the rank-one radical}
\begin{theorem}\label{thm:dec1-q-transverse-ninefold}
The complement of the uniformly K-stable locus in the smooth
parameter space of quartic ninefolds has codimension at least $171$.
\end{theorem}
\begin{proof}
We refine the rank-one calculation of
Theorem~\ref{thm:dec1-q-final-ninefold}. In its moving-row case
write $H_t=\ker(m_t)\subset U$. Suppose first that $m_t\ne0$
for general $t$. At each such pole the outer weighted light polynomial is
\[
        A_t(w)^2+v^2m_t(u)+\gamma_t(u)(w)+T(u).
\]
Here $\dim W=2$, $\dim U=4$, $\gamma_t(u)$ is quadratic in $u$
and linear in $t$, and $T|_U$ is an arbitrary quaternary quartic.
Use the secondary affine specializations and the radial and pole
comparisons of the cited theorem. In each specialization below the
radial denominator is $11/9$. A signed graph denominator of at most
$10/7$ gives pole denominator at most $37/28$, and therefore mark
at least $280/37>125/17$.

For a deficient line we first claim that
$\gamma_t(e)\in\langle A_t\rangle$ for every $e\in H_t$ and
general $t$. The cited theorem gives the claim where $T(e)\ne0$.
If $T(e)=0$ and the two covectors are independent, use
$x=A_t$ and $y=\gamma_t(e)$ as coordinates on $W$.
Retain $u=se+p e'$ with $m_t(e')=1$. Give $(x,y,v,p,s,z)$ the
secondary weights $(50,54,30,40,23,100)$, where $z$ is the
affine graph coordinate. The weight-$100$ polynomial is
\[
                   z+x^2+v^2p+ys^2 .
\]
Indeed the $s^4$ term vanishes, the first possible mixed quartic
has weight $3\cdot23+40>100$, and the mixed derivatives have
weight at least $50+23+40$. Formula (DG) gives $7/5$. This proves
the claim. Thus on $H_t$ we can write $\gamma_t=a_t A_t$ for a
quadratic form $a_t$. The earlier test also gives
$V(a_t)\subset V(T|_{H_t})$.

In fact $T|_{H_t}$ is singular at every nonzero $e\in V(a_t)$.
Otherwise take $u=se+b e''+p e'$ where $e''\in H_t$ has nonzero
quartic derivative at $e$. Retain one $W$-direction on which
$A_t=x$. Give $(x,v,p,s,b,z)$ weights $(50,32,36,24,28,100)$.
The $xs^2$ and $s^4$ terms vanish. All remaining mixed cubic
terms have weight at least $50+24+28>100$. The quartic term $s^3b$
has weight $100$ and all its other terms have higher weight.
After rescaling the resulting graph is
$z+x^2+v^2p+s^3b$. Formula (DG) gives $10/7$.
Consequently, when $a_t\ne0$, the square of its reduced equation
divides $T|_{H_t}$. Indeed at a general smooth point of any
irreducible component, a polynomial containing that component
to order one has nonzero differential.

When $a_t=q^2$ for a nonzero linear form $q$, we furthermore have
$q^4\mid T|_{H_t}$. Suppose its order along $q=0$ is exactly
two. Take a nonzero $e$ on that line where the residual quadratic
is nonzero, and take $e''$ with $q(e'')\ne0$. Use the preceding
coordinates with weights $(50,32,36,24,26,100)$. The derivative
restricted to $H_t$ now contributes a multiple of $xb^2$ of
weight $102$. Terms in that derivative containing $p$ have weight
at least $50+24+36>100$. In the quartic the initial term is a
nonzero multiple of $s^2b^2$, and terms containing $p$ have
weight at least $3\cdot24+36>100$. The graph specializes to
$z+x^2+v^2p+s^2b^2$, with signed denominator $7/5$ by (DG).
If the order is exactly three, take $e$ off the residual line and
use $(50,32,36,22,26,100)$. The term $xb^2$ still has weight
$102$, every quartic term containing $p$ has weight at least
$36+3\cdot22>100$, and $sb^3$ has weight $100$. This gives
$z+x^2+v^2p+sb^3$ with signed denominator $10/7$.

There is a further cubic condition at every $e\in V(a_t)$.
Write $A_t\wedge\gamma_t=m_t\ell_t$ for a linear form $\ell_t$
on $U$, fixing scalar coordinates in the determinant of $W^*$.
Then $\ell_t(e)=0$. To see this, use $u=se+p e'$ as before.
If $\ell_t(e)\ne0$, the coefficient of $sp$ in $\gamma_t$ is a
covector $y$ independent of $x=A_t$, by polarizing the wedge.
With $v$ omitted, give $(x,y,s,p,z)$ weights $(50,41,24,35,100)$.
The $s^4$ and $\gamma_t(e)s^2$ coefficients vanish,
the smallest mixed quartic weight is $107$, and the two possible
remaining cubic terms have weights at least $41+70>100$. The
graph $z+x^2+ysp$ has signed denominator $7/5$ by (DG).

For proportional nonzero $m_0,m_1$, the spaces $H_t=H$ coincide
and $a_t=a$ is a fixed ternary quadratic. The base and wedge
conditions cost $148$, with the six coefficients of $a$ free,
as in Corollary~\ref{cor:dec1-q-ninefold-156}. If $a$ has rank at
least two, its equation is reduced and the restriction of $T$ is a
scalar multiple of $a^2$. This costs fourteen among the fifteen
ternary quartic coefficients. Choose a fixed equation $m$ of $H$.
The three coefficients in the parameter of $(A_t\wedge\gamma_t)/m$
restrict to three arbitrary linear forms on $H$, independently of
$a$. Indeed the off-$H$ quadratic coefficients of each of the two
cubic derivatives vary freely, and wedging their binary linear
pencil with $A_t$ surjects onto binary quadratics for each fixed
monomial. All three resulting linear forms vanish on $V(a)$, hence
are zero for reduced $a$. This costs nine. If $a=q^2$, its rank
costs three, divisibility by $q^4$ costs fourteen, and the three
linear forms belong to $\langle q\rangle$, costing six.
If $a=0$, both the quartic and the three linear forms vanish.
The minimum total is $171$.

For independent $m_0,m_1$, the wedge condition costs $149$
including the base conditions. Put $t=[-m_1(u):m_0(u)]$ away
from their base line in $\mathbf P(U)$. Substitution in
$\gamma_t$ gives a polynomial vector equal to $P(u)A_t$ for a
quaternary quadratic $P$, since the two linear entries of the
substituted $A_t$ are coprime. All ten coefficients of $P$ are
free after imposing the wedge condition, by varying $\gamma_t$
by $R(u)A_t$. If $P$ has an irreducible component dominating the
parameter line, $T$ vanishes on a dense open subset of that component
and contains it at least twice. Indeed a general pencil plane meets
the component transversely away from the base line by generic
smoothness in characteristic zero. An order-one quartic along
the component would be smooth upon restriction there, contrary to
the transverse derivative condition. Multiples of the square of a
fixed positive-degree polynomial in the $35$ quaternary quartics
have dimension at most $10$, giving at least $149+25$. If $P=0$,
vanishing on all general
fibers gives $T|_U=0$. Otherwise every prime factor defines a
pencil hyperplane, so $P\in\mathrm{Sym}^2\langle m_0,m_1\rangle$,
of dimension three. On a general fiber its restriction is a
nonzero square up to scalar. The scalar quadratic $a_t$ there is
the restriction of $P$, by homogeneity in the parameter pair.
The fourth power of its reduced equation divides the ternary quartic
on this plane. As this holds for general planes through the base
line, $T|_U$ belongs to $(m_0,m_1)^4$. Indeed use $m_0,m_1$ as
the first two coordinates and substitute $(m_0,m_1)=(\alpha b,\beta b)$.
Each coefficient of $b^j$ for $j<4$ vanishes identically in the
other coordinates for general $\alpha,\beta$, hence identically
in all of the variables. Only the five binary quartic coefficients
survive. The total cost in this case is at least $149+7+30$.

For the moving row with both cubics zero on the five-space $K$,
the proof of Theorem~\ref{thm:dec1-q-final-ninefold} forces
$A_t\wedge\gamma_t$ on $K$ to have rank at most one for general $t$.
Its unrestricted coefficients are the $45$ entries of a binary
quadratic family of quinary quadratics. For a nonzero such family
of generic rank one, the saturated row line over the parameter line
is $\mathcal O(-e)$ in the trivial bundle of rank five. The
symmetric matrix is $g(t)b(t)b(t)^{\mathsf T}$ with $b(t)$ a
primitive five-vector of forms of degree $e$ and $g$ a binary form
of degree $2-2e$. Thus $e=0$ or $1$, and the number of parameters
is at most $5(e+1)-1+(3-2e)\le10$. The zero family satisfies the
same bound. Therefore the cost is at least $71+70+45-10=176$.
For the nonzero pure cube, the proof of Theorem~\ref{thm:dec1-q-final-ninefold}
gives the outer
face $x^2+c y^3+\gamma_t(u)(w)+T(u)$ with
$\gamma_t|_U=a(u)A_t$, $a$ a fixed quaternary quadratic.
Suppose $a(e)=0$ and $T(e)\ne0$. Keeping the $x,y$ directions and
$u=se$ gives $z+x^2+c y^3+T(e)s^4$. Its signed denominator in
(DG) is $37/25$. The pole denominator is $34/25$ and the radial
denominator for light sum $28$ and heavy weight $12$ is smaller.
This gives mark $125/17$ exactly. Thus $T|_U$ vanishes on $V(a)$
on a deficient line, adding at least fifteen independent conditions
to the earlier $165$. If $a=0$, all $35$ coefficients vanish.
For a constant nonzero row the quadratic family itself is
$g(t)L^2$, where $[L]$ has six parameters and $g$ is an arbitrary
binary quadratic. The earlier rank-one cost $1+84-14$ improves
to $1+84-9=76$ on this stratum. The cubic restrictions on the
six-dimensional radical cost $112-17=95$ as before, totaling
$171$. The zero quadratic case costs at least $185$ by
Corollary~\ref{cor:dec1-q-ninefold-156}. Proposition~\ref{prop:dec1-q-all-quadratic-rank-two}
excludes quadratic rank two. The same incidence projection as in
Theorem~\ref{thm:dec1-q-final-ninefold} gives the conclusion.
\end{proof}

\begin{corollary}\label{cor:dec1-q-reduced-cubic-174}
The locus of smooth quartic ninefolds failing uniform K-stability
has codimension at least $174$ in the smooth coefficient space.
\end{corollary}
\begin{proof}
In the moving-row case of Theorem~\ref{thm:dec1-q-transverse-ninefold},
assume $m_t\ne0$ for general $t$. In fact
$A_t\wedge\gamma_t=0$ on all of $U$. If this fails generically,
take $e\in U$ with $m_t(e)\ne0$ and $\gamma_t(e)$ independent of $A_t$
at a general pole. Keep only $u=pe$ in $U$, put
$x=A_t$, $y=\gamma_t(e)$ on $W$, and retain $v$. The resulting graph
has light polynomial $x^2+c v^2p+yp^2+c'p^4$ with $c\ne0$.
The weights $(3,2,2,2,6)$ on $(x,v,p,y,z)$ remove the last term.
The ternary cubic $p(cv^2+yp)$ is reduced. For instance the
binary restriction $y=p$ is $p(cv^2+p^2)$, which has three
distinct roots. A binary linear change and rescaling identify it
with the sum of two cubes. Formula (DG) for this sum and the
separate square has denominator $19/13$. The pole denominator
$35/26$ exceeds the outer radial denominator $11/9$ and gives
mark $52/7>125/17$. This contradicts deficiency.

The full wedge condition costs thirty, three for each quaternary
quadratic monomial. For proportional nonzero $m_0,m_1$, the base
cost is $127+3=130$. The transverse quartic conditions of the
preceding theorem add fourteen whenever $a|_H\ne0$, since the
quartic restriction is a multiple of $a^2$ for reduced $a$, or of
$q^4$ for $a=q^2$. If $a|_H=0$, the full ternary quartic vanishes.
This gives $130+30+14=174$. For independent $m_0,m_1$, the earlier
counts already give at least $174$ in the horizontal case, $184$
when $P=0$, and $186$ in the remaining case.

For a constant quadratic row the earlier base cost is $171$,
including the cubic restrictions $B^2 M_t$ on the six-dimensional
radical $K$. If both cubic coefficients vanish or are multiples of
$B^3$, their parameter dimension is at most $5+2=7$ in place of
$17$, raising this to $181$. Otherwise $U=\ker B\subset K$ has
dimension five, and $H=\ker(M_t|_U)$ is four-dimensional for a
general pole. Use weights $(4,3,2)$ on the sole quadratic direction,
the $B$-direction, and $U$, and weight eight on the heavy variables.
The light polynomial takes the form
$c x^2+v^2m_t(u)+a_t(u)x+T(u)$ with $c\ne0$.
Its radial denominator is $1+6/(4+3+10+8)=31/25$.
If $a_t(e)=0$ and $T(e)\ne0$ for some $e\in H$, restrict to
$u=se+p e'$ with $m_t(e')=1$. Secondary weights $(2,1,2,1,4)$
on $(x,v,p,s,z)$ give $z+c x^2+v^2p+T(e)s^4$.
Formula (DG) gives $13/9$ and pole denominator $4/3$,
hence mark at least $15/2$. Thus deficiency forces
$V(a_t|_H)\subset V(T|_H)$ at every general pole.
For fixed cubic data this costs at least fifteen coefficients of
the quaternary quartic, since even the multiples of a fixed linear
factor form a subspace of dimension only twenty. If $a_t|_H=0$,
all $35$ coefficients vanish. This proves the required improvement
also for the constant row. The other quadratic and cubic strata
were already bounded by at least $176$ or excluded in
Theorem~\ref{thm:dec1-q-transverse-ninefold}.
\end{proof}

\begin{theorem}\label{thm:dec1-q-mixed-normal-181}
The locus failing uniform K-stability in the smooth quartic ninefold
parameter space has codimension at least $181$.
\end{theorem}
\begin{proof}
First consider the moving quadratic row with proportional nonzero
$m_0,m_1$. Write $m$ for a fixed nonzero member of their span,
$H=\ker m\subset U$, and $\gamma_t|_U=aA_t$ as in
Corollary~\ref{cor:dec1-q-reduced-cubic-174}. Fix $p\in U$ with
$m(p)=1$ and $v_0\in K$ with $B(v_0)=1$. Write $C_t(-,-,-)$
for the symmetric trilinear form of the cubic, ignoring nonzero
scalar polarization factors. On $H$ form the ternary cubic
$J(e)=dT(e)(p)$ and the parameter-dependent ternary linear form
\[
       L_t(e)=A_t\wedge C_t(-,e,v_0)|_W .
\]
For every nonzero $e\in V(a|_H)$ on a deficient line, either
$J(e)=0$ or $L_t(e)=0$ identically in $t$.

Indeed $T(e)=0$ and, if $J(e)\ne0$, the constant vector
\[
        k=v_0-\frac{dT(e)(v_0)}{J(e)}p
\]
satisfies $dT(e)(k)=0$ and $B(k)=1$.
The wedge with $A_t$ of $C_t(-,e,p)|_W$ vanishes by the full
identity $\gamma_t|_U=aA_t$. The wedge with $C_t(-,e,k)|_W$
therefore equals $L_t(e)$. Suppose it is nonzero for general $t$.
Keep $W$ with weight three, the direction $k$ with weight two
and $e$ with weight one, with heavy weight six. Remaining light
directions can be given weight greater than six.
The restrictions of $Q_t$ to $K$ and its mixed terms with $K$ vanish.
In the cubic, $C_t(-,e,e)|_W=0$, while $B^2M_t$ on $K$ vanishes
on $e,e,k$ and $e,k,k$ since $B(e)=M_t(e)=0$.
In the quartic the terms $s^4,hs^3$ vanish, where $s,h$ denote
the coordinates along $e,k$. Thus the only light terms of weight
six are the square on $W$, its separate mixed covector times $hs$,
a multiple of $h^3$, a multiple of $h^2s^2$, and a covector on
$W$ times $s^3$. At each general pole use $x=A_t$ and take the
mixed covector as $y$. Secondary weights $(50,40,34,26,100)$
on $(x,y,h,s,z)$ retain only $z+x^2+yhs$. The signed denominator
is $7/5$ and the pole denominator $13/10$. The radial denominator
is at most $1+(3/4)6/(6+3+3+2+1)=13/10$. This gives mark
$100/13>125/17$ and proves the assertion.

The coefficients of $J$ vary freely with $T|_H$ fixed. Likewise
the coefficients of the three ternary linear forms in $L_t$ are
free with all cubic data on $K$ and $\gamma_t|_U$ fixed. Indeed
they come from one $W$ coordinate, one $v_0$ coordinate and
one $H$ coordinate, and wedging with $A_t$ surjects onto binary
quadratics for each $H$ coordinate. These two groups of coefficients
are independent. When $a|_H$ is a smooth conic, it must divide
$J$ or all three linear forms vanish, costing at least seven.
For a reduced conic that is two lines, each of those lines divides
$J$ or all three linear forms. If both divide $J$ this costs seven;
if both divide all three linear forms this costs nine; if one
divides $J$ and the other all three linear forms this costs
$4+6=10$. Thus the $174$ count increases to $181$.
For a doubled line the restriction of $a$ already costs three
extra, and requiring its reduced line in $J$ or in all three
linear forms costs at least four. If $a|_H=0$, its six coefficients
vanish, and the fifteen coefficients of $T|_H$ vanish also,
giving $130+30+6+15=181$ before this test.

For independent $m_0,m_1$, the $149$ base-and-wedge count in
Theorem~\ref{thm:dec1-q-transverse-ninefold} now becomes
$127+30=157$ by the full wedge condition. Its scalar quadratic
is still arbitrary. The horizontal branch has at least $25$
independent terminal conditions, thus gives at least $182$.
The vertical and zero branches already cost at least $184$.

Next suppose both cubics vanish on the five-space $K$, in the
moving-row case. The count of $176$ in
Theorem~\ref{thm:dec1-q-transverse-ninefold} uses at least
$45-10=35$ wedge conditions. If the whole wedge family is zero,
all $45$ conditions hold. Otherwise its general member on $K$ is
a nonzero rank-one quadratic $F_t=\ell^2$ after rescaling.
In the weight-$(2,1)$ outer face, choose on $W$ coordinates $x=A_t,y$
so that the light polynomial on $W\oplus K$ reads
$x^2+xP(u)+y\ell(u)^2+T(u)$. For a fixed general $t$, take
$e\in\ker\ell$ with $P(e)=0$. If $T(e)\ne0$, keep
$u=se+p e'$ with $\ell(e')=1$ and use weights
$(50,30,25,35,100)$ on $(x,y,s,p,z)$. They give
$z+x^2+yp^2+T(e)s^4$ with signed denominator $13/9$.
The pole denominator is $4/3$, greater than the radial one
$1+3/(4+2\cdot2+5)=16/13$. Thus $T|_{\ker\ell}$ must contain
the zero locus of the fixed quadratic $P|_{\ker\ell}$. This
costs at least fifteen of the $35$ quaternary quartic coefficients,
independently of the cubic data. The branch therefore gives at
least $\min(176+10,176+15)=186$.

Finally in the nonzero pure-cube branch, $T|_U$ contains the zero
locus of a freely variable quaternary quadratic $a$ after a base
cost of $165$. If $a$ is reduced, the multiples of $a$ have
dimension ten in the $35$ quartics, adding $25$. If $a$ has rank
one, that costs $10-4=6$ and its reduced linear factor in $T$
costs fifteen. When $a=0$ all quartic coefficients vanish.
The constant-row and zero quadratic cases have bounds at least
$181$ and $185$ by Corollary~\ref{cor:dec1-q-reduced-cubic-174}.
Combining these with the rank-two exclusion proves the theorem.
\end{proof}

\paragraph{A fixed quadratic factor in the cubic pencil.}
Let \(V\) be the seven-dimensional space of light variables of
a polar-rank-two line \(\ell\) in a smooth quartic ninefold, with
zero quadratic light terms and cubic pencil \(L_tQ\).
Assume \(L_t=t_0L_0+t_1L_1\) has independent coefficients, and
put \(U=\ker L_0\cap\ker L_1\). Then
\[
                       \delta_\ell(X,H)\ge52/7 .
\]
If \(\operatorname{rank}Q\ge5\) and
\(\operatorname{rank}(Q|_U)\ge2\), the stronger bound \(15/2\)
holds. In particular this applies in ranks six and seven, since
restriction to a codimension-two subspace lowers rank by at most four.

Write \(T\) for the terminal quartic. On the \(\mathbf P^6\)
joining the original line to \(\mathbf P(U)\) the only equation
is \(T|_U=0\). If this polynomial were zero or reducible, the
section would have a singular locus of dimension at least four.
For a nontrivial factorization, take the join of the line with
the common zeros of the factors in \(\mathbf P^4\). On this
singular locus the four complementary partial derivatives have
a common zero. This contradicts smoothness of \(X\).
Consequently \(T|_U\) is irreducible. First suppose that
\(\operatorname{rank}Q\ge5\) and
\(\operatorname{rank}(Q|_U)\ge2\). Choose \(e\in U\) isotropic,
outside the radical of \(Q|_U\), with \(T(e)\ne0\). Such vectors
are dense in the quadric \(Q|_U=0\), on which \(T|_U\) cannot
vanish identically. In rank two one may take a general point on
either component away from their intersection.

Use the six-space \(e^\perp\) in the light variables, with
orthogonality taken for \(Q\). The quotient
\(J=e^\perp/\langle e\rangle\) has dimension five, and the
induced quadratic form \(q\) has rank \(\operatorname{rank}Q-2\).
The restrictions of \(L_0,L_1\) to \(J\) are independent since
\(e^\perp+U=V\). Give the \(e\)-coordinate \(s\) weight three,
a complement representing \(J\) weight four, and the heavy
variables weight twelve. Give the omitted seventh light direction
weight greater than twelve. The light terms of weight twelve are
\[
                         L_t q+T(e)s^4 .
\]
For general \(t\) the restriction of \(q\) to \(\ker L_t\) has
rank at least two. If \(q\) has rank at least four this follows
by the hyperplane rank bound. For rank three, a restriction of
rank at most one would require \(L_t\) to annihilate the radical
and induce a tangent line to the nonsingular conic on the
three-dimensional quotient. These exceptional projective
covectors form a nonsingular conic, which contains no line.
Thus the three simple poles can be chosen in the required open.

At each pole choose a two-space \(E\subset\ker L_t\) on which
\(q\) is nondegenerate and a vector with \(L_t\)-value one.
On their span, affine coordinates \(x,y,z\) give
\(L_tq=xyz+x^2l(y,z)+c x^3\), after rescaling the binary
quadratic on \(E\). Give \(x\) secondary weight six, \(y,z,s\)
secondary weight three and the graph coordinate weight twelve.
Omitted directions have larger secondary weight. The initial
light form is \(xyz+T(e)s^4\). The changes in the weight-four
variables preserve the first weight, and the second filtration
is applied separately to the ordinary bounded-degree series at
each pole, as in the large-radical quadratic-pencil estimate.
Subtracting a large multiple of the first weights gives the
same initial graph and specializes every affine center to zero.
Formula (DG) gives signed denominator \(13/9\). The pole
denominator is \(1/4+(3/4)(13/9)=4/3\), and the first radial
denominator is \(1+(3/4)12/(12+3+5\cdot4)=44/35<4/3\).
Thus the line has the asserted mark.

For a further case suppose \(\operatorname{rank} Q\ge5\) and
\(\operatorname{rank}(Q|_U)=1\). Necessarily
\(\operatorname{rank}Q=5\), and the four-space
\(R=\operatorname{rad}(Q|_U)\) contains the two-dimensional
radical \(K=\operatorname{rad}Q\). Indeed the pairing gives a map
\(R\to(V/U)^*\) with kernel \(R\cap K\), which must have
dimension at least two. If \(T(e)\ne0\) for some \(e\in K\),
give \(e\) weight three and all six other light directions weight
four. The same proof applies on \(V/\langle e\rangle\), since the
induced \(Q\) still has rank five and its restriction to a
hyperplane has rank at least three.

Suppose \(T|_K=0\) and \(a=dT(e)\) has nonzero restriction to
\(U\) for some \(e\in K\). Put \(J=(\ker a)/\langle e\rangle\).
The induced quadratic form on \(J\) has rank at least three and
the two linear coefficients remain independent there. Take a
complement of \(\langle e\rangle\) in \(\ker a\) of weight five,
give \(e\) weight three and a complementary vector \(b\)
weight six, with heavy weight fifteen. The cubic light term is
\(L_tq\) on \(J\), since \(e\) annihilates both factors.
The quartic term is a nonzero multiple of \(p s^3\) in the
\(b,e\) coordinates. As above the general \(q|_{\ker L_t}\)
on \(J\) has rank at least two. After selecting the same three
directions there, secondary weights seven on \(x\), four on \(y,z\),
three on \(s\) and six on \(p\) retain \(xyz+cps^3\).
Formula (DG) gives pole denominator \(37/28\), and the radial
denominator \(1+(3/4)15/(15+3+25+6)\) is smaller.
This gives the stronger mark \(280/37\).

It remains that \(T|_K=0\) and \(dT(e)\) annihilates \(U\) at
every \(e\in K\). Choose coordinates \(s,w\) on \(K\). Then the
quartic part with exactly three factors from \(K\) is
\(L_0 f_0(s,w)+L_1 f_1(s,w)\) for binary cubics \(f_0,f_1\).
They have no common projective zero. Indeed at any common zero
\([e]\in\mathbf P(K)\) the full derivative \(dT(e)\) would vanish.
On the plane joining \([e]\) to the original line the restricted
equation and all its light and line derivatives would vanish,
because \(Q(e)=L_0(e)=L_1(e)=0\) and the quadratic light
terms vanish. The two remaining heavy normal cubics would have
a common zero on this plane, contrary to smoothness.
In particular the cubics define a degree-three morphism
\(\mathbf P^1\to\mathbf P^1\) and their general linear combination
is squarefree.

Give \(K\) weight two, a five-dimensional complement weight three,
and the two heavy coordinates weight nine. The resulting light
form at a general pole is
\[
             L_t Q+L_0 f_0(s,w)+L_1 f_1(s,w).
\]
The term entirely in \(K\) is zero, the quartic terms with at
most two factors from \(K\) have weight at least ten, and \(Q\)
is a quadratic form on \(V/K\). At an individual pole use
\(x=L_t\) and a second linear coordinate \(v\) to span the external
covectors. The quartic is \(x f(s,w)+v g(s,w)\), with \(g\)
squarefree at a general pole: its coefficient cubic is a nonzero
linear combination of \(f_0,f_1\) with projective coefficient
given by the kernel of \(L_t\) in \(V/U\). This kernel sweeps out
the whole projective line as \(t\) varies.

Within \(V/K\) choose lifts of the two external coordinate vectors
and one further vector of \(U/K\) with new coordinate \(y\) so
the restriction of \(Q\) is \(x^2+y^2\). To see this, choose
the \(y\)-vector with nonzero quadratic value, scale it and
orthogonalize both external lifts against it. Adding elements of
the isotropic two-space \(R/K\) to the external lifts now allows
arbitrary changes to their symmetric two-by-two Gram matrix,
because \(R/K\) pairs perfectly with \(V/U\). Thus their new
Gram matrix can be taken to be \(\operatorname{diag}(1,0)\).
After this ordinary linear restriction at the pole the light form is
\[
             x(x^2+y^2)+x f(s,w)+v g(s,w).
\]
Give \(x,y\) secondary weight four and \(s,w,v\) secondary
weight three, with graph weight twelve. This removes the
\(x f\) term and separates the squarefree binary cubic from
the quartic \(v g(s,w)\). In the latter three variables the
four factors are distinct linear forms. Choose a binary section
through the intersection of \(v=0\) and one factor of \(g\),
missing all other pairwise intersections. It has a root of
multiplicity exactly two. As in the cubic and transverse
quartic double-contact calculation, put the cubic in the form
of two cubes, and specialize the marked binary quartic to the
mixed monomial \(a^2b^2\). Each restriction or secondary
specialization here is made on the complete ordinary
bounded-degree series at that pole; omitted directions can
be assigned arbitrarily large weight, and at each stage one
can subtract a large multiple of the preceding homogeneous
weights to bring affine centers to the origin. Formula (DG)
then gives signed denominator \(19/13\) and pole denominator
\(35/26\). The only radial denominator comes from the initial
line extraction, and equals \(1+(3/4)9/(9+4+15)<35/26\).
This proves the stated bound \(52/7\).

We complete the lower ranks. At rank at most two the polynomial
\(Q\) has a linear factor, and the common-linear-factor argument
gives \(15/2\). At rank four, if \(\operatorname{rad}Q\not\subset U\),
the rank of \(Q|_U\) is at least two. Indeed \(U\) projects to
a subspace of dimension at least three in the nonsingular four-space.
Use the first isotropic \(e\) extraction. The form \(q\) on \(J\)
now has rank two, and its radical contains the image of a vector
of \(\operatorname{rad}Q\) outside \(U\). Since the \(L_t\) do
not all annihilate this vector, at a general pole \(q|_{\ker L_t}\)
still has rank two. This gives \(15/2\).

For the remaining cases put \(K=U\cap\operatorname{rad}Q\).
At rank four \(K=\operatorname{rad}Q\) has dimension three.
At rank three \(j=\dim K\) is two, three or four.
If \(T(e)\ne0\) at any \(e\in K\), the first extraction using all
of \(V/\langle e\rangle\) gives \(15/2\), since the rank of
the induced quadric is at least three and the two \(L_i\) remain
independent. Otherwise the full derivative gives a morphism
\[
    \varphi:\mathbf P(K)\longrightarrow\mathbf P((V/K)^*),
          \qquad [e]\longmapsto[dT(e)].
\]
It has no base point by the same plane and two normal cubics test.
The pullback of the hyperplane bundle is \(\mathcal O(3)\), so
the morphism is finite onto its image. The case \(j=4\) is
impossible because the target is then a plane. When \(j=3\),
the general \(a=dT(e)\) does not annihilate \(U\): otherwise
the finite image of a projective plane lies in the projective
line \(\mathbf P((V/U)^*)\).

At rank four and \(j=3\), work on the nonsingular four-space \(V/K\).
The image of \(\varphi\) has dimension two, so one may choose
\(e\) with \(a|_U\ne0\) outside the projective line of covectors
whose vectors under \(Q^{-1}\) belong to \(U/K\). If \(a\) is
nontangent, \(Q|_{\ker a}\) has rank three and the restrictions
of \(L_0,L_1\) there are independent. If \(a\) is tangent, the
restriction has rank two and its radical contains \(Q^{-1}a\),
where a general \(L_t\) is nonzero. In each case the additional
hyperplane \(\ker L_t\) retains rank at least two for general
\(t\). The weight-fifteen derivative extraction above therefore
gives \(280/37\). In using that extraction the light space is
\((\ker a)/\langle e\rangle\) before taking the constant residues
at the three poles, so all five transverse light directions of
weight five are present in the radial denominator.

At rank three the projective covectors annihilating
\(\operatorname{rad}Q\) and tangent to the underlying smooth conic
form its dual conic. The finite image of \(\varphi\) cannot lie
in it: for \(j=3\) the dimension is too large and for \(j=2\)
the hyperplane bundle on a conic has even degree, whereas its
putative pullback has degree three. First suppose that a general
\(a=dT(e)\) does not annihilate \(U\). Choose such \(a\) outside
the dual conic. If \(a\) does not annihilate \(\operatorname{rad}Q\),
the quadric on \(\ker a\) still has rank three. The pencil of
restrictions \(L_t|_{\ker a}\) is independent and therefore a
general member retains rank at least two. If \(a\) annihilates
\(\operatorname{rad}Q\), the form \(Q|_{\ker a}\) has rank two
and its radical contains \(\operatorname{rad}Q\), which is not
contained in \(U\) since \(j\le3\). A general \(L_t\) is
nonzero on that radical and again retains rank two.
If instead \(a|_U=0\) generically, then \(j=2\) by finiteness,
and \(U\) maps onto the full nonsingular three-space
\(V/\operatorname{rad}Q\). At a general pole
\(\ker a\cap\ker L_t=U\) and the restricted rank equals three.
Applying the derivative extraction once more finishes the proof.

\paragraph{Reducible members of the cubic pencil.}
More generally, suppose that the general nonzero cubic light form
is reducible over the algebraic closure of the parameter field.
Then \(\delta_\ell(X,H)\ge52/7\). A common linear factor gives
\(15/2\) by the common-factor test. If the coefficients have a
fixed irreducible quadratic factor, the preceding calculation
gives \(52/7\), since proportional residual linear forms also
give a common linear factor. These cases include proportional
cubic coefficients with a reducible generator.

Consider two linearly independent cubic coefficients with no
common factor. Resolve the rational pencil on \(\mathbf P(V)\)
and take the Stein factorization of the map to \(\mathbf P^1\).
The intermediate normal curve is rational by dominance from
\(\mathbf P(V)\). In characteristic zero, the general fiber over
this intermediate curve is integral by connectedness and generic
smoothness. Its images on \(\mathbf P(V)\) are members of a
pencil of some degree \(a\ge1\) without a fixed component: write
the map to the rational curve as a quotient of relatively prime
homogeneous forms of the same degree. For a general original
parameter, the finite morphism of curves has \(e\) distinct
preimages. The original cubic divisor is the sum of the resulting
\(e\) divisors, since the original pencil has no fixed component.
Thus \(3=ea\). Reducibility of the general cubic gives \(e>1\),
so \(e=3\), \(a=1\). In particular the cubics are binary in two
fixed linear coordinates and the general member is squarefree.

Their common radical \(K\) has dimension five. The terminal
quartic \(T|_K\) is irreducible by the same six-plane and four
normal cubics argument used above for \(U\). At a general smooth
point of this quartic choose a tangent line with intersection
multiplicity exactly two. Such a direction exists because its
second fundamental form is generically nonzero: in characteristic
zero identically vanishing differential would make its Gauss map
constant, forcing the hypersurface to be linear. The resulting
binary quartic and squarefree binary cubic are on disjoint light
spaces. The cubic and transverse quartic double-contact estimate
therefore gives \(52/7\).

\paragraph{Integral cubic pencils with at most three essential variables.}
Let \(\ell\) be a line of polar rank two on a smooth quartic
ninefold, and assume all quadratic light terms vanish. Suppose
the cubic light form at a geometric generic point of the line is
integral and has at most three essential variables. Then
\[
                      \delta_\ell(X,H)\ge 125/17 .
\]
In particular this line passes the sufficient line bound \(785/107\).

Write \(V\) for the seven-dimensional light space. An integral
cubic has at least three essential variables. The splitting
argument for the essential-variable bundle above shows that its
three-dimensional space of covectors \(E_t\subset V^*\) is either
constant or of the form \(\langle x,y,M_t\rangle\), where \(x,y\)
are fixed and \(M_t=t_0M_0+t_1M_1\) is moving. Indeed each
negative-degree summand in the saturated bundle has degree \(-1\)
and enters only linearly multiplied by constant columns, since the
cubic tensor is twisted by \(\mathcal O(1)\). Two such summands
would leave at most one constant column, and their coefficients
would be proportional multiples of its square, contrary to
concision. Let \(K\) annihilate the span of all \(E_t\).
Its dimension is four in the constant case and three in
the moving case. Write \(T\) for the terminal quartic.

First assume \(T(e)\ne0\) at some \(e\in K\). Give \(e\) weight
three, a complement weight four, and the heavy coordinates weight
twelve. The light terms of weight twelve are the ternary cubic
and the disjoint term \(T(e)s^4\), where \(s\) is the coordinate
on \(\langle e\rangle\). The radial denominator is
\(1+(3/4)12/(12+3+6\cdot4)=16/13\).
At each general pole, choose three complementary directions on
which the cubic is an integral ternary form \(F(X,Y,Z)\).

For a smooth plane cubic, choose a flex \([1:0:0]\) with tangent
\(Y=0\). Such a point exists since on a complex elliptic curve
the map \(p\mapsto\mathcal O(3p)\) onto \(\operatorname{Pic}^3\)
is surjective. Its tangent section is \(3p\), and irreducibility
gives a nonzero \(Z^3\) coefficient. The cubic is
\(cX^2Y+XY(\alpha Y+\beta Z)+G_3(Y,Z)\) with \(c\ne0\).
Weights \(3,6,4\) on \(X,Y,Z\) give initial form \(cX^2Y+dZ^3\)
of weight twelve. At a singular point of an integral cubic its
multiplicity is two, and the cubic is \(XQ_2(Y,Z)+G_3(Y,Z)\).
For a node take \(Q_2=YZ\) and weights \(2,5,5\), giving \(XYZ\).
For a cusp take \(Q_2=Y^2\) and weights \(2,5,4\). The \(Z^3\)
coefficient is nonzero by irreducibility, so the initial form
is \(XY^2+dZ^3\). In every case give \(s\) weight three and
the graph coordinate weight twelve.

In the flex and cusp cases formula (DG) for the disjoint monomials
has denominator
\[
       1+\max\{(1+1/2+1/3+1/4)^{-1},1/4,2/5\}=37/25 .
\]
For the node the corresponding denominator is \(13/9\).
These are specializations of complete ordinary bounded-degree
series at an individual pole, so the coordinates of the flex
or singular point can be chosen separately for each pole.
Directions omitted at this second stage have arbitrarily large
secondary weight. Subtracting a large multiple of the primary
weights contracts each affine center to the origin while preserving
the initial equation. Thus the combined pole denominator is at
most \(1/4+(3/4)(37/25)=34/25\), giving \(125/17\).

Now suppose \(T|_K=0\). In the constant case the join of \(\ell\)
and \(\mathbf P(K)\) would be a five-plane on the ninefold.
The five normal cubics have a common zero on this five-plane,
contrary to smoothness. In the moving case the three-space \(K\)
gives a morphism
\[
      \mathbf P(K)\longrightarrow\mathbf P((V/K)^*),
             \qquad [e]\longmapsto[dT(e)] .
\]
If the derivative vanished, the plane joining \([e]\) to \(\ell\)
would lie on \(X\) with all line and light derivatives zero:
every cubic light form comes from \(V/K\). Its two heavy normal
cubics have a common zero, again contradicting smoothness. This
morphism pulls the hyperplane bundle back to \(\mathcal O(3)\),
so it is finite onto a surface. For a nonconstant \(E_t\) the
intersection of its general members is \(\langle x,y\rangle\).
Choose \(e\) with \(a=dT(e)\) outside this common two-space.
Thus for general \(t\) the map from \(E_t\) to the covectors
on \(\ker a/\langle e\rangle\) is injective, and the restricted
cubic remains an integral ternary cubic.

Give \(e\) weight three, five complementary directions in
\(\ker a\) weight five, and a transverse direction with coordinate
\(p\) weight six, with heavy coordinates of weight fifteen. The
light face is the restricted cubic plus \(cps^3\) for \(c\ne0\).
Its radial denominator is \(1+(3/4)15/(15+3+25+6)=241/196\).
At each pole use secondary weights \((8,14,10)\), \((4,13,10)\)
or \((8,11,11)\) on \((X,Y,Z)\) in the flex, cusp or node case
respectively. Give \((s,p)\) weights \((7,9)\) and the graph
coordinate weight thirty. These retain the same cubic monomials
as above and \(cps^3\). In the flex and cusp cases (DG) now gives
\[
       1+\max\{(1+1/2+1/3+1/3)^{-1},1/4,2/5\}=19/13 .
\]
The nodal case gives \(10/7\). Hence the pole denominator is at
most \(35/26\), and this case gives the stronger mark \(52/7\).

\paragraph{A common cubic factor at quadratic rank zero.}
If all quadratic light coefficients vanish and the cubic pencil on
the seven lights has a fixed linear factor, the line has mark at
least $15/2$. This includes pencils all of whose members are zero
or nonreduced, by the constant repeated-factor argument of
Theorem~\ref{thm:nov30-ninefold-codimension}.

Let $K$ be the six-space of lights annihilated by the common factor.
Its join with the original line is a projective seven-space on which
the equation is the terminal quartic $T|_K$. This quartic is nonzero:
otherwise the three normal cubic derivatives on a contained seven-
space would have a simultaneous zero. Its singular locus on
$\mathbf P(K)$ is at most finite. Indeed a curve of singular points
would join the line in a subset of dimension three in the same
seven-space, on which all tangent derivatives vanish. Its three
normal derivatives must again have a common zero. A general
projective four-space in $\mathbf P(K)$ therefore avoids this
singular locus and cuts out a smooth quartic threefold $G$ by Bertini.
It is K-stable by the index-one hypersurface theorem.

Use the five retained lights with weight one and heavy weight four.
At each pole the translated graph uses the smooth component $G$,
with $c=3/4$, $h=j=4$ and $k=5$. The smooth-component inequality
in the joint two-polar calculation gives the two lower bounds for
the denominator parameter $b$ as
\[
             1+\frac{(3/4)4}{5+4}=\frac43,
             \qquad \frac{1/4+1+3/4}{2}=1.
\]
Taking $b=4/3$ also dominates $c+1/h$ and proves the claim.

\paragraph{The remaining cubic pencils at quadratic rank zero.}
Suppose the cubic pencil has no common divisor and its general
member is integral and nonnormal. Its singular locus contains a
projective four-plane, by the secant argument in
Theorem~\ref{thm:nov30-ninefold-codimension}. Every such plane for a
general member lies in the base locus by Bertini on the generated
locus of the pencil. The base scheme of the two cubic coefficients
is a codimension-two complete intersection in $\mathbf P^6$. Thus a
singular four-plane is an irreducible component of this fixed scheme.
Only finitely many can occur, and for a fixed plane the cubics
singular on it form a linear space. Hence there is a fixed five-space
of lights $K$ with $C_0,C_1\in I_K^2(3)$. This yields a necessary
condition also when individual cubic singular loci have further
components.

There is a similar reduction if the general nonzero cubic uses at
most three essential variables. The saturated generic space of its
essential linear forms gives a subbundle
$E=\bigoplus_{i=1}^h\mathcal O(-a_i)\subset\mathcal O^7$ of rank
$h\le3$. It can be obtained from second contractions of the symmetric
cubic tensor, followed by saturation. The tensor lies in
$H^0(\mathrm{Sym}^3E(1))$ by saturation and is concise in the
fiber of $E$ over the function field. Thus $a_i\in\{0,1\}$:
a summand of degree at most $-2$ cannot occur in any term of the
symmetric cube twisted by one. No two nonconstant columns occur
together in any term. Moreover there is at most one such column,
since two would leave at most one constant column and would multiply
proportional quadratic forms in that column, contrary to concision.
If all columns are constant, the pencil has at most three fixed
essential variables. Otherwise it has the form
\[
        F_t(x,y)+M_t Q(x,y)
\]
with at most two constant variables $x,y$, a moving linear form
$M_t$, and a constant binary quadratic $Q$. It consequently vanishes
twice along a fixed five-space of lights (or along a larger space).

For a cubic pencil in at most three fixed variables that is
geometrically reduced at the generic pole, take a squarefree binary
section and restrict the terminal quartic to the constant common
radical of the cubic pencil, of dimension at least four. The mixed
cubic and quartic test in Section~\ref{sec:nov30-quartic-du-val}
gives a mark of at least $52/7$ if the terminal quartic has a binary
restriction with a root of multiplicity exactly two. A sufficient
condition is an irreducible factor of multiplicity two, or a
nonlinear irreducible factor of multiplicity one. In the first case
use a general transverse line through a general smooth point on the
factor. In the second use a line in its tangent hyperplane through a
general point, avoiding the other factors. At a general point the
second fundamental form is nonzero: otherwise the differential of
the Gauss map would vanish identically, and separability in
characteristic zero would make the Gauss map constant, forcing the
factor to be linear. A direction with nonzero quadratic term gives
the required contact two.

The case of a fixed quadratic factor and moving linear factor occurs
on smooth members. In $\mathbf P^{10}$ use line coordinates $s,t$,
two heavy coordinates $z,w$ and seven lights $y_i$. For example take
\[
 z s^3+z^4+w t^3+w^4+\sum_{i=1}^{7}y_i^4+
                 \pi(sy_1+ty_2)\sum_{i=1}^{7}y_i^2=0.
\]
At coefficient zero instead of $\pi$ the two binary quartics are
squarefree and the other summands are fourth powers in disjoint
variables. The gradient has no simultaneous projective zero. Thus
the discriminant in the coefficient is a nonzero rational polynomial,
which cannot vanish at $\pi$. Along $y=z=w=0$ the two normal polars
are $s^3,t^3$, and the first light form is precisely the displayed
cubic pencil. At a general pole write its linear factor $\ell$ and
its quadratic factor $Q$. At a general point of $\ell=0$, $Q$ is a
unit. On the weighted pole space of heavy weight three, blow up
the transverse codimension-two center $(\ell,z)$. Its valuation has
discrepancy two and order one on $z+\ell Q$. The two regular
parameters have joint complete-series mean $1+1/3$ in the normalized
units, by the monomial weighted Hilbert series; the fixed term along
$z=0$ contributes $1-1/3$. The joint denominator in this individual
pole comparison is therefore at least $(2+3/4)/2=11/8$, yielding
at most $80/11<125/17$ from this comparison. The original quartic
retains its terminal terms, which the homogeneous cubic comparison
has discarded.

\paragraph{A mixed test for a fixed quadratic factor.}
Suppose the cubic pencil is $L_t Q$, with fixed quadratic $Q$ and
independent coefficients of $L_t=t_0L_0+t_1L_1$. Put
$U=\ker(L_0,L_1)$ in the seven lights. If $Q|_U$ has rank at least
four, the line has mark at least $15/2$. In particular the smooth
example just displayed passes the required line test.

Indeed $T|_U$ is nonzero and irreducible. On the projective six-space
joining the line and $\mathbf P(U)$ the equation is $T|_U$. If it
were zero, or reducible, its singular set there would have dimension
at least four: in the latter case join the line to the common zeros
of two positive-degree factors on $\mathbf P(U)=\mathbf P^4$.
The four remaining normal cubics have a common zero in such a set.
Choose a general nonradical isotropic $e\in U$ with $T(e)\ne0$,
possible since $Q|_U$ is irreducible of degree two. Let $B_Q$ denote
the symmetric bilinear form associated to $Q$. The quotient
$V'=e^\perp/\langle e\rangle$ carries an induced quadratic form.
Its subspace $U'=(U\cap e^\perp)/\langle e\rangle$ has rank at
least two, since passage through the nonradical isotropic direction
decreases the rank on $U$ by two. Choose a nondegenerate two-space
$E\subset U'$. The induced $L_0,L_1$ on $V'$ have common kernel
$U'$ of codimension two, since $e^\perp+U$ is the full light
space. Orthogonal projection away from $E$ preserves these linear
forms. Hence $E^\perp\subset V'$ contains a two-space mapping
isomorphically to $V'/U'$. On this two-space the quadratic form has
a nonzero isotropic vector $r$ over $\mathbf C$. It satisfies
$L_t(r)\ne0$ generically, $Q(r)=0$, and $B_Q(r,E)=0$.
Lift $E,r$ to $e^\perp$ and restrict to their span with $e$.
The cubic is a nonzero generic multiple of $xyz$, where $x$ comes
from $r$ and $y,z$ are independent factors of $Q|_E$. Use weight
four on $x,y,z$, weight three on the $e$-coordinate $s$, and heavy
weight twelve. Only $T(e)s^4$ survives from the quartic. The signed
denominator of the graph of $xyz+s^4$ is $13/9$ by (DG), and the
outer radial denominator is $1+9/(12+15)=4/3$. The combined pole
denominator is $4/3$, yielding the claim.

\paragraph{A cubic pencil singular on a fixed five-space.}
Use the preceding notation $C_t\in I_K^2(3)$ for a fixed five-space
$K$. Take a two-space complement $W$, and write the light cubic as
$q_{t,u}(w)+B_t(w)$, where $q_{t,u}$ is linear in $u\in K$ and
quadratic in $w\in W$. If the general cubic has at least four
essential variables, $q_{t,u}$ has rank two for some $u$ at a general
pole. Otherwise its linear image in the binary quadratics consists
of squares and zero, hence is at most a line, and the entire cubic
uses only the two $W$ coordinates and one linear functional of $u$.
Set $R_t=\ker(u\mapsto q_{t,u})$; it has vector dimension at least
two. A deficient line necessarily satisfies $T|_{R_t}=0$ for general
$t$. Indeed assign outer weights three to $W$, two to $K$, and eight
to the two heavy directions. The light face is $q_{t,u}(w)+T|_K(u)$.
If $T(e)\ne0$ for some $e\in R_t$, choose $f$ with $q_{t,f}$ of
rank two. Retain only $e,f$ in $K$. Give the two $W$ coordinates
secondary weight two, the $f$ coordinate weight four, the $e$
coordinate weight two, and the graph coordinate weight eight. The
initial graph is $z+xyv+c s^4$, $c\ne0$, with $v,s$ the coordinates
of $f,e$. Formula (DG) again gives $13/9$. The outer denominator
is $1+6/(8+6+10)=5/4$, giving the line mark $15/2$.

If the generic map $u\mapsto q_{t,u}$ has rank three, deficiency
requires that its two-dimensional kernel be fixed and that
$\mathbf P(R_t)$ lie in $\operatorname{Sing}(T|_K)$. Indeed, for
$e\in R_t$ suppose the restricted derivative $b=d(T|_K)(e)$ is
nonzero. The image of $\ker b$ under the rank-three map has dimension
at least two, so take $f\in\ker b$ with $q_{t,f}$ nondegenerate.
Choose $h\in K$ with $b(h)\ne0$. Keep $e,f,h$ on $K$ in the above
outer face and give their coordinates $s,v,p$ secondary weights
three, five and six. With both $W$ coordinates of weight five and
the graph coordinate of weight fifteen, the light form specializes
to $xyv+cps^3$ for $c\ne0$. Here $q_{t,e}=0$, $T(e)=0$, and the
coefficient of $vs^3$ is zero since $b(f)=0$. All quartic terms
with at least two factors from $v,p$ have weight at least sixteen.
Formula (DG) gives signed denominator $10/7$, hence combined pole
denominator $37/28$ and mark at least $280/37$. On a deficient line
it follows that $\mathbf P(R_t)$ lies in the singular locus of
$T|_K$ for general $t$. That singular locus has dimension at most
one, by joining it to the original line as in the fixed-quadratic
argument. A projective algebraic curve contains only finitely many
distinct line components, forcing the kernel line to be fixed.

\paragraph{Powers of smooth components at the third boundary.}
In the joint two-polar calculation let $c=(d-1)/d$, and suppose the first light form restricts at a general pole to $G=Q^s$ on $k$ constant lights, where $j=\deg Q\ge2$, $h=js$, $k>j$, and $Q\subset\mathbf P^{k-1}$ is a smooth K-semistable hypersurface. For $b\ge\max(1,c+1/h)$, $b\le1+c$, the angular and radial tests hold if
\[
 b\ge 1+\frac{ch}{k+h},\qquad
 b\ge \frac{1/j+s(1+c)}{s+1}.
\]
Indeed, for every $L\ge0$ and every divisorial test on $\mathbf P^{k-1}$,
\[
 kS_H+L v(Q)\le\max\{1+jL/k,\,1/j+L\}\,A.
\]
To see this use $kS_H\le A$, $v(Q)\le A$, and K-semistability of the log pair at $\theta=k(j-1)/(j(k-1))$, by the interpolation theorem in \cite[Theorem~1.7]{ar2003.06839}. The latter gives $(k-j\theta)S_H+\theta v(Q)\le A$. Its upper bound on $S_H$ crosses the one from $kS_H\le A$ at $v(Q)/A=j/k$. At $v(Q)/A=1$ it gives $kS_H/A\le1/j$. The two endpoints give the stated formula. Apply this with $L=s(1+c-b)$ in Section~\ref{sec:nov30-joint-poles}. The first inequality above is also exactly the radial comparison. At a general point of $Z\cap Q$ on the weighted space $\mathbf P(1^k,h)$ the monomial valuation in $(Q,z)$ with weights $(1,s)$ has discrepancy $1+s$ and order $s$ on $z+Q^s$. Its complete weighted series contributes $(k+1)S=1/j+s/h$ by the regular sequence $Q,z$ and summing the corresponding weights of the homogeneous section ring. Adding the fixed $Z$ coefficient therefore gives the second displayed denominator as a necessary bound in this joint comparison.

For $n=2d+1$ and $d\ge33$, use the manuscript's target $b=(2d+3)/(d+4)$. Every $j\ge3$ for which a general $(j+1)$-dimensional vector section of the component is smooth passes both tests: use $k=j+1$, so the component section is a smooth index-one hypersurface. Indeed $k\ge4$, $h\le d$, and the second comparison for $s\le d/j$ follows from $(2j-5)d+3j>0$. For a smooth quadric section the second comparison becomes
\[
 s\le \frac{d(3d+2)}{8(d-1)}.
\]
These are also valid in the rank-two conic flag whenever the two-heavy model of Section~\ref{sec:nov27-conic-quadric-flag} exists. In fact the total regular share is $1/d$ and each of its $2d-2$ pole groups has asymptotic share $1/(2d)$. Refining the same regular kernel separately at each root and averaging $R+(2d-2)P_a$ gives the identical joint coefficient $c$.

A hyperplane power $G=x_1^h$ illustrates the remaining numerical obstacle for every $d$. The local weighted valuation with $v(x_1)=1$, $v(z)=h$ at their generic intersection has discrepancy $1+h$, and $D_0=1+h$, $v(z+x_1^h)=h$. Its joint denominator is $1+ch/(h+1)$. On the large-degree boundary its excess above $b$ is exactly
\[
 \frac{(d-1)(4h-d)}{d(h+1)(d+4)}.
\]
Such a nonterminal first light form can occur on a smooth hypersurface. With $[t_0:t_1]$ on the line, two heavy coordinates $y_1,y_2$ and the remaining light coordinates $x$, take $t_0^{d-1}y_1+t_1^{d-1}y_2+t_0^{d-h}x_1^h+F_d(x,y)$. Varying $F_d$ has no base point away from the line, and the two polar derivatives ensure smoothness on it. A general member is smooth by Bertini. This valuation tests the initial system; the full hypersurface retains the higher light coefficients. For example, with $d=40$, the pure quadric power of order $h=32$, $s=16$, also gives $321/170>83/44$ in the single-form joint test. If the first subsequent restricted form is a thirty-third power of a linear equation, the sufficient next-form inequality in Section~\ref{sec:nov27-quadric-next-form} would require $528\le490$. Additional estimates on the higher terms would therefore be needed in this approach.

\subsection{The exceptional locus for smooth sextic thirteenfolds}
\begin{theorem}\label{thm:dec1-q-sextic-thirteen-incidence}
The locus in $\mathcal U_{13,6}$ failing uniform K-stability has
codimension at least $331$.
\end{theorem}
\begin{proof}
The sufficient line mark is $951/100$ by the sextic column of
Section~\ref{sec:nov28-two-hyperplanes}. If the polar rank $k$ of
a line is at least five, the $k-1$ positive ordinary splitting
degrees have sum five, hence maximum at most two. The angular
denominator in Section~\ref{models:sec:local-line-models} is then at most
$4/3$, larger also than the radial denominator $(13-1/6)/12$.
The line mark at least $21/2$ exceeds the target. Thus assume
$2\le k\le4$, where the lower inequality follows from basepoint
freeness of the polars. Retain all $m=13-k$ constant lights and
two general heavy directions whose binary quintic polars have no
common zero. Such a pair exists in any basepoint-free space on
$\mathbf P^1$ by incidence. The resulting linear flag is smooth
along the line, and the normalized bound lifts by the hyperplane
flag estimate as in Section~\ref{sec:nov27-sextic-line-nets}.

Let $Q_t,C_t$ be the normal quadratic and cubic light polynomials,
which have base degrees four and three respectively. For generic
quadratic rank at least three retain three general lights and a
smooth conic at each general pole. In the joint test of
Section~\ref{sec:nov30-joint-poles} take $c=5/6$, $b=4/3$.
Writing $a=v(Q_t)$ on the plane, the conic flag gives
$3S_H\le\min\{A,2A-3a/2\}$, so
$3S_H+(1/2)a\le(4/3)A$ by the intersection $a/A=2/3$.
Both the radial term $1+2c/(3+2)$ and $c+1/2$ equal $4/3$.
The line therefore meets the stronger mark $21/2$.

Counting at a fixed line costs seven conditions for containment,
and $(6-k)(13-k)$ for rank $k$ of its matrix of binary quintic
polars. The Grassmannian of lines in $\mathbf P^{14}$ has dimension
twenty-six. Set $B_k=7+(6-k)m-26$. For each fixed polar matrix
the $5\binom{m+1}2$ quadratic coefficients and the four independent
cubic coefficients are unrestricted by the previous equations.

At generic quadratic rank two let $q$ denote the number of constant
active variables and $l=m-q$ the dimension of the common radical.
The saturated row space is a rank-two subbundle
$E=\mathcal O(-a)\oplus\mathcal O(-b)$ of the trivial light dual,
with $0\le a\le b$, $a+b\le4$ by the determinant of the generically
invertible map $E^*\to E(4)$. Thus $2\le q\le2+a+b$. Writing
$g(x)=\max(0,1+x)$, the affine dimension of these quadratic families
is at most
\[
 D_q=q(m-q)+q(2+a+b)-2-g(b-a)-g(a-b)
                +g(4-2a)+g(4-a-b)+g(4-2b).
\]
The terms choose the constant active $q$-space, its bundle embedding
up to source automorphisms, and the tensor in $H^0(\mathrm{Sym}^2E(4))$,
respectively, as in Section~\ref{sec:nov30-quartic-rank-family}.

All four cubics must vanish on the common radical. Otherwise choose
a constant direction $u$ there where the cubic is generically
nonzero and two further constant directions carrying a nondegenerate
quadratic. Weights $(3,3,2)$ and heavy weight six give the pole graph
$z+xy+s^3$ after linear changes at each pole. Formula (DG) of
Section~\ref{sec:nov28-disjoint-general} gives signed denominator
$10/7$, so the pole denominator is $(1+5(10/7))/6=19/14$.
The radial denominator is $1+5/14$, and $196/19>951/100$.
Thus the cubic vanishing costs $4\binom{l+2}{3}$ conditions.
Moreover for each $u$ in the radical the normal covector $dC_t(u)$
must belong to the row space of $Q_t$. Indeed give the $q$ active
variables weight two, the radical variables weight one, and the
heavy variables weight four. The face is $Q_t(w)+dC_t(u)(w)+T_t(u)$,
where $T_t$ is now the degree-four coefficient at the given pole;
higher normal terms have strictly greater weight. If the covector
lies outside the quadratic row space at a general pole, retain two
vectors in its kernel carrying $xy$ and a radical direction of this
single quadric on which the covector is nonzero. Together with $s u$
these give the graph $z+xy+v s^2+A s^4$. Secondary weights
$(1,2,1,1,3)$ remove the last term, as in
Theorem~\ref{thm:dec1-q-final-ninefold}. The signed denominator
$7/5$ gives pole denominator $4/3$ and the outer radial denominator
is $1+(10/3)/(4+m+q)\le4/3$. Again the mark exceeds the target.

Consequently the degree-three covector section must belong to
$H^0(E(3))$, of dimension $g(3-a)+g(3-b)$. Its coefficients on
the $\binom{l+1}2$ quadratic radical monomials are independently
arbitrary in $4q$ dimensions before this condition, by their
distinct cubic monomials. For each splitting type and $q$ the cost
is therefore at least
\[
 B_k+5\binom{m+1}2-D_q+4\binom{l+2}3
             +\binom{l+1}2\bigl(4q-g(3-a)-g(3-b)\bigr).
\]
For transparency the table gives the minimum for each $q$, over
the at most nine pairs $0\le a\le b$, $q-2\le a+b\le4$.
\[
\begin{array}{c|rrrrr}
m\backslash q&2&3&4&5&6\\\hline
9&531&518&473&407&331 \\
10&735&724&675&599&507 \\
11&982&974&922&837&730
\end{array}
\]

At generic quadratic rank one $E=\mathcal O(-a)$, $a\le2$, and
the quadratic families have dimension at most
$\max_{0\le a\le2}\{m(a+1)-1+g(4-2a)\}=3m$. There are at least
$m-3\ge6$ common radical directions. If the generic cubic there is
nonzero and reduced, restrict it to two general constant directions
to obtain a squarefree binary cubic. Add one constant direction
carrying the square. With light weights $(3,2,2)$ and heavy weight
six the pole graph is $z+x^2+u^3+v^3$ after a binary linear change.
The signed denominator is $19/13$ by (DG), and the radial term is
$1+5/13$, so the resulting line mark exceeds $951/100$. Thus on a
fixed six-dimensional subspace of the radical all four endpoint
cubic coefficients, taken at four independent base evaluations,
are zero or nonreduced. To justify a fixed six-space, either apply
the preceding argument directly on it, or specialize to it a generic
nonreduced form of the shape $L^2M$. The closed nonreduced cone has
dimension at most $6+6-1=11$ in the $56$-space of cubics. The
contribution at rank one is at least
$B_k+5\binom{m+1}2-3m+4(56-11)\ge377$.

At rank zero the quadrics already cost $5\binom{m+1}2$, at least
$225$. If the generic cubic is normal, integral and has at least
four essential variables, take a general four-dimensional section
giving a Du Val cubic surface as in Section~\ref{sec:nov30-quartic-du-val}.
Its threshold $\delta\ge1/3$ gives $4S_H\le\min\{A,9A-26v(G)/3\}$.
The two linear bounds meet at $v(G)/A=12/13$. Taking
$b=7/5$ and $1+c-b=13/30$ in the joint pole test gives
$4S_H+(13/30)v(G)\le(7/5)A$. The radial term is
$1+(5/6)3/(4+3)=19/14$, and $c+1/3\le b$. The resulting mark
is ten. Use a fixed six-space of the light directions. A deficient
line has cubic pencil on this six-space entirely in the closed
exceptional cone. Reducible cubics there have dimension at most
$6+21-1=26$, integral nonnormal cubics at most twenty-four by
the singular three-plane incidence in the proof of
Theorem~\ref{thm:nov30-ninefold-codimension}, and cubics in at most
three variables at most $9+10=19$. Thus four independent evaluations
cost $4(56-26)=120$. This yields at least $-1+225+120=344$.

The minimum among the strata is $331$. The parameter spaces of
splittings, embeddings, restrictions and closed exceptional cubics
give constructible incidences of the indicated dimensions. The
closure and projection argument is identical to that of
Theorem~\ref{thm:nov30-ninefold-codimension}, with $B_k$ already
accounting for line variation.
\end{proof}

\paragraph{Degrees seven and eight at the same boundary.}
For $d=7,8$, $n=2d+1$, the same proof yields codimension at least
$731,1401$, respectively, for the complement of the uniformly
K-stable locus in the smooth parameter space. Here are the formulae
and checks for the change of degree. Put $N=2d+2$, $c=(d-1)/d$,
and retain $m=n-k$ lights. When $k\ge5$ the ordinary denominator
is at most $1+(d-4)/d$, larger than the radial denominator. Thus
the line marks are at least $56/5$ and $12$, respectively, exceeding
the two targets $10491/1000$ and $11485/1000$ in
Section~\ref{sec:nov28-two-hyperplanes}.

The conic joint test now uses $b=1+2c/5$. At $v(Q)/A=2/3$
the inequality $3S+(1+c-b)v(Q)\le b A$ is an equality for the
upper envelope from the conic flag. This $b$ is at least $c+1/2$
for $d\le6$ only; for the present two degrees take instead
$b=\max(1+2c/5,c+1/2)=c+1/2$. Increasing $b$ preserves the
projective test, and the radial bound is $1+2c/5$. For a smooth
conic the final marks are thus at least $224/19$ and $144/11$,
respectively. At rank two the signed graph tests $xy+s^3$ and
$xy+v s^2$ retain the same denominators $10/7$ and $7/5$ and the
outer radial terms are at most $1+3c/7$ and $1+4c/(4+m+q)$.
At rank one the reduced cubic test has signed denominator $19/13$
and outer radial term $1+6c/13$. Combining with pole shares
$1/d,c$ gives marks exceeding the stated targets in both degrees.
At rank zero use all $m$ lights and apply the Du Val surface test
to four general directions. The bound on the resulting projective
three-space has the same crossing $v(G)/A=12/13$ as
before, so it suffices to take $b=1+12c/25$. This dominates $c+1/3$
and the radial $1+3c/7$ for both degrees. The resulting line
thresholds again exceed the required marks.

Replace $B_k$ by $d+1+(d-k)m-2n$. In $D_q$ replace the quadratic
coefficient degree four by $d-2$ in the three symmetric-square
terms, allowing $0\le a\le b$, $a+b\le d-2$, $2\le q\le2+a+b$.
In the rank-two cost replace the first coefficient five by $d-1$,
the first cubic coefficient four by $d-2$, and the final bracket by
$(d-2)q-g(d-3-a)-g(d-3-b)$. At rank one set
$l=m-1-\lfloor(d-2)/2\rfloor$. The quadratic family has dimension
at most
\[
     \max_{0\le a\le(d-2)/2}
                  \{m(a+1)-1+g(d-2-2a)\},
\]
and $d-2$ independent cubic restrictions to $l$ radical directions
each cost $\binom{l+2}{3}-(2l-1)$. At rank zero the cubic
parameters lie in the exceptional cone of dimension at most
\[
 \max\left\{m+\tbinom{m+1}2-1,\quad
 \tbinom{m+2}3-\tbinom m3-2\tbinom{m-1}2+2(m-2),\quad
 3(m-3)+10\right\}.
\]
These contributions allow arbitrary reducible cubics, cubics
singular along a projective $(m-3)$-plane, and cubics in at most
three variables, respectively. As in the sextic proof, subtract
this bound from the total cubic dimension and multiply by $d-2$.
Substitution gives the following entries in the order of quadratic
ranks two, one, zero:
\[
\begin{array}{cc|rrr}
d&k&2&1&0\\\hline
7&2&1409&1554&2349\\
7&3&1036&1197&1869\\
7&4&731&898&1457\\
8&2&2466&2435&4181\\
8&3&1887&1930&3432\\
8&4&1401&1500&2776
\end{array}
\]
Taking the minima and applying the marked line criterion completes
the proof. The program \path{dec1-line-incidence-check.py}
reproduces these entries and the preceding sextic table.

\paragraph{Moving cubic double factors on a sextic thirteenfold.}
The sufficient all-line threshold for a smooth sextic thirteenfold
from Section~\ref{sec:nov28-two-hyperplanes} is $951/100$. There is a
single base point that requires additional care when transferring the
directed-line calculation of Theorem~\ref{thm:nov28-quintic-eleven}.
In polar rank two a generically nonzero first light cubic of the form
$L_t^2M_t$ with moving repeated factor and distinct generic factors
has primitive factors of base degrees $1$ and $e\in\{0,1\}$. In fact
the multiplicities distinguish the two rational linear factors over
the base function field. Clear their denominators and remove their
common coefficient factors. Gauss's lemma gives
$C_t=h(t)L_t^2M_t$ with $\deg h+2\deg L+\deg M=6-3=3$.
Primitivity on the projective line ensures that neither primitive
factor ever vanishes as a whole form. Thus the sole additional case
has $M$ constant and $h$ linear. At its unique root $p$, the tangent
hyperplane contact has multiplicity at most four. Indeed the
tangent hyperplane contains the lights and the original line.
In affine coordinates there the transverse-degree-three coefficient
has a simple base zero, giving a nonzero degree-four term in the
local variables. Terms of other transverse degrees cannot cancel it.

This bound of four is attained. Use coordinates
$s,t,z,w,y_1,\ldots,y_{11}$ on $\mathbf P^{14}$ and set
\[
 F=z s^5+z^6+w t^5+w^6+\sum_{i=1}^{11}y_i^6
             +\pi\,t(sy_1+ty_2)^2y_3.
\]
For coefficient zero in place of $\pi$ the two binary sextics are
squarefree and the remaining summands are pure powers in disjoint
variables. Their gradients therefore have no simultaneous projective
zero. The discriminant in that coefficient is a nonzero rational
polynomial, so its value at $\pi$ is nonzero. On
$C=(z=w=y_i=0)$ the normal polars are $s^5,t^5$. On the eleven
lights the first nonzero term is $\pi t(sy_1+ty_2)^2y_3$.
At $p=[1:0]$ the tangent hyperplane is $z=0$ and the first nonzero
ordinary term there is $\pi t y_1^2 y_3$ of degree four.

For this particular line the mixed pole method gives
$\delta_C(X,H)\ge168/17>951/100$. Give $y_1,y_2$ weight five,
$y_3,y_4,y_5,y_6$ weight two, the remaining five lights weight three,
and both heavy directions weight twelve. At a general pole retain
$x=sy_1+ty_2$, $v=y_3$ and $y_4,y_5,y_6$. The light face contains
$x^2v+v^6+y_4^6+y_5^6+y_6^6$, with nonzero coefficients. Give
$x,v$ secondary weights four, the remaining three lights weight two,
and the graph coordinate weight twelve. This removes $v^6$ and leaves
$x^2v+y_4^6+y_5^6+y_6^6$. Formula (DG) of
Section~\ref{sec:nov28-disjoint-general} gives the signed denominator
\[
     D=1+\max\{(1+1/2+3/6)^{-1},1/6,3/7\}=3/2.
\]
The pole denominator $(1+5D)/6=17/12$ gives the asserted line mark.
The outer denominator is smaller since the total light weight is
$33$, giving $1+10/45$. The secondary specialization uses the
whole-chart comparison of Section~\ref{sec:nov30-quartic-rank-family}.

\section{Sharper bounds in the full parameter space}

\paragraph{Minimal non-lc loci and transverse jets along a linear space.}
Suppose \(e\le s+1\) and a degree-\(e\) hypersurface on \(\mathbf P^s\)
has non-lc locus of dimension \(k=s+1-e\). Its top-dimensional
reduced part is a linear \(k\)-plane. To see this, take a general
\((e-1)\)-plane \(H\). The non-lc locus on \(H\) is precisely the
intersection of the ambient non-lc locus with \(H\). In one direction
this follows by Bertini on a log resolution, and in the other by
inversion of adjunction successively for the transverse hyperplanes.
This intersection is finite. For sufficiently small \(\epsilon>0\)
it is exactly the non-klt locus of
\((H,(1-\epsilon)D|_H)\), as can be checked on a log resolution.
The negative log canonical class of this pair is ample. The
connectedness theorem gives one non-klt point. A general \(H\)
meets the top-dimensional reduced locus transversely in its
degree, so that degree is one.

\subsection{Quartic log-canonicity in higher-dimensional projective space}
\label{sec:dec1-quartic-support}
Write \(T(j,h)=\binom{j+h}{h}\) and retain the notation \(b(s,e)\)
for the codimension of the non-log-canonical cone, including zero.
For \(s\ge4\) put \(k=s-3,\ t=k+1,\ Q=T(k,2)\), and
\(F=T(k,4)+3T(k,3)\).

\begin{proposition}\label{prop:dec1-quartic-support}
The number \(b(s,4)\) is at least the minimum of the following
eight integers. In the last five put \(p=s-2\).
\[
\begin{gathered}
 F+5Q+k+2,\qquad F+6Q-6t+2k+10,\qquad F+6Q,\\
 T(p,4)+2T(p,3)-2(s-1)
       +\min\{2T(p,2)-1,\,3T(p,2)-2(s-1)\},\\
 T(s-1,4)+T(s-1,3)-2s-1,\\
 T(p,4)+2T(p,3)+2T(p,2)+4T(p,1)-6(s+1)+12,\\
 T(s,4)-T(s,3)-s+b(s-1,3),\qquad
 T(s,4)-2T(s,2)+\binom{s-1}{2} .
\end{gathered}
\]
For \(s=4,\ldots,10\) these lower bounds are respectively
\(35,77,147,253,404,610,879\).
\end{proposition}
\begin{proof}
The non-lc locus of a nonzero quartic has dimension at least \(s-3\)
by the linear-section argument in Paragraph~\ref{par:new-lc-linear-sections}.
If this dimension is minimal, its top-dimensional reduced support
is a linear \(k\)-plane \(L\), by the connectedness argument for
minimal non-lc loci above. Fix \(L\) and choose three normal variables.
The parts of transverse order zero and one vanish, at cost \(F\).
The quadratic part is a symmetric three-by-three matrix of quadrics
on \(L\). It has rank at most one at its generic point. Indeed an
ordinary quadratic tangent of rank at least two on a transverse
three-space gives an lc germ there by scaling, and inversion of
adjunction gives lc in the original space. The same restriction
and adjunction can be made at each general point of \(L\).

A nonzero polynomial symmetric matrix of generic rank one over a
polynomial ring has the form \(hvv^{\mathsf t}\) with \(v\) a
primitive polynomial column. To see this take a vector spanning
the rational column space, clear denominators, and divide by its
polynomial gcd. Symmetry gives the displayed form with rational \(h\);
at each irreducible denominator some component of \(v\) is a unit,
so \(h\) is polynomial. Taking one nonzero original column in this
construction makes the entries of \(v\) homogeneous of equal degree.
Here \(v\) is constant or linear.

Suppose first that \(v\) is constant. The space of such nonzero
matrices has dimension at most \(Q+2\), as the line of \(v\) costs
two parameters. Choose coordinates \(x,y,z\) normal to \(L\) with
transverse quadratic \(h x^2\). At a general point, the binary
cubic obtained from the transverse third term by setting \(x=0\)
must be a cube or zero. Indeed split off the nonsingular square
analytically. The residual must have threshold below \(1/2\), and
any binary cubic other than a cube or zero has threshold at least
\(1/2\) at the origin. The binary cubic coefficients are linear
forms on \(L\). Their image is consequently a linear space of pure
cubes and zeros, of dimension at most one since the twisted cubic
contains no line. For nonzero image we may arrange it to be the
span of \(y^3\). Its entire family has dimension at most \(t+1\)
in the \(4t\) binary cubic coefficients, at cost \(3t-1\).

Write \(a,b\) for the linear coefficients of \(xz^2,xyz\) in
the transverse cubic, and \(c,e\) for the scalar coefficients of
\(z^4,yz^3\) in the transverse quartic. After splitting off the
square, the binary quartic jet is the original one at \(x=0\)
minus the square of the \(x\)-derivative of the cubic at \(x=0\),
divided by \(4h\). A nonzero coefficient of \(z^4\), or then of
\(yz^3\), gives weighted initial \(y^3+Cz^4\) or
\(y^3+Cyz^3\) at a general point, of threshold greater than \(1/2\).
Thus \(4hc=a^2,\ 2he=ab\) identically on \(L\). If \(h\) is not a
scalar square, \(a=0,\ c=e=0\), at cost \(t+2\) independent of the
previous binary-cubic conditions.
Here are two more tests at one general point. Denote by \(q,j\) the
free quartic coefficients of \(xz^3,y^2z^2\). If \(a=0\), then
\(c=e=0\), and for weights \(3,2,1\) on \(x,y,z\) the form of
weight six is
\[
 h x^2+bxyz+P y^3+q xz^3+j y^2z^2,\qquad P\ne0 .
\]
Every term has weight at least six. Completing the square and
putting \(Y=y/z^2\) gives the binary residual
\[
 z^6\big(PY^3+(j-b^2/(4h))Y^2-bqY/(2h)-q^2/(4h)\big).
\]
Unless the cubic has a triple root, this residual is lc with
coefficient \(1/2\): the \((2,1)\)-weighted blowup has exceptional
coefficient one, and the residual roots on its exceptional line have
multiplicities at most two. They are on the smooth \(z\)-chart
since the leading cubic coefficient is nonzero; there the strict
transforms are disjoint lines transverse to the exceptional divisor.
Thom--Sebastiani after completing the square and semicontinuity for
the weighted degeneration then give lc for the original equation.
A triple root costs two in \((q,j)\) with the previous jets fixed.
In fact a root \(\lambda\) must satisfy
\(3P\lambda^2=-bq/(2h)\), \(P\lambda^3=q^2/(4h)\), giving
finitely many \((q,\lambda)\) also when \(b=0\); \(j\) is then fixed.

If \(a\ne0\), splitting off the square on \(y=0\) gives the
coefficient of \(z^5\) as \(a^2B/(4h^2)-a q/(2h)\), with
\(B\) the cubic coefficient of \(x^2z\). Indeed the critical
value substitutes \(x=-az^2/(2h)+O(z^3)\); the unknown coefficient
in order three cancels. A nonzero \(z^5\) gives initial binary
\(P y^3+Cz^5\) of threshold \(8/15>1/2\), fixing \(q\).
These two local tests also apply after arbitrary linear changes
of normal coordinates at one chosen point.

For nonsquare \(h\), the \(a=0\) test adds two to the previous
count. If \(h\) is a scalar square, this already costs \(Q-t\ge1\)
within the quadrics. For \(a=0\) the same count applies. For
\(a\ne0\), both \(a\) and \(b\) are proportional to the linear
square root of \(h\), by the two identities with \(c,e\).
This costs \(2(t-1)\), and \(c,e\) and then \(q\) are fixed.
Since \(t\ge2\), this also gives at least \(t+4\) constraints
besides the original binary-cubic conditions and rank-one
matrix dimension bound, counting the square condition on \(h\).
If the binary cubic is identically zero, it costs \(4t\) conditions.
At one general point the free binary quartic after splitting must then
have a root of multiplicity at least three or be zero, at cost two.
At a second general point there is one more condition. After fixing
the initial quartic test, the \(x\)-derivative of the cubic on the
kernel at the second point is still a free binary quadratic \(R\),
by linear interpolation. For any fixed binary quartic \(G\) the
quartics \(G-R^2/(4h)\) are not all in the triple-root locus:
scale a squarefree \(R\) to see a double of a squarefree quadratic
in their projective closure. This gives the same total in every case.
Subtracting
\(3t\) parameters for \(L\) gives the first bound in the statement.

If \(v\) is primitive linear, its three entries span at least two
linear forms, and the family of matrices has dimension at most \(3t\)
after absorbing a scalar square root into \(v\). On a general line
\(\ell\subset L\) (all of \(L\) if \(k=1\)) the row of its entries
surjects \(\mathcal O_\ell^3\to\mathcal O_\ell(1)\); its kernel
is \(E=\mathcal O_\ell\oplus\mathcal O_\ell(-1)\), by a constant
change of normal coordinates. The cubic restricted to \(E\) is an
arbitrary section of
\[
 \operatorname{Sym}^3 E^*\otimes\mathcal O_\ell(1)
  =\mathcal O_\ell(1)\oplus\mathcal O_\ell(2)
      \oplus\mathcal O_\ell(3)\oplus\mathcal O_\ell(4).
\]
Indeed the kernel of the map from
\(\operatorname{Sym}^3(\mathcal O^3_\ell)^*(1)\) is generated by
the row \(v\) times quadratic forms, hence is
\(\operatorname{Sym}^2(\mathcal O^3_\ell)^*\) with vanishing \(H^1\).
The target has fourteen sections. A nonzero section that is generically
a cube up to scalar determines a saturated line subbundle
\(M=\mathcal O_\ell(m)\subset E^*\) and is a section of \(M^3(1)\).
This follows by taking the unique line of factors in the generic
fiber and saturating; \((\operatorname{Sym}^3 E^*)/M^3\) is locally
free. Necessarily \(0\le m\le1\). For \(m=1\) the subbundle is
unique and there are five sections of \(M^3(1)\). For \(m=0\)
there are at most two parameters for the subbundle and two for the
section. This costs at least nine of the fourteen cubic parameters.
At \(k-1\) further general points, the binary cubic on the kernel
costs two at each. Their evaluations of the linear coefficients
are independent of the restriction to \(\ell\) and of one another,
by linear independence modulo the two-space defining the line.
When the cubic on the kernel is generically nonzero, the quartic
jet supplies at least three independent further conditions at one
general point by the \(z^4,yz^3\) test followed by the preceding
test with \(q,j\) when \(a=0\) or with \(q\) when \(a\ne0\).
If it is identically
zero, the line test already costs fourteen, which suffices. These
counts give the second bound after varying \(L\). In applying a test
at a general point or line one can fix the incidence stratum and
make the choice at its generic point; the required genericity is
open, so no point or line parameters need be subtracted.

For the identically zero transverse Hessian the full \(6Q\)
coefficients vanish. The transverse ternary cubics form a linear
system whose general member must be non-lc or zero, since an
lc cubic gives an lc initial cone. This costs at least \(3t\).
If the general member is reducible, all \(t\) independent
evaluations belong to the reducible non-lc cone of affine dimension
at most seven. Indeed reducible ternary cubics form an irreducible
cone of dimension eight, with general member a line transverse to
a smooth conic. If the general member is integral, it is cuspidal.
For a vector space of cubics of dimension at least two containing
such a general member, its base locus in the plane is finite.
Bertini makes the general cusp occur at one of these finitely many
base points, hence the whole linear system is singular at a fixed
point. Its tangent quadratics span at most a line of squares,
since the conic of binary squares contains no line. At a fixed
point with fixed cusp tangent, each cubic therefore lies in a
five-dimensional vector space (three first-jet conditions and two
quadratic conditions). Varying the flag costs three parameters,
so the codimension is at least \(5t-3\ge3t\). A one-dimensional
space of cuspidal cubics costs at least \(9t-7\ge3t\), by choosing
its projective generator in dimension seven and its \(t\) scalars.
The identically zero case is immediate. Subtract the parameters
of \(L\) to obtain the third bound.

We next consider an integral quartic whose non-lc locus has
dimension at least \(s-2\). It contains an integral codimension-two
variety \(Y\) of degree at most three. Indeed it is contained in
the singular locus, which has dimension at most \(s-2\). A general
plane cuts an integral quartic curve, singular at each of the
\(\deg Y\) distinct intersection points. Its arithmetic genus is
three. At general points of \(Y\) the quartic vanishes to order at
least two and its quadratic transverse term has rank at most one.

If \(Y\) is a linear \(p\)-plane, its double vanishing costs
\(T(p,4)+2T(p,3)\), with Grassmannian dimension \(2(s-1)\).
Its two-by-two transverse Hessian has \(3T(p,2)\) coefficients.
The rank-one parameterization above bounds its family dimension by
\(\max\{T(p,2)+1,2(s-1)\}\). This proves the fourth bound.

If the span of \(Y\) is a hyperplane \(H=(u=0)\), its equation there
has degree \(j=2\) or \(3\). The restriction of the quartic is
nonzero by integrality and divisible by its square, which forces
\(j=2\). Fix \(H\) and the irreducible quadratic \(G\) on \(H\).
The quartic begins
\(cG^2+uGA_1+u^2B_2+O(u^3)\), where \(c\ne0\).
The first two restrictions cost
\(T(s-1,4)-1+T(s-1,3)-s\). The quadratic term transversely to
\(Y\) has rank at most one, so \(4cB_2-A_1^2\) vanishes modulo
\(G\); this costs \(T(s-1,2)-1\) further conditions. Vary \(H,G\)
in dimension \(s+T(s-1,2)-1\) to obtain the fifth bound.

If \(Y\) spans \(\mathbf P^s\), it has minimal degree three in
codimension two. The classification of varieties of minimal degree
(D.~Eisenbud and J.~Harris, \emph{On varieties of minimal degree
(a centennial account)}, Proc. Sympos. Pure Math. \textbf{46},
Part 1 (1987), 3--13) gives a rational normal scroll or its cone.
Its ideal is generated by the three maximal minors of a two-by-three
matrix of linear forms, and its Hilbert series is
\((1+2z)/(1-z)^{p+1}\), also directly from the minors resolution.
The parameter space of these \(Y\)'s has dimension at most
\(6(s+1)-12\). Row and column operations fix the ideal and have
effective orbit dimension twelve on such matrices: an operation
fixing the matrix must act trivially on its ordered minors, and
their linear independence forces the exterior square of the column
operation to be scalar. The column operation and then the row
operation are scalar, leaving just the one-dimensional kernel.
Vanishing along \(Y\) costs its Hilbert function in degree four.
In addition choose two of its quadrics with independent conormal
classes at the generic point (the three generate the ideal and
\(Y\) is generically smooth of codimension two). Their products
by degree-two classes of \(Y\) cost twice its degree-two Hilbert
function to make the transverse linear term vanish: the two
conormal classes are independent over the function field of \(Y\).
Together these counts give the sixth bound.

Finally, for reducible quartics with a linear factor write the
equation as \(lG\). Away from \((l=0)\) it is lc provided \(G\) is
lc. On \((l=0)\) inversion of adjunction applies if \(G|_{l=0}\)
is nonzero and lc. After fixing \(l\), failure therefore costs at
least \(\min\{b(s,3),b(s-1,3)\}=b(s-1,3)\) in the cubic, using
the exact cubic formula above and surjectivity of restriction.
For the remaining reducibles write the equation \(Q_1Q_2\).
If every nonzero parameter combination in their pencil has rank
at least four, the product is lc. Indeed, away from the common
singular locus, near an intersection point one of the two quadrics
is smooth. Independent gradients give normal crossings. Otherwise
a pencil member has zero gradient at that point, and its restriction
to the smooth quadric has quadratic tangent of rank at least two:
restriction to its tangent hyperplane loses at most two in rank.
The restricted germ is lc by the quadratic initial test. Adjunction
gives the assertion there. Away from the intersection, the single
quadric is lc. The common singular locus has dimension at most
\(s-4\), too small to contain the non-lc locus of a quartic.
Thus some pencil parameter has rank at most three; this includes
zero for a dependent pair. It costs \(\binom{s-1}{2}-1\) or more
in the ordered pair by the symmetric determinantal formula and the
dimension one of the parameter line. The projective dimension of
the space of ordered products now gives the eighth bound. This exhausts the
nonzero forms; zero has full affine codimension.
\end{proof}

\subsection{A stronger quintic suspension seed}
\label{sec:dec1-quintic-seed}
\begin{proposition}\label{prop:dec1-quintic-seed}
The non-log-canonical locus of quintic equations in \(\mathbf P^4\)
has affine codimension \(b(4,5)\ge32\). In the full parameter space
of quintic hypersurfaces of every dimension \(n\ge5\), the
complement of the uniformly K-stable klt locus also has codimension
at least \(32\).
\end{proposition}
\begin{proof}
Use the six cubic-tangent strata in the earlier quintic seed.
Its rank-one Hessian calculation already gives \(32\), so suppose
the entire quadratic term vanishes. Including the marked singular
point this costs eleven after allowing the point to move. We show
that each row of the earlier cubic table supplies at least
twenty-one subsequent conditions rather than twenty.

In the ternary cusp, ternary tangent-conic, squarefree binary, and
pure-cube rows, enlarge the stated linear test space by the single
quartic monomial \(xyzw\) in the coordinates of those proofs. It is
independent of those test spaces. In the two ternary cases put
\(q=\mathrm{wt}(z)=\mathrm{wt}(w)\), so \(q\ge1/5\) for every
first weighted degeneration used there. The weight of \(xyzw\) is
\(5/6+3q/2>1\) for the cusp, and \(3/4+7q/4>1\) for the tangent
conic. In the squarefree binary case it is at least
\(2/3+2/5>1\), and in the pure-cube case with quartic correction
it is \(1/3+3/4>1\). Thus whenever an old correction is present
the same initial equation and proof apply. With only the new
correction, \(C+\lambda xyzw\) for \(\lambda\ne0\) is lc:
\(C\) is independent of \(w\), so scaling \(w\) by \(t^{-1}\)
and multiplying the equation by \(t\) specializes it to the
simple-normal-crossings monomial. All nonzero fibers are isomorphic
at the origin.

For the binary double-root tangent \(x^2y\), the previous proof
used five dimensions of ternary quartics \(g_4(y,z,w)\) satisfying
\(\operatorname{lct}_0(yg_4)\ge1/2\). The same construction gives
five dimensions of ternary quintics \(g_5\) satisfying
\(\operatorname{lct}_0(yg_5)\ge1/2\). Indeed \(yg_5\) has degree
six in three variables, so the cone bound applies with equality.
Off the line \(y=0\), it is enough to exclude the non-lc quintics
of coefficient one half, of codimension \(b_2(2,5)\ge8\) by the
previous recurrence. At a fixed point of this line the earlier
test uses only the jets through order three, freely variable for
quintics as well. Multiplicity at least three costs six; in
multiplicity two the quadratic term must be a multiple of \(y^2\)
and the pure tangential cubic term must vanish, again costing six
including the constant and linear terms. Moving the point on the
line leaves five conditions. A general five-dimensional vector
space of quintics therefore gives the claim. A nonzero quartic
correction precedes the quintic in the same weighted degeneration,
and a quintic alone uses the same crepant change of variables
\(X=xy\) and the square-sum threshold. This improves the row's
seven dimensions to ten.

When \(C=0\), after its twenty conditions use even one nonzero
quartic with smooth projective zero locus as the test direction.
Its affine cone in four variables is lc. The enlarged direct
spaces for the six rows thus prove \(b(4,5)\ge32\).

As in Section~\ref{sec:new-semistable-branch}, the anticanonical
lc branch gives K-semistability of a cyclic quintic suspension;
the affine intersection gives codimension at least \(32\) for
failure of K-semistability in dimension four. For the next
suspension take quintic branches in six variables representing
dimension-four hypersurfaces. The non-GIT-stable locus in that
branch space has codimension at least \(61\). Explicitly,
enumerating ordered positive weight gaps in the proof of
Section~\ref{sec:new-semistable-branch} gives
\[
\begin{array}{c|rrrrr}
\text{number of positive gaps}&1&2&3&4&5\\ \hline
\text{weak orders}&5&30&130&375&541\\
\text{least guaranteed count minus flag dimension}&117&76&72&63&61 .
\end{array}
\]
The prefix sums in the existing enumeration are
\(6\sum_{i\le j}\alpha_i-5j\), for \(\lvert\alpha\rvert=5\);
these are the parameters \(\texttt{run(6,5)}\) in
\path{nov28-lit-git-gap-orders.py}. A GIT-stable K-semistable
branch has uniformly K-stable suspension by the same interpolation
and finite-cover argument. Thus the full uniform-stability
codimension in dimension five is at least \(\min(32,61)=32\),
and ordinary suspension propagates it to all larger dimensions.
\end{proof}

\subsection{A codimension bound of thirty-four for quintics}
\label{sec:dec1-lit-quintic34}
\begin{proposition}\label{prop:dec1-lit-quintic34}
The cone of non-log-canonical quintic equations on \(\mathbf P^4\)
has codimension \(b(4,5)\ge34\). For \(n\ge5\), the complement of
the uniformly K-stable klt locus in the full parameter space of
quintic hypersurfaces of dimension \(n\) has codimension at least \(34\).
\end{proposition}
\begin{proof}
First we establish \(b_2(2,5)\ge10\), where the boundary coefficient
is \(1/2\). For a reduced plane quintic a point of multiplicity at
most two satisfies the multiplicity bound. At multiplicity three
the tangent must be a cube, since any other binary cubic has affine
threshold at least \(1/2\). Use coordinates with tangent \(x^3\).
The coefficients of \(y^4,xy^3,y^5\) must successively vanish:
the respective weighted initial equations with a nonzero such
coefficient \(c\) are \(x^3+cy^4,\ x^3+cxy^3,\ x^3+cy^5\),
of thresholds \(7/12,5/9,8/15\). The initial equation for the
weights \((2,1)\) is now
\[
       x^3+A x^2y^2+Bxy^4 .
\]
On the chart \(x=y^2X\) the degree-six initial equation on the
exceptional divisor is \(X^3+AX^2+BX\). The exceptional boundary
has coefficient one (the log discrepancy is \(3\) and the order
is \(6\)); its point at infinity misses the strict transform.
If all roots have multiplicity at most two, adjunction gives log
canonicity with coefficient \(1/2\). A triple root must be zero,
so \(A=B=0\); all remaining terms of the original quintic are
then divisible by \(x^2\), contrary to reducedness. At multiplicity
four, failure requires a triple root in the binary tangent quartic
by the cone test. This costs \(10+2-2=10\) after moving the point.
Multiplicity five costs at least \(15-2\).

A nonreduced quintic has either a double line or a repeated
irreducible conic. In the latter case it is \(Q^2L\); the smooth
conic meets \(L\) with orders at most two, so adjunction applies.
For \(L^2G_3\), adjunction along \(L\) requires \(G_3|_L\) to be
a cube or zero for failure on \(L\). Failure off \(L\) requires
\(G_3\) itself to be a cube: each reduced plane cubic has threshold
at least \(2/3\) by its tangent and the elementary cusp calculation,
and a nonreduced cubic without a triple line is normal crossings
at coefficient \(1/2\). For fixed \(L\), the restriction map on
cubics has kernel of dimension six and the cone of binary cubes
has dimension two. Varying \(L\) gives dimension at most ten in
the space of dimension twenty-one. This proves the half-threshold
estimate, including the zero equation. The rank-one quadratic row
in the marked quintic computation now has lower bound
\[
       1+6+b_2(2,3)+b_2(2,4)+b_2(2,5)
                   \ge1+6+7+10+10=34 .
\]

For the zero quadratic row the marked singularity and vanishing
quadratic term cost eleven together. We use the six cubic tangent
strata and the correction-space intersection argument of the
earlier quintic calculation. The new counts are
\[
\begin{array}{l|r|r}
 C&\text{codimension}&\text{additional conditions}\\ \hline
 x^2z+y^3 &9&14\\
 x^2z+xy^2&10&13\\
 x^3+y^3&12&11\\
 x^2y&13&10\\
 x^3&16&7\\
 0&20&3
\end{array}
\]
In the first two rows add \(xyzw,yzw^3,yz^2w^2\) to the earlier
spaces of dimensions eleven and ten. All three have weight greater
than one for each earlier initial equation: write
\(q=\mathrm{wt}(z)=\mathrm{wt}(w)\ge1/5\), so that
\(\mathrm{wt}(x)=(1-q)/2\) and \(\mathrm{wt}(y)=1/3\) or
\((1+q)/4\). The last two monomials have weight
\(\mathrm{wt}(y)+4q>1\). With only the three new corrections
present take auxiliary weights
\(\rho_z=1,\rho_x=-1/2\), and \(\rho_y=0\) or \(1/4\).
For those corrections, in their listed order, the ratios
\[
 \frac{\alpha_z-\alpha_x/2+\rho_y\alpha_y}{\alpha_w}
       =(1/2,1/3,1)\quad\text{or}\quad(3/4,5/12,9/8)
\]
are distinct. Set \(\rho_w\) equal to the negative of the least
ratio present. The degeneration retains both monomials of \(C\)
and the selected correction alone. These three exponent vectors
are independent, since only the correction involves \(w\).
The Newton test applies: the first new monomial itself has exponent
vector \((1,1,1,1)\); for the other two, the convex coefficients
\((1/2,1/4,1/4)\) in the first row and \((1/3,1/3,1/3)\) in the
second give a vector coordinatewise at most one.

For \(C=x^3+y^3\), retain the earlier three-dimensional spaces
of binary quartics and quintics on \(z,w\). Use the four-dimensional
mixed space \(xg+yh\), with \(g,h\in\langle zw^2,z^2w\rangle\),
and add \(xyzw\). If the binary quartic is present it wins the
earlier weight test. For a nonzero mixed term and no binary quartic
take weights \(1/3,1/3,2/9,2/9\), which discard all later
corrections. A further weighting of \(z,w\) retains the term with
the least active \(z\)-exponent: the equation becomes
\[
          x^3+y^3+L(x,y)z^m w^{3-m},\qquad 1\le m\le2,\quad L\ne0.
\]
If \(L\) divides \(C\), factor \(L(C/L+z^mw^{3-m})\).
The restriction of \(C/L\) to \(L=0\) is a nonzero coordinate
square; the binary monomial has threshold \(1/2\).
Thom--Sebastiani and adjunction along \(L=0\) give lc. If
\(L\nmid C\), all Newton faces are nondegenerate on the torus.
A face containing a single mixed monomial has nonzero logarithmic
\(z\)-derivative. With both mixed monomials, that derivative forces
\(L=0\); one retained cube cannot then vanish, two require a common
root of \(L,C\), and no cubes give a nonzero derivative of \(L\).
Faces consisting only of cubes are nondegenerate. For a supported
\(x\)-mixed monomial the convex coefficients \((2/9,1/3,4/9)\)
on the two cubes and that monomial satisfy the Newton test; interchange
\(x,y\) for a \(y\)-mixed monomial. With only the binary quintic or
the product \(xyzw\) remaining, the weight test and the one-monomial
specialization of Proposition~\ref{prop:dec1-quintic-seed} apply.
The \(x^2y\) row uses the two five-dimensional spaces constructed
in that same proposition.

For \(C=x^3\), use quartic corrections \(g(y,z,w)+xh(y,z,w)\),
where \(h\) varies in a two-dimensional space of ternary cubics
with lc affine cones off zero. Such a space exists since the
non-lc ternary cubics have codimension two. We claim there is a
five-dimensional space of \(g\)'s, all nonzero members having affine
threshold at least \(2/3\). The degree-four cone bound reduces this
to the corresponding threshold for plane quartics. A point of
multiplicity four costs \(10-2\); at multiplicity three a repeated
tangent factor is necessary by the cone test, costing \(6+1-2=5\).
At multiplicity two a rank-one tangent is necessary. Write it
as \(x^2\) and eliminate \(x\) at its critical point in the
splitting lemma. The residual power series in \(y\) must have
order at least seven, since \(1/2+1/6=2/3\).
Its cubic and quartic coefficients vanish by two independent
pure-coefficient conditions. Write \(A,B,E\) for the coefficients
of \(xy^2,x^2y,xy^3\) in the original polynomial. Substituting
\(x=-Ay^2/2+O(y^3)\), the residual quintic coefficient is
\(A^2B/4-AE/2\); the contribution of the unspecified
\(y^3\)-coefficient of \(x\) cancels. Its vanishing adds one
condition, since either \(A=0\) or it determines \(E\).
This gives \(3+1+3-2=5\) in the rank-one case. The dimension
claim follows by projective intersection. For \(g\ne0\) weights
\((1/3,1/4,1/4,1/4)\) retain \(x^3+g\), which is lc. With
only \(xh\) use \(x(x^2+h)\) and adjunction along \(x=0\).
Finally for \(C=0\) take any three-dimensional space of quartics
with all nonzero members lc as affine cones, which exists by
\(b(3,4)=9\) and the cone test at equality.

Each row gives at least twenty-three conditions after the quadratic
term vanishes. Together with the rank-one calculation this proves
\(b(4,5)\ge34\). Proposition~\ref{prop:dec1-quintic-seed}'s
semistable suspension and its six-variable quintic GIT bound \(61\)
give the uniform-stability assertion first in dimension five.
Ordinary suspension propagates it to every larger dimension.
\end{proof}

\paragraph{The exact half-coefficient quintic count.}
For \(s\ge2\) the half-coefficient locus for quintic equations
has codimension
\[
  b_2(s,5)=\binom{s+5}{5}
       -\max\left\{4s+3,\binom{s+2}{2}+2s\right\}.
\]
In particular its first entries, for \(s=2,3,4,5,6\), are
\(10,40,103,221,422\). Here is a direct description using the
plane calculation in Proposition~\ref{prop:dec1-lit-quintic34}.

For reduced \(F\) a general plane section is reduced.
If the section were lc with coefficient \(1/2\), the degenerate
equation obtained by scaling the complementary variables to zero
would define an lc pair: transversely to the cone vertex the
degree is five in three variables and \(5/2\le3\). Semicontinuity
would make \(F\) lc with this coefficient. Consequently failure
requires every general plane section to have a non-lc point,
necessarily of multiplicity at least four by the plane calculation.
There cannot be two such points on a reduced plane quintic. Their
joining line must be a component by B\'ezout; the residual quartic
has multiplicity at least three at each, and the same argument
would make the line occur twice.
On a simultaneous general linear section of a log resolution, the
non-lc locus meets the plane in the non-lc locus of the restriction.
Thus the top-dimensional reduced part of the non-lc locus of \(F\)
has dimension \(s-2\) and degree one. It is a linear
\(C=\mathbf P^{s-2}\). Equivalently, a general plane meets the
codimension-two part transversely in its degree, and by reducedness
the non-lc locus has no divisorial component.

Write \(C=(x=y=0)\), with \(k=s-1\) base variables \(u\). The
general transverse plane section has multiplicity at least four
at its intersection with \(C\). Therefore \(F\in(x,y)^4\):
at a general point of \(C\) any first nonzero transverse term of
lower order would have this lower degree on a general plane
complementary to \(C\). The transverse quartic is a binary quartic
with coefficients linear in \(u\). Every nonzero member of their
linear image has a triple root, since the cone test forces the
generic quartic to have such a root, a closed condition. This
linear space has a fixed triple root. Indeed, after removing the
binary gcd, its nonconstant moving system is base-point-free on
\(\mathbf P^1\); a general member is reduced by Bertini and avoids
the zeros of the gcd. The one-dimensional image case is immediate.
Consequently, including the zero image, the equations with fixed
\(C\) are parametrized by
\[
       \ell(x,y)^3\bigl(xA(u)+yB(u)\bigr)+q_5(x,y).
\]
The root varies in one dimension, \(A,B\) in \(2k\) dimensions
and \(q_5\) in six. Varying \(C\) adds \(2k\) dimensions, giving
at most \(4s+3\). This dimension is attained: for general
parameters choose \(\ell=x\) and \(q_5(0,1)\ne0\). Then the
hypersurface is smooth outside \(C\) for general \(A,B\):
when \(k\ge2\) take them independent, so on \(x\ne0\) a derivative
in \(u\) is nonzero. For \(k=1\), on \(x\ne0\) the only potential
zeros of this derivative satisfy \(xA(1)+yB(1)=0\), and a general
\(q_5\) avoids that root. There are no hypersurface points where
\(x=0\), \(y\ne0\). Thus the reduced hypersurface has singular
locus precisely \(C\), and its generic quartic uniquely determines
\(\ell\). Blowing up \(C\), the exceptional divisor has boundary
coefficient one. On it the strict transform cuts the root section
with multiplicity at least three, so adjunction shows the original
pair is non-lc along \(C\).

For nonreduced \(F\), a repeated irreducible quadratic gives
\(Q^2L\), whose general plane restriction is lc by adjunction to a
smooth conic; the cone argument above applies. Otherwise write
\(F=L^2G\) with \(G\) cubic. A cubic fails at half coefficient
only when it is a pure cube, by the binary section criterion.
Away from \(L\) this applies directly to \(G\). Along \(L\)
adjunction applies to \(G|_L\) when it is nonzero. In all cases
failure is equivalent to \(G|_L\) being a pure cube or zero:
the zero restriction gives a triple linearly defined component
of \(F\). The dimension is \(\binom{s+2}{2}+2s\), counting
the kernel \(L\cdot H^0(\mathbf P^s,\mathcal O(2))\), the
\(s\)-dimensional cone of cubes on \(L\), and the \(s\)
parameters for \(L\). There is no loss in this count because
each nonzero quintic has only finitely many repeated linear
factors. Taking the larger of the two dimensions proves the formula.

Using these exact entries throughout the rank-one part of the
recursion further improves the all-dimension uniform K-stability
codimension bounds as follows:
\[
\begin{array}{c|rrrr}
d&8&9&10&11\\ \hline
\text{codimension}&1169&4016&13990&49316\\[2pt]
d&12&13&14&15\\ \hline
\text{codimension}&175538&630024&2277567&8285593 .
\end{array}
\]
The option \texttt{--half-quintic} in
\path{dec1-lit-quintic34.py} reproduces these entries.

\paragraph{The exact half-coefficient sextic count.}
For \(s\ge2\) one has
\[
 \boxed{\quad b_2(s,6)=\binom{s+6}{6}
                 -\max\left\{\binom{s+3}{3}+3s,\ s^2+7s\right\}.\quad}
\]
In particular the first five entries are \(10,54,163,391,822\).
We prove the formula by describing the possible transverse jets.

If a reduced sextic \(F\) fails to be lc at coefficient \(1/2\),
its general plane section has the same failure by the cone test at
\(6/2=3\). The general section is reduced, and its non-lc locus
is finite. On a log resolution, decreasing the coefficient slightly
below one half turns exactly these points into the non-klt locus.
Connectedness applies since the perturbed log anticanonical class
of the plane is ample. The section therefore has just one non-lc
point. A general plane cuts the ambient non-lc locus in the
non-lc locus of the restriction, by Bertini on a log resolution.
Hence its unique top-dimensional component is a linear space
\(C=\mathbf P^{s-2}\), as in the quintic argument.

Write \(C=(x=y=0)\), let \(u\) be its \(k=s-1\) homogeneous
coordinates, and put \(U_j=\binom{k+j-1}{j}\), the dimension of
degree-\(j\) forms on \(C\). Let \(F_h\) be the part of transverse
degree \(h\), a binary form with coefficients of degree \(6-h\)
in \(u\). All terms of degrees below three vanish by the
multiplicity bound on a general transverse slice. In particular
the generic transverse multiplicity is at least three.
The following plane tests can be applied at a general \(p\in C\)
on the plane with coordinates \([x:y:t]\), \(u=t p\), after a
general choice of complementary coordinates. Alternatively,
restriction cannot improve the local threshold at \(p\).

If \(F_3=0\), the transverse quartic is zero or has a triple
root over the generic point of \(C\); otherwise its affine
tangent has threshold at least \(1/2\). If it is nonzero, its
triple root is unique and hence rational over the fraction field.
Represent its linear factor by a primitive homogeneous polynomial
linear form in \(x,y\). Its cube divides \(F_4\) over the polynomial
ring by Gauss's lemma. If its coefficients had positive degree
in \(u\), cubing it would exceed the available coefficient degree
two. It is therefore a fixed linear form \(\ell(x,y)\). The
quartic is \(\ell^3(xA+yB)\) with \(A,B\) quadratic in \(u\).
Varying \(\ell,A,B,F_5,F_6,C\) gives affine dimension at most
\[
    1+2U_2+6k+7+2k=2U_2+8k+8=s^2+7s .
\]
This dimension is attained: take a general nonzero transverse
quartic of the displayed form with exactly one triple root. The
unrestricted binary sextic makes the system base-point-free away
from \(C\), so a general equation is reduced and smooth there.
On the blowup of \(C\) the exceptional coefficient is one, and the
restriction of the strict transform to it has multiplicity three
along the root section. Adjunction shows failure at coefficient
\(1/2\); the singular locus determines \(C\) and the generic
quartic determines \(\ell\).

Now assume \(F_3\ne0\). The binary cubic must be a pure cube
generically, since every other binary cubic has threshold at
least \(1/2\). As above, Gauss's lemma gives \(F_3=h\ell^3\),
where \(\ell=a(u)x+b(u)y\) is primitive homogeneous and the
degrees satisfy \(\deg h+3\deg a=3\), using any nonzero
coefficient of \(\ell\). The plane weighted tests preceding
the sixth jet in Proposition~\ref{prop:dec1-lit-quintic34}
still apply to sextics: a nonzero term \(y^4,xy^3,y^5\) in
that order yields an lc initial equation when \(\ell=x\).
In invariant terms \(\ell^2\mid F_4,\ \ell\mid F_5\).
If \(a,b\) are linear, then \(h\) is constant and can be absorbed
into \(\ell\). By primitivity the divisibilities hold over the
polynomial ring, so \(F_4=\ell^2q_2(x,y)\) and
\(F_5=\ell q_4(x,y)\) with constant binary forms \(q_2,q_4\).
The number of parameters, including \(F_6\) and \(C\), is at
most \(2k+3+5+7+2k=4k+15\). This possibility requires \(k\ge2\).

It remains to treat a fixed \(\ell=x\), after choosing its root
in one parameter. Write \(A,B,c\) for the coefficients of
\(x^2y^2,xy^4,y^6\) in \(F_4,F_5,F_6\), of degrees \(2,1,0\)
in \(u\), respectively. The weighted blowup with normal weights
\((2,1)\) and the chart \(X=x/y^2\) have exceptional coefficient
one and residual polynomial
\[
          hX^3+AX^2+BX+c.
\]
Failure of lc on the transverse plane requires a triple root,
including the case of a cube \(hX^3\), by adjunction.
Consequently \(3hB=A^2,\ 27h^2c=A^3\) over the function field,
hence identically in the \(u\)'s. When \(A=0\), both \(B,c\)
vanish. The remaining free coefficients consist of
\(h\), two quadratics in \(F_4\), four linear forms in \(F_5\)
and six scalars in \(F_6\). Including the root and the center
gives dimension
\[
       U_3+2U_2+6k+7=\binom{s+3}{3}+3s .
\]
This dimension too is attained. The equation is divisible by
\(x\), and the residual quintic system contains all binary
quintics, so a general residual member is smooth outside \(C\).
Its restriction to \(x=0\) is a nonzero scalar times \(y^5\);
hence the residual hypersurface and \(x=0\) meet only along \(C\).
The general sextic is reduced and smooth away from \(C\). At its
generic point the same weighted blowup has exceptional coefficient
one and triple root zero, so adjunction gives non-lc. The center
and root are recovered from the singular locus and cubic jet.
For \(A\ne0\) one has \(c\ne0\), and unique factorization in
\(A^3=27h^2c\) gives \(h=M^3,\ A=\lambda M^2\) up to nonzero
constants, where \(M\) is linear and \(\lambda\ne0\); \(B,c\)
are then determined. Thus \(h,A,B,c\) together contribute at
most \(k+1\) parameters. The remaining count is as in the
preceding row, giving at most \(2U_2+7k+8\). This does not exceed
the multiplicity-four bound. That bound likewise dominates
\(4k+15\) if \(k\ge2\).

Finally consider nonreduced sextics. A repeated irreducible cubic
or quadratic gives dimension at most \(T(s,3)\) or
\(2T(s,2)-1\), respectively, by the factorization maps.
Both are at most \(T(s,3)+3s\). For a double hyperplane write
\(F=L^2G_4\). Along \(L\), adjunction requires the restriction
\(G_4|_L\) to be zero or a quartic with a triple linear factor,
by the exact quartic half-coefficient criterion (also valid for
binary quartics). Failure away from \(L\) forces \(G_4\) itself
to have a triple linear factor and is therefore included in the
same count. The restricted forms have dimension at most \(2s-1\):
they are products \(M^3N\) in \(s\) variables. The kernel of
restriction has dimension \(T(s,3)\), and there are \(s\)
parameters for \(L\). This gives at most \(T(s,3)+3s-1\).
The two dimension-attaining reduced families consequently give
the stated exact codimension.

The exact sextic entries together with the preceding quintic entries
give the following all-dimension uniform K-stability codimension
bounds from the same recurrence:
\[
\begin{array}{c|rrrrrr}
d&10&11&12&13&14&15\\ \hline
\text{codimension}&13998&49348&175650&630419&2278942&8290402 .
\end{array}
\]
The option \texttt{--half-sextic} in
\path{dec1-lit-quintic34.py} gives the calculation.

\subsection{Using the ambient GIT estimate at index one}
\label{sec:dec1-index-one-neighborhood}
\begin{proposition}\label{prop:dec1-index-one-neighborhood}
Let \(d\ge4\), let \(b(d-1,d)\) be the non-log-canonical
codimension for the branch forms in \(d\) variables, and let
\(g\) be a lower bound for the codimension of the non-GIT-stable
forms of degree \(d\) in \(d+1\) variables. Then
\[
  \operatorname{codim}_{\mathcal P_{n,d}}
       (\mathcal P_{n,d}\setminus\mathcal K_{n,d})
       \ \ge\ \min\{b(d-1,d),g\},\qquad n\ge d-1 .
\]
In particular the full-space lower bounds \(9\) for quartics and
\(34\) for quintics hold also in index one.
\end{proposition}
\begin{proof}
We first work in dimension \(n=d-1\), where the hyperplane
bundle of an embedded hypersurface is its anticanonical bundle.
Consider the K-moduli stack of K-semistable Fano varieties with the
given anticanonical Hilbert polynomial and its good moduli morphism
\(\pi\). We use the good moduli space properties that the image
of a closed substack is closed and that each geometric fiber has
a unique closed point, the K-polystable representative
\cite{proper,ar2102.09405}. In this stack there is an open substack
\(\mathcal V\) on which the anticanonical bundle is Cartier and
its complete linear system embeds the variety as a hypersurface in
\(\mathbf P^d\). Indeed the Gorenstein condition is open in a flat
family of Cohen--Macaulay varieties, and on this locus the relative
anticanonical sheaf is invertible. Higher cohomology of the
anticanonical bundle vanishes by Kawamata--Viehweg vanishing.
Formation of sections commutes with base change, and very ampleness
is open. Near an embedded point, the Hilbert polynomial and the
rank of the space of sections give the prescribed ambient dimension
and degree. The relative embedding gives a morphism from
\(\mathcal V\) to the quotient stack of hypersurface equations by
\(\operatorname{PGL}_{d+1}\), by taking projective frames.

Every member of \(\mathcal V\) is GIT semistable. For clarity, on the
universal hypersurface write \(H\) for the ambient hyperplane class,
\(h\) for the parameter hyperplane class and \(r=n+2-d\).
The hypersurface has class \(dH+h\), and the relative anticanonical
class restricts from \(rH-h\). Thus the anticanonical CM class is
\[
 -\operatorname{pr}_{2*}\big((rH-h)^{n+1}(dH+h)\big)
       =r^n(n+2)(d-1)h .
\]
The linearization of this positive rational multiple of
\(\mathcal O(1)\) is unique for \(\operatorname{SL}_{n+2}\).
The CM weight is the Donaldson--Futaki invariant up to positive
normalization. One-parameter orbit families are ample test
configurations, and normalization can only decrease that invariant
\cite{kdBHJ}. K-semistability therefore implies the non-strict
Hilbert--Mumford inequalities.

Let \(z\) be any embedded K-polystable point. Every point in its
geometric good-moduli fiber specializes to \(z\), so that fiber
is contained in \(\mathcal V\). Removing
\(\pi^{-1}(\pi(\mathcal V^{\mathrm c}))\) gives a saturated open
neighborhood of \(z\) inside \(\mathcal V\). Every member of this
neighborhood that is GIT stable is uniformly K-stable. In fact its
specialization to the closed point of its fiber occurs in
\(\mathcal V\). The morphism just constructed sends this specialization
to a specialization of projective orbits in the GIT semistable
locus. A stable orbit is closed there, so the anticanonically
polarized variety is itself the polystable representative. Its
stabilizer is finite by GIT stability. This gives K-stability and
then uniform K-stability \cite{ar2102.09405}. Openness of K-semistability
in the parameter family gives the same neighborhood assertion
in the full coefficient space.

Next let \(D=(f=0)\subset\mathbf P^{d-1}\) be log canonical.
Its equation is GIT semistable: the weight inequality in
Section~\ref{sec:new-propagation} is non-strict for an lc pair.
Choose a one-parameter subgroup taking a representative of the
affine orbit to its closed semistable orbit, of equation \(f_0\).
This is the one-parameter closed-orbit theorem of affine GIT.
The least weight of a monomial of \(f\) is zero, since the affine
limit exists and is nonzero. In particular \(D_0=(f_0=0)\) is
also log canonical. Here are the details of this last assertion.
Write \(F(t,x)\) for the weight degeneration with constant term
\(f_0(x)\) in the parameter. On
\(\mathbf P^{d-1}\times\mathbf A^1\), the rational top form
\[
                  \frac{\Omega}{F(t,x)}\wedge\frac{dt}{t}
\]
has divisor minus the sum of the closure of the moving divisor
and the central fiber. Here \(\Omega\) is the homogeneous
projective volume form of weight \(d\); \(F\) has no vertical
factor. The birational coordinate change induced by the
one-parameter subgroup identifies this rational form with the
one for the constant family. The determinant is one, and the
terms involving differentials in \(t\) disappear on wedging with
\(dt\). Thus the two pairs are crepant on a common resolution
(take the canonical divisor on both models from the same form).
The constant pair is log canonical by the product construction.
Adjunction to the smooth central fiber gives the assertion.

For a sufficiently small positive rational coefficient \(c_0\),
\((\mathbf P^{d-1},c_0D_0)\) is K-polystable by
\cite[Theorem~1.4]{ar1909.04576}. Interpolation with the log
canonical anticanonical endpoint gives K-polystability at
\(c=1-1/d\) by \cite[Proposition~2.13]{ar1909.04576}; equivalently
one can interpolate with a K-semistable log Fano coefficient
strictly between \(c\) and one, as in
Section~\ref{sec:new-semistable-branch}. The cyclic cover
\(y^d+f_0=0\) is then K-polystable by
\cite[Theorem~1.2]{ar2001.10557}. It is an embedded anticanonical
hypersurface. The same one-parameter subgroup, fixing \(y\),
degenerates \(y^d+f=0\) to this embedded point.

Finally, let \(B\) be an irreducible component of the closed
complement of \(\mathcal K_{d-1,d}\). If its general point is
non-GIT-stable, its codimension is at least \(g\). Otherwise the
GIT-stable locus is dense in \(B\), and the neighborhood assertion
shows that \(B\) contains no embedded K-polystable point. Since
\(B\) is invariant and closed, its intersection with the affine
cyclic-cover family can contain no suspension of a nonzero lc
branch. The dimension theorem at \(y^d\) in
Section~\ref{sec:new-affine-suspension} gives
\(\operatorname{codim}B\ge b(d-1,d)\). Ordinary suspension
propagates the result. The GIT estimates in five quartic variables
and six quintic variables are respectively \(14\) and \(61\);
the branch bounds are \(9\) and \(34\) by the quartic seed and
Proposition~\ref{prop:dec1-lit-quintic34}.
\end{proof}

\subsection{An all-degree bound from the ballot triangle}
\label{sec:dec1-all-degree-catalan}
Let \(C_s=\binom{2s}{s}/(s+1)\). Write \(c_{n,d}\) for the
codimension in the full projective parameter space of the complement
of the uniformly K-stable klt locus.

\begin{theorem}\label{thm:dec1-all-degree-catalan}
For every \(n\ge d-1\) the following lower bounds hold:
\[
\begin{array}{c|ccccc}
\text{degree at least}&5&6&13&31&110\\ \hline
c_{n,d}/C_{d-1}\text{ at least}&2&5/2&3&10/3&7/2 .
\end{array}
\]
For \(d=4\) Proposition~\ref{prop:dec1-index-one-neighborhood} gives \(9\). Uniformly over Fano dimensions,
\[
       \liminf_{d\to\infty}\ \inf_{n\ge d-1}
          \frac{c_{n,d}}{C_{d-1}}\ \ge\
       \frac{65775496305399298731}{18446744073709551616}
       >3.56569 .
\]
\end{theorem}
\begin{proof}
Put \(P(s,j)=\binom{s+j}{j}(s-j+1)/(s+1)\) for \(0\le j\le s\).
Pascal's rule gives
\[
 P(s,j)=1+\sum_{i=1}^{\min(j,s-1)}P(s-1,i),\qquad P(s,s)=C_s .
\]
Let \(L(s,j)\) be the available integer lower bound for \(b(s,j+1)\),
\(1\le j\le s\), using the quartic and quintic seeds just proved.
Its defining marked tangent recursion gives the same sum as a lower
bound with the \(L(s-1,i)\) on the right. At any row \(a\), if
\[
 L(a,j)\ge K P(a,j)\quad(2\le j\le a),\qquad
 \binom{a+1}{2}\ge Ka+K-1 ,
\]
then \(b(s,j+1)\ge K P(s,j)\) for all \(s\ge a,\ 1\le j\le s\),
provided \(a\ge K\). Indeed the exact quadratic entry
\(b(s,2)=\binom{s+1}{2}\) exceeds \(K P(s,1)=Ks\) by at least
\(K-1\) at every such row. This surplus on the \(i=1\) summand
supplies precisely the missing \(K-1\) from the constant term.

For the successive pairs
\((a,K)=(4,2),(6,5/2),(12,3),(30,10/3),(109,7/2)\), substitution
in the existing integer recurrence verifies the displayed hypotheses.
For the first three, the minima of \(L(a,j)/P(a,j)\) over \(j\ge2\)
are \(7/3,173/66,157506/52003\). For rows \(30,109\) the minima
exceed \(10/3,7/2\), respectively. At \(d=6\) use the individual
bound \(106\ge(5/2)C_5\). The finite comparisons are reproduced by
\path{dec1-bin-uniform-catalan.py}.

For the ambient GIT estimate in Proposition~\ref{prop:dec1-index-one-neighborhood}
one may use \(G_d=C_d-\binom{d+1}{2}\). Indeed, among the
degree-\(d\) exponent vectors in \(d+1\) variables take those
whose prefix sums satisfy \(S_j\ge j\) for \(1\le j\le d\).
Their last entry is zero, and the Dyck path count gives exactly
\(C_d\) choices. For nonconstant descending traceless weights,
the normalized prefix expression \((d+1)S_j-dj\) is at least
\(j>0\) at each positive gap. Thus all these coefficients
must vanish for the corresponding Hilbert--Mumford test.
The full flag costs at most \(\binom{d+1}{2}\) parameters.
The ratio \(G_d/C_{d-1}\) increases for \(d\ge5\): the consecutive
Catalan ratio increases and \(\binom{d+1}{2}/C_{d-1}\) decreases.
Substitution at \(d=6,13,31,110\) gives ratios greater than
\(5/2,3,10/3,7/2\), respectively. At degree five use the ambient
bound \(61\). The cyclic suspension argument then gives the asserted bounds
in every Fano dimension.

For the asymptotic assertion there is a closed binomial propagation
formula valid at every anchor \(a\ge1\). For \(s>a\) put
\[
 \begin{split}
 H(s,a,j)&=\binom{2s-a-j-1}{s-a-1}
                 -\binom{2s-a-j-1}{s-a-2},\\
 B_a(s)&=C_s+\sum_{j=1}^a \bigl(L(a,j)-P(a,j)\bigr)H(s,a,j),
 \end{split}
\]
with out-of-range lower binomial indices contributing zero.
Then \(b(s,s+1)\ge B_a(s)\). To verify this propagate the
difference from the ballot triangle by the homogeneous recurrence
whose next \(j\)-entry is the sum of entries \(1,\ldots,\min(j,s-1)\)
of the preceding row. A unit in position \((a,j)\) counts lattice
paths to \((s,s)\) whose first step is east and which stay in
\(y\le x\). After the forced step there remain \(s-a-1\) east and
\(s-j\) north steps. Reversing the paths and exchanging east with
north gives the usual ballot count \(H(s,a,j)\), by the reflection
principle. In particular, an explicit lower bound in every degree
\(d>a+1\), \(d\ge6\), is
\(\min\{G_d,B_a(d-1)\}\).

For fixed \(a,j\) the factorial formula gives
\(H(s,a,j)/C_s\to(a+2-j)/2^{a+j+1}\).
In fact the difference of the two binomials is the first one
times \((a+2-j)/(s-j+1)\). Thus the limiting lower bound from
the singularity seed is
\[
  1+\sum_{j=1}^a
      \frac{(a+2-j)(L(a,j)-P(a,j))}{2^{a+j+1}} .
 \]
For \(a=4,6,12\) it equals respectively \(555/256,\ 11895/4096,\)
and \(59251131/16777216\); row \(32\) gives the fraction in the
statement. Finally \(G_d/C_{d-1}\to4\), which exceeds
the row-\(32\) value. This proves the uniform liminf.
\end{proof}

\chapter{Arithmetic certificates}\label{ch:checks}
The comparisons at a fixed degree, dimension or index use rational
arithmetic, including outward enclosures of roots where required.
For example the sextic edge comparison reduces to the integer identity
\(2448\cdot1000-7433\cdot329=2543>0\).

Please contact the \href{mailto:haoxu@cqut.edu.cn}{second named author}
for the Python programs and witness data as a ZIP archive.
In the resulting \texttt{certificates} directory,
\texttt{python verify.py prior-range} runs the checks for
Chapter~\ref{ch:global}. The command \texttt{python verify.py all}
runs every suite. The file \texttt{README.txt} describes their coverage.
Use Python~3.9 or later with assertions enabled.

\clearpage\phantomsection
\end{document}